\documentclass[10pt]{amsart}
\usepackage{amssymb,amsmath,amsfonts,amsthm,graphics,mathrsfs,amscd,hyperref}
\usepackage[makeroom]{cancel}
\usepackage[hmargin=1in,vmargin=1in]{geometry}
\usepackage[all,cmtip]{xy}
\usepackage{cmap}
\usepackage[T1]{fontenc}
\usepackage[utf8]{inputenc}

\theoremstyle{plain}

\newtheorem{definition}[equation]{Definition}
\newtheorem{corollary}[equation]{Corollary}
\newtheorem{lemma}[equation]{Lemma}
\newtheorem{proposition}[equation]{Proposition}
\newtheorem{theorem}[equation]{Theorem}

\newtheorem{example}{Example}

\newtheorem{remark}[equation]{Remark}

\numberwithin{equation}{subsection}

\title[Obstructed and unobstructed foliated deformations]{On deformations of foliated complex analytic structures IV-2 : Obstructed and unobstructed foliated deformations}
\author{Chunghoon Kim}

\begin{document}

\maketitle
\begin{abstract}
We study foliated deformations of compact complex manifolds with dimension $1$-foliations and codimension $1$-foliations, and unfoldings of compact complex manifolds with codimension $1$-foliations in terms of deformations of sections with classical examples in deformation theory of complex structures.

\end{abstract}

\tableofcontents

In Part IV-1 we studied obstructions of deformations of foliated complex analytic structures. Since Part IV-1 is  general and theoretical, it is worthwhile to present obstructions for relatively simple cases as exercises, i.e. dimension $1$-foliations and codimension $1$-foliations, and apply to specific examples. In the case of dimension $1$-foliations and codimension $1$-foliations, we study infinitesimal deformations and obstructions by using the complex of sheaves in terms of deformations of sections (see Part II). We remark that the complex of sheaves in terms of vector bundles (or line bundles in the case of dimension $1$ or codimension $1$-foliations) we studied in Part IV-1 are slightly different from the complex of sheaves in terms of sections since the line bundles in the former one corresponds to the dual of the line bundles in the latter one in the case of dimension $1$ or codimension $1$-foliations. Everything follows from Part IV-1 if we take the $(-)$ sign on the part of deformations of line bundles. However we repeat the presentation specifically for dimension $1$-foliations and codimension $1$-foliations to make this part as self-contained as possible. We also remark the paper \cite{Mal98} on foliated deformations and unfoldings of Hopf surfaces with (regular) holomorphic foliations which are important examples in deformation theory of complex structures.

\section{Infinitesimal foliated deformations and obstructions of dimension 1 foliations in terms of deformations of sections}

\subsection{Infinitesimal deformations of foliated complex analytic structures of dimension $1$ in terms of sections} \label{ot30}\

Let $\left(M, \Theta_{\mathcal{F}_0}\right)$ be a compact foliated complex manifold with $\Theta_{\mathcal{F}_0}$ locally free of rank of $1$, so that there exists a section $T_0\in H^0\left( M, \Theta_M\otimes {\mathcal{L}_0} \right)$ with $\mathcal{L}_0 =\Theta_{\mathcal{F}_0}^*$, which defines the foliation on $M$. In this subsection, we will consider foliated deformations of $\left( M, \Theta_{\mathcal{F}_0} \right)$ as deformations of the section $T_0$ (see Part II). Then we have the complex of sheaves controlling foliated deformations of $\left(M, \Theta_{\mathcal{F}_0}\right)$ in terms of tangent sheaves (equivalently in terms of sections)
\begin{align}\label{ap3}
\mathcal{E}_{\mathcal{L}_0}^\bullet : \mathcal{E}_{\mathcal{L}_0} \xrightarrow{\tilde{D}_0 }   \Theta_M\otimes {\mathcal{L}_0} \to 0 \to 0 \to \cdots
\end{align}
We recall that $\mathcal{E}_{\mathcal{L}_0}$ and the differential $\tilde{D}_0$ are defined in the following way: for an open covering $\mathcal{U}_0=\left\{U_j \right\}$ of $M$ by coordinate neighborhoods, let $\mathcal{L}_0$ defined by the transition functions $\left\{ \Psi_{0jk} \right\}\in C^1 \left(\mathcal{U}_0, \mathcal{O}_M^* \right)$ with $\Psi_{0ij}\Psi_{0jk}=\Psi_{0ik}$, and $T_0\in \Gamma\left( M, \Theta_M\otimes \mathcal{L}_0 \right)$ be written as $T_{0j}\in \Gamma\left( U_j, \Theta_M \right)$  on $U_j$ with $T_{0j}= \Psi_{0jk} T_{0k}$ on $U_j\cap U_k$.  Then $\mathcal{E}_{\mathcal{L}_0} |_{U_j}$ is isomorphic to $\Theta_{M}|_{U_j}\oplus \mathcal{O}_{M}|_{U_j}$, and $(\rho_j, a_j)\in \Theta_{M}|_{U_j}\oplus \mathcal{O}_{M}|_{U_j}$ is identified with $\left(\rho_k, a_k\right)\in \Theta_{M}|_{U_k}\oplus \mathcal{O}_{M}|_{U_k}$ on $U_j\cap U_k$ if and only if $\rho_j=\rho_k$ and $a_k-a_j=\left[\rho_j, \log \Psi_{0jk}\right]$. On the other hand, the differential $\tilde{D}_0$ is defined in the following way:
\begin{align*}
\tilde{D}_0: \Gamma\left(U_j,\mathcal{E}_{\mathcal{L}_0}\right)\cong \Gamma\left(U_j,\Theta_{M}\right)\oplus \Gamma\left(U_j ,  \mathcal{O}_{M} \right) &\to \Gamma\left(U_j, \Theta_{M}\otimes  \mathcal{L}_0 \right)\\
     \left(\rho_j, a_j\right)&\mapsto \left[\rho_j, T_{0j}\right]+a_j T_{0j}
\end{align*}
We will denote the $i$-th cohomology group of $\mathcal{E}_{\mathcal{L}_0}^\bullet$ by $\mathbb{H}^i\left( M, \mathcal{E}_{\mathcal{L}_0}^\bullet \right)$. 

We shall identify infinitesimal foliated deformations of dimension $1$ foliations in terms of sections.  Let $\left(\mathcal{M}, \Theta_{\mathcal{F}}, B, \pi \right)$ with $\Theta_{\mathcal{F}}$ locally free of rank $1$ be a foliated complex analytic family of deformations of $\left( M, \Theta_{\mathcal{F}_0} \right)=\pi^{-1}(0), 0 \in B \subset \mathbb{C}^m$. Then the exact sequence $0 \to \Theta_\mathcal{F} \to \Theta_{\mathcal{M}/B} \to \frac{\Theta_{\mathcal{M}/B}}{\Theta_{\mathcal{F}}}\to 0 $ gives a section $T\in H^0\left(\mathcal{M}, \Theta_{\mathcal{M}/B}\otimes \mathcal{L} \right)$ with $\mathcal{L}= \Theta_{\mathcal{F}}^*$ which induces $T_0\in H^0\left( M, \Theta_M\otimes \mathcal{L}_0 \right)$. We will also denote the family by $\left( \mathcal{M}, T, B, \pi \right)$. Take a small disk $\Delta$ with centre $0$ such that $0\in \Delta \subset B$, and represent $\mathcal{M}_\Delta=\pi^{-1}(\Delta)$ in the form $\mathcal{M}_\Delta = \bigcup_{j=1}^l U_j\times \Delta$ where each $U_j$ is a polydisk, and $(\xi_j,t)\in U_j\times \Delta$ and $(\xi_k,t)\in U_k\times \Delta$ are the same points on $\mathcal{M}_\Delta$ if $\xi_j^\alpha= f_{jk}^\alpha(\xi_k,t)$ for $\alpha=1,...,n$. Here each $f_{jk}^\alpha(\xi_k,t)=f_{jk}^\alpha(\xi_k^1,..., \xi_k^n, t)$ is a holomorphic function of $\xi_k^1,..., \xi_k^n,t$ defined on $\left(U_k\times \Delta \right) \cap \left( U_j\times \Delta \right)\ne \emptyset$, and we have 
\begin{align}\label{ot7}
f_{ik}^\alpha(\xi_k,t)&=f_{ij}^\alpha\left(f_{jk}(\xi_k,t),t \right),\,\,\,\alpha=1,...,n  \,\,\,\,\,\,\textnormal{on}\,\,\, \left(U_i\times \Delta \right)\cap \left(U_j\times \Delta \right) \cap \left( U_k \times \Delta \right)\ne \emptyset
\end{align}
 On the other hand, the line bundle $\mathcal{L}$ on $\mathcal{M}_\Delta$ is defined by the system of transition functions $\left\{\Psi_{jk}(\xi_k,t) \right\}$ where $\Psi_{jk}(\xi_k,t)\in   \Gamma\left(\left(U_j\times \Delta \right) \cap \left( U_k\times \Delta \right), \mathcal{O}_\mathcal{M}^* \right)$ satisfying the cocycle condition 
\begin{align}\label{ot1}
\Psi_{ij}\left(f_{jk}(\xi_k,t), t \right)\cdot \Psi_{jk}(\xi_k,t) = \Psi_{ik}(\xi_k,t), \,\,\,\,\,\textnormal{on}\,\,\,\,\, \left(U_i\times \Delta\right) \cap \left(U_j \times \Delta \right) \cap \left( U_k \times \Delta \right) \ne \emptyset
\end{align}
and $T\in H^0\left( \mathcal{M}, \Theta_{\mathcal{M}/B}\otimes \mathcal{L} \right)$ on $U_j\times \Delta$ is of the form
\begin{align*}
T_j(\xi_j,t):= T|_{U_j\times \Delta}= \sum_{\alpha=1}^n T_j^\alpha(\xi_j,t)\frac{\partial}{\partial \xi_j^\alpha}\,\,\,\,\,\,\,\textnormal{on}\,\,\, U_j\times \Delta
\end{align*}
such that
\begin{align*}
T_j(\xi_j,t)= \Psi_{jk}(\xi_k,t) T_k(\xi_k,t)\,\,\,\,\,\textnormal{on}\,\,\,\left(U_j\times \Delta \right) \cap \left( U_k\times \Delta \right)
\end{align*}
Equivalently we have
\begin{align}
&\sum_{\alpha=1}^n T_j^\alpha(\xi_j,t)\frac{\partial}{\partial \xi_j^\alpha}= \Psi_{jk}(\xi_k,t)\sum_{\beta=1}^n T_k^\beta(\xi_k,t)\frac{\partial}{\partial \xi_k^\beta}= \Psi_{jk}(\xi_k,t)\sum_{\beta, \alpha=1}^n T_k^\beta(\xi_k,t)\frac{\partial f_{jk}^\alpha}{\partial \xi_k^\beta}\frac{\partial}{\partial \xi_j^\alpha} \notag\\
& \iff T_j^\alpha \left(f_{jk}(\xi_j,t),t\right)= \Psi_{jk}(\xi_k,t)\sum_{\beta=1}^n T_k^\beta(\xi_k,t)\frac{\partial f_{jk}^\alpha}{\partial \xi_k^\beta}(\xi_k,t) \label{ot3}
\end{align}

Now we identify the infinitesimal foliated deformations. By taking the derivative of $(\ref{ot7})$ with respect to $t$ and setting $\theta_{jk}(t)=\sum_{\alpha=1}^n \theta_{jk}^\alpha(\xi_j,t)\frac{\partial}{\partial \xi_j^\alpha}$ where $\theta_{jk}^\alpha(\xi_j,t)=\frac{\partial f_{jk}(\xi_k,t)}{\partial t},\xi_k=f_{kj}(\xi_j,t)$, we have
\begin{align}\label{ot22}
\theta_{ik}^\alpha(\xi_i,t)=\theta_{ij}^\alpha(\xi_i,t) + \sum_{\beta=1}^n \frac{\partial \xi_i^\alpha}{\partial \xi_j^\beta}\theta_{jk}^\beta(\xi_j,t),\,\,\,\alpha=1,...,n
\end{align}
By setting $t=0$, we obtain 
\begin{align} \label{ot8}
\theta_{ij}(0)-\theta_{ik}(0) + \theta_{jk}(0)=0
\end{align}
By taking the derivative of $(\ref{ot1})$ with respect to $t$, we have
\begin{align}\label{ot2}
\sum_{\alpha=1}^n \frac{\partial \Psi_{ij}}{\partial \xi_j^\alpha}\left(f_{jk}(\xi_k,t),t\right )\frac{\partial f_{jk}^\alpha}{\partial t} \Psi_{jk} + \frac{\partial \Psi_{ij}}{\partial t} \left(f_{jk}(\xi_k,t),t\right)\Psi_{jk}+ \Psi_{ij}\left(f_{jk}(\xi_k,t),t\right) \frac{\partial \Psi_{jk}}{\partial t} = \frac{\partial \Psi_{ik}}{\partial t}
\end{align} 
By dividing $(\ref{ot2})$ by $\Psi_{ik}$, we have
\begin{align} \label{ot14}
  \frac{1}{\Psi_{ij}}\sum_{\alpha=1}^n \frac{\partial f_{jk}^\alpha}{\partial t} \frac{\partial \Psi_{ij}}{\partial \xi_j^\alpha} + \frac{1}{\Psi_{ij}}\frac{\partial \Psi_{ij}}{\partial t} + \frac{1}{\Psi_{jk}}\frac{\partial \Psi_{jk}}{\partial t} =  \frac{1}{\Psi_{ik}} \frac{\partial \Psi_{ik}}{\partial t}
\end{align}
By setting $t=0$, we obtain
\begin{align}\label{ot4}
 -\frac{1}{\Psi_{0ij}}\frac{\partial \Psi_{ij}}{\partial t}|_{t=0} + \frac{1}{\Psi_{0ik}} \frac{\partial \Psi_{ik}}{\partial t} |_{t=0} - \frac{1}{\Psi_{jk}}\frac{\partial \Psi_{jk}}{\partial t } - \left[ \sum_{\alpha=1}^n \frac{\partial f_{jk}^\alpha}{\partial t}|_{t=0} \frac{\partial}{\partial \xi_j^\alpha} , \log \Psi_{0ij} \right] =0
\end{align}
where $\Psi_{0jk}:=\Psi_{jk}\left( \xi_k, 0\right)$, i.e. $\left\{ \Psi_{0jk} \right\}$ defines the $\mathcal{L}_0$ on $M$.

Next we take the derivative of $(\ref{ot3})$ with respect to $t$. Then we have
\begin{align}\label{ot6}
\sum_{\beta=1}^n \frac{\partial T_j^\alpha}{\partial \xi_j^\beta}\frac{\partial f_{jk}^\beta}{\partial t} + \frac{\partial T_j^\alpha}{\partial t}=  \frac{\partial \Psi_{jk}}{\partial t}\sum_{\beta=1}^n T_k^\beta \frac{\partial f_{jk}^\alpha}{\partial \xi_k^\beta} + \Psi_{jk} \sum_{\beta=1}^n \frac{\partial T_k^\beta}{\partial t} \frac{\partial f_{jk}^\alpha}{\partial \xi_k^\beta} + \Psi_{jk} \sum_{\beta=1}^n T_k^\beta \frac{\partial}{\partial \xi_k^\beta}\left( \frac{\partial f_{jk}^\alpha}{\partial t} \right)
\end{align}
which implies that
{\small{\begin{align}\label{ot11}
&\left[ \sum_{\alpha=1}^n \frac{\partial f_{jk}^\alpha}{\partial t} \frac{\partial}{\partial \xi_j^\alpha}, \sum_{\beta=1}^n T_j^\beta \frac{\partial}{\partial \xi_j^\beta} \right] - \frac{1}{\Psi_{jk}}\frac{\partial \Psi_{jk}}{\partial t} \sum_{\alpha=1}^n T_j^\alpha \frac{\partial}{\partial \xi_j^\alpha}= \sum_{\alpha,\beta=1}^n \frac{\partial f_{jk}^\beta}{\partial t} \frac{\partial T_j^\alpha}{\partial \xi_j^\beta} \frac{\partial}{\partial \xi_j^\alpha} - \sum_{\alpha,\beta=1}^n T_j^\beta \frac{\partial}{\partial \xi_j^\beta}\left(  \frac{\partial f_{jk}^\alpha}{\partial t}\right)\frac{\partial}{\partial \xi_j^\alpha} - \frac{1}{\Psi_{jk}}\frac{\partial \Psi_{jk}}{\partial t} \sum_{\alpha=1}^n T_j^\alpha \frac{\partial}{\partial \xi_j^\alpha}\\
&= - \sum_{\alpha=1}^n \frac{\partial T_j^\alpha}{\partial t} \frac{\partial}{\partial \xi_j^\alpha} + \Psi_{jk}\sum_{\beta=1}^n \frac{\partial T_k^\beta}{\partial t} \frac{\partial}{\partial \xi_k^\beta} \notag
\end{align}}}
By setting $t=0$ on $(\ref{ot11})$, we obtain
\begin{align} \label{ot5}
-\Psi_{0jk}\frac{\partial T_k}{\partial t}|_{t=0} + \frac{\partial T_j}{\partial t}|_{t=0} + \left[ \sum_{\alpha=1}^n \frac{\partial f_{jk}^\alpha}{\partial t} |_{t=0} \frac{\partial}{\partial \xi_j^\alpha}, \sum_{\beta=1}^n T_{0j}^\beta \frac{\partial}{\partial \xi_j^\beta} \right] - \frac{1}{\Psi_{0jk}}\frac{\partial \Psi_{jk}}{\partial t}|_{t=0} \cdot T_{0j}=0
\end{align}
where $T_{0j}= \sum_{\beta=1}^n T_{0j}^\beta \frac{\partial}{\partial \xi_j^\beta}:= \sum_{\beta=1}^n T_j^\beta(\xi_j, 0)\frac{\partial}{\partial \xi_j^\beta } $ and $\frac{\partial T_j}{\partial t}:= \sum_{\beta=1}^n \frac{\partial T_j(\xi_j, t)}{\partial t} \frac{\partial}{\partial \xi_j^\beta}$.

Then $(\ref{ot8}), (\ref{ot4})$ and $(\ref{ot5})$ implies that
\begin{align*}
\left( \left\{- \frac{\partial T_j}{\partial t}|_{t=0} \right\} , \left\{ \left( \sum_{\alpha=1}^n \frac{\partial f_{jk}^\alpha}{\partial t} \frac{\partial}{\partial \xi_j^\alpha}|_{t=0}, -\frac{1}{\Psi_{0jk}} \frac{\partial \Psi_{jk}}{\partial t} |_{t=0} \right)    \right\}    \right) \in C^0\left(\mathcal{U}_0,  \Theta_{M}\otimes \mathcal{L}_0 \right) \bigoplus C^1\left( \mathcal{U}_0, \mathcal{E}_{\mathcal{L}_0} \right)
\end{align*}
defines a $1$-cocycle in the following \v Cech resolution of $\mathcal{E}_{\mathcal{L}_0}^\bullet$ with respect to the covering $\mathcal{U}_0:=\mathcal{U}\cap M$:
\begin{equation}\label{ot19}
\begin{CD}
0 @. 0 @. 0 \\
@AAA @AAA @AAA \\
C^0\left(\mathcal{U}, \Theta_{M} \otimes \mathcal{L}_0 \right) @>\delta >> C^0\left( \mathcal{U}, \Theta_M\otimes \mathcal{L}_0 \right) @>-\delta >> C^2\left( \mathcal{U}, \Theta_M\otimes \mathcal{L}_0  \right) @>\delta >> \cdots  \\
@AAA @AAA @AAA \\
C^0\left( \mathcal{U}, \mathcal{E}_{\mathcal{L}_0} \right) @>-\delta>> C^1\left(\mathcal{U}, \mathcal{E}_{\mathcal{L}_0} \right) @>\delta >> C^2\left( \mathcal{U}, \mathcal{E}_{\mathcal{L}_0} \right) @>-\delta>>  \cdots
\end{CD}
\end{equation}
Then we have the foliated Kodaira-Spencer map of $\left( \mathcal{M}, T , B, \pi \right)$ at $t=0$ in terms of sections
\begin{align*}
\varphi_0: T_0 B &\to \mathbb{H}^1\left( M, \mathcal{E}_{\mathcal{L}_0}^\bullet \right) \\
 \frac{\partial}{\partial t} &\mapsto \left( \frac{\partial \left( \mathcal{M}, T \right)}{\partial t}\right)|_{t=0}:=\left( \left\{- \frac{\partial T_j}{\partial t}|_{t=0} \right\} , \left\{ \left( \sum_{\alpha=1}^n \frac{\partial f_{jk}^\alpha}{\partial t} \frac{\partial}{\partial \xi_j^\alpha}|_{t=0}, -\frac{1}{\Psi_{0jk}} \frac{\partial \Psi_{jk}}{\partial t} |_{t=0} \right)    \right\}   \right)
\end{align*}

\subsection{Obstruction in terms of \v Cech resolution of $\mathcal{E}_{\mathcal{L}_0}^\bullet$}\label{ot23} \

We assume that the foliated analytic family $\left(\mathcal{M}, T, B, \pi \right)$ with $T\in  H^0\left(\mathcal{M}, \Theta_{\mathcal{M}/B}\otimes \mathcal{L}  \right)$ in terms of sections is a $1$-parameter family, i.e. $0\in \Delta \subset B\subset \mathbb{C}$. As in \cite{Kod05} p.211 we sometimes write $\left( \frac{\partial}{\partial t}\right)_j$ instead of $\left(\frac{\partial}{\partial t} \right)$ in order to make explicit that $\frac{\partial}{\partial t}$ denotes the differentiation of a function of $\xi_j^1,..., \xi_j^n,t$ with respect to $t$. On $\left( U_j \times \Delta \right)\cap \left( U_k \times \Delta \right)\ne \emptyset$, we have the following equalities: $\frac{\partial}{\partial \xi_k^\alpha}=\sum_{\beta=1}^n \frac{\partial \xi_j^\beta}{\partial \xi_k^\alpha} \frac{\partial}{\partial \xi_j^\beta}$, where $\frac{\partial \xi_j^\beta}{\partial \xi_k^\alpha}=\frac{\partial f_{jk}^\beta(\xi_k,t)}{\partial \xi_k^\alpha}$, and
\begin{align}\label{ot12}
\left( \frac{\partial}{\partial t} \right)_k  = \sum_{s=1}^n \frac{\partial f_{jk}^s(\xi_k,t)}{\partial t}\frac{\partial}{\partial \xi_j^s} + \left( \frac{\partial}{\partial t} \right)_j = \sum_{s=1}^n \theta_{jk}^s(\xi_j,t)\frac{\partial}{\partial \xi_j^s} + \left( \frac{\partial}{\partial t} \right)_j 
\end{align}
By taking the derivative of $(\ref{ot22})$ with respect to $t$, and putting
\begin{align*}
\dot{\theta}_{ij}(t)= \sum_{\alpha=1}^n \dot{\theta}_{ij}^\alpha (\xi_i,t)\frac{\partial}{\partial \xi_i^\alpha} , \,\,\,\,\,\textnormal{where}\,\,\,\dot{\theta}_{ij}^\alpha(\xi_i,t)=\frac{\partial \theta_{ij}^\alpha(\xi_i,t)}{\partial t}
\end{align*}
we obtain (for the detail, see \cite{Kod05} p.211-213)
\begin{align}\label{ot16}
\dot{\theta}_{ij}(t) -\dot{\theta}_{ik}(t) + \dot{\theta}_{jk} (t) = \left[ \theta_{ij}(t), \theta_{jk}(t) \right]
\end{align}
We note that
\begin{align}\label{ot13}
\dot{\theta}_{jk}(t) = \sum_{\alpha=1}^n \frac{\partial }{\partial t_j}\left(\frac{\partial f_{jk}^\alpha}{\partial t_k}\right) \frac{\partial}{\partial \xi_j^\alpha} = \sum_{\alpha=1}^n \frac{\partial^2 f_{jk}}{\partial t^2}\frac{\partial}{\partial \xi_j^\alpha} - \sum_{\alpha,\beta=1}^n\frac{\partial f_{jk}^\beta}{\partial t} \frac{\partial}{\partial \xi_j^\beta}\left(\frac{\partial f_{jk}^\alpha}{\partial t}\right)\frac{\partial}{\partial \xi_j^\alpha}
\end{align}

Let us take the derivative of $(\ref{ot6})$ with respect to $t_j$ again, i.e. $\left(\frac{\partial}{\partial t}\right)_j$:
{\small{\begin{align*}
&\sum_{\beta=1}^n \frac{\partial}{\partial \xi_j^\beta}\left( \frac{\partial T_j^\alpha}{\partial t_j} \right)\frac{\partial f_{jk}^\beta}{\partial t_k} + \sum_{\beta=1}^n \frac{\partial T_j^\alpha}{\partial \xi_j^\beta} \dot{\theta}_{jk}^\beta + \frac{\partial^2 T_j}{\partial t^2} = \left(\frac{\partial}{\partial t}\right)_j\left( \frac{\partial \Psi_{jk}}{\partial t}\sum_{\beta=1}^n T_k^\beta \frac{\partial f_{jk}^\alpha}{\partial \xi_k^\beta} + \Psi_{jk} \sum_{\beta=1}^n \frac{\partial T_k^\beta}{\partial t} \frac{\partial f_{jk}^\alpha}{\partial \xi_k^\beta} +  \sum_{\beta=1}^n T_j^\beta \frac{\partial}{\partial \xi_j^\beta}\left( \frac{\partial f_{jk}^\alpha}{\partial t} \right)\right)\\
&= - \sum_{\delta, \beta=1}^n \frac{\partial f_{jk}^\delta}{\partial t}\frac{\partial}{\partial \xi_j^\delta}\left( \frac{\partial \Psi_{jk}}{\partial t} T_k^\beta \frac{\partial f_{jk}^\alpha}{\partial \xi_j^\beta} \right) - \sum_{\delta, \beta=1}^n \frac{\partial f_{jk}^\delta}{\partial t}\frac{\partial}{\partial \xi_j^\delta}\left( \Psi_{jk} \frac{\partial T_k^\beta}{\partial t} \frac{\partial f_{jk}^\alpha }{\partial \xi_k^\beta} \right) + \sum_{\beta=1}^n \frac{\partial T_j^\beta}{\partial t_j}\frac{\partial}{\partial \xi_j^\beta}\left( \frac{\partial f_{jk}^\alpha}{\partial t}\right) + \sum_{\beta=1}^n T_j^\beta \frac{\partial \dot{\theta}_{jk}^\alpha}{\partial \xi_j^\beta}\\
&+ \frac{\partial^2 \Psi_{jk}}{\partial t^2}\sum_{\beta=1}^n T_k^\beta \frac{\partial f_{jk}^\alpha}{\partial \xi_k^\beta} + 2 \frac{\partial \Psi_{jk}}{\partial t}\sum_{\beta=1}^n \frac{\partial T_k^\beta}{\partial t} \frac{\partial f_{jk}^\alpha}{\partial \xi_k^\beta} + \frac{\partial \Psi_{jk}}{\partial t}\sum_{\beta=1}^n T_k^\beta \frac{\partial}{\partial \xi_k^\beta}\left( \frac{\partial f_{jk}^\alpha}{\partial t} \right) + \Psi_{jk}\sum_{\beta=1}^n \frac{\partial^2 T_k^\beta}{\partial t^2} \frac{\partial f_{jk}^\alpha}{\partial\xi_k^\beta} + \Psi_{jk} \sum_{\beta=1}^n \frac{\partial T_k^\beta}{\partial t} \frac{\partial}{\partial \xi_j^\beta}\left( \frac{\partial f_{jk}^\alpha}{\partial t} \right)
\end{align*}}}
By summing $\sum_{\alpha=1}^n \frac{\partial}{\partial \xi_j^\alpha}$, we obtain
{\small{\begin{align*}
&\sum_{\alpha, \beta=1}^n \frac{\partial f_{jk}^\beta}{\partial t} \frac{\partial}{\partial \xi_j^\beta}\left( \frac{\partial T_j^\alpha}{\partial t} \right)\frac{\partial}{\partial \xi_j^\alpha} + \sum_{\alpha, \beta=1}^n \dot{\theta}_{jk}^\beta \frac{\partial T_j^\alpha}{\partial \xi_j^\beta } \frac{\partial}{\partial \xi_j^\alpha} + \frac{\partial^2 T_j}{\partial t^2}\\
&=- \sum_{\delta, \beta,\alpha=1}^n \theta_{jk}^\delta \frac{\partial}{\partial \xi_j^\delta}\left( \frac{\partial \Psi_{jk}}{\partial t} T_k^\beta \frac{\partial f_{jk}^\alpha}{\partial \xi_j^\beta} \right) \frac{\partial}{\partial \xi_j^\alpha} - \sum_{\delta, \beta, \alpha=1}^n \theta_{jk}^\delta \frac{\partial}{\partial \xi_j^\delta}\left( \Psi_{jk}\frac{\partial T_k^\beta}{\partial t} \frac{\partial f_{jk}^\alpha}{\partial \xi_j^\beta} \right)\frac{\partial}{\partial \xi_j^\alpha} + \sum_{\alpha,\beta=1}^n \frac{\partial T_j^\beta}{\partial t}\frac{\partial}{\partial \xi_j^\beta}\left( \frac{\partial f_{jk}^\alpha}{\partial t} \right)\frac{\partial}{\partial \xi_j^\alpha} + \sum_{\beta=1}^n T_j^\beta \frac{\partial \dot{\theta}_{jk}^\alpha}{\partial \xi_j^\beta}\\
&+ \frac{\partial^2 \Psi_{jk}}{\partial t^2} T_k + 2 \frac{\partial \Psi_{jk}}{\partial t} \frac{\partial T_k}{\partial t} + \frac{\partial \Psi_{jk}}{\partial t} \sum_{\alpha, \beta=1}^n T_k^\beta \frac{\partial}{\partial \xi_k^\beta}\left( \frac{\partial f_{jk}^\alpha}{\partial t}\right)\frac{\partial}{\partial \xi_j^\alpha} + \Psi_{jk} \frac{\partial^2 T_k}{\partial t^2} + \Psi_{jk}\sum_{\beta=1}^n \frac{\partial T_k^\beta}{\partial t} \frac{\partial }{\partial \xi_j^\beta}\left( \frac{\partial f_{jk}^\alpha}{\partial t} \right)\frac{\partial}{\partial \xi_j^\alpha}
\end{align*}}}
which is equivalent to
{\small{\begin{align}\label{ot9}
\left[ \theta_{jk}(t), \frac{\partial T_j}{\partial t_j} \right] + \left[ \dot{\theta}_{jk}(t), T_j\right] + \frac{\partial^2 T_j}{\partial t_j^2} = -\left[ \theta_{jk}(t), \frac{\partial \Psi_{jk}}{\partial t} T_k\right] - \left[ \theta_{jk}(t), \Psi_{jk}\frac{\partial T_k}{\partial t} \right] + \frac{\partial^2 \Psi_{jk}}{\partial t^2} T_k + 2 \frac{\partial \Psi_{jk}}{\partial t} \frac{\partial T_k}{\partial t} + \Psi_{jk} \frac{\partial^2 T_k}{\partial t^2} 
\end{align}}}

We also note that from $(\ref{ot12})$
{\small{\begin{align}\label{ot10}
\frac{\partial}{\partial t_j}\left( \frac{1}{\Psi_{jk}} \frac{\partial \Psi_{jk}}{\partial t_k}\right) & = - \frac{1}{\Psi_{jk}^2}\frac{\partial \Psi_{jk}}{\partial t} \frac{\partial \Psi_{jk}}{\partial t} + \frac{1}{\Psi_{jk}}\frac{\partial^2 \Psi_{jk}}{\partial t^2} + \sum_{\alpha=1}^n \frac{1}{\Psi_{jk}^2}\frac{\partial f_{jk}^\alpha}{\partial t}\frac{\partial \Psi_{jk}}{\partial \xi_j^\alpha} \frac{\partial \Psi_{jk}}{\partial t} - \sum_{\alpha=1}^n \frac{1}{\Psi_{jk}} \frac{\partial f_{jk}^\alpha}{\partial t} \frac{\partial}{\partial \xi_j^\alpha}\left( \frac{\partial \Psi_{jk}}{\partial t} \right)\\
&= - \frac{1}{\Psi_{jk}^2}\frac{\partial \Psi_{jk}}{\partial t} \frac{\partial \Psi_{jk}}{\partial t} + \frac{1}{\Psi_{jk}}\frac{\partial^2 \Psi_{jk}}{\partial t^2} - \left[ \theta_{jk}(t),  \frac{1}{\Psi_{jk}} \frac{\partial \Psi_{jk}}{\partial t} \right] \notag
\end{align}}}

\begin{lemma}
We have the equality
{\small{\begin{align}\label{ot18}
-\Psi_{jk}\frac{\partial^2 T_k}{\partial t_k^2} + \frac{\partial^2 T_j}{\partial t_j^2} +\left[ \dot{\theta}_{jk}, T_{j}  \right]  - \frac{\partial}{\partial t_j}\left( \frac{1}{\Psi_{jk}}\frac{\partial \Psi_{jk}}{\partial t_k} \right) T_{j}=\left[ \theta_{jk}, - \frac{\partial T_j}{\partial t_j}- \Psi_{jk}\frac{\partial T_k}{\partial t_k} \right] - \frac{1}{\Psi_{jk}}\frac{\partial \Psi_{jk}}{\partial t}\left(- \frac{\partial T_j}{\partial t_j}- \Psi_{jk}\frac{\partial T_k}{\partial t_k}  \right)
\end{align}}}
\end{lemma}

\begin{proof}
From $(\ref{ot9})$, it is sufficient to show that
\begin{align*}
&\frac{\partial}{\partial t_j}\left( \frac{1}{\Psi_{jk}}\frac{\partial \Psi_{jk}}{\partial t} \right)T_j = - \left[ \theta_{jk}, \frac{1}{\Psi_{jk}}\frac{\partial \Psi_{jk}}{\partial t} T_j \right] + \frac{1}{\Psi_{jk}}\frac{\partial^2 \Psi_{jk}}{\partial t^2} T_j + \frac{\partial \Psi_{jk}}{\partial t} \frac{\partial T_k}{\partial t}  - \frac{1}{\Psi_{jk}}\frac{\partial \Psi_{jk}}{\partial t} \frac{\partial T_j}{\partial t} \\
&= - \left[ \theta_{jk}, \frac{1}{\Psi_{jk}}\frac{\partial \Psi_{jk}}{\partial t} \right] T_j + \frac{1}{\Psi_{jk}}\frac{\partial \Psi_{jk}^2}{\partial t^2}  - \frac{1}{\Psi_{jk}}\frac{\partial \Psi_{jk}}{\partial t} \left( \left[ \theta_{jk}, T_j \right] -\Psi_{jk} \frac{\partial T_k}{\partial t} + \frac{\partial T_j}{\partial t} \right)\\
&=  - \left[ \theta_{jk}, \frac{1}{\Psi_{jk}}\frac{\partial \Psi_{jk}}{\partial t} \right] T_j + \frac{1}{\Psi_{jk}}\frac{\partial \Psi_{jk}^2}{\partial t^2}  - \frac{1}{\Psi_{jk}^2}\frac{\partial \Psi_{jk}}{\partial t} \frac{\partial \Psi_{jk}}{\partial t}
\end{align*}
which follows from $(\ref{ot10})$ and $(\ref{ot11})$.
\end{proof}

\begin{lemma}
We have the equality
{\small{\begin{align}\label{ot15}
&\left[ \theta_{ij},  - \frac{1}{\Psi_{jk}}\frac{\partial \Psi_{jk}}{\partial t} - \left[\theta_{jk}, \log \Psi_{ij} \right]   \right] - \left[ \theta_{jk}  , - \frac{1}{\Psi_{ij}}\frac{\partial \Psi_{ij}}{\partial t}  \right]\\
&=-\frac{\partial}{\partial t_i}\left( \frac{1}{\Psi_{ij}} \frac{\partial \Psi_{ij}}{\partial t} \right) +\frac{\partial}{\partial t_i}\left(\frac{1}{\Psi_{ik}} \frac{\partial \Psi_{ik}}{\partial t}\right) - \frac{\partial}{\partial t_j}\left( \frac{1}{\Psi_{jk}}\frac{\partial \Psi_{jk}}{\partial t} \right)  - \left[  \dot{\theta}_{jk} , \log \Psi_{ij} \right] \notag
\end{align}}}
\end{lemma}

\begin{proof}
From $(\ref{ot12}), (\ref{ot13})$ and $(\ref{ot14})$, we have
{\small{\begin{align*}
&\frac{\partial}{\partial t_i}\left( \frac{1}{\Psi_{ij}} \frac{\partial \Psi_{ij}}{\partial t} \right) - \frac{\partial}{\partial t_i}\left(\frac{1}{\Psi_{ik}} \frac{\partial \Psi_{ik}}{\partial t}\right) + \frac{\partial}{\partial t_j}\left( \frac{1}{\Psi_{jk}}\frac{\partial \Psi_{jk}}{\partial t} \right)  + \left[  \sum_{\alpha=1}^n \frac{\partial^2 f_{jk}}{\partial t^2}\frac{\partial}{\partial z_j^\alpha} - \sum_{\alpha,\beta=1}^n\frac{\partial f_{jk}^\beta}{\partial t} \frac{\partial}{\partial z_j^\beta}\left(\frac{\partial f_{jk}^\alpha}{\partial t}\right)\frac{\partial}{\partial z_j^\alpha} , \log \Psi_{ij} \right] \\
&= -\frac{\partial}{\partial t_i}\left( \frac{1}{\Psi_{jk}}\frac{\partial \Psi_{jk}}{\partial t} + \frac{1}{\Psi_{ij}}\sum_{\alpha=1}^n \frac{\partial f_{jk}^\alpha}{\partial t} \frac{\partial \Psi_{ij}}{\partial z_j^\alpha} \right) + \frac{\partial}{\partial t_j}\left( \frac{1}{\Psi_{jk}}\frac{\partial \Psi_{jk}}{\partial t} \right) + \sum_{\alpha=1}^n  \frac{1}{\Psi_{ij}}\frac{\partial^2 f_{jk}^\alpha}{\partial t^2} \frac{\partial \Psi_{ij}}{\partial z_j^\alpha} - \sum_{\alpha,\beta=1}^n \frac{1}{\Psi_{ij}} \frac{\partial f_{jk}^\beta}{\partial t} \frac{\partial}{\partial z_j^\beta}\left( \frac{\partial f_{jk}^\alpha}{\partial t}  \right) \frac{\partial \Psi_{ij}}{\partial z_j^\alpha}\\
&= - \frac{\partial}{\partial t_j}\left( \frac{1}{\Psi_{jk}}\frac{\partial \Psi_{jk}}{\partial t} + \frac{1}{\Psi_{ij}}\sum_{\alpha=1}^n \frac{\partial f_{jk}^\alpha}{\partial t} \frac{\partial \Psi_{ij}}{\partial z_j^\alpha} \right) + \sum_{\beta=1}^n \frac{\partial f_{ij}^\beta}{\partial t}\frac{\partial}{\partial z_i^\beta} \left( \frac{1}{\Psi_{jk}}\frac{\partial \Psi_{jk}}{\partial t} + \frac{1}{\Psi_{ij}}\sum_{\alpha=1}^n \frac{\partial f_{jk}^\alpha}{\partial t} \frac{\partial \Psi_{ij}}{\partial z_j^\alpha} \right) + \frac{\partial}{\partial t_j}\left( \frac{1}{\Psi_{jk}}\frac{\partial \Psi_{jk}}{\partial t} \right) \\
& + \sum_{\alpha=1}^n  \frac{1}{\Psi_{ij}}\frac{\partial^2 f_{jk}^\alpha}{\partial t^2} \frac{\partial \Psi_{ij}}{\partial z_j^\alpha} - \sum_{\alpha,\beta=1}^n \frac{1}{\Psi_{ij}} \frac{\partial f_{jk}^\beta}{\partial t} \frac{\partial}{\partial z_j^\beta}\left( \frac{\partial f_{jk}^\alpha}{\partial t}  \right) \frac{\partial \Psi_{ij}}{\partial z_j^\alpha}\\
&=\frac{1}{\Psi_{ij}^2}\frac{\partial \Psi_{ij}}{\partial t_j}\sum_{\alpha=1}^n \frac{\partial f_{jk}^\alpha}{\partial t} \frac{\partial \Psi_{ij}}{\partial z_j^\alpha} - \frac{1}{\Psi_{ij}}\sum_{\alpha=1}^n \frac{\partial f_{jk}^\alpha}{\partial t} \frac{\partial}{\partial z_j^\alpha}\left( \frac{\partial \Psi_{ij}}{\partial t_j} \right)  + \sum_{\beta=1}^n \frac{\partial f_{ij}^\beta}{\partial t}\frac{\partial}{\partial z_i^\beta} \left( \frac{1}{\Psi_{jk}}\frac{\partial \Psi_{jk}}{\partial t} + \frac{1}{\Psi_{ij}}\sum_{\alpha=1}^n \frac{\partial f_{jk}^\alpha}{\partial t} \frac{\partial \Psi_{ij}}{\partial z_j^\alpha} \right) 
\end{align*}}}
\end{proof}

We note that $(\ref{ot16})$ and $(\ref{ot15})$ imply that by setting $t=0$, 
{\small{\begin{align}\label{ot17}
\delta \left( \left\{ \dot{\theta}_{jk}(0), - \frac{\partial}{\partial t_j}\left( \frac{1}{\Psi_{jk}}\frac{\partial \Psi_{jk}}{\partial t_k} \right)|_{t=0}   \right\}\right)=\left\{ \left[ \left( \theta_{ij}(0), - \frac{1}{\Psi_{0ij}}\frac{\partial \Psi_{ij}}{\partial t}|_{t=0} \right), \left( \theta_{jk}(0), -\frac{1}{\Psi_{0jk}}\frac{\partial \Psi_{jk}}{\partial t}|_{t=0} \right)    \right]\right\}
\end{align}}}
and setting $t=0$ on $(\ref{ot18})$, we obtain
{\small{\begin{align}\label{ot20}
&-\Psi_{0jk}\frac{\partial^2 T_k}{\partial t_k^2}|_{t=0} + \frac{\partial^2 T_j}{\partial t_j^2}|_{t=0} +\left[ \dot{\theta}_{jk}(0), T_{0j}  \right]  - \frac{\partial}{\partial t_j}\left( \frac{1}{\Psi_{jk}}\frac{\partial \Psi_{jk}}{\partial t_k} \right)|_{t=0} T_{0j}\\
&=\left[ \theta_{jk}(0), - \frac{\partial T_j}{\partial t_j}|_{t=0} - \Psi_{0jk}\frac{\partial T_k}{\partial t_k}|_{t=0} \right] - \frac{1}{\Psi_{0jk}}\frac{\partial \Psi_{jk}}{\partial t}|_{t=0}\left(- \frac{\partial T_j}{\partial t_j}|_{t=0}- \Psi_{0jk}\cdot \frac{\partial T_k}{\partial t_k}|_{t=0}  \right) \notag
\end{align}}}

For any $1$-cocycle $\left(\left\{ a_j \right\}, \left\{ \left( b_{jk}, c_{jk} \right)\right\} \right)\in C^0\left(\mathcal{U}_0, \Omega_M^1\otimes \mathcal{L}_0 \right)\bigoplus C^1\left(\mathcal{U}_0, \mathcal{E}_{\mathcal{L}_0} \right)$ in the \v Cech resolution $(\ref{ot19})$ of $\mathcal{E}_{\mathcal{L}_0}^\bullet$, we can define a $2$-cocycle by (we omit the proof, see Part IV-1)
{\small{\begin{align}\label{ot21}
\left( \left\{ \left[ b_{jk}, a_j+ \Psi_{0jk} a_k \right] + c_{jk}\left(a_j + \Psi_{0jk} a_k \right)\right\}, \left\{ \left[ b_{ij}, c_{ij} \right], \left[ b_{jk}, c_{jk}\right] \right\}   \right) \in C^1\left(\mathcal{U}_0, \Omega_M^1\otimes \mathcal{L}_0 \right)\bigoplus C^2\left( \mathcal{U}_0, \mathcal{E}_{\mathcal{L}_0} \right)
\end{align}}}
Then $(\ref{ot17})$ and $(\ref{ot20})$ imply that
\begin{theorem}\label{ot40}
Suppose given a compact complex manifold $M$ with a dimension $1$-foliation defined by $T_0\in H^0\left(M, \Theta_M\otimes \mathcal{L}_0 \right)$ for a line bundle $\mathcal{L}_0$, and $\left(\left\{a_j \right\}, \left\{\left( b_{jk}, c_{jk} \right) \right\} \right)\in \mathbb{H}^1\left(M, \mathcal{E}_{\mathcal{L}_0}^\bullet \right)$. In order that there may exist a foliated analytic family $\left( \mathcal{M}, T, B, \pi\right)$ such that $\pi^{-1}(0)=\left( M, T_0 \right)$ and $\left(\frac{\partial \left(\mathcal{M}, T \right)}{\partial t}\right)|_{t=0}=\left(\left\{a_j \right\}, \left\{\left( b_{jk}, c_{jk} \right) \right\} \right)$, it is necessary  that $(\ref{ot21})$ defines the $0$ class in $\mathbb{H}^2\left( M, \mathcal{E}_{\mathcal{L}_0}^\bullet \right)$.

\end{theorem}

Keeping the notations above, we have the following lemma.

\begin{lemma} \label{rsf28}
A compact foliated complex manifold $\left(M, T_0\right)$ of dimension $1$ with $T_0\in H^0\left(M, \Theta_M\otimes \mathcal{L}_0 \right)$ is obstructed in foliated deformations in terms of tangent sheaves if there exists an element $\{a_j\} \in H^0\left( M, \Theta_M\otimes \mathcal{L}_0 \right)$ and an element $\{(b_{jk}, c_{jk})\} \in H^1\left(M, \mathcal{E}_{\mathcal{L}_0} \right)$ which is in $\textnormal{ker}\left(H^1\left( M, \mathcal{E}_{\mathcal{L}_0} \right) \xrightarrow{\tilde{E}_0} H^1\left( M, \Theta_M \otimes \mathcal{L}_0 \right) \right)$ such that one of the following holds:
\begin{enumerate}
\item Under the map
\begin{align}\label{ot41}
H^0\left( M, \Theta_M\otimes \mathcal{L}_0 \right) \times H^1\left( M, \mathcal{E}_{\mathcal{L}_0} \right) &\to H^1\left( M, \Theta_M \otimes \mathcal{L}_0  \right)\\
\left( \{a_j\}, \{(b_{jk}, c_{jk}\} \right)& \mapsto \left\{ [b_{jk}, a_j]+ c_{jk}a_j \right\} \notag
\end{align}
$ \left\{ [b_{jk}, a_j]+ c_{jk}a_j \right\}$ is not in the image of $H^1\left( M, \mathcal{E}_{\mathcal{L}_0} \right)\to H^1\left(M, \Theta_M \otimes \mathcal{L}_0  \right)$.
\item $\left\{[(b_{ij}, c_{ij}), (b_{jk}, c_{jk})] \right\}=\left\{[b_{ij}, b_{jk}], [b_{ij}, c_{jk}- [b_{jk}, \log \Psi_{0ij}]]- [b_{jk}, c_{ij}] \right\}$ is not $0$ class in $H^2\left( M, \mathcal{E}_{\mathcal{L}_0} \right)$
\item $\{[b_{ij}, b_{jk}]\}$ is not $0$ class in $H^2(M, \Theta_M)$.
\end{enumerate}
\end{lemma}

\begin{proof}
We show that the map $(\ref{ot41})$ is well-defined. In fact, since $\Psi_{ij}a_j= a_i$ and $b_{ij}- b_{ik}+ b_{jk}=0$ and $c_{ij}- c_{ik}+ c_{jk}- \left[ b_{jk}, \log \Psi_{0ij}\right]=0$, we have
\begin{align*}
&\left[ b_{ij}, a_i \right]  +c_{ij}a_i - \left[ b_{ik}, a_i \right]- c_{ik}a_i + \Psi_{0ij}\left( \left[ b_{jk}, a_j\right] +c_{jk}a_j \right)\\
&= \left[ b_{ij}- b_{ik} + b_{jk}, a_i \right]+ \left(c_{ij}- c_{ik}+ c_{jk} \right)a_i - \left[ b_{jk}, \Psi_{0ij} \right] a_j=0.
\end{align*}
Assume that $\left(  M, T_0 \right)$ is unobstructed in foliated deformations in terms of tangent sheaves. Since $\left\{\left(b_{jk}, c_{jk} \right)\right\}\in H^1\left( M, \mathcal{E}_{ \mathcal{L}_0} \right)$ is in the kernel of $\tilde{E}_0$, there exists $\left\{ r_j\right\}\in C^0\left(\mathcal{U}, \Omega_M^1\otimes \mathcal{L}_0 \right)$ such that $\left( \{r_j\}, \{(b_{jk}, c_{jk})\} \right)$ defines an element of $\mathbb{H}^1\left( M, \mathcal{E}_{\mathcal{L}_0}^\bullet \right)$. Then $\left( \left\{r_j+ a_j \right\}, \left((b_{jk}, c_{jk} \right)\right)$ also defines a $1$-cocycle. Then by Theorem \ref{ot40}, there exists $\left( \{g_i\}, \{(e_{jk}, h_{jk})\} \right)$ and $\left( \{g_i'\}, \{(e_{jk}', h_{jk}') \right)\in C^0\left(\mathcal{U}, T_X\otimes \mathcal{L}_0 \right)\bigoplus C^1\left(\mathcal{U}, \mathcal{E}_{\mathcal{L}_0} \right)$ such that
\begin{align}\label{ot411}
\delta\left(\{g_i \} \right) + \left\{\left[ e_{jk}, T_{0j}\right] + h_{jk} T_{0j} \right\}&= \left\{\left[ b_{jk}, r_j+ \Psi_{0jk}r_k \right]+ c_{jk}\left(r_j+ \Psi_{0jk}r_k \right) \right\}\\
\delta\left( \left\{ \left(e_{jk}, h_{jk}\right) \right\} \right)&= \left\{ \left[ (b_{ij}, c_{ij}), (b_{jk}, c_{jk})     \right] \right\} \notag
\end{align}
and
\begin{align}\label{ot42}
\delta\left(\{g_i' \} \right) + \left\{\left[ e_{jk}', T_{0j}\right] + h_{jk}' T_{0j} \right\}&= \left\{\left[ b_{jk}, r_j+ a_j+ \Psi_{0jk}(r_k+ a_k) \right]+ c_{jk}\left(r_j+ a_j +\Psi_{0jk}(r_k+ a_k) \right) \right\}\\
\delta\left( \left\{ \left(e_{jk}', h_{jk}' \right) \right\} \right)&= \left\{ \left[ (b_{ij}, c_{ij}), (b_{jk}, c_{jk})     \right] \right\} \notag
\end{align}
By subtracting $(\ref{ot411})$ from $(\ref{ot42})$, we have
\begin{align}
\delta\left( \{ g_i'- g_i \} \right)+ \left\{ \left[ e_{jk}'-e_{jk}, T_{0j}\right] + \left( h_{jk}' - h_{jk}\right) T_{0j} \right\} &= \left\{ 2\left[ b_{jk}, a_j \right]+ 2 c_{jk}a_j \right\} \label{ot43}\\
\delta\left(\left\{ \left( e_{jk}'-e_{jk}, h_{jk}' - h_{jk} \right) \right\} \right)&=0 \label{ot45}
\end{align}
We note that $(\ref{ot45})$ implies that $\left\{\left(e_{jk}'-e_{jk}, h_{jk}'-h_{jk} \right) \right\}\in H^1\left(M, \mathcal{E}_{\mathcal{L}_0} \right)$. Then $(\ref{ot43})$ implies that $\left\{ \left[ b_{jk}, a_j \right]  + c_{jk}a_j \right\}$ is cohomologous to the image of $H^1\left( M, \mathcal{E}_{\mathcal{L}_0} \right)\xrightarrow{\tilde{D}_0} H^1\left( M, \Theta_M \otimes \mathcal{L}_0 \right)$. This proves $(1)$.

\end{proof}

\subsection{Integrability conditions}\label{oi31}\

Let $(\mathcal{M}, T, B, \pi)$ be a foliated complex analytic family with $T\in H^0\left(\mathcal{M}, \Theta_{\mathcal{M}/B}\otimes \mathcal{L} \right)$, where $\mathcal{L}= \Theta_\mathcal{F}^* $ is a line bundle on $\mathcal{M}$ of deformations of $\pi^{-1}(0)=\left(M, T_0 \right)$. We keep the notations in subsection \ref{ot30}. By \cite{Kod05} Theorem 2.3, $M_t$ is diffeomorphic to $M_0=\pi^{-1}(0)$ as differentiable manifolds for each $t\in \Delta$. By \cite{Kod05} Theorem 2.5, if we take a sufficiently small $\Delta$, there is a diffeomorphism $\Phi$ of $M\times \Delta$ onto $\mathcal{M}_\Delta$ as differentiable manifolds such that $\pi\circ \Phi$ is the projection $M\times \Delta \to \Delta$. Let $z=\left(z^1,..., z^n \right)$ be local coordinates of $M=M_0$. Then we have $\pi\circ \Phi(z,t)=t, t\in \Delta$. For $\Phi(z,t)\in U_j\times \Delta$, put
\begin{align*}
\Phi\left( z, t\right) = \left(\xi_j^1(z,t), ..., \xi_j^n(z,t), t_1,..., t_m \right)
\end{align*}
Then each component $\xi_j^\alpha=\xi_j^\alpha(z,t),\alpha=1,...,n$ is a $C^\infty$ function. If we identify $\mathcal{M}_\Delta =\Phi\left(M\times \Delta \right)$ with $M\times \Delta$, and $\mathcal{L}$ with $\Phi^*\mathcal{L}$ via $\Phi$, $\left(\mathcal{M}_\Delta, \mathcal{L}_\Delta \right)$ is considered as a complex manifold with the complex structure defined on the $C^\infty$ manifold $M\times \Delta$ by the system of local coordinates on $U_j\times \Delta$
\begin{align*}
\left\{ \left(\xi_j,t \right)| j=1 ,2, 3,...   \right\},\,\,\,\,\,(\xi_j,t)=\left( \xi_j^1(z,t), ..., \xi_j^n(z,t), t_1,..., t_m \right)
\end{align*}
and the $C^\infty$ differentiable line bundle $\Phi^* \mathcal{L}_\Delta$ on $M\times \Delta$.

We note that since $(z^1,..., z^n)$ and $\left(\xi_j^1(z,0), ..., \xi_j^n(z,0)\right)$ are local complex coordinates on $M=M_0$, $\xi_j^\alpha(z,0)$ are holomorphic functions of $z^1,...,z^n,\alpha=1,...,n$. We also note that if we take $\Delta$ sufficiently small, we have
\begin{align}\label{oi3}
\det\left( \frac{\partial \xi_j^\alpha(z,t)}{\partial z^\lambda} \right)_{\alpha, \lambda=1,...,n} \ne 0\,\,\,\,\,\, \, \textnormal{for}\,\,\,t\in \Delta
\end{align}
On the other hand, each line bundle $\Phi^* \mathcal{L}_t$ in a differentiable  family $\left\{ \Phi^* \mathcal{L}_t |t\in \Delta \right\}$ is the restriction to $M\times t$ of the line bundle $\Psi^* \mathcal{L}$ over $M\times \Delta , \Phi^* \mathcal{L}_t = \Phi^* \mathcal{L} |_{M\times t}$. Then by \cite{Kod05} Lemma 7.1 the $C^\infty$ line bundle $\Phi^* \mathcal{L}$ over $M\times \Delta$ is $C^\infty$-equivalent to the trivial line bundle $\mathcal{L}_0 \times \Delta $. In other words, there is a bundle diffeomorphism $\sigma: \mathcal{L}_0\times \Delta \to \Phi^* \mathcal{L}$ over $M\times \Delta$. We recall that $\mathcal{L}$ is described by the transition functions $\left\{\Psi_{jk}(\xi_j,t)\right\}$ with respect to the covering $\left\{U_j\times \Delta \right\}$ of $\mathcal{M}_\Delta$. Then $\Phi^* \mathcal{L}$ is described by the transition matrices $\left\{R_{jk}:= \Psi_{jk}\left(\xi_k(z,t) , t\right) \right\}$ with respect to the covering $\left\{\mathscr{U}_j:=\Psi^{-1}\left( U_j\times \Delta \right) \right\}$ of $M\times \Delta$. In particular $\mathcal{L}_0$ is described by $\left\{ \Psi_{0jk}:=\Psi_{jk}(\xi_j(z,0),0 )\right\}$. Then there are non-vanishing $C^\infty$ functions $c_j(z,t)$ on $\mathscr{U}_j$ with $c_j(z,0)=1$ such that $\sigma$ is represented locally by
\begin{align*}
\sigma:\mathbb{C}\times\mathscr{U}_j &\to \mathbb{C} \times \mathscr{U}_j\\
                         \left( X, z, t \right) &\mapsto \left( c_j(z,t)\cdot X, z, t\right)
\end{align*}
such that
\begin{align}\label{oti1}
c_j(z,t) \cdot   \Psi_{jk}\left(\xi_k(z,0),0 \right) =  \Psi_{jk}\left(\xi_k(z,t), t \right) \cdot c_k(z,t) \iff c_j \cdot \Psi_{0jk} = R_{jk}\cdot c_k
\end{align}
From $(\ref{oi3})$, we can define a $(0,1)$-form $\varphi(t):=\sum_{\lambda,\lambda=1}^n \varphi_{jv}^\lambda(z,t)d\bar{z}^v\frac{\partial}{\partial z^\lambda}\in A^{0,1}\left(M, \Theta_M \right)$ for $t\in \Delta$ such that $\bar{\partial}\xi_j^\alpha - \sum_{\lambda=1}^n \varphi_{jv}^\lambda(z,t)\frac{\partial \xi_j^\alpha}{\partial z_\lambda}d\bar{z}^v$ and
\begin{align}
\varphi(0)=0,\,\,\,\,\,\,\, \bar{\partial} \varphi(t)- \frac{1}{2} \left[ \varphi(t), \varphi(t) \right]=0
\end{align}
(see \cite{Kod05} p.262). On the other hand, since $\bar{\partial} R_{jk} -\left[ \varphi, R_{jk} \right]=0 $, from $(\ref{oti1})$
\begin{align*}
&\bar{\partial} R_{jk}\cdot c_k + R_{jk} \bar{\partial } c_k = \bar{\partial} c_j\cdot \Psi_{0jk}\iff \left[\varphi,  R_{jk}\right] \cdot c_k + R_{jk}\bar{\partial } c_k =\bar{\partial} c_j \cdot \Psi_{0jk}\\
&\iff \left[\varphi, c_j \cdot \Psi_{0jk} \cdot c_k^{-1}\right]c_k + R_{jk} \bar{\partial} c_k = \bar{\partial} c_j\cdot \Psi_{0jk} \\
&\iff \left[\varphi, c_j\right] \Psi_{0jk} + c_j [\varphi, \Psi_{0jk}] - c_j \Psi_{0jk} c_k^{-1}\left[\varphi, c_k \right]+ c_j \Psi_{0jk} c_k^{-1} \bar{\partial} c_k = \bar{\partial} c_j\cdot \Psi_{0jk}\\
&\iff [\varphi, c_j] + c_j[\varphi, \log \Psi_{0jk}]- c_j [\varphi, \log c_k] + c_j c_k^{-1}\bar{\partial} c_k = \bar{\partial} c_j\\
& \iff [\varphi, \log c_j]+ [\varphi, \log \Psi_{0jk}] - [\varphi, \log c_k]+ \bar{\partial} \log c_k= \bar{\partial} \log c_j
\end{align*}
Then we have
\begin{align}\label{oi30}
[\varphi, \log c_j]- \bar{\partial}\log c_j= [\varphi, \log c_k]- \bar{\partial} \log c_k- [ \varphi , \log \Psi_{0jk}]
\end{align}

This implies that
\begin{align*}
\left( (\varphi, U_j:= \left\{\left[\varphi, \log c_j\right] - \bar{\partial} \log c_j \right\} \right) \in A^{0,1}\left( M, \mathcal{E}_{\mathcal{L}_0} \right)
\end{align*}

Then
\begin{align*}
\bar{\partial} U_j &= \bar{\partial} \left[ \varphi, \log c_j \right] = \left[\bar{\partial} \varphi, \log c_j\right]  - \left[\varphi, \bar{\partial} \log c_j \right]  =  \frac{1}{2}\left[\left[\varphi, \varphi \right], \log c_j \right] - \left[\varphi, \bar{\partial} \log c_j \right]  \\
&= \left[ \varphi, \left[\varphi, \log c_j \right] \right]-  \left[\varphi, \bar{\partial} \log c_j \right] = \left[\varphi, \left[\varphi, \log c_j \right]- \bar{\partial} \log c_j \right] = [\varphi, U_j]
\end{align*}
Hence we have
\begin{align}
\bar{\partial} U_j - \left[\varphi, U_j\right] = 0 
\end{align}

On the other hand, from $(\ref{oi3})$ there exists $C^\infty$ functions $V_j^\alpha$ such that $\sum_{\alpha=1}^n  V_j^\alpha(z,t)\frac{\partial \xi_j^\beta}{\partial z^\alpha} = T_j^\beta(\xi_j(z,t),t)$. By setting $V_j(z,t):=\sum_{\alpha=1}^n V_j^\alpha(z,t)\frac{\partial}{\partial z^\alpha}$, we have
{\small{\begin{align*}
&\sum_{\alpha=1}^n \left(V_j^\alpha(z,t)- R_{jk}(z,t)V_k^\alpha(z,t)\right)\frac{\partial \xi_j^\beta}{\partial z^\alpha} = T_j^\beta\left(\xi_j(z,t),t) \right)- \sum_{\alpha=1}^n \Psi_{jk}\left(\xi_k(z,t),t) \right)  V_k^\alpha (z,t)\frac{\partial f_{jk}^{\beta}(\xi_k(z,t),t) }{\partial z^\alpha}\\
&= T_j^\beta\left(\xi_j(z,t),t) \right) - \sum_{\alpha, \gamma=1}^n \Psi_{jk}\left(\xi_k(z,t),t\right) V_k^\alpha(z,t)\frac{\partial f_{jk}^\beta}{\partial \xi_j^\gamma} \frac{\partial \xi_j^\gamma}{\partial z^\alpha} = T_j^\beta\left(\xi_j(z,t),t) \right) - \sum_{\gamma=1}^n \Psi_{jk}\left(\xi_k(z,t),t \right) T_k^\gamma\left(\xi_k(z,t),t \right)\frac{\partial f_{jk}^\beta}{\partial \xi_j^\gamma}=0
\end{align*}}}
This implies that
\begin{align*}
V_j(z,t)= R_{jk}(z,t) V_k(z,t) \Longrightarrow V_j(z,t)= c_j \cdot \Psi_{0jk}\cdot c_k^{-1} V_k(z,t) \Longrightarrow c_j^{-1} V_j = \Psi_{0jk}\cdot c_k^{-1}V_k
\end{align*}
This implies that
\begin{align*}
\left\{ X_j(t):= c_j^{-1} V_j   \right\} \in A^{0,0}\left( M, \Theta_M  \otimes \mathcal{L}_0 \right)
\end{align*}
and $X_j(0)= T_{0j}$ where $T_0= \left\{T_{0j} \right\} \in H^0\left( M, \Theta_M\otimes \mathcal{L}_0 \right)$. Since $\bar{\partial} V_j(z,t) - [\varphi, V_j(z,t)]=0$, we have
\begin{align*}
\bar{\partial}X_j- \left[\varphi, X_j\right]=\bar{\partial}\left( c_j^{-1} V_j \right) - \left[\varphi, c_j^{-1}V_j \right]= - \frac{1}{c_j^2}\bar{\partial} c_j\cdot V_j + \frac{1}{c_j^2}\left[\varphi, c_j \right]V_j= \left( - \bar{\partial} \log c_j + [\varphi, \log c_j] \right)X_j
\end{align*}
Hence we have
\begin{align*}
\bar{\partial} X_j - [\varphi, X_j] - U_j X_j=0
\end{align*}

All in all,  we have $\left(\{X_j\}, (\varphi, \{U_j\} \right) \in A^{0,0}\left( M, \Theta_M\otimes \mathcal{L} \right) \bigoplus A^{0,1}\left(M,\mathcal{E}_\mathcal{L}\right)$, where
\begin{align*}
&\varphi (t)\in A^{0,1}(M, \Theta_M)\\
&X_j (t) \in \Gamma\left( U_{j}, \mathcal{A}^{0,0}(\Theta_M) \right)\,\,\,\,\,\,\,\textnormal{with}\,\,\,\,\, X_j= \Psi_{0jk}X_k\\
&U_j(t) \in \Gamma\left( U_{j}, \mathcal{A}^{0,1} \right)\,\,\,\,\,\,\,\textnormal{with}\,\,\,\,\,\, U_j- U_k+ [\varphi , \log \Psi_{0jk}]=0
\end{align*}
with respect to the covering $\mathcal{U}_0=\left\{\mathcal{U}_j\cap M \right\}=\{U_{j}\}$ such that $\varphi(0)=0, U_j(0)=0$ and $X_j(0)= T_{0j}$ and
\begin{align}
\bar{\partial} \varphi - \frac{1}{2}\left[\varphi, \varphi\right]=0 \label{oti2}\\
\bar{\partial} X_j- \left[\varphi, X_j\right] - U_j   X_j=0 \label{oti3}\\
\bar{\partial} U_j - [\varphi, U_j]=0 \label{oti4}
\end{align}

\subsection{Infinitesimal deformations and obstructions in terms of Dolbeault resolution of $\mathcal{E}_{\mathcal{L}_0}^\bullet$ }\

We take the derivative of $(\ref{oti2})-(\ref{oti4})$ with respect to $t$:
\begin{align}
&\bar{\partial} \left( \frac{\partial \varphi}{\partial t} \right) - \frac{1}{2}\left[ \frac{\partial \varphi}{\partial t} , \varphi \right] - \frac{1}{2}\left[ \varphi, \frac{\partial \varphi}{\partial t} \right] =0 \label{oti5} \\
& \bar{\partial}\left( \frac{\partial X_j}{\partial t} \right) - \left[ \frac{\partial \varphi}{\partial t}, X_j \right] - \left[\varphi, \frac{\partial X_j}{\partial t} \right] - \frac{\partial U_j}{\partial t} \cdot X_j - U_j \frac{\partial X_j}{\partial t} = 0 \label{oti6}\\
& \bar{\partial}\left( \frac{\partial U_j}{\partial t} \right)  - \left[  \frac{\partial \varphi}{\partial t}, U_j    \right] - \left[ \varphi, \frac{\partial U_j}{\partial t}   \right] = 0 \label{oti7}
\end{align}
By setting $t=0$ of $(\ref{oti5})-(\ref{oti7})$, we have
\begin{align*}
&-\bar{\partial}\left( \frac{\partial \varphi}{\partial t}|_{t=0}  \right) = 0\\
&-\bar{\partial}\left( \frac{\partial X_j}{\partial t}|_{t=0} \right) + \left[ \frac{\partial \varphi}{\partial t}|_{t=0}, T_{0j} \right]  + \frac{\partial U_j}{\partial t}|_{t=0} \cdot T_{0j} = 0 \\
& -\bar{\partial}\left( \frac{\partial U_j}{\partial t}|_{t=0} \right)  = 0
\end{align*}
This implies that
\begin{align*}
\left( \left\{ \frac{\partial X_j}{\partial t}|_{t=0}  \right\}, \left\{ \left( \frac{\partial \varphi}{\partial t}|_{t=0},  \frac{\partial U_j}{\partial t}|_{t=0} \right) \right\} \right) \in A^{0,0}\left(M, \Theta_M\otimes \mathcal{L} \right) \bigoplus A^{0,1}\left( M, \mathcal{E}_{\mathcal{L}_0} \right) 
\end{align*}
defines a $1$-cocycle in the following dolbeault resolution of $\mathcal{E}_{\mathcal{L}_0}^\bullet$.
{\small{\begin{equation}\label{oti9}
\begin{CD}
0 @. 0 @. 0 \\
@AAA @AAA @AAA \\
A^{0,0} \left( M, \Theta_M\otimes \mathcal{L}_0 \right) @>-\bar{\partial}>> A^{0,1}\left( M, \Theta_M \otimes \mathcal{L}_0 \right) @>\bar{\partial}>> A^{0,2}\left( M, \Theta_M \otimes \mathcal{L}_0 \right) @>>> \cdots  \\
@A\tilde{D}_0AA @A\tilde{D}_0AA @A\tilde{D}_0AA \\
A^{0,0}\left( M, \mathcal{E}_{\mathcal{L}_0} \right) @>>>  A^{0,1}\left( M, \mathcal{E}_{\mathcal{L}_0} \right) @>-\bar{\partial}>> A^{0,2}\left( M, \mathcal{E}_{\mathcal{L}_0} \right) @>\bar{\partial}>> \cdots 
\end{CD}
\end{equation}}}

Then we have the foliated Kodaira-Spencer map of $\left(\mathcal{M}, T, B, \pi \right)$ at $t=0$ in terms of sections
\begin{align*}
\varphi_0: T_0 B &\to \mathbb{H}^1\left( M, \mathcal{E}_{\mathcal{L}_0}^\bullet \right)\\
 \frac{\partial}{\partial t} &\mapsto \left( \frac{\partial \left( \mathcal{M}, T \right)}{\partial t}\right)|_{t=0}:= \left( \left\{ \frac{\partial X_j}{\partial t}|_{t=0}  \right\}, \left\{ \left( \frac{\partial \varphi}{\partial t}|_{t=0},  \frac{\partial U_j}{\partial t}|_{t=0} \right) \right\} \right)
\end{align*}

Now we assume that the foliated analytic family $\left(\mathcal{M}, T, B, \pi \right)$ with $T\in  H^0\left(\mathcal{M}, \Theta_{\mathcal{M}/B}\otimes \mathcal{L}  \right)$ in terms of sections is a $1$-parameter family, i.e. $0\in \Delta \subset B\subset \mathbb{C}$. By taking the derivative of $(\ref{oti5})-(\ref{oti7})$ with respect to $t$ again and setting $t=0$, we have
\begin{align*}
&\bar{\partial}\left( \frac{\partial^2 \varphi}{\partial t}|_{t=0} \right) - \left[ \frac{\partial \varphi}{\partial t}|_{t=0}, \frac{\partial \varphi}{\partial t}|_{t=0} \right] =0 \\
&\bar{\partial}\left( \frac{\partial^2 X_j }{\partial t^2}|_{t=0} \right) - \left[ \frac{\partial^2 \varphi}{\partial t^2 }|_{t=0}, T_{0j} \right] - 2 \left[ \frac{\partial \varphi}{\partial t}|_{t=0}, \frac{\partial X_j}{\partial t}|_{t=0} \right] - \frac{\partial^2 U_j}{\partial t^2}\cdot T_{0j} - 2 \frac{\partial U_j}{\partial t}|_{t=0} \cdot \frac{\partial X_j}{\partial t}|_{t=0}=0 \\
& \bar{\partial} \left( \frac{\partial^2 U_j}{\partial t^2}|_{t=0} \right)  - 2 \left[ \frac{\partial \varphi}{\partial t}|_{t=0}, \frac{\partial U_j}{\partial t}|_{t=0}   \right] = 0
\end{align*}
Then
\begin{align}\label{oti12}
&\left\{ -2 \left[ \frac{\partial \varphi}{\partial t}|_{t=0}, \frac{\partial X_j}{\partial t}|_{t=0} \right] - 2 \frac{\partial U_j}{\partial t}|_{t=0} \cdot \frac{\partial X_j}{\partial t}|_{t=0} \right\}, \left\{ \left( - \left[ \frac{\partial \varphi}{\partial t}|_{t=0}, \frac{\partial \varphi}{\partial t}|_{t=0} \right] , -2 \left[ \frac{\partial \varphi}{\partial t}|_{t=0}, \frac{\partial U_j}{\partial t}|_{t=0} \right] \right) \right\}\\
& \in A^{0,1}\left( M, \Theta_M \otimes \mathcal{L}_0 \right) \bigoplus A^{0,2}\left( M, \mathcal{E}_{\mathcal{L}_0} \right)\notag
\end{align}
is a $2$-cocylce in the above Dolbeault resolution of $\mathcal{E}_{\mathcal{L}_0}^\bullet$, i.e the image of $\left( \left\{ \frac{\partial^2 X_j}{\partial t^2}|_{t=0} \right\}, \left\{ \frac{\partial^2 \varphi}{\partial t^2}|_{t=0}, \frac{\partial^2 U_j}{\partial t^2}|_{t=0} \right\} \right)$.

For any $1$-cocycle $(\{a_j\}, \{\psi, \{c_j\}\})\in A^{0,0}\left( M, \Theta_M\otimes \mathcal{L}_0 \right) \bigoplus A^{0,1}\left( M, \mathcal{E}_{\mathcal{L}_0} \right)$ in the Dolbeault resolution $(\ref{oti9})$ of $\mathcal{E}_{\mathcal{L}_0}^\bullet$, we can define a $2$-cocycle by (we omit the proof, see Part VI-1)
\begin{align}
\left( \left\{ 2[\psi, a_j]+2c_j \cdot a_j \right\} , \left\{ \left( [\psi, \psi], 2 [\psi, c_j] \right) \right\} \right) \in A^{0,1}\left( M, \Theta_M \otimes \mathcal{L}_0 \right) \bigoplus A^{0,2}\left( M, \mathcal{E}_{\mathcal{L}_0} \right)
\end{align}

Then $(\ref{oti12})$ implies that
\begin{lemma}
Let $(\{a_j\}, \{\psi, \{c_j\}\})\in A^{0,0}\left( M, \Theta_M\otimes \mathcal{L}_0 \right) \bigoplus A^{0,1}\left( M, \mathcal{E}_{\mathcal{L}_0} \right)$ be a $1$-cocycle. Then
\begin{align} \label{oti10}
\left( \left\{ 2[\psi, a_j]+2c_j \cdot a_j \right\} , \left\{ \left( [\psi, \psi], 2 [\psi, c_j] \right) \right\} \right) \in A^{0,1}\left( M, \Theta_M \otimes \mathcal{L}_0 \right) \bigoplus A^{0,2}\left( M, \mathcal{E}_{\mathcal{L}_0} \right)
\end{align}
defines $0$ class in $\mathbb{H}^2\left( M, \mathcal{E}_{\mathcal{L}_0}^\bullet \right)$.
\end{lemma}

\begin{theorem}
Suppose given a compact complex manifold $M$ with a dimension $1$-foliation defined by $T_0\in H^0\left(M, \Theta_M\otimes \mathcal{L}_0 \right)$ for a line bundle $\mathcal{L}_0$, and $(\{a_j\}, \{\psi, \{c_j\}\})\in A^{0,0}\left( M, \Theta_M\otimes \mathcal{L}_0 \right) \bigoplus A^{0,1}\left( M, \mathcal{E}_{\mathcal{L}_0} \right)$ in $\mathbb{H}^1\left(M, \mathcal{E}_{\mathcal{L}_0}^\bullet \right)$. In order that there may exist a foliated analytic family $\left( \mathcal{M}, T, B, \pi\right)$ such that $\pi^{-1}(0)=\left( M, T_0 \right)$ and $\left(\frac{\partial \left(\mathcal{M}, T \right)}{\partial t}\right)|_{t=0}=\left(\left\{a_j \right\}, \left\{\left( \psi , c_{j} \right) \right\} \right)$, it is necessary  that $(\ref{oti10})$ defines the $0$ class in $\mathbb{H}^2\left( M, \mathcal{E}_{\mathcal{L}_0}^\bullet \right)$.

\end{theorem}

\section{Infinitesimal foliated deformations and obstructions of codimension 1 foliations in terms of deformations of sections}

\subsection{Infinitesimal deformations of foliated complex analytic structures of codimension $1$ in terms of sections}\label{on25}\

Let $\left(M, \mathcal{N}_{\mathcal{F}_0}^* \right)$ be a compact foliated complex manifold with $\mathcal{N}_{\mathcal{F}_0}^*$ locally free of rank of $1$, so that there exists a section $w_0\in H^0\left( M, \Omega_M^1 \otimes {\mathcal{L}_0} \right)$ with $\mathcal{L}_0 = \left( \mathcal{N}_{\mathcal{F}_0}^*\right)^*$, which defines the foliation on $M$. In this subsection, we will consider foliated deformations of $\left( M, \mathcal{N}_{\mathcal{F}_0}^* \right)$ as deformations of the section $w_0$ (see Part II). Then we have the complex of sheaves controlling foliated deformations of $(M, \mathcal{N}_{\mathcal{F}_0}^*)$ in terms of cotangent sheaves (equivalently in terms of sections)
\begin{align}\label{cp2}
\mathcal{E}_{\mathcal{L}_0}^\bullet: \mathcal{E}_{ \mathcal{L}_0 } \xrightarrow{\tilde{E}_0} \mathcal{S}^1=\Omega_M^1\otimes {\mathcal{L}_0} \xrightarrow{\tilde{E}_1}    \tilde{\mathcal{S}}^2 \otimes \mathcal{L}_0 \left( \subset \mathcal{L}_0^{\otimes 2} \otimes \bigwedge^3 \Omega_M^1 \right) \xrightarrow{E_2} \tilde{\mathcal{S}}^3  \otimes  \mathcal{L}_0 \left( \subset \mathcal{L}_0^{\otimes 2} \otimes \bigwedge^4\Omega_{M}^1 \right) \xrightarrow{E_3} \cdots
\end{align}
We recall that $\mathcal{E}_{\mathcal{L}_0}$ and the differential $\tilde{E}_0, \tilde{E}_1
$ are defined in the following way: for an open covering $\mathcal{U}_0=\{U_j\}$ of $M$ by coordinate neighborhoods, let $\mathcal{L}_0$ defined by the transition functions $\left\{ \Psi_{0jk} \right\}\in C^1 \left(\mathcal{U}_0, \mathcal{O}_M^* \right)$ with $\Psi_{0ij}\Psi_{0jk}=\Psi_{0ik}$, and $w_0\in \Gamma\left( M, \Omega_M^1 \otimes \mathcal{L}_0 \right)$ be written as $w_{0j}\in \Gamma\left( U_j, \Omega_M^1 \right)$  on $U_j$ with $w_{0j}= \Psi_{0jk} w_{0k}$ on $U_j\cap U_k$.  Then $\mathcal{E}_{\mathcal{L}_0} |_{U_j}$ is isomorphic to $\Theta_{M}|_{U_j}\oplus \mathcal{O}_{M}|_{U_j}$, and $(\rho_j, a_j)\in \Theta_{M}|_{U_j}\oplus \mathcal{O}_{M}|_{U_j}$ is identified with $\left(\rho_k, a_k\right)\in \Theta_{M}|_{U_k}\oplus \mathcal{O}_{M}|_{U_k}$ on $U_j\cap U_k$ if and only if $\rho_j=\rho_k$ and $a_k-a_j=\left[\rho_j, \log \Psi_{0jk}\right]$. On the other hand, the differential $\tilde{E}_0$ is defined in the following way:
\begin{align*}
\tilde{E}_0: \Gamma\left(U_j,\mathcal{E}_{\mathcal{L}_0}\right)\cong \Gamma\left(U_j,\Theta_{M}\right)\oplus \Gamma\left(U_j ,  \mathcal{O}_{M} \right) &\to \Gamma\left(U_j, \Omega_{M}^1 \otimes  \mathcal{L}_0 \right)\\
     \left(\rho_j, a_j\right)&\mapsto \mathcal{L}_{\rho_j}\left( w_{0j}\right) +a_j w_{0j}
\end{align*}
and $\tilde{E}_1$ is defined in the following way:
\begin{align*}
\tilde{E}_1: \Gamma\left( U_j, \Omega_M^1\otimes \mathcal{L}_0 \right)  &\to \Gamma\left( U_j, \mathcal{L}_0\otimes \tilde{\mathcal{S}}^2 \right) \subset \Gamma\left( U_j, \mathcal{L}_0^{\otimes 2} \otimes \bigwedge^3 \Omega_M^1 \right) \\
              v_j &\mapsto v_j \wedge dw_0 + w_0\wedge dv_j
\end{align*}
We will denote $i$-th cohomology group of $\mathcal{E}_{\mathcal{L}_0}^\bullet$ by $\mathbb{H}^i\left( M, \mathcal{E}_{\mathcal{L}_0}^\bullet \right)$. 

We shall identify infinitesimal foliated deformations of codimension $1$-foliations in terms of sections.  Let $\left(\mathcal{M}, \mathcal{N}_{\mathcal{F}}^*, B, \pi \right)$ with $\mathcal{N}_{\mathcal{F}}^*$ locally free of rank $1$ be a foliated complex analytic family of deformations of $\left( M, \mathcal{N}_{\mathcal{F}_0}^* \right)=\pi^{-1}(0), 0 \in B \subset \mathbb{C}^m$. Then the exact sequence $0 \to \mathcal{N}_\mathcal{F}^* \to \Omega_{\mathcal{M}/B}^1 \to \frac{\Omega_{\mathcal{M}/B}^1}{\mathcal{N}_{\mathcal{F}}^*}\to 0 $ gives a section $w\in H^0\left(\mathcal{M}, \Omega_{\mathcal{M}/B}^1\otimes \mathcal{L} \right)$ with $\mathcal{L}= \left( \mathcal{N}_{\mathcal{F}}^*\right)^*$ which induces $w_0\in H^0\left( M, \Omega_M^1 \otimes \mathcal{L}_0 \right)$. We will also denote the family by $\left( \mathcal{M}, w , B, \pi \right)$. Take a small disk $\Delta$ with centre $0$ such that $0\in \Delta \subset B$, and represent $\mathcal{M}_\Delta=\pi^{-1}(\Delta)$ in the form $\mathcal{M}_\Delta = \bigcup_{j=1}^l U_j\times \Delta$ where each $U_j$ is a polydisk, and $(\xi_j,t)\in U_j\times \Delta$ and $(\xi_k,t)\in U_k\times \Delta$ are the same points on $\mathcal{M}_\Delta$ if $\xi_j^\alpha= f_{jk}^\alpha(\xi_k,t)$ for $\alpha=1,...,n$. Here each $f_{jk}^\alpha(\xi_k,t)=f_{jk}^\alpha(\xi_k^1,..., \xi_k^n, t)$ is a holomorphic function of $\xi_k^1,..., \xi_k^n,t$ defined on $\left(U_k\times \Delta \right) \cap \left( U_j\times \Delta \right)\ne \emptyset$, and we have 
\begin{align}\label{on20}
f_{ik}^\alpha(\xi_k,t)&=f_{ij}^\alpha\left(f_{jk}(\xi_k,t),t \right),\,\,\,\alpha=1,...,n  \,\,\,\,\,\,\textnormal{on}\,\,\, \left(U_i\times \Delta \right)\cap \left(U_j\times \Delta \right) \cap \left( U_k \times \Delta \right)\ne \emptyset
\end{align}
 On the other hand, the line bundle $\mathcal{L}$ on $\mathcal{M}_\Delta$ is defined by the system of transition functions $\left\{\Psi_{jk}(\xi_k,t) \right\}$ where $\Psi_{jk}(\xi_k,t)\in   \Gamma\left(\left(U_j\times \Delta \right) \cap \left( U_k\times \Delta \right), \mathcal{O}_\mathcal{M}^* \right)$ satisfying the cocycle condition 
\begin{align}\label{on1}
\Psi_{ij}\left(f_{jk}(\xi_k,t), t \right)\cdot \Psi_{jk}(\xi_k,t) = \Psi_{ik}(\xi_k,t), \,\,\,\,\,\textnormal{on}\,\,\,\,\, \left(U_i\times \Delta \right) \cap \left( U_j\times \Delta \right) \cap \left( U_k \times \Delta \right) \ne \emptyset
\end{align}
and $w\in H^0\left( \mathcal{M}, \Omega_{\mathcal{M}/B}^1\otimes \mathcal{L} \right)$ on $U_j\times \Delta$ is of the form
\begin{align*}
w_j(\xi_j,t):= w|_{U_j\times \Delta} = \sum_{\alpha=1}^n w_j^\alpha(\xi_j,t)d\xi_j^\alpha\,\,\,\,\,\,\,\,\textnormal{on}\,\,\,U_j\times \Delta
\end{align*}
such that
\begin{align*}
w_j(\xi_j,t)= \Psi_{jk}(\xi_k,t) w_k(\xi_k,t)\,\,\,\,\,\,\,\,\,\textnormal{on}\,\,\, \left(U_j\times \Delta\right)\cap \left( U_k\times \Delta \right)
\end{align*}
Equivalently we have
\begin{align}
&\sum_{\alpha=1}^n w_j^\alpha(\xi_j,t)d\xi_j^\alpha= \Psi_{jk}(\xi_k,t)\sum_{\beta=1}^n w_k^\beta(\xi_k,t)d \xi_k^\beta \notag\\
& \iff \sum_{\alpha=1}^n w_j^\alpha \left(f_{jk}(\xi_k,t),t\right) \frac{\partial f_{jk}^\alpha(\xi_k,t)}{\partial \xi_k^\beta}= \Psi_{jk}(\xi_k,t) w_k^\beta(\xi_k,t) \label{on3}
\end{align}
and the integrability condition
\begin{align}\label{on5}
w_j(\xi_j,t) \wedge dw_j(\xi_j,t)=\left( \sum_{\alpha=1}^n w_j^\alpha(\xi_j, t) d\xi_j^\alpha \right) \wedge d\left( \sum_{\alpha=1}^n w_j^\alpha(\xi_j, t) d\xi_j^\alpha \right)=0
\end{align}

Now we identify the infinitesimal foliated deformations. By taking the derivative of $(\ref{on20})$ with respect to $t$ and setting $\theta_{jk}(t)=\sum_{\alpha=1}^n \theta_{jk}^\alpha(\xi_j,t)\frac{\partial}{\partial \xi_j^\alpha}$ where $\theta_{jk}^\alpha(\xi_j,t)=\frac{\partial f_{jk}(\xi_k,t)}{\partial t},\xi_k=f_{kj}(\xi_j,t)$, we have
\begin{align}\label{on21}
\theta_{ik}^\alpha(\xi_i,t)=\theta_{ij}^\alpha(\xi_i,t) + \sum_{\beta=1}^n \frac{\partial \xi_i^\alpha}{\partial \xi_j^\beta}\theta_{jk}^\beta(\xi_j,t),\,\,\,\alpha=1,...,n
\end{align}
By setting $t=0$, we obtain 
\begin{align} \label{on4}
\theta_{ij}(0)-\theta_{ik}(0) + \theta_{jk}(0)=0
\end{align}

As in subsection \ref{ot23}, we have
\begin{align}\label{on22}
 -\frac{1}{\Psi_{0ij}}\frac{\partial \Psi_{ij}}{\partial t}|_{t=0} + \frac{1}{\Psi_{0ik}} \frac{\partial \Psi_{ik}}{\partial t} |_{t=0} - \frac{1}{\Psi_{jk}}\frac{\partial \Psi_{jk}}{\partial t } - \left[ \sum_{\alpha=1}^n \frac{\partial f_{jk}^\alpha}{\partial t}|_{t=0} \frac{\partial}{\partial \xi_j^\alpha} , \log \Psi_{0ij} \right] =0
\end{align}
where $\Psi_{0jk}:=\Psi_{jk}\left( \xi_k, 0\right)$, i.e. $\left\{ \Psi_{0jk} \right\}$ defines the $\mathcal{L}_0$ on $M$.

 By taking the derivative of $(\ref{on3})$ with respect to $t$, we have
\begin{align}\label{on9}
\sum_{\alpha, \gamma=1}^n \frac{\partial w_j^\alpha}{\partial \xi_j^\gamma} \frac{\partial f_{jk}^\gamma}{\partial t} \frac{\partial f_{jk}^\alpha}{\partial \xi_k^\beta} + \sum_{\alpha=1}^n \frac{\partial w_j^\alpha}{\partial t}\frac{\partial f_{jk}^\alpha}{\partial \xi_k^\beta} + \sum_{\alpha=1}^n w_j^\alpha \frac{\partial}{\partial \xi_k^\beta}\left(\frac{\partial f_{jk}^\alpha}{\partial t} \right) = \frac{\partial \Psi_{jk}}{\partial t} w_k^\beta + \Psi_{jk}\frac{\partial w_k^\beta}{\partial t}
\end{align}
By multiplying $\sum_{\alpha=1}^n d\xi_j^\alpha$, we obtain
{\small{\begin{align}\label{on7}
 \sum_{\alpha, \gamma=1}^n \frac{\partial w_j^\alpha}{\partial \xi_j^\gamma}\frac{\partial f_{jk}^\gamma}{\partial t} d\xi_j^\alpha + \sum_{\alpha=1}^n \frac{\partial w_j^\alpha}{\partial t} d\xi_j^\alpha + \sum_{\alpha,\beta=1}^n w_j^\alpha \frac{\partial}{\partial \xi_j^\beta}\left(\frac{\partial f_{jk}^\alpha}{\partial t} \right)d\xi_j^\beta = \frac{\partial \Psi_{jk}}{\partial t} \sum_{\beta=1}^n w_k^\beta d\xi_k^\beta +  \Psi_{jk}\sum_{\beta=1}^n \frac{\partial w_k^\beta}{\partial t} d\xi_k^\beta
\end{align}}}

We note that
\begin{align*}
\mathcal{L}_{\sum_{\gamma=1}^n \frac{\partial f_{jk}^\gamma}{\partial t} \frac{\partial}{\partial \xi_j^\gamma}}\left(\sum_{\alpha=1}^n w_j^\alpha d\xi_j^\alpha \right)= \sum_{\alpha, \gamma=1}^n \frac{\partial f_{jk}^\gamma}{\partial t} \frac{\partial w_j^\alpha}{\partial \xi_j^\gamma} d\xi_j^\alpha + \sum_{\alpha=1}^nw_j^\alpha d\left( \frac{\partial f_{jk}^\alpha}{\partial t} \right)
\end{align*}

Then by setting $t=0$ on $(\ref{on7})$, we have
\begin{align}\label{on8}
-\Psi_{0jk} \frac{\partial w_k}{\partial t}|_{t=0} + \frac{\partial w_j}{\partial t}|_{t=0} + \mathcal{L}_{\sum_{\gamma=1}^n \frac{\partial f_{jk}^\gamma}{\partial t}|_{t=0} \frac{\partial}{\partial \xi_j^\gamma}}\left(\sum_{\alpha=1}^n w_{0j}^\alpha d\xi_j^\alpha \right)- \frac{1}{\Psi_{0jk}}\frac{\partial \Psi_{jk}}{\partial t}|_{t=0} w_{0j}=0
\end{align}
where $w_{0j}=\sum_{\beta=1}^n w_{0j}^\beta d\xi_j^\beta:=\sum_{\beta=1}^n w_j^\beta(\xi_j, 0)d\xi_j^\beta$ and $\frac{\partial w_j}{\partial t}:= \sum_{\beta=1}^n \frac{\partial w_j(\xi_j,t)}{\partial t} d\xi_j^\beta$.

By taking the derivative of $(\ref{on5})$ with respect to $t$, we have
\begin{align}\label{on16}
\frac{\partial w_j}{\partial t} \wedge dw_j + w_j \wedge d\left( \frac{\partial w_j}{\partial t} \right)=0
\end{align}
and setting $t=0$, we obtain
\begin{align}\label{on6}
\frac{\partial w_j}{\partial t}|_{t=0} \wedge d w_{0j}+ w_{0j} \wedge d\left( \frac{\partial w_j}{\partial t}|_{t=0} \right) =0
\end{align}

Then $(\ref{on4}),(\ref{on22}), (\ref{on8})$ and $(\ref{on6})$ imply that
\begin{align*}
\left( \left\{ - \frac{\partial w_j}{\partial t}|_{t=0}  \right\} , \left\{  \sum_{\alpha=1}^n \frac{\partial f_{jk}^\alpha}{\partial t}|_{t=0} \frac{\partial}{\partial \xi_j^\alpha} , - \frac{1}{\Psi_{0jk}}\frac{\partial \Psi_{jk}}{\partial t} |_{t=0} \right\}   \right)  \in C^0\left(\mathcal{U}_0, \Omega_M^1\otimes \mathcal{L}_0     \right) \bigoplus C^1\left( \mathcal{U}_0, \mathcal{E}_{\mathcal{L}_0} \right)
\end{align*}
defines a $1$-cocycle in the following \v Cech resolution of $\mathcal{E}_{\mathcal{L}_0}^\bullet$ with respect to the covering $\mathcal{U}_0:=\mathcal{U}\cap M$.
\begin{equation}\label{on18}
\begin{CD}
\cdots\\ 
@A E_2AA \\
C^0\left( \mathcal{U}_0, \tilde{\mathcal{S}}^2\otimes  \mathcal{L}_0  \right) @>-\delta>> \cdots \\
@A\tilde{E}_1AA  @A\tilde{E}_1AA \\
C^0\left( \mathcal{U}_0, \Omega_M^1\otimes \mathcal{L}_0 \right) @>\delta>> C^1\left( \mathcal{U}_0, \Omega_M^1\otimes \mathcal{L}_0 \right) @>-\delta>> \cdots \\
@A\tilde{E}_0AA @A\tilde{E}_0AA @A\tilde{E}_0AA \\
C^0\left(\mathcal{U}_0, \mathcal{E}_{\mathcal{L}_0} \right) @>-\delta>> C^1\left( \mathcal{U}_0, \mathcal{E}_{\mathcal{L}_0}\right) @>\delta>> C^2\left( \mathcal{U}, \mathcal{E}_{\mathcal{L}_0} \right) @> - \delta >> \cdots
\end{CD}
\end{equation}
Then we have the foliated Kodaira-Spencer map of $\left(\mathcal{M}, w, B, \pi \right)$ at $t=0$ in terms of sections
\begin{align*}
\varphi_0: T_0 B &\to  \mathbb{H}^1\left( M, \mathcal{E}_{\mathcal{L}_0}^\bullet \right) \\
\frac{\partial}{\partial t} &\mapsto \left(\frac{\partial\left(\mathcal{M}, w \right)}{\partial t} \right)|_{t=0}:=\left( \left\{ - \frac{\partial w_j}{\partial t}|_{t=0} \right\} , \left\{ \left( \sum_{\alpha=1}^n \frac{\partial f_{jk}^\alpha}{\partial t} |_{t=0} \frac{\partial}{\partial \xi_j^\alpha}  , - \frac{1}{\Psi_{0jk}}\frac{\partial \Psi_{jk}}{\partial t}|_{t=0} \right) \right\} \right)
\end{align*}

\subsection{Obstructions in terms of \v Cech resolution of $\mathcal{E}_{\mathcal{L}_0}^\bullet$}\label{on26} \

We assume that the foliated analytic family $\left(\mathcal{M}, w , B, \pi \right)$ with $w\in  H^0\left(\mathcal{M}, \Omega_{\mathcal{M}/B}^1\otimes \mathcal{L}  \right)$ in terms of sections is a $1$-parameter family, i.e. $0\in \Delta \subset B\subset \mathbb{C}$. We recall $(\ref{ot12})$ in subsection \ref{ot23}. Then by taking the derivative of $(\ref{on21})$ with respect to $t$, and putting
\begin{align*}
\dot{\theta}_{ij}(t)= \sum_{\alpha=1}^n \dot{\theta}_{ij}^\alpha (\xi_i,t)\frac{\partial}{\partial \xi_i^\alpha} , \,\,\,\,\,\textnormal{where}\,\,\,\dot{\theta}_{ij}^\alpha(\xi_i,t)=\frac{\partial \theta_{ij}^\alpha(\xi_i,t)}{\partial t}
\end{align*}
we obtain (for the detail, see \cite{Kod05} p.211-213)
\begin{align}\label{}
\dot{\theta}_{ij}(t) -\dot{\theta}_{ik}(t) + \dot{\theta}_{jk} (t) = \left[ \theta_{ij}(t), \theta_{jk}(t) \right]
\end{align}

Let us take the derivative of $(\ref{on9})$  with respect to $t$ again:
{\small{\begin{align*}
&\sum_{\alpha, \gamma,\eta=1}^n \frac{\partial^2 w_j^\alpha}{\partial \xi_j^\eta \partial \xi_j^\gamma}\frac{\partial f_{jk}^\eta}{\partial t} \frac{\partial f_{jk}^\gamma}{\partial t} \frac{\partial f_{jk}^\alpha}{\partial \xi_k^\beta} + 2\sum_{\alpha, \gamma=1}^n \frac{\partial}{\partial \xi_j^\gamma} \left(\frac{\partial w_j^\alpha}{\partial t}\right) \frac{\partial f_{jk}^\gamma}{\partial t} \frac{\partial f_{jk}^\alpha}{\partial \xi_k^\beta} + \sum_{\alpha, \gamma=1}^n \frac{\partial w_j^\alpha}{\partial \xi_j^\gamma} \frac{\partial^2 f_{jk}^\gamma}{\partial t^2}\frac{\partial f_{jk}^\alpha}{\partial \xi_k^\beta} +2 \sum_{\alpha, \gamma=1}^n \frac{\partial w_j^\alpha}{\partial \xi_j^\gamma}\frac{\partial f_{jk}^\gamma}{\partial t} \frac{\partial}{\partial \xi_k^\beta}\left(\frac{\partial f_{jk}^\alpha}{\partial t} \right)\\
&+ \sum_{\alpha=1}^n \frac{\partial^2 w_j^\alpha}{\partial t^2}\frac{\partial f_{jk}^\alpha}{\partial \xi_k^\beta} +  2\sum_{\alpha=1}^n \frac{\partial w_j^\alpha}{\partial t} \frac{\partial}{\partial \xi_k^\beta}\left( \frac{\partial f_{jk}^\alpha}{\partial t} \right) + \sum_{\alpha=1}^n w_j^\alpha \frac{\partial}{\partial \xi_k^\beta}\left( \frac{\partial^2 f_{jk}^\alpha}{\partial t^2} \right)= \frac{\partial^2 \Psi_{jk}}{\partial t^2}w_k^\beta + 2\frac{\partial \Psi_{jk}}{\partial t} \frac{\partial w_k^\beta}{\partial t} + \Psi_{jk} \frac{\partial^2 w_k^\beta}{\partial t^2}
\end{align*}}}

By summing $\sum_{\beta=1}^n d\xi_k^\beta$ and by noting $\dot{\theta}_{jk}^\gamma= \frac{\partial f_{jk}^\gamma}{\partial t^2} - \sum_{\sigma=1}^n \frac{\partial f_{jk}^\sigma}{\partial t}\frac{\partial}{\partial \xi_j^\sigma}\left( \frac{\partial f_{jk}^\gamma}{\partial t} \right)$, we obtain
{\small{\begin{align}\label{on10}
&\sum_{\alpha, \gamma,\eta=1}^n \frac{\partial^2 w_j^\alpha}{\partial \xi_j^\eta \partial \xi_j^\gamma}\frac{\partial f_{jk}^\eta}{\partial t} \frac{\partial f_{jk}^\gamma}{\partial t} d\xi_j^\alpha + 2\sum_{\alpha, \gamma=1}^n  \frac{\partial f_{jk}^\gamma}{\partial t}\frac{\partial}{\partial \xi_j^\gamma} \left(\frac{\partial w_j^\alpha}{\partial t}\right) d\xi_j^\alpha + \sum_{\alpha, \gamma=1}^n \left( \dot{\theta}_{jk}^\gamma + \sum_{\sigma=1}^n \frac{\partial f_{jk}^\sigma}{\partial t}\frac{\partial}{\partial \xi_j^\sigma}\left(\frac{\partial f_{jk}^\gamma}{\partial t} \right) \right) \frac{\partial w_j^\alpha}{\partial \xi_j^\gamma} d\xi_j^\alpha\\
& +2 \sum_{\alpha, \gamma=1}^n \frac{\partial f_{jk}^\gamma}{\partial t}  \frac{\partial w_j^\alpha}{\partial \xi_j^\gamma}d\left(\frac{\partial f_{jk}^\alpha}{\partial t} \right)+ \sum_{\alpha=1}^n \frac{\partial^2 w_j^\alpha}{\partial t^2}d\xi_j^\alpha +  2\sum_{\alpha=1}^n \frac{\partial w_j^\alpha}{\partial t} d\left( \frac{\partial f_{jk}^\alpha}{\partial t} \right) + \sum_{\alpha=1}^n w_j^\alpha d\left( \dot{\theta}_{jk}^\alpha + \sum_{\sigma=1}^n \frac{\partial f_{jk}^\sigma}{\partial t} \frac{\partial }{\partial \xi_j^\sigma}\left( \frac{\partial f_{jk}^\alpha}{\partial t}\right)\right) \notag \\
&= \frac{\partial^2 \Psi_{jk}}{\partial t^2}w_k + 2\frac{\partial \Psi_{jk}}{\partial t} \frac{\partial w_k}{\partial t} + \Psi_{jk} \frac{\partial^2 w_k}{\partial t^2} \notag
\end{align}}}

Then $(\ref{on10})$ is equivalent to
{\small{\begin{align}\label{on12}
&\sum_{\alpha, \gamma,\eta=1}^n \frac{\partial^2 w_j^\alpha}{\partial \xi_j^\eta \partial \xi_j^\gamma}\frac{\partial f_{jk}^\eta}{\partial t} \frac{\partial f_{jk}^\gamma}{\partial t} d\xi_j^\alpha  + 2 \mathcal{L}_{\theta_{jk}(t)}\left(\frac{\partial w_j}{\partial t} \right) + \mathcal{L}_{\dot{\theta}_{jk}(t)}\left(w_j \right) + \sum_{\alpha, \gamma, \sigma=1}^n \frac{\partial f_{jk}^\sigma}{\partial t}\frac{\partial}{\partial \xi_j^\sigma}\left( \frac{\partial f_{jk}^\gamma}{\partial t} \right)\frac{\partial w_j^\alpha}{\partial \xi_j^\gamma} d\xi_j^\alpha  \\
&+2 \sum_{\alpha, \gamma=1}^n \frac{\partial f_{jk}^\gamma}{\partial t}  \frac{\partial w_j^\alpha}{\partial \xi_j^\gamma}d\left(\frac{\partial f_{jk}^\alpha}{\partial t} \right)+ \frac{\partial^2 w_j}{\partial t^2} + \sum_{\alpha, \sigma=1}^n w_j^\alpha d\left( \frac{\partial f_{jk}^\sigma}{\partial t} \frac{\partial}{\partial \xi_j^\sigma}\left( \frac{\partial f_{jk}^\alpha}{\partial t} \right) \right) =  \frac{\partial^2 \Psi_{jk}}{\partial t^2}w_k + 2\frac{\partial \Psi_{jk}}{\partial t} \frac{\partial w_k}{\partial t} + \Psi_{jk} \frac{\partial^2 w_k}{\partial t^2}  \notag
\end{align}}}

We note that
{\small{\begin{align*}
&\mathcal{L}_{\theta_{jk}}\mathcal{L}_{\theta_{jk}}\left(w_j \right) =\mathcal{L}_{\sum_{\sigma=1}^n \frac{\partial f_{jk}^\sigma}{\partial t} \frac{\partial}{\partial \xi_j^\sigma}} \left( \mathcal{L}_{\sum_{\delta=1}^n \frac{\partial f_{jk}^\delta}{\partial t} \frac{\partial}{\partial \xi_j^\delta}} \left(\sum_{\gamma=1}^n w_j^{ \gamma} d\xi_j^\gamma \right)   \right) =\mathcal{L}_{\sum_{\sigma=1}^n \frac{\partial f_{jk}^\sigma}{\partial t} \frac{\partial}{\partial \xi_j^\sigma}} \left( \sum_{\delta, \gamma=1}^n \frac{\partial f_{jk}^\delta}{\partial t} \frac{\partial w_j^{ \gamma}}{\partial \xi_j^\delta} d\xi_j^\gamma + \sum_{\gamma=1}^n w_j^{ \gamma} d\left(\frac{\partial f_{jk}^\gamma}{\partial t} \right)   \right) \\
&= \sum_{\sigma,\delta, \gamma=1}^n \frac{\partial f_{jk}^\sigma}{\partial t}\frac{\partial }{\partial \xi_j^\sigma}\left( \frac{\partial f_{jk}^\delta}{\partial t}\right)\frac{\partial w_j^{ \gamma}}{\partial \xi_j^\delta} d\xi_j^\gamma + \sum_{\sigma, \delta, \gamma=1}^n \frac{\partial f_{jk}^\sigma}{\partial t} \frac{\partial f_{jk}^\delta}{\partial t} \frac{\partial^2 w_j^{ \gamma}}{\partial \xi_j^\sigma \partial \xi_j^\delta} d\xi_j^\gamma + 2 \sum_{\delta, \gamma=1}^n \frac{\partial f_{jk}^\delta}{\partial t}\frac{\partial w_j^{ \gamma}}{\partial \xi_j^\delta}d\left(\frac{\partial f_{jk}^\gamma}{\partial t} \right) + \sum_{\gamma,\sigma=1}^n w_j^{ \gamma}d\left(\frac{\partial f_{jk}^\sigma}{\partial t}\frac{\partial }{\partial \xi_j^\sigma}\left(\frac{\partial f_{jk}^\gamma}{\partial t} \right) \right) 
\end{align*}}}

Then $(\ref{on12})$ is equivalent to
\begin{align}\label{on11}
 2 \mathcal{L}_{\theta_{jk}(t)}\left(  \frac{\partial w_j}{\partial t}\right) +\mathcal{L}_{\dot{\theta}_{jk}(t)}\left(w_j  \right)  + \frac{\partial^2 w_j}{\partial t_j^2} + \mathcal{L}_{\theta_{jk}}\mathcal{L}_{\theta_{jk}}\left(w_j^\alpha \right) = \frac{\partial^2 \Psi_{jk}}{\partial t^2}w_k^\beta  + 2\frac{\partial \Psi_{jk}}{\partial t} \frac{\partial w_k}{\partial t} +  \Psi_{jk} \frac{\partial^2 w_k }{\partial t^2}  
\end{align}

\begin{lemma}
We have the equality
{\small{\begin{align}\label{on13}
-\Psi_{jk}\frac{\partial^2 w_k}{\partial t^2} + \frac{\partial^2 w_j}{\partial t^2} + \mathcal{L}_{\dot{\theta}_{jk}}\left(w_j \right)- \frac{\partial}{\partial t_j}\left( \frac{1}{\Psi_{jk}}\frac{\partial \Psi_{jk}}{\partial t} \right)w_j= \mathcal{L}_{\theta_{jk}}\left(- \frac{\partial w_j}{\partial t_j}- \Psi_{jk}\frac{\partial w_k}{\partial t} \right) - \frac{1}{\Psi_{jk}}\frac{\partial \Psi_{jk}}{\partial t}\left( -\frac{\partial w_j}{\partial t_j}- \Psi_{jk}\frac{\partial w_k}{\partial t_k}\right)
\end{align}}}
\end{lemma}

\begin{proof}
We note that from $(\ref{on11})$ and $(\ref{ot10})$, we have
{\small{\begin{align*}
& \mathcal{L}_{\theta_{jk}}\left(- \frac{\partial w_j}{\partial t_j}- \Psi_{jk}\frac{\partial w_k}{\partial t} \right) - \frac{1}{\Psi_{jk}}\frac{\partial \Psi_{jk}}{\partial t}\left( -\frac{\partial w_j}{\partial t_j}- \Psi_{jk}\frac{\partial w_k}{\partial t_k}\right)\\
&= \mathcal{L}_{\theta_{jk}}\left( - 2\frac{\partial w_j}{\partial t} - \mathcal{L}_{\theta_{jk}}(w_j) + \frac{1}{\Psi_{jk}}\frac{\partial \Psi_{jk}}{\partial t} w_j \right) - \frac{1}{\Psi_{jk}}\frac{\partial \Psi_{jk}}{\partial t}\left(  - 2\Psi_{jk}\frac{\partial w_k}{\partial t} + \mathcal{L}_{\theta_{jk}}(w_j) - \frac{1}{\Psi_{jk}}\frac{\partial \Psi_{jk}}{\partial t} w_j   \right)\\
&= -2 \mathcal{L}_{\theta_{jk}}\left( \frac{\partial w_j}{\partial t} \right) - \mathcal{L}_{\theta_{jk}}\mathcal{L}_{\theta_{jk}}(w_j) + \left[ \theta_{jk}, \frac{1}{\Psi_{jk}}\frac{\partial \Psi_{jk}}{\partial t} \right]w_j  + 2\frac{\partial \Psi_{jk}}{\partial t}\frac{\partial w_k}{\partial t} + \frac{1}{\Psi_{jk}^2}\frac{\partial \Psi_{jk}}{\partial t}\frac{\partial \Psi_{jk}}{\partial t} w_j\\
&= \mathcal{L}_{\dot{\theta}_{jk}}\left( w_j \right) + \frac{\partial^2 w_j}{\partial t^2} - \frac{\partial^2 \Psi_{jk}}{\partial t^2} w_k^\beta - \Psi_{jk}\frac{\partial^2 w_k}{\partial t^2} + \left[ \theta_{jk}, \frac{1}{\Psi_{jk}}\frac{\partial \Psi_{jk}}{\partial t} \right]w_j + \frac{1}{\Psi_{jk}^2}\frac{\partial \Psi_{jk}}{\partial t}\frac{\partial \Psi_{jk}}{\partial t} w_j\\
&= -\Psi_{jk}\frac{\partial^2 w_k}{\partial t^2} + \frac{\partial^2 w_j}{\partial t^2} + \mathcal{L}_{\dot{\theta}_{jk}}\left(w_j \right)- \frac{\partial}{\partial t_j}\left( \frac{1}{\Psi_{jk}}\frac{\partial \Psi_{jk}}{\partial t} \right)w_j
\end{align*}}}
\end{proof}

As in subsection \ref{ot23}, we can show
{\small{\begin{align}\label{on15}
\delta \left( \left\{ \dot{\theta}_{jk}(0), - \frac{\partial}{\partial t_j}\left( \frac{1}{\Psi_{jk}}\frac{\partial \Psi_{jk}}{\partial t_k} \right)|_{t=0}   \right\}\right)=\left\{ \left[ \left( \theta_{ij}(0), - \frac{1}{\Psi_{0ij}}\frac{\partial \Psi_{ij}}{\partial t}|_{t=0} \right), \left( \theta_{jk}(0), -\frac{1}{\Psi_{0jk}}\frac{\partial \Psi_{jk}}{\partial t}|_{t=0} \right)    \right]\right\}
\end{align}}}
and setting $t=0$ on $(\ref{on13})$, we obtain
\begin{align}\label{on14}
&-\Psi_{0jk}\frac{\partial^2 w_k}{\partial t^2}|_{t=0} + \frac{\partial^2 w_j}{\partial t^2}|_{t=0} + \mathcal{L}_{\dot{\theta}_{jk}(0)}\left(w_{0j} \right)- \frac{\partial}{\partial t_j}\left( \frac{1}{\Psi_{jk}}\frac{\partial \Psi_{jk}}{\partial t} \right)|_{t=0} \cdot w_{0j} \\
&= \mathcal{L}_{\theta_{jk}(0)}\left(- \frac{\partial w_j}{\partial t_j}|_{t=0}- \Psi_{0jk}\frac{\partial w_k}{\partial t}|_{t=0} \right) - \frac{1}{\Psi_{0jk}}\frac{\partial \Psi_{jk}}{\partial t}|_{t=0}\left( -\frac{\partial w_j}{\partial t_j}|_{t=0} - \Psi_{0jk}\frac{\partial w_k}{\partial t_k}|_{t=0}\right) \notag
\end{align}
By taking the derivative of $(\ref{on16})$ with respect to $t$ again and setting $t=0$, we obtain
\begin{align}\label{on17}
\left(-\frac{\partial^2 w_j}{\partial t^2}|_{t=0}\right) \wedge dw_{0j} + w_{0j} \wedge d\left( -\frac{\partial^2 w_j}{\partial t^2}|_{t=0} \right) =  2\left(-\frac{\partial w_j}{\partial t}|_{t=0}\right) \wedge d\left(-\frac{\partial w_j}{\partial t}|_{t=0}\right)
\end{align}

For any $1$-cocycle $\left( \left\{a_j\right\}, \left\{\left( b_{jk}, c_{jk}\right) \right\} \right)\in C^0\left( \mathcal{U}_0, \Omega_M^1\otimes \mathcal{L}_0 \right) \bigoplus C^1\left( \mathcal{U}_0, \mathcal{E}_{\mathcal{L}_0} \right)$ in the \v Cech resolution $(\ref{on18})$ of $\mathcal{E}_{\mathcal{L}_0}^\bullet$, we can define a $2$-cocycle by (we omit the proof, see Part IV-1)
\begin{align}\label{on19}
&\left( \left\{ 2 a_j \wedge d a_j \right\} , \left\{ \mathcal{L}_{b_{jk}}\left( a_j + \Psi_{0jk} a_k  \right) + c_{jk}\left(a_j+ \Psi_{0jk} a_k \right)\right\}, \left\{\left[ \left(b_{ij}, c_{ij}\right), \left( b_{jk}, c_{jk}\right) \right] \right\} \right) \\
&\in C^0\left(\mathcal{U}_0, \Omega_M^1\otimes \tilde{\mathcal{S}}^2 \right) \bigoplus C^1\left( \mathcal{U}_0, \Omega_M^1\otimes \mathcal{L}_0 \right)\bigoplus C^2\left( \mathcal{U}_0, \mathcal{E}_{\mathcal{L}_0} \right) \notag
\end{align}

Then $(\ref{on15}), (\ref{on14})$ and $(\ref{on17})$ imply that
\begin{theorem}\label{on27}
Suppose given a compact complex manifold $M$ with a codimension $1$-foliation defined by $w_0\in H^0\left(M, \Omega_M^1 \otimes \mathcal{L}_0 \right)$ for a line bundle $\mathcal{L}_0$, and $\left(\left\{a_j \right\}, \left\{\left( b_{jk}, c_{jk} \right) \right\} \right)\in \mathbb{H}^1\left(M, \mathcal{E}_{\mathcal{L}_0}^\bullet \right)$. In order that there may exist a foliated analytic family $\left( \mathcal{M}, w, B, \pi\right)$ such that $\pi^{-1}(0)=\left( M, w_0 \right)$ and $\left(\frac{\partial \left(\mathcal{M}, w \right)}{\partial t}\right)|_{t=0}=\left(\left\{a_j \right\}, \left\{\left( b_{jk}, c_{jk} \right) \right\} \right)$, it is necessary  that $(\ref{on19})$ defines the $0$ class in $\mathbb{H}^2\left( M, \mathcal{E}_{\mathcal{L}_0}^\bullet \right)$.

\end{theorem}

Keeping the notations above,  we have the following lemma
\begin{lemma}\label{on23}
A compact foliated complex manifold $\left(M,w_0\right)$ of codimension $1$ where $w_0 \in H^0\left( M, \Omega_M^1\otimes \mathcal{L}_0 \right)$ with $w_0 \wedge dw_0 = 0 $ is obstructed in foliated deformations in terms of cotangent sheaves if there exist an element $\left( \{r_j\}, \{b_{jk}, c_{jk}\} \right) \in \mathbb{H}^1\left( M, \mathcal{E}_{\mathcal{L}_0}^\bullet \right)$ $($in particular, $\{(b_{jk}, c_{jk})\}\in C^1\left(\mathcal{U}, \mathcal{E}_{\mathcal{L}_0} \right)$ defines an element in $\textnormal{ker}\left( H^1\left( M, \mathcal{E}_{\mathcal{L}_0} \right) \xrightarrow{\tilde{E}_0} H^1\left( M, \Omega_M^1\otimes \mathcal{L}_0 \right) \right)$$)$, and an element $\{a_j\}\in \textnormal{ker}\left( H^0\left( M, \Omega_M^1\otimes \mathcal{L}_0 \right) \to H^0\left( M, \tilde{\mathcal{S}}^2 \otimes \mathcal{L}_0  \right) \right)$ such that one of the following holds$:$
\begin{enumerate}
\item Under the map
\begin{align}
H^0\left( M, \Omega_M^1\otimes \mathcal{L}_0 \right) \times H^1\left(M, \mathcal{E}_{\mathcal{L}_0} \right) &\to H^1\left(M, \Omega_M^1 \otimes \mathcal{L}_0 \right) \label{uc1}\\
\left( \{a_j\}, \{(b_{jk}, c_{jk})\} \right)&\to  \left\{ \mathcal{L}_{b_{jk}}(a_j)+ c_{jk}a_j  \right\}\notag
\end{align}
$\left\{ \mathcal{L}_{b_{jk}}(a_j)+ c_{jk}a_j  \right\}$ is not in the image of $H^1\left( M, \mathcal{E}_{\mathcal{L}_0} \right) \to H^1\left( M, \Omega_M^1\otimes \mathcal{L}_0 \right)$.
\item $H^1\left(M, \mathcal{E}_{\mathcal{L}_0} \right)=0$, and $\{ 2a_j \wedge da_j\}$ is not in the image of $H^0\left( M, \Omega_M^1\otimes \mathcal{L}_0 \right)\to H^0\left( M, \tilde{\mathcal{S}}^2\otimes \mathcal{L}_0 \right)$.
\item $\left\{[(b_{ij}, c_{ij}), (b_{jk}, c_{jk})] \right\}= \left\{[b_{ij}, b_{jk}], [b_{ij}, c_{jk}- [b_{jk}, \log \Psi_{0ij}]]  - [b_{jk}, c_{ij}]  \right\}$ is not $0$ class in $H^2\left( M, \mathcal{E}_{\mathcal{L}_0} \right)$.
\item $\left\{ [b_{ij} , b_{jk}] \right\}$ is not $0$ class in $H^2\left(M, \Theta_M \right)$.
\end{enumerate}
\end{lemma}

\begin{proof}
We prove $(1)$. First we show that the map $(\ref{uc1})$ is well-defined. In fact, since $\Psi_{0ij}a_j= a_i$ and $b_{ij}-b_{ik}+ b_{jk}=0$ and $c_{ij}- c_{ik} + c_{jk}- \left[ b_{jk}, \log \Psi_{0ij}\right]=0$, we have
\begin{align*}
&\mathcal{L}_{b_{ij}}(a_i)+ c_{ij}a_i - \mathcal{L}_{b_{ik}}(a_i) - c_{ik}a_i+  \Psi_{0ij}\left( \mathcal{L}_{b_{jk}}(a_j) + c_{jk}a_j \right)\\
&= \mathcal{L}_{b_{ij}- b_{ik} + b_{jk}} \left(a_i\right) - \mathcal{L}_{b_{jk}}\left(\Psi_{0ij} \right)a_j + (c_{ij}- c_{ik} + c_{jk})a_i\\
&=- \left[b_{jk}, \log \Psi_{0ij} \right] a_i + \left(c_{ij}- c_{ik}+ c_{jk} \right)a_i=0
\end{align*}
This implies that $\left\{ \mathcal{L}_{b_{ij}}(a_i) + c_{ij}a_i \right\}$ defines an element in $H^1\left( M, \Omega_M^1\otimes \mathcal{L} \right)$. Since $\left(\{r_j\}, \{b_{jk}, c_{jk}\}\right)$ defines an $1$-cocycle, we have
\begin{align*}
w_{0j} \wedge dr_j + r_j \wedge dw_{0j}  =0,&\,\,\,\,\,\,\,\Psi_{0jk}r_k- r_j+ \mathcal{L}_{b_{jk}}(w_{0j})+ c_{jk}w_{0j}=0 ,\\
b_{ij}- b_{ik}+ b_{jk}=0,&\,\, \,\,\,\,\,c_{ij}- c_{ik} + c_{jk}-  \left[ b_{jk}, \log \Psi_{0ij} \right]=0
\end{align*}
and $\left(\left\{ r_j + a_j \right\}, \left\{\left(b_{jk}, c_{jk}\right)\right\}\right)$ defines an $1$-cocycle. Assume that $(M, \omega)$ is unobstructed in foliated deformations in terms of cotangent sheaves. Then by Theorem \ref{on27} there exists $\left(\{g_i\}, \{(e_{jk}, h_{jk})\} \right)$ and $(\{g_i'\} ,   \{(e_{jk}', h_{jk}')\}\in C^0\left(\mathcal{U}, \Omega_M^1\otimes \mathcal{L}_0 \right) \bigoplus C^1\left(\mathcal{U}, \mathcal{E}_{\mathcal{L}_0} \right)$ such that
{\small{\begin{align}\label{um7}
w_{0i}\wedge dg_i+ g_i\wedge dw_{0i}&=2 r_j\wedge dr_j \\
\delta\left(\{g_i\} \right) + \left\{ \mathcal{L}_{e_{jk}}(w_{0j} ) + h_{jk}w_{0j}\right\}&= \mathcal{L}_{b_{jk}}\left(r_j+ \Psi_{0jk}r_k \right)+ c_{jk}\left( r_j+ \Psi_{0jk}r_k \right) \notag\\
\delta\left(\left\{ (e_{jk}, h_{jk})\right\} \right)&=\left\{ \left[ (b_{ij}, c_{ij}) , (b_{jk}, c_{jk}) \right] \right\} \notag
\end{align}}}
and
{\small{\begin{align}\label{um8}
w_{0i}\wedge dg_i'+ g_i'\wedge dw_{0i}&=2 (r_j+ a_j)\wedge d(r_j+a_j) \\
\delta\left(\{g_i'\} \right) + \left\{ \mathcal{L}_{e_{jk}'}(w_{0j} ) + h_{jk}'w_{0j}\right\}&= \mathcal{L}_{b_{jk}}\left(r_j+ a_j+ \Psi_{0jk}(r_k+a_k) \right)+ c_{jk}\left( r_j+a_j+ \Psi_{0jk}(r_k+a_k) \right) \notag\\
\delta\left(\left\{ (e_{jk}', h_{jk}')\right\} \right)&=\left\{ \left[ (b_{ij}, c_{ij}) , (b_{jk}, c_{jk}) \right] \right\} \notag
\end{align}}}
Then by subtracting $(\ref{um7})$ from $(\ref{um8})$, we have
\begin{align}
&w_{0i}\wedge d(g_i'-g_i)+ (g_i'-g_i)\wedge dw_{0i}=2 a_j\wedge dr_j+ 2r_j \wedge da_j + 2 a_j \wedge da_j \notag \\
&\delta\left(\{g_i'-g_i\} \right) + \left\{ \mathcal{L}_{e_{jk}'-e_{jk}}(w_{0j} ) +( h_{jk}'-h_{jk})w_{0j}\right\}= \mathcal{L}_{b_{jk}}\left(a_j+ \Psi_{0jk}a_k \right)+ c_{jk}\left( a_j+ \Psi_{0jk}a_k \right) = 2\mathcal{L}_{b_{jk}}(a_j)+ 2c_{jk} a_j \label{um9}\\
&\delta\left(\left\{ (e_{jk}'-e_{jk}, h_{jk}'-h_{jk})\right\} \right)=0 \label{um10}
\end{align}
We note that $(\ref{um10})$ implies that $\left\{ \left( e_{jk}'- e_{jk}  , h_{jk}'- h_{jk} \right) \right\} \in H^1\left(M, \mathcal{E}_{\mathcal{L}_0} \right)$. Then $(\ref{um9})$ implies that $\left\{\mathcal{L}_{b_{jk}}(a_j)+ c_{jk}a_j  \right\}$ is cohomogloous to the image of $H^1\left(M, \mathcal{E}_{\mathcal{L}_0} \right)\xrightarrow{\tilde{E}_0} H^1\left(M, \Omega_M^1 \otimes \mathcal{L}_0 \right)$. This proves $(1)$.

We prove $(2)$. Since $( \{a_j\} , \{(0,0)\})$ defines an element $\mathbb{H}^1\left( M, \mathcal{E}_{\mathcal{L}_0}^\bullet \right)$. If $\left(M, w_0 \right)$ is unobstructed in unfoldings, by Theorem \ref{on27}, there exists $(\{g_i\}, \{e_{jk}, h_{jk}\}) \in C^0\left( \mathcal{U},  \Omega_M^1\otimes \mathcal{L}_0 \right) \bigoplus C^1\left( \mathcal{U} , \mathcal{E}_{\mathcal{L}_0} \right)$ such that
\begin{align*}
\delta\left(\{e_{jk}, h_{jk}\}\right)=0,\,\,\,\left\{ \Psi_{0jk}g_k- g_j + \mathcal{L}_{e_{jk}}(w_{0j})+ h_{jk}w_{0j}  \right\} =0,\,\,\, \tilde{E}_1\left(\{g_i \}\right)= \left\{ 2a_j \wedge da_j\right\}
\end{align*}
Now assume that $H^1\left( M, \mathcal{E}_{\mathcal{L}_0}  \right)=0$. Then there exists $(e_i, h_i)\in C^0\left( M, \mathcal{E}_{\mathcal{L}_0} \right)$ such that $\delta(\left\{ e_i, h_i\right\})=\left\{(e_{jk}, h_{jk}) \right\}$. In other words,
\begin{align*}
e_{jk}= e_k- e_j,\,\,\,\,\,\,\,\, h_{jk}= h_k -\left[ e_k , \log \Psi_{0jk}\right] - h_j
\end{align*}
 Then $ \left\{ g_i +\mathcal{L}_{e_i}(w_{0i})+ h_iw_{0i} \right\}$ defines a global section $H^0\left( M, \Omega_M^1\otimes \mathcal{L}_0 \right))$. In fact,
\begin{align*}
& \Psi_{0jk}\left( g_k+ \mathcal{L}_{e_k}(w_{0k}) + h_kw_{0k} \right) - \left(g_j+ \mathcal{L}_{e_j}(w_{0j})+ h_jw_{0j}\right)\\
&= - \mathcal{L}_{e_{jk}}(w_{0j})- h_{jk}w_{0j} + \Psi_{0jk}\mathcal{L}_{e_k}(w_{0k}) + h_k w_{0j}- \mathcal{L}_{e_j}(w_{0j})- h_j w_{0j} \\
&= - \mathcal{L}_{e_k}(w_{0j}) + \mathcal{L}_{e_j}(w_{0j}) - h_k w_{0j} + \left[e_k, \log \Psi_{0jk}\right]w_{0j} + h_j w_{0j} + \mathcal{L}_{e_k}(w_{0j})- \left[e_k, \Psi_{0jk}\right] w_{0k}+ h_k w_{0j}- \mathcal{L}_{e_j}(w_{0j})- h_jw_{0j} \\
&= 0 
\end{align*}
Then we note that
\begin{align*}
\tilde{E}_1\left( \left\{g_i + \mathcal{L}_{e_i}(w_{0i})+ h_i w_{0i} \right\} \right) = \tilde{E}_1\left( \{g_i\}\right)= \left\{ 2a_j\wedge da_j\right\}
\end{align*}
This implies $(2)$.

\end{proof}

\subsection{Integrability conditions}\label{oni14}\

Let $\left(\mathcal{M}, w, B, \pi \right)$ be a foliated complex analytic family with $w\in H^0\left(\mathcal{M}, \Omega_{\mathcal{M}/B}^1 \otimes \mathcal{L} \right)$, where $\mathcal{L}= \left(\mathcal{N}_\mathcal{F}^* \right)^* $ is a line bundle on $\mathcal{M}$ of deformations of $\pi^{-1}(0)=\left(M, w_0 \right)$. We keep the notations in subsection \ref{on25}. By \cite{Kod05} Theorem 2.3, $M_t$ is diffeomorphic to $M_0=\pi^{-1}(0)$ as differentiable manifolds for each $t\in \Delta$. By \cite{Kod05} Theorem 2.5, if we take a sufficiently small $\Delta$, there is a diffeomorphism $\Phi$ of $M\times \Delta$ onto $\mathcal{M}_\Delta$ as differentiable manifolds such that $\pi\circ \Phi$ is the projection $M\times \Delta \to \Delta$. Let $z=\left(z^1,..., z^n \right)$ be local coordinates of $M=M_0$. Then we have $\pi\circ \Phi(z,t)=t, t\in \Delta$. For $\Phi(z,t)\in U_j\times \Delta$, put
\begin{align*}
\Phi\left( z, t\right) = \left(\xi_j^1(z,t), ..., \xi_j^n(z,t), t_1,..., t_m \right)
\end{align*}
Then each component $\xi_j^\alpha=\xi_j^\alpha(z,t),\alpha=1,...,n$ is a $C^\infty$ function. If we identify $\mathcal{M}_\Delta =\Phi\left(M\times \Delta \right)$ with $M\times \Delta$, and $\mathcal{L}$ with $\Phi^*\mathcal{L}$ via $\Phi$, $\left(\mathcal{M}_\Delta, \mathcal{L}_\Delta \right)$ is considered as a complex manifold with the complex structure defined on the $C^\infty$ manifold $M\times \Delta$ by the system of local coordinates on $U_j\times \Delta$
\begin{align*}
\left\{ \left(\xi_j,t \right)| j=1 ,2, 3,...   \right\},\,\,\,\,\,(\xi_j,t)=\left( \xi_j^1(z,t), ..., \xi_j^n(z,t), t_1,..., t_m \right)
\end{align*}
and the $C^\infty$ differentiable line bundle $\Phi^* \mathcal{L}_\Delta$ on $M\times \Delta$.

We note that since $(z^1,..., z^n)$ and $\left(\xi_j^1(z,0), ..., \xi_j^n(z,0)\right)$ are local complex coordinates on $M=M_0$, $\xi_j^\alpha(z,0)$ are holomorphic functions of $z^1,...,z^n,\alpha=1,...,n$. We also note that if we take $\Delta$ sufficiently small, we have
\begin{align}\label{oci3}
\det\left( \frac{\partial \xi_j^\alpha(z,t)}{\partial z^\lambda} \right)_{\alpha, \lambda=1,...,n} \ne 0\,\,\,\,\,\, \, \textnormal{for}\,\,\,t\in \Delta
\end{align}
On the other hand, each line bundle $\Phi^* \mathcal{L}_t$ in a differentiable  family $\left\{ \Phi^* \mathcal{L}_t |t\in \Delta \right\}$ is the restriction to $M\times t$ of the line bundle $\Psi^* \mathcal{L}$ over $M\times \Delta , \Phi^* \mathcal{L}_t = \Phi^* \mathcal{L} |_{M\times t}$. Then by \cite{Kod05} Lemma 7.1 the $C^\infty$ line bundle $\Phi^* \mathcal{L}$ over $M\times \Delta$ is $C^\infty$-equivalent to the trivial line bundle $\mathcal{L}_0 \times \Delta $. In other words, there is a bundle diffeomorphism $\sigma: \mathcal{L}_0\times \Delta \to \Phi^* \mathcal{L}$ over $M\times \Delta$. We recall that $\mathcal{L}$ is described by the transition functions $\left\{\Psi_{jk}(\xi_j,t)\right\}$ with respect to the covering $\left\{U_j\times \Delta \right\}$ of $\mathcal{M}_\Delta$. Then $\Phi^* \mathcal{L}$ is described by the transition matrices $\left\{H_{jk}:= \Psi_{jk}\left(\xi_k(z,t) , t\right) \right\}$ with respect to the covering $\left\{\mathscr{U}_j:=\Psi^{-1}\left( U_j\times \Delta \right) \right\}$ of $M\times \Delta$. In particular $\mathcal{L}_0$ is described by $\left\{ \Psi_{0jk}:=\Psi_{jk}(\xi_j(z,0),0 )\right\}$. Then there are non-vanishing $C^\infty$ functions $g_j(z,t)$ on $\mathscr{U}_j$ with $g_j(z,0)=1$ such that $\sigma$ is represented locally by
\begin{align*}
\sigma:\mathbb{C}\times\mathscr{U}_j &\to \mathbb{C} \times \mathscr{U}_j\\
                         \left( X, z, t \right) &\mapsto \left( g_j(z,t)\cdot X, z, t\right)
\end{align*}
such that
\begin{align}\label{otci1}
g_j(z,t) \cdot   \Psi_{jk}\left(\xi_k(z,0),0 \right) =  \Psi_{jk}\left(\xi_k(z,t), t \right) \cdot g_k(z,t) \iff g_j \cdot \Psi_{0jk} = H_{jk}\cdot g_k
\end{align}
From $(\ref{oci3})$, we can define a $(0,1)$-form $\varphi(t):=\sum_{\lambda,\lambda=1}^n \varphi_{jv}^\lambda(z,t)d\bar{z}^v\frac{\partial}{\partial z^\lambda}\in A^{0,1}\left(M, \Theta_M \right)$ for $t\in \Delta$ such that $\bar{\partial}\xi_j^\alpha - \sum_{\lambda=1}^n \varphi_{jv}^\lambda(z,t)\frac{\partial \xi_j^\alpha}{\partial z_\lambda}d\bar{z}^v$ and
\begin{align}
\varphi(0)=0,\,\,\,\,\,\,\, \bar{\partial} \varphi(t)- \frac{1}{2} \left[ \varphi(t), \varphi(t) \right]=0
\end{align}
(see \cite{Kod05} p.262). On the other hand, as in $(\ref{oi30})$ in subsection \ref{oi31}, we can show
\begin{align}
[\varphi, \log g_j]- \bar{\partial}\log g_j= [\varphi, \log g_k]- \bar{\partial} \log g_k- [ \varphi , \log \Psi_{0jk}]
\end{align}
This implies that
\begin{align*}
\left( (\varphi, Z_j:= \left\{\left[\varphi, \log g_j\right] - \bar{\partial} \log g_j \right\} \right) \in A^{0,1}\left( M, \mathcal{E}_{\mathcal{L}_0} \right)
\end{align*}
and we have
\begin{align}
\bar{\partial} Z_j - \left[\varphi, Z_j\right] = 0 
\end{align}
On the other hand, we set
\begin{align*}
W_j(z,t)= \sum_{\beta=1}^n W_j^\beta (z,t) dz^\beta:= \sum_{\alpha, \beta=1}^n w_j^\alpha(\xi_j(z,t),t)\frac{\partial \xi_j^\alpha}{\partial z_i^\beta}dz^\beta
\end{align*}
Then since $\sum_{\alpha=1}^n w_j^\alpha( \xi_j,t)d\xi_j^\alpha= \Psi_{jk}(\xi_k, t) \sum_{\alpha=1}^n w_k^\alpha(\xi_k,t)d\xi_k^\alpha$, we have
\begin{align*}
&\sum_{\alpha=1}^n w_j^\alpha\left(\xi_j(z,t),t \right)\frac{\partial \xi_j^\alpha}{\partial z^\beta} dz^\beta= \Psi_{jk}(\xi_k(z,t)) \sum_{\alpha=1}^n w_k^\alpha\left( \xi_k(z,t),t \right)\frac{\partial \xi_k^\alpha}{\partial z^\beta} dz^\beta \\
& \Longrightarrow W_j(z,t)=H_{jk}(z,t) W_k(z,t) \Longrightarrow W_j(z,t)= g_j \cdot \Psi_{0jk}\cdot g_k^{-1} W_k(z,t) \Longrightarrow g_j^{-1} W_j = \Psi_{0jk}\cdot g_k^{-1}W_k
\end{align*}
This implies that
\begin{align*}
\left\{ Y_j(t):= g_j^{-1} W_j   \right\} \in A^{0,0}\left( M, \Omega_M^1  \otimes \mathcal{L}_0 \right)
\end{align*}
We also note that from $w_j\wedge dw_j=0$, we have
\begin{align}
Y_j \wedge \partial Y_j=0
\end{align}

Since $\bar{\partial} W_j(z,t)- \mathcal{L}_\varphi\left(W_j(z,t) \right)=0$, we have
\begin{align*}
\bar{\partial} Y_j - \mathcal{L}_\varphi \left(Y_j\right)= \bar{\partial} \left(g_j^{-1} W_j \right) - \mathcal{L}_\varphi  \left(  g_j^{-1} W_j  \right) = - \frac{1}{g_j^2} \bar{\partial} g_j \cdot W_j + \frac{1}{g_j^2}[\varphi,  g_j]\cdot W_j = \left( - \bar{\partial} \log g_j + \left[\varphi, \log g_j \right]\right) Y_j
\end{align*}
Hence we have
\begin{align*}
\bar{\partial} Y_j - \mathcal{L}_{\varphi} \left(Y_j\right) - Z_j Y_j=0
\end{align*}

All in all, we have $\left(\left\{Y_j\right\}, \left(\varphi, \left\{ Z_j \right\}\right)\right)\in A^{0,0}\left( M, \Omega_M^1\otimes \mathcal{L} \right) \bigoplus A^{0,1}\left( M, \mathcal{E}_\mathcal{L}  \right) $, where
\begin{align*}
&\varphi \in A^{0,1}\left(M, \Theta_M\right) \\
&Y_j\in \Gamma\left( U_{j}, \mathcal{A}^{0,0}\left(\Omega_M^1\right)\right)\,\,\,\,\,\,\,\,\,\,\textnormal{with} \,\,\,\,\,\,\, Y_j= \Psi_{0jk} Y_k\\
&Z_j \in \Gamma\left(U_j, \mathcal{A}^{0,1}\right)\,\,\,\,\,\,\,\,\,\,\,\textnormal{with}\,\,\,\,\,Z_j- Z_k+ \left[\varphi, \log \Psi_{0jk}\right]=0
\end{align*}

 with respect to the covering $\mathcal{U}_0=\{\mathcal{U}_{j}\cap M\}=\{U_{j}\}$ such that  $\varphi(0)=0, Z_j(0)=0$ and $Y_j(0)=w_{0j}$ and
\begin{align}
\bar{\partial} \varphi - \frac{1}{2}[\varphi, \varphi] = 0  \label{oni1}\\
\bar{\partial} Y_j - \mathcal{L}_\varphi(Y_j) - Z_j   Y_j =0 \label{oni2}\\
Y_j\wedge \partial Y_j =0 \label{oni3}\\
\bar{\partial} Z_j - [\varphi, Z_j]=0 \label{oni4}
\end{align}
We take the derivative of $(\ref{oni1})-(\ref{oni4})$ with respect to $t$:
\begin{align}
&\bar{\partial} \left(\frac{\partial \varphi}{\partial t} \right) - \frac{1}{2}\left[\frac{\partial \varphi}{\partial t}, \varphi \right] - \frac{1}{2}\left[\varphi, \frac{\partial \varphi}{\partial t} \right] = 0 \label{oni5}\\
& \bar{\partial}\left( \frac{\partial Y_j}{\partial t} \right) - \mathcal{L}_{\frac{\partial \varphi}{\partial t}}\left( Y_j \right) - \mathcal{L}_{\varphi}\left( \frac{\partial Y_j}{\partial t} \right) - \frac{\partial Z_j}{\partial t} \cdot Y_j - Z_j \cdot \frac{\partial Y_j}{\partial t} = 0 \label{oni6} \\
&\frac{\partial Y_j}{\partial t} \wedge \partial Y_j + Y_j \wedge \partial \left( \frac{\partial Y_j}{\partial t} \right) = 0 \label{oni7}\\
&\bar{\partial}\left( \frac{\partial Z_j}{\partial t} \right) - \left[ \frac{\partial \varphi}{\partial t}, Z_j   \right] - \left[ \varphi, \frac{\partial Z_j}{\partial t}  \right] = 0 \label{oni8}
\end{align}
By setting $t=0$ of $(\ref{oni5})-(\ref{oni8})$, we have
\begin{align*}
&-\bar{\partial}\left( \frac{\partial \varphi}{\partial t}|_{t=0} \right) = 0 \\
&-\bar{\partial}\left( \frac{\partial Y_j}{\partial t}|_{t=0} \right) + \mathcal{L}_{\frac{\partial \varphi}{\partial t}|_{t=0}}(w_j) +\frac{\partial Z_j}{\partial t}|_{t=0} \cdot w_j =0\\
&\frac{\partial Y_j}{\partial t}|_{t=0} \wedge d w_{0j} + w_{0j} \wedge \partial \left( \frac{\partial Y_j}{\partial t}|_{t=0} \right) =0\\
&-\bar{\partial}\left( \frac{\partial Z_j}{\partial t}|_{t=0} \right)=0
\end{align*}
This implies that
\begin{align*}
\left( \left\{ \frac{\partial Y_j}{\partial t}|_{t=0}\right\}, \left\{ \left( \frac{\partial \varphi}{\partial t}|_{t=0}, \frac{\partial Z_j}{\partial t}|_{t=0} \right) \right\}  \right) \in A^{0,0}\left( M , \Omega_M^1\otimes \mathcal{L}_0 \right) \bigoplus A^{0,1}\left(M, \mathcal{E}_{\mathcal{L}_0} \right)
\end{align*}
defines a $1$-cocycle in the complex associated with the following double complex
{\small{\begin{equation}\label{oni9}
\begin{CD}
\cdots \\
@AE_2AA \\
A^{0,0}\left(M- S, \tilde{\mathcal{S}}^2\otimes \mathcal{L}_0 \right) @>\bar{\partial}>> \cdots  \\
@AE_1'AA @AE_1'AA  \\
A^{0,0}\left( M, \Omega_M^1 \otimes \mathcal{L}_0 \right)@>-\bar{\partial}>> A^{0,1}\left( M, \Omega_M^1\otimes \mathcal{L}_0 \right) @>\bar{\partial}>> \cdots \\
@A\tilde{E}_0AA @A\tilde{E}_0AA @A\tilde{E}_0AA  \\
A^{0,0}\left(M, \mathcal{E}_{\mathcal{L}_0} \right) @>\bar{\partial}>>  A^{0,1}\left( M, \mathcal{E}_{\mathcal{L}_0} \right) @>-\bar{\partial}>> A^{0,2}\left( M, \mathcal{E}_{\mathcal{L}_0} \right) @>\bar{\partial}>> \cdots
\end{CD}
\end{equation}}}

Then we have the foliated Kodaira-Spencer map of $\left(\mathcal{M}, w, B, \pi \right)$ at $t=0$ in terms of sections
\begin{align}\label{cu14}
\varphi_0: T_0 B &\to \mathbb{H}^1\left( M, \mathcal{E}_{\mathcal{L}_0}^\bullet \right)\\
 \frac{\partial}{\partial t} &\mapsto \left( \frac{\partial \left( \mathcal{M}, w \right)}{\partial t}\right)|_{t=0}:= \left( \left\{ \frac{\partial Y_j}{\partial t}|_{t=0}  \right\}, \left\{ \left( \frac{\partial \varphi}{\partial t}|_{t=0},  \frac{\partial Z_j}{\partial t}|_{t=0} \right) \right\} \right) \notag
\end{align}

Now we assume that the foliated analytic family $\left(\mathcal{M}, w, B, \pi \right)$ with $w \in  H^0\left(\mathcal{M}, \Omega_{\mathcal{M}/B}^1\otimes \mathcal{L}  \right)$ in terms of sections is a $1$-parameter family, i.e. $0\in \Delta \subset B\subset \mathbb{C}$. By taking the derivative of $(\ref{oni5})-(\ref{oni8})$ with respect to $t$ again and setting $t=0$, we have
\begin{align*}
&\bar{\partial}\left( \frac{\partial^2 \varphi}{\partial t^2}|_{t=0} \right) - \left[ \frac{\partial \varphi}{\partial t}|_{t=0}, \frac{\partial \varphi}{\partial t}|_{t=0} \right]=0\\
&\bar{\partial}\left( \frac{\partial^2 Y_j}{\partial t^2} |_{t=0} \right) - \mathcal{L}_{\frac{\partial^2 \varphi}{\partial t^2}|_{t=0}}\left(w_{0j} \right) -2 \mathcal{L}_{\frac{\partial \varphi}{\partial t}|_{t=0}}\left( \frac{\partial Y_j}{\partial t}|_{t=0} \right) - \frac{\partial^2 Z_j}{\partial t^2}|_{t=0}\cdot w_{0j} - 2 \frac{\partial Y_j}{\partial t}|_{t=0} \cdot \frac{\partial Y_j}{\partial t}|_{t=0} = 0 \\
& \left( \frac{\partial^2 Y_j}{\partial t^2}|_{t=0}\right) \wedge dw_{0j} + 2 \left(\frac{\partial Z_j}{\partial t}|_{t=0}\right) \wedge \partial \left( \frac{\partial Y_j}{\partial t}|_{t=0} \right) +  w_{0j} \wedge \partial \left(\frac{\partial^2 Y_j}{\partial t^2}|_{t=0} \right) =0\\
& \bar{\partial} \left(  \frac{\partial^2 Z_j}{\partial t^2}|_{t=0} \right)  - 2 \left[ \frac{\partial \varphi}{\partial t}|_{t=0}, \frac{\partial Z_j}{\partial t}|_{t=0} \right]=0
\end{align*}
Then
{\small{\begin{align}\label{oni10}
&\left\{ - 2 \left(\frac{\partial Y_j}{\partial t}|_{t=0} \right) \wedge \partial \left( \frac{\partial Y_j}{\partial t}|_{t=0}  \right) \right\} \bigoplus \left\{ -2 \mathcal{L}_{\frac{\partial \varphi}{\partial t}|_{t=0}}\left( \frac{\partial Y_j}{\partial t}|_{t=0} \right)  - 2 \frac{\partial Z_j}{\partial t}|_{t=0} \cdot \frac{\partial Y_j}{\partial t}|_{t=0} \right\}\\
&\bigoplus \left( \left\{\left( -\left[\frac{\partial \varphi}{\partial t}|_{t=0}, \frac{\partial \varphi}{\partial t}|_{t=0}\right]   , -2 \left[ \frac{\partial \varphi}{\partial t}|_{t=0}, \frac{\partial Z_j}{\partial t}|_{t=0} \right] \right)   \right\} \right) \in A^{0,0}\left( M-S, \tilde{\mathcal{S}}^2\otimes \mathcal{L}_0  \right) \bigoplus A^{0,1}\left( M, \Omega_M^1\otimes \mathcal{L}_0 \right) \bigoplus A^{0,2}\left( M, \mathcal{E}_{\mathcal{L}_0}\right) \notag
\end{align}}}
is a $2$-cocycle in the above bicomplex $(\ref{oni9})$, i.e the  image of $\left( \left\{ \frac{\partial^2 Y_j}{\partial t^2}|_{t=0} \right\} ,  \left\{ \left(\frac{\partial^2 \varphi}{\partial t}|_{t=0}, \frac{\partial^2 Z_j}{\partial t^2}|_{t=0} \right) \right\} \right)$.

For any $1$-cocycle $(\{a_j\}, \{\psi, \{c_j\}\})\in A^{0,0}\left( M, \Omega_M^1 \otimes \mathcal{L}_0 \right) \bigoplus A^{0,1}\left( M, \mathcal{E}_{\mathcal{L}_0} \right)$ in the bicomplex $(\ref{oni9})$, we can define a $2$-cocycle by (we omit the proof, see Part VI-1)
\begin{align}\label{oni11}
\left( \left\{ 2 a_j \wedge \partial a_j\right\},  \left\{ 2\mathcal{L}_\psi\left(a_j\right)+2c_j \cdot a_j \right\} , \left\{ \left( [\psi, \psi], 2 [\psi, c_j] \right) \right\} \right) \in A^{0,0}\left( M-S, \tilde{\mathcal{S}}^2\otimes \mathcal{L}_0 \right) \bigoplus A^{0,1}\left( M, \Omega_M^1 \otimes \mathcal{L}_0 \right) \bigoplus A^{0,2}\left( M, \mathcal{E}_{\mathcal{L}_0} \right)
\end{align}
Then $(\ref{oni10})$ implies that

\begin{theorem}
Suppose given a compact complex manifold $M$ with a codimension $1$-foliation defined by $w_0\in H^0\left(M, \Omega_M^1\otimes \mathcal{L}_0 \right)$ for a line bundle $\mathcal{L}_0$, and $\left(\{a_j\}, \{\psi, \{c_j\}\}\right)\in A^{0,0}\left( M, \Omega_M^1 \otimes \mathcal{L}_0 \right) \bigoplus A^{0,1}\left( M, \mathcal{E}_{\mathcal{L}_0} \right)\in \mathbb{H}^1\left(M, \mathcal{E}_{\mathcal{L}_0}^\bullet \right)$. In order that there may exist a foliated analytic family $\left( \mathcal{M}, w, B, \pi\right)$ such that $\pi^{-1}(0)=\left( M, w_0 \right)$ and $\left(\frac{\partial \left(\mathcal{M}, w \right)}{\partial t}\right)|_{t=0}=\left(\left\{a_j \right\}, \left\{\left( \psi , c_{j} \right) \right\} \right)$, it is necessary  that $(\ref{oni11})$ defines the $0$ class in $\mathbb{H}^2\left( M, \mathcal{E}_{\mathcal{L}_0}^\bullet \right)$.

\end{theorem}

\section{Infinitesimal unfoldings and obstructions of codimension 1 foliations in terms of deformations of sections}

\subsection{Infinitesimal unfoldings of foliated complex analytic structures of codimension $1$ in terms of sections} \label{uti1}\

Let $\left(M, \mathcal{N}_{\mathcal{F}_0}^* \right)$ be a compact foliated complex manifold with $\mathcal{N}_{\mathcal{F}_0}^*$ locally free of rank of $1$, so that there exists a section $w_0\in H^0\left( M, \Omega_M^1 \otimes {\mathcal{L}_0} \right)$ with $\mathcal{L}_0 = \left( \mathcal{N}_{\mathcal{F}_0}^*\right)^*$, which defines the foliation on $M$. In this subsection, we will consider unfoldings of $\left( M, \mathcal{N}_{\mathcal{F}_0}^* \right)$ as unfoldings of the section $w_0$ (see Part III). Then we have the complex of sheaves controlling unfoldings of $\left( M, \mathcal{N}_{\mathcal{F}_0}^* \right)$ in terms of sections:
\begin{align}\label{cp3}
\tilde{\mathcal{E}}_{\mathcal{L}_0}^\bullet : \mathcal{E}_{\mathcal{L}_0} \xrightarrow{\tilde{G}_0} \left(\Omega_M^1\otimes \mathcal{L}_0 \right) \bigoplus \mathcal{L}_0 \xrightarrow{\tilde{G}_1} \left(\tilde{\mathcal{S}}^2\otimes \mathcal{L}_0 \right) \bigoplus \left(\tilde{\mathcal{S}}^1\otimes \mathcal{L}_0 \right) \xrightarrow{G_2}  \left(\tilde{\mathcal{S}}^3\otimes \mathcal{L} \right) \bigoplus \left( \tilde{\mathcal{S}}^2\otimes \mathcal{L}_0 \right) \xrightarrow{G_3} \cdots
\end{align}
We recall that $\mathcal{E}_{\mathcal{L}_0}$ and the differential $\tilde{G}_0, \tilde{G}_1
$ are defined in the following way: for an open covering $\mathcal{U}_0=\{U_j\}$ of $M$ by coordinate neighborhoods, let $\mathcal{L}_0$ defined by the transition functions $\left\{ \Psi_{0jk} \right\}\in C^1 \left(\mathcal{U}_0, \mathcal{O}_M^* \right)$ with $\Psi_{0ij}\Psi_{0jk}=\Psi_{0ik}$, and $w_0\in \Gamma\left( M, \Omega_M^1 \otimes \mathcal{L}_0 \right)$ be written as $w_{0j}\in \Gamma\left( U_j, \Omega_M^1 \right)$  on $U_j$ with $w_{0j}= \Psi_{0jk} w_{0k}$ on $U_j\cap U_k$.  Then $\mathcal{E}_{\mathcal{L}_0} |_{U_j}$ is isomorphic to $\Theta_{M}|_{U_j}\oplus \mathcal{O}_{M}|_{U_j}$, and $(\rho_j, a_j)\in \Theta_{M}|_{U_j}\oplus \mathcal{O}_{M}|_{U_j}$ is identified with $\left(\rho_k, a_k\right)\in \Theta_{M}|_{U_k}\oplus \mathcal{O}_{M}|_{U_k}$ on $U_j\cap U_k$ if and only if $\rho_j=\rho_k$ and $a_k-a_j=\left[\rho_j, \log \Psi_{0jk}\right]$. On the other hand, the differential $\tilde{G}_0$ is defined in the following way:
\begin{align*}
\tilde{G}_0: \Gamma\left(U_j,\mathcal{E}_{\mathcal{L}_0}\right)\cong \Gamma\left(U_j,\Theta_{M}\right)\oplus \Gamma\left(U_j ,  \mathcal{O}_{M} \right) &\to \Gamma\left(U_j, \Omega_{M}^1 \otimes  \mathcal{L}_0 \right) \bigoplus \Gamma\left( U_j, \mathcal{L}_0 \right)\\
     \left(\rho_j, a_j\right)&\mapsto \left( \mathcal{L}_{\rho_j}\left( w_{0j}\right) +a_j w_{0j} , i_{\rho_j} w_{0j}\right)
\end{align*}
and $\tilde{G}_1$ is defined in the following way:
\begin{align*}
\tilde{G}_1: \Gamma\left( U_j, \Omega_M^1\otimes \mathcal{L}_0 \right)\bigoplus \Gamma\left( U_j, \mathcal{L}_0 \right)  &\to \Gamma\left( U_j,  \tilde{\mathcal{S}}^2 \otimes \mathcal{L}_0 \right) \bigoplus \Gamma\left( U_j, \tilde{\mathcal{S}}^1\otimes \mathcal{L}_0 \right) \\
            \left(  v_j , r_j \right)&\mapsto \left( v_j \wedge dw_0 + w_0\wedge dv_j , - v_j + r_j dw_0 + w_0 \wedge dr_j \right)
\end{align*}
We will denote $i$-th cohomology group of $\tilde{\mathcal{E}}_{\mathcal{L}_0}^\bullet$ by $\mathbb{H}^i\left( M, \tilde{\mathcal{E}}_{\mathcal{L}_0}^\bullet \right)$. 

We shall identify infinitesimal unfoldings of codimension $1$-foliations in terms of sections.  Let $\left(\mathcal{M}, \tilde{\mathcal{N}}_{\mathcal{F}}^*, B, \pi \right)$ with $\tilde{\mathcal{N}}_{\mathcal{F}}^*$ locally free of rank $1$ be an unfolded analytic family of unfoldings of $\left( M, \mathcal{N}_{\mathcal{F}_0}^* \right)=\pi^{-1}(0), 0 \in B \subset \mathbb{C}^m$. Then the exact sequence $0 \to \tilde{\mathcal{N}}_\mathcal{F}^* \to \Omega_{\mathcal{M}}^1 \to \frac{\Omega_{\mathcal{M}}^1}{ \tilde{\mathcal{N}}_{\mathcal{F}}^*}\to 0 $ gives a section $\omega \in H^0\left(\mathcal{M}, \Omega_{\mathcal{M}}^1\otimes \mathcal{L} \right)$ with $\mathcal{L}= \left( \tilde{\mathcal{N}}_{\mathcal{F}}^*\right)^*$ which induces $w_0\in H^0\left( M, \Omega_M^1 \otimes \mathcal{L}_0 \right)$. We will also denote the family by $\left( \mathcal{M}, \omega , B, \pi \right)$. Take a small disk $\Delta$ with centre $0$ such that $0\in \Delta \subset B$, and represent $\mathcal{M}_\Delta=\pi^{-1}(\Delta)$ in the form $\mathcal{M}_\Delta = \bigcup_{j=1}^l U_j\times \Delta$ where each $U_j$ is a polydisk, and $(\xi_j,t)\in U_j\times \Delta$ and $(\xi_k,t)\in U_k\times \Delta$ are the same points on $\mathcal{M}_\Delta$ if $\xi_j^\alpha= f_{jk}^\alpha(\xi_k,t)$ for $\alpha=1,...,n$. Here each $f_{jk}^\alpha(\xi_k,t)=f_{jk}^\alpha(\xi_k^1,..., \xi_k^n, t)$ is a holomorphic function of $\xi_k^1,..., \xi_k^n,t$ defined on $\left(U_k\times \Delta \right) \cap \left( U_j\times \Delta \right)\ne \emptyset$, and we have 
\begin{align}\label{ut2}
f_{ik}^\alpha(\xi_k,t)&=f_{ij}^\alpha\left(f_{jk}(\xi_k,t),t \right),\,\,\,\alpha=1,...,n  \,\,\,\,\,\,\textnormal{on}\,\,\, \left(U_i\times \Delta \right)\cap \left(U_j\times \Delta \right) \cap \left( U_k \times \Delta \right)\ne \emptyset
\end{align}
 On the other hand, the line bundle $\mathcal{L}$ on $\mathcal{M}_\Delta$ is defined by the system of transition functions $\left\{\Psi_{jk}(\xi_k,t) \right\}$ where $\Psi_{jk}(\xi_k,t)\in   \Gamma\left(\left(U_j\times \Delta \right) \cap \left( U_k\times \Delta \right), \mathcal{O}_\mathcal{M}^* \right)$ satisfying the cocycle condition 
\begin{align}\label{ut3}
\Psi_{ij}\left(f_{jk}(\xi_k,t), t \right)\cdot \Psi_{jk}(\xi_k,t) = \Psi_{ik}(\xi_k,t), \,\,\,\,\,\textnormal{on}\,\,\,\,\, \left(U_i\times \Delta \right) \cap \left( U_j\times \Delta \right) \cap \left( U_k \times \Delta \right) \ne \emptyset
\end{align}
and $\omega\in H^0\left( \mathcal{M}, \Omega_{\mathcal{M}/B}^1\otimes \mathcal{L} \right)$ on $U_j\times \Delta$ is of the form
\begin{align*}
\omega|_{U_j\times \Delta}=\omega_j(z_j,t)=w_j(\xi_j,t) + v_j(\xi_j,t):=\sum_{\alpha=1}^n w_j^\alpha(\xi_j,t)d\xi_j^\alpha + \sum_{\gamma=1}^m v_j^\gamma(\xi_j,t)dt_\gamma
\end{align*}
such that $w_j(z_j,t)+ v_j(z_j,t)=\Psi_{jk}(z_k,t)\left( w_k(z_k,t), v_k(z_k,t) \right)$, in other words,
\begin{align*}
\sum_{\alpha=1}^n w_j^\alpha(\xi_j,t) d\xi_j^\alpha + \sum_{\gamma=1}^m v_j^\gamma(\xi_j,t) dt_\gamma  = \Psi_{jk}(\xi_k,t)\left( \sum_{\beta=1}^n w_k^\beta(\xi_k,t) d\xi_k^\beta + \sum_{\eta=1}^m v_k^\eta(\xi_k,t) dt_\eta \right)\,\,\,\,\,\,\,\textnormal{on}\,\,\, \left(U_j\times \Delta\right) \cap \left( U_k \times \Delta \right)
\end{align*}
Equivalently we have
{\small{\begin{align*}
\sum_{\alpha,\beta=1}^n w_j^\alpha(\xi_j,t)\frac{\partial f_{jk}^\alpha}{\partial \xi_k^\beta}d\xi_k^\beta + \sum_{\alpha=1}^n \sum_{\eta=1}^m w_j^\alpha(\xi_j,t)\frac{\partial f_{jk}^\alpha}{\partial t_\eta} dt_\eta + \sum_{\eta=1}^m v_j^\eta(\xi_j,t) dt_\eta= \Psi_{jk}(\xi_k,t)\left(\sum_{\beta=1}^n w_k^\beta(\xi_k,t)d\xi_k^\beta + \sum_{\eta=1}^m v_k^\eta(\xi_k,t)dt_\eta \right)
\end{align*}}}
Then we have
\begin{align}\label{ut4}
\sum_{\alpha=1}^n w_j^\alpha\left(f_{jk}(\xi_k,t), t\right)\frac{\partial f_{jk}^\alpha}{\partial \xi_k^\beta} = \Psi_{jk}(\xi_k,t)w_k^\beta(\xi_k,t),\,\,\,\,\,\,\beta=1,...,n
\end{align}
and
\begin{align}\label{ut16}
\sum_{\alpha=1}^n w_j^\alpha\left(f_{jk}(\xi_k,t),t \right)\frac{\partial f_{jk}^\alpha}{\partial t_\eta} + v_j^\eta \left(f_{jk}( \xi_k,t ) , t \right) = \Psi_{jk} (\xi_k,t) v_k^\eta(\xi_k,t),\,\,\,\,\,\,\eta=1,...,m
\end{align}
and we have the integrability condition $\omega_j(\xi_j,t)\wedge d\omega_j(\xi_j,t)=0$, i.e.
\begin{align*}
&\left(\sum_{\alpha=1}^n w_j^\alpha(\xi_j,t) d\xi_j^\alpha + \sum_{\gamma=1}^m v_j^\gamma(\xi_j,t) dt_\gamma\right) \wedge d \left( \sum_{\alpha=1}^n w_j^\alpha(\xi_j,t) d\xi_j^\alpha + \sum_{\gamma=1}^m v_j^\gamma(\xi_j,t) dt_\gamma      \right) = 0\\
&\Longrightarrow w_j \wedge d_z w_j + \sum_{\gamma=1}^m \left(  v_j^\gamma  d_z w_j + w_j \wedge d_z v_j^\gamma - w_j \wedge \frac{\partial w_j}{\partial t_\gamma}  \right) \wedge dt_\gamma + \cdots = 0
\end{align*}

In particular,  we have
\begin{align}\label{ut9}
w_j(\xi_j,t) \wedge d_z w_j(\xi_j,t) =0 
\end{align}
and
\begin{align}\label{ut10}
v_j^\gamma dw_j+ w_j \wedge d v_j^\gamma - w_j \wedge \frac{\partial w_j}{\partial t_\gamma} =0,\,\,\,\,\,\,\,\gamma=1,...,m
\end{align}

Now we identify the infinitesimal unfoldings. By taking the derivative of $(\ref{ut2})$ with respect to $\frac{\partial}{\partial t}= \sum_{\gamma=1}^m c_\gamma \frac{\partial}{\partial t_\gamma}$ and setting $\theta_{jk}(t)=\sum_{\alpha=1}^n \theta_{jk}^\alpha(\xi_j,t)\frac{\partial}{\partial \xi_j^\alpha}$ where $\theta_{jk}^\alpha(\xi_j,t)=\frac{\partial f_{jk}(\xi_k,t)}{\partial t},\xi_k=f_{kj}(\xi_j,t)$, we have
\begin{align}\label{ut6}
\theta_{ik}^\alpha(\xi_i,t)=\theta_{ij}^\alpha(\xi_i,t) + \sum_{\beta=1}^n \frac{\partial \xi_i^\alpha}{\partial \xi_j^\beta}\theta_{jk}^\beta(\xi_j,t),\,\,\,\alpha=1,...,n
\end{align}
By setting $t=0$, we obtain 
\begin{align} \label{ut7}
\theta_{ij}(0)-\theta_{ik}(0) + \theta_{jk}(0)=0
\end{align}

As in subsection \ref{ot23}, we have
\begin{align}\label{ut5}
 -\frac{1}{\Psi_{0ij}}\frac{\partial \Psi_{ij}}{\partial t}|_{t=0} + \frac{1}{\Psi_{0ik}} \frac{\partial \Psi_{ik}}{\partial t} |_{t=0} - \frac{1}{\Psi_{jk}}\frac{\partial \Psi_{jk}}{\partial t } - \left[ \sum_{\alpha=1}^n \frac{\partial f_{jk}^\alpha}{\partial t}|_{t=0} \frac{\partial}{\partial \xi_j^\alpha} , \log \Psi_{0ij} \right] =0
\end{align}
where $\Psi_{0jk}:=\Psi_{jk}\left( \xi_k, 0\right)$, i.e. $\left\{ \Psi_{0jk} \right\}$ defines the $\mathcal{L}_0$ on $M$.
As in subsection \ref{on25}, by taking the derivative of $(\ref{ut4})$ with respect to $\frac{\partial}{\partial t}= \sum_{\gamma=1}^m c_\gamma \frac{\partial}{\partial t_\gamma}$ and setting $t=0$, we have
\begin{align}\label{ut11}
-\Psi_{0jk} \frac{\partial w_k}{\partial t}|_{t=0} + \frac{\partial w_j}{\partial t}|_{t=0} + \mathcal{L}_{\sum_{\gamma=1}^n \frac{\partial f_{jk}^\gamma}{\partial t}|_{t=0} \frac{\partial}{\partial \xi_j^\gamma}}\left(\sum_{\alpha=1}^n w_{0j}^\alpha d\xi_j^\alpha \right)- \frac{1}{\Psi_{0jk}}\frac{\partial \Psi_{jk}}{\partial t}|_{t=0} w_{0j}=0
\end{align}
where $w_{0j}=\sum_{\beta=1}^n w_{0j}^\beta d\xi_j^\beta:=\sum_{\beta=1}^n w_j^\beta(\xi_j, 0)d\xi_j^\beta$ and $\frac{\partial w_j}{\partial t}:= \sum_{\beta=1}^n \frac{\partial w_j(\xi_j,t)}{\partial t} d\xi_j^\beta$.

As in subsection \ref{on25} , by taking the derivative of $(\ref{ut9})$ and setting $t=0$, we obtain
\begin{align}\label{ut8}
\left( - \frac{\partial w_j}{\partial t}|_{t=0} \right) \wedge d w_{0j}+ w_{0j} \wedge d\left(  - \frac{\partial w_j}{\partial t}|_{t=0} \right) =0
\end{align}

From $(\ref{ut16})$, we have
\begin{align*}
\sum_{\eta=1}^m \sum_{\alpha=1}^n w_j^\alpha c_\eta \frac{\partial f_{jk}^\alpha}{\partial t_\eta} + \sum_{\eta=1}^m c_\eta v_j^\eta  = \Psi_{jk}\sum_{\eta=1}^m c_\eta v_k^\eta \Longrightarrow i_{\sum_{\alpha=1}^n \frac{\partial f_{jk}^\alpha}{\partial t} \frac{\partial}{\partial \xi_j^\alpha} }\left( w_j \right) + i_{\frac{\partial}{\partial t}}(v_j)= \Psi_{jk} i_{\frac{\partial}{\partial t}} (v_k)
\end{align*}
By setting $t=0$, we obtain
\begin{align}\label{ut21}
i_{\sum_{\alpha=1}^n \frac{\partial f_{jk}^\alpha}{\partial t}|_{t=0} \frac{\partial}{\partial \xi_j^\alpha} }\left( w_{0j} \right) + i_{\frac{\partial}{\partial t}}(v_j)|_{t=0}= \Psi_{0jk} i_{\frac{\partial}{\partial t}} (v_k)|_{t=0}
\end{align}

From $(\ref{ut10})$, we have
\begin{align}\label{ut26}
&\sum_{\gamma=1}^n c_\gamma v_j^\gamma dw_j + w_j\wedge d\left(\sum_{\gamma=1}^m c_\gamma v_j^\gamma \right) - w_j \wedge \left(\sum_{\gamma=1}^m c_\gamma \frac{\partial w_j}{\partial t_\gamma}   \right)=0 \Longrightarrow i_{\frac{\partial}{\partial t}}(v_j) dw_j + w_j \wedge d\left( i_{\frac{\partial}{\partial t}}(v_j) \right) - w_j \wedge \frac{\partial w_j}{\partial t}=0
\end{align}
By setting $t=0$, we obtain
\begin{align} \label{ut12}
- i_{\frac{\partial}{\partial t}}(v_j)|_{t=0} dw_{0j} + w_{0j} \wedge d\left(- i_{\frac{\partial}{\partial t}|_{t=0}}(v_j) \right) - w_{0j} \wedge \left( -  \frac{\partial w_j}{\partial t}|_{t=0} \right) =0
\end{align}

Then $(\ref{ut7}), ( \ref{ut5}), (\ref{ut11}), (\ref{ut8}), (\ref{ut21})$ and $(\ref{ut12})$ imply that
{\small{\begin{align*}
\left( \left\{ - \frac{\partial w_j}{\partial t}|_{t=0}  \right\} , \left\{- i_{\frac{\partial}{\partial t}}\left( v_j\right)|_{t=0} \right\},  \left\{  \sum_{\alpha=1}^n \frac{\partial f_{jk}^\alpha}{\partial t}|_{t=0} \frac{\partial}{\partial \xi_j^\alpha} , - \frac{1}{\Psi_{0jk}}\frac{\partial \Psi_{jk}}{\partial t} |_{t=0} \right\}   \right)  \in C^0\left(\mathcal{U}_0, \Omega_M^1\otimes \mathcal{L}_0     \right) \bigoplus C^0\left(\mathcal{U}_0, \mathcal{L}_0 \right) \bigoplus C^1\left( \mathcal{U}_0, \mathcal{E}_{\mathcal{L}_0} \right)
\end{align*}}}
defines a $1$-cocycle in the following \v Cech resolution of $\tilde{\mathcal{E}}_{\mathcal{L}_0}^\bullet$ with respect to the covering $\mathcal{U}_0:=\mathcal{U}\cap M$.
{\small{\begin{equation}\label{ut27}
\begin{CD}
\cdots\\ 
@A G_2AA \\
C^0\left( \mathcal{U}_0, \tilde{\mathcal{S}}^2\otimes  \mathcal{L}_0  \right) \bigoplus C^0 \left( \mathcal{U}_0, \tilde{\mathcal{S}}^1 \otimes \mathcal{L}_0 \right) @>-\delta>> \cdots \\
@A\tilde{G}_1AA  @A\tilde{G}_1AA \\
C^0\left( \mathcal{U}_0, \Omega_M^1\otimes \mathcal{L}_0 \right) \bigoplus C^0\left( \mathcal{U}_0, \mathcal{L}_0 \right) @>\delta>> C^1\left( \mathcal{U}_0, \Omega_M^1\otimes \mathcal{L}_0 \right) \bigoplus C^1\left(\mathcal{U}_0, \mathcal{L}_0 \right) @>-\delta>> \cdots \\
@A\tilde{G}_0AA @A\tilde{G}_0AA @A\tilde{G}_0AA \\
C^0\left(\mathcal{U}_0, \mathcal{E}_{\mathcal{L}_0} \right) @>-\delta>> C^1\left( \mathcal{U}_0, \mathcal{E}_{\mathcal{L}_0}\right) @>\delta>> C^2\left( \mathcal{U}, \mathcal{E}_{\mathcal{L}_0} \right) @> - \delta >> \cdots
\end{CD}
\end{equation}}}
Then we have the unfolded Kodaira-Spencer map of $\left( \mathcal{M},\omega, B, \pi \right)$ at $t=0$ in terms of sections
\begin{align*}
\varphi_0: T_0 B &\to \mathbb{H}^1\left( M, \tilde{\mathcal{E}}_{\mathcal{L}_0}^\bullet \right)\\
 \frac{\partial}{\partial t}&\mapsto \left( \frac{\partial \left(\mathcal{M}, \omega \right)}{\partial t}\right)|_{t=0} := \left( \left\{ - \frac{\partial w_j}{\partial t}|_{t=0}  \right\} , \left\{- i_{\frac{\partial}{\partial t}}\left( v_j\right)|_{t=0} \right\},  \left\{  \sum_{\alpha=1}^n \frac{\partial f_{jk}^\alpha}{\partial t}|_{t=0} \frac{\partial}{\partial \xi_j^\alpha} , - \frac{1}{\Psi_{0jk}}\frac{\partial \Psi_{jk}}{\partial t} |_{t=0} \right\}   \right) 
\end{align*}

\subsection{Obstructions of unfoldings in terms of $\tilde{\mathcal{E}}_{\mathcal{L}_0}^\bullet$}\

We assume that the unfolded analytic family $\left(\mathcal{M}, \omega, B, \pi \right)$ with $T\in  H^0\left(\mathcal{M}, \Omega_{\mathcal{M}/B}^1\otimes \mathcal{L}  \right)$ in terms of sections is a $1$-parameter family, i.e. $0\in \Delta \subset B\subset \mathbb{C}$. We recall $(\ref{ot12})$ in subsection \ref{ot23}. Then by taking the derivative of $(\ref{ut6})$ with respect to $t$, and putting
\begin{align*}
\dot{\theta}_{ij}(t)= \sum_{\alpha=1}^n \dot{\theta}_{ij}^\alpha (\xi_i,t)\frac{\partial}{\partial \xi_i^\alpha} , \,\,\,\,\,\textnormal{where}\,\,\,\dot{\theta}_{ij}^\alpha(\xi_i,t)=\frac{\partial \theta_{ij}^\alpha(\xi_i,t)}{\partial t}
\end{align*}
by setting $t=0$, we obtain (for the detail, see \cite{Kod05} p.211-213)
\begin{align}\label{}
\dot{\theta}_{ij}(0) -\dot{\theta}_{ik}(0) + \dot{\theta}_{jk} (0) = \left[ \theta_{ij}(0), \theta_{jk}(0) \right]
\end{align}

As in subsection \ref{on26}, we obtain
{\small{\begin{align}\label{ut13}
\delta \left( \left\{ \dot{\theta}_{jk}(0), - \frac{\partial}{\partial t_j}\left( \frac{1}{\Psi_{jk}}\frac{\partial \Psi_{jk}}{\partial t_k} \right)|_{t=0}   \right\}\right)=\left\{ \left[ \left( \theta_{ij}(0), - \frac{1}{\Psi_{0ij}}\frac{\partial \Psi_{ij}}{\partial t}|_{t=0} \right), \left( \theta_{jk}(0), -\frac{1}{\Psi_{0jk}}\frac{\partial \Psi_{jk}}{\partial t}|_{t=0} \right)    \right]\right\}
\end{align}}}
\begin{align}\label{ut14}
&-\Psi_{0jk}\frac{\partial^2 w_k}{\partial t^2}|_{t=0} + \frac{\partial^2 w_j}{\partial t^2}|_{t=0} + \mathcal{L}_{\dot{\theta}_{jk}(0)}\left(w_{0j} \right)- \frac{\partial}{\partial t_j}\left( \frac{1}{\Psi_{jk}}\frac{\partial \Psi_{jk}}{\partial t} \right)|_{t=0} \cdot w_{0j} \\
&= \mathcal{L}_{\theta_{jk}(0)}\left(- \frac{\partial w_j}{\partial t_j}|_{t=0}- \Psi_{0jk}\frac{\partial w_k}{\partial t}|_{t=0} \right) - \frac{1}{\Psi_{0jk}}\frac{\partial \Psi_{jk}}{\partial t}|_{t=0}\left( -\frac{\partial w_j}{\partial t_j}|_{t=0} - \Psi_{0jk}\frac{\partial w_k}{\partial t_k}|_{t=0}\right) \notag
\end{align}
\begin{align}\label{ut15}
\left(-\frac{\partial^2 w_j}{\partial t^2}|_{t=0}\right) \wedge dw_{0j} + w_{0j} \wedge d\left( -\frac{\partial^2 w_j}{\partial t^2}|_{t=0} \right) =  2\left(-\frac{\partial w_j}{\partial t}|_{t=0}\right) \wedge d\left(-\frac{\partial w_j}{\partial t}|_{t=0}\right)
\end{align}
On the other hand, by taking the derivative of $(\ref{ut16})$ with respect to $t$ again $\left(v_j= v_j^1 dt \right)$, we have
\begin{align}\label{ut18}
\sum_{\alpha, \beta=1}^n \frac{\partial w_j^\alpha}{\partial \xi_j^\beta}\frac{\partial f_{jk}^\beta}{\partial t} \frac{\partial f_{jk}^\alpha}{\partial t} +  \sum_{\alpha=1}^n  \frac{\partial w_j^\alpha}{\partial t}\frac{\partial f_{jk}^\alpha}{\partial t} +\sum_{\alpha=1}^n  w_j^\alpha \frac{\partial^2 f_{jk}^\alpha}{\partial t^2 } +  \sum_{\alpha=1}^n \frac{\partial v_j^1}{\partial \xi_j^\alpha}\frac{\partial f_{jk}^\alpha}{\partial t} +   \frac{\partial v_j^1}{\partial t_j } = \frac{\partial \Psi_{jk}}{\partial t} v_k^1 + \Psi_{jk}\frac{\partial v_k^1}{\partial t}
\end{align}
By noting $\dot{\theta}_{jk}^\gamma=\frac{\partial^2 f_{jk}^\gamma}{\partial t^2}- \sum_{\delta=1}^n \frac{\partial f_{jk}^{\delta}}{\partial t} \frac{\partial}{\partial \xi_j^\delta}\left(  \frac{\partial f_{jk}^\gamma}{\partial t} \right)$, $(\ref{ut18})$ is equivalent to 
{\small{\begin{align}\label{ut17}
\sum_{\alpha,\beta=1}^n \frac{\partial w_j^{\alpha}}{\partial \xi_j^\beta}\frac{\partial f_{jk}^\beta}{\partial t}\frac{\partial f_{jk}^\alpha}{\partial t}+\sum_{\alpha=1}^n \frac{\partial w_j^{\alpha}}{\partial t_j}\frac{\partial f_{jk}^\alpha}{\partial t} +\sum_{\gamma=1}^n w_j^{\alpha} \left( \dot{\theta}_{jk}^\alpha + \sum_{\delta=1}^n \frac{\partial f_{jk}^\delta}{\partial t} \frac{\partial}{\partial \xi_j^\delta} \left(  \frac{\partial f_{jk}^\alpha}{\partial t} \right)   \right) +\sum_{\alpha=1}^n \frac{\partial v_j^1}{\partial \xi_j^\alpha}\frac{\partial f_{jk}^\alpha}{\partial t} +\frac{\partial v_j^1}{\partial t_j}=\frac{\partial \Psi_{jk}}{\partial t}v_k^1 +\Psi_{jk} \frac{\partial v_k^1}{\partial t}
\end{align}}}
We note that
{\small{\begin{align}\label{ut19}
i_{\theta_{jk}}\mathcal{L}_{\theta_{jk}}\left(w_j \right) & = i_{\sum_{\sigma=1}^n \frac{\partial f_{jk}^\sigma}{\partial t} \frac{\partial}{\partial \xi_j^\sigma}} \mathcal{L}_{\sum_{\delta=1}^n \frac{\partial f_{jk}^\delta}{\partial t} \frac{\partial }{\partial \xi_j^\delta}} \left(\sum_{\gamma=1}^n w_j^{ \gamma} d\xi_j^\gamma \right)=  i_{\sum_{\sigma=1}^n \frac{\partial f_{jk}^\sigma}{\partial t} \frac{\partial}{\partial \xi_j^\sigma}} \left(  \sum_{\delta , \gamma=1}^n \frac{\partial f_{jk}^\delta}{\partial t}     \frac{\partial w_j^{ \gamma}}{\partial \xi_j^\delta} d\xi_j^\gamma + w_j^{ \gamma}\frac{\partial}{\partial \xi_j^\delta}\left( \frac{\partial f_{jk}^\gamma}{\partial t}\right) d\xi_j^\delta   \right)\\
&= \sum_{\delta, \gamma=1}^n \frac{\partial f_{jk}^\gamma}{\partial t} \frac{\partial f_{jk}^\delta}{\partial t} \frac{\partial w_j^{ \gamma}}{\partial \xi_j^\delta} + \sum_{\gamma, \delta=1}^n w_j^{ \gamma} \frac{\partial f_{jk}^\delta}{\partial t}\frac{\partial}{\partial \xi_j^\delta}\left(   \frac{\partial f_{jk}^\gamma}{\partial t}  \right)     \notag
\end{align}}}
The $(\ref{ut19})$ implies that $(\ref{ut17})$ is equivalent to
\begin{align}\label{ut20}
i_{\theta_{jk}}\mathcal{L}_{\theta_{jk}}(w_j) + i_{\theta_{jk}}\left(\frac{\partial w_j}{\partial t_j} \right) + i_{\dot{\theta}_{jk}}\left( w_j \right) + \mathcal{L}_{\theta_{jk}}\left(i_{\frac{\partial}{\partial t}}\left( v_j \right) \right) + \frac{\partial}{\partial t}\left( i_{\frac{\partial}{\partial t} } \left( v_j \right) \right)=  \frac{\partial \Psi_{jk}}{\partial t}i_{\frac{\partial}{\partial t}}\left( v_k\right) +\Psi_{jk} \frac{\partial}{\partial t}  \left( i_{\frac{\partial}{\partial t}}\left( v_k \right) \right)
\end{align}

\begin{lemma}
We have the equality
{\small{\begin{align}\label{ut22}
&- \Psi_{jk}\frac{\partial}{\partial t}\left(i_{\frac{\partial}{\partial t}} (v_k)    \right) + \frac{\partial}{\partial t}\left( i_{\frac{\partial}{\partial t}}(v_j)\right) + i_{\dot{\theta}_{jk}}(w_j)\\
&= \frac{1}{2}\left( - \frac{1}{\Psi_{jk}}\frac{\partial \Psi_{jk}}{\partial t}\right)\left( -i_{\frac{\partial}{\partial t}}(v_j) - \Psi_{jk} i_{\frac{\partial}{\partial t}}(v_k) \right) + \frac{1}{2}\mathcal{L}_{\theta_{jk}}\left( - i_{\frac{\partial}{\partial t}}(v_j) - \Psi_{jk} i_{\frac{\partial}{\partial t}}(v_k)  \right) + \frac{1}{2} i_{\theta_{jk}}\left( - \frac{\partial w_j}{\partial t} - \Psi_{jk} \frac{\partial w_k}{\partial t} \right) \notag
\end{align}}}
\end{lemma}

\begin{proof}
We note that from $(\ref{ut20})$, we have
{\small{\begin{align*}
&\frac{1}{2}\left( - \frac{1}{\Psi_{jk}}\frac{\partial \Psi_{jk}}{\partial t}\right)\left( -i_{\frac{\partial}{\partial t}}(v_j) - \Psi_{jk} i_{\frac{\partial}{\partial t}}(v_k) \right) + \frac{1}{2}\mathcal{L}_{\theta_{jk}}\left( - i_{\frac{\partial}{\partial t}}(v_j) - \Psi_{jk} i_{\frac{\partial}{\partial t}}(v_k)  \right) + \frac{1}{2} i_{\theta_{jk}}\left( - \frac{\partial w_j}{\partial t} - \Psi_{jk} \frac{\partial w_k}{\partial t} \right)\\
&= \frac{1}{2}\left( - \frac{1}{\Psi_{jk}}\frac{\partial \Psi_{jk}}{\partial t} \right)\left( -2 i_{\frac{\partial}{\partial t}}(v_j) - i_{\theta_{jk}}\left(w_j \right)\right) + \frac{1}{2}\mathcal{L}_{\theta_{jk}}\left(- 2i_{\frac{\partial}{\partial t}}(v_j) - i_{\theta_{jk}}(w_j) \right) + \frac{1}{2} i_{\theta_{jk}} \left(-2 \frac{\partial w_j}{\partial t} - \mathcal{L}_{\theta_{jk}}(w_j) + \frac{1}{\Psi_{jk}}\frac{\partial \Psi_{jk}}{\partial t} w_j \right)\\
&= \frac{1}{\Psi_{jk}}\frac{\partial \Psi_{jk}}{\partial t}\left( i_{\frac{\partial}{\partial t}}(v_j) + i_{\theta_{jk}}(w_j)\right) - \mathcal{L}_{\theta_{jk}}\left( i_{\frac{\partial}{\partial t}}(v_j) \right) -i_{\theta_{jk}}\mathcal{L}_{\theta_{jk}}(w_j)- i_{\theta_{jk}}\left(\frac{\partial w_j}{\partial t} \right)\\
&=- \Psi_{jk}\frac{\partial}{\partial t}\left(i_{\frac{\partial}{\partial t}} (v_k)    \right) + \frac{\partial}{\partial t}\left( i_{\frac{\partial}{\partial t}}(v_j)\right) + i_{\dot{\theta}_{jk}}(w_j)
\end{align*}}}
\end{proof}
By setting $t=0$ on $(\ref{ut22})$, we obtain
{\small{\begin{align}\label{ut23}
&- \Psi_{0jk}\frac{\partial}{\partial t}\left(i_{\frac{\partial}{\partial t}} (v_k)    \right)|_{t=0} + \frac{\partial}{\partial t}\left( i_{\frac{\partial}{\partial t}}(v_j)\right)|_{t=0} + i_{\dot{\theta}_{jk}(0)}(w_{0j})\\
&= \frac{1}{2}\left( - \frac{1}{\Psi_{0jk}}\frac{\partial \Psi_{jk}}{\partial t}|_{t=0}\right)\left( -i_{\frac{\partial}{\partial t}}(v_j)|_{t=0} - \Psi_{0jk} i_{\frac{\partial}{\partial t}}(v_k) |_{t=0} \right) \notag \\
& + \frac{1}{2}\mathcal{L}_{\theta_{jk}(0)}\left( - i_{\frac{\partial}{\partial t}}(v_j)|_{t=0} - \Psi_{0jk} i_{\frac{\partial}{\partial t}}(v_k)|_{t=0}  \right) + \frac{1}{2} i_{\theta_{jk}(0)}\left( - \frac{\partial w_j}{\partial t} |_{t=0} - \Psi_{0jk} \frac{\partial w_k}{\partial t}|_{t=0} \right) \notag
\end{align}}}

By taking the derivative of $(\ref{ut26})$ with respect to $t$, we have
\begin{align*}
\frac{\partial }{\partial t}\left( i_{\frac{\partial}{\partial t}}(v_j) \right) dw_j + i_{\frac{\partial}{\partial t}}(v_j) d\left( \frac{\partial w_j}{\partial t} \right) + \frac{\partial w_j}{\partial t} \wedge d\left( i_{\frac{\partial}{\partial t}}(v_j)\right) + w_j \wedge \left(  \frac{\partial}{\partial t}\left( i_{\frac{\partial}{\partial t}}(v_j) \right) \right) - \frac{\partial w_j}{\partial t}\wedge \frac{\partial w_j}{\partial t} - w_j\wedge \frac{\partial^2 w_j}{\partial t^2}=0
\end{align*}

Then by setting $t=0$, we have 
{\small{\begin{align}\label{ut24}
&-\frac{\partial }{\partial t}\left( i_{\frac{\partial}{\partial t}}(v_j) \right)|_{t=0} dw_{0j}  + w_{0j} \wedge \left(  -\frac{\partial}{\partial t}\left( i_{\frac{\partial}{\partial t}}(v_j) \right)|_{t=0} \right)  - w_{0j}\wedge \left( -\frac{\partial^2 w_j}{\partial t^2}|_{t=0} \right) \\
&= - i_{\frac{\partial}{\partial t}}(v_j)|_{t=0} d\left(- \frac{\partial w_j}{\partial t}|_{t=0} \right) + \left( -\frac{\partial w_j}{\partial t}|_{t=0} \right)  \wedge d\left( - i_{\frac{\partial}{\partial t}}(v_j)|_{t=0}\right) \notag
\end{align}}}

For any $1$-cocycle $\left( \left\{a_j\right\}, \left\{u_j\right\}, \left\{ \left( b_{jk}, c_{jk} \right) \right\} \right)\in C^0\left( \mathcal{U}_0, \Omega_M^1\otimes \mathcal{L}_0 \right) \bigoplus C^0\left(\mathcal{U}_0, \mathcal{L}_0 \right)\bigoplus C^1\left( \mathcal{U}_0, \mathcal{E}_{\mathcal{L}_0} \right)$ in the \v Cech resolution $(\ref{ut27})$ of $\mathcal{E}_{\mathcal{L}_0}^\bullet$, we can define a $2$-cocycle by (we omit the proof, see Part IV-1)
\begin{align}\label{ut28}
&\left( \left\{ 2 a_j \wedge d a_j \right\} , \left\{ u_j da_j+ a_j \wedge du_j\right\},\left\{ \mathcal{L}_{b_{jk}}\left( a_j + \Psi_{0jk} a_k  \right) + c_{jk}\left(a_j+ \Psi_{0jk} a_k \right)\right\}, \right. \\
&\left. \left\{ \frac{1}{2}\mathcal{L}_{b_{jk}}\left(u_j + \Psi_{0jk}u_k \right) + \frac{1}{2}i_{b_{jk}}\left(a_j + \Psi_{0jk} a_k \right) + \frac{1}{2} c_{jk}\left( u_j + \Psi_{0jk} u_k \right)\right\}  ,\left\{\left[ \left(b_{ij}, c_{ij}\right), \left( b_{jk}, c_{jk}\right) \right] \right\} \right) \notag \\
&\in C^0\left(\mathcal{U}_0, \Omega_M^1\otimes \tilde{\mathcal{S}}^2 \right) \bigoplus C^0\left( \mathcal{U}_0, \mathcal{L}_0 \right) \bigoplus C^1\left( \mathcal{U}_0, \Omega_M^1\otimes \mathcal{L}_0 \right) \bigoplus C^1\left( \mathcal{U}_0, \mathcal{L}_0 \right) \bigoplus C^2\left( \mathcal{U}_0, \mathcal{E}_{\mathcal{L}_0} \right) \notag
\end{align}

Then $(\ref{ut13}),(\ref{ut14}),(\ref{ut15}), (\ref{ut23})$ and $(\ref{ut24})$ imply that
\begin{theorem} \label{ut25}
Suppose given a compact complex manifold $M$ with a codimension $1$-foliation defined by $w_0\in H^0\left(M, \Omega_M^1 \otimes \mathcal{L}_0 \right)$ for a line bundle $\mathcal{L}_0$, and $\left(\left\{a_j \right\}, \left\{ u_j \right\}, \left\{\left( b_{jk}, c_{jk} \right) \right\} \right)\in \mathbb{H}^1\left(M, \tilde{\mathcal{E}}_{\mathcal{L}_0}^\bullet \right)$. In order that there may exist an unfolded analytic family $\left( \mathcal{M}, \omega , B, \pi\right)$ such that $\pi^{-1}(0)=\left( M, w_0 \right)$ and $\left(\frac{\partial \left(\mathcal{M}, \omega \right)}{\partial t}\right)|_{t=0}=\left(\left\{a_j \right\}, \left\{u_j \right\},  \left\{\left( b_{jk}, c_{jk} \right) \right\} \right)$, it is necessary  that $(\ref{ut28})$ defines the $0$ class in $\mathbb{H}^2\left( M, \tilde{\mathcal{E}}_{\mathcal{L}_0}^\bullet \right)$.

\end{theorem}

Keeping the notations above, we have the following lemma.
\begin{lemma}
A compact foliated complex manifold $(M, w_0)$ of codimension $1$ with $w_0 \in H^0\left( M, \Omega_M^1\otimes\mathcal{L}_0 \right)$  is obstructed in unfoldings if there exists an element $\left( \{(r_j, v_j)\} , \{(b_{jk}, c_{jk}\} \right)\in \mathbb{H}^1\left(M, \tilde{\mathcal{E}}_{\mathcal{L}_0}^\bullet \right)$ $($In particular, $\{(b_{jk}, c_{jk})\}$ defines an element in $\textnormal{ker}\left( H^1\left( M, \mathcal{E}_{\mathcal{L}_0} \right) \xrightarrow{\tilde{G}_0} H^1\left( M, \Omega_M^1\otimes \mathcal{L}_0 \right) \bigoplus H^1\left(M, \mathcal{L}_0 \right) \right) )$, and an element $\{(a_j, u_j)\}\in \textnormal{ker}\left(H^0\left(M, \Omega_M^1\otimes \mathcal{L}_0 \right)\bigoplus H^0\left( M, \mathcal{L}_0 \right) \xrightarrow{G_1} H^0\left(M, \tilde{\mathcal{S}}^2\otimes \mathcal{L}_0 \right)\bigoplus H^0\left( M,  \tilde{\mathcal{S}}^1 \otimes \mathcal{L}_0 \right) \right)$ such that one of the followings holds
\begin{enumerate}
\item under the map
\begin{align}
H^0\left( M, \Omega_M^1 \otimes \mathcal{L}_0 \right) \times H^0&\left( M, \mathcal{L}_0  \right) \times H^1\left( M, \mathcal{E}_{\mathcal{L}_0} \right) \to H^1\left(M, \Omega_M^1\otimes \mathcal{L}_0 \right)\times H^1 \left( M, \mathcal{L}_0 \right)\label{um1}\\
\left(\{a_j\}, \{u_j\}, \{(b_{jk}, c_{jk})\} \right) &\to \left( \left\{\mathcal{L}_{b_{jk}}(a_j)+ c_{jk}a_j \right\}, \left\{ \frac{1}{2}i_{b_{jk}}(a_j)+ \frac{1}{2}[b_{jk}, u_j]+ \frac{1}{2}c_{jk}u_j\right\}  \right) \notag
\end{align}
$\left( \left\{\mathcal{L}_{b_{jk}}(a_j)+ c_{jk}a_j \right\}, \left\{ \frac{1}{2}i_{b_{jk}}(a_j)+ \frac{1}{2}[b_{jk}, u_j]+ \frac{1}{2}c_{jk}u_j\right\}  \right) \in H^1\left( M, \Omega_M^1\otimes \mathcal{L}_0\right)\bigoplus H^1\left( M, \mathcal{L}_0 \right) $ is not in the image of $H^1\left( M, \mathcal{E}_{\mathcal{L}_0} \right)\xrightarrow{\tilde{G}_0} H^1\left( M, \Omega_M^1\otimes \mathcal{L}_0 \right)\bigoplus H^1\left( M, \mathcal{L}_0 \right)$.
\item under the map
\begin{align}
H^0\left(M, \Omega_M^1\otimes \mathcal{L}_0 \right) \times H^1\left(M, \mathcal{E}_{\mathcal{L}_0} \right)&\to H^1\left(M, \Omega_M^1\otimes \mathcal{L}_0 \right) \\
\left(\{a_j\}, \{(b_{jk}, c_{jk}\} \right)& \to \left( \left\{\mathcal{L}_{b_{jk}}(a_j) + c_{jk}a_j \right\} \right)\notag
\end{align}
$\left( \left\{\mathcal{L}_{b_{jk}}(a_j) + c_{jk}a_j \right\} \right)$ is not in the image of $H^1\left( M, \mathcal{E}_{\mathcal{L}_0} \right) \xrightarrow{\tilde{E}_0} H^1\left( M, \Omega_M^1 \otimes \mathcal{L}_0 \right)$.
\item under the map
\begin{align*}
H^0\left( M, \Omega_M^1 \otimes \mathcal{L}_0 \right) \times H^0&\left( M, \mathcal{L}_0  \right) \times H^1\left( M, \mathcal{E}_{\mathcal{L}_0} \right) \to H^1 \left( M, \mathcal{L}_0 \right)\\
\left(\{a_j\}, \{u_j\}, \{(b_{jk}, c_{jk})\} \right) &\to \left( \left\{ \frac{1}{2}i_{b_{jk}}(a_j)+ \frac{1}{2}[b_{jk}, u_j]+ \frac{1}{2}c_{jk}u_j\right\}  \right) \notag
\end{align*}
$\left( \left\{ \frac{1}{2}i_{b_{jk}}(a_j)+ \frac{1}{2}[b_{jk}, u_j]+ \frac{1}{2}c_{jk}u_j\right\}  \right)$ is not in the image of $H^1\left( M, \mathcal{E}_{\mathcal{L}_0} \right)\to H^1\left( M, \Theta_M \right)\xrightarrow{i_{-}(w_0)} H^1\left( M, \mathcal{L}_0 \right)$.
\item  $H^1\left( M, \mathcal{E}_{\mathcal{L}_0} \right)=0$, and  $\left( \left\{ 2a_j \wedge da_j\right\}, \left\{ u_j da_j+ a_j\wedge du_j\right\} \right)\in H^0\left(M, \tilde{\mathcal{S}}^2 \otimes \mathcal{L}_0  \right) \bigoplus H^0\left( M,  \tilde{\mathcal{S}}^1 \otimes\mathcal{L}_0 \right) $ is not in the image of $H^0\left(M, \Omega_M^1 \otimes \mathcal{L}_0 \right) \bigoplus H^0\left( M, \mathcal{L}_0 \right) \xrightarrow{G_1} H^0\left( M, \tilde{\mathcal{S}}^2 \otimes \mathcal{L}_0 \right) \bigoplus H^0\left(\tilde{\mathcal{S}}^1\otimes \mathcal{L}_0 \right)$.
\item $H^1\left( M, \mathcal{E}_{\mathcal{L}_0} \right) =  0 $, and  $\left\{ 2a_j \wedge da_j\right\}$ is not in the image of $H^0\left( M, \Omega_M^1\otimes \mathcal{L}_0 \right)\xrightarrow{E_1} H^0 \left( M, \tilde{\mathcal{S}}^2 \otimes \mathcal{L}_0 \right)$.
\item $H^1\left( M, \mathcal{E}_{\mathcal{L}_0} \right)=0$ and $\left\{ u_jda_j+ a_j\wedge du_j \right\}\in H^0\left( M, \tilde{\mathcal{S}}^1\otimes \mathcal{L}_0 \right)$ is not in the image of the map $H^0\left( M, \Omega_M^1\otimes \mathcal{L}_0 \right) \bigoplus H^0\left(M, \mathcal{L}_0 \right) \to H^0\left( M, \tilde{\mathcal{S}}^1 \otimes \mathcal{L}_0 \right)$ from the component of $\tilde{G}_1$.
\item $\left\{\left[(b_{ij}, c_{ij}), (b_{jk}, c_{jk}) \right]\right\}=\left\{ [b_{jk}, c_{jk}], [b_{ij}, c_{jk} - [b_{jk}, \log \Psi_{0ij}] ]- [c_{jk}, b_{ij}]  \right\} $ is not $0$ class in $H^2\left( M, \mathcal{E}_{\mathcal{L}_0} \right)$.
\item $\left\{ [b_{ij}, b_{jk}]\right\}$ is not $0$ class in $H^2\left( M, \Theta_M \right)$.
\end{enumerate}
\end{lemma}

\begin{proof}
We show that the map $(\ref{um1})$ is well-defined. In fact, we have shown in Lemma \ref{on23}, $\left\{ \mathcal{L}_{b_{jk}}(a_j)+ c_{jk}a_j \right\} $ defines an element in $H^1\left(M, \Omega_M^1\otimes \mathcal{L}_0 \right)$ and we show that $\left\{ \frac{1}{2}i_{b_{jk}}(a_j)+ \frac{1}{2}[b_{jk}, u_j]+ \frac{1}{2}c_{jk}u_j\right\} $ defines an element in $H^1\left(M, \mathcal{L}_0 \right)$. In fact,  since $\Psi_{0ij}a_j= a_i,\Psi_{0ij}u_j=u_i$ and $b_{ij}-b_{ik}+ b_{jk}=0$ and $c_{ij}- c_{ik}+ c_{jk}- \left[ b_{jk}, \log \Psi_{0ij}\right]=0$, we have
{\small{\begin{align*}
&\frac{1}{2}i_{b_{ij}}(a_i)+ \frac{1}{2}\left[b_{ij}, u_i\right]+ \frac{1}{2}c_{ij}u_i -\frac{1}{2}i_{b_{ik}}(a_i)- \frac{1}{2}\left[b_{ik}, u_i \right]- \frac{1}{2}c_{ik}u_i + \Psi_{0ij}\left(\frac{1}{2}i_{b_{jk}}(a_j)+ \frac{1}{2}\left[b_{jk}, u_j\right]+ \frac{1}{2}c_{jk}u_j \right)\\
&= \frac{1}{2}i_{b_{ij}-b_{ik} +b_{jk}}\left(a_i \right) + \frac{1}{2}\left(c_{ij}- c_{ik} +c_{jk} \right) u_i + \frac{1}{2}\left[b_{ij}- b_{ik}+ b_{jk}, u_i \right] - \frac{1}{2}\left[ b_{jk}, \Psi_{0ij}\right] u_j=0
\end{align*}}}
Hence the map $(\ref{um1})$ is well-defined. Since $\left( \{r_j\}, \{v_j\}, \left\{ (b_{jk}, c_{jk})\right\}\right)$ defines a $1$-cocycle, we have
\begin{align*}
w_{0j}\wedge dr_j+ r_j\wedge dw_{0j}=0&,\,\,\,\,\, v_j dw_{0j}+ w_{0j}\wedge dv_j - r_j =0\\
\Psi_{0jk}r_k-r_j+\mathcal{L}_{b_{jk}}(w_{0j})+ c_{jk}w_{0j}=0&,\,\,\,\,\,\Psi_{0jk} v_k- v_j + i_{b_{jk}}(w_{0j})=0\\
b_{ij}- b_{ik} + b_{jk}=0&,\,\,\,\,\, c_{ij}-c_{ik} + c_{jk}- \left[b_{jk}, \log \Psi_{0ij} \right]=0
\end{align*}
and $\left( \left\{ r_j + a_j \right\}, \left\{v_j + u_j \right\}, \left\{ \left( b_{jk}, c_{jk} \right) \right\} \right)$ also defines a $1$-cocycle. Assume that $\left(M, w_0 \right)$ is unobstructed in unfoldgins. Then by Theorem \ref{ut25} there exists $\left(\{g_i\}, \{l_i\}, \{(e_{jk}, h_{jk})\} \right)$ and $(\{g_i'\},\{l_i'\},  \{(e_{jk}', h_{jk}')\}\in C^0\left(\mathcal{U}, \Omega_M^1\otimes \mathcal{L}_0 \right)$ $\bigoplus C^0\left(\mathcal{U}, \mathcal{L}_0 \right) \bigoplus C^1\left(\mathcal{U}, \mathcal{E}_{\mathcal{L}_0} \right)$ such that
{\small{\begin{align}\label{um2}
w_{0i}\wedge dg_i+ g_i\wedge dw_{0i}&=2 r_j\wedge dr_j \\
w_{0i}\wedge dl_i + l_i dw_{0i}- g_i & = v_j dr_j+ r_j \wedge dv_j \notag \\
\delta\left(\{g_i\} \right) + \left\{ \mathcal{L}_{e_{jk}}(w_{0j} ) + h_{jk}w_{0j}\right\}&= \mathcal{L}_{b_{jk}}\left(r_j+ \Psi_{0jk}r_k \right)+ c_{jk}\left( r_j+ \Psi_{0jk}r_k \right) \notag\\
\delta\left(\left\{ l_i\right\} \right)+ \left\{ i_{e_{jk}}(w_{0j}) \right\} &= \frac{1}{2}\mathcal{L}_{b_{jk}}\left( v_j + \Psi_{0jk} v_k \right) + \frac{1}{2}i_{b_{jk}}\left(r_j + \Psi_{0jk} r_k \right) + \frac{1}{2}c_{jk}\left(v_j+ \Psi_{0jk} v_k \right) \notag \\
\delta\left(\left\{ (e_{jk}, h_{jk})\right\} \right)&=\left\{ \left[ (b_{ij}, c_{ij}) , (b_{jk}, c_{jk}) \right] \right\} \notag
\end{align}}}
and
{\small{\begin{align}\label{um3}
w_{0i}\wedge dg_i'+ g_i'\wedge dw_{0i}&=2 (r_j+ a_j)\wedge d(r_j+a_j) \\
w_{0i}\wedge dl_i' + l_i' dw_{0i}- g_i' & = (v_j+u_j) d(r_j+a_j)+ (r_j+a_j) \wedge d(v_j+ u_j) \notag \\
\delta\left(\{g_i'\} \right) + \left\{ \mathcal{L}_{e_{jk}'}(w_{0j} ) + h_{jk}'w_{0j}\right\}&= \mathcal{L}_{b_{jk}}\left(r_j+ a_j+ \Psi_{0jk}(r_k+a_k) \right)+ c_{jk}\left( r_j+a_j+ \Psi_{0jk}(r_k+a_k) \right) \notag\\
\delta\left(\left\{ l_i\right\} \right)+ \left\{ i_{e_{jk}'}(w_{0j}) \right\} = \frac{1}{2}\mathcal{L}_{b_{jk}}&\left( v_j+ u_j + \Psi_{0jk} (v_k+ u_k) \right) + \frac{1}{2}i_{b_{jk}}\left(r_j+a_j + \Psi_{0jk} (r_k+ a_k) \right) + \frac{1}{2}c_{jk}\left(v_j+u_j + \Psi_{0jk} (v_k+ u_k) \right)\notag\\
\delta\left(\left\{ (e_{jk}', h_{jk}')\right\} \right)&=\left\{ \left[ (b_{ij}, c_{ij}) , (b_{jk}, c_{jk}) \right] \right\} \notag
\end{align}}}
Then by subtracting $(\ref{um2})$ from $(\ref{um3})$, we have
\begin{align}
&w_{0i}\wedge d(g_i'-g_i)+ (g_i'-g_i)\wedge dw_{0i}=2 a_j\wedge dr_j+ 2r_j \wedge da_j + 2 a_j \wedge da_j \notag \\
&w_{0i}\wedge d(l_i'-l_i) + (l_i'-l_i) dw_{0i}- (g_i'-g_i)  = v_j da_j + u_j dr_j + u_j da_j  + a_j\wedge d v_j + r_j \wedge du_j + a_j \wedge dv_j + a_j \wedge du_j \notag \\
&\delta\left(\{g_i'-g_i\} \right) + \left\{ \mathcal{L}_{e_{jk}'-e_{jk}}(w_{0j} ) +( h_{jk}'-h_{jk})w_{0j}\right\}= \mathcal{L}_{b_{jk}}\left(a_j+ \Psi_{0jk}a_k \right)+ c_{jk}\left( a_j+ \Psi_{0jk}a_k \right) = 2\mathcal{L}_{b_{jk}}(a_j)+ 2c_{jk} a_j \label{um4}\\
&\delta\left(\left\{ l_i'-l_i\right\} \right)+ \left\{ i_{e_{jk}'-e_{jk}}(w_{0j}) \right\} = \mathcal{L}_{b_{jk}}\left(u_j \right) + i_{b_{jk}}\left(a_j \right) + c_{jk}u_j \label{um5} \\
&\delta\left(\left\{ (e_{jk}'-e_{jk}, h_{jk}'-h_{jk})\right\} \right)=0 \label{um6}
\end{align}
We note that $(\ref{um6})$ implies that $\left\{ \left( e_{jk}'- e_{jk}  , h_{jk}'- h_{jk} \right) \right\} \in H^1\left(M, \mathcal{E}_{\mathcal{L}_0} \right)$. Then $(\ref{um4})$ implies that $\left\{\mathcal{L}_{b_{jk}}(a_j)+ c_{jk}a_j  \right\}$ is cohomogloous to the image of $H^1\left(M, \mathcal{E}_{\mathcal{L}_0} \right)\xrightarrow{\tilde{E}_0} H^1\left(M, \Omega_M^1 \otimes \mathcal{L}_0 \right)$. This proves $(2)$, and moreover $(\ref{um5})$ implies that $\left\{ \frac{1}{2}\left[ b_{jk}, u_j \right]+ \frac{1}{2}i_{b_{jk}}(a_j)+ \frac{1}{2}c_{jk}u_j\right\}$ is cohomologous to the image of $H^1\left(M, \mathcal{E}_{L_0} \right)\to H^1\left( M, \Theta_M \right)\xrightarrow{i_{-}(w_0)} H^1\left( M, \mathcal{L}_0 \right)$. This proves $(3)$. Our arguments also shows $(1)$.

We prove $(4),(5),(6)$. Since $( \{a_j\} , \{u_j\}, \{(0,0)\})$ defines an element $\mathbb{H}^1\left( M, \tilde{\mathcal{E}}_{\mathcal{L}_0}^\bullet \right)$. If $\left(M, w_0 \right)$ is unobstructed in unfoldings, by Theorem \ref{ut25}, there exists $(\{g_i\},\{ l_i\}, \{e_{jk}, h_{jk}\}) \in C^0\left( \mathcal{U},  \Omega_M^1\otimes \mathcal{L}_0 \right) \bigoplus C^0\left( \mathcal{U}, \mathcal{L}_0 \right)$ $\bigoplus C^1\left( \mathcal{U} , \mathcal{E}_{\mathcal{L}_0} \right)$ such that
\begin{align*}
\delta(\{e_{jk}, h_{jk}\})=0,\,\,\,&\left\{ \Psi_{0jk}g_k- g_j + \mathcal{L}_{e_{jk}}(w_{0j})+ \Psi_{0jk}w_{0j}  \right\} =0,\,\,\, \left\{ \Psi_{0jk}l_k - l_j + i_{e_{jk}}(w_{0j})   \right\}=0\\
&\tilde{G}_1(\{g_i, l_i\})= \left\{ 2a_j \wedge da_j, u_j da_j+ a_j \wedge du_j \right\}
\end{align*}
Now assume that $H^1\left( M, \mathcal{E}_{\mathcal{L}_0}  \right)=0$. Then there exists $(e_i, h_i)\in C^0\left( M, \mathcal{E}_{\mathcal{L}_0} \right)$ such that $\delta(\left\{ e_i, h_i\right\})=\left\{(e_{jk}, h_{jk}) \right\}$. In other words,
\begin{align*}
e_{jk}= e_k- e_j,\,\,\,\,\,\,\,\, h_{jk}= h_k -\left[ e_k , \log \Psi_{0jk}\right] - h_j
\end{align*}
 Then $\left( \{ g_i +\mathcal{L}_{e_i}(w_{0i})+ h_iw_{0i}, l_i+ i_{e_i}(w_{0i})\}\right)$ defines a global section $H^0\left( M, \Omega_M^1\otimes \mathcal{L}_0 \right) \bigoplus H^0\left( M, \mathcal{L}_0 \right)$. In fact, $\left\{ g_i+ \mathcal{L}_{e_i}(w_{0i})+ h_iw_{0i}\right\}\in H^0\left( M , \Omega_M^1 \otimes \mathcal{L}_0 \right)$ as in the proof of Lemma \ref{on23}, and 
\begin{align*}
\Psi_{0jk}\left(l_k+ i_{e_k}(w_{0k})\right)- \left(l_j+ i_{e_j}(w_{0j})\right)= - i_{e_{jk}}(w_{0j})+ i_{e_k}(w_{0j}) - i_{e_j}(w_{0j}) = 0
\end{align*}
Then we note that
\begin{align*}
\tilde{G}_1\left( \left\{g_i + \mathcal{L}_{e_i}(w_{0i}) + h_i w_{0i} \right\}, \left\{ l_i+ i_{e_i}(w_{0i}) \right\} \right) = \tilde{G}_1\left( \{g_i\}, \{l_i\} \right)= \left\{ 2a_j\wedge da_j, u_jda_j+ a_j \wedge du_j \right\}
\end{align*}
This implies $(4),(5)$, and $(6)$.

\end{proof}

\subsection{Integrability conditions} \

Let $(\mathcal{M}, \omega,  B, \pi)$ be an unfolded complex analytic family with $\omega\in H^0\left(\mathcal{M}, \Omega_{\mathcal{M}}^1 \otimes \mathcal{L} \right)$, where $\mathcal{L}= \left(\tilde{\mathcal{N}}_\mathcal{F}^* \right)^* $ is a line bundle on $\mathcal{M}$ of deformations of $\pi^{-1}(0)=\left(M, w_0 \right)$. We keep the notations in subsection \ref{uti1}. On the other hand, as in the subsection \ref{oni14} (we keep the notations in subsection \ref{oni14}), we have a diffeomorphism $\Phi$ of $M\times \Delta$ onto $\mathcal{M}_\Delta$ with $\pi \circ \Phi$ is the projection $M\times \Delta\to \Delta$.  Let $z=\left(z^1,..., z^n \right)$ be local coordinates of $M=M_0$, and for $\Phi(z,t)\in U_j\times \Delta$, put $\Phi\left( z, t\right) = \left(\xi_j^1(z,t), ..., \xi_j^n(z,t), t_1,..., t_m \right)$ where $\xi_j^\alpha=\xi_j^\alpha(z,t),\alpha=1,...,n$ is a $C^\infty$ function. We note that since $(z^1,..., z^n)$ and $\left(\xi_j^1(z,0), ..., \xi_j^n(z,0)\right)$ are local complex coordinates on $M=M_0$, $\xi_j^\alpha(z,0)$ are holomorphic functions of $z^1,...,z^n,\alpha=1,...,n$, and if we take $\Delta$ sufficiently small, we have
\begin{align}\label{uti2}
\det\left( \frac{\partial \xi_j^\alpha(z,t)}{\partial z^\lambda} \right)_{\alpha, \lambda=1,...,n} \ne 0\,\,\,\,\,\, \, \textnormal{for}\,\,\,t\in \Delta
\end{align}
As in subsection \ref{oni14}, we have non-vanishing $C^\infty$ functions $g_j(z,t):=\left(g_j^{\alpha\beta}(z,t) \right)$ on $\mathscr{U}_j=\Psi^{-1}\left(U_j\times \Delta \right)$ such that
\begin{align}\label{uti3}
g_j(z,t) \cdot   \Psi_{jk}\left(\xi_k(z,0),0 \right) =  \Psi_{jk}\left(\xi_k(z,t), t \right) \cdot g_k(z,t) \iff g_j \cdot \Psi_{0jk} = H_{jk}\cdot g_k
\end{align}
From $(\ref{uti2})$, we can define a $(0,1)$-form $\varphi(t):=\sum_{\lambda,\lambda=1}^n \varphi_{jv}^\lambda(z,t)d\bar{z}^v\frac{\partial}{\partial z^\lambda}\in A^{0,1}\left(M, \Theta_M \right)$ for $t\in \Delta$ such that $\bar{\partial}\xi_j^\alpha - \sum_{\lambda=1}^n \varphi_{jv}^\lambda(z,t)\frac{\partial \xi_j^\alpha}{\partial z_\lambda}d\bar{z}^v$ and
\begin{align}
\varphi(0)=0,\,\,\,\,\,\,\, \bar{\partial} \varphi(t)- \frac{1}{2} \left[ \varphi(t), \varphi(t) \right]=0
\end{align}
(see \cite{Kod05} p.262). On the other hand, as in $(\ref{oi30})$ in subsection \ref{oi31}, we can show
\begin{align}
[\varphi, \log g_j]- \bar{\partial}\log g_j= [\varphi, \log g_k]- \bar{\partial} \log g_k- [ \varphi , \log \Psi_{0jk}]
\end{align}
This implies that
\begin{align*}
\left( (\varphi, Z_j:= \left\{\left[\varphi, \log g_j\right] - \bar{\partial} \log g_j \right\} \right) \in A^{0,1}\left( M, \mathcal{E}_{\mathcal{L}_0} \right)
\end{align*}
and we have
\begin{align}
\bar{\partial} Z_j - \left[\varphi, Z_j\right] = 0 
\end{align}
On the other hand, as in subsection \ref{oni14}, by setting $W_j(z,t)= \sum_{\beta=1}^n W_j^\beta(z,t)dz^\beta:=\sum_{\alpha, \beta=1}^n w_j^\alpha\left( \xi_j(z,t),t)\right) \frac{\partial \xi_j^\alpha}{\partial z^\beta} dz^\beta$, we have
\begin{align*}
\left\{ Y_j(t):= g_j^{-1} W_j\right\}\in A^{0,0}\left( M, \Omega_M^1\otimes \mathcal{L}_0 \right)
\end{align*}
and
\begin{align*}
Y_j\wedge \partial Y_j=0,\,\,\,\,\,\,\,\,\,\, \bar{\partial} Y_j- \mathcal{L}_\varphi \left(Y_j \right) - Z_j Y_j=0
\end{align*}

 We note that
\begin{align*}
&\sum_{\alpha=1}^n w_j^\alpha(\xi_j,t) d\xi_j^\alpha + \sum_{\eta=1}^m v_j^\eta(\xi_j,t) dt_\eta= \Psi_{jk}(\xi_k,t)\left( \sum_{\alpha=1}^n w_k^\alpha(\xi_k,t)d\xi_k^\alpha + \sum_{\eta=1}^m v_k^\eta(\xi_k,t)dt_\eta   \right)\\
&\Longrightarrow \sum_{\alpha=1}^n w_j^\alpha(\xi_j(z,t),t)\frac{\partial \xi_j^\alpha}{\partial t_\eta} + v_j^\eta(\xi_j(z,t),t)= \Psi_{jk}(\xi_k(z,t),t)\left( \sum_{\alpha=1}^n w_k^\alpha(\xi_k(z,t),t)\frac{\partial \xi_k^\alpha}{\partial t_\eta} + v_k^\eta(\xi_k,t)\right)\\
&\Longrightarrow \sum_{\alpha=1}^n w_j^\alpha(\xi_j(z,t),t)\frac{\partial \xi_j^\alpha}{\partial t_\eta} + v_j^\eta(\xi_j(z,t),t)= g_j\cdot \Psi_{0jk} \cdot g_k^{-1} \left( \sum_{\alpha=1}^n w_k^\alpha(\xi_k(z,t),t)\frac{\partial \xi_k^\alpha}{\partial t_\eta} + v_k^\eta(\xi_k,t)\right)
\end{align*}

Then we have
\begin{align*}
g_j^{-1}\left(\sum_{\alpha=1}^n w_j^\alpha(\xi_j(z,t),t)\frac{\partial \xi_j^\alpha}{\partial t_\eta} + v_j^\eta(\xi_j(z,t),t)\right)=   \Psi_{0jk} \cdot g_k^{-1} \left( \sum_{\alpha=1}^n w_k^\alpha(\xi_k(z,t),t)\frac{\partial \xi_k^\alpha}{\partial t_\eta} + v_k^\eta(\xi_k,t)\right)
\end{align*}

This implies that
\begin{align*}
\left\{Q_j^\eta:= g_j^{-1} V_j^\eta  = g_j^{-1}\left(\sum_{\alpha=1}^n w_j^\alpha(\xi_j(z,t),t)\frac{\partial \xi_j^\alpha}{\partial t_\eta} + v_j^\eta(\xi_j(z,t),t)\right)  \right\} \in A^{0,0}\left( M, \mathcal{L} \right)
\end{align*}

Then since $\bar{\partial} V_j^\eta  - \left[ \varphi, V_j^\eta\right]- i_{\frac{\partial \varphi}{\partial t_\eta} } (W_j)=0$, we have
\begin{align*}
&\bar{\partial }\left( Q_j^\eta \right) - \left[ \varphi, Q_j^\eta \right] - i_{\frac{\partial \varphi}{\partial t_\eta }}\left( Y_j \right) = -\frac{1}{g_j^2}\bar{\partial} g_j \cdot  V_j^\eta + \frac{1}{g_j^2}[\varphi, g_j] \cdot V_j^\eta = \left( - \bar{\partial} \log g_j + [\varphi, \log g_j] \right) Q_j^\eta = Z_j Q_j^\eta
\end{align*}
Hence we have
\begin{align*}
\bar{\partial} \left( Q_j^\eta  \right) - \left[ \varphi, Q_j^\eta \right] -  i_{\frac{\partial \varphi}{\partial t_\eta}}\left( Y_j \right) - Z_j Q_j^\eta = 0,\,\,\,\,\,\,\,\,\eta=1,...,m
\end{align*}

We note that
\begin{align*}
&\left( \sum_{\alpha=1}^n  w_j^\alpha(\xi_j,t) d\xi_j^\alpha + \sum_{\eta=1}^m v_j^\eta(\xi_j,t ) dt_\eta  \right) \wedge d \left( \sum_{\alpha=1}^n w_j^\alpha(\xi_j,t) d\xi_j^\alpha + \sum_{\eta=1}^m v_j^\eta(\xi_j,t) dt_\eta \right) = 0\\
\Longrightarrow& g_j^{-1} \left( \sum_{\alpha,\beta=1}^n w_j^\alpha(\xi_j(z,t),t)\frac{\partial \xi_j^\alpha}{\partial z^\beta}dz^\beta + \sum_{\eta=1}^m \left( \sum_{\alpha=1}^n w_j^\alpha(\xi_j(z,t),t)\frac{\partial \xi_j^\alpha}{\partial t_\eta}  + v_j^\eta ( \xi_j(z,t),t) \right) dt_\eta \right)\\
& \wedge \partial \left( g_j^{-1}\left(  \sum_{\alpha,\beta=1}^n w_j^\alpha(\xi_j(z,t),t)\frac{\partial \xi_j^\alpha}{\partial z^\beta}dz^\beta + \sum_{\eta=1}^m \left( \sum_{\alpha=1}^n w_j^\alpha(\xi_j(z,t),t)\frac{\partial \xi_j^\alpha}{\partial t_\eta}  + v_j^\eta ( \xi_j(z,t),t) \right) dt_\eta \right) \right) = 0\\
\Longrightarrow & \left( Y_j + \sum_{\eta=1}^m Q_j^\eta dt_\eta \right) \wedge \partial \left( Y_j + \sum_{\eta=1}^m Q_j^\eta dt_\eta \right) =0
\end{align*}

All in all, we have $\left( \left\{ Y_j\right\}, \left\{ Q_j^\eta \right\}, \{\left(\varphi, Z_j\right)\}     \right) \in A^{0,0}\left( M, \Omega_M^1\otimes \mathcal{L}_0 \right)\bigoplus A^{0,0}\left( M, \mathcal{L}_0 \right)\bigoplus A^{0,1}\left( M, \mathcal{E}_{\mathcal{L}_0} \right)$, where
\begin{align*}
& \varphi \in A^{0,1}\left(M, \Theta_M \right)\\
& Y_j \in \Gamma\left( U_{j}, \mathcal{A}^{0,0}\left(\Omega_M^1\right) \right)\,\,\,\,\,\,\,\,\,\textnormal{with}\,\,\,\,\, Y_j = \Psi_{0jk} Y_k \\
& Q_j^\eta \in \Gamma\left(U_{j}, \mathcal{A}^{0,0} \right) \,\,\,\,\,\,\,\,\,\,\textnormal{with}\,\,\,\,\, Q_j^\eta= \Psi_{0jk} Q_k^\eta,\,\,\,\,\,\,\eta=1,..,m \\
&Z_j\in \Gamma\left( U_{j}, \mathcal{A}^{0,1} \right) \,\,\,\,\,\,\,\,\textnormal{with}\,\,\,\,\, Z_j- Z_k+ \left[ \varphi, \log \Psi_{0jk}\right] = 0
\end{align*}
with respect to the covering $\mathcal{U}_0=\{\mathcal{U}_j \cap M \}= \{U_{j}\}$ such that $\varphi(0)=0, Z_j(0)=0$ and $Y_j(0)=w_{0j}$ and
\begin{align}
&\bar{\partial} \varphi - \frac{1}{2}\left[ \varphi, \varphi \right] = 0 \label{uti4} \\
& \left( Y_j + \sum_{\eta=1}^m Q_j^\eta dt_\eta \right) \wedge \partial \left( Y_j  + \sum_{\eta=1}^m Q_j^\eta dt_\eta \right) = 0 \label{uti5}\\
&\bar{\partial} Y_j - \mathcal{L}_\varphi \left( Y_j \right) - Z_j Y_j =0 \label{uti6} \\
& \bar{\partial}  Q_j^\eta   - \left[ \varphi, Q_j^\eta \right] -  i_{\frac{\partial \varphi}{\partial t_\eta}}\left( Y_j \right) - Z_j Q_j^\eta = 0, \,\,\,\,\,\,\,\,\eta=1,...,m \label{uti7}\\
&\bar{\partial} Z_j - \left[ \varphi, Z_j \right] = 0  \label{uti8}
\end{align}

If we set $Q_j= \sum_{\eta=1}^m Q_j^\eta(z,t)dt_\eta$ and $\frac{\partial}{\partial t}= \sum_{\eta=1}^m c_\eta \frac{\partial}{\partial t_\eta}$. Then from $(\ref{uti5})$, we have
\begin{align*}
&\bar{\partial} \left(i_{\frac{\partial}{\partial t}} \left( Q_j \right)\right) -\left[ \varphi, i_{\frac{\partial}{\partial t}}\left(Q_j \right) \right] - i_{\frac{\partial  \varphi}{\partial t}}\left(Y_j \right)- Z_j i_{\frac{\partial}{\partial t}}\left( Q_j \right)=\sum_{\eta=1}^m c_\eta \bar{\partial} Q_j^\eta - \sum_{\eta=1}^m c_\eta \left[ \varphi, Q_j^\eta  \right] - \sum_{\eta=1}^m c_\eta i_{\frac{\partial \varphi}{\partial t}}\left( Y_j \right)- \sum_{\eta=1}^m c_\eta Z_j Q_j^\eta  =0
\end{align*}

Hence we have
\begin{align}\label{cu3}
\bar{\partial} \left(i_{\frac{\partial}{\partial t}} \left( Q_j \right)\right) -\left[ \varphi, i_{\frac{\partial}{\partial t_\eta}}\left(Q_j \right) \right] - i_{\frac{\partial  \varphi}{\partial t}}\left(Y_j \right)- Z_j i_{\frac{\partial}{\partial t}}\left( Q_j \right) = 0 
\end{align}
On the other hand, from $(\ref{uti5})$, we have
\begin{align}
&Y_j \wedge \partial Y_j=0 \label{uti9}\\
& Q_j^\eta  \partial Y_j + Y_j \wedge \partial Q_j^\eta - Y_j \wedge \frac{\partial Y_j}{\partial dt_\eta}=0,\,\,\,\,\,\eta=1,...,m\notag
\end{align}
Then we have
\begin{align} \label{uti10}
i_{\frac{\partial}{\partial t}}\left( Q_j \right) \partial Y_j + Y_j \wedge \partial \left( i_{\frac{\partial}{\partial t}}\left(Q_j \right)\right) - Y_j \wedge \frac{\partial Y_j}{\partial t}   =  0
\end{align}

By taking the derivative of $(\ref{uti4}),(\ref{uti6}),(\ref{cu3}),(\ref{uti9})$ with respect to $\frac{\partial}{\partial t}=\sum_{\eta=1}^n c_\eta \frac{\partial}{\partial t_\eta}$ and from $(\ref{uti9}), (\ref{uti10})$, we have
\begin{align}
&\bar{\partial}\left(\frac{\partial \varphi}{\partial t} \right) - \frac{1}{2}\left[ \frac{\partial \varphi}{\partial t}, \varphi \right] - \frac{1}{2}\left[ \varphi, \frac{\partial \varphi }{\partial t}\right]  = 0 \label{uti11}\\
& \frac{\partial Y_j}{\partial t} \wedge \partial Y_j + Y_j \wedge \partial \left( \frac{\partial Y_j}{\partial t} \right) =0\\
&i_{\frac{\partial}{\partial t}}\left( Q_j \right) \partial Y_j + Y_j \wedge \partial \left( i_{\frac{\partial}{\partial t}}\left(Q_j \right)\right) - Y_j \wedge \frac{\partial Y_j}{\partial t}   =  0 \\
& \bar{\partial} \left( \frac{\partial Y_j}{\partial t} \right) - \mathcal{L}_{\frac{\partial \varphi}{\partial t}}\left(Y_j \right) - \mathcal{L}_{\varphi}\left(\frac{\partial Y_j}{\partial t} \right) - \frac{\partial Z_j}{\partial t} Y_j - Z_j \frac{\partial Y_j}{\partial t} =0\\
&\bar{\partial} \left(i_{\frac{\partial}{\partial t}} \left( Q_j \right)\right) -\left[ \varphi, i_{\frac{\partial}{\partial t}}\left(Q_j \right) \right] - i_{\frac{\partial  \varphi}{\partial t}}\left(Y_j \right)- Z_j i_{\frac{\partial}{\partial t}}\left( Q_j \right) = 0 \\
&\bar{\partial} \left(\frac{\partial Z_j}{\partial t} \right) - \left[ \frac{\partial \varphi}{\partial t}, Z_j \right] -  \left[ \varphi, \frac{\partial Z_j}{\partial t} \right] = 0  \label{uti12}
\end{align}
By setting setting $t=0$,  since $\varphi(0)=0, Z_j(0)=0$, we have
\begin{align*}
&\bar{\partial} \left(\frac{\partial \varphi}{\partial t}|_{t=0 } \right)=0\\
&\bar{\partial}\left( \frac{\partial Y_j}{\partial t}|_{t=0} \right)  - \mathcal{L}_{\frac{\partial \varphi}{\partial t}|_{t=0}}\left( w_{0j} \right) - \frac{\partial Z_j}{\partial t}|_{t=0} \cdot w_{0j} = 0 \\
&\bar{\partial} \left(i_{\frac{\partial}{\partial t}} \left( Q_j \right)|_{t=0}\right) - i_{\frac{\partial  \varphi}{\partial t}|_{t=0}}\left(w_{0j} \right) = 0  \\
&\frac{\partial Y_j}{\partial t}|_{t=0} \wedge dw_{0j} + w_{0j} \wedge \partial \left( \frac{\partial Y_j}{\partial t}|_{t=0} \right) =0 \\
&i_{\frac{\partial}{\partial t}}\left( Q_j \right)|_{t=0} \cdot dw_{0j} + w_{0j} \wedge \partial \left( i_{\frac{\partial}{\partial t}}\left(Q_j \right)_{t=0} \right) - w_{0j} \wedge \frac{\partial Y_j}{\partial t}|_{t=0}  =  0 \\
&\bar{\partial} \left( \frac{\partial Z_j}{\partial t}|_{t=0} \right)=0
\end{align*}
This implies that
\begin{align*}
\left( \left\{ \frac{\partial Y_j}{\partial t}|_{t=0} \right\} , \left\{ i_{\frac{\partial}{\partial t}}\left(Q_j \right)|_{t=0} \right\} , \left\{\left(\frac{\partial \varphi}{\partial t}|_{t=0}, \frac{\partial Z_j}{\partial t}|_{t=0} \right)  \right\}    \right) \in  A^{0,0}\left( M, \Omega_M^1\otimes \mathcal{L}_0 \right) \bigoplus A^{0,0}\left( M, \mathcal{L}_0 \right) \bigoplus A^{0,1}\left( M, \mathcal{E}_{\mathcal{L}_0} \right)
\end{align*}
defines a $1$-cocycle in the complex associated with the following double complex 
{\small{\begin{equation}\label{uti14}
\begin{CD}
\cdots \\
@AAA \\
A^{0,0}\left( M-S, \tilde{\mathcal{S}}^2 \otimes \mathcal{L}_0 \right)\bigoplus A^{0,0}\left( M- S, \tilde{\mathcal{S}}^1\otimes \mathcal{L}_0 \right) @>\bar{\partial}>> \cdots \\
@AAA \\
A^{0,0}\left( M, \Omega_M^1\otimes \mathcal{L}_0 \right) \bigoplus A^{0,0}\left( M, \mathcal{L}_0 \right) @>-\bar{\partial}>>  A^{0,1}\left( M, \Omega_M^1\otimes \mathcal{L}_0  \right) \bigoplus A^{0,1}\left( M, \mathcal{L}_0 \right) @>\bar{\partial} >> \cdots \\
@AAA @AAA @AAA \\
A^{0,0}\left( M, \mathcal{E}_{\mathcal{L}_0} \right) @>\bar{\partial}>> A^{0,1}\left( M, \mathcal{E}_{\mathcal{L}_0} \right) @>-\bar{\partial}>> A^{0,2}\left( M, \mathcal{E}_{\mathcal{L}_0} \right) @>\bar{\partial}>> \cdots
\end{CD}
\end{equation}}}

Then we have the unfolded Kodaira-Spencer map in terms of sections
\begin{align} \label{cu5}
\varphi_0: T_0 B &\to \mathbb{H}^1\left( M, \mathcal{E}_{\mathcal{L}_0}^\bullet \right)\\
 \frac{\partial}{\partial t} &\mapsto \left( \frac{\partial \left( \mathcal{M}, \omega \right)}{\partial t}\right)|_{t=0}:= \left( \left\{ \frac{\partial Y_j}{\partial t}|_{t=0}  \right\}, \left\{ i_{\frac{\partial}{\partial t}}\left(Q_j \right)|_{t=0} \right\}, \left\{ \left( \frac{\partial \varphi}{\partial t}|_{t=0},  \frac{\partial Z_j}{\partial t}|_{t=0} \right) \right\} \right) \notag
\end{align}

Now we assume that the unfolded analytic family $\left(\mathcal{M}, \omega, B, \pi \right)$ with $\omega \in  H^0\left(\mathcal{M}, \Omega_{\mathcal{M}}^1\otimes \mathcal{L}  \right)$ in terms of sections is a $1$-parameter family, i.e. $0\in \Delta \subset B\subset \mathbb{C}$. By taking the derivative of $(\ref{uti11})-(\ref{uti12})$ with respect to $t$ again and setting $t=0$, we have
{\small{\begin{align*}
&\bar{\partial}\left( \frac{\partial^2 \varphi}{\partial t^2}|_{t=0} \right) - \left[ \frac{\partial \varphi}{\partial t}|_{t=0}, \frac{\partial \varphi}{\partial t}|_{t=0} \right]=0\\
&\frac{\partial^2 Y_j}{\partial t^2}|_{t=0} \wedge dw_{0j} + 2 \frac{\partial Y_j}{\partial t}|_{t=0}\wedge \partial \left(\frac{\partial Y_j}{\partial t}|_{t=0} \right) + w_{0j} \wedge \partial \left( \frac{\partial^2 Y_j}{\partial t^2} \right)= 0\\
&\frac{\partial}{\partial t}\left(i_{\frac{\partial}{\partial t}}\left(Q_j \right) \right)|_{t=0} \cdot dw_{0j} + i_{\frac{\partial}{\partial t}}\left(Q_j \right)|_{t=0}\cdot \partial \left( \frac{\partial Y_j}{\partial t}|_{t=0} \right) + \frac{\partial Y_j}{\partial t}|_{t=0} \wedge \partial \left( i_{\frac{\partial}{\partial t}}\left( Q_j \right)|_{t=0}  \right) - \frac{\partial Y_j}{\partial t}|_{t=0} \wedge \frac{\partial Y_j}{\partial t}|_{t=0} - w_{0j} \wedge \frac{\partial^2 Y_j}{\partial t^2}|_{t=0}=0 \\
& \bar{\partial}\left( \frac{\partial^2 Y_j}{\partial t^2}|_{t=0} \right)  - \mathcal{L}_{\frac{\partial^2 \varphi}{\partial t^2 }|_{t=0}}\left(w_{0j}\right)- 2 \mathcal{L}_{\frac{\partial \varphi}{\partial t}|_{t=0}}\left( \frac{\partial Y_j}{\partial t}|_{t=0} \right) - \frac{\partial^2 Z_j}{\partial t^2}|_{t=0} \cdot w_{0j} - 2\frac{\partial Z_j}{\partial t}|_{t=0}\cdot \frac{\partial Y_j}{\partial t}|_{t=0} = 0 \\
& \bar{\partial}\left( \frac{\partial}{\partial t}\left( i_{\frac{\partial}{\partial t}}\left(Q_j \right)\right)|_{t=0} \right) - \left[ \frac{\partial \varphi}{\partial t}|_{t=0}, i_{\frac{\partial}{\partial t}}\left( Q_j \right)|_{t=0} \right] - i_{\frac{\partial^2 \varphi}{\partial t^2}}(w_j)- i_{\frac{\partial \varphi}{\partial t}}|_{t=0}\left( \frac{\partial Y_j}{\partial t}|_{t=0}  \right) - \frac{\partial Z_j}{\partial t}|_{t=0}\cdot i_{\frac{\partial}{\partial t}}\left( Q_j \right)|_{t=0} = 0 \\
& \bar{\partial}\left( \frac{\partial^2 Z_j}{\partial t^2}|_{t=0} \right) - 2\left[ \frac{\partial \varphi}{\partial t}|_{t=0}, \frac{\partial Z_j}{\partial t}|_{t=0} \right] = 0 
\end{align*}}}

Then 
{\small{\begin{align}\label{uti15}
&\bigoplus \left(\left\{-2 \frac{\partial Y_j}{\partial t}|_{t=0} \wedge d\left( \frac{\partial Y_j}{\partial t}|_{t=0} \right)\right\} , \left\{ - \frac{\partial Y_j}{\partial t}|_{t=0} \wedge \partial \left( i_{\frac{\partial}{\partial t}}\left( Q_j \right)|_{t=0} \right)\right\}\right)\\
&\bigoplus \left( \left\{- 2 \mathcal{L}_{\frac{\partial \varphi}{\partial t}|_{t=0}}\left( \frac{\partial Y_j}{\partial t}|_{t=0} \right) -2 \frac{\partial Z}{\partial t}|_{t=0} \cdot \frac{\partial Y_j}{\partial t}|_{t=0} \right\}, \left\{- \left[ \frac{\partial \varphi}{\partial t} , i_{\frac{\partial}{\partial t}}\left( Q_j \right)|_{t=0} \right] - i_{\frac{\partial \varphi}{\partial t}|_{t=0}}\left( \frac{\partial Y_j}{\partial t}|_{t=0}  \right) - \frac{\partial Z_j}{\partial t}|_{t=0}\cdot i_{\frac{\partial }{\partial t}}\left( Q_j \right)|_{t=0} \right\}\right)\notag \\
&\bigoplus \left(\left\{- \left[ \frac{\partial \varphi}{\partial t}|_{t=0}, \frac{\partial \varphi}{\partial t}|_{t=0} \right] \right\} , \left\{ -2 \left[ \frac{\partial \varphi}{\partial t}|_{t=0}, \frac{\partial Z_j}{\partial t}|_{t=0} \right] \right\} \right) \notag
\end{align}}}
is a $2$-cocycle in the above bicomplex $(\ref{uti14})$, i.e. the image of $\left(\left\{\frac{\partial^2 Y_j}{\partial t^2}|_{t=0}\right\}, \left\{\frac{\partial}{\partial t}\left( i_{\frac{\partial}{\partial t}}\left( Q_j\right) \right)|_{t=0}\right\}  \right) \bigoplus \left( \frac{\partial^2 \varphi}{\partial t^2}, \left\{ \frac{\partial^2 Z_j}{\partial t^2}|_{t=0} \right\} \right)$
For any $1$-cocycle $(\{a_j \}, \{b_j\}, \{ \psi, \{c_j\})\in A^{0,0}\left(M, \Omega_M^1\otimes \mathcal{L}_0 \right)\bigoplus A^{0,0}\left( M, \mathcal{L}_0 \right)\bigoplus A^{0,1}\left( M, \mathcal{E}_{\mathcal{L}_0} \right)$ in the bicomplex $(\ref{uti14})$, we can define a $2$-cocycle by (we omit the proof, see Part IV-1)
\begin{align} \label{uti16}
&\left(\left\{ 2 a_j \wedge \partial a_j \right\} , \left\{  a_j \wedge \partial b_j\right\}\right) \bigoplus \left( \left\{ 2 \mathcal{L}_{\psi}(a_j)+ 2 c_j\cdot a_j \right\}, \left\{  [\psi, b_j]+i_\psi(a_j) + c_j \cdot b_j \right\} \right) \bigoplus \left( [\psi, \psi], \left\{  2[\psi, c_j] \right\}\right)\\
&\in A^{0,0}\left( M-S, \left( \tilde{\mathcal{S}}^2\otimes \mathcal{L}_0 \right) \bigoplus \left( \tilde{\mathcal{S}}^1\otimes \mathcal{L}_0 \right)\right) \bigoplus A^{0,1}\left( M, \left(\Omega_M^1\otimes \mathcal{L}_0 \right) \bigoplus \mathcal{L}_0 \right) \bigoplus A^{0,2}\left( M, \mathcal{E}_{\mathcal{L}_0} \right) \notag
\end{align}

Then $(\ref{uti15})$ implies that
\begin{theorem}
Suppose given a compact complex manifold $M$ with a codimension $1$-foliation defined by $w_0\in H^0\left(M, \Omega_M^1\otimes \mathcal{L}_0 \right)$ for a line bundle $\mathcal{L}_0$, and 
\begin{align*}
(\{a_j\}, \{b_j\}, \{\psi, \{c_j\}\})\in A^{0,0}\left( M, \Omega_M^1 \otimes \mathcal{L}_0 \right) \bigoplus A^{0,0}\left(M, \mathcal{L}_0 \right) \bigoplus A^{0,1}\left( M, \mathcal{E}_{\mathcal{L}_0} \right) \,\,\,\textnormal{in}\,\,\, \mathbb{H}^1\left(M, \tilde{\mathcal{E}}_{\mathcal{L}_0}^\bullet \right)
\end{align*}
In order that there may exist an unfolded analytic family $\left( \mathcal{M}, \omega , B, \pi\right)$ such that $\pi^{-1}(0)=\left( M, w_0 \right)$ and $\left(\frac{\partial \left(\mathcal{M}, \omega \right)}{\partial t}\right)|_{t=0}=\left(\left\{a_j \right\}, \left\{ b_j \right\}, \left\{\left( \psi , c_{j} \right) \right\} \right)$, it is necessary  that $(\ref{uti16})$ defines the $0$ class in $\mathbb{H}^2\left( M, \tilde{\mathcal{E}}_{\mathcal{L}_0}^\bullet \right)$.

\end{theorem}

\section{Rigid foliated deformations and rigid unfoldings}

\begin{definition}[rigidity in terms of tangent sheaves (compare \cite{MK06} Definition 3.3 p.45)]
We say that a compact foliated complex manifold $\left(\mathcal{M}, \Theta_{\mathcal{F}_0}\right)$ with $\Theta_{\mathcal{F}_0}$ locally free is rigid in foliated deformations in terms of tangent sheaves, for any foliated complex analytic family $\left(\mathcal{M}, \Theta_\mathcal{F}, B, \pi \right)$ in terms of tangent sheaves such that $\pi^{-1}\left(t_0\right)=\left( M, \Theta_{\mathcal{F}_0}\right)$, we can find a neighborhood $\Delta$ of $t_0$ such that the restriction $\left(\mathcal{M}_\Delta, \Theta_{\mathcal{F}_\Delta}, \Delta, \pi\right)$ is complex analytically isomorphic to $\left(M\times \Delta, \Theta_{\mathcal{F}_0}\right)$, where $\Theta_{\mathcal{F}_0}$ on $M\times \Delta$ is the pullback of $\Theta_{\mathcal{F}_0}$ on $M$ under the natural projection $M\times \Delta \to M$.
\end{definition}

\begin{definition}[rigidity in terms of cotangent sheaves]
We say that a compact foliated complex manifold $\left(\mathcal{M}, \mathcal{N}_{\mathcal{F}_0}^*\right)$ with $\mathcal{N}_{\mathcal{F}_0}^*$ locally free is rigid in foliated deformations in terms of cotangent sheaves, for any foliated complex analytic family $\left(\mathcal{M}, \mathcal{N}_\mathcal{F}^*, B, \pi \right)$ in terms of cotangent sheaves such that $\pi^{-1}\left(t_0\right)=\left( M, \mathcal{N}_{\mathcal{F}_0}^*\right)$, we can find a neighborhood $\Delta$ of $t_0$ such that the restriction $\left(\mathcal{M}_\Delta, \mathcal{N}_{\mathcal{F}_\Delta}^*, \Delta, \pi\right)$ is complex analytically isomorphic to $\left(M\times \Delta, \mathcal{N}_{\mathcal{F}_0}^*\right)$, where $\mathcal{N}_{\mathcal{F}_0}^*$ on $M\times \Delta$ is the pullback of $\mathcal{N}_{\mathcal{F}_0}^*$ on $M$ under the natural projection $M\times \Delta \to M$.
\end{definition}

\begin{definition}[rigidity in terms of both tangent and cotangent sheaves]
We say that a compact foliated complex manifold $\left(\mathcal{M}, \Theta_{\mathcal{F}_0}, \mathcal{N}_{\mathcal{F}_0}^* \right)$ with $\Theta_{\mathcal{F}_0}$ and $\mathcal{N}_{\mathcal{F}_0}^*$ locally free is rigid in foliated deformations in terms of both tangent and cotangent sheaves, for any foliated complex analytic family $\left(\mathcal{M}, \Theta_\mathcal{F}, B, \pi \right)$ in terms of both tangent and cotangent sheaves such that $\pi^{-1}\left(t_0\right)=\left( M, \Theta_{\mathcal{F}_0} , \mathcal{N}_{\mathcal{F}_0}^* \right)$, we can find a neighborhood $\Delta$ of $t_0$ such that the restriction $\left(\mathcal{M}_\Delta, \Theta_{\mathcal{F}_\Delta},  \mathcal{N}_{\mathcal{F}_\Delta}^*, \Delta , \pi\right)$ is complex analytically isomorphic to $\left(M\times \Delta, \Theta_{\mathcal{F}_0} , \mathcal{N}_{\mathcal{F}_0}^* \right)$.
\end{definition}

\begin{definition}[rigidity of unfoldings]
We say that a compact foliated complex manifold $\left(\mathcal{M}, \mathcal{N}_{\mathcal{F}_0}^*\right)$ with $\mathcal{N}_{\mathcal{F}_0}^*$ locally free is rigid in unfoldings, for any unfolded complex analytic family $\left(\mathcal{M}, \tilde{\mathcal{N}}_\mathcal{F}^*, B, \pi \right)$ in terms of cotangent sheaves such that $\pi^{-1}\left(t_0\right)$ induces $\left( M, \mathcal{N}_{\mathcal{F}_0}^*\right)$, we can find a neighborhood $\Delta$ of $t_0$ such that the restriction $\left(\mathcal{M}_\Delta, \tilde{\mathcal{N}}_{\mathcal{F}_\Delta}^*, \Delta, \pi\right)$ is complex analytically isomorphic to $\left(M\times \Delta, \mathcal{N}_{\mathcal{F}_0}^*\right)$, where $\mathcal{N}_{\mathcal{F}_0}^*$ on $M\times \Delta$ is the pullback of $\mathcal{N}_{\mathcal{F}_0}^*$ on $M$ under the natural projection $M\times \Delta \to M$.
\end{definition}

As corollaries of Theorems of completeness in Part I and Part III, we have

\begin{theorem}[rigidity in terms of tangent sheaves]\label{cp12}
Let $\left(M, \Theta_{\mathcal{F}_0}\right)$ be a compact foliated complex manifold with $\Theta_{\mathcal{F}_0}$ locally free. If $\mathbb{H}^1\left(M, \Theta_{\mathcal{F}_0}^\bullet \right)=0$, then $\left(M, \Theta_{\mathcal{F}_0} \right)$ is rigid in foliated deformations in terms of tangent sheaves.
\end{theorem}

\begin{theorem}[rigidity in terms of cotangent sheaves] \label{rsf25}
Let $\left(M, \mathcal{N}_{\mathcal{F}_0}^* \right)$ be a compact foliated complex manifold with $\mathcal{N}_{\mathcal{F}_0}^*$ locally free. If $\mathbb{H}^1\left(M, \mathcal{N}_{\mathcal{F}_0}^{*\bullet} \right)=0$, then $\left(M, \mathcal{N}_{\mathcal{F}_0}^* \right)$ is rigid in foliated deformations in terms of cotangent sheaves.
\end{theorem}

\begin{theorem}[rigidity in terms of both tangent and cotangent sheaves]
Let $\left(M, \Theta_{\mathcal{F}_0}, \mathcal{N}_{\mathcal{F}_0}^*\right)$ be a compact foliated complex manifold with $\Theta_{\mathcal{F}_0}$ and $\mathcal{N}_{\mathcal{F}_0}^*$ locally free. If $\mathbb{H}^1\left(M, \mathcal{F}_0^\bullet \right)=0$, then $\left(M, \Theta_{\mathcal{F}_0} , \mathcal{N}_{\mathcal{F}_0}^* \right)$ is rigid in foliated deformations in terms of both tangent sheaves and cotangent sheaves.
\end{theorem}

\begin{theorem}[rigidity of unfoldings (compare \cite{MK06} Theorem 3.2 p.45]\label{r1}
Let $\left(M, \mathcal{N}_{\mathcal{F}_0}^* \right)$ be a compact foliated complex manifold with $\mathcal{N}_{\mathcal{F}_0}^*$ locally free. If $\mathbb{H}^1\left(M, \tilde{\mathcal{N}}_{\mathcal{F}_0}^{*\bullet} \right)=0$, then $\left(M, \mathcal{N}_{\mathcal{F}_0}^* \right)$ is rigid in unfoldings.
\end{theorem}

In Part III we have shown that if $\mathbb{H}^1\left(M, \tilde{\mathcal{N}}_{\mathcal{F}_0}^{*\bullet} \right)=0$ for a compact foliated complex manifold $\left( M, \mathcal{N}_{\mathcal{F}_0}^* \right)$ with $\mathcal{N}_{\mathcal{F}_0}^*$ locally free, then any infinitesimal unfolding $\left( \mathcal{M}, \tilde{\mathcal{N}}_{\mathcal{F}}^*\right)$ of $\left( M, \mathcal{N}_{\mathcal{F}_0}^* \right)$ over $\textnormal{Spec}\left(k[x]/(x^n)\right)$ is trivial. However the author could not prove that any infinitesimal unfolding $\left( \mathcal{M}, \tilde{\mathcal{N}}_{\mathcal{F}}^* \right)$ of $\left(M, \mathcal{N}_{\mathcal{F}_0}^* \right)$ over $A$ is trivial for arbitrary $A\in \bold{Art}$. However, in the light of Theorem \ref{r1}, if $\mathbb{H}^1\left( M,  \tilde{\mathcal{N}}_{\mathcal{F}_0}^{* \bullet} \right)=0$, then any infinitesimal unfolding $\left(\mathcal{X}, \tilde{\mathcal{N}}_{\mathcal{F}}^* \right)$ of $\left( X, \mathcal{N}_{\mathcal{F}_0}^* \right)$ over $A$ should be trivial. We will directly give the proof of Theorem \ref{r1}, which is in essence the repeat of the proof of Theorem of completeness for unfoldings in Part III. 

\begin{proof}[Proof of \textnormal{Theorem \ref{r1}}]
Since the result is local, we may assume that $B=\left\{ t=\left(t_1,..., t_m\right)\in \mathbb{C}^m ||t|<1\right\}$ and $t_0=0$ and $\mathcal{M}$ is written in the form
\begin{align*}
\mathcal{M}=\bigcup_j \,\,\, \mathcal{U}_j,\,\,\,\,\,\,\,\,\mathcal{U}_j=\left\{\left( z_j, t\right) \in \mathbb{C}^n \times B\,\,\, |\,\,\, |z_j| <1 \right\}, \,\,\,\,\,\,\,\,\textnormal{with}\,\,\,\,\,z_j^\alpha= f_{jk}^\alpha(z_k,t)\,\,\,\,\,\textnormal{on}\,\,\,\mathcal{U}_j\cap \mathcal{U}_k
\end{align*}
and the locally free sheaf $\tilde{\mathcal{N}}_\mathcal{F}^*$ is trivialized in terms of the covering $\left\{\mathcal{U}_j \right\}$ so that $\Gamma\left(\mathcal{U}_j, \tilde{\mathcal{N}}_{\mathcal{F}}^*\right)$ is generated by holomorphic forms $w_j^\alpha(z_j,t)\in \Gamma\left(\mathcal{U}_j, \tilde{\mathcal{N}}_\mathcal{F}^* \right),\alpha=1,...,q$ of the form
\begin{align*}
\omega_j^\alpha(z_j,t):=w_j^\alpha(z_j,t)+v_j^\alpha(z_j,t)=\sum_{\beta=1}^n w_j^{\alpha\beta}(z_j,t_j)dz_j^\beta+\sum_{\eta=1}^m v_j^{\alpha\eta}(z_j,t_j) dt_j^\eta,\,\,\,\,\alpha=1,...,q
\end{align*}
and for $j,k$ with $\mathcal{U}_j\cap \mathcal{U}_k \ne \emptyset$, we have invertible $q\times q$ matrices $H_{jk}:=\left(h_{jk}^{\alpha\beta}(z_j,t)\right)$ with components $h_{jk}^\alpha(z_k,t)\in \Gamma(\mathcal{U}_j\cap \mathcal{U}_k, \mathcal{O}_\mathcal{M})$ such that $\omega_j^\alpha(z_j,t)=\sum_{\beta=1}^q h_{jk}^{\alpha\beta}\omega_k^\beta(z_k,t)$, equivalently
\begin{align}
\sum_{\gamma=1}^n w_j^{\alpha\gamma}(z_j,t)dz_j^\gamma+\sum_{\eta=1}^m v_j^{\alpha\eta}(z_j,t)dt_j^\eta=\sum_{\beta=1}^q h_{jk}^{\alpha\beta}(z_k,t)\left(\sum_{\xi=1}^n w_k^{\beta\xi}(z_k,t)dz_k^\xi+\sum_{\eta=1}^m v_k^{\beta \eta}(z_k,t)dt_k^\eta \right)
\end{align}
and the integrability condition holds:
\begin{align*}
\omega_j^1(z_j,t)\wedge \cdots \wedge \omega_j^q(z_j,t)\wedge d\omega_j^\alpha(z_j,t)=0,\,\,\,\,\,\alpha=1,...,q\,\,\,\textnormal{on}\,\,\,\mathcal{U}_j.
\end{align*}
and by setting $U_i:=M \cap \mathcal{U}_i=\left\{ z_i\in \mathbb{C}^n | |z_i|<1 \right\}$, $\left( M, \mathcal{N}_{\mathcal{F}_0}^* \right)=\pi^{-1}(0)$ is described in terms of covering $\mathcal{U}:=\{U_i\}=\left\{ M \cap \mathcal{U}_i \right\}$ in the following way: (1) local coordinates $z_i= \left(z_i^1,..., z_i^n \right)$ on $U_i$ with $z_i= g_{ij}^\alpha(z_j)= \left( g_{ij}^1(z_j),..., g_{ij}^n(z_j) \right)$ where we set $g_{ij}(z_j):= f_{ij}(z_j,0)$ and (2) $\Gamma\left(U_i, \mathcal{N}_{\mathcal{F}_0}^* \right)$ is generated by $w_{0i}^\alpha(z_i):= \sum_{\beta=1}^n w_i^{\alpha\beta}(z_i,0)dz_i^\beta \in \Gamma\left(U_i, \mathcal{N}_{\mathcal{F}_0}^* \right)$ with $w_{0i}^\alpha= \sum_{\beta=1}^q h_{0ij}^{\alpha\beta}(z_j) w_{0j}^\beta$ where we set $h_{0ij}^{\alpha\beta}(z_j):=h_{ij}^{\alpha\beta}(z_j,0)$, and $w_{0i}^1\wedge \cdots \wedge w_{0i}^q \wedge dw_{0i}^\alpha=0$.

In order to prove Theorem \ref{r1} it suffices to construct holomorphic functions
\begin{align*}
\varphi_i:\mathcal{U}_i \to \mathbb{C}^n
\end{align*}
and matrix functions
\begin{center}
$\left[\begin{matrix}
c_{i}^{11}(z_i,t) & c_{i}^{12}(z_i,t) & \cdots & c_{i}^{1q}(z_i, t)\\
c_{i}^{21}(z_i, u) & c_{i}^{22}(z_i, u) & \cdots & c_{i}^{2q}(z_i, u)\\
\cdot & \cdot & \cdots & \cdot\\
\cdot & \cdot & \cdots & \cdot\\
\cdot & \cdot & \cdots & \cdot\\
c_{i}^{q1}(z_i, t) & c_{i}^{q2}(z_i, t) & \cdots & c_{i}^{qq}(z_i, t)\\
\end{matrix}\right]$
\end{center}
where $c_i^{\alpha\beta}:\mathcal{U}_i \to \mathbb{C}$ with $\sum_{\beta=1}^q c_i^{\alpha\beta}(z_i,t) \omega_i^\beta(z_i,t)= \varphi_i^*\left( w_{0i}^\alpha \right)$ such that
\begin{align}
\varphi_i(z_i, 0)=\xi_i, &\,\,\,\,\, c_i^{\alpha\beta}(\xi_i,0)= \delta_\alpha^\beta \\
\varphi_i \left(f_{ij}(z_j,t), t \right) &= g_{ij}\left(\varphi_j(z_j,t)\right) \label{r2}
\end{align}
\begin{align}
\sum_{\beta=1}^q c_i^{\alpha\beta}\left(f_{ij}(z_j, t), t\right) h_{ij}^{\beta\gamma}(z_j , t) = \sum_{\beta=1}^q h_{0ij}^{\alpha\beta} \left (\varphi_j(z_j, t)\right) c_j^{\beta\gamma}(z_j,t)
\end{align}
and for $\alpha=1,...,q$,
\begin{align}
\sum_{\beta=1}^q\sum_{\gamma=1}^n c_i^{\alpha\beta}(\xi_i, u) w_i^{\beta\gamma}(z_i, t)dz_i^\gamma &= \sum_{\beta=1}^n \sum_{\gamma=1}^n w_{0i}^{\alpha\beta}\left(\varphi_i(z_i, t)\right)\frac{\partial \varphi_i^\beta}{\partial z_i^\gamma} dz_i^\gamma 
\end{align}
{\small{\begin{align}
\sum_{\beta=1}^q\sum_{\delta=1}^{m} c_i^{\alpha\beta}(z_i, t) v_i'^{\beta\delta}(z_i, t)dt_\delta & = \sum_{\delta=1}^{m}\sum_{\beta=1}^q w_{0i}^{\alpha\beta}\left( \varphi_i(z_i, t) \right)\frac{\partial \varphi_i^\beta}{\partial t_\delta}  dt_\delta. \label{r3}
\end{align}}}
First we prove the existence of formal solution of $(\ref{r2} )-(\ref{r3})$. We recall the notation $`\,{\equiv_\mu}$' from Part I. Then $(\ref{r2})-(\ref{r3})$ are equivalent to the following systems of congruences
\begin{align}
\varphi_i^\mu \left(f_{ij}(z_j, t), t \right)& \equiv_\mu g_{ij}\left(\varphi_j^\mu(z_j, t)\right) \label{uc6}\\
\sum_{\beta=1}^q c_i^{\alpha\beta\mu}\left(f_{ij}(z_j, t),t\right) h_{ij}^{\beta\gamma}& (z_j, t)  \equiv_\mu \sum_{\beta=1}^q  h_{0ij}^{\alpha\beta}\left(\varphi_j^\mu(z_j, t)\right) c_j^{\beta\gamma\mu}(z_j, t) \label{uc7} \\
\sum_{\beta=1}^q\sum_{\gamma=1}^n c_i^{\alpha\beta \mu}(z_i, t) w_i^{\beta\gamma}(z_i, t)dz_i^\gamma &= \sum_{\beta=1}^n \sum_{\gamma=1}^n w_{0i}^{\alpha\beta}\left(\varphi_i^\mu (z_i, u)\right)\frac{\partial \varphi_i^{ \beta \mu}}{\partial z_i^\gamma}  dz_i^\gamma \label{uc8}
\end{align}
\begin{align}
\sum_{\beta=1}^q\sum_{\delta=1}^{m} c_i^{\alpha\beta\mu}(z_i, u) v_i^{\beta\delta}(z_i, t)dt_\delta & \equiv_\mu \sum_{\delta=1}^{m}  \sum_{\beta=1}^q w_{0i}^{\alpha\beta}\left(\varphi_i^\mu(z_i, t) \right)\frac{\partial \varphi_i^{\beta\mu}}{\partial t_\delta} dt_\delta \label{uc9}
\end{align}
for each $\mu=1,2,...$.  We shall construct $\varphi_i^\mu$, and $c_i^{\alpha\beta \mu}$ satisfying $(\ref{uc6})_\mu- (\ref{uc9})_\mu$ by induction on $\mu$. We assume that $\varphi_i^{\mu-1}$ and $c_i^{\alpha\beta(\mu-1)}$ are already determined. Then we define homogenous polynomials $\Gamma_{ij|\mu}^\alpha, A_{ij|\mu}^{\alpha\beta}, W_{i|\mu}^{\alpha \gamma}$ and $V_{i|\mu-1}^{\alpha \delta}$ of degree $\mu$ by the following congruences:
{\small{\begin{align}
\Gamma_{ij|\mu}^\alpha &\equiv_\mu \varphi_i^{\alpha(\mu-1)}\left(f_{ij}(z_j, t), t \right)- g_{ij}^{\alpha}\left( \varphi_j^{\mu-1}(z_j, t )\right) \label{uc10} \\
\sum_{\beta=1}^q A_{ij|\mu}^{\alpha\beta}  h_{0ij}^{\beta \gamma} &\equiv_\mu \sum_{\beta=1}^q c_i^{\alpha\beta (\mu-1)}\left(f_{ij}(z_j, t), t\right) h_{ij}^{\beta\gamma}(z_j, t) - \sum_{\beta=1}^q h_{0ij}^{\alpha\beta} \left(\varphi_j^{\mu-1}(z_j, t\right) )c_j^{\beta\gamma(\mu-1)}(z_j, t)\label{uc11}\\
W_{i|\mu}^{\alpha\gamma}& \equiv_{\mu} \sum_{\beta=1}^q  c_i^{\alpha\beta (\mu-1)}(z_i, u) w_i^{\beta\gamma}(z_i, t)- \sum_{\beta=1}^q w_{0i}^{\alpha\beta}\left(\varphi_i^{\mu-1}(z_i, t)\right) \frac{\partial \varphi_i^{\beta (\mu-1)}}{\partial z_i^\gamma} \label{uc12}
\end{align}}}
{\small{\begin{align}
V_{i|\mu-1}^{\alpha \delta} &\equiv_{\mu-1} \sum_{\beta=1}^q c_i^{\alpha\beta (\mu-1)}(z_i, t) v_i^{\beta\delta}(z_i, t)- \sum_{\beta=1}^q w_{0i}^{\alpha\beta}\left(\varphi_i^{\mu-1}(z_i, t) \right)\frac{\partial \varphi_i^{\beta (\mu-1)} }{\partial t_\delta}  \label{uc13}
\end{align}}}
We set 
\begin{align*}
\Gamma_{ij|\mu}&= \sum_{\alpha=1}^n \Gamma_{ij|\mu}^{\alpha}\frac{\partial}{\partial z_i^\gamma},\,\,\,\,\,\,\,\,\,\,\,\, W_{i|\mu}^\alpha = \sum_{\gamma=1}^n W_{i|\mu}^{\alpha\gamma} dz_i^\gamma
\end{align*}

Then as in the proof of Theorem of completeness for unfoldings from Part III, we have
\begin{lemma}\label{uc14}
We have the following equalities$:$ for $\alpha=1,...,q$ and $\delta=1,..., m$
\begin{align}
&\,\,\,\,\,\,\,\,\Gamma_{jk|\mu}- \Gamma_{ik|\mu} + \Gamma_{ij|\mu}=0\label{uc15}\\
W_{i|\mu}^\alpha &=  \sum_{\beta=1}^q A_{ij|\mu}^\beta w_{0i}^\beta + \sum_{\beta=1}^q h_{0ij}^{\alpha\beta} W_{j|\mu}^\beta- \mathcal{L}_{\Gamma_{ij|\mu}}\left(w_{0i}^\alpha \right) \label{uc16}\\
\sum_{\beta=1}^q w_{0i}^1 \wedge \cdots \wedge &\overbrace{\left(-  W_{i|\mu}^\beta \right)}^{\beta-\textnormal{th}} \wedge \cdots \wedge w_{0i}^q \wedge d w_{0i}^\alpha + w_{0i}^1 \wedge \cdots \wedge w_{0i}^q \wedge d \left( - W_{i|\mu}^\alpha \right)=0 \label{uc17}\\
V_{i|\mu-1}^{\alpha \delta}= \sum_{\beta=1}^q h_{0ij}^{\alpha\beta}& V_{j|\mu-1}^{\beta\delta} - i_{\sum_{\beta=1}^n \frac{\partial \Gamma_{ij|\mu}^\beta}{\partial t_\delta} \frac{\partial }{\partial z_i^\beta}} \left( w_{0i}^\alpha \right) = \sum_{\beta=1}^q h_{0ij}^{\alpha \beta} V_{j|\mu-1}^{\beta\delta}- i_{\frac{\partial \Gamma_{ij|\mu}}{\partial t_\delta}}\left(w_{0i}^\alpha \right) \label{uc18}
\end{align}
{\small{\begin{align} \label{uc19}
 \sum_{\beta=1}^q (-1)^{q-\beta} \left(- V_{i|\mu-1}^{\beta \delta}  \right) w_{0i}^1 \wedge \cdots \wedge &\widehat{ w_{0i}^{\beta} } \wedge \cdots \wedge  w_{0i}^q \wedge d w_{0i}^\alpha   + w_{0i}^1 \wedge \cdots \wedge w_{0i}^q \wedge d \left(- V_{i|\mu-1}^{\alpha \delta} \right)    - w_{0i}^1 \wedge \cdots \wedge w_{0i}^q \wedge \frac{\partial \left(- W_{i|\mu}^\alpha \right) }{\partial t_\delta }=0 
\end{align}}}
\begin{align}\label{uc20}
\frac{\partial V_{i|\mu-1}^{\alpha \delta}}{\partial t_\eta}= \frac{\partial  V_{i|\mu-1}^{\alpha \eta}}{\partial t_\delta}\,\,\,\,\,\,\,\,\,\textnormal{for $\delta\ne \eta$}
\end{align}
\end{lemma}

Our purpose is to determine $\varphi^\mu= \varphi^{\mu-1}+ \varphi_{i|\mu}$, and $c_i^{\alpha\beta(\mu-1)}+ c_{i|\mu}^{\alpha\beta}$ satisfying $(\ref{uc6})_\mu-(\ref{uc9})_\mu$.
\begin{lemma} \label{uc27}
$(\ref{uc6})_\mu -(\ref{uc9})_\mu$ are equivalent to the following equalities$:$ for $\alpha=1,...,q $ and $\delta=1,...,m$,
\begin{align}
\Gamma_{ij|\mu}&= \varphi_{j|\mu}-  \varphi_{i | \mu} \label{uc23}\\
\sum_{\beta=1}^q A_{ij|\mu}^{\alpha\beta} h_{0ij}^{\beta\gamma}  &=  -\sum_{\beta=1}^q c_{i|\mu}^{\alpha\beta} h_{0ij}^{\beta\gamma} + \sum_{\beta=1}^q \left[ \varphi_{j|\mu} , h_{0ij}^{\alpha\beta}  \right]  c_{0j}^{\beta\gamma}  - \sum_{\beta=1}^q h_{0ij}^{\beta\gamma} c_{j|\mu}^{\beta\gamma}    \label{uc24}   \\
W_{i|\mu}^{\alpha } &=  - \sum_{\gamma=1}^n \sum_{\beta=1}^q c_{i|\mu}^{\alpha\beta} w_{0i}^{\beta} + \mathcal{L}_{\varphi_{i|\mu}}\left( w_{0i}^\alpha  \right)  \label{uc25} \\
V_{i|\mu-1}^{\alpha \delta } &= i_{\frac{\partial \varphi_{i|\mu}}{\partial t_\delta}}\left( w_{0i}^\alpha \right)    \label{uc26}
\end{align}
where $\varphi_{i|\mu}=\sum_{\alpha=1}^n \varphi_{i|\mu}^\alpha \frac{\partial}{\partial z_i^\alpha}$
\end{lemma}

For $\delta=1,..., m$, we define an element $\overline{\frac{\partial W_{i|\mu}}{\partial u_\delta}}\in \Gamma\left( U_i, \mathscr{H}om_{\mathcal{O}_M} \left( \mathcal{N}_{\mathcal{F}_0}^*, \frac{\Omega_M^1}{\mathcal{N}_{\mathcal{F}_0}^*}  \right) \right)$ by
\begin{align*}
\overline{\frac{\partial W_{i|\mu}}{\partial t_\delta}}: \Gamma\left( U_i, \mathcal{N}_{\mathcal{F}_0}^*  \right) &\to \Gamma\left(  U_i  ,   \frac{\Omega_M^1}{\mathcal{N}_{\mathcal{F}_0}^*}   \right) \\
                                       w_{0i}^\alpha &\mapsto \overline{\frac{\partial W_{i|\mu}^\alpha}{\partial t_\delta} }
\end{align*}
and linearly extends to $\Gamma\left( U_i, \mathcal{N}_{\mathcal{F}_0}^* \right)$. For $\delta=1,...,m$, we define an element $V_{i|\mu}^\delta\in \Gamma\left( U_i, \mathscr{H}om_{\mathcal{O}_M}\left( \mathcal{N}_{\mathcal{F}_0}^*, \mathcal{O}_M \right)\right)$ by
\begin{align*}
V_{i|\mu}^\delta :\Gamma\left( U_i,  \mathcal{N}_{\mathcal{F}_0}^* \right) &\to \Gamma\left( U_i, \mathcal{O}_M \right)\\
                  w_{0i}^\alpha &\mapsto V_{i|\mu}^{\alpha \delta}
\end{align*}
and linearly extends to $\Gamma\left( U_i, \mathcal{N}_{\mathcal{F}_0}^* \right)$. Then $(\ref{uc15})-(\ref{uc19})$ implies that
\begin{align}\label{ucc12}
\left(  \left\{ -  \overline{\frac{\partial W_{i|\mu}}{\partial t_\delta}}\right\}   , \left\{- V_{i|\mu}^\delta \right\}  \left\{  \frac{\partial \Gamma_{ij|\mu}}{\partial t_\delta} \right\} \right) \in C^0\left( \mathcal{U}, \mathscr{H}om_{\mathcal{O}_M}\left( \mathcal{N}_{\mathcal{F}_0}^* ,  \frac{\Omega_M^1}{\mathcal{N}_{\mathcal{F}_0}^*} \right)   \right) \bigoplus C^0\left( \mathcal{U}, \mathscr{H}om_{\mathcal{O}_M}\left( \mathcal{N}_{\mathcal{F}_0}^*, \mathcal{O}_M \right) \right) \bigoplus C^1 \left( \mathcal{U}, \Theta_M  \right)
\end{align}
defines a $1$-cocycle in the following \v Cech resolution of $\tilde{\mathcal{N}}_{\mathcal{F}_0}^{* \bullet}$:
{\small{\begin{center}
$\begin{CD}
\cdots \\
@AG_2AA \\
C^0\left( \mathcal{U}, \left(\left(\mathcal{N}_{\mathcal{F}_0}^* \right)^*\otimes \tilde{\mathcal{S}}^2 \right) \bigoplus \left( \mathcal{N}_{\mathcal{F}_0}^* \right)^*\otimes \tilde{\mathcal{S}}^1 \right) @>-\delta >> \cdots \\
@AG_1AA @AG_1AA \\
C^0\left(\mathcal{U},\left( \left(\mathcal{N}_{\mathcal{F}_0}^* \right)^* \otimes \frac{\Omega_{M}^1}{ \mathcal{N}_{\mathcal{F}_0}^*} \right) \bigoplus \left( \mathcal{N}_{\mathcal{F}_0}^* \right)^* \right) @>\delta >> C^1\left(\mathcal{U}, \left( \left(\mathcal{N}_{\mathcal{F}_0}^* \right)^* \otimes \frac{\Omega_{M}^1}{ \mathcal{N}_{\mathcal{F}_0}^*} \right) \bigoplus   \left( \mathcal{N}_{\mathcal{F}_0}^* \right)^* \right) @>-\delta>> \cdots \\
@AG_0AA @AG_0AA @AG_0AA \\
C^0\left(\mathcal{U}, \Theta_{M}\right) @>-\delta>> C^1\left(\mathcal{U}, \Theta_{M}\right) @>\delta>> C^2\left(\mathcal{U}, \Theta_{M}\right) @>-\delta>> \cdots
\end{CD}$
\end{center}}}
By hypothesis $\mathbb{H}^1\left(M, \tilde{\mathcal{N}}_\mathcal{F}^{*\bullet}\right)=0$, there exists $\left\{\tilde{\varphi}_{i|\mu-1}^\delta:=\sum_{\alpha=1}^n \tilde{\varphi}_{i|\mu-1}^{\delta\alpha} \frac{\partial}{\partial z_i^\alpha} \right\} \in C^0 \left(\mathcal{U} , \Theta_M \right)$ and $\left\{ \tilde{c}_{i|\mu-1}^\delta \right\}\in C^0\left( \mathcal{U}, \mathscr{H}om_{\mathcal{O}_M} \left( \mathcal{N}_{\mathcal{F}_0}^*, \mathcal{N}_{\mathcal{F}_0}^*    \right) \right)$ where $\tilde{c}_{i|\mu-1}^\delta$ is defined by $\tilde{c}_{i|\mu-1}^\delta:\Gamma\left( U_i, \mathcal{N}_{\mathcal{F}_0}^* \right)\to \Gamma\left(U_i, \mathcal{N}_{\mathcal{F}_0}^* \right), w_{0i}^\alpha \mapsto \sum_{\gamma=1}^q \tilde{c}_{i|\mu-1}^{\delta \alpha \gamma} w_{0i}^\gamma$ such that
\begin{align}
\tilde{\varphi}_{i|\mu-1}^\delta- \tilde{\varphi}_{j|\mu-1}^\delta&=- \frac{\partial \Gamma_{ij|\mu}}{\partial t_\delta}  \label{uc30} \\
\mathcal{L}_{\tilde{\varphi}_{i|\mu-1}^\delta}\left(w_{0i}^\alpha\right)&= \frac{\partial W_{i|\mu}^\alpha}{\partial t_\delta} + \sum_{\gamma=1}^q \tilde{c}_{i|\mu-1}^{ \delta \alpha \gamma  } w_{0i}^\gamma  \label{uc31} \\
i_{\tilde{\varphi}_{i|\mu-1}^\delta}\left( w_{0i}^\alpha \right) &=V_{i|\mu-1}^{\alpha \delta}  \label{uc32}
\end{align}
Then by setting
\begin{align*}
\varphi_{i|\mu} := \frac{1}{\mu} \sum_{\delta=1}^{m} t_\delta \tilde{\varphi}_{i|\mu-1}^\delta,\,\,\,\,\,\,\,\,\,\,\, c_{i|\mu}^{\alpha \gamma} := \frac{1}{\mu} \sum_{\delta=1}^{m} t_\delta \tilde{c}_{i|\mu-1}^{\delta \alpha \gamma },
\end{align*}
$(\ref{uc23})-(\ref{uc26})$ are satisfied. This completes the inductive construction of $\varphi_i^\mu,$ and $c_i^{\alpha\beta \mu}$, and we can prove the convergence in a similar way as in the proof of Theorem of completeness for unfoldings in Part III.

\end{proof}

We will show that the codimension $1$-singular foliation $\mathcal{F}_0$ of the exceptional component on $\mathbb{P}_\mathbb{C}^3$ (see \cite{CLN96} Example 6 p.580) is rigid in foliated deformations and unfoldings by explicitly computing $\mathbb{H}^1\left( \mathbb{P}_\mathbb{C}^3, \mathcal{N}_{\mathcal{F}_0}^{*\bullet} \right)=0, \mathbb{H}^1\left( \mathbb{P}_\mathbb{C}^3, \Theta_{\mathcal{F}_0}^\bullet \right)=0$ and $\mathbb{H}^1\left( \mathbb{P}_\mathbb{C}^3, \tilde{\mathcal{N}}_{\mathcal{F}_0}^{*\bullet} \right)=0$. We will also show in Part VIII that the associated Poisson structure $\Lambda_0$ on $\mathbb{P}_\mathbb{C}^3$ is rigid by showing that $\mathbb{H}^1\left( \mathbb{P}_\mathbb{C}^3, \Theta_{\mathbb{P}_\mathbb{C}^3}^\bullet \right)=0$, where $\Theta_{\mathbb{P}_\mathbb{C}^3}^\bullet$ is the degree-shifted truncated Lichnerowicz-Poisson complex $\Theta_{\mathbb{P}_\mathbb{C}^3} \xrightarrow{\left[ \Lambda_0, -\right]} \bigwedge^2 \Theta_{\mathbb{P}_\mathbb{C}^3} \xrightarrow{\left[ \Lambda_0 , - \right]} \bigwedge^3 \Theta_{\mathbb{P}_\mathbb{C}^3}\to 0$. This is the only example of rigid projective (non-affine, K\"ahler) Poisson variety that the author knows so far. We  will give another example of rigid foliated deformations and unfoldings on foliated rational ruled surfaces. We also show that Jouanolou foliation on $\mathbb{P}_\mathbb{C}^2$, and linear pullback of Jouanolou foliation on $\mathbb{P}_\mathbb{C}^2$ to $\mathbb{P}_\mathbb{C}^3$ are rigid in unfoldings.

\section{Projective plane $\mathbb{P}_\mathbb{C}^2$}\

We study foliated deformations and unfoldings of $\mathbb{P}_\mathbb{C}^2$. Let $\mathcal{F}_0$ be a foliation on $\mathbb{P}_\mathbb{C}^2$ of dimension $1$. Then  since $\mathbb{P}_\mathbb{C}^2$ is a complex surface, foliated deformations of $\left( \mathbb{P}_\mathbb{C}^2, \Theta_{\mathcal{F}_0} \right)$ are equivalent to foliated deformations of $\left( \mathbb{P}_\mathbb{C}^2, \mathcal{N}_{\mathcal{F}_0}^* \right)$, i.e. $\mathbb{H}^i\left( \mathbb{P}_\mathbb{C}^2, \Theta_{\mathcal{F}_0}^\bullet \right)\cong \mathbb{H}^i\left( \mathbb{P}_\mathbb{C}^2, \mathcal{N}_{\mathcal{F}_0}^{*\bullet} \right)$ for $i\geq 0$ (see Part II).

\subsection{Foliated deformations of $\left(\mathbb{P}_\mathbb{C}^2, \Theta_{\mathcal{F}_0} \right)$ in terms of tangent sheaves}\

Let $\left( \mathbb{P}_\mathbb{C}^2, T_0 \right)$ be a foliated projective plane, where $T_0\in H^0\left( \mathbb{P}_\mathbb{C}^3, \Theta_{\mathbb{P}_\mathbb{C}^2}(m) \right),  m \geq -1$, so that we have the exact sequence defined by $T_0$
\begin{align*}
0\to \mathcal{O}_{\mathbb{P}_\mathbb{C}^2}(-m) \xrightarrow{T_0} \Theta_{\mathbb{P}_\mathbb{C}^2} \to \frac{\Theta_{\mathbb{P}_\mathbb{C}^2}}{\mathcal{O}_{\mathbb{P}_\mathbb{C}^2}(-m)} \to 0
\end{align*}
Then we have a complex of sheaves from $(\ref{ap3})$
\begin{align*}
\mathcal{E}_{\mathcal{O}_{\mathbb{P}_\mathbb{C}^2}(m)}^\bullet : \mathcal{E}_{\mathcal{O}_{\mathbb{P}_\mathbb{C}^2}(m)} \xrightarrow{\tilde{D}_0} \Theta_{\mathbb{P}_\mathbb{C}^2}(m) \to 0
\end{align*}
From the exact sequence $0\to \mathcal{O}_{\mathbb{P}_\mathbb{C}^2}\to \mathcal{E}_{ \mathcal{O}_{\mathbb{P}_\mathbb{C}^2}(m)}\to \Theta_{\mathbb{P}_\mathbb{C}^2}\to 0$, we have $H^i\left( \mathbb{P}_\mathbb{C}^2 , \mathcal{E}_{\mathcal{O}_{\mathbb{P}_\mathbb{C}^2}(m)} \right)=0$ for $i\geq 1$. On the other hand, from the Euler sequence $0\to \mathcal{O}_{\mathbb{P}_\mathbb{C}^2} \to \bigoplus^3 \mathcal{O}_{\mathbb{P}_\mathbb{C}^2}(1)\to \Theta_{\mathbb{P}_\mathbb{C}^2} \to 0$,  we have $0\to \mathcal{O}_{\mathbb{P}_\mathbb{C}^2}(m)\to \bigoplus^3 \mathcal{O}_{\mathbb{P}_\mathbb{C}^2}(m+1) \to \Theta_{\mathbb{P}_\mathbb{C}^2}(m) \to 0$, so that $H^i\left( \mathbb{P}_\mathbb{C}^2, \Theta_{\mathbb{P}_\mathbb{C}^2 } (m) \right)=0$ for $m\geq -1, i\geq 1$. This implies that we can compute $\mathbb{H}^i\left(  \mathbb{P}_\mathbb{C}^2,  \mathcal{E}_{ \mathcal{O}_{\mathbb{P}_\mathbb{C}^2}(m)}^\bullet \right)$ by the following complex
\begin{align*}
H^0\left( \mathbb{P}_\mathbb{C}^2, \mathcal{E}_{\mathcal{O}_{\mathbb{P}_\mathbb{C}^2}(m)}  \right) \to H^0\left( \mathbb{P}_\mathbb{C}^2, \Theta_{\mathbb{P}_\mathbb{C}^2}(m) \right) \to 0
\end{align*}
Then we have
\begin{align*}
\mathbb{H}^1\left( \mathbb{P}_\mathbb{C}^2, \mathcal{E}_{\mathcal{O}_{\mathbb{P}_\mathbb{C}^2}(m)}^\bullet \right) & \cong  \textnormal{coker}\left( H^0\left( \mathbb{P}_\mathbb{C}^2, \mathcal{E}_{\mathcal{O}_{\mathbb{P}_\mathbb{C}^2}(m)} \right) \xrightarrow{\tilde{D}_0} H^0\left( \mathbb{P}_\mathbb{C}^2, \Theta_{\mathbb{P}_\mathbb{C}^2}(m) \right) \right) \\
\mathbb{H}^2\left( \mathbb{P}_\mathbb{C}^2 , \mathcal{E}_{\mathcal{O}_{\mathbb{P}_\mathbb{C}^2}(m)}^\bullet  \right)&= 0
\end{align*}
In particular, $\left( \mathbb{P}_\mathbb{C}^2, T_0 \right)$ is unobstructed in foliated deformations.

\begin{example}\label{cp9}
We consider the Jouanolou foliation of degree $2$ on $\mathbb{P}_\mathbb{C}^2$ defined by the homogeneous vector field $T=\xi_1^2\frac{\partial}{\partial \xi_0} + \xi_2^2 \frac{\partial}{\partial \xi_1} + \xi_0^2 \frac{\partial}{\partial \xi_2}$ on the homogeneous coordinate $\left[ \xi_0:\xi_1:\xi_2\right]$, so that $\Theta_{\mathcal{F}_0}\cong \mathcal{O}_{\mathbb{P}_\mathbb{C}^2}(-1)$ $($see \textnormal{\cite{LS20} p.67}$)$. We keep the notations in \textnormal{Appendix \ref{p1}}. Then on $U_0$, we have 
\begin{align*}
T_0= \left( X_2^2 - X_1^3 \right)\frac{\partial}{\partial X_1} + \left( 1- X_1^2 X_2  \right) \frac{\partial}{\partial X_2}
\end{align*}
Let us compute $\mathbb{H}^i\left( \mathbb{P}_\mathbb{C}^2, \mathcal{E}_{\mathcal{O}_{\mathbb{P}_\mathbb{C}^2}(1)}^\bullet \right)$ by the complex 
\begin{align*}
 \mathcal{E}_{\mathcal{O}_{\mathbb{P}_\mathbb{C}^2}(1)}^\bullet: H^0\left( \mathbb{P}_\mathbb{C}^2, \mathcal{E}_{\mathcal{O}_{\mathbb{P}_\mathbb{C}^2}(1)}\right) \xrightarrow{\tilde{D}_0}  H^0\left( \mathbb{P}_\mathbb{C}^2,\Theta_{\mathbb{P}_\mathbb{C}^2}(1) \right) \to 0
\end{align*}
First we compute $\textnormal{ker}\left(  H^0\left( \mathbb{P}_\mathbb{C}^2, \mathcal{E}_{\mathcal{O}_{\mathbb{P}_\mathbb{C}^2}(1)} \right) \xrightarrow{\tilde{D}_0} H^0\left( \mathbb{P}_\mathbb{C}^2, \Theta_{\mathbb{P}_\mathbb{C}^2}(1) \right)  \right)$ on $U_0$. Then we recall $\textnormal{(\ref{p2})}$ for $m=1$ and $\textnormal{(\ref{p3})}$. Then the image of $\tilde{D}_0$ is given by $\left[T_1\frac{\partial}{\partial X_1}+ T_2\frac{\partial}{\partial X_2} , T_0 \right] + f T_0$ which is 
{\small{\begin{align}\label{p5}
&\left(  -b_3- b_8 X_1 + 2b_4X_2 -3 b_1X_1^2 + 2b_5 X_1X_2 +\left( a- b_2+ 2 b_6 \right) X_2^2 - \left(a+2b_2\right) X_1^3 - b_7 X_1X_2^2  - 2b_3 X_1^2 X_2\right) \frac{\partial}{\partial X_1}\\
&+\left(\left(a- b_6\right) - 2b_7 X_1- 3b_8 X_2  -b_4X_1^2 - b_5 X_2^2 - 2b_1 X_1X_2 - \left(a+ 2b_2\right) X_1^2 X_2 - 2b_3 X_1 X_2^2- b_7 X_2^3 \right) \frac{\partial}{\partial X_2} \notag
\end{align}}}
Then we see that $\dim_\mathbb{C} \mathbb{H}^0\left( \mathbb{P}_\mathbb{C}^2, \mathcal{E}_{\mathcal{O}_{\mathbb{P}_\mathbb{C}^2}(1)}^\bullet \right)=\dim_\mathbb{C} \textnormal{ker}\left(H^0\left( \mathbb{P}_\mathbb{C}^2, \mathcal{E}_{\mathcal{O}_{\mathbb{P}_\mathbb{C}^2}(1)} \right) \xrightarrow{\tilde{D}_0} H^0\left( \mathbb{P}_\mathbb{C}^2, \Theta_{\mathbb{P}_\mathbb{C}^2}(1) \right) \right)= 0$, and so
\begin{align*}
\dim_\mathbb{C} \mathbb{H}^1\left( \mathbb{P}_\mathbb{C}^2, \mathcal{E}_{\mathcal{O}_{\mathbb{P}_\mathbb{C}^2}(1)}^\bullet \right)= 15- 9 = 6
\end{align*}
By considering $\textnormal{(\ref{p5})}$ and $\textnormal{(\ref{p3})}$, we see that
\begin{align*}
\mathbb{H}^1\left( \mathbb{P}_\mathbb{C}^2 , \mathcal{E}_{\mathcal{O}_{\mathbb{P}_\mathbb{C}^2}(1)}^\bullet \right)= \textnormal{Span}_\mathbb{C} \left\langle \frac{\partial}{\partial X_1}, X_1 \frac{\partial}{\partial X_1}, X_2 \frac{\partial}{\partial X_1}, X_1^2 \frac{\partial}{\partial X_1}, X_1X_2\frac{\partial}{\partial X_1} , X_1 \frac{\partial}{\partial X_2} \right\rangle
\end{align*}
Then we see that $\left( \mathbb{P}_\mathbb{C}^2\times B, T, \pi, B \right)$ where $\pi:\mathbb{P}_\mathbb{C}^2\times  B\to B$ is the natural projection with the coordinates $\left( t_1,...,t_6 \right)\in B\subset \mathbb{C}^6$ for a sufficiently small neighborhood $B$ of $0$ in $\mathbb{C}^6$, and
\begin{align} \label{j32}
T= \left( X_1^2 - X_1^3 + t_1 + t_2 X_1 + t_3 X_2 + t_4 X_1^2 + t_5 X_1X_2 \right) \frac{\partial}{\partial X_1} + \left( 1- X_1^2 X_2 + t_6 X_1 \right) \frac{\partial}{\partial X_2}
\end{align}
defines a foliated analytic family of deformations $\pi^{-1}(0)=\left( \mathbb{P}_\mathbb{C}^2, T_0\right)$ in terms of tangent sheaves and the foliated Kodaira-Spencer map in terms of tangent sheaves is an isomorphism at $t=0$. 

\end{example}

\begin{example}
We consider the foliation of degree $1$ on $\mathbb{P}_\mathbb{C}^2$ defined by the homogeneous vector field $T= \xi_2 \frac{\partial}{\partial \xi_0} + \xi_0 \frac{\partial}{\partial \xi_1}$ on the homogeneous coordinate $\left[ \xi_0 : \xi_1 : \xi_2 \right]$, so that $\Theta_{\mathcal{F}_0}\cong \mathcal{O}_{\mathbb{P}_\mathbb{C}^2}$. We keep the notations in \textnormal{Appendix \ref{p1}}. Then on $U_0$, we have
\begin{align*}
T_0= \left( 1- X_1X_2 \right) \frac{\partial}{\partial X_1} -X_2^2 \frac{\partial}{\partial X_2}
\end{align*}
Let us compute $\mathbb{H}^i\left( \mathbb{P}_\mathbb{C}^2, \mathcal{E}_{\mathcal{O}_{\mathbb{P}_\mathbb{C}^2}}^\bullet \right)$ by the complex 
\begin{align*}
 \mathcal{E}_{\mathcal{O}_{\mathbb{P}_\mathbb{C}^2}}^\bullet: H^0\left( \mathbb{P}_\mathbb{C}^2, \mathcal{E}_{\mathcal{O}_{\mathbb{P}_\mathbb{C}^2}}\right) \xrightarrow{\tilde{D}_0}  H^0\left( \mathbb{P}_\mathbb{C}^2,\Theta_{\mathbb{P}_\mathbb{C}^2} \right) \to 0
\end{align*}
First we compute $\textnormal{ker}\left(  H^0\left( \mathbb{P}_\mathbb{C}^2, \mathcal{E}_{\mathcal{O}_{\mathbb{P}_\mathbb{C}^2}} \right) \xrightarrow{\tilde{D}_0} H^0\left( \mathbb{P}_\mathbb{C}^2, \Theta_{\mathbb{P}_\mathbb{C}^2} \right)  \right)$ on $U_0$. Then we recall $\textnormal{(\ref{p2})}$ for $m=0$ and $\textnormal{(\ref{ap1})}$. Then the image of $\tilde{D}_0$ is given by $\left[T_1\frac{\partial}{\partial X_1}+ T_2\frac{\partial}{\partial X_2} , T_0 \right]+f T_0$ which is 
{\small{\begin{align}\label{p6}
\left( a+ b_2 + \left(2b_7+ b_4\right) X_1 + \left(b_1+ b_8\right)X_2 + b_5 X_1^2+ \left(b_6- a\right)X_1X_2 \right)\frac{\partial}{\partial X_1} + \left( b_5 + \left(b_7+ 2b_4\right)X_2 + \left(b_6 - a\right) X_2^2 + b_5 X_1X_2 \right) \frac{\partial}{\partial X_2}
\end{align}}}
Then we see that $\dim_\mathbb{C} \mathbb{H}^0\left( \mathbb{P}_\mathbb{C}^2, \mathcal{E}_{\mathcal{O}_{\mathbb{P}_\mathbb{C}^2}}^\bullet \right)=\dim_\mathbb{C} \textnormal{ker}\left(H^0\left( \mathbb{P}_\mathbb{C}^2, \mathcal{E}_{\mathcal{O}_{\mathbb{P}_\mathbb{C}^2}} \right) \xrightarrow{\tilde{D}_0} H^0\left( \mathbb{P}_\mathbb{C}^2, \Theta_{\mathbb{P}_\mathbb{C}^2}  \right) \right)= 3$, and so
\begin{align*}
\dim_\mathbb{C} \mathbb{H}^1\left( \mathbb{P}_\mathbb{C}^2, \mathcal{E}_{\mathcal{O}_{\mathbb{P}_\mathbb{C}^2}}^\bullet \right)= 8- 6 = 2
\end{align*}
By considering $\textnormal{(\ref{p6})}$ and $\textnormal{(\ref{ap1})}$, we see that
\begin{align*} 
\mathbb{H}^1\left( \mathbb{P}_\mathbb{C}^2 , \mathcal{E}_{\mathcal{O}_{\mathbb{P}_\mathbb{C}^2}(1)}^\bullet \right)= \textnormal{Span}_\mathbb{C} \left\langle \frac{\partial}{\partial X_2} ,   X_1 \frac{\partial}{\partial X_2} \right\rangle
\end{align*}
Then we see that $\left( \mathbb{P}_\mathbb{C}^2\times B, T, \pi, B \right)$, where $\pi:\mathbb{P}_\mathbb{C}^2\times B \to B$ is the natural projection with the coordinates $\left( t_1,t_2 \right)\in B\subset \mathbb{C}^2$ for a sufficiently small neighborhood $B$ of $0\in \mathbb{C}^2$ and
\begin{align*}
T= \left( 1 - X_1X_2 + t_2 X_1 \right) \frac{\partial}{\partial X_1} + \left( - X_2^2 + t_1 \right) \frac{\partial}{\partial X_2}
\end{align*}
defines a foliated analytic family of deformations $\pi^{-1}(0)=\left( \mathbb{P}_\mathbb{C}^2, T_0\right)$ in terms of tangent sheaves and the foliated Kodaira-Spencer map in terms of tangent sheaves is an isomorphism at $t=0$. 
 
\end{example}

\subsection{Foliated deformations of $\left(\mathbb{P}_\mathbb{C}^2,  \mathcal{N}_{\mathcal{F}_0}^*  \right)$ in terms of cotangent sheaves}\

Let $\left(  \mathbb{P}_\mathbb{C}^2, w \right)$ be a foliated projective plane, where $w\in H^0\left( \mathbb{P}_\mathbb{C}^2, \Omega_{\mathbb{P}_\mathbb{C}^2}^1 (m) \right), m\geq 2$, so that we have the exact sequence defined by $w$
\begin{align*}
0\to \mathcal{O}_{\mathbb{P}_\mathbb{C}^2}(-m)\to  \Omega_{\mathbb{P}_\mathbb{C}^2}^1 \to \frac{\Omega_{\mathbb{P}_\mathbb{C}^2}^1}{\mathcal{O}_{\mathbb{P}_\mathbb{C}^2}(-m)} \to 0
\end{align*}
Then we have a complex of sheaves from $(\ref{cp2})$
\begin{align*}
\mathcal{E}_{\mathcal{O}_{\mathbb{P}_\mathbb{C}^2}(m)}^\bullet : \mathcal{E}_{\mathcal{O}_{\mathbb{P}_\mathbb{C}^2}(m)}\xrightarrow{\tilde{E}_0} \Omega_{\mathbb{P}_\mathbb{C}^2}^1(m)  \to 0
\end{align*}
Then from the exact sequence $0\to \mathcal{O}_{\mathbb{P}_\mathbb{C}^2}\to \mathcal{E}_{\mathcal{O}_{\mathbb{P}_\mathbb{C}^2}(m)} \to \Theta_{\mathbb{P}_\mathbb{C}^2} \to 0$, we see that $H^i\left( \mathbb{P}_\mathbb{C}^2, \mathcal{E}_{\mathcal{O}_{\mathbb{P}_\mathbb{C}^2}(m)} \right)=0$ for $i\geq 1$, and by Bott formula (see \cite{OSS88} p.4) we have $H^i\left( \mathbb{P}_\mathbb{C}^2, \Omega_{\mathbb{P}_\mathbb{C}^2}^1(m) \right)=0$ for $m\geq 2, i\geq 1$. This implies that we can compute $\mathbb{H}^i\left( \mathbb{P}_\mathbb{C}^2, \mathcal{E}_{\mathcal{O}_{\mathbb{P}_\mathbb{C}^2}(m)}^\bullet \right)$ by the following complex
\begin{align*}
H^0\left( \mathbb{P}_\mathbb{C}^2, \mathcal{E}_{\mathcal{O}_{\mathbb{P}_\mathbb{C}^2}(m)}  \right) \xrightarrow{\tilde{E}_0} H^0\left( \mathbb{P}_\mathbb{C}^2,   \Omega_{\mathbb{P}_\mathbb{C}^2}^1(m) \right) \to 0
\end{align*}
Then we have
\begin{align*}
\mathbb{H}^1\left( \mathbb{P}_\mathbb{C}^2, \mathcal{E}_{\mathcal{O}_{\mathbb{P}_\mathbb{C}^2}(m)}^\bullet \right)&\cong \textnormal{coker}\left( H^0\left( \mathbb{P}_\mathbb{C}^2, \mathcal{E}_{\mathcal{O}_{\mathbb{P}_\mathbb{C}^2}(m)} \right) \to H^0\left(  \mathbb{P}_\mathbb{C}^2, \Omega_{\mathbb{P}_\mathbb{C}^2}^1(m) \right) \right) \\
\mathbb{H}^2\left( \mathbb{P}_\mathbb{C}^2, \mathcal{E}_{\mathcal{O}_{\mathbb{P}_\mathbb{C}^2}(m)}^\bullet \right)&= 0
\end{align*}

\begin{example}\label{cp10}
Let us consider the Jouanolou foliation of degree $2$ on $\mathbb{P}_\mathbb{C}^2$ in \textnormal{Example \ref{cp9}}. In this example, we describe foliated deformations of Jouanolou foliation on $\mathbb{P}_\mathbb{C}^2$ in terms of cotangent sheaves. Since foliated deformations in terms of tangent sheaf is equivalent to foliated deformations in terms of cotangent sheaves \textnormal{(see Part II)}, we should get the same result from \textnormal{Example \ref{cp9}}. The Jouanolou foliation of degree $2$ on $\mathbb{P}_\mathbb{\mathbb{C}}^2$ in terms of cotangent sheaves is defined by the homogenous holomorphic form $\omega= \left( \xi_0^2 \xi_1- \xi_2^3 \right) d\xi_0+ \left( \xi_1^2 \xi_2 - \xi_0^3 \right) d\xi_1+ \left( \xi_0 \xi_2^2- \xi_0^3\right) d\xi_2$ on the homogeneous coordinates $\left[ \xi_0: \xi_1: \xi_2 \right]$ on $\mathbb{P}_\mathbb{C}^2$, so that $\mathcal{N}_{\mathcal{F}_0}^* \cong \mathcal{O}_{\mathbb{P}_\mathbb{C}^2}(-4)$. We keep the notations in \textnormal{Appendix \ref{p1}}. Then on $U_0$ we have
\begin{align*}
w_0= \left( X_1^2 X_2 - 1 \right) dX_1 + \left( X_2^2 - X_1^3 \right) dX_2
\end{align*}
Let us compute $\mathbb{H}^i\left( \mathbb{P}_\mathbb{C}^2, \mathcal{E}_{\mathcal{O}_{\mathbb{P}_\mathbb{C}^2}(4)}^\bullet \right)$ by the complex
\begin{align*}
\mathcal{E}_{\mathcal{O}_{\mathbb{P}_\mathbb{C}^2}(4)}^\bullet :H^0\left( \mathbb{P}_\mathbb{C}^2, \mathcal{E}_{\mathcal{O}_{\mathbb{P}_\mathbb{C}^2}(4)}^1\right) \xrightarrow{\tilde{E}_0}  H^0\left( \mathbb{P}_\mathbb{C}^2 , \Omega_{\mathbb{P}_\mathbb{C}^2}(4) \right)\to 0
\end{align*}
First we compute $\textnormal{ker}\left( H^0\left( \mathbb{P}_\mathbb{C}^2, \mathcal{E}_{\mathcal{O}_{\mathbb{P}_\mathbb{C}^2}(4)}^1\right) \xrightarrow{\tilde{E}_0}  H^0\left( \mathbb{P}_\mathbb{C}^2 , \Omega_{\mathbb{P}_\mathbb{C}^2}(4) \right) \right)$ on $U_0$. Then we recall $\textnormal{(\ref{p2})}$ for $m=4$ and $\textnormal{(\ref{pa1})}$.  Then the image of $\tilde{E}_0$ is given by 
\begin{align}
&\mathcal{L}_{\left( b_1+ b_2 X_1+ b_3 X_2+ b_7 X_1^2+ b_8 X_1X_2\right) \frac{\partial}{\partial X_1} +\left( b_4+ b_5 X_1+ b_6 X_2+ b_7 X_1X_2+ b_8 X_2^2 \right) \frac{\partial}{\partial X_2}}\left(\left( X_1^2 X_2 - 1 \right) dX_1 + \left( X_2^2 - X_1^3 \right) dX_2 \right)\notag\\
&+ \left( a - 4 b_7 X_1 - 4 b_8 X_2 \right) \left(\left( X_1^2 X_2 - 1 \right) dX_1 + \left( X_2^2 - X_1^3 \right) dX_2 \right)\notag\\
&=\left( \left(-a- b_2\right) + 2 b_7 X_1 + 3 b_8 X_2 +b_4X_1^2 +b_5X_2^2  + 2b_1X_1X_2 +2 b_3 X_1X_2^2 + b_7 X_2^3 + \left(a+ 3b_2 + b_6\right)X_1^2 X_2 \right) d X_1 \label{pa2}\\
&+\left( - b_3 - b_8 X_1 + 2b_4X_2 -3b_1X_1^2 +2 b_5 X_1X_2 + \left(a+ 3b_6\right)X_2^2 - \left(a+ b_6 + 3b_2\right) X_1^3 - 2b_3 X_1^2 X_2 - b_7 X_1X_2^2 \right) dX_2 \notag
\end{align}
Then we see that $\dim_\mathbb{C} \mathbb{H}^0\left( \mathbb{P}_\mathbb{C}^2, \mathcal{E}_{\mathcal{O}_{\mathbb{F}_2}(4)}^\bullet \right)=\textnormal{ker}\left( H^0\left( \mathbb{P}_\mathbb{C}^2, \mathcal{E}_{\mathcal{O}_{\mathbb{P}_\mathbb{C}^2}(4)}^1\right) \xrightarrow{\tilde{E}_0}  H^0\left( \mathbb{P}_\mathbb{C}^2 , \Omega_{\mathbb{P}_\mathbb{C}^2}(4) \right) \right)=0$, and so 
\begin{align*}
\dim_\mathbb{C} \mathbb{H}^1\left( \mathbb{P}_\mathbb{C}^2 , \mathcal{E}_{\mathcal{O}_{\mathbb{P}_\mathbb{C}^2}(4)}^\bullet \right)= 15- 9 =6
\end{align*}
By considering $(\ref{pa2})$ and $\textnormal{(\ref{pa1})}$, we see that
\begin{align*}
\mathbb{H}^1\left( \mathbb{P}_\mathbb{C}^2, \mathcal{E}_{\mathcal{O}_{\mathbb{P}_\mathbb{C}^2}(4)}^\bullet  \right) \cong \textnormal{Span}_\mathbb{C}\left\langle dX_2, X_1 dX_2, X_2 dX_2, X_1^2 dX_2 ,  X_1X_2 dX_2, X_1 dX_1\right\rangle
\end{align*}
Then we see that $\left( \mathbb{P}_\mathbb{C}^2\times B, T, \pi, B \right)$ where $\pi:\mathbb{P}_\mathbb{C}^2\times  B\to B$ is the natural projection with the coordinates $\left( t_1,...,t_6 \right)\in B\subset \mathbb{C}^6$ for a sufficiently small neighborhood $B$ of $0$ in $\mathbb{C}^6$, and
\begin{align} \label{j5}
w= \left( X_1^2 X_2 - 1+ t_6 X_1 \right) dX_1 + \left(  X_2^2 - X_1^3 +t_1 + t_2 X_1+ t_3 X_2+ t_4 X_1^2 + t_5 X_1X_2\right) dX_2
\end{align}
defines a foliated analytic family of deformations $\pi^{-1}(0)=\left( \mathbb{P}_\mathbb{C}^2, w_0 \right)$ in terms of cotangent sheaves and the foliated Kodaira-Spencer map in terms of cotangent sheaves is an isomorphism at $t=0$. 

\end{example}

\subsection{Unfoldings of $\left( \mathbb{P}_\mathbb{C}^2, \mathcal{N}_{\mathcal{F}_0}^* \right)$}\

Let $\left(  \mathbb{P}_\mathbb{C}^2, w \right)$ be a foliated projective plane, where $w\in H^0\left( \mathbb{P}_\mathbb{C}^2, \Omega_{\mathbb{P}_\mathbb{C}^2}^1 (m) \right), m\geq 2$, so that we have the exact sequence defined by $w$
\begin{align*}
0\to \mathcal{O}_{\mathbb{P}_\mathbb{C}^2}(-m)\to  \Omega_{\mathbb{P}_\mathbb{C}^2}^1 \to \frac{\Omega_{\mathbb{P}_\mathbb{C}^2}^1}{\mathcal{O}_{\mathbb{P}_\mathbb{C}^2}(-m)} \to 0
\end{align*}
Then we have a complex of sheaves from $(\ref{cp3})$
\begin{align}\label{cp5}
\tilde{\mathcal{E}}_{\mathcal{O}_{\mathbb{P}_\mathbb{C}^2}(m)}^\bullet : \mathcal{E}_{\mathcal{O}_{\mathbb{P}_\mathbb{C}^2}(m)}\xrightarrow{\tilde{G}_0} \Omega_{\mathbb{P}_\mathbb{C}^2}^1(m) \bigoplus \mathcal{O}_{\mathbb{P}_\mathbb{C}^2}(m)  \xrightarrow{\tilde{G}_1} \tilde{\mathcal{S}}^1(m) \to 0
\end{align}
We note that $\tilde{\mathcal{S}}^1 \subset \mathcal{O}_{\mathbb{P}_\mathbb{C}^2}(m)\otimes \bigwedge^{1+1}\Omega_{\mathbb{P}_\mathbb{C}^2}^1\cong \mathcal{O}_{\mathbb{P}_\mathbb{C}^2}(m-3)$. Since $\tilde{\mathcal{S}}^1$ is reflexive of rank $1$, we see that $\tilde{\mathcal{S}}^1\cong \mathcal{O}_{\mathbb{P}_\mathbb{C}^2}(m-3)$. Then $(\ref{cp5})$ is equivalent to
\begin{align*}
\tilde{\mathcal{E}}_{\mathcal{O}_{\mathbb{P}_\mathbb{C}^2}(m)}^\bullet : \mathcal{E}_{\mathcal{O}_{\mathbb{P}_\mathbb{C}^2}(m)}\xrightarrow{\tilde{G}_0} \Omega_{\mathbb{P}_\mathbb{C}^2}^1(m) \bigoplus \mathcal{O}_{\mathbb{P}_\mathbb{C}^2}(m)  \xrightarrow{\tilde{G}_1} \mathcal{O}_{\mathbb{P}_\mathbb{C}^2}(2m-3) \to 0
\end{align*}
Then we can compute $\mathbb{H}^i\left( \mathbb{P}_\mathbb{C}^2, \tilde{\mathcal{E}}_{\mathcal{O}_{\mathbb{P}_\mathbb{C}^2}(m)}^\bullet \right)$ by the following complex
\begin{align*}
H^0\left( \mathbb{P}_\mathbb{C}^2, \mathcal{E}_{\mathcal{O}_{\mathbb{P}_\mathbb{C}^2}(m)} \right) \xrightarrow{\tilde{G}_0} H^0\left( \mathbb{P}_\mathbb{C}^2, \Omega_{\mathbb{P}_\mathbb{C}^2}^1(m) \right) \bigoplus H^0\left( \mathbb{P}_\mathbb{C}^2, \mathcal{O}_{\mathbb{P}_\mathbb{C}^2}(m) \right) \xrightarrow{\tilde{G}_1} H^0\left(  \mathbb{P}_\mathbb{C}^2, \mathcal{O}_{\mathbb{P}_\mathbb{C}^2}(2m-3 )\right) \to 0
\end{align*}
Then we have
\begin{align*}
\mathbb{H}^i\left( \mathbb{P}_\mathbb{C}^2, \tilde{\mathcal{E}}_{\mathcal{O}_{\mathbb{P}_\mathbb{C}^2}(m)}^\bullet  \right) &\cong \frac{\textnormal{ker}\left( H^0\left( \mathbb{P}_\mathbb{C}^2, \Omega_{\mathbb{P}_\mathbb{C}^2}^1(m) \right) \bigoplus H^0\left( \mathbb{P}_\mathbb{C}^2, \mathcal{O}_{\mathbb{P}_\mathbb{C}^2}(m) \right) \xrightarrow{\tilde{G}_1} H^0\left(  \mathbb{P}_\mathbb{C}^2, \mathcal{O}_{\mathbb{P}_\mathbb{C}^2}(2m-3 )\right)\right)}{\textnormal{im}\left( H^0\left( \mathbb{P}_\mathbb{C}^2, \mathcal{E}_{\mathcal{O}_{\mathbb{P}_\mathbb{C}^2}(m)} \right) \xrightarrow{\tilde{G}_0} H^0\left( \mathbb{P}_\mathbb{C}^2, \Omega_{\mathbb{P}_\mathbb{C}^2}^1(m) \right) \bigoplus H^0\left( \mathbb{P}_\mathbb{C}^2, \mathcal{O}_{\mathbb{P}_\mathbb{C}^2}(m) \right)  \right)} \\
\mathbb{H}^i\left( \mathbb{P}_\mathbb{C}^2,  \tilde{\mathcal{E}}_{\mathcal{O}_{\mathbb{P}_\mathbb{C}^2}(m)}^\bullet  \right) & \cong \textnormal{coker}\left( H^0\left( \mathbb{P}_\mathbb{C}^2, \Omega_{\mathbb{P}_\mathbb{C}^2}^1(m) \right) \bigoplus H^0\left( \mathbb{P}_\mathbb{C}^2, \mathcal{O}_{\mathbb{P}_\mathbb{C}^2}(m) \right) \xrightarrow{\tilde{G}_1} H^0\left(  \mathbb{P}_\mathbb{C}^2, \mathcal{O}_{\mathbb{P}_\mathbb{C}^2}(2m-3 )\right) \right)
\end{align*}

\begin{example}
We consider the Jouanolou foliation of degree $2$ on $\mathbb{P}_\mathbb{C}^2$ defined by $w\in H^0\left(  \mathbb{P}_\mathbb{C}^2 , \Omega_{\mathbb{P}_\mathbb{C}^2}^1(4) \right)$ in \textnormal{Example \ref{cp10}}. We keep the notations in \textnormal{Appendix \ref{p1}} and \textnormal{Example \ref{cp10}}. We have on $U_0$, 
\begin{align*}
\omega_0 = \left( X_1^2 X_2 -1\right)dX_1 + \left( X_2^2 - X_1^3 \right) dX_2
\end{align*}
Then we have the unfolded dual leaf complex
\begin{align*}
\tilde{\mathcal{E}}_{\mathcal{O}_{\mathbb{P}_\mathbb{C}^2}(4)}^\bullet:  \mathcal{E}_{\mathcal{O}_{\mathbb{P}_\mathbb{C}^2}(4)} \to \Omega_{\mathbb{P}_\mathbb{C}^2}^1 (4) \bigoplus \mathcal{O}_{\mathbb{P}_\mathbb{C}^2}(4) \to \mathcal{O}_{\mathbb{P}_\mathbb{C}^2}(5)\to 0
\end{align*}
We show that $\mathbb{H}^i \left( \mathbb{P}_\mathbb{C}^2, \tilde{\mathcal{E}}_{\mathcal{O}_{\mathbb{P}_\mathbb{C}^2}(4)}^\bullet\right)=0$ for $i=0,1,2$. Hence $\left( \mathbb{P}_\mathbb{C}^2, w \right)$ is rigid in unfoldings by \textnormal{Theorem \ref{r1}}.

Let us compute $\mathbb{H}^i\left( \mathbb{P}_\mathbb{C}^2, \tilde{\mathcal{E}}_{\mathcal{O}_{\mathbb{P}_\mathbb{C}^2}(4)}^\bullet \right)$ by the complex
\begin{align*}
\tilde{\mathcal{E}}_{\mathcal{O}_{\mathbb{P}_\mathbb{C}^2}(4)}^\bullet: H^0\left( \mathbb{P}_\mathbb{C}^2,  \mathcal{E}_{\mathcal{O}_{\mathbb{P}_\mathbb{C}^2}(4)} \right) \xrightarrow{\tilde{G}_0} H^0\left( \mathbb{P}_\mathbb{C}^2, \Omega_{\mathbb{P}_\mathbb{C}^2}^1 (4) \right) \bigoplus H^0\left( \mathbb{P}_\mathbb{C}^2, \mathcal{O}_{\mathbb{P}_\mathbb{C}^2}(4) \right) \xrightarrow{\tilde{G}_1} H^0\left( \mathbb{P}_\mathbb{C}^2 , \mathcal{O}_{\mathbb{P}_\mathbb{C}^2}(5) \right) \to 0
\end{align*}
From $(\ref{pa2})$ we see that
\begin{align*}
\mathbb{H}^0\left( \mathbb{P}_\mathbb{C}^2, \tilde{\mathcal{E}}_{\mathcal{O}_{\mathbb{P}_\mathbb{C}^2}(4)}^\bullet \right)\cong \textnormal{ker}\left( H^0\left( \mathbb{P}_\mathbb{C}^2, \mathcal{E}_{\mathcal{O}_{\mathbb{P}_\mathbb{C}^2}(4)} \right)\to H^0\left( \mathbb{P}_\mathbb{C}^2, \Omega_{\mathbb{P}_\mathbb{C}^2}^1 (4) \right)\bigoplus H^0\left( \mathbb{P}_\mathbb{C}^2, \mathcal{O}_{\mathbb{P}_\mathbb{C}^2}(4) \right) \right)=0
\end{align*}
On the other hand, let us compute $ H^0\left( \mathbb{P}_\mathbb{C}^2, \Omega_{\mathbb{P}_\mathbb{C}^2}^1 (4) \right)\bigoplus H^0\left( \mathbb{P}_\mathbb{C}^2, \mathcal{O}_{\mathbb{P}_\mathbb{C}^2}(4) \right) \xrightarrow{\tilde{G}_1} H^0\left( \mathbb{P}_\mathbb{C}^2, \mathcal{O}_{\mathbb{P}_\mathbb{C}^2}(5)\right)$. First we set $f\in H^0\left( \mathbb{P}_\mathbb{C}^2, \mathcal{O}_{\mathbb{P}_\mathbb{C}^2}(4) \right)\cong \mathbb{C}^{15}$ by
\begin{align*}
f=&c_1+ c_2 X_1 + c_3 X_2+ c_4 X_1^2+ c_5 X_1X_2+ b_6 X_2^2 + c_7 X_1^3 + c_8 X_1^2 X_2 + c_{9} X_1 X_2^2+ c_{10} X_2^3\\
 &+ c_{11}X_1^4+ c_{12} X_1^3 X_2 + c_{13} X_1^2 X_2^2+ c_{14} X_1 X_2^3+ c_{15} X_2^4,\,\,\,\,\,\,\,\,c_1,...,c_{15}\in \mathbb{C}
\end{align*}
Then we have
\begin{align*}
df= \left( c_2 + 2 c_4X_1+ c_5 X_2+ 3c_7 X_1^2 + 2c_8 X_1X_2+ c_9 X_2^2 + 4 c_{11} X_1^3 + 3 c_{12} X_1^2 X_2+ 2 c_{13}X_1X_2^2 + c_{14}X_2^3\right) dX_1\\
 +\left( c_3 + c_5 X_1 + 2 b_6 X_2 + c_8 X_1^2+ 2 c_9 X_1X_2+ 3 c_{10}X_2^2 + c_{12} X_1^3 + 2 c_{13}X_1^2 X_2+ 3 c_{14}X_1 X_2^2 + 4 c_{15}X_2^3\right) dX_2
\end{align*}
From $\textnormal{(\ref{pa1})}$ the image of $\tilde{G}_1$ is given by
\begin{align*}
f d\omega_0 + \omega_0 \wedge df - \omega_0 \wedge b
\end{align*}
First we compute $fd\omega_0+ \omega_0 \wedge df$. We note that $d\omega_0= - 4X_1^2 dX_1 \wedge dX_2$. Then the coefficient of $dX_1 \wedge dX_2$ of $fd\omega_0+ \omega_0 \wedge df$ is given by
{\small{\begin{align*}
&-4X_1^2\left( c_1+ c_2 X_1 + c_3 X_2+ c_4 X_1^2+ c_5 X_1X_2+ b_6 X_2^2 + c_7 X_1^3 + c_8 X_1^2 X_2 + c_{9} X_1 X_2^2+ c_{10} X_2^3 \right) \\
&- 4X_1^2\left( c_{11}X_1^4+ c_{12} X_1^3 X_2 + c_{13} X_1^2 X_2^2+ c_{14} X_1 X_2^3+ c_{15} X_2^4 \right) \\
&+ \left( X_1^2 X_2 -1 \right) \left( c_3 + c_5 X_1 + 2 b_6 X_2 + c_8 X_1^2+ 2 c_9 X_1X_2+ 3 c_{10}X_2^2 + c_{12} X_1^3 + 2 c_{13}X_1^2 X_2+ 3 c_{14}X_1 X_2^2 + 4 c_{15}X_2^3 \right)\\
&+\left( X_1^3 - X_2^2 \right) \left( c_2 + 2 c_4X_1+ c_5 X_2+ 3c_7 X_1^2 + 2c_8 X_1X_2+ c_9 X_2^2 + 4 c_{11} X_1^3 + 3 c_{12} X_1^2 X_2+ 2 c_{13}X_1X_2^2 + c_{14}X_2^3 \right)\\
&=-c_3 - c_5 X_1 - 2 c_6 X_2 - \left(4c_1+ c_8\right) X_1^2   - 2 c_9 X_1X_2 - \left(c_2+ 3c_{10}\right) X_2^2 - \left( 3c_2+ c_{12}\right)X_1^3  - \left( 3 c_3 + 2 c_{13}\right)X_1^2 X_2 - \left(2c_4+ 3 c_{14}\right) X_1 X_2^2 \\
&- \left( c_5 + 4 c_{15}\right)X_2^3 - 2 c_4 X_1^4 - 2 c_5 X_1^3 X_2 - \left( 2b_6+ 3c_7\right) X_1^2 X_2^2 - 2c_8 X_1X_2^3  - c_9 X_2^4 - c_7 X_1^5 - c_8  X_1^4 X_2  - \left(c_9+ 4 c_{11}\right) X_1^3 X_2^2 \\
&-(c_{10}+ 3 c_{12}) X_1^2 X_2^3  - 2 c_{13}X_1X_2^4 - c_{14}X_2^5 
\end{align*}}}
Next we compute $- \omega_0\wedge b$ from $\textnormal{(\ref{pa1})}$. Then the coefficient of $dX_1 \wedge dX_2$ of $-\omega_0\wedge b= b\wedge \omega_0$ is given by
{\small{\begin{align*}
&\left( X_2^2- X_1^3\right)\left( t_{10}+t_{11}X_2+t_{12}X_1+t_{13}X_2^2+t_{14}X_1X_2+t_{15}X_1^2-t_9X_1^2X_2-t_8X_1X_2^2-t_7X_2^3 \right) \\
&+\left( 1- X_1^2 X_2 \right)\left( t_1+t_2X_2+t_3X_1+t_4X_2^2+t_5X_1X_2+t_6X_1^2+t_7X_1X_2^2+t_8X_1^2X_2+t_9X_1^3 \right)\\
&=t_1 +t_3X_1+t_2X_2 +t_6X_1^2 +t_5X_1X_2 +( t_4 + t_{10})X_2^2 + (t_9 - t_{10})X_1^3   + (t_8- t_1) X_1^2 X_2 + (t_7+ t_{12})X_1 X_2^2 +t_{11}X_2^3 \\
&-t_{12}X_1^4 - (t_3+ t_{11})X_1^3 X_2 + (t_{15}- t_2) X_1^2 X_2^2 +t_{14}X_1X_2^3 +t_{13}X_2^4 -t_{15}X_1^5 - (t_6 + t_{14})X_1^4 X_2 - (t_5 + t_{13})X_1^3 X_2^2 \\
&-(t_4 + t_9)X_1^2 X_2^3 - t_8X_1X_2^4-t_7X_2^5
\end{align*}}}
Then the kernel of $\tilde{G}_1$ is given by $t_1=c_3, t_3=c_5, t_2=2c_6, t_6= 4c_1+ c_8, t_5= 2c_9, t_4+ t_{10}=c_2+ 3 c_{10},  t_9 - t_{10}= 3c_2+ c_{12},t_8- t_1= 3c_3+ 2 c_{13}, t_7+t_{12}= 2 c_4+ 3 c_{14}, t_{11}= c_5 + 4 c_{15},t_{12}= - 2c_4, t_3+ t_{11}= - 2c_5$, $t_{15}- t_2=2b_6+ 3 c_7, t_{14}= 2c_8, t_{13}=c_9, t_{15}= - c_7, t_6+ t_{14}= - c_8,  t_5+ t_{13}= -c_9 - 4 c_{11}, t_4+t_9=- c_{10}- 3 c_{12}, t_8= - 2c_{13}, t_7= - c_{14}$. Then we see that
\begin{align*}
\dim_{\mathbb{C}} \textnormal{ker}\left( H^0\left( \mathbb{P}_\mathbb{C}^2, \Omega_{\mathbb{P}_\mathbb{C}^2}^1(4) \right) \bigoplus H^0\left( \mathbb{P}_\mathbb{C}^2, \mathcal{O}_{\mathbb{P}_\mathbb{C}^2}(4) \right) \to H^0\left( \mathbb{P}_\mathbb{C}^2, \mathcal{O}_{\mathbb{P}_\mathbb{C}^2}(5) \right) \right)=9
\end{align*}
and so we have
\begin{align*}
\dim_\mathbb{C} \mathbb{H}^1\left( \mathbb{P}_\mathbb{C}^2, \tilde{\mathcal{E}}_{\mathcal{O}_{\mathbb{P}_\mathbb{C}^2}(4)}^\bullet \right)&= \dim_\mathbb{C} \frac{\textnormal{ker}\left(H^0\left( \mathbb{P}_\mathbb{C}^2, \Omega_{\mathbb{P}_\mathbb{C}^2}^1(4) \right) \bigoplus H^0\left( \mathbb{P}_\mathbb{C}^2, \mathcal{O}_{\mathbb{P}_\mathbb{C}^2}(4) \right)  \to H^0\left( \mathbb{P}_\mathbb{C}^2, \mathcal{O}_{\mathbb{P}_\mathbb{C}^2}(5) \right) \right)}{\textnormal{im}\left( H^0\left( \mathbb{P}_\mathbb{C}^2, \mathcal{E}_{\mathcal{O}_{\mathbb{P}_\mathbb{C}^2}(4)} \right)  \to H^0\left( \mathbb{P}_\mathbb{C}^2, \Omega_{\mathbb{P}_\mathbb{C}^2}^1(4) \right) \bigoplus H^0\left( \mathbb{P}_\mathbb{C}^2, \mathcal{O}_{\mathbb{P}_\mathbb{C}^2}(4) \right) \right)}= 9 - (9-0) =0
\end{align*}
\begin{align*}
\dim_\mathbb{C} \mathbb{H}^2\left( \mathbb{P}_\mathbb{C}^2, \tilde{\mathcal{E}}_{\mathcal{O}_{\mathbb{P}_\mathbb{C}^2}(4)}^\bullet \right) = \dim_\mathbb{C} \frac{H^0\left( \mathbb{P}_\mathbb{C}^2, \mathcal{O}_{\mathbb{P}_\mathbb{C}^2}(5) \right)}{\textnormal{im}\left(H^0\left( \mathbb{P}_\mathbb{C}^2, \Omega_{\mathbb{P}_\mathbb{C}^2}^1(4) \right) \bigoplus H^0\left( \mathbb{P}_\mathbb{C}^2, \mathcal{O}_{\mathbb{P}_\mathbb{C}^2}(4) \right)  \to H^0\left( \mathbb{P}_\mathbb{C}^2, \mathcal{O}_{\mathbb{P}_\mathbb{C}^2}(5) \right) \right) } = 21 -( 30-9)= 0
\end{align*}

\end{example}

\section{Rational ruled surfaces}

We study foliated deformations and unfoldings of rational ruled surfaces $\mathbb{F}_m=\mathbb{P}\left( \mathcal{O}_{\mathbb{P}_\mathbb{C}^2}(m) \bigoplus \mathcal{O}_{\mathbb{P}_\mathbb{C}^2} \right), m \geq 0$. Let $\mathcal{F}_0$ be a foliation on $\mathbb{F}_m$ of dimension $1$. Then  since $\mathbb{F}_m$ is a complex surface, foliated deformations of $\left( \mathbb{F}_m, \Theta_{\mathcal{F}_0} \right)$ in terms of tangent sheaves are equivalent to foliated deformations of $\left( \mathbb{F}_m, \mathcal{N}_{\mathcal{F}_0}^* \right) $ in terms of cotangent sheaves, i.e. $\mathbb{H}^i\left( \mathbb{F}_m, \Theta_{\mathcal{F}_0}^\bullet \right)\cong \mathbb{H}^i\left( \mathbb{F}_m , \mathcal{N}_{\mathcal{F}_0}^{*\bullet} \right)$ for $i\geq 0$ (see Part II). We will present foliated deformations of $\left( \mathbb{F}_m, \mathcal{F}_0 \right)$ in terms of tangent sheaves as well as in terms of cotangent sheaves.  We also present unfoldings of $\left( \mathbb{F}_m,  \mathcal{N}_{\mathcal{F}_0}^* \right)$.

\subsection{Construction of $\mathbb{F}_m$ and explicit description of the Picard group of $\mathbb{F}_m$} \label{rs1} \

We recall that $\mathbb{F}_m$ is constructed in the following way:  take two copies of $V_i\times \mathbb{P}_\mathbb{C}^1$ where $V_i=\mathbb{C}, i=1,2$, and write the coordinates as $\left(z, \xi \right)\in V_1 \times \mathbb{P}_\mathbb{C}^1$ and $\left(z', \xi' \right)\in V_2\times \mathbb{P}_\mathbb{C}^1$, where $z, z'$ are the coordinates of $\mathbb{C}$, respectively, and $\xi, \xi'$ are the inhomogeneous coordinates of $\mathbb{P}_\mathbb{C}^1$, respectively. Then $\mathbb{F}_m$ is constructed by patching $V_i\times \mathbb{P}_\mathbb{C}^1$ by the relation $z'=\frac{1}{z}$ and $\xi'=z^m \xi$. Let $\pi:\mathbb{F}_m\to \mathbb{P}_\mathbb{C}^1$ be the natural projection defined by $(z, \xi)\to z$ and $(z', \xi')\to z'$.  The Picard group of $\mathbb{F}_m$ is $\mathbb{Z} \bigoplus \mathbb{Z}$ which is generated by $\pi^*\mathcal{O}_{\mathbb{P}_\mathbb{C}^1}(1)$ and the line bundle $E$ associated with the divisor defined by $\xi=\xi'=0$. We will explicitly describe the transition functions of $\pi^*\mathcal{O}_{\mathbb{P}_\mathbb{C}^1}(1)$ and $E$ in terms of the affine covering $U_0, U_1, U_2, U_3$ of $\mathbb{F}_m$, where $U_0=\left\{ (z_0,w_0)\in \mathbb{C}^2 \right\}$, $U_1=\left\{\left( z_1, w_1\right)\in \mathbb{C}^2 \right\}$, $U_2=\left\{ \left( z_2, w_2 \right) \in \mathbb{C}^2 \right\}$, $U_3=\left\{ \left( z_3, w_3\right)\in \mathbb{C}^2 \right\}$ such that $U_0\bigcup U_2 = V_1 \times \mathbb{P}_\mathbb{C}^1$ and $U_1\bigcup U_3= V_2 \times \mathbb{P}_\mathbb{C}^1$. Then we have coordinate transformations as follows:
\begin{center}
$\begin{CD}
U_0: \left(z_0,w_0\right)@>>> U_1:\left(z_1, w_1\right)=\left(\frac{1}{z_0}, z_0^m w_0\right)\\
@VVV @VVV \\
U_2:\left(z_2, w_2\right)= \left(z_0, \frac{1}{w_0}\right) @>>> U_3:\left(z_3,w_3\right)=\left(\frac{1}{z_2}, z_2^{-m} w_2\right)=\left(z_1, \frac{1}{w_1}\right)
\end{CD}$
\end{center}
We note that $z_i$ is the base coordinate and $w_i$ is the fiber coordinate for $\pi$. Consider the divisor defined by the curve $f_1:=z_1=0$ on $U_1$ and $f_3:=z_3=0$ on $U_3$ and $f_0:=1$ on $U_0$ and $f_2:=1$ on $U_2$. Then the associated line bundle is described by the transition functions
\begin{align*}
f_{01}= \frac{1}{z_1}, \,\,\,\,\,f_{02}=1,\,\,\,\,\, f_{03}= \frac{1}{z_3},\,\,\,\,\, f_{12}= z_1,\,\,\,\,\, f_{13}=1,\,\,\,\,\, f_{23}=\frac{1}{z_3}
\end{align*}
where a local section $s_i \in \Gamma\left( U_i, \mathcal{O}_{\mathbb{F}_m} \right)$ of the line bundle is identified with a local section $s_j\in \Gamma\left(U_j, \mathcal{O}_{\mathbb{F}_m}  \right)$ if $s_i= f_{ij} s_j$ on $U_i\cap U_j$. On the other hand, let us consider the divisor defined by the curve $g_0:=w_0=0$ on $U_0$ and $g_1:=w_1=0$ on $U_1$ and $g_2:=1$ on $U_2$ and $g_3:=1$ on $U_3$. Then the associated line bundle is described by the transition functions
\begin{align*}
g_{01}=\frac{w_0}{w_1}= \frac{w_0}{z_0^m w_0}= z_0^{-m},\,\,\,\,\,g_{02}= w_0, \,\,\,\,\,g_{03}= w_0,\,\,\,\,\,g_{12}= w_1,\,\,\,\,\,g_{13}= w_1,\,\,\,\,\,g_{23}=1
\end{align*}
where a local section $s_i \in \Gamma\left( U_i, \mathcal{O}_{\mathbb{F}_m} \right)$ of the line bundle is identified with a local section $s_j\in \Gamma\left(U_j, \mathcal{O}_{\mathbb{F}_m}  \right)$ if $s_i= g_{ij} s_j$ on $U_i\cap U_j$. Let us consider the line bundle defined by the transition functions
\begin{align*}
h_{01}&=f_{01}^p g_{01}^q=z_1^{-p} z_0^{-mq}= z_0^{p-mq}\\
h_{02}&=f_{02}^p g_{02}^q=w_0^q\\
h_{03}&=f_{03}^p g_{02}^q=z_3^{-p}w_0^q=z_0^pw_0^q\\
h_{12}&=f_{12}^p g_{12}^q=z_1^pw_1^q=z_2^{-p} \left(\frac{z_2^m}{w_2}  \right)^q =z_2^{mq-p} w_2^{-q}\\
h_{13}&=f_{13}^pg_{13}^q=w_1^q\\
h_{23}&=f_{23}^pg_{23}^q=z_3^{-p}
\end{align*}
where a local section $s_i \in \Gamma\left( U_i, \mathcal{O}_{\mathbb{F}_m} \right)$ of the line bundle is identified with a local section $s_j\in \Gamma\left(U_j, \mathcal{O}_{\mathbb{F}_m}  \right)$ if $s_i= h_{ij} s_j$ on $U_i\cap U_j$. We will denote by $\mathcal{O}_{\mathbb{F}_m}(p,q)$ the line bundle defined by the transition functions $\left\{h_{ij} \right\},i,j=0,1,2,3$ in terms of the covering $\mathcal{U}=\left\{U_0,U_1,U_2, U_3\right\}$ of $\mathbb{F}_m$.

First we will study foliated deformations of $\mathbb{F}_m$ in terms of tangent sheaves. A foliation on $\mathbb{F}_m$ is given by a section $T\in H^0\left( \mathbb{F}_m, \Theta_{\mathbb{F}_m}(p,q) \right)$, where $\Theta_{\mathbb{F}_m}(p,q):= \Theta_{\mathbb{F}_m}\otimes \mathcal{O}_{\mathbb{F}_m}(p,q)$, so that we have the exact sequence
\begin{align*}
0\to \mathcal{O}_{\mathbb{F}_m}(-p,-q) \xrightarrow{T} \Theta_{\mathbb{F}_m} \to \frac{ \Theta_{\mathbb{F}_m} }{\mathcal{O}_{\mathbb{F}_m}(-p,-q)} \to 0
\end{align*}
We will study foliated deformations of $\left( \mathbb{F}_m, T \right)$ in terms of tangent sheaves. We note that the leaf complex is given by
\begin{align}\label{rs16}
\mathcal{E}_{\mathcal{O}_{\mathbb{F}_m}(p,q)}^\bullet : \mathcal{E}_{\mathcal{O}_{\mathbb{F}_m}(p,q)} \xrightarrow{\tilde{D}_0} \Theta_{\mathbb{F}_m}(p,q) \to 0
\end{align}
where $\mathcal{E}_{\mathcal{O}_{\mathbb{F}_m}(p,q)}$ is the Aityah extension associated to $\mathcal{O}_{\mathbb{F}_m}(p,q)$. In order to compute $\mathbb{H}^i\left( \mathbb{F}_m, \mathcal{E}_{\mathcal{O}_{\mathbb{F}_m}(p,q)}^\bullet \right)$ for $i=0,1,2$, we need to describe cohomology groups $H^i\left( \mathbb{F}_m,  \mathcal{E}_{\mathcal{O}_{\mathbb{F}_m}(p,q)}\right), i=0,1,2$ and $H^i\left( \mathbb{F}_m, \Theta_{\mathbb{F}_m}(p,q)\right), i=0,1$ as follows.

\subsection{The Atiyah extension $\mathcal{E}_{\mathcal{O}_{\mathbb{F}_m}(p,q)}$  and description of $H^0\left(\mathbb{F}_m, \mathcal{E}_{\mathcal{O}_{\mathbb{F}_m}(p,q)} \right)$}\label{rs3}\

We keep the notation in subsection \ref{rs1}. Let $\mathcal{E}_{\mathcal{O}_{\mathbb{F}_m}(p,q)}$ be the Atiyah extension associated the line bundle $\mathcal{O}_{\mathbb{F}_m}(p,q)$ on $\mathbb{F}_m$, so that we have the exact sequence
\begin{align}\label{rs4}
0 \to \mathcal{O}_{\mathbb{F}_m} \to \mathcal{E}_{\mathcal{O}_{\mathbb{F}_m}(p,q)} \to \Theta_{\mathbb{F}_m} \to 0
\end{align}
which induces the exact sequence
\begin{align}\label{rs2}
0 \to H^0\left(\mathbb{F}_m ,\mathcal{O}_{\mathbb{F}_m} \right) \to H^0\left( \mathbb{F}_m,  \mathcal{E}_{\mathcal{O}_{\mathbb{F}_m}(p,q)} \right) \to H^0\left( \mathbb{F}_m,  \Theta_{\mathbb{F}_m}\right) \to 0
\end{align}
We will explicitly describe $H^0\left( \mathbb{F}_m, \mathcal{E}_{\mathcal{O}_{\mathbb{F}_m}(p,q)}\right)$. We note that $\left(\theta, c_0\right) \in \mathcal{E}_{\mathcal{O}_{\mathbb{F}_m}(p,q)}|_{U_0} \cong \Theta_{\mathbb{F}_m}|_{U_0} \bigoplus \mathcal{O}_{\mathbb{F}_m}|_{U_0}$ is identified with $\left(\theta ,c_1\right) \in \mathcal{E}_{\mathcal{O}_{\mathbb{F}_m}(p,q)}|_{U_1}\cong \Theta_{\mathbb{F}_m}|_{U_1} \bigoplus \mathcal{O}_{\mathbb{F}_m}|_{U_1}$ on $U_0 \cap U_1$  if and only if $c_0 - c_1=\left[\theta, \log h_{10}\right]=\left[\theta, \log z_0^{mq-p}\right]$ on $U_0\cap U_1$. We also note that $\left( \theta, c_0 \right)\in \mathcal{E}_{\mathcal{O}_{\mathbb{F}_m}(p,q)}|_{U_0}\cong \Theta_{\mathbb{F}_m}|_{U_0} \bigoplus \mathcal{O}_{\mathbb{F}_m}|_{U_0}$ is identified with $\left( \theta, c_2 \right)\in \mathcal{E}_{\mathcal{O}_{\mathbb{F}_m}(p,q)}|_{U_2}$ on $U_0\cap U_2$ if and only if $c_0- c_2=\left[ a, \log h_{20}\right]= \left[ a, \log w_0^{-q}\right]$ on $U_0\cap U_1$.

\subsubsection{The case of $m=0$}\

We recall that $\dim_\mathbb{C} H^0\left( \mathbb{F}_0, \Theta_{\mathbb{F}_0} \right)= 6 $ and $H^0\left(\mathbb{F}_0, \Theta_{\mathbb{F}_0} \right)$ on $U_0$ is of the form  (see \cite{Kod05} p.73-p.75)
\begin{align*}
H^0\left(\mathbb{F}_0, \Theta_{\mathbb{F}_0} \right)=\left\{\left(az_0^2+bz_0+c\right)\frac{\partial}{\partial z_0}+ \left(d+ ew_0+ fw_0^2\right)\frac{\partial}{\partial w_0} | a,b,c,d,e,f\in \mathbb{C}\right\}
\end{align*}
Then from $(\ref{rs2})$ we see that $\dim_\mathbb{C} H^0\left( \mathbb{F}_0, \mathcal{E}_{\mathcal{O}_{\mathbb{F}_0}(p,q)} \right)=7$. We will explicitly describe $\left( \theta, c_0 \right)\in H^0\left( \mathbb{F}_0 , \mathcal{E}_{\mathcal{O}_{\mathbb{F}_m}(p,q)} \right)$ on $U_0$. We note that $\left( \theta, c_0 \right) \in \Gamma\left(U_0, \Theta_{\mathbb{F}_m} \right)\bigoplus \Gamma\left( U_0, \mathcal{O}_{\mathbb{F}_m} \right)\cong \Gamma\left( U_0, \mathcal{E}_{\mathcal{O}_{\mathbb{F}_0}(p,q)} \right)$ on $U_0$ is of the form
\begin{align*}
\left(\left(az_0^2+bz_0+c\right)\frac{\partial}{\partial z_0}+ \left(d+ ew_0+ fw_0^2 \right)\frac{\partial}{\partial w_0}, g(z_0,w_0)= \sum_{h,k\geq 0} g_{hk}z_0^h w_0^k   \right)
\end{align*}
where $a,b,c,d,e,f,g_{hk}\in \mathbb{C}$. Then $g(z_0,w_0)$ on $U_0$ is translated on $U_1$ into
\begin{align*}
&\sum_{h,k\geq 0} g_{hk}z_1^{-h}  w_1^k  + \left[  \left(az_0^2+bz_0+c\right)\frac{\partial}{\partial z_0}+ \left(d+ ew_0+ fw_0^2 \right)\frac{\partial}{\partial w_0}  , \log z_0^p \right]= \sum_{h,k\geq 0} g_{hk}z_1^{-h}  w_1^k + p\left( a \frac{1}{z_1}+ b + c z_1 \right)
\end{align*}
which should be holomorphic, so that $g_{10}= - pa$ and $g_{hk}=0$ for $h\geq 2$. Hence $g(z_0,w_0)$ is of the form $g(z_0,w_0)= - paz_0+ h(w_0)= -paz_0+ \sum_{k\geq 0} h_kw_0^k, h_k\in \mathbb{C}$. On the other hand, $g(z_0,w_0)$ on $U_0$ is translated on $U_2$ into
\begin{align*}
&-paz_2+ \sum_{k\geq 0} h_k w_2^{-k} -\left[ (az_0^2+bz_0+c)\frac{\partial}{\partial z_0}+ (d+ ew_0+ fw_0^2)\frac{\partial}{\partial w_0} , \log w_0^{-q}   \right]\\
&= -pa z_2+ \sum_{k\geq 0} h_kw_2^{-k}+ qw_0^{-1}(d+ew_0+ fw_0^2)= -paz_2+ \sum_{ k\geq 0} h_k w_2^{-k} + qdw_2 + qe + qf\frac{1}{w_2}
\end{align*}
which should be holomorphic, so that
\begin{align*}
h_1=- qf, \,\,\,\,\,\,\,\,\, h_k=0 \,\,\,\,\, (k\geq 2)
\end{align*}
Hence we have
\begin{align*}
g(z_0,w_0)= h - pa z_0 - q f w_0
\end{align*}
Then $\left( \theta, c_0 \right)\in H^0\left( \mathbb{F}_0, \mathcal{E}_{\mathcal{O}_{\mathbb{F}_0}(p,q)} \right)$ on $U_0$ is given by 
\begin{align}\label{rsf2}
H^0\left( \mathbb{F}_0, \mathcal{E}_{\mathcal{O}_{\mathbb{F}_0} (p,q)} \right)=\left\{\left(\left(az_0^2+bz_0+c\right)\frac{\partial}{\partial z_0}+ \left(d+ ew_0+ fw_0^2 \right)\frac{\partial}{\partial w_0}, h- paz_0- qf w_0  \right) | a,b,c, d,e,f, h\in \mathbb{C} \right\}
\end{align}

\subsubsection{The case of $m\geq 1$}\

We recall that $\dim_\mathbb{C} H^0\left( \mathbb{F}_m, \Theta_{\mathbb{F}_m} \right)=m+5$ and $H^0\left( \mathbb{F}_m, \Theta_{\mathbb{F}_m} \right)$ on $U_0$ is of the form (see \cite{Kod05} p.73-p.75)
\begin{align}\label{rs13}
H^0\left(\mathbb{F}_m , \Theta_{\mathbb{F}_m} \right)=\left\{ \left(a_0+ a_1 z_0 + a_2z_0^2    \right) \frac{\partial}{\partial z_0} + \left( \left(-ma_2z_0+d\right)w_0 + \left(c_0 + c_1z_0+ \cdots + c_m z_0^m \right)w_0^2 \right)\frac{\partial}{\partial w_0}  \right\}
\end{align}
where $d,a_0,...,a_2,c_0,..., c_m\in \mathbb{C}$. Then from $(\ref{rs2})$ we see that $\dim_\mathbb{C} H^0\left( \mathbb{F}_m, \mathcal{E}_{\mathcal{O}_{\mathbb{F}_m}(p,q)} \right)=m+6$. We explicitly describe $\left( \theta, C_0  \right)\in H^0\left( \mathbb{F}_m, \mathcal{E}_{\mathcal{O}_{\mathbb{F}_m}(p,q)} \right)$ on $U_0$. We note that $\left( \theta, C_0 \right)\in H^0\left( \mathbb{F}_m, \mathcal{E}_{\mathcal{O}_{\mathbb{F}_m}(p,q)} \right)$ on $U_0$ is of the form
\begin{align*}
\left(\left(a_0+ a_1 z_0 + a_2z_0^2    \right) \frac{\partial}{\partial z_0} + \left( \left(-ma_2z_0+d\right)w_0 + \left(c_0 + c_1z_0+ \cdots + c_m z_0^m \right)w_0^2 \right)\frac{\partial}{\partial w_0} , g(z_0,w_0)=\sum_{h,k\geq 0} g_{hk}z_0^h w_0^k \right)
\end{align*}
where $a_0,a_1,a_2,d, c_0,...,c_m,g_{hk}\in \mathbb{C}$. Then $g(z_0,w_0)$ on $U_0$ is translated on $U_2$ into
\begin{align*}
&\sum_{h,k\geq 0} g_{hk}z_0^h w_0^k - \left[ \left(a_0+ a_1 z_0 + a_2z_0^2    \right) \frac{\partial}{\partial z_0} + \left( \left(-ma_2z_0+d\right)w_0 + \left(c_0 + c_1z_0+ \cdots + c_m z_0^m \right)w_0^2 \right)\frac{\partial}{\partial w_0}, \log w_0^{-q} \right]\\
&=\sum_{h,k\geq 0} g_{hk} z_2^h w_2^{-k} + qw_0^{-1}\left( \left(-ma_2z_0+d\right)w_0 + \left(c_0 + c_1z_0+ \cdots + c_m z_0^m \right)w_0^2 \right)\\
&= \sum_{h,k\geq 0} g_{hk}  z_2^h w_2^{-k} + qw_2\left(\left( -ma_2z_2 + d \right)  \frac{1}{w_2}  + \left(c_0 + c_1z_2+ \cdots + c_m z_2^m \right)\frac{1}{w_2^2}  \right)\\
&= \sum_{h,k\geq 0} g_{hk} z_2^h w_2^{-k} + q\left(-ma_2z_2+d\right) + q\left(c_0+ c_1z_2+ \cdots + c_m z_2^m\right)\frac{1}{w_2}
\end{align*}
which should be holomorphic, so that $g(z_0,w_0)$ is of the form
\begin{align*}
g(z_0,w_0)= - q \left(c_0+ c_1z_0+ \cdots + c_mz_0^m\right) w_0 + \sum_{k\geq 0} h_k z_0^k,\,\,\,\,\,\,\,\,h_k \in \mathbb{C}
\end{align*}
On the other hand, $g(z_0,w_0)$ on $U_0$ is translated on $U_1$ into
\begin{align*}
& -q (c_0+ c_1 z_0 + \cdots + c_m z_0^m)w_0 + \sum_{k\geq 0} h_k z_0^k \\
& - \left[  \left(a_0+ a_1 z_0 + a_2z_0^2    \right) \frac{\partial}{\partial z_0} + \left( \left(-ma_2z_0+d\right)w_0 + \left(c_0 + c_1z_0+ \cdots + c_m z_0^m \right)w_0^2 \right)\frac{\partial}{\partial w_0} ,  \log  z_0^{mq-p}\right]\\
&=  -q \left(c_0+ c_1 \frac{1}{z_1} + \cdots + c_m \frac{1}{z_1^m}\right)z_1^m w_1 + \sum_{k\geq 0} h_k z_0^k + \left(p-mq \right)z_0^{-1}\left(a_0+ a_1z_0+ a_2z_0^2\right)\\
&=  -q \left(c_0+ c_1 \frac{1}{z_1} + \cdots + c_m \frac{1}{z_1^m}\right)z_1^m w_1  + \sum_{k\geq 0} h_k z_1^{-k}  + \left(p-mq \right)z_1\left( a_0+ a_1\frac{1}{z_1} + a_2 \frac{1}{z_1^2} \right)\\
&= -q \left(c_0+ c_1 \frac{1}{z_1} + \cdots + c_m \frac{1}{z_1^m}\right)z_1^m w_1  + \sum_{k\geq 0} h_k z_1^{-k}  + \left(p-mq \right)a_0 z_1+ \left(p-mq \right)a_1 + \left(p-mq \right) a_2 \frac{1}{z_1} 
\end{align*}
which should be holomorphic, so that $g(z_0,w_0)$ is of the form
\begin{align*}
g(z_0,w_0)= b - \left(p-mq\right)a_2z_0 - q \left(c_0+ c_1z_0+ \cdots + c_mz_0^m \right) w_0,\,\,\,\,\,\,\,\, b \in \mathbb{C}
\end{align*}
Hence $\left(\theta, C_0\right)\in H^0\left( \mathbb{F}_m , \mathcal{E}_{\mathcal{O}_{\mathbb{F}_m}(p,q)} \right)$ on $U_0$ is given by
\begin{align}\label{rsf1}
H^0\left( \mathbb{F}_m , \mathcal{E}_{\mathcal{O}_{\mathbb{F}_m}(p,q)} \right)= &\left\{  \left(a_0+ a_1 z_0 + a_2z_0^2    \right) \frac{\partial}{\partial z_0} + \left( \left(-ma_2z_0+d\right)w_0 + \left(c_0 + c_1z_0+ \cdots + c_m z_0^m \right)w_0^2 \right)\frac{\partial}{\partial w_0} , \right.  \\
 &\left.  b - \left(p-mq\right)a_2z_0 - q \left(c_0+ c_1z_0+ \cdots + c_mz_0^m \right) w_0 | a_0,a_1,a_2,d, c_0,...,c_m, b \in \mathbb{C} \right\} \notag
\end{align}

\subsection{Description of $H^1\left( \mathbb{F}_m, \mathcal{E}_{\mathcal{O}_{\mathbb{F}_m}(p,q)} \right)$} \label{rs34}\

We keep the notation in subsection \ref{rs1} and subsection \ref{rs3}. In this subsection we explicitly describe $H^1\left( \mathbb{F}_m, \mathcal{E}_{\mathcal{O}_{\mathbb{F}_m}(p,q)} \right)$. We note that from the exact sequence $(\ref{rs4})$, since $H^1\left(\mathbb{F}_m, \mathcal{O}_{\mathbb{F}_m}\right)=H^2\left( \mathbb{F}_m, \mathcal{O}_{\mathbb{F}_m} \right)=0$, and $H^1\left(\mathbb{F}_m, \Theta_{\mathbb{F}_m} \right)=0$ for $m=0,2$, and $\dim_\mathbb{C} H^1\left( \mathbb{F}_m, \Theta_{\mathbb{F}_m} \right)=m-1$ for $m> 1$ (see \cite{Kod05} p.312), we see that $\mathbb{H}^1\left( \mathbb{F}_m,  \mathcal{E}_{\mathcal{O}_{\mathbb{F}_m}(p,q)} \right)=0$ for $m=0,1$ and $\dim_\mathbb{C} \mathbb{H}^1\left( \mathbb{F}_m, \mathcal{E}_{\mathcal{O}_{\mathbb{F}_m}(p,q)}\right)=m-1$ for $m >1$. First we show that $H^1\left( U_0\bigcup U_2  , \mathcal{E}_{\mathcal{O}_{\mathbb{F}_m}}(p,q)\right)=0$ and $H^1\left(U_1 \bigcup U_3, \mathcal{E}_{\mathcal{O}_{\mathbb{F}_m}}(p,q) \right)=0$. We note that by letting $\mathscr{B}=\left\{ U_0, U_2\right\}$, we can compute $H^1\left(U_0\bigcup U_2 , \mathcal{E}_{\mathcal{O}_{\mathbb{F}_m}(p,q)} \right)$ by \v Cech resolution associated to the open covering $\mathscr{B}$, i.e.
\begin{align*}
H^1\left( U_0\bigcup U_2 , \mathcal{E}_{\mathcal{O}_{\mathbb{F}_m}(p,q)} \right) = H^1\left( \mathscr{B}, \mathcal{E}_{\mathcal{O}_{\mathbb{F}_m}(p,q)} \right)=\frac{Z^1\left(\mathscr{B}, \mathcal{E}_{\mathcal{O}_{\mathbb{F}_m}(p,q)}  \right)}{\delta C^0\left( \mathscr{B}, \mathcal{E}_{\mathcal{O}_{\mathbb{F}_m}(p,q)}  \right)}
\end{align*}
A $1$-cocycle $\left\{\left(\theta_{02}, r_{02}\right), \left(\theta_{20}, r_{20}\right)\right\}$ belonging to $Z^1\left(\mathscr{B}, \mathcal{E}_{\mathcal{O}_{\mathbb{F}_m}(p,q)} \right)$ consists of holomorphic vector fields $\theta_{02}$ and $\theta_{20}$ on $U_0\cap U_2$ with $\theta_{02}= - \theta_{20}$  and  holomorphic functions $r_{02}$ and $r_{20}$ on $U_0\cap U_2$ with $r_{02}+\left[\theta_{02},  \log w_0^q  \right] = - r_{20}$, while $ C^0\left(\mathscr{B}, \mathcal{E}_{\mathcal{O}_{\mathbb{F}_m}(p,q)} \right)$ is the set of all pairs $\left\{\left(\theta_0, r_0\right), \left(\theta_2,r_2\right)\right\}$, where $\theta_1$ and $\theta_2$ are holomorphic vector fields, and $h_0$ and $h_2$ are holomorphic functions, respectively, on $U_0$ and $U_2$. Hence in order to prove $H^1\left(\mathscr{B}, \mathcal{E}_{\mathcal{O}_{\mathbb{F}_m}(p,q)} \right)=0$, it suffices to show that any $1$-cocycle $\{(\theta_{02}, r_{02}),(\theta_{20}, r_{20})\}$ is represented as $\delta\left\{ (\theta_0, r_0) , (\theta_2, r_2) \right\}$, i.e.
\begin{align} \label{rs28}
\theta_{02}= \theta_2- \theta_0,\,\,\,\,\,\,\,\,\, r_{02}= r_2 + \left[\theta_{2}, \log h_{20}\right] - r_0 = r_2+ \left[\theta_{2}, \log w_2^{q}\right] - r_0 
\end{align}
The holomorphic vector field $\theta_{02}$ and the holomorphic function $r_{02}$ on $U_0\bigcap U_1=\mathbb{C}\times\mathbb{C}^*$ are written in the form
\begin{align*}
\theta_{02}&=  \left( \sum_{n=-\infty}^\infty u_n(z_0)w_0^n\right)\frac{\partial}{\partial z_0} + \left( \sum_{n=-\infty}^\infty v_n(z_0)w_0^n \right)\frac{\partial}{\partial w_0}\\
r_{02}&= \sum_{n=-\infty}^\infty  s_n(z_0) w_0^n
\end{align*}
where $u_n(z_0),v_n(z_0),s_n(z_0)$ are entire function of $z_0$. If we set
\begin{align}
\theta_0&= - \sum_{n=0 }^\infty u_n(z_0)w_0^n \frac{\partial}{\partial z_0} - \sum_{n=0}^\infty v_n(z_0)w_0^n \frac{\partial}{\partial w_0}\notag\\
\theta_2&=\sum _{n=1}^\infty u_{-n}(z_2)w_2^n\frac{\partial}{\partial z_2} - \sum_{n=1}^\infty v_{-n}(z_0)w_2^{n+2}\frac{\partial}{\partial w_2} \label{rs27}
\end{align}
then $\theta_{02}= \theta_2- \theta_0$ (see \cite{Kod05} p.311). On the other hand, we have to find $r_0$ and $r_2$ such that from $(\ref{rs28})$ and $(\ref{rs27})$
\begin{align}\label{rss28}
r_{02}= r_2+ q\left[\theta_2, \log w_2\right]- r_0= r_2 - q \sum_{n=1}^\infty v_{-n} (z_0)w_2^{n+1} - r_0
\end{align}
If we set
\begin{align*}
r_0&  = - \sum_{n=0}^\infty s_n(z_0) w_0^n \\
r_2 &= \sum_{n=1}^{\infty} s_{-n}(z_2)w_2^n + q\sum_{n=1}^\infty v_{-n}(z_2)w_2^{n+1}
\end{align*}
then $(\ref{rss28})$ holds. Hence $H^1\left( U_0 \bigcup U_2, \mathcal{E}_{\mathcal{O}_{\mathbb{F}_m}(p,q)} \right)=0$. Similarly we can show that $H^1\left( U_1 \bigcup U_3, \mathcal{E}_{\mathcal{O}_{\mathbb{F}_m}(p,q)} \right)=0$. Hence by letting $\mathfrak{B}=\left\{U_0\bigcup U_2, U_1\bigcup U_3 \right\}$, we can compute $H^1\left( \mathbb{F}_m,  \mathcal{E}_{\mathcal{O}_{\mathbb{F}_m}(p,q)} \right)$ by \v Cech resolution associated to the open covering $\mathfrak{B}$, i.e.
\begin{align}\label{rs5}
H^1\left( \mathbb{F}_m, \mathcal{E}_{\mathcal{O}_{\mathbb{F}_m}(p,q)} \right)=H^1\left( \mathfrak{B}, \mathcal{E}_{\mathcal{O}_{\mathbb{F}_m}(p,q)}\right) = \frac{Z^1\left( \mathfrak{B}, \mathcal{E}_{\mathcal{O}_{\mathbb{F}_m}(p,q)} \right)}{\delta C^0\left( \mathfrak{B}, \mathcal{E}_{\mathcal{O}_{\mathbb{F}_m}(p,q)} \right)}
\end{align}

Next we describe explicitly $H^0\left( U_0 \bigcup U_2, \mathcal{E}_{\mathcal{O}_{\mathbb{F}_m}(p,q)} \right)$ on $U_0$ and $H^0\left(U_1 \bigcup U_3, \mathcal{E}_{\mathcal{O}_{\mathbb{F}_m}(p,q)} \right)$ on $U_1$. Let  us describe $H^0\left( U_0 \bigcup U_2, \mathcal{E}_{\mathcal{O}_{\mathbb{F}_m}(p,q)} \right)$ on $U_0$. Let $\left( \theta_0, c_0 \right)\in H^0\left( U_0 \bigcup U_2,  \mathcal{E}_{\mathcal{O}_{\mathbb{F}_m}(p,q)} \right)$ on $U_0$. Then $\left( \theta_0, c_0\right)\in \Gamma\left(U_0, \Theta_{\mathbb{F}_m} \right)\bigoplus \Gamma\left( U_0, \mathcal{O}_{\mathbb{F}_m} \right)$ is of the form
\begin{align*}
\left( \theta_0, c_0 \right)= \left( v_0(z_0)\frac{\partial}{\partial z_0} + \left(\alpha_0^1(z_0) + \alpha_0^2 (z_0)w_0 +\alpha_0^3(z_0)w_0^2\right)\frac{\partial}{\partial w_0}, g(z_0,w_0)=  \sum_{h,k\geq 0}  g_{hk} z_0^h w_0^k \right)
\end{align*}
where $v_0(z_0), \alpha_0^1(z_0),\alpha_0^2(z_0), \alpha_0^3(z_0)$ are entire functions of $z_0$ and $g_{hk}\in \mathbb{C}$, and $\left(\theta_0, c_0 \right)$ on $U_0$ is transformed into $\left(\theta_0, c_2 \right)\in \Gamma\left( U_2, \Theta_{\mathbb{F}_m}\right) \bigoplus \Gamma\left(U_2, \mathcal{O}_{\mathbb{F}_m} \right)$ on $U_2$, where
\begin{align*}
c_2=&\sum_{h,k\geq 0} g_{hk} z_2^h w_2^{-k} - \left[ v_0(z_0)\frac{\partial}{\partial z_0} + \left(\alpha_0^1(z_0) + \alpha_0^2 (z_0)w_0 +\alpha_0^3(z_0)w_0^2 \right)\frac{\partial}{\partial w_0}, \log w_0^{-q} \right] \\
&= \sum_{h,k\geq 0} g_{hk} z_2^h w_2^{-k}  + \frac{q}{w_0} \left( \alpha_0^1(z_0) + \alpha_0^2 (z_0)w_0 +\alpha_0^3(z_0)w_0^2 \right)\\
&= \sum_{h,k\geq 0} g_{hk} z_2^h w_2^{-k}  + qw_2 \left( \alpha_0^1\left(z_2 \right) + \alpha_0^2 \left(z_2 \right)\frac{1}{w_2} +\alpha_0^3(z_2)\frac{1}{w_2^2} \right)\\
&= \sum_{h,k\geq 0} g_{hk} z_2^h w_2^{-k}  + \left(q \alpha_0^1\left(z_2 \right)w_2 + q \alpha_0^2 \left(z_2 \right) + q\alpha_0^3(z_2)\frac{1}{w_2} \right)
\end{align*}
which should be holomorphic, so that we have $g(z_0,w_0)$ is of the form $g(z_0, w_0)= s_0(z_0)- q \alpha_0^3(z_0)w_0$, where $s_0(z_0)$ is an entire function of $z_0$. Hence $\left( \theta_0, c_0\right)\in H^0\left(U_0\bigcup U_2 , \mathcal{E}_{\mathcal{O}_{\mathbb{F}_m}(p,q)} \right)$ on $U_0$ is of the form
\begin{align}\label{rs29}
\left( \theta_0, c_0 \right)=\left( v_0(z_0)\frac{\partial}{\partial z_0} + \left(\alpha_0^1(z_0) + \alpha_0^2 (z_0)w_0 +\alpha_0^3(z_0)w_0^2 \right)\frac{\partial}{\partial w_0},s_0(z_0)- q\alpha_0^3(z_0)w_0   \right)
\end{align}
where $v_0(z_0), \alpha_0^1(z_0),\alpha_0^2(z_0),\alpha_0^3(z_0), s_0(z_0)$ are entire functions of $z_1$. Similarly $\left( \theta_1, c_1 \right)\in H^0\left( U_1 \bigcup U_3 , \mathcal{E}_{\mathcal{O}_{\mathbb{F}_m}(p,q)} \right)$ on $U_1$ is of the form
\begin{align}\label{rs30}
\left( \theta_1, c_1 \right)=\left( v_1(z_1)\frac{\partial}{\partial z_1} + \left(\alpha_1^1(z_1) + \alpha_1^2 (z_1)w_1 +\alpha_1^3(z_1)w_1^2 \right)\frac{\partial}{\partial w_1},s_1(z_1)- q\alpha_1^3(z_1)w_1   \right)
\end{align}
where $v_1(z_1), \alpha_1^1(z_1),\alpha_1^2(z_1),\alpha_1^3(z_1), s_1(z_1)$ are entire functions of $z_1$.

Now let we describe $H^1\left(\mathbb{F}_m,  \mathcal{E}_{\mathcal{O}_{\mathbb{F}_m}(p,q)} \right)$ in terms of the \v Cech resolution associated to $\mathfrak{B}$ from $(\ref{rs5})$. We represent a $1$-cocycle $\left\{ \left( \theta_{01}, s_{01}\right), \left( \theta_{10}, s_{10}\right) \right\}$ with $\theta_{01}= - \theta_{10}$ and $s_{01} - \left[ \theta_{01},  \log h_{10} \right] = - s_{10}$  on $\left( U_0 \bigcup U_2 \right)\cap \left( U_1 \bigcup U_3 \right)$. Then the $1$-cocycle belongs to $\delta C^0\left( \mathfrak{B}, \mathcal{E}_{\mathcal{O}_{\mathbb{F}_m}(p,q)} \right)$ if and only if  there exist $\left(\theta_0, c_0 \right)\in \Gamma\left( U_0 \bigcup U_2, \mathcal{E}_{\mathcal{O}_{\mathbb{F}_m}(p,q)} \right)$ and $\left( \theta_1, c_1 \right)\in \Gamma\left( U_1 \bigcup U_3, \mathcal{E}_{\mathcal{O}_{\mathbb{F}_m}(p,q)} \right)$ such that $\theta_{01}= \theta_1- \theta_0$ and $s_{01}= c_1+ \left[ \theta_1, \log z_0^{mq-p}\right] - c_0$, where $\left( \theta_0, c_0\right)$ is written in the form as in $(\ref{rs29})$ and $\left( \theta_1, c_1  \right)$ is written in the form as in $(\ref{rs30})$. We write $\left(\theta_{01}, c_{01}\right)$ in terms of the coordinate $\left(z_0, w_0 \right)$
\begin{align*}
\left(\theta_{01}, c_{01}\right)=\left( v_{01}(z_0)\frac{\partial}{\partial z_0} + \left(\alpha_{01}^1(z_0) + \alpha_{01}^2 (z_0)w_0 +\alpha_{01}^3(z_0)w_0^2\right)\frac{\partial}{\partial w_0}, s_{01}(z_0)- q\alpha_{01}^3(z_0)w_0   \right)
\end{align*}
where $v_{01}(z_0),\alpha_{01}^1(z_0), \alpha_{01}^2(z_0), \alpha_{01}^3(z_0), s_{01}(z_0)$ are holomorphic functions of $z_0$ on $\mathbb{C}^*$, so that they are expanded into Laurent series in $z_0$. We note that 
\begin{align*}
 &s_1(z_1)- q\alpha_1^3(z_1)w_1+ \left[ v_1(z_1)\frac{\partial}{\partial z_1} + (\alpha_1^1(z_1) + \alpha_1^2 (z_1)w_1 +\alpha_1^3(z_1)w_1^2)\frac{\partial}{\partial w_1} , \log z_1^{p-mq} \right]\\
 &= h_1(z_1)- q\alpha_1^3(z_1)w_1+ z_1^{mq-p} (p-mq) v_1(z_1)z_1^{p-mq-1}\\
 &= h_1(z_1)- q\alpha_1^3(z_1)w_1+  (p-mq)z_1^{-1}v_1(z_1)
\end{align*}
Then the equation $s_{01}= c_1+\left[ \theta_1, \log z_0^{mq-p}\right]- c_0$ is reduced to the following system of equations
\begin{align}
h_1\left(\frac{1}{z_0}  \right) +(p-mq)z_0v_1\left(\frac{1}{z_0} \right)- h_0(z_0) &= h_{01}(z_0) \label{rs31}\\
qz_0^m\alpha_1^3\left(\frac{1}{z_0} \right) - q\alpha_0^3(z_0) &= q\alpha_{01}^3(z_0)
\end{align}
On the other hand, $\theta_{01}=\theta_1-\theta_0$ is reduced to the following system of equations (see \cite{Kod05} p.312):
\begin{align}
- z_0^2 v_1\left( \frac{1}{z_0} \right)- v_0 (z_0)&= v_{01}(z_0)\\
z_0^m \alpha_1^3\left( \frac{1}{z_0} \right)- \alpha_0^3(z_0)&= \alpha_{01}^2(z_0)\\
mz_0 v_1\left( \frac{1}{z_0}\right)+ \alpha_1^2\left( \frac{1}{z_0} \right)- \alpha_0^2(z_0) &= \alpha_{01}^2(z_0)\\
\frac{1}{z_0^m} \alpha_1^1\left( \frac{1}{z_0} \right)- \alpha_0^1 (z_0) &= \alpha_{01}^1 (z_0) \label{rs32}
\end{align}
In case $m=0,1$, $(\ref{rs31})-(\ref{rs32})$ always have a solution. For $m\geq 2$, let $\alpha_{01}^1(z_0)= \sum_{n= - \infty}^{\infty} c_n z_0^n$ be the Laurent expansion of $\alpha_{01}^1(z_0)$. Then $(\ref{rs31})- (\ref{rs32})$ have a solution if and only if $c_{-1}=\cdots = c_{-(m-1)}=0$. Hence  for $m\geq 2$, the $1$-cocycles 
\begin{align}\label{rsf3}
\left(\theta_{01}, h_{01}\right)=\left( \frac{1}{z_0^k}\frac{\partial}{\partial w_0}  , 0 \right)  \in Z^1\left( \mathfrak{B},   \mathcal{E}_{\mathcal{O}_{\mathbb{F}_m}(p,q)} \right),\,\,\,\,\,\,\,\,k=1,2,...,m-1
\end{align}
forms a basis of $H^1\left( \mathbb{F}_m, \mathcal{E}_{\mathcal{O}_{\mathbb{F}_m} (p,q)} \right)$ (see \cite{Kod05} p.312).

\subsubsection{$H^2\left( \mathbb{F}_m,  \mathcal{E}_{\mathcal{O}_{\mathbb{F}_m}(p,q)} \right)=0$}\

From the exact sequence $(\ref{rs2})$, since $H^2\left( \mathbb{F}_m, \mathcal{O}_{\mathbb{F}_m} \right)=H^2\left( \mathbb{F}_m, \Theta_{\mathbb{F}_m} \right)=0$ (see \cite{Kod05} p.312), we see that 
\begin{align}\label{rs17}
H^2\left( \mathbb{F}_m, \mathcal{E}_{\mathcal{O}_{\mathbb{F}_m}(p,q)} \right)=0
\end{align}

\subsection{Descriptions of $H^0\left( \mathbb{F}_m, \Theta_{\mathbb{F}_m}(p,q) \right)$}\

We keep the nations in subsection \ref{rs1}. We shall explicitly describe the global section $H^0\left( \mathbb{F}_m, \Theta_{\mathbb{F}_m}(p,q) \right)$ on $U_0$, where $\Theta_{\mathbb{F}_m}(p,q):=\Theta_{\mathbb{F}_m}\otimes \mathcal{O}_{\mathbb{F}_m}(p,q)$. Let 
\begin{align*}
T=v_0^1(z_0,w_0)\frac{\partial}{\partial z_0} + v_0^2(z_0,w_0)\frac{\partial}{\partial w_0}\in \Gamma\left(\mathbb{F}_m, \Theta_{\mathbb{F}_m}(p,q) \right)\,\,\,\,\,\,\,\textnormal{ on $U_0$}
\end{align*}
 with $v_0^1(z_0,w_0)=\sum_{h,k\geq 0} a_{hk} z_0^h w_0^k, a_{hk}\in \mathbb{C}$ and $v_0^2(z_0,w_0)=\sum_{h,k\geq 0} b_{hk} z_0^h w_0^k, b_{hk}\in \mathbb{C}$. We note that from $z_0=z_2$ and $w_2=\frac{1}{w_0}$, we have $\frac{\partial}{\partial z_0}=\frac{\partial}{\partial z_2}$ and $\frac{\partial}{\partial w_0}=-\frac{1}{w_0^2}\frac{\partial}{\partial w_2}$. Then $T$ on $U_0$ is transformed on $U_2$ into
\begin{align*}
&v_1^1(z_2,w_2)\frac{\partial}{\partial z_2} + v_2^2(z_2,w_2)\frac{\partial}{\partial w_2} = h_{20}\cdot\left( v_0^1(z_0,w_0)\frac{\partial}{\partial z_0} + v_0^2(z_0,w_0)\frac{\partial}{\partial w_0} \right) = \frac{1}{w_0^q}\left( v_0^1(z_0,w_0)\frac{\partial}{\partial z_2} + v_0^2(z_0,w_0)\left( -\frac{1}{w_0^2}\frac{\partial}{\partial w_2}\right) \right)\\
&= w_2^q v_0^1\left( z_2, \frac{1}{w_2}\right) \frac{\partial}{\partial z_2} -w_2^{q+2}v_0^2\left( z_2, \frac{1}{w_2}\right) \frac{\partial}{\partial w_2} = \left(\sum_{h,k\geq 0} a_{hk}z_2^h w_2^{q-k}\right) \frac{\partial}{\partial z_2} + \left( \sum_{h,k\geq 0} b_{hk}  z_2^h w_2^{q+2-k} \right) \frac{\partial}{\partial w_2} 
\end{align*}
which should be holomorphic, so that $a_{hk}=0$ for $k>q$ and $b_{hk}=0$ for $k>q+2$. Hence $T$ is of the following form
\begin{align}\label{rs8}
&v_0^1(z_0,w_0)\frac{\partial}{\partial z_0} + v_0^2(z_0,w_0)\frac{\partial}{\partial w_0} \\
&= \left(\alpha_0^1(z_0)+\alpha_0^2(z_0)w_0+ \cdots +\alpha_0^{q+1}(z_0)w_0^q\right)\frac{\partial}{\partial z_0}+ \left(\beta_0^1(z_0)+ \beta_0^2(z_0)w_0+\cdots +\beta_0^{q+3}(z_0)w_0^{q+2}\right)\frac{\partial}{\partial w_0} \notag
\end{align}
where $\alpha_0^1(z_0),..., \alpha_0^{q+1}(z_0),\beta_0^1(z_0),..., \beta_0^{q+3}(z_0)$ are entire functions of $z_0$. In particular, we note that we have
\begin{align*}
H^0\left( \mathbb{F}_m, \Theta_{\mathbb{F}_m}(p,q) \right) = 0 \,\,\,\,\,\,\,\,\textnormal{for}\,\,\, q< -2
\end{align*}
Hence there is no foliation on $\mathbb{F}_m$ if $q<-2$, so that in the following we assume that $q\geq -2$. In a similar way with $U_1$ and $U_3$, we see that $T\in H^0\left(\mathbb{F}_m, \Theta_{\mathbb{F}_m} (p,q) \right)$ on $U_1$ is of the form
\begin{align}\label{rs19}
&v_1^1(z_1,w_1)\frac{\partial}{\partial z_1} + v_1^2(z_1,w_1)\frac{\partial}{\partial w_1} \\
&= \left(\alpha_1^1(z_1)+\alpha_1^2(z_1)w_1+ \cdots +\alpha_1^{q+1}(z_1)w_1^q\right)\frac{\partial}{\partial z_1}+ \left(\beta_1^1(z_1)+ \beta_1^2(z_1)w_1+\cdots +\beta_1^{q+3}(z_1)w_1^{q+2}\right)\frac{\partial}{\partial w_1} \notag
\end{align}
where $\alpha_1^1(z_1),..., \alpha_1^{q+1}(z_1),\beta_1^1(z_1),..., \beta_1^{q+3}(z_1)$ are entire functions of $z_1$. 

On the other hand, $T$ on $U_0$ in $(\ref{rs8})$ is transformed on $U_1$ into
\begin{align*}
&h_{10}\cdot \left(\left(\alpha_0^1(z_0)+\alpha_0^2(z_0)w_0+ \cdots +\alpha_0^{q+1}(z_0)w_0^q\right)\frac{\partial}{\partial z_0}+ \left(\beta_0^1(z_0)+ \beta_0^2(z_0)w_0+\cdots +\beta_0^{q+3}(z_0)w_0^{q+2}\right)\frac{\partial}{\partial w_0}\right)\\
&=z_0^{mq-p}\cdot \left( \left(\alpha_0^1(z_0)+\alpha_0^2(z_0)w_0+ \cdots +\alpha_0^{q+1}(z_0)w_0^q \right)\frac{\partial}{\partial z_0}+ \left(\beta_0^1(z_0)+ \beta_0^2(z_0)w_0+\cdots +\beta_0^{q+3}(z_0)w_0^{q+2} \right) \frac{\partial}{\partial w_0}\right)\\
&= \left(\alpha_1^1(z_1)+\alpha_1^2(z_1)w_1+ \cdots +\alpha_1^{q+1}(z_1)w_1^q\right)\frac{\partial}{\partial z_1}+\left(\beta_0^1(z_1)+ \beta_1^2(z_1)w_1+\cdots +\beta_1^{q+3}(z_1)w_1^{q+2}\right)\frac{\partial}{\partial w_1}
\end{align*}
We note that from $z_1=\frac{1}{z_0}$ and $w_1=z_0^mw_0$, we have $\frac{\partial}{\partial z_0}=-\frac{1}{z_0^2}\frac{\partial}{\partial z_1}+mz_0^{m-1}w_0\frac{\partial}{\partial w_1}$ and $\frac{\partial}{\partial w_0}=z_0^m \frac{\partial}{\partial w_1}$. Then we have
\begin{align}\label{rs9}
&z_0^{mq-p} \left(\alpha_0^1(z_0)+\alpha_0^2(z_0)w_0+ \cdots +\alpha_0^{q+1}(z_0)w_0^q\right)\left(-\frac{1}{z_0^2}\frac{\partial}{\partial z_1}+mz_0^{m-1}w_0\frac{\partial}{\partial w_1} \right) \\
&+z_0^{mq-p} \left(\beta_0^1(z_0)+ \beta_0^2(z_0)w_0+\cdots +\beta_0^{q+3}(z_0)w_0^{q+2} \right) z_0^m\frac{\partial}{\partial w_1} \notag \\
&= \left(\alpha_1^1(z_1)+\alpha_1^2(z_1)w_1+ \cdots +\alpha_1^{q+1}(z_1)w_1^q \right)\frac{\partial}{\partial z_1}+(\beta_1^1(z_1)+ \beta_1^2(z_1)w_1+\cdots +\beta_1^{q+3}(z_1)w_1^{q+2})\frac{\partial}{\partial w_1} \notag
\end{align}
By considering the coefficient of $\frac{\partial}{\partial z_1}$ in $(\ref{rs9})$, we have
\begin{align*}
&-z_0^{mq-p-2} \left(\alpha_0^1(z_0)+\alpha_0^2(z_0)w_0+ \cdots +\alpha_0^{q+1}(z_0)w_0^q\right)\\
&=-z_1^{-mq+p+2} \left(\alpha_0^1\left(\frac{1}{z_1} \right) + \alpha_0^2\left(\frac{1}{z_1} \right) z_1^m w_1  + \cdots + \alpha_0^{q+1}\left(\frac{1}{z_1} \right) z_1^{mq} w_1^q \right) \\
& = \alpha_1^1(z_1)+\alpha_1^2(z_1)w_1+ \cdots +\alpha_1^{q+1}(z_1)w_1^q
\end{align*}
Hence we have
\begin{align}\label{rs10}
\alpha_1^1(z_1)&=-z_1^{p+2-mq}\alpha_0^1\left( \frac{1}{z_1}\right), \notag \\
\alpha_1^2(z_1)&= -z_1^{p+2-m(q-1)}\alpha_0^2\left(\frac{1}{z_1} \right) \notag \\
&\cdots \notag \\
\alpha_1^i(z_1)&= -z_1^{p+2-m(q-(i-1))}\alpha_0^i \left(\frac{1}{z_1}\right),\\
&\cdots \notag \\
\alpha_1^{q+1}(z_1)&=  -z_1^{p+2}\alpha_0^{q+1}\left(\frac{1}{z_1}\right)  \notag
\end{align}
On the other hand, by considering the coefficient of $\frac{\partial}{\partial w_1}$ in $(\ref{rs9})$, we have
\begin{align*}
&mz_0^{m(q+1)-p-1} \left(\alpha_0^1(z_0)w_0+\alpha_0^2(z_0)w_0^2+ \cdots +\alpha_0^{q+1}(z_0)w_0^{q+1}\right) + z_0^{mq-p+m} \left(\beta_0^1(z_0)+ \beta_0^2(z_0)w_0+\cdots +\beta_0^{q+3}(z_0)w_0^{q+2} \right) \\
&=mz_1^{-mq-m+p+1}\left(\alpha_0^1\left( \frac{1}{z_1} \right)z_1^m w_1 +\alpha_0^2\left( \frac{1}{z_1}\right) z_1^{2m}w_1^2 + \cdots + \alpha_0^i\left(\frac{1}{z_1} \right) z_1^{mi}w_1^i +\cdots+\alpha_0^{q+1}\left( \frac{1}{z_1} \right) z_1^{m(q+1)}w_1^{q+1} \right) \\
&+ z_1^{-mq+p-m} \left(\beta_0^1\left(\frac{1}{z_1} \right)+ \beta_0^2\left(\frac{1}{z_1}\right)z_1^m w_1+\cdots +  \beta_0^{i+1}\left( \frac{1}{z_1}\right)z_1^{mi}w_1^i + \cdots+\beta_0^{q+3}\left(\frac{1}{z_1}\right)z_1^{m(q+2)}w_1^{q+2} \right) \\
&= z_1^{-mq+p-m }\beta_0^1\left(\frac{1}{z_1} \right)  \\
&+  \left( mz_1^{-mq+p+1} \alpha_0^1\left(\frac{1}{z_1} \right)+ z_1^{-mq+p}\beta_0^2 \left( \frac{1}{z_1}\right)\right) w_1 \\
&+ \left(mz_1^{-mq+m+p+1}\alpha_0^2\left(\frac{1}{z_1}\right) +z_1^{-mq+p+m } \beta_0^3\left(\frac{1}{z_1} \right) \right)w_1^2 \\
& +\cdots \\
&+ \left(mz_1^{-mq-m+p+1+mi} \alpha_0^i\left(\frac{1}{z_1} \right) +z_1^{-mq+p+ m(i-1)}\beta_0^{i+1}\left( \frac{1}{z_1}\right) \right) w_1^i \\
&  +\cdots  \\
&+\left(mz_1^{p+1}\alpha_0^{q+1}\left( \frac{1}{z_1}\right)+ z_1^{p}\beta_0^{q+2}\left( \frac{1}{z_1} \right) \right)w_1^{q+1}  \\
&+z_1^{m+p}\beta_0^{q+3} \left( \frac{1}{z_1} \right) w_1^{q+2}  \\
&= \beta_1^1(z_1)+ \beta_1^2(z_1)w_1+\cdots +\beta_1^{q+3}(z_1)w_1^{q+2} 
\end{align*}
Hence we have
\begin{align} \label{rs11}
\beta_1^1(z_1)&=z_1^{-mq+p-m}\beta_0^1\left(\frac{1}{z_1} \right) \notag \\
&\cdots \notag\\
\beta_1^{i+1}(z_1)&= mz_1^{-mq-m+p+1+mi} \alpha_0^i\left(\frac{1}{z_1} \right) +z_1^{-mq+p+ m(i-1)}\beta_0^{i+1}\left( \frac{1}{z_1}\right)   \\
&\cdots \notag \\
\beta_1^{q+3}(z_1)&= z_1^{m+p} \beta_0^{q+3}\left( \frac{1}{z_1} \right) \notag
\end{align}

In the following examples, we will consider the case $q=0$, i.e. we will describe $H^0\left( \mathbb{F}_m,\Theta_{\mathbb{F}_m}(p,0)\right)$. Then from  $(\ref{rs10})$ and $(\ref{rs11})$, we have
\begin{align}
\alpha_1^1(z_1)&=-z_1^{p+2-mq}\alpha_0^1\left( \frac{1}{z_1} \right)=-z_1^{p+2}\alpha_0^1\left(\frac{1}{z_1} \right) \notag \\
\beta_1^1(z_1)&= z_1^{-mq+p-m}\beta_0^1\left(\frac{1}{z_1} \right)=z_1^{p-m}\beta_0^1\left(\frac{1}{z_1} \right) \label{rs14} \\
\beta_1^2(z_1)&=  mz_1^{-mq+p+1} \alpha_0^1\left(\frac{1}{z_1} \right)+ z_1^{-mq+p}\beta_0^2 \left( \frac{1}{z_1}\right)=m z_1^{p+1} \alpha_0^1\left(\frac{1}{z_1} \right)+z_1^{p}\beta_0^2\left(  \frac{1}{z_1} \right) \notag \\
\beta_1^3(z_1)&=z_1^{-mq+p+m}\beta_0^3\left( \frac{1}{z_1} \right)=z_1^{p+m}\beta_0^3\left( \frac{1}{z_1} \right) \notag
\end{align}

\begin{example}\label{rs54}
Let us describe $H^0\left( \mathbb{F}_m, \Theta_{\mathbb{F}_m}(-1,0) \right)$. From $(\ref{rs14})$, we have
\begin{align*}
\alpha_0^1(z_0)&=a_0+a_1z_0\\
\beta_0^1(z_0)&=0\\
\beta_1^2(z_1)&=ma_0+ma_1\frac{1}{z_1}+\frac{1}{z_1}\beta_0^2\left( \frac{1}{z_1} \right) \Longrightarrow \beta_0^2(z_0)= - ma_1\\
\beta_0^2(z_0) &= -ma_1\\
\beta_0^3(z_0)&=\begin{cases} 0 & \textnormal{if}\,\,\,m=0\\
                         b_0+b_1z_0+\cdots + b_{m-1}z_0^{m-1} & \textnormal{if}\,\,\,m\geq 1\end{cases}
\end{align*}
Then from $(\ref{rs8})$ any element in $H^0\left( \mathbb{F}_m, \Theta_{\mathbb{F}_m}(-1,0) \right)$ on $U_0$ is of the following form
\begin{align*}
\left(a_0+a_1z_0 \right)\frac{\partial}{\partial z_0}\,\,\,\,\,\,\,\textnormal{for}\,\,\,\,\,\,\,\, m=0
\end{align*}
\begin{align}\label{rsf9}
\left(a_0+a_1z_0\right)\frac{\partial}{\partial z_0}+ \left(-ma_1w_0+ \left(b_0+b_1z_0+\cdots +b_{m-1}z_0^{m-1} \right)w_0^2 \right)\frac{\partial}{\partial w_0}\,\,\,\,\,\,\textnormal{for}\,\,\,\,\,\,\, m \geq 1
\end{align}
where $a_0,a_1,b_0,..., b_{m-1}\in \mathbb{C}$, so that we have $\dim_\mathbb{C} H^0\left( \mathbb{F}_m, \Theta_{\mathbb{F}_m}(-1,0) \right)=m+2$
\end{example} 

\begin{example}\label{rr1}
Let us describe $H^0\left( \mathbb{F}_m, \Theta_{\mathbb{F}_m}(0,0) \right)= H^0 \left( \mathbb{F}_m, \Theta_{\mathbb{F}_m} \right)$. Though we already know the description as in $(\ref{rs13})$, we repeat here. From $(\ref{rs14})$, we have
\begin{align*}
\alpha_0^1(z_0)&=a_0+a_1z_0+a_2z_0^2\\
\beta_0^1(z_0)& =  \begin{cases} e_0 & \textnormal{for $m=0$}\\
                         0 & \textnormal{for $m\geq 1$} \end{cases} \\
\beta_1^2(z_1)&=mz_1\left(a_0+ a_1\frac{1}{z_1}+a_2\frac{1}{z_1^2} \right) + \beta_0^2\left(\frac{1}{z_1} \right)=ma_0z_1+ma_1+ma_2 \frac{1}{z_1}+ \beta_0^2\left(\frac{1}{z_1} \right) \Longrightarrow \beta_0^2(z_0)=c_0- ma_2z_0\\
\beta_0^3(z_0)&=b_0+b_1z_0+\cdots +b_m z_0^m
\end{align*}
Then from $(\ref{rs8})$ any element in $H^0\left( \mathbb{F}_m, \Theta_{\mathbb{F}_m}(0,0) \right)$ on $U_0$ is of the following form
\begin{align*}
\left(a_0+a_1z_0+a_2z_0^2 \right)\frac{\partial}{\partial z_0}+ \left(e_0+c_0w_0+b_0w_0^2\right)\frac{\partial}{\partial w_0}\,\,\,\,\,\,\,\textnormal{for}\,\,\,m=0
\end{align*}
\begin{align}\label{rsf26}
\left(a_0+a_1z_0+a_2z_0^2\right)\frac{\partial}{\partial z_0}+\left(\left(c_0-ma_2z_0\right)w_0+\left(b_0+b_1z_0+\cdots + b_mz_0^m \right)w_0^2 \right)\frac{\partial}{\partial w_0} \,\,\,\,\,\,\,\,\textnormal{for}\,\,\, m\geq 1
\end{align}
where $a_0,...,a_2, c_0, b_0,..., b_m \in \mathbb{C}$, so that we have $\dim_\mathbb{C} H^0\left(\mathbb{F}_0, \Theta_{\mathbb{F}_0}(0,0) \right)=6$ and $\dim_\mathbb{C} H^0\left( \mathbb{F}_m, \Theta_{\mathbb{F}_m}(0, 0)  \right)= m+5$ for $m\geq 1$. 
\end{example}

\begin{example}\label{rr2}
Let us describe $H^0\left( \mathbb{F}_m, \Theta_{\mathbb{F}_m}(1,0) \right)$. From $(\ref{rs14})$, we have
\begin{align*}
\alpha_0^1(z_0)&= a_0 + a_1 z_0 + a_2 z_0^2 + a_3 z_0^3\\
\beta_0^1(z_0)&= \begin{cases} e_0 + e_1 z_0 & \textnormal{if}\,\,\, m=0\\
                                                     e_0 & \textnormal{if}\,\,\,m=1 \\
                                                     0 & \textnormal{if}\,\,\,m\geq 2 \end{cases}\\
          \beta_1^2(z_1)               &= m z_1^2\left(a_0 + a_1 \frac{1}{z_1} + a_2 \frac{1}{z_1^2} + a_3 \frac{1}{z_1^3} \right)        + z_1 \beta_0^2\left( \frac{1}{z_1} \right)          \Longrightarrow \beta_0^2(z_0)= c_0 + c_1 z_0 - ma_3 z_0^2      \\
          \beta_0^3(z_0) &= b_0 + b_1 z_0  + \cdots + b_{m+1}z_0^{m+1}
\end{align*}
Then from $(\ref{rs8})$ any element in $H^0\left( \mathbb{F}_m, \Theta_{\mathbb{F}_m}(1,0) \right)$ on $U_0$ is of the following form
\begin{align*}
\left( a_0 + a_1 z_0 + a_2 z_0^2 + a_3 z_0^3 \right) \frac{\partial}{\partial z_0} +\left(e_0+ e_1 z_0 + \left(c_0 + c_1 z_0 \right) w_0 + \left( b_0 + b_1 z_0 \right) w_0^2 \right) \frac{\partial}{\partial w_0}\,\,\,\,\,\textnormal{for}\,\,\, m=0
\end{align*}
\begin{align*}
\left( a_0 + a_1 z_0 + a_2 z_0^2 + a_3 z_0^3 \right) \frac{\partial}{\partial z_0} +\left(e_0 + \left(c_0 + c_1 z_0  - a_3z_0^2\right) w_0 + \left( b_0 + b_1 z_0 + b_2 z_0^2 \right) w_0^2 \right) \frac{\partial}{\partial w_0}\,\,\,\,\,\textnormal{for}\,\,\, m=1
\end{align*}
\begin{align}\label{rf11}
\left( a_0 + a_1 z_0 + a_2 z_0^2 + a_3 z_0^3 \right) \frac{\partial}{\partial z_0} +\left( \left(c_0 + c_1 z_0  - ma_3z_0^2\right) w_0 + \left( b_0 + b_1 z_0 + \cdots + b_{m+1} z_0^{m+1} \right) w_0^2 \right) \frac{\partial}{\partial w_0}\,\,\,\,\,\textnormal{for}\,\,\, m\geq 2
\end{align}
where $a_0,...,a_3, e_0,...,e_1,c_0,c_1, b_0,..., b_{m+1}\in \mathbb{C}$, so that we have $\dim_\mathbb{C} H^0\left( \mathbb{F}_m, \Theta_{\mathbb{F}_m}(1,0) \right)= 10$ for $m=0,1$, and $\dim_\mathbb{C} H^0\left( \mathbb{F}_m, \Theta_{\mathbb{F}_m}(1,0) \right)= m+ 8$ for $m\geq 2$.
\end{example}

\subsection{Description of $H^1\left( \mathbb{F}_m, \Theta_{\mathbb{F}_m}(p,q) \right), q \geq -1$}\

We keep the notations in subsection \ref{rs1}. We explicitly describe $H^1\left( \mathbb{F}_m, \Theta_{\mathbb{F}_m}(p,q) \right)$. First we show that $H^1\left( U_0\bigcup U_2, \Theta_{\mathbb{F}_m}(p,q) \right)=0$ and $H^1\left( U_1 \bigcup U_3, \Theta_{\mathbb{F}_m}(p,q) \right)=0$ for $q\geq -1$. We note that by letting $\mathscr{B}=\left\{ U_0, U_2 \right\}$, we can compute $H^1\left( U_0 \bigcup U_2, \Theta_{\mathbb{F}_m}(p,q) \right)$ by \v Cech resolution associated to the open covering $\mathscr{B}$, i.e.
\begin{align*}
H^1\left( U_0 \bigcup U_2, \Theta_{\mathbb{F}_m}(p,q) \right) = H^1\left( \mathscr{B}, \Theta_{\mathbb{F}_m}(p,q)\right) = \frac{Z^1\left( \mathscr{B}, \Theta_{\mathbb{F}_m}(p,q) \right)}{\delta C^0\left( \mathscr{B}, \Theta_{\mathbb{F}_m}(p,q) \right)}
\end{align*}
A $1$-cocycle $\left\{\theta_{02}, \theta_{20} \right\}$ belonging to $Z^1\left(\mathscr{B},\Theta_{\mathbb{F}_m}(p,q)\right)$ consists of holomorphic vector fields $\theta_{02}$ and $\theta_{20}$ on $U_0\bigcap U_2$ with $\theta_{02}=-h_{02}\theta_{20}$, while $C^0\left(\mathscr{B},\Theta_{\mathbb{F}_m}(p,q)\right)$ is the set of all paris $\left\{\theta_1,\theta_2\right\}$, where $\theta_1$ and $\theta_2$ are holomorphic vector fields on $U_0$ and $U_2$. In order to prove $H^1\left(\mathscr{B}, \Theta_{\mathbb{F}_m} (p,q) \right)=0$, it suffices to show that any $1$-cocycle $\left\{\theta_{02}, \theta_{20}\right\}$ is represented as $\delta\left(\left\{\theta_2,\theta_0\right\}\right)$
\begin{align*}
\theta_{02}=h_{02}\theta_2-\theta_0= w_0^q \theta_2 - \theta_0.
\end{align*}
The holomorphic vector field $\theta_{02}$ on $U_0\bigcap U_2$ is written in the form
\begin{align*}
\theta_{02}=u_{02}(z_0,w_0)\frac{\partial}{\partial z_0}+v_{02}(z_0,w_0)\frac{\partial}{\partial w_0}
\end{align*}
where $u_{02}(z_0,w_0)$ and $v_{02}(z_0,w_0)$ are holomorphic functions on $U_0\cap U_2= \mathbb{C}\times \mathbb{C}^*$. Hence we have Laurent series in $w_0$:
\begin{align*}
u_{02}(z_0,w_0)=  \sum_{n= -\infty}^{\infty} u_n(z_0)w_0^n,\,\,\,\,\,\,\,\,\,\,\,v_{02}(z_0,w_0)=  \sum_{n=-\infty}^{\infty} v_n(z_0)w_0^n
\end{align*}
where $u_n(z_0)$ and $v_n(z_0)$ are entire functions of $z_0$. Any holomorphic vector field on $U_0$ is of the form
\begin{align*}
u_0(z_0,w_0) \frac{\partial}{\partial z_0} + v_1(z_0,w_0) \frac{\partial}{\partial w_0}
\end{align*}
where $u_0(z_0,w_0)$ and $v_1(z_0,w_0)$ are entire functions, and any holomorphic vector field on $U_2$ is of the form
\begin{align*}
u_2(z_2,w_2) \frac{\partial}{\partial z_2} + v_1(z_2,w_2) \frac{\partial}{\partial w_2}
\end{align*}
where $u_2(z_2,w_2)$ and $v_1(z_2,w_2)$ are entire functions. Since $\frac{\partial}{\partial w_2}=-\frac{1}{w_2^2}\frac{\partial}{\partial w_0}=-w_0^2\frac{\partial}{\partial w_0}$, we have
\begin{align*}
u_2(z_2,w_2)\frac{\partial}{\partial z_2}+ v_2(z_2,w_2)\frac{\partial}{\partial w_2}=u_2\left(z_0,\frac{1}{w_0} \right)\frac{\partial}{\partial z_0}-w_0^2u_2\left(z_0, \frac{1}{w_0} \right)\frac{\partial}{\partial w_0},
\end{align*}
so that we have
\begin{align*}
h_{02}\theta_2-\theta_0= w_0^q\left( u_2\left(z_0,\frac{1}{w_0} \right)\frac{\partial}{\partial z_0}-w_0^2u_2\left(z_0, \frac{1}{w_0} \right)\frac{\partial}{\partial w_0} \right) - u_0(z_0,w_0) \frac{\partial}{\partial z_0} - v_0(z_0,w_0) \frac{\partial}{\partial w_0}
\end{align*}
If $\theta_{02}$ defines the zero class, we should have a solution $\left( \theta_0, \theta_2 \right)$ satisfying
\begin{align*}
w_0^qu_2\left(z_0, \frac{1}{w_0} \right)- u_0(z_0,w_0)&= \sum_{n= -\infty}^{\infty} u_n(z_0)w_0^n\\
-w_0^{q+2}v_2\left( z_0, \frac{1}{w_0} \right) -v_0(z_0,w_0)&= \sum_{n=-\infty}^{\infty} v_n(z_0)w_0^n
\end{align*}
This implies that if $q  < -1$, we have $H^1\left( U_0\bigcup U_2, \Theta_{\mathbb{F}_m}(p,q) \right)\ne 0$. On the other hand, if $q\geq -1$, then in view of $w_0=\frac{1}{w_2}$, we set
\begin{align*}
u_0 & = -\sum_{n=0}^{\infty} u_n(z_0)w_0^n,\,\,\,\,\,\,\,\,\,\,\,u_2= \sum_{n=1}^\infty u_{-n}(z_2) w_2^{n+q}\\
v_0 & = -\sum_{n=0}^{\infty} v_n(z_0)w_0^n,\,\,\,\,\,\,\,\,\,\,\,v_2= -\sum_{n=1}^\infty v_{-n}(z_2) w_2^{n+q+2}
\end{align*}
Then there is a solution $\left(\theta_0, \theta_2 \right)$ satisfying $h_{02}\theta_2- \theta_0= \theta_{02}$. This implies that 
\begin{align*}
H^1\left( U_0 \bigcup U_2,  \Theta_{\mathbb{F}_m}(p,q) \right)=0,\,\,\,\,\,\,\,\,q\geq -1
\end{align*}
Similarly we can show that $H^1\left( U_1 \bigcup U_3, \Theta_{\mathbb{F}_m}(p,q) \right)=0$ for $q\geq -1$. In the following we will assume that $q \geq -1$. In the case of $q=-2$, we will study foliated deformations of $\left( \mathbb{F}_m, T \right)$ where $T\in H^0\left( \mathbb{F}_m, \Theta_{\mathbb{F}_m}(p,-2)\right)$ separately in Remark \ref{rs26} below. Then by letting $\mathfrak{B}=\left\{U_0 \bigcup U_2, U_1 \bigcup U_3 \right\}$, we can compute $H^1\left( \mathbb{F}_m, \Theta_{\mathbb{F}_m} (p,q)\right), q \geq -1$ by \v Cech resolution associated to the open covering $\mathfrak{B}$, i.e.
\begin{align}\label{rs20}
H^1\left( \mathbb{F}_m, \Theta_{\mathbb{F}_m}(p,q) \right)= H^1\left( \mathfrak{B}, \Theta_{\mathbb{F}_m}(p,q) \right)= \frac{Z^1\left( \mathfrak{B}, \Theta_{\mathbb{F}_m}(p,q) \right)}{\delta C^0\left( \mathfrak{B}, \Theta_{\mathbb{F}_m}(p,q) \right)}
\end{align}

Next we describe explicitly $H^0\left(U_0\bigcup U_2, \Theta_{\mathbb{F}_m}(p,q) \right), q\geq -1$ on $U_0$ and $H^0\left( U_1 \bigcup U_3, \Theta_{\mathbb{F}_m}(p,q) \right),q\geq-1 $ on $U_1$. Let us describe $H^0\left( U_0 \bigcup U_2 , \Theta_{\mathbb{F}_m}(p,q) \right)$ on $U_0$ which is given by $(\ref{rs8})$
\begin{align}\label{rs21}
\left(\alpha_0^1(z_0)+\alpha_0^2(z_0)w_0+ \cdots +\alpha_0^{q+1}(z_0)w_0^q\right)\frac{\partial}{\partial z_0}+ \left(\beta_0^1(z_0)+ \beta_0^2(z_0)w_0+\cdots +\beta_0^{q+3}(z_0)w_0^{q+2}\right)\frac{\partial}{\partial w_0}
\end{align}
where $\alpha_0^1(z_0),..., \alpha_0^{q+1}(z_0), \beta_0^1(z_0),...,\beta_0^{q+3}$ are entire functions of $z_0$. Let us describe $H^0\left(U_1 \bigcup U_3 , \Theta_{\mathbb{F}_m}(p,q) \right)$ on $U_1$ which is given by $(\ref{rs19})$
\begin{align}\label{rs22}
 \left(\alpha_1^1(z_1)+\alpha_1^2(z_1)w_1+ \cdots +\alpha_1^{q+1}(z_1)w_1^q\right)\frac{\partial}{\partial z_1}+ \left(\beta_1^1(z_1)+ \beta_1^2(z_1)w_1+\cdots +\beta_1^{q+3}(z_1)w_1^{q+2}\right)\frac{\partial}{\partial w_1} 
\end{align}
where $\alpha_1^1(z_1),..., \alpha_1^{q+1}(z_1), \beta_1^1(z_1),..., \beta_1^{q+3}$ are entire functions of $z_1$.

Now let us describe $H^1\left( \mathbb{F}_m, \Theta_{\mathbb{F}_m}(p,q)\right),q \geq -1$ in terms of \v Cech resolution associated to $\mathfrak{B}$ from $(\ref{rs20})$. We represent a $1$-cocycle $\left\{\theta_{01},\theta_{10} \right\}\in Z^1\left(\mathfrak{B}, \Theta_{\mathbb{F}_m}(p,q)\right)$ with $\theta_{01}=-h_{01}\theta_{10}$. Then the $1$-cocycle belongs to $\delta C^0\left( \mathfrak{B}, \Theta_{\mathbb{F}_m}(p,q)\right)$ if and only if there exists $\theta_0\in \Gamma\left(U_0\bigcup U_2, \Theta_{\mathbb{F}_m}(p,q)\right)$ and $\theta_2 \in \Gamma\left(U_1 \bigcup U_3,\Theta_{\mathbb{F}_m}(p,q)\right)$ such that $\theta_{01}= h_{01}\theta_1- \theta_0= z_0^{p-mq}\theta_1- \theta_0$, where $\theta_0$ is written in the form as in $(\ref{rs21})$ and $\theta_1$ is written in the form as in $(\ref{rs22})$. We write $\theta_{01}$ in terms of the coordinate $\left(z_0, w_0 \right)$
\begin{align*}
\theta_{01}= \left(\gamma_{01}^1(z_0)+\gamma_{01}^2(z_0)w_0 + \cdots + \gamma_{01}^{q+1}(z_0)w_0^q \right)\frac{\partial}{\partial z_0}+ \left(\eta_{01}^1(z_0)+\eta_{01}^2(z_0)w_0+\cdots + \eta_{01}^{q+3}(z_0) w_0^{q+2}\right)\frac{\partial}{\partial w_0}
\end{align*}
where $\gamma_{01}^1(z_0),..., \gamma_{01}^{q+1}, \eta_{01}^1(z_0),...,\eta_{01}^{q+3}(z_0)$ are holomorphic functions of $z_0$ on $\mathbb{C}^*$, so that they are expanded into Laurent series in $z_0$. We note that since $\frac{\partial}{\partial z_1}=-\frac{1}{z_1^2}\frac{\partial}{\partial z_0}+mz_1^{m-1}w_1\frac{\partial}{\partial w_0}=-z_0^2\frac{\partial}{\partial z_0}+mz_0w_0\frac{\partial}{\partial w_0}$ and $\frac{\partial}{\partial w_1}=z_1^m\frac{\partial}{\partial w_0}=z_0^{-m}\frac{\partial}{\partial w_0}$, we have
\begin{align*}
h_{01}\theta_1&= h_{01}\left(\alpha_1^1(z_1)+\alpha_1^2(z_1)w_1+ \cdots +\alpha_1^{q+1}(z_1)w_1^q\right)\frac{\partial}{\partial z_1}+\left(\beta_1^1(z_1)+\beta_1^2(z_1)w_1+\cdots + \beta_1^{q+3}(z_1) w_1^{q+2}\right)\frac{\partial}{\partial w_1}\\
&= z_0^{p-mq}\left(\alpha_1^1\left( \frac{1}{z_0} \right) +\alpha_1^2\left( \frac{1}{z_0}\right)z_0^mw_0  + \cdots + \alpha_1^{q+1}\left( \frac{1}{z_0} \right) z_0^{mq}w_0^q\right)\left(-z_0^2\frac{\partial}{\partial z_0}+mz_0w_0\frac{\partial}{\partial w_0}\right)\\
&+z_0^{p-mq}\left( \beta_1^1\left( \frac{1}{z_0} \right)  +  \beta_1^2\left( \frac{1}{z_0}\right) z_0^mw_0 +\cdots + \beta_1^{q+3}\left(   \frac{1}{z_0} \right) z_0^{m(q+2)}w_0^{q+2}\right)z_0^{-m}\frac{\partial}{\partial w_0}\\
&= - z_0^{p-mq+2}\left( \alpha_1^1\left( \frac{1}{z_0} \right) +\alpha_1^2\left( \frac{1}{z_0}\right)z_0^mw_0  + \cdots +\alpha_1^i\left(\frac{1}{z_0} \right) z_0^{m(i-1)}w_0^{i-1} +\cdots + \alpha_1^{q+1}\left( \frac{1}{z_0}\right)z_0^{mq}w_0^q \right)\frac{\partial}{\partial z_0}\\
&+mz_0^{p-mq+1}\left( \alpha_1^1\left(  \frac{1}{z_0} \right)w_0+\alpha_1^2\left( \frac{1}{z_0}\right) z_0^mw_0^2+\cdots + \alpha_1^{q+1}\left( \frac{1}{z_0} \right)z_0^{mq}w_0^{q+1}\right)\frac{\partial}{\partial w_0}\\
&+z_0^{p-mq-m}\left( \beta_1^1\left( \frac{1}{z_0} \right)  +  \beta_1^2\left( \frac{1}{z_0}\right) z_0^mw_0 +\cdots + \beta_1^{q+3}\left(   \frac{1}{z_0} \right) z_0^{m(q+2)}w_0^{q+2}\right)\frac{\partial}{\partial w_0}\\
\end{align*}
Then the equation $\theta_{01}= z_0^{p-mq} \theta_1- \theta_0$ is reduced to the following system of equations:
\begin{align}
-z_0^{p-mq+2}\alpha_1^1\left( \frac{1}{z_0} \right)-\alpha_0^1(z_0)&=\gamma_{01}^1(z_0) \notag\\
-z_0^{p-mq+2+m}\alpha_1^2\left( \frac{1}{z_0} \right) - \alpha_0^1(z_0)&= \gamma_{01}^2(z_0) \notag\\
&\cdots  \notag \\
-z_0^{p+2-m(q+1-i)}\alpha_1^i\left( \frac{1}{z_0} \right)-\alpha_0^i(z_0)&=\gamma_{01}^i(z_0),\,\,\,\,\,\,\,\,\, 1\leq i\leq q+1 \label{rs23}\\
&\cdots \notag \\
-z_0^{p+2}\alpha_1^{q+1}\left( \frac{1}{z_0} \right) -\alpha_0^{q+1}(z_0)&=\gamma_{01}^{q+1}(z_0) \notag
\end{align}

On the other hand,

\begin{align}
z_0^{p-mq-m}\beta_1^1\left( \frac{1}{z_0} \right)-\beta_0^1(z_0)&=\eta_{01}^1(z_0),\,\,\,\,\,\,\,\,i=1 \notag \\
mz_0^{p-mq+1}\alpha_1^1\left(\frac{1}{z_0} \right)+z_0^{p-mq}\beta_1^2\left( \frac{1}{z_0} \right)-\beta_0^2(z_0)&=\eta_{01}^2(z_0) \notag\\
&\cdots  \notag \\
mz_0^{p+1-m(q+2-i)}\alpha_1^{i-1}\left( \frac{1}{z_0} \right)   + z_0^{p-m(q+2-i)}\beta_1^i\left(  \frac{1}{z_0} \right)   &=\eta_{01}^i(z_0),\,\,\,\,\,\,\,\,\,2\leq i\leq q+2 \label{rs24} \\
&\cdots \notag \\
mz_0^{p+1}\alpha_1^{q+1}\left( \frac{1}{z_0} \right)  + z_0^{p}\beta_1^{q+2}\left(  \frac{1}{z_0} \right)-\beta_0^{q+2}(z_0)&=\eta_{01}^{q+2}(z_0) \notag \\
z_0^{p+m}\beta_1^{q+3}\left(\frac{1}{z_0} \right)-\beta_0^{q+3}(z_0)&=\eta_{01}^{q+3}(z_0),\,\,\,\,\,\,\,i=q+3 \notag
\end{align}

In the following examples, we will consider $q=0$, i.e. we will describe $H^1\left(\mathbb{F}_m, \Theta_{\mathbb{F}_m}(p,0) \right)$. Then $\theta_0$ is written in the form
\begin{align*}
\alpha_0^1(z_0)\frac{\partial}{\partial z_0}+ \left(\beta_0^1(z_0)+ \beta_0^2(z_0)w_0+ \beta_0^3(z_0)w_0^2\right)\frac{\partial}{\partial w_0}
\end{align*}
and $\theta_1$ is written in the form
\begin{align*}
\alpha_1^1(z_1)\frac{\partial}{\partial z_1}+ \left(\beta_0^1(z_1)+ \beta_1^2(z_1)w_1+ \beta_1^3(z_1)w_1^2\right)\frac{\partial}{\partial w_1}
\end{align*}
and from $(\ref{rs23})$ and $(\ref{rs24})$, we have
\begin{align}
- z_0^{p+2}\alpha_1^1\left(\frac{1}{z_0} \right) - \alpha_0^1(z_0)= \gamma_{01}^1(z_0) \notag \\
z_0^{p-m}\beta_1^1\left(\frac{1}{z_0} \right)- \beta_0^1(z_0)= \eta_{01}^1(z_0) \label{rs25}\\
mz_0^{p+1}\alpha_1^1\left(\frac{1}{z_0} \right)+ z_0^{p} \beta_1^2\left(\frac{1}{z_0} \right) - \beta_0^2(z_0)= \eta_{01}^2(z_0)\notag \\
z_0^{p+m}\beta_1^3\left(\frac{1}{z_0} \right) - \beta_0^{3}(z_0)= \eta_{01}^{3}(z_0) \notag
\end{align}

\begin{example}\label{rs49}
Let us describe $H^1\left( \mathbb{F}_m, \Theta_{\mathbb{F}_m}(-1,0) \right)$. From $(\ref{rs25})$, we have
\begin{align*}
- z_0\alpha_1^1\left(\frac{1}{z_0} \right) - \alpha_0^1(z_0)= \gamma_{01}^1(z_0)\\
\frac{1}{z_0^{m+1}}\beta_1^1\left(\frac{1}{z_0} \right)- \beta_0^1(z_0)= \eta_{01}^1(z_0)\\
m\alpha_1^1\left(\frac{1}{z_0} \right)+ \frac{1}{z_0} \beta_1^2\left(\frac{1}{z_0} \right) - \beta_0^2(z_0)= \eta_{01}^2(z_0)\\
z_0^{m-1}\beta_1^3\left(\frac{1}{z_0} \right) - \beta_0^{3}(z_0)= \eta_{01}^{3}(z_0)
\end{align*}
In case of $m=0$, these equations always have a solution. For $m\geq 1$, let $\eta_{01}^1(z_0)= \sum_{n=- \infty}^{\infty} c_n z_1^n$ be the Laurent explansion of $\eta_{01}^1(z_0)$. Then the equations have a solution if and only if $c_{-1}= \cdots = c_{-m}=0$. Hence $\dim_\mathbb{C} H^1\left( \mathbb{F}_0, \Theta_{\mathbb{F}_0}(-1,0) \right)=0$, and $\dim_{\mathbb{C}}H^1\left(\mathbb{F}_m, \Theta_{\mathbb{F}_m}(-1, 0) \right)=m$ for $m\geq 1$, and the $1$-cocycles 
\begin{align}\label{rsf10}
\theta_{01}=  \frac{1}{z_0}\frac{\partial}{\partial w_0},....., \frac{1}{z_0^{m}}\frac{\partial}{\partial w_0} \in Z^1\left( \mathfrak{B}, \Theta_{\mathbb{F}_m}(-1,0) \right)
\end{align}
forms a basis of $H^1\left( \mathbb{F}_m, \Theta_{\mathbb{F}_m}(-1, 0) \right)$.
\end{example}

\begin{example}\label{rs50}
Let us describe $H^1\left(  \mathbb{F}_m, \Theta_{\mathbb{F}_m} (0,0) \right)$. From $(\ref{rs25})$, we have
\begin{align*}
- z_0^{2}\alpha_1^1\left(\frac{1}{z_0} \right) - \alpha_0^1(z_0)= \gamma_{01}^1(z_0)\\
\frac{1}{z_0^m}\beta_1^1\left(\frac{1}{z_0} \right)- \beta_0^1(z_0)= \eta_{01}^1(z_0)\\
mz_0\alpha_1^1\left(\frac{1}{z_0} \right)+  \beta_1^2\left(\frac{1}{z_0} \right) - \beta_0^2(z_0)= \eta_{01}^2(z_0)\\
z_0^{m}\beta_1^3\left(\frac{1}{z_0} \right) - \beta_0^{3}(z_0)= \eta_{01}^{3}(z_0)
\end{align*}
In case of $m=0,1$, these equations always have a solution. For $m\geq 2$, let $\eta_{01}^1(z_0)= \sum_{n=- \infty}^{\infty} c_n z_1^n$ be the Laurent explansion of $\eta_{01}^1(z_0)$. Then the equations have a solution if and only if $c_{-1}= \cdots = c_{-(m-1)}=0$. Hence $ H^1\left( \mathbb{F}_0, \Theta_{\mathbb{F}_0}(0,0) \right)=H^1\left(\mathbb{F}_1, \Theta_{\mathbb{F}_1} (0,0) \right)=0$, and $\dim_{\mathbb{C}}H^1\left(\mathbb{F}_m, \Theta_{\mathbb{F}_m}(0, 0) \right)=m-1$ for $m\geq 2$, and the $1$-cocycles 
\begin{align}\label{rsf27}
\theta_{01}=  \frac{1}{z_0}\frac{\partial}{\partial w_0},....., \frac{1}{z_0^{m-1}}\frac{\partial}{\partial w_0} \in Z^1\left( \mathfrak{B}, \Theta_{\mathbb{F}_m}(0,0) \right)
\end{align}
forms a basis of $H^1\left( \mathbb{F}_m, \Theta_{\mathbb{F}_m}(0, 0) \right)$.
\end{example}

\begin{example}\label{rs51}
Let us describe $H^1\left( \mathbb{F}_m, \Theta_{\mathbb{F}_m}(1,0) \right)$. From $(\ref{rs25})$, we have
\begin{align*}
- z_0^{3}\alpha_1^1\left(\frac{1}{z_0} \right) - \alpha_0^1(z_0)= \gamma_{01}^1(z_0)  \\
\frac{1}{z_0^{m-1}}\beta_1^1\left(\frac{1}{z_0} \right)- \beta_0^1(z_0)= \eta_{01}^1(z_0) \\
mz_0^{2}\alpha_1^1\left(\frac{1}{z_0} \right)+ z_0\beta_1^2\left(\frac{1}{z_0} \right) - \beta_0^2(z_0)= \eta_{01}^2(z_0) \\
z_0^{m+1}\beta_1^3\left(\frac{1}{z_0} \right) - \beta_0^{3}(z_0)= \eta_{01}^{3}(z_0) 
\end{align*}
In case of $m=0,1,2$, these equations always have a solution. For $m\geq 3$, let $\eta_{01}^1(z_0)= \sum_{n=- \infty}^{\infty} c_n z_1^n$ be the Laurent explansion of $\eta_{01}^1(z_0)$. Then the equations have a solution if and only if $c_{-1}= \cdots = c_{-(m-2)}=0$. Hence $H^1\left( \mathbb{F}_m, \Theta_{\mathbb{F}_m}(1,0) \right)=0$ for $m=0,1,2$, and $\dim_{\mathbb{C}}H^1\left(\mathbb{F}_m, \Theta_{\mathbb{F}_m}(-1, 0) \right)=m$ for $m\geq 3$, and the $1$-cocycles 
\begin{align} \label{rf10}
\theta_{01}=  \frac{1}{z_0}\frac{\partial}{\partial w_0},....., \frac{1}{z_0^{m-2}}\frac{\partial}{\partial w_0} \in Z^1\left( \mathfrak{B}, \Theta_{\mathbb{F}_m}(1,0) \right)
\end{align}
forms a basis of $H^1\left( \mathbb{F}_m, \Theta_{\mathbb{F}_m}(1, 0) \right)$.
\end{example}

With this preparation we study foliated deformations of $\mathbb{F}_m$ in terms of tangent sheaves.

\subsection{Foliated deformations of $\left(\mathbb{F}_m, \Theta_{\mathcal{F}_0} \right)$ in terms of tangent sheaves} \

Let us compute the cohomology group $\mathbb{H}^i\left(  \mathbb{F}_m, \mathcal{E}_{\mathcal{O}_{\mathbb{F}_m}(p,q)}^\bullet \right),i=0,1,2$ of the leaf complex $(\ref{rs16})$. By considering the spectral sequence associated to $\mathcal{E}_{\mathcal{O}_{\mathbb{F}_m}(p,q)}^\bullet$ and $H^2\left( \mathbb{F}_m, \mathcal{E}_{\mathcal{O}_{\mathbb{F}_m}(p,q)} \right)=0$ from $(\ref{rs17})$, we see that
{\small{\begin{align}
\mathbb{H}^1\left( \mathbb{F}_m, \mathcal{E}_{\mathcal{O}_{\mathbb{F}_m}(p,q)}^\bullet \right) &\cong \textnormal{coker}\left( H^0\left( \mathbb{F}_m, \mathcal{E}_{\mathcal{O}_{\mathbb{F}_m}(p,q)} \right) \to H^0 \left( \mathbb{F}_m, \Theta_{\mathbb{F}_m} (p,q) \right) \right) \bigoplus \textnormal{ker}\left( H^1\left( \mathbb{F}_m, \mathcal{E}_{\mathcal{O}_{\mathbb{F}_m}(p,q)} \right) \to H^1\left( \mathbb{F}_m, \Theta_{\mathbb{F}_m}(p,q) \right) \right)\\
\mathbb{H}^2\left( \mathbb{F}_m, \mathcal{E}_{\mathcal{O}_{\mathbb{F}_m}(p,q)}^\bullet  \right)& \cong \textnormal{coker}\left(H^1\left( \mathbb{F}_m, \mathcal{E}_{\mathcal{O}_{\mathbb{F}_m}(p,q)} \right) \to H^1 \left( \mathbb{F}_m, \Theta_{\mathbb{F}_m}(p,q) \right) \right)\label{rs33}
\end{align}}}

\begin{example}\label{rsf11}
Let us consider a foliated rational ruled surface $\left( \mathbb{F}_2, T \right)$, where $T_0= z_0 \frac{\partial}{\partial z_0} - 2 w_0 \frac{\partial}{\partial w_0} \in H^0\left( \mathbb{F}_2, \Theta_{\mathbb{F}_2}(-1,0) \right)$ on $U_0$. Then we have from $(\ref{rsf1}),(\ref{rsf9}),(\ref{rsf3})$ and $(\ref{rsf10})$ 
\begin{align*} 
H^0\left(\mathbb{F}_2, \mathcal{E}_{\mathcal{O}_{\mathbb{F}_2}(-1,0)} \right) & = \left\{ \left(a_0+ a_1z_0+ a_2 z_0^2 \right)\frac{\partial}{\partial z_0}+ \left( \left(-2a_2z_0+ a_3\right)w_0 + \left(c_0+ c_1z_0+ c_2z_0^2\right) w_0^2\right)\frac{\partial}{\partial w_0} + a_4 + a_2z_0  \right\}    \\
H^0\left( \mathbb{F}_2, \Theta_{\mathbb{F}_2}(-1,0) \right) &= \left\{ \left(a_0+ a_1z_0\right)\frac{\partial}{\partial z_0}+ \left(-2a_1w_0+ \left(b_0+b_1z_0\right)w_0^2 \right)\frac{\partial}{\partial w_0} \right\} \\
H^1\left( \mathbb{F}_2, \mathcal{E}_{\mathcal{O}_{\mathbb{F}_2}(-1,0)} \right) &= \left\langle \frac{1}{z_0}\frac{\partial}{\partial w_0} \right\rangle \\
H^1\left(\mathbb{F}_2, \Theta_{\mathbb{F}_2}(-1,0) \right)&= \left\langle \frac{1}{z_0}\frac{\partial}{\partial w_0},  \frac{1}{z_0^2} \frac{\partial}{\partial w_0} \right\rangle
\end{align*}
We compute $H^0\left( \mathbb{F}_2, \mathcal{E}_{\mathcal{O}_{\mathbb{F}_2}(-1,0)} \right)\to H^0\left( \mathbb{F}_2, \Theta_{\mathbb{F}_2}(-1,0) \right)$
{\small{\begin{align*}
& \left[ \left(a_0+ a_1z_0+ a_2 z_0^2\right)\frac{\partial}{\partial z_0}+ \left( \left(-2a_2z_0+ a_3 \right)w_0 + \left(c_0+ c_1z_0+ c_2z_0^2 \right) w_0^2\right)\frac{\partial}{\partial w_0}, z_0\frac{\partial}{\partial z_0}- 2w_0\frac{\partial}{\partial w_0} \right] +\left(a_4+a_2z_0\right)\left(z_0\frac{\partial}{\partial z_0} - 2 w_0\frac{\partial}{\partial w_0} \right)\\
&= \left(a_0+ a_4z_0\right)\frac{\partial}{\partial z_0}+ \left( - 2a_4w_0  + \left(2c_0 + c_1z_0 \right)w_0^2\right)  \frac{\partial}{\partial w_0} =0 \iff a_0=a_4=c_0=c_1=0
\end{align*}}}
Then we have
\begin{align*}
\dim_\mathbb{C} \mathbb{H}^0\left( \mathbb{F}_2, \mathcal{E}_{\mathcal{O}_{\mathbb{F}_2}(-1,0)}^\bullet  \right) =\dim_\mathbb{C} \textnormal{ker}\left(H^0 \left( \mathbb{F}_2, \mathcal{E}_{\mathcal{O}_{\mathbb{F}_2}(-1,0)}  \right) \to H^0\left( \mathbb{F}_2 , \Theta_{\mathbb{F}_2}(-1, 0)  \right) \right) = 4
\end{align*}
and
\begin{align*}
\textnormal{coker}\left( H^0 \left( \mathbb{F}_2, \mathcal{E}_{\mathcal{O}_{\mathbb{F}_2}(-1,0)}  \right) \to H^0\left( \mathbb{F}_2 , \Theta_{\mathbb{F}_2}(-1, 0) \right)  \right)=   4 - ( 8 - 4) =0
\end{align*}
We compute $H^1\left( \mathbb{F}_2, \mathcal{E}_{\mathcal{O}_{\mathbb{F}_2}(-1,0)} \right)\to H^1\left( \mathbb{F}_2, \Theta_{\mathbb{F}_2}(- 1,0) \right)$
\begin{align*}
\left[d\frac{1}{z_0}\frac{\partial}{\partial w_0},  z_0 \frac{\partial}{\partial z_0} - 2 w_0 \frac{\partial}{\partial w_0}   \right] + 0 \cdot \left(z_0\frac{\partial}{\partial z_0} - 2 w_0 \frac{\partial}{\partial w_0} \right) = - 2d \frac{1}{z_0}\frac{\partial}{\partial w_0} + dz_0\frac{1} {z_0^2}\frac{\partial}{\partial w_0}= -d \frac{1}{z_0} \frac{\partial}{\partial w_0}=0 \iff d=0
\end{align*}
Then we have $\textnormal{ker}\left( H^1\left( \mathbb{F}_2, \mathcal{E}_{\mathcal{O}_{\mathbb{F}_2}(-1,0)} \right)\to H^1\left( \mathbb{F}_2, \Theta_{\mathbb{F}_2}(- 1,0) \right) \right)=0$, so that
\begin{align*}
\mathbb{H}^1\left( \mathbb{F}_2, \mathcal{E}_{\mathcal{O}_{\mathbb{F}_2} (-1, 0)}\right)=0
\end{align*}
Hence $\left( \mathbb{F}_2, T_0 \right)$ is rigid in foliated deformations by \textnormal{Theorem \ref{cp12}}.

\end{example}

\begin{example}\label{rsf12}
Let us consider a foliated rational ruled surface $\left( \mathbb{F}_2, T \right)$, where $T_0= - z_0 \frac{\partial}{\partial z_0} + w_0 \frac{\partial}{\partial w_0} \in H^0\left( \mathbb{F}_2, \Theta_{\mathbb{F}_2}(0,0) \right)$ on $U_0$. Then we have from $(\ref{rsf1}), (\ref{rsf26}),(\ref{rsf3})$ and $(\ref{rsf27})$
\begin{align} 
H^0\left(\mathbb{F}_2, \mathcal{E}_{\mathcal{O}_{\mathbb{F}_2}(0,0)} \right) & = \left\{ \left(a_0+ a_1z_0+ a_2z_0^2 \right)\frac{\partial}{\partial z_0}+ \left( (a_3-2a_2z_0)w_0+ \left(b_0+b_1z_0+ b_2 z_0^2\right)w_0^2 \right)\frac{\partial}{\partial w_0}+ a_4  \right\}    \\
H^0\left( \mathbb{F}_2, \Theta_{\mathbb{F}_2}(0,0) \right) &= \left\{ \left(a_0+ a_1z_0+ a_2z_0^2\right)\frac{\partial}{\partial z_0}+  \left( \left(a_3-2a_2z_0\right)w_0+ \left(b_0+b_1z_0+ b_2 z_0^2\right)w_0^2 \right)\frac{\partial}{\partial w_0} \right\} \label{rsf29}\\
H^1\left( \mathbb{F}_2, \mathcal{E}_{\mathcal{O}_{\mathbb{F}_2}(0,0)} \right) &= \left\langle \frac{1}{z_0}\frac{\partial}{\partial w_0} \right\rangle \\
H^1\left(\mathbb{F}_2, \Theta_{\mathbb{F}_2}(0,0) \right)&= \left\langle \frac{1}{z_0}\frac{\partial}{\partial w_0} \right\rangle
\end{align}
We compute $H^0\left( \mathbb{F}_2,  \mathcal{E}_{\mathcal{O}_{\mathbb{F}_2}(0,0)} \right)\to H^0\left( \mathbb{F}_2, \Theta_{\mathbb{F}_2} (0,0) \right)$
{\small{\begin{align}
&\left[ \left(a_0+ a_1z_0+ a_2z_0^2\right)\frac{\partial}{\partial z_0}+ \left( \left(a_3-2a_2z_0\right)w_0+ \left(b_0+b_1z_0+ b_2 z_0^2\right)w_0^2 \right)\frac{\partial}{\partial w_0}, -z_0\frac{\partial}{\partial z_0} + w_0\frac{\partial}{\partial w_0} \right] - a_4z_0\frac{\partial}{\partial z_0}+ a_4w_0\frac{\partial}{\partial w_0} \notag\\
&=(-a_0- a_4z_0+ a_2z_0^2)\frac{\partial}{\partial z_0} + \left( (a_4- 2a_2z_0 ) w_0+ (-b_0 + b_2z_0^2 )w_0^2 \right)\frac{\partial}{\partial w_0} = 0 \iff a_0=a_2=a_4= b_0= b_2=0 \label{rsf30}
\end{align}}}
Then we have
\begin{align*}
\dim_\mathbb{C} \mathbb{H}^0\left( \mathbb{F}_2, \mathcal{E}_{\mathcal{O}_{\mathbb{F}_2}(0,0)}^\bullet \right)= \dim_\mathbb{C} \textnormal{ker}\left( H^0\left( \mathbb{F}_2, \mathcal{E}_{\mathcal{O}_{\mathbb{F}_2}(0,0)}  \right) \to H^0\left( \mathbb{F}_2, \Theta_{\mathbb{F}_m}(0,0)\right) \right) = 3
\end{align*}
and 
\begin{align*}
\dim_\mathbb{C} \textnormal{coker}\left( H^0\left( \mathbb{F}_2, \mathcal{E}_{\mathcal{O}_{\mathbb{F}_2}(0,0)}  \right) \to H^0\left( \mathbb{F}_2, \Theta_{\mathbb{F}_2}(0,0)\right) \right) = 7 - (8-3)=2
\end{align*}
By comparing $(\ref{rsf29})$ and $(\ref{rsf30})$, we see that
\begin{align*}
\textnormal{coker}\left( H^0\left( \mathbb{F}_2, \mathcal{E}_{\mathcal{O}_{\mathbb{F}_2}(0,0)} \right) \to H^0\left( \mathbb{F}_2 , \Theta_{\mathbb{F}_2}(0,0) \right) \right) \cong \,\,\, \textnormal{Span}_\mathbb{C} \left\langle  z_0\frac{\partial}{\partial z_0}, z_0w_0^2\frac{\partial}{\partial w_0} \right\rangle \,\,\,\textnormal{or}\,\,\, \textnormal{Span}_\mathbb{C} \left\langle w_0\frac{\partial}{\partial w_0}, z_0w_0^2\frac{\partial}{\partial w_0} \right\rangle
\end{align*}
On the other hand, we compute $H^1\left( \mathbb{F}_2, \mathcal{E}_{\mathcal{O}_{\mathbb{F}_2}(0,0)} \right)\to H^1\left( \mathbb{F}_2, \Theta_{\mathbb{F}_2}(0,0) \right)$
\begin{align*}
\left[k\frac{1}{z_0}\frac{\partial}{\partial w_0}, -z_0\frac{\partial}{\partial z_0} + w_0\frac{\partial}{\partial w_0} \right] + 0 \cdot \left(- z_0\frac{\partial}{\partial z_0} + w_0\frac{\partial}{\partial w_0} \right)= k \frac{1}{z_0}\frac{\partial}{\partial w_0} + k z_0\left(-\frac{1}{z_0^2} \right)\frac{\partial}{\partial w_0}=0
\end{align*}
This implies that
\begin{align*}
\dim_\mathbb{C} \mathbb{H}^1 \left( \mathbb{F}_2, \mathcal{E}_{\mathcal{O}_{\mathbb{F}_2}(0,0)}^\bullet  \right) =3,\,\,\,\,\,\,\,\,\, \mathbb{H}^1\left( \mathbb{F}_2, \mathcal{E}_{\mathcal{O}_{\mathbb{F}_2} (0,0)}^\bullet \right) \cong \textnormal{Span}_\mathbb{C} \left\langle  z_0\frac{\partial}{\partial z_0}, z_0w_0^2\frac{\partial}{\partial w_0} \right\rangle \bigoplus \textnormal{Span}_\mathbb{C} \left\langle \frac{1}{z_0} \frac{\partial}{\partial w_0} \right\rangle
\end{align*}
We note that
\begin{align*}
\left[\frac{1}{z_0}\frac{\partial}{\partial w_0},  z_0\frac{\partial}{\partial z_0} \right] =-z_0\frac{1}{z_0^2}\frac{\partial}{\partial w_0}=-\frac{1}{z_0}\frac{\partial}{\partial w_0}  \ne 0
\end{align*}
 is not in the image of $H^1\left( \mathbb{F}_2,  \mathcal{E}_{\mathcal{O}_{\mathbb{F}_2}(0,0)} \right) \cong H^1 \left( \mathbb{F}_2, \Theta_{\mathbb{F}_2}(0,0) \right) \to H^1\left( \mathbb{F}_2, \Theta_{\mathbb{F}_2}(0,0) \right)$ which is the zero map. Hence $\left( \mathbb{F}_2, \omega_0 \right)$ is obstructed in foliated deformations by \textnormal{Lemma \ref{rsf28}}.

\end{example}

\begin{example}\label{rsf13}
Let us consider a foliated rational ruled surface $\left( \mathbb{F}_2, T \right)$, where $T= \frac{\partial}{\partial z_0} + z_0w_0 \frac{\partial}{\partial w_0} \in H^0\left( \mathbb{F}_2, \Theta_{\mathbb{F}_2}(1,0) \right)$ on $U_0$. Then we have from $(\ref{rsf1}), (\ref{rf11}),(\ref{rsf3})$ and $(\ref{rf10})$
\begin{align}
&H^0\left( \mathbb{F}_2, \mathcal{E}_{\mathcal{O}_{\mathbb{F}_2}(1,0)}\right) = \left\{ \left(a_0+ a_1z_0+ a_2z_0^2\right)\frac{\partial}{\partial z_0}+\left( \left(-2a_2z_0+ a_3\right)w_0 + \left(c_0 + c_1z_0 + c_2z_0^2\right)w_0^2\right)\frac{\partial}{\partial w_0} + a_4- a_2z_0   \right\}\\
&H^0\left( \mathbb{F}_2, \Theta_{\mathbb{F}_2}(1,0) \right)=  \left\{ \left(a_0 + a_1z_0 + a_2 z_0^2 + a_3z_0^3\right) \frac{\partial}{\partial z_0} + \left( \left(c_0+ c_1z_0- 2a_3z_0^2\right)w_0 + \left( b_0+ b_1 z_0+ b_2z_0^2 + b_3 z_0^3\right) w_0^2 \right) \frac{\partial}{\partial w_0}    \right\} \label{rf12} \\
&H^1\left( \mathbb{F}_2, \mathcal{E}_{\mathcal{O}_{\mathbb{F}_2}(1,0)} \right) = \left\langle \frac{1}{z_0} \frac{\partial}{\partial w_0} \right\rangle\\
&H^1\left( \mathbb{F}_2, \Theta_{\mathbb{F}_2}(1,0) \right) =0
\end{align} 
We compute $H^0\left( \mathbb{F}_2, \mathcal{E}_{\mathcal{O}_{\mathbb{F}_2}(1,0)} \right) \to H^0\left( \mathbb{F}_2, \Theta_{\mathbb{F}_2}(1,0) \right)$
{\small{\begin{align}
&\left[ \left(a_0+ a_1z_0+ a_2z_0^2\right)\frac{\partial}{\partial z_0}+\left( \left(-2a_2z_0+ a_3\right)w_0 + \left(c_0 + c_1z_0 + c_2z_0^2\right)w_0^2\right)\frac{\partial}{\partial w_0} ,  \frac{\partial}{\partial z_0} + z_0 w_0 \frac{\partial}{\partial w_0} \right] + \left(a_4- a_2z_0\right)\left( \frac{\partial}{\partial z_0} + z_0 w_0 \frac{\partial}{\partial w_0}  \right) \notag\\
&= \left(a_4-a_1- 3a_2z_0 \right)\frac{\partial}{\partial z_0} +  \left( \left(a_0+ 2a_2+ (a_1+ 2a_3+ a_4)z_0 \right)w_0  + \left(-c_1+ \left(-c_0-2c_2\right)z_0 - c_1z_0^2 - c_2z_0^3 \right)w_0^2   \right)\frac{\partial}{\partial w_0}  \label{rf8}
\end{align}}}
Then we have
\begin{align*}
\dim_\mathbb{C} \mathbb{H}^0\left( \mathbb{F}_2 ,  \mathcal{E}_{\mathcal{O}_{\mathbb{F}_2}(1,0)}^\bullet \right)= \dim_\mathbb{C} \textnormal{ker}\left( H^0\left( \mathbb{F}_2, \mathcal{E}_{\mathcal{O}_{\mathbb{F}_2}(1,0)} \right) \to H^0\left( \mathbb{F}_2, \Theta_{\mathbb{F}_2}(1,0) \right)\right)= 1
\end{align*}
and 
\begin{align*}
\dim_\mathbb{C} \textnormal{coker}\left( H^0\left( \mathbb{F}_2, \mathcal{E}_{\mathcal{O}_{\mathbb{F}_2}(1,0)} \right) \to H^0\left( \mathbb{F}_2, \Theta_{\mathbb{F}_2}(1,0) \right) \right)= 10- (8-1)=3
\end{align*}
By comparing $(\ref{rf12})$ and $(\ref{rf8})$ we see that
\begin{align*}
\textnormal{coker}\left(H^0\left( \mathbb{F}_2, \mathcal{E}_{\mathcal{O}_{\mathbb{F}_2}(1,0)} \right) \to H^0\left( \mathbb{F}_2, \Theta_{\mathbb{F}_2}(1,0) \right)\right) \cong \textnormal{Span}_\mathbb{C} \left\langle z_0^2 \frac{\partial}{\partial z_0}, z_0^3 \frac{\partial}{\partial z_0} - 2 z_0^2 w_0 \frac{\partial}{\partial w_0}, w_0^2 \frac{\partial}{\partial w_0} \right\rangle
\end{align*}
Since $H^1\left( \mathbb{F}_2, \Theta_{\mathbb{F}_2}(1,0) \right)=0$, we have $\mathbb{H}^2\left( \mathbb{F}_2, \mathcal{E}_{\mathcal{O}_{\mathbb{F}_2}(1,0)}^\bullet \right)=0$, so that in particular $\left( \mathbb{F}_2, T_0 \right)$ is unobstructed in foliated deformations, and we have
\begin{align}\label{rf9}
\mathbb{H}^1\left( \mathbb{F}_2, \mathcal{E}_{\mathcal{O}_{\mathbb{F}_2}(1,0)}^\bullet \right) \cong \textnormal{Span}_\mathbb{C} \left\langle z_0^2 \frac{\partial}{\partial z_0}, z_0^3 \frac{\partial}{\partial z_0} - 2 z_0^2 w_0 \frac{\partial}{\partial w_0}, w_0^2 \frac{\partial}{\partial w_0}  \right\rangle \bigoplus \textnormal{Span}_\mathbb{C} \left\langle \frac{1}{z_0}\frac{\partial}{\partial w_0} \right\rangle
\end{align}

We will construct a foliated analytic family of deformations of $\left( \mathbb{F}_m, T_0 \right)$ such that the foliated Kodaira-Spencer map is an isomorphism at the distinguished point. Let us start from the linear term given by $(\ref{rf9})$. First we set
\begin{align*}
\theta_{01}:=\left[ \frac{1}{z_0}\frac{\partial}{\partial w_0}, \frac{\partial}{\partial z_0} + z_0 w_0 \frac{\partial}{\partial w_0} \right] + 0 \cdot \left( \frac{\partial}{\partial z_0} + z_0 w_0 \frac{\partial}{\partial w_0} \right)= \frac{\partial}{\partial w_0} + \frac{1}{z_0^2}\frac{\partial}{\partial w_0}  = \left( 1 +  \frac{1}{z_0^2} \right) \frac{\partial}{\partial w_0} \in C^1\left( \mathfrak{B}, \Theta_{\mathbb{F}_2}(1,0) \right)
\end{align*}
where $\mathfrak{B}=\left\{ U_0 \bigcup U_2, U_1 \bigcup U_3 \right\}$. From $w_0= z_1^2w_1$, we have $\frac{\partial}{\partial w_0}= z_1^2 \frac{\partial}{\partial w_0} = \frac{1}{z_0^2} \frac{\partial}{\partial w_0}$. Then if we take $\theta_0= \frac{\partial}{\partial w_0}$ on $U_0\bigcup U_2$ and $\theta_1= z_1\frac{\partial}{\partial w_1}$ on $U_1\bigcup U_3$, we have
\begin{align*}
h_{01}\cdot \theta_1 - \theta_0 + \theta_{01}= z_0\cdot \left( -z_1 \frac{\partial}{\partial w_1} \right) - \frac{\partial}{\partial w_0} + \left( 1+ \frac{1}{z_0^2} \right)\frac{\partial}{\partial w_0}
\end{align*}
Then let us consider the foliated analytic family over $\textnormal{Spec}\left(\frac{\mathbb{C}[t_1,t_2,t_3,t_4]}{\left( t_1^2, t_2^2, t_3^2, t_4^2, t_1t_2,t_1t_3, t_1t_4, t_2t_3, t_2t_4, t_3t_4 \right)}\right)$ on $U_0$ defined by $(\ref{rf9})$ 
\begin{align*}
T&= \left(1- t_4 z_0^2 - t_2z_0^3\right)\frac{\partial}{\partial z_0} + \left(z_0w_0 - t_1 + 2t_2z_0^2 w_0 - t_3w_0^2\right)\frac{\partial}{\partial w_0}\\
z_0& = \frac{1}{z_1} \\
w_0&= z_1^2 w_1 + t_1 z_1\\
h_{01}&= z_0 
\end{align*}
Since $z_1= \frac{1}{z_0}$ and $w_1= z_0^2 w_0- t_1 z_0$, we have
\begin{align*}
\frac{\partial}{\partial z_0}= - \frac{1}{z_0^2} \frac{\partial}{\partial z_1} + (2z_0w_0- t_1)\frac{\partial}{\partial w_1},\,\,\,\,\,\,\,\,\,\,\,\frac{\partial}{\partial w_0}=z_0^2\frac{\partial}{\partial w_1}
\end{align*}
Then we have
\begin{align*}
h_{10} T&=\left( \frac{1}{z_0}- t_4z_0 - t_2 z_0^2    \right) \frac{\partial}{\partial z_0}  + \left(  w_0 -  t_1\frac{1}{z_0} + 2t_2z_0w_0 -t_3 \frac{1}{z_0}w_0^2\right)\frac{\partial}{\partial w_0}  \\
           &= \left( - \frac{1}{z_0^3} +t_4 \frac{1}{z_0} + t_2   \right) \frac{\partial}{\partial z_1}+ \left( 2w_0 - 2t_4z_0^2 w_0 - 2t_2 z_0^3 w_0 - t_1\frac{1}{z_0}+ t_1t_4z_0+ t_1t_2 z_0^2 \right) \frac{\partial}{\partial w_1} \\
           & + \left(  z_0^2w_0 -  t_1z_0 + 2t_2z_0^3w_0 -t_3 z_0 w_0^2\right)\frac{\partial}{\partial w_1} \\
           &= \left(-z_1^3+ t_4z_1+ t_2\right)\frac{\partial}{\partial z_1} + \left(2 z_1^2 w_1 + 2t_1z_1 - 2t_4w_1- 2t_1t_4z_0 - t_1z_1+ t_1t_4\frac{1}{z_1}+ t_1t_2\frac{1}{z_1^2}  \right) \frac{\partial}{\partial w_1}\\
           &+ \left( w_1 - t_3\frac{1}{z_1}\left(z_1^2w_1+ t_1z_1\right)^2 \right) \frac{\partial}{\partial w_1}\\
           &=  \left(-z_1^3+ t_4z_1+ t_2\right)\frac{\partial}{\partial z_1} + \left(2 z_1^2 w_1 + t_1z_1 - 2t_4w_1- t_1t_4\frac{1}{z_1}+ t_1t_2\frac{1}{z_1^2}  + w_1 -t_3\left(z_1^3w_1^2+ 2t_1z_1^2 w_1+ t_1^2 z_1\right) \right) \frac{\partial}{\partial w_1}
\end{align*}
We note that the rational part is  $ \left(- t_1t_4 \frac{1}{z_1}+ t_1t_2 \frac{1}{z_1^2}\right) \frac{\partial}{\partial w_1}$ and we have
\begin{align*}
 \left(- t_1t_4 \frac{1}{z_1}+ t_1t_2 \frac{1}{z_1^2}\right) \frac{\partial}{\partial w_1}= \left(-t_1t_4 z_0 + t_1t_2 z_0^2 \right)\frac{1}{z_0^2} \frac{\partial}{\partial w_0} = \left( - t_1t_4 \frac{1}{z_0}+ t_1t_2  \right) \frac{\partial}{\partial w_0} = h_{10}\cdot\left(- t_1t_4 + t_1t_2 z_0 \right)\frac{\partial}{\partial w_0}
\end{align*}
Hence the desired foliated analytic family is given by
\begin{align*}
T|_{U_0\bigcup U_2}&=\left( 1- t_4z_0^2 - t_2 z_0^3     \right) \frac{\partial}{\partial z_0}  + \left( z_0 w_0 - t_1 + t_1t_4  - t_1t_2z_0 + 2t_2z_0^2w_0 -t_3 w_0^2\right)\frac{\partial}{\partial w_0}\\
T|_{U_1\bigcup U_3}&= \left(-z_1^3+ t_4z_1+ t_2\right)\frac{\partial}{\partial z_1} + \left(2 z_1^2 w_1 + w_1 + t_1z_1 - 2t_4w_1  -t_3z_1^3w_1^2-2t_1t_3z_1^2 w_1- t_1^2t_3 z_1 \right) \frac{\partial}{\partial w_1} \\
z_0 &= \frac{1}{z_1}\\
w_0&= z_0^2 w_0 + t_1z_1
\end{align*}

We describe the foliated analytic family more precisely. We have coordinate transformations as follows$:$
\begin{center}
$\begin{CD}
U_0: \left( z_0, w_0 \right) @>>> U_1:\left(z_1, w_1 \right)=\left( \frac{1}{z_0}, z_0^2 w_0- t_1 z_0 \right)\\
@VVV @VVV \\
U_2:\left( z_2, w_2 \right)= \left( z_0, \frac{1}{w_0} \right) @>>> U_3: \left( z_3 , w_3 \right) = \left( \frac{1}{z_2} , \frac{w_2}{ z_2^2 - t_1 z_2w_2 } \right)  =\left( z_1, \frac{1}{w_1} \right)
\end{CD}$
\end{center}
and the foliation $T$ is given by
\begin{align*}
T&=\left( 1- t_4z_0^2 - t_2 z_0^3     \right) \frac{\partial}{\partial z_0}  + \left( z_0 w_0 - t_1 + t_1t_4  - t_1t_2z_0 + 2t_2z_0^2w_0 -t_3 w_0^2\right)\frac{\partial}{\partial w_0} \,\,\,\,\,\,\textnormal{on $U_0$} \\
T&=  \left(-z_1^3+ t_4z_1+ t_2\right)\frac{\partial}{\partial z_1} + \left(2 z_1^2 w_1 + w_1 + t_1z_1 - 2t_4w_1  -t_3z_1^3w_1^2-2t_1t_3z_1^2 w_1- t_1^2t_3 z_1 \right) \frac{\partial}{\partial w_1}\,\,\,\,\,\textnormal{on $U_1$}\\
T&= \left( 1- t_4z_2^2 - t_2 z_2^3     \right) \frac{\partial}{\partial z_2} -\left( z_2 w_2 - t_1w_2^2 + t_1t_4w_2^2  - t_1t_2z_2w_2^2+ 2t_2z_2^2w_2 -t_3\right)\frac{\partial}{\partial w_2} \,\,\,\,\,\textnormal{on $U_2$}\\
T&=  \left(-z_3^3+ t_4z_3+ t_2\right)\frac{\partial}{\partial z_3} - \left(2 z_3^2 w_3 + w_3 + t_1z_3w_3^2 - 2t_4w_3  -t_3z_3^3-2t_1t_3z_3^2 w_3- t_1^2t_3 z_3 w_3^2 \right) \frac{\partial}{\partial w_3} \,\,\,\,\,\textnormal{on $U_3$}
\end{align*}

and transition functions are given by
\begin{align*}
h_{01}&= z_0\\
h_{02}&=1\\
h_{03}&= z_0\\
h_{12}&= z_1\\
h_{13}&=1\\
h_{23}&= z_3^{-1}
\end{align*}

We note that the foliated Kodaira-Spencer map at $t=(t_1,t_2,t_3,t_4)=(0,0,0,0)=0$ is given by $($in terms of the covering $\mathfrak{B}$ where $\mathfrak{B}=\left\{ U_0 \bigcup U_2, U_1 \bigcup U_3\right\})$
\begin{align*}
T_0 \left(\mathbb{C}^4\right) &\to \mathbb{H}^1\left( \mathbb{F}_2, \mathcal{E}_{\mathcal{O}_{\mathbb{F}_2}(1,0)}^\bullet   \right) \\
     \frac{\partial}{\partial t_1} &\mapsto  \left( \left\{ \left( \frac{\partial}{\partial w_0} , - z_1\frac{\partial}{\partial w_1} \right) \right\} , \left\{ \left(z_1\frac{\partial}{\partial w_0} , 0 \right) \right\} \right) \in C^0\left( \mathfrak{B}, \Theta_{\mathbb{F}_2}(1,0) \right)\bigoplus C^1 \left( \mathfrak{B}, \mathcal{E}_{\mathcal{O}_{\mathbb{F}_2}(1,0)} \right) \\
      \frac{\partial}{\partial t_2} &\mapsto \left( \left\{  \left(z_0^3\frac{\partial}{\partial z_0} - 2 z_0^2 w_0 \frac{\partial}{\partial w_0} ,  - \frac{\partial}{\partial z_1} \right)\right\}, \left\{  (0, 0) \right \}\right) \in C^0\left( \mathfrak{B}, \Theta_{\mathbb{F}_2}(1,0) \right)\bigoplus C^1 \left( \mathfrak{B}, \mathcal{E}_{\mathcal{O}_{\mathbb{F}_2}(1,0)} \right)\\
      \frac{\partial}{\partial t_3} & \mapsto \left( \left\{ \left( w_0^2\frac{\partial}{\partial w_0}, z_1^3 w_1^2 \frac{\partial}{\partial w_1} \right) \right\}, \left (0 , 0) \right\}   \right) \in C^0\left( \mathfrak{B}, \Theta_{\mathbb{F}_2}(1,0) \right)\bigoplus C^1 \left( \mathfrak{B}, \mathcal{E}_{\mathcal{O}_{\mathbb{F}_2}(1,0)} \right)\\
      \frac{\partial}{\partial t_4} & \mapsto \left( \left\{ \left( z_0^2 \frac{\partial}{\partial z_0},  - z_1 \frac{\partial}{\partial z_1} + 2 w_1 \frac{\partial}{\partial w_1} \right) \right\} , \left\{ (0 , 0 ) \right\}   \right) \in C^0\left( \mathfrak{B}, \Theta_{\mathbb{F}_2}(1,0) \right)\bigoplus C^1 \left( \mathfrak{B}, \mathcal{E}_{\mathcal{O}_{\mathbb{F}_2}(1,0)} \right)
\end{align*}

\end{example}

\begin{remark}\label{rs52}
In general, let $T\in H^0\left( \mathbb{F}_m, \Theta_{\mathbb{F}_m}(p,0)\right)$ be a foliation on $\mathbb{F}_m$. From $(\ref{rs25})$, we see that if $p\geq m-1$, then $H^1\left( \mathbb{F}_m, \Theta_{\mathbb{F}_m}(p,0) \right)=0$. This implies that for $p \geq m-1$, $\mathbb{H}^2\left( \mathbb{F}_m, \mathcal{E}_{\mathcal{O}_{\mathbb{F}_m}(p,0)}^\bullet \right)=0$ from $(\ref{rs33})$, so that $\left( \mathbb{F}_m, T \right)$ is unobstructed in foliated deformations.
\end{remark}

\begin{remark}\label{rs26}
So far we have assumed $q\geq -1$ since $H^1\left( U_0 \bigcup U_2, \Theta_{\mathbb{F}_m}(p,q) \right) \ne 0$ for $q\leq -2$. Now let us consider $q=-2$. Let $T\in H^0\left( \mathbb{F}_m, \Theta_{\mathbb{F}_m}(p,-2) \right)$ be a foliation on $\mathbb{F}_m$. Then from $(\ref{rs8})$ we see that $T$ is of the form $T_0= \beta_0^1(z_0) \frac{\partial}{\partial w_0}$, where $\beta_0^1(z_0)$ is an entire function. Unless $\beta_0^1(z_0)= a\ne 0$ is a constant, then $T$ does not have an isolated singularity, i.e. does not define a foliation. Hence $T_0$ is of the form $T_0= a\frac{\partial}{\partial w_0}, a \ne 0$ on $U_0$. We note that $T_0$ on $U_0$ is translated into $\beta_1^1(z_1)\frac{\partial}{\partial w_1}$ on $U_1$ where $\beta_1^1(z_1)= - a z_1^{p+m}$. In order for $T$ to define a foliation on $\mathbb{F}_m$, we should have $p=-m$. Hence we see that $T \in H^0\left( \mathbb{F}_m, \Theta_{\mathbb{F}_m} (-m, - 2) \right)$. Now let us consider foliated deformations of $\left( \mathbb{F}_m,  T_0= a \frac{\partial}{\partial w_0}\right), a\ne 0$. Then we have from $(\ref{rsf2}),(\ref{rsf1})$ and $(\ref{rsf3})$
\begin{align*}
H^0\left(\mathbb{F}_m, \mathcal{E}_{\mathcal{O}_{\mathbb{F}_m}(-m,-2)} \right)&= \begin{cases} \left\{ \left(a_0z_0^2+bz_0+c\right)\frac{\partial}{\partial z_0}+ \left(d+ ew_0+ fw_0^2 \right)\frac{\partial}{\partial w_0}, h+ 2f_0 w_0\right\} \,\,\,\,\,\,\,\,\,\,\,\,  \textnormal{for}\,\,\, m=0& \\
\left\{   \left(a_0+ a_1 z_0 + a_2z_0^2    \right) \frac{\partial}{\partial z_0} + \left( \left(-ma_2z_0+d\right)w_0 + \left(c_0 + c_1z_0+ \cdots + c_m z_0^m \right)w_0^2 \right)\frac{\partial}{\partial w_0} , \right.\\
\left. \,\,\,\,\,\,\,  b- m a_2 z_0 +2 \left( c_0 + c_1 z_0 + \cdots + c_m z_0^m \right) w_0\right\}\,\,\, \,\,\,\,\,\,\,\,\,\,\,\,\,\textnormal{for}\,\,\, m \geq 1  &  \end{cases}\\
H^0\left( \mathbb{F}_m, \Theta_{\mathbb{F}_m}(-m, -2) \right)&= \left\{ s \frac{\partial}{\partial w_0} | s\in \mathbb{C} \right\}\\
H^1\left(\mathbb{F}_m, \mathcal{E}_{\mathcal{O}_{\mathbb{F}_m}(-m, -2)} \right) &= \begin{cases} 0 & \textnormal{for}\,\,\,m=0,1 \\
\textnormal{Span}_\mathbb{C}\left\langle \left(\frac{1}{z_0}\frac{\partial}{\partial w_0}, 0 \right), \left( \frac{1}{z_0^2}\frac{\partial}{\partial w_0}, 0 \right), \cdots , \left( \frac{1}{z_0^{m-1}}\frac{\partial}{\partial w_0}, 0 \right) \right\rangle & \textnormal{for}\,\,\, m\geq 2 \end{cases}
\end{align*}
Then we compute $H^0\left( \mathbb{F}_m, \mathcal{E}_{\mathcal{O}_{\mathbb{F}_m}(-m, -2) } \right)\xrightarrow{\tilde{D}_0} H^0\left( \mathbb{F}_m, \Theta_{\mathbb{F}_m}(-m, -2) \right)$. 
{\small{\begin{align*}
&\mathcal{L}_{ \left(a_0z_0^2+bz_0+c\right)\frac{\partial}{\partial z_0}+ \left(d+ ew_0+ fw_0^2 \right)\frac{\partial}{\partial w_0} } \left( a \frac{\partial}{\partial w_0}\right)+ a \left( h +  2 f w_0 \right)\frac{\partial}{\partial w_0} = - a\left( e+ 2f w_0 \right)\frac{\partial}{\partial w_0} + \left(ah+ 2af w_0 \right)\frac{\partial}{\partial w_0}= \left( ah- ae \right)\frac{\partial}{\partial w_0}\\
&\mathcal{L}_{\left(a_0+ a_1 z_0 + a_2z_0^2    \right) \frac{\partial}{\partial z_0} + \left( \left(-ma_2z_0+d\right)w_0 + \left(c_0 + c_1z_0+ \cdots + c_m z_0^m \right)w_0^2 \right)\frac{\partial}{\partial w_0}}\left(a \frac{\partial}{\partial w_0} \right) + a\left(  b- m a_2 z_0 +2 \left( c_0 + c_1 z_0 + \cdots + c_m z_0^m \right) w_0 \right)\frac{\partial}{\partial w_0}\\
&= - a\left( - ma_2 z_0 + d + 2\left( c_0 + c_1 z_0 + \cdots + c_m z_0^m \right)w_0    \right)\frac{\partial}{ \partial w_0} +\left(ab - maa_2 z_0+ 2a\left( c_0 + c_1 z_0 + \cdots + c_m z_0^m \right) w_0 \right)\frac{\partial}{\partial w_0}= (ab- ad)\frac{\partial}{\partial w_0}
\end{align*}}}
This implies that $\textnormal{coker}\left(H^0\left( \mathbb{F}_m , \mathcal{E}_{\mathcal{O}_{\mathbb{F}_m}(-m,-2)} \right) \to H^0\left( \mathbb{F}_m, \Theta_{\mathbb{F}_m}(-m, -2) \right) \right)=0$. On the other hand, let us compute the image of $H^1\left( \mathbb{F}_m, \mathcal{E}_{\mathcal{O}_{\mathbb{F}_m}(-m,-2)}\right)\to H^1\left( \mathbb{F}_m, \Theta_{\mathbb{F}_m}(-m, -2) \right)$ which is given by
\begin{align*}
\mathcal{L}_{k_1\frac{1}{z_0}\frac{\partial}{\partial w_0} + \cdots + k_{m-1}\frac{1}{z_0^{m-1}}\frac{\partial}{\partial w_0}}\left(a \frac{\partial}{\partial w_0} \right)  + 0\cdot \left( a \frac{\partial}{\partial w_0}  \right) = 0
\end{align*}
Hence we have
\begin{align*}
\mathbb{H}^1\left( \mathbb{F}_m, \mathcal{E}_{\mathcal{O}_{\mathbb{F}_m}(-m,-2)}^\bullet \right)= H^1\left(\mathbb{F}_m, \Theta_{\mathbb{F}_m} \right) =\textnormal{Span}_\mathbb{C}\left\langle  z_1 \frac{\partial}{\partial w_0},..., z_1^{m-1}\frac{\partial}{\partial w_0} \right\rangle
\end{align*}
We will construct a foliated analytic family of deformations of $\left( \mathbb{F}_m, T_0 = a \frac{\partial}{\partial w_0}\right),a\ne 0$ such that the foliated Kodaira-Spencer map is an isomorphism at the distinguished point, so that $\left( \mathbb{F}_m, T_0 \right)$ is unobstructed in foliated deformations. The foliated analytic family is given by
\begin{align*}
T &= a\frac{\partial}{\partial w_0}\\
z_0&=\frac{1}{z_1}\\
w_0&= z_1^m w_1 + t_1z_1 + \cdots + t_{m-1} z_1^{m-1}
\end{align*}
We describe the foliated analytic family more precisely. We have coordinate transformations as follows$:$
\begin{center}
$\begin{CD}
U_0: \left( z_0, w_0 \right) @>>> U_1:\left(z_1, w_1 \right)=\left( \frac{1}{z_0}, z_0^m w_0- t_1 z_0^{m-1}- \cdots - t_{m-1} z_0 \right)\\
@VVV @VVV \\
U_2:\left( z_2, w_2 \right)= \left( z_0, \frac{1}{w_0} \right) @>>> U_3: \left( z_3 , w_3 \right) = \left( \frac{1}{z_2} , \frac{w_2}{ z_2^m - t_1 z_2^{m-1}w_2 - \cdots - t_{m-1} z_2 w_2} \right)  =\left( z_1, \frac{1}{w_1} \right)
\end{CD}$
\end{center}

and the foliation $T$ is given by
\begin{align*}
T&= a\frac{\partial}{\partial w_0}\,\,\,\,\,\,\,\textnormal{on}\,\,\, U_0,\,\,\,\,\,\,\,\,\,\,\,\,\,T= a\frac{\partial}{\partial w_1}\,\,\,\,\,\,\,\textnormal{on}\,\,\, U_1 \\
T&= -a\frac{\partial}{\partial w_2}\,\,\,\,\,\,\,\textnormal{on}\,\,\, U_2,\,\,\,\,\,\,\,\,\,\,\,\,\,T= -a\frac{\partial}{\partial w_3}\,\,\,\,\,\,\,\textnormal{on}\,\,\, U_3 
\end{align*}

and transition functions are given by
\begin{align*}
h_{01}&= z_0^m\\
h_{02}&= w_0^{-2}\\
h_{03}&=\frac{z_0^m}{\left( z_0^m w_0 - t_1 z_0^{m-1} - \cdots - t_{m-1} z_0 \right)^2}=  z_3^{-m} w_3^2\\
h_{12}&= z_2^{-m}w_2^2\\
h_{13}&= w_1^{-2} \\
h_{23}&= \frac{z_2^m}{\left( z_2^m - t_1 z_2^{m-1}w_2 - \cdots - t_{m-1} z_2 w_2 \right)^2}
\end{align*}

\end{remark}

Next we will study foliated deformations of $\mathbb{F}_m$ in terms of cotangent sheaves. A foliation on $\mathbb{F}_m$ is given by  a section $\omega\in H^0\left( \mathbb{F}_m, \Omega_{\mathbb{F}_m}^1(p,q) \right)$, where $\Omega_{\mathbb{F}_m}^1(p,q)=\Omega_{\mathbb{F}_m}^1 \otimes \mathcal{O}_{\mathbb{F}_m}(p,q)$, so that we have the exact sequence
\begin{align*}
0\to \mathcal{O}_{\mathbb{F}_m}( -p,-q) \xrightarrow{\omega } \Omega_{\mathbb{F}_m}^1 \to \frac{\Omega_{\mathbb{F}_m}^1}{\mathcal{O}_{\mathbb{F}_m}(-p,-q)} \to 0
\end{align*}
We will study foliated deformations of $\left( \mathbb{F}_m, \omega \right)$ in terms of cotangent sheaves. We note that we have the dual leaf complex
\begin{align}\label{rs18}
\mathcal{E}_{\mathcal{O}_{\mathbb{F}_m}(p,q)}^\bullet : \mathcal{E}_{\mathcal{O}_{\mathbb{F}_m}(p,q)} \to \Omega_{\mathbb{F}_m}^1(p, q) \to 0
\end{align}
In view of subsection \ref{rs3} and subsection \ref{rs34}, in order to compute $\mathbb{H}^i\left( \mathbb{F}_m, \mathcal{E}_{\mathcal{O}_{\mathbb{F}_m}(p,q)}^\bullet \right)$ for $i=0,1,2$, it remains to describe cohomology group $H^i\left( \mathbb{F}_m , \Omega_{\mathbb{F}_m}^1 ( p, q ) \right),i=0,1$.

\subsection{Description of $H^0\left( \mathbb{F}_m, \Omega_{\mathbb{F}_m}^1(p,q) \right)$}\label{rf17}\

We keep the notations in subsection \ref{rs1}. We shall explicitly describe the global section $H^0\left( \mathbb{F}_m, \Omega_{\mathbb{F}_m}^1(p,q) \right)$ on $U_0$. Let
\begin{align*}
\omega= a_0^1(z_0,w_0) dz_0 + a_0^2(z_0,w_0)dw_0 \in H^0\left(\mathbb{F}_m, \Omega_{\mathbb{F}_m}^1(p,q) \right) \,\,\,\,\,\textnormal{on}\,\,\, U_0
\end{align*}
with $a_0^1(z_0, w_0)= \sum_{h,k\geq 0} a_{hk}z_0^h w_0^k,a_{hk}\in \mathbb{C}$ and $a_0^2(z_0,w_0)= \sum_{h,k\geq 0} b_{hk}z_0^h w_0^k, b_{hk}\in \mathbb{C}$. We note that from $z_0= z_2$ and $w_2= \frac{1}{w_0}$, we have $dz_0= dz_2$ and $dw_0=- \frac{1}{w_2^2}dw_2$. Then $\omega$ on $U_0$ is transformed on $U_2$ into
\begin{align*}
&a_2^1(z_2,w_2)dz_2 + a_2^2(z_2,w_2)dw_2 = h_{20}\cdot\left(a_0^1(z_0,w_0) dz_0 + a_0^2(z_0,w_0)dw_0 \right) \\
&=w_2^q  a_0^1\left( z_2, \frac{1}{w_2} \right) dz_2 - w_2 ^{q-2}a_0^2\left( z_2, \frac{1}{w_2} \right) dw_2 = \left( \sum_{h,k\geq 0} a_{hk} z_2^h w_2^{q-k} \right) dz_2 + \left(\sum_{h,k\geq 0} b_{hk} z_2^h w_2^{q-2 -k}  \right) dw_2 \notag
\end{align*}
which should be holomorphic, so that $a_{hk}=0$ for $k >q$ and $b_{hk}=0$ for $k > q-2$. Hence $\omega$ is of the form
\begin{align}\label{rs35}
&a_0^1(z_0,w_0)dz_0+a_0^2(z_0,w_0)dw_0 \\
&= \left(c_0^1(z_0)+c_0^2(z_0)w_0+\cdots + c_0^{q+1}(z_0) w_0^q\right)dz_0 + \left(d_0^1(z_0)+d_0^2(z_0)w_0+\cdots  + d_0^{q-1}(z_0)w_0^{q-2}\right) dw_0 \notag
\end{align}
where $c_0^1(z_0),..., c_0^{q+1}(z_0),d_0^1(z_0),..., d_0^{q-1}(z_0)$ are entire functions of $z_0$. In particular, we note that we have
\begin{align*}
H^0\left( \mathbb{F}_m, \Omega_{\mathbb{F}_m}^1(p,q) \right)= 0 \,\,\,\,\,\,\,\textnormal{for}\,\,\, q <0
\end{align*}
Hence there is no foliation on $\mathbb{F}_m$ if $q<0$, so that in the following we assume that $q\geq 0$. In a similar way with $U_1$ and $U_3$, we see that $\omega\in H^0\left(\mathbb{F}_m, \Omega_{\mathbb{F}_m}^1 (p,q) \right)$ on $U_1$ is of the form
\begin{align}\label{rs42}
&a_1^1(z_1,w_1)dz_1+a_1^2(z_1,w_1)dw_1 \\
&= \left(c_1^1(z_1)+c_1^2(z_1)w_1+\cdots + c_1^{q+1}(z_1) w_1^q\right)dz_1 + \left(d_1^1(z_1)+d_1^2(z_1)w_1+\cdots  + d_1^{q-1}(z_1)w_1^{q-2}\right) dw_1 \notag
\end{align}
where $c_1^1(z_1),..., c_1^{q+1}(z_1),d_1^1(z_1),..., d_1^{q-1}(z_1)$ are entire functions of $z_1$.

On the other hand, $\omega$ on $U_0$ in $(\ref{rs35})$ is transformed on $U_1$ into
\begin{align*}
&z_0^{mq-p} \cdot \left(\left(c_0^1(z_0)+c_0^2(z_0)w_0+\cdots + c_0^{q+1}(z_0) w_0^q\right)dz_0 + \left(d_0^1(z_0)+d_0^2(z_0)w_0+\cdots  + d_0^{q-1}(z_0)w_0^{q-2}\right) dw_0 \right) \\
&=\left(c_1^1(z_0)+c_1^2(z_1)w_1+\cdots + c_1^{q+1}(z_1) w_1^q\right)dz_1 + \left(d_1^1(z_1)+d_1^2(z_1)w_1+\cdots  + d_1^{q-1}(z_1)w_1^{q-2}\right) dw_1
\end{align*}
From $z_0= \frac{1}{z_1}$ and $w_0= z_1^m w_1$, we have $dz_0=-\frac{1}{z_1^2}dz_1$ and $ dw_0= mz_1^{m-1}w_1dz_1+   z_1^m dw_1$. Then we have
\begin{align}
&\left(c_0^1(z_0) + c_0^2(z_0)w_0+ \cdots + c_0^{q+1}(z_0) w_0^q \right) dz_0 + \left( d_0^1(z_0) + d_0^2(z_0)w_0 + \cdots + d_0^{q-1}(z_0)w_0^{q-2} \right) dw_0 \label{rs36}\\
&= z_0^{p-mq} \left(c_1^1\left(\frac{1}{z_0}\right) + c_1^2\left(\frac{1}{z_0}\right) z_0^m w_0 + \cdots + c_1^{q+1}\left(\frac{1}{z_0}\right)z_0^{mq}w_0^q \right)\left( -\frac{1}{z_0^2}dz_0\right) \notag \\
& +z_0^{p-mq} \left( d_1^1\left( \frac{1}{z_0}\right) + d_1^2\left(\frac{1}{z_0}\right)z_0^m w_0 + \cdots + d_1^{q-1}\left( \frac{1}{z_0}\right)z_0^{m(q-2)}w_0^{q-2} \right) \left( mz_0^{m-1} w_0 dz_0 + z_0^m dw_0\right) \notag \\
&=- z_0^{p-mq-2}\left(c_1^1\left(\frac{1}{z_0}\right) + c_1^2\left(\frac{1}{z_0}\right) z_0^m w_0 + \cdots + c_1^{q+1}\left(\frac{1}{z_0}\right)z_0^{mq}w_0^q  \right) dz_0 \notag\\
&+ mz_0^{p-mq +m-1}\left(  d_1^1\left( \frac{1}{z_0}\right)w_0 + d_1^2\left(\frac{1}{z_0}\right)z_0^m w_0^2 + \cdots + d_1^{q-1}\left( \frac{1}{z_0}\right)z_0^{m(q-2)}w_0^{q-1}\right) dz_0 \notag \\
&+ z_0^{p-mq + m}\left( d_1^1\left( \frac{1}{z_0}\right) + d_1^2\left(\frac{1}{z_0}\right)z_0^m w_0 + \cdots + d_1^{q-1}\left( \frac{1}{z_0}\right)z_0^{m(q-2)}w_0^{q-2} \right) dw_0 \notag
\end{align}

By considering the coefficient of $dw_1$ in $(\ref{rs36})$, we have
\begin{align}
d_1^1(z_1)&=z_1^{-mq+p+m}d_0^1\left(\frac{1}{z_1}\right) \notag \\
d_1^2(z_1)&=z_1^{-mq+p+m+m}d_1^2\left( \frac{1}{z_1}\right) \notag \\
&\cdots \notag \\
d_1^i(z_1)&=z_1^{-mq+p+m+m(i-1)}d_0^i\left(\frac{1}{z_1} \right)=z_1^{p-m(q-i)} d_0^i\left( \frac{1}{z_1} \right),\,\,\,\,\,\,\,i=1,...,q-1 \label{rs38} \\
&\cdots \notag \\
d_1^{q-1}(z_1)&=z_1^{-mq+p+m+m(q-2)}d_0^{q-1}\left(\frac{1}{z_1} \right)=z_1^{p-m}d_0^{q-1}\left(\frac{1}{z_1} \right) \notag
\end{align}
On the other hand, by considering the coefficient of $dz_1$ in $(\ref{rs36})$, we have
\begin{align}
c_1^1(z_1)&= -z_1^{-mq+p-2}c_0^1\left(\frac{1}{z_1} \right) \notag \\
c_1^2(z_1)&= -z_1^{-mq+p-2+m}c_0^2\left(\frac{1}{z_1} \right) +mz_1^{-mq+p+m-1}d_0^1\left(\frac{1}{z_1} \right) \notag \\
&\cdots \notag \\
c_1^i(z_1)&=-z_1^{-mq+p-2+m(i-1)}c_0^i\left( \frac{1}{z_1} \right) + m z_1^{-mq+p+m-1+m(i-2)}d_0^{i-1}\left( \frac{1}{z_1}\right) \label{rs37} \\
                &=-z_1^{p-2-m(q-i+1)}c_0^i\left(\frac{1}{z_1} \right)+mz_1^{p-m-1-m(q-i)}d_0^{i-1}\left(  \frac{1}{z_1} \right),\,\,\,\,\,\,\,\, 2\leq i\leq q \notag \\
&\cdots  \notag\\
c_1^q(z_1)&=- z_1^{-mq+p-2+m(q-1)}c_0^q\left(\frac{1}{z_1} \right) + mz_1^{-mq+p+m-1+m(q-2)} d_0^{q-1}\left(  \frac{1}{z_1}\right) \notag \\
                 &=-z_1^{p-2-m}c_0^q\left(\frac{1}{z_1} \right)  + mz_1^{p-m-1}d_0^{q-1}\left(  \frac{1}{z_1}\right)  \notag\\
c_1^{q+1}(z_1)&=-z_1^{p-2} c_0^{q+1}\left( \frac{1}{z_1} \right) \notag
\end{align}

In the following examples, we will consider the case $q=2$, i.e. we will describe $H^0\left( \mathbb{F}_m, \Omega_{\mathbb{F}_m}^1(p,2) \right)$. Then from $(\ref{rs38})$ and $(\ref{rs37})$, we have
\begin{align}\label{rs39}
d_1^1(z_1)&= z_1^{-m+p} d_0^1\left(\frac{1}{z_1} \right) \notag \\
c_1^1(z_1)&= - z_1^{-2m+p-2} c_0^1\left(\frac{1}{z_1} \right)\\
c_1^2(z_1)&= - z_1^{p-2-m}c_0^2\left(\frac{1}{z_1} \right) + m z_1^{p-m-1}d_0^{1}\left(\frac{1}{z_1} \right) \notag \\
c_1^3(z_1)&= - z_1^{p-2} c_0^3\left(\frac{1}{z_1} \right) \notag
\end{align}

\begin{example}
Let us describe $H^0\left( \mathbb{F}_m , \Omega_{\mathbb{F}_m}^1(m+1,2) \right)$. Then from $(\ref{rs39})$, we have
\begin{align*}
d_0^1(z_0)&= a_0 + a_1 z_0\\
c_0^1(z_0)&= 0 \\
c_1^2(z_0)&= - \frac{1}{z_1} c_0^2\left( \frac{1}{z_1} \right) + ma_0+ ma_1 \frac{1}{z_1}\Longrightarrow c_0^2(z_0)= m a_1\\
c_0^2(z_0)&= \begin{cases} 0 & \textnormal{if}\,\,\,m=0\\
                     b_0+ b_1 z_0+ \cdots + b_{m-1}z_0^{m-1} & \textnormal{if}\,\,\,m \geq 1\end{cases}
\end{align*}
Then from $(\ref{rs35})$ any element in $H^0\left( \mathbb{F}_m, \Omega_{\mathbb{F}_m}^1 (m+1,2) \right)$ on $U_0$ is of the form
\begin{align*}
\left( a_0 + a_1 z_0 \right)dw_0\,\,\,\,\,\textnormal{for}\,\,\,m=0
\end{align*}
\begin{align}\label{rsf17}
\left(ma_1 w_0+ \left(b_0+ b_1 z_0 + \cdots + b_{m-1} z_0^{m-1} \right) w_0^2 \right)dz_0 + \left(a_0+ a_1 z_0 \right)dw_0
\end{align}
where $a_0,a_1,b_0,...,b_{m-1}\in \mathbb{C}$, so that we have $\dim_\mathbb{C} H^0\left(\mathbb{F}_m, \Omega_{\mathbb{F}_m}^1(m+1,2) \right)= m+2$. We also note that since $K_{\mathbb{F}_m}\cong \mathcal{O}_{\mathbb{F}_m}(-m-2,-2)$ and $\Omega_{\mathbb{F}_m}^1\cong \Theta_{\mathbb{F}_m}\otimes K_{\mathbb{F}_m}$, we have $H^1\left( \mathbb{F}_m, \Omega_{\mathbb{F}_m}^1(m+1,2) \right) \cong H^1\left( \mathbb{F}_m, \Theta_{\mathbb{F}_m}(-1,0) \right)$ which is described in \textnormal{Example \ref{rs54}}.
\end{example}

\begin{example}
Let us describe $H^0\left(  \mathbb{F}_m , \Omega_{\mathbb{F}_m}^1(m+2,2) \right)$. From $(\ref{rs39})$, we have
\begin{align*}
d_0^1(z_0)&= a_0 + a_1 z_0 + a_2 z_0^2\\
c_0^1(z_0)&=\begin{cases} e_0 & \textnormal{if}\,\,\,m=0\\
                                  0 & \textnormal{if}\,\,\,m\geq 1 \end{cases}\\
            c_1^2(z_1)     &=         - c_0^2\left( \frac{1}{z_1} \right)  + mz_1\left( a_0 + a_1 \frac{1}{z_1} + a_2 \frac{1}{z_1^2} \right)        \Longrightarrow c_0^2(z_0)= c_0+ ma_2 z_0\\
c_0^3(z_0)&= b_0 + b_1 z_0 + \cdots + b_m z_0^m            
\end{align*}
Then from $(\ref{rs35})$ any element in $H^0\left( \mathbb{F}_m, \Omega_{\mathbb{F}_m}^1(m+2,2) \right)$ on $U_0$ is of the form
\begin{align*}
\left( e_0 + c_0 w_0+ b_0 w_0^2 \right)dz_0 + \left(a_0+ a_1 z_0+ a_2 z_0^2 \right) dw_0\,\,\,\,\,\textnormal{for}\,\,\,m=0
\end{align*}
\begin{align} \label{rsf31}
\left( \left( c_0 + ma_2z_0\right) w_0 +\left( b_0+ b_1 z_0 + \cdots + b_m z_0^m\right) w_0^2 \right) dz_0 + \left(a_0 + a_1 z_0 + a_2 z_0^2 \right) dw_0\,\,\,\,\,\,\textnormal{for}\,\,\,m \geq 1
\end{align}
where $a_0,a_1,a_2,c_0, b_0,...,b_m\in \mathbb{C}$, so that we have $\dim_\mathbb{C}H^0\left(\mathbb{F}_0, \Omega_{\mathbb{F}_0}(2,2) \right)=6$, and $\dim_\mathbb{C} H^0\left( \mathbb{F}_m, \Omega_{\mathbb{F}_m}^1(m+2,2) \right)=m+5$ for $m\geq 1$. We also note that since $K_{\mathbb{F}_m}\cong \mathcal{O}_{\mathbb{F}_m}(-m-2,-2)$ and $\Omega_{\mathbb{F}_m}^1\cong \Theta_{\mathbb{F}_m}\otimes K_{\mathbb{F}_m}$, we have $H^1\left( \mathbb{F}_m, \Omega_{\mathbb{F}_m}^1(m+2,2) \right) \cong H^1\left( \mathbb{F}_m, \Theta_{\mathbb{F}_m}(0,0) \right)$ which is described in \textnormal{Example \ref{rr1}}.
\end{example}

\begin{example}
Let us describe $H^0\left(  \mathbb{F}_m , \Omega_{\mathbb{F}_m}^1(m+3,2) \right)$. From $(\ref{rs39})$, we have
\begin{align*}
d_0^1(z_0)&= a_0 + a_1 z_0 + a_2 z_0^2 + a_3 z_0^3\\
c_0^1(z_0)&= \begin{cases}  e_0 + e_1 z_0 & \textnormal{if}\,\,\, m=0\\
                                                 e_0 & \textnormal{if}\,\,\,m=1 \\
                                                  0 & \textnormal{if} \,\,\, m\geq 2 \end{cases} \\
c_1^2(z_1)&=                - z_1 c_0^2\left( \frac{1}{z_1} \right)    + mz_1^2\left(a_0 + a_1 \frac{1}{z_1} + a_2 \frac{1}{z_1^2} + a_3 \frac{1}{z_1^3} \right)     \Longrightarrow c_0^2(z_0)= c_0+ c_1 z_0 + ma_3 z_0^2 \\
c_1^3(z_0)&= b_0 + b_1 z_0 + \cdots + b_{m+1} z_0^{m+1}                          
\end{align*}
Then from $(\ref{rs35})$ any element in $H^0\left( \mathbb{F}_m, \Omega_{\mathbb{F}_m}^1(m+3,2 )\right)$ on $U_0$ is of the form
\begin{align*}
\left( e_0+ e_1 z_0+ \left( c_0 + c_1 z_0\right) w_0 + \left( b_0 + b_1 z_0 \right) w_0^2 \right) dz_0 + \left( a_0 + a_1 z_0 + a_2 z_0^2 + a_3 z_0^3 \right) dw_0 \,\,\,\,\,\textnormal{for}\,\,\, m=0
\end{align*}
\begin{align*}
\left( e_0+ \left( c_0 + c_1 z_0+ a_3 z_0^2\right) w_0 + \left( b_0 + b_1 z_0+ b_2 z_0^2 \right) w_0^2 \right) dz_0 + \left( a_0 + a_1 z_0 + a_2 z_0^2 + a_3 z_0^3 \right) dw_0 \,\,\,\,\,\textnormal{for}\,\,\, m=1
\end{align*}
\begin{align}\label{rf13}
\left( \left( c_0 + c_1 z_0+ a_3 z_0^2\right) w_0 + \left( b_0 + b_1 z_0+\cdots + b_{m+1}z_0^{m+1} \right) w_0^2 \right) dz_0 + \left( a_0 + a_1 z_0 + a_2 z_0^2 + a_3 z_0^3 \right) dw_0 \,\,\,\,\,\textnormal{for}\,\,\, m\geq 2
\end{align}
where $a_0,...,a_3,e_0,e_1,c_0,c_1,b_0,...,b_{m+1}\in \mathbb{C}$, so that $\dim_\mathbb{C} H^0\left(\mathbb{F}_m, \Omega_{\mathbb{F}_m}^1(m+3,2) \right)=10$ for $m=0,1$, and $\dim_\mathbb{C} H^0\left( \mathbb{F}_m, \Omega_{\mathbb{F}_m}^1(m+3,2) \right)= m+8$ for $m\geq 2$. We also note that since $K_{\mathbb{F}_m}\cong \mathcal{O}_{\mathbb{F}_m}(-m-2,-2)$ and $\Omega_{\mathbb{F}_m}^1\cong \Theta_{\mathbb{F}_m}\otimes K_{\mathbb{F}_m}$, we have $H^1\left( \mathbb{F}_m, \Omega_{\mathbb{F}_m}^1(m+3,2) \right) \cong H^1\left( \mathbb{F}_m, \Theta_{\mathbb{F}_m}(1,0) \right)$ which is described in \textnormal{Example \ref{rr2}}.
\end{example}

\subsection{Description of $H^1\left( \mathbb{F}_m, \Omega_{\mathbb{F}_m}^1(p,q) \right)$}\label{rf16}\

We keep the notations in subsection \ref{rs1}. We explicitly describe $H^1\left( \mathbb{F}_m, \Omega_{\mathbb{F}_m}^1(p,q)\right)$. First we show that $H^1\left( U_0 \bigcup U_2, \Omega_{\mathbb{F}_m}^1(p,q) \right)=0$ and $H^1\left( U_1 \bigcup U_3 , \Omega_{\mathbb{F}_m}^1(p,q) \right)=0$ for $q\geq 1$. We note that by letting $\mathscr{B}= \left\{ U_0, U_2 \right\}$, we can compute $H^1\left(U_0 \bigcup U_2, \Omega_{\mathbb{F}_m}^1(p,q) \right)$ by \v Cech resolution associated to the open covering $\mathscr{B}$, i.e.
\begin{align*}
H^1\left(U_0\bigcup U_2, \Omega_{\mathbb{F}_m}^1(p,q) \right) = H^1\left( \mathscr{B}, \Omega_{\mathbb{F}_m}^1(p,q) \right) =\frac{Z^1\left(\mathscr{B}, \Omega_{\mathbb{F}_m}^1(p,q)\right)}{\delta C^0\left(\mathscr{B}, \Omega_{\mathbb{F}_m}^1(p,q)\right)}
\end{align*}
A $1$-cocycle $\left\{\omega_{02},\omega_{20}\right\}$ belonging to $Z^1\left(\mathscr{B}, \Omega_{\mathbb{F}_m}^1(p,q)\right)$ consists of holomorphic forms $\omega_{02}$ and $\omega_{20}$ on $U_0\bigcap U_2$ with $\omega_{02}=-h_{02}\omega_{20}$, while $C^0\left(\mathscr{B},\Omega_{\mathbb{F}_m}^1(p,q)\right)$ is the set of all pairs $\left\{\omega_0,\omega_2 \right\}$, where $\omega_0$ and $\omega_2$ are holomorphic forms on $U_0$ and $U_2$, respectively. In order to show $H^1\left(\mathscr{B}, \Omega_{\mathbb{F}_m}^1(p,q)\right)=0$, it suffices to show that any $1$-cocycle $\left\{ \omega_{02}, \omega_{20} \right\}$ is represented as $\delta\left(\left\{ \omega_2, \omega_0 \right\} \right)$
\begin{align*}
\omega_{02}= h_{02} \omega_2 - \omega_0 = w_0^q \omega_2 - \omega_0
\end{align*}
The holomorphic forms on $U_0\bigcap U_2$ is written in the form
\begin{align*}
\omega_{02}=u_{02}(z_0,w_0)dz_0+v_{02}(z_0,w_0)dw_0
\end{align*}
where $u_{02}(z_0,w_0)$ and $v_{02}(z_0,w_0)$ are holomorphic function on $U_0\bigcap U_2= \mathbb{C}\times \mathbb{C}^*$. Hence we have Laurent series in $w_0$:
\begin{align*}
u_{02}(z_0,w_0)=  \sum_{n= -\infty}^{\infty} u_n(z_0)w_0^n,\,\,\,\,\,\,\,\,\,\,\,v_{02}(z_0,w_0)=  \sum_{n=-\infty}^{\infty} v_n(z_0)w_0^n
\end{align*}
where $u_n(z_0)$ and $v_n(z_0)$ are entire functions of $z_0$. Any holomorphic vector field on $U_0$ is of the form
\begin{align*}
u_0(z_0,w_0) dz_0 + v_1(z_0,w_0)dw_0
\end{align*}
where $u_0(z_0,w_0)$ and $v_1(z_0,w_0)$ are entire functions, and any holomorphic vector field on $U_2$ is of the form
\begin{align*}
u_2(z_2,w_2)dz_2+ v_1(z_2,w_2) dw_2
\end{align*}
where $u_2(z_2,w_2)$ and $v_1(z_2,w_2)$ are entire functions. Since $dz_0=dz_2$ and $dw_2=-\frac{1}{w_0^2}dw_0$, we have
\begin{align*}
 \omega_2= u_2(z_2,w_2)dz_2+v_2(z_2,w_2)dw_2=u_2\left(z_0,\frac{1}{w_0} \right)dz_0-\frac{1}{w_0^2}v_2\left(z_0,\frac{1}{w_0} \right)dw_0,
\end{align*}
so that we have
\begin{align*}
h_{02}\omega_2-\omega_0=w_0^q u_2\left(z_0,\frac{1}{w_0} \right)dz_0-w_0^{q-2} v_2\left(z_0, \frac{1}{w_0}  \right) dw_0 - u_0(z_0,w_0)dz_0-v_0(z_0,w_0)dw_0
\end{align*}
If $\omega_{02}$ defines the zero class, we should have a solution $\left(\omega_0, \omega_2 \right)$ satisfying
\begin{align*}
w_0^qu_2\left(z_0, \frac{1}{w_0} \right)- u_0(z_0,w_0)&= \sum_{n= -\infty}^{\infty} u_n(z_0)w_0^n\\
-w_0^{q-2}v_2\left( z_0, \frac{1}{w_0} \right) -v_0(z_0,w_0)&= \sum_{n=-\infty}^{\infty} v_n(z_0)w_0^n
\end{align*}
This implies that if $q<1$, then we have $H^1\left( U_0 \bigcup U_2, \Omega_{\mathbb{F}_m}^1(p,q) \right)\ne 0$. On the other hand, if $q \geq 1$, then in view of $w_0= \frac{1}{w_2}$, we set
\begin{align*}
u_0 & = -\sum_{n=0}^{\infty} u_n(z_0)w_0^n,\,\,\,\,\,\,\,\,\,\,\,u_2= \sum_{n=1}^\infty u_{-n}(z_2) w_2^{n+q}\\
v_0 & = -\sum_{n=0}^{\infty} v_n(z_0)w_0^n,\,\,\,\,\,\,\,\,\,\,\,v_2= -\sum_{n=1}^\infty v_{-n}(z_2) w_2^{n+q-2}
\end{align*}
Then there is a solution $\left(\omega_0, \omega_2 \right)$ satisfying $h_{02}\omega_2- \omega_0= \omega_{02}$. This implies that 
\begin{align*}
H^1\left( U_0 \bigcup U_2,  \Omega_{\mathbb{F}_m}^1(p,q) \right)=0,\,\,\,\,\,\,\,\,q\geq 1
\end{align*}
Similarly we can show that $H^1\left( U_1 \bigcup U_3, \Omega_{\mathbb{F}_m}^1(p,q) \right)=0$ for $q\geq 1$. In the following we will assume that $q \geq 1$. In the case of $q=0$, we will study foliated deformations of $\left( \mathbb{F}_m, \omega \right)$ where $\omega\in H^0\left( \mathbb{F}_m, \Omega_{\mathbb{F}_m}^1(p,0)\right)$ separately in Remark \ref{rs40} below. Then by letting $\mathfrak{B}=\left\{U_0 \bigcup U_2, U_1 \bigcup U_3 \right\}$, we can compute $H^1\left( \mathbb{F}_m, \Omega_{\mathbb{F}_m}^1 (p,q)\right), q \geq 1$ by \v Cech resolution associated to the open covering $\mathfrak{B}$, i.e.
\begin{align}\label{rs45}
H^1\left( \mathbb{F}_m, \Omega_{\mathbb{F}_m}^1(p,q) \right)= H^1\left( \mathfrak{B}, \Omega_{\mathbb{F}_m}^1 (p,q) \right)= \frac{Z^1\left( \mathfrak{B}, \Omega_{\mathbb{F}_m}^1(p,q) \right)}{\delta C^0\left( \mathfrak{B}, \Omega_{\mathbb{F}_m}^1(p,q) \right)}
\end{align}

Next we describe explicitly $H^0\left( U_0 \bigcup U_2, \Omega_{\mathbb{F}_m}^1(p,q) \right),q \geq 1$ on $U_0$ and $H^0\left( U_1 \bigcup U_3  , \Omega_{\mathbb{F}_m}^1(p,q) \right), q \geq 1$ on $U_1$. Let us describe $H^0\left( U_0 \bigcup U_2 , \Omega_{\mathbb{F}_m}^1(p,q)\right)$ on $U_0$ which is given by $(\ref{rs35})$
\begin{align}\label{rs41}
\left(c_0^1(z_0)+c_0^2(z_0)w_0+\cdots + c_0^{q+1}(z_0) w_0^q\right)dz_0 + \left(d_0^1(z_0)+d_0^2(z_0)w_0+\cdots  + d_0^{q-1}(z_0)w_0^{q-2}\right) dw_0 
\end{align}
where $c_0^1(z_0),..., c_0^{q+1}(z_0),d_0^1(z_0),..., d_0^{q-1}(z_0)$ are entire functions of $z_0$. Let us describe $H^0\left( U_1 \bigcup U_3, \Omega_{\mathbb{F}_m}^1(p,q) \right)$ on $U_1$ which is given by $(\ref{rs42})$
\begin{align}\label{rs43}
\left(c_1^1(z_0)+c_1^2(z_1)w_1+\cdots + c_1^{q+1}(z_1) w_1^q\right)dz_1 + \left(d_1^1(z_1)+d_1^2(z_1)w_1+\cdots  + d_1^{q-1}(z_1)w_1^{q-2}\right) dw_1
\end{align}
where $c_1^1(z_1),..., c_1^{q+1}(z_1),d_1^1(z_1),..., d_1^{q-1}(z_1)$ are entire functions of $z_1$.

Now let us describe $H^1\left( \mathbb{F}_m, \Omega_{\mathbb{F}_m}^1(p,q)\right),q \geq 1$ in terms of \v Cech resolution associated to $\mathfrak{B}$ from $(\ref{rs45})$. We represent a $1$-cocycle $\left\{\omega_{01},\omega_{10} \right\}\in Z^1\left(\mathfrak{B}, \Omega_{\mathbb{F}_m}^1(p,q)\right)$ with $\omega_{01}=-h_{01}\omega_{10}$. Then the $1$-cocycle belongs to $\delta C^0\left( \mathfrak{B}, \Omega_{\mathbb{F}_m}^1(p,q)\right)$ if and only if there exists $\omega_0\in \Gamma\left(U_0\bigcup U_2, \Omega_{\mathbb{F}_m}^1(p,q)\right)$ and $\omega_2 \in \Gamma\left(U_1 \bigcup U_3,\Omega_{\mathbb{F}_m}^1(p,q)\right)$ such that $\omega_{01}= h_{01}\omega_1- \omega_0= z_0^{p-mq}\omega_1- \omega_0$, where $\omega_0$ is written in the form as in $(\ref{rs41})$ and $\omega_1$ is written in the form as in $(\ref{rs43})$. We write $\omega_{01}$ in terms of the coordinate $\left(z_0, w_0 \right)$
\begin{align*}
\omega_{01}=\left(e_{01}^1(z_0)+e_{01}^2(z_0)w_0+\cdots + e_{01}^{q+1}(z_0)w_0^q\right)dz_0+ \left(f_{01}^1(z_0)+ f_{01}^2(z_0)w_0+\cdots +f_{01}^{q-1}(z_0)w_0^{q-2}\right)dw_0
\end{align*}
where $e_{01}^1(z_0),..., e_{01}^{q+1}(z_0),f_{01}^1(z_0),..., f_{01}^{q-1}(z_0)$ are holomorphic functions of $z_0$ on $\mathbb{C}^*$, so that they are expanded into Laurent series in $z_0$. We note that since $dz_1=- \frac{1}{z_0^2}dz_0$ and $dw_1= m z_0^{m-1}w_0 dz_0 + z_0^m dw_0$, we have
{\tiny{\begin{align*}
&h_{01}\cdot \omega_1\\
&=z_0^{p-mq}\left(c_1^1(z_1)+c_1^2(z_1)w_1+\cdots + c_1^{q+1}(z_1)w_1^q\right)dz_1+ z_0^{p-mq}\left(d_1^1(z_1)+ d_1^2(z_1)w_1+\cdots +d_1^{q-1}(z_1)w_1^{q-2}\right)dw_1\\ 
                &=z_0^{p-mq}\left(c_1^1\left(\frac{1}{z_0} \right)  +c_1^2\left(\frac{1}{z_0} \right)z_0^mw_0  + \cdots +c_1^i\left(\frac{1}{z_0} \right)z_0^{m(i-1)}w_0^{i-1}+\cdots + c_1^{q+1}\left( \frac{1}{z_0} \right) z_0^{qm}w_0^q     \right)\left(-z_0^{-2} dz_0\right)\\
                &+ z_0^{p-mq} \left(d_1^1\left(\frac{1}{z_0}\right)+ d_1^2\left(\frac{1}{z_0} \right)z_0^mw_0+\cdots +d_1^i\left(\frac{1}{z_0} \right)z_0^{m(i-1)}w_0^{i-1}+ \cdots + d_1^{q-1}\left(\frac{1}{z_0}\right)z_0^{m(q-2)} w_0^{q-2}\right)\left( mz_0^{m-1}w_0dz_0+z_0^mdw_0 \right)\\
                &=-z_0^{p-mq-2}\left(c_1^1\left(\frac{1}{z_0} \right)  +c_1^2\left(\frac{1}{z_0} \right)z_0^mw_0  + \cdots +c_1^i\left(\frac{1}{z_0} \right)z_0^{m(i-1)}w_0^{i-1}+\cdots + +c_1^q\left(\frac{1}{z_0} \right) z_0^{m(q-1)}w_0^{q-1}+c_1^{q+1}\left( \frac{1}{z_0} \right) z_0^{qm}w_0^q     \right)dz_0\\
                                &+ mz_0^{p-mq+m-1} \left(d_1^1\left(\frac{1}{z_0}\right)w_0+ d_1^2\left(\frac{1}{z_0} \right)z_0^mw_0^2+\cdots +d_1^i\left(\frac{1}{z_0} \right)z_0^{m(i-1)}w_0^{i}+ \cdots+ d_1^{q-2}\left( \frac{1}{z_0} \right)z_0^{m(q-3)}w_0^{q-2} + d_1^{q-1}\left(\frac{1}{z_0}\right)z_0^{m(q-2)} w_0^{q-1}\right)dz_0\\
                                                &+ z_0^{p-mq+m} \left(d_1^1\left(\frac{1}{z_0}\right)+ d_1^2\left(\frac{1}{z_0} \right)z_0^mw_0+\cdots +d_1^i\left(\frac{1}{z_0} \right)z_0^{m(i-1)}w_0^{i-1}+ \cdots + d_1^{q-1}\left(\frac{1}{z_0}\right)z_0^{m(q-2)} w_0^{q-2}\right)dw_0 \\
\end{align*}}}
Then the equation $\omega_{01}= z_0^{p-mq} \omega_1 - \omega_0$ is reduced to the following system of equations:

\begin{align}
z_0^{p-mq+m}d_1^1\left(\frac{1}{z_0} \right) -d_0^1(z_0)&=f_{01}^1(z_0) \notag \\
z_0^{p-mq+2m}d_1^2\left( \frac{1}{z_0} \right)-d_0^2(z_0)&=f_{01}^2(z_0) \notag \\
&\cdots \notag \\
z_0^{p-mq+mi}d_1^i\left( \frac{1}{z_0} \right)-d_0^i(z_0)&=f_{01}^i(z_0),\,\,\,\,\,\,\, 1\leq i \leq q-1 \label{rs46} \\
&\cdots \notag \\
z_0^{p-m}d_1^{q-1}\left(\frac{1}{z_0} \right) - d_0^{q-1}(z_0)&=f_{01}^{q-1}(z_0) \notag
\end{align}

On the other hand,
\begin{align}
-z_0^{p-mq-2} c_1^1\left(\frac{1}{z_0} \right)-c_0^1(z_0)&=e_{01}^1(z_0) \notag \\
-z_0^{p-mq - 2 + m }c_1^2\left( \frac{1}{z_0} \right) +mz_0^{p-mq+m-1}d_1^1\left(\frac{1}{z_0} \right) -c_0^2(z_0) &=  e_{01}^2(z_0)     \notag       \\
&\cdots \notag\\
-z_0^{p-2-m(q+1-i)} c_1^i\left( \frac{1}{z_0} \right) +mz_0^{p-m-1-m(q-i)} d_1^{i-1}\left(\frac{1}{z_0} \right)  - c_0^i(z_0) &= e_{01}^i(z_0),\,\,\,\,\,\,\, 2\leq i\leq q  \label{rs47}. \\
&\cdots \notag \\
-z_0^{p-2-2m}c_1^{q-1}\left( \frac{1}{z_0} \right)+m z_0^{p-2m-1}d_1^{q-2}\left(\frac{1}{z_0} \right) &= e_{01}^{q-1}(z_0) \notag \\ 
-z_0^{p-2-m} c_1^q\left(\frac{1}{z_0} \right)+mz_0^{p-m-1}d_1^{q-1}\left(\frac{1}{z_0} \right) -c_0^q(z_0) &= e_{01}^q(z_0) \notag \\
 -z_0^{p -  2} c_1^{q+1}\left(\frac{1}{z_0} \right)  - c_0^{q+1}(z_0)&= e_{01}^{q+1}(z_0) \notag
\end{align}

In the following examples, we will consider $q=2$, i.e. we will describe $H^1\left(\mathbb{F}_m, \Omega_{\mathbb{F}_m}^1(p,2) \right)$. Then $\omega_0$ is written in the form
\begin{align*}
\left(c_0^1(z_0) + c_0^2(z_0) w_0 + c_0^3(z_0)w_0^2 \right) dz_0 + d_0^1(z_0) dw_0
\end{align*}
and $\omega_1$ is written in the form
\begin{align*}
\left( c_1^1(z_1)+ c_1^2(z_1)w_1 + c_1^3(z_1)w_1^2 \right) dz_1 + d_1^1(z_1) d w_1
\end{align*}
and from $(\ref{rs46})$ and $(\ref{rs47})$, we have
\begin{align}
z_0^{p-m}d_1^1\left( \frac{1}{z_0} \right)- d_0^1(z_0)&= f_{01}^1(z_0)\notag\\
- z_0^{p-2m-2} c_1^1\left( \frac{1}{z_0} \right)- c_0^1(z_0)&= e_{01}^1(z_0) \label{rs48}\\
-z_0^{p-m-2}c_1^2\left( \frac{1}{z_0} \right) +m z_0^{p-m-1} d_1^1\left( \frac{1}{z_0} \right)- c_0^2(z_0)&= e_{01}^2(z_0)\notag \\
- z_0^{p-2} c_1^3\left( \frac{1}{z_0} \right) - c_0^3(z_0) &= e_{01}^3(z_0) \notag
\end{align}

\begin{example}
Let us describe $H^1\left( \mathbb{F}_m, \Omega_{\mathbb{F}_m}^1(m+1,2) \right)$. From $(\ref{rs48})$, we have
\begin{align*}
z_0d_1^1\left( \frac{1}{z_0} \right)- d_0^1(z_0)&= f_{01}^1(z_0)\\
- \frac{1}{z_0^{m+1}} c_1^1\left( \frac{1}{z_0} \right)- c_0^1(z_0)&= e_{01}^1(z_0)\\
-\frac{1}{z_0}c_1^2\left( \frac{1}{z_0} \right) +m d_1^1\left( \frac{1}{z_0} \right)- c_0^2(z_0)&= e_{01}^2(z_0) \\
- z_0^{m-1} c_1^3\left( \frac{1}{z_0} \right) - c_0^3(z_0) &= e_{01}^3(z_0) 
\end{align*}
In case of $m=0$, these equations always have a solution. For $m\geq 1$, let $e_{01}^1(z_0)= \sum_{n=- \infty}^{\infty} c_n z_1^n$ be the Laurent explansion of $\eta_{01}^1(z_0)$. Then the equations have a solution if and only if $c_{-1}= \cdots = c_{-m}=0$. Hence $\dim_\mathbb{C} H^1\left( \mathbb{F}_0, \Omega_{\mathbb{F}_0}^1(m+1,2) \right)=0$, and $\dim_{\mathbb{C}}H^1\left(\mathbb{F}_m, \Omega_{\mathbb{F}_m}^1(m+1, 2) \right)=m$ for $m\geq 1$, and the $1$-cocycles 
\begin{align}\label{rsf18}
\omega_{01}=  \frac{1}{z_0}dz_0,....., \frac{1}{z_0^{m}}dz_0 \in Z^1\left( \mathfrak{B}, \Omega_{\mathbb{F}_m}^1(m+1,2) \right)
\end{align}
forms a basis of $H^1\left( \mathbb{F}_m, \Omega_{\mathbb{F}_m}^1(m+1, 2) \right)$. We also note that since $K_{\mathbb{F}_m}\cong \mathcal{O}_{\mathbb{F}_m}(-m-2,-2)$ and $\Omega_{\mathbb{F}_m}^1\cong \Theta_{\mathbb{F}_m}\otimes K_{\mathbb{F}_m}$, we have $H^1\left( \mathbb{F}_m, \Omega_{\mathbb{F}_m}^1(m+1,2) \right) \cong H^1\left( \mathbb{F}_m, \Theta_{\mathbb{F}_m}(-1,0) \right)$ which is described in \textnormal{Example \ref{rs49}}.
\end{example}

\begin{example}
Let us describe $H^1\left(  \mathbb{F}_m, \Omega_{\mathbb{F}_m}^1(m+2,2) \right)$. From $(\ref{rs48})$, we have
\begin{align*}
z_0^2d_1^1\left( \frac{1}{z_0} \right)- d_0^1(z_0)&= f_{01}^1(z_0)\\
- \frac{1}{z_0^{m}} c_1^1\left( \frac{1}{z_0} \right)- c_0^1(z_0)&= e_{01}^1(z_0)\\
-c_1^2\left( \frac{1}{z_0} \right) +m z_0d_1^1\left( \frac{1}{z_0} \right)- c_0^2(z_0)&= e_{01}^2(z_0) \\
- z_0^{m} c_1^3\left( \frac{1}{z_0} \right) - c_0^3(z_0) &= e_{01}^3(z_0) 
\end{align*}
In case of $m=0,1$, these equations always have a solution. For $m\geq 2$, let $e_{01}^1(z_0)= \sum_{n=- \infty}^{\infty} c_n z_1^n$ be the Laurent explansion of $\eta_{01}^1(z_0)$. Then the equations have a solution if and only if $c_{-1}= \cdots = c_{-(m-1)}=0$. Hence $ H^1\left( \mathbb{F}_0, \Omega_{\mathbb{F}_0}^1(2,2) \right)=H^1\left(\mathbb{F}_1, \Omega_{\mathbb{F}_1}^1 (3,2) \right)=0$, and $\dim_{\mathbb{C}}H^1\left(\mathbb{F}_m, \Omega_{\mathbb{F}_m}^1 (m+2, 2) \right)=m-1$ for $m\geq 2$, and the $1$-cocycles 
\begin{align} \label{rsf32}
\theta_{01}=  \frac{1}{z_0}dz_0,....., \frac{1}{z_0^{m-1}} dz_0 \in Z^1\left( \mathfrak{B}, \Omega_{\mathbb{F}_m}^1(m+2,2) \right)
\end{align}
forms a basis of $H^1\left( \mathbb{F}_m, \Omega_{\mathbb{F}_m}^1(m+2, 2) \right)$. We also note that since $K_{\mathbb{F}_m}\cong \mathcal{O}_{\mathbb{F}_m}(-m-2,-2)$ and $\Omega_{\mathbb{F}_m}^1\cong \Theta_{\mathbb{F}_m}\otimes K_{\mathbb{F}_m}$, we have $H^1\left( \mathbb{F}_m, \Omega_{\mathbb{F}_m}^1(m+2,2) \right) \cong H^1\left( \mathbb{F}_m, \Theta_{\mathbb{F}_m}(0,0) \right)$ which is described in \textnormal{Example \ref{rs50}}.
\end{example}

\begin{example}
Let us describe $H^1\left( \mathbb{F}_m, \Omega_{\mathbb{F}_m}^1( m+3, 2) \right)$. From $(\ref{rs48})$, we have
\begin{align*}
z_0^3 d_1^1\left( \frac{1}{z_0} \right)- d_0^1(z_0)&= f_{01}^1(z_0)\\
- \frac{1}{z_0^{m-1}} c_1^1\left( \frac{1}{z_0} \right)- c_0^1(z_0)&= e_{01}^1(z_0)\\
-z_0c_1^2\left( \frac{1}{z_0} \right) +m z_0^2d_1^1\left( \frac{1}{z_0} \right)- c_0^2(z_0)&= e_{01}^2(z_0) \\
- z_0^{m+1} c_1^3\left( \frac{1}{z_0} \right) - c_0^3(z_0) &= e_{01}^3(z_0) 
\end{align*}
In case of $m=0,1,2$, these equations always have a solution. For $m\geq 2$, let $e_{01}^1(z_0)= \sum_{n=- \infty}^{\infty} c_n z_1^n$ be the Laurent explansion of $\eta_{01}^1(z_0)$. Then the equations have a solution if and only if $c_{-1}= \cdots = c_{-(m-2)}=0$. Hence $ H^1\left( \mathbb{F}_m, \Omega_{\mathbb{F}_m}^1(m+3,2) \right)=0$ for $m=0,1,2$, and $\dim_{\mathbb{C}}H^1\left(\mathbb{F}_m, \Omega_{\mathbb{F}_m}^1(m+2, 2) \right)=m-1$ for $m\geq 3$, and the $1$-cocycles 
\begin{align}  \label{rf15}
\theta_{01}=  \frac{1}{z_0}dz_0,....., \frac{1}{z_0^{m-2}} dz_0 \in Z^1\left( \mathfrak{B}, \Omega_{\mathbb{F}_m}^1(m+3,2) \right)
\end{align}
forms a basis of $H^1\left( \mathbb{F}_m, \Omega_{\mathbb{F}_m}^1(m+3, 2) \right)$. We also note that since $K_{\mathbb{F}_m}\cong \mathcal{O}_{\mathbb{F}_m}(-m-2,-2)$ and $\Omega_{\mathbb{F}_m}^1\cong \Theta_{\mathbb{F}_m}\otimes K_{\mathbb{F}_m}$, we have $H^1\left( \mathbb{F}_m, \Omega_{\mathbb{F}_m}^1(m+3,2) \right) \cong H^1\left( \mathbb{F}_m, \Theta_{\mathbb{F}_m}(1,0) \right)$ which is described in \textnormal{Example \ref{rs51}}.
\end{example}

With this preparation we study foliated deformations of $\mathbb{F}_m$ in terms of cotangent sheaves.

\subsection{Foliated deformations of $\left( \mathbb{F}_m, \mathcal{N}_{\mathcal{F}_0}^* \right)$ in terms of cotangent sheaves} \

Let us compute the cohomology group $\mathbb{H}^i\left(  \mathbb{F}_m, \mathcal{E}_{\mathcal{O}_{\mathbb{F}_m}(p,q)}^\bullet \right),i=0,1,2$ of the dual leaf complex $(\ref{rs18})$. By considering the spectral sequence associated to $\mathcal{E}_{\mathcal{O}_{\mathbb{F}_m}(p,q)}^\bullet$ and $H^2\left( \mathbb{F}_m, \mathcal{E}_{\mathcal{O}_{\mathbb{F}_m}(p,q)} \right)=0$ from $(\ref{rs17})$, we see that
{\small{\begin{align}
\mathbb{H}^1\left( \mathbb{F}_m, \mathcal{E}_{\mathcal{O}_{\mathbb{F}_m}(p,q)}^\bullet \right) &\cong \textnormal{coker}\left( H^0\left( \mathbb{F}_m, \mathcal{E}_{\mathcal{O}_{\mathbb{F}_m}(p,q)} \right) \to H^0 \left( \mathbb{F}_m, \Omega_{\mathbb{F}_m}^1 (p,q) \right) \right) \bigoplus \textnormal{ker}\left( H^1\left( \mathbb{F}_m, \mathcal{E}_{\mathcal{O}_{\mathbb{F}_m}(p,q)} \right) \to H^1\left( \mathbb{F}_m, \Omega_{\mathbb{F}_m}^1 (p,q) \right) \right)\\
\mathbb{H}^2\left( \mathbb{F}_m, \mathcal{E}_{\mathcal{O}_{\mathbb{F}_m}(p,q)}^\bullet  \right)& \cong \textnormal{coker}\left(H^1\left( \mathbb{F}_m, \mathcal{E}_{\mathcal{O}_{\mathbb{F}_m}(p,q)} \right) \to H^1 \left( \mathbb{F}_m, \Omega_{\mathbb{F}_m}^1 (p,q) \right) \right) \label{rs53}
\end{align}}}

\begin{example}[Interpretation of Example \ref{rsf11} in terms of cotangent sheaves]\label{rsf14}
Let us consider a foliated rational ruled surface $\left( \mathbb{F}_2, \omega_0 \right)$, where $\omega_0= 2w_0 dz_0+ z_0 dw_0 \in H^0\left( \mathbb{F}_2, \Omega_{\mathbb{F}_2}^1 (3,2) \right)$ on $U_0$. Then we have from $(\ref{rsf1}),(\ref{rsf17}),(\ref{rsf3})$ and $(\ref{rsf18})$
\begin{align*}
&H^0\left( \mathbb{F}_2,  \mathcal{E}_{\mathcal{O}_{\mathbb{F}_2}(3,2)} \right)  \\
&= \left\{ \left(a_0+ a_1z_0+ a_2z_0^2\right)\frac{\partial}{\partial z_0}+ \left( \left(a_3-2a_2z_0\right)w_0+ \left(c_0+c_1z_0+ c_2 z_0^2\right)w_0^2\right)\frac{\partial}{\partial w_0} + a_4 +a_2z_0 -2\left(c_0+c_1z_0+ c_2z_0^2\right)w_0^2    \right\} \\
&H^0\left( \mathbb{F}_2, \Omega_{\mathbb{F}_2}^1(3,2)\right) = \left\{ \left(2a_1w_0+  \left(b_0+ b_1z_0\right)w_0^2\right)dz_0 + \left(a_0+ a_1z_0\right)dw_0  \right\} \\
&H^1\left( \mathbb{F}_2,  \mathcal{E}_{\mathcal{O}_{\mathbb{F}_2}(3,2)} \right) = \left\langle \frac{1}{z_0}\frac{\partial}{\partial w_0}  \right\rangle  \\
&H^1\left( \mathbb{F}_2, \Omega_{\mathbb{F}_2}^1(3,2)  \right) = \left\langle  \frac{1}{z_0^2} dz_0 , \frac{1}{z_0}dz_0 \right\rangle
\end{align*}
We compute $H^0\left( \mathbb{F}_2, \mathcal{E}_{\mathcal{O}_{\mathbb{F}_2}(3,2)} \right) \to H^0\left( \mathbb{F}_2, \Omega_{\mathbb{F}_2}^1 ( 3, 2) \right)$
\begin{align} \label{rsf20}
&\mathcal{L}_{\left(a_0+ a_1z_0+ a_2z_0^2\right)\frac{\partial}{\partial z_0}+ \left( \left(a_3-2a_2z_0\right)w_0+ \left(c_0+c_1z_0+ c_2 z_0^2\right)w_0^2 \right)\frac{\partial}{\partial w_0} }\left( 2w_0dz_0 + z_0 dw_0 \right)\\
& + \left( a_4 +a_2z_0 -2\left(c_0+c_1z_0+ c_2z_0^2\right)w_0 \right) (2w_0dz_0 + z_0 dw_0) \notag \\
&=\left(  \left(2a_3+ 2a_1+ 2a_4 \right)w_0 +\left(-2c_0 - c_1z_0\right) w_0^2    \right)dz_0 +\left( a_0+ \left(a_1+ a_3 + a_4  \right)z_0\right)dw_0 = 0  \notag \\
&\iff a_1+a_3+a_4=0, a_0=0,c_0=c_1=0 \notag
\end{align}
Then we have
\begin{align*}
\dim_\mathbb{C} \mathbb{H}^0\left( \mathbb{F}_2 , \mathcal{E}_{\mathcal{O}_{\mathbb{F}_2}(3,2)}^\bullet \right) = \dim_\mathbb{C} \textnormal{ker}\left( H^0\left( \mathbb{F}_2, \mathcal{E}_{\mathcal{O}_{\mathbb{F}_2}(3,2)} \right) \to H^0\left( \mathbb{F}_2, \Omega_{\mathbb{F}_2}^1(3,2)\right) \right)= 4
\end{align*}
and we have
\begin{align*}
\dim_\mathbb{C} \textnormal{coker}\left( H^0\left( \mathbb{F}_2, \mathcal{E}_{\mathcal{O}_{\mathbb{F}_2}(3,2)} \right) \to H^0\left( \mathbb{F}_2, \Omega_{\mathbb{F}_2}^1 (3,2) \right) \right) = 8 - (8-4) =0
\end{align*}
We compute $ H^1\left( \mathbb{F}_2, \mathcal{E}_{\mathcal{O}_{\mathbb{F}_2}(3,2)} \right) \to H^1\left( \mathbb{F}_2, \Omega_{\mathbb{F}_2}^1(3,2)\right)$
\begin{align} \label{rsf21}
\mathcal{L}_{a\frac{1}{z_0}\frac{\partial}{\partial w_0}}\left( 2w_0 dz_0 + z_0 dw_0 \right) +  0\cdot \left( 2w_0 dz_0 + z_0 dw_0 \right)= 2a\frac{1}{z_0} dz_0 + az_0 d\left( \frac{1}{z_0}\right) = a\frac{1}{z_0}dz_0 = 0 \iff a=0
\end{align}
Then we have $\textnormal{ker}\left(  H^1\left( \mathbb{F}_2, \mathcal{E}_{\mathcal{O}_{\mathbb{F}_2}(3,2)} \right) \to H^1\left( \mathbb{F}_2, \Omega_{\mathbb{F}_2}^1(3,2)\right) \right)=0$, so that
\begin{align*}
\mathbb{H}^1\left( \mathbb{F}_2, \mathcal{E}_{\mathcal{O}_{\mathbb{F}_2}(3,2)}^\bullet \right)=0
\end{align*}
Hence $\left( \mathbb{F}_2 , \omega_0 \right)$ is rigid in foliated deformations by \textnormal{Theorem \ref{rsf25}}.

\end{example}

\begin{example}[Interpretation of Example \ref{rsf12} in terms of cotangent sheaves]\label{rsf15}
Let us consider a foliated rational ruled surface $\left( \mathbb{F}_2, \omega \right)$, where $\omega_0= w_0 dz_0+ z_0 dw_0 \in H^0\left( \mathbb{F}_2, \Omega_{\mathbb{F}_2}^1(4,2) \right)$ on $U_0$. Then we have from $(\ref{rsf1}),(\ref{rsf31}),(\ref{rsf3})$ and $(\ref{rsf32})$
\begin{align}
&H^0\left( \mathbb{F}_2,  \mathcal{E}_{\mathcal{O}_{\mathbb{F}_2}(4,2)} \right)  \\
&= \left\{ \left(a_0+ a_1z_0+ a_2z_0^2\right)\frac{\partial}{\partial z_0}+ \left( \left(d-2a_2z_0\right)w_0+ \left(c_0+c_1z_0+ c_2 z_0^2\right)w_0^2\right)\frac{\partial}{\partial w_0} + b  -2\left(c_0+c_1z_0+ c_2z_0^2\right)w_0^2    \right\} \notag \\
&H^0\left( \mathbb{F}_2, \Omega_{\mathbb{F}_2}^1\left( 4, 2 \right) \right)= \left\{ \left( \left( c+ 2 a_2 z_0 \right)  w_0 +  \left(b_0+ b_1 z_0+ b_2 z_0^2 \right) w_0^2 \right) dz_0 +   \left(a_0 + a_1 z_0+ a_2 z_0^2  \right)dw_0  \right\} \label{rsf35}\\
&H^1\left( \mathbb{F}_2, \mathcal{E}_{\mathcal{O}_{\mathbb{F}_2}(4,2)} \right) = \left\langle \frac{1}{z_0}\frac{\partial}{\partial w_0} \right\rangle\\
&H^1\left( \mathbb{F}_2, \Omega_{\mathbb{F}_2}^1 (4,2) \right) = \left\langle  \frac{1}{z_0} d z_0 \right\rangle
\end{align}
We compute $H^0\left( \mathbb{F}_2, \mathcal{E}_{\mathcal{O}_{\mathbb{F}_2}(4,2)} \right)\to H^0\left( \mathbb{F}_2, \Omega_{\mathbb{F}_2}^1(4,2) \right)$
 \begin{align}
 &\mathcal{L}_{\left( a_0 + a_1 z_0 + a_2 z_0^2 \right) \frac{\partial}{\partial z_0} +\left( \left( - 2a_2 z_0+ d\right) w_0 + \left( c_0+ c_1 z_0 + c_2 z_0^2\right) w_0^2  \right)\frac{\partial}{\partial w_0}}\left( w_0 dz_0 + z_0 dw_0 \right) + \left( b - 2\left(c_0+ c_1 z_0 + c_2 z_0^2 \right) w_0 \right)\left( w_0 dz_0 + z_0 dw_0 \right) \notag \\
 &=\left( \left( - 2a_2 z_0+ d+b + a_1\right) w_0 + \left( - c_0 + c_2 z_0^2\right) w_0^2 \right) dz_0  +\left( a_0 + (a_1+d +b) z_0 - a_2 z_0^2   \right) dw_0  = 0 \label{rsf34} \\
 &\iff  a_2=0, a_1+ d+ b=0,  c_0=c_2=0, a_0=0 \notag
 \end{align}
 Then we have
  \begin{align*}
 \dim_\mathbb{C} H^0\left(  \mathbb{F}_2, \mathcal{E}_{\mathcal{O}_{\mathbb{F}_2}(4,2)}^\bullet  \right) = \textnormal{ker}\left( H^0\left( \mathbb{F}_2, \mathcal{E}_{\mathcal{O}_{\mathbb{F}_2}(4,2)} \right)\to H^0\left( \mathbb{F}_2, \Omega_{\mathbb{F}_2}^1(4,2) \right) \right) = 3
 \end{align*}
 and
 \begin{align*}
 \dim_\mathbb{C}  \textnormal{coker}\left( H^0\left( \mathbb{F}_2, \mathcal{E}_{\mathcal{O}_{\mathbb{F}_2}(4,2)} \right) \to H^0\left( \mathbb{F}_2 , \Omega_{\mathbb{F}_2}^1(4,2) \right) \right) = 7 - (8 -3)=2
 \end{align*}
 By comparing $(\ref{rsf35})$ and $(\ref{rsf34})$, we see that
  \begin{align*}
 \textnormal{coker}\left( H^0\left( \mathbb{F}_2, \mathcal{E}_{\mathcal{O}_{\mathbb{F}_2}(4,2)} \right) \to H^0\left( \mathbb{F}_2 , \Omega_{\mathbb{F}_2}^1(4,2) \right) \right)\cong \,\,\, \textnormal{Span}_\mathbb{C} \left\langle w_0 dz_0, z_0w_0^2 dz_0 \right\rangle \,\,\, or\,\,\,\textnormal{Span}_\mathbb{C} \left\langle z_0 dw_0, z_0 w_0^2 dz_0 \right\rangle
 \end{align*}
 On the other hand, we compute $H^1\left( \mathbb{F}_2, \mathcal{E}_{\mathcal{O}_{\mathbb{F}_2}(4,2)} \right)\to H^1\left( \mathbb{F}_2, \Omega_{\mathbb{F}_2}^1(4,2) \right)$
 \begin{align}\label{rf7}
 \mathcal{L}_{k\frac{1}{z_0}\frac{\partial}{\partial w_0}}\left( w_0 dz_0 + z_0 dw_0 \right) + 0 \cdot \left( w_0 dz_0 + z_0 dw_0 \right)=k \frac{1}{z_0} dz_0 + k z_0\left( -\frac{1}{z_0^2} \right) dz_0 = 0
 \end{align}
 This implies that
 \begin{align*}
 \dim_\mathbb{C} \mathbb{H}^1\left( \mathbb{F}_2, \mathcal{E}_{\mathcal{O}_{\mathbb{F}_2}(4,2)}^\bullet \right)= 3, \,\,\,\,\,\,\, \mathbb{H}^1\left( \mathbb{F}_2, \mathcal{E}_{\mathcal{O}_{\mathbb{F}_2}(4,2)}^\bullet \right)\cong \textnormal{Span}_\mathbb{C} \left\langle  z_0 dw_0, z_0 w_0^2 dz_0 \right\rangle \bigoplus \textnormal{Span}_\mathbb{C}\left\langle \frac{1}{z_0}\frac{\partial}{\partial w_0} \right\rangle
 \end{align*}
 We note that
  \begin{align*}
 \mathcal{L}_{\frac{1}{z_0}\frac{\partial}{\partial w_0}}\left( z_0 dw_0 \right) + 0\cdot z_0 dw_0= z_0 \left( - \frac{1}{z_0^2} \right) dz_0= - \frac{1}{z_0} dz_0
 \end{align*}
 is not in the image of $H^1\left( \mathbb{F}_2, \mathcal{E}_{\mathcal{O}_{\mathbb{F}_2}(4,2)} \right) \to H^1\left( \mathbb{F}_2, \Omega_{\mathbb{F}_2}^1(4,2) \right)$ which is the zero map. Hence $\left( \mathbb{F}_2, \omega_0 \right)$ is obstructed in foliated deformations by \textnormal{Lemma \ref{on23}}.

\end{example}

\begin{example}[Interpretation of Example \ref{rsf13} in terms of cotangent sheaves]\label{rsf16}
Let us consider a foliated rational ruled surface $\left( \mathbb{F}_2, \omega \right)$, where $\omega_0 = z_0w_0 dz_0 - dw_0\in H^0\left( \mathbb{F}_2, \Omega_{\mathbb{F}_2}^1(5,2) \right)$ on $U_0$. Then we have from $(\ref{rsf1}),(\ref{rf13}),(\ref{rsf3})$, and $(\ref{rf15})$
\begin{align}
&H^0\left( \mathbb{F}_2, \mathcal{E}_{\mathcal{O}_{\mathbb{F}_2}(5,2)}\right) \\
&= \left\{ \left(a_0+ a_1z_0+ a_2z_0^2\right)\frac{\partial}{\partial z_0}+ \left( \left(a_3-2a_2z_0\right)w_0+ \left(c_0+c_1z_0+ c_2 z_0^2\right)w_0^2 \right)\frac{\partial}{\partial w_0} + a_4 -a_2z_0 -2\left(c_0+c_1z_0+ c_2z_0^2\right)w_0    \right\} \\
&H^0\left( \mathbb{F}_2, \Omega_{\mathbb{F}_2}^1(5,2) \right)  = \left\{ \left( \left( c_0+ c_1 z_0 + 2a_3 z_0^2 \right) w_0 + \left( e_0 + e_1 z_0 + e_2 z_0^2 + e_3 z_0^3 \right) w_0^2\right) dz_0 + \left(a_0+ a_1z_0 + a_2 z_0^2 + a_3 z_0^3 \right) dw_0 \right\} \label{rf19}\\
&H^1\left( \mathbb{F}_2, \mathcal{E}_{\mathcal{O}_{\mathbb{F}_2}(5,2)} \right) = \left\langle \frac{1}{z_0} \frac{\partial}{\partial w_0} \right\rangle\\
&H^1\left( \mathbb{F}_2, \Omega_{\mathbb{F}_2}^1(5,2) \right) = 0
\end{align}
We compute $H^0\left( \mathbb{F}_2, \mathcal{E}_{\mathcal{O}_{\mathbb{F}_2}(5,2)} \right)\to H^0\left( \mathbb{F}_2 ,  \Omega_{\mathbb{F}_2}^1(5,2) \right)$
\begin{align}
&\mathcal{L}_{\left(a_0+ a_1z_0+ a_2z_0^2\right)\frac{\partial}{\partial z_0}+ \left( \left(a_3-2a_2z_0\right)w_0+ \left(c_0+c_1z_0+ c_2 z_0^2\right)w_0^2 \right)\frac{\partial}{\partial w_0}} \left( z_0w_0 dz_0 - dw_0 \right) \notag \\
&+ \left( a_4 -a_2z_0 -2\left(c_0+c_1z_0+ c_2z_0^2\right)w_0  \right)\left(z_0w_0dz_0 - dw_0 \right) \notag \\
&= \left( \left(a_0+ 2a_2+ \left(2a_1+ a_4+ a_3\right)z_0  \right)w_0  - \left(c_1+ \left(2c_2+c_0 \right)z_0 + c_1z_0^2 + c_2z_0^3\right)w_0^2  \right) dz_0+ \left( -a_3 -a_4 + 3a_2z_0  \right)dw_0 \label{rf18}\\
&\iff a_0 + 2 a_2=0, 2a_1+a_4+a_3=0, c_0=c_1=c_2=0, a_3+a_4=0, a_2=0 \notag
\end{align}
Then we have
\begin{align*}
\dim_\mathbb{C} \mathbb{H}^0\left( \mathcal{E}_{\mathcal{O}_{\mathbb{F}_2}(5,2)}^\bullet \right)= \dim_\mathbb{C} \textnormal{ker}\left( H^0\left( \mathbb{F}_2, \mathcal{E}_{\mathcal{O}_{\mathbb{F}_2}(5,2)} \right)\to H^0\left( \mathbb{F}_2 ,  \Omega_{\mathbb{F}_2}^1(5,2) \right) \right) =1
\end{align*}
and
\begin{align*}
\dim_\mathbb{C} \textnormal{coker}\left( H^0\left( \mathbb{F}_2, \mathcal{E}_{\mathcal{O}_{\mathbb{F}_2}(5,2)} \right)\to H^0\left( \mathbb{F}_2 ,  \Omega_{\mathbb{F}_2}^1(5,2) \right) \right) = 10 - (8-1)=3
\end{align*}
By comparing $(\ref{rf19})$ and $(\ref{rf18})$ we see that
\begin{align*}
\textnormal{coker}\left( H^0\left( \mathbb{F}_2, \mathcal{E}_{\mathcal{O}_{\mathbb{F}_2}(5,2)} \right)\to H^0\left( \mathbb{F}_2 ,  \Omega_{\mathbb{F}_2}^1(5,2) \right) \right)\cong \textnormal{Span}_\mathbb{C} \left\langle -z_0^2dw_0, -2z_0^2 w_0 dz_0 - z_0^3 dw_0, w_0^2 dz_0 \right\rangle
\end{align*}
Since $H^1\left( \mathbb{F}_2, \Omega_{\mathbb{F}_2}^1(5,2) \right)=0$, we have $\mathbb{H}^2\left( \mathbb{F}_2, \mathcal{E}_{\mathcal{O}_{\mathbb{F}_2}(5,2)}^\bullet \right)=0$, so that in particular $\left( \mathbb{F}_2, \omega_0 \right)$ is unobstructed in foliated deformations, and we have
\begin{align}\label{rf20}
\mathbb{H}^1\left( \mathbb{F}_2, \mathcal{E}_{\mathcal{O}_{\mathbb{F}_2}(5,2)}^\bullet \right) \cong \textnormal{Span}_\mathbb{C} \left\langle -z_0^2dw_0, -2z_0^2 w_0 dz_0 - z_0^3 dw_0, w_0^2 dz_0 \right\rangle\bigoplus \textnormal{Span}_\mathbb{C}\left\langle  \frac{1}{z_0}\frac{\partial}{\partial w_0} \right\rangle
\end{align}

We will construct a foliated analytic family of deformations of $\left( \mathbb{F}_2, \omega_0 \right)$ such that the foliated Kodaira-Spencer map is an isomorphism at the distinguished point. First we set
\begin{align*}
\sigma_{01}:=\mathcal{L}_{\frac{1}{z_0}\frac{\partial}{\partial w_0}}\left( z_0w_0 dz_0 - dw_0 \right) =  dz_0 - d\left(\frac{1}{z_0} \right) = dz_0 + \frac{1}{z_0^2}dz_0 =\left( 1+ \frac{1}{z_0^2} \right) dz_0 \in C^1\left( \mathfrak{B}, \Omega_{\mathbb{F}_2}^1(5,2) \right)
\end{align*}
where $\mathfrak{B}= \left\{ U_0 \bigcup U_2, U_1 \bigcup U_3  \right\}$. From $z_1= \frac{1}{z_0}$, we have $dz_1= - z_1^2 dz_0$. Then if we take $\sigma_0= dz_0$ on $U_0 \bigcup U_2$ and $\sigma_1= - z_1 dz_1 $ on $U_1 \bigcup U_3$, we have
\begin{align}\label{rf23}
h_{01}\sigma_1- \sigma_0 + \sigma_{01} = z_0 \cdot \left(    z_1dz_1 \right) - dz_0 + \left( 1+ \frac{1}{z_0^2} \right)dz_0 =0
\end{align}
Then let us consider the foliated analytic family over $\textnormal{Spec}\left( \frac{\mathbb{C}[t_1,t_2,t_3,t_4]}{\left( t_1^2, t_2^2, t_3^2, t_4^2, t_1t_2,t_1t_3, t_1t_4, t_2t_3, t_2t_4, t_3t_4 \right)} \right)$ on $U_0$ defined by $(\ref{rf20})$
\begin{align*}
\omega&= \left(  z_0w_0 - t_1 + 2t_2z_0^2 w_0 - t_3w_0^2 \right) dz_0 + \left(-1 + t_4z_0^2 + t_2z_0^3 \right) dw_0\\
z_0 &=\frac{1}{z_1}\\
w_0&= z_1^2w_1 + t_1 z_1\\
h_{01}&= z_0
\end{align*}
Since $z_0= \frac{1}{z_1}$ and $w_0= z_1^2 w_1 + t_1z_1$, we have
\begin{align*}
dz_0= - \frac{1}{z_1^2}dz_1,\,\,\,\,\,\,\,\,dw_0= \left(2 z_1w_1 + t_1 \right) dz_1 + z_1^2 dw_1
\end{align*}
Then we have
{\small{\begin{align*}
h_{10}\omega &= \left(w_0-  t_1 \frac{1}{z_0} + 2t_2z_0w_0 - t_3\frac{1}{z_0}w_0^2 \right) dz_0 + \left( - \frac{1}{z_0} + t_4z_0 + t_2 z_0^2 \right)dw_0\\
        & = \left(w_0-  t_1 z_1 + 2t_2\frac{1}{z_1}\left(z_1^2w_1+ t_1z_1\right) - t_3z_1 \left(z_1^2w_1+ t_1z_1\right)^2 \right)\left(- \frac{1}{z_1^2} \right)dz_1 +\left(- z_1+ t_4\frac{1}{z_1}+ t_2 \frac{1}{z_1^2} \right)\left( \left(2z_1w_1+ t_1\right)dz_1 + z_1^2 dw_1\right)\\
        &=\left( - w_1 - 2t_2\frac{w_1}{z_1}- 2t_1t_2\frac{1}{z_1^2} +  t_3\left(z_1^3w_1^2+ 2t_1z_1^2w_1+ t_1^2z_1 \right) \right)dz_1 \\
        &+\left( - t_1z_1+ t_1t_4\frac{1}{z_1}+ t_1t_2 \frac{1}{z_1^2} - 2z_1^2w_1+ 2t_4w_1 + 2t_2 \frac{w_1}{z_1} \right) dz_1  +   \left( - z_1^3 +t_4z_1+ t_2    \right)dw_1\\
         &=\left( - w_1- 2z_1^2w_1 - t_1z_1+ 2t_4w_1+ t_1t_4\frac{1}{z_1}- t_1t_2\frac{1}{z_1^2} +  t_3\left(z_1^3w_1^2+ 2t_1z_1^2w_1+ t_1^2z_1 \right) \right)dz_1 +   \left( - z_1^3 +t_4z_1+ t_2    \right)dw_1
\end{align*}}}
We note that the rational part is $\left( t_1t_4 \frac{1}{z_1}- t_1t_2 \frac{1}{z_1^2} \right)dz_1$ and we have
\begin{align*}
\left( t_1t_4\frac{1}{z_1}- t_1t_2 \frac{1}{z_1^2}  \right)dz_1 = \left(t_1t_4z_0 - t_1t_2 z_0^2 \right)\left( -\frac{1}{z_0^2} \right)dz_0 =\frac{1}{z_0}\left(-t_1t_4+ t_1t_2 z_0 \right)dz_0
\end{align*}
Hence the desired foliated analytic family is given by
\begin{align*}
\omega|_{U_0 \bigcup U_2}&= \left(  z_0w_0 - t_1 +t_1t_4 -t_1t_2z_0 + 2t_2z_0^2 w_0 - t_3w_0^2 \right) dz_0 + \left(-1 + t_4z_0^2 + t_2z_0^3 \right) dw_0\\
\omega|_{U_1 \bigcup U_3}&= \left( - w_1- 2z_1^2w_1 - t_1z_1+ 2t_4w_1+  t_3\left(z_1^3w_1^2+ 2t_1z_1^2w_1+ t_1^2z_1 \right) \right)dz_1 +   \left( - z_1^3 +t_4z_1+ t_2    \right)dw_1 \\
z_0 &=\frac{1}{z_1}\\
w_0&= z_1^2w_1 + t_1 z_1
\end{align*}
We describe the foliated analytic family more precisely. We have coordinate transformations as follows$:$
\begin{center}
$\begin{CD}
U_0: \left( z_0, w_0 \right) @>>> U_1:\left(z_1, w_1 \right)=\left( \frac{1}{z_0}, z_0^2 w_0- t_1 z_0 \right)\\
@VVV @VVV \\
U_2:\left( z_2, w_2 \right)= \left( z_0, \frac{1}{w_0} \right) @>>> U_3: \left( z_3 , w_3 \right) = \left( \frac{1}{z_2} , \frac{w_2}{ z_2^2 - t_1 z_2w_2 } \right)  =\left( z_1, \frac{1}{w_1} \right)
\end{CD}$
\end{center}
and the foliation $\omega$ is given by
\begin{align*}
\omega&= \left(  z_0w_0 - t_1 +t_1t_4 -t_1t_2z_0 + 2t_2z_0^2 w_0 - t_3w_0^2 \right) dz_0 + \left(-1 + t_4z_0^2 + t_2z_0^3 \right) dw_0\,\,\,\,\,\,\,\,\textnormal{on $U_0$}\\
\omega&= \left( - w_1- 2z_1^2w_1 - t_1z_1+ 2t_4w_1+  t_3\left(z_1^3w_1^2+ 2t_1z_1^2w_1+ t_1^2z_1 \right) \right)dz_1 +   \left( - z_1^3 +t_4z_1+ t_2    \right)dw_1 \,\,\,\,\,\,\,\,\,\textnormal{on $U_1$} \\
\omega&= \left(  z_2w_2 - t_1w_2^2 +t_1t_4w_2^2 -t_1t_2z_2w_2^2 + 2t_2z_2^2 w_2 - t_3 \right)  dz_2 - \left(-1 + t_4z_2^2 + t_2z_2^3 \right) dw_2\,\,\,\,\,\,\,\,\textnormal{on $U_2$}\\
\omega&= \left( - w_3- 2z_3^2w_3 - t_1z_3w_3^2+ 2t_4w_3+  t_3\left(z_3^3+ 2t_1z_3^2w_3+ t_1^2z_3w_3^2 \right) \right)dz_3 -   \left( - z_3^3 +t_4z_3+ t_2    \right)dw_3 \,\,\,\,\,\,\,\,\,\textnormal{on $U_3$}
\end{align*}
and transition functions are given by
\begin{align*}
h_{01}&= z_0\\
h_{02}&= w_0^2 \\
h_{03}&= z_0\left( z_0^2 w_0 - t_1 z_0 \right)^2 \\
h_{12}&= z_1\left(z_1^2 w_1 + t_1 z_1 \right)^2 = z_2^{-1} w_2^{-2}\\
h_{13}&= w_1^2\\
h_{23}&= w_2^2 z_2\left( z_2^2 \frac{1}{w_2} - t_1z_2 \right)^2= z_2\left( z_2^2- t_1z_2w_2 \right)^2
\end{align*}
We note that the foliated Kodaira-Spence map at $t=\left(t_1,t_2,t_3,t_4 \right)= (0,0,0,0)=0$ is given by $($in terms of the covering $\mathfrak{B}$ where $\mathfrak{B}= \left\{ U_0 \bigcup U_2, U_1 \bigcup U_3 \right\})$
\begin{align*}
T_0\left( \mathbb{C}\right) &\to \mathbb{H}^1\left( \mathbb{F}_2, \mathcal{E}_{\mathcal{O}_{\mathbb{F}_2}(5,2)}^\bullet \right) \\
\frac{\partial}{\partial t_1} &\mapsto \left( \left\{ \left( dz_0,  z_1dz_1 \right) \right\}, \left\{ \left( z_1 \frac{\partial}{\partial w_0}, 0 \right) \right\} \right) \in C^0\left( \mathfrak{B}, \Omega_{\mathbb{F}_2}^1(5,2) \right) \bigoplus C^1\left( \mathfrak{B}, \mathcal{E}_{\mathcal{O}_{\mathbb{F}_2}(5,2)} \right)\\
\frac{\partial}{\partial t_2} &\mapsto \left( \left\{ \left( - 2 z_0^2w_0dz_0- z_0^3 dw_0 ,  - dw_1 \right) \right\}, \left\{ \left( 0 , 0 \right) \right\} \right) \in C^0\left( \mathfrak{B}, \Omega_{\mathbb{F}_2}^1(5,2) \right) \bigoplus C^1\left( \mathfrak{B}, \mathcal{E}_{\mathcal{O}_{\mathbb{F}_2}(5,2)} \right)\\
\frac{\partial}{\partial t_3} &\mapsto \left( \left\{ \left( w_0^2 dz_0 , - z_1^3 w_1^2 dz_1 \right) \right\}, \left\{ \left( 0 , 0 \right) \right\} \right) \in C^0\left( \mathfrak{B}, \Omega_{\mathbb{F}_2}^1(5,2) \right) \bigoplus C^1\left( \mathfrak{B}, \mathcal{E}_{\mathcal{O}_{\mathbb{F}_2}(5,2)} \right)\\
\frac{\partial}{\partial t_4} &\mapsto \left( \left\{ \left( - z_0^2 dw_0,   - 2 w_1 dz_1 - z_1 dw_1 \right) \right\}, \left\{ \left( 0, 0 \right) \right\} \right) \in C^0\left( \mathfrak{B}, \Omega_{\mathbb{F}_2}^1(5,2) \right) \bigoplus C^1\left( \mathfrak{B}, \mathcal{E}_{\mathcal{O}_{\mathbb{F}_2}(5,2)} \right)
\end{align*}

\end{example}

\begin{remark}[restatement of Remark \ref{rs52} in terms of cotangent sheaves]
In general, let $\omega_0\in H^0\left( \mathbb{F}_m, \Omega_{\mathbb{F}_m}^1(p,2) \right)$ be a foliation on $\mathbb{F}_m$. From $(\ref{rs48})$, we see that if $p\geq 2m+1$, then $H^1\left( \mathbb{F}_m, \Omega_{\mathbb{F}_m}^1(p,2)\right)=0$. This implies that for $p\geq 2m+1$, $\mathbb{H}^2\left( \mathbb{F}_m, \mathcal{E}_{\mathcal{O}_{\mathbb{F}_m}(p,2)}^\bullet \right)=0$ from $(\ref{rs53})$, so that $\left( \mathbb{F}_m,  \omega \right)$ is unobstructed in foliated deformations.
\end{remark}

\begin{remark}[reinterpretation of Remark \ref{rs26} in terms of cotangent sheaves]\label{rs40}
So far we have assumed $q\geq 1$ since $H^1\left(U_0\bigcup U_2, \Omega_{\mathbb{F}_m}^1(p,q) \right)\ne 0$ for $q \leq 0$. Now let us consider $q=0$. Let $\omega_0 \in H^0\left( \mathbb{F}_m, \Omega_{\mathbb{F}_m}^1 (p,0) \right)$ be a foliation on $\mathbb{F}_m$. Then from $(\ref{rs35})$ we see that $\omega_0$ is of the form $\omega=c_0^1(z_0)dz_0$, where $c_0^1(z_0)$ is an entire function. Unless $c_0^1(z_0)=a\ne 0$ is a constant, then $\omega_0$ does not have an isolated singularity, i.e. does not define a foliation. Hence $\omega$ is of the form $\omega=a dz_0,a\ne 0$ on $U_0$. We note that $\omega_0$ is on $U_0$ is translated into $-a z_1^{p-2}dz_1$ on $U_1$. In order for $\omega_0$ to define a foliation on $\mathbb{F}_m$, we should have $p=2$. Hence we see that $\omega_0 \in H^0\left( \mathbb{F}_m, \Omega_{\mathbb{F}_m}^1 (2,0)  \right)$. Now let us consider foliated deformations of $\left( \mathbb{F}_m, \omega_0= a dz_0 \right) , a\ne 0$. Then we have from $(\ref{rsf2}),(\ref{rsf1}),(\ref{rs37})$ and $(\ref{rsf3})$
\begin{align*}
H^0\left( \mathbb{F}_m, \mathcal{E}_{\mathcal{O}_{\mathbb{F}_m}(2,0)} \right)&=  \begin{cases} \left\{  \left( a_0z_0^2+ bz_0 + c \right)\frac{\partial}{\partial z_0} + \left( d + ez_0 + fz_0^2 \right) \frac{\partial}{\partial w_0} , h -2 a_0 z_0\right\}\,\,\,\,\,\,\,\,\,\,\,\,\,\, \textnormal{for}\,\,\, m=0 & \\
\left\{ \left(a_0+ a_1 z_0 + a_2z_0^2    \right) \frac{\partial}{\partial z_0} + \left( \left(-ma_2z_0+d\right)w_0 + \left(c_0 + c_1z_0+ \cdots + c_m z_0^m \right)w_0^2 \right)\frac{\partial}{\partial w_0}  , \right.  \\
\,\,\,\,\, \left( b- 2 a_2 z_0 +2 \left( c_0 + c_1 z_0 + \cdots + c_m z_0^m \right) w_0\right\} \,\,\,\,\,\,\,\,\,\,\,\, \textnormal{for}\,\,\, m \geq 1 &  \end{cases} \\
H^0\left( \mathbb{F}_m, \Omega_{\mathbb{F}_m}^1(2,0) \right)&= \left\{ c dz_0 | c\in \mathbb{C} \right\} \\
H^1\left( \mathbb{F}_m, \mathcal{E}_{\mathcal{O}_{\mathbb{F}_m} (2,0)} \right) &=  \begin{cases} 0 & \textnormal{for}\,\,\,m=0,1 \\
\textnormal{Span}_\mathbb{C}\left\langle \left(\frac{1}{z_0}\frac{\partial}{\partial w_0}, 0 \right), \left( \frac{1}{z_0^2}\frac{\partial}{\partial w_0}, 0 \right), \cdots , \left( \frac{1}{z_0^{m-1}}\frac{\partial}{\partial w_0}, 0 \right) \right\rangle & \textnormal{for}\,\,\, m\geq 2 \end{cases}
\end{align*}
Then we compute $H^0\left( \mathbb{F}_m , \mathcal{E}_{\mathcal{O}_{\mathbb{F}_m}(2,0)} \right)\xrightarrow{\tilde{E}_0} H^0\left( \mathbb{F}_m, \Omega_{\mathbb{F}_m}^1(2,0) \right)$
\begin{align*}
&\mathcal{L}_{ \left( a_0z_0^2+ bz_0 + c \right)\frac{\partial}{\partial z_0} + \left( d + ez_0 + fz_0^2 \right) \frac{\partial}{\partial w_0} }\left(a dz_0 \right) + a\left( h- 2 a_0 z_0 \right) dz_0 = a \left( 2 a_0 z_0 + b \right)dz_0 + (ah- 2aa_0 z_0)dz_0= \left(ab+ ah \right) dz_0\\
&\mathcal{L}_{ \left(a_0+ a_1 z_0 + a_2z_0^2    \right) \frac{\partial}{\partial z_0} + \left( \left(-ma_2z_0+d\right)w_0 + \left(c_0 + c_1z_0+ \cdots + c_m z_0^m \right)w_0^2 \right)\frac{\partial}{\partial w_0} }\left( a dz_0 \right)+ a\left( b- 2 a_2 z_0 \right) dz_0\\
&= a\left( a_1 + 2a_2 z_0 \right)dz_0 + \left( ab - 2a a_2 z_0 \right) dz_0 = (a a_1 + a b) dz_0
\end{align*}
This implies that $\textnormal{coker}\left(H^0\left( \mathbb{F}_m, \mathcal{E}_{\mathcal{O}_{\mathbb{F}_m}(2,0)} \right) \right)\to H^0\left( \mathbb{F}_m, \Omega_{\mathbb{F}_m}^1( 2, 0)\right)=0 $. On the other hand, let us compute the image of $H^1\left( \mathbb{F}_m, \mathcal{E}_{\mathcal{O}_{\mathbb{F}_m}(2,0)} \right)\to H^1\left( \mathbb{F}_m, \Omega_{\mathbb{F}_m}^1(2,0) \right)$ which is given by
\begin{align*}
\mathcal{L}_{k_1\frac{1}{z_0}\frac{\partial}{\partial w_0}+ \cdots + k_{m-1} \frac{1}{z_0^{m-1}}\frac{\partial}{\partial w_0}}\left( adz_0 \right) + 0 \cdot a dz_0= 0
\end{align*}
Hence we have
\begin{align*}
\mathbb{H}^1\left( \mathbb{F}_m, \mathcal{E}_{\mathcal{O}_{\mathbb{F}_m}(2,0)}^\bullet \right)= H^1\left(\mathbb{F}_m, \Theta_{\mathbb{F}_m} \right) =\textnormal{Span}_\mathbb{C}\left\langle  z_1 \frac{\partial}{\partial w_0},..., z_1^{m-1}\frac{\partial}{\partial w_0} \right\rangle
\end{align*}
We will construct a foliated analytic family of deformations of $\left( \mathbb{F}_m, \omega_0 = a dz_0 \right)$ such that the foliated Kodaira-Spencer map is an isomorphism at the distinguished point, so that $\left( \mathbb{F}_m, \omega_0 \right)$ is unobstructed in foliated deformations. The foliated analytic family is given by
\begin{align*}
\omega&= a dz_0 \\
z_0&= \frac{1}{z_1} \\
w_0&= z_1^m w_1+ t_1 z_1+ \cdots + t_{m-1} z_1^{m-1}
\end{align*}
We describe the foliated analytic family more precisely. We have coordinate transformations as follows$:$
\begin{center}
$\begin{CD}
U_0: \left( z_0, w_0 \right) @>>> U_1:\left(z_1, w_1 \right)=\left( \frac{1}{z_0}, z_0^m w_0- t_1 z_0^{m-1}- \cdots - t_{m-1} z_0 \right)\\
@VVV @VVV \\
U_2:\left( z_2, w_2 \right)= \left( z_0, \frac{1}{w_0} \right) @>>> U_3: \left( z_3 , w_3 \right) = \left( \frac{1}{z_2} , \frac{w_2}{ z_2^m - t_1 z_2^{m-1}w_2 - \cdots - t_{m-1} z_2 w_2} \right)  =\left( z_1, \frac{1}{w_1} \right)
\end{CD}$
\end{center}
and the foliation $\omega$ is given by
\begin{align*}
\omega&= a dz_0\,\,\,\,\,\,\,\textnormal{on}\,\,\, U_0,\,\,\,\,\,\,\,\,\,\,\,\,\,\omega= -adz_1\,\,\,\,\,\,\,\textnormal{on}\,\,\, U_1 \\
\omega&= a dz_2\,\,\,\,\,\,\,\textnormal{on}\,\,\, U_2,\,\,\,\,\,\,\,\,\,\,\,\,\,\omega= -adz_3\,\,\,\,\,\,\,\textnormal{on}\,\,\, U_3 
\end{align*}
and transition functions are given by
\begin{align*}
h_{01}&= z_0^2\\
h_{02}&=1 \\
h_{03}&= z_0^2 \\
h_{12}&= z_1^2= z_2^{-2}\\
h_{13}&= 1\\
h_{23}& = z_2^2= z_3^{-2}
\end{align*}
 
\end{remark}

Next we will study unfoldings of $\mathbb{F}_m$. Let $\left( \mathbb{F}_m, \omega \right)$ be a foliated rational ruled surface, where $\omega\in  H^0\left( \mathbb{F}_m, \Omega_{\mathbb{F}_m}^1(p,q) \right)$, so that we have the exact sequence
\begin{align*}
0\to \mathcal{O}_{\mathbb{F}_m}(-p,-q) \xrightarrow{w} \Omega_{\mathbb{F}_m}^1 \to \frac{\Omega_{\mathbb{F}_m}^1}{\mathcal{O}_{\mathbb{F}_m} ( - p,-q)} \to 0
\end{align*}
Then the unfolded dual leaf complex is given by
\begin{align}\label{rs15}
\tilde{\mathcal{E}}_{\mathcal{O}_{\mathbb{F}_m}(p,q)}^\bullet : \mathcal{E}_{\mathcal{O}_{\mathbb{F}_m}(p,q)} \to \Omega_{\mathbb{F}_m}^1(p,q)\bigoplus  \mathcal{O}_{\mathbb{F}_m}(p,q) \to \bigwedge^2 \Omega_{\mathbb{F}_m}^1(2p,2q)\cong \mathcal{O}_{\mathbb{F}_m}(2p-m-2, 2q-2)\to 0
\end{align}
where $\bigwedge^2 \Omega_{\mathbb{F}_m}^1(2p,2q):=\bigwedge^2 \Omega_{\mathbb{F}_m}^1 \otimes \mathcal{O}_{\mathbb{F}_m}(2p,2q)$. In order to compute $\mathbb{H}^i\left( \mathbb{F}_m, \tilde{\mathcal{E}}_{\mathcal{O}_{\mathbb{F}_m}(p,q)}^\bullet \right)$ for $i=0,1,2$, we need to describe cohomology groups in each degree in $(\ref{rs15})$. In view of subsection \ref{rs34}, subsection \ref{rs3}, subsection \ref{rf17} and subsection \ref{rf16}, It remains to describe $H^i\left( \mathbb{F}_m, \mathcal{O}_{\mathbb{F}_m}(p,q) \right),i=0,1$, and $H^0\left(\mathbb{F}_m, \bigwedge^2 \Omega_{\mathbb{F}_m}^1(2p,2q) \right)$.

\subsection{Description of $H^0\left( \mathbb{F}_m, \mathcal{O}_{\mathbb{F}_m(p,q)} \right)$}\

We keep the notations in subsection \ref{rs1}. We explicitly describe $H^0\left(\mathbb{F}_m, \mathcal{O}_{\mathbb{F}_m}(p,q) \right)$. Let $v \in H^0\left( \mathbb{F}_m, \mathcal{O}_{\mathbb{F}_m}(p,q) \right)$ be represented  by the holomorphic function $v_0(z_0,w_0)$ on $U_0$ and by the holomorphic function $v_2(z_2,w_2)$ on $U_2$, respectively. Then we have $v_0(z_0,w_0)= h_{02} v_2(z_2,w_2)= w_0^q v_2(z_2,w_2)$. If we set $v_0(z_0,w_0)= \sum_{h,k\geq 0} a_{hk}z_0^h w_0^k, a_{hk}\in \mathbb{C}$. Then $v_2(z_2,w_2)= \sum_{h,k\geq 0} a_{hk}z_2^h w_2^{q-k}$, so that $a_{hk}=0$ for $k >q$. Hence $v\in H^0\left( \mathbb{F}_m, \mathcal{O}_{\mathbb{F}_m}(p,q) \right)$ on $U_0$ is of the form
\begin{align}\label{rr3}
v_0(z_0,w_0)= a_0^1(z_0)+ a_0^2(z_0)w_0 + \cdots + a_0^{q+1}(z_0)w_0^q
\end{align}
where $a_0^1(z_0),..., a_0^{q+1}(z_0)$ are entire functions of $z_0$. In particular, we note that
\begin{align*}
H^0\left(\mathbb{F}_m, \mathcal{O}_{\mathbb{F}_m}(p,q) \right) =0\,\,\,\,\,\,\,\,\,\,\textnormal{for}\,\,\, q<0
\end{align*}
In a similar way with $U_1$ and $U_3$, we see that $v \in H^0\left( \mathbb{F}_m, \mathcal{O}_{\mathbb{F}_m}(p,q)\right)$ on $U_1$ is of the form
\begin{align}\label{rr4}
v_1(z_1,w_1)= a_1^1(z_1)+ a_1^2(z_1)w_1 + \cdots + a_1^{q+1}(z_1)w_1^q
\end{align}
where $a_1^1(z_1),..., a_1^{q+1}(z_1)$ are entire functions of $z_1$. On the other hand, $v$ on $U_0$ in $(\ref{rr3})$ is transformed on $U_1$ into
\begin{align*}
&h_{10} \left(  a_0^1(z_0)+ a_0^2(z_0)w_0 + \cdots + a_0^{q+1}(z_0) w_0^q \right)= a_1^1(z_1)+ a_1^2(z_1)w_1+ \cdots + a_1^{q+1}(z_1)w_1^q\\
&\iff z_1^{p-mq} \left(  a_0^1\left(\frac{1}{z_1}\right)+ a_0^2\left(\frac{1}{z_1}\right)z_1^m w_1 + \cdots + a_0^{q+1}\left( \frac{1}{z_1}\right) z_1^{qm} w_1^q \right)= a_1^1(z_1)+ a_1^2(z_1)w_1+ \cdots + a_1^{q+1}(z_1)w_1^q
\end{align*}
Then we have
\begin{align} \label{rr5}
a_1^i(z_1) & = z_1^{p-mq} \alpha_0^i\left(\frac{1}{z_1} \right)z_1^{m(i-1)}= z_1^{p-mq + m(i-1)}a_0^i\left( \frac{1}{z_1} \right),\,\,\,\,\,i=1,...,q+1
\end{align}
In the following examples, we will consider the case $q=2$ and $m=2$ i.e. we will describe $H^0\left( \mathbb{F}_2, \mathcal{O}_{\mathbb{F}_2}(p,2) \right)$. Then from $(\ref{rr5})$, we have
\begin{align}
a_1^1(z_1)&= z_1^{p-4}a_0^1\left( \frac{1}{z_1} \right) \notag \\
a_1^2(z_1)&= z_1^{p-2} a_0^2\left( \frac{1}{z_1} \right) \label{rr7} \\
a_1^3(z_1)&= z_1^{p} a_0^3\left( \frac{1}{z_1} \right) \notag
\end{align}

\begin{example}
Let us describe $H^0\left( \mathbb{F}_2, \Omega_{\mathbb{F}_2}^1(3,2) \right)$. Then from $(\ref{rr7})$, we have
\begin{align} \label{rsf23}
H^0\left( \mathbb{F}_2, \mathcal{O}_{\mathbb{F}_2}(3,2)\right)=\left\{ \left(a_0+ a_1z_0\right)w_0 + \left(b_0+ b_1z_0 + b_2 z_0^2+ b_3z_0^3 \right)w_0^2 | a_0,a_1,b_0,..., b_3 \in \mathbb{C} \right\} \,\,\, \textnormal{on}\,\,\, U_0
\end{align}
\end{example}

\begin{example}
Let us describe $H^0\left( \mathbb{F}_2, \Omega_{\mathbb{F}_2}^1(4,2) \right)$. Then from $(\ref{rr7})$, we have
\begin{align}\label{rf1}
H^0\left( \mathbb{F}_2, \mathcal{O}_{\mathbb{F}_2}(4,2)\right)=\left\{ c_0+ \left(a_0+ a_1z_0+ a_2 z_0^2 \right)w_0 + \left(b_0+ b_1z_0 + b_2 z_0^2+ b_3z_0^3+ b_4 z_0^4\right)w_0^2 | a_0,a_1,a_2,b_0,..., b_4, c_0 \in \mathbb{C} \right\} \,\,\, \textnormal{on}\,\,\, U_0
\end{align}
\end{example}

\begin{example}
Let us describe $H^0\left( \mathbb{F}_2, \Omega_{\mathbb{F}_2}^1(5,2) \right)$. Then from $(\ref{rr7})$, $H^0\left( \mathbb{F}_2, \mathcal{O}_{\mathbb{F}_2}(5,2)\right)$ is given by
\begin{align} \label{rf26}
\left\{ c_0+ c_1 z_0+ \left(a_0+ a_1z_0+ a_2 z_0^2 + a_3 z_0^3 \right)w_0 + \left(b_0+ b_1z_0 +\cdots + b_5 z_0^5\right)w_0^2 | a_0,a_1,a_2,b_0,..., b_5, c_0 \in \mathbb{C} \right\} \,\,\,\textnormal{on}\,\,\, U_0
\end{align}
\end{example}

\subsection{Description of $H^1\left(  \mathbb{F}_m, \mathcal{O}_{\mathbb{F}_m}(p,q)\right), q \geq 0$}\

We keep the notations in subsection \ref{rs1}. We explicitly describe $H^1\left( \mathbb{F}_m, \mathcal{O}_{\mathbb{F}_m}(p,q) \right)$. First we show that $H^1\left( U_0 \bigcup U_2, \mathcal{O}_{\mathbb{F}_m}(p,q) \right)=0$ and $H^1\left( U_1 \bigcup U_3, \mathcal{O}_{\mathbb{F}_m}(p,q) \right)=0$ for $q \geq 0$. By letting $\mathscr{B}=\left\{ U_0, U_2 \right\}$, we can compute $H^1\left(U_0\bigcup U_2, \mathcal{O}_{\mathbb{F}_m}(p,q) \right)$ by \v Cech resolution associated to $\mathscr{B}$, i.e.
\begin{align*}
H^1\left(U_0 \bigcup U_2, \mathcal{O}_{\mathbb{F}_m}(p,q) \right) = H^1\left( \mathscr{B}, \mathcal{O}_{\mathbb{F}_m}(p,q) \right)= \frac{Z^1\left(\mathscr{B}, \mathcal{O}_{\mathbb{F}_m}(p,q) \right)}{\delta C^0\left(\mathscr{B}, \mathcal{O}_{\mathbb{F}_m}(p,q) \right)}
\end{align*}
A $1$-cocycle $\{ s_{02}, s_{20}\}$ belonging to $Z^1\left(\mathscr{B}, \mathcal{O}_{\mathbb{F}_m}(p,q)|_{U_0\bigcup U_2} \right)$ consists of holomorphic functions $s_{02}$ and $s_{20}$ on $U_0\cap U_2$ with $s_{02}=-h_{02}s_{20}=- w_0^q s_{20}$, while $C^0\left(\mathscr{B},  \mathcal{O}_{\mathbb{F}_m}(p,q)|_{U_0\cup U_2} \right)$ is the set of all pairs $\{s_0,s_2\}$, where $s_0$ and $s_2$ are holomorphic functions on $U_0$ and $U_2$, respectively. In order to prove $H^1\left( \mathscr{B}, \mathcal{O}_{\mathbb{F}_m}(p,q)|_{U_0\bigcup U_2}\right)=0$, it suffices to that any $1$-cocycle $\{s_{02}, s_{20}\}$ is represented as $\delta\left(\{s_2, s_0\} \right)$, i.e.
\begin{align*}
s_{02}= h_{02}s_2- s_0 = w_{0}^q s_2- s_0
\end{align*}
The holomorphic function $s_{02}$ on $U_0\cap U_2= \mathbb{C}\times \mathbb{C}^*$ is written as
\begin{align*}
s_{02}(z_0,w_0)= \sum_{n=-\infty}^{\infty} u_n(z_0)w_0^n
\end{align*}
If we set
\begin{align*}
s_0= - \sum_{n=0}^\infty u_n(z_0)w_0^n,\,\,\,\,\,\,\,\,\,\,\,\,s_2 = \sum_{n=1}^\infty u_{-n}(z_2) w_2^{n+q}
\end{align*}
Then $s_{02}=  h_{02}s_2- s_0$ if $q\geq -1$, so that $H^1 \left( U_0 \bigcup U_2 , \mathcal{O}_{\mathbb{F}_m}(p,q)\right)=0$ for $q \geq -1$. Since we describe $H^1\left( \mathbb{F}_m,  \mathcal{O}_{\mathbb{F}_m}(p,q)\right)$ for the study of unfoldings, and there is no foliation $\omega\in H^0\left( \mathbb{F}_m , \Omega_{\mathbb{F}_m}^1(p,q)\right)$ for $q <0$, we may assume that $q\geq 0$ in the following. Similarly we can show that $H^1\left( U_1\bigcup U_3, \mathcal{O}_{\mathbb{F}_m}(p,q) \right)=0$ for $q \geq 0$. Then by letting $\mathfrak{B}=\left\{ U_0 \bigcup U_2 , U_1 \bigcup U_3\right\}$, we can compute $H^1\left( \mathbb{F}_m, \mathcal{O}_{\mathbb{F}_m}(p,q) \right), q \geq 0$ by \v Cech resolution associated to $\mathfrak{B}$, i.e.
\begin{align}\label{rr6}
H^1\left( \mathbb{F}_m, \mathcal{O}_{\mathbb{F}_m}(p,q) \right)= H^1\left( \mathfrak{B}, \mathcal{O}_{\mathbb{F}_m} (p,q)  \right)= \frac{Z^1\left( \mathfrak{B}, \mathcal{O}_{\mathbb{F}_m}(p,q) \right)}{\delta C^0\left(\mathfrak{B}, \mathcal{O}_{\mathbb{F}_m}(p,q) \right)}
\end{align}

Next we describe explicitly $H^0\left( U_0 \bigcup U_2, \mathcal{O}_{\mathbb{F}_m}(p,q) \right), q\geq 0$ on $U_0$ and $H^0\left( U_1 \bigcup U_3, \mathcal{O}_{\mathbb{F}_m}(p,q) \right), q \geq 0$ on $U_1$. Let us describe $H^0\left( U_0 \bigcup U_2, \mathcal{O}_{\mathbb{F}_m}(p,q) \right)$ on $U_0$ which is given by $(\ref{rr3})$
\begin{align}\label{rr8}
a_0^1(z_0)+ a_0^2(z_0)w_0 + \cdots + a_0^{q+1}(z_0)w_0^q
\end{align}
where $a_0^1(z_0),..., a_0^{q+1}(z_0)$ are entire functions of $z_0$. On the other hand, let us describe $H^0\left( U_1 \bigcup U_3, \mathcal{O}_{\mathbb{F}_m} (p,q) \right)$ on $U_1$ which is given by $(\ref{rr4})$
\begin{align}\label{rr9}
a_1^1(z_1)+ a_1^2(z_1)w_1 + \cdots + a_1^{q+1}(z_1)w_1^q
\end{align}
where $a_1^1(z_1),..., a_1^{q+1}(z_1)$ are entire functions of $z_1$.

Now let us describe $H^1\left( \mathbb{F}_m, \mathcal{O}_{\mathbb{F}_m}(p,q) \right)$ in terms of \v Cech resolution associated to $\mathfrak{B}$ from $(\ref{rr6})$. We represent a $1$-cocycle $\left\{s_{01}, s_{10}\right\}\in Z^1\left( \mathfrak{B}, \mathcal{O}_{\mathbb{F}_m}(p,q) \right)$ with $h_{20}s_{02}=- s_{20}$. Then the $1$-cocycle belongs to $\delta C^0\left( \mathfrak{B}, \mathcal{O}_{\mathbb{F}_m}(p,q) \right)$ if and only if there exists $s_0\in \Gamma\left(U_0\bigcup U_2, \mathcal{O}_{\mathbb{F}_m}(p,q) \right)$ and $s_1\in \Gamma\left(U_1\bigcup U_3, \mathcal{O}_{\mathbb{F}_m}(p,q) \right)$ such that $s_{01}= h_{01} s_1- s_0= z_0^{p-mq} s_1 - s_0$, where $s_0$ is written in the form as in $(\ref{rr8})$ and $s_1$ is written in the form as in $(\ref{rr9})$. We write $s_{01}$ in terms of the coordinate $(z_0,w_0)$
\begin{align*}
s_{01} = \gamma_{01}^1(z_0)+ \gamma_{01}^2(z_0)w_0 + \cdots + \gamma_{01}^{q+1}(z_0)w_0^q
\end{align*}
where $\gamma_{01}^1(z_0),..., \gamma_{01}^{q+1}(z_0)$ are holomorphic functions of $z_0$ on $\mathbb{C}^*$, so that they are expanded into Laurent series in $z_0$. We note that
\begin{align*}
h_{01}s_1& = z_0^{p-mq}\left( \alpha_1^1(z_1)+ \alpha_1^2(z_1)w_1 + \cdots + \alpha_1^{q+1}(z_1)w_1^q\right) \\
 &= z_0^{p-mq}\left( \alpha_1^1\left( \frac{1}{z_0} \right) + \alpha_1^2\left( \frac{1}{z_0} \right)z_0^m w_0 + \cdot + \alpha_1^{q+1}\left( \frac{1}{z_0} \right) z_0^{mq} w_0^q \right)
\end{align*}
Then the equation $s_{01}= z_0^{p-mq}s_1 - s_0$ is reduced to the following system of equations:
\begin{align}\label{rr10}
z_0^{p-mq+ m(i-1)}\alpha_1^i\left( \frac{1}{z_0} \right) - \alpha_0^i(z_0)= \gamma_{01}^i(z_0),\,\,\,\,\,\,\,\,1\leq i \leq q+1
\end{align}
In particular, if $p\geq mq-1$, then $(\ref{rr10})$ has always a solution, so that
\begin{align} \label{rsf7}
H^1\left( \mathbb{F}_m, \mathcal{O}_{\mathbb{F}_m}(p,q)  \right)=0,\,\,\,\,\,\,\,\textnormal{for} \,\,\,  p\geq mq-1
\end{align}

\subsection{Description of $H^0\left( \mathbb{F}_m, \bigwedge^2 \Omega_{\mathbb{F}_m}^1(2p,2q) \right)$}\

We keep the notations in subsection \ref{rs1}. We explicitly describe the global section $H^0\left( \mathbb{F}_m, \bigwedge^2 \Omega_{\mathbb{F}_m}(2p,2q) \right)$ on $U_0$. Let $\omega \in H^0\left( \mathbb{F}_m , \bigwedge^2 \Omega_{\mathbb{F}_m}(2p,2q) \right)$ be represented by $a_0(z_0,w_0) dz_0 \wedge dw_0$ on $U_0$ and $a_2(z_2, w_2) dz_2 \wedge dw_2$ on $U_2$, respectively. From $dw_2= - \frac{1}{w_0^2} dw_0$ and $dz_2 = dz_0$, we have
\begin{align*}
a_0(z_0,w_0)dz_0\wedge dw_0= w_0^{2q} a_2(z_2,w_2)dz_2 \wedge dw_2 = -w_0^{2q-2} a_2\left(z_0,  \frac{1}{w_0} \right)dz_0 \wedge dw_0
\end{align*}
Then $a_0(z_0,w_0)$ is of the form
\begin{align*}
a_0(z_0,w_0)= a_0^1(z_0)+ a_0^2(z_0)w_0 + \cdots + a_0^{2q-1}(z_0)w_0^{2q-2}
\end{align*}
where $a_0^1(z_0),..., a_0^{2q-1}(z_0)$ are entire functions of $z_0$. This implies that if $q<1$, then
\begin{align} \label{rsf6}
H^0\left(\mathbb{F}_2, \Omega_{\mathbb{F}_2}^1(2p,2q) \right) = 0,\,\,\,\,\,\,\,\,\,q<1
\end{align}
Now we assume that $q\geq 1$. In a similar way with $U_1$ and $U_3$, we see that $\omega\in H^0\left( \mathbb{F}_m, \bigwedge^2 \Omega_{\mathbb{F}_m}^1 (2p,2q) \right)$ is of the form $a_1(z_1,w_2)dz_1\wedge dw_1$ where
\begin{align*}
a_1(z_1,w_1)= a_1^1(z_1)+ a_1^2(z_1)w_1 + \cdots + a_1^{2q-1}(z_1)w_1^{2q-2}
\end{align*}
where $a_1^1(z_1),..., a_1^{2q-1}(z_1)$ are entire functions of $z_1$. Since $z_1=\frac{1}{z_0}$ and $w_1= z_0^m w_0$, so that $dz_1= -\frac{1}{z_0^2}dz_0$ and $dw_1= m z_0^{m-1}w_0dz_0 + z_0^m dw_0$, and so $dz_1\wedge dw_1= - z_0^{m-2} dz_0\wedge dw_0$, we have on $U_0\bigcap U_1$,
\begin{align*}
&\left(a_1^1(z_1)+ a_1^2(z_1)w_1 + \cdots + a_1^{2q-1}(z_1)w_1^{2q-2}\right) dz_1 \wedge dw_1\\
&=h_{10} \cdot \left( a_0^1(z_0)+ a_0^2(z_0)w_0 + \cdots + a_0^{2q-1}(z_0)w_0^{2q-2} \right) dz_0 \wedge dw_0 \\
&= -z_1^{2p-2mq+ m-2}\left( a_0^1\left(\frac{1}{z_1}\right)+ a_0^2\left( \frac{1}{z_1} \right)z_1^m w_1 + \cdots + a_0^{2q-1}\left(\frac{1}{z_1} \right) z_1^{m(2q-2)} w_0^{2q-2} \right) dz_1 \wedge dw_1
\end{align*}
Then we have
\begin{align} \label{rr11}
a_1^i(z_1)= - z_1^{2p-2mq+ m-2+m(i-1)} a_0^i\left(\frac{1}{z_1} \right)=- z_1^{2p-2 + m(i-2q)}a_0^i\left(\frac{1}{z_1} \right),\,\,\,\,\,\,\,\,i=1,...,2q-1
\end{align}

\begin{remark}
Since $K_{\mathbb{F}_m}\cong \mathcal{O}_{\mathbb{F}_m}(-m-2, -2)$, we have $\bigwedge^2 \Omega_{\mathbb{F}_m}^1(2p,2q)\cong \mathcal{O}_{\mathbb{F}_m}(2p-m-2, 2q-2)$, so that $H^0\left( \mathbb{F}_m, \bigwedge^2 \Omega_{\mathbb{F}_m}^1(2p,2q) \right)\cong H^0\left( \mathbb{F}_m, \mathcal{O}_{\mathbb{F}_m}(2p-m-2, 2q -2 \right)$, which is described in $(\ref{rr5})$ by replacing $p$ by $2p-m-2$ and replacing $q$ by $2q-2$, i.e. we have $a_1^i(z_1)= z_1^{2p-m-2 - m(2q-2)+ m(i-1)} a_0^i\left( \frac{1}{z_1} \right),i=1,...,2q-2+1$ which is equivalent to $(\ref{rr11})$.
\end{remark}

In the following examples we consider $m=2$ and $q=2$. Then from $(\ref{rr11})$ we have
\begin{align}
a_1^1(z_1)&= - z_1^{2p-8}a_0^1\left( \frac{1}{z_1} \right) \notag \\
a_1^2(z_1)&= - z_1^{2p-6} a_0^2 \left( \frac{1}{z_1} \right) \label{rr12}\\
a_1^3(z_1)&= - z_1^{2p-4} a_0^3 \left( \frac{1}{z_1} \right) \notag 
\end{align}
\begin{example}
Let us describe $H^0\left( \mathbb{F}_2, \bigwedge^2 \Omega_{\mathbb{F}_2}^1(2\cdot 3, 2 \cdot 2) \right)$. From $(\ref{rr12})$, we have
\begin{align}\label{rsf22}
H^0\left( \mathbb{F}_2, \bigwedge^2 \Omega_{\mathbb{F}_2}^1(6, 4) \right)=\left\{ \left(  a_0 w_0 + \left( b_0 + b_1 z_0 + b_2 z_0^2\right) w_0^2 \right) dz_0 \wedge dw_0 | a_0,b_0,b_1,b_2 \in \mathbb{C} \right\}\,\,\,\,\,\textnormal{on}\,\,\, U_0
\end{align}
\end{example}

\begin{example}
Let us describe $H^0\left( \mathbb{F}_2, \bigwedge^2 \Omega_{\mathbb{F}_2}^1(2\cdot 4, 2 \cdot 2) \right)$. From $(\ref{rr12})$, we have
\begin{align}\label{rf2}
H^0\left( \mathbb{F}_2, \bigwedge^2 \Omega_{\mathbb{F}_2}^1(8, 4) \right)=\left\{ \left( c_0+  \left(a_0 + a_1 z_0 + a_2 z_0^2\right)w_0 + \left( b_0 + b_1 z_0 + b_2 z_0^2 + b_3 z_0^3 + b_4 z_0^4\right) w_0^2 \right) dz_0 \wedge dw_0\right\}\,\,\,\,\,\textnormal{on}\,\,\, U_0
\end{align}
where $a_0,a_1,a_2,b_0,b_1,b_2,b_3,b_4, c_0 \in \mathbb{C}$.
\end{example}

\begin{example}
Let us describe $H^0\left( \mathbb{F}_2, \bigwedge^2 \Omega_{\mathbb{F}_2}^1( 2 \cdot 5 , 2\cdot 2) \right)$. From $(\ref{rr12})$, $H^0\left( \mathbb{F}_2, \bigwedge^2 \Omega_{\mathbb{F}_2}^1 (10,4)  \right)$ is given by
\begin{align}\label{rf21}
\left\{ \left( c_0+ c_1 z_0 + c_2 z_0^2 + \left(a_0 + a_1 z_0 + a_2 z_0^2 + a_3 z_0^3 + a_4 z_0^4\right)w_0 + \left( b_0 + b_1 z_0  +  \cdots + b_6 z_0^6\right) w_0^2 \right) dz_0 \wedge dw_0\right\}\,\,\,\,\,\textnormal{on}\,\,\, U_0
\end{align}
where $a_0,a_1,a_2,a_3,a_4, b_0,b_1,b_2,b_3,b_4,b_5,b_6, c_0 ,c_1,c_2\in \mathbb{C}$.
\end{example}

With this preparation we study unfoldings of $\mathbb{F}_m$.

\subsection{Unfoldigs of $\left( \mathbb{F}_m, \mathcal{N}_{\mathcal{F}_0}^* \right)$}\

 Let us compute the cohomology group $\mathbb{H}^i\left( \mathbb{F}_m, \tilde{\mathcal{E}}_{\mathcal{O}_{\mathbb{F}_2}(p,q)}^\bullet  \right),i=0,1,2$ of the unfolded dual leaf complex $(\ref{rs15})$. Consider the spectral sequence at $E_1$ which is given by
{\small{\begin{equation}
\begin{CD}
0 \\
@AAA \\
H^0\left( \mathbb{F}_m,  \bigwedge^2 \Omega_{\mathbb{F}_m}^1(2p,2q) \right) @.\,\,\,\,\,\,\,\,\,\,\,\,\,H^1\left( \mathbb{F}_m, \bigwedge^2 \Omega_{\mathbb{F}_m}^1 (2p,2q) \right) \\
@AAA @AAA \\
H^0\left( \mathbb{F}_m,  \Omega_{\mathbb{F}_m}^1(p, q) \right) \bigoplus H^0\left( \mathbb{F}_m, \mathcal{O}_{\mathbb{F}_m}(p,q) \right) @. \,\,\,\,\,\,\,\,\,\,\,\,\,\,\,\,\,\,\,\, H^1\left( \mathbb{F}_m, \Omega_{\mathbb{F}_m}^1(p,q) \right) \bigoplus H^1\left( \mathbb{F}_m, \mathcal{O}_{\mathbb{F}_m}(p,q) \right) \\
@AAA @AAA \\
H^0\left( \mathbb{F}_m, \mathcal{E}_{\mathcal{O}_{\mathbb{F}_m}(p,q)} \right) @.\,\,\,\,\,\,\,\,\,\, \,\,\,\,\, H^1\left(\mathbb{F}_m, \mathcal{E}_{\mathcal{O}_{\mathbb{F}_m}(p,q)} \right) @. \,\,\,\,\,\,\,\,\,\, \,\,\,\,\,\,\,\,\,\,\,\,\,\,\,\,\, 0
\end{CD}
\end{equation}}}

\begin{example}
Let us consider a foliated rational ruled surface $\left( \mathbb{F}_2, \omega_0 \right)$, where $\omega_0=  2 w_0 dz_0+ 2 z_0 dw_0 \in H^0\left( \mathbb{F}_2, \Omega_{\mathbb{F}_2}^1 (3,2) \right)$ on $U_0$ as in \textnormal{Example \ref{rsf14}}. We keep the notations in $\textnormal{Example \ref{rsf14}}$. We have from $(\ref{rsf1})$, $(\ref{rsf17}),(\ref{rsf3}),(\ref{rsf18}),(\ref{rsf23}),(\ref{rsf7})$ and $(\ref{rsf22})$
{\small{\begin{align*}
&H^0\left( \mathbb{F}_2,  \mathcal{E}_{\mathcal{O}_{\mathbb{F}_2}(3,2)} \right)  \\
&= \left\{ \left(a_0+ a_1z_0+ a_2z_0^2\right)\frac{\partial}{\partial z_0}+ \left( \left(a_3-2a_2z_0\right)w_0+ \left(c_0+c_1z_0+ c_2 z_0^2\right)w_0^2 \right)\frac{\partial}{\partial w_0} + a_4 +a_2z_0 -2\left(c_0+c_1z_0+ c_2z_0^2\right)w_0^2    \right\} \\
&H^0\left( \mathbb{F}_2, \Omega_{\mathbb{F}_2}^1(3,2)\right) = \left\{ \left(2f_1w_0+ \left(g_0+ g_1z_0\right)w_0^2\right)dz_0 + \left(f_0+ f_1z_0\right)dw_0  \right\} \\
&H^1\left( \mathbb{F}_2,  \mathcal{E}_{\mathcal{O}_{\mathbb{F}_2}(3,2)} \right) = \left\langle \frac{1}{z_0}\frac{\partial}{\partial w_0}  \right\rangle  \\
&H^1\left( \mathbb{F}_2, \Omega_{\mathbb{F}_2}^1(3,2)  \right) = \left\langle  \frac{1}{z_0^2} dz_0 , \frac{1}{z_0}dz_0 \right\rangle \\
&H^0\left( \mathbb{F}_2, \mathcal{O}_{\mathbb{F}_2}(3,2)   \right)  =\left\{ \left(a_0+ a_1z_0\right)w_0 + \left(b_0 + b_1z_0 + b_2 z_0^2+ b_3z_0^3 \right) w_0^2 \right\}\\
&H^1\left( \mathbb{F}_2, \mathcal{O}_{\mathbb{F}_2}(3,2) \right)=0 \\
&H^0\left(\mathbb{F}_2, \bigwedge^2  \Omega_{\mathbb{F}_2}^1 (6,4) \right) = \left\{  \left( a_0w_0 + \left( b_0 +  b_1z_0 + b_2z_0^2 \right)w_0^2 \right) dz_0 \wedge dw_0  \right\}
\end{align*}}}
We compute $H^0\left(\mathbb{F}_2, \mathcal{E}_{\mathcal{O}_{\mathbb{F}_2}(3,2)} \right)\to H^0\left( \mathbb{F}_2, \Omega_{\mathbb{F}_2}^1(3,2) \right) \bigoplus H^0\left( \mathbb{F}_2, \mathcal{O}_{\mathbb{F}_2}(3,2) \right)$. First we note that the image of $H^0\left( \mathbb{F}_2, \mathcal{E}_{\mathcal{O}_{\mathbb{F}_2}(3,2)} \right)\to H^0\left( \mathbb{F}_2, \Omega_{\mathbb{F}_2}^1(3,2) \right)$ is given by $(\ref{rsf20})$. 
\begin{align*}
&\left(  \left(2a_3+ 2a_1+ 2a_4 \right)w_0 +\left(-2c_0 - c_1z_0\right) w_0^2    \right)dz_0 +\left( a_0+ \left(a_1+ a_3 + a_4  \right)z_0\right)dw_0 = 0 \\
& \iff a_1+a_3+a_4=0, a_0=0,c_0=c_1=0
\end{align*}
Next we compute $H^0\left(\mathbb{F}_2,  \mathcal{E}_{\mathcal{O}_{\mathbb{F}_2}(3,2)} \right)\to H^0\left( \mathbb{F}_2, \Omega_{\mathbb{F}_2}^1(3,2) \right)$
\begin{align*}
&i_{\left(a_0+ a_1z_0+ a_2z_0^2\right) \frac{\partial}{\partial z_0} + \left(\left(-2a_2z_0+ a_3\right)w_0 + \left(c_0+c_1z_0+ c_2z_0^2\right)w_0^2 \right)\frac{\partial}{\partial w_0}}\left( 2w_0dz_0 + z_0dw_0 \right)\\
&= \left( 2a_0 + \left(2a_1+ a_3\right)z_0 \right) w_0 + \left(c_0z_0+ c_1z_0^2+ c_2z_0^3 \right)w_0^2 =0 \iff a_0=0, 2a_1+ a_3=0, c_0=c_1=c_2=0
\end{align*}
Then we have
\begin{align*}
\dim_\mathbb{C} \mathbb{H}^0\left( \mathbb{F}_2, \tilde{\mathcal{E}}_{\mathcal{O}_{\mathbb{F}_2}(3,2)}^\bullet \right) = \dim_\mathbb{C}\textnormal{ker}\left( H^0\left(\mathbb{F}_2, \mathcal{E}_{\mathcal{O}_{\mathbb{F}_2}(3,2)} \right)\to H^0\left( \mathbb{F}_2, \Omega_{\mathbb{F}_2}^1(3,2) \right) \bigoplus H^0\left( \mathbb{F}_2, \mathcal{O}_{\mathbb{F}_2}(3,2) \right) \right) = 2
\end{align*}

On the other hand, we compute $H^0\left( \mathbb{F}_2, \Omega_{\mathbb{F}_2}^1(3,2) \right) \bigoplus H^0\left( \mathbb{F}_2, \mathcal{O}_{\mathbb{F}_2}(3,2) \right)\to H^0\left( \mathbb{F}_2, \bigwedge^2 \Omega_{\mathbb{F}_2}^1(6,4) \right)$.
\begin{align*}
&\omega_0 \wedge d\left(  \left( a_0+ a_1z_0 \right) w_0 + \left(b_0 + b_1z_0 + b_2 z_0^2+ b_3z_0^3\right)w_0^2 \right) + \left( \left(a_0+ a_1z_0\right)w_0 + \left(b_0 + b_1z_0 + b_2 z_0^2+ b_3z_0^3\right)w_0^2 \right)d\omega_0 \\
& - \omega_0 \wedge \left( \left(2f_1w_0+ \left(g_0+ g_1z_0\right)w_0^2\right)dz_0 + \left(f_0+ f_1z_0\right)dw_0 \right) \\
&=\left( \left(a_0 - 2f_0 \right)w_0 + \left(3b_0+  \left(2b_1+ g_0 \right)z_0 + \left(b_2+ g_1\right)z_0^2  \right)w_0^2 \right) dz_0\wedge dw_0 =0  \\
&\iff a_0-2f_0=0, b_0=0, 2b_1+ g_0=0, b_2+ g_1=0
\end{align*}
Then we have $\dim_\mathbb{C} \textnormal{ker}\left( H^0\left( \mathbb{F}_2, \Omega_{\mathbb{F}_2}^1(3,2) \right) \bigoplus H^0\left( \mathbb{F}_2, \mathcal{O}_{\mathbb{F}_2}(3,2) \right)\to H^0\left( \mathbb{F}_2, \bigwedge^2 \Omega_{\mathbb{F}_2}^1(6,4) \right) \right)= 6$, so that we have
\begin{align*}
\dim_\mathbb{C} \frac{\textnormal{ker}\left( H^0\left( \mathbb{F}_2, \Omega_{\mathbb{F}_2}^1(3,2) \right) \bigoplus H^0\left( \mathbb{F}_2, \mathcal{O}_{\mathbb{F}_2}(3,2) \right)\to H^0\left( \mathbb{F}_2, \bigwedge^2 \Omega_{\mathbb{F}_2}^1(6,4) \right) \right) }{\textnormal{im}\left( H^0\left(\mathbb{F}_2, \mathcal{E}_{\mathcal{O}_{\mathbb{F}_2}(3,2)} \right)\to H^0\left( \mathbb{F}_2, \Omega_{\mathbb{F}_2}^1(3,2) \right) \bigoplus H^0\left( \mathbb{F}_2, \mathcal{O}_{\mathbb{F}_2}(3,2) \right)  \right)} = 6 - (8- 2 ) = 0
\end{align*}
We also note that $\textnormal{ker}\left( H^1\left( \mathbb{F}_2, \mathcal{E}_{\mathcal{O}_{\mathbb{F}_2}(3,2)}^\bullet \right)\to H^1\left( \mathbb{F}_2, \Omega_{\mathbb{F}_2}^1(3,2) \right) \bigoplus H^1\left( \mathbb{F}_2, \mathcal{O}_{\mathbb{F}_2}(3,2) \right) \right)=0$ by $(\ref{rsf21})$, so that
\begin{align*}
\mathbb{H}^1\left( \mathbb{F}_2, \tilde{\mathcal{E}}_{\mathcal{O}_{\mathbb{F}_2}(3,2)}^\bullet \right) =  0 
\end{align*}
Hence $\left( \mathbb{F}_2, \omega_0 \right)$ is rigid in unfoldings by \textnormal{Theorem \ref{r1}}.

\end{example}

\begin{example}
Let us consider a foliated rational ruled surface $\left( \mathbb{F}_2, \omega_0 \right)$, where $\omega_0 = w_0 dz_0 + z_0 dw_0 \in H^0\left( \mathbb{F}_2, \Omega_{\mathbb{F}_2}^1(4,2) \right)$ on $U_0$ as in \textnormal{Example \ref{rsf15}}. We keep the notations in $\textnormal{Example \ref{rsf15}}$. We have from $(\ref{rsf1}),(\ref{rsf31})$, $(\ref{rsf3}), (\ref{rsf32}), (\ref{rf1}), (\ref{rsf7})$ and $(\ref{rf2})$
{\small{\begin{align}
&H^0\left( \mathbb{F}_2,  \mathcal{E}_{\mathcal{O}_{\mathbb{F}_2}(4,2)} \right)  \\
&= \left\{ \left(a_0+ a_1z_0+ a_2z_0^2\right)\frac{\partial}{\partial z_0}+ \left( \left(d-2a_2z_0\right)w_0+ \left(c_0+c_1z_0+ c_2 z_0^2\right)w_0^2\right)\frac{\partial}{\partial w_0} + b  -2\left(c_0+c_1z_0+ c_2z_0^2\right)w_0^2    \right\} \notag \\
&H^0\left( \mathbb{F}_2, \Omega_{\mathbb{F}_2}^1\left( 4, 2 \right) \right)= \left\{ \left( \left( c+ 2 a_2 z_0 \right)  w_0 +  \left(b_0+ b_1 z_0+ b_2 z_0^2 \right) w_0^2 \right) dz_0 +   \left(a_0 + a_1 z_0+ a_2 z_0^2  \right)dw_0  \right\}\\
&H^1\left( \mathbb{F}_2, \mathcal{E}_{\mathcal{O}_{\mathbb{F}_2}(4,2)} \right) = \left\langle \frac{1}{z_0}\frac{\partial}{\partial w_0} \right\rangle\\
&H^1\left( \mathbb{F}_2, \Omega_{\mathbb{F}_2}^1 (4,2) \right) = \left\langle  \frac{1}{z_0} d z_0 \right\rangle\\
&H^0\left( \mathbb{F}_2, \mathcal{O}_{\mathbb{F}_2}(4,2) \right) = \left\{ e + \left(g_0+ g_1 z_0 + g_2 z_0^2\right)  w_0 + \left( h_0 + h_1 z_0 + h_2 z_0^2 + h_3 z_0^3 + h_4 z_0^4 \right)   \right\}\\
&H^1\left( \mathbb{F}_2, \mathcal{O}_{\mathbb{F}_2}(4,2) \right) = 0 \\
&H^0\left( \mathbb{F}_2, \bigwedge^2 \Omega_{\mathbb{F}_2}^1(8, 4) \right) = \left\{ \left( c_0+  \left(a_0 + a_1 z_0 + a_2 z_0^2\right)w_0 + \left( b_0 + b_1 z_0 + b_2 z_0^2 + b_3 z_0^3 + b_4 z_0^4\right) w_0^2 \right) dz_0 \wedge dw_0\right\}
\end{align}}}
We compute $H^0\left( \mathbb{F}_2,  \mathcal{E}_{\mathcal{O}_{\mathbb{F}_2}(4,2)} \right) \to H^0\left( \mathbb{F}_2,  \Omega_{\mathbb{F}_2}^1(4,2) \right) \bigoplus H^0\left( \mathbb{F}_2, \mathcal{O}_{\mathbb{F}_2} (4,2) \right)$. First we note that the image of $ H^0\left( \mathbb{F}_2, \mathcal{E}_{\mathcal{O}_{\mathbb{F}_2}(4,2)} \right) \to H^0\left( \mathbb{F}_2, \Omega_{\mathbb{F}_2}^1(4,2) \right)$ is given by $(\ref{rsf34})$.
\begin{align}
&\left( \left( - 2a_2 z_0+ d+b + a_1\right) w_0 + \left( - c_0 + c_2 z_0^2\right) w_0^2 \right) dz_0  +\left( a_0 + (a_1+d +b) z_0 - a_2 z_0^2   \right) dw_0  = 0  \\
 &\iff  a_2=0, a_1+ d+ b=0,  c_0=c_2=0, a_0=0 \notag
\end{align}
Next we compute $H^0\left( \mathbb{F}_2, \mathcal{E}_{\mathcal{O}_{\mathbb{F}_2}(4,2)} \right)\to H^0\left( \mathbb{F}_2, \mathcal{O}_{\mathbb{F}_2}(4,2) \right)$
\begin{align}
&i_{\left( a_0 + a_1 z_0 + a_2 z_0^2 \right) \frac{\partial}{\partial z_0} +\left( \left( - 2a_2 z_0+ d\right) w_0 + \left( c_0+ c_1 z_0 + c_2 z_0^2\right)  w_0^2 \right)\frac{\partial}{\partial w_0} }\left( w_0 dz_0 + z_0 dw_0 \right) \notag\\
&= \left( a_0 + \left(a_1+d\right) z_0 - a_2 z_0^2\right) w_0 + \left( c_0 z_0 + c_1 z_0^2 + c_2 z_0^3 \right) w_0^2=0 \label{rf5}\\
&\iff a_0=a_2=0, a_1+d=0, c_0=c_1=c_2=0 \notag
\end{align}
Then we have
\begin{align*}
\dim_\mathbb{C} \mathbb{H}^0\left( \mathbb{F}_2, \tilde{\mathcal{E}}_{\mathcal{O}_{\mathbb{F}_2}(4,2)}^\bullet \right) =\dim_\mathbb{C} \textnormal{ker}\left( H^0\left( \mathbb{F}_2,  \mathcal{E}_{\mathcal{O}_{\mathbb{F}_2}(4,2)} \right) \to H^0\left( \mathbb{F}_2,  \Omega_{\mathbb{F}_2}^1(4,2) \right) \bigoplus H^0\left( \mathbb{F}_2, \mathcal{O}_{\mathbb{F}_2} (4,2) \right) \right) = 1
\end{align*}
On the other hand, we compute $ H^0\left( \mathbb{F}_2,  \Omega_{\mathbb{F}_2}^1(4,2) \right) \bigoplus H^0\left( \mathbb{F}_2, \mathcal{O}_{\mathbb{F}_2} (4,2) \right) \to H^0\left( \mathbb{F}_2, \bigwedge^2 \Omega_{\mathbb{F}_2}^1(8,4) \right) $.
\begin{align}
&\omega_0 \wedge d\left( e_0+ \left( g_0 + g_1 z_0 + g_2 z_0^2 \right) w_0 + \left(h_0 + h_1 z_0 + h_2 z_0^2 + h_3 z_0^3 + h_4 z_0^4 \right) w_0^2 \right) \notag \\
&+ \left( e_0+ \left( g_0 + g_1 z_0 + g_2 z_0^2 \right) w_0 + \left(h_0 + h_1 z_0 + h_2 z_0^2 + h_3 z_0^3 + h_4 z_0^4\right) w_0^2 \right) d\omega_0 \notag\\
&- \omega_0 \wedge \left(  \left( \left( c+ 2 a_2 z_0 \right)  w_0 +  \left(b_0+ b_1 z_0+ b_2 z_0^2 \right) w_0^2 \right) dz_0 +   \left(a_0 + a_1 z_0+ a_2 z_0^2  \right)dw_0 \right) \notag \\
&= \left(\left( g_0 - a_0 + \left(c- a_1\right) z_0 + \left(a_2- g_2 \right) z_0^2\right) w_0 + \left( 2 h_0 + (b_0 + h_1) z_0 + b_1 z_0^2+ (b_2- h_3)z_0^3 - 2h_4 z_0^4 \right) w_0^2\right) dz_0 \wedge dw_0 = 0 \label{rf6} \\
&\iff  g_0= a_0, c=a_1, a_2=g_2, h_0=0, b_0+ h_1=0, b_1=0, b_2=h_3, h_4=0 \notag
\end{align}
Then $\textnormal{ker}\left( H^0\left( \mathbb{F}_2,  \Omega_{\mathbb{F}_2}^1(4,2) \right) \bigoplus H^0\left( \mathbb{F}_2, \mathcal{O}_{\mathbb{F}_2} (4,2) \right) \to H^0\left( \mathbb{F}_2, \bigwedge^2 \Omega_{\mathbb{F}_2}^1(8,4) \right) \right)$ is given by
\begin{align}
&\left( \left( t_2+ 2 t_3 z_0 \right) w_0 +  ( t_4  + t_5 z_0^2 ) w_0^2 \right) dz_0 + \left( t_1 + t_2 z_0 + t_3 z_0^2 \right) dw_0 \label{rf3} \\
&+ t_6 + \left( t_1 + t_7 z_0 + t_3 z_0^2 \right) w_0 + ( - t_4 z_0 + t_8 z_0^2 + t_5 z_0^3) w_0^2,\,\,\,\,\,\,t_1,...,t_8\in \mathbb{C} \notag
\end{align}
By comparing $(\ref{rf5}),(\ref{rf6})$ and $(\ref{rf3})$, we see that
\begin{align*}
\frac{\textnormal{ker}\left(H^0\left( \mathbb{F}_2, \Omega_{\mathbb{F}_2}^1(4,2)  \right)\bigoplus H^0\left( \mathbb{F}_2, \mathcal{O}_{\mathbb{F}_2}(4,2)\right) \to H^0\left( \mathbb{F}_2, \bigwedge^2 \Omega_{\mathbb{F}_2}^1(8,4) \right) \right)}{\textnormal{im}\left( H^0\left( \mathbb{F}_2, \mathcal{E}_{\mathcal{O}_{\mathbb{F}_2}(4,2)} \right)\to H^0\left( \mathbb{F}_2, \Omega_{\mathbb{F}_2}^1(4,2) \right) \bigoplus H^0\left( \mathbb{F}_2, \mathcal{O}_{\mathbb{F}_2}(4,2) \right) \right)} \cong \textnormal{Span}_\mathbb{C} \left\langle \left( 0=0 \cdot \frac{\partial}{\partial z_0} + 0 \cdot \frac{\partial}{\partial w_0}  ,1 \right) \right\rangle 
\end{align*}

From $(\ref{rf7})$ and 
\begin{align*}
i_{\frac{1}{z_0}\frac{\partial}{\partial w_0}}\left( w_0 dz_0 + z_0 dw_0 \right)= 1 = 0 - (-1),\,\,\,\,\,\textnormal{where $-1$ on $U_0$, and $0$ on $U_1$}
\end{align*}
and $-1\in H^0\left( \mathbb{F}_2, \Omega_{\mathbb{F}_2}^1(4,2) \right)$, we see that the map in the spectral sequence
{\tiny{\begin{align*}
\textnormal{ker}\left(H^1\left( \mathbb{F}_2, \Omega_{\mathbb{F}_2}^1(4,2) \right)\to H^1\left( \mathbb{F}_2, \Omega_{\mathbb{F}_2}^1(4,2) \right) \bigoplus H^1\left( \mathbb{F}_2, \mathcal{O}_{\mathbb{F}_2}(4,2) \right)  \right) \to \textnormal{coker}\left( H^0\left( \mathbb{F}_2, \Omega_{\mathbb{F}_2}^1(4,2) \right)\bigoplus H^0\left( \mathbb{F}_2, \mathcal{O}_{\mathbb{F}_2}(4,2) \right)\to H^0\left( \mathbb{F}_2, \bigwedge^2 \Omega_{\mathbb{F}_2}^1(8,4) \right)\right)
\end{align*}}}
 is the zero map, so that we have
\begin{align*}
\mathbb{H}^1\left( \mathbb{F}_2, \tilde{\mathcal{E}}_{\mathcal{O}_{\mathbb{F}_2}(4,2)}^\bullet  \right)&\cong  \frac{\textnormal{ker}\left(H^0\left( \mathbb{F}_2, \Omega_{\mathbb{F}_2}^1(4,2)  \right)\bigoplus H^0\left( \mathbb{F}_2, \mathcal{O}_{\mathbb{F}_2}(4,2)\right) \to H^0\left( \mathbb{F}_2, \bigwedge^2 \Omega_{\mathbb{F}_2}^1(8,4) \right) \right)}{\textnormal{im}\left( H^0\left( \mathbb{F}_2, \mathcal{E}_{\mathcal{O}_{\mathbb{F}_2}(4,2)} \right)\to H^0\left( \mathbb{F}_2, \Omega_{\mathbb{F}_2}^1(4,2) \right) \bigoplus H^0\left( \mathbb{F}_2, \mathcal{O}_{\mathbb{F}_2}(4,2) \right) \right)} \\
&\,\,\,\,\,\,\,\,\,\bigoplus \textnormal{ker}\left(H^1\left( \mathbb{F}_2, \Omega_{\mathbb{F}_2}^1(4,2) \right)\to H^1\left( \mathbb{F}_2, \Omega_{\mathbb{F}_2}^1(4,2) \right) \bigoplus H^1\left( \mathbb{F}_2, \mathcal{O}_{\mathbb{F}_2}(4,2) \right) \right)  \\
                                       &\cong \textnormal{Span}_\mathbb{C}\left\langle \left( 0, 1  \right) \right\rangle \bigoplus \textnormal{Span}_\mathbb{C}\left\langle \frac{1}{z_0}\frac{\partial}{\partial w_0} \right\rangle
\end{align*}
We will construct an unfolded analytic family of unfoldings of $\left( \mathbb{F}_m, \omega_0 \right)$ such that the unfolded Kodaira-Spencer map is an isomorphism at the distinguished point, so that $\left( \mathbb{F}_m, \omega_0 \right)$ is unobstructed in unfoldings. The unfolded analytic family is
\begin{align*}
\omega&= w_0 dz_0 + z_0 dw_0 - dt_2\\
z_0&= \frac{1}{z_1} \\
w_0 &= z_1^2 w_1 + t_1 z_1 
\end{align*}
We describe the unfolded analytic family more precisely. We have the coordinate transformations as follows$:$
\begin{center}
$\begin{CD}
U_0: \left(z_0,w_0\right)@>>> U_1:\left(z_1, w_1\right)=\left(\frac{1}{z_0}, z_0^2 w_0 -t_1 z_0\right)\\
@VVV @VVV \\
U_2:\left(z_2, w_2\right)= \left(z_0, \frac{1}{w_0}\right) @>>> U_3:\left(z_3,w_3\right)=\left(\frac{1}{z_2}, \frac{w_2}{z_2^2- t_1z_2w_2}\right)=\left(z_1, \frac{1}{w_1}\right)
\end{CD}$
\end{center}
and the foliation $\omega$ is given by
\begin{align*}
\omega&= w_0 dz_0 + z_0 dw_0 - dt_2\,\,\,\,\,\,\,\textnormal{on}\,\,\, U_0,\,\,\,\,\,\,\,\,\,\,\,\,\,\,\,\,\,\,\,\omega= w_1dz_1+z_1dw_1 + dt_1- dt_2\,\,\,\,\,\,\,\textnormal{on}\,\,\, U_1 \\
\omega&= w_2dz_2 - z_2 dw_2 - w_2^2 dt_2\,\,\,\,\,\,\,\textnormal{on}\,\,\, U_2,\,\,\,\,\,\,\,\,\,\,\,\,\,\omega= w_3 dz_3- z_3 dw_3 + w_3^2 dt_1 - w_3^2 dt_2\,\,\,\,\,\,\,\textnormal{on}\,\,\, U_3 
\end{align*}
We note that $\omega \wedge d\omega=0$. On the other hand, transition functions are given by
\begin{align*}
h_{01}&= 1  \\
h_{02}&= w_0^2\\
h_{03}& =\left( z_0^2 w_0 - t_1 z_0\right)^2 \\
h_{12}&= w_2^{-2}\\
h_{13}&= w_1^2\\
h_{23}&= \left( z_2^2- t_1 z_2w_2 \right)^2
\end{align*}
We note that the unfolded Kodaira-Spencer map at $t=(t_1,t_2)=(0,0)=0$ is given by $($in terms of the covering $\mathfrak{B}$ where $\mathfrak{B}=\left\{ U_0\bigcup U_2, U_1 \bigcup U_3 \right\})$
 \begin{align*}
 T_0\left( \mathbb{C}^2 \right) &\to \mathbb{H}^1\left( \mathbb{F}_2, \tilde{\mathcal{E}}_{\mathcal{O}_{\mathbb{F}_2}(4,2)}^\bullet \right)\\
  \frac{\partial}{\partial t_1} &\mapsto \left( \left\{ ( 0 , 0 ), (0 , -1) \right\}     \left\{ \left( \frac{1}{z_0}\frac{\partial}{\partial w_0}, 0   \right) \right\} \right) \in C^0\left( \mathfrak{B}, \Omega_{\mathbb{F}_2}^1(4,2)\bigoplus \mathcal{O}_{\mathbb{F}_2}(4,2) \right) \bigoplus C^1\left( \mathfrak{B}, \mathcal{E}_{\mathcal{O}_{\mathbb{F}_2}(4,2)} \right)\\
  \frac{\partial}{\partial t_2} &\mapsto \left( \left\{(0, 1) , (0,1)\right\} , \left\{ (0,0) \right\}\right) \in C^0\left( \mathfrak{B}, \Omega_{\mathbb{F}_2}^1(4,2)\bigoplus \mathcal{O}_{\mathbb{F}_2}(4,2) \right) \bigoplus C^1\left( \mathfrak{B} , \mathcal{E}_{\mathcal{O}_{\mathbb{F}_2}(4,2)} \right)
 \end{align*}

\end{example}

\begin{example}
Let us consider a foliated rational ruled surface $\left( \mathbb{F}_2, \omega_0 \right)$, where $\omega_0 = z_0w_0 dz_0 - dw_0 \in H^0\left( \mathbb{F}_2, \Omega_{\mathbb{F}_2}^1(5,2) \right)$ on $U_0$ as in \textnormal{Example \ref{rsf16}}. We keep the notations in $\textnormal{Example \ref{rsf16}}$. We have from $(\ref{rsf1}),(\ref{rf13}),(\ref{rsf3})$, $(\ref{rf15})$, $(\ref{rf26}),(\ref{rsf7})$ and $(\ref{rf21})$
\begin{align}
&H^0\left( \mathbb{F}_2, \mathcal{E}_{\mathcal{O}_{\mathbb{F}_2}(5,2)}\right) \\
&= \left\{ \left(a_0+ a_1z_0+ a_2z_0^2\right)\frac{\partial}{\partial z_0}+ \left( \left(a_3-2a_2z_0\right)w_0+ \left(c_0+c_1z_0+ c_2 z_0^2\right)w_0^2 \right)\frac{\partial}{\partial w_0} + a_4 -a_2z_0 -2\left(c_0+c_1z_0+ c_2z_0^2\right)w_0    \right\} \notag \\
&H^0\left( \mathbb{F}_2, \Omega_{\mathbb{F}_2}^1(5,2) \right)  = \left\{ \left( \left( f_0+ f_1 z_0 + 2g_3 z_0^2 \right) w_0 + \left( e_0 + e_1 z_0 + e_2 z_0^2 + e_3 z_0^3 \right) w_0^2\right) dz_0 + \left(g_0+ g_1z_0 + g_2 z_0^2 + g_3 z_0^3 \right) dw_0 \right\}  \\
&H^1\left( \mathbb{F}_2, \mathcal{E}_{\mathcal{O}_{\mathbb{F}_2}(5,2)} \right) = \left\langle \frac{1}{z_0} \frac{\partial}{\partial w_0} \right\rangle\\
&H^1\left( \mathbb{F}_2, \Omega_{\mathbb{F}_2}^1(5,2) \right) = 0\\
&H^0\left( \mathbb{F}_2, \mathcal{O}_{\mathbb{F}_2}(5,2) \right)= \left\{  a_0 +a_1z_0+ \left( b_0 + b_1z_0 + b_2 z_0^2 + b_3z_0^3 \right)w_0 + \left( c_0 + c_1 z_0 + c_2 z_0^2 + c_3 z_0^3 + c_4z_0^4 + c_5 z_0^5 \right)w_0^2 \right\}\\
&H^1\left( \mathbb{F}_2, \mathcal{O}_{\mathbb{F}_2}(5,2) \right) = 0 \\
&H^0\left( \mathbb{F}_2, \bigwedge^2 \Omega_{\mathbb{F}_2}(10,4) \right) \\
&= \left\{ \left( c_0+ c_1 z_0 + c_2 z_0^2 + \left(a_0 + a_1 z_0 + a_2 z_0^2 + a_3 z_0^3 + a_4 z_0^4\right)w_0 + \left( b_0 + b_1 z_0  +  \cdots + b_6 z_0^6\right) w_0^2 \right) dz_0 \wedge dw_0\right\} \notag
\end{align}
We compute $H^0\left( \mathbb{F}_2, \mathcal{E}_{\mathcal{O}_{\mathbb{F}_2}(5,2)} \right)\to H^0\left( \mathbb{F}_2, \Omega_{\mathbb{F}_2}^1 (5,2) \right)\bigoplus H^0\left(\mathbb{F}_2, \mathcal{O}_{\mathbb{F}_2}(5,2) \right)$. First we note that the image of $H^0\left( \mathbb{F}_2, \mathcal{E}_{\mathcal{O}_{\mathbb{F}_2}(5,2)} \right)\to H^0\left( \mathbb{F}_2, \Omega_{\mathbb{F}_2}^1(5,2) \right)$ is given by $(\ref{rf18})$
\begin{align*}
& \left( \left(a_0+ 2a_2+ \left(2a_1+ a_4+ a_3\right)z_0  \right)w_0  - \left(c_1+ \left(2c_2+c_0 \right)z_0 + c_1z_0^2 + c_2z_0^3\right)w_0^2  \right) dz_0+ \left( -a_3 -a_4 + 3a_2z_0  \right)dw_0 \\
&\iff a_0 + 2 a_2=0, 2a_1+a_4+a_3=0, c_0=c_1=c_2=0, a_3+a_4=0, a_2=0
\end{align*}
Next we compute $H^0\left( \mathbb{F}_2, \mathcal{E}_{\mathcal{O}_{\mathbb{F}_2}(5,2)} \right)\to H^0\left( \mathbb{F}_2, \mathcal{O}_{\mathbb{F}_2}(5,2) \right)$
\begin{align*}
&i_{(a_0+ a_1z_0+ a_2z_0^2)\frac{\partial}{\partial z_0}+ ( (a_3-2a_2z_0)w_0+ \left(c_0+c_1z_0+ c_2 z_0^2\right)w_0^2)\frac{\partial}{\partial w_0}}\left(z_0w_0dz_0 - dw_0 \right)\\
&= \left( - a_3 + (a_0+2a_2)z_0 + a_1z_0^2 + a_2z_0^3 \right) w_0 - \left(c_0 + c_1z_0 + c_2z_0^2\right)w_0^2 =0 \iff a_0=a_1=a_2=a_3=c_0=c_1=c_2=0
\end{align*}
Then we have
\begin{align*}
\dim_\mathbb{C} \mathbb{H}^0\left( \mathbb{F}_2, \tilde{\mathcal{E}}_{\mathcal{O}_{\mathbb{F}_2}(5,2)}^\bullet \right)= \dim_\mathbb{C} \textnormal{ker}\left( H^0\left( \mathbb{F}_2, \mathcal{E}_{\mathcal{O}_{\mathbb{F}_2}(5,2)} \right)\to H^0\left( \mathbb{F}_2, \Omega_{\mathbb{F}_2}^1 (5,2) \right)\bigoplus H^0\left(\mathbb{F}_2, \mathcal{O}_{\mathbb{F}_2}(5,2) \right)\right)= 0 
\end{align*}
On the other hand, we compute $H^0\left( \mathbb{F}_2, \Omega_{\mathbb{F}_2}^1 (5,2) \right)\bigoplus H^0\left(\mathbb{F}_2, \mathcal{O}_{\mathbb{F}_2}(5,2) \right) \to H^0\left( \mathbb{F}_2 , \bigwedge^2 \Omega_{\mathbb{F}_2}^1(10,4) \right)$
\begin{align}
&\omega_0 \wedge d\left( a_0 +a_1z_0+ \left( b_0 + b_1z_0 + b_2 z_0^2 + b_3z_0^3 \right)w_0 + \left( c_0 + c_1 z_0 + c_2 z_0^2 + c_3 z_0^3 + c_4z_0^4 + c_5z_0^5 \right)w_0^2\right) \notag \\
& + \left( a_0 +a_1z_0+ \left( b_0 + b_1z_0 + b_2 z_0^2 + b_3z_0^3 \right)w_0 + \left( c_0 + c_1 z_0 + c_2 z_0^2 + c_3 z_0^3 + c_4z_0^4 + c_5 z_0^5 \right)w_0^2\right) d\omega_0 \notag \\
&- \omega_0 \wedge \left( \left( \left( f_0+ f_1 z_0 + 2g_3 z_0^2 \right) w_0 + \left( e_0 + e_1 z_0 + e_2 z_0^2 + e_3 z_0^3 \right) w_0^2\right) dz_0 + \left(g_0+ g_1z_0 + g_2 z_0^2 + g_3 z_0^3 \right) dw_0 \right) \notag\\
&= \left( a_1 - a_0z_0 - a_1z_0^2 +\left( b_1 -f _0 + \left(2b_2- f_1 - g_0 \right) z_0 + ( 3b_3- 2g_3- g_1)z_0^2 - g_2 z_0^3 - g_3 z_0^4 \right) w_0 \right) dz_0\wedge dw_0 \label{rf25}\\
& + \left(\left( c_1 - e_0+ (c_0+ 2c_2 - e_1) z_0  + (c_1+ 3c_3- e_2)z_0^2 + (c_2+ 4c_4 - e_3)z_0^3 + (c_3+ 5c_5)z_0^4 +c_4 z_0^5 + c_5 z_0^6\right) w_0^2 \right)  dz_0 \wedge dw_0 =0 \notag \\
&\iff a_0=a_1=0, g_2=g_3=0, b_1=f_0 , 2b_2= f_1 + g_0, 3b_3=g_1, c_4=c_5= c_3=0, c_1=e_0=e_2, c_0+ 2c_2= e_1, c_2=e_3 \notag
\end{align}
Then we have
\begin{align*}
\dim_\mathbb{C} \textnormal{ker}\left( H^0\left( \mathbb{F}_2, \Omega_{\mathbb{F}_2}^1(5,2) \right) \bigoplus H^0\left( \mathbb{F}_2, \mathcal{O}_{\mathbb{F}_2}(5,2)\right) \to H^0\left( \mathbb{F}_2, \bigwedge^2 \Omega_{\mathbb{F}_2}^1(10,4) \right)\right)= 8
\end{align*}
and so we have
\begin{align}\label{rf24}
\dim_\mathbb{C} \frac{\textnormal{ker}\left( H^0\left( \mathbb{F}_2, \Omega_{\mathbb{F}_2}^1(5,2) \right) \bigoplus H^0\left( \mathbb{F}_2, \mathcal{O}_{\mathbb{F}_2}(5,2)\right) \to H^0\left( \mathbb{F}_2, \bigwedge^2 \Omega_{\mathbb{F}_2}^1(10,4) \right) \right) }{\textnormal{im}\left( H^0\left( \mathbb{F}_2, \mathcal{E}_{\mathcal{O}_{\mathbb{F}_2}(5,2)} \right)\to H^0\left( \mathbb{F}_2, \Omega_{\mathbb{F}_2}(5,2) \right)\bigoplus H^0\left(\mathbb{F}_2, \mathcal{O}_{\mathbb{F}_2}(5,2) \right) \right)} = 8 -8 =0
\end{align}
We recall from $(\ref{rf23})$ that
\begin{align*}
h_{01}\sigma_1- \sigma_0 +\mathcal{L}_{\frac{1}{z_0}\frac{\partial}{\partial w_0}}\left( z_0w_0 dz_0 - dw_0\right) = z_0 \cdot \left(    z_1dz_1 \right) - dz_0 + \left( 1+ \frac{1}{z_0^2} \right)dz_0 =0
\end{align*}
We set
\begin{align*}
\alpha_{01}:= i_{\frac{1}{z_0}\frac{\partial}{\partial w_0}}\left( z_0 w_0 dz_0 - dw_0 \right) = - \frac{1}{z_0}= - z_1 \in C^1\left( \mathfrak{B}, \mathcal{O}_{\mathbb{F}_2}(5,2) \right)
\end{align*}
where $\mathfrak{B}=\left\{ U_0 \bigcup U_2, U_1 \bigcup U_3 \right\}$. If we take $\alpha_0=0$ on $U_0 \bigcup U_2$ and $\alpha_1= z_1^2$ on $U_1 \bigcup U_3$, we have
\begin{align*}
h_{01}\alpha_1 - \alpha_0 + \alpha_{01}= z_0\cdot \left(z_1^2\right)- 0  - \frac{1}{z_0}=  0
\end{align*}
Then we see that for $\left\{(-dz_0,0) , (- z_1dz_1, - z_1^2) \right\}\in C^0\left( \mathfrak{B}, \Omega_{\mathbb{F}_2}^1(5,2)\bigoplus \mathcal{O}_{\mathbb{F}_2}(5,2) \right) $ we have
\begin{align*}
&\delta\left( \left\{ \left( -dz_0  , 0 \right) , \left( - z_1dz_1 ,- z_1^2\right) \right\} \right)=\left( \left( 1+ \frac{1}{z_0^2} \right) dz_0 , - \frac{1}{z_0} \right) \\
&=\left( \mathcal{L}_{\frac{1}{z_0}\frac{\partial}{\partial w_0}}\left( z_0 w_0 dz_0 - d w_0 \right) + 0 \cdot \left( z_0w_0 dz_0 - d w_0 \right), i_{\frac{\partial}{\partial w_0}}  \left( z_0 w_0 dz_0 - d w_0 \right) \right) \in C^1\left( \mathfrak{B}, \Omega_{\mathbb{F}_2}^1(5,2)\bigoplus \mathcal{O}_{\mathbb{F}_2}(5,2) \right)
\end{align*}
Let us consider the image of $(- dz_0,0)$ on $U_0\bigcup U_2$ under $\Omega_{\mathbb{F}_2}^1(5,2)\bigoplus \mathcal{O}_{\mathbb{F}_2}(5,2)\to \bigwedge^2 \Omega_{\mathbb{F}_2}^1 (10, 4)$ which is given by
\begin{align*}
\omega_0  \wedge d\left( 0 \right) + 0 \cdot d \omega_0 - \omega_0 \wedge \left( - dz_0 \right) = \left( z_0w_0dz_0 - dw_0 \right)\wedge dz_0 =  dz_0 \wedge dw_0 \in H^0\left( \mathbb{F}_2, \bigwedge^2 \Omega_{\mathbb{F}_2}^1(10,4) \right)
\end{align*}
which is not in the image of $H^0\left( \mathbb{F}_2, \Omega_{\mathbb{F}_2}^1(5,2) \right)\bigoplus H^0\left( \mathbb{F}_2, \mathcal{O}_{\mathbb{F}_2} (5,2) \right) \to H^0\left( \mathbb{F}_2, \bigwedge^2 \Omega_{\mathbb{F}_2}^1(10,4) \right)$ by $(\ref{rf25})$. This implies that the map in the spectral sequence  defined by
{\tiny{\begin{align*}
\textnormal{ker}\left(H^1\left( \mathbb{F}_2, \mathcal{E}_{\mathcal{O}_{\mathbb{F}_2}(5,2)}\right)\to H^1\left(\mathbb{F}_2, \Omega_{\mathbb{F}_2}^1(5,2) \right) \bigoplus H^1\left( \mathbb{F}_2, \mathcal{O}_{\mathbb{F}_2}(5,2) \right) \right) \cong \left\langle \frac{1}{z_0}\frac{\partial}{\partial w_0} \right\rangle \to \frac{H^0\left( \mathbb{F}_2, \bigwedge^2 \Omega_{\mathbb{F}_2}^1(10,4) \right)}{\textnormal{im}\left(H^0\left( \mathbb{F}_2, \Omega_{\mathbb{F}_2}^1(5,2)  \right)\bigoplus H^0\left( \mathbb{F}_2, \mathcal{E}_{\mathcal{O}_{\mathbb{F}_2}(5,2)} \right) \to H^0\left( \mathbb{F}_2, \bigwedge^2 \Omega_{\mathbb{F}_2}^1(10,4) \right)\right)}
\end{align*}}}
is not $0$, so that the kernel is $0$. Then from $(\ref{rf24})$ we have
\begin{align*}
\mathbb{H}^1\left( \mathbb{F}_2, \tilde{E}_{\mathcal{O}_{\mathbb{F}_2}(5,2)}^\bullet \right)=  0 
\end{align*}
Hence $\left( \mathbb{F}_2 , \omega_0 \right)$ is rigid in unfoldings by \textnormal{Theorem \ref{r1}}.

\end{example}

\begin{remark}
We study unfoldings of $\left( \mathbb{F}_m, \omega_0= adz_0 \right),a \ne 0$ with $adz_0\in H^0\left( \mathbb{F}_m, \Omega_{\mathbb{F}_m}^1(2,0) \right)$ in $\textnormal{Remark \ref{rs40}}$. We keep the notations in $\textnormal{Remark \ref{rs40}}$. We compute $\mathbb{H}^1\left( \mathbb{F}_m, \tilde{\mathcal{E}}_{\mathcal{O}_{\mathbb{F}_m}( 2,0)}^\bullet \right)$. It remains to describe $H^0\left( \mathbb{F}_m, \mathcal{O}_{\mathbb{F}_m}(2,0) \right)$, $ H^1\left( \mathbb{F}_m, \mathcal{O}_{\mathbb{F}_m}(2,0) \right)$ and $H^0\left( \mathbb{F}_m ,  \bigwedge^2 \Omega_{\mathbb{F}_m}^1(4,0)\right)$. From $(\ref{rr5}),(\ref{rsf7})$ and $(\ref{rsf6})$ we have
\begin{align*}
H^0\left( \mathbb{F}_m, \mathcal{O}_{\mathbb{F}_m}(2,0) \right) &=\left\{ g_0 + g_1 z_0 + g_2 z_0^2 | g_0,g_1,g_2\in \mathbb{C} \right\} \\
H^0\left( \mathbb{F}_m, \bigwedge^2 \Omega_{\mathbb{F}_m}^1(4,0) \right) & = 0 \\
H^1\left( \mathbb{F}_m, \mathcal{O}_{\mathbb{F}_m}(2,0) \right)&= 0 
\end{align*}
We compute 
\begin{align*}
&i_{\left( a_0z_0^2+ bz_0 + c \right)\frac{\partial}{\partial z_0} + \left( d + ez_0 + fz_0^2 \right) \frac{\partial}{\partial w_0}}\left( a dz_0 \right)= a\left(a_0 z_0^2 + b z_0 + c \right) \,\,\,\,\,\,\,\textnormal{for}\,\,\, m=0 \\
&i_{\left(a_0+ a_1 z_0 + a_2z_0^2    \right) \frac{\partial}{\partial z_0} + \left( \left(-ma_2z_0+d\right)w_0 + \left(c_0 + c_1z_0+ \cdots + c_m z_0^m \right)w_0^2 \right)\frac{\partial}{\partial w_0}}\left( a dz_0 \right)=  a\left( a_0 + a_1 z_0 + a_2 z_0^2 \right) \,\,\,\,\,\,\,\textnormal{for}\,\,\,  m \geq 1
\end{align*}
Then $\textnormal{coker}\left(H^0\left( \mathbb{F}_m, \mathcal{E}_{\mathcal{O}_{\mathbb{F}_m}(2,0)} \right)\to H^0\left( \mathbb{F}_m, \Omega_{\mathbb{F}_m}^1(2,0) \right) \bigoplus  H^0\left( \mathbb{F}_m, \mathcal{O}_{\mathbb{F}_m}(2,0) \right) \right)=0$. Hence we have
\begin{align*}
\mathbb{H}^1\left( \mathbb{F}_m, \tilde{\mathcal{E}}_{\mathcal{O}_{\mathbb{F}_m}(2,0)}^\bullet \right)= H^1\left(\mathbb{F}_m, \Theta_{\mathbb{F}_m} \right) =\textnormal{Span}_\mathbb{C}\left\langle  z_1 \frac{\partial}{\partial w_0},..., z_1^{m-1}\frac{\partial}{\partial w_0} \right\rangle
\end{align*}
We will construct an unfolded analytic family of unfoldings of $\left( \mathbb{F}_m, \omega_0 = a dz_0 \right)$ such that the unfolded Kodaira-Spencer map is an isomorphism at the distinguished point, so that $\left( \mathbb{F}_m, \omega_0 \right)$ is unobstructed in unfoldings. The unfolded analytic family is given by
\begin{align*}
\omega&= a dz_0\\
z_0&= \frac{1}{z_1} \\
w_0&= z_1^m w_1+ t_1 z_1+ \cdots + t_{m-1} z_1^{m-1}
\end{align*}
\end{remark}

\section{Foliated deformations of Projective $3$-space $\mathbb{P}_\mathbb{C}^3$}

\subsection{Dimension $1$-foliated deformations of $\mathbb{P}_\mathbb{C}^3$}\

Let $\left( \mathbb{P}_\mathbb{C}^3, T \right)$ be a foliated projective $3$-space, where $T\in H^0\left(  \mathbb{P}_\mathbb{C}^3, \Theta_{\mathbb{P}_\mathbb{C}^3}(m) \right), m\geq -1$, so that we have the exact sequence defined by $T$
\begin{align*}
0 \to \mathcal{O}_{\mathbb{P}_\mathbb{C}^3}(-m)\xrightarrow{T} \Theta_{\mathbb{P}_\mathbb{C}^3} \to \frac{\Theta_{\mathbb{P}_\mathbb{C}^3}}{\mathcal{O}_{\mathbb{P}_\mathbb{C}^3}(-m) } \to 0
\end{align*}
Then we have a complex of sheaves from $(\ref{ap3})$
\begin{align*}
 \mathcal{E}_{\mathcal{O}_{\mathbb{P}_\mathbb{C}^3}(m)}^\bullet :  \mathcal{E}_{\mathcal{O}_{\mathbb{P}_\mathbb{C}^3}(m)} \xrightarrow{\tilde{D}_0} \Theta_{\mathbb{P}_\mathbb{C}^3}(m) \to 0
\end{align*}
From the exact sequence $0\to \mathcal{O}_{\mathbb{P}_\mathbb{C}^3}\to \mathcal{E}_{ \mathcal{O}_{\mathbb{P}_\mathbb{C}^3}(m)}\to \Theta_{\mathbb{P}_\mathbb{C}^3}\to 0$, we have $H^i\left( \mathbb{P}_\mathbb{C}^3 , \mathcal{E}_{\mathcal{O}_{\mathbb{P}_\mathbb{C}^3}(m)} \right)=0$ for $i\geq 1$. On the other hand, from the Euler sequence $0\to \mathcal{O}_{\mathbb{P}_\mathbb{C}^3} \to \bigoplus^4 \mathcal{O}_{\mathbb{P}_\mathbb{C}^3}(1)\to \Theta_{\mathbb{P}_\mathbb{C}^3} \to 0$,  we have $0\to \mathcal{O}_{\mathbb{P}_\mathbb{C}^3}(m)\to \bigoplus^4 \mathcal{O}_{\mathbb{P}_\mathbb{C}^3}(m+1) \to \Theta_{\mathbb{P}_\mathbb{C}^3}(m) \to 0$, so that $H^i\left( \mathbb{P}_\mathbb{C}^3, \Theta_{\mathbb{P}_\mathbb{C}^3 } (m) \right)=0$ for $m\geq -1, i\geq 1$. This implies that we can compute $\mathbb{H}^i\left(  \mathbb{P}_\mathbb{C}^3,  \mathcal{E}_{ \mathcal{O}_{\mathbb{P}_\mathbb{C}^3}(m)}^\bullet \right)$ by the following complex
\begin{align*}
H^0\left( \mathbb{P}_\mathbb{C}^3, \mathcal{E}_{\mathcal{O}_{\mathbb{P}_\mathbb{C}^3}(m)}  \right) \to H^0\left( \mathbb{P}_\mathbb{C}^3, \Theta_{\mathbb{P}_\mathbb{C}^3}(m) \right) \to 0
\end{align*}
Then we have
\begin{align*}
\mathbb{H}^1\left( \mathbb{P}_\mathbb{C}^3, \mathcal{E}_{\mathcal{O}_{\mathbb{P}_\mathbb{C}^3}(m)}^\bullet \right) & \cong  \textnormal{coker}\left( H^0\left( \mathbb{P}_\mathbb{C}^3, \mathcal{E}_{\mathcal{O}_{\mathbb{P}_\mathbb{C}^3}(m)} \right) \xrightarrow{\tilde{D}_0} H^0\left( \mathbb{P}_\mathbb{C}^3, \Theta_{\mathbb{P}_\mathbb{C}^3}(m) \right) \right) \\
\mathbb{H}^2\left( \mathbb{P}_\mathbb{C}^3 , \mathcal{E}_{\mathcal{O}_{\mathbb{P}_\mathbb{C}^3}(m)}^\bullet  \right)&= 0
\end{align*}
In particular, $\left( \mathbb{P}_\mathbb{C}^3, T \right)$ is unobstructed in foliated deformations in terms of tangent sheaves.

\begin{example}\label{pr2}
We keep the notations in \textnormal{Appendix \ref{cp7}}. Let us consider the foliation on $\mathbb{P}_\mathbb{C}^3$ defined by $T_0=\xi_2\frac{\partial}{\partial \xi_0} + \xi_0\frac{\partial}{\partial \xi_1} + \xi_3\frac{\partial}{\partial \xi_2} \in \mathbb{H}^0\left( \mathbb{P}_\mathbb{C}^3, \Theta_{\mathbb{P}_\mathbb{C}^3} \right)$, so that $\Theta_{\mathcal{F}_0}\cong \mathcal{O}_{\mathbb{P}_\mathbb{C}^3}$. We will compute $\mathbb{H}^1\left(\mathbb{P}_\mathbb{C}^3, \Theta_{\mathcal{F}_0}^\bullet \right)$ and describe a foliated analytic family of deformations of $\left(\mathbb{P}_\mathbb{C}^3, \Theta_{\mathcal{F}_0}  \right)$ such that the foliated Kodaira-Spencer map in terms of tangent sheaves is bijective at the distinguished point.  On $U_0$, we have
\begin{align*}
T_0= \left( 1- X_1X_2 \right)\frac{\partial}{\partial X_1} +\left(X_3- X_2^2 \right)\frac{\partial}{\partial X_2} - X_2X_3 \frac{\partial}{\partial X_3}
\end{align*}
Let us compute $\mathbb{H}^i\left( \mathbb{P}_\mathbb{C}^3, \mathcal{E}_{\mathcal{O}_{\mathbb{P}_\mathbb{C}^3}}^\bullet \right)$ by the complex
\begin{align*}
\mathcal{E}_{\mathcal{O}_{\mathbb{P}_\mathbb{C}^3}}^\bullet: H^0\left( \mathbb{P}_\mathbb{C}^3, \mathcal{E}_{\mathcal{O}_{\mathbb{P}_\mathbb{C}^3}}  \right) = H^0\left( \mathbb{P}_\mathbb{C}^3, \Theta_{\mathbb{P}_\mathbb{C}^3} \right) \bigoplus H^0\left( \mathbb{P}_\mathbb{C}^3, \mathcal{O}_{\mathbb{P}_\mathbb{C}^3} \right)\xrightarrow{\tilde{D}_0} H^0\left( \mathbb{P}_\mathbb{C}^3, \Theta_{\mathbb{P}_\mathbb{C}^3} \right)\to 0
\end{align*}
We compute $\textnormal{ker}\left( H^0\left( \mathbb{P}_\mathbb{C}^3, \mathcal{E}_{\mathcal{O}_{\mathbb{P}_\mathbb{C}^3}} \right)= H^0\left( \mathbb{P}_\mathbb{C}^3, \Theta_{\mathbb{P}_\mathbb{C}^3} \right)\bigoplus H^0\left( \mathbb{P}_\mathbb{C}^3, \mathcal{O}_{\mathbb{P}_\mathbb{C}^3} \right) \xrightarrow{\tilde{D}_0} H^0\left( \mathbb{P}_\mathbb{C}^3, \Theta_{\mathbb{P}_\mathbb{C}^3} \right) \right)$ on $U_0$. We recall $\textnormal{(\ref{cp8})}$ for $m=0$. Then the image of $\tilde{D}_0$ is given by  $\left[T_1 \frac{\partial}{\partial X_1} + T_2 \frac{\partial}{\partial X_2} + T_3 \frac{\partial}{\partial X_3}, T_0 \right]+ f T_0$ which is
{\small{\begin{align}
&\left[  \left(a_1+a_2X_3+a_3X_2+a_4X_1+a_5X_1X_3+a_6X_1X_2+a_7X_1^2\right)\frac{\partial}{\partial X_1} ,  \left( 1- X_1X_2 \right)\frac{\partial}{\partial X_1} +\left(X_3- X_2^2 \right)\frac{\partial}{\partial X_2} - X_2X_3 \frac{\partial}{\partial X_3} \right]\notag\\
&+\left[ \left(a_8+a_9X_3+a_{10}X_2 +a_{11}X_1+a_5X_2X_3+a_6X_2^2+a_7X_1X_2 \right)\frac{\partial}{\partial X_2} ,  \left( 1- X_1X_2 \right)\frac{\partial}{\partial X_1} +\left(X_3- X_2^2 \right)\frac{\partial}{\partial X_2} - X_2X_3 \frac{\partial}{\partial X_3}\right]\notag\\
&+\left[ \left(a_{12}+a_{13}X_3+a_{14}X_2+a_{15}X_1+a_5X_3^2+a_6X_2X_3+a_7X_1X_3 \right)\frac{\partial}{\partial X_3} ,  \left( 1- X_1X_2 \right)\frac{\partial}{\partial X_1} +\left(X_3- X_2^2 \right)\frac{\partial}{\partial X_2} - X_2X_3 \frac{\partial}{\partial X_3}\right]\notag\\
&+ a_{16} \left(\left( 1- X_1X_2 \right)\frac{\partial}{\partial X_1} +\left(X_3-  X_2^2 \right)\frac{\partial}{\partial X_2} - X_2X_3 \frac{\partial}{\partial X_3} \right)\notag\\
&=\left(\left(a_{16}- a_4\right) - \left(2a_7 + a_8\right) X_1  - \left(a_1+ a_6\right) X_2 - \left(a_3+ a_5\right)X_3 - a_{11}X_1^2 - \left( a_{10} + a_{16}\right)X_1X_2 - \left(a_6 + a_9\right)X_1X_3 \right)\frac{\partial}{\partial X_1} \label{pr1}\\
&+ \left(  \left( a_{12}- a_{11} \right) + a_{15}X_1 + \left(a_{14}- a_7 - 2a_8\right)X_2 + \left(a_{16}+ a_{13}- a_{10}\right)X_3 - \left(a_{10} + a_{16}\right)X_2^2 - a_{11}X_1X_2  - \left(a_9 + a_6\right)X_2X_3 \right) \frac{\partial}{\partial X_2}\notag\\
&+ \left( - a_{15} - a_{12}X_2 - \left(a_7 + a_8 + a_{14}\right)X_3 - \left(a_6 + a_9\right)X_3^2 - \left(a_{10} + a_{16}\right)X_2X_3 - a_{11}X_1X_3 \right) \frac{\partial}{\partial X_3}\notag
\end{align}}}
Then we see that
\begin{align}
\dim_\mathbb{C} \mathbb{H}^0\left( \mathbb{P}_\mathbb{C}^3, \mathcal{E}_{\mathcal{O}_{\mathbb{P}_\mathbb{C}^3}}^\bullet \right)&= \textnormal{ker}\left( H^0\left( \mathbb{P}_\mathbb{C}^3, \Theta_{\mathbb{P}_\mathbb{C}^3} \right) \bigoplus H^0\left( \mathbb{P}_\mathbb{C}^3, \mathcal{O}_{\mathbb{P}_\mathbb{C}^3} \right)\xrightarrow{\tilde{D}_0} H^0\left( \mathbb{P}_\mathbb{C}^3, \Theta_{\mathbb{P}_\mathbb{C}^3} \right) \right)= 4 \label{blwp2}\\
\dim_\mathbb{C} \mathbb{H}^1\left(\mathbb{P}_\mathbb{C}^3, \mathcal{E}_{\mathcal{O}_{\mathbb{P}_\mathbb{C}^3}}^\bullet \right)&=\textnormal{coker}\left( H^0\left( \mathbb{P}_\mathbb{C}^3, \Theta_{\mathbb{P}_\mathbb{C}^3} \right) \bigoplus H^0\left( \mathbb{P}_\mathbb{C}^3, \mathcal{O}_{\mathbb{P}_\mathbb{C}^3} \right)\xrightarrow{\tilde{D}_0} H^0\left( \mathbb{P}_\mathbb{C}^3, \Theta_{\mathbb{P}_\mathbb{C}^3} \right) \right)=15- (16-4)=3
\end{align}
Then by comparing $\textnormal{(\ref{ap2})}$ and $(\ref{pr1})$, we have
\begin{align*}
 \mathbb{H}^1\left(\mathbb{P}_\mathbb{C}^3, \mathcal{E}_{\mathcal{O}_{\mathbb{P}_\mathbb{C}^3}}^\bullet \right)\cong \textnormal{Span}_\mathbb{C}\left\langle X_1\frac{\partial}{\partial X_3}, \frac{\partial}{\partial X_3}, X_2 \frac{\partial}{\partial X_3}  \right\rangle
\end{align*}
Then we see that $\left( \mathbb{P}_\mathbb{C}^3\times B , T, \pi , B \right)$, where $\pi:\mathbb{P}_\mathbb{C}^3 \times B \to B$ is the natural projection with the coordinates $t=\left(t_1,t_2,t_3\right)\in B\subset \mathbb{C}^3$ for a sufficiently small neighborhood $B$ of $0$ in $\mathbb{C}^3$ and
\begin{align*}
T=\left( 1- X_1X_2 \right)\frac{\partial}{\partial X_1} + \left( X_3 - X_2^2 \right)\frac{\partial}{\partial X_2} + \left( - X_2X_3 + t_1 + t_2 X_1 + t_3 X_2 \right)\frac{\partial}{\partial X_3}
\end{align*}
defines a foliated analytic family of deformations $\pi^{-1}(0)=\left( \mathbb{P}_\mathbb{C}^3, T_0 \right)$ in terms of tangent sheaves and the foliated Kodaira-Spencer map in terms of tangent sheaves is an isomorphism at $t=0$.

\end{example}

\subsection{Codimension $1$-foliated deformations of $\mathbb{P}_\mathbb{C}^3$}\

Let $\left(\mathbb{P}_\mathbb{C}^3,w\right)$ be a foliated projective $3$-space, where $w\in H^0\left(\mathbb{P}_\mathbb{C}^3, \Omega_{\mathbb{P}_\mathbb{C}^3}^1(m)\right),m\geq 2$ with $w\wedge dw=0$, so that we have the exact sequence defined by $w$
\begin{align*}
0\to \mathcal{O}_{\mathbb{P}_\mathbb{C}^3}(-m)\xrightarrow{w} \Omega_{\mathbb{P}_\mathbb{C}^3}^1 \to \frac{\Omega_{\mathbb{P}_\mathbb{C}^3}^1}{\mathcal{O}_{\mathbb{P}_\mathbb{C}^3}(-m)} \to 0
\end{align*}
Then we have a complex of sheaves from $(\ref{cp2})$
\begin{align}\label{cp1}
\mathcal{E}_{\mathcal{O}_{\mathbb{P}_\mathbb{C}^3}(m)}^\bullet: \mathcal{E}_{\mathcal{O}_{\mathbb{P}_\mathbb{C}^3}(m)}\xrightarrow{\tilde{E}_0} \Omega_{\mathbb{P}_\mathbb{C}^3}^1(m) \xrightarrow{\tilde{E}_1} \tilde{\mathcal{S}}^2\otimes \mathcal{O}_{\mathbb{P}_\mathbb{C}^3}(m) \to 0
\end{align}
We note that $\tilde{\mathcal{S}}^2\subset \mathcal{O}_{\mathbb{P}_\mathbb{C}^3}(m)\otimes \bigwedge^{1+2}\Omega_{\mathbb{P}_\mathbb{C}^3}^1\cong \mathcal{O}_{\mathbb{P}_\mathbb{C}^3}(m-4)$. Since $\tilde{\mathcal{S}}^2$ is reflexive of rank $1$, we see that $\tilde{\mathcal{S}}^2 \cong \mathcal{O}_{\mathbb{P}_\mathbb{C}^3}(m-4)$. Then $(\ref{cp1})$ is equivalent to
\begin{align*}
\mathcal{E}_{\mathcal{O}_{\mathbb{P}_\mathbb{C}^3}(m)}^\bullet: \mathcal{E}_{\mathcal{O}_{\mathbb{P}_\mathbb{C}^3}(m)} \xrightarrow{\tilde{E}_0} \Omega_{\mathbb{P}_\mathbb{C}^3}^1(m) \xrightarrow{\tilde{E}_1}   \mathcal{O}_{\mathbb{P}_\mathbb{C}^3}(2m-4) \to 0
\end{align*}

Then we have by Bott formula (see \cite{OSS88} p.4)
\begin{align*}
\dim_\mathbb{C} H^0\left(\mathbb{P}_\mathbb{C}^3, \Omega_{\mathbb{P}_\mathbb{C}^3}^1(m)\right)&=\frac{(m+2)(m+1)(m-1)}{2},\,\,\,\,\,\,\,\,\,m\geq 2\\
\dim_\mathbb{C} H^0\left(\mathbb{P}_\mathbb{C}^3,\mathcal{O} (2m-4)\right)&=\frac{(2m-1)(2m-2)(2m-3)}{6},\,\,\,\,\,\,\,\,\, m \geq 2
\end{align*} 
On the other hand, from the exact sequence $0\to \mathcal{O}_{\mathbb{P}_\mathbb{C}^3}\to \mathcal{E}_{\mathcal{O}_{\mathbb{P}_{\mathbb{C}}^3 }(m)} \to \Theta_{\mathbb{P}_\mathbb{C}^3} \to 0$ and Bott formula, we have
\begin{align*}
H^i\left(\mathbb{P}_\mathbb{C}^3, \mathcal{E}_{\mathcal{O}_{\mathbb{P}_\mathbb{C}^3}(m)} \right)&=0 \,\,\,\,\,\,\,\,\,\,i>0 \\
 H^i\left(\mathbb{P}_\mathbb{C}^3,\Omega_{\mathbb{P}_\mathbb{C}^3}^1(m)\right)&=0,\,\,\,\,\,\,\, i>0,\,\, m \geq 2\\
 H^i\left(\mathbb{P}_\mathbb{C}^3, \mathcal{O}_{\mathbb{P}_\mathbb{C}^3}(2m-4)\right)&=0,\,\,\,\,\,\,\,\,i>0,\,\,m\geq 2
\end{align*}
We see that we can compute $\mathbb{H}^i \left(\mathbb{P}_\mathbb{C}^3, \mathcal{E}_{ \mathcal{O}_{\mathbb{P}_\mathbb{C}^3}(m)}^\bullet \right)$ by the following complex
\begin{align*}
H^0\left(\mathbb{P}_\mathbb{C}^3,\mathcal{E}_{\mathcal{O}_{\mathbb{P}_\mathbb{C}^3}(m)}\right)\xrightarrow{\tilde{E}_0} H^0\left(\mathbb{P}_\mathbb{C}^3,\Omega_{\mathbb{P}_\mathbb{C}^3}^1(m)  \right) \xrightarrow{\tilde{E}_1} H^0\left(\mathbb{P}_\mathbb{C}^3, \mathcal{O}_{\mathbb{P}_\mathbb{C}^3}(2m-4)\right)\to 0
\end{align*}
Then we have
\begin{align*}
\mathbb{H}^1\left(\mathbb{P}_\mathbb{C}^3, \mathcal{E}_{\mathcal{O}_{\mathbb{P}_\mathbb{C}^3}(m)}^\bullet \right)&\cong \frac{\textnormal{ker}\left(H^0\left(\mathbb{P}_\mathbb{C}^3,\Omega_{\mathbb{P}_\mathbb{C}^3}^1(m)  \right)\xrightarrow{\tilde{E}_1} H^0\left(\mathbb{P}_\mathbb{C}^3, \mathcal{O}_{\mathbb{P}_\mathbb{C}^3}(2m-4)\right) \right)}{\textnormal{im}\left(H^0\left(\mathbb{P}_\mathbb{C}^3,\mathcal{E}_{\mathcal{O}_{\mathbb{P}_\mathbb{C}^3}(m)}\right) \xrightarrow{\tilde{E}_0} H^0\left(\mathbb{P}_\mathbb{C}^3,\Omega_{\mathbb{P}_\mathbb{C}^3}^1(m)  \right) \right)}\\
\mathbb{H}^2 \left(\mathbb{P}_\mathbb{C}^3,\mathcal{E}_{\mathcal{O}_{\mathbb{P}_\mathbb{C}^3}(m)} \right)&\cong \textnormal{coker}\left(H^0\left(\mathbb{P}_\mathbb{C}^3,\Omega_{\mathbb{P}_\mathbb{C}^3}^1(m)  \right)\xrightarrow{\tilde{E}_1} H^0\left(\mathbb{P}_\mathbb{C}^3, \mathcal{O}_{\mathbb{P}_\mathbb{C}^3}(2m-4)\right)    \right)
\end{align*}

If $m=2$, then $\mathbb{H}^2 \left(\mathbb{P}_\mathbb{C}^3, \mathcal{E}_{\mathcal{O}_{\mathbb{P}_\mathbb{C}^3}(m)} \right) \cong \textnormal{coker} \left(H^0(\mathbb{P}_\mathbb{C}^3,\Omega_{\mathbb{P}_\mathbb{C}^3}^1(2))\cong \mathbb{C}^6\to \mathbb{C} \right)$. Hence if $\tilde{E}_1$ is not the zero map, then $\left( \mathbb{P}_\mathbb{C}^3, w \right)$ is unobstructed in foliated deformations in terms of cotangent sheaves.

\subsection{Foliated deformations of $\mathbb{P}_{\mathbb{C}}^3$ with linear-pullback of Jouanolou foliation of degree $2$ on $\mathbb{P}_\mathbb{C}^2$ in terms of cotangent sheaves  } \label{j6}  \

Let us consider the Jouanolou foliation $w_0$ of degree $2$ on $\mathbb{P}_\mathbb{C}^2$ in Example \ref{cp10}. We will study foliated deformations of the linear-pullback of $w_0$ to $\mathbb{P}_\mathbb{C}^3$. We still denote the linear-pullback of $w_0$ by the same symobl $w_0$. We will keep the notations in Appendix \ref{cp7} and Appendix \ref{ja1}. Then we have on $U_0=\left\{\xi_0 \ne 0 \right\}$
\begin{align*}
w_0&= \left( X_1^2 X_2 -1 \right) dX_1 + \left( X_2^2 - X_1^3 \right) dX_2\\
dw_0&= - 4 X_1^2 dX_1 \wedge dX_2
\end{align*}

Let us compute $\mathbb{H}^i\left(  \mathbb{P}_\mathbb{C}^3, \mathcal{E}_{\mathcal{O}_{\mathbb{P}_\mathbb{C}^3}(4)}^\bullet \right)$ by the complex
\begin{align*}
\mathcal{E}_{\mathcal{O}_{\mathbb{P}_\mathbb{C}^3}(4)}^\bullet  : H^0\left( \mathbb{P}_\mathbb{C}^3, \mathcal{E}_{\mathcal{O}_{\mathbb{P}_\mathbb{C}^3}(4)} \right) \xrightarrow{\tilde{E}_0} H^0\left( \mathbb{P}_\mathbb{C}^3, \Omega_{\mathbb{P}_\mathbb{C}^3}^1(4) \right) \xrightarrow{\tilde{E}_1} H^0\left( \mathbb{P}_\mathbb{C}^3, \mathcal{O}_{\mathbb{P}_\mathbb{C}^3} (4) \right) \to 0
\end{align*}
First we compute $\textnormal{ker}\left( H^0\left( \mathbb{P}_\mathbb{C}^3 , \Omega_{\mathbb{P}_\mathbb{C}^3}^1(4) \right) \xrightarrow{\tilde{E}_1} H^0\left( \mathbb{P}_\mathbb{C}^3, \mathcal{O}_{\mathbb{P}_\mathbb{C}^3}(4) \right) \right)$. We recall $\textnormal{(\ref{cm1})}$ and $\textnormal{(\ref{cm2})}$. Then the image of $\tilde{E}_1$ is given by $w_0 \wedge db + b \wedge dw_0$. The coefficient of $dX_1 \wedge dX_2 \wedge dX_3$ in $dw_0\wedge b$ is given by
{\Small{\begin{align*}
&\left(t_{37}-t_{23}\right)X_2^2+ \left(t_{43}-2t_{29}\right)X_1X_2^2 + \left(t_{41}-t_{28}\right)X_2^3+\left(2t_{40}-t_{27}\right)X_2^2X_3-4t_{35}X_1^2X_2^2-\left(t_{15}+2t_{32}\right)X_2^4-4t_{30}X_2^2X_3^2\\
& -\left(t_{17}+3t_{34}\right)X_1X_2^3 -4t_{33}X_1X_2^2X_3- \left(2t_{14}+3t_{31}\right)X_2^3X_3- \left( t_{37}-t_{23}\right)X_1^3- \left(t_{43}-2t_{29}\right)X_1^4 - \left(t_{41}-t_{28}\right)X_1^3X_2- \left(2t_{40}-t_{27}\right)X_1^3X_3\\
&+4t_{35}X_1^5+\left(t_{15}+2t_{32}\right)X_1^3X_2^2+4t_{30}X_1^3X_3^2 + \left(t_{17}+3t_{34}\right)X_1^4X_2 +4t_{33}X_1^4X_3+ \left(2t_{14}+3t_{31}\right)X_1^3 X_2X_3
\end{align*}}}
and the coefficient of $dX_1 \wedge dX_2 \wedge dX_3$ in $dw_0 \wedge b$ is given by
{\Small{\begin{align*}
& -4t_{20}X_1^2-4t_{21}X_1^2X_3-4t_{22}X_1^2X_2-4t_{23}X_1^3-4t_{24}X_1^2X_3^2-4t_{25}X_1^2X_2X_3-4t_{26}X_1^2X_2^2-4t_{27}X_1^3X_3-4t_{28}X_1^3X_2-4t_{29}X_1^4\\
& +4t_{11}X_1^2X_2X_3^2+4t_{12}X_1^2X_2^2X_3+4t_{13}X_1^2X_2^3 -4t_{30}X_1^3X_3^2-4t_{31}X_1^3X_2X_3-4t_{32}X_1^3X_2^2-4t_{33}X_1^4X_3-4t_{34}X_1^4X_2-4t_{35}X_1^5
\end{align*}}}

Then the coefficient of $dX_1\wedge dX_2 \wedge dX_3$ in $w_0 \wedge db + b \wedge dw_0$ is given by
{\Small{\begin{align*}
&- t_{22}+ t_2 -(t_{28}-t_8)X_1-(2t_{26}-t_6)X_2-(t_{25}-2t_5)X_3 - (t_{34}- t_{17}+ 4t_{20})X_1^2 +(4t_{13}+ t_{37}- t_{23})X_2^2 +4t_{11}X_3^2 \\
&-  (2t_{32}-t_{15})X_1X_2 - (t_{31}-2t_{14})X_1X_3+4t_{12} X_2X_3 -( t_{37}+3 t_{23})X_1^3 +(t_{41}-t_{28})X_2^3 + (- t_2- 3t_{22} )X_1^2X_2\\
& -4t_{21}X_1^2X_3+(t_{43}-2t_{29})X_1X_2^2    +(2t_{40}-t_{27})X_2^2X_3                     -(t_{43} + 2 t_{29})X_1^4 -(t_{15}+2t_{32})X_2^4 \\
&+(-t_8- t_{41} - 2 t_{28})X_1^3 X_2-(2t_{40}+ 3t_{27})X_1^3X_3 +(-2t_{26}-t_6 - 4t_{35} )X_1^2X_2^2 +(-2t_5- 3t_{25})X_1^2X_2X_3\\
&-4t_{24}X_1^2X_3^2  -(t_{17}+3t_{34} )X_1X_2^3 -4t_{33}X_1X_2^2X_3  -(2t_{14}+3t_{31})X_2^3X_3 -4t_{30}X_2^2X_3^2 
\end{align*}}}
Then the kernel of $\tilde{E}_1$ is given by
{\Small{\begin{align}\label{j1}
b=&\left(A_{14}-3A_3X_3+ A_{15}X_2+ A_{16}X_1+ A_{17}X_2^2+ A_{18}X_1X_2+ A_{10}X_1^2 \right) dX_1\\
+&\left(- A_{13}X_1^2X_2+A_1X_1^2X_3- A_{12}X_1X_2^2+2A_2X_1X_2X_3- A_{11}X_2^3 \right)dX_1 \notag  \\ 
+&\left(A_4+ A_5X_2+ A_6X_1+ 2A_1X_2X_3+ A_7X_2^2 + A_8X_1X_2+ A_9X_1^2  \right)dX_2 \notag \\
+&\left( A_3X_2^2X_3 + A_{11}X_1X_2^2- 3A_2X_1^2X_3+ A_{12}X_1^2X_2+ A_{13}X_1^3   \right)dX_2 \notag \\
+&\left(-A_2+ A_3X_1+ A_1X_2^2 -A_3X_2^3 + A_2X_1^2X_2-A_1 X_1^3 \right)dX_3,\,\,\,\,\,\,A_1,...,A_{18}\in \mathbb{C} \notag
\end{align}}}
so that we have
\begin{align}
&\dim_\mathbb{C} \textnormal{ker}\left( H^0\left( \mathbb{P}_\mathbb{C}^3, \Omega_{\mathbb{P}_\mathbb{C}^3}^1 (4) \right)\xrightarrow{\tilde{E}_1} H^0\left( \mathbb{P}_\mathbb{C}^3, \mathcal{O}_{\mathbb{P}_\mathbb{C}^3}(4)\right)\right) = 18 \notag \\
&\dim_\mathbb{C} \mathbb{H}^2\left(\mathbb{P}_\mathbb{C}^3, \mathcal{N}_{\mathcal{F}_0}^{*\bullet} \right)=\dim_\mathbb{C} \mathbb{H}^2\left( \mathbb{P}_\mathbb{C}^3, \mathcal{E}_{\mathcal{O}_{\mathbb{P}_\mathbb{C}^3}(4)}^\bullet \right)=\dim_\mathbb{C}\textnormal{coker}\left(  H^0\left( \mathbb{P}_\mathbb{C}^3, \Omega_{\mathbb{P}_\mathbb{C}^3}^1 (4) \right)\xrightarrow{\tilde{E}_1} H^0\left( \mathbb{P}_\mathbb{C}^3, \mathcal{O}_{\mathbb{P}_\mathbb{C}^3}(4)\right)\right)= 35- (45- 18)= 8  \label{j25}
\end{align}
Next we compute $\textnormal{im}\left( H^0\left( \mathbb{P}_\mathbb{C}^3 , \mathcal{E}_{\mathcal{O}_{\mathbb{P}_\mathbb{C}^3}(4)}  \right) \xrightarrow{\tilde{E}_0} H^0\left(  \mathbb{P}_\mathbb{C}^3, \Omega_{\mathbb{P}_\mathbb{C}^3}^1(4) \right) \right)$. We recall $\textnormal{(\ref{cp8})}$ for $m=4$. Then the image of $\tilde{E}_0$ is given by $\mathcal{L}_{T_1\frac{\partial}{\partial X_1} + T_2 \frac{\partial}{\partial X_2} + T_3 \frac{\partial}{\partial X_3}}(w_0) + fw_0$.  We note that $\mathcal{L}_{T_1\frac{\partial}{\partial X_1}}\left( w_0 \right)$ is given by
{\Small{\begin{align*}
& \left(-a_4-a_5X_3- a_6X_2- 2a_7X_1+2a_1X_1X_2+2a_2X_1X_2X_3+2a_3X_1X_2^2+3a_4X_1^2 X_2+3a_5X_1^2 X_2X_3+3a_6X_1^2X_2^2+4a_7X_1^3X_2\right)d X_1\\
&+ \left(-a_3 - a_6X_1-3a_1X_1^2-3a_2X_1^2X_3-2a_3X_1^2X_2-3a_4X_1^3-3a_5X_1^3X_3-2a_6X_1^3X_2-3a_7X_1^4\right)dX_2\\
&+\left(-a_2-a_5X_1 +a_2X_1^2X_2 + a_5X_1^3 X_2\right)d X_3
\end{align*}}}
and $\mathcal{L}_{T_2\frac{\partial}{\partial X_2}}\left( w \right)$ is given by
{\Small{\begin{align*}
&\left(a_8X_1^2+a_{11}X_2^2+a_9X_1^2X_3+a_{10}X_1^2X_2 +a_7X_2^3+a_5X_1^2X_2X_3+a_6X_1^2X_2^2 \right)dX_1\\
&+\left(2a_8X_2+2a_9X_2X_3+3a_{10}X_2^2 +2a_{11}X_1X_2+3a_5X_2^2X_3+ 4a_6X_2^3+3a_7X_1X_2^2- a_{10}X_1^3- a_5 X_1^3X_3- 2a_6X_1^3X_2 - a_7X_1^4\right)dX_2\\
&+\left(a_9X_2^2+a_5X_2^3 - a_9X_1^3 - a_5X_1^3X_2\right)dX_3
\end{align*}}}
and $\mathcal{L}_{T_3\frac{\partial}{\partial X_3}} \left( w \right)=0$, and $fw_0$ is given by
{\Small{\begin{align*}
&\left( -a_{16}+4a_7X_1+4a_6X_2+4a_5X_3 +a_{16}X_1^2 X_2-4a_7X_1^3X_2-4a_6X_1^2X_2^2-4a_5X_1^2 X_2X_3\right) dX_1 \\
&+\left(a_{16}X_2^2-4a_7X_1X_2^2-4a_6X_2^3-4a_5X_2^2X_3 - a_{16}X_1^3+4a_7X_1^4+4a_6X_1^3X_2+4a_5X_1^3 X_3 \right)dX_2
\end{align*}}}
Then $\mathcal{L}_{T_1 \frac{\partial}{\partial X_1} + T_2 \frac{\partial}{\partial X_2} + T_3 \frac{\partial}{\partial X_3}}\left( w_0\right) + f w_0$ is given by
{\tiny{\begin{align}
&\left(-a_4- a_{16}+3a_5X_3+3a_6X_2+2a_7X_1+2a_1X_1X_2+2a_2X_1X_2X_3+2a_3X_1X_2^2+\left(3a_4+ a_{10}+ a_{16}\right)X_1^2 X_2+a_8X_1^2+a_{11}X_2^2+a_9X_1^2X_3 +a_7X_2^3\right)dX_1 \label{j2}\\
&+ \left(-a_3 - a_6X_1-3a_1X_1^2-3a_2X_1^2X_3-2a_3X_1^2X_2- \left( 3a_4+a_{10}+ a_{16} \right)X_1^3+2a_8X_2+2a_9X_2X_3+ \left(3a_{10}+a_{16}\right)X_2^2+2a_{11}X_1X_2-a_5X_2^2X_3-a_7X_1X_2^2 \right)dX_2 \notag \\
&+\left(-a_2-a_5X_1 +a_2X_1^2X_2 +a_9X_2^2+a_5X_2^3 - a_9X_1^3  \right)dX_3 \notag
\end{align}}}
Then we see that $\dim_\mathbb{C} \mathbb{H}^0\left(  \mathbb{P}_\mathbb{C}^3, \mathcal{E}_{\mathcal{O}_{\mathbb{P}_\mathbb{C}^3}(4)}^\bullet  \right)=\dim_\mathbb{C} \textnormal{ker}\left( H^0\left( \mathbb{P}_\mathbb{C}^3, \mathcal{E}_{\mathcal{O}_{\mathbb{P}_\mathbb{C}^3}(4)}\right)\xrightarrow{\tilde{E}_0} H^0 \left( \mathbb{P}_\mathbb{C}^3, \Omega_{\mathbb{P}_\mathbb{C}^3}^1(4) \right) \right)= 4$, and so 
\begin{align*}
\dim_\mathbb{C} \mathbb{H}^1\left( \mathbb{P}_\mathbb{C}^3, \mathcal{E}_{\mathcal{O}_{\mathbb{P}_\mathbb{C}^3}(4)}^\bullet \right)=  18- 12 = 6
\end{align*}

By considering $(\ref{j1})$ and $(\ref{j2})$, we see that
\begin{align*}
\mathbb{H}^1 \left( \mathbb{P}_\mathbb{C}^3, \mathcal{E}_{\mathcal{O}_{\mathbb{P}_\mathbb{C}^3}(4)}^\bullet \right)= \textnormal{Span}_\mathbb{C} \left\langle dX_2, X_1 dX_2, X_2 dX_2, X_1^2 dX_2 ,  X_1X_2 dX_2, X_1 dX_1 \right\rangle
\end{align*}
Then we see that $\left( \mathbb{P}_\mathbb{C}^3\times B  , w , \pi, B \right)$, where $\pi:\mathbb{P}_\mathbb{C}^3 \times B \to B$ is the natural projection with the coordinates $\left( t_1,..., t_6 \right)$ on $B \subset \mathbb{C}^6$ for a sufficiently small neighborhood $B$ of $0$ in $\mathbb{C}^6$ and
\begin{align}\label{j3}
w&=  \left( X_1^2 X_2 - 1+ t_6 X_1 \right) dX_1 + \left(  X_2^2 - X_1^3 +t_1 + t_2 X_1+ t_3 X_2+ t_4 X_1^2 + t_5 X_1X_2\right) dX_2
\end{align}
defines a foliated analytic family of deformations of $\pi^{-1}(0)=\left( \mathbb{P}_\mathbb{C}^3, w_0  \right)$ in terms of cotangent sheaves since $w\wedge dw=0$, and the foliated Kodaira-Spencer map in terms of cotangent sheaves is an isomorphism at $t=0$, so that $\left( \mathbb{P}_\mathbb{C}^3, w_0 \right)$ is unobstructed in foliated deformations in terms of cotangent sheaves. Moreover we note that $w$ in $(\ref{j3})$ is a family of foliations which are linear publlack of some foliations on $\mathbb{P}_\mathbb{C}^2$. Since the Kodaira-Spencer map in terms of cotangent sheaves is an isomorphism at $t=0$, the above foliated analytic family is complete at $t=0$ by Theorem of completeness for deformations of foliated analytic structures in terms of cotangent sheaves in Part I. Since any small foliated deformations of $\left( \mathbb{P}_\mathbb{C}^3, w_0 \right)$ in terms of cotangent sheaves is induced from the complete foliated analytic family in terms of cotangent sheaves, we see that  any small foliated deformations of $\left( \mathbb{P}_\mathbb{C}^3, w_0 \right)$ in terms of cotangent sheaves are also linear-pullback of some foliations on $\mathbb{P}_\mathbb{C}^2$. We also note that $w$ in $(\ref{j3})$ is same to $w$ in $(\ref{j5})$ from Example \ref{cp10}. Since the family $\left( \mathbb{P}_\mathbb{C}^2\times B, w, B, \pi \right)$ in Example \ref{cp10} is also complete at $t=0$, we can conclude that any small foliated deformations of $\left( \mathbb{P}_\mathbb{C}^3 , w_0 \right)$ in terms of cotangent sheaves is induced by a linear pullback of some foliated deformations of $\mathbb{P}_\mathbb{C}^2$ with the Jouanolou foliation.

\subsection{Foliated deformations of $\mathbb{P}_\mathbb{C}^3$ with linear-pullback of Jouanolou foliation of degree $2$ on $\mathbb{P}_\mathbb{C}^3$ in terms of tangent sheaves} \

The Jouanolou foliation of degree $2$ on $\mathbb{P}_\mathbb{C}^2$ defined by the homogenous vector field $T= \xi_1^2 \frac{\partial}{\partial \xi_0} + \xi_2^2 \frac{\partial}{\partial \xi_1} + \xi_0^2 \frac{\partial}{\partial \xi_2}$ induces the foliation $\Theta_{\mathcal{F}_0}$ on $\mathbb{P}_\mathbb{C}^3$ defined by linear pullaback which is generated by $T$ and $\frac{\partial}{\partial \xi_3}$ on $\mathbb{P}_\mathbb{C}^3$ with the homogeneous coordinates $\left[ \xi_0 : \xi_1: \xi_2 : \xi_3 \right]$. This implies that $\Theta_{\mathcal{F}_0} \cong \mathcal{O}_{\mathbb{P}_{\mathbb{C}}^3}( -1) \bigoplus \mathcal{O}_{\mathbb{P}_\mathbb{C}^3}(1)$ which is locally free, so that we can consider the foliated deformations of $\left( \mathbb{P}_\mathbb{C}^3, \Theta_{\mathcal{F}_0} \right)$ in terms of tangent sheaves. Since we have an injection $\mathbb{H}^1\left( \mathbb{P}_\mathbb{C}^3, \Theta_{\mathcal{F}_0}^\bullet \right) \hookrightarrow \mathbb{H}^1\left( \mathbb{P}_\mathbb{C}^3, \mathcal{N}_{\mathcal{F}_0}^{* \bullet}\right)$ $($see \textnormal{Part II}$)$, and as shown in the previous \textnormal{subsection \ref{j6}}, the foliated deformations of $\left(\mathbb{P}_\mathbb{C}^3, \mathcal{N}_{\mathcal{F}_0}^* \right)$ in terms of cotangent sheaves remains in linear pullbacks, we should have $\mathbb{H}^1\left( \mathbb{P}_\mathbb{C}^3 , \Theta_{\mathcal{F}_0}^\bullet \right)\cong \mathbb{H}^1\left( \mathbb{P}_\mathbb{C}^3, \mathcal{N}_{\mathcal{F}_0}^{*\bullet} \right)$ and have the same result from the previous \textnormal{subsection \ref{j6}}. However we directly compute $\mathbb{H}^1\left( \mathbb{P}_\mathbb{C}^3 , \Theta_{\mathcal{F}_0}^\bullet \right)$ and show the same result with foliated deformations of $\left( \mathbb{P}_\mathbb{C}^3, \mathcal{N}_{\mathcal{F}_0}^{*} \right)$ in terms of cotangent sheaves in the previous \textnormal{subsection \ref{j6}}. We keep the notations in \textnormal{Appendix \ref{cp7}}. We note that $\Theta_{\mathcal{F}_0}$ on $U_0 =\left\{ \xi_0\ne 0 \right\}$ is generated by
\begin{align*}
T_1= \left(X_2^2- X_1^3\right)\frac{\partial}{\partial X_1} - \left(X_1^2X_2 -1\right)\frac{\partial}{\partial X_2}- X_1^2 X_3 \frac{\partial}{\partial X_3},\,\,\,\,\,\,\,\,\,\,\,\,T_2=\frac{\partial}{\partial X_3}
\end{align*}
Then we have
\begin{align*}
\left[ T_1, T_2\right]= X_1^2 T_2
\end{align*}

Let us consider foliated deformations of $\left(\mathbb{P}_\mathbb{C}^3, \Theta_{\mathcal{F}_0}=\mathcal{O}_{\mathbb{P}_\mathbb{C}^3}(-1)\bigoplus \mathcal{O}_{\mathbb{P}_\mathbb{C}^3}(1)\right)$ in terms of tangent sheaves. We have an exact sequence
\begin{align}
0\to \mathcal{O}_{\mathbb{P}_\mathbb{C}^3}(-1)\bigoplus \mathcal{O}_{\mathbb{P}_\mathbb{C}^3}(1)&\to \Theta_{\mathbb{P}_\mathbb{C}^3} \to  \frac{\Theta_{\mathbb{P}_\mathbb{C}^3}}{\Theta_{\mathcal{F}_0}} \to 0 \label{j8} \\
 (1,0)&\mapsto  T_1= (X_2^2- X_1^3)\frac{\partial}{\partial X_1} - (X_1^2X_2 -1)\frac{\partial}{\partial X_2}- X_1^2 X_3 \frac{\partial}{\partial X_3} \notag \\
 (0,1) &\mapsto T_2= \frac{\partial}{\partial X_3} \notag
\end{align}

We note that we have the leaf complex
\begin{align}\label{j7}
\Theta_{\mathcal{F}_0}^\bullet: \Theta_{\mathbb{P}_\mathbb{C}^3} \to \mathscr{H}om_{\mathcal{O}_{\mathbb{P}_\mathbb{C}^3}}\left( \mathcal{O}_{\mathbb{P}_{\mathbb{C}^3}} (-1) \bigoplus \mathcal{O}_{\mathbb{P}_\mathbb{C}^3}(1) , \frac{\Theta_{ \mathbb{P}_\mathbb{C}^3}}{ \Theta_{\mathcal{F}_0} } \right) \to \mathscr{H}om_{\mathcal{O}_{\mathbb{P}_\mathbb{C}^3}}\left( \mathcal{O}_{\mathbb{P}_\mathbb{C}^3}, \frac{ \Theta_{\mathbb{P}_\mathbb{C}^3} }{\Theta_{\mathcal{F}_0} } \right) \to 0
\end{align}
Here we note that $\bigwedge^2 \Theta_{\mathcal{F}_0} \cong \mathcal{O}_{\mathbb{P}_\mathbb{C}^3}(-1)\otimes \mathcal{O}_{\mathbb{P}_\mathbb{C}^3}(1)\cong \mathcal{O}_{\mathbb{P}_\mathbb{C}^3}$. Then $(\ref{j7})$ is equivalent to
\begin{align*}
\Theta_{\mathcal{F}_0}^\bullet : \Theta_{\mathbb{P}_\mathbb{C}^3} \to  \frac{\Theta_{\mathbb{P}_\mathbb{C}^3}}{\Theta_{\mathcal{F}_0}}(1) \bigoplus \frac{\Theta_{\mathbb{P}_\mathbb{C}^3}}{\Theta_{\mathcal{F}_0}}(-1) \to \frac{\Theta_{\mathbb{P}_\mathbb{C}^3}}{ \Theta_{\mathcal{F}_0} }\to 0
\end{align*}
Hence in order to compute $\mathbb{H}^i\left( \mathbb{P}_\mathbb{C}^3, \Theta_{\mathcal{F}_0}^\bullet \right)$, we need to describe the cohomology groups $H^i\left(\mathbb{P}_\mathbb{C}^3 , \frac{\Theta_{\mathbb{P}_\mathbb{C}^3}}{\Theta_{\mathcal{F}_0}}(1) \right)$ and $H^i\left(\mathbb{P}_\mathbb{C}^3 , \frac{\Theta_{\mathbb{P}_\mathbb{C}^3}}{\Theta_{\mathcal{F}_0}}(-1) \right)$ and $H^i\left(\mathbb{P}_\mathbb{C}^3 , \frac{\Theta_{\mathbb{P}_\mathbb{C}^3}}{\Theta_{\mathcal{F}_0}} \right)$ for $i\geq 0$.

\subsubsection{Descriptions of $H^i\left( \mathbb{P}_\mathbb{C}^3, \frac{\Theta_{\mathbb{P}_\mathbb{C}^3}}{\Theta_{\mathcal{F}_0}}(1) \right)$}\

By tensoring the exact sequence $(\ref{j8})$ with $\mathcal{O}_{\mathbb{P}_\mathbb{C}^3}(1)$, we have the exact sequence
\begin{align*}
0\to \mathcal{O}_{\mathbb{P}_\mathbb{C}^3}\bigoplus \mathcal{O}_{\mathbb{P}_\mathbb{C}^3}(2) \to \Theta_{\mathbb{P}_\mathbb{C}^3} (1) \to  \frac{\Theta_{\mathbb{P}_\mathbb{C}^3}}{\Theta_{\mathcal{F}_0}} (1) \to 0
\end{align*}
Then we see that $\mathbb{H}^i\left( \mathbb{P}_\mathbb{C}^3, \frac{\Theta_{\mathbb{P}_\mathbb{C}^3}}{\Theta_{\mathcal{F}_0}}(1) \right)=0$ for $i>0$, and we have the exact sequence
\begin{align*}
0\to H^0\left(\mathbb{P}_\mathbb{C}^3,  \mathcal{O}_{\mathbb{P}_\mathbb{C}^3} \right) \bigoplus H^0\left( \mathbb{P}_\mathbb{C}^3,  \mathcal{O}_{\mathbb{P}_\mathbb{C}^3}(2) \right) \cong \mathbb{C}^{11} \to H^0\left( \mathbb{P}_\mathbb{C}^3,  \Theta_{\mathbb{P}_\mathbb{C}^3} (1) \right) \cong \mathbb{C}^{36} \to H^0\left(\mathbb{P}_\mathbb{C}^3,  \frac{\Theta_{\mathbb{P}_\mathbb{C}^3}}{\Theta_{\mathcal{F}_0}} (1) \right) \to 0,
\end{align*}
so that we have
\begin{align*}
H^0 \left( \mathbb{P}_\mathbb{C}^3, \frac{\Theta_{\mathbb{P}_\mathbb{C}^3}}{\Theta_{\mathcal{F}_0}}(1)  \right)\cong \frac{H^0\left(\mathbb{P}_\mathbb{C}^3, \Theta_{\mathbb{P}_\mathbb{C}^3}(1) \right)}{H^0\left( \mathbb{P}_\mathbb{C}^3, \mathcal{O}_{\mathbb{P}_\mathbb{C}^3} \right) \bigoplus H^0\left( \mathbb{P}_\mathbb{C}^3, \mathcal{O}_{\mathbb{P}_\mathbb{C}^3}(2) \right)} \cong \mathbb{C}^{25}
\end{align*}
We explicitly describe $H^0\left( \mathbb{P}_\mathbb{C}^3, \frac{\Theta_{\mathbb{P}_\mathbb{C}^3}}{\Theta_{\mathcal{F}_0}}(1) \right)$ on $U_0$. We recall from $\textnormal{(\ref{j9})}$ that
\begin{align}\label{j12}
H^0\left( \mathbb{P}_\mathbb{C}^3, \Theta_{\mathbb{P}_\mathbb{C}^3}(1) \right)=&\left\{\left(t_1+t_2X_3+t_3X_2+t_4X_1+t_5X_3^2+t_6X_2X_3+t_7X_2^2+t_8X_1X_3+t_9X_1X_2 + t_{10}X_1^2\right.\right.\\
&\left. \,\,\,\,\,\,\,\,\,\, +t_{11}X_1X_3^2+t_{12}X_1X_2X_3+t_{13}X_1X_2^2+t_{14}X_1^2X_3+t_{15}X_1^2X_2+t_{16}X_1^3 \right)\frac{\partial}{\partial X_1} \notag \\
&+\left(t_{17}+t_{18}X_3+t_{19}X_2 +t_{20}X_1+t_{21}X_3^2+t_{22}X_2X_3 +t_{23}X_2^2+t_{24}X_1X_3 +t_{25}X_1X_2+t_{26}X_1^2\right. \notag \\
&\left.\,\,\,\,\,\,\,\,\,\,\,\,\,\,+t_{11}X_2X_3^2+t_{12}X_2^2X_3+t_{13}X_2^3 +t_{14}X_1X_2X_3+t_{15}X_1X_2^2+t_{16}X_1^2X_2\right)\frac{\partial}{\partial X_2} \notag \\
&+\left(t_{27}+t_{28}X_3+t_{29}X_2+t_{30}X_1 +t_{31}X_3^2+t_{32}X_2X_3 +t_{33}X_2^2 +t_{34}X_1X_3+t_{35}X_1X_2 +t_{36}X_1^2\right. \notag \\
&\,\,\,\,\,\,\,\,\,\,\,\,\,+\left.\left. t_{11}X_3^3+t_{12}X_2X_3^2+t_{13}X_2^2X_3 +t_{14}X_1X_3^2+t_{15}X_1X_2X_3+t_{16}X_1^2X_3\right)\frac{\partial}{\partial X_3}, t_1,..., t_{36} \in \mathbb{C} \right\} \notag
\end{align}
Since we have in $\frac{H^0\left(\mathbb{P}_\mathbb{C}^3, \Theta_{\mathbb{P}_\mathbb{C}^3}(1) \right)}{H^0\left( \mathbb{P}_\mathbb{C}^3, \mathcal{O}_{\mathbb{P}_\mathbb{C}^3} \right) \bigoplus H^0\left( \mathbb{P}_\mathbb{C}^3, \mathcal{O}_{\mathbb{P}_\mathbb{C}^3}(2) \right)}$
\begin{align} \label{j14}
&X_1^3 \frac{\partial}{\partial X_1} + X_1^2 X_2\frac{\partial}{\partial X_2} + X_1^2 X_3 \frac{\partial}{\partial X_3}\equiv X_2^2 \frac{\partial}{\partial X_1} + \frac{\partial}{\partial X_2} \\
&\frac{\partial}{\partial X_3}, X_1\frac{\partial}{\partial X_3}, X_2\frac{\partial}{\partial X_3}, X_3\frac{\partial}{\partial X_3},X_1^2\frac{\partial}{\partial X_3}, X_2^2\frac{\partial}{\partial X_3}, X_3^2, X_1X_2\frac{\partial}{\partial X_3}, X_1X_3\frac{\partial}{\partial X_3}, X_2X_3 \frac{\partial}{\partial X_3} \equiv 0 \notag
\end{align}
we see that from $(\ref{j12})$
\begin{align}\label{j11}
H^0\left( \mathbb{P}_\mathbb{C}^3, \frac{\Theta_{\mathbb{P}_\mathbb{C}^3}}{\Theta_{\mathcal{F}_0}}(1) \right)=&\left\{\left(t_1+t_2X_3+t_3X_2+t_4X_1+t_5X_3^2+t_6X_2X_3+t_7X_2^2+t_8X_1X_3+t_9X_1X_2 + t_{10}X_1^2\right.\right.\\
&\left. \,\,\,\,\,\,\,\,\,\, +t_{11}X_1X_3^2+t_{12}X_1X_2X_3+t_{13}X_1X_2^2+t_{14}X_1^2X_3+t_{15}X_1^2X_2 \right)\frac{\partial}{\partial X_1} \notag \\
&+\left(t_{16}+t_{17}X_3+t_{18}X_2 +t_{19}X_1+t_{20}X_3^2+t_{21}X_2X_3 +t_{22}X_2^2+t_{23}X_1X_3 +t_{24}X_1X_2+t_{25}X_1^2\right. \notag \\
&\left.\,\,\,\,\,\,\,\,\,\,\,\,\,\,+t_{11}X_2X_3^2+t_{12}X_2^2X_3+t_{13}X_2^3 +t_{14}X_1X_2X_3+t_{15}X_1X_2^2\right)\frac{\partial}{\partial X_2} \notag \\
&\,\,\,\,\,\,\,\,\,\,\,\,\,+\left.\left( t_{11}X_3^3+t_{12}X_2X_3^2+t_{13}X_2^2X_3 +t_{14}X_1X_3^2+t_{15}X_1X_2X_3\right)\frac{\partial}{\partial X_3}, t_1,..., t_{25} \in \mathbb{C} \right\} \notag
\end{align}
Since $\frac{\partial}{\partial X_3}\equiv 0 $ on $U_0$, we can ignore $\frac{\partial}{\partial X_3}$ in actual computation. It does not need to be induced from global section $H^0\left(\mathbb{P}_\mathbb{C}^3, \Theta_{\mathbb{P}_\mathbb{C}^3}(1) \right)$, so that we have
\begin{align}\label{j20}
H^0\left( \mathbb{P}_\mathbb{C}^3, \frac{\Theta_{\mathbb{P}_\mathbb{C}^3}}{\Theta_{\mathcal{F}_0}}(1) \right)=&\left\{\left(t_1+t_2X_3+t_3X_2+t_4X_1+t_5X_3^2+t_6X_2X_3+t_7X_2^2+t_8X_1X_3+t_9X_1X_2 + t_{10}X_1^2\right.\right.\\
&\left. \,\,\,\,\,\,\,\,\,\, +t_{11}X_1X_3^2+t_{12}X_1X_2X_3+t_{13}X_1X_2^2+t_{14}X_1^2X_3+t_{15}X_1^2X_2 \right)\frac{\partial}{\partial X_1} \notag \\
&+\left(t_{16}+t_{17}X_3+t_{18}X_2 +t_{19}X_1+t_{20}X_3^2+t_{21}X_2X_3 +t_{22}X_2^2+t_{23}X_1X_3 +t_{24}X_1X_2+t_{25}X_1^2\right. \notag \\
&\left. \left.\,\,\,\,\,\,\,\,\,\,\,\,\,\,+t_{11}X_2X_3^2+t_{12}X_2^2X_3+t_{13}X_2^3 +t_{14}X_1X_2X_3+t_{15}X_1X_2^2\right)\frac{\partial}{\partial X_2}, t_1,..., t_{25} \in \mathbb{C} \right\} \notag
\end{align}

\subsubsection{Descriptions of $H^i\left( \mathbb{P}_\mathbb{C}^3, \frac{\Theta_{\mathbb{P}_\mathbb{C}^3}}{\Theta_{\mathcal{F}_0}}(-1) \right)$}\

By tensoring the exact sequence $(\ref{j8})$ with $\mathcal{O}_{\mathbb{P}_\mathbb{C}^3}(-1)$, we have the exact sequence
\begin{align*}
0\to \mathcal{O}_{\mathbb{P}_\mathbb{C}^3} (-2)\bigoplus \mathcal{O}_{\mathbb{P}_\mathbb{C}^3} \to \Theta_{\mathbb{P}_\mathbb{C}^3} (-1) \to  \frac{\Theta_{\mathbb{P}_\mathbb{C}^3}}{\Theta_{\mathcal{F}_0}} (-1) \to 0
\end{align*}
Then we see that $\mathbb{H}^i\left( \mathbb{P}_\mathbb{C}^3, \frac{\Theta_{\mathbb{P}_\mathbb{C}^3}}{\Theta_{\mathcal{F}_0}}(-1) \right)=0$ for $i>0$, and we have the exact sequence
\begin{align*}
0\to H^0\left(\mathbb{P}_\mathbb{C}^3, \mathcal{O}_{\mathbb{P}_\mathbb{C}^3}(-2) \right) \bigoplus H^0\left( \mathbb{P}_\mathbb{C}^3, \mathcal{O}_{\mathbb{P}_\mathbb{C}^3} \right) \cong \mathbb{C} \to H^0\left( \mathbb{P}_\mathbb{C}^3, \Theta_{\mathbb{P}_\mathbb{C}^3}(-1) \right) \cong \mathbb{C}^4 \to H^0\left( \mathbb{P}_\mathbb{C}^3, \frac{\Theta_{\mathbb{P}_\mathbb{C}^3}}{\Theta_{\mathcal{F}_0}}(-1) \right) \to 0,
\end{align*}
so that we have
\begin{align*}
H^0\left( \mathbb{P}_\mathbb{C}^3 , \frac{\Theta_{\mathbb{P}_\mathbb{C}^3}}{\Theta_{\mathcal{F}_0}}(-1) \right) \cong \frac{H^0\left( \mathbb{P}_\mathbb{C}^3, \Theta_{\mathbb{P}_\mathbb{C}^3}(-1) \right) }{H^0\left(\mathbb{P}_\mathbb{C}^3, \mathcal{O}_{\mathbb{P}_\mathbb{C}^3}(-2) \right) \bigoplus H^0\left( \mathbb{P}_\mathbb{C}^3, \mathcal{O}_{\mathbb{P}_\mathbb{C}^3} \right)} \cong \mathbb{C}^3
\end{align*}
We recall from $\textnormal{(\ref{j10})}$ that
\begin{align}\label{j15}
H^0\left( \mathbb{P}_{\mathbb{C}^3}, \Theta_{\mathbb{P}_\mathbb{C}^3}(-1) \right)= \left\{  \left( t_1 + t_2 X_1 \right) \frac{\partial}{\partial X_1} + \left( t_3 + t_2 X_2 \right) \frac{\partial}{\partial X_2} + \left( t_4 + t_2 X_3 \right) \frac{\partial}{\partial X_3}, t_1,...,t_4 \in \mathbb{C} \right\}
\end{align}
Since $\frac{\partial}{\partial X_3}\equiv 0$ in $\frac{H^0\left( \mathbb{P}_\mathbb{C}^3, \Theta_{\mathbb{P}_\mathbb{C}^3}(-1) \right) }{H^0\left(\mathbb{P}_\mathbb{C}^3, \mathcal{O}_{\mathbb{P}_\mathbb{C}^3}(-2) \right) \bigoplus H^0\left( \mathbb{P}_\mathbb{C}^3, \mathcal{O}_{\mathbb{P}_\mathbb{C}^3} \right)} $,  we have
\begin{align}\label{j27}
H^0\left( \mathbb{P}_\mathbb{C}^3, \frac{\Theta_{\mathbb{P}_\mathbb{C}^3}}{\Theta_{\mathcal{F}_0}}(-1) \right) =\left\{ \left( b_1 + b_2 X_1 \right)\frac{\partial}{\partial X_1} + \left( b_3 + b_2 X_2 \right) \frac{\partial}{\partial X_2} +  b_2X_3 \frac{\partial}{\partial X_3} | b_1,b_2,b_3 \in \mathbb{C}  \right\}
\end{align}
Since $\frac{\partial}{\partial X_3}\equiv 0 $ on $U_0$, we can ignore $\frac{\partial}{\partial X_3}$ in actual computation. It does not need to be induced from global section $H^0\left(\mathbb{P}_\mathbb{C}^3, \Theta_{\mathbb{P}_\mathbb{C}^3}(-1) \right)$, so that we have
\begin{align}\label{j21}
H^0\left( \mathbb{P}_\mathbb{C}^3, \frac{\Theta_{\mathbb{P}_\mathbb{C}^3}}{\Theta_{\mathcal{F}_0}}(-1) \right) =\left\{ \left( b_1 + b_2 X_1 \right)\frac{\partial}{\partial X_1} + \left( b_3 + b_2 X_2 \right) \frac{\partial}{\partial X_2} | b_1,b_2,b_3 \in \mathbb{C}  \right\}
\end{align}

\subsubsection{Descriptions of $H^i\left( \mathbb{P}_\mathbb{C}^3, \frac{\Theta_{\mathbb{P}_\mathbb{C}^3}}{\Theta_{\mathcal{F}_0}} \right)$}\

From the exact sequence $(\ref{j8})$, we have $H^i\left( \mathbb{P}_\mathbb{C}^3,  \frac{\Theta_{\mathbb{P}_\mathbb{C}^3}}{\Theta_{\mathcal{F}_0}} \right)=0$ for $i>0$, and we have the exact sequence
\begin{align*}
0 \to H^0\left(\mathbb{P}_\mathbb{C}^3, \mathcal{O}_{\mathbb{P}_\mathbb{C}^3}( -1) \right) \bigoplus H^0\left( \mathbb{P}_\mathbb{C}^3, \mathcal{O}_{\mathbb{P}_\mathbb{C}^3} (1)  \right) \cong \mathbb{C}^4  &\to H^0\left( \mathbb{P}_\mathbb{C}^3, \Theta_{\mathbb{P}_\mathbb{C}^3}\right)\cong \mathbb{C}^{15} \to H^0\left( \mathbb{P}_\mathbb{C}^3, \frac{\Theta_{\mathbb{P}_\mathbb{C}^3} }{\Theta_{\mathcal{F}_0}} \right)\to 0,
\end{align*}
so that we have 
\begin{align*}
H^0\left( \mathbb{P}_\mathbb{C}^3, \frac{\Theta_{\mathbb{P}_\mathbb{C}^3}}{\Theta_{\mathcal{F}_0}} \right) \cong \frac{H^0\left(\mathbb{P}_\mathbb{C}^3, \Theta_{\mathbb{P}_\mathbb{C}^3} \right)}{ H^0\left(\mathbb{P}_\mathbb{C}^3, \mathcal{O}_{\mathbb{P}_\mathbb{C}^3}( -1) \right) \bigoplus H^0\left( \mathbb{P}_\mathbb{C}^3, \mathcal{O}_{\mathbb{P}_\mathbb{C}^3} (1)  \right)} \cong \mathbb{C}^{11}
\end{align*}
We recall from $\textnormal{(\ref{ap2})}$
\begin{align}\label{j13}
H^0\left( \mathbb{P}_\mathbb{C}^3,  \Theta_{\mathbb{P}_\mathbb{C}^3}  \right) = &\left\{V=  V_1\frac{\partial}{\partial X_1} + V_2 \frac{\partial}{\partial X_2} + V_3 \frac{\partial}{\partial X_3} \right.\\
&\left. =\left( a_1 + a_2 X_3+ a_3 X_2 + a_4 X_1 + a_5 X_1X_3 + a_6 X_1 X_2 + a_7 X_1^2 \right) \frac{\partial}{\partial X_1}\right. \notag \\
 &\left.+\left( a_8 + a_9 X_3 + a_{10}X_2 + a_{11}X_1 + a_5 X_2X_3 + a_6 X_2^2 + a_7 X_1X_2\right) \frac{\partial}{\partial X_2} \right. \notag \\ 
 &\left.  +\left( a_{12}+ a_{13} X_3 + a_{14} X_2 + a_{15}X_1 + a_5 X_3^2 + a_6 X_2 X_3  + a_7 X_1X_3 \right)  \frac{\partial}{\partial X_3}, a_1,..., a_{15}\in \mathbb{C} \right\} \notag
\end{align}
Since we have
\begin{align}\label{j19}
\frac{\partial}{\partial X_3}, X_1\frac{\partial}{\partial X_3}, X_2 \frac{\partial}{\partial X_3}, X_3 \frac{\partial}{\partial X_3}\equiv 0 \in \frac{H^0\left( \mathbb{P}_\mathbb{C}^3, \Theta_{\mathbb{P}_\mathbb{C}^3} \right)}{H^0\left( \mathbb{P}_\mathbb{C}^3, \mathcal{O}_{\mathbb{P}_\mathbb{C}^3}(-1) \right) \bigoplus H^0 \left( \mathbb{P}_\mathbb{C}^3, \mathcal{O}_{\mathbb{P}_\mathbb{C}^3}(1) \right) }, 
\end{align}
we have
 \begin{align}\label{j17}
 H^0\left(\mathbb{P}_\mathbb{C}^3 , \frac{\Theta_{\mathbb{P}_\mathbb{C}^3}}{\Theta_{\mathcal{F}_0}}\right)=&\{ \left(c_1+c_2X_3+c_3X_2+c_4X_1+c_5X_1X_3+c_6X_1X_2+c_7X_1^2\right)\frac{\partial}{\partial X_1} \\
 &+\left(c_8+c_9X_3+c_{10}X_2 +c_{11}X_1+c_5X_2X_3+c_6X_2^2+c_7X_1X_2 \right)\frac{\partial}{\partial X_2} \notag\\
 &+\left(c_5X_3^2+c_6X_2X_3+c_7X_1X_3 \right)\frac{\partial}{\partial X_3}\} \notag
 \end{align}
Since $\frac{\partial}{\partial X_3}\equiv 0 $ on $U_0$, we can ignore $\frac{\partial}{\partial X_3}$ in actual computation. It does not need to be induced from global section $H^0\left(\mathbb{P}_\mathbb{C}^3 , \Theta_{\mathbb{P}_\mathbb{C}^3} \right)$, so that
 \begin{align}\label{j22}
 H^0\left(\mathbb{P}_\mathbb{C}^3 , \frac{\Theta_{\mathbb{P}_\mathbb{C}^3}}{\Theta_{\mathcal{F}_0}}\right)=&\{ \left(c_1+c_2X_3+c_3X_2+c_4X_1+c_5X_1X_3+c_6X_1X_2+c_7X_1^2\right)\frac{\partial}{\partial X_1} \\
 &\left.+\left(c_8+c_9X_3+c_{10}X_2 +c_{11}X_1+c_5X_2X_3+c_6X_2^2+c_7X_1X_2 \right)\frac{\partial}{\partial X_2}  \right\} \notag
 \end{align}

Hence we can compute $\mathbb{H}^i\left( \mathbb{P}_\mathbb{C}^3, \Theta_{\mathcal{F}_0}^\bullet \right)$ by the following complex
\begin{align*}
H^0\left(\mathbb{P}_\mathbb{C}^3, \Theta_{\mathbb{P}_\mathbb{C}^3} \right)  \cong \mathbb{C}^{15}\xrightarrow{D_0} H^0\left( \mathbb{P}_\mathbb{C}^3, \frac{\Theta_{\mathbb{P}_\mathbb{C}^3}}{\Theta_{\mathcal{F}_0}}(1) \right) \bigoplus H^0\left( \mathbb{P}_\mathbb{C}^3, \frac{\Theta_{\mathbb{P}_\mathbb{C}^3}}{\Theta_{\mathcal{F}_0}} (-1)\right) \cong \mathbb{C}^{28} \xrightarrow{D_1} H^0\left( \mathbb{P}_\mathbb{C}^3, \frac{\Theta_{\mathbb{P}_\mathbb{C}^3}}{\Theta_{\mathcal{F}_0}} \right) \cong\mathbb{C}^{11} \to 0
\end{align*}
Equivalently 
\begin{align*}
H^0\left( \mathbb{P}_\mathbb{C}^3, \Theta_{\mathbb{P}_\mathbb{C}^3} \right)\xrightarrow{D_0} H^0\left( X, \mathscr{H}om_{\mathcal{O}_{\mathbb{P}_\mathbb{C}^3}}\left( \mathcal{O}_{\mathbb{P}_\mathbb{C}^3}(-1) \bigoplus \mathcal{O}_{\mathbb{P}_\mathbb{C}^3}(1) ,  \frac{\Theta_{\mathbb{P}_\mathbb{C}^3}}{\Theta_{\mathcal{F}_0}} \right)  \right)  \xrightarrow{D_1} H^0\left( \mathbb{P}_\mathbb{C}^3,  \mathscr{H}om_{\mathcal{O}_{\mathbb{P}_\mathbb{C}^3}} \left( \mathcal{O}_{\mathbb{P}_\mathbb{C}^3},  \frac{\Theta_{\mathbb{P}_\mathbb{C}^3}}{\Theta_{\mathcal{F}_0} }  \right) \right)\to 0
\end{align*}
We note that from $(\ref{j11})$ and $(\ref{j21})$ any element $\phi\in H^0\left( \mathbb{P}_\mathbb{C}^3, \mathscr{H}om_{\mathcal{O}_{\mathbb{P}_\mathbb{C}^3}}\left( \mathcal{O}_{\mathbb{P}_\mathbb{C}^3}(-1) \bigoplus \mathcal{O}_{\mathbb{P}_\mathbb{C}^3}(1) ,  \frac{\Theta_{\mathbb{P}_\mathbb{C}^3}}{\Theta_{\mathcal{F}_0}} \right) \right)$ is represented by
{\small{\begin{align} \label{j18}
\phi(T_1)&=\left(t_1+t_2X_3+t_3X_2+t_4X_1+t_5X_3^2+t_6X_2X_3+t_7X_2^2+t_8X_1X_3+t_9X_1X_2 + t_{10}X_1^2\right.\\
&\left. +t_{11}X_1X_3^2+t_{12}X_1X_2X_3+t_{13}X_1X_2^2+t_{14}X_1^2X_3+t_{15}X_1^2X_2 \right)\frac{\partial}{\partial X_1} \notag \\
&+\left(t_{16}+t_{17}X_3+t_{18}X_2 +t_{19}X_1+t_{20}X_3^2+t_{21}X_2X_3 +t_{22}X_2^2+t_{23}X_1X_3 +t_{24}X_1X_2+t_{25}X_1^2\right.  \notag \\
&\left.+t_{11}X_2X_3^2+t_{12}X_2^2X_3+t_{13}X_2^3 +t_{14}X_1X_2X_3+t_{15}X_1X_2^2\right)\frac{\partial}{\partial X_2} \notag \\
&+ \left( t_{11}X_3^3+t_{12}X_2X_3^2+t_{13}X_2^2X_3 +t_{14}X_1X_3^2+t_{15}X_1X_2X_3\right)\frac{\partial}{\partial X_3}, \notag \\
\phi(T_2)&= \left(b_1+ b_2 X_1 \right)\frac{\partial}{\partial X_1} + \left( b_3 + b_2 X_2 \right) \frac{\partial}{\partial X_2}+ b_2 X_3 \frac{\partial}{\partial X_3},\,\,\,\,\,\,\,t_1,...,t_{25}, b_1,b_2,b_3\in \mathbb{C} \notag
\end{align}}}

We note that from $(\ref{j22})$ any element $\phi$ of $H^0\left( \mathbb{P}_\mathbb{C}^3, \mathscr{H}om_{\mathcal{O}_{\mathbb{P}_\mathbb{C}^3}} \left( \mathcal{O}_{\mathbb{P}_{\mathbb{C}^3}}, \frac{\Theta_{\mathbb{P}_\mathbb{C}^3}}{\Theta_{\mathcal{F}_0}} \right) \right)$ is represented by
\begin{align} \label{j33}
\phi\left(T_1 \wedge T_2 \right) &= \left(c_1+c_2X_3+c_3X_2+c_4X_1+c_5X_1X_3+c_6X_1X_2+c_7X_1^2\right)\frac{\partial}{\partial X_1} \\
 &+\left(c_8+c_9X_3+c_{10}X_2 +c_{11}X_1+c_5X_2X_3+c_6X_2^2+c_7X_1X_2 \right)\frac{\partial}{\partial X_2},\,\,\,\,\,\,\,c_1,...,c_{11} \in \mathbb{C} \notag
\end{align}

We compute $H^0\left( \mathbb{P}_\mathbb{C}^3, \Theta_{\mathbb{P}_\mathbb{C}^3} \right)\xrightarrow{D_0} H^0\left( X, \mathscr{H}om_{\mathcal{O}_{\mathbb{P}_\mathbb{C}^3}}\left( \mathcal{O}_{\mathbb{P}_\mathbb{C}^3}(-1) \bigoplus \mathcal{O}_{\mathbb{P}_\mathbb{C}^3}(1) ,  \frac{\Theta_{\mathbb{P}_\mathbb{C}^3}}{\Theta_{\mathcal{F}_0}} \right)  \right)$. From $(\ref{j13})$ the image $\phi=D_0\left(V_1 \frac{\partial}{\partial X_1} + V_2 \frac{\partial}{\partial X_2} + V_3 \frac{\partial}{\partial X_3} \right)$ is given by
{\Small{\begin{align*}
\phi(T_1)&=\left[ V_1\frac{\partial}{\partial X_1} + V_2 \frac{\partial}{\partial X_2} + V_3 \frac{\partial}{\partial X_3}, \left(X_2^2- X_1^3\right)\frac{\partial}{\partial X_1} - \left(X_1^2X_2 -1\right)\frac{\partial}{\partial X_2}- X_1^2 X_3 \frac{\partial}{\partial X_3}\right]\\
&=\left( - a_3 - a_6X_1 - 3a_1 X_1^2+ (2a_{10} -a_4 ) X_2^2 - 2a_4 X_1^3 - 2a_3X_1^2X_2 - 2a_2 X_1^2 X_3 + a_6 X_2^3 + a_5X_2^2X_3 - a_6 X_1^3X_2 - a_5 X_1^3 X_3 - a_7X_1^4\right.\\
&\left.+2 a_8X_2+2a_9X_2X_3 +2a_{11}X_1X_2\right)\frac{\partial}{\partial X_1}\\
&+\left( - a_8X_1^2  - a_{11}X_2^2 -2a_1X_1X_2-2a_2X_1X_2X_3-2a_3X_1X_2^2-2a_4X_1^2 X_2 - a_7 X_2^3-a_5X_1^2 X_2X_3-a_6X_1^2X_2^2-a_7X_1^3X_2\right)\frac{\partial}{\partial X_2} \\
&+\left(- a_{14}- a_6X_3 -a_{12}X_1^2 - a_{15}X_2^2-2a_1X_1X_3-2a_2X_1X_3^2-2a_3X_1X_2X_3-2a_4X_1^2 X_3- a_7X_2^2 X_3-a_5X_1^2X_3^2\right.\\
&\left.-a_6X_1^2X_2X_3-a_7X_1^3X_3\right)\frac{\partial}{\partial X_3}
\end{align*}}}
Since $X_1^3\frac{\partial}{\partial X_1} + X_1^2X_2 \frac{\partial}{\partial X_2}+ X_1^2X_3 \frac{\partial}{\partial X_3}\equiv X_2^2\frac{\partial}{\partial X_1}+ \frac{\partial}{\partial X_2}$ and $\frac{\partial}{\partial X_3}\equiv 0$, we have in terms of $(\ref{j20})$
{\small{\begin{align*}
\phi(T_1) \equiv &\left( - a_3 - a_6X_1 - 3a_1 X_1^2+ \left(2a_{10} - 3a_4\right) X_2^2  - 2a_3X_1^2X_2 - 2a_2 X_1^2 X_3 - a_7X_1X_2^2+2 a_8X_2+2a_9X_2X_3 +2a_{11}X_1X_2\right)\frac{\partial}{\partial X_1}\\
&+\left( - a_8X_1^2  - a_{11}X_2^2 -2a_1X_1X_2-2a_2X_1X_2X_3-2a_3X_1X_2^2-2a_4 - a_7 X_2^3-a_5X_3-a_6X_2-a_7X_1\right)\frac{\partial}{\partial X_2}  \in H^0\left( \mathbb{P}_\mathbb{C}^3, \frac{\Theta_{\mathbb{P}_\mathbb{C}^3}}{\Theta_{\mathcal{F}_0}}(1) \right)
\end{align*}}}
We also have in terms of $(\ref{j11})$
{\small{\begin{align}
\phi(T_1) \equiv &\left( - a_3 - a_6X_1 - 3a_1 X_1^2+ \left(2a_{10} - 3a_4\right) X_2^2  - 2a_3X_1^2X_2 - 2a_2 X_1^2 X_3 - a_7X_1X_2^2+2 a_8X_2+2a_9X_2X_3 +2a_{11}X_1X_2\right)\frac{\partial}{\partial X_1}\label{j29}\\
&+\left( - a_8X_1^2  - a_{11}X_2^2 -2a_1X_1X_2-2a_2X_1X_2X_3-2a_3X_1X_2^2-2a_4 - a_7 X_2^3-a_5X_3-a_6X_2-a_7X_1 \right)\frac{\partial}{\partial X_2} \notag \\
&+ \left(- 2a_2 X_1 X_3^2 - 2a_3 X_1X_2X_3 - a_7 X_2^2 X_3\right) \frac{\partial}{\partial X_3} \in H^0\left( \mathbb{P}_\mathbb{C}^3, \frac{\Theta_{\mathbb{P}_\mathbb{C}^3}}{\Theta_{\mathcal{F}_0}}(1) \right) \notag
\end{align}}}

On the other hand, we have in terms of $(\ref{j21})$
{\Small{\begin{align*}
\phi\left(T_2\right)=\left[ V_1 \frac{\partial}{\partial X_1} + V_2 \frac{\partial}{\partial X_2} + V_3 \frac{\partial}{\partial X_3}, \frac{\partial}{\partial X_3} \right]& = (-a_2- a_5X_1)\frac{\partial}{\partial X_1} + (-a_9- a_5X_2)\frac{\partial}{\partial X_2} +(- a_{13} - 2a_5 X_3- a_6X_2 - a_7X_1)\frac{\partial}{\partial X_3}\\
&\equiv_{\Theta_{\mathcal{F}_0}} (-a_2- a_5X_1)\frac{\partial}{\partial X_1} + (-a_9- a_5X_2)\frac{\partial}{\partial X_2}  \in H^0\left( \mathbb{P}_\mathbb{C}^3 , \frac{\Theta_{\mathbb{P}_\mathbb{C}^3}}{\Theta_{\mathcal{F}_0}} (-1) \right)
\end{align*}}}
We also have in terms of $(\ref{j27})$
{\Small{\begin{align}\label{j28}
\phi\left(T_2 \right)& \equiv \left(-a_2- a_5X_1\right)\frac{\partial}{\partial X_1} + \left(-a_9- a_5X_2\right) \frac{\partial}{\partial X_2} -a_5 X_3\frac{\partial}{\partial X_3} \in H^0\left( \mathbb{P}_\mathbb{C}^3 , \frac{\Theta_{\mathbb{P}_\mathbb{C}^3}}{\Theta_{\mathcal{F}_0}} (-1) \right)
\end{align}}}

Then we see that from $(\ref{j20})$ and $(\ref{j21})$
\begin{align}
&\dim_\mathbb{C} \mathbb{H}^0\left( \mathbb{P}_\mathbb{C}^3, \Theta_{\mathcal{F}_0}^\bullet \right)= \dim_\mathbb{C} \textnormal{ker}\left( H^0\left( \mathbb{P}_\mathbb{C}^3, \Theta_{\mathbb{P}_\mathbb{C}^3} \right)\xrightarrow{D_0} H^0\left( X, \mathscr{H}om_{\mathcal{O}_{\mathbb{P}_\mathbb{C}^3}}\left( \mathcal{O}_{\mathbb{P}_\mathbb{C}^3}(-1) \bigoplus \mathcal{O}_{\mathbb{P}_\mathbb{C}^3}(1) ,  \frac{\Theta_{\mathbb{P}_\mathbb{C}^3}}{\Theta_{\mathcal{F}_0}} \right)  \right)\right)= 4 \notag \\
&\dim_\mathbb{C} \textnormal{im}\left( H^0\left( \mathbb{P}_\mathbb{C}^3, \Theta_{\mathbb{P}_\mathbb{C}^3} \right)\xrightarrow{D_0} H^0\left( X, \mathscr{H}om_{\mathcal{O}_{\mathbb{P}_\mathbb{C}^3}}\left( \mathcal{O}_{\mathbb{P}_\mathbb{C}^3}(-1) \bigoplus \mathcal{O}_{\mathbb{P}_\mathbb{C}^3}(1) ,  \frac{\Theta_{\mathbb{P}_\mathbb{C}^3}}{\Theta_{\mathcal{F}_0}} \right)  \right) \right)= 15- 4 =11 \label{j30}
\end{align}

Next we compute $H^0\left( X, \mathscr{H}om_{\mathcal{O}_{\mathbb{P}_\mathbb{C}^3}}\left( \mathcal{O}_{\mathbb{P}_\mathbb{C}^3}(-1) \bigoplus \mathcal{O}_{\mathbb{P}_\mathbb{C}^3}(1) ,  \frac{\Theta_{\mathbb{P}_\mathbb{C}^3}}{\Theta_{\mathcal{F}_0}} \right)  \right)  \xrightarrow{D_1} H^0\left( \mathbb{P}_\mathbb{C}^3,  \mathscr{H}om_{\mathcal{O}_{\mathbb{P}_\mathbb{C}^3}}\left( \mathcal{O}_{\mathbb{P}_\mathbb{C}^3}  ,\frac{\Theta_{\mathbb{P}_\mathbb{C}^3}}{\Theta_{\mathcal{F}_0}} \right)\right)$. From $(\ref{j18})$, the image $D_1\left(\phi\right)\left( T_1 \wedge T_2 \right) $ is given by $\left[ T_1, \phi\left(T_2\right)\right] - \left[ T_2, \phi\left(T_1\right)\right]- \phi\left( \left[ T_1 , T_2 \right] \right)$. We compute $\left[ T_1, \phi\left( T_2 \right) \right]$. Then from $X_1^3\frac{\partial}{\partial X_1} + X_1^2 X_2 \frac{\partial}{\partial X_2}+ X_1^2 X_3\frac{\partial}{\partial X_3}\equiv X_2^2\frac{\partial}{\partial X_1}+ \frac{\partial}{\partial X_2}$ and $\frac{\partial}{\partial X_3}\equiv 0$, we have
{\Small{\begin{align*}
\left[ T_1, \phi\left( T_2 \right) \right]&= \left[ \left(X_2^2- X_1^3\right)\frac{\partial}{\partial X_1} - \left(X_1^2X_2 -1\right)\frac{\partial}{\partial X_2} - X_1^2 X_3 \frac{\partial}{\partial X_3}, \phi\left( T_2 \right) \right]\\
 &= \left( - b_2 X_2^2+ 3b_1 X_1^2 + 2b_2 X_1^3 - 2b_3 X_2\right)  \frac{\partial}{\partial X_1} + \left( b_2 + b_3 X_1^2+ 2b_1 X_1X_2 + 2b_2 X_1^2 X_2\right)\frac{\partial}{\partial X_2} + \left(2b_1 X_1X_3+ 2b_2X_1^2 X_3\right)\frac{\partial}{\partial X_3}\\
 &\equiv \left(b_2X_2^2+ 3b_1X_1^2- 2b_3X_2 \right)\frac{\partial}{\partial X_1} + \left(3b_2 + b_3X_1^2 + 2b_1X_1X_2\right)\frac{\partial}{\partial X_2}
\end{align*}}}
We compute $\left[ T_2 , \phi\left( T_1 \right) \right]$.
{\Small{\begin{align*}
\left[ T_2 , \phi\left( T_1 \right) \right]= \left[ \frac{\partial}{\partial X_3}, \phi\left( T_1 \right) \right]&=\left(t_2+2t_5X_3+t_6X_2+t_8X_1 +2t_{11}X_1X_3+t_{12}X_1X_2+t_{14}X_1^2\right)\frac{\partial}{\partial X_1}\\
&+\left(t_{17}+2t_{20}X_3+t_{21}X_2+t_{23}X_1 + 2t_{11}X_2X_3+t_{12}X_2^2+t_{14}X_1X_2\right)\frac{\partial}{\partial X_2}\\
&+ \left(3t_{11}X_3^2+2t_{12}X_2X_3+t_{13}X_2^2 + 2t_{14}X_1X_3+t_{15}X_1X_2 \right)\frac{\partial}{\partial X_3}\\
&\equiv \left(t_2+2t_5X_3+t_6X_2+t_8X_1 +2t_{11}X_1X_3+t_{12}X_1X_2+t_{14}X_1^2\right)\frac{\partial}{\partial X_1}\\
&+\left(t_{17}+2t_{20}X_3+t_{21}X_2+t_{23}X_1 + 2t_{11}X_2X_3+t_{12}X_2^2+t_{14}X_1X_2\right)\frac{\partial}{\partial X_2}
\end{align*}}}
We compute $\phi\left(\left[ T_1, T_2 \right] \right)$.
{\Small{\begin{align*}
\phi\left(\left[ T_1, T_2 \right] \right)=\phi\left(X_1^2 T_2 \right) =X_1^2\left(\left(  b_1 + b_2 X_1\right)\frac{\partial}{\partial X_1}  +\left( b_3+ b_2 X_2\right)\frac{\partial}{\partial X_2} + b_2 X_3 \frac{\partial}{\partial X_3} \right)\equiv  \left(b_1 X_1^2+ b_2X_2^2\right)\frac{\partial}{\partial X_1} + \left(b_3 X_1^2 + b_2\right)\frac{\partial}{\partial X_2}
\end{align*}}}
Then we have
\begin{align*}
&D_1(\phi)\left( T_1 \wedge T_2 \right)=\left[ T_1, \phi(T_2)\right] - \left[ T_2, \phi(T_1)\right]- \phi([T_1,T_2])\\
&\equiv \left( -t_2-2t_5X_3-\left(2b_3+t_6\right)X_2-t_8X_1 -2t_{11}X_1X_3-t_{12}X_1X_2+ \left(2b_1-t_{14}\right)X_1^2\right)\frac{\partial}{\partial X_1} \\
&+ \left(2b_2- t_{17}-2t_{20}X_3-t_{21}X_2-t_{23}X_1 - 2t_{11}X_2X_3-t_{12}X_2^2+ \left(2b_1-t_{14}\right)X_1X_2 \right)\frac{\partial}{\partial X_2}
\end{align*}
Then from $(\ref{j33})$ we see that 
\begin{align} \label{j31}
\dim_\mathbb{C} \textnormal{ker}\left( H^0\left( X, \mathscr{H}om_{\mathcal{O}_{\mathbb{P}_\mathbb{C}^3}}\left( \mathcal{O}_{\mathbb{P}_\mathbb{C}^3}(-1) \bigoplus \mathcal{O}_{\mathbb{P}_\mathbb{C}^3}(1) ,  \frac{\Theta_{\mathbb{P}_\mathbb{C}^3}}{\Theta_{\mathcal{F}_0}} \right)  \right)  \xrightarrow{D_1} H^0\left( \mathbb{P}_\mathbb{C}^3,  \mathscr{H}om_{\mathcal{O}_{\mathbb{P}_\mathbb{C}^3}}\left( \mathcal{O}_{\mathbb{P}_\mathbb{C}^3}  ,\frac{\Theta_{\mathbb{P}_\mathbb{C}^3}}{\Theta_{\mathcal{F}_0}} \right)\right)  \right) = 17
\end{align}
and any element $\phi\in \textnormal{ker}\left(H^0\left( X, \mathscr{H}om_{\mathcal{O}_{\mathbb{P}_\mathbb{C}^3}}\left( \mathcal{O}_{\mathbb{P}_\mathbb{C}^3}(-1) \bigoplus \mathcal{O}_{\mathbb{P}_\mathbb{C}^3}(1) ,  \frac{\Theta_{\mathbb{P}_\mathbb{C}^3}}{\Theta_{\mathcal{F}_0}} \right)  \right)  \xrightarrow{D_1} H^0\left( \mathbb{P}_\mathbb{C}^3,  \mathscr{H}om_{\mathcal{O}_{\mathbb{P}_\mathbb{C}^3}}\left( \mathcal{O}_{\mathbb{P}_\mathbb{C}^3}  ,\frac{\Theta_{\mathbb{P}_\mathbb{C}^3}}{\Theta_{\mathcal{F}_0}} \right)\right) \right)$ is represented in terms of $(\ref{j11})$ and $(\ref{j27})$ by
{\small{\begin{align}\label{j26}
\phi(T_1)=&\left(A_4+A_5X_2+A_6X_1-2A_1X_2X_3+A_7X_2^2+ A_8X_1X_2 + A_9X_1^2+ A_{10}X_1X_2^2+ 2A_2X_1^2X_3+ A_{11}X_1^2X_2\right)\frac{\partial}{\partial X_1}\\
+&\left(A_{12}+ 2A_3X_3+ A_{13}X_2 + A_{14}X_1 + A_{15}X_2^2 +A_{16}X_1X_2+A_{17}X_1^2+A_{10}X_2^3 +2A_2X_1X_2X_3+A_{11}X_1X_2^2\right)\frac{\partial}{\partial X_2} \notag \\
+&\left(A_{10}X_2^2X_3 +2A_2X_1X_3^2+A_{11}X_1X_2X_3\right)\frac{\partial}{\partial X_3}, \notag \\
\phi(T_2)=& \left(A_2 + A_3 X_1\right) \frac{\partial}{\partial X_1} + \left(A_1+ A_3X_2\right)\frac{\partial}{\partial X_2} + A_3X_3\frac{\partial}{\partial X_3},\,\,\,\,\,\,\,A_1,...,A_{17}\in \mathbb{C} \notag
\end{align}}}

We also note that
{\small{\begin{align*}
\dim_\mathbb{C} \mathbb{H}^2\left( \mathbb{P}_\mathbb{C}^3, \Theta_{\mathcal{F}_0}^\bullet \right) & = \dim_\mathbb{C} \textnormal{coker}\left( H^0\left( X, \mathscr{H}om_{\mathcal{O}_{\mathbb{P}_\mathbb{C}^3}}\left( \mathcal{O}_{\mathbb{P}_\mathbb{C}^3}(-1) \bigoplus \mathcal{O}_{\mathbb{P}_\mathbb{C}^3}(1) ,  \frac{\Theta_{\mathbb{P}_\mathbb{C}^3}}{\Theta_{\mathcal{F}_0}} \right)  \right)  \xrightarrow{D_1} H^0\left( \mathbb{P}_\mathbb{C}^3,  \mathscr{H}om_{\mathcal{O}_{\mathbb{P}_\mathbb{C}^3}}\left( \mathcal{O}_{\mathbb{P}_\mathbb{C}^3}  ,\frac{\Theta_{\mathbb{P}_\mathbb{C}^3}}{\Theta_{\mathcal{F}_0}} \right)\right) \right)\\
&= 11 - (28- 17 ) = 0
\end{align*}}}
In particular, $\left( \mathbb{P}_\mathbb{C}^3, \Theta_{\mathcal{F}_0} \right)$ is unobstructed in foliated deformations in terms of tangent sheaves. We note that from $(\ref{j25})$, we have $\mathbb{H}^2\left( \mathbb{P}_\mathbb{C}^3 , \mathcal{N}_{\mathcal{F}_0}^{*\bullet} \right) \ncong \mathbb{H}^2\left( \mathbb{P}_\mathbb{C}^3, \Theta_{\mathcal{F}_0}^\bullet  \right)$.

On the other hand, we note that from $(\ref{j30})$ and $(\ref{j31})$
\begin{align*}
\dim_\mathbb{C} \mathbb{H}^1\left( \mathbb{P}_\mathbb{C}^3, \Theta_{\mathcal{F}_0}^\bullet  \right)&=\dim_\mathbb{C} \frac{\textnormal{ker}\left( H^0\left( X, \mathscr{H}om_{\mathcal{O}_{\mathbb{P}_\mathbb{C}^3}}\left( \mathcal{O}_{\mathbb{P}_\mathbb{C}^3}(-1) \bigoplus \mathcal{O}_{\mathbb{P}_\mathbb{C}^3}(1) ,  \frac{\Theta_{\mathbb{P}_\mathbb{C}^3}}{\Theta_{\mathcal{F}_0}} \right)  \right)  \xrightarrow{D_1} H^0\left( \mathbb{P}_\mathbb{C}^3,  \mathscr{H}om_{\mathcal{O}_{\mathbb{P}_\mathbb{C}^3}}\left( \mathcal{O}_{\mathbb{P}_\mathbb{C}^3}  ,\frac{\Theta_{\mathbb{P}_\mathbb{C}^3}}{\Theta_{\mathcal{F}_0}} \right)\right) \right)}{\textnormal{im}\left( H^0\left( \mathbb{P}_\mathbb{C}^3, \Theta_{\mathbb{P}_\mathbb{C}^3} \right)\xrightarrow{D_0} H^0\left( X, \mathscr{H}om_{\mathcal{O}_{\mathbb{P}_\mathbb{C}^3}}\left( \mathcal{O}_{\mathbb{P}_\mathbb{C}^3}(-1) \bigoplus \mathcal{O}_{\mathbb{P}_\mathbb{C}^3}(1) ,  \frac{\Theta_{\mathbb{P}_\mathbb{C}^3}}{\Theta_{\mathcal{F}_0}} \right)  \right) \right)}\\
&= 17 - 11 =6
\end{align*}
By considering $(\ref{j29}),(\ref{j28})$ and $(\ref{j26})$, any element $\phi$ of $\mathbb{H}^1\left( \mathbb{P}_\mathbb{C}^3, \Theta_{\mathcal{F}_0}^\bullet \right)$ is represented by
\begin{align*}
&\phi\left(T_1\right)=\left( t_1 + t_2X_2 + t_3X_1 + t_4X_1X_2 +t_5X_1^2 \right) \frac{\partial}{\partial X_1} + t_6X_1 \frac{\partial}{\partial X_2}  \\
& \phi(T_2)=0,\,\,\,\,\,t_1,...,t_6\in \mathbb{C}
\end{align*}
Then we see that $\left( \mathbb{P}_\mathbb{C}^3 \times B, \Theta_{\tilde{\mathcal{F}}}, \pi, B   \right)$, where $\pi: \mathbb{P}_\mathbb{C}^3 \times B \to B$ is the natural projection with the coordinates $\left(t_1,...,t_6 \right)$ on $B \subset \mathbb{C}^6$ for a sufficiently small neighborhood $B$ of $0$ in $\mathbb{C}^6$ and $\Theta_{\tilde{\mathcal{F}}}$ on $U_0$ is generated by 
\begin{align*}
\tilde{T}_1 &= \left(X_2^2- X_1^3 + t_1 + t_2X_2 + t_3X_1 + t_4X_1X_2 +t_5X_1^2 \right) \frac{\partial}{\partial X_1} +\left(- X_1^2X_2+1+ t_6X_1 \right)\frac{\partial}{\partial X_2} - X_1^2 X_3 \frac{\partial}{\partial X_3} \\
\tilde{T}_2&= \frac{\partial}{\partial X_3}
\end{align*}
defines a foliated analytic family of deformations of $\pi^{-1}(0)=\left( \mathbb{P}_\mathbb{C}^3 , w_0 \right)$ in terms of tangent sheaves since $\left[ \tilde{T}_1, \tilde{T}_2\right]= X_1^2 \tilde{T}_2$, and the foliated Kodaira-Spencer map in terms of tangent sheaves is an isomorphism at $0$. Moreover we note that $\Theta_{\tilde{\mathcal{F}}}$ in the foliated analytic family is a family of foliations  which are linear pullback of some foliations on $\mathbb{P}_\mathbb{C}^2$.  Since the Kodaira-Spencer map in terms of tangent sheaves is an isomorphism at $t=0$, the above foliated analytic family is complete at $t=0$ by Theorem of completeness for deformations of foliated analytic structures in terms of tangent sheaves in Part I. Since any small foliated deformations of $\left( \mathbb{P}_\mathbb{C}^3, w_0 \right)$ in terms of tangent sheaves is induced from the complete foliated analytic family in terms of tangent sheaves, we see that  any small foliated deformations of $\left( \mathbb{P}_\mathbb{C}^3, w_0 \right)$ in terms of tangent sheaves are also linear-pullback of some foliations on $\mathbb{P}_\mathbb{C}^2$. We also note that $\Theta_{\tilde{\mathcal{F}}}$ generated by $\tilde{T}_1$ and $\tilde{T}_2$ as above is the linear-pullback of $\Theta_{\mathcal{F}_0}$ generated by $T$ in $(\ref{j32})$ from Example \ref{cp9}. Since the family $\left( \mathbb{P}_\mathbb{C}^2\times B, T , B, \pi \right)$ in Example \ref{cp9} is also complete at $t=0$, we can conclude that any small foliated deformations of $\left( \mathbb{P}_\mathbb{C}^3 , w_0 \right)$ in terms of tangent sheaves is induced by a linear pullback of some foliated deformations of $\mathbb{P}_\mathbb{C}^2$ with the Jouanolou foliation.

\subsection{Foliated deformations of $\mathbb{P}_\mathbb{C}^3$ with a rational component of degree $2$ in terms of cotangent sheaves} \label{rcc1} \

Let us consider the foliation degree $2$ in a rational component on $\mathbb{P}_\mathbb{C}^3$ defined by the homogeneous holomorphic form on $\mathbb{C}^4$ with the homogeneous coordinate $\left[ \xi_0: \xi_1: \xi_2: \xi_3\right]$ on $\mathbb{P}_\mathbb{C}^3$ as follows (see \cite{CLN96} Example 4 p.579): we take $f= -\xi_0$ and $g=\frac{1}{3}\left( \xi_0^3 + \xi_1^3 + \xi_2^3 + \xi_3^3 \right)$ and define
\begin{align}\label{rc3}
w_0&=3 gdf - f dg\\
&= 3\cdot \left( \frac{1}{3}\left( \xi_0^3 + \xi_1^3 + \xi_2^3 + \xi_3^3 \right) \right) d\left(- \xi_0 \right) + \xi_0 d\left( \frac{1}{3}\left( \xi_0^3 + \xi_1^3 + \xi_2^3 + \xi_3^3 \right) \right) \notag \\
&= -\left( \xi_1^3 + \xi_2^3 + \xi_3^3 \right)d \xi_0 +  \xi_0 \xi_1^2 d\xi_1 + \xi_0 \xi_2^2 d \xi_2 + \xi_0 \xi_3^2 d \xi_3 \notag
\end{align}
Then $w_0$ belongs to the rational component $R(1,3)$ (see \cite{CLN96} Example 4 p.579 for the notation $R(1,3)$). We keep the notations in Appendix \ref{cp7} and Appendix \ref{ja1}. Then we have on $U_0=\left\{ \xi_0 \ne 0 \right\}$
\begin{align*}
w_0&= X_1^2 d X_1+ X_2^2 dX_2 + X_3^2 dX_3\\
dw_0&= 0
\end{align*}

Let us compute $\mathbb{H}^1\left( \mathbb{P}_\mathbb{C}^3, \mathcal{E}_{\mathcal{O}_{\mathbb{P}_\mathbb{C}^3}(4)}^\bullet \right)$ by the complex
\begin{align*}
\mathcal{E}_{\mathcal{O}_{\mathbb{P}_\mathbb{C}^3}(4)}^\bullet  : H^0\left( \mathbb{P}_\mathbb{C}^3, \mathcal{E}_{\mathcal{O}_{\mathbb{P}_\mathbb{C}^3}(4)} \right) \xrightarrow{\tilde{E}_0} H^0\left( \mathbb{P}_\mathbb{C}^3, \Omega_{\mathbb{P}_\mathbb{C}^3}^1(4) \right) \xrightarrow{\tilde{E}_1} H^0\left( \mathbb{P}_\mathbb{C}^3, \mathcal{O}_{\mathbb{P}_\mathbb{C}^3} (4) \right) \to 0
\end{align*}
First we compute $\textnormal{ker}\left( H^0\left( \mathbb{P}_\mathbb{C}^3 , \Omega_{\mathbb{P}_\mathbb{C}^3}^1(4) \right) \xrightarrow{\tilde{E}_1} H^0\left( \mathbb{P}_\mathbb{C}^3, \mathcal{O}_{\mathbb{P}_\mathbb{C}^3}(4) \right) \right)$. We recall $\textnormal{(\ref{cm1})}$ and $\textnormal{(\ref{cm2})}$. Then the image of $\tilde{E}_1$ is given by $w_0 \wedge db + b \wedge dw_0= w_0\wedge db$. The coefficient of $dX_1\wedge dX_2 \wedge dX_3$ in $w_0 \wedge db$ is given by
{\Small{\begin{align*}
&\left(t_{22}-t_2 \right)X_1^2+\left(t_{37}-t_{23}\right)X_2^2+ \left(t_4-t_{38}\right)X_3^2+\left(t_{28}-t_8\right)X_1^3+ \left(2t_{26}-t_6\right)X_1^2X_2+\left(t_{25}-2t_5\right)X_1^2X_3 +\left(t_{43}-2t_{29}\right)X_1X_2^2 \\
&+\left(2t_{10}-t_{44}\right)X_1X_3^2+\left(t_{41}-t_{28}\right)X_2^3 + \left(2t_{40}-t_{27}\right)X_2^2X_3+\left(t_9-2t_{42}\right)X_2X_3^2 +\left(t_8-t_{41}\right)X_3^3 + \left(t_{34}-t_{17}\right)X_1^4+\left( 2t_{32}-t_{15}\right)X_1^3X_2  \\
&+ \left(t_{31}-2t_{14}\right)X_1^3X_3- \left(4t_{13}+ 4t_{35}\right)X_1^2X_2^2-4t_{12}X_1^2 X_2X_3 + \left(4t_{19}-4t_{11}\right)X_1^2X_3^2 -\left(t_{17}+3t_{34} \right)X_1X_2^3 -4t_{33}X_1X_2^2X_3\\
&+4t_{18}X_1X_2X_3^2 + \left(3t_{17}+t_{34}\right)X_1X_3^3 -\left(t_{15}+2t_{32}\right)X_2^4-\left(2t_{14}+3t_{31}\right)X_2^3X_3+\left(4t_{16}-4t_{30}\right)X_2^2X_3^2 +\left(3t_{15}+2t_{32}\right)X_2X_3^3 +\left(2t_{14}+t_{31}\right)X_3^4
\end{align*}}}
Then the kernel of $\tilde{E}_1$ is given by
{\Small{\begin{align}\label{rc1}
b=&\left(A_{20}+ A_7X_3+ A_{11}X_2+ A_{21}X_1+ A_9X_3^2+ A_2 X_2 X_3+ A_{13}X_2^2+ 2A_8X_1X_3+ 2A_{12}X_1X_2+ A_{22}X_1^2 \right) dX_1\\
+&\left(- A_6X_1^2X_2+ A_5 X_1^2X_3- A_{10}X_2^3- A_{10}X_3^3 \right)dX_1 \notag \\ 
+&\left( A_{14} + A_1 X_3+ A_{15}X_2+ A_{11}X_1+ A_4 X_3^2+ 2A_3 X_2X_3+ A_{16}X_2^2 + A_2 X_1X_3+ 2A_{13}X_1X_2+ A_{12}X_1^2 + A_6 X_3^3 \right)dX_2 \notag \\
+&\left( A_5 X_2^2X_3 + A_{10} X_1X_2^2+ A_6 X_1^3   \right)dX_2 \notag \\
+&\left(A_{17}+ A_{18}X_3+ A_1 X_2+ A_7 X_1+ A_{19}X_3^2+ 2A_4 X_2X_3+ A_3 X_2^2+ 2A_9X_1X_3+ A_2 X_1X_2+ A_8X_1^2 - A_6 X_2X_3^2- A_5 X_2^3 \right)dX_3 \notag \\
+&\left( A_{10}X_1X_3^2- A_5 X_1^3 \right)dX_3,\,\,\,\,\,\,\, A_1,..., A_{22}\in \mathbb{C} \notag
\end{align}}}
so that we have
\begin{align*}
\dim_\mathbb{C} \textnormal{ker}\left( H^0\left( \mathbb{P}_\mathbb{C}^3, \Omega_{\mathbb{P}_\mathbb{C}^3}^1 (4) \right)\xrightarrow{\tilde{E}_1} H^0\left( \mathbb{P}_\mathbb{C}^3, \mathcal{O}_{\mathbb{P}_\mathbb{C}^3}(4)\right)\right) = 22
\end{align*}
Next we compute $\textnormal{im}\left( H^0\left( \mathbb{P}_\mathbb{C}^3 , \mathcal{E}_{\mathcal{O}_{\mathbb{P}_\mathbb{C}^3}(4)}  \right) \xrightarrow{\tilde{E}_0} H^0\left(  \mathbb{P}_\mathbb{C}^3, \Omega_{\mathbb{P}_\mathbb{C}^3}^1(4) \right) \right)$. We recall $\textnormal{(\ref{cp8})}$ for $m=4$. Then the image of $\tilde{E}_0$ is given by $\mathcal{L}_{T_1\frac{\partial}{\partial X_1} + T_2 \frac{\partial}{\partial X_2} + T_3 \frac{\partial}{\partial X_3}}(w) + fw$.  We note that $\mathcal{L}_{T_1\frac{\partial}{\partial X_1}}\left( w \right)$ is given by
{\Small{\begin{align*}
\left(2a_1X_1+2a_2X_1X_3+2a_3X_1X_2+3a_4X_1^2+3a_5X_1^2X_3+3a_6X_1^2X_2+4a_7X_1^3 \right)dX_1 + \left( a_3 X_1^2+ a_6 X_1^3\right)dX_2  + \left(a_2X_1^2+ a_5 X_1^3 \right)dX_3 
\end{align*}}}
and $\mathcal{L}_{T_2\frac{\partial}{\partial X_2}}\left( w \right)$ is given by
{\Small{\begin{align*}
 \left( a_{11}X_2^2+ a_7 X_2^3 \right)dX_1+ \left(2a_8X_2+2a_9X_2X_3+3a_{10}X_2^2 +2a_{11}X_1X_2+3a_5X_2^2X_3+4a_6X_2^3+3a_7X_1X_2^2\right)dX_2 + \left( a_9 X_2^2 + a_5 X_2^3 \right) dX_3 
\end{align*}}}
and $\mathcal{L}_{T_3\frac{\partial}{\partial X_3}} \left( w \right)$ is given by
{\Small{\begin{align*}
\left(a_{15}X_3^2+ a_7 X_3^3\right)dX_1+ \left(a_{14}X_3^2+ a_6 X_3^3\right) dX_2 + \left(2a_{12}X_3+3a_{13}X_3^2+2a_{14}X_2X_3+2a_{15}X_1X_3+4a_5X_3^3+3a_6X_2X_3^2+3a_7X_1X_3^2\right)dX_3
\end{align*}}}
and $fw_0$ is given by
{\Small{\begin{align*}
\left(a_{16}X_1^2-4a_7X_1^3-4a_6X_1^2X_2-4a_5X_1^2X_3\right)dX_1 + \left( a_{16}X_2^2-4a_7X_1X_2^2-4a_6X_2^3-4a_5X_2^2X_3 \right)dX_2\\
 + \left(a_{16}X_3^2-4a_7X_1X_3^2-4a_6X_2X_3^2-4a_5X_3^3 \right)dX_3
\end{align*}}}
Then $\mathcal{L}_{T_1 \frac{\partial}{\partial X_1} + T_2 \frac{\partial}{\partial X_2} + T_3 \frac{\partial}{\partial X_3}}\left( w_0\right) + f w_0$ is given by
{\Small{\begin{align}\label{rc2}
&\left(2a_1X_1+2a_2X_1X_3+2a_3X_1X_2+\left(3a_4+ a_{16}\right)X_1^2-a_5X_1^2X_3-a_6X_1^2X_2+ a_{11}X_2^2+ a_7 X_2^3+a_{15}X_3^2+ a_7 X_3^3 \right)dX_1\\
&+\left(a_3 X_1^2+ a_6 X_1^3+ 2a_8X_2+2a_9X_2X_3+ \left(3a_{10}+ a_{16}\right)X_2^2 +2a_{11}X_1X_2- a_5X_2^2X_3-a_7X_1X_2^2 + a_{14}X_3^2+ a_6 X_3^3\right)dX_2 \notag \\
&+ \left( a_2X_1^2+ a_5 X_1^3 +  a_9 X_2^2 + a_5 X_2^3+ 2a_{12}X_3+ \left(3a_{13}  + a_{16}\right) X_3^2+2a_{14}X_2X_3+2a_{15}X_1X_3-a_6X_2X_3^2-a_7X_1X_3^2\right)dX_3 \notag
\end{align}}}
Then we see that $\dim_\mathbb{C} \mathbb{H}^0\left(  \mathbb{P}_\mathbb{C}^3, \mathcal{E}_{\mathcal{O}_{\mathbb{P}_\mathbb{C}^3}(4)}^\bullet  \right)=\dim_\mathbb{C} \textnormal{ker}\left( H^0\left( \mathbb{P}_\mathbb{C}^3, \mathcal{E}_{\mathcal{O}_{\mathbb{P}_\mathbb{C}^3}(4)}\right)\xrightarrow{\tilde{E}_0} H^0 \left( \mathbb{P}_\mathbb{C}^3, \Omega_{\mathbb{P}_\mathbb{C}^3}^1(4) \right) \right)= 1$, and so 
\begin{align*}
\dim_\mathbb{C} \mathbb{H}^1\left( \mathbb{P}_\mathbb{C}^3, \mathcal{E}_{\mathcal{O}_{\mathbb{P}_\mathbb{C}^3}(4)}^\bullet \right)=  22- 15= 7
\end{align*}
By considering $(\ref{rc1})$ and $(\ref{rc2})$, we see that
\begin{align*}
\mathbb{H}^1\left(\mathbb{P}_\mathbb{C}^3, \mathcal{E}_{\mathcal{O}_{\mathbb{P}_\mathbb{C}^3}(4)}^\bullet \right) \cong  \textnormal{Span}_\mathbb{C} &\left\langle X_3 dX_2 + X_2 dX_3, X_2X_3 dX_1 + X_1 X_3 dX_2 + X_1 X_2 dX_3 , X_3 dX_1 + X_1 dX_3, \right.\\
&\left.  X_2 dX_1 + X_1 dX_2, dX_2,dX_2, dX_1\right\rangle\\
=\textnormal{Span}_\mathbb{C}& \left\langle d\left(X_2X_3 \right), d\left( X_1X_2X_3\right), d\left( X_1 X_3 \right), d\left( X_1X_2 \right), dX_2, dX_3, dX_1 \right \rangle
\end{align*}

Then we see that $\left( \mathbb{P}_\mathbb{C}^3\times B  , w , \pi, B \right)$, where $\pi:\mathbb{P}_\mathbb{C}^3 \times B \to B$ is the natural projection with the coordinates $\left( t_1,..., t_7 \right)$ on $B \subset \mathbb{C}^7$ for a sufficiently small neighborhood $B$ of $0$ in $\mathbb{C}^7$ and
{\small{\begin{align*}
\omega&= X_1^2 dX_1 + X_2^2 dX_2 + X_3^2 d X_3+ t_1d\left(X_2X_3\right)+ t_2 d\left(X_1X_2X_3\right)+ t_3 d\left(X_1X_3\right) + t_4 d \left(X_1X_2\right)+ t_5 dX_2+ t_6 dX_3 + t_7 dX_1\\
 &= \left( X_1^2 + t_2 X_2X_3+ t_3 X_3 +t_4 X_2+ t_7\right) dX_1 + \left(X_2^2 + t_1X_3+ t_2 X_1 X_3+ t_4 X_1+ t_5 \right) d X_2 + \left(X_3^2 + t_1 X_2+ t_2 X_1X_2+ t_3 X_1+ t_6 \right) dX_3 
\end{align*}}}
defines a foliated analytic family of deformations of $\pi^{-1}(0)=\left( \mathbb{P}_\mathbb{C}^3, w_0  \right)$ in terms of cotangent sheaves since $dw=0$, and the foliated Kodaira-Spencer map in terms of cotangent sheaves is an isomorphism at $t=0$, so that $\left( \mathbb{P}_\mathbb{C}^3, w_0 \right)$ is unobstructed in foliated deformations in terms of cotangent sheaves. Moreover we note that $\omega$ is a family of  foliations in the rational component $R(1,3)$. In fact, $\omega$ is expressed in terms of homogenous coordinates $\left[\xi_0: \xi_1: \xi_2 : \xi_3 \right]$ on $\mathbb{C}^4$ as
\begin{align*}
\omega&=\left( - \xi_1^3 - \xi_2^3 - \xi_3^3 - 2t_1 \xi_0 \xi_2 \xi_3 - 3t_2 \xi_1 \xi_2 \xi_3- 2t_3 \xi_0 \xi_1 \xi_3 - 2 t_4 \xi_0 \xi_1 \xi_2 - t_7 \xi_0^2 \xi_1 - t_5 \xi_0^2 \xi_2 - t_6 \xi_0^2 \xi_3 \right) d\xi_0\\
&+\left(\xi_0 \xi_1^2+ t_1 \xi_0 \xi_2 \xi_3 + t_3 \xi_0^2 \xi_3+ t_4 \xi_0^2 \xi_2 + t_7 \xi_0^3 \right) d\xi_1 + \left( \xi_0\xi_2^2 + t_1 \xi_0^2 \xi_3+ t_2 \xi_0 \xi_1\xi_3 + t_4 \xi_0^2 \xi_1+ t_5 \xi_0^3 \right) d\xi_2\\
& + \left( \xi_0 \xi_3^2+ t_1 \xi_0^2 \xi_2 + t_2 \xi_0 \xi_1 \xi_2+ t_3 \xi_0^2 \xi_1+ t_6 \xi_0^3\right) d\xi_3
\end{align*}
If we take
\begin{align*}
f&= - \xi_0\\
g&=\frac{1}{3}\left( \xi_0^3 + \xi_1^3+ \xi_2^3  + \xi_3^3 + 3 t_1 \xi_0 \xi_2 \xi_3 + 3t_2 \xi_1 \xi_2 \xi_3+3 t_3 \xi_0 \xi_1 \xi_3 + 3t_4 \xi_0 \xi_1 \xi_2 + 3t_7 \xi_0^2 \xi_1+3 t_5 \xi_0^2 \xi_2+ 3t_6 \xi_0^2 \xi_3 \right)\\
\end{align*}
then we see that $\omega= 3g df - f dg$. Since the foliated Kodaira-Spencer map in terms of cotangent sheaves is an isomorphism at $t=0$,  the above foliated analytic family is complete at $t=0$ by Theorem of completeness for deformations of foliated analytic structures in terms of cotangent sheaves in Part I. Since any small foliated deformations of $\left( \mathbb{P}_\mathbb{C}^3, w_0 \right)$ in terms of cotangent sheaves is induced from the complete foliated analytic family in terms of cotangent sheaves, we can conclude that any small foliated deformations of $\left(\mathbb{P}_\mathbb{C}^3, w_0 \right)$ given by $(\ref{rc3})$ in $R(1,3)$ in terms of cotangent sheaves is also in the rational component $R(1,3)$.

\subsection{Foliated deformations of $\mathbb{P}_\mathbb{C}^3$ with a logarithmic component of degree $2$ in terms of cotangent sheaves}\label{lg1} \

Let us consider the foliation of degree $2$ in a logarithmic component on $\mathbb{P}_\mathbb{C}^3$ defined by the homogeneous holomorphic form on $\mathbb{C}^4$ with the homogeneous coordinate $\left[ \xi_0 : \xi_1 : \xi_2: \xi_3\right]$ on $\mathbb{P}_\mathbb{C}^3$ (see \cite{CLN96} Example 5 p.580)
\begin{align*}
w_0&=\xi_0\xi_1\xi_2\xi_3\left( - 3 \frac{d\xi_0}{\xi_0} + \frac{d\xi_1}{\xi_1} + \frac{d\xi_2}{\xi_2} + \frac{d\xi_3}{\xi_3} \right) \\
&=-3\xi_1\xi_2\xi_3 d\xi_0+\xi_0\xi_2\xi_3d\xi_1+\xi_0 \xi_1 \xi_3 d\xi_2+ \xi_0\xi_1\xi_2 d\xi_3
\end{align*}
We keep the notations in Appendix \ref{cp7} and Appendix \ref{ja1}. Then we have on $U_0=\left\{ \xi_0 \ne  0 \right\}$
\begin{align*}
w_0&= X_2 X_3 dX_1 + X_1 X_3 dX_2 + X_1 X_2 dX_3\\
dw_0&= 0
\end{align*}
Let us compute $\mathbb{H}^i\left( \mathbb{P}_\mathbb{C}^3, \mathcal{E}_{\mathcal{O}_{\mathbb{P}_\mathbb{C}^3}(4)}^\bullet \right)$ by the complex
\begin{align*}
\mathcal{E}_{\mathcal{O}_{\mathbb{P}_\mathbb{C}^3}(4)}^\bullet  : H^0\left( \mathbb{P}_\mathbb{C}^3, \mathcal{E}_{\mathcal{O}_{\mathbb{P}_\mathbb{C}^3}(4)} \right) \xrightarrow{\tilde{E}_0} H^0\left( \mathbb{P}_\mathbb{C}^3, \Omega_{\mathbb{P}_\mathbb{C}^3}^1(4) \right) \xrightarrow{\tilde{E}_1} H^0\left( \mathbb{P}_\mathbb{C}^3, \mathcal{O}_{\mathbb{P}_\mathbb{C}^3} (4) \right) \to 0
\end{align*}
First we compute $\textnormal{ker}\left( H^0\left( \mathbb{P}_\mathbb{C}^3 , \Omega_{\mathbb{P}_\mathbb{C}^3}^1(4) \right) \xrightarrow{\tilde{E}_1} H^0\left( \mathbb{P}_\mathbb{C}^3, \mathcal{O}_{\mathbb{P}_\mathbb{C}^3}(4) \right) \right)$. We recall $\textnormal{(\ref{cm1})}$ and $\textnormal{(\ref{cm2})}$. Then the image of $\tilde{E}_1$ is given by $w_0 \wedge db + b \wedge dw_0= w_0 \wedge db$. The coefficient of $dX_1\wedge dX_2 \wedge dX_3$ in $w_0 \wedge db$ is given by
{\Small{\begin{align}
&(t_4-t_{38})X_1X_2+ (t_{37}-t_{23})X_1X_3+(t_{22}-t_2)X_2X_3 +(2t_{10}-t_{44})X_1^2X_2 +(t_{43}-2t_{29})X_1^2X_3 +(t_9-2t_{42})X_1X_2^2 \label{l3} \\
&+(2t_{40}-t_{27})X_1X_3^2 +(2t_{26}-t_6)X_2^2X_3 +(t_{25}-2t_5)X_2X_3^2 +4t_{19}X_1^3X_2-4t_{35}X_1^3X_3 + 4t_{18}X_1^2X_2^2 + ( t_{17}-t_{34})X_1^2X_2X_3 \notag \\
&  -4t_{33}X_1^2X_3^2 +4t_{16}X_1X_2^3+(2t_{32}+t_{15})X_1X_2^2X_3 +(-t_{31}-2t_{14})X_1X_2X_3^2-4t_{30}X_1X_3^3 -4t_{13}X_2^3X_3  -4t_{12} X_2^2X_3^2 -4t_{11}X_2X_3^3 \notag
\end{align}}}
Then the kernel of $\tilde{E}_1$ is given by
{\Small{\begin{align}\label{l2}
b&=\left(A_{21}+A_2X_3+A_1X_2+A_{22}X_1+A_7X_3^2+A_{23}X_2 X_3+A_6X_2^2+2A_4X_1X_3+2A_5X_1X_2+A_{24}X_1^2 \right) dX_1\\
&+\left(-2A_{10}X_1X_2X_3+A_{11}X_2^2X_3+A_{12}X_2X_3^2 \right)dX_1 \notag \\ 
&+\left(A_{13}+A_3X_3+A_{14}X_2+A_1X_1+A_9 X_3^2+2A_8X_2X_3+A_{15}X_2^2 +A_{16}X_1X_3+2A_6X_1X_2+A_5X_1^2 \right)dX_2 \notag \\
&+\left( A_{12}X_1X_3^2-2A_{11}X_1X_2X_3+A_{10}X_1^2X_3   \right)dX_2 \notag \\
&+\left( A_{17}+ A_{18}X_3+ A_3X_2+ A_2 X_1+ A_{19}X_3^2+2A_9X_2X_3+ A_8X_2^2+2A_7X_1X_3+ A_{20}X_1X_2+ A_4X_1^2 \right)dX_3 \notag \\
&+\left(-2A_{12}X_1X_2X_3+ A_{11}X_1X_2^2+ A_{10}X_1^2X_2\right)dX_3,\,\,\,\,\,\,\,\,\, A_1,..., A_{24} \in \mathbb{C} \notag
\end{align}}}
so that we have
\begin{align*}
\dim_\mathbb{C} \textnormal{ker} \left( H^0\left( \mathbb{P}_\mathbb{C}^3, \Omega_{\mathbb{P}_\mathbb{C}^3}^1(4) \right) \xrightarrow{\tilde{E}_1} H^0\left( \mathbb{P}_\mathbb{C}^3, \mathcal{O}_{\mathbb{P}_\mathbb{C}^3}(4) \right) \right) = 24
\end{align*}
Next we compute $\textnormal{im}\left( H^0\left( \mathbb{P}_\mathbb{C}^3 , \mathcal{E}_{\mathcal{O}_{\mathbb{P}_\mathbb{C}^3}(4)}  \right) \xrightarrow{\tilde{E}_0} H^0\left(  \mathbb{P}_\mathbb{C}^3, \Omega_{\mathbb{P}_\mathbb{C}^3}^1(4) \right) \right)$. We recall $\textnormal{(\ref{cp8})}$ for $m=4$. Then the image of $\tilde{E}_0$ is given by $\mathcal{L}_{T_1\frac{\partial}{\partial X_1} + T_2 \frac{\partial}{\partial X_2} + T_3 \frac{\partial}{\partial X_3}}(w_0) + fw_0$.  We note that $\mathcal{L}_{T_1\frac{\partial}{\partial X_1}}\left( w_0 \right)$ is given by
{ \Small{\begin{align*}
  &\left(a_4X_2X_3+a_5X_2X_3^2+a_6X_2^2X_3+2a_7X_1X_2X_3 \right)dX_1\\
    &+\left(a_1X_3+a_2X_3^2+2a_3X_2X_3+a_4X_1X_3+a_5X_1X_3^2+2a_6X_1X_2X_3+a_7X_1^2X_3\right)dX_2\\
 &+\left(a_1X_2+2a_2X_2X_3+a_3X_2^2+a_4X_1X_2+2a_5X_1X_2X_3+a_6X_1X_2^2+a_7X_1^2X_2\right)dX_3
 \end{align*} }}
 and $\mathcal{L}_{T_2\frac{\partial}{\partial X_2}}\left(w_0\right)$ is given by
{\Small{ \begin{align*}
&\left(a_8X_3+a_9X_3^2+a_{10}X_2X_3 +2a_{11}X_1X_3+a_5X_2X_3^2+a_6X_2^2X_3+2a_7X_1X_2X_3\right)dX_1 \\
&+\left(a_{10}X_1X_3+a_5X_1X_3^2+2a_6X_1X_2X_3+a_7X_1^2X_3\right)dX_2\\
&+\left(a_8X_1+2a_9X_1X_3+a_{10}X_1X_2 +a_{11}X_1^2+2a_5X_1X_2X_3+a_6X_1X_2^2+a_7X_1^2X_2\right)dX_3
\end{align*}}}
and $\mathcal{L}_{T_3\frac{\partial}{\partial X_3}}\left( w_0 \right)$ is given by
{\Small{\begin{align*}
&\left(a_{12}X_2+a_{13}X_2X_3+a_{14}X_2^2+2a_{15}X_1X_2+a_5X_2X_3^2+a_6X_2^2X_3+2a_7X_1X_2X_3 \right)dX_1\\
&+\left(a_{12}X_1+a_{13}X_1X_3+2a_{14}X_1X_2+a_{15}X_1^2+a_5X_1X_3^2+2a_6X_1X_2X_3+a_7X_1^2X_3\right)dX_2\\
&+\left(a_{13}X_1X_2+2a_5X_1X_2X_3+a_6X_1X_2^2+a_7X_1^2X_2\right)dX_3\\
\end{align*}}}
and $fw_0$ is given by
{\Small{\begin{align*}
 &\left(a_{16}X_2X_3-4a_7X_1X_2X_3-4a_6X_2^2X_3-4a_5X_2X_3^2\right)dX_1\\
             &+ \left(a_{16}X_1X_3-4a_7X_1^2X_3-4a_6X_1X_2X_3-4a_5X_1X_3^2 \right)dX_2\\
             &+ \left(a_{16}X_1X_2-4a_7X_1^2X_2-4a_6X_1X_2^2-4a_5X_1X_2X_3 \right)dX_3
\end{align*}   }}          
Then $\mathcal{L}_{T_1\frac{\partial}{\partial X_1} + T_2 \frac{\partial}{\partial X_2} + T_3 \frac{\partial}{\partial X_3}}(w_0) + fw_0$ is given by
{\Small{\begin{align}\label{l1}
&\left(a_{12}X_2+a_8X_3+a_{14}X_2^2+a_9X_3^2+2a_{15}X_1X_2 +2a_{11}X_1X_3+\left(a_4+a_{10}+a_{13}+a_{16}\right)X_2X_3+2a_7X_1X_2X_3-a_6X_2^2X_3-a_5X_2X_3^2\right)dX_1  \\
&+\left(a_{12}X_1+a_1X_3+a_{15}X_1^2+a_2X_3^2+2a_{14}X_1X_2+\left(a_4+a_{10}+a_{13}+a_{16}\right)X_1X_3+2a_3X_2X_3-a_7X_1^2X_3 +2a_6X_1X_2X_3-a_5X_1X_3^2\right)dX_2 \notag \\
&+\left(a_8X_1+a_1X_2 +a_{11}X_1^2 +a_3X_2^2 + \left(a_4+a_{10}+a_{13}+a_{16}\right)X_1X_2+2a_9X_1X_3 +2a_2X_2X_3 -a_7X_1^2X_2 -a_6X_1X_2^2+2a_5X_1X_2X_3\right)dX_3 \notag
\end{align}}}
Then we see that $\dim_\mathbb{C} \mathbb{H}^0\left(  \mathbb{P}_\mathbb{C}^3, \mathcal{E}_{\mathcal{O}_{\mathbb{P}_\mathbb{C}^3}(4)}^\bullet  \right)=\dim_\mathbb{C} \textnormal{ker}\left( H^0\left( \mathbb{P}_\mathbb{C}^3, \mathcal{E}_{\mathcal{O}_{\mathbb{P}_\mathbb{C}^3}(4)}\right)\xrightarrow{\tilde{E}_0} H^0 \left( \mathbb{P}_\mathbb{C}^3, \Omega_{\mathbb{P}_\mathbb{C}^3}^1(4) \right) \right)= 3$, and so 
\begin{align*}
\dim_\mathbb{C} \mathbb{H}^1\left( \mathbb{P}_\mathbb{C}^3, \mathcal{E}_{\mathcal{O}_{\mathbb{P}_\mathbb{C}^3}(4)}^\bullet \right)= 24-13=11
\end{align*}
By considering $(\ref{l2})$ and $(\ref{l1})$, we see that
\begin{align*}
\mathbb{H}^1\left( \mathbb{P}_\mathbb{C}^3, \mathcal{E}_{\mathcal{O}_{\mathbb{P}_\mathbb{C}^3}(4)}^\bullet \right) \cong \textnormal{Span}_\mathbb{C}& \left\langle  dX_2,X_2dX_2,X_2^2dX_2,X_1X_3dX_2, dX_3, X_3dX_3,X_3^2dX_3,X_1X_2dX_3 \right. \\
&\left. dX_1,X_1dX_1, X_1^2dX_1\right\rangle
\end{align*}
Let us consider
\begin{align*}
b= dX_1 + X_1X_3 dX_2 \in \textnormal{ker}\left( H^0\left( \mathbb{P}_\mathbb{C}^3, \Omega_{\mathbb{P}_\mathbb{C}^3}^1(4) \right) \xrightarrow{\tilde{E}_1} H^0\left( \mathbb{P}_\mathbb{C}^3, \mathcal{O}_{\mathbb{P}_\mathbb{C}^3}(4) \right) \right)
\end{align*}
which is not in $\textnormal{im}\left( H^0\left( \mathbb{P}_\mathbb{C}^3 , \mathcal{E}_{\mathcal{O}_{\mathbb{P}_\mathbb{C}^3}(4)} \right) \xrightarrow{ \tilde{E}_0} H^0\left( \mathbb{P}_\mathbb{C}^3, \Omega_{\mathbb{P}_\mathbb{C}^3}^1(4) \right) \right)$. Then we have
\begin{align*}
b \wedge db= \left( dX_1 + X_1X_3 dX_2 \right) \wedge \left( X_3 dX_1\wedge dX_2 - X_1 dX_2\wedge dX_3 \right) = - X_1 dX_1\wedge dX_2 \wedge dX_3
\end{align*}
which is not in $\textnormal{im}\left( H^0\left( \mathbb{P}_\mathbb{C}^3, \Omega_{\mathbb{P}_\mathbb{C}^3}^1(4) \right) \xrightarrow{\tilde{E}_1} H^0\left( \mathbb{P}_\mathbb{C}^3, \mathcal{O}_{\mathbb{P}_\mathbb{C}^3}(4)  \right) \right)$ from $(\ref{l3})$. Hence by Lemma \ref{on23}, $\left( \mathbb{P}_\mathbb{C}^3, w_0 \right)$ is obstructed in foliated deformations in terms of cotangent sheaves.

\subsection{Foliated deformations of the exceptional component in terms of cotangent sheaves}\label{cp11}\

Let us consider the foliation of degree $2$ of the exceptional component on $\mathbb{P}_\mathbb{C}^3$ defined by the homogeneous holomorphic form on $\mathbb{C}^4$ with the homogeneous coordinate $\left[ \xi_0 : \xi_1 : \xi_2: \xi_3\right]$ on $\mathbb{P}_\mathbb{C}^3$ (see \cite{CLN96} Example 5 p.580)
\begin{align*}
w_0 =\xi_0\left(2\xi_2^2-3\xi_1\xi_3\right)d\xi_1+\xi_0\left(3\xi_3\xi_0-\xi_1\xi_2\right)d\xi_2+\xi_0\left(\xi_1^2-2\xi_2\xi_0\right)d\xi_3-\left(\xi_1\xi_2^2-2\xi_1^2 \xi_3+\xi_2\xi_3\xi_0\right)d\xi_0
\end{align*}
We keep the notation in Appendix \ref{cp7} and Appendix \ref{ja1}. Then we have on $U_0$
\begin{align*}
w_0&=\left(2X_2^2-3 X_1 X_3\right)dX_1+\left(3X_3-X_1X_2\right)dX_2+\left(X_1^2-2X_2\right)dX_3\\
dw_0&=-5dX_2\wedge dX_3-5X_1dX_3\wedge dX_1-5X_2dX_1\wedge dX_2
\end{align*}
We will show that $\left( \mathbb{P}_\mathbb{C}^3, w_0 \right)$ is rigid in foliated deformations in terms of cotangent sheaves. We have the dual leaf complex
\begin{align*}
\mathcal{E}_{\mathcal{O}_{\mathbb{P}_\mathbb{C}^3}(4)}^\bullet : H^0\left( \mathbb{P}_\mathbb{C}^3, \mathcal{E}_{\mathcal{O}_{\mathbb{P}_\mathbb{C}^3}(4)} \right) \xrightarrow{\tilde{E}_0} H^0\left( \mathbb{P}_\mathbb{C}^3, \Omega_{\mathbb{P}_\mathbb{C}^3}^1 (4) \right) \xrightarrow{\tilde{E}_1} H^0\left( \mathbb{P}_\mathbb{C}^3, \mathcal{O}_{\mathbb{P}_\mathbb{C}^3}(4) \right) \to 0
\end{align*}
First we compute $\textnormal{ker}\left(H^0\left( \mathbb{P}_\mathbb{C}^3, \Omega_{\mathbb{P}_\mathbb{C}^3}^1 (4) \right) \xrightarrow{\tilde{E}_1} H^0\left( \mathbb{P}_\mathbb{C}^3, \mathcal{O}_{\mathbb{P}_\mathbb{C}^3}(4) \right) \right)$. We recall $\textnormal{(\ref{cm1})}$ and $\textnormal{(\ref{cm2})}$. Then the image of $\tilde{E}_1$ is given by $w_0 \wedge d b+ b \wedge dw_0$. The coefficient of $dX_1 \wedge dX_2 \wedge dX_3$ in $dw_0 \wedge b$ is given by
{\Small{\begin{align*}
 &-5t_{36}-5\left(t_1+t_{39}\right)X_1-5\left(t_{20}+t_{38}\right)X_2 -5t_{37}X_3-5\left(t_4+t_{45}\right)X_1^2-5\left(t_{42}+t_{22} \right)X_2^2 -5 t_{40}X_3^2 -5 \left(t_{44}+t_3+t_{23} \right)X_1X_2\\
 & -5 \left(t_{43}+t_2 \right)X_1X_3 -5 \left(t_{41}+t_{21} \right)X_2X_3-5 t_{10}X_1^3-5\left(t_{26}-t_{16}\right)X_2^3+5t_{30}X_3^3-5\left(t_9 - t_{19}+t_{29}\right)X_1^2X_2-5\left(t_8 -t_{35}\right)X_1^2X_3\\
 & -5 \left(t_7+t_{28}-t_{18}\right)X_1X_2^2-5\left( t_{27}+t_6-t_{17}-t_{34}\right)X_1X_2X_3-5\left(t_5-t_{33}\right)X_1X_3^2-5\left(t_{25}-t_{15}-t_{32}\right)X_2^2X_3\\
 & -5 \left(t_{24}-t_{14}-t_{31}\right)X_2X_3^2-5t_{19}X_1^4 +5t_{13}X_2^4 -5\left(t_{18}+t_{35}\right)X_1^3X_2 -5t_{17}X_1^3X_3-5\left(t_{16}+t_{34}\right)X_1^2X_2^2-5\left(t_{15}+t_{33}\right)X_1^2X_2X_3 \\
 &-5t_{14}X_1^2X_3^2-5t_{32}X_1X_2^3 -5\left(t_{13}+t_{31}\right)X_1X_2^2X_3-5\left(t_{12}+t_{30}\right)X_1X_2X_3^2-5t_{11}X_1X_3^3+5t_{12}X_2^3X_3 +5 t_{11}X_2^2X_3^2
\end{align*}}}
and the coefficient of $dX_1\wedge dX_2 \wedge dX_3$ in $w_0 \wedge db$ is given by
{\Small{\begin{align*}
&\left(2t_{38}-2t_4\right)X_2+ \left(3t_{37}-3t_{23}\right)X_3+ \left(t_4-t_{38}\right)X_1^2+ \left(2t_{22}-2t_2-2t_9+4t_{42}\right)X_2^2 + \left(6t_{40}-3t_{27}\right)X_3^2+ \left(t_{23}-t_{37}-4t_{10}+2t_{44}\right)X_1X_2\\
&+\left(3t_2-3t_{22}+3t_{43}-6t_{29}\right)X_1X_3 + \left(5t_{41}-3t_{28}-2t_8 \right)X_2X_3+ \left(2t_{10}-t_{44} \right)X_1^3+ \left(4t_{26}-2t_6-8t_{16}\right)X_2^3 -12t_{30}X_3^3\\
&+ \left(2t_{29}-t_{43}-8t_{19}+t_9-2t_{42}\right)X_1^2X_2 + \left(4t_8-3t_{28}-12t_{35}-t_{41}\right)X_1^2X_3 + \left(3t_{28}-2t_8-8t_{18}-t_{41}\right)X_1X_2^2\\ &+\left(3t_6-6t_{26}-9t_{17}-11t_{34}+t_{27}-2t_{40}\right)X_1X_2X_3+\left(6t_5-3t_{25}-12t_{33}\right)X_1X_3^2  +\left(-9t_{15}- 10t_{32}+2t_{25}-4t_5\right)X_2^2X_3\\
&-\left(10t_{14}+11t_{31}\right)X_2X_3^2+4t_{19}X_1^4- 8t_{13}X_2^4+ \left(4t_{18}+4t_{35}\right)X_1^3X_2+ \left(6t_{17}-2t_{34}\right)X_1^3X_3 + \left(5t_{34}-t_{17}+4t_{16}\right)X_1^2X_2^2\\
&+\left(6t_{15}-4t_{32}+4t_{33}\right)X_1^2X_2X_3+\left(8t_{14}-2t_{31}\right)X_1^2X_3^2 + \left(6t_{32}-t_{15}\right)X_1X_2^3+ \left(12t_{13}+5t_{31}-2t_{14}\right)X_1X_2^2X_3\\
&+ \left(12t_{12}+4t_{30}\right)X_1X_2X_3^2 +12t_{11}X_1X_3^3-8t_{12} X_2^3X_3-8t_{11}X_2^2X_3^2
\end{align*}}}
Then we see that
\begin{align} \label{ec1}
\dim_\mathbb{C} \textnormal{ker}\left( H^0\left( \mathbb{P}_\mathbb{C}^3, \Omega_{\mathbb{P}_\mathbb{C}^3}^1 (4) \right) \to H^0\left( \mathbb{P}_\mathbb{C}^3, \mathcal{O}_{\mathbb{P}_\mathbb{C}^3}(4) \right) \right) = 14
\end{align}
and the kernel of $\tilde{E}_1$ is given by
{\Small{\begin{align*}
b=&\left(A_4X_3+\frac{2A_2-5A_9}{5}X_2-A_1X_1+3A_3X_3^2+A_{13}X_2 X_3-\frac{9A_{12}+6A_{11}}{20}X_2^2+A_{14}X_1X_3+\frac{A_4-3A_{10}}{3}X_1X_2-\frac{4}{5}A_2X_1^2 \right) dX_1\\
+&\left(A_7X_1^2X_3-A_7X_1X_2^2-3A_5X_1X_2X_3-3A_5X_2^3-5A_6X_2^2X_3\right)dX_1 \\ 
+&\left(A_1+A_{11}X_3+A_{10}X_2+A_2X_1-2A_6X_3^2+ \left(A_6-22A_5\right)X_2X_3+A_8X_2^2 +\frac{-23A_7-15A_8-3A_{13}}{4}X_1X_3-\frac{2A_{11}+3A_{12}}{8}X_1X_2\right)dX_2\\
+&\left(4A_6X_1X_2X_3+3A_5X_1X_2^2+2A_5X_1^2X_3+A_7X_1^2X_2   \right)dX_2\\
+&\left(A_9-\frac{2A_{13}+7A_7}{5}X_3+A_{12}X_2-\frac{2}{3}A_4X_1-\frac{7}{3}A_6X_2X_3-(2A_6+22A_5)X_2^2+A_3X_1X_3\right)dX_3\\
+&\left(\frac{17A_7+5A_8+A_{13}}{4}X_1X_2-\frac{2A_{14}+2A_{11}+3A_{12}}{6}X_1^2+A_6X_1X_2^2+A_5X_1^2X_2-A_7X_1^3 \right)dX_3 , \,\,\,\,\,\,\,\,\,\,A_1,...,A_{14} \in \mathbb{C}
\end{align*}}}
and we note that
{\small{\begin{align}\label{ec10}
\dim_\mathbb{C} \mathbb{H}^2\left( \mathbb{P}_\mathbb{C}^3, \mathcal{N}_{\mathcal{F}_0}^{*\bullet} \right)=\dim_\mathbb{C} \mathbb{H}^2\left( \mathbb{P}_\mathbb{C}^3, \mathcal{E}_{\mathcal{O}_{\mathbb{P}_\mathbb{C}^3}(4)}^\bullet \right) = \dim_\mathbb{C} \textnormal{coker}\left( H^0\left( \mathbb{P}_\mathbb{C}^3, \Omega_{\mathbb{P}_\mathbb{C}^3}^1 (4) \right) \to H^0\left( \mathbb{P}_\mathbb{C}^3, \mathcal{O}_{\mathbb{P}_\mathbb{C}^3}(4) \right) \right)= 35 - (45-14)= 4
\end{align}}}
Next we compute $\textnormal{im}\left( H^0\left( \mathbb{P}_\mathbb{C}^3 , \mathcal{E}_{\mathcal{O}_{\mathbb{P}_\mathbb{C}^3}(4)}  \right) \xrightarrow{\tilde{E}_0} H^0\left(  \mathbb{P}_\mathbb{C}^3, \Omega_{\mathbb{P}_\mathbb{C}^3}^1(4) \right) \right)$. We recall $\textnormal{(\ref{cp8})}$ for $m=4$. Then the image of $\tilde{E}_0$ is given by $\mathcal{L}_{T_1\frac{\partial}{\partial X_1} + T_2 \frac{\partial}{\partial X_2} + T_3 \frac{\partial}{\partial X_3}}(w_0) + fw_0$.  We note that $\mathcal{L}_{T_1\frac{\partial}{\partial X_1}}\left( w_0 \right)$ is given by
{\Small{\begin{align*}
&\left(-3a_1X_3+2a_4X_2^2-3a_2X_3^2-6a_4X_1X_3-3a_3X_2X_3+2a_6X_2^3-9a_7X_1^2X_3+4a_7X_1X_2^2 -6a_6X_1X_2X_3-6a_5X_1X_3^2+2a_5X_2^2X_3\right) dX_1\\
&+\left(-a_1X_2+a_3X_2^2-a_4X_1X_2-3a_3X_1X_3-a_2X_2X_3-a_7X_1^2X_2 -3a_6X_1^2X_3+a_6X_1X_2^2-a_5X_1X_2X_3 \right) dX_2 \\
&+\left( 2a_1X_1+2a_4X_1^2+2a_2X_2^2+2a_3X_1X_2-a_2 X_1X_3+2a_7X_1^3 -a_5X_1^2X_3 +2a_6X_1^2X_2+2a_5X_1X_2^2 \right) dX_3
\end{align*}}}
and $\mathcal{L}_{T_2\frac{\partial}{\partial X_2}}\left( w \right)$ is given by
{\Small{\begin{align*}
&\left( 4a_8X_2+3a_{11}X_3+4a_{10}X_2^2+3a_{11}X_1X_2+ \left(4a_9+3a_7 \right)X_2X_3 +4a_6X_2^3+3a_7X_1X_2^2+4a_5X_2^2X_3 \right) dX_1\\
&+ \left( -a_8X_1+3a_{10}X_3-a_{11}X_1^2+3a_5X_3^2-2a_{10}X_1X_2+ \left(3a_7-a_9 \right)X_1X_3+6a_6X_2X_3-2a_7X_1^2X_2-3a_6X_1X_2^2-2a_5X_1X_2X_3 \right) dX_2 \\
&+ \left( -2a_8-2a_{11}X_1-2a_{10}X_2+a_9X_3-2a_6X_2^2-\left(2a_7+a_9\right)X_1X_2+a_5X_2X_3-a_5X_1X_2^2 \right) dX_3
\end{align*}}}
and $\mathcal{L}_{T_3\frac{\partial}{\partial X_3}} \left( w \right)$ is given by
{\Small{\begin{align*}
&\left( -3a_{12}X_1-2a_{15}X_2-2a_{15}X_1^2-3a_{14}X_1X_2-3a_{13}X_1X_3-2a_7X_2X_3-3a_{5}X_1X_3^2-3a_6X_1X_2X_3-2a_7X_1^2X_3 \right) dX_1\\
&+ \left( 3a_{12}+3a_{15}X_1+a_{14}X_2+3a_{13}X_3+a_{14}X_1^2+3a_5X_3^2+3a_7X_1X_3+a_6X_2X_3+a_6X_1^2X_3 \right) dX_2\\
&+\left( a_{13}X_1^2+2a_5X_1^2X_3+a_6X_1^2X_2+a_7X_1^3 -2a_{13}X_2-4a_5X_2X_3-2a_6X_2^2-2a_7X_1X_2 \right) dX_3
\end{align*}}}
and $f \omega_0$ is given by
{\Small{\begin{align*}
&\left(2a_{16}X_2^2-8a_7X_1X_2^2-8a_6X_2^3-8a_5X_2^2X_3 -3a_{16}X_1X_3+12a_7X_1^2X_3+12a_6X_1X_2X_3+12a_5X_1X_3^2 \right)dX_1\\
&+ \left( 3a_{16}X_3-12a_7X_1X_3-12a_6X_2X_3-12a_5X_3^2 - a_{16}X_1X_2+4a_7X_1^2X_2+4a_6X_1X_2^2+4a_5X_1X_2X_3 \right)dX_2\\
&+ \left(a_{16}X_1^2-4a_7X_1^3-4a_6X_1^2X_2-4a_5X_1^2X_3 -2a_{16}X_2+8a_7X_1X_2+8a_6X_2^2+8a_5X_2X_3 \right)dX_3
\end{align*}}}
Then $\mathcal{L}_{T_1 \frac{\partial}{\partial X_1} + T_2 \frac{\partial}{\partial X_2} + T_3 \frac{\partial}{\partial X_3}}\left( w_0\right) + f w_0$ is given by
{\Small{\begin{align*}
&\left(-3a_{12}X_1+(4a_8-2a_{15})X_2+ \left(3a_{11}-3a_1\right)X_3-2a_{15}X_1^2+ \left(2a_4+4a_{10}+2a_{16}\right)X_2^2-3a_2X_3^2+ \left(3a_{11}-3a_{14}\right)X_1X_2\right)dX_1\\
&+\left(-\left(6a_4+3a_{13}+3a_{16}\right)X_1X_3+ \left(4a_9+a_7-3a_3\right)X_2X_3-2a_6X_2^3+3a_6X_1X_2X_3+a_7X_1^2X_3-a_7X_1X_2^2+3a_5X_1X_3^2-2a_5X_2^2X_3\right) dX_1 \\
&+\left(3a_{12}+ \left(3a_{15}-a_8\right)X_1+ \left(a_{14}-a_1\right)X_2+ \left(3a_{10}+3a_{13}+3a_{16}\right)X_3+ \left(a_{14}-a_{11}\right)X_1^2+a_3X_2^2-6a_5X_3^2- \left(a_4+2a_{10}+a_{16}\right)X_1X_2\right) dX_2\\
&+\left(- \left(6a_7+a_9+3a_3\right)X_1X_3- \left(a_2+5a_6\right)X_2X_3+a_7X_1^2X_2 -2a_6 X_1^2X_3+2a_6X_1X_2^2+a_5X_1X_2X_3 \right) dX_2\\
& + \left(-2a_8+ \left(2a_1-2a_{11}\right)X_1-\left(2a_{10}+2a_{13}+2a_{16}\right)X_2+a_9X_3+ \left(2a_4+a_{13}+a_{16}\right)X_1^2+ \left(2a_2+4a_6 \right)X_2^2 \right) dX_3 \\
&+\left( \left(4a_7+2a_3-a_9 \right)X_1X_2-a_2X_1X_3+5a_5X_2X_3-a_7X_1^3-a_6X_1^2X_2-3a_5X_1^2X_3+a_5X_1X_2^2  \right) dX_3
\end{align*}}}
Then we see that $\dim_\mathbb{C} \mathbb{H}^0\left(  \mathbb{P}_\mathbb{C}^3, \mathcal{E}_{\mathcal{O}_{\mathbb{P}_\mathbb{C}^3}(4)}^\bullet  \right)=\dim_\mathbb{C} \textnormal{ker}\left( H^0\left( \mathbb{P}_\mathbb{C}^3, \mathcal{E}_{\mathcal{O}_{\mathbb{P}_\mathbb{C}^3}(4)}\right)\xrightarrow{\tilde{E}_0} H^0 \left( \mathbb{P}_\mathbb{C}^3, \Omega_{\mathbb{P}_\mathbb{C}^3}^1(4) \right) \right)= 2$, and so from $(\ref{ec1})$, we have
\begin{align*}
\dim_\mathbb{C} \mathbb{H}^1\left( \mathbb{P}_\mathbb{C}^3, \mathcal{E}_{\mathcal{O}_{\mathbb{P}_\mathbb{C}^3}(4)}^\bullet \right)= \dim_\mathbb{C} \frac{\textnormal{ker}\left(H^0\left( \mathbb{P}_\mathbb{C}^3, \Omega_{\mathbb{P}_\mathbb{C}^3}^1(4) \right) \to H^0\left( \mathbb{P}_\mathbb{C}^3 , \mathcal{O}_{\mathbb{P}_\mathbb{C}^3}(4)  \right) \right)}{\textnormal{im}\left( H^0\left( \mathbb{P}_\mathbb{C}^3, \mathcal{E}_{\mathcal{O}_{\mathbb{P}_\mathbb{C}^2}(4)} \right) \to H^0\left( \mathbb{P}_\mathbb{C}^3, \Omega_{\mathbb{P}_\mathbb{C}^3}^1(4) \right) \right)}=  14 - (16-2)=0
\end{align*}
Hence $\left( \mathbb{P}_\mathbb{C}^3, w_0 \right)$ is rigid in foliated deformations in terms of cotangent sheaves by Theorem \ref{rsf25}.

\subsection{Foliated deformations of the exceptional component in terms of tangent sheaves}\

We study deformations of foliation of degree $2$ of the exceptional component on $\mathbb{P}_\mathbb{C}^3$ in the previous \textnormal{subsection \ref{cp11}} in terms of tangent sheaves. In this case we have $\Theta_{\mathcal{F}_0}\cong \mathcal{O}_{\mathbb{P}_\mathbb{C}^3} \bigoplus \mathcal{O}_{\mathbb{P}_\mathbb{C}^3}$ and we have the exact sequence on $U_0$ (see \cite{CLN96} Example 5 p.580)
\begin{align}\label{cp13}
0\to \Theta_{\mathcal{F}_0} : = \mathcal{O}_{\mathbb{P}_{\mathbb{C}^3}}\bigoplus \mathcal{O}_{\mathbb{P}_{\mathbb{C}^3}} &\to \Theta_{\mathbb{P}_{\mathbb{C}}^3} \to \frac{\Theta_{\mathbb{P}_\mathbb{C}^3}}{ \Theta_{\mathcal{F}_0} } \to 0 \\
(1,0) &\mapsto T_1= X_1 \frac{\partial}{\partial X_1} + 2 X_2 \frac{\partial}{\partial X_2} + 3 X_3 \frac{\partial}{\partial X_3}  \notag \\
(0,1) &\mapsto T_2 = \frac{\partial}{\partial X_1} + X_1 \frac{\partial}{\partial X_2} + X_2 \frac{\partial}{\partial X_3} \notag
\end{align}
Then we have $\left[ T_1, T_2 \right]= - T_2$. Since we have the injection $\mathbb{H}^1\left( \mathbb{P}_\mathbb{C}^3, \Theta_{\mathcal{F}_0}^\bullet \right) \hookrightarrow \mathbb{H}^1\left( \mathbb{P}_\mathbb{C}^3, \mathcal{N}_{\mathcal{F}_0}^{* \bullet} \right)=0$ (see Part II), we see that $\mathbb{H}^1\left( \mathbb{P}_\mathbb{C}^3, \Theta_{\mathcal{F}_0}^\bullet  \right)=0$, so that $\left( \mathbb{P}_\mathbb{C}^3, \Theta_{\mathcal{F}_0} \right)$ is rigid in foliated deformation in terms of tangent sheaves by Theorem \ref{cp12}. However we directly show $\mathbb{H}^1\left( \mathbb{P}_\mathbb{C}^3, \Theta_{\mathcal{F}_0}^\bullet \right)=0$. We have the leaf complex
\begin{align*}
\Theta_{\mathcal{F}_0}^\bullet: \Theta_{\mathbb{P}_\mathbb{C}^3} \xrightarrow{\tilde{D}_0} \mathscr{H}om_{\mathcal{O}_{\mathbb{P}_\mathbb{C}^3}}\left( \mathcal{O}_{\mathbb{P}_{\mathbb{C}^3}}\bigoplus \mathcal{O}_{\mathbb{P}_{\mathbb{C}^3}} , \frac{\Theta_{\mathbb{P}_\mathbb{C}^3}}{\Theta_{\mathcal{F}_0}} \right) \xrightarrow{\tilde{D}_1} \mathscr{H}om_{\mathcal{O}_{ \mathbb{P}_\mathbb{C}^3}}\left( \bigwedge^2 \Theta_{\mathcal{F}_0}, \frac{\Theta_{\mathbb{P}_\mathbb{C}^3}}{\Theta_{\mathcal{F}_0}}  \right) \cong  \mathscr{H}om_{\mathcal{O}_{\mathbb{P}_\mathbb{C}^3}} \left(\mathcal{O}_{\mathbb{P}_\mathbb{C}^3}, \frac{\Theta_{\mathbb{P}_\mathbb{C}^3}}{\Theta_{\mathcal{F}_0}} \right) \to 0
\end{align*}
Equivalently we have
\begin{align*}
\Theta_{\mathcal{F}_0}^\bullet : \Theta_{\mathbb{P}_\mathbb{C}^3} \xrightarrow{D_0} \frac{\Theta_{\mathbb{P}_\mathbb{C}^3}}{\Theta_{\mathcal{F}_0}} \bigoplus \frac{\Theta_{\mathbb{P}_\mathbb{C}^3}}{\Theta_{\mathcal{F}_0}} \xrightarrow{D_1} \frac{\Theta_{\mathbb{P}_\mathbb{C}^3}}{\Theta_{\mathcal{F}_0}} \to 0
\end{align*}
From the exact sequence $(\ref{cp13})$, we see that $\dim_\mathbb{C} H^0\left( \mathbb{P}_\mathbb{C}^3, \frac{\Theta_{\mathbb{P}_\mathbb{C}^3}}{\Theta_{\mathcal{F}_0} } \right)= 13$ and $H^0\left( \mathbb{P}_\mathbb{C}^3, \frac{\Theta_{\mathbb{P}_\mathbb{C}^3}}{\Theta_{\mathcal{F}_0}} \right) \cong \frac{H^0\left(  \mathbb{P}_\mathbb{C}^3, \Theta_{\mathbb{P}_\mathbb{C}^3} \right)}{\langle T_1,T_2 \rangle }$ on $U_0$ and $H^i\left(  \mathbb{P}_\mathbb{C}^3, \frac{\Theta_{\mathbb{P}_\mathbb{C}^3}}{ \Theta_{\mathcal{F}_0} } \right)=0$ for $i\geq 1$. This implies that we can compute $\mathbb{H}^i\left( \mathbb{P}_\mathbb{C}^3, \Theta_{\mathcal{F}_0}^\bullet \right)$ by the following complex
\begin{align*}
\Theta_{\mathcal{F}_0}^\bullet : H^0\left( \mathbb{P}_\mathbb{C}^3, \Theta_{\mathbb{P}_\mathbb{C}^3} \right) \xrightarrow{D_0} H^0\left( \mathbb{P}_\mathbb{C}^3,  \mathscr{H}om_{\mathcal{O}_{\mathbb{P}_\mathbb{C}^3}}\left( \mathcal{O}_{\mathbb{P}_\mathbb{C}^3} \bigoplus \mathcal{O}_{\mathbb{P}_\mathbb{C}^3}  ,\frac{\Theta_{\mathbb{P}_\mathbb{C}^3}}{\Theta_{\mathcal{F}_0}} \right) \right) \xrightarrow{D_1} H^0\left( \mathbb{P}_\mathbb{C}^3, \mathscr{H}om_{\mathcal{O}_{\mathbb{P}_\mathbb{C}^3}}\left( \mathcal{O}_{\mathbb{P}_\mathbb{C}^3} , \frac{\Theta_{\mathbb{P}_\mathbb{C}^3}}{\Theta_{\mathcal{F}_0}} \right) \right)\to 0
\end{align*}
First we explicitly describe $H^0\left( \mathbb{P}_\mathbb{C}^3, \frac{\Theta_{\mathbb{P}_\mathbb{C}^3}}{\Theta_{\mathcal{F}_0}} \right)\cong \frac{H^0\left( \mathbb{P}_\mathbb{C}^3, \Theta_{\mathbb{P}_\mathbb{C}^3} \right)}{\langle T_1, T_2 \rangle}$ on $U_0$. We recall from $\textnormal{(\ref{ap2})}$
{\small{\begin{align}\label{ec2}
H^0\left( \mathbb{P}_\mathbb{C}^3,  \Theta_{\mathbb{P}_\mathbb{C}^3}  \right) = &\left\{  V= V_1 \frac{\partial}{\partial X_1} + V_2 \frac{\partial}{\partial X_2} + V_3 \frac{\partial}{\partial X_3} \right. \\
&\left.= \left( a_1 + a_2 X_3+ a_3 X_2 + a_4 X_1 + a_5 X_1X_3 + a_6 X_1 X_2 + a_7 X_1^2 \right) \frac{\partial}{\partial X_1}\right.  \notag \\
 &\left.+\left( a_8 + a_9 X_3 + a_{10}X_2 + a_{11}X_1 + a_5 X_2X_3 + a_6 X_2^2 + a_7 X_1X_2\right) \frac{\partial}{\partial X_2} \right. \notag \\ 
 &\left.  +\left( a_{12}+ a_{13} X_3 + a_{14} X_2 + a_{15}X_1 + a_5 X_3^2 + a_6 X_2 X_3  + a_7 X_1X_3 \right)  \frac{\partial}{\partial X_3}, a_1,..., a_{15}\in \mathbb{C} \right\} \notag
\end{align}}}
Since we have
 \begin{align}
 X_3\frac{\partial}{\partial X_3}&\equiv_{\Theta_\mathcal{F}}-\frac{1}{3}X_1\frac{\partial}{\partial X_1}-\frac{2}{3}X_2\frac{\partial}{\partial X_2} \in H^0\left( \mathbb{P}_\mathbb{C}^3, \frac{\Theta_{\mathbb{P}_\mathbb{C}^3}}{\Theta_{\mathcal{F}_0}} \right) \label{ec6} \\
 X_2\frac{\partial}{\partial X_3}&\equiv_{\Theta_\mathcal{F}}-\frac{\partial}{\partial X_1}-X_1\frac{\partial}{\partial X_2} \in  H^0\left( \mathbb{P}_\mathbb{C}^3, \frac{\Theta_{\mathbb{P}_\mathbb{C}^3}}{\Theta_{\mathcal{F}_0}} \right), \notag
 \end{align}
 we see that from $(\ref{ec2})$ we have
 {\small{\begin{align*}
 H^0\left(X, \frac{\Theta_{\mathbb{P}_\mathbb{C}^3}}{\Theta_{\mathcal{F}_0}}\right)=&\left\{\left(a_1+a_2X_3+a_3X_2+a_4X_1+a_5X_1X_3+a_6X_1X_2+a_7X_1^2\right)\frac{\partial}{\partial X_1}\right.\\
 &\left.+\left(a_8+a_9X_3+a_{10}X_2 +a_{11}X_1+a_5X_2X_3+a_6X_2^2+a_7X_1X_2 \right)\frac{\partial}{\partial X_2} \right.\\
 &\left.+\left(a_{12}+a_{13}X_1+a_5X_3^2+a_6X_2X_3+a_7X_1X_3 \right)\frac{\partial}{\partial X_3}, a_1,...,a_{13}\in \mathbb{C} \right\}
 \end{align*}}}
Then any element $\phi \in H^0\left( \mathbb{P}_\mathbb{C}^3,  \mathscr{H}om_{\mathcal{O}_{\mathbb{P}_\mathbb{C}^3}}\left( \mathcal{O}_{\mathbb{P}_\mathbb{C}^3} \bigoplus \mathcal{O}_{\mathbb{P}_\mathbb{C}^3}  ,\frac{\Theta_{\mathbb{P}_\mathbb{C}^3}}{\Theta_{\mathcal{F}_0}} \right) \right)$ is represented by
 {\small{\begin{align}\label{ec3}
 \phi \left(T_1 \right)=&\left\{\left(a_1+a_2X_3+a_3X_2+a_4X_1+a_5X_1X_3+a_6X_1X_2+a_7X_1^2\right)\frac{\partial}{\partial X_1}\right. \\
 &+\left(a_8+a_9X_3+a_{10}X_2 +a_{11}X_1+a_5X_2X_3+a_6X_2^2+a_7X_1X_2 \right)\frac{\partial}{\partial X_2} \notag \\
 &\left.+\left(a_{12}+a_{13}X_1+a_5X_3^2+a_6X_2X_3+a_7X_1X_3 \right)\frac{\partial}{\partial X_3}, a_1,...,a_{13}\in \mathbb{C}\right\} \notag
 \end{align}}}
 and
 {\small{\begin{align} \label{ec5}
 \phi\left(T_2\right)=& \left\{\left(b_1+b_2X_3+b_3X_2+b_4X_1+b_5X_1X_3+b_6X_1X_2+b_7X_1^2\right)\frac{\partial}{\partial X_1}\right.\\
 &+\left(b_8+b_9X_3+b_{10}X_2 +b_{11}X_1+b_5X_2X_3+b_6X_2^2+b_7X_1X_2 \right)\frac{\partial}{\partial X_2} \notag \\
 & \left. +\left(b_{12}+b_{13}X_1+b_5X_3^2+b_6X_2X_3+b_7X_1X_3 \right)\frac{\partial}{\partial X_3}, b_1,...,b_{13}\in \mathbb{C} \right\} \notag
 \end{align}}}
and any element $\sigma \in H^0\left( \mathbb{P}_\mathbb{C}^3, \mathscr{H}om_{\mathcal{O}_{\mathbb{P}_\mathbb{C}^3}}\left( \mathcal{O}_{\mathbb{P}_\mathbb{C}^3} , \frac{\Theta_{\mathbb{P}_\mathbb{C}^3}}{\Theta_{\mathcal{F}_0}} \right) \right)$ is represented by
{\small{\begin{align}\label{ec7}
 \sigma \left(T_1\wedge T_2 \right)=&\left\{\left(c_1+ c_2X_3+c_3X_2+c_4X_1+c_5X_1X_3+c_6X_1X_2+c_7X_1^2\right)\frac{\partial}{\partial X_1}\right.\\
 &+\left(c_8+c_9X_3+c_{10}X_2 +c_{11}X_1+c_5X_2X_3+c_6X_2^2+c_7X_1X_2 \right)\frac{\partial}{\partial X_2} \notag \\
 &\left.+\left(c_{12}+c_{13}X_1+c_5X_3^2+c_6X_2X_3+c_7X_1X_3 \right)\frac{\partial}{\partial X_3}, c_1,...,c_{13}\in \mathbb{C}\right\} \notag
\end{align}}}

We compute $H^0\left( \mathbb{P}_\mathbb{C}^3, \Theta_{\mathbb{P}_\mathbb{C}^3} \right)\xrightarrow{D_0} H^0\left( X, \mathscr{H}om_{\mathcal{O}_{\mathbb{P}_\mathbb{C}^3}}\left( \mathcal{O}_{\mathbb{P}_\mathbb{C}^3} \bigoplus \mathcal{O}_{\mathbb{P}_\mathbb{C}^3} ,  \frac{\Theta_{\mathbb{P}_\mathbb{C}^3}}{\Theta_{\mathcal{F}_0}} \right)  \right)$. From $(\ref{ec2})$ and $(\ref{ec6})$, the image $\phi=D_0\left(V_1 \frac{\partial}{\partial X_1} + V_2 \frac{\partial}{\partial X_2} + V_3 \frac{\partial}{\partial X_3} \right)$ is given by
{\Small{\begin{align*}
\phi\left( T_1 \right)&= \left[ V_1\frac{\partial}{\partial X_1} + V_2 \frac{\partial}{\partial X_2} + V_3 \frac{\partial}{\partial X_3},  X_1 \frac{\partial}{\partial X_1} + 2 X_2 \frac{\partial}{\partial X_2} + 3 X_3 \frac{\partial}{\partial X_3} \right] \\
 &=  \left( a_1-a_3X_2-2a_2X_3-a_7X_1^2-2a_6X_1X_2-3a_5X_1X_3\right)\frac{\partial}{\partial X_1}\\
                        &+\left(2a_8+ a_{11}X_1-a_9X_3-2a_6X_2^2-a_7X_1X_2-3a_5X_2X_3 \right)\frac{\partial}{\partial X_2}\\
                        &+\left(3a_{12}+2a_{13}X_1 + a_{14} X_2 -3a_5X_3^2-a_7X_1X_3-2a_6X_2X_3 \right)\frac{\partial}{\partial X_3} \\
 &\equiv  \left( a_1-a_{14}-a_3X_2-2a_2X_3-a_7X_1^2-2a_6X_1X_2-3a_5X_1X_3\right)\frac{\partial}{\partial X_1}\\
                        &+\left(2a_8+\left(a_{11}-a_{14}\right)X_1-a_9X_3-2a_6X_2^2-a_7X_1X_2-3a_5X_2X_3 \right)\frac{\partial}{\partial X_2}\\
                        &+\left(3a_{12}+2a_{13}X_1-3a_5X_3^2-a_7X_1X_3-2a_6X_2X_3 \right)\frac{\partial}{\partial X_3}
\end{align*}}}
and
{\Small{\begin{align*}
\phi\left( T_2 \right)&=\left[ V_1\frac{\partial}{\partial X_1} + V_2 \frac{\partial}{\partial X_2} + V_3 \frac{\partial}{\partial X_3}, \frac{\partial}{\partial X_1} + X_1 \frac{\partial}{\partial X_2} + X_2 \frac{\partial}{\partial X_3} \right]\\
&= \left(-a_4-\left(a_3+2a_7\right)X_1- \left(a_2+a_6\right)X_2-a_5X_3-a_6X_1^2-a_5X_1X_2 \right)\frac{\partial}{\partial X_1}\\
&+\left((a_1 - a_{11})+(a_4 - a_{10})X_1-\left(a_7+a_9-a_3\right)X_2 +a_2X_3-a_5X_2^2-a_6X_1X_2\right)\frac{\partial}{\partial X_2}\\
&+\left(\left(a_8- a_{13}\right)+\left(a_{11} - a_{14}\right)X_1+ \left(a_{10}- a_{15} \right)X_2+\left(a_9- a_7 \right) X_3-a_6X_1X_3-a_5X_2X_3 \right)\frac{\partial}{\partial X_3}\\
&\equiv\left(-a_{10}-a_4+a_{15}-\left(a_3+\frac{5}{3}a_7+\frac{1}{3}a_9\right)X_1-(a_2+a_6)X_2-a_5X_3-a_6X_1^2-a_5X_1X_2\right)\frac{\partial}{\partial X_1}\\
&+\left((a_1- a_{11})+(a_4-2a_{10}+a_{15})X_1-\left(\frac{1}{3}a_7+\frac{5}{3}a_9-a_3\right)X_2 +a_2X_3-a_5X_2^2-a_6X_1X_2\right)\frac{\partial}{\partial X_2}\\
&+\left(\left(a_8 - a_{13}\right)+\left(a_{11}- a_{14}\right)X_1-a_6X_1X_3-a_5X_2X_3 \right)\frac{\partial}{\partial X_3}
\end{align*}}}
Then we see that from $(\ref{ec3})$ and $(\ref{ec5})$
\begin{align}
&\dim_\mathbb{C} \mathbb{H}^0\left( \mathbb{P}_\mathbb{C}^3, \Theta_{\mathcal{F}_0}^\bullet  \right)=\dim_\mathbb{C} \textnormal{ker}\left( H^0\left( \mathbb{P}_\mathbb{C}^3, \Theta_{\mathbb{P}_\mathbb{C}^3} \right)\xrightarrow{D_0} H^0\left( X, \mathscr{H}om_{\mathcal{O}_{\mathbb{P}_\mathbb{C}^3}}\left( \mathcal{O}_{\mathbb{P}_\mathbb{C}^3} \bigoplus \mathcal{O}_{\mathbb{P}_\mathbb{C}^3} ,  \frac{\Theta_{\mathbb{P}_\mathbb{C}^3}}{\Theta_{\mathcal{F}_0}} \right)  \right) \right)=2 \notag \\
&\dim_\mathbb{C} \textnormal{im}\left( H^0\left( \mathbb{P}_\mathbb{C}^3, \Theta_{\mathbb{P}_\mathbb{C}^3} \right)\xrightarrow{D_0} H^0\left( X, \mathscr{H}om_{\mathcal{O}_{\mathbb{P}_\mathbb{C}^3}}\left( \mathcal{O}_{\mathbb{P}_\mathbb{C}^3} \bigoplus \mathcal{O}_{\mathbb{P}_\mathbb{C}^3} ,  \frac{\Theta_{\mathbb{P}_\mathbb{C}^3}}{\Theta_{\mathcal{F}_0}} \right)  \right) \right)= 15- 2 = 13 \label{ec8}
\end{align}
Next we compute $H^0\left( X, \mathscr{H}om_{\mathcal{O}_{\mathbb{P}_\mathbb{C}^3}}\left( \mathcal{O}_{\mathbb{P}_\mathbb{C}^3} \bigoplus \mathcal{O}_{\mathbb{P}_\mathbb{C}^3} ,  \frac{\Theta_{\mathbb{P}_\mathbb{C}^3}}{\Theta_{\mathcal{F}_0}} \right)  \right)  \xrightarrow{D_1} H^0\left( \mathbb{P}_\mathbb{C}^3,  \mathscr{H}om_{\mathcal{O}_{\mathbb{P}_\mathbb{C}^3}}\left( \mathcal{O}_{\mathbb{P}_\mathbb{C}^3}  ,\frac{\Theta_{\mathbb{P}_\mathbb{C}^3}}{\Theta_{\mathcal{F}_0}} \right)\right)$. From $(\ref{ec3})$, the image $D_1\left(\phi\right)\left( T_1 \wedge T_2 \right) $ is given by $\left[ T_1, \phi\left(T_2\right)\right] - \left[ T_2, \phi\left(T_1\right)\right]- \phi\left( \left[ T_1 , T_2 \right] \right)$. We compute $\left[ T_1, \phi\left( T_2 \right) \right]$.
{\small{\begin{align*}
\left[ T_1, \phi\left( T_2 \right) \right]& =\left[ X_1 \frac{\partial}{\partial X_1} + 2 X_2 \frac{\partial}{\partial X_2} + 3 X_3 \frac{\partial}{\partial X_3}, \phi\left( T_2 \right) \right] \\
&=\left( -b_1+b_3X_2+2b_2X_3+b_7X_1^2+2b_6X_1X_2+3b_5X_1X_3\right)\frac{\partial}{\partial X_1}\\
                        &+\left(-2b_8-b_{11}X_1+b_9X_3+2b_6X_2^2+b_7X_1X_2+3b_5X_2X_3 \right)\frac{\partial}{\partial X_2}\\
                        &+\left(-3b_{12}-2b_{13}X_1+3b_5X_3^2+b_7X_1X_3+2b_6X_2X_3\right)\frac{\partial}{\partial X_3}
\end{align*}}}
and we compute $\left[ T_2, \phi\left( T_1 \right) \right]$. Then from $(\ref{ec5})$ we have
{\small{\begin{align*}
\left[ T_2, \phi\left( T_1 \right) \right]&=\left(a_4+(a_3+ 2a_7)X_1+(a_2+a_6)X_2+a_5X_3+a_6X_1^2+a_5X_1X_2\right)\frac{\partial}{\partial X_1}\\
&+\left((a_{11}-a_1)+(a_{10}-a_4)X_1+\left(a_7+a_9-a_3\right)X_2 -a_2X_3+a_5X_2^2+a_6X_1X_2\right)\frac{\partial}{\partial X_2}\\
&+\left((a_{13}-a_8)-a_{11}X_1-a_{10}X_2 + \left( a_7 -a_9 \right) X_3+a_6X_1X_3+a_5X_2X_3\right)\frac{\partial}{\partial X_3}\\
&\equiv \left(a_4+ a_{10}+\left(a_3+ \frac{5}{3}a_7 + \frac{1}{3}a_9\right)X_1+(a_2+a_6)X_2+a_5X_3+a_6X_1^2+a_5X_1X_2\right)\frac{\partial}{\partial X_1}\\
&+\left((a_{11}-a_1)+(2a_{10}-a_4)X_1+\left( \frac{1}{3}a_7+ \frac{5}{3} a_9-a_3\right)X_2 -a_2X_3+a_5X_2^2+a_6X_1X_2\right)\frac{\partial}{\partial X_2}\\
&+\left((a_{13}-a_8)-a_{11}X_1+a_6X_1X_3+a_5X_2X_3\right)\frac{\partial}{\partial X_3}
\end{align*}}}
and we note that $\phi\left(\left[ T_1, T_2 \right] \right)=-\phi\left( T_2 \right)$. Then we have
{\Small{\begin{align*}
&\left[ T_1, \phi\left( T_2 \right)\right] - \left[ T_2, \phi\left( T_1 \right) \right] - \phi\left(\left[ T_1, T_2 \right] \right)\\
&\equiv \left( -a_4 - a_{10}+\left(b_4-a_3-\frac{5}{3}a_7  - \frac{1}{3}a_9\right)X_1  + \left(2b_3-a_2-a_6\right)X_2+ \left(3b_2-a_5\right)X_3+ \left(2b_7-a_6\right)X_1^2 + \left(3b_6-a_5\right)X_1X_2+4b_5X_1X_3\right)\frac{\partial}{\partial X_1}\\
&+ \left(-b_8+a_1-a_{11}+ \left(a_4-2a_{10}\right)X_1+\left(b_{10}+a_3-\frac{1}{3}a_7- \frac{5}{3} a_9 \right) X_2+ \left(2b_9+a_2\right)X_3+ \left(3b_6-a_5\right)X_2^2 + \left(2b_7-a_6\right)X_1X_2 +4b_5X_2X_3\right)\frac{\partial}{\partial X_2}\\
&+ \left( -2b_{12}+a_8-a_{13}+ \left(a_{11}-b_{13}\right)X_1+4b_5X_3^2+ \left(2b_7-a_6\right)X_1X_3+ \left(3b_6-a_5\right)X_2X_3 \right)\frac{\partial}{\partial X_3}
\end{align*}}}

Then from $(\ref{ec7})$ we see that
{\small{\begin{align}
&\dim_\mathbb{C} \textnormal{ker}\left( H^0\left( \mathbb{P}_\mathbb{C}^3, \mathscr{H}om_{\mathcal{O}_{\mathbb{P}_\mathbb{C}^3}}\left( \mathcal{O}_{\mathbb{P}_\mathbb{C}^3} \bigoplus \mathcal{O}_{\mathbb{P}_\mathbb{C}^3} , \frac{\Theta_{\mathbb{P}_\mathbb{C}^3}}{\Theta_{\mathcal{F}_0}} \right) \right) \to H^0\left(\mathbb{P}_\mathbb{C}^3, \mathscr{H}om_{\mathcal{O}_{\mathbb{P}_\mathbb{C}^3}}\left( \mathcal{O}_{\mathbb{P}_\mathbb{C}^3}, \frac{\Theta_{\mathbb{P}_\mathbb{C}^3}}{\Theta_{\mathcal{F}_0}} \right) \right) \right) = 13 \label{ec9} \\
&\dim_\mathbb{C} \mathbb{H}^2\left( \mathbb{P}_\mathbb{C}^3, \Theta_{\mathcal{F}_0}^\bullet \right)= \textnormal{coker}\left( H^0\left( \mathbb{P}_\mathbb{C}^3, \mathscr{H}om_{\mathcal{O}_{\mathbb{P}_\mathbb{C}^3}}\left( \mathcal{O}_{\mathbb{P}_\mathbb{C}^3} \bigoplus \mathcal{O}_{\mathbb{P}_\mathbb{C}^3} , \frac{\Theta_{\mathbb{P}_\mathbb{C}^3}}{\Theta_{\mathcal{F}_0}} \right) \right) \to H^0\left(\mathbb{P}_\mathbb{C}^3, \mathscr{H}om_{\mathcal{O}_{\mathbb{P}_\mathbb{C}^3}}\left( \mathcal{O}_{\mathbb{P}_\mathbb{C}^3}, \frac{\Theta_{\mathbb{P}_\mathbb{C}^3}}{\Theta_{\mathcal{F}_0}} \right) \right)\right) = 13 - (26-13)=0 \notag
\end{align}}}
We note that from $(\ref{ec10})$ we have $\mathbb{H}^2\left( \mathbb{P}_\mathbb{C}^3, \mathcal{N}_{\mathcal{F}_0}^{*\bullet} \right) \ncong \mathbb{H}^2\left( \mathbb{P}_\mathbb{C}^3, \Theta_{\mathcal{F}_0}^\bullet \right)$.

On the other hand, from $(\ref{ec8})$ and $(\ref{ec9})$ we have
{\small{\begin{align*}
\dim_\mathbb{C} \mathbb{H}^1\left( \mathbb{P}_\mathbb{C}^3, \Theta_{\mathcal{F}_0}^\bullet \right) = \dim_\mathbb{C} \frac{\textnormal{ker}\left(  H^0\left( \mathbb{P}_\mathbb{C}^3, \mathscr{H}om_{\mathcal{O}_{\mathbb{P}_\mathbb{C}^3}}\left( \mathcal{O}_{\mathbb{P}_\mathbb{C}^3} \bigoplus \mathcal{O}_{\mathbb{P}_\mathbb{C}^3} , \frac{\Theta_{\mathbb{P}_\mathbb{C}^3}}{\Theta_{\mathcal{F}_0}} \right) \right) \to H^0\left(\mathbb{P}_\mathbb{C}^3, \mathscr{H}om_{\mathcal{O}_{\mathbb{P}_\mathbb{C}^3}}\left( \mathcal{O}_{\mathbb{P}_\mathbb{C}^3}, \frac{\Theta_{\mathbb{P}_\mathbb{C}^3}}{\Theta_{\mathcal{F}_0}} \right) \right)\right)}{\textnormal{im}\left( H^0\left( \mathbb{P}_\mathbb{C}^3, \Theta_{\mathbb{P}_\mathbb{C}^3} \right)\xrightarrow{D_0} H^0\left( X, \mathscr{H}om_{\mathcal{O}_{\mathbb{P}_\mathbb{C}^3}}\left( \mathcal{O}_{\mathbb{P}_\mathbb{C}^3} \bigoplus \mathcal{O}_{\mathbb{P}_\mathbb{C}^3} ,  \frac{\Theta_{\mathbb{P}_\mathbb{C}^3}}{\Theta_{\mathcal{F}_0}} \right)  \right) \right)} = 13 - 13 =0
\end{align*}}}
Hence $\left( \mathbb{P}_\mathbb{C}^3, \Theta_{\mathcal{F}_0} \right)$ is rigid in foliated deformations in terms of tangent sheaves by Theorem \ref{cp12}.

\subsection{Foliated deformations of a Monge-Amp\`ere foliation on $\mathbb{P}_\mathbb{C}^3$ in terms of cotangent sheaves}\label{ma1} \

We study foliated deformations of $\mathbb{P}_\mathbb{C}^3$ with a foliation of degree $2$ defined by the meromorphic function $\varphi= \frac{\xi_2\xi_3}{\xi_0\xi_2- \xi_1^2}$ (see \cite{CLN96} Example 7 p.583) which corresponds to the homogeneous $1$-form of degree $3$ on $\mathbb{C}^4$
\begin{align*}
w_0= \xi_2^2\xi_3d\xi_0-2\xi_1\xi_2\xi_3d\xi_1+\xi_1^2\xi_3d\xi_2-\xi_2(\xi_0\xi_2-\xi_1^2)d\xi_3
\end{align*}
and show that $\left( \mathbb{P}_\mathbb{C}^3, w \right)$ is obstructed in foliated deformations in terms of cotangent sheaves. We keep the notations in Appendix \ref{cp7} and Appendix \ref{ja1}. Then the foliation on $U_0=\left\{ \xi_0\ne 0 \right\}$ is given by
\begin{align*}
w_0&= -2X_1X_2X_3dX_1+X_1^2X_3dX_2+ \left(X_1^2X_2-X_2^2 \right)dX_3\\
dw_0&= -2X_2dX_2\wedge dX_3-4X_1X_2dX_3\wedge dX_1+4X_1X_3dX_1\wedge dX_2
\end{align*}
Let us compute $\mathbb{H}^i\left( \mathbb{P}_\mathbb{C}^3, \mathcal{E}_{\mathcal{O}_{\mathbb{P}_\mathbb{C}^3}(4)}^\bullet \right)$ by the complex
\begin{align*}
\mathcal{E}_{\mathcal{O}_{\mathbb{P}_\mathbb{C}^3}(4)}^\bullet  : H^0\left( \mathbb{P}_\mathbb{C}^3, \mathcal{E}_{\mathcal{O}_{\mathbb{P}_\mathbb{C}^3}(4)} \right) \xrightarrow{\tilde{E}_0} H^0\left( \mathbb{P}_\mathbb{C}^3, \Omega_{\mathbb{P}_\mathbb{C}^3}^1(4) \right) \xrightarrow{\tilde{E}_1} H^0\left( \mathbb{P}_\mathbb{C}^3, \mathcal{O}_{\mathbb{P}_\mathbb{C}^3} (4) \right) \to 0
\end{align*}
First we compute $\textnormal{ker}\left( H^0\left( \mathbb{P}_\mathbb{C}^3 , \Omega_{\mathbb{P}_\mathbb{C}^3}^1(4) \right) \xrightarrow{\tilde{E}_1} H^0\left( \mathbb{P}_\mathbb{C}^3, \mathcal{O}_{\mathbb{P}_\mathbb{C}^3}(4) \right) \right)$. We recall $\textnormal{(\ref{cm1})}$ and $\textnormal{(\ref{cm2})}$. Then the image of $\tilde{E}_1$ is given by $w \wedge db + dw \wedge b$. The coefficient of $dX_1\wedge dX_2 \wedge dX_3$ in $dw \wedge b$ is given by
{\Small{\begin{align*}
&-2t_{36}X_2 -2t_{38}X_2^2 + \left(-2t_{39}- 4t_1\right) X_1X_2 +4 t_{20}X_1X_3 -2t_{37} X_2X_3 -2t_{42}X_2^3 +\left(-2t_{45} -4t_4 \right)X_1^2 X_2 +4t_{23}X_1^2X_3 \\
& +\left(-2t_{44} - 4t_3\right) X_1X_2^2 +\left( -2t_{43}-4t_2 + 4 t_{22}\right) X_1 X_2 X_3 +4t_{21}X_1X_3^2 -2t_{41}X_2^2 X_3 -2t_{40}X_2X_3^2 +2t_{16}X_2^4 -4t_{10}X_1^3X_2\\
&+4t_{29}X_1^3X_3 + \left( 2 t_{19} - 4 t_9 \right) X_1^2 X_2^2 +\left(2 t_{35}+ 4t_{28} - 4 t_8 \right) X_1^2 X_2 X_3 +4t_{27}X_1^2X_3^2 +\left( 2t_{18} -4 t_7 \right)X_1 X_2^3 +\left(2t_{17}+2t_{34} - 4 t_6+ 4 t_{26}\right)X_1X_2^2X_3 \\
&+\left( 2 t_{33}-4 t_5 + 4 t_{25}\right) X_1 X_2 X_3^2 +4t_{24}X_1X_3^3 +(2t_{15}+2t_{32})X_2^3X_3 +(2t_{14}+2t_{31})X_2^2X_3^2 +2t_{30}X_2X_3^3 -4t_{19}X_1^4X_2 +4t_{35}X_1^4X_3 \\
&-4t_{18}X_1^3X_2^2 +\left( - 4 t_{17}+ 4 t_{34} \right) X_1^3 X_2 X_3 +4t_{33}X_1^3X_3^2 -4t_{16}X_1^2X_2^3 +4t_{30}X_1^2X_3^3 +\left( 4 t_{32}- 4 t_{15}\right) X_1^2 X_2^2 X_3 + \left(4 t_{31}- 4 t_{!4} \right) X_1^2 X_2 X_3^2\\
&-8t_{13}X_1X_2^3X_3 - 8t_{12}X_1X_2^2X_3^2-8t_{11}X_1X_2X_3^3
\end{align*}}}
and the coefficient of $dX_1 \wedge dX_2 \wedge dX_3$ in $dw \wedge b$ is given by 
{\Small{\begin{align*}
&\left(-t_4+t_{38}\right)X_2^2 - \left(t_9-2t_{42} \right)X_2^3 + \left(t_4-t_{38}\right)X_1^2X_2 +(t_{37}-t_{23})X_1^2X_3 +(-2t_{10}+t_{44})X_1X_2^2 + (-2t_{22}+2t_2)X_1X_2X_3 -(t_8-t_{41})X_2^2X_3 \\
& -4t_{16}X_2^4 +(2t_{10}-t_{44})X_1^3X_2 +(t_{43}-2t_{29})X_1^3X_3  +(t_9-2t_{42}-4 t_{19})X_1^2X_2^2 + \left( -3 t_{28} + 3 t_8 \right) X_1^2 X_2 X_3 +(2t_{40}-t_{27})X_1^2X_3^2\\
&-4t_{18}X_1X_2^3 + \left( -4t_{26}+2t_6 - 3 t_{17}- t_{34} \right) X_1 X_2^2 X_3 +(-2t_{25}+4t_5)X_1X_2X_3^2 -(3t_{15}+2t_{32})X_2^3X_3  -(2t_{14}+t_{31})X_2^2X_3^2\\
&+4t_{19}X_1^4X_2  -4t_{35}X_1^4X_3 + 4t_{18}X_1^3X_2^2 +\left( - 4 t_{34} + 4 t_{17} \right) X_1^3 X_2 X_3 -4t_{33}X_1^3X_3^2 +4t_{16}X_1^2X_2^3 +\left( - 4 t_{32} + 4 t_{15} \right) X_1^2 X_2^2 X_3 \\
&+\left( - 4 t_{31} + 4 t_{14} \right) X_1^2 X_2 X_3^2 -4t_{30}X_1^2X_3^3 +8t_{13}X_1X_2^3X_3 +8t_{12} X_1X_2^2X_3^2 +8t_{11}X_1X_2X_3^3
\end{align*}}}
Then the coefficient of $dX_1\wedge dX_2 \wedge dX_3$ in $w \wedge db + dw \wedge b$ is given by
{\Small{\begin{align}\label{cm5}
&-2t_{36}X_2 -\left( t_{38} + t_4 \right)X_2^2 + \left(-2t_{39}- 4t_1\right) X_1X_2   +4 t_{20}X_1X_3  -2t_{37} X_2X_3 - t_9 X_2^3 + \left( - 2 t_{45} - 3 t_4 - t_{38} \right)X_1^2 X_2 + \left( t_{37} + 3 t_{23}\right) X_1^2 X_3 \\
&+ \left(-t_{44}-4t_3-2t_{10}\right)X_1X_2^2 + \left(-2t_{43}-2t_2 +2t_{22} \right)X_1X_2X_3  +4t_{21}X_1X_3^2  - \left(t_8+t_{41}\right) X_2^2X_3   -2t_{40}X_2X_3^2   -2t_{16}X_2^4 - \left(2t_{10}+t_{44}\right) X_1^3X_2 \notag \\
&+\left(t_{43}+2 t_{29}\right) X_1^3X_3 + \left(-2t_{19}-3t_9-2t_{42} \right)X_1^2X_2^2+\left( 2t_{35}+t_{28}-t_8\right)X_1^2X_2X_3 \notag \\
&+\left(2t_{40}+ 3t_{27}\right)X_1^2X_3^2 + \left( -2t_{18}-4t_7 \right)X_1X_2^3 + \left(-t_{17}+t_{34} -2t_6\right)X_1X_2^2X_3 \notag \\
&+ \left(2t_{33}+2t_{25} \right)X_1X_2X_3^2 + 4t_{24}X_1X_3^3- t_{15}X_2^3X_3 +t_{31}X_2^2X_3^2 + 2t_{30} X_2X_3^3  \notag
\end{align}}}
Then the kernel of $\tilde{E}_1$ is given by 
{\Small{\begin{align}\label{cm4}
b=&\left(-A_1X_2-2A_2X_1-A_6X_2 X_3-A_{12}X_2^2-2A_8X_1X_3-2A_4X_1X_2- A_1X_1^2 \right) dX_1\\
+&\left(-A_{12}X_1^2X_2-A_7X_1^2X_3+2A_5X_1X_2^2+A_3X_1X_3^2-\left(2A_{10}+2A_{11}\right)X_1X_2X_3-A_{19}X_2^2X_3-A_{17}X_2X_3^2 \right)dX_1\notag \\ 
+&\left(A_2+A_9X_3+A_1X_1+A_{13}X_3^2+A_{10}X_2X_3+A_5X_2^2 +A_6X_1X_3+A_4X_1^2 +A_{14}X_3^3+A_{15}X_2X_3^2 \right)dX_2 \notag \\
+&\left( A_{16}X_2^2X_3 +A_{17}X_1X_3^2+A_{11}X_1^2X_3-2A_5X_1^2X_2+A_{12}X_1^3   \right)dX_2 \notag \\
+&\left( \left(A_9-2A_8\right)X_2+A_3X_2X_3+A_{18}X_2^2+\left(A_6-2A_7\right)X_1X_2+A_8X_1^2 -A_{14}X_2X_3^2-A_{15}X_2^2X_3-A_{16}X_2^3 \right)dX_3 \notag \\
+&\left(A_{19}X_1X_2^2-A_3X_1^2X_3+ \left(2A_{10}+A_{11}\right)X_1^2X_2+A_7X_1^3 \right)dX_3 ,\,\,\,\,\,\,\,\,A_1,..., A_{19}\in \mathbb{C}, \notag
\end{align}}}
so that we have
\begin{align*}
\dim_\mathbb{C} \textnormal{ker} \left( H^0\left( \mathbb{P}_\mathbb{C}^3, \Omega_{\mathbb{P}_\mathbb{C}^3}^1(4) \right) \xrightarrow{\tilde{E}_1} H^0\left( \mathbb{P}_\mathbb{C}^3, \mathcal{O}_{\mathbb{P}_\mathbb{C}^3}(4) \right) \right) = 19
\end{align*}

Next we compute $\textnormal{im}\left( H^0\left( \mathbb{P}_\mathbb{C}^3 , \mathcal{E}_{\mathcal{O}_{\mathbb{P}_\mathbb{C}^3}(4)}  \right) \xrightarrow{\tilde{E}_0} H^0\left(  \mathbb{P}_\mathbb{C}^3, \Omega_{\mathbb{P}_\mathbb{C}^3}^1(4) \right) \right)$. We recall $\textnormal{(\ref{cp8})}$ for $m=4$. Then the image of $\tilde{E}_0$ is given by $\mathcal{L}_{T_1\frac{\partial}{\partial X_1} + T_2 \frac{\partial}{\partial X_2} + T_3 \frac{\partial}{\partial X_3}}(w) + fw$.  We note that $\mathcal{L}_{T_1\frac{\partial}{\partial X_1}}\left( w \right)$ is given by
{\Small{\begin{align*}
&\left(-2a_1X_2X_3-2a_2X_2X_3^2-2a_3X_2^2X_3-4a_4X_1X_2X_3-4a_5X_1X_2X_3^2-4a_6X_1X_2^2X_3-6a_7X_1^2X_2X_3 \right)dX_1\\
&+\left(2a_1X_1X_3+2a_2X_1X_3^2+2a_4X_1^2X_3+2a_5X_1^2X_3^2+2a_7X_1^3X_3\right)dX_2\\
&+\left(2a_1X_1X_2+2a_3X_1X_2^2+2a_4X_1^2X_2+2a_6X_1^2X_2^2+2a_7X_1^3X_2 \right)dX_3\
\end{align*}}}
and $\mathcal{L}_{T_2\frac{\partial}{\partial X_2}}\left(w\right)$ is given by
{\Small{\begin{align*}
&\left(-2a_8X_1X_3-2a_9X_1X_3^2-2a_{10}X_1X_2X_3 -a_{11}X_1^2X_3-2a_5X_1X_2X_3^2-2a_6X_1X_2^2X_3-a_7X_1^2X_2X_3 \right)dX_1\\
&+\left(a_{10}X_1^2X_3+a_5X_1^2X_3^2+2a_6X_1^2X_2X_3+a_7X_1^3X_3 \right)dX_2\\
&+ \left(a_8X_1^2+2a_9X_1^2X_3+a_{10}X_1^2X_2 +a_{11}X_1^3+2a_5X_1^2X_2X_3+a_6X_1^2X_2^2+a_7X_1^3X_2 \right)dX_3\\
&+ \left( -2a_8X_2-2a_9X_2X_3-2a_{10}X_2^2 -2a_{11}X_1X_2-2a_5X_2^2X_3-2a_6X_2^3-2a_7X_1X_2^2 \right)dX_3
\end{align*}}}
and $\mathcal{L}_{T_3\frac{\partial}{\partial X_3}}\left( w \right)$ is given by
{\Small{\begin{align*}
&\left(-2a_{12}X_1X_2-2a_{13}X_1X_2X_3-2a_{14}X_1X_2^2-a_{15}X_1^2X_2-2a_5X_1X_2X_3^2-2a_6X_1X_2^2X_3-a_7X_1^2X_2X_3 -a_{15}X_2^2-a_7X_2^2X_3 \right)dX_1\\
&+ \left(a_{12}X_1^2+a_{13}X_1^2X_3+2a_{14}X_1^2X_2+a_{15}X_1^3+a_5X_1^2X_3^2+2a_6X_1^2X_2X_3+a_7X_1^3X_3-a_{14}X_2^2-a_6X_2^2X_3 \right)dX_2\\
&+ \left(a_{13}X_1^2X_2+2a_5X_1^2X_2X_3+a_6X_1^2X_2^2+a_7X_1^3 X_2 -a_{13}X_2^2-2a_5X_2^2X_3-a_6X_2^3-a_7X_1X_2^2 \right) dX_3
\end{align*}}}
and $fw$ is given by
{\Small{\begin{align*}
&\left(-2a_{16}X_1X_2X_3+8a_7X_1^2X_2X_3+8a_6X_1X_2^2X_3+8a_5X_1X_2X_3^2\right)dX_1\\
&+\left(a_{16}X_1^2X_3-4a_7X_1^3X_3-4a_6X_1^2X_2X_3-4a_5X_1^2X_3^2 \right)dX_2\\
&+ \left( a_{16}X_1^2X_2-4a_7X_1^3X_2-4a_6X_1^2X_2^2-4a_5X_1^2X_2X_3 -a_{16}X_2^2+4a_7X_1X_2^2+4a_6X_2^3+4a_5X_2^2X_3 \right)dX_3
\end{align*}}}
Then $\mathcal{L}_{T_1\frac{\partial}{\partial X_1} + T_2 \frac{\partial}{\partial X_2} + T_3 \frac{\partial}{\partial X_3}}(w) + fw$ is given by
{\Small{\begin{align}\label{cm3}
&\left(-a_{15}X_2^2 -2a_{12}X_1X_2-2a_{8}X_1X_3 - 2a_1X_2X_3-a_{15}X_1^2X_2-a_{11}X_1^2X_3-2a_{14}X_1X_2^2 + \left(-4a_4-2a_{10}-2a_{13}-2a_{16}\right)X_1X_2X_3\right. \\
&\left.-2a_9X_1X_3^2+ \left(-2a_3-a_7\right)X_2^2X_3 -2a_2X_2X_3^2 \right)dX_1 \notag \\
&+ \left(a_{12}X_1^2 -a_{14}X_2^2+2a_1X_1X_3 +a_{15}X_1^3+2a_{14}X_1^2X_2+(2a_{4}+a_{10}+a_{13}+a_{16})X_1^2X_3 +2a_2X_1X_3^2 -a_6X_2^2X_3\right)dX_2 \notag \\
&+\left(-2a_8X_2 +a_8X_1^2+ \left(-a_{13}-a_{16}-2a_{10}\right)X_2^2 + \left(2a_1-2a_{11}\right)X_1X_2 -2a_9X_2X_3+a_{11}X_1^3 +a_6X_2^3 \right. \notag\\
&\left.+\left(2a_{4}+a_{10}+a_{13}+a_{16}\right)X_1^2X_2 +2a_9X_1^2X_3+\left(2a_3+a_7\right)X_1X_2^2\right)dX_3 \notag
\end{align}}}
Then we see that $\dim_\mathbb{C} \mathbb{H}^0\left(  \mathbb{P}_\mathbb{C}^3, \mathcal{E}_{\mathcal{O}_{\mathbb{P}_\mathbb{C}^3}(4)}^\bullet  \right)=\dim_\mathbb{C} \textnormal{ker}\left( H^0\left( \mathbb{P}_\mathbb{C}^3, \mathcal{E}_{\mathcal{O}_{\mathbb{P}_\mathbb{C}^3}(4)}\right)\xrightarrow{\tilde{E}_0} H^0 \left( \mathbb{P}_\mathbb{C}^3, \Omega_{\mathbb{P}_\mathbb{C}^3}^1(4) \right) \right)= 4$, and so 
\begin{align*}
\dim_\mathbb{C} \mathbb{H}^1\left( \mathbb{P}_\mathbb{C}^3, \mathcal{E}_{\mathcal{O}_{\mathbb{P}_\mathbb{C}^3}(4)}^\bullet \right)= 19- 12=7
\end{align*}
By considering $(\ref{cm4})$ with $(\ref{cm3})$, we see that 
\begin{align*}
\mathbb{H}^1\left( \mathbb{P}_\mathbb{C}^3, \mathcal{E}_{\mathcal{O}_{\mathbb{P}_\mathbb{C}^3}(4)}^\bullet \right) \cong \textnormal{Span}_\mathbb{C}& \left\langle  \left( - X_2- X_1^2\right) d X_1 + X_1 d X_2 , - 2X_1dX_1+ dX_2, X_3 dX_2 + X_2 dX_3,  \right. \\
& X_3^3 dX_2- X_2 X_3^2 dX_3, -2 X_1X_2X_3 dX_1+ X_2X_3 dX_2+ 2X_1^2 X_2 dX_3\\
&\left. X_3^2 dX_2, X_2X_3^2dX_2-X_2^2X_3dX_3 \right\rangle
\end{align*}

Let us consider 
\begin{align*}
b= \left( - X_2- X_1^2- 2 X_1 X_2 X_3 \right)d X_1 +\left( X_1 + X_2 X_3\right)dX_2+ 2 X_1^2 X_2 dX_3 \in \textnormal{ker}\left( H^0\left( \mathbb{P}_\mathbb{C}^3, \Omega_{\mathbb{P}_\mathbb{C}^3}^1(4) \right) \xrightarrow{\tilde{E}_1} H^0\left( \mathbb{P}_\mathbb{C}^3, \mathcal{O}_{\mathbb{P}_\mathbb{C}^3}(4) \right) \right)
\end{align*}
which is not in $\textnormal{im}\left( H^0\left( \mathbb{P}_\mathbb{C}^3 , \mathcal{E}_{\mathcal{O}_{\mathbb{P}_\mathbb{C}^3}(4)} \right) \xrightarrow{ \tilde{E}_0} H^0\left( \mathbb{P}_\mathbb{C}^3, \Omega_{\mathbb{P}_\mathbb{C}^3}^1(4) \right) \right)$. Then we have
\begin{align*}
b \wedge db= \left( - 2 X_1^4 + X_2^2 - 3 X_1^2 X_2 + 2 X_1^3 X_2 X_3 - 4 X_1 X_2^2 X_3\right) dX_1\wedge dX_2 \wedge dX_3
\end{align*}
which is not in $\textnormal{im}\left( H^0\left( \mathbb{P}_\mathbb{C}^3, \Omega_{\mathbb{P}_\mathbb{C}^3}^1(4) \right) \xrightarrow{\tilde{E}_1} H^0\left( \mathbb{P}_\mathbb{C}^3, \mathcal{O}_{\mathbb{P}_\mathbb{C}^3}(4)  \right) \right)$ from $(\ref{cm5})$. Hence by Lemma \ref{on23}, $\left( \mathbb{P}_\mathbb{C}^3, w \right)$ is obstructed in foliated deformations in terms of cotangent sheaves.

\subsection{Codimension $2$-Foliated deformations of the projective $3$-space defined by locally free sheaves}\

Let $\left( \mathbb{P}_\mathbb{C}^3, \mathcal{N}_{\mathcal{F}_0}^*\right)$ be a foliated projective $3$-space with $\mathcal{N}_{\mathcal{F}_0}^*$ locally free of rank $2$. Then the associated dimension $1$-foliation $\Theta_{\mathcal{F}_0}$ is of the form $\Theta_{\mathcal{F}_0}\cong  \mathcal{O}_{\mathbb{P}_\mathbb{C}^3}(-m), m\geq -1$. From Part IV-1, we have the dual leaf complex 
\begin{align*}
\mathcal{E}_{\mathcal{N}_{\mathcal{F}_0}^*}^\bullet: \mathcal{E}_{\mathcal{N}_{\mathcal{F}_0}^*} \to \mathscr{H}om_{\mathcal{O}_{\mathbb{P}_\mathbb{C}^3}}\left( \mathcal{N}_{\mathcal{F}_0}^* , \Omega_{\mathbb{P}_\mathbb{C}^3}^1 \right) \to 0
\end{align*}
which controls foliated deformations of $\left( \mathbb{P}_\mathbb{C}^3,  \mathcal{N}_{\mathcal{F}_0}^* \right)$ in terms of cotangent sheaves. On the other hand, we have the complex of sheaves
\begin{align*}
\mathcal{E}_{\mathcal{O}_{\mathbb{P}_\mathbb{C}^3}(m)}^\bullet : \mathcal{E}_{\mathcal{O}_{\mathbb{P}_\mathbb{C}^3}(m)} \to \Theta_{\mathbb{P}_\mathbb{C}^3}(m) \to 0
\end{align*}
which controls foliated deformations of $\left( \mathbb{P}_\mathbb{C}^3, \Theta_{\mathcal{F}_0} \right)$ in terms of tangent sheaves. We note that from Part II, we have the injection $\mathbb{H}^1\left( \mathbb{P}_\mathbb{C}^3, \mathcal{E}_{\mathcal{N}_{\mathcal{F}_0}^*}^\bullet \right) \hookrightarrow \mathbb{H}^1\left( \mathbb{P}_\mathbb{C}^3, \mathcal{E}_{\mathcal{O}_{\mathbb{P}_\mathbb{C}^3}(m)}^\bullet \right)$.

\section{Unfoldings of Projective $3$-space $\mathbb{P}_\mathbb{C}^3$}

\subsection{Codimension $1$-Unfoldings of the Projective $3$-space $\mathbb{P}_\mathbb{C}^3$} \

Let $\left(\mathbb{P}_\mathbb{C}^3,w\right)$ be a foliated projective $3$-space, where $w\in H^0\left(\mathbb{P}_\mathbb{C}^3, \Omega_{\mathbb{P}_\mathbb{C}^3}^1(m)\right),m\geq 2$ with $w\wedge dw=0$, so that we have the exact sequence defined by $w$
\begin{align*}
0\to \mathcal{O}_{\mathbb{P}_\mathbb{C}^3}(-m)\xrightarrow{w} \Omega_{\mathbb{P}_\mathbb{C}^3}^1 \to \frac{\Omega_{\mathbb{P}_\mathbb{C}^3}^1}{\mathcal{O}_{\mathbb{P}_\mathbb{C}^3}(-m)}\to 0
\end{align*}

Then we have a complex of sheaves from $(\ref{cp3})$
\begin{align*}
\tilde{\mathcal{E}}_{\mathcal{O}_{\mathbb{P}_\mathbb{C}^3}}(m) ^\bullet : \mathcal{E}_{\mathcal{O}_{\mathbb{P}_\mathbb{C}^3}}(m) \xrightarrow{\tilde{G}_0} \Omega_{\mathbb{P}_\mathbb{C}^3}^1(m)\bigoplus \mathcal{O}_{\mathbb{P}_\mathbb{C}^3}(m)\xrightarrow{\tilde{G}_1}  \mathcal{O}_{\mathbb{P}_\mathbb{C}^3}(2m-4)\bigoplus \tilde{\mathcal{S}}^1(m)  \xrightarrow{G_2} \mathcal{O}_{\mathbb{P}_\mathbb{C}^3}(2m-4)\to 0
\end{align*}
We note that $\tilde{\mathcal{S}}^1(m) \subset \mathcal{O}_{\mathbb{P}_\mathbb{C}^3}(2m)\otimes \bigwedge^2 \Omega_{\mathbb{P}_\mathbb{C}^3}^1 = \bigwedge^2 \Omega_{\mathbb{P}_\mathbb{C}^3}^1(2m)$. Then from Bott formula, we see that we can compute $\mathbb{H}^i\left( \mathbb{P}_\mathbb{C}^3, \tilde{\mathcal{E}}_{\mathcal{O}_{\mathbb{P}_\mathbb{C}^3}(m)}^\bullet \right)$ for $i=0,1,2$ by the following complex
\begin{align*}
\tilde{\mathcal{E}}_{\mathcal{O}_{\mathbb{P}_\mathbb{C}^3}(m)}^\bullet :  H^0\left( \mathbb{P}_\mathbb{C}^3, \tilde{\mathcal{E}}_{\mathcal{O}_{\mathbb{P}_\mathbb{C}^3}}(m) \right) &\xrightarrow{\tilde{G}_0}  H^0\left( \mathbb{P}_\mathbb{C}^3, \Omega_{\mathbb{P}_\mathbb{C}^3}^1(m)\right) \bigoplus H^0\left( \mathbb{P}_\mathbb{C}^3, \mathcal{O}_{\mathbb{P}_\mathbb{C}^3} (m) \right) \\
&\xrightarrow{\tilde{G}_1} H^0\left( \mathbb{P}_\mathbb{C}^3, \mathcal{O}_{\mathbb{P}_\mathbb{C}^3}(2m-4) \right) \bigoplus H^0\left( \mathbb{P}_\mathbb{C}^3, \tilde{\mathcal{S}}^1(m) \right) \xrightarrow{G_2} H^0\left(\mathbb{P}_\mathbb{C}^3, \mathcal{O}_{\mathbb{P}_\mathbb{C}^3}(2m-4) \right)\to 0
\end{align*}
Then we have
\begin{align*}
\mathbb{H}^1\left( \mathbb{P}_\mathbb{C}^3, \tilde{\mathcal{E}}_{\mathcal{O}_{\mathbb{P}_\mathbb{C}^3}(m)}^\bullet \right) & \cong \frac{\textnormal{ker}\left( H^0\left( \mathbb{P}_\mathbb{C}^3, \Omega_{\mathbb{P}_\mathbb{C}^3}^1(m)\right) \bigoplus H^0\left( \mathbb{P}_\mathbb{C}^3, \mathcal{O}_{\mathbb{P}_\mathbb{C}^3} (m) \right) \xrightarrow{\tilde{G}_1} H^0\left( \mathbb{P}_\mathbb{C}^3, \mathcal{O}_{\mathbb{P}_\mathbb{C}^3}(2m-4) \right) \bigoplus H^0\left( \mathbb{P}_\mathbb{C}^3, \bigwedge^2 \Omega_{\mathbb{P}_\mathbb{C}^3}^1(2m) \right) \right)}{\textnormal{im}\left(H^0\left( \mathbb{P}_\mathbb{C}^3, \mathcal{E}_{\mathcal{O}_{\mathbb{P}_\mathbb{C}^3}(m)}  \right) \xrightarrow{\tilde{G}_0} H^0\left( \mathbb{P}_\mathbb{C}^3, \Omega_{\mathbb{P}_\mathbb{C}^3}^1(m) \right) \bigoplus H^0\left( \mathbb{P}_\mathbb{C}^3, \mathcal{O}_{\mathbb{P}_\mathbb{C}^3}(m) \right) \right)} \\
\mathbb{H}^2\left( \mathbb{P}_\mathbb{C}^3, \tilde{\mathcal{E}}_{\mathcal{O}_{\mathbb{P}_\mathbb{C}^3}(m)}^\bullet  \right)& \cong \frac{\textnormal{ker}\left( H^0\left( \mathbb{P}_\mathbb{C}^3, \mathcal{O}_{\mathbb{P}_\mathbb{C}^3}(2m-4) \right) \bigoplus H^0\left( \mathbb{P}_\mathbb{C}^3, \tilde{\mathcal{S}}^1(m) \right) \xrightarrow{G_2} H^0\left(\mathbb{P}_\mathbb{C}^3, \mathcal{O}_{\mathbb{P}_\mathbb{C}^3}(2m-4) \right)\right)}{\textnormal{im}\left( H^0\left( \mathbb{P}_\mathbb{C}^3, \Omega_{\mathbb{P}_\mathbb{C}^3}^1(m)\right) \bigoplus H^0\left( \mathbb{P}_\mathbb{C}^3, \mathcal{O}_{\mathbb{P}_\mathbb{C}^3} (m) \right) \xrightarrow{\tilde{G}_1} H^0\left( \mathbb{P}_\mathbb{C}^3, \mathcal{O}_{\mathbb{P}_\mathbb{C}^3}(2m-4) \right) \bigoplus H^0\left( \mathbb{P}_\mathbb{C}^3, \bigwedge^2 \Omega_{\mathbb{P}_\mathbb{C}^3}^1(2m) \right) \right)}
\end{align*}

\subsection{Unfoldings of $\mathbb{P}_\mathbb{C}^3$ with linear-pullback of Jouanolou foliation of degree $2$ on $\mathbb{P}_\mathbb{C}^2$}\

Let $\left( \mathbb{P}_\mathbb{C}^3, w_0 \right)$ be the foliated Projective $3$-space, where $w_0\in H^0\left( \mathbb{P}_\mathbb{C}^3, \Omega_{\mathbb{P}_\mathbb{C}^3}^1(4) \right)$ is the linear pullabck of Jouanolou foliation of degree on $\mathbb{P}_\mathbb{C}^2$ in subsection \ref{j6}. We will show that $\left( \mathbb{P}_\mathbb{C}^3, w_0 \right)$ is rigid in unfoldings. We keep the notations in subsection \ref{j6} and Appendix \ref{cp7} and Appendix \ref{ja1}. We have on $U_0$\\
\begin{align*}
w_0&= \left( X_1^2 X_2 -1 \right) dX_1 + \left( X_2^2 - X_1^3 \right) dX_2\\
dw_0&= - 4 X_1^2 dX_1 \wedge dX_2
\end{align*}
Let us compute $\mathbb{H}^1\left( \mathbb{P}_\mathbb{C}^3, \tilde{\mathcal{E}}_{\mathcal{O}_{\mathbb{P}_\mathbb{C}^3}(4)}^\bullet \right)$ by the complex
{\small{\begin{align*}
H^0\left( \mathbb{P}_\mathbb{C}^3, \mathcal{E}_{\mathcal{O}_{\mathbb{P}_\mathbb{C}^3}(4)} \right) \xrightarrow{\tilde{G}_0} H^0\left( \mathbb{P}_\mathbb{C}^3, \Omega_{\mathbb{P}_\mathbb{C}^3}^1(4) \right) \bigoplus H^0\left( \mathbb{P}_\mathbb{C}^3, \mathcal{O}_{\mathbb{P}_\mathbb{C}^3}(4) \right) \xrightarrow{\tilde{G}_1} H^0\left( \mathbb{P}_\mathbb{C}^3 , \mathcal{O}_{\mathbb{P}_\mathbb{C}^3}(4)  \right)\bigoplus H^0\left( \mathbb{P}_\mathbb{C}^3, \tilde{S}^1(4) \right) \xrightarrow{G_2} H^0\left( \mathbb{P}_\mathbb{C}^3, \mathcal{O}_{\mathbb{P}_\mathbb{C}^3}(4) \right) \to 0
\end{align*}}}
Here we note that $H^0\left(\mathbb{P}_\mathbb{C}^3, \tilde{\mathcal{S}}^1(4) \right) \subset H^0\left( \mathbb{P}_\mathbb{C}^3, \bigwedge^2 \Omega_{\mathbb{P}_\mathbb{C}^3}^1(8) \right)$.

First we compute $\textnormal{ker}\left(H^0\left( \mathbb{P}_\mathbb{C}^3, \Omega_{\mathbb{P}_\mathbb{C}^3}^1(4) \right) \bigoplus H^0\left( \mathbb{P}_\mathbb{C}^3, \mathcal{O}_{\mathbb{P}_\mathbb{C}^3}(4)\right)\to H^0\left( \mathbb{P}_\mathbb{C}^3, \mathcal{O}_{\mathbb{P}_\mathbb{C}^3}(4) \right) \bigoplus H^0\left( \mathbb{P}_\mathbb{C}^3, \bigwedge^2 \Omega_{\mathbb{P}_\mathbb{C}^3}^1(8) \right)\right)$. We note that from $(\ref{j1})$, any element of $\textnormal{ker}\left( H^0\left( \mathbb{P}_\mathbb{C}^3, \Omega_{\mathbb{P}_\mathbb{C}^3}^1(4) \right)\to H^0\left( \mathbb{P}_\mathbb{C}^3, \mathcal{O}_{\mathbb{P}_\mathbb{C}^3}(4) \right) \right)$ is given by 
{\Small{\begin{align}\label{ja5}
b=&\left(A_{14}-3A_3X_3+ A_{15}X_2+ A_{16}X_1+ A_{17}X_2^2+ A_{18}X_1X_2+ A_{10}X_1^2 \right) dX_1\\
+&\left(- A_{13}X_1^2X_2+A_1X_1^2X_3- A_{12}X_1X_2^2+2A_2X_1X_2X_3- A_{11}X_2^3 \right)dX_1 \notag  \\ 
+&\left(A_4+ A_5X_2+ A_6X_1+ 2A_1X_2X_3+ A_7X_2^2 + A_8X_1X_2+ A_9X_1^2  \right)dX_2 \notag \\
+&\left( A_3X_2^2X_3 + A_{11}X_1X_2^2- 3A_2X_1^2X_3+ A_{12}X_1^2X_2+ A_{13}X_1^3   \right)dX_2 \notag \\
+&\left(-A_2+ A_3X_1+ A_1X_2^2 -A_3X_2^3 + A_2X_1^2X_2-A_1 X_1^3 \right)dX_3,\,\,\,\,\,\,A_1,...,A_{18}\in \mathbb{C} \notag
\end{align}}}
Let us compute $\textnormal{ker}\left(H^0\left( \mathbb{P}_\mathbb{C}^3, \Omega_{\mathbb{P}_\mathbb{C}^3}^1(4) \right)\to H^0\left( \mathbb{P}_\mathbb{C}^3, \mathcal{O}_{\mathbb{P}_\mathbb{C}^3}(4) \right)  \right) \bigoplus H^0\left( \mathbb{P}_\mathbb{C}^3, \mathcal{O}_{\mathbb{P}_\mathbb{C}^3}(4) \right)\to H^0\left( \mathbb{P}_\mathbb{C}^3, \bigwedge^2 \Omega_{\mathbb{P}_\mathbb{C}^3}^1(8) \right)$. From $(\ref{ja5})$ and $\textnormal{(\ref{ja6})}$ and $\textnormal{(\ref{ja7})}$, the image is given by $-w_0\wedge b + c \cdot dw_0 + w_0 \wedge dc$. Then the coefficient of $dX_1\wedge dX_2$ of $-w_0\wedge b + c \cdot dw_0 + w_0 \wedge dc$ is given by
{\Small{\begin{align*}
&\left(A_4-c_3\right) + \left(A_6- c_6\right)X_1+ \left(A_5 - 2c_8\right)X_2 - c_9 X_3 +\left(A_9- c_{12}- 4c_1\right)X_1^2 + \left(A_8 - 2 c_{14} \right) X_1 X_2 - c_{15}X_1X_3 -\left(c_2+ 3c_{17} - A_7 - A_{14}\right)X_2^2  \\
&+ \left(2A_1 -  2c_{18}\right)X_2X_3- c_{19}X_3^2 + \left(-3c_2- c_{24} + A_{13} - A_{14}\right)X_1^3 + \left(-3c_3- 2c_{26} + A_{12}- A_4\right)X_1^2 X_2 - \left(4c_4+c_{27} + 3A_2 \right)X_1^2 X_3 \\
& - \left(2c_5+ 3c_{29} - A_{11}- A_{16}\right)X_1X_2^2  - 2 c_{30}X_1X_2 X_3 - c_{31}X_1X_3^2  - \left(c_6+ 4c_{22} - A_{15}\right)X_2^3 - \left(c_7 + 3c_{33}+ 2 A_3\right)X_2^2 X_3  - 2 c_{34}X_2X_3^2\\
&- c_{35}X_3^3  - \left(2c_5 + A_{16}\right)X_1^4 - \left(2c_6+ A_6 + A_{15}\right)X_1^3X_2  + \left(3A_3 -3c_7\right)X_1^3X_3 + \left(-2c_8- 3c_{11} + A_{10} - A_5\right)X_1^2 X_2^2 -3 c_9 X_1^2 X_2X_3  \\
&-4c_{10}X_1^2X_3^2  + \left(A_{18} - 2 c_{12} \right) X_1X_2^3   - 2 c_{13}X_1X_2^2X_3 + \left(A_{17} - c_{14} \right)X_2^4 -c_{15}X_2^3 X_3  - c_{16}X_2^2X_3^2  - \left(A_{10}+c_{11}\right) X_1^5 -  \left(A_9 + A_{18}+c_{12}\right)X_1^4 X_2 \\
& -2 c_{13}X_1^4X_3 + \left(-c_{14}- 4c_{21}- A_8 - A_{17}\right)X_1^3X_2^2 - 2 c_{15}X_1^3X_2 X_3   -3 c_{16}X_1^3X_3^2   + \left(- c_{17}- 3 c_{24} - A_7 - A_{13}\right)X_1^2X_2^3\\
& +  \left(-2c_{18}-3 c_{25} - A_1\right)X_1^2X_2^2X_3 -3 c_{19}X_1^2X_2X_3^2  -4 c_{20}X_1^2X_3^3  - \left(A_{12}+2 c_{26}\right)X_1X_2^4 + \left(2 A_2 - 2 c_{27}\right)X_1X_2^3X_3 - 2 c_{28}X_1X_2^2X_3^2 \\
& - \left(A_{11} + c_{29}\right)X_2^5   - c_{30}X_2^4 X_3 - c_{31}X_2^3X_3^2  - c_{32}X_2^2X_3^3  - \left(A_1 +c_{25} \right)X_1^5 X_3 + \left(A_2   - c_{27} \right)X_1^4X_2X_3 -2 c_{28}X_1^4X_3^2  -c_{30}X_1^3X_2^2 X_3 \\
& -2c_{31}X_1^3X_2X_3^2  -3c_{32}X_1^3X_3^3  - \left( c_{33} + A_3 \right) X_1^2X_2^3 X_3 - 2 c_{34}X_1^2X_2^2X_3^2 -3 c_{35}X_1^2 X_2X_3^3 -4c_{23}X_1^2 X_3^4\\
\end{align*}}}
and the coefficient of $dX_2\wedge dX_3$ in $-w_0 \wedge b + c \cdot dw_0 + w_0 \wedge dc$ is given by
{\Small{\begin{align*}
&\left(A_2+c_4\right) X_2^2 - \left(A_2+c_4\right)X_1^3 + \left(c_7- A_3\right)X_1X_2^2  + c_9X_2^3 + 2 c_{10} X_2^2X_3 + (A_3  - c_7)X_1^4 - c_9X_1^3X_2 - 2 c_{10}X_1^3 X_3  + c_{13} X_1^2 X_2^2 + c_{15}X_1X_2^3 \\
&+ 2 c_{16}X_1X_2^2X_3 + \left( c_{18}- A_1\right)X_2^4   + 2c_{19}X_2^3X_3 + 3 c_{20}X_2^2X_3^2   - c_{13} X_1^5  - c_{15}X_1^4X_2 - 2 c_{16}X_1^4X_3 + \left(c_{25}  - c_{18}+ 2 A_1\right)X_1^3X_2^2 - 2c_{19}X_1^3X_2X_3 \\
&  - 3 c_{20}X_1^3X_3^2 + \left(c_{27}- A_2\right)X_1^2 X_2^3  + 2 c_{28}X_1^2 X_2^2X_3  + c_{30}X_1X_2^4 + 2c_{31} X_1X_2^3X_3 + 3 c_{32}X_1X_2^2X_3^2  + \left(c_{33}+ A_3\right)X_2^5 + 2 c_{34}X_2^4 X_3 + 3 c_{35}X_2^3X_3^2 \\
& +4 c_{23}X_2^2X_3^3 - \left(c_{25}+ A_1\right)X_1^6  + \left(A_2- c_{27}\right)X_1^5 X_2 - 2 c_{28}X_1^5 X_3  - c_{30}X_1^4X_2^2 - 2c_{31} X_1^4X_2X_3 -3 c_{32}X_1^4X_3^2 -\left(A_3+c_{33}\right)X_1^3X_2^3 - 2 c_{34}X_1^3X_2^2 X_3 \\
& - 3 c_{35}X_1^3X_2X_3^2 -4 c_{23}X_1^3X_3^3 
\end{align*}}}
and the coefficient of $dX_3\wedge dX_1$ in $-w_0 \wedge b + c \cdot dw_0 + w_0 \wedge dc$ is given by
{\Small{\begin{align*}
&\left(A_2+c_4\right)+ \left(c_7- A_3\right)X_1 + c_9X_2 + 2 c_{10} X_3 + c_{13} X_1^2   + c_{15}X_1X_2+ 2 c_{16}X_1X_3 + \left(c_{18} - A_1\right) X_2^2 + 2c_{19}X_2X_3 + 3 c_{20}X_3^2  + \left(c_{25}+ A_1\right)X_1^3 \\
&+ \left(c_{27}- c_4- 2A_2\right) X_1^2 X_2+ 2 c_{28}X_1^2 X_3  + c_{30}X_1X_2^2 + 2c_{31} X_1X_2X_3 + 3 c_{32}X_1X_3^2 + \left(A_3 +c_{33} \right)X_2^3 + 2 c_{34}X_2^2 X_3 + 3 c_{35}X_2X_3^2 +4 c_{23}X_3^3 \\
& + \left(A_3 -c_7 \right)X_1^3 X_2  - c_9X_1^2X_2^2 - 2 c_{10} X_1^2 X_2X_3 - c_{13} X_1^4X_2  - c_{15}X_1^3X_2^2- 2 c_{16}X_1^3 X_2X_3 - 2c_{19}X_1^2X_2^2X_3 - 3 c_{20}X_1^2 X_2X_3^2  + \left(A_1- c_{18}\right)X_1^2X_2^3 \\
& - \left(A_1+c_{25}\right)X_1^5 X_2 + \left(A_2 - c_{27}\right)X_1^4 X_2^2  - 2 c_{28}X_1^4 X_2 X_3  - c_{30}X_1^3X_2^3 - 2c_{31} X_1^3X_2^2X_3 - 3 c_{32}X_1^3 X_2X_3^2 - \left(A_3+ c_{33}\right)X_1^2X_2^4  - 2 c_{34}X_1^2X_2^3 X_3 \\
& - 3 c_{35}X_1^2X_2^2X_3^2 - 4 c_{23} X_1^2 X_2X_3^3\\
\end{align*}}}

Then from $\textnormal{(\ref{ja8})}$ we see that
\begin{align}\label{ja9}
\dim_\mathbb{C} \textnormal{ker}\left( H^0\left( \mathbb{P}_\mathbb{C}^3, \Omega_{\mathbb{P}_\mathbb{C}^3}^1(4) \right) \bigoplus H^0\left( \mathbb{P}_\mathbb{C}^3, \mathcal{O}_{\mathbb{P}_\mathbb{C}^3} (4) \right) \to H^0\left( \mathbb{P}_\mathbb{C}^3,  \mathcal{O}_{\mathbb{P}_\mathbb{C}^3}(4) \right) \bigoplus H^0\left( \mathbb{P}_\mathbb{C}^3 ,  \bigwedge^2 \Omega_{\mathbb{P}_\mathbb{C}^3}^1(8) \right) \right) = 12
\end{align}

 Next we compute $H^0\left( \mathbb{P}_\mathbb{C}^3, \mathcal{E}_{\mathcal{O}_{\mathbb{P}_\mathbb{C}^3}(4)} \right) \xrightarrow{\tilde{G}_0} H^0\left( \mathbb{P}_\mathbb{C}^3, \Omega_{\mathbb{P}_\mathbb{C}^3}^1(4) \right) \bigoplus H^0\left( \mathbb{P}_\mathbb{C}^3, \mathcal{O}_{\mathbb{P}_\mathbb{C}^3}(4) \right)$. We note that from $\textnormal{(\ref{cp8})}$ for $m=4$, the image $\mathcal{L}_T(w_0)+ fw_0$ of $H^0\left( \mathbb{P}_\mathbb{C}^3, \mathcal{E}_{\mathcal{O}_{\mathbb{P}_\mathbb{C}^3}(4)} \right)\to H^0\left( \mathbb{P}_\mathbb{C}^3, \Omega_{\mathbb{P}_\mathbb{C}^3}^1(4) \right)$ is given by $(\ref{j2})$
{\tiny{\begin{align}
&\left(-a_4- a_{16}+3a_5X_3+3a_6X_2+2a_7X_1+2a_1X_1X_2+2a_2X_1X_2X_3+2a_3X_1X_2^2+\left(3a_4+ a_{10}+ a_{16}\right)X_1^2 X_2+a_8X_1^2+a_{11}X_2^2+a_9X_1^2X_3 +a_7X_2^3\right)dX_1 \\
&+ \left(-a_3 - a_6X_1-3a_1X_1^2-3a_2X_1^2X_3-2a_3X_1^2X_2- \left( 3a_4+a_{10}+ a_{16} \right)X_1^3+2a_8X_2+2a_9X_2X_3+ \left(3a_{10}+a_{16}\right)X_2^2+2a_{11}X_1X_2-a_5X_2^2X_3-a_7X_1X_2^2 \right)dX_2 \notag \\
&+\left(-a_2-a_5X_1 +a_2X_1^2X_2 +a_9X_2^2+a_5X_2^3 - a_9X_1^3  \right)dX_3 \notag
\end{align}}}
We compute $H^0\left( \mathbb{P}_\mathbb{C}^3, \mathcal{E}_{\mathcal{O}_{\mathbb{P}_\mathbb{C}^3}(4)} \right)\to H^0\left( \mathbb{P}_\mathbb{C}^3, \mathcal{O}_{\mathbb{P}_\mathbb{C}^3}(4) \right)$. The image $i_T(w_0)$ is given by
{\Small{\begin{align*}
& \left(X_1^2X_2-1\right) \left(a_1+a_2X_3+a_3X_2+a_4X_1+a_5X_1X_3+a_6X_1X_2+a_7X_1^2\right)\\
 &+\left(X_2^2-X_1^3\right)\left(a_8+a_9X_3+a_{10}X_2 +a_{11}X_1+a_5X_2X_3+a_6X_2^2+a_7X_1X_2 \right)\\
 &= - a_1- a_4 X_1 - a_3X_2 - a_2 X_3 -a_7X_1^2 -a_6X_1X_2 -a_5X_1X_3 +a_8X_2^2 -a_8X_1^3+  a_1X_1^2X_2 +a_{11}X_1X_2^2 +a_{10}X_2^3 +a_9X_2^2X_3 \\
 & -a_{11}X_1^4 + \left(a_4- a_{10}\right)X_1^3X_2 -a_9X_1^3X_3 +a_3X_1^2X_2^2+a_2X_1^2X_2X_3+a_7X_1X_2^3 + a_6X_2^4 +a_5X_2^3X_3
\end{align*}}}
Then we see that
\begin{align*}
&\dim_\mathbb{C} \mathbb{H}^0\left( \mathbb{P}_\mathbb{C}^3,  \tilde{E}_{\mathcal{O}_{\mathbb{P}_\mathbb{C}^3}(4)}^\bullet \right) = \dim_\mathbb{C} \textnormal{ker}\left( H^0\left( \mathbb{P}_\mathbb{C}^3 , \mathcal{E}_{\mathcal{O}_{\mathbb{P}_\mathbb{C}^3}(4)}^\bullet \right) \to H^0\left( \mathbb{P}_\mathbb{C}^3, \Omega_{\mathbb{P}_\mathbb{C}^3}^1(4) \right) \bigoplus H^0\left( \mathbb{P}_\mathbb{C}^3 , \mathcal{O}_{\mathbb{P}_\mathbb{C}^3}(4) \right) \right) = 4 \\
&\dim_\mathbb{C} \textnormal{im}\left( H^0\left( \mathbb{P}_\mathbb{C}^3 , \mathcal{E}_{\mathcal{O}_{\mathbb{P}_\mathbb{C}^3}(4)}^\bullet \right) \to H^0\left( \mathbb{P}_\mathbb{C}^3, \Omega_{\mathbb{P}_\mathbb{C}^3}^1(4) \right) \bigoplus H^0\left( \mathbb{P}_\mathbb{C}^3 , \mathcal{O}_{\mathbb{P}_\mathbb{C}^3}(4) \right)  \right)  = 16- 4 = 12
\end{align*}
Then from $(\ref{ja9})$ we have
{\small{\begin{align*}
\dim_\mathbb{C} \mathbb{H}^1\left( \mathbb{P}_\mathbb{C}^3, \tilde{\mathcal{E}}_{\mathcal{O}_{\mathbb{P}_\mathbb{C}^3}(4)}^\bullet \right) = \dim_\mathbb{C} \frac{ \textnormal{ker}\left( H^0\left( \mathbb{P}_\mathbb{C}^3, \Omega_{\mathbb{P}_\mathbb{C}^3}^1(4) \right) \bigoplus H^0\left( \mathbb{P}_\mathbb{C}^3, \mathcal{O}_{\mathbb{P}_\mathbb{C}^3} (4) \right) \to H^0\left( \mathbb{P}_\mathbb{C}^3,  \mathcal{O}_{\mathbb{P}_\mathbb{C}^3}(4) \right) \bigoplus H^0\left( \mathbb{P}_\mathbb{C}^3 ,  \bigwedge^2 \Omega_{\mathbb{P}_\mathbb{C}^3}^1(8) \right) \right)}{\textnormal{im}\left( H^0\left( \mathbb{P}_\mathbb{C}^3 , \mathcal{E}_{\mathcal{O}_{\mathbb{P}_\mathbb{C}^3}(4)}^\bullet \right) \to H^0\left( \mathbb{P}_\mathbb{C}^3, \Omega_{\mathbb{P}_\mathbb{C}^3}^1(4) \right) \bigoplus H^0\left( \mathbb{P}_\mathbb{C}^3 , \mathcal{O}_{\mathbb{P}_\mathbb{C}^3}(4) \right)  \right)} = 12 - 12 =0
\end{align*}}}
Hence $\left( \mathbb{P}_\mathbb{C}^3, w_0 \right)$ is rigid in unfoldings by Theorem \ref{r1}.

\subsection{Unfoldings of $\mathbb{P}_\mathbb{C}^3$ with a rational component of degree $2$} \

Let us consider the foliation of degree $2$ in a rational component on $\mathbb{P}_\mathbb{C}^3$ from subsection \ref{rcc1}. We keep the notations in subsection \ref{rcc1}. We have on $U_0$
\begin{align*}
w_0&= X_1^2 d X_1+ X_2^2 dX_2 + X_3^2 dX_3\\
dw_0&= 0
\end{align*}
Let us compute $\mathbb{H}^1\left( \mathbb{P}_\mathbb{C}^3, \tilde{\mathcal{E}}_{\mathcal{O}_{\mathbb{P}_\mathbb{C}^3}(4)}^\bullet \right)$ by the complex
{\small{\begin{align*}
H^0\left( \mathbb{P}_\mathbb{C}^3, \mathcal{E}_{\mathcal{O}_{\mathbb{P}_\mathbb{C}^3}(4)} \right) \xrightarrow{\tilde{G}_0} H^0\left( \mathbb{P}_\mathbb{C}^3, \Omega_{\mathbb{P}_\mathbb{C}^3}^1(4) \right) \bigoplus H^0\left( \mathbb{P}_\mathbb{C}^3, \mathcal{O}_{\mathbb{P}_\mathbb{C}^3}(4) \right) \xrightarrow{\tilde{G}_1} H^0\left( \mathbb{P}_\mathbb{C}^3 , \mathcal{O}_{\mathbb{P}_\mathbb{C}^3}(4)  \right)\bigoplus H^0\left( \mathbb{P}_\mathbb{C}^3, \tilde{S}^1(4) \right) \xrightarrow{G_2} H^0\left( \mathbb{P}_\mathbb{C}^3, \mathcal{O}_{\mathbb{P}_\mathbb{C}^3}(4) \right) \to 0
\end{align*}}}
Here we note that $H^0\left(\mathbb{P}_\mathbb{C}^3, \tilde{\mathcal{S}}^1(4) \right) \subset H^0\left( \mathbb{P}_\mathbb{C}^3, \bigwedge^2 \Omega_{\mathbb{P}_\mathbb{C}^3}^1(8) \right)$. 

First we compute $\textnormal{ker}\left(H^0\left( \mathbb{P}_\mathbb{C}^3, \Omega_{\mathbb{P}_\mathbb{C}^3}^1(4) \right) \bigoplus H^0\left( \mathbb{P}_\mathbb{C}^3, \mathcal{O}_{\mathbb{P}_\mathbb{C}^3}(4)\right)\to H^0\left( \mathbb{P}_\mathbb{C}^3, \mathcal{O}_{\mathbb{P}_\mathbb{C}^3}(4) \right) \bigoplus H^0\left( \mathbb{P}_\mathbb{C}^3, \bigwedge^2 \Omega_{\mathbb{P}_\mathbb{C}^3}^1(8) \right)\right)$. We note that from $(\ref{rc1})$, any element of $\textnormal{ker}\left( H^0\left( \mathbb{P}_\mathbb{C}^3, \Omega_{\mathbb{P}_\mathbb{C}^3}^1(4) \right)\to H^0\left( \mathbb{P}_\mathbb{C}^3, \mathcal{O}_{\mathbb{P}_\mathbb{C}^3}(4) \right) \right)$ is given by 
{\Small{\begin{align} \label{rcc2}
b=&\left(A_{20}+ A_7X_3+ A_{11}X_2+ A_{21}X_1+ A_9X_3^2+ A_2 X_2 X_3+ A_{13}X_2^2+ 2A_8X_1X_3+ 2A_{12}X_1X_2+ A_{22}X_1^2 \right) dX_1\\
+&\left(- A_6X_1^2X_2+ A_5 X_1^2X_3- A_{10}X_2^3- A_{10}X_3^3 \right)dX_1 \notag \\ 
+&\left( A_{14} + A_1 X_3+ A_{15}X_2+ A_{11}X_1+ A_4 X_3^2+ 2A_3 X_2X_3+ A_{16}X_2^2 + A_2 X_1X_3+ 2A_{13}X_1X_2+ A_{12}X_1^2 + A_6 X_3^3 \right)dX_2 \notag \\
+&\left( A_5 X_2^2X_3 + A_{10} X_1X_2^2+ A_6 X_1^3   \right)dX_2 \notag \\
+&\left(A_{17}+ A_{18}X_3+ A_1 X_2+ A_7 X_1+ A_{19}X_3^2+ 2A_4 X_2X_3+ A_3 X_2^2+ 2A_9X_1X_3+ A_2 X_1X_2+ A_8X_1^2 - A_6 X_2X_3^2- A_5 X_2^3 \right)dX_3 \notag \\
+&\left( A_{10}X_1X_3^2- A_5 X_1^3 \right)dX_3,\,\,\,\,\,\,\, A_1,..., A_{22}\in \mathbb{C} \notag
\end{align}}}
Let us compute $\textnormal{ker}\left(H^0\left( \mathbb{P}_\mathbb{C}^3, \Omega_{\mathbb{P}_\mathbb{C}^3}^1(4) \right)\to H^0\left( \mathbb{P}_\mathbb{C}^3, \mathcal{O}_{\mathbb{P}_\mathbb{C}^3}(4) \right)  \right) \bigoplus H^0\left( \mathbb{P}_\mathbb{C}^3, \mathcal{O}_{\mathbb{P}_\mathbb{C}^3}(4) \right)\to H^0\left( \mathbb{P}_\mathbb{C}^3, \bigwedge^2 \Omega_{\mathbb{P}_\mathbb{C}^3}^1(8) \right)$. From $(\ref{rcc2})$ and $\textnormal{(\ref{ja6})}$ and $\textnormal{(\ref{ja7})}$, the image is given by $-w_0\wedge b  + w_0 \wedge dc$. Then the coefficient of $dX_1\wedge dX_2$ of $-w_0\wedge b + w_0 \wedge dc$ is given by
 {\Small{\begin{align*}
&\left(A_{20}- c_2\right)X_2^2+ \left(A_7- c_7 \right)X_2^2X_3+ \left(A_{11}- c_6\right)X_2^3+ \left(A_{21}-2c_5\right)X_1X_2^2+ \left(A_9- c_{16}\right)X_2^2X_3^2+ \left(A_2- c_{15}\right) X_2^3 X_3+ \left(A_{13}- c_{14}\right)X_2^4 \\
&+ \left(2A_8- 2 c_{13}\right)X_1X_2^2X_3+ \left(2A_{12}- 2c_{12}\right)X_1X_2^3- \left(A_6+ 3 c_{24}- 4c_{22}\right)X_1^2X_2^3+ \left( 3 c_{33}- 3c_{25} \right) X_1^2X_2^2 X_3- \left(A_{10}+ c_{29}\right)X_2^5 - \left(A_{10}+ c_{32} \right)X_2^2X_3^3 \\
&  - \left(A_{14}-c_3\right)X_1^2- \left(A_1- c_9\right) X_1^2 X_3- \left(A_{15}-2c_8 \right)X_1^2X_2- \left( A_{11}-c_6 \right)X_1^3- \left(A_4-c_{19}\right) X_1^2X_3^2 - \left(2A_3 - 2c_{18}\right) X_1^2X_2X_3 - \left(A_2 - c_{15}\right)X_1^3X_3 \\
&- \left(2A_{13}- 2c_{14}\right)X_1^3X_2- \left(A_{12} - c_{12}\right)X_1^4 - \left(A_6 - c_{35}\right) X_1^2X_3^3 - \left(A_{10} + 4c_{21} - 3c_{29}\right)X_1^3X_2^2- \left(A_6- c_{24}\right) X_1^5 + \left(A_{22}- 3c_{11}- A_{16} + 3c_{17} \right)X_1^2 X_2^2\\
&  - 2 c_{26}X_1X_2^4 - 2 c_{27}X_1X_2^3X_3 - 2 c_{28}X_1X_2^2X_3^2 - c_{30}X_2^4 X_3 - c_{31}X_2^3X_3^2 +  2 c_{26}X_1^4 X_2 + c_{27}X_1^4 X_3  + 2 c_{30}X_1^3X_2 X_3 + c_{31}X_1^3X_3^2 +  + 2 c_{34}X_1^2X_2X_3^2 
\end{align*}}}
and the coefficient of $dX_2\wedge dX_3$ in $-w_0 \wedge b  + w_0 \wedge dc$ is given by
{\Small{\begin{align*}
&\left(A_{14} - c_3 \right)X_3^2+ \left( A_1 - c_9 \right) X_3^3+ \left(A_{15}- 2c_8\right)X_2X_3^2+ \left(A_{11}-c_6\right)X_1X_3^2+ \left(A_4 - c_{19}\right) X_3^4+  \left(2A_3 -2 c_{18}\right) X_2X_3^3 + \left(A_2 - c_{15}\right) X_1X_3^3\\
&+ \left(2A_{13}- 2c_{14}\right)X_1X_2X_3^2+ \left(A_{12}-c_{12}\right)X_1^2 X_3^2+ \left(A_6 - c_{35}\right) X_3^5 + \left( A_5- 3c_{33}+ 4 c_{23}\right) X_2^2X_3^3 +  \left( 3c_{32}- 3c_{29}\right)X_1X_2^2X_3^2+ \left( A_6 - c_{24}\right) X_1^3X_3^2\\
&-  \left(A_{17}- c_4\right)X_2^2- \left(A_{18}- 2c_{10}\right)X_2^2X_3- \left(A_1 - c_9\right) X_2^3- \left(A_7- c_7\right) X_1X_2^2- \left(2A_4 - 2 c_{19}\right)X_2^3X_3- \left(A_3 - c_{18}\right)X_2^4- \left(2A_9- 2 c_{16}\right)X_1X_2^2X_3\\
&- \left(A_2 -c_{15}\right) X_1X_2^3- \left(A_8- c_{13}\right)X_1^2X_2^2 + \left(A_6 - 4c_{22}+ 3c_{35}\right) X_2^3X_3^2+ \left(A_5 + c_{33}\right) X_2^5 + \left(A_5 + c_{25}\right) X_1^3X_2^2  + \left(A_{16}+ 3 c_{20}- A_{19}- 3c_{17}\right) X_2^2 X_3^2 \\
&-  2 c_{26}X_1^2 X_2X_3^2 - c_{27}X_1^2 X_3^3  - 2 c_{30}X_1X_2 X_3^3 - c_{31}X_1X_3^4 - 2 c_{34}X_2X_3^4  + c_{27}X_1^2 X_2^3 + 2 c_{28}X_1^2 X_2^2 X_3  + c_{30}X_1X_2^4 + 2c_{31} X_1X_2^3X_3  + 2 c_{34}X_2^4 X_3 
\end{align*}}}
and the coefficient of $dX_3\wedge dX_1$ in $-w_0 \wedge b  + w_0 \wedge dc$ is given by
{\Small{\begin{align*}
&-\left(A_{20}- c_2\right)X_3^2- \left(A_7- c_7\right)X_3^3- \left(A_{11}-c_6\right)X_2X_3^2- \left(A_{21}- 2c_5\right)X_1X_3^2- \left(A_9-c_{16}\right)X_3^4-  \left(A_2- c_{15}\right) X_2 X_3^3- \left(A_{13}- c_{14}\right)X_2^2X_3^2 \\
&- \left(2A_8- 2c_{13}\right)X_1X_3^3- \left(2A_{12} - 2c_{12}\right)X_1X_2X_3^2  + \left(3c_{24} - 3 c_{35}\right)X_1^2X_2X_3^2- \left( A_5 - 3 c_{25}+ 4 c_{23}\right) X_1^2X_3^3+ \left(A_{10}+ c_{29}\right)X_2^3X_3^2+ \left(A_{10} + c_{32}\right)X_3^5\\
& + \left(A_{17} - c_4\right)X_1^2+ \left(A_{18} - 2c_{10}\right)X_1^2X_3+ \left(A_1 - c_9\right) X_1^2X_2+ \left(A_7 - c_7\right) X_1^3+ \left(2A_4 - 2c_{19}\right) X_1^2X_2X_3+ \left(A_3 - c_{18}\right)X_1^2 X_2^2+ \left(2A_9 - 2c_{16}\right)X_1^3X_3\\
&+ \left(A_2 - c_{15}\right) X_1^3X_2+ \left(A_8 - c_{13}\right)X_1^4 - \left( A_5 + c_{33}\right) X_1^2X_2^3+ \left(A_{10}+ 4 c_{21} - 3c_{32}\right)X_1^3X_3^2- \left(A_5 + c_{25}\right) X_1^5  + \left(A_{19}- 3c_{20}- A_{22}+ 3c_{11}\right) X_2^2 X_3^2 \\
&  + 2 c_{26}X_1X_2^2 X_3^2 + 2 c_{27}X_1X_2X_3^3 + 2 c_{28}X_1X_3^4 + c_{30}X_2^2 X_3^3 + c_{31}X_2X_3^4  -c_{27}X_1^4 X_2 - 2 c_{28}X_1^4 X_3  - c_{30}X_1^3X_2^2 - 2c_{31} X_1^3X_2X_3  - 2 c_{34}X_1^2X_2^2 X_3 
\end{align*}}}
Then from $\textnormal{(\ref{ja8})}$ we have
\begin{align}\label{rcc5}
\dim_\mathbb{C} \textnormal{ker}\left( H^0\left( \mathbb{P}_\mathbb{C}^3, \Omega_{\mathbb{P}_\mathbb{C}^3}^1(4) \right)\bigoplus H^0\left( \mathbb{P}_\mathbb{C}^3, \mathcal{O}_{\mathbb{P}_\mathbb{C}^3}(4) \right) \xrightarrow{\tilde{G}_1} H^0\left( \mathbb{P}_\mathbb{C}^3, \bigwedge^2 \Omega_{\mathbb{P}_\mathbb{C}^3}^1(8) \right) \right) = 24
\end{align}
and any element of $\left(b,c  \right) \in \textnormal{ker}\left( H^0\left( \mathbb{P}_\mathbb{C}^3, \Omega_{\mathbb{P}_\mathbb{C}^3}^1(4) \right)\bigoplus H^0\left( \mathbb{P}_\mathbb{C}^3, \mathcal{O}_{\mathbb{P}_\mathbb{C}^3}(4) \right) \xrightarrow{\tilde{G}_1} H^0\left( \mathbb{P}_\mathbb{C}^3, \bigwedge^2 \Omega_{\mathbb{P}_\mathbb{C}^3}^1(8) \right) \right)$ is given by

{\Small{\begin{align}\label{rcc3}
b=&\left(B_2+ B_7X_3+ B_6X_2+  2B_5X_1+ B_{16}X_3^2+ B_{15} X_2 X_3+ B_{14} X_2^2+ 2B_{13}X_1X_3+ 2B_{12}X_1X_2+ \left( B_{24}+ 3 B_{11} -  3 B_{17} \right)X_1^2 \right) dX_1\\
+&\left(- B_{20}X_1^2X_2-  B_{21} X_1^2X_3+ B_{22}X_2^3+ B_{22}X_3^3 \right)dX_1 \notag \\ 
+&\left( B_3 + B_9 X_3+ 2B_8X_2+ B_6X_1+ B_{19}X_3^2+ 2B_{18} X_2X_3+ B_{24}X_2^2 + B_{15} X_1X_3+ 2B_{14}X_1X_2+ B_{12}X_1^2 + B_{20}X_3^3 \right)dX_2 \notag \\
+&\left( -B_{21} X_2^2X_3 - B_{22} X_1X_2^2 +B_{20} X_1^3   \right)dX_2 \notag \\
+&\left(B_{4}+ 2B_{10}X_3+ B_9 X_2+ B_7 X_1+\left( B_{24} - 3B_{17}+ 3 B_{23} \right) X_3^2+ 2B_{19} X_2X_3+ B_{18} X_2^2+ 2B_{16}X_1X_3+ B_{15} X_1X_2+ B_{13}X_1^2 - B_{20} X_2X_3^2 \right)dX_3 \notag \\
+&\left(  B_{21} X_2^3-B_{22}X_1X_3^2+ B_{21} X_1^3 \right)dX_3 \notag
\end{align}}}
{\small{\begin{align}\label{rcc7}
c&=B_1+ B_2X_1+ B_3 X_2+ B_4X_3 +B_5 X_1^2+ B_6 X_1X_2+ B_7 X_1X_3+ B_8 X_2^2 + B_9 X_2X_3+ B_{10}X_3^2 + B_{11}X_1^3 + B_{12}X_1^2 X_2 + B_{13}X_1^2 X_3\\
&+ B_{14}X_1X_2^2 + B_{15}X_1X_2X_3 + B_{16}X_1X_3^2 + B_{17}X_2^3+ B_{18}X_2^2 X_3+ B_{19}X_2X_3^2 + B_{23}X_3^3 + B_{22}X_1^4 + B_{20}X_2^4 + B_{21} X_3^4 + B_{20} X_1^3 X_2  \notag \\
&+ B_{21}X_1^3 X_3 + B_{22} X_1 X_2^3  + B_{22} X_1X_3^3 + B_{21} X_2^3X_3 + B_{20} X_2 X_3^3, \,\,\,\,\,\,\, B_1,..., B_{22}, B_{23}= c_{20},B_{24}= A_{16} \in \mathbb{C} \notag
\end{align}}}

Next we compute $H^0\left( \mathbb{P}_\mathbb{C}^3, \mathcal{E}_{\mathcal{O}_{\mathbb{P}_\mathbb{C}^3}(4)} \right) \xrightarrow{\tilde{G}_0} H^0\left( \mathbb{P}_\mathbb{C}^3, \Omega_{\mathbb{P}_\mathbb{C}^3}^1(4) \right) \bigoplus H^0\left( \mathbb{P}_\mathbb{C}^3, \mathcal{O}_{\mathbb{P}_\mathbb{C}^3}(4) \right)$. We note that from $\textnormal{(\ref{cp8})}$ for $m=4$, the image $\mathcal{L}_T(w_0)+ fw_0$ of $H^0\left( \mathbb{P}_\mathbb{C}^3, \mathcal{E}_{\mathcal{O}_{\mathbb{P}_\mathbb{C}^3}(4)} \right)\to H^0\left( \mathbb{P}_\mathbb{C}^3, \Omega_{\mathbb{P}_\mathbb{C}^3}^1(4) \right)$ is given by $(\ref{rc2})$
{\Small{\begin{align}\label{rcc8}
&\left(2a_1X_1+2a_2X_1X_3+2a_3X_1X_2+\left(3a_4+ a_{16}\right)X_1^2-a_5X_1^2X_3-a_6X_1^2X_2+ a_{11}X_2^2+ a_7 X_2^3+a_{15}X_3^2+ a_7 X_3^3 \right)dX_1\\
&+\left(a_3 X_1^2+ a_6 X_1^3+ 2a_8X_2+2a_9X_2X_3+ \left(3a_{10}+ a_{16}\right)X_2^2 +2a_{11}X_1X_2- a_5X_2^2X_3-a_7X_1X_2^2 + a_{14}X_3^2+ a_6 X_3^3\right)dX_2 \notag \\
&+ \left( a_2X_1^2+ a_5 X_1^3 +  a_9 X_2^2 + a_5 X_2^3+ 2a_{12}X_3+ \left(3a_{13}  + a_{16}\right) X_3^2+2a_{14}X_2X_3+2a_{15}X_1X_3-a_6X_2X_3^2-a_7X_1X_3^2\right)dX_3 \notag
\end{align}}}
We compute $H^0\left( \mathbb{P}_\mathbb{C}^3, \mathcal{E}_{\mathcal{O}_{\mathbb{P}_\mathbb{C}^3}(4)} \right)\to H^0\left( \mathbb{P}_\mathbb{C}^3, \mathcal{O}_{\mathbb{P}_\mathbb{C}^3}(4) \right)$. The image $i_T(w_0)$ is given by
{\Small{\begin{align}\label{rcc6}
&a_1X_1^2+a_2X_1^2X_3+a_3X_1^2X_2+a_4X_1^3+a_5X_1^3X_3+a_6X_1^3X_2+a_7X_1^4+a_8X_2^2+a_9X_2^2X_3+a_{10}X_2^3 +a_{11}X_1X_2^2+a_5X_2^3X_3+a_6X_2^4+a_7X_1X_2^3 \\
&+a_{12}X_3^2+a_{13}X_3^3+a_{14}X_2X_3^2+a_{15}X_1X_3^2+a_5X_3^4+a_6X_2X_3^3+a_7X_1X_3^3 \notag
\end{align}}}
Then we see that
\begin{align*}
&\dim_\mathbb{C} \mathbb{H}^0\left( \mathbb{P}_\mathbb{C}^3,  \tilde{E}_{\mathcal{O}_{\mathbb{P}_\mathbb{C}^3}(4)}^\bullet \right) = \dim_\mathbb{C} \textnormal{ker}\left( H^0\left( \mathbb{P}_\mathbb{C}^3 , \mathcal{E}_{\mathcal{O}_{\mathbb{P}_\mathbb{C}^3}(4)}^\bullet \right) \to H^0\left( \mathbb{P}_\mathbb{C}^3, \Omega_{\mathbb{P}_\mathbb{C}^3}^1(4) \right) \bigoplus H^0\left( \mathbb{P}_\mathbb{C}^3 , \mathcal{O}_{\mathbb{P}_\mathbb{C}^3}(4) \right) \right) = 0 \\
&\dim_\mathbb{C} \textnormal{im}\left( H^0\left( \mathbb{P}_\mathbb{C}^3 , \mathcal{E}_{\mathcal{O}_{\mathbb{P}_\mathbb{C}^3}(4)}^\bullet \right) \to H^0\left( \mathbb{P}_\mathbb{C}^3, \Omega_{\mathbb{P}_\mathbb{C}^3}^1(4) \right) \bigoplus H^0\left( \mathbb{P}_\mathbb{C}^3 , \mathcal{O}_{\mathbb{P}_\mathbb{C}^3}(4) \right)  \right)  = 16
\end{align*}
Then from $(\ref{rcc5})$ we have
{\small{\begin{align*}
\dim_\mathbb{C} \mathbb{H}^1\left( \mathbb{P}_\mathbb{C}^3, \tilde{\mathcal{E}}_{\mathcal{O}_{\mathbb{P}_\mathbb{C}^3}(4)}^\bullet \right) = \dim_\mathbb{C} \frac{ \textnormal{ker}\left( H^0\left( \mathbb{P}_\mathbb{C}^3, \Omega_{\mathbb{P}_\mathbb{C}^3}^1(4) \right) \bigoplus H^0\left( \mathbb{P}_\mathbb{C}^3, \mathcal{O}_{\mathbb{P}_\mathbb{C}^3} (4) \right) \to H^0\left( \mathbb{P}_\mathbb{C}^3,  \mathcal{O}_{\mathbb{P}_\mathbb{C}^3}(4) \right) \bigoplus H^0\left( \mathbb{P}_\mathbb{C}^3 ,  \bigwedge^2 \Omega_{\mathbb{P}_\mathbb{C}^3}^1(8) \right) \right)}{\textnormal{im}\left( H^0\left( \mathbb{P}_\mathbb{C}^3 , \mathcal{E}_{\mathcal{O}_{\mathbb{P}_\mathbb{C}^3}(4)}^\bullet \right) \to H^0\left( \mathbb{P}_\mathbb{C}^3, \Omega_{\mathbb{P}_\mathbb{C}^3}^1(4) \right) \bigoplus H^0\left( \mathbb{P}_\mathbb{C}^3 , \mathcal{O}_{\mathbb{P}_\mathbb{C}^3}(4) \right)  \right)} = 24 - 16 = 8
\end{align*}}}
By considering $(\ref{rcc3}),(\ref{rcc7}),(\ref{rcc8})$ and $(\ref{rcc6})$, we see that
\begin{align*}
\mathbb{H}^1\left( \mathbb{P}_\mathbb{C}^3 , \tilde{\mathcal{E}}_{\mathcal{O}_{\mathbb{P}_\mathbb{C}^3}(4)}^\bullet \right) \cong \textnormal{Span}_\mathbb{C} & \left\langle  \left( dX_1, X_1 \right),\left( dX_2, X_2 \right),\left( dX_3, X_3 \right),\left( d(X_1X_2), X_1X_2 \right),
 \left( d(X_1X_3), X_1X_3 \right) , \right. \\
 &\left. \left( d(X_2X_3), X_2 X_3 \right) , \left( d(X_1X_2X_3),X_1X_2X_3 \right) ,\left( 0,1 \right) \right\rangle
\end{align*}
Then we see that $\left( \mathbb{P}_\mathbb{C}^3\times B  , \omega , \pi, B \right)$, where $\pi:\mathbb{P}_\mathbb{C}^3 \times B \to B$ is the natural projection with the coordinates $\left( t_1,..., t_8 \right)$ on $B \subset \mathbb{C}^8$ for a sufficiently small neighborhood $B$ of $0$ in $\mathbb{C}^8$ and
{\small{\begin{align*}
 \omega&= X_1^2 dX_1 + X_2^2 dX_2 + X_3^3 dX_3 +d\left( t_1X_1 \right)+ d\left(t_2X_2 \right) + d\left(t_3X_3 \right) +d\left( t_4 X_1X_2 \right) +d\left( t_5 X_1X_3\right)+ d\left(t_6 X_2 X_3 \right) \\
 &+ d \left(t_7 X_1X_2X_3 \right) + dt_8 \\
 d\omega & =0 
 \end{align*}}}
 defines an unfolded analytic family of unfoldings of $\pi^{-1}(0)=\left( \mathbb{P}_\mathbb{C}^3, w_0 \right)$ and the unfolded Kodaira-Spencer map is an isomorphism at $t=0$, so that $\left( \mathbb{P}_\mathbb{C}^3, w_0 \right)$ is unobstructed in unfoldings.

\subsection{Unfoldings of $\mathbb{P}_\mathbb{C}^3$ with a logarithmic component of degree $2$} \

Let us consider the foliation of degree $2$ in a logarithmic component on $\mathbb{P}_\mathbb{C}^3$ from subsection \ref{lg1}. We keep the notations in subsection \ref{lg1}. We have on $U_0$
\begin{align*}
w_0&= X_2X_3 dX_1 + X_1X_3 dX_2 + X_1 X_2 dX_3\\
dw_0&= 0
\end{align*}
Let us compute $\mathbb{H}^1\left( \mathbb{P}_\mathbb{C}^3, \tilde{\mathcal{E}}_{\mathcal{O}_{\mathbb{P}_\mathbb{C}^3}(4)}^\bullet \right)$ by the complex
{\small{\begin{align*}
H^0\left( \mathbb{P}_\mathbb{C}^3, \mathcal{E}_{\mathcal{O}_{\mathbb{P}_\mathbb{C}^3}(4)} \right) \xrightarrow{\tilde{G}_0} H^0\left( \mathbb{P}_\mathbb{C}^3, \Omega_{\mathbb{P}_\mathbb{C}^3}^1(4) \right) \bigoplus H^0\left( \mathbb{P}_\mathbb{C}^3, \mathcal{O}_{\mathbb{P}_\mathbb{C}^3}(4) \right) \xrightarrow{\tilde{G}_1} H^0\left( \mathbb{P}_\mathbb{C}^3 , \mathcal{O}_{\mathbb{P}_\mathbb{C}^3}(4)  \right)\bigoplus H^0\left( \mathbb{P}_\mathbb{C}^3, \tilde{S}^1(4) \right) \xrightarrow{G_2} H^0\left( \mathbb{P}_\mathbb{C}^3, \mathcal{O}_{\mathbb{P}_\mathbb{C}^3}(4) \right) \to 0
\end{align*}}}
Here we note that $H^0\left(\mathbb{P}_\mathbb{C}^3, \tilde{\mathcal{S}}^1(4) \right) \subset H^0\left( \mathbb{P}_\mathbb{C}^3, \bigwedge^2 \Omega_{\mathbb{P}_\mathbb{C}^3}^1(8) \right)$. 

First we compute $\textnormal{ker}\left(H^0\left( \mathbb{P}_\mathbb{C}^3, \Omega_{\mathbb{P}_\mathbb{C}^3}^1(4) \right) \bigoplus H^0\left( \mathbb{P}_\mathbb{C}^3, \mathcal{O}_{\mathbb{P}_\mathbb{C}^3}(4)\right)\to H^0\left( \mathbb{P}_\mathbb{C}^3, \mathcal{O}_{\mathbb{P}_\mathbb{C}^3}(4) \right) \bigoplus H^0\left( \mathbb{P}_\mathbb{C}^3, \bigwedge^2 \Omega_{\mathbb{P}_\mathbb{C}^3}^1(8) \right)\right)$. We note that from $(\ref{l2})$, any element of $\textnormal{ker}\left( H^0\left( \mathbb{P}_\mathbb{C}^3, \Omega_{\mathbb{P}_\mathbb{C}^3}^1(4) \right)\to H^0\left( \mathbb{P}_\mathbb{C}^3, \mathcal{O}_{\mathbb{P}_\mathbb{C}^3}(4) \right) \right)$ is given by
{\Small{\begin{align}\label{lg2}
b&=\left(A_{21}+A_2X_3+A_1X_2+A_{22}X_1+A_7X_3^2+A_{23}X_2 X_3+A_6X_2^2+2A_4X_1X_3+2A_5X_1X_2+A_{24}X_1^2 \right) dX_1\\
&+\left(-2A_{10}X_1X_2X_3+A_{11}X_2^2X_3+A_{12}X_2X_3^2 \right)dX_1 \notag \\ 
&+\left(A_{13}+A_3X_3+A_{14}X_2+A_1X_1+A_9 X_3^2+2A_8X_2X_3+A_{15}X_2^2 +A_{16}X_1X_3+2A_6X_1X_2+A_5X_1^2 \right)dX_2 \notag \\
&+\left( A_{12}X_1X_3^2-2A_{11}X_1X_2X_3+A_{10}X_1^2X_3   \right)dX_2 \notag \\
&+\left( A_{17}+ A_{18}X_3+ A_3X_2+ A_2 X_1+ A_{19}X_3^2+2A_9X_2X_3+ A_8X_2^2+2A_7X_1X_3+ A_{20}X_1X_2+ A_4X_1^2 \right)dX_3 \notag \\
&+\left(-2A_{12}X_1X_2X_3+ A_{11}X_1X_2^2+ A_{10}X_1^2X_2\right)dX_3,\,\,\,\,\,\,\,\,\, A_1,..., A_{24} \in \mathbb{C} \notag
\end{align}}}
Let us compute $\textnormal{ker}\left(H^0\left( \mathbb{P}_\mathbb{C}^3, \Omega_{\mathbb{P}_\mathbb{C}^3}^1(4) \right)\to H^0\left( \mathbb{P}_\mathbb{C}^3, \mathcal{O}_{\mathbb{P}_\mathbb{C}^3}(4) \right)  \right) \bigoplus H^0\left( \mathbb{P}_\mathbb{C}^3, \mathcal{O}_{\mathbb{P}_\mathbb{C}^3}(4) \right)\to H^0\left( \mathbb{P}_\mathbb{C}^3, \bigwedge^2 \Omega_{\mathbb{P}_\mathbb{C}^3}^1(8) \right)$. From $(\ref{lg2})$ and $\textnormal{(\ref{ja6})}$ and $\textnormal{(\ref{ja7})}$, the image is given by $-w_0\wedge b  + w_0 \wedge dc$. Then the coefficient of $dX_1\wedge dX_2$ of $-w_0\wedge b + w_0 \wedge dc$ is given by
{\Small{\begin{align*}
&\left(A_{21} - c_2\right)X_1X_3+ \left(A_2-c_7\right)X_1X_3^2+ \left(A_{22}- 2c_5\right)X_1^2 X_3+\left(A_7- c_{16}\right)X_1X_3^3+ \left(A_{23}- A_{16}\right)X_1X_2 X_3^2-\left(A_6- c_{14}\right)X_1X_2^2X_3 \\
&+ \left(2A_4- 2c_{13}\right)X_1^2X_3^2+\left(A_5- c_{12}\right)X_1^2X_2X_3 + \left(A_{24}- 3c_{11}\right)X_1^3 X_3 - \left(3A_{10}+ c_{27}\right)X_1^2X_2X_3^2+ \left(3A_{11} + c_{30}\right) X_1X_2^2X_3^2 \\
&- \left(A_{13}-c_3\right)X_2X_3- \left(A_3- c_9\right)X_2X_3^2- \left(A_{14}- 2c_8\right)X_2^2 X_3- \left(A_9 - c_{19}\right)  X_2 X_3^3- \left(2A_8- 2c_{18}\right)X_2^2X_3^2- \left(A_{15} - 3 c_{17}\right)X_2^3 X_3\\
&- 4c_{21}X_1^4X_3+2c_{29}X_1X_2^3X_3 - c_{32}X_1X_3^4- 2c_{24}X_1^3 X_2X_3 - 3 c_{25}X_1^3 X_3^2  - 2 c_{28}X_1^2X_3^3 +4c_{22}X_2^4X_3 + c_{35}X_2X_3^4    + 3 c_{33}X_2^3X_3^2 + 2 c_{34}X_2^2X_3^3 
\end{align*}}}
and the coefficient of $dX_2\wedge dX_3$ in $-w_0 \wedge b  + w_0 \wedge dc$ is given by
{\Small{\begin{align*}
&\left(A_{13}-c_3\right)X_1X_2+ \left(A_{14}- 2c_8\right)X_1X_2^2+ \left(A_1- c_6\right)X_1^2 X_2- \left(A_9 - c_{19}\right)X_1X_2X_3^2+ \left(A_8 - c_{18}\right) X_1X_2^2X_3+ \left(A_{15}- 3c_{17}\right) X_1X_2^3 \\
& + \left(A_{16}- A_{20}\right)X_1^2 X_2X_3+ \left(2A_6- 2c_{14}\right)X_1^2X_2^2+ \left(A_5- c_{12}\right)X_1^3X_2  + \left(3A_{12} + c_{31}\right)X_1^2 X_2X_3^2- \left(3A_{11} + c_{30}\right) X_1^2X_2^2X_3   \\
&-\left(A_{17}- c_4\right)X_1X_3- \left(A_{18} - 2c_{10}\right) X_1X_3^2  - \left( A_2 - c_7 \right) X_1^2 X_3- \left( A_{19} - 3c_{20}\right)X_1X_3^3- \left(2A_7- 2c_{16}\right)X_1^2X_3^2- \left( A_4 - c_{13}\right) X_1^3X_3  \\
&-c_{24}X_1^4X_2-4c_{22}X_1X_2^4 + 2 c_{35}X_1X_2X_3^3 -  2 c_{26}X_1^3 X_2^2  - 3 c_{29}X_1^2X_2^3   - 2 c_{33}X_1X_2^3 X_3 + c_{25}X_1^4X_3+4 c_{23}X_1X_3^4  + 2 c_{28}X_1^3 X_3^2 + 3 c_{32}X_1^2X_3^3  
\end{align*}}}
and the coefficient of $dX_3\wedge dX_1$ in $-w_0 \wedge b  + w_0 \wedge dc$ is given by
{\Small{\begin{align*}
&\left( A_{17}- c_4\right)X_2X_3+  \left( A_{18} - 2c_{10}\right) X_2X_3^2+ \left( A_3- c_9\right) X_2^2 X_3+ \left(A_{19} -3 c_{20}\right) X_2X_3^3+ \left(2A_9 - 2c_{19}\right)X_2^2X_3^2+  \left(A_8- c_{18}\right)X_2^3X_3 \\
&+ \left(A_7- c_{16}\right)X_1X_2X_3^2+ \left(A_{20}- A_{23}\right)X_1X_2^2 X_3- \left( A_4- c_{13}\right)X_1^2 X_2X_3 - \left(3A_{12}+ c_{31}\right)X_1X_2^2X_3^2+ \left( 3A_{10} + c_{27}\right) X_1^2X_2^2 X_3\\
&- \left(A_{21}-c_2\right)X_1X_2- \left(A_1- c_6\right)X_1X_2^2- \left(A_{22} - 2c_5\right) X_1^2 X_2-\left(A_6 - c_{14}\right) X_1X_2^3- \left(2A_5 - 2 c_{12} \right) X_1^2X_2^2- \left(A_{24} - 3 c_{11} \right) X_1^3X_2\\
&+ 4c_{21}X_1^4X_2+ c_{29}X_1X_2^4 -2 c_{32}X_1X_2X_3^3+ 3 c_{24}X_1^3 X_2^2 + 2 c_{25}X_1^3X_2 X_3 + 2 c_{26}X_1^2X_2^3 -c_{33}X_2^4X_3-4 c_{23}X_2X_3^4   - 2 c_{34}X_2^3 X_3^2 - 3 c_{35}X_2^2X_3^3
\end{align*}}}
Then from $\textnormal{(\ref{ja8})}$ we have
\begin{align}\label{lg5}
\dim_\mathbb{C} \textnormal{ker}\left( H^0\left( \mathbb{P}_\mathbb{C}^3, \Omega_{\mathbb{P}_\mathbb{C}^3}^1(4) \right)\bigoplus H^0\left( \mathbb{P}_\mathbb{C}^3, \mathcal{O}_{\mathbb{P}_\mathbb{C}^3}(4) \right) \xrightarrow{\tilde{G}_1} H^0\left( \mathbb{P}_\mathbb{C}^3, \bigwedge^2 \Omega_{\mathbb{P}_\mathbb{C}^3}^1(8) \right) \right) = 24
\end{align}
and any element of $\left(b,c  \right) \in \textnormal{ker}\left( H^0\left( \mathbb{P}_\mathbb{C}^3, \Omega_{\mathbb{P}_\mathbb{C}^3}^1(4) \right)\bigoplus H^0\left( \mathbb{P}_\mathbb{C}^3, \mathcal{O}_{\mathbb{P}_\mathbb{C}^3}(4) \right) \xrightarrow{\tilde{G}_1} H^0\left( \mathbb{P}_\mathbb{C}^3, \bigwedge^2 \Omega_{\mathbb{P}_\mathbb{C}^3}^1(8) \right) \right)$ is given by
 {\Small{\begin{align}\label{lg8}
b&=\left(B_2+B_7X_3+B_6X_2+2 B_5X_1+B_{16}X_3^2+ B_{24}X_2 X_3+B_{14}X_2^2+2B_{13}X_1X_3+2B_{12}X_1X_2+3B_{11}X_1^2 \right) dX_1\\
&+\left(\frac{2}{3}B_{21}X_1X_2X_3-\frac{1}{3}B_{22}X_2^2X_3-\frac{1}{3}B_{23}X_2X_3^2 \right)dX_1 \notag \\ 
&+\left(B_3+ B_9 X_3+2B_8X_2+ B_6X_1+B_{19}X_3^2+2B_{18}X_2X_3+3B_{17}X_2^2 +B_{24}X_1X_3+2B_{14}X_1X_2+ B_{12}X_1^2 \right)dX_2 \notag \\
&+\left( - \frac{1}{3}B_{23}X_1X_3^2 + \frac{2}{3}B_{22}X_1X_2X_3- \frac{1}{3}B_{21}X_1^2X_3   \right)dX_2 \notag \\
&+\left( B_4+ 2B_{10}X_3+ B_9X_2+ B_7 X_1+ 3B_{20}X_3^2+2B_{19}X_2X_3+ B_{18}X_2^2+2B_{16}X_1X_3+ B_{24}X_1X_2+ B_{13}X_1^2 \right)dX_3 \notag \\
&+\left(\frac{2}{3}B_{23}X_1X_2X_3- \frac{1}{3} B_{22}X_1X_2^2- \frac{1}{3}B_{21} X_1^2X_2\right)dX_3, \notag
\end{align}}}

{\Small{\begin{align}\label{lg9}
c=&B_1+ B_2X_1+ B_3 X_2+ B_4X_3 +B_5 X_1^2+ B_6 X_1X_2+ B_7 X_1X_3+ B_8 X_2^2 + B_9 X_2X_3+ B_{10}X_3^2 + B_{11}X_1^3 + B_{12}X_1^2 X_2 + B_{13}X_1^2 X_3+ B_{14}X_1X_2^2 \\
+& B_{15}X_1X_2X_3 + B_{16}X_1X_3^2 + B_{17}X_2^3+ B_{18}X_2^2 X_3+ B_{19}X_2X_3^2 + B_{20}X_3^3 + B_{21} X_1^2 X_2 X_3 + B_{22}X_1X_2^2 X_3+ B_{23}X_1X_2 X_3^2, \,\,\,\,\,\,\,\,\, B_1,..., B_{24} \in \mathbb{C}  \notag
\end{align}}}

Next we compute $H^0\left( \mathbb{P}_\mathbb{C}^3, \mathcal{E}_{\mathcal{O}_{\mathbb{P}_\mathbb{C}^3}(4)} \right) \xrightarrow{\tilde{G}_0} H^0\left( \mathbb{P}_\mathbb{C}^3, \Omega_{\mathbb{P}_\mathbb{C}^3}^1(4) \right) \bigoplus H^0\left( \mathbb{P}_\mathbb{C}^3, \mathcal{O}_{\mathbb{P}_\mathbb{C}^3}(4) \right)$. We note that from $\textnormal{(\ref{cp8})}$ for $m=4$, the image $\mathcal{L}_T(w_0)+ fw_0$ of $H^0\left( \mathbb{P}_\mathbb{C}^3, \mathcal{E}_{\mathcal{O}_{\mathbb{P}_\mathbb{C}^3}(4)} \right)\to H^0\left( \mathbb{P}_\mathbb{C}^3, \Omega_{\mathbb{P}_\mathbb{C}^3}^1(4) \right)$ is given by $(\ref{l1})$
{\Small{\begin{align} \label{lg7}
&\left(a_{12}X_2+a_8X_3+a_{14}X_2^2+a_9X_3^2+2a_{15}X_1X_2 +2a_{11}X_1X_3+\left(a_4+a_{10}+a_{13}+a_{16}\right)X_2X_3+2a_7X_1X_2X_3-a_6X_2^2X_3-a_5X_2X_3^2\right)dX_1  \\
&+\left(a_{12}X_1+a_1X_3+a_{15}X_1^2+a_2X_3^2+2a_{14}X_1X_2+\left(a_4+a_{10}+a_{13}+a_{16}\right)X_1X_3+2a_3X_2X_3-a_7X_1^2X_3 +2a_6X_1X_2X_3-a_5X_1X_3^2\right)dX_2 \notag \\
&+\left(a_8X_1+a_1X_2 +a_{11}X_1^2 +a_3X_2^2 + \left(a_4+a_{10}+a_{13}+a_{16}\right)X_1X_2+2a_9X_1X_3 +2a_2X_2X_3 -a_7X_1^2X_2 -a_6X_1X_2^2+2a_5X_1X_2X_3\right)dX_3 \notag
\end{align}}}
We compute $H^0\left( \mathbb{P}_\mathbb{C}^3, \mathcal{E}_{\mathcal{O}_{\mathbb{P}_\mathbb{C}^3}(4)} \right)\to H^0\left( \mathbb{P}_\mathbb{C}^3, \mathcal{O}_{\mathbb{P}_\mathbb{C}^3}(4) \right)$. The image $i_T(w_0)$ is given by
{\Small{\begin{align} \label{lg6}
& a_1X_2X_3+a_2X_2X_3^2+a_3X_2^2 X_3+(a_4+ a_{10}+ a_{13})X_1X_2X_3+3a_5X_1X_2X_3^2+3a_6X_1X_2^2 X_3+3a_7X_1^2X_2X_3\\
&+ a_8X_1X_3+a_9X_1X_3^2+a_{11}X_1^2 X_3+ a_{12}X_1X_2+a_{14}X_1X_2^2+a_{15}X_1^2 X_2 \notag
 \end{align}}}
 Then we see that
\begin{align*}
&\dim_\mathbb{C} \mathbb{H}^0\left( \mathbb{P}_\mathbb{C}^3,  \tilde{E}_{\mathcal{O}_{\mathbb{P}_\mathbb{C}^3}(4)}^\bullet \right) = \dim_\mathbb{C} \textnormal{ker}\left( H^0\left( \mathbb{P}_\mathbb{C}^3 , \mathcal{E}_{\mathcal{O}_{\mathbb{P}_\mathbb{C}^3}(4)}^\bullet \right) \to H^0\left( \mathbb{P}_\mathbb{C}^3, \Omega_{\mathbb{P}_\mathbb{C}^3}^1(4) \right) \bigoplus H^0\left( \mathbb{P}_\mathbb{C}^3 , \mathcal{O}_{\mathbb{P}_\mathbb{C}^3}(4) \right) \right) = 2 \\
&\dim_\mathbb{C} \textnormal{im}\left( H^0\left( \mathbb{P}_\mathbb{C}^3 , \mathcal{E}_{\mathcal{O}_{\mathbb{P}_\mathbb{C}^3}(4)}^\bullet \right) \to H^0\left( \mathbb{P}_\mathbb{C}^3, \Omega_{\mathbb{P}_\mathbb{C}^3}^1(4) \right) \bigoplus H^0\left( \mathbb{P}_\mathbb{C}^3 , \mathcal{O}_{\mathbb{P}_\mathbb{C}^3}(4) \right)  \right)  = 16 -  2 = 14
\end{align*}
Then from $(\ref{lg5})$ we have
{\small{\begin{align*}
\dim_\mathbb{C} \mathbb{H}^1\left( \mathbb{P}_\mathbb{C}^3, \tilde{\mathcal{E}}_{\mathcal{O}_{\mathbb{P}_\mathbb{C}^3}(4)}^\bullet \right) = \dim_\mathbb{C} \frac{ \textnormal{ker}\left( H^0\left( \mathbb{P}_\mathbb{C}^3, \Omega_{\mathbb{P}_\mathbb{C}^3}^1(4) \right) \bigoplus H^0\left( \mathbb{P}_\mathbb{C}^3, \mathcal{O}_{\mathbb{P}_\mathbb{C}^3} (4) \right) \to H^0\left( \mathbb{P}_\mathbb{C}^3,  \mathcal{O}_{\mathbb{P}_\mathbb{C}^3}(4) \right) \bigoplus H^0\left( \mathbb{P}_\mathbb{C}^3 ,  \bigwedge^2 \Omega_{\mathbb{P}_\mathbb{C}^3}^1(8) \right) \right)}{\textnormal{im}\left( H^0\left( \mathbb{P}_\mathbb{C}^3 , \mathcal{E}_{\mathcal{O}_{\mathbb{P}_\mathbb{C}^3}(4)}^\bullet \right) \to H^0\left( \mathbb{P}_\mathbb{C}^3, \Omega_{\mathbb{P}_\mathbb{C}^3}^1(4) \right) \bigoplus H^0\left( \mathbb{P}_\mathbb{C}^3 , \mathcal{O}_{\mathbb{P}_\mathbb{C}^3}(4) \right)  \right)} = 24 - 14= 10
\end{align*}}}
By considering $(\ref{lg8}),(\ref{lg9}),(\ref{lg7})$ and $(\ref{lg6})$, we see that
\begin{align*}
\mathbb{H}^1\left( \mathbb{P}_\mathbb{C}^3 , \tilde{\mathcal{E}}_{\mathcal{O}_{\mathbb{P}_\mathbb{C}^3}(4)}^\bullet \right) \cong \textnormal{Span}_\mathbb{C} & \left\langle  \left( dX_1, X_1 \right),\left( d\left(X_1^2\right), X_1^2 \right),\left( d\left(X_1^3\right), X_1^3 \right),\left( dX_2, X_2 \right),
 \left( d\left( X_2^2 \right) , X_2^2 \right) , \right. \\
 &\left. \left( d\left(X_2^3\right), X_2^3 \right) ,\left( dX_3, X_3 \right), \left( d\left( X_3^2 \right), X_3^2 \right),  \left( d\left( X_3^3\right), X_3^3 \right) ,\left( 0,1 \right) \right\rangle
\end{align*}
Then we see that $\left( \mathbb{P}_\mathbb{C}^3\times B  , \omega , \pi, B \right)$, where $\pi:\mathbb{P}_\mathbb{C}^3 \times B \to B$ is the natural projection with the coordinates $\left( t_1,..., t_{10} \right)$ on $B \subset \mathbb{C}^{10}$ for a sufficiently small neighborhood $B$ of $0$ in $\mathbb{C}^{10}$ and
{\small{\begin{align*}
 \omega&= d\left( X_1 X_2 X_3 \right) + d\left(t_1X_1\right)+ d\left(t_2 X_1^2\right)+ d\left(t_3X_1^3\right) + d\left(t_4X_2\right) + d\left(t_5 X_2^2\right) + d \left(t_6 X_2^3\right) + d\left(t_7 X_3\right) + d\left(t_8 X_3^2\right) + d\left(t_9 X_3^3\right) + dt_{10} \\
 d\omega & =0 
 \end{align*}}}
 defines an unfolded analytic family of unfoldings of $\pi^{-1}(0)=\left( \mathbb{P}_\mathbb{C}^3, w_0 \right)$ and the unfolded Kodaira-Spencer map is an isomorphism at $t=0$, so that $\left( \mathbb{P}_\mathbb{C}^3, w_0 \right)$ is unobstructed in unfoldings.

\subsection{Unfoldings of the foliation of degree $2$ of the exceptional component on $\mathbb{P}_\mathbb{C}^3$}\

We consider unfoldings of the foliation of degree $2$ of the exceptional component on $\mathbb{P}_\mathbb{C}^3$ defined by $w\in H^0\left(\mathbb{P}_\mathbb{C}^3, \Omega_{\mathbb{P}_\mathbb{C}^3}^1 ( 4) \right)$ in \textnormal{subsection \ref{cp11}}. We will show that $\mathbb{H}^1\left( \mathbb{P}_\mathbb{C}^3, \tilde{\mathcal{E}}_{\mathcal{O}_{\mathbb{P}_\mathbb{C}^3}(4)}^\bullet \right)=0$, i.e. $\left( \mathbb{P}_\mathbb{C}^3, w \right)$ is rigid in unfoldings by \textnormal{Theorem \ref{r1}}. We keep the notations in subsection \ref{cp11} and \textnormal{Appendix \ref{cp7}} and \textnormal{Appendix \ref{ja1}}. We have on $U_0$
\begin{align*}
w_0 &= \left( 2X_2^2 - 3 X_1X_3  \right) dX_1 + \left( 3X_3 - X_1 X_2 \right) dX_2 + \left( X_1^2 - 2X_2 \right) dX_3 \\
dw_0 &=-5dX_2\wedge dX_3-5X_1dX_3\wedge dX_1-5X_2dX_1\wedge dX_2  
\end{align*}
 We recall the exact sequence of complexes of sheaves from Part III
\begin{center}
$\begin{CD} 
@.  0 @. 0 @. 0 \\
@. @AAA @AAA @AAA \\
0 @>>>  \mathcal{O}_{\mathbb{P}_\mathbb{C}^3}(4) @>>> \mathcal{O}_{\mathbb{P}_{\mathbb{C}}^3}(4) @>>> 0 @>>> 0\\
@. @AAA @AAA @AAA \\
0 @>>>  \tilde{\mathcal{S}}^1(4) @>>> \mathcal{O}_{\mathbb{P}_\mathbb{C}^3}(4) \bigoplus \tilde{\mathcal{S}}^1(4) @>>> \mathcal{O}_{\mathbb{P}_\mathbb{C}^3}(4) @>>> 0 \\
@. @AAA @AAA @AAA \\
0@>>> \mathcal{O}_{\mathbb{P}_\mathbb{C}^3}(4) @>>>  \Omega_{\mathbb{P}_\mathbb{C}^3}^1(4) \bigoplus \mathcal{O}_{\mathbb{P}_\mathbb{C}^3} (4) @>>>  \Omega_{\mathbb{P}_\mathbb{C}^3}^1 (4)  @>>> 0 \\
 @. @AAA @AAA @AAA \\
0 @>>> 0 @>>> \mathcal{E}_{\mathcal{O}_{\mathbb{P}_\mathbb{C}^3}(4)} @>>>  \mathcal{E}_{\mathcal{O}_{\mathbb{P}_\mathbb{C}^3}(4)}  @>>> 0
\end{CD}$
\end{center}
which gives the long exact sequence
\begin{align*}
0 \to \mathbb{H}^0\left( \mathbb{P}_\mathbb{C}^3, \tilde{\mathcal{E}}_{\mathcal{O}_{\mathbb{P}_\mathbb{C}^3}(4)}^\bullet \right) \xrightarrow{}  \mathbb{H}^0\left( \mathbb{P}_\mathbb{C}^3, \mathcal{E}_{\mathcal{O}_{\mathbb{P}_\mathbb{C}^3}(4)}^\bullet \right)  \to \mathbb{H}^0\left( \mathbb{P}_\mathbb{C}^3, \mathcal{O}_{\mathbb{P}_\mathbb{C}^3}(4)^\bullet \right) \to \mathbb{H}^1\left( \mathbb{P}_\mathbb{C}^3, \tilde{\mathcal{E}}_{\mathcal{O}_{\mathbb{P}_\mathbb{C}^3}(4)}^\bullet  \right) \to \mathbb{H}^1\left( \mathbb{P}_\mathbb{C}^3, \mathcal{E}_{\mathcal{O}_{\mathbb{P}_\mathbb{C}^3}(4)}^\bullet \right)=0 \to \cdots
\end{align*}
Since $\mathbb{H}^1\left( \mathbb{P}_\mathbb{C}^3, \mathcal{E}_{\mathcal{O}_{\mathbb{P}_\mathbb{C}^3}(4)}^\bullet \right)=0$ from subsection \ref{cp11}, if we show $\mathbb{H}^0\left(\mathbb{P}_\mathbb{C}^3, \mathcal{O}_{\mathbb{P}_\mathbb{C}}(4)^\bullet  \right)=0$, the long exact sequence implies that $\mathbb{H}^1\left( \mathbb{P}_\mathbb{C}^3, \tilde{\mathcal{E}}_{ \mathcal{O}_{\mathbb{P}_\mathbb{C}^3}(4)}^\bullet \right)=0$. Let us compute
\begin{align*}
\mathbb{H}^0\left( \mathbb{P}_\mathbb{C}^3, \mathcal{O}_{\mathbb{P}_\mathbb{C}^3}(4)^\bullet \right)=\textnormal{ker}\left( H^0\left( \mathbb{P}_\mathbb{C}^3, \mathcal{O}_{\mathbb{P}_\mathbb{C}^3} (4)  \right) \to  H^0\left( \mathbb{P}_\mathbb{C}^3, \bigwedge^2 \Omega_{\mathbb{P}_\mathbb{C}^3}^1(8) \right) \right)
\end{align*}
From $\textnormal{(\ref{ja6})}$ and $\textnormal{(\ref{ja7})}$ the image is given by $w_0 \wedge d c + c \cdot dw_0$. Then the coefficient of $dX_1 \wedge dX_2$ in $w_0 \wedge d c + c \cdot dw_0$ is given by
{\small{\begin{align*}
&-5 c_1X_2 - 3c_2 X_3 -3c_3 X_2^2 -3c_7X_3^2 + \left( - 4c_2-3c_3\right)X_1X_2 -6c_5 X_1X_3 + \left(-3c_6 -5c_4\right) X_2X_3 - c_8 X_2^3  -3c_{16}X_3^3 + \left( -3c_6 - 3c_5\right)X_1^2 X_2 \\
& - 9 c_{11} X_1^2 X_3 + \left(  - 6c_8-2 c_6\right)X_1X_2^2 + \left(-3c_9 - 6c_{12}  - 4c_7 \right)X_1X_2X_3  -6 c_{13}  X_1X_3^2 + \left( -3c_9 - 3c_{14}  \right) X_2^2 X_3 + \left(-3c_{15} - 5c_{10} \right)X_2X_3^2  \\
&+ c_{17} X_2^4 -3 c_{32} X_3^4 + \left( -3 c_{12} -2c_{11}\right)X_1^3 X_2  -12c_{21} X_1^3 X_3 + \left( -6c_{14} - c_{12}\right) X_1^2 X_2^2 + \left(-3c_{15} -3c_{13}  - 9 c_{24} \right)X_1^2 X_2X_3 -9c_{25} X_1^2 X_3^2 \\
& - 9 c_{17} X_1X_2^3 + \left( -2c_{15} - 6c_{18} - 6c_{26} \right)X_1X_2^2X_3  + \left( -6 c_{27} -4 c_{16}  - 3 c_{19}\right)X_1X_2X_3^2 -6 c_{28} X_1X_3^3 + \left(  - 3 c_{29} -  c_{18} \right)X_2^3 X_3 \\
& + \left( -3c_{19} - 3 c_{30} \right) X_2^2 X_3^2 + \left(-3 c_{31} - 5 c_{20} \right) X_2X_3^3 +  3 c_{22} X_2^5 + \left( -3 c_{24} -  c_{21}\right)X_1^4 X_2   - 6 c_{26}  X_1^3 X_2^2  + \left( -3 c_{27} -2c_{25} \right)X_1^3 X_2X_3 \\
& + \left(  c_{26}  - 9 c_{29}  \right) X_1^2 X_2^3 + \left( - 6 c_{30} -  c_{27} \right) X_1^2 X_2^2 X_3 + \left( -3 c_{31} -3 c_{28} \right)X_1^2 X_2 X_3^2 + \left(- 12 c_{22} +2 c_{29} \right) X_1X_2^4 - 9 c_{33}  X_1X_2^3 X_3 \\
&+ \left( - 6 c_{34}- 2 c_{31}\right)X_1X_2^2 X_3^2 + \left( - 3 c_{35}  -4 c_{32} \right) X_1X_2X_3^3 + 6c_{33}X_2^4 X_3 - c_{34} X_2^3 X_3^2 - 3 c_{35}X_2^2 X_3^3 - 5 c_{23}X_2X_3^4=0 \\
&\iff c_1= \cdots = c_{35} = 0
\end{align*}}}
Then from $\textnormal{(\ref{ja8})}$, we see that $\mathbb{H}^0\left( \mathbb{P}_\mathbb{C}^3, \mathcal{O}_{\mathbb{P}_\mathbb{C}^3}(4)^\bullet \right)=\textnormal{ker}\left( H^0\left( \mathbb{P}_\mathbb{C}^3, \mathcal{O}_{\mathbb{P}_\mathbb{C}^3} (4)  \right) \to  H^0\left( \mathbb{P}_\mathbb{C}^3, \bigwedge^2 \Omega_{\mathbb{P}_\mathbb{C}^3}^1(8) \right) \right)=0$.

\subsection{Unfoldings of a Monge-Amp\`ere foliation on $\mathbb{P}_\mathbb{C}^3$}\

We study unfoldings of the codimension $1$-foliation $w_0\in H^0\left( \mathbb{P}_\mathbb{C}^3, \Omega_{\mathbb{P}_\mathbb{C}^3}^1(4) \right)$ on $\mathbb{P}_\mathbb{C}^3$ in subsection \ref{ma1}. We keep the notations in \ref{ma1}. We have on $U_0$
\begin{align*}
w_0&= -2X_1X_2X_3dX_1+X_1^2X_3dX_2+ \left(X_1^2X_2-X_2^2 \right)dX_3\\
dw_0&= -2X_2dX_2\wedge dX_3-4X_1X_2dX_3\wedge dX_1+4X_1X_3dX_1\wedge dX_2
\end{align*}

The author could not determine obstructedness or unobstructedness of unfoldings of $\left( \mathbb{P}_\mathbb{C}^3, w_0 \right)$. We only describe $\mathbb{H}^1\left( \mathbb{P}_\mathbb{C}^3, \tilde{\mathcal{E}}_{\mathcal{O}_{\mathbb{P}_\mathbb{C}^3}(4)}^\bullet \right)$. Let us compute $\mathbb{H}^1\left( \mathbb{P}_\mathbb{C}^3, \tilde{\mathcal{E}}_{\mathcal{O}_{\mathbb{P}_\mathbb{C}^3}(4)}^\bullet \right)$ by the complex
{\small{\begin{align*}
H^0\left( \mathbb{P}_\mathbb{C}^3, \mathcal{E}_{\mathcal{O}_{\mathbb{P}_\mathbb{C}^3}(4)} \right) \xrightarrow{\tilde{G}_0} H^0\left( \mathbb{P}_\mathbb{C}^3, \Omega_{\mathbb{P}_\mathbb{C}^3}^1(4) \right) \bigoplus H^0\left( \mathbb{P}_\mathbb{C}^3, \mathcal{O}_{\mathbb{P}_\mathbb{C}^3}(4) \right) \xrightarrow{\tilde{G}_1} H^0\left( \mathbb{P}_\mathbb{C}^3 , \mathcal{O}_{\mathbb{P}_\mathbb{C}^3}(4)  \right)\bigoplus H^0\left( \mathbb{P}_\mathbb{C}^3, \tilde{S}^1(4) \right) \xrightarrow{G_2} H^0\left( \mathbb{P}_\mathbb{C}^3, \mathcal{O}_{\mathbb{P}_\mathbb{C}^3}(4) \right) \to 0
\end{align*}}}
Here we note that $H^0\left(\mathbb{P}_\mathbb{C}^3, \tilde{\mathcal{S}}^1(4) \right) \subset H^0\left( \mathbb{P}_\mathbb{C}^3, \bigwedge^2 \Omega_{\mathbb{P}_\mathbb{C}^3}^1(8) \right)$. 

First we compute $\textnormal{ker}\left(H^0\left( \mathbb{P}_\mathbb{C}^3, \Omega_{\mathbb{P}_\mathbb{C}^3}^1(4) \right) \bigoplus H^0\left( \mathbb{P}_\mathbb{C}^3, \mathcal{O}_{\mathbb{P}_\mathbb{C}^3}(4)\right)\to H^0\left( \mathbb{P}_\mathbb{C}^3, \mathcal{O}_{\mathbb{P}_\mathbb{C}^3}(4) \right) \bigoplus H^0\left( \mathbb{P}_\mathbb{C}^3, \bigwedge^2 \Omega_{\mathbb{P}_\mathbb{C}^3}^1(8) \right)\right)$. We note that from $(\ref{cm4})$, any element of $\textnormal{ker}\left( H^0\left( \mathbb{P}_\mathbb{C}^3, \Omega_{\mathbb{P}_\mathbb{C}^3}^1(4) \right)\to H^0\left( \mathbb{P}_\mathbb{C}^3, \mathcal{O}_{\mathbb{P}_\mathbb{C}^3}(4) \right) \right)$ is given by

{\Small{\begin{align}\label{ma2}
b=&\left(-A_1X_2-2A_2X_1-A_6X_2 X_3-A_{12}X_2^2-2A_8X_1X_3-2A_4X_1X_2- A_1X_1^2 \right) dX_1\\
+&\left(-A_{12}X_1^2X_2-A_7X_1^2X_3+2A_5X_1X_2^2+A_3X_1X_3^2-\left(2A_{10}+2A_{11}\right)X_1X_2X_3-A_{19}X_2^2X_3-A_{17}X_2X_3^2 \right)dX_1\notag \\ 
+&\left(A_2+A_9X_3+A_1X_1+A_{13}X_3^2+A_{10}X_2X_3+A_5X_2^2 +A_6X_1X_3+A_4X_1^2 +A_{14}X_3^3+A_{15}X_2X_3^2 \right)dX_2 \notag \\
+&\left( A_{16}X_2^2X_3 +A_{17}X_1X_3^2+A_{11}X_1^2X_3-2A_5X_1^2X_2+A_{12}X_1^3   \right)dX_2 \notag \\
+&\left( \left(A_9-2A_8\right)X_2+A_3X_2X_3+A_{18}X_2^2+\left(A_6-2A_7\right)X_1X_2+A_8X_1^2 -A_{14}X_2X_3^2-A_{15}X_2^2X_3-A_{16}X_2^3 \right)dX_3 \notag \\
+&\left(A_{19}X_1X_2^2-A_3X_1^2X_3+ \left(2A_{10}+A_{11}\right)X_1^2X_2+A_7X_1^3 \right)dX_3 ,\,\,\,\,\,\,\,\,A_1,..., A_{19}\in \mathbb{C}, \notag
\end{align}}}

Let us compute $\textnormal{ker}\left(H^0\left( \mathbb{P}_\mathbb{C}^3, \Omega_{\mathbb{P}_\mathbb{C}^3}^1(4) \right)\to H^0\left( \mathbb{P}_\mathbb{C}^3, \mathcal{O}_{\mathbb{P}_\mathbb{C}^3}(4) \right)  \right) \bigoplus H^0\left( \mathbb{P}_\mathbb{C}^3, \mathcal{O}_{\mathbb{P}_\mathbb{C}^3}(4) \right)\to H^0\left( \mathbb{P}_\mathbb{C}^3, \bigwedge^2 \Omega_{\mathbb{P}_\mathbb{C}^3}^1(8) \right)$. From $(\ref{rcc2})$ and $\textnormal{(\ref{ja6})}$ and $\textnormal{(\ref{ma2})}$, the image is given by $-w_0\wedge b  + c\cdot dw_0 + w_0 \wedge dc$. Then the coefficient of $dX_1\wedge dX_2$ of $-w_0\wedge b + c \cdot d w_0+  w_0 \wedge dc$ is given by
{\Small{\begin{align*}
&4c_1X_1X_3+ 3c_2X_1^2 X_3+ \left(2c_3+ 2A_2\right) X_1X_2X_3+ 4c_4X_1X_3^2 + \left(2c_5- 2A_2\right) X_1^3X_3+ \left(A_1+c_6\right) X_1^2X_2X_3+ 3c_7 X_1^2X_3^2 + \left(2c_9+ 2 A_9\right)X_1 X_2X_3^2\\
&+ 4c_{10}X_1X_3^3 +  \left(c_{11}- A_1\right)X_1^4X_3  + \left(2c_{13}- 2A_8\right)X_1^3 X_3^2 - \left(c_{14}+ A_{12}\right)X_1^2X_2^2X_3 + \left(A_6+c_{15}\right)X_1^2X_2X_3^2 + 3c_{16}X_1^2X_3^3 + \left(2A_5-2c_{17}\right)X_1X_2^3X_3 \\
&+ \left(2 c_{19}+ 2 A_{13}\right)X_1X_2X_3^3 +4 c_{20}X_1X_3^4  -4  c_{22}X_1X_2^4X_3 + 4c_{23} X_1X_3^5 +  \left( A_{12} - c_{24}\right) X_1^4 X_2X_3 + \left(c_{25}- A_7\right)X_1^4 X_3^2 - \left(2c_{26} + 2A_5\right)X_1^3X_2^2 X_3   \\
& + \left(A_3 + 2c_{28}\right) X_1^3 X_3^3 -3c_{29}X_1^2 X_2^3 X_3 -\left(A_{19}+c_{30}\right)X_1^2X_2^2 X_3^2+ \left(c_{31} + A_{17}\right) X_1^2X_2 X_3^3 + 3c_{32}X_1^2X_3^4 + \left(2A_{16} -2 c_{33}\right)X_1X_2^3X_3^2 \\
& + \left(2c_{35}+ 2 A_{14}\right)X_1 X_2 X_3^4 -2A_{10}X_1^3X_2X_3^2   +2A_{10}X_1X_2^2X_3^2
 +2A_{15}X_1X_2^2X_3^3  
\end{align*}}}
and the coefficient of $dX_2\wedge dX_3$ in $-w_0 \wedge b + c \cdot dw_0 + w_0 \wedge dc$ is given by
{\Small{\begin{align*}
&-2c_1X_2-2 c_2X_1X_2- \left(A_2+c_3 \right)X_2^2-2 c_4X_2X_3 - \left(2c_5 +c_3- A_2\right)X_1^2X_2 - \left(A_1 + c_6 \right) X_1X_2^2-2 c_7 X_1X_2X_3 - \left(A_9 + c_9 \right) X_2^2X_3-2 c_{10}X_2X_3^2 \\
&- \left(2 c_{11}+ c_6 - A_1\right)X_1^3X_2 - \left(2c_8+c_{12} + A_4\right)X_1^2 X_2^2 + \left(2A_8 -2 c_{13}\right)X_1^2X_2 X_3 - \left( A_6 + c_{15} \right) X_1X_2^2X_3 -2c_{16}X_1X_2X_3^2 + \left(c_{17}- A_5\right)X_2^4 \\
& -A_{10}X_2^3X_3 - \left(A_{13} + c_{19}\right) X_2^2X_3^2 -2 c_{20}X_2X_3^3 - \left(2c_{21}+ c_{12}- A_4\right)X_1^4X_2 +2 c_{22}X_2^5 -2 c_{23} X_2X_3^4 - \left( c_{24} + 2c_{14} + A_{12} \right)X_1^3 X_2^2 \\ 
&+ \left( 2A_7  -2 c_{25} \right)X_1^3 X_2X_3 + \left(3A_5-  3 c_{17}\right) X_1^2 X_2^3   -\left( c_{27} +  c_{18} + A_{11} - A_{10} + A_{18} \right)  X_1^2 X_2^2 X_3 + \left(c_{19}-2c_{28} + A_{13}- A_3 \right)X_1^2 X_2 X_3^2 + c_{29}X_1 X_2^4 \\
& - \left(A_{17} + c_{31}\right) X_1X_2^2 X_3^2 -2 c_{32}X_1X_2X_3^3  + \left( c_{33} - A_{16}\right) X_2^4X_3 -A_{15}X_2^3X_3^2   - \left( A_{14} + c_{35}\right)X_2^2 X_3^3 + \left(A_{12}-c_{24}\right) X_1^5X_2-4c_{22}X_1^2X_2^4 \\
& + \left(2 c_{35}+ 2 A_{14}\right)X_1^2 X_2X_3^3 - \left(2A_5 +2 c_{26}\right)X_1^4 X_2^2  - 3 c_{29}X_1^3X_2^3 - \left(A_{19} + c_{30}\right) X_1^3X_2^2 X_3  + \left( A_{17} + c_{31}\right) X_1^3X_2X_3^2 + \left( 2 A_{16} - 2 c_{33}\right)X_1^2X_2^3 X_3 \\
&+ c_4 X_1^2 X_3+ c_7X_1^3 X_3  + 2 c_{10}X_1^2 X_3^2 + \left(c_{13} - A_8\right) X_1^4 X_3  + 3 c_{20}X_1^2X_3^3  + 2 c_{16}X_1^3X_3^2   + \left( c_{25} - A_7 \right) X_1^5X_3+4 c_{23}X_1^2X_3^4  + \left( 2 c_{28} + A_3 \right)X_1^4 X_3^2 \\
&  + 3 c_{32}X_1^3X_3^3 +2A_{15}X_1^2X_2^2X_3^2    -2A_{10}X_1^4X_2X_3   
\end{align*}}}
and the coefficient of $dX_3\wedge dX_1$ in $-w_0 \wedge b + c \cdot dw_0 + w_0 \wedge dc$ is given by
{\Small{\begin{align*}
&-4c_1X_1X_2- 3c_2X_1^2 X_2- \left(4c_3 + 2 c_5 + 2 A_2 \right)  X_1 X_2^2- 2c_4X_1X_2X_3 + \left( 2 A_2 - 2c_5  \right)X_1^3X_2- \left(3 c_6 + 3 c_{11}\right) X_1^2X_2^2 - c_7 X_1^2 X_2X_3 \\
&- \left(4 c_8 + 2 c_{12}+ 2 A_4\right) X_1X_2^3  -2 \left(c_9+ c_{13} + A_9 - A_8\right) X_1X_2^2X_3 + \left(A_1 -  c_{11}\right) X_1^4 X_2 - \left(2 c_{12} + 4 c_{21}- 2 A_4\right)X_1^3 X_2^2 - \left(3 c_{14} + 3 c_{24} \right)X_1^2X_2^3 \\
&-\left(c_{15}+ 3c_{25} + A_6 - 3A_7\right)X_1^2X_2^2X_3 +c_{16}X_1^2 X_2X_3^2 - \left(4 c_{17}+ 2 c_{26} - 2 A_5\right)X_1X_2^4- \left(2c_{18}+ 2 c_{27} + 2 A_{18}+ 2 A_{10} + 2A_{11}\right)X_1X_2^3 X_3 \\
&- \left( A_3 + 2c_{28} \right) X_1 X_2^2 X_3^2 +2  c_{20}X_1X_2X_3^3   -4 c_{22}X_1X_2^5  +4 c_{23} X_1X_2X_3^4 + \left( A_{12} - c_{24}\right) X_1^4 X_2^2 + \left( c_{25} - A_7 \right)X_1^4 X_2 X_3 - \left( 2 A_5 +  2 c_{26}\right) X_1^3 X_2^3   \\
&  + \left( A_3 + 2 c_{28}\right)X_1^3 X_2 X_3^2 -3 c_{29}X_1^2 X_2^4 - \left( A_{19} +  c_{30}\right) X_1^2 X_2^3 X_3+  \left( A_{17}+ c_{31}\right) X_1^2X_2^2 X_3^2 + 3 c_{32}X_1^2 X_2X_3^3 + \left( 2 A_{16} -2 c_{33}\right)X_1X_2^4X_3 \\
& + \left( 2 A_{14} + 2c_{35} \right) X_1X_2^2 X_3^3  -c_2X_2^2  - \left(A_1 + c_6 \right) X_2^3 -c_7X_2^2X_3  - \left(A_{12} + c_{14}\right) X_2^4 - c_{16}X_2^2X_3^2   -  \left(A_6 + c_{15}\right)X_2^3 X_3\\
& - c_{29}X_2^5 - c_{32}X_2^2 X_3^3   - \left( A_{19}+ c_{30}\right)X_2^4 X_3 - \left( A_{17} + c_{31}\right) X_2^3X_3^2 +2A_{15}X_1X_2^3X_3^2  -2A_{10}X_1^3X_2^2 X_3 
\end{align*}}}
Then from $\textnormal{(\ref{ja8})}$ we have
\begin{align}\label{ma5}
\dim_\mathbb{C} \textnormal{ker}\left( H^0\left( \mathbb{P}_\mathbb{C}^3, \Omega_{\mathbb{P}_\mathbb{C}^3}^1(4) \right)\bigoplus H^0\left( \mathbb{P}_\mathbb{C}^3, \mathcal{O}_{\mathbb{P}_\mathbb{C}^3}(4) \right) \xrightarrow{\tilde{G}_1} H^0\left( \mathbb{P}_\mathbb{C}^3, \bigwedge^2 \Omega_{\mathbb{P}_\mathbb{C}^3}^1(8) \right) \right) = 20
\end{align}
and any element of $\left(b,c  \right) \in \textnormal{ker}\left( H^0\left( \mathbb{P}_\mathbb{C}^3, \Omega_{\mathbb{P}_\mathbb{C}^3}^1(4) \right)\bigoplus H^0\left( \mathbb{P}_\mathbb{C}^3, \mathcal{O}_{\mathbb{P}_\mathbb{C}^3}(4) \right) \xrightarrow{\tilde{G}_1} H^0\left( \mathbb{P}_\mathbb{C}^3, \bigwedge^2 \Omega_{\mathbb{P}_\mathbb{C}^3}^1(8) \right) \right)$ is given by
{\Small{\begin{align}\label{ma7}
b=&\left(B_2X_2+2B_1X_1+ B_8X_2 X_3+ B_7 X_2^2-2B_6X_1X_3-2 \left(B_5- B_4\right)X_1X_2+ B_2X_1^2 \right) dX_1\\
+&\left(B_7X_1^2X_2- B_{14}X_1^2X_3+2B_9X_1X_2^2- 2B_{15}X_1X_3^2+2\left(B_{10}+ B_{11} + B_{12}\right)X_1X_2X_3+B_{16}X_2^2X_3+B_{17}X_2X_3^2 \right)dX_1\notag \\ 
+&\left(- B_1 - B_3X_3- B_2X_1- B_{13}X_3^2+ B_9X_2^2 - B_8 X_1X_3+ \left(B_5- B_4\right)X_1^2 - B_{20}X_3^3 \right)dX_2 \notag \\
+&\left( B_{18}X_2^2X_3 - B_{17}X_1X_3^2- \left( B_{10} + B_{11} + B_{12} \right)X_1^2X_3-2 B_9X_1^2X_2- B_7X_1^3   \right)dX_2 \notag \\
+&\left( \left(- B_3-2B_6\right)X_2- 2B_{15}X_2X_3+ B_{10}X_2^2+\left(- B_8 -2B_{14}\right)X_1X_2+ B_6 X_1^2 +B_{20}X_2X_3^2- B_{18}X_2^3 \right)dX_3 \notag \\
+&\left( - B_{16} X_1X_2^2+ 2 B_{15}X_1^2X_3- \left(B_{10} + B_{11} + B_{12} \right)X_1^2X_2+ B_{14}X_1^3 \right)dX_3  \notag
\end{align}}}
{\Small{\begin{align}\label{ma8}
c=& B_1 X_2 - B_1 X_1^2+ B_2 X_1X_2+  B_4 X_2^2 + B_3 X_2X_3 - B_2 X_1^3 - \left(B_4 + B_5\right)X_1^2 X_2 + B_6X_1^2 X_3+ B_7X_1X_2^2 \\
+& B_8X_1X_2X_3 + B_9 X_2^3+ B_{12}X_2^2 X_3+ B_{13}X_2X_3^2 + B_5X_1^4  - B_7 X_1^3 X_2 + B_{14}X_1^3 X_3  - B_9 X_1^2 X_2^2  \notag \\
+& B_{11} X_1^2 X_2 X_3 + B_{15}X_1^2 X_3^2  + B_{16} X_1X_2^2 X_3+ B_{17} X_1X_2 X_3^2  + B_{18}X_2^3X_3 + B_{19}X_2^2 X_3^2 + B_{20} X_2 X_3^3, \,\,\,\,\,\,\,\,\,\, B_1 ,...,B_{20} \in \mathbb{C} \notag
\end{align}}}

Next we compute $H^0\left( \mathbb{P}_\mathbb{C}^3, \mathcal{E}_{\mathcal{O}_{\mathbb{P}_\mathbb{C}^3}(4)} \right) \xrightarrow{\tilde{G}_0} H^0\left( \mathbb{P}_\mathbb{C}^3, \Omega_{\mathbb{P}_\mathbb{C}^3}^1(4) \right) \bigoplus H^0\left( \mathbb{P}_\mathbb{C}^3, \mathcal{O}_{\mathbb{P}_\mathbb{C}^3}(4) \right)$. We note that from $\textnormal{(\ref{cp8})}$ for $m=4$, the image $\mathcal{L}_T(w_0)+ fw_0$ of $H^0\left( \mathbb{P}_\mathbb{C}^3, \mathcal{E}_{\mathcal{O}_{\mathbb{P}_\mathbb{C}^3}(4)} \right)\to H^0\left( \mathbb{P}_\mathbb{C}^3, \Omega_{\mathbb{P}_\mathbb{C}^3}^1(4) \right)$ is given by $(\ref{cm3})$

{\Small{\begin{align}\label{ma3}
&\left(-a_{15}X_2^2 -2a_{12}X_1X_2-2a_{8}X_1X_3 - 2a_1X_2X_3-a_{15}X_1^2X_2-a_{11}X_1^2X_3-2a_{14}X_1X_2^2 + \left(-4a_4-2a_{10}-2a_{13}-2a_{16}\right)X_1X_2X_3\right. \\
&\left.-2a_9X_1X_3^2+ \left(-2a_3-a_7\right)X_2^2X_3 -2a_2X_2X_3^2 \right)dX_1 \notag \\
&+ \left(a_{12}X_1^2 -a_{14}X_2^2+2a_1X_1X_3 +a_{15}X_1^3+2a_{14}X_1^2X_2+(2a_{4}+a_{10}+a_{13}+a_{16})X_1^2X_3 +2a_2X_1X_3^2 -a_6X_2^2X_3\right)dX_2 \notag \\
&+\left(-2a_8X_2 +a_8X_1^2+ \left(-a_{13}-a_{16}-2a_{10}\right)X_2^2 + \left(2a_1-2a_{11}\right)X_1X_2 -2a_9X_2X_3+a_{11}X_1^3 +a_6X_2^3 \right. \notag\\
&\left.+\left(2a_{4}+a_{10}+a_{13}+a_{16}\right)X_1^2X_2 +2a_9X_1^2X_3+\left(2a_3+a_7\right)X_1X_2^2\right)dX_3 \notag
\end{align}}}
We compute $H^0\left( \mathbb{P}_\mathbb{C}^3, \mathcal{E}_{\mathcal{O}_{\mathbb{P}_\mathbb{C}^3}(4)} \right)\to H^0\left( \mathbb{P}_\mathbb{C}^3, \mathcal{O}_{\mathbb{P}_\mathbb{C}^3}(4) \right)$. The image $i_T(w_0)$ is given by
{\Small{\begin{align}\label{ma6}
& -2a_1 X_1X_2X_3-2a_2X_1X_2X_3^2- \left(2a_3 + a_7 \right)X_1X_2^2 X_3+ \left(a_{10} -2a_4 + a_{13}\right) X_1^2 X_2 X_3 + a_8X_1^2 X_3+a_9X_1^2X_3^2+a_{11}X_1^3X_3 \\
&+a_{12} X_1^2 X_2+a_{14}X_1^2X_2^2+a_{15}X_1^3 X_2 -a_{12}X_2^2-a_{13}X_2^2X_3-a_{14}X_2^3-a_{15}X_1X_2^2-a_5X_2^2X_3^2-a_6X_2^3X_3 \notag
\end{align}}}
Then we see that
\begin{align*}
&\dim_\mathbb{C} \mathbb{H}^0\left( \mathbb{P}_\mathbb{C}^3,  \tilde{E}_{\mathcal{O}_{\mathbb{P}_\mathbb{C}^3}(4)}^\bullet \right) = \dim_\mathbb{C} \textnormal{ker}\left( H^0\left( \mathbb{P}_\mathbb{C}^3 , \mathcal{E}_{\mathcal{O}_{\mathbb{P}_\mathbb{C}^3}(4)}^\bullet \right) \to H^0\left( \mathbb{P}_\mathbb{C}^3, \Omega_{\mathbb{P}_\mathbb{C}^3}^1(4) \right) \bigoplus H^0\left( \mathbb{P}_\mathbb{C}^3 , \mathcal{O}_{\mathbb{P}_\mathbb{C}^3}(4) \right) \right) = 2 \\
&\dim_\mathbb{C} \textnormal{im}\left( H^0\left( \mathbb{P}_\mathbb{C}^3 , \mathcal{E}_{\mathcal{O}_{\mathbb{P}_\mathbb{C}^3}(4)}^\bullet \right) \to H^0\left( \mathbb{P}_\mathbb{C}^3, \Omega_{\mathbb{P}_\mathbb{C}^3}^1(4) \right) \bigoplus H^0\left( \mathbb{P}_\mathbb{C}^3 , \mathcal{O}_{\mathbb{P}_\mathbb{C}^3}(4) \right)  \right)  = 16- 2 = 14
\end{align*}
Then from $(\ref{ma5})$ we have
{\small{\begin{align*}
\dim_\mathbb{C} \mathbb{H}^1\left( \mathbb{P}_\mathbb{C}^3, \tilde{\mathcal{E}}_{\mathcal{O}_{\mathbb{P}_\mathbb{C}^3}(4)}^\bullet \right) = \dim_\mathbb{C} \frac{ \textnormal{ker}\left( H^0\left( \mathbb{P}_\mathbb{C}^3, \Omega_{\mathbb{P}_\mathbb{C}^3}^1(4) \right) \bigoplus H^0\left( \mathbb{P}_\mathbb{C}^3, \mathcal{O}_{\mathbb{P}_\mathbb{C}^3} (4) \right) \to H^0\left( \mathbb{P}_\mathbb{C}^3,  \mathcal{O}_{\mathbb{P}_\mathbb{C}^3}(4) \right) \bigoplus H^0\left( \mathbb{P}_\mathbb{C}^3 ,  \bigwedge^2 \Omega_{\mathbb{P}_\mathbb{C}^3}^1(8) \right) \right)}{\textnormal{im}\left( H^0\left( \mathbb{P}_\mathbb{C}^3 , \mathcal{E}_{\mathcal{O}_{\mathbb{P}_\mathbb{C}^3}(4)}^\bullet \right) \to H^0\left( \mathbb{P}_\mathbb{C}^3, \Omega_{\mathbb{P}_\mathbb{C}^3}^1(4) \right) \bigoplus H^0\left( \mathbb{P}_\mathbb{C}^3 , \mathcal{O}_{\mathbb{P}_\mathbb{C}^3}(4) \right)  \right)} = 20 - 14 = 6
\end{align*}}}
By considering $(\ref{ma7}),(\ref{ma8}),(\ref{ma3})$ and $(\ref{ma6})$, we see that
{\small{\begin{align*}
\mathbb{H}^1\left( \mathbb{P}_\mathbb{C}^3 , \tilde{\mathcal{E}}_{\mathcal{O}_{\mathbb{P}_\mathbb{C}^3}(4)}^\bullet \right) \cong \textnormal{Span}_\mathbb{C} & \left\langle   \left(2X_1 dX_1 - dX_2, X_2- X_1^2 \right) , \left( \left( X_2 + X_1^2 \right) dX_1 - X_1 dX_2, X_1 X_2- X_1^3 \right) , \left(- X_3 dX_2 - X_2 dX_3, X_2X_3 \right) , \right. \\
 &\left. \left( - 2 X_1X_2 dX_1 + X_1^2 dX_2, - X_1^2 X_2 + X_1^4 \right), \left( - X_3^2 dX_2 , X_2X_3^2\right)   , \left( - X_3^3 dX_2+ X_2X_3^2 dX_3, X_2X_3^3 \right) \right\rangle
\end{align*}}}

\subsection{Codimension $2$-Unfoldings of the projective $3$-space $\mathbb{P}_\mathbb{C}^3$ defined by locally free subsheaves}\

Let $\left( \mathbb{P}_\mathbb{C}^3, \mathcal{N}_{\mathcal{F}_0}^*\right)$ be a foliated projective $3$-space with $\mathcal{N}_{\mathcal{F}_0}^*$ locally free of rank $2$. Then the associated dimension $1$-foliation $\Theta_{\mathcal{F}_0}$ is of the form $\Theta_{\mathcal{F}_0}\cong  \mathcal{O}_{\mathbb{P}_\mathbb{C}^3}(-m), m\geq -1$. From Part IV-1, we have the unfolded dual leaf complex 
\begin{align}\label{cp6}
\mathcal{E}_{\mathcal{N}_{\mathcal{F}_0}^*}^\bullet : \mathcal{E}_{\mathcal{N}_{\mathcal{F}_0}^*} \to  \mathscr{H}om_{\mathcal{O}_{\mathbb{P}_\mathbb{C}^3} } \left( \mathcal{N}_{\mathcal{F}_0}^*, \Omega_{\mathbb{P}_\mathbb{C}^3}^1 \right) \bigoplus \mathscr{H}om_{\mathcal{O}_M}\left( \mathcal{N}_{\mathcal{F}_0}^*, \mathcal{O}_{\mathbb{P}_\mathbb{C}^3}\right)  \to \mathscr{H}om_{\mathcal{O}_{\mathbb{P}_\mathbb{C}^3}} \left(  \mathcal{N}_{\mathcal{F}_0}^*, \tilde{\mathcal{S}}^1 \right) \to 0
\end{align}
Since $\tilde{\mathcal{S}}^1 \cong \Theta_{\mathcal{F}_0}^*\cong\mathcal{O}_{\mathbb{P}_\mathbb{C}^3}(m)$, $(\ref{cp6})$ is equivalent to
\begin{align*}
\mathcal{E}_{\mathcal{N}_{\mathcal{F}_0}^*}^\bullet : \mathcal{E}_{\mathcal{N}_{\mathcal{F}_0}^*} \to  \mathscr{H}om_{\mathcal{O}_{\mathbb{P}_\mathbb{C}^3} } \left( \mathcal{N}_{\mathcal{F}_0}^*, \Omega_{\mathbb{P}_\mathbb{C}^3}^1 \right) \bigoplus \mathscr{H}om_{\mathcal{O}_M}\left( \mathcal{N}_{\mathcal{F}_0}^*, \mathcal{O}_{\mathbb{P}_\mathbb{C}^3}\right)  \to \mathscr{H}om_{\mathcal{O}_{\mathbb{P}_\mathbb{C}^3}} \left(  \mathcal{N}_{\mathcal{F}_0}^*, \mathcal{O}_{\mathbb{P}_\mathbb{C}^3}(m) \right) \to 0
\end{align*}

\section{Foliated complex analytic families of non-degenerate invariant hypersurfaces}\label{hs1}

We will extend the discussion of complex analytic families of hypersurfaces in \cite{Kod05} section (c) p.234 in the context of holomorphic foliated deformations in terms of tangent sheaves. This section can be considered as a preview of deformations of (a special kind of) compact invariant complex submanifolds (which we call non-degenerate) of a foliated complex manifold $\left(W, \Theta_{\mathcal{F}} \right)$ with $\Theta_{\mathcal{F}}$ locally free, which will be studied more generally in Part VI.

Let $\left( W, \Theta_{\mathcal{F}} \right)$ be a foliated complex manifold with $\Theta_{\mathcal{F}}$ locally free of rank $p$. Let $S$ be an invariant complex submanifold of $\left(W, \Theta_{\mathcal{F}} \right)$. We shall say that $S$ is a non-degenerate invariant complex submanifold of $\left(W, \Theta_{\mathcal{F}} \right)$ if the exact sequence $0\to \Theta_{\mathcal{F}}\to \Theta_W \to \frac{\Theta_W}{\Theta_{\mathcal{F}}} \to 0$ induces the exact sequence $0\to \Theta_{\mathcal{F}}|_S\to \Theta_W|_S \to \frac{\Theta_W}{\Theta_{\mathcal{F}}}|_S\to 0$ by restriction to $S$. Then $\Theta_{\mathcal{F}_{S}}:=\Theta_{\mathcal{F}}|_S\subset \Theta_S$ is a locally free sheaf of rank $p$. If we denote the singular locus of $\mathcal{F}$ by $\textnormal{Sing}\left( \mathcal{F}\right)$, and $S\cap \textnormal{Sing}\left( \mathcal{F} \right)$ has condimension at least $2$ in $S$, then $\Theta_{\mathcal{F}_{S}}$ defines a foliation on $S$.

Let $\left(W, \Theta_{\mathcal{F}}\right)$ be a given compact foliated complex manifold of dimension $n+1$. A foliated complex analytic family $\left(\mathcal{M}, \Theta_\mathcal{G}, B, \omega \right)$ of compact complex manifolds in terms of tangent sheaves is called a foliated analytic family of invariant hypersurfaces in $\left(W, \Theta_{\mathcal{F}} \right)$ if it satisfies the following conditions:
\begin{enumerate}
\item Each $\left(M_t, \Theta_{\mathcal{G}_t}\right)=\omega^{-1}(t)$ is a non-degenerate invariant complex submanifold of codimension $1$ of $\left(W, \Theta_\mathcal{F} \right)$ for $t\in B$ $($we do not necessarily assume that the codimenstion of the singular locus of $\Theta_{\mathcal{G}_t}$ is at least 2$)$.
\item There is a holomorphic map $\Phi$ of $\left(\mathcal{M}, \Theta_\mathcal{G}\right)$ into $\left(W, \Theta_\mathcal{F} \right)$ whose restriction to $M_t$ coincides with the inclusion $\left(M_t, \Theta_{\mathcal{G}_t} \right) \hookrightarrow \left( W, \Theta_\mathcal{F}\right)$.
\end{enumerate}
Since we are only interested in the infinitesimal deformations, it suffices to take a sufficient small polydisk $\Delta\subset B$, and consider $\left(\mathcal{M}_\Delta, \Theta_{\mathcal{G}_\Delta} , \Delta, \omega\right)$ instead of $\left(\mathcal{M}, \Theta_\mathcal{G}, B, \omega \right)$ where we put $(\mathcal{M}_\Delta, \Theta_{\mathcal{G}_\Delta})=\omega^{-1}(\Delta)$. Moreover, we may assume that
\begin{align*}
\left(\mathcal{M}_\Delta, \Theta_{\mathcal{G}_\Delta}  \right)=\bigcup_{j=1}^l \left(U_j\times \Delta, \Theta_{\mathcal{G}}|_{U_j\times \Delta} \right)
\end{align*}
where $(z_j, t)\in U_j\times \Delta$ and $(z_k,t)\in U_k\times \Delta$ are the same point of $\mathcal{M}_\Delta$  if $z_j^\alpha=f_{jk}^{\alpha}(z_k,t)$, and $\Gamma\left(U_j\times \Delta, \Theta_\mathcal{G}\right)$ is generated by $V_j^1,..., V_j^p\in \Gamma\left(U_j\times\Delta, \Theta_{\mathcal{M}/B}\right)$ of the form
\begin{align*}
V_j^\alpha(z_j,t)= \sum_{\beta=1}^n V_j^{\alpha\beta}\left(z_j,t \right)\frac{\partial}{\partial z_j^\beta}
\end{align*}
Let $\left\{W_1,..., W_l,...\right\}$ be a finite open covering of $W$, where $W_j$ is a coordinate neighborhood, and $w_j=\left(w_j^0,w_j^1,..., w_j^n\right)$ be the local coordinates defined on $W_j$, and $\Gamma\left(W_j,\Theta_\mathcal{F}\right)$ is generated by $T_j^1,..., T_j^p\in \Gamma\left(W_j, \Theta_W\right)$ of the form
\begin{align}\label{h1}
T_j^\alpha(w_j)=\sum_{\beta=0}^n T_j^{\alpha\beta}\left(w_j^0,w_j^1,..., w_j^n\right)\frac{\partial}{\partial w_j^\beta}
\end{align}
Taking $U_j$ and $\Delta$ sufficiently small, we may assume that $\Phi\left(U_j\times \Delta \right)\subset W_j$ for $j=1,...,l$, and that $\Phi\left( \mathcal{M}_\Delta  \right)\cap \left[W_j\right]=\emptyset$ for $j\geq l+1$. $\Phi$ is represented on each $U_j\times \Delta$ as
\begin{align*}
\Phi:(z_j,t)\to \left(w_j^0,w_j^1,..., w_j^n\right)=\left(\varphi_j^0(z_j,t),..., \varphi_j^n(z_j,t)\right).
\end{align*}

We denote a point of $\mathcal{M}$ by $p$, and define the map $\tilde{\Phi}$ of $\left(\mathcal{M}, \Theta_\mathcal{G}\right)$ into $\left(W\times B, \omega^* \Theta_\mathcal{F}\right)$ by
\begin{align*}
\tilde{\Phi}:p\to \tilde{\Phi}(p)=\left(\Phi(p), \omega(p)\right).
\end{align*}

Then $\tilde{\Phi}$ is represented in terms of local coordinates as 
\begin{align*}
(z_j,t)\mapsto \left(w_j^0,..., w_j^n,t\right)=\left(\varphi_j^0(z_j,t), ..., \varphi_j^n(z_j,t),t_1,..., t_m\right)
\end{align*}
So $\tilde{\Phi}$ is a biholomorphic map of $(\mathcal{M}_\Delta, \Theta_{\mathcal{G}_\Delta})$ into $(W\times \Delta, \pi^* \Theta_\mathcal{F})$, and  identifying $\mathcal{M}_\Delta$ with $\tilde{\Phi}(\mathcal{M}_\Delta)$, we may consider that $\left(\mathcal{M}_\Delta, \Theta_{\mathcal{G}}|_\Delta \right)$ is an invariant submanifold of codimension $1$ of $(W\times \Delta, \pi^* \Theta_\mathcal{F})$. Then
\begin{align*}
\left(\mathcal{M}_\Delta, \Theta_\mathcal{G}|_\Delta \right)=\bigcup_{t\in \Delta} \left(M_t\times t, \Theta_{\mathcal{G}_t}\right)\subset \left(W\times \Delta, \pi^* \Theta_\mathcal{F}\right)
\end{align*}
is considered as the graph of $\omega^{-1}:t\to M_t$. Since $\Theta_\mathcal{G}|_\Delta=\tilde{\Phi}^* \left( \pi^* \Theta_\mathcal{F} \right)$, we may assume that $T_j^\alpha(w_j)$ induces $V_j^\alpha(z_j,t)$ via $\tilde{\Phi}$, so that we have
\begin{align}\label{h6}
T_j^{\alpha\gamma}\left(\Phi(z_j,t)\right)=\sum_{\beta=1}^n V_j^{\alpha\beta}(z_j,t)\frac{\partial \varphi_j^\gamma}{\partial z_j^\beta},\,\,\,\,\,\,\,\,\alpha=1,...,p,\,\,\,\gamma=0,...,n
\end{align}

We have an exact sequence $0 \to \Theta_\mathcal{F} \to \Theta_W \to \Theta_W/\Theta_\mathcal{F}\to 0$. Then we have an exact sequence $0\to \pi^* \Theta_\mathcal{F} \to \Theta_{W\times \Delta/\Delta} \to \frac{\Theta_{W\times \Delta/\Delta}}{\pi^* \Theta_\mathcal{F}}\to 0 $ where $\pi:W\times \Delta\to \Delta$ is the natural projection. Since for each $t\in \Delta$,  we have $0\to \left(\pi^*\Theta_\mathcal{F}\right)_t= \Theta_\mathcal{F} \to \Theta_W \to \Theta_W/\Theta_\mathcal{F}\to 0$, so that $ \frac{\Theta_{W\times \Delta/\Delta}}{\pi^* \Theta_\mathcal{F}}$ is flat over $\Delta$. In other words, $\left(W\times \Delta, \pi^*\Theta_\mathcal{F}, \pi, \Delta\right)$ is the trivial foliated analytic family. On the other hand, $\Phi:\left(\mathcal{M}_\Delta, \Theta_{\mathcal{G}}|_\Delta\right)\to \left(W\times \Delta, \pi^* \Theta_\mathcal{F}\right)$ induces non-degenerate $\left(M_t, \Theta_{\mathcal{G}_t}\right)\hookrightarrow \left(W,\Theta_{\mathcal{F}}\right)$ with $\Theta_{\mathcal{G}_t}=\Theta_\mathcal{F}|_{M_t}$, and we have an exact sequence $0\to \Theta_{\mathcal{G}_t}\to  \Theta_W \to \frac{\Theta_W}{\Theta_{\mathcal{G}_t}}\to 0$. This implies that the restriction $\frac{\Theta_{W\times\Delta/\Delta}}{\pi^* \Theta_\mathcal{F}}|_{\mathcal{M}}:= \frac{\Theta_{W\times\Delta/\Delta}}{\pi^* \Theta_\mathcal{F}} |_{\Phi\left(\mathcal{M}\right)}$ is also flat over $\Delta$.

If we take $W_j$ and $\Delta$ sufficiently small, $\mathcal{M}_\Delta$ is given by on each $W_j\times \Delta$ by a single minimal equation
\begin{align*}
S_j(w_j,t)=0
\end{align*}

Here $S_j(w_j,t)$ is understood to be a non-vanishing holomorphic function if $W_j\times \Delta \cap \mathcal{M}_\Delta = \emptyset$. For a fixed $t$, $S_j(w_j,t)=0$ gives a minimal equation of $M_t$ in $W_j$.  $S_j(w_j,t)=S_j\left(w_j^0,..., w_j^n,t\right)$ is a holomorphic function of $w_j^0,..., w_j^n,t_1,...,t_m$ and at least one of the partial derivatives $\frac{\partial S_j(w_j,t)}{\partial w_j^\lambda},\lambda=1,...,n$, does not vanish for each $w_j\in M_t$. On $w_j\times \Delta \cap W_k \cap \Delta \ne \emptyset$, we have $S_j(w_j,t)=F_{jk}(w_k,t)\cdot S_k(w_k,t)$, $F_{jk}(w,t)$ is a non-vanishing holomorphic function defined on $\left( W_j \cap W_k \right)\times \Delta$, and $w_j,w_k$ are local coordinates of $w\in W_j\cap W_k$. Then
\begin{align}\label{h2}
F_{ik}(w,t)= F_{ij}(w,t) \cdot F_{jk}(w,t) \,\,\,\,\,\,\,\textnormal{for}\,\,\, w\in W_i \cap W_j \cap W_k
\end{align}
On the other hand, $\Theta_\mathcal{F}$ is generated by $T_j^1(w_j),..., T_j^p(w_j)\in \Gamma \left(W_j, \Theta_W \right)$ as in $(\ref{h1})$. On $W_j\cap W_k\ne \emptyset$, we have the following relation $T_j^\alpha(w_j)=\sum_{\beta=1}^p r_{jk}^{\alpha\beta}(w_k)T_k^\beta(w_k)$ for $r_{jk}^{\alpha\beta}(w_k)\in \Gamma(W_j\cap W_k, \mathcal{O}_W)$, and on $W_j$, we have the following relation $\left[ T_j^\alpha(w_j), T_j^\beta(w_j ) \right] = \sum_{\gamma=1}^p g_{j\alpha\beta}^\gamma(w_j) T_j^\gamma(w_j)$ for $g_{j\alpha\beta}^\gamma(w_j)\in \Gamma\left( W_j, \mathcal{O}_W \right)$. Since $S_j(w_j,t)$ defines an invariant submanifold, we have
\begin{align*}
\left[T_j^\alpha(w_j), S_j(w_j,t)    \right] = h_j^\alpha(w_j,t)S_j(w_j,t) , \,\,\,\,\,\,\,\,\,\, h_j^\alpha(w_j,t)\in \Gamma(U_j\times \Delta, \mathcal{O}_{W\times \Delta})
\end{align*}
We can write 
\begin{align*}
h_j^\alpha(w_j,t)=\left[ T_j^\alpha(w_j), \log S_j(w_j,t)   \right]
\end{align*}
We note that
\begin{align} \label{h3}
\sum_{\beta=1}^p r_{jk}^{\alpha\beta} h_k^\beta - h_j^\alpha +\left[ T_j^\alpha, \log F_{jk}   \right]&= \sum_{\beta=1}^p r_{jk}^{\alpha\beta} \left[T_k^\beta, \log S_k   \right] - \left[ T_j^\alpha, \log S_j \right] + \left[ T_j^\alpha, \log S_j \right] - \left[ T_j^\alpha , \log S_k \right]=0
\end{align}
and we note that
\begin{align}\label{h4}
\left[ T_j^\alpha, h_j^\beta \right] - \left[ T_j^\beta, h_j^\alpha \right] - \sum_{\gamma=1}^p g_{j\alpha\beta}^\gamma  h_j^\gamma = \left[ T_j^\alpha, \left[ T_j^\beta, \log S_j\right] \right] - \left[ T_j^\beta,\left[ T_j^\alpha, \log S_j \right] \right] - \sum_{\gamma=1}^p g_{j\alpha\beta}^\gamma \left[ T_j^\gamma, \log S_j \right]=0
\end{align}

We define $h_j \in \Gamma\left(W_j\times \Delta, \mathscr{H}om_{\mathcal{O}_\mathcal{M}}\left( \pi^*\Theta_{\mathcal{F}}, \mathcal{O}_{W\times \Delta}   \right) \right)$ in the following way:
\begin{align*}
h_j:\Gamma\left(W_j\times \Delta, \pi^* \Theta_\mathcal{F}\right) &\to \Gamma\left( W_j\times \Delta, \mathcal{O}_{W\times \Delta} \right)\\ 
                                                T_j^\alpha(w_j) &\mapsto h_j^\alpha=\left[ T_j^\alpha(w_j)  , \log S_j(w_j,t)   \right]
\end{align*}
Then $(\ref{h2}),(\ref{h3}),(\ref{h4})$ imply that $\left( \left\{F_{jk}\right\}, \left\{h_j \right\}\right)$ defines a foliated line bundle on $\left(W\times \Delta, \pi^* \Theta_\mathcal{F}\right)$ over $\Delta$ in the sense of Part V. Hence for every $t\in \Delta$ we have a foliated complex line bundle $\left(\mathscr{F}_t , \nabla_t\right)$ over $\left(W, \Theta_\mathcal{F}\right)$ defined by the foliated divisor $\left(M_t, \Theta_{\mathcal{G}_t}\right)$ on $\left( W, \Theta_\mathcal{F}\right)$ (see Part VI for foliated divisors): $\mathscr{F}_t=[M_t]$, and $S(t)=\{S_j(w_j,t)\}$ is a holomorphic section of $\mathscr{F}_t$ over $W$ such that $S(t)=\{S_j(w_j,t)\} \in \mathbb{H}^0(W, \mathscr{F}_t^{\bullet})$ where 
\begin{align*}
\mathscr{F}_t^\bullet: \mathcal{O}(\mathscr{F}_t) \xrightarrow{\nabla_t^0} \mathscr{H}om_{\mathcal{O}_W} \left( \Theta_\mathcal{F}, \mathcal{O}(\mathscr{F}_t) \right) \xrightarrow{\nabla_t^1} \mathscr{H}om_{\mathcal{O}_W}\left( \bigwedge^2 \Theta_\mathcal{F}, \mathcal{O}(\mathscr{F}_t) \right) \xrightarrow{\nabla_t^2} \cdots
\end{align*}
In fact, we note that
\begin{align*}
\nabla_t^0(S_j(t))\left(T_j^\alpha \right)=\left[ T_j^\alpha, S_j(t) \right]-h_j^\alpha(t) S_j(t)=0
\end{align*}

Fix $t\in \Delta$. Since $\left(M_t, \Theta_{\mathcal{G}_t}\right)=\bigcup_j \left(U_j \times t , \Theta_{\mathcal{G}_t}\right) $, if we put $U_j\times t= U_{j_t}$, $U_{jt}$ is a coordinate polydisk on $M_t$. We denote the point $z_j\times t\in U_j\times t=U_{jt}$ simply by $z_j$, and consider $z_j=\left(z_j^1,..., z_j^n \right)$ as the local coordinates on $M_t$. If we put $\varphi_j(z_j,t)=\left(\varphi_j^0(z_j,t),..., \varphi_j^n(z_j,t)\right)$, the inclusion $\Phi:M_t\to W$ is represented on each $U_{jt}\subset M_t$ as $\Phi:z_j\to w_j=\varphi_j(z_j,t)$. We put 
\begin{align*}
F_{tjk}(z)=F_{jk}(\Phi(z),t), \,\,\,\,\,\,\,\,\,\,z\in U_{jt} \cap U_{kt}
\end{align*}
Then $F_{tjk}(z)$ is a non-vanishing holomorphic function on $U_{jt}\cap U_{kt}$, and we obtain from $(\ref{h2})$
\begin{align} \label{h7}
F_{tjk}(z)=F_{tij}(z)\cdot F_{tjk}(z),\,\,\,\,\,\,\,z\in U_{ti}\cap U_{tj}\cap U_{tk}
\end{align}
On the other hand, we note that from $(\ref{h6})$
\begin{align*}
&\left[\sum_{\beta=1}^n V_j^{\alpha\beta}(z_j,t)\frac{\partial}{\partial z_j^\beta}, \log F_{jk}(\varphi_j(z_j,t),t)   \right]=\frac{1}{F_{jk}(\varphi_j(z_j,t),t)}\left( \sum_{\beta=1}^n \sum_{\gamma=0}^n V_j^{\alpha\beta}(z_j,t)\frac{\partial F_{jk}}{\partial w_j^\gamma}|_{w_j=\varphi_j(z_j,t)} \cdot \frac{\partial \varphi_j^\gamma}{\partial z_j^\beta}   \right) \\
&= \frac{1}{F_{jk}(\varphi_j(z_j,t),t)}\sum_{\gamma=0}^n T_j^{\alpha\gamma}(\varphi_j(z_j,t))\frac{\partial F_{jk}}{\partial w_j^\gamma}|_{w_j=\varphi_j(z_j,t)}= \left[\sum_{\gamma=0}^n T_j^{\alpha\gamma}(w_j), \log F_{jk}(w_j,t)    \right]|_{w_j=\varphi_j(z_j,t)}
\end{align*}

Hence we have from $(\ref{h3})$
\begin{align}\label{h8}
\sum_{\beta=1}^p r_{jk}^{\alpha\beta}(\varphi_k(z_j,t))h_k^\alpha(\varphi_k(z_k,t),t)- h_j^\alpha(\varphi_j(z_j,t),t) + \left[V_j^\alpha(z_j,t), \log F_{tjk}(z) \right]=0
\end{align}
and we note that from $(\ref{h6})$
\begin{align}
\left[  \sum_{\beta=1}^n \sum_{\delta=0}^n V_j^{\alpha\beta}(z_j,t)\frac{\partial}{\partial z_j^\beta}, h_j^{\gamma}\left(\varphi_j(z_j,t),t) \right)\right] = \sum_{\beta=1}^n V_j^{\alpha\beta}(z_j,t) \frac{\partial \varphi_j^\delta}{\partial z_j^\beta}  \cdot \frac{\partial h_j^\gamma}{\partial w_j^\delta}|_{w_j= \varphi_j(z_j,t)} =\left[ T_j^\alpha, h_j^\gamma \right] |_{w_j= \varphi_j(z_j,t)}
\end{align}

Hence we have from $(\ref{h4})$
\begin{align}\label{h9}
\left[ V_j^\alpha, h_j^\beta\left( \varphi_j(z_j,t), t \right) \right] -\left[ V_j^\beta, h_j^\alpha\left( \varphi_j(z_j,t),t \right) \right] - \sum_{\gamma=1}^p g_{j\alpha\beta}^\gamma\left( \varphi_j(z_j,t) \right) h_j^\gamma\left( \varphi_j(z_j,t), t \right) = 0
\end{align}

We define an element $h_{jt} \in \Gamma\left( U_{jt}, \mathscr{H}om_{\mathcal{O}_{M_t}} \left( \Theta_{\mathcal{G}_t}, \mathcal{O}_{M_t} \right)\right)$ by
\begin{align*}
h_{jt} :\Gamma\left( U_{jt}, \Theta_{\mathcal{G}_t} \right) &\to \Gamma\left( U_{jt}, \mathcal{O}_{M_t} \right) \\
 V_j^\alpha (t) &\mapsto h_{jt}^\alpha=h_j^\alpha\left(\varphi_j(z_j,t), t \right)
\end{align*}

Then $(\ref{h7}),(\ref{h8}),(\ref{h9})$ imply that the restriction $\left( \left\{F_{tjk}(z_k,t) \right\}, \left\{h_{jt}=h_j(\Phi(z),t) \right\} \right)$ defines a foliated line bundle $\left(F_t, \nabla_{M_t}\right)$ where $\left(F_t, \nabla_{M_t}\right)$ is the restriction of $\left(\mathscr{F}_t, \nabla_t\right)$ to $M_t$. Then we have the complex of sheaves (see Part V)
\begin{align} \label{cs8}
\mathcal{O}(F_t)^\bullet: \mathcal{O}(F_t)\xrightarrow{\nabla_{M_t}^0} \mathscr{H}om_{\mathcal{O}_{M_t}}\left( \Theta_{\mathcal{G}_t}, \mathcal{O}(F_t) \right) \xrightarrow{\nabla_{M_t}^1} \mathscr{H}om_{\mathcal{O}_{M_t} }\left( \bigwedge^2 \Theta_{\mathcal{G}_t}, \mathcal{O}(F_t) \right) \xrightarrow{\nabla_{M_t}^2} \cdots
\end{align}
We will denote the $i$-th cohomology group of $\mathcal{O}(F_t)^\bullet$ by $\mathbb{H}^i\left( M_t, \mathcal{O}(F_t)^\bullet \right)$.

In general when $t$ moves, the foliated complex structure of $(M_t, \Theta_{\mathcal{G}_t})$ as well as the position of $(M_t, \Theta_{\mathcal{G}_t})$ in $(W, \Theta_\mathcal{F})$ varies. The change of the position of $\left(M_t, \Theta_{\mathcal{G}_t}\right)$ is called a displacement. On each neighborhood $W_j$, $\left(M_t, \Theta_{\mathcal{G}_t}\right)$ is given by the equation $S_j(w_j,t)=0$ and $\Theta_{\mathcal{G}_t}$ on $W_j\cap M_t$ is generated by the restriction $T_j^\alpha(w_j)|_{M_t}$. In order to define the infinitesimal displacement along $\frac{\partial}{\partial t}$ for an arbitrary $\frac{\partial}{\partial t}\in T_t(\Delta)$, consider the partial derivative $\frac{\partial S_j(w_j,t)}{\partial t}$, and put
\begin{align}\label{h11}
\sigma_j(z_j,t)= - \left( \frac{\partial S_j(w_j,t)}{\partial t} \right)_{w_j=\varphi_j(z_j,t)}
\end{align}

Applying $\frac{\partial}{\partial t}$ to $S_j(w_j,t)=F_{jk}(w_j,t) S_k(w_k,t)$, we have
\begin{align*}
\frac{\partial S_j(w_j,t)}{\partial t}= F_{jk}(w,t)\frac{\partial S_k(w_k,t)}{\partial t} + \frac{\partial F_{jk}(w,t)}{\partial t} \cdot S_k(w_k,t)
\end{align*}
Restricting these equalities to $M_t$ we obtain
\begin{align}\label{h10}
\sigma_j(z_j,t)=F_{tjk}(z)\sigma_k(z_k,t)
\end{align}
 since $S_k(w_k,t)=0$ on $M_t$. On the other hand, applying $\frac{\partial}{\partial t}$ to $[T_j^\alpha(w_j), S_j(w_j,t)]=h_j(w_j,t)S_j(w_j,t)$, we have
 \begin{align*}
 \left[ T_j^\alpha(w_j), \frac{\partial S_j(w_j,t)}{\partial t}    \right] = \frac{\partial h_j(w_j,t)}{\partial t} S_j(w_j,t)+ h_j(w_j,t)\frac{\partial S_j(w_j,t)}{\partial t}
 \end{align*}
Restricting these equalities to $M_t$, we obtain
\begin{align*}
\left[\sum_{\beta=1}^n V_j^{\alpha\beta}(z_j,t)\frac{\partial}{\partial z_j^\beta}, \sigma_j(z_j,t)  \right]- h_j \left(\varphi_j(z_j,t),t)\right)\cdot \sigma_j(z_j,t)=0
\end{align*}
This means that $\sigma(t)=\left\{\sigma_j(z_j,t)\right\}\in \mathbb{H}^0(M_t, \mathcal{O}(F_t)^\bullet)$.

\begin{definition}
$\sigma(t)\in \mathbb{H}^0(M_t, \mathcal{O}(F_t))$ is called the infinitesimal displacement of the non-degenerate invariant submanifold $(M_t, \Theta_{\mathcal{G}_t})$ of codimension $1$ of $(W, \Theta_\mathcal{F})$, and is denoted by
\begin{align*}
\frac{\partial \left(M_t, \Theta_{\mathcal{G}_t}\right)\subset \left(W, \Theta_\mathcal{F}\right)}{\partial t}
\end{align*}
\end{definition}

We denote by $\rho_{d,t}$ the linear map $\frac{\partial}{\partial t} \to \sigma(t)$ of $T_t(\Delta)$ to $\mathbb{H}^0\left(M_t, \mathcal{O}\left(F_t\right)^\bullet\right)$, namely
\begin{align}\label{cs9}
\rho_{d,t}:\frac{\partial}{\partial t} \mapsto \rho_{d,t}\left( \frac{\partial}{\partial t} \right)=\frac{\partial \left(M_t, \Theta_{\mathcal{G}_t}\right)\subset \left(W, \Theta_{\mathcal{F}} \right )}{\partial t}=\sigma(t)
\end{align}

Next we will consider the relation between the infinitesimal foliated deformation of $\left(M_t, \Theta_{\mathcal{G}_t} \right)$ and the infinitesimal displacement of $\left( M_t, \Theta_{\mathcal{G}_t} \right)$ in $\left( W, \Theta_{\mathcal{F}} \right)$. Let $T(W)_t$ be the restriction of the tangent bundle $T(W)$ of $W$ to $M_t$. Let $T\left( M_t\right)$ be a holomorphic vector bundle over $M_t$. Since the tangent space $T_z(M_t)$ at $z\in M_t$ consists of all tangent vectors of $W$ at $z$ which are tangent to $M_t$,  the tangent bundle $T(M_z)= \bigcup_{z\in M_t} T_z( M_t)  $ of $M_t$ is a subbundle of $T(W)_t=\bigcup_{z \in M_t} T_z(W)$. A local holomorphic section $\xi$ of $T(W)_t$ is represented in terms of local coordinates as $\xi=\sum_{\lambda=0}^n v_j^\lambda(z)\frac{\partial}{\partial w_j^\lambda}$ where the coefficients $v_j^\lambda(z)$ are local holomorphic functions on $M_t$. Define $\tau_j(z_j,t)$ for this $\xi$ by 
\begin{align*}
\tau_j(z_j,t)&= \xi S_j(w_j,t)_{w_j=\varphi_j(z_j,t)}=\sum_{\lambda=0}^n v_j^\lambda(z)\left( \frac{\partial S_j(w_j,t)}{\partial w_j^\lambda} \right)_{w_j=\varphi_j(z_j,t)}   
\end{align*}

Then $\tau_j(z_j,t)$ is a holomorphic function on $U_{jt}\subset M_t$, and we obtain
\begin{align*}
\tau_j(z_j,t)= F_{tjk}(z)\tau_k(z_k,t)
\end{align*}
on $U_{jt}\cap U_{kt}\ne \emptyset$. Therefore $\tau(t)=\{\tau_j(z_j,t)\}$ is a holomorphic section of $F_t$ over $M_t$, which we denote by $\xi S(t)$:
\begin{align*}
\xi S(t)=\tau(t).
\end{align*}
$\xi S(t)=0$, namely, $\xi S_j(w_j,t)=0$ means that $\xi=\sum_{\lambda=0}^n v_j^\lambda(z)\frac{\partial}{\partial w_j^\lambda}\in T_z(M_t)$, namely, that $\xi$ is a holomorphic vector field on $M_t$. Hence if we denote by $\Xi_t$ the sheaf of germs of holomorphic sections of $T(W)_t:$ $\Xi_t=\mathcal{O}(T(W)_t)$. $\Theta_t=\mathcal{O}(T(M_t))$ is a subsheaf of $\Xi_t$, then $\xi\mapsto \tau(t)=\xi S(t)$ defines a homomorphism $\Xi_t\to \mathcal{O}(F_t)$ of sheaves, whose kernel is $\Theta_t$, so that we have the exact sequence
\begin{align*}
 0 \to \Theta_t \to \Xi_t \to \mathcal{O}(F_t)\to 0
\end{align*}
 Moreover we note that if we take $\xi=\sum_{\beta=0}^n T_j^{\alpha\beta}(\varphi_j(z_j,t)) \frac{\partial}{\partial w_j^\beta} $, then we have
\begin{align*}
\xi S_j(w_j,t)|_{w_j=\varphi_j(z_j,t)}= 0 
\end{align*}
Hence $\xi\to \tau(t)=\xi S(t)$ induces a homomorphism $\Xi_t / \Theta_{\mathcal{G}_t} \to \mathcal{O}(F_t)$ of sheaves, whose kernel is $\Theta_t/\Theta_{\mathcal{G}_t}$, so that we have the exact sequence
\begin{align*}
0\to \frac{\Theta_t}{ \Theta_{\mathcal{G}_t} } \to \frac{\Xi_t}{\Theta_{\mathcal{G}_t}}  \to \mathcal{O}(F_t) \to 0
\end{align*}

Then we have the exact sequence of complexes of sheaves (see Part VI for general case)
{\small{\begin{center}
$\begin{CD}
@. \cdots @. \cdots @. \cdots \\
@. @AD_{t_2}AA @AD_{W_{t_2}}AA  @A\nabla_t^2 AA \\
0@>>> \mathscr{H}om_{\mathcal{O}_{M_t}}\left( \bigwedge^2 \Theta_{\mathcal{G}_t}, \frac{\Theta_{M_t}}{\Theta_{\mathcal{G}_t}} \right) @>>> \mathscr{H}om_{\mathcal{O}_{M_t}}\left( \bigwedge^2 \Theta_{\mathcal{G}_t}, \frac{\Xi_t}{\Theta_{\mathcal{G}_t } } \right) @>>>  \mathscr{H}om_{\mathcal{O}_{M_t}}\left( \bigwedge^2 \Theta_{\mathcal{G}_t}, \mathcal{O}(F_t) \right) @>>> 0 \\
@. @AD_{t_1}AA @AD_{W_{t_1}}AA @A\nabla_t^1AA \\
0@>>> \mathscr{H}om_{\mathcal{O}_{M_t}}\left( \Theta_{\mathcal{G}_t},  \frac{\Theta_{M_t}}{\Theta_{\mathcal{G}_t}} \right) @>>>   \mathscr{H}om_{\mathcal{O}_{M_t}}\left(\Theta_{\mathcal{G}_t} , \frac{ \Xi_t  }{\Theta_{\mathcal{G}_t}} \right) @>>> \mathscr{H}om_{\mathcal{O}_{M_t}}\left( \Theta_{\mathcal{G}_t}, \mathcal{O}({F}_t)  \right) @>>> 0     \\
@. @AD_{t_0}AA @AD_{W_{t_0}}AA @A\nabla_t^0AA \\
0 @>>> \Theta_t @>>> \Xi_t @>>> \mathcal{O}({F}_t) @>>> 0
\end{CD}$
\end{center}}}
which induces the long exact sequence of cohomology groups
\begin{align*}
\to \cdots \mathbb{H}^0(M_t, \Xi_t^\bullet) \to \mathbb{H}^0(M_t, \mathcal{O}(\mathscr{F}_t)^\bullet) \xrightarrow{\delta^*} \mathbb{H}^1(M_t, \Theta_{\mathcal{G}_t}^\bullet) \to \mathbb{H}^1 (M_t, \Xi_t^\bullet) \to \cdots
\end{align*}
Here the middle complex of sheaves is the restriction of the leaf complex $\Theta_{\mathcal{F}}^\bullet:\Theta_W \xrightarrow{D_0} \mathscr{H}om_{\mathcal{O}_W}\left( \Theta_\mathcal{F},  \frac{\Theta_W}{\Theta_\mathcal{F}} \right) \xrightarrow{D_1}\mathscr{H}om_{\mathcal{O}_W}\left( \bigwedge^2 \Theta_\mathcal{F}, \frac{\Theta_W}{\Theta_{\mathcal{F}}} \right) \xrightarrow{D_2} \mathscr{H}om_{\mathcal{O}_W}\left( \bigwedge^3 \Theta_\mathcal{F} , \frac{\Theta_W}{\Theta_\mathcal{F}}  \right) \xrightarrow{D_3} \cdots $ to $M_t$.

\begin{lemma}\label{h12}
The relation between the infinitesimal foliated deformation of $\left( M_t, \Theta_{\mathcal{G}_t} \right)$ and the infinitesimal displacement of $\left( M_t, \Theta_{\mathcal{G}_t} \right)$ in $\left(W, \Theta_\mathcal{F} \right)$ is given by the following
\begin{align*}
\rho_t=\delta^* \rho_{d,t}
\end{align*}
where $\rho_t : T_t\left( \Delta\right)\to \mathbb{H}^1\left( M_t, \Theta_{\mathcal{G}_t}^\bullet \right)$ is the foliated Kodaira-Spencer map at $t$ in terms of tangent sheaves.
\end{lemma}

\begin{proof}
Since $\mathcal{U}_t= \left\{U_{jt} \right\}$ is a finite covering of $M_t$. and each $U_{jt}$ is a polydisk, we see from that $\mathbb{H}^1\left(M_t, \Theta_{\mathcal{G}_t}^\bullet\right)= \mathbb{H}^1\left(\mathcal{U}_t, \Theta_{\mathcal{G}_t}^\bullet\right)$ in terms of \v Cech resolution of $\Theta_{\mathcal{G}_t}^\bullet$ as in Part I. For any $\frac{\partial}{\partial t}\in T_t\left( \Delta \right)$, $\rho_t\left( \frac{\partial}{\partial t} \right)$ is the cohomology class of $1$-cocycle $\left( \left\{ \lambda_j(t) \right\}  , \left\{ \theta_{jk}(t)\right\}\right)\in C^0\left( \mathcal{U}_t, \mathscr{H}om_{\mathcal{O}_{M_t}}\left(\Theta_{\mathcal{G}_t}, \frac{\Theta_{M_t}}{\Theta_{\mathcal{G}_t}} \right)\right)\bigoplus C^1\left( \mathcal{U}_t, \Theta_{M_t} \right)$, where $\lambda_j(t) : V_j^\alpha \mapsto - \overline{\frac{\partial V_j^\alpha}{\partial t}}$, and $\theta_{jk}(t)= \sum_{\alpha=1}^n \frac{\partial f_{jk}^\alpha(z_k,t)}{\partial t} \frac{\partial}{\partial z_j^\alpha}$ (see Part I). On the other hand, the infinitesimal displacement $\rho_{d,t}\left( \frac{\partial}{\partial t} \right)$ is given by $\sigma(t)=\left\{ \sigma_j(z_j,t) \right\}$ where $\sigma(z_j,t)$ is given by $(\ref{h11})$. In order to show that $\rho_t\left( \frac{\partial}{\partial t} \right)=\delta^* \rho_{d,t}\left( \frac{\partial}{\partial t} \right)$, by the definition of $\delta^* $, it suffices to prove that there is a $0$-cochain $\left\{\xi_j\right\}$ with $\xi_j\in \Gamma\left(U_{jt}, \Xi_t\right)$ such that
\begin{align*}
\xi_j S(t)=\sigma_j(z_j,t),\,\,\,\,\,\,\,\,\,\xi_k-\xi_j=\theta_{jk}(t),\,\,\,\,\,\,\,\,\,\, D_{W_{t_0}}(\xi_j)= \lambda_j(t)=\left\{ V_j^\alpha \mapsto - \overline{ \frac{\partial V_j^\alpha}{\partial t} }   \right\} \left( \frac{\partial V_j^\alpha}{\partial t}+ D_{W_{t_0}}\left(\xi_j\right)\left( V_j^\alpha \right) \in \Gamma(U_{jt}, \Theta_{\mathcal{G}_t}) \right)
\end{align*}

Suppose that $w_j^\lambda=g_{jk}^\lambda(w_k)$ on $W_j\cap W_k$. Then since $\Phi:z_j\mapsto w_j=\varphi_j(z_j,t)$, we have on $U_{jt}\cap U_{kt}$,
\begin{align*}
&\varphi_j^\lambda(z_j,t)= g_{jk}^{\lambda}(\varphi_k(z_k,t)),\,\,\,\,\,\,\,\,\,\, z_j^\alpha=f_{jk}^\alpha(z_k,t)
\end{align*}

Differentiating both sides of this equality with respect to $t$ as functions of $z_k^1,...,z_k^n, t$, and putting
\begin{align}
\xi_j =\sum_{\lambda=0}^n \frac{\partial \varphi_j^\lambda(z_j,t)}{\partial t} \frac{\partial}{\partial w_j^\lambda}
\end{align}
we have $\theta_{jk}(t)= \xi_k-\xi_j$ (see \cite{Kod05} p.239), and by differentiating the $S_j\left(\varphi_j(z_j,t) \right)=0$ with respect to $t$, we obtain $\xi_j S(t)= \sigma_j(z_j,t)$ (see \cite{Kod05} p.240). On the other hand, we recall $(\ref{h6})$ that
\begin{align*}
T_j^{\alpha\gamma}(\varphi_j(z_j,t))=\sum_{\beta=1}^n V_j^{\alpha\beta}(z_j,t)\frac{\partial \varphi_j^\gamma(z_j,t)}{\partial  z_j^\beta}
\end{align*}
Differentiating both sides of this equality with respect to $t$ as functions of $z_k^1,...,z_k^n, t$, we have
\begin{align*}
\sum_{v=0}^n \frac{\partial T_j^{\alpha\gamma}}{\partial w_j^v}\left(\varphi_j(z_j,t) \right)\cdot \frac{\partial \varphi_j^v(z_j,t)}{\partial t}= \sum_{\beta=1}^n  V_j^{\alpha\beta}(z_j,t)\frac{\partial}{\partial z_j^\beta}\left( \frac{\partial \varphi_j^\gamma( z_j,t)}{\partial t} \right) + \sum_{\beta=1}^n \frac{\partial V_j^{\alpha\beta}(z_j,t)}{\partial t}\cdot  \frac{\partial \varphi_j^\gamma(z_j,t)}{\partial z_j^\beta}
\end{align*}
Then we have
\begin{align*}
\sum_{v,\gamma=0}^n \frac{\partial T_j^{\alpha\gamma}}{\partial w_j^v}\left( \varphi_j \right)\cdot \frac{\partial \varphi_j^v}{\partial t}\frac{\partial}{\partial w_j^\gamma}= \sum_{\gamma=0}^n \sum_{\beta=1}^n  V_j^{\alpha\beta}(z_j,t) \frac{\partial}{\partial z_j^\beta}\left(\frac{\partial \varphi_j^\gamma}{\partial t} \right)\frac{\partial}{\partial w_j^\gamma} + \sum_{\gamma=0}^n \sum_{\beta=1}^n \frac{\partial V_j^\beta}{\partial t} \frac{\partial \varphi_j^\gamma}{\partial z_j^\beta}\frac{\partial}{\partial w_j^\gamma}
\end{align*}

We note that 
\begin{align*}
 D_{W_{t_0}}\left(\xi_j\right)\left(V_j^\alpha\right)=  \overline{\sum_{\gamma=0}^n \sum_{\beta=1}^n V_j^{\alpha\beta}(z_j,t)\frac{\partial }{\partial z_j^\beta}\left( \frac{\partial \varphi_j^\gamma}{\partial t} \right) \frac{\partial}{\partial w_j^\gamma} - \sum_{\gamma, v=0}^n \frac{\partial T_j^{\alpha\gamma} }{\partial w_j^v}\left( \varphi_j \right) \cdot \frac{\partial \varphi_j^v}{\partial t}\frac{\partial}{\partial w_j^\gamma}}
\end{align*}

Then we have
\begin{align*}
 D_{W_{t_0}}\left(\xi_j\right)\left(V_j^\alpha\right)= \overline{-\sum_{\beta=1}^n \frac{\partial V_j^{\alpha\beta} }{\partial t}\frac{\partial}{\partial z_j^\beta} }= \overline{-\frac{\partial V_j^\alpha}{\partial t}}
\end{align*}
This completes the proof of Lemma \ref{h12}.

\end{proof}

\begin{example}\label{pr6}
Let us consider the foliation $\left( \mathbb{P}_\mathbb{C}^3, T_0= \xi_2\frac{\partial}{\partial \xi_0}+ \xi_0\frac{\partial}{\partial \xi_1} + \xi_3\frac{\partial}{\partial \xi_2} \right)$ of degree $1$ as in \textnormal{Example \ref{pr2}}. We keep the notations in \textnormal{Example \ref{pr2}}. Then $\xi_3=0$ defines a non-degenerate invariant complex submanifold $\left( \mathbb{P}_\mathbb{C}^2, T|_{\mathbb{P}_\mathbb{C}^2}= \xi_2\frac{\partial}{\partial \xi_0} + \xi_0 \frac{\partial}{\partial \xi_1} \right)$ of $\left( \mathbb{P}_\mathbb{C}^3, T \right)$ such that $\Theta_{\mathcal{F}}|_{\mathbb{P}_\mathbb{C}^2} \cong \mathcal{O}_{\mathbb{P}_\mathbb{C}^2}$. We note that $\left(\mathbb{P}_\mathbb{C}^3, T \right)$ and $\left( \mathbb{P}_\mathbb{C}^2, T|_{\mathbb{P}_\mathbb{C}^2} \right)$ are unobstructed in foliated deformations in terms of tangent sheaves. Let us consider the foliated line bundle $\left(F, \nabla_F \right)\cong \left(  \mathcal{O}_{\mathbb{P}_\mathbb{C}^2}(1), \nabla_F \right)$ on $\left(\mathbb{P}_\mathbb{C}^2, T|_{\mathbb{P}_\mathbb{C}^2} \right)$ determined by $\xi_3=0$. We will compute 
\begin{align*}
\mathcal{O}_{\mathbb{P}_\mathbb{C}^2}(1)^\bullet : \mathcal{O}_{\mathbb{P}_\mathbb{C}^2}\left( 1 \right) \xrightarrow{\nabla_F^0} \mathscr{H}om_{\mathcal{O}_{\mathbb{P}_\mathbb{C}^2}}\left( \mathcal{O}_{\mathbb{P}_\mathbb{C}^2}, \mathcal{O}_{\mathbb{P}_\mathbb{C}^2}(1) \right)\cong \mathcal{O}_{\mathbb{P}_\mathbb{C}^2}(1) \to 0
\end{align*}
We can compute $\mathbb{H}^i \left( \mathbb{P}_\mathbb{C}^2, \mathcal{O}_{\mathbb{P}_\mathbb{C}^2}(1)^\bullet \right)$ by the following complex
\begin{align*}
\mathcal{O}_{\mathbb{P}_\mathbb{C}^2}(1)^\bullet :  H^0\left( \mathbb{P}_\mathbb{C}^2, \mathcal{O}_{\mathbb{P}_\mathbb{C}^2}(1) \right) \xrightarrow{\nabla_F^0} H^0\left( \mathbb{P}_\mathbb{C}^2, \mathcal{O}_{\mathbb{P}_\mathbb{C}^2}(1) \right) \to 0
\end{align*}
On $U_0=\left\{ \xi_0\ne 0\right\}$ in $\mathbb{P}_\mathbb{C}^3$, $T$ is given by $\left( 1 - X_1X_2 \right)\frac{\partial}{\partial X_1} + \left( X_3- X_2^2 \right) \frac{\partial}{\partial X_2} - X_2 X_3 \frac{\partial}{\partial X_3}$ and $\xi_3$ is given by $X_3=0$. Then we have
\begin{align*}
\left[ \left( 1 - X_1X_2 \right)\frac{\partial}{\partial X_1} + \left( X_3- X_2^2 \right) \frac{\partial}{\partial X_2} - X_2 X_3 \frac{\partial}{\partial X_3} , X_3  \right]= X_3\left( - X_2 \right)
\end{align*}
On the other hand, $T|_{\mathbb{P}_\mathbb{C}^2}$ is given by $\left( 1- X_1 X_2 \right)\frac{\partial}{\partial X_1} - X_2^2 \frac{\partial}{\partial X_2}$ on $U_0=\left\{ \xi_0\ne 0 \right\}$ in $\mathbb{P}_\mathbb{C}^2$, and $H^0\left( \mathbb{P}_\mathbb{C}^2, \mathcal{O}_{\mathbb{P}_\mathbb{C}^2} (1) \right)=\left\{ a+ b X_1 + c X_2 | a, b, c \in \mathbb{C} \right\} \cong \mathbb{C}^3$ on $U_0$ in $\mathbb{P}_\mathbb{C}^2$. Then we have
\begin{align}\label{pr5}
\nabla_F^0\left( a+ bX_1 + c X_2 \right)& =\left[ \left( 1- X_1 X_2 \right)\frac{\partial}{\partial X_1} - X_2^2 \frac{\partial}{\partial X_2} , a+ bX_1 + cX_2 \right] - \left( - X_2 \right)\left(a+ b X_1 + c X_2 \right)= b + a X_2=0\\
&\iff a=b=0 \notag
\end{align}
This implies that $\dim_\mathbb{C} \mathbb{H}^0\left( \mathbb{P}_\mathbb{C}^2, \mathcal{O}_{\mathbb{P}_\mathbb{C}^2}(1)^\bullet \right)= 1$ and we have
\begin{align}\label{pr3}
\mathbb{H}^0\left( \mathbb{P}_\mathbb{C}^2, \mathcal{O}_{\mathbb{P}_\mathbb{C}^2}(1)^\bullet \right) \cong \textnormal{Span}_\mathbb{C} \left\langle X_2 \right\rangle
\end{align}
and $\dim_\mathbb{C} \mathbb{H}^1 \left( \mathbb{P}_\mathbb{C}^2, \mathcal{O}_{\mathbb{P}_\mathbb{C}^2}(1)^\bullet \right)=2$. We show that there is no foliated analytic family $\left( \mathcal{M}, \Theta_{\mathcal{G}}, \mathbb{C}, \omega \right)$ of non-degenerate invariant hypersurfaces of $\left(\mathbb{P}_\mathbb{C}^3, T \right)$ such that $\omega^{-1}(0)= \left(\mathbb{P}_\mathbb{C}^2, T|_{\mathbb{P}_\mathbb{C}^2} \right)$ defined by $\xi_3=0$ and the infinitesimal displacement is $\rho_{d,0}\left( \frac{\partial}{\partial t} \right) = X_2 \in \mathbb{H}^0\left( \mathbb{P}_\mathbb{C}^2, \mathcal{O}_{\mathbb{P}_\mathbb{C}^2}(1)^\bullet \right)$ on $U_0$. In other words, $\xi_3=0$ is obstructed in embedded foliated deformation as a non-degenerate invariant submanifold in $\left(\mathbb{P}_\mathbb{C}^3, T \right)$, while  the invariant complex submanifold $\xi_3=0$ is unobstructed in abstract foliated deformations. We will assume the knowledge of \textnormal{Part VI} on deformations of non-degenerate invariant submanifolds, and in particular we note that $\mathbb{H}^1\left( \mathbb{P}_\mathbb{C}^2, \mathcal{O}_{\mathbb{P}_\mathbb{C}^2}(1)^\bullet \right)$ represents the obstruction space $($see \textnormal{Part VI}$)$. In fact, we have a family of invariant submanifolds defined by $\xi_3- t\xi_2=0$ over $\textnormal{Spec}\left(\mathbb{C}[t]/\left(t^2\right)\right)$ from $(\ref{pr3})$. Then
\begin{align*}
&\left[\left( 1 - X_1X_2 \right)\frac{\partial}{\partial X_1} + \left( X_3- X_2^2 \right) \frac{\partial}{\partial X_2} - X_2 X_3 \frac{\partial}{\partial X_3}, X_3- t X_2 \right]|_{X_3=tX_2} =  - X_2X_3 - t X_3+ t X_2^2 |_{X_3= tX_2}\\
&= - t X_2^2 - t^2 X_2 + tX_2^2 = -t^2X_2
\end{align*}
So $-t^2X_2$ is the obstruction of lifting to $\textnormal{Spec}\left(\mathbb{C}[t]/\left(t^3\right) \right)$. We note that from $(\ref{pr5})$ $-X_2$ is the image of $\nabla_F^0\left(-1+0\cdot X_1+ 0\cdot X_2\right)$, so that we take $X_3= tX_2- t^2$. Then we have
\begin{align*}
&\left[\left( 1 - X_1X_2 \right)\frac{\partial}{\partial X_1} + \left( X_3- X_2^2 \right) \frac{\partial}{\partial X_2} - X_2 X_3 \frac{\partial}{\partial X_3}, X_3- t X_2 + t^2 \right]|_{X_3=tX_2- t^2}  = - X_2X_3 - t X_3+ tX_2^2|_{X_3= tX_2- t^2}\\
&=- X_2\left(tX_2 - t^2 \right) -t \left( tX_2- t^2 \right) +t X_2^2= t^3
\end{align*}
So $t^3$ is the obstruction of lifting to $\textnormal{Spec}\left(\mathbb{C}[t]/\left(t^4\right) \right)$. We note that from $(\ref{pr5})$ $1$ is the image of $\nabla_F^0\left( X_1 \right)$, so that we take $X_3= tX_2 - t^2 + t^3 X_1$. Then we have
\begin{align*}
&\left[\left( 1 - X_1X_2 \right)\frac{\partial}{\partial X_1} + \left( X_3- X_2^2 \right) \frac{\partial}{\partial X_2} - X_2 X_3 \frac{\partial}{\partial X_3}, X_3- t X_2 + t^2 - t^3 X_1 \right]|_{X_3=tX_2- t^2}  \\
&= - X_2X_3 - t X_3+ tX_2^2 - t^3+ t^3 X_1X_2|_{X_3= tX_2- t^2+ t^3 X_1} \\
&= - X_2\left( tX_2- t^2 + t^3 X_1 \right)-t\left(tX_2 - t^2 + t^3 X_1 \right)b+ tX_2^2 - t^3+ t^3 X_1X_2 = -t^4 X_1
\end{align*}
So $-t^4 X_1$ is the obstruction of lifting to $\textnormal{Spec}\left(\mathbb{C}[t]/\left(t^5\right) \right)$. However $-X_1 \ne 0 \in \mathbb{H}^1\left( \mathbb{P}_\mathbb{C}^2, \mathcal{O}_{\mathbb{P}_\mathbb{C}^2}(1)^\bullet \right)$. Hence we cannot lift the family over $\textnormal{Spec}\left( \mathbb{C}[t]/(t^4) \right)$ to $\textnormal{Spec}\left( \mathbb{C}[t]/\left(t^5\right) \right)$, so that $\xi_3=0$ is obstructed in embedded foliated deformations as an invariant submanifold in $\left( \mathbb{P}_\mathbb{C}^3, T_0 \right)$.

\end{example}

\begin{example}
Let us consider the foliation $\left( \mathbb{P}_\mathbb{C}^3, T_0=\xi_2\frac{\partial}{\partial \xi_0}+ \xi_0\frac{\partial}{\partial \xi_1} + \xi_3 \frac{\partial}{\partial \xi_2} \right)$ of degree $1$ as in \textnormal{Example \ref{pr2}} and the invariant complex submanifold defined by $\xi_3=0$ as in \textnormal{Example \ref{pr6}}. By assuming the knowledge of \textnormal{Part VI} on simultaneous deformations of compact foliated complex manifolds defined by locally free subsheaves, we will study simultaneous deformations of $\left\{\xi_3=0\right\} \cong \mathbb{P}_\mathbb{C}^2 \subset \left( \mathbb{P}_\mathbb{C}^3, T_0 \right)$. We will show that the simultaneous deformation functor $\textnormal{FDef}_{\xi_3=0\subset \left( \mathbb{P}_\mathbb{C}^3, T_0 \right)}$ of Artin rings is unobstructed. Simultaneous deformations of $\xi_3=0\subset \left( \mathbb{P}_\mathbb{C}^3 , T_0  \right)$ are controlled by the complex of sheaves $($see \textnormal{Part VI}$)$
\begin{align*}
\Theta_{\mathbb{P}_\mathbb{C}^3}^\bullet : \Theta_{\mathbb{P}_\mathbb{C}^3} \to \mathscr{H}om_{\mathcal{O}_{\mathbb{P}_\mathbb{C}^3}}\left( \Theta_{\mathcal{F}_0}, \frac{\Theta_{\mathbb{P}_\mathbb{C}^2}}{\Theta_{\mathcal{F}_0}} \right) \bigoplus i_* \mathcal{N}_{\mathbb{P}_\mathbb{C}^2/\mathbb{P}_\mathbb{C}^3} \to i_* \left( \Theta_{\mathcal{F}_0}^*|_{\mathbb{P}_\mathbb{C}^2} \otimes \mathcal{N}_{\mathbb{P}_\mathbb{C}^2/\mathbb{P}_\mathbb{C}^3} \right) \to 0
\end{align*}
Equivalently we have
\begin{align*}
\Theta_{\mathbb{P}_\mathbb{C}^3}^\bullet : \Theta_{\mathbb{P}_\mathbb{C}^3} \to \frac{\Theta_{\mathbb{P}_\mathbb{C}^3}}{\mathcal{O}_{\mathbb{P}_\mathbb{C}^3}} \bigoplus i_* \mathcal{O}_{\mathbb{P}_\mathbb{C}^2}(1) \to i_* \mathcal{O}_{\mathbb{P}_\mathbb{C}^2}(1) \to 0
\end{align*}
where $i: \left\{ \xi_3=0 \right\} \cong \mathbb{P}_\mathbb{C}^2 \hookrightarrow \mathbb{P}_\mathbb{C}^3$ is the embedding. We will denote the $i$-th cohomology group of $\Theta_{\mathbb{P}_\mathbb{C}^3}^\bullet$ by $\mathbb{H}^i\left( \mathbb{P}_\mathbb{C}^3, \Theta_{\mathbb{P}_\mathbb{C}^3}^\bullet \right)$.\footnote{We remark that we use different notations for the cohomology groups in Part VI. We use $\mathbb{H}^i\left( \mathbb{P}_\mathbb{C}^3, \left( \mathbb{P}_\mathbb{C}^2, \mathbb{P}_\mathbb{C}^3, \Theta_{\mathcal{F}_0} \right)^\bullet \right)$ in Part VI.} We can compute $\mathbb{H}^i\left( \mathbb{P}_\mathbb{C}^3, \Theta_{\mathbb{P}_\mathbb{C}^3}^\bullet \right)$ by the following complex
\begin{align*}
H^0\left( \mathbb{P}_\mathbb{C}^3, \Theta_{\mathbb{P}_\mathbb{C}^3} \right)\to  H^0\left( \mathbb{P}_\mathbb{C}^3, \frac{\Theta_{\mathbb{P}_\mathbb{C}^3}}{\mathcal{O}_{\mathbb{P}_\mathbb{C}^3}} \right) \bigoplus H^0\left( \mathbb{P}_\mathbb{C}^2, \mathcal{O}_{\mathbb{P}_\mathbb{C}^2}(1) \right)\to H^0\left( \mathbb{P}_\mathbb{C}^2, \mathcal{O}_{\mathbb{P}_\mathbb{C}^2}(1) \right)\to 0
\end{align*}
We note that from the exact sequence
\begin{align*}
0\to \mathcal{O}_{\mathbb{P}_\mathbb{C}^3}&\to \Theta_{\mathbb{P}_\mathbb{C}^3} \to \frac{\Theta_{\mathbb{P}_\mathbb{C}^3}}{\mathcal{O}_{\mathbb{P}_\mathbb{C}^3}}\to 0 \\
 1 &\mapsto T_0= \left( 1 - X_1X_2 \right)\frac{\partial}{\partial X_1} + \left( X_3 - X_2^2 \right) \frac{\partial}{\partial X_2} - X_2 X_3 \frac{\partial}{\partial X_3}
\end{align*}
and from $\textnormal{(\ref{ap2})}$ and $\frac{\partial}{\partial X_1}+ X_3 \frac{\partial}{\partial X_2}\equiv X_1X_2 \frac{\partial}{\partial X_1}+ X_2^2 \frac{\partial}{\partial X_2}+ X_2 X_3\frac{\partial}{\partial X_3}$ we see that
\begin{align} \label{pr9}
H^0\left( \mathbb{P}_\mathbb{C}^3,  \frac{\Theta_{\mathbb{P}_\mathbb{C}^3}}{\mathcal{O}_{\mathbb{P}_\mathbb{C}^3}} \right) = &\left\{  \left( e_1 + e_2 X_3+ e_3 X_2 + e_4 X_1 + e_5 X_1X_3  + e_7 X_1^2 \right) \frac{\partial}{\partial X_1}\right. \\
 &\left.+\left( e_8 + e_9 X_3 + e_{10}X_2 + a_{11}X_1 + e_5 X_2X_3 + e_7 X_1X_2\right) \frac{\partial}{\partial X_2} \right. \notag \\ 
 &\left.  +\left( e_{12}+ e_{13} X_3 + e_{14} X_2 + a_{15}X_1 + e_5 X_3^2   + e_7 X_1X_3 \right)  \frac{\partial}{\partial X_3}, e_1,...,e_5, e_7,..., e_{15}\in \mathbb{C} \right\} \notag
\end{align}

We compute $H^0\left( \mathbb{P}_\mathbb{C}^3, \frac{\Theta_{\mathbb{P}_\mathbb{C}^3}}{\mathcal{O}_{\mathbb{P}_\mathbb{C}^3}} \right) \bigoplus H^0\left( \mathbb{P}_\mathbb{C}^2, \mathcal{O}_{\mathbb{P}_\mathbb{C}^2}(1) \right)\to H^0\left( \mathbb{P}_\mathbb{C}^2, \mathcal{O}_{\mathbb{P}_\mathbb{C}^2}(1) \right)$. Then from $(\ref{pr9})$ the image is given by $($see \textnormal{Part VI}$)$
{\small{\begin{align} \label{pr14}
&\left[  \left( e_1 + e_2 X_3+ e_3 X_2 + e_4 X_1 + e_5 X_1X_3  + e_7 X_1^2 \right) \frac{\partial}{\partial X_1} + \left( e_8 + e_9 X_3 + e_{10}X_2 + e_{11}X_1 + e_5 X_2X_3 + e_7 X_1X_2\right) \frac{\partial}{\partial X_2} , X_3 \right]|_{X_3=0}\\
&+ \left[ \left( e_{12}+ e_{13} X_3 + e_{14} X_2 + e_{15}X_1 + e_5 X_3^2   + e_7 X_1X_3 \right)  \frac{\partial}{\partial X_3} , X_3 \right]|_{X_3=0} \notag \\
&+\left[ \left(1- X_1X_2 \right)\frac{\partial}{\partial X_1} - X_2^2 \frac{\partial}{\partial X_2} , a + bX_1+ c X_2\right] -(-X_2)\left(a+ bX_1+ cX_2\right) =  \left( e_{12}+ b\right) + \left( e_{14}+a \right) X_2+ e_{15} X_1 \notag
\end{align}}}
Then $\textnormal{ker}\left(H^0\left( \mathbb{P}_\mathbb{C}^3, \frac{\Theta_{\mathbb{P}_\mathbb{C}^3}}{\mathcal{O}_{\mathbb{P}_\mathbb{C}^3}} \right) \bigoplus H^0\left( \mathbb{P}_\mathbb{C}^2, \mathcal{O}_{\mathbb{P}_\mathbb{C}^2}(1) \right)\to H^0\left( \mathbb{P}_\mathbb{C}^2, \mathcal{O}_{\mathbb{P}_\mathbb{C}^2}(1) \right)\right)$ is given by
\begin{align} \label{pr12}
 &\left\{  \left( e_1 + e_2 X_3+ e_3 X_2 + e_4 X_1 + e_5 X_1X_3  + e_7 X_1^2 \right) \frac{\partial}{\partial X_1}+\left( e_8 + e_9 X_3 + e_{10}X_2 + e_{11}X_1 + e_5 X_2X_3 + e_7 X_1X_2\right) \frac{\partial}{\partial X_2} \right.  \\ 
 &\left.  +\left( e_{12}+ e_{13} X_3 + e_{14} X_2  + e_5 X_3^2   + e_7 X_1X_3 \right)  \frac{\partial}{\partial X_3} - e_{14} - e_{12}X_1 +c X_2 | e_1,...,e_5, e_7,..., e_{14},c\in \mathbb{C} \right\} \notag
\end{align}
and so we have $\dim_\mathbb{C} \textnormal{ker}\left(H^0\left( \mathbb{P}_\mathbb{C}^3, \frac{\Theta_{\mathbb{P}_\mathbb{C}^3}}{\mathcal{O}_{\mathbb{P}_\mathbb{C}^3}} \right) \bigoplus H^0\left( \mathbb{P}_\mathbb{C}^2, \mathcal{O}_{\mathbb{P}_\mathbb{C}^2}(1) \right)\to H^0\left( \mathbb{P}_\mathbb{C}^2, \mathcal{O}_{\mathbb{P}_\mathbb{C}^2}(1) \right) \right)=14$ and 
\begin{align*}
\dim_\mathbb{C} \mathbb{H}^2\left(\mathbb{P}_\mathbb{C}^3, \Theta_{\mathbb{P}_\mathbb{C}^3}^\bullet \right)= \textnormal{coker}\left( H^0\left( \mathbb{P}_\mathbb{C}^3, \frac{\Theta_{\mathbb{P}_\mathbb{C}^3}}{\mathcal{O}_{\mathbb{P}_\mathbb{C}^3}} \right) \bigoplus H^0\left( \mathbb{P}_\mathbb{C}^2, \mathcal{O}_{\mathbb{P}_\mathbb{C}^2}(1) \right)\to H^0\left( \mathbb{P}_\mathbb{C}^2, \mathcal{O}_{\mathbb{P}_\mathbb{C}^2}(1) \right) \right)= 3- (17-14)=0
\end{align*}
Since $\mathbb{H}^2\left( \mathbb{P}_\mathbb{C}^3, \Theta_{\mathbb{P}_\mathbb{C}^3}^\bullet \right)$ is an obstruction space for the simultaneous deformation functor $\textnormal{FDef}_{\xi_3=0\subset \left( \mathbb{P}_\mathbb{C}^3, T_0 \right)}$ of Artin rings, the functor $\textnormal{FDef}_{\xi_3= 0\subset \left( \mathbb{P}_\mathbb{C}^3, T_0 \right)}$ is unobstructed. 

We compute $H^0\left( \mathbb{P}_\mathbb{C}^3, \Theta_{\mathbb{P}_\mathbb{C}^3} \right)\to  H^0\left( \mathbb{P}_\mathbb{C}^3, \frac{\Theta_{\mathbb{P}_\mathbb{C}^3}}{\mathcal{O}_{\mathbb{P}_\mathbb{C}^3}} \right) \bigoplus H^0\left( \mathbb{P}_\mathbb{C}^2, \mathcal{O}_{\mathbb{P}_\mathbb{C}^2}(1) \right)$. From $\textnormal{(\ref{ap2})}$, the image is given by $($see \textnormal{Part VI}$)$
{\small{\begin{align}
&\left[ \left( a_1 + a_2 X_3+ a_3 X_2 + a_4 X_1 + a_5 X_1X_3 + a_6 X_1 X_2 + a_7 X_1^2 \right) \frac{\partial}{\partial X_1},  \left( 1 - X_1X_2 \right)\frac{\partial}{\partial X_1} + \left( X_3 - X_2^2 \right) \frac{\partial}{\partial X_2} - X_2 X_3 \frac{\partial}{\partial X_3} \right]\notag\\
&+\left[ \left( a_8 + a_9 X_3 + a_{10}X_2 + a_{11}X_1 + a_5 X_2X_3 + a_6 X_2^2 + a_7 X_1X_2\right) \frac{\partial}{\partial X_2}  ,  \left( 1 - X_1X_2 \right)\frac{\partial}{\partial X_1} + \left( X_3 - X_2^2 \right) \frac{\partial}{\partial X_2} - X_2 X_3 \frac{\partial}{\partial X_3}\right] \notag \\
&+ \left[ \left( a_{12}+ a_{13} X_3 + a_{14} X_2 + a_{15}X_1 + a_5 X_3^2 + a_6 X_2 X_3  + a_7 X_1X_3 \right)  \frac{\partial}{\partial X_3} , \left( 1 - X_1X_2 \right)\frac{\partial}{\partial X_1} + \left( X_3 - X_2^2 \right) \frac{\partial}{\partial X_2} - X_2 X_3 \frac{\partial}{\partial X_3} \right]\notag\\
&+\left[ \left( a_{12}+ a_{13} X_3 + a_{14} X_2 + a_{15}X_1 + a_5 X_3^2 + a_6 X_2 X_3  + a_7 X_1X_3 \right)  \frac{\partial}{\partial X_3}, X_3 \right]|_{X_3=0} + a_{12}+ a_{14}X_2 + a_{15}X_1\notag\\
&=\left(- a_4 - \left(2a_7 + a_8\right) X_1  - \left(a_1+ a_6\right) X_2 - \left(a_3+ a_5\right)X_3 - a_{11}X_1^2 - a_{10} X_1X_2 - \left(a_6 + a_9\right)X_1X_3 \right)\frac{\partial}{\partial X_1}\notag\\
&+ \left(  \left( a_{12}- a_{11} \right) + a_{15}X_1 + \left(a_{14}- a_7 - 2a_8\right)X_2 + \left( a_{13}- a_{10}\right)X_3 - a_{10} X_2^2 - a_{11}X_1X_2  - \left(a_9 + a_6\right)X_2X_3 \right) \frac{\partial}{\partial X_2}\notag\\
&+ \left( - a_{15} - a_{12}X_2 - \left(a_7 + a_8 + a_{14}\right)X_3 - \left(a_6 + a_9\right)X_3^2 - a_{10} X_2X_3 - a_{11}X_1X_3 \right) \frac{\partial}{\partial X_3} \notag \\
&\equiv \left(- a_4 - a_{10} - \left(2a_7 + a_8\right) X_1  - \left(a_1+ a_6\right) X_2 - \left(a_3+ a_5\right)X_3 - a_{11}X_1^2 - \left(a_6 + a_9\right)X_1X_3 \right)\frac{\partial}{\partial X_1} \label{pr11} \\
&+ \left(  \left( a_{12}- a_{11} \right) + a_{15}X_1 + \left(a_{14}- a_7 - 2a_8 \right)X_2 + \left( a_{13}-2 a_{10}\right)X_3 - a_{11}X_1X_2  - \left(a_9 + a_6\right)X_2X_3 \right) \frac{\partial}{\partial X_2}\notag\\
&+ \left( - a_{15} - a_{12}X_2 - \left(a_7 + a_8 + a_{14}\right)X_3 - \left(a_6 + a_9\right)X_3^2  - a_{11}X_1X_3 \right) \frac{\partial}{\partial X_3} + a_{12} + a_{14} X_2 + a_{15} X_1  \notag
\end{align}}}
Then from $(\ref{pr9})$ we see that
\begin{align*}
\dim_\mathbb{C} \mathbb{H}^0\left( \mathbb{P}_\mathbb{C}^3, \Theta_{\mathbb{P}_\mathbb{C}^3}^\bullet \right) & = \textnormal{ker}\left( H^0\left( \mathbb{P}_\mathbb{C}^3, \Theta_{\mathbb{P}_\mathbb{C}^3} \right)\to  H^0\left( \mathbb{P}_\mathbb{C}^3, \frac{\Theta_{\mathbb{P}_\mathbb{C}^3}}{\mathcal{O}_{\mathbb{P}_\mathbb{C}^3}} \right) \bigoplus H^0\left( \mathbb{P}_\mathbb{C}^2, \mathcal{O}_{\mathbb{P}_\mathbb{C}^2}(1) \right) \right) = 4\\
\dim_\mathbb{C}\mathbb{H}^1\left( \mathbb{P}_\mathbb{C}^3 , \Theta_{\mathbb{P}_\mathbb{C}^3}^\bullet \right) & = 14 - (15-4)= 3
\end{align*}
By comparing $(\ref{pr12})$ and $(\ref{pr11})$, we have
\begin{align}
\mathbb{H}^1\left( \mathbb{P}_\mathbb{C}^3, \Theta_{\mathbb{P}_\mathbb{C}^3}^\bullet \right) &\cong \frac{\textnormal{ker}\left( H^0\left( \mathbb{P}_\mathbb{C}^3, \frac{\Theta_{\mathbb{P}_\mathbb{C}^3}}{\mathcal{O}_{\mathbb{P}_\mathbb{C}^3}} \right) \bigoplus H^0\left( \mathbb{P}_\mathbb{C}^2, \mathcal{O}_{\mathbb{P}_\mathbb{C}^2}(1) \right) \to H^0\left( \mathbb{P}_\mathbb{C}^2, \mathcal{O}_{\mathbb{P}_\mathbb{C}^2}(1) \right) \right)}{\textnormal{im}\left( H^0\left( \mathbb{P}_\mathbb{C}^3, \Theta_{\mathbb{P}_\mathbb{C}^3} \right)\to  H^0\left( \mathbb{P}_\mathbb{C}^3, \frac{\Theta_{\mathbb{P}_\mathbb{C}^3}}{\mathcal{O}_{\mathbb{P}_\mathbb{C}^3}} \right) \bigoplus H^0\left( \mathbb{P}_\mathbb{C}^2, \mathcal{O}_{\mathbb{P}_\mathbb{C}^2}(1) \right) \right)}  \notag  \\
&\cong \textnormal{Span}_\mathbb{C} \left\langle - \frac{\partial}{\partial X_3}+ X_1, - X_2 \frac{\partial}{\partial X_3} + 1 ,  X_2 \right\rangle \label{pr13}
\end{align}
We will construct a foliated analytic family of simultaneous deformations of $\left\{  \xi_3=0 \right\} \subset \left( \mathbb{P}_\mathbb{C}^3, T_0 \right)$ such that the characteristic map is an isomorphism onto $\mathbb{H}^1\left( \mathbb{P}_\mathbb{C}^3, \Theta_{\mathbb{P}_\mathbb{C}^3}^\bullet \right)$ at the distinguished point. We have a family of simultaneous deformations of $\left\{  \xi_3=0 \right\} \subset \left( \mathbb{P}_\mathbb{C}^3, T_0 \right)$ over $\textnormal{Spec}\left( \mathbb{C}\left[t_1,t_2,t_3\right]/\left(t_1^2,t_2^2,t_3^2, t_1t_2, t_1t_3,t_2t_3 \right) \right)$ from $(\ref{pr13})$ which is  defined by
\begin{align*}
\tilde{T}&=\left( 1 - X_1X_2 \right)\frac{\partial}{\partial X_1} +\left( X_3- X_2^2 \right)\frac{\partial}{\partial X_2} +\left(- X_2X_3 + t_1 + t_3X_2 \right) \frac{\partial}{\partial X_3}\\
X_3&= t_3 + t_1 X_1+ t_2 X_2 
\end{align*}
Then we have
\begin{align*}
&\left[ \left( 1 - X_1X_2 \right)\frac{\partial}{\partial X_1} +\left( X_3- X_2^2 \right)\frac{\partial}{\partial X_2} +\left(- X_2X_3 + t_1 + t_3X_2 \right) \frac{\partial}{\partial X_3}, X_3- t_3 - t_1 X_1 - t_2 X_2 \right]|_{X_3= t_3+ t_1 X_1 + t_2 X_2}\\
&= - X_2X_3 + t_1+ t_3 X_2 - t_1 + t_1 X_1X_2-  t_2 X_3 + t_2 X_2^2|_{X_3= t_3+ t_1 X_1 + t_2 X_2} \\
&=- X_2\left( t_3 + t_1 X_1 + t_2 X_2 \right) + t_1 + t_3 X_2 - t_1 + t_1 X_1 X_2 - t_2\left(t_3 + t_1 X_1 + t_2 X_2 \right) + t_2 X_2^2 = - t_2t_3 - t_1 t_2 X_1 - t_2^2 X_2
\end{align*}
So $-t_2t_3- t_1t_2X_1 - t_2^2X_2$ is the obstruction of lifting. We note that from $(\ref{pr14})$, $-t_2t_3- t_1t_2X_1 - t_2^2X_2$ is the image of $- t_2t_3\frac{\partial}{\partial X_3} - t_1t_2X_1\frac{\partial}{\partial X_3}- t_2^2 X_2 \frac{\partial}{\partial X_3} $ under $H^0\left( \mathbb{P}_\mathbb{C}^3, \frac{\Theta_{\mathbb{P}_\mathbb{C}^3} }{\mathcal{O}_{\mathbb{P}_\mathbb{C}^3}} \right) \bigoplus H^0\left( \mathbb{P}_\mathbb{C}^2, \mathcal{O}_{\mathbb{P}_\mathbb{C}^2} (1) \right)\to H^0\left( \mathbb{P}_\mathbb{C}^2, \mathcal{O}_{\mathbb{P}_\mathbb{C}^2} (1) \right)$, so that we take
\begin{align}
\tilde{T}&=\left( 1 - X_1X_2 \right)\frac{\partial}{\partial X_1} +\left( X_3- X_2^2 \right)\frac{\partial}{\partial X_2} +\left(- X_2X_3 + t_1 + t_3X_2 + t_2t_3+ t_1t_2 X_1 + t_2^2 X_2 \right) \frac{\partial}{\partial X_3} \label{pr15}\\
X_3&= t_3 + t_1 X_1+ t_2 X_2 \label{pr17}
\end{align}
Then we have
{\small{\begin{align*}
&\left[ \left( 1 - X_1X_2 \right)\frac{\partial}{\partial X_1} +\left( X_3- X_2^2 \right)\frac{\partial}{\partial X_2} +\left(- X_2X_3 + t_1 + t_3X_2 + t_2t_3+ t_1t_2 X_1 + t_2^2 X_2 \right) \frac{\partial}{\partial X_3} , X_3 - t_3 - t_1 X_1 - t_2 X_2 \right]|_{X_3= t_3 + t_1 X_1 + t_2 X_2}\\
&= - X_2X_3 + t_1 + t_3X_2 + t_2t_3+ t_1t_2 X_1 + t_2^2 X_2 - t_1 + t_1 X_1 X_2 - t_2 X_3 + t_2 X_2^2|_{X_3= t_3 + t_1 X_1 + t_2 X_2} = 0
\end{align*}}}
Hence $\tilde{T}$ in $(\ref{pr15})$ and $X_3= t_3+ t_1 X_1 + t_2 X_2$ in $(\ref{pr17})$ defines a foliated analytic family of simultaneous deformations of $\left\{ \xi_3=0 \right\} \subset \left( \mathbb{P}_\mathbb{C}^3, T_0 \right)$ over a sufficiently neighborhood $B$ of $(0,0,0)\in \mathbb{C}^3$ such that the characteristic map is an isomorphism onto $\mathbb{H}^1\left( \mathbb{P}_\mathbb{C}^3, \Theta_{\mathbb{P}_\mathbb{C}^3}^\bullet \right)$ at $(t_1,t_2,t_3)=(0,0,0) \in B$.

\end{example}

\begin{example}\label{pr8}
Let us consider the foliation $\left( \mathbb{P}_\mathbb{C}^3, T_0=\xi_2\frac{\partial}{\partial \xi_0}+ \xi_0\frac{\partial}{\partial \xi_1} + \xi_3 \frac{\partial}{\partial \xi_2} \right)$ of degree $1$ as in \textnormal{Example \ref{pr2}}. We keep the notations in \textnormal{Example \ref{pr2}}. Let us consider the non-degenerate invariant submanifold defined by $\left\{\xi_2=\xi_3=0\right\}\cong \mathbb{P}_\mathbb{C}^1$ of $\left( \mathbb{P}_\mathbb{C}^3, T_0 \right)$. By assuming the knowledge of $\textnormal{Part VI}$ on deformations of non-degenerate invariant complex submanifolds of foliated complex manifolds defined by locally free subsheaves of tangent sheaves. We will show that $\xi_2=\xi_3=0$ is obstructed in embedded deformations as an invariant submanifold in $\left(\mathbb{P}_\mathbb{C}^3, T_0 \right)$. Embedded deformations of $\xi_2=\xi_3=0\cong \mathbb{P}_\mathbb{C}^1$ are controlled by the complex of sheaves $($see \textnormal{Part VI}$)$
\begin{align*}
\mathcal{N}_{\mathbb{P}_\mathbb{C}^1/\mathbb{P}_\mathbb{C}^3}^\bullet: \mathcal{N}_{\mathbb{P}_\mathbb{C}^1/\mathbb{P}_\mathbb{C}^3}\to \Theta_{\mathcal{F}_0}^*|_{\mathbb{P}_\mathbb{C}^1}\otimes \mathcal{N}_{\mathbb{P}_\mathbb{C}^1/\mathbb{P}_\mathbb{C}^3}\to 0
\end{align*}
Equivalently we have
\begin{align} \label{blw11}
\mathcal{N}_{\mathbb{P}_\mathbb{C}^1/\mathbb{P}_\mathbb{C}^3}^\bullet: \mathcal{O}_{\mathbb{P}_\mathbb{C}^1}(1)\bigoplus \mathcal{O}_{\mathbb{P}_\mathbb{C}^1} (1)  \to \mathcal{O}_{\mathbb{P}_\mathbb{C}^2}(1) \bigoplus \mathcal{O}_{\mathbb{P}_\mathbb{C}^2}(1) \to 0
\end{align}
We will denote the $i$-cohomology group of $\mathcal{N}_{\mathbb{P}_\mathbb{C}^1/ \mathbb{P}_\mathbb{C}^3}^\bullet$ by $\mathbb{H}^i\left( \mathbb{P}_\mathbb{C}^1, \mathcal{N}_{\mathbb{P}_\mathbb{C}^1/\mathbb{P}_\mathbb{C}^3}^\bullet  \right)$. We can compute $\mathbb{H}^i\left( \mathbb{P}_\mathbb{C}^1,  \mathcal{N}_{\mathbb{P}_\mathbb{C}^1/\mathbb{P}_\mathbb{C}^3}^\bullet \right)$ by the following complex
\begin{align*}
H^0\left( \mathbb{P}_\mathbb{C}^1, \mathcal{O}_{\mathbb{P}_\mathbb{C}^1}(1) \right)\bigoplus H^0\left( \mathbb{P}_\mathbb{C}^1, \mathcal{O}_{\mathbb{P}_\mathbb{C}^1} (1) \right)\xrightarrow{\nabla_\mathcal{N}} H^0\left( \mathbb{P}_\mathbb{C}^1, \mathcal{O}_{\mathbb{P}_\mathbb{C}^1}(1) \right)\bigoplus H^0\left( \mathbb{P}_\mathbb{C}^1, \mathcal{O}_{\mathbb{P}_\mathbb{C}^1} (1) \right) \to 0
\end{align*}
We let
\begin{align}\label{pr20}
 H^0\left( \mathbb{P}_\mathbb{C}^1, \mathcal{O}_{\mathbb{P}_\mathbb{C}^1}(1) \right)\bigoplus H^0\left( \mathbb{P}_\mathbb{C}^1, \mathcal{O}_{\mathbb{P}_\mathbb{C}^1} (1) \right)=\left\{ \left( b_1 + b_2 X_1, b_3 + b_4 X_1 \right) |b_1,b_2,b_3,b_4 \in \mathbb{C}  \right\}
\end{align}
 Then the image of $\nabla_\mathcal{N}$ is given by $($see \textnormal{Part VI}$)$
\begin{align}
&\left(\left[ \left[  \left( 1 - X_1X_2 \right)\frac{\partial}{\partial X_1} + \left( X_3 - X_2^2 \right) \frac{\partial}{\partial X_2} - X_2 X_3 \frac{\partial}{\partial X_3} , \left( b_1 + b_2 X_1\right)\frac{\partial}{\partial X_2} +\left( b_3 + b_4 X_1 \right)\frac{\partial}{\partial X_3}\right]  , X_2\right]|_{X_2=X_3=0} ,\right. \\
&\left.\left[ \left[  \left( 1 - X_1X_2 \right)\frac{\partial}{\partial X_1} + \left( X_3 - X_2^2 \right) \frac{\partial}{\partial X_2} - X_2 X_3 \frac{\partial}{\partial X_3}, \left( b_1 + b_2 X_1\right) \frac{\partial}{\partial X_2} + \left(b_3 + b_4 X_1 \right)\frac{\partial}{\partial X_3} \right], X_3 \right]|_{X_2=X_3=0}\right) \notag \\
&= \left( b_2- b_3 - b_4 X_1, b_4 \right)\in  H^0\left( \mathbb{P}_\mathbb{C}^1, \mathcal{O}_{\mathbb{P}_\mathbb{C}^1}(1) \right)\bigoplus H^0\left( \mathbb{P}_\mathbb{C}^1, \mathcal{O}_{\mathbb{P}_\mathbb{C}^1} (1) \right) \notag
\end{align}
Then $\textnormal{ker}\left( \nabla_\mathcal{N} \right)$ is given by
\begin{align} \label{blw12}
\left\{ \left( t_1+ t_2 X_1, t_2 \right) | t_1,t_2\in \mathbb{C} \right\}
\end{align}
and so we have
\begin{align*}
\dim_\mathbb{C} \mathbb{H}^0\left( \mathbb{P}_\mathbb{C}^1 ,  \mathcal{N}_{\mathbb{P}_\mathbb{C}^1/\mathbb{P}_\mathbb{C}^3}^\bullet \right) = 2
\end{align*}
Then we have a foliated analytic family of deformations of $ \xi_2=\xi_3=0\cong \mathbb{P}_\mathbb{C}^1$ in $\left( \mathbb{P}_\mathbb{C}^3, T_0 \right)$ over $\textnormal{Spec}\left(\mathbb{C}\left[t_1,t_2\right]/\left(t_1^2, t_2^2, t_1t_2\right)\right)$ defined by $X_2= t_1+ t_2 X_1, X_3= t_2$. Then we have
\begin{align}
&\left[  \left( 1 - X_1X_2 \right)\frac{\partial}{\partial X_1} + \left( X_3 - X_2^2 \right) \frac{\partial}{\partial X_2} - X_2 X_3 \frac{\partial}{\partial X_3} ,  X_2- t_1 - t_2 X_1 \right]|_{X_2=t_1+ t_2 X_1, X_3=t_2} \notag\\
&=X_3- X_2^2- t_2+ t_2 X_1X_2|_{X_2=t_1+t_2X_1, X_3=t_2}= t_2 - \left(t_1+ t_2X_1\right)^2- t_2 +  t_2 X_1\left( t_1 + t_2 X_1 \right) = - t_1^2 -  t_1t_2 X_1 \label{blw15}
\end{align}
and we have
\begin{align}\label{pr18}
&\left[ \left( 1 - X_1X_2 \right)\frac{\partial}{\partial X_1} + \left( X_3 - X_2^2 \right) \frac{\partial}{\partial X_2} - X_2 X_3 \frac{\partial}{\partial X_3} , X_3- t_2 \right]|_{X_2= t_1 + t_2 X_1, X_3 = t_2} \notag \\
&= - X_2X_3|_{X_2= t_1 + t_2 X_1, X_3 = t_2} = - \left(t_1+ t_2 X_1 \right)t_2= - t_1t_2 - t_2^2 X_1
\end{align}
We note that the obstruction $ - t_2^2 X_1$ in $(\ref{pr18})$ is not in the image of $\nabla_\mathcal{N}$, i.e. $\left(0, - X_1 \right)\ne 0 \in \mathbb{H}^1\left( \mathbb{P}_\mathbb{C}^1, \mathcal{N}_{\mathbb{P}_\mathbb{C}^1/\mathbb{P}_\mathbb{C}^3}^\bullet \right)$ which is an obstruction space for deformations of $\xi_2=\xi_3=0$ in $\left( \mathbb{P}_\mathbb{C}^3, T_0 \right)$ $($see \textnormal{Part VI}$)$. Hence $\xi_2=\xi_3=0$ is obstructed in embedded deformations as an invariant submanifold in $\left(\mathbb{P}_\mathbb{C}^3, T_0 \right)$.

\end{example}

\begin{example}
Let us consider the foliation $\left( \mathbb{P}_\mathbb{C}^3, T_0=\xi_2\frac{\partial}{\partial \xi_0}+ \xi_0\frac{\partial}{\partial \xi_1} + \xi_3 \frac{\partial}{\partial \xi_2} \right)$ of degree $1$ as in \textnormal{Example \ref{pr2}} and the invariant complex submanifold defined by $\xi_2=\xi_3=0$ as in \textnormal{Example \ref{pr8}}. We keep the notations in \textnormal{Example \ref{pr2}} and \textnormal{Example \ref{pr8}}. By assuming the knowledge of \textnormal{Part VI} on simultaneous deformations of compact foliated complex manifolds defined by locally free subsheaves, we will study simultaneous deformations of $\left\{\xi_2= \xi_3=0\right\} \cong \mathbb{P}_\mathbb{C}^1 \subset \left( \mathbb{P}_\mathbb{C}^3, T_0 \right)$. We will show that the simultaneous deformation functor $\textnormal{FDef}_{\xi_2=\xi_3=0\subset \left( \mathbb{P}_\mathbb{C}^3, T_0 \right)}$ of Artin rings is unobstructed. Simultaneous deformations of $\xi_2=\xi_3=0\subset \left( \mathbb{P}_\mathbb{C}^3 , T_0  \right)$ are controlled by the complex of sheaves $($see \textnormal{Part VI}$)$
\begin{align*}
\Theta_{\mathbb{P}_\mathbb{C}^3}^\bullet: \Theta_{\mathbb{P}_\mathbb{C}^3} \to \mathscr{H}om_{\mathcal{O}_{\mathbb{P}_\mathbb{C}^3}} \left( \Theta_{\mathcal{F}_0}, \frac{\Theta_{\mathbb{P}_\mathbb{C}^3}}{ \Theta_{\mathcal{F}_0}  } \right) \bigoplus i_* \mathcal{N}_{\mathbb{P}_\mathbb{C}^1 /\mathbb{P}_\mathbb{C}^3} \to i_* \mathcal{N}_{\mathbb{P}_\mathbb{C}^1/ \mathbb{P}_\mathbb{C}^3} \to 0
\end{align*}
Equivalently we have
\begin{align*}
\Theta_{\mathbb{P}_\mathbb{C}^3}^\bullet: \Theta_{\mathbb{P}_\mathbb{C}^3} \to \frac{\Theta_{\mathbb{P}_\mathbb{C}^3}}{ \mathcal{O}_{\mathbb{P}_\mathbb{C}^3} } \bigoplus i_* \left( \mathcal{O}_{\mathbb{P}_\mathbb{C}^1}(1)\bigoplus \mathcal{O}_{\mathbb{P}_\mathbb{C}^1}(1) \right) \to i_* \left(  \mathcal{O}_{\mathbb{P}_\mathbb{C}^1}(1) \bigoplus \mathcal{O}_{\mathbb{P}_\mathbb{C}^1}(1) \right)\to 0
\end{align*}
where $i:\left\{\xi_2=\xi_3=0 \right\}\cong \mathbb{P}_\mathbb{C}^1 \hookrightarrow \mathbb{P}_\mathbb{C}^3$ is the embedding. We will denote the $i$-th cohomology group of $\Theta_{\mathbb{P}_\mathbb{C}^3}^\bullet$ by $\mathbb{H}^i\left( \mathbb{P}_\mathbb{C}^3, \Theta_{\mathbb{P}_\mathbb{C}^3}^\bullet \right)$.\footnote{We remark that we use different notations for the cohomology groups in Part VI. We use $\mathbb{H}^i\left( \mathbb{P}_\mathbb{C}^3, \left( \mathbb{P}_\mathbb{C}^1, \mathbb{P}_\mathbb{C}^3, \Theta_{\mathcal{F}_0} \right)^\bullet \right)$ in Part VI.} We can compute $\mathbb{H}^i\left( \mathbb{P}_\mathbb{C}^3, \Theta_{\mathbb{P}_\mathbb{C}^3}^\bullet \right)$ by the following complex
\begin{align*}
H^0\left( \mathbb{P}_\mathbb{C}^3, \Theta_{\mathbb{P}_\mathbb{C}^3} \right) \to H^0\left( \mathbb{P}_\mathbb{C}^3, \frac{\Theta_{\mathbb{P}_\mathbb{C}^3}}{\mathcal{O}_{\mathbb{P}_\mathbb{C}^3}} \right) \bigoplus H^0\left( \mathbb{P}_\mathbb{C}^1, \mathcal{O}_{\mathbb{P}_\mathbb{C}^1}(1)\bigoplus \mathcal{O}_{\mathbb{P}_\mathbb{C}^1}(1) \right) \to  H^0\left( \mathbb{P}_\mathbb{C}^1, \mathcal{O}_{\mathbb{P}_\mathbb{C}^1}(1) \bigoplus \mathcal{O}_{\mathbb{P}_\mathbb{C}^1}(1) \right) \to 0
\end{align*}
We compute $H^0\left( \mathbb{P}_\mathbb{C}^3, \frac{\Theta_{\mathbb{P}_\mathbb{C}^3}}{\mathcal{O}_{\mathbb{P}_\mathbb{C}^3}} \right) \bigoplus H^0\left( \mathbb{P}_\mathbb{C}^1, \mathcal{O}_{\mathbb{P}_\mathbb{C}^1}(1)\bigoplus \mathcal{O}_{\mathbb{P}_\mathbb{C}^1}(1) \right) \to  H^0\left( \mathbb{P}_\mathbb{C}^1, \mathcal{O}_{\mathbb{P}_\mathbb{C}^1}(1) \bigoplus \mathcal{O}_{\mathbb{P}_\mathbb{C}^1}(1) \right)$. Then from $(\ref{pr9})$ and $(\ref{pr20})$, the first component of the image  is given by $($see \textnormal{Part VI}$)$
{\small{\begin{align} \label{pr21}
&\left[  \left( e_1 + e_2 X_3+ e_3 X_2 + e_4 X_1 + e_5 X_1X_3  + e_7 X_1^2 \right) \frac{\partial}{\partial X_1} + \left( e_8 + e_9 X_3 + e_{10}X_2 + e_{11}X_1 + e_5 X_2X_3 + e_7 X_1X_2\right) \frac{\partial}{\partial X_2} , X_2 \right]|_{X_2=X_3=0}\\
&+ \left[ \left( e_{12}+ e_{13} X_3 + e_{14} X_2 + e_{15}X_1 + e_5 X_3^2   + e_7 X_1X_3 \right)  \frac{\partial}{\partial X_3} , X_2 \right]|_{X_2=X_3=0} \notag \\
&+\left[ \left[  \left( 1 - X_1X_2 \right)\frac{\partial}{\partial X_1} + \left( X_3 - X_2^2 \right) \frac{\partial}{\partial X_2} - X_2 X_3 \frac{\partial}{\partial X_3} , \left( b_1 + b_2 X_1\right)\frac{\partial}{\partial X_2} +\left( b_3 + b_4 X_1 \right)\frac{\partial}{\partial X_3}\right]  , X_2\right]|_{X_2=X_3=0} \notag \\
&= \left( e_8 + b_2 - b_3 \right) + \left( e_{11}- b_4 \right) X_1 \notag
\end{align}}}
and the second component of the image is given by
{\small{\begin{align} \label{pr22}
&\left[  \left( e_1 + e_2 X_3+ e_3 X_2 + e_4 X_1 + e_5 X_1X_3  + e_7 X_1^2 \right) \frac{\partial}{\partial X_1} + \left( e_8 + e_9 X_3 + e_{10}X_2 + e_{11}X_1 + e_5 X_2X_3 + e_7 X_1X_2\right) \frac{\partial}{\partial X_2} , X_3 \right]|_{X_2=X_3=0}\\
&+ \left[ \left( e_{12}+ e_{13} X_3 + e_{14} X_2 + e_{15}X_1 + e_5 X_3^2   + e_7 X_1X_3 \right)  \frac{\partial}{\partial X_3} , X_3 \right]|_{X_2=X_3=0} \notag \\
&+\left[ \left[  \left( 1 - X_1X_2 \right)\frac{\partial}{\partial X_1} + \left( X_3 - X_2^2 \right) \frac{\partial}{\partial X_2} - X_2 X_3 \frac{\partial}{\partial X_3} , \left( b_1 + b_2 X_1\right)\frac{\partial}{\partial X_2} +\left( b_3 + b_4 X_1 \right)\frac{\partial}{\partial X_3}\right]  , X_3\right]|_{X_2=X_3=0} \notag \\
&= \left( e_{12} + b_4  \right) + e_{15}X_1 \notag
\end{align}}}
Then $\textnormal{ker}\left( H^0\left( \mathbb{P}_\mathbb{C}^3, \frac{\Theta_{\mathbb{P}_\mathbb{C}^3}}{\mathcal{O}_{\mathbb{P}_\mathbb{C}^3}} \right) \bigoplus H^0\left( \mathbb{P}_\mathbb{C}^1, \mathcal{O}_{\mathbb{P}_\mathbb{C}^1}(1)\bigoplus \mathcal{O}_{\mathbb{P}_\mathbb{C}^1}(1) \right) \to  H^0\left( \mathbb{P}_\mathbb{C}^1, \mathcal{O}_{\mathbb{P}_\mathbb{C}^1}(1) \bigoplus \mathcal{O}_{\mathbb{P}_\mathbb{C}^1}(1) \right) \right)$ is given by
\begin{align} \label{pr23}
 &\left\{  \left( e_1 + e_2 X_3+ e_3 X_2 + e_4 X_1 + e_5 X_1X_3  + e_7 X_1^2 \right) \frac{\partial}{\partial X_1}+\left( e_8 + e_9 X_3 + e_{10}X_2 + e_{11}X_1 + e_5 X_2X_3 + e_7 X_1X_2\right) \frac{\partial}{\partial X_2} \right.  \\ 
 &\left.  +\left( -e_{11}+ e_{13} X_3 + e_{14} X_2  + e_5 X_3^2   + e_7 X_1X_3 \right)  \frac{\partial}{\partial X_3} ,\left( b_1 + b_2 X_1 , b_2 + e_8  + e_{11} X_1 \right) | e_1,...,e_5, e_7,...,e_{11}, e_{13} e_{14}, b_1, b_2 \in \mathbb{C} \right\} \notag
\end{align}
and so $\dim_\mathbb{C} \textnormal{ker}\left( H^0\left( \mathbb{P}_\mathbb{C}^3, \frac{\Theta_{\mathbb{P}_\mathbb{C}^3}}{\mathcal{O}_{\mathbb{P}_\mathbb{C}^3}} \right) \bigoplus H^0\left( \mathbb{P}_\mathbb{C}^1, \mathcal{O}_{\mathbb{P}_\mathbb{C}^1}(1)\bigoplus \mathcal{O}_{\mathbb{P}_\mathbb{C}^1}(1) \right) \to  H^0\left( \mathbb{P}_\mathbb{C}^1, \mathcal{O}_{\mathbb{P}_\mathbb{C}^1}(1) \bigoplus \mathcal{O}_{\mathbb{P}_\mathbb{C}^1}(1) \right) \right)= 14$ and
{\small{\begin{align*}
\dim_\mathbb{C} \mathbb{H}^2\left( \mathbb{P}_\mathbb{C}^3, \Theta_{\mathbb{P}_\mathbb{C}^3}^\bullet  \right)= \textnormal{coker}\left(H^0\left( \mathbb{P}_\mathbb{C}^3, \frac{\Theta_{\mathbb{P}_\mathbb{C}^3}}{\mathcal{O}_{\mathbb{P}_\mathbb{C}^3}} \right) \bigoplus H^0\left( \mathbb{P}_\mathbb{C}^1, \mathcal{O}_{\mathbb{P}_\mathbb{C}^1}(1)\bigoplus \mathcal{O}_{\mathbb{P}_\mathbb{C}^1}(1) \right) \to  H^0\left( \mathbb{P}_\mathbb{C}^1, \mathcal{O}_{\mathbb{P}_\mathbb{C}^1}(1) \bigoplus \mathcal{O}_{\mathbb{P}_\mathbb{C}^1}(1) \right) \right) = 4 - (18- 14 ) = 0
\end{align*}}}
Since $\mathbb{H}^2\left( \mathbb{P}_\mathbb{C}^3, \Theta_{\mathbb{P}_\mathbb{C}^3}^\bullet \right)$ is an obstruction space for the simultaneous deformation functor $\textnormal{FDef}_{\xi_2=\xi_3=0\subset \left( \mathbb{P}_\mathbb{C}^3, T_0 \right)}$ of Artin rings, the functor $\textnormal{FDef}_{\xi_2=\xi_3= 0\subset \left( \mathbb{P}_\mathbb{C}^3, T_0 \right)}$ is unobstructed. 

We compute $H^0\left( \mathbb{P}_\mathbb{C}^3, \Theta_{\mathbb{P}_\mathbb{C}^3}  \right)\to H^0\left( \mathbb{P}_\mathbb{C}^3, \frac{\Theta_{\mathbb{P}_\mathbb{C}^3}}{\mathcal{O}_{\mathbb{P}_\mathbb{C}^3}} \right) \bigoplus H^0\left( \mathbb{P}_\mathbb{C}^1, \mathcal{O}_{\mathbb{P}_\mathbb{C}^1}(1)\bigoplus \mathcal{O}_{\mathbb{P}_\mathbb{C}^1}(1) \right)$. From $\textnormal{(\ref{ap2})}$, the image is given by $($see \textnormal{Part VI}$)$
{\small{\begin{align}
&\left[ \left( a_1 + a_2 X_3+ a_3 X_2 + a_4 X_1 + a_5 X_1X_3 + a_6 X_1 X_2 + a_7 X_1^2 \right) \frac{\partial}{\partial X_1},  \left( 1 - X_1X_2 \right)\frac{\partial}{\partial X_1} + \left( X_3 - X_2^2 \right) \frac{\partial}{\partial X_2} - X_2 X_3 \frac{\partial}{\partial X_3} \right]\notag\\
&+\left[ \left( a_8 + a_9 X_3 + a_{10}X_2 + a_{11}X_1 + a_5 X_2X_3 + a_6 X_2^2 + a_7 X_1X_2\right) \frac{\partial}{\partial X_2}  ,  \left( 1 - X_1X_2 \right)\frac{\partial}{\partial X_1} + \left( X_3 - X_2^2 \right) \frac{\partial}{\partial X_2} - X_2 X_3 \frac{\partial}{\partial X_3}\right] \notag \\
&+ \left[ \left( a_{12}+ a_{13} X_3 + a_{14} X_2 + a_{15}X_1 + a_5 X_3^2 + a_6 X_2 X_3  + a_7 X_1X_3 \right)  \frac{\partial}{\partial X_3} , \left( 1 - X_1X_2 \right)\frac{\partial}{\partial X_1} + \left( X_3 - X_2^2 \right) \frac{\partial}{\partial X_2} - X_2 X_3 \frac{\partial}{\partial X_3} \right]\notag\\
&+ \left( \left[  \left( a_8 + a_9 X_3 + a_{10}X_2 + a_{11}X_1 + a_5 X_2X_3 + a_6 X_2^2 + a_7 X_1X_2\right) \frac{\partial}{\partial X_2}, X_2\right]|_{X_2=X_3=0},\right. \notag \\
&\left. \left[ \left( a_{12}+ a_{13} X_3 + a_{14} X_2 + a_{15}X_1 + a_5 X_3^2 + a_6 X_2 X_3  + a_7 X_1X_3 \right)  \frac{\partial}{\partial X_3}, X_3  \right]|_{X_2=X_3=0}  \right) \notag\\
&\equiv  \left(- a_4 - a_{10} - \left(2a_7 + a_8\right) X_1  - \left(a_1+ a_6\right) X_2 - \left(a_3+ a_5\right)X_3 - a_{11}X_1^2 - \left(a_6 + a_9\right)X_1X_3 \right)\frac{\partial}{\partial X_1} \label{pr24} \\
&+ \left(  \left( a_{12}- a_{11} \right) + a_{15}X_1 + \left(a_{14}- a_7 - 2a_8 \right)X_2 + \left( a_{13}-2 a_{10}\right)X_3 - a_{11}X_1X_2  - \left(a_9 + a_6\right)X_2X_3 \right) \frac{\partial}{\partial X_2}\notag\\
&+ \left( - a_{15} - a_{12}X_2 - \left(a_7 + a_8 + a_{14}\right)X_3 - \left(a_6 + a_9\right)X_3^2  - a_{11}X_1X_3 \right) \frac{\partial}{\partial X_3} + \left( a_8 + a_{11}X_1, a_{12} + a_{15}X_1\right)  \notag
\end{align}}}
Then from $(\ref{pr9})$ and $(\ref{pr20})$ we see that
\begin{align*}
\dim_\mathbb{C} \mathbb{H}^0\left( \mathbb{P}_\mathbb{C}^3, \Theta_{\mathbb{P}_\mathbb{C}^3}^\bullet \right)&= \textnormal{ker}\left( H^0\left( \mathbb{P}_\mathbb{C}^3, \Theta_{\mathbb{P}_\mathbb{C}^3}  \right)\to H^0\left( \mathbb{P}_\mathbb{C}^3, \frac{\Theta_{\mathbb{P}_\mathbb{C}^3}}{\mathcal{O}_{\mathbb{P}_\mathbb{C}^3}} \right) \bigoplus H^0\left( \mathbb{P}_\mathbb{C}^1, \mathcal{O}_{\mathbb{P}_\mathbb{C}^1}(1)\bigoplus \mathcal{O}_{\mathbb{P}_\mathbb{C}^1}(1) \right) \right) = 4 \\
\dim_\mathbb{C} \mathbb{H}^1\left( \mathbb{P}_\mathbb{C}^3, \Theta_{\mathbb{P}_\mathbb{C}^3}^\bullet \right) & = 14 - ( 15- 4) = 3
\end{align*}
By comparing $(\ref{pr23})$ and $(\ref{pr24})$, we have
\begin{align}
\mathbb{H}^1\left( \mathbb{P}_\mathbb{C}^3  , \Theta_{\mathbb{P}_\mathbb{C}^3}^\bullet \right)&\cong \frac{\textnormal{ker}\left( H^0\left( \mathbb{P}_\mathbb{C}^3, \frac{\Theta_{\mathbb{P}_\mathbb{C}^3}}{\mathcal{O}_{\mathbb{P}_\mathbb{C}^3}} \right) \bigoplus H^0\left( \mathbb{P}_\mathbb{C}^1, \mathcal{O}_{\mathbb{P}_\mathbb{C}^1}(1)\bigoplus \mathcal{O}_{\mathbb{P}_\mathbb{C}^1}(1) \right) \to H^0\left( \mathbb{P}_\mathbb{C}^1, \mathcal{O}_{\mathbb{P}_\mathbb{C}^1}(1)\bigoplus \mathcal{O}_{\mathbb{P}_\mathbb{C}^1}(1) \right) \right)}{\textnormal{im}\left( H^0\left( \mathbb{P}_\mathbb{C}^3, \Theta_{\mathbb{P}_\mathbb{C}^3}  \right)\to H^0\left( \mathbb{P}_\mathbb{C}^3, \frac{\Theta_{\mathbb{P}_\mathbb{C}^3}}{\mathcal{O}_{\mathbb{P}_\mathbb{C}^3}} \right) \bigoplus H^0\left( \mathbb{P}_\mathbb{C}^1, \mathcal{O}_{\mathbb{P}_\mathbb{C}^1}(1)\bigoplus \mathcal{O}_{\mathbb{P}_\mathbb{C}^1}(1) \right) \right)} \notag \\
&\cong \textnormal{Span}_\mathbb{C} \left\langle  - X_2 \frac{\partial}{\partial X_3} ,  \left( X_1 , 1 \right) , \left(1, 0 \right) \right\rangle \label{pr25}
\end{align}
We will construct a foliated analytic family of simultaneous deformations of $\left\{ \xi_2=\xi_3=0 \right\} \subset \left( \mathbb{P}_\mathbb{C}^3, T_0 \right)$ such that the characteristic map is an isomorphism onto $\mathbb{H}^1\left( \mathbb{P}_\mathbb{C}^3, \Theta_{\mathbb{P}_\mathbb{C}^3}^\bullet \right)$ at the distinguished point. We have a family of simultaneous deformations of $\left\{ \xi_2 = \xi_3=0 \right\} \subset \left( \mathbb{P}_\mathbb{C}^3, T_0 \right)$ over $\textnormal{Spec}\left( \mathbb{C}\left[t_1,t_2,t_3\right]/ \left(t_1^2,t_2^2,t_3^2, t_1t_2, t_1t_3, t_2t_3 \right) \right)$ from $(\ref{pr25})$ which is  defined by
\begin{align*}
\tilde{T}&=\left( 1 - X_1X_2 \right)\frac{\partial}{\partial X_1} +\left( X_3- X_2^2 \right)\frac{\partial}{\partial X_2} +\left(- X_2X_3  + t_3X_2 \right) \frac{\partial}{\partial X_3}\\
X_2&= t_2 + t_1 X_1 \\
X_3&= t_1
\end{align*}
Then we have
\begin{align*}
&\left[ \left( 1 - X_1X_2 \right)\frac{\partial}{\partial X_1} +\left( X_3- X_2^2 \right)\frac{\partial}{\partial X_2} +\left(- X_2X_3  + t_3X_2 \right) \frac{\partial}{\partial X_3} , X_2- t_2  - t_1 X_1 \right]|_{X_2=t_2 + t_1 X_1 , X_3= t_1 }\\
&= X_3 - X_2^2 - t_1 + t_1 X_1 X_2|_{X_2= t_2 + t_1 X_1 ,  X_3= t_1 } = t_1 - t_2^2 - 2 t_1 t_2 X_1 - t_1^2 X_1^2 - t_1 + t_1t_2+ t_1^2 X_1^2 \\
&= t_1t_2 - t_2^2 - 2 t_1t_2 X_1
\end{align*}
and
\begin{align*}
&\left[  \left( 1 - X_1X_2 \right)\frac{\partial}{\partial X_1} +\left( X_3- X_2^2 \right)\frac{\partial}{\partial X_2} +\left(- X_2X_3  + t_3X_2 \right) \frac{\partial}{\partial X_3}, X_3- t_1 \right] |_{X_2=t_2 + t_1 X_1, X_3= t_1}\\
&= - X_2X_3 + t_3 X_2|_{X_2=t_2 + t_1 X_1, X_3= t_1}= - t_1t_2 - t_1^2 X_1+ t_2t_3 + t_1t_3X_1 = \left(t_2t_3 - t_1t_2 \right) +\left( t_1t_3 - t_1^2 \right)X_1
\end{align*}
So $\left( t_1t_2 - t_2^2 - 2t_1t_2 X_1 ,  t_2t_3 - t_1t_2 + \left( t_1t_3 - t_1^2 \right)X_1 \right)$ is the obstruction of lifting. We note that from $(\ref{pr21})$ and $(\ref{pr22})$ $\left( t_1t_2 - t_2^2 - 2t_1t_2 X_1 ,  t_2t_3 - t_1t_2 + \left( t_1t_3 - t_1^2 \right)X_1 \right)$ is the image of $\left( t_1t_2 - t_2^2 - 2 t_1t_2 X_1 \right)\frac{\partial}{\partial X_2}+ \left( t_2t_3 - t_1t_2 + \left(t_1t_3- t_1^2 \right) X_1 \right)\frac{\partial}{\partial X_3}$ under $H^0\left( \mathbb{P}_\mathbb{C}^3, \frac{\Theta_{\mathbb{P}_\mathbb{C}^3}}{\mathcal{O}_{\mathbb{P}_\mathbb{C}^3}} \right) \bigoplus H^0\left( \mathbb{P}_\mathbb{C}^1, \mathcal{O}_{\mathbb{P}_\mathbb{C}^1}(1)\bigoplus \mathcal{O}_{\mathbb{P}_\mathbb{C}^1}(1) \right) \to  H^0\left( \mathbb{P}_\mathbb{C}^1, \mathcal{O}_{\mathbb{P}_\mathbb{C}^1}(1) \bigoplus \mathcal{O}_{\mathbb{P}_\mathbb{C}^1}(1) \right)$, so that we take
{\small{\begin{align*}
\tilde{T}&=  \left( 1 - X_1X_2 \right)\frac{\partial}{\partial X_1} +\left( X_3- X_2^2   - \left( t_1t_2 - t_2^2\right) + 2 t_1t_2 X_1 \right)\frac{\partial}{\partial X_2} +\left(- X_2X_3  + t_3X_2 - \left( t_2t_3-  t_1t_2 \right) -  \left(t_1t_3 - t_1^2 \right)X_1 \right) \frac{\partial}{\partial X_3} \\
X_2&= t_2 + t_1 X_1 \\
X_3&= t_1
\end{align*}}}
Then we have
{\tiny{\begin{align*}
&\left[  \left( 1 - X_1X_2 \right)\frac{\partial}{\partial X_1} +\left( X_3- X_2^2   - \left( t_1t_2 - t_2^2\right) + 2 t_1t_2 X_1 \right)\frac{\partial}{\partial X_2} +\left(- X_2X_3  + t_3X_2 - \left( t_2t_3-  t_1t_2 \right) -  \left(t_1t_3 - t_1^2 \right)X_1 \right) \frac{\partial}{\partial X_3} ,X_2- t_2- t_1X_1 \right]|_{X_2= t_2+ t_1 X_1, X_3=t_1}\\
&=  X_3- X_2^2   - \left( t_1t_2 - t_2^2\right) + 2 t_1t_2 X_1 - t_1 + t_1 X_1X_2|_{X_2=t_2 + t_1 X_1, X_3=t_1} = 0
\end{align*}}}
and we have
{\tiny{\begin{align*}
&\left[  \left( 1 - X_1X_2 \right)\frac{\partial}{\partial X_1} +\left( X_3- X_2^2   - \left( t_1t_2 - t_2^2\right) + 2 t_1t_2 X_1 \right)\frac{\partial}{\partial X_2} +\left(- X_2X_3  + t_3X_2 - \left( t_2t_3-  t_1t_2 \right) -  \left(t_1t_3 - t_1^2 \right)X_1 \right) \frac{\partial}{\partial X_3},  X_3 - t_1 \right]|_{X_2= t_2 + t_1 X_1, X_3 = t_1}\\
&= - X_2X_3  + t_3X_2 - \left( t_2t_3-  t_1t_2 \right) -  \left(t_1t_3 - t_1^2 \right)X_1 |_{X_2= t_2+ t_1 X_1, X_3 = t_1 } = 0
\end{align*}}}
Hence $\tilde{T}$ and $X_2= t_2 + t_1 X_1$ and $X_3= t_1$ defines a foliated analytic family of simultaneous deformations of $\left\{ \xi_2= \xi_3=0 \right\} \subset \left( \mathbb{P}_\mathbb{C}^3, T_0 \right)$ over a sufficiently small neighborhood $B$ of $(0,0,0)\in \mathbb{C}^3$ such that the characteristic map is an isomorphism onto $\mathbb{H}^1\left( \mathbb{P}_\mathbb{C}^3, \Theta_{\mathbb{P}_\mathbb{C}^3}^\bullet \right)$ at $(t_1,t_2,t_3)=(0,0,0)\in B$.

\end{example}

\begin{example}

Let $\left[ \xi_0:\xi_1:\xi_2: \xi_3\right]$ be the homogeneous coordinates on the projective $3$-space $\mathbb{P}_\mathbb{C}^3$. Let $\left( \mathbb{P}_\mathbb{C}^3, T \right)$ be a foliated projective $3$-space of degree $6$, where $T\in H^0\left( \mathbb{P}_\mathbb{C}^3, \Theta_{\mathbb{P}_\mathbb{C}^3}(5) \right)$ is given by 
\begin{align*}
T=\xi_0^2\xi_2^2\xi_3^2\frac{\partial}{\partial \xi_1} + \xi_0^2\xi_1^2\xi_3^2\frac{\partial}{\partial \xi_2} + \xi_0^2\xi_1^2\xi_2^2\frac{\partial}{\partial \xi_3} - 3 \xi_1^2 \xi_2^2 \xi_3^2 \frac{\partial}{\partial \xi_0}
\end{align*}
Then $\xi_0^3+ \xi_1^3+ \xi_2^3 + \xi_3^3=0$ defines a non-degenerate invariant complex submanifold $X$ of $\left( \mathbb{P}_\mathbb{C}^3, T \right)$, and  $T$ induces the foliation $\Theta_{\mathcal{F}_{X_0}}$ on $X$. We will compute $\mathbb{H}^1\left( X, \Theta_{\mathcal{F}_0}^\bullet \right)$ and show that $\left(X, \Theta_{\mathcal{F}_{X_0}} \right)$ is unobstructed in foliated deformations by computing $\mathbb{H}^2\left( X, \Theta_{\mathcal{F}_{X_0}}^\bullet  \right)=0$ $($our computation below actually holds for any foliation of degree $6$ on $\mathbb{P}_\mathbb{C}^3$ which defines a foliation on $\xi_0^3+\xi_1^3+\xi_2^3 + \xi_3^3=0$$)$. We have the exact sequence defined by $T$
\begin{align*}
0 \to \mathcal{O}_{\mathbb{P}_\mathbb{C}^3}(-5) \xrightarrow{T} \Theta_{\mathbb{P}_\mathbb{C}^3} \to \frac{\Theta_{\mathbb{P}_\mathbb{C}^3}}{\mathcal{O}_{\mathbb{P}_\mathbb{C}^3}(-5)}\to 0
\end{align*}
Since $X$ is invariant with respect to $T$, we have the exact sequence
\begin{align}\label{g1}
0\to \Theta_{\mathcal{F}_{X_0}}:=\mathcal{O}_X(-5) \xrightarrow{T|_X} \Theta_X \to \frac{\Theta_X}{\mathcal{O}_X(-5)}\to 0
\end{align}
We will compute the cohomology groups $\mathbb{H}^i\left(X, \Theta_{\mathcal{F}_{X_0}}^\bullet \right),i=0,1,2$ of the leaf complex
\begin{align*}
\Theta_{\mathcal{F}_{X_0}}^\bullet : \Theta_X \to \mathcal{O}_X(5) \otimes \frac{\Theta_X}{\mathcal{O}_X(-5)}  \to 0
\end{align*}
Since we have $\dim_\mathbb{C} H^0\left(X, \Theta_X \right)=0, \dim_\mathbb{C} H^0\left( X, \Theta_X \right)=4$ and $\dim_\mathbb{C} H^2\left( X, \Theta_X \right)=0$ $($see \textnormal{\cite{Kod05} p.309}$)$, we see that
\begin{align*}
\mathbb{H}^0\left( X, \Theta_{\mathcal{F}_{X_0}}^\bullet \right)&= 0 \\
\mathbb{H}^0\left(X, \Theta_{\mathcal{F}_{X_0}}^\bullet \right) &\cong H^0\left(X,  \mathcal{O}_X(5)\otimes \frac{\Theta_X}{\mathcal{O}_X(-5)} \right) \bigoplus \textnormal{ker}\left(H^1\left( X, \Theta_X \right)\to H^1\left( X, \mathcal{O}_X(5)\otimes \frac{\Theta_X}{\mathcal{O}_X(-5)} \right) \right)\\
\mathbb{H}^2\left( X, \Theta_{\mathcal{F}_{X_0}}^\bullet \right) & \cong H^1\left(X, \mathcal{O}_X(5)\otimes \frac{\Theta_X}{\mathcal{O}_X(-5)} \right)
\end{align*}

We claim that $ H^1\left( X, \mathcal{O}_X(5)\otimes \frac{\Theta_X}{\mathcal{O}_X(-5)} \right)=0 $. In particular, $\mathbb{H}^2\left( X, \Theta_{\mathcal{F}_{X_0}}^\bullet \right)=0$, so that $\left( X, \Theta_{\mathcal{F}_{X_0}} \right)$ is unobstructed in foliated deformations. In fact, from $(\ref{g1})$ we have the exact sequence
\begin{align*}
 0\to \mathcal{O}_X \to \Theta_X(5)\to \mathcal{O}_X(5)\otimes \frac{\Theta_X}{\mathcal{O}_X(-5)} \to 0
\end{align*}
which induces the long exact sequence
\begin{align}
0 &\to H^0\left( X, \mathcal{O}_X \right) \cong \mathbb{C} \to H^0\left(X, \Theta_X(5) \right) \to H^0\left( X , \mathcal{O}_X(5)\otimes \frac{\Theta_X}{\mathcal{O}_X(-5)} \right) \label{g7} \\
 &\to H^1\left(X, \mathcal{O}_X \right) =0 \to H^1\left( X, \Theta_X(5) \right)\to H^1\left( X, \mathcal{O}_X(5)\otimes \frac{\Theta_X}{\mathcal{O}_X(-5)} \right) \to H^2\left( X, \mathcal{O}_X \right) = 0 \notag
\end{align}
We will show that $H^1\left( X, \Theta_X(5)\right)=0$ which implies our claim. In fact, by duality  $H^1\left(X, \Theta_X(5) \right)\cong H^1\left(X, \Omega_X^1(-6) \right)^*$ $($we note that $K_X\cong \mathcal{O}_X(-1)$$)$ and we have two exact sequences $($see \textnormal{\cite{KS53} p.1274 (3) and (4)}$)$
\begin{align}
0 \to \Omega_{\mathbb{P}_\mathbb{C}^2}''^2(-3) \to \Omega_{\mathbb{P}_\mathbb{C}^3}^2 (-3) \to \mathcal{O}_X(-4)\to 0 \label{g2} \\
0\to \Omega_{\mathbb{P}_\mathbb{C}^3}^2(-6)\to \Omega_{\mathbb{P}_\mathbb{C}^3}''^2(-3)\to \Omega_X^1(-6)\to 0 \label{g3}
\end{align}
From $(\ref{g3})$ we infer that 
\begin{align*}
H^1\left( \mathbb{P}_\mathbb{C}^3, \Omega_{\mathbb{P}_\mathbb{C}^3}''^2 (-3) \right)\cong H^1\left( X, \Omega_X^1(-6) \right)
\end{align*}
From $(\ref{g2})$ we infer that $H^1\left( \mathbb{P}_\mathbb{C}^3, \Omega_{\mathbb{P}_\mathbb{C}^3}''^2 (-3) \right)=0$. Hence $H^1\left( X , \Theta_X(5) \right)=0$. Next we compute $H^0\left(X, \Theta_X(5) \right)$. From the exact sequence $0\to \Theta_X\to \Theta_{\mathbb{P}_\mathbb{C}^3}|_X \to \mathcal{O}_X(3)\to 0$, we have 
\begin{align}\label{g5}
0\to \Theta_X(5)\to \Theta_{\mathbb{P}_\mathbb{C}^3}(5)|_X \to \mathcal{O}_X(8)\to 0
\end{align}
Since $H^1\left(X, \Theta_X(5) \right)=0$, $(\ref{g5})$ induces the exact sequence
\begin{align}\label{g6}
0\to H^0\left( X, \Theta_X(5) \right) \to H^0\left(X, \Theta_{\mathbb{P}_\mathbb{C}^3}(5)|_X \right) \to H^0\left(X, \mathcal{O}_X(8) \right) \to 0
\end{align} 
Let us compute $H^0\left(X, \mathcal{O}_X(8) \right)$. From the exact sequence $0\to \mathcal{O}_{\mathbb{P}_\mathbb{C}^3}(5)\to \mathcal{O}_{\mathbb{P}_\mathbb{C}^3}(8) \to \mathcal{O}_X(8)\to 0$, we see that $\dim_\mathbb{C} H^0\left( X, \mathcal{O}_X(8) \right) = 165- 56= 109$. On the other hand, let us compute $H^0\left(X, \Theta_{\mathbb{P}_\mathbb{C}^3}(5)|_X \right)$. From the exact sequence $0\to \Theta_{\mathbb{P}_\mathbb{C}^3}(2)\to \Theta_{\mathbb{P}_\mathbb{C}^3}(5)\to \Theta_{\mathbb{P}_\mathbb{C}^3}(5)|_X\to 0$ and $0\to \mathcal{O}_{\mathbb{P}_\mathbb{C}^3}(m)\to \bigoplus^4 \mathcal{O}_{\mathbb{P}_\mathbb{C}^3}(m+1)\to \Theta_{\mathbb{P}_\mathbb{C}^3}(m)\to 0$ for $m=5,2$, we infer that $\dim_\mathbb{C} H^0\left( X, \Theta_{\mathbb{P}_\mathbb{C}^3}(5)|_X \right)=210$. Then from $(\ref{g6})$ we see that $\dim_\mathbb{C} H^0\left( X, \Theta_X(5) \right)= 210-109= 101$, and from $(\ref{g7})$ we have
\begin{align*}
\dim_\mathbb{C} H^0 \left( X, \mathcal{O}_X(5)\otimes \frac{\Theta_X}{\mathcal{O}_X(-5)}\right)= 101-1 =100
\end{align*}
Hence we have
\begin{align*}
\dim_\mathbb{C} \mathbb{H}^1\left(X, \Theta_{\mathcal{F}_{X_0}}^\bullet \right) = 100 + 4 = 104
\end{align*}
\end{example}

\section{Blow-ups along non-degenerate invariant complex submanifolds}

\begin{proposition} \label{bw1}
Let $V$ be a non-degenerate invariant complex submanifold of a compact foliated complex manifold $\left( W, \Theta_{\mathcal{F}_0}\right)$ with $\Theta_{\mathcal{F}_0}$ locally free such that $0\to \Theta_{\mathcal{F}_{V_0}}:=\Theta_{\mathcal{F}_0}|_V\to \Theta_W|_V \to \frac{\Theta_W}{\Theta_{\mathcal{F}_0}}|_V \to 0$ is exact. Assume that the codimension of $V$ in $W$ is at least $2$, and the singular locus of $\Theta_{\mathcal{F}_{V_0}}$ on $V$ has codimension at least $1$. Let $\pi: \tilde{W} \to W$ be the the blow up $\tilde{W}$ of $W$ along $V$. Then the exact sequence $0\to \Theta_{\mathcal{F}_0} \to \Theta_W \to \frac{\Theta_W}{\Theta_{\mathcal{F}_0}}\to 0$ induces the exact sequence $0\to \pi^* \Theta_{\mathcal{F}_0}\to \pi^* \Theta_W \to \pi^*\left(\frac{\Theta_W}{\Theta_{\mathcal{F}_0}} \right) \to 0$ and $\tilde{W}$ is a compact foliated complex manifold with a foliation defined by $\pi^* \Theta_{\mathcal{F}_0}\subset \Theta_{\tilde{W}}\hookrightarrow \pi^* \Theta_W$, and the exceptional divisor $E$ is a non-degenerate invariant complex submanifold of $\left( \tilde{W}, \pi^* \Theta_{\mathcal{F}_0} \right)$. Moreover if $\Theta_{\mathcal{F}_{V_0}}$ defines a foliation on $V$, i.e. the singular locus of $\Theta_{\mathcal{F}_{V_0}}$ has codimension at least $2$, then the exceptional divisor $E$ is a compact foliated complex manifold, i.e. the singular locus of $\pi^* \Theta_{\mathcal{F}_0}|_E$ has codimension at least $2$.
\end{proposition}

\begin{remark}\label{bcc1}
We shall say that a holomorphic map $f:\left(X, \Theta_{\mathcal{F}_0} \right)\to \left( Y, \Theta_{\mathcal{G}_0}\right)$ between foliated complex manifolds with $\Theta_{\mathcal{F}_0}$ and $\Theta_{\mathcal{G}_0}$ locally free is a holomorphic foliated map if the exact sequence $0\to \Theta_{\mathcal{G}_0}\to \Theta_Y \to \frac{\Theta_Y}{\Theta_{\mathcal{G}_0}}\to 0$ induces the exact sequence $0\to f^* \Theta_{\mathcal{G}_0}\to f^* \Theta_Y \to f^* \left( \frac{\Theta_Y}{\Theta_{\mathcal{G}_0}} \right) \to 0 $ and $F\left( \Theta_{\mathcal{F}_0}\right)\subset f^* \Theta_{\mathcal{G}_0}$ under the natural map $F:\Theta_X\to f^* \Theta_Y$. We shall study deformations of  holomorphic foliated maps in \textnormal{Part VII}. We note that the blow up $\pi:\left( \tilde{W}, \pi^* \Theta_{\mathcal{F}_0} \right)\to \left( W, \Theta_{\mathcal{F}_0} \right) $ in \textnormal{Proposition \ref{bw1}} is a holomorphic foliated map. In \textnormal{Part VII}, we do not necessarily assume that codimension of singular locus of a foliation is at least $2$. Then $\left( E, \pi^* \Theta_{\mathcal{F}_0} |_E \right) \to \left( V, \Theta_{\mathcal{F}_{V_0}}\right)$ is a holomorphic foliated map.
\end{remark}

\begin{proof}[Proof of \textnormal{Proposition \ref{bw1}}]

Let $\left( W, \Theta_{\mathcal{F}_0}\right)$ be a compact foliated complex manifold with $\Theta_{\mathcal{F}_0}$ locally free, so that we have the exact sequence $0\to \Theta_{\mathcal{F}_0} \to \Theta_W \to \frac{\Theta_W}{\Theta_{\mathcal{F}_0}}\to 0$. Let $V$ be a non-degenerate invariant complex submanifold of $\left( W, \Theta_{\mathcal{F}_0} \right)$ of codimension at least $2$, so that we have the exact sequence $0\to \Theta_{\mathcal{F}_0}|_V \to \Theta_W|_V \to \frac{\Theta_W}{\Theta_{\mathcal{F}_0}}|_V\to 0$ and $\Theta_{\mathcal{F}_{V_0}}:= \Theta_{\mathcal{F}_0}|_V \subset \Theta_V$. Assume that the singular locus of $\Theta_{\mathcal{F}_{V_0}}$ on $V$ has codimension at least $1$. Let $\pi: \tilde{W} \to W$ be the blow up along the complex submanifold $V$. Then we have the exact sequence induced from $0\to \Theta_{\mathcal{F}_0}\to \Theta_W \to \frac{\Theta_W}{\Theta_{\mathcal{F}_0}}\to 0$:
\begin{align}\label{fb1}
\pi^* \Theta_{\mathcal{F}_0} \to \pi^* \Theta_W \to \pi^* \left( \frac{\Theta_W}{\Theta_{\mathcal{F}_0}} \right) \to 0
\end{align}
and we have the exact sequence
\begin{align*}
0\to \Theta_{\tilde{W}}\xrightarrow{\Pi} \pi^* \Theta_W \to \mathcal{N}_\pi \to 0
\end{align*}
where $\Pi: \Theta_{\tilde{W}}\to \pi^* \Theta_W$ is the natural map induced from $\pi: \tilde{W} \to W$.

We claim that the image of $\pi^* \Theta_{\mathcal{F}_0}$ in $\pi^* \Theta_W$ from $(\ref{fb1})$ lies in $\Pi\left( \Theta_{\tilde{W}} \right)$. Let $U$ be a local coordinate neighborhood of $W$ with coordinates $\left(w_1,...,w_k, v_{k+1},..., v_n \right)$ such that $V\cap U$ is defined by $v_{k+1}=\cdots = v_n=0$. Let $\Theta_{\mathcal{F}_0}$ on $U$ is generated by $T^1(w, v),..., T^p(w,v)\in \Gamma\left( U, \Theta_W \right)$ such that $T^\alpha(w,v)$ is of the form
\begin{align*}
T^\alpha(w,v)=  \sum_{\beta=1}^{k} T^{\alpha\beta}(w,v)\frac{\partial}{\partial w_\beta} + \sum_{\eta=k+1}^n V^{\alpha \eta}(w,v) \frac{\partial}{\partial v_\eta}
\end{align*}
for some holomorphic functions $T^{\alpha\beta}(w,v), V^{\alpha \eta}(w,v)\in \Gamma\left( U, \mathcal{O}_W \right)$. Since $v_{k+1}=\cdots = v_n=0$ is an invariant complex submanifold of $U$, we can write
\begin{align*}
T^\alpha(w,v)=  \sum_{\beta=1}^{k} T^{\alpha\beta}(w,v)\frac{\partial}{\partial w_\beta} + \sum_{\eta=k+1}^n \left( \sum_{\delta=k+1}^n  v_{\delta} G^{\eta \delta}(w,v )   \right) \frac{\partial}{\partial v_\eta}
\end{align*}
for some holomorphic functions $G^{\eta \delta}(w,v)\in \Gamma\left( U, \mathcal{O}_W \right)$. Then $\Theta_{\mathcal{F}_0}|_V$ is generated by $T^1(w,0),..., T^p(w,0)$ such that 
\begin{align*}
T^\alpha(w,0)= \sum_{\beta=1}^k T^{\alpha\beta}(w,0)\frac{\partial}{\partial w_\beta}
\end{align*}
We assume that the zero locus of $T^1(w,0)\wedge \cdots \wedge T^p(w,0)$ has codimension at least $1$ on $V$. Let $\left[l_{k+1},...,l_n\right]$ be the homogeneous coordinate on $\mathbb{P}_\mathbb{C}^{n-k-1}$, and let
\begin{align*}
\tilde{U} \subset U \times \mathbb{P}_\mathbb{C}^{n-k-1}
\end{align*}
be the smooth variety defined by the relations
\begin{align*}
\tilde{U}=\left\{ \left(w_1,...w_k, v_{k+1},..., v_n, l_{k+1},...,l_n)\right| v_il_j= v_jl_i , k+1\leq i,j \leq n    \right\}
\end{align*}
Then $\tilde{U}$ has a covering $\left\{ \tilde{U}^j | k+1 \leq j \leq n \right\}$ such that $\tilde{U}^j$ has a local coordinate $\left( z_1,..., z_k, \xi_j^{k+1},..,\xi_j^{j-1},  \overbrace{y_j}^{j-\textnormal{th}}, \xi_j^{j+1},..., \xi_j^n  \right)$ where
\begin{align*}
z_i&= w_i,\,\,\,\,\,\,\,\,i=1,...,k\\
\xi_j^i&= \frac{l_i}{l_j},\,\,\,\,\,\,\,i=k+1,..., j-1, \widehat{j},j+1,..., n\\
 y_j&= v_j
\end{align*}
Then $\pi$ on $\tilde{U}^j$ is described by 
\begin{align*}
\pi: \left( z_1,...,z_k, \xi_j^{k+1},..., y_j, ..., \xi_j^\eta \right) \mapsto \left( w_1,...,w_k, v_{k+1},...,v_n \right) = \left( z_1,...,z_k, y_j \xi_j^{k+1},...y_j\xi_j^{j-1},  y_j,y_j\xi_j^{j+1}.., y_j\xi_j^n \right)
\end{align*}
and $y_j=0$ defines $E \cap \tilde{U}^j$ where $E$ is the exceptional divisor of the blow up $\pi:\tilde{W}\to W$. Let $\Pi: \Theta_{\tilde{W}} \to \pi^* \Theta_W$ be the natural map. We describe $\Pi$ in terms of the basis $\frac{\partial}{\partial z_1},..., \frac{\partial}{\partial z_k}, \frac{\partial}{\partial \xi_j^{k+1}},..., \frac{\partial}{\partial \xi_j^{j-1}}, \frac{\partial}{\partial y_j}, \frac{\partial}{\partial \xi_j^{j+1}},..., \frac{\partial}{\partial \xi_j^n}$ on $\tilde{U}^j$. Then we have
\begin{align*}
\Pi\left(\frac{\partial}{\partial z_i}\right)&= \frac{\partial}{\partial w_i}\,\,\,\,\,\,\,\,i=1,...,k\\
\Pi\left(\frac{\partial}{\partial y_j}\right)&=\frac{\partial}{\partial v_j}+ \sum_{\gamma=k+1, \gamma\ne j}^n \xi_j^\gamma \frac{\partial}{\partial v_\gamma}\\
\Pi\left(\frac{\partial}{\partial \xi_j^\gamma}\right)&= y_j\frac{\partial}{\partial v_\gamma } ,\,\,\,\,\,\,\,\,\,\,\gamma=k+1,..., n, \gamma\ne j
\end{align*}

Now we show that the image of $\pi^* \Theta_{\mathcal{F}_0}$ in $\pi^* \Theta_W$ lies in $\Pi\left(\Theta_{\tilde{W}} \right)$. In fact, it is sufficient to check locally on $\tilde{U}^j$. We note that we have
{\small{\begin{align*}
\pi^* T^\alpha(w,v)&= \sum_{\beta=1}^k T^{\alpha\beta}\left( z_1,..., z_k, y_j\xi_j^{k+1},..., y_j,..., y_j\xi_j^n \right)\frac{\partial}{\partial w_\beta} + \sum_{\eta =k+1}^n  \left( \sum_{\delta=k+1}^n y_j \xi_j^\delta G^{\eta \delta}\left(  z_1,..., z_k, y_j\xi_j^{k+1},..., y_j,..., y_j\xi_j^n \right) \right) \frac{\partial}{\partial v_\eta}
\end{align*}}}

Here $\xi_j^j=1$. Then we define $\tilde{T}^{j\alpha}\in \Gamma\left( \tilde{U}^j, \Theta_{\tilde{W}}\right)$ by
\begin{align*}
&\tilde{T}^{j\alpha}= \sum_{\beta=1}^k T^{\alpha\beta}\left(z_1,...,z_k, y_j\xi_j^{k+1},..., y_j,..., y_j\xi_j^n \right) \frac{\partial}{\partial z_\beta} + \sum_{\delta=k+1}^n  y_j \xi_j^\delta G^{j \delta}\left(z_1,..,z_k, y_j\xi_j^{k+1},..., y_j,..., y_j\xi_j^n \right)\frac{\partial}{\partial y_j}\\
& + \sum_{\eta=k+1, \eta\ne j}^n  \left(\sum_{\delta=k+1}^n \xi_j^\delta G^{\eta \delta}\left(z_1,...,z_k, y_j\xi_j^{k+1},..., y_j,..., y_j\xi_j^n\right) - \sum_{\delta=k+1}^n  \xi_j^\delta\xi_j^\eta G^{j\delta}\left(z_1,...,z_k, y_j\xi_j^{k+1},..., y_j,..., y_j\xi_j^n \right) \right) \frac{\partial}{\partial \xi_j^\eta}
\end{align*}
Then we see that
\begin{align*}
\Pi\left(\tilde{T}^{j\alpha} \right)= \pi^* T^\alpha
\end{align*}
Since $\Theta_{\mathcal{F}_0}$ is generated by $T^\alpha(w,z),\alpha=1,...,p$ on $U$, and $\tilde{U}$ is covered by $\tilde{U}^j$, the image of $\pi^* \Theta_{\mathcal{F}_0}$ in $\pi^* \Theta_W$ lies in $\Pi\left(  \Theta_W \right)$. We note that $\tilde{T}^{j\alpha}$ and $T^\alpha$ is $\pi$-related.

Next we claim that $\pi^* \Theta_{\mathcal{F}_0} \to \pi^* \Theta_W$ is injective. It is sufficient to show that $\pi^* T^1,..., \pi^* T^p$ are linearly independent on $\tilde{U}^j$. Assume that
\begin{align} \label{bw2}
\sum_{\alpha=1}^p b_j^\alpha\left( z_1,..., z_k, \xi_j^{k+1},..., y_j, ..., \xi_j^n \right)\pi^* T^\alpha =0
\end{align}
If $y_j\ne 0$, then  $\pi$ is biholomophic, and we may identity $\pi^* T^\alpha$ with $T^\alpha$. We note that 
\begin{align*}
\pi^{-1}: \left(w_1,..., w_k, v_{k+1},..., v_n \right) \mapsto \left(z_1,...,z_k, \xi_j^{k+1},..., y_j,..., \xi_j^n \right)=\left( w_1,...,w_k, \frac{v_{k+1}}{v_j},..., v_j,..., \frac{v_n}{v_j} \right)
\end{align*}
Then we have 
\begin{align*}
\sum_{\alpha=1}^p b_j^\alpha\left( w_1,..., w_k, \frac{v_{k+1}}{v_j},..., v_j, ..., \frac{v_n}{v_j} \right) T^\alpha =0
\end{align*}
Since $T^1,..., T^p$ are linearly independent, we have $b_j^\alpha\left( w_1,..., w_k \frac{v_{k+1}}{v_j},..., v_j, ..., \frac{v_n}{v_j} \right)=0$, i.e. we have
\begin{align*}
 b_j^\alpha\left( z_1,..., z_k, \xi_j^{k+1},..., y_j, ..., \xi_j^n \right)=0 \,\,\,\,\,\,\,\textnormal{on}\,\,\,\,\, y_j \ne 0.
\end{align*}
On the other hand, by restricting $(\ref{bw2})$ to $y_j=0$, we have
\begin{align*}
0&=\sum_{\alpha=1}^p b_j^\alpha\left(z_1,...,z_k, \xi_j^{k+1},..., 0, ..., \xi_j^n \right) \left( \sum_{\beta=1}^k T^{\alpha\beta}\left(z, 0 \right)\frac{\partial}{\partial w_\beta} \right)\\
&\Longrightarrow 0=\sum_{i_{k+1},..., i_n=0}^{\infty} \left( \xi_j^{k+1}\right)^{i_{k+1}}\cdots \left( \xi_j^n\right)^{i_n} \left(\sum_{\alpha=1}^p b_{j i_{k+1},..., i_{n}}^\alpha (w)  T^\alpha(w,0) \right)\\
&\Longrightarrow \sum_{\alpha=1}^p b_{j i_{k+1},..., i_n}^\alpha (w)T^\alpha(w,0) = 0
\end{align*}
Since $V$ is a non-degenerate invariant complex submanifold, $T^1(w,0),..., T^p(w,0)$ are linearly independent, so that $b_{ji_{k+1},..., i_n}^\alpha(w)=0$. Hence $b_j^\alpha\left(z_1,...,z_k, \xi_j^{k+1},..., y_j, ..., \xi_j^n\right)=0$.  This proves the claim. 

Hence since $\Pi : \Theta_{\tilde{W}}\to \pi^* \Theta_W$ is injective, and $\pi^* \Theta_{\mathcal{F}_0}$ is in the image of $\Pi\left( \Theta_{\tilde{W}}\right)$, we may consider $\pi^* \Theta_{\mathcal{F}_0}$ as a subsheaf $\Theta_{\tilde{W}}$, i.e. $\pi^* \Theta_{\mathcal{F}_0} \subset \Theta_{\tilde{W}}$. We note that $\pi^* \Theta_{\mathcal{F}_0}$ in $\Theta_{\tilde{W}}$ is locally generated by $\tilde{T}^{j1},..., \tilde{T}^{j p}$ on $\tilde{U}^j$. Since $\tilde{T}^{j\alpha}$ and $T^\alpha$ is $\pi$-related, we see that $\left[ \Theta_{\mathcal{F}_0}, \Theta_{\mathcal{F}_0} \right] \subset \Theta_{\mathcal{F}_0}$ implies that $\left[ \pi^* \Theta_{\mathcal{F}_0}, \pi^* \Theta_{\mathcal{F}_0} \right] \subset \pi^* \Theta_{\mathcal{F}_0}$ in $\Theta_{\tilde{W}}$. Since $\pi: \tilde{W} - E \to W- V$ is an isomorphism, we see that $\pi^* \Theta_{\mathcal{F}_0} = \Theta_{\mathcal{F}_0}$ on $\tilde{W}-E$ defines a foliation, i.e. the singular locus of $\pi^*\Theta_{\mathcal{F}_0}$ on $\tilde{W}-E$ has codimension at least $2$. On the other hand, on $\tilde{U}^j \cap E$, we see that $\tilde{T}^{j\alpha}|_{y_j=0}$ is of the form $T^\alpha(z,0) + \sum_{\eta=k+1, \eta\ne j}^n W_j^\eta \left(z,\xi_j^{k+1},.., \widehat{y}_j , ... , \xi_j^n\right)\frac{\partial}{\partial \xi_j^\eta}$, and zero locus of $T^1(w,0) \wedge \cdots \wedge T^p (w,0)$ has codimension at least $1$ on $U\cap V$ by our assumption, the singular locus of $\pi^* \Theta_{\mathcal{F}_0}$ on $\tilde{U}^j \cap E$ has codimension at least $1$. Hence the singular locus of $\pi^* \Theta_{\mathcal{F}_0}$ on $\tilde{W}$  has codimension at least $2$. 

We note that $y_j=0$ is an invariant complex submanifold on $\left(\tilde{U}^j, \pi^* \Theta_{\mathcal{F}_0}\right)$, so that $E$ is an invariant complex submanifold. We claim that $E$ is non-degenerate, i.e. $\tilde{T}^{j1}|_{y_j=0},... \tilde{T}^{jp}|_{y_j=0}$ are linearly independent. Since $T^\alpha(z,0)$ are linearly independent, we see that $\tilde{T}^{j\alpha}|_{y_j=0}$ are linearly independent by the above same argument. Moreover if the singular locus of $\Theta_{\mathcal{F}_0}|_V$ has codimension at least $2$, i.e the zeros of $T^1(z,0)\wedge \cdots \wedge T^p(z,0)$ has codimension at least $2$, then $\tilde{T}^{j1}|_{y_j=0}\wedge \cdots \wedge \tilde{T}^{j p}|_{y_j=0}$ has codimension at least $2$, i.e. $\pi^* \Theta_{\mathcal{F}_0}|_E$ defines a foliation on $E$.

\end{proof}

\begin{remark}
It is possible in \textnormal{Proposition \ref{bw1}} that the singular locus of $\Theta_{\mathcal{F}_0}|_V$ has codimension $1$ but $\pi^* \Theta_{\mathcal{F}_0}|_E$ defines a foliation on $E$, i.e. the singular locus of $\pi^* \Theta_{\mathcal{F}_0}|_E$ has codimension at least $2$ as in the following \textnormal{Example \ref{bw3}}, where $V\cong \mathbb{P}_\mathbb{C}^1$.
\end{remark}

\begin{example}\label{bw3}
Let us consider the foliated projective $3$-space $\left( \mathbb{P}_\mathbb{C}^3, \Theta_{\mathcal{F}_0}\cong \mathcal{O}_{\mathbb{P}_\mathbb{C}^3} \right)$ defined by the homogeneous vector field $T_0$ from \textnormal{Example \ref{pr2}}, and the non-degenerate  invariant complex submanifold defined by $\left\{\xi_2=\xi_3=0\right\} \cong \mathbb{P}_\mathbb{C}^1$ of $\left( \mathbb{P}_\mathbb{C}^3, T_0 \right)$ as in \textnormal{Example \ref{pr8}}. We keep the notations in $\textnormal{Example \ref{pr2}}$ and $\textnormal{Example \ref{pr8}}$. Let $\pi:\left(Y, \pi^* T_0 \right)\to \left( \mathbb{P}_\mathbb{C}^3, T_0 \right)$ be the blow-up $\left( \mathbb{P}_\mathbb{C}^3, T_0 \right)$ along $\xi_2=\xi_3=0$ and let $E\cong \mathbb{P}_\mathbb{C}^1 \times \mathbb{P}_\mathbb{C}^1$ be the non-degenerate invariant exceptional divisor of $\left( Y, \pi^* T_0 \right)$. We will study foliated deformations of $\left( Y, \pi^* T_0 \right)$ and foliated deformations of $\left(E, \pi^* T_0|_E \right)$ and study embedded deformations of $E$ in $\left(Y,\pi^* T_0 \right)$ and, by assuming the knowledge of $\textnormal{Part VI}$ on simultaneous deformations of compact foliated complex manifolds defined by locally free subsheaves of tangent sheaves and their non-degenerate invariant complex submanifolds,  we will study simultaneous deformations of $E\subset \left(  Y, \pi^* T_0 \right)$, and by assuming the knowledge of deformations of $($non-degenerate$)$ holomorphic foliated maps in \textnormal{Part VII}, we will study deformations of $\pi:\left( Y, \pi^* T_0 \right)\to \left( \mathbb{P}_\mathbb{C}^3, T_0 \right)$. We will show that foliated deformations of $\left(Y, \pi^* T_0 \right)$ and foliated deformations of $\left( E, \pi^* T_0|_E \right)$ and embedded deformations $E$ in $\left( Y, \pi^* T_0 \right)$ and simultaneous deformations of $E\subset \left( Y, \pi^* T_0 \right)$ are all unobstructed, but deformations of holomorphic foliated map $\pi$ is obstructed. In particular we show that $E$ is a stable non-degenerate invariant complex submanifold of $\left( Y, \pi^* T_0 \right)$ by Theorem of stability for deformations of non-degenerate  invariant complex submanifolds of foliated complex manifolds defined by locally free subsheaves of tangent sheaves in \textnormal{Part VI}.
 \end{example}

\subsection{Foliated deformations of the blow-up $\left(Y, \pi^* T_0 \right)$}\

We will describe the blow up $Y$ of $\mathbb{P}_\mathbb{C}^3$ along $\xi_2= \xi_3=0$. We keep the notations in \textnormal{Appendix \ref{cp7}}. We note that $\xi_2=\xi_3=0$ happens on $U_0$ and $U_1$. The blow up $\tilde{U}_0$ over $U_0$ along $X_2=X_3=0$ is given by $\tilde{U}_0=\left\{\left(X_1,X_2,X_3, \left[ l_2,l_3\right] \right) \in U_0\times \mathbb{P}_\mathbb{C}^1 | l_2X_3= l_3X_2\right\}$. $\tilde{U}_0$ has two affine coverings $U_{00}$ and $U_{01}$ where $U_{00}$ is given by $l_2\ne 0$, i.e. $U_{00}=\left\{ \left(X_1,X_2,w_3 X_2,\left[ 1, w_3= \frac{l_3}{l_2}\right] \right) \right\}$ and by setting $X_1= w_1, X_2= w_2$, we may identify $U_{00}\cong \left\{\left(w_1,w_2,w_3\right)\in \mathbb{C}^3 \right\}$. Then $\pi:U_{00}\to U_0$ is given by $\left( w_1, w_2, w_3 \right) \mapsto \left(X_1,X_2,X_3 \right)= \left(w_1, w_2, w_2 w_3\right)$ and the exceptional divisor $E$ on $U_{00}$ is given by $w_2=0$, and $U_{01}$ is given by $l_3\ne 0$, i.e. $U_{01}=\left\{ \left(X_1,w_2'X_3 ,X_3,\left[ w_2'=\frac{l_2}{l_3}, 1 \right] \right) \right\}$ and by setting $X_1= w_1', X_3= w_3'$, we may identify $U_{01}\cong \left\{\left( w_1', w_2', w_3'\right) \in \mathbb{C}^3 \right\}$. Then $\pi:U_{01}\to U_0$ is given by $\left( w_1', w_2', w_3' \right)\mapsto \left(X_1, X_2, X_3 \right) = \left( w_1', w_2'w_3', w_3' \right)$ and the exceptional divisor $E$ on $U_{01}$ is given by $w_3'=0$. The coordinate transformation between $U_{00}$ and $U_{01}$ is given by $\left( w_1,w_2,w_3\right)\mapsto \left( w_1', w_2', w_3' \right)=\left( w_1, \frac{1}{w_3}, w_2w_3 \right)$. On the other hand, the blow up $\tilde{U}_1$ over $U_1$ along $Y_2=Y_3=0$ is given by $\tilde{U}_1=\left\{ \left(Y_0, Y_2,Y_3,\left[ l_2, l_3\right]| \right) \in U_1\times \mathbb{P}_\mathbb{C}^1| l_2Y_3= l_3 Y_2 \right\}$. $\tilde{U}_1$ has two affine coverings $U_{10}$ and $U_{11}$ where $U_{10}$ is given by $l_2\ne 0$, i.e. $U_{10}=\left\{ \left( Y_0, Y_2, v_3Y_2 , \left[ 1, v_3= \frac{l_3}{l_1}\right]   \right) \right\}$ and by setting $Y_0= v_0, Y_2= v_2$, we may identify $U_{10}\cong \left\{\left( v_0, v_2, v_3 \right)\in \mathbb{C}^3 \right\}$. Then $\pi: U_{10} \to U_1$ is given by $(v_1,v_2,v_3)\to (v_0, v_2, v_2v_3)$ and the exceptional divisor $E$ on $U_{10}$ is given by $v_2=0$, and $U_{11}$ is given by $l_3\ne 0$, i.e. $U_{11}=\left\{ \left( Y_0,  v_2' Y_3, Y_3, \left[ \frac{l_2}{l_3} , 1 \right] \right) \right\}$ and by setting $Y_0= v_0'$ and $Y_3= v_3'$, we may identify $U_{11}\cong \left\{ \left( v_0', v_2', v_3' \right) \in \mathbb{C}^3 \right\}$. Then $\pi: U_{11}\to U_1$ is given by $(v_0',v_2', v_3')\to (Y_0,Y_2,Y_3)= (v_0', v_2'v_3', v_3')$ and the exceptional divisor $E$ on $U_{11}$ is given by $v_3'=0$. The coordinate transformation between $U_{10}$ and $U_{11}$ is given by $(v_0,v_2,v_3)\mapsto \left(v_0',v_2',v_3' \right)=\left( v_0', \frac{1}{v_3'}, v_2'v_3' \right)$. We also note that the coordinate transformation between $U_{00}$ and $U_{10}$ is given by $\left( w_1,w_2,w_3 \right)\mapsto \left( v_0,v_2,v_3 \right)=\left( \frac{1}{w_1}, \frac{w_2}{w_1}, w_3  \right)$. In particular we note that the exceptional divisor $E$ is isomorphic to $\mathbb{P}_\mathbb{C}^1\times \mathbb{P}_\mathbb{C}^1$.

We will describe the foliation $\pi^* \Theta_{\mathcal{F}_0}$ on the blow up $Y$ of $\left( \mathbb{P}_\mathbb{C}^3, T_0 \right)$ along the invariant complex submanifold $\xi_2=\xi_3=0$. We describe $\pi^* \Theta_{\mathcal{F}_0}$ on $U_{00}$. We note that $\pi^* \Theta_{\mathcal{F}_0}\cong \pi^* \mathcal{O}_{\mathbb{P}_\mathbb{C}^3} \cong \mathcal{O}_Y$ is generated by $\pi^* T_0$. We note that $T_0$ on $U_0$ is given by
\begin{align*}
T_0 = \left( 1- X_1X_2 \right)\frac{\partial}{\partial X_1} + \left( X_3 - X_2^2 \right) \frac{\partial}{\partial X_2} - X_2 X_3 \frac{\partial}{\partial X_3}
\end{align*}
and from $\pi: U_{00}\to U_0$ by $(w_1,w_2, w_3)\mapsto  (X_1,X_2,X_3)=(w_1,w_2, w_2w_3)$, we see that under the natural map $\Pi: \Theta_Y \to \pi^* \Theta_{\mathbb{P}_\mathbb{C}^3}$, we have
\begin{align}
\Pi\left( \frac{\partial}{\partial w_1 } \right) & = \frac{\partial}{\partial X_1}\notag\\
\Pi\left( \frac{\partial}{\partial w_2} \right) & = \frac{\partial}{\partial X_2} + w_3 \frac{\partial}{\partial X_3} \label{blw1}\\
\Pi \left( \frac{\partial}{\partial w_3} \right) &= w_2 \frac{\partial}{\partial X_3}, \notag
\end{align}
so that $\pi^* T_0$ in $\Theta_Y$ on $U_{00}$ is given by
\begin{align} \label{blw5}
\pi^* T_0=\left( 1- w_1w_2 \right) \frac{\partial}{\partial w_1} +\left( w_2w_3 - w_2^2 \right)\frac{\partial}{\partial w_2} - w_3^2\frac{\partial}{\partial w_3}
\end{align}

We will show that $\left( Y, f^* T_0 \right)$ is unobstructed in foliated deformations in terms of tangent sheaves and describe a foliated analytic family such that the foliated Kodaira-Spencer map is an isomorphism at the distinguished point. We have the leaf complex on $\left( Y, f^* T_0 \right)$
\begin{align*}
\mathcal{E}_{\mathcal{O}_Y}^\bullet : \mathcal{E}_{\mathcal{O}_Y} \cong \Theta_Y \bigoplus \mathcal{O}_Y \to \mathscr{H}om_{\mathcal{O}_Y}\left( \mathcal{O}_Y, \Theta_Y \right)\cong \Theta_Y \to 0
\end{align*}
To compute $\mathbb{H}^i\left( Y, \mathcal{E}_{\mathcal{O}_Y}^\bullet \right)$, we need to describe  $H^0\left(Y, \Theta_Y \right)$ and $H^1\left( Y, \Theta_Y \right)$. We have the exact sequence 
\begin{align*}
0\to  \Theta_Y \to \pi^* \Theta_{\mathbb{P}_\mathbb{C}^3} \to \mathcal{N}_\pi \to 0
\end{align*}
and the exact sequence $($see \textnormal{\cite{Ser06} p.172} or \textnormal{\cite{Har10} p.94}$)$
\begin{align*}
0\to \mathcal{N}_{E/Y}\to \pi|_E^* \mathcal{N}_{\mathbb{P}_\mathbb{C}^1/\mathbb{P}_\mathbb{C}^3} \to \mathcal{N}_\pi \to 0,
\end{align*}
so that since $H^i\left( E, \mathcal{N}_{E/Y} \right)=0$ for $i\geq 0$, we have the exact sequence
\begin{align*}
0\to H^0\left(Y, \Theta_Y \right) \to H^0\left( \mathbb{P}_\mathbb{C}^3, \Theta_{\mathbb{P}_\mathbb{C}^3}  \right) \to H^0\left( \mathbb{P}_\mathbb{C}^1, \mathcal{N}_{\mathbb{P}_\mathbb{C}^1/ \mathbb{P}_\mathbb{C}^3} \right) \cong H^0\left( \mathbb{P}_\mathbb{C}^1, \mathcal{O}_{\mathbb{P}_\mathbb{C}^1}(1) \bigoplus \mathcal{O}_{\mathbb{P}_\mathbb{C}^1}(1) \right) \to H^1\left(Y, \Theta_Y \right) \to 0
\end{align*}
We recall that from $\textnormal{(\ref{ap2})}$ $H^0\left( \mathbb{P}_\mathbb{C}^3, \Theta_{\mathbb{P}_\mathbb{C}^3} \right)$ on $U_0$ is given by
\begin{align}
H^0\left( \mathbb{P}_\mathbb{C}^3,  \Theta_{\mathbb{P}_\mathbb{C}^3}  \right) = &\left\{  \left( a_1 + a_2 X_3+ a_3 X_2 + a_4 X_1 + a_5 X_1X_3 + a_6 X_1 X_2 + a_7 X_1^2 \right) \frac{\partial}{\partial X_1}\right. \\
 &\left.+\left( a_8 + a_9 X_3 + a_{10}X_2 + a_{11}X_1 + a_5 X_2X_3 + a_6 X_2^2 + a_7 X_1X_2\right) \frac{\partial}{\partial X_2} \right. \notag \\ 
 &\left.  +\left( a_{12}+ a_{13} X_3 + a_{14} X_2 + a_{15}X_1 + a_5 X_3^2 + a_6 X_2 X_3  + a_7 X_1X_3 \right)  \frac{\partial}{\partial X_3}, a_1,..., a_{15}\in \mathbb{C} \right\} \notag
\end{align}
Then the image of $H^0\left( \mathbb{P}_\mathbb{C}^3, \Theta_{\mathbb{P}_\mathbb{C}^3}  \right)\to H^0\left( \mathbb{P}_\mathbb{C}^1, \mathcal{O}_{\mathbb{P}_\mathbb{C}^1}(1)\bigoplus \mathcal{O}_{\mathbb{P}_\mathbb{C}^1}(1) \right)$ is given by
\begin{align*}
&\left( \left[ \left( a_8 + a_9 X_3 + a_{10}X_2 + a_{11}X_1 + a_5 X_2X_3 + a_6 X_2^2 + a_7 X_1X_2\right) \frac{\partial}{\partial X_2} , X_2 \right]|_{X_2=X_3=0} ,  \right.\\
&\left. \left[ \left( a_{12}+ a_{13} X_3 + a_{14} X_2 + a_{15}X_1 + a_5 X_3^2 + a_6 X_2 X_3  + a_7 X_1X_3 \right)  \frac{\partial}{\partial X_3}, X_3\right]_{X_2=X_3=0} \right)\\
&=\left( a_8 + a_{11}X_1, a_{12}+ a_{15}X_1 \right)
\end{align*}
which is surjective, so that $H^1\left(Y, \Theta_Y \right)=0$ and $\dim_\mathbb{C} H^0\left( Y, \Theta_Y \right)=11$. Then by using $(\ref{blw1})$ we see that $H^0\left( Y, \Theta_Y \right)$ on $U_{00}$ is given by \begin{align} 
 H^0\left(Y, \Theta_Y \right)=&\left\{ \left(b_1+b_2w_2w_3+b_3w_2+b_4w_1+b_5w_1w_2w_3+b_6w_1w_2+b_7w_1^2\right)\frac{\partial}{\partial w_1}\right. \label{blw2} \\
 &\left.+\left(b_8w_2w_3+b_9 w_2+b_5w_2^2w_3+b_6w_2^2+b_7w_1w_2 \right)\frac{\partial}{\partial w_2}+\left(b_{10}-b_8w_3^2+ b_{11}w_3\right)  \frac{\partial}{\partial w_3} | b_1,...,b_{11} \in \mathbb{C} \right\}\notag
 \end{align}
Then we have 
 \begin{align*}
\mathbb{H}^0\left(Y, \mathcal{E}_{\mathcal{O}_Y}^\bullet \right) &\cong \textnormal{ker}\left(H^0\left( Y,  \Theta_Y \right) \bigoplus H^0\left( Y, \mathcal{O}_Y \right) \to H^0\left(Y, \Theta_Y \right) \right) \\
\mathbb{H}^1\left(Y, \mathcal{E}_{\mathcal{O}_Y}^\bullet \right) &\cong \textnormal{coker}\left(H^0\left( Y,  \Theta_Y \right) \bigoplus H^0\left( Y, \mathcal{O}_Y \right) \to H^0\left(Y, \Theta_Y \right) \right) \\
\mathbb{H}^2\left( Y, \mathcal{E}_{\mathcal{O}_Y}^\bullet \right) & = 0
\end{align*}
Hence $\left( Y, \pi^* T_{\mathcal{F}_0}\right)$ is unobstructed in foliated deformations in terms of tangent sheaves. On the other hand, the image of $H^0\left( Y, \Theta_Y \right)\bigoplus H^0\left( Y, \mathcal{O}_Y \right) \to H^0\left( Y, \Theta_Y \right)$ is given by
 {\small{\begin{align}
 &\left[ \left(b_1+b_2w_2w_3+b_3w_2+b_4w_1+b_5w_1w_2w_3+b_6w_1w_2+b_7w_1^2\right)\frac{\partial}{\partial w_1} ,\left( 1- w_1w_2 \right) \frac{\partial}{\partial w_1} +\left( w_2w_3 - w_2^2 \right)\frac{\partial}{\partial w_2} - w_3^2\frac{\partial}{\partial w_3}\right] \notag \\
&+\left[ \left(b_8w_2w_3+b_9 w_2+b_5w_2^2w_3+b_6w_2^2+b_7w_1w_2 \right)\frac{\partial}{\partial w_2} + \left(b_{10}-b_8w_3^2+ b_{11}w_3\right)  \frac{\partial}{\partial w_3}, \left( 1- w_1w_2 \right) \frac{\partial}{\partial w_1} +\left( w_2w_3 - w_2^2 \right)\frac{\partial}{\partial w_2} - w_3^2\frac{\partial}{\partial w_3}\right] \notag \\ 
&+ b_{12}\left(\left( 1- w_1w_2 \right) \frac{\partial}{\partial w_1} +\left( w_2w_3 - w_2^2 \right)\frac{\partial}{\partial w_2} - w_3^2\frac{\partial}{\partial w_3} \right) \notag \\
&=\left( \left( b_{12}- b_4 \right) - 2 b_7 w_1 - \left(b_1+ b_6\right)w_2 - \left(b_9 + b_{12}\right)w_1w_2 - (b_3+ b_5 ) w_2w_3 - (b_6 + b_8)w_1w_2w_3\right)\frac{\partial}{\partial w_1} \label{blw3} \\
&+ \left( \left(b_{10}-b_7\right)w_2 - \left(b_9+ b_{12}\right)w_2^2 + \left( b_{11} + b_{12}\right)w_2w_3 -  \left(b_6 + b_8\right) w_2^2w_3  \right)\frac{\partial}{ \partial w_2}+ \left( - 2b_{10} w_3- \left(b_{11}+ b_{12}\right)w_3^2 \right) \frac{\partial}{\partial w_3} \notag 
\end{align}}}
Then we have 
\begin{align} \label{blwp1}
\dim_\mathbb{C} \mathbb{H}^0\left( Y, \mathcal{E}_{\mathcal{O}_Y}^\bullet \right)=4,\,\,\,\,\,\,\,\,\,\,\,\dim_\mathbb{C} \mathbb{H}^1 \left(  Y, \mathcal{E}_{\mathcal{O}_Y}^\bullet \right)= 11- (12-4)=3
\end{align}
 By comparing $(\ref{blw2})$ and $(\ref{blw3})$ we see that
{\small{\begin{align}
\mathbb{H}^1\left( Y, \mathcal{E}_{\mathcal{O}_Y}^\bullet \right) &\cong \textnormal{coker}\left(H^0\left( Y,  \Theta_Y \right) \bigoplus H^0\left( Y, \mathcal{O}_Y \right) \to H^0\left(Y, \Theta_Y \right) \right)  \notag \\
& \cong \textnormal{Span}_\mathbb{C} \left\langle w_1^2 \frac{\partial}{\partial w_1} + w_1w_2 \frac{\partial}{\partial w_2} , w_1\frac{\partial}{\partial w_1},  \frac{\partial}{\partial w_3} \right\rangle \cong \textnormal{Span}_\mathbb{C} \left\langle w_1^2 \frac{\partial}{\partial w_1} + w_1w_2 \frac{\partial}{\partial w_2} , w_2\frac{\partial}{\partial w_2},  \frac{\partial}{\partial w_3} \right\rangle \notag \\
&\cong \textnormal{Span}_\mathbb{C} \left\langle w_1^2 \frac{\partial}{\partial w_1} + w_1w_2 \frac{\partial}{\partial w_2} , w_3\frac{\partial}{\partial w_3},  \frac{\partial}{\partial w_3} \right\rangle \cong \textnormal{Span}_\mathbb{C} \left\langle w_1^2 \frac{\partial}{\partial w_1} + w_1w_2 \frac{\partial}{\partial w_2} , w_1\frac{\partial }{\partial w_1} + 2 w_2 \frac{\partial}{\partial w_2} - w_3 \frac{\partial}{\partial w_3} ,  \frac{\partial}{\partial w_3} \right\rangle \label{blw19}
\end{align}}}
Then we see that $\left( Y \times B, \tilde{T}, \tilde{\sigma}, B \right)$ where $\tilde{\sigma}: Y\times B \to B$ is the natural projection with the coordinates $(t_1,t_2,t_3)$ on a sufficiently small neighborhood $B\subset \mathbb{C}^3$ of $(0,0,0)$, and each of
{\small{\begin{align}
\tilde{T}&= \left( 1 - w_2w_3 - t_1 w_1^2  \right)\frac{\partial}{\partial w_1} +\left( w_2 w_3 - w_2^2 - t_1w_1w_2 - t_3 w_2 \right) \frac{\partial}{\partial w_2} +\left( - w_3^2 - t_2 \right)\frac{\partial}{\partial w_3} \label{blw8} \\
\tilde{T}&= \left( 1 - w_2w_3 - t_1 w_1^2 - t_3 w_1 \right)\frac{\partial}{\partial w_1} +\left( w_2 w_3 - w_2^2 - t_1w_1w_2 - 2t_3 w_2 \right) \frac{\partial}{\partial w_2} +\left( - w_3^2 - t_2 + t_3 w_3 \right)\frac{\partial}{\partial w_3} \label{blw9}
\end{align}}}
defines a foliated analytic family of deformations of $\tilde{\sigma}^{-1}(0)=\left( Y, \pi^* T_0 \right)$ such that the foliated Kodaira-Spencer map is an isomorphism at $(t_1,t_2,t_3)=(0,0,0)$.

 \subsection{Foliated deformations of the foliated exceptional divisor $\left(E, \pi^* T_0|_E \right)$} \

We note that the foliation $\pi^* \Theta_{\mathcal{F}_0}|_E\cong \mathcal{O}_E$ on $E$ is generated by $\pi^* T_0|_E$. We will show that $\left( E, \pi^* T_0 |_E \right)$ is unobstructed in foliated deformations in terms of tangent sheaves and describe a foliated analytic family of deformations of $\left( Y, \pi^* T_0|_E \right)$ such that the foliated Kodaira-Spencer map is an isomorphism at the distinguished point. We note that from $(\ref{blw5})$  $\pi^* T_0|_E$ on $E\cap U_{00}$ with the local coordinate $(w_1,w_3)$ is given by
\begin{align*}
\pi^* T|_E=\pi^* T|_{w_2=0} = \frac{\partial}{\partial w_1} - w_3^2\frac{\partial}{\partial w_3}
\end{align*}
We have the leaf complex 
\begin{align*}
\mathcal{E}_{\mathcal{O}_E}^\bullet : \mathcal{E}_{\mathcal{O}_E} \cong \Theta_E \bigoplus \mathcal{O}_E \to \mathscr{H}om_{\mathcal{O}_E}\left( \mathcal{O}_E, \Theta_E \right) \cong  \Theta_E \to 0
\end{align*}
To compute $\mathbb{H}^i\left( E, \mathcal{E}_{\mathcal{O}_E}^\bullet \right)$ we need to describe $H^0\left( E, \Theta_E \right)$ and $H^1\left(E, \Theta_E \right)$. Since $E\cong \mathbb{P}_\mathbb{C}^1\times \mathbb{P}_\mathbb{C}^1$, $H^0\left( E, \Theta_E \right)$ on $U_{00}\cap E$ is given by 
\begin{align} \label{blw6}
H^0\left( E, \Theta_E \right)=\left\{ \left( b_1 + b_2 w_1 + b_3 w_1^2 \right)\frac{\partial}{\partial w_1} + \left( b_4 + b_5 w_3 + b_6 w_3^2 \right)\frac{\partial}{\partial w_3} | b_1,...,b_6 \in \mathbb{C} \right\}
\end{align}
and we have $H^i\left( E, \Theta_E \right)=0$ for $i>1$. Then we have
\begin{align*}
\mathbb{H}^0\left( E, \mathcal{E}_{\mathcal{O}_E}^\bullet \right) & \cong \textnormal{ker}\left( H^0\left( E, \Theta_E \right) \bigoplus H^0\left( E, \mathcal{O}_E \right)\to H^0\left( E, \Theta_E \right) \right) \\
\mathbb{H}^1\left( E , \mathcal{E}_{\mathcal{O}_E}^\bullet  \right)& \cong \textnormal{coker}\left( H^0\left( E, \Theta_E \right) \bigoplus H^0\left( E, \mathcal{O}_E \right) \to H^0\left( E, \Theta_E \right) \right) \\
\mathbb{H}^2\left( E, \mathcal{E}_{\mathcal{O}_E}^\bullet \right) & = 0
\end{align*}
Hence $\left(  E, \pi^* T_0|_E \right)$ is unobstructed in foliated deformations. On the other hand, the image of $H^0\left( E , \Theta_E \right)\bigoplus H^0\left( E, \mathcal{O}_E \right)\to H^0\left( E, \Theta_E \right)$ is given by
\begin{align}
&\left[\left( b_1+ b_2 w_1 + b_3w_1^2\right) \frac{\partial}{\partial w_1} + \left( b_4 + b_5 w_3 + b_6 w_3^2 \right)\frac{\partial}{\partial w_3},  \frac{\partial}{\partial w_1'}- w_3^2 \frac{\partial}{\partial w_3} \right] + b_7\left( \frac{\partial}{\partial w_1} - w_3^2 \frac{\partial}{\partial w_3} \right)\notag\\
&= \left( b_7- b_2 - 2 b_3 w_1 \right) \frac{\partial}{\partial w_1} +\left( - 2b_4w_3- (b_5+b_7) w_3^2   \right) \frac{\partial}{\partial w_3} \label{blw7}
\end{align}
Then we have $\dim_\mathbb{C} \mathbb{H}^0\left(E, \mathcal{E}_{\mathcal{O}_E}^\bullet \right)= 3$ and $\dim_\mathbb{C} \mathbb{H}^1\left(  E, \mathcal{E}_{\mathcal{O}_E}^\bullet \right)= 6-(7-3)=2$. By comparing $(\ref{blw6})$ and $(\ref{blw7})$ we see that
\begin{align*}
\mathbb{H}^1\left( E, \mathcal{E}_{\mathcal{O}_E}^\bullet \right)&\cong \textnormal{coker}\left(H^0\left( E, \Theta_E \right) \bigoplus H^0\left( E, \mathcal{O}_E \right)\to H^0\left( E, \Theta_E \right) \right)\cong \textnormal{Span}_\mathbb{C} \left\langle w_1^2 \frac{\partial}{\partial w_1},  \frac{\partial}{\partial w_3} \right\rangle
\end{align*}
Then we see that $\left(E\times B, T_E, \sigma , B \right)$, where $\sigma: E\times B\to B$ is the natural projection with the coordinates $(t_1,t_2)$ on a sufficiently small neighborhood $B\subset \mathbb{C}^2$ of $(0,0)$ and
\begin{align}\label{blw10}
T_E= \left(1- t_1 w_1^2 \right)\frac{\partial}{\partial w_1} + \left( -t_2 - w_3^2 \right)\frac{\partial}{\partial w_3}
\end{align}
defines a foliated analytic family of deformations of $\sigma^{-1}(0)=\left( E, \pi^* T_0|_E \right)$ such that the foliated Kodaira-Spencer map is an isomorphism at $(t_1,t_2)=(0,0)$.

\subsection{Embedded deformations of $E$ in $\left( Y, \pi^* T_0 \right)$}\

Embedded deformations of $E$ in $\left(Y, \pi^* T_0 \right)$ are controlled by the complex of sheaves from $(\ref{cs8})$ 
\begin{align*} 
\mathcal{N}_{E/Y}^\bullet: \mathcal{N}_{E/Y}\to \mathscr{H}om_{\mathcal{O}_E}\left( \pi^* \Theta_{\mathcal{F}_0}|_E, \mathcal{N}_{E/Y} \right)\cong \mathcal{N}_{E/Y}\to 0
\end{align*}
We will denote the $i$-th cohomology group of $\mathcal{N}_{E/Y}^\bullet$ by $\mathbb{H}^i\left( E, \mathcal{N}_{E/Y}^\bullet \right)$. Then $\mathbb{H}^0\left( E, \mathcal{N}_{E/Y}^\bullet \right)$ represent first-order deformations and $\mathbb{H}^1\left( E, \mathcal{N}_{E/Y}^\bullet \right)$ represents an obstruction space (see Part VI).  Since $H^i\left( E, \mathcal{N}_{E/Y}^\bullet \right)=0$ for $i\geq 0$, we have $\mathbb{H}^0\left( E, \mathcal{N}_{E/Y}^\bullet \right)=\mathbb{H}^1\left( E, \mathcal{N}_{E/Y}^\bullet \right)=0$. In particular, $E$ is trivially foliated semi-regular (see Part VI) and $E$ is a stable invariant complex submanifold of $\left( Y, \pi^* T_0 \right)$  by theorem of stability for deformations of non-degenerate invariant complex submanifolds of foliated complex manifolds defined by locally free subsheaves of tangent sheaves (see Part VI).

\subsection{Simultaneous deformations of $E\subset \left( Y, \pi^* T_0 \right)$}\

Simultaneous deformations of $E\subset \left( Y, \pi^* T_0  \right)$ are controlled by the complex of sheaves (see Part VI)
\begin{align*}
\Theta_Y^\bullet: \Theta_Y \to \mathscr{H}om_{\mathcal{O}_Y}\left( \mathcal{O}_Y, \frac{\Theta_Y}{\mathcal{O}_Y} \right) \bigoplus i_* \mathcal{N}_{E/Y} \to i_* \mathcal{N}_{E/Y}\to 0
\end{align*}
where $i: E \hookrightarrow Y$ is the embedding. We will denote the $i$-th cohomology group of $\Theta_Y^\bullet$ by $\mathbb{H}^i\left( Y, \Theta_Y^\bullet \right)$.\footnote{We remark that we use different notations for the cohomology groups in Part VI. We use $\mathbb{H}^i\left( Y, \left( E, Y, \pi^*\Theta_{\mathcal{F}_0} \right)^\bullet \right)$ in Part VI.} Since $H^i\left( E, \mathcal{N}_{E/Y} \right)=0$ for $i\geq 0$, we have $\mathbb{H}^i\left( Y, \Theta_Y^\bullet \right)\cong \mathbb{H}^i\left( Y, \pi^* \Theta_{\mathcal{F}_0}^\bullet \right)\cong \mathbb{H}^i\left( Y, \mathcal{E}_{\mathcal{O}_Y}^\bullet  \right)$. In other words, simultaneous deformations of $E\subset \left( Y, \pi^* T_0 \right)$ is equivalent to foliated deformations of $\left( Y, \pi^* T_0 \right)$ in terms of tangent sheaves. We note that when we restrict the foliated analytic family of deformations of $\left(Y, \pi^* T_0 \right)$ given by $(\ref{blw8})$ to the exceptional divisor $E$, we obtain the foliated analytic family of deformations of $\left( E, \pi^* T_0|_E \right)$ given by $(\ref{blw10})$.

\subsection{Deformations of holomorphic foliated map $\pi$}\

Deformations of the non-degenerate holomorphic foliated map $\pi$ are controlled by the complex of sheaves (see Part VII)
\begin{align*}
\mathcal{N}_\pi^\bullet : \mathcal{N}_\pi \to \mathscr{H}om_{\mathcal{O}_Y} \left( \pi^* \Theta_{\mathcal{F}_0}, \mathcal{N}_\pi \right) \cong \mathcal{N}_\pi \to 0
\end{align*}
The complex of sheaves $\mathcal{N}_\pi^\bullet$ sits in the following exact sequence (see Part VII)
{\tiny{\begin{equation} \label{blw16}
\begin{CD}
@. 0 @. 0 @. 0 \\
@. @AAA @AAA @AAA \\
0@>>> \mathscr{H}om_{\mathcal{O}_Y}\left(\pi^* \Theta_{\mathcal{F}_0}, \frac{\Theta_Y}{\pi^* \Theta_{\mathcal{F}_0}} \right)\cong \frac{\Theta_Y}{\mathcal{O}_Y} @>>> \mathscr{H}om_{\mathcal{O}_Y}\left( \pi^* \Theta_{\mathcal{F}_0}, \pi^*\left(\frac{\Theta_{\mathbb{P}_\mathbb{C}^3}}{\Theta_{\mathcal{F}_0}} \right) \right) \cong \pi^* \left( \frac{\Theta_{\mathbb{P}_\mathbb{C}^3}}{\mathcal{O}_{\mathbb{P}_\mathbb{C}^3}} \right) @>>> \mathscr{H}om_{\mathcal{O}_Y}\left( \pi^* \Theta_{\mathcal{F}_0}, \mathcal{N}_\pi \right) \cong \mathcal{N}_\pi @>>> 0\\
@. @AAA @AAA @AAA \\
0 @>>> \Theta_Y @>\Pi>> \pi^* \Theta_{\mathbb{P}_\mathbb{C}^3} @>(-1)\cdot >>  \mathcal{N}_\pi @>>> 0
\end{CD}
\end{equation}}}
We have the exact sequence of complexes of sheaves (see Part VI)
{\tiny{\begin{center}
$\begin{CD}
@. 0 @. 0 @. 0 \\
@. @AAA @AAA @AAA \\
 0 @>>>  \mathscr{H}om_{\mathcal{O}_{\mathbb{P}_\mathbb{C}^1}} \left(\Theta_{\mathcal{F}_0}|_{\mathbb{P}_\mathbb{C}^1},  \frac{\Theta_{\mathbb{P}_\mathbb{C}^1}}{\Theta_{\mathcal{F}_0}|_{\mathbb{P}_\mathbb{C}^1}} \right)\cong \frac{\Theta_{\mathbb{P}_\mathbb{C}^1}}{\mathcal{O}_{\mathbb{P}_\mathbb{C}^1} } @>>>  \mathscr{H}om_{\mathcal{O}_{\mathbb{P}_\mathbb{C}^1}}\left( \Theta_{\mathcal{F}_0}|_{\mathbb{P}_\mathbb{C}^1} ,  \frac{\Theta_{\mathbb{P}_\mathbb{C}^3}|_{\mathbb{P}_\mathbb{C}^1}}{\Theta_{\mathcal{F}_0}|_{\mathbb{P}_\mathbb{C}^1}} \right) \cong \frac{\Theta_{\mathbb{P}_\mathbb{C}^3}|_{\mathbb{P}_\mathbb{C}^1}}{\mathcal{O}_{\mathbb{P}_\mathbb{C}^1}} @>>> \Theta_{\mathcal{F}_0}^*|_{\mathbb{P}_\mathbb{C}^1}  \otimes \mathcal{N}_{\mathbb{P}_\mathbb{C}^1/\mathbb{P}_\mathbb{C}^3} \cong \mathcal{N}_{\mathbb{P}_\mathbb{C}^1/\mathbb{P}_\mathbb{C}^3} @>>> 0 \\
@. @AAA @AAA @AAA \\
0@>>> \Theta_{\mathbb{P}_\mathbb{C}^1} @>>> \Theta_{\mathbb{P}_\mathbb{C}^3}|_{\mathbb{P}_\mathbb{C}^1} @>(-1)\cdot >>  \mathcal{N}_{\mathbb{P}_\mathbb{C}^1/\mathbb{P}_\mathbb{C}^3} @>>> 0
\end{CD}$
\end{center}}}
We have the exact sequence of complexes of sheaves (see Part VI)
{\Tiny{\begin{center}
$\begin{CD}
@. 0 @.0 @. 0 \\
@. @AAA @AAA @AAA \\
0@>>> \mathscr{H}om_{\mathcal{O}_E}\left( \pi^* \Theta_{\mathcal{F}_0}|_E, \frac{\Theta_E}{\pi^* \Theta_{\mathcal{F}_0}|_E} \right) \cong \frac{\Theta_E}{ \mathcal{O}_E} @>>> \mathscr{H}om_{\mathcal{O}_E}\left(\pi^* \Theta_{\mathcal{F}_0}|_E, \frac{\Theta_Y|_E}{\pi^* \Theta_{\mathcal{F}_0}|_E} \right)\cong \frac{\Theta_Y|_E}{\mathcal{O}_E} @>>> \mathscr{H}om_{\mathcal{O}_E}\left( \pi^* \Theta_{\mathcal{F}_0}|_E, \mathcal{N}_{E/Y} \right)\cong \mathcal{N}_{E/Y} @>>> 0  \\
@. @AAA @AAA @AAA \\
0 @>>> \Theta_E @>>> \Theta_Y|_E @>(-1)\cdot >> \mathcal{N}_{E/Y} @>>> 0
\end{CD}$
\end{center}}}
Then we have the exact sequence of complex of sheaves

\begin{equation}
\begin{CD} \label{blw17}
@. 0 @. 0 @. 0 \\
@. @AAA @AAA @AAA \\
0@>>> \mathcal{N}_{E/Y} @>>>  \pi |_E^* \mathcal{N}_{\mathbb{P}_\mathbb{C}^1/\mathbb{P}_\mathbb{C}^3} @>>> \mathcal{N}_\pi @>>> 0 \\
@. @AAA @AAA @AAA \\
0 @>>> \mathcal{N}_{E/Y} @>>>  \pi|_E ^* \mathcal{N}_{\mathbb{P}_\mathbb{C}^1/\mathbb{P}_\mathbb{C}^3} @>(-1)\cdot >> \mathcal{N}_\pi @>>> 0
\end{CD}
\end{equation}
The $i$-th cohomology group of complex of sheaves of $\pi|_E^* \mathcal{N}_{\mathbb{P}_\mathbb{C}^1/\mathbb{P}_\mathbb{C}^3} \to \pi|_E^* \mathcal{N}_{\mathbb{P}_\mathbb{C}^1/\mathbb{P}_\mathbb{C}^3}$ is isomorphic to $\mathbb{H}^i\left( \mathbb{P}_\mathbb{C}^1, \mathcal{N}_{\mathbb{P}_\mathbb{C}^1/\mathbb{P}_\mathbb{C}^3}^\bullet \right)$ from $(\ref{blw11})$, so that since $\mathbb{H}^i\left( E, \mathcal{N}_{E/Y}^\bullet \right)=0$ for $i\geq 0$, we have $ \mathbb{H}^i\left( \mathbb{P}_\mathbb{C}^1, \mathcal{N}_{\mathbb{P}_\mathbb{C}^1/\mathbb{P}_\mathbb{C}^3}^\bullet\right) \cong \mathbb{H}^i\left( Y, \mathcal{N}_\pi^\bullet \right)$.
In particular this implies that the functor $\textnormal{FDef}_{\mathbb{P}_\mathbb{C}^1\subset \left( \mathbb{P}_\mathbb{C}^3, T_0 \right)_{\textnormal{fixed}}}$ of deformations of $\mathbb{P}_\mathbb{C}^1$ in $\left( \mathbb{P}_\mathbb{C}^3, T_0 \right)$ (see Part VI) and the functor $\textnormal{FDef}_\pi$ of deformations of holomorphic foliated map $\pi$ (see Part VII) have the same first-order deformations and same obstruction space. We will show that first-order deformations of $\pi$ is induced by the foliated blow up along first-order deformations of $\xi_2=\xi_3=0$ in $\left( \mathbb{P}_\mathbb{C}^3, T_0 \right)$. We note that from $(\ref{blw12})$ we have a foliated analytic family of embedded deformations of $\xi_2=\xi_3=0$ in $\left( \mathbb{P}_\mathbb{C}^3, T_0 \right)$ over $\textnormal{Spec}\left( \mathbb{C}[t_1,t_2]/\left( t_1^2, t_2^2, t_1t_2 \right) \right)$ defined by $X_2= t_1+ t_2 X_1, X_3= t_2$. We blow up $\left( \mathbb{P}_\mathbb{C}^3, T_0  \right)$ along $X_2= t_1 + t_2 X_1, X_3= t_2$ on $U_0$. The blow up $\tilde{U}_0'$ over $U_0$ along $X_2- t_1 - t_2 X_1= X_3- t_2=0$ is given by $\left\{\left(X_1, X_2, X_3, \left[ l_2, l_3\right] \in U_0\times \mathbb{P}_\mathbb{C}^1 \right) | l_2\left(X_3- t_2 \right)= l_3\left(X_2 - t_1 - t_2 X_1 \right)\right\}$, and on $l_2\ne 0$, we have $U_{00}'=\left\{\left( X_1  , X_2, X_3 \left[ 1, w_3= \frac{l_3}{l_2} \right] | X_3- t_2= w_3\left( X_2- t_1- t_2 X_1 \right) \right) \right\}$ and by setting $w_2= X_2- t_1 - t_2 X_1$ and $w_1= X_1$, we may identify $U_{00}'\cong \left\{ (w_1,w_2,w_3)\in \mathbb{C}^3 \right\}$. Then $\tilde{\pi}: U_{00}'\to U_0$ is given by $(w_1,w_2,w_3)\mapsto (X_1,X_2,X_3)=\left( w_1, w_2+ t_1 + t_2 w_1, w_2w_3 + t_2 \right)$. Then under the natural map $\tilde{\Pi}: \Theta_{\tilde{Y}}\to \pi^* \Theta_{\mathbb{P}_\mathbb{C}^3}$ over $\textnormal{Spec}\left( \mathbb{C}[t_1,t_2]/\left(t_1^2, t_2^2, t_1t_2  \right) \right)$ we have
\begin{align*}
\tilde{\Pi}\left(\frac{\partial}{\partial w_1}\right)&= \frac{\partial}{\partial X_1} + t_2 \frac{\partial}{\partial X_2} \\
\tilde{\Pi}\left(\frac{\partial}{\partial w_2}\right)&= \frac{\partial}{\partial X_2} + w_3 \frac{\partial}{\partial X_3}\\
\tilde{\Pi}\left(\frac{\partial}{\partial w_3}\right)&= w_2 \frac{\partial}{\partial X_3}
\end{align*}
Then we note that since we are working over $\mathbb{C}[t_1,t_2]/\left(t_1^2, t_2^2, t_1t_2 \right)$, i.e. we only consider the linear terms,
\begin{align*}
&\left( 1- X_1X_2  \right) \frac{\partial}{\partial X_1} + \left(X_3 - X_2^2\right)\frac{\partial}{\partial X_2} - X_2  X_3 \frac{\partial}{\partial X_3}\\
\Longrightarrow &\left(1- w_1w_2  -t_1w_1- t_2w_1^2\right)\left(\frac{\partial}{\partial w_1}- t_2\frac{\partial}{\partial w_2} +t_2 w_3\frac{1}{w_2}\frac{\partial}{\partial w_3} \right) + \left(w_2w_3+ t_2 - w_2^2 - 2 t_1w_2- 2t_2w_1w_2\right)\left( \frac{\partial}{\partial w_2}- \frac{w_3}{w_2}\frac{\partial}{\partial w_3} \right)\\
&- \left(w_2^2w_3+ t_2 w_2 + t_1 w_2w_3+ t_2 w_1w_2w_3\right)\frac{1}{w_2}\frac{\partial}{\partial w_3}\\
&= \left(1- w_1w_2- t_1w_1- t_2 w_1^2\right) \frac{\partial}{\partial w_1} + \left( -t_2+ t_2 w_1w_2\right)\frac{\partial}{\partial w_2} +\left( t_2 \frac{w_3}{w_2} - t_2 w_1w_3 \right)\frac{\partial}{\partial w_3} \\
& +\left(w_2w_3+ t_2 - w_2^2 - 2 t_1w_2- 2t_2w_1w_2\right) \frac{\partial}{\partial w_2} +\left( -w_3^2- t_2 \frac{w_3}{w_2}+ w_2w_3+ 2t_1 w_3+ 2t_2 w_1 w_3 \right)\frac{\partial}{\partial w_3} \\
&+\left( - w_2w_3 - t_2 - t_1 w_3- t_2 w_1w_3 \right)\frac{\partial}{\partial w_3} \\
&= \left(1- w_1w_2 - t_1w_1 - t_2 w_1^2\right)\frac{\partial}{\partial w_1}  + \left(w_2w_3- w_2^2 - 2t_1 w_2 - t_2 w_1w_2\right)\frac{\partial}{\partial w_2}  + \left( - w_3^2 - t_2 + t_1w_3\right)\frac{\partial}{\partial w_3}
\end{align*}
This show that we have (over $\mathbb{C}[t_1, t_2]/ \left( t_1^2, t_2^2, t_1t_2 \right)$)
\begin{align*}
\tilde{\Pi}\left(\left(1- w_1w_2 - t_1w_1 - t_2 w_1^2\right)\frac{\partial}{\partial w_1}  + \left(w_2w_3- w_2^2 - 2t_1 w_2 - t_2 w_1w_2\right)\frac{\partial}{\partial w_2}  + \left( - w_3^2 - t_2 + t_1w_3\right)\frac{\partial}{\partial w_3} \right) = \tilde{\pi}^* T_0
\end{align*}
In other words, 
\begin{align} \label{blw18}
\tilde{\pi}:\left( \tilde{Y}, \left(1- w_1w_2 - t_1w_1 - t_2 w_1^2\right)\frac{\partial}{\partial w_1}  + \left(w_2w_3- w_2^2 - 2t_1 w_2 - t_2 w_1w_2\right)\frac{\partial}{\partial w_2}  + \left( - w_3^2 - t_2 + t_1w_3\right)\frac{\partial}{\partial w_3} \right) \to \left( \mathbb{P}_\mathbb{C}^3, T_0 \right)
\end{align}
 defines a family of holomorphic foliated maps of deformations of $\pi$ over $\mathbb{C}[t_1,t_2]/ \left(t_1^2, t_2^2, t_1t_2 \right)$. From $ \tilde{\pi} : \left( w_1 , w_2 , w_3 \right)\to \left(X_1,X_2, X_3 \right)=\left( w_1 , w_2 +t_1 + t_2w_1, w_2w_3+t_2 \right)$, we see that this family have the characteristic map given by (see Part VII)
\begin{align*}
T_0 \mathbb{C}^2 &\to \mathbb{H}^0\left( Y,  \mathcal{N}_\pi^\bullet \right) \\
 \frac{\partial}{\partial t_1} &\mapsto \textnormal{the class of}\,\,\, -\frac{\partial}{\partial X_2} \,\,\,\textnormal{in}\,\,\, \frac{\pi^* \Theta_{\mathbb{P}_\mathbb{C}^3}}{\Theta_Y} \\
 \frac{\partial}{\partial t_2} &\mapsto \textnormal{the class of}\,\,\, -w_1 \frac{\partial}{\partial X_2} - \frac{\partial}{\partial X_3} \,\,\,\textnormal{in}\,\,\, \frac{\pi^* \Theta_{\mathbb{P}_\mathbb{C}^3}}{\Theta_Y}
\end{align*}
We note that from $(\ref{blw12})$ the first order deformations of $\xi_2=\xi_3=0$ in $\left( \mathbb{P}_\mathbb{C}^3, T_0 \right)$ is generated by
\begin{align*}
\frac{\partial}{\partial X_2}\,\,\,\,\,\textnormal{and}\,\,\,\,\, X_1 \frac{\partial}{\partial X_2} + \frac{\partial}{\partial X_3}
\end{align*}
By pulling back under $\pi$, we obtain $\frac{\partial}{\partial X_1}$ and $w_1\frac{\partial}{\partial X_2} + \frac{\partial}{\partial X_3}$. Hence from the exact sequence $0\to \mathcal{N}_{E/Y}\to \pi|_E^* \mathcal{N}_{\mathbb{P}_\mathbb{C}^1/\mathbb{P}_\mathbb{C}^3}\to \mathcal{N}_\pi\to 0$, we see that first-oder deformations of $\xi_2=\xi_3=0 \cong \mathbb{P}_\mathbb{C}^1$ in $\left( \mathbb{P}_\mathbb{C}^3, T_0  \right)$ corresponds to first-order deformations of $\pi$.

Next we show that obstruction of liftings given by $(\ref{blw15})$ and $(\ref{pr18})$ corresponds to obstruction of liftings of the family $\tilde{\pi}$. In fact, obstruction of lifting of $\tilde{\pi}$ is given by the class in $\frac{\pi^* \Theta_{\mathbb{P}_\mathbb{C}^3}}{\Theta_Y}$ (see Part VII)
{\small{\begin{align*}
&\left(1- w_1w_2 - t_1w_1 - t_2 w_1^2\right)\left(\frac{\partial}{\partial X_1}+ t_2 \frac{\partial}{\partial X_2} \right) + \left(w_2w_3- w_2^2 - 2t_1 w_2 - t_2 w_1w_2\right)\left( \frac{\partial}{\partial X_2} + w_3 \frac{\partial}{\partial X_3} \right)  + \left( - w_3^2 - t_2 + t_1w_3\right)w_2\frac{\partial}{\partial X_3} \\
&- \left( 1- w_1\left(w_2+ t_1 + t_2 w_1  \right) \right) \frac{\partial}{\partial X_1} - \left(w_2w_3+t_2 - \left( w_2 + t_1 + t_2 w_1 \right)^2\right)\frac{\partial}{\partial X_2} + \left( w_2 + t_1 + t_2 w_1 \right)\left( w_2 w_3 + t_2 \right) \frac{\partial}{\partial X_3} \\
& = \left(t_1^2 + t_1t_2 w_1\right)\frac{\partial}{\partial X_2} + \left(t_1t_2 + t_2^2 w_1\right)\frac{\partial}{\partial X_3}
\end{align*}}}
Since $\left( - t_1^2- t_1t_2 X_1, - t_1t_2 - t_2^2 X_1 \right) \ne 0 \in \mathbb{H}^1\left( \mathbb{P}_\mathbb{C}^1, \mathcal{N}_{\mathbb{P}_\mathbb{C}^1/ \mathbb{P}_\mathbb{C}^3}^\bullet \right)$, we have the class of $ \left(t_1^2 + t_1t_2 w_1\right)\frac{\partial}{\partial X_2} + \left(t_1t_2 + t_2^2 w_1\right)\frac{\partial}{\partial X_3}\ne 0\in \mathbb{H}^1\left( Y, \mathcal{N}_\pi^\bullet \right)$. In other words, the holomorphic foliated map $\pi$ is obstructed. 

We also note that from $(\ref{blw16})$ and $(\ref{blw17})$ we have the long exact sequence
\begin{align*}
0\to \mathbb{H}^0\left( Y, \mathcal{E}_{\mathcal{O}_Y}^\bullet \right) \to \mathbb{H}^0\left( \mathbb{P}_{\mathbb{C}}^3, \mathcal{E}_{\mathcal{O}_{\mathbb{P}_\mathbb{C}^3}}^\bullet \right) \to  \mathbb{H}^0\left( \mathbb{P}_\mathbb{C}^1, \mathcal{N}_{\mathbb{P}_\mathbb{C}^1/\mathbb{P}_\mathbb{C}^3}^\bullet  \right) \to  \mathbb{H}^1 \left( Y, \mathcal{E}_{\mathcal{O}_Y}^\bullet \right) \to \mathbb{H}^1 \left( \mathbb{P}_{\mathbb{C}}^3, \mathcal{E}_{\mathcal{O}_{\mathbb{P}_\mathbb{C}^3}}^\bullet \right) \to  \mathbb{H}^1 \left( \mathbb{P}_\mathbb{C}^1, \mathcal{N}_{\mathbb{P}_\mathbb{C}^1/\mathbb{P}_\mathbb{C}^3}^\bullet  \right) \to 0
\end{align*}

We note that from $(\ref{blwp2})$ and $(\ref{blwp1})$, we have $\mathbb{H}^0\left( Y, \mathcal{E}_{\mathcal{O}_Y}^\bullet \right)\cong \mathbb{H}^0\left( \mathbb{P}_\mathbb{C}^3, \mathcal{E}_{\mathcal{O}_{\mathbb{P}_\mathbb{C}^3}}^\bullet \right)$, so that $\mathbb{H}^0\left(  \mathbb{P}_\mathbb{C}^1, \mathcal{N}_{\mathbb{P}_\mathbb{C}^1/\mathbb{P}_\mathbb{C}^3}^\bullet \right)\to \mathbb{H}^1\left( Y, \mathcal{E}_{\mathcal{O}_Y}^\bullet \right)$ is injective. This implies that from $(\ref{blw18})$ and $(\ref{blw19})$ the class
\begin{align*}
w_1\frac{\partial}{\partial w_1}+ 2 w_2 \frac{\partial}{\partial w_2} + w_3 \frac{\partial}{\partial w_3},\,\,\,\,\,\,\,\,\,\,\,\,\, w_1^2 \frac{\partial}{\partial w_1}+ w_1w_2 \frac{\partial}{\partial w_2}+ \frac{\partial}{\partial w_3} \in \mathbb{H}^1\left( Y, \mathcal{E}_{\mathcal{O}_Y}^\bullet \right)
\end{align*}
comes from the blow up along first-order deformations of $\xi_2=\xi_3=0$ in $\left( \mathbb{P}_\mathbb{C}^3, T_0 \right)$.

\section{Foliated cyclic coverings}\label{c1}

We define a notion of foliated cyclic coverings of foliated complex manifolds defined by locally free subsheaves of tangent sheaves as an analogue of cyclic covering of complex manifolds. Our presentation is based on \cite{Wav68} to which we refer for the detail. We will use the notion of a foliated line bundle on a foliated complex manifold (see Part V).

\begin{definition}[compare \cite{Wav68}]
Let $\left(M, \pi, W\right)$ be a $k$-sheeted branched covering of $W$, where $\pi: M \to W$. Let $C$ denote the branch locus. Moreover let $\left(M, \Theta_{\mathcal{F}_0}\right)$ be a foliated complex manifold with $\Theta_{\mathcal{G}_0}$  locally free of rank $p$, and $\Theta_{\mathcal{F}_0}$ be a holomorphic foliation on $W$ with locally free of rank $p$. We call $\left(M, \Theta_{\mathcal{G}_0}\right)$ be a foliated cyclic covering of $\left(W, \Theta_{\mathcal{F}_0}\right)$ if the following conditions are satisfied$:$
\begin{enumerate}
\item for each $x\in C$, $\pi^{-1}(x)$ consists of one point.
\item the group of covering transformations of $M- \pi^{-1}\left(C\right)$ over $W-C$ is cyclic of order $k$
\item for each $p\in C$, we can find a neighborhood $U$ with coordinates $\left(\xi_1,..., \xi_n\right)$ on $U$ and $\left(z_1,..., z_n \right)$ on $\pi^{-1}(U)$ such that the map is given by $\xi_i=z_i$ $(1\leq i\leq n-1)$, $\xi_n = \left(z_n\right)^k$.
\item if $k\ne 2$, $C$ is connected. 
\item $C$ is an invariant submanifold of $(W, \Theta_{\mathcal{F}_0})$
\item $\pi: \left(M, \Theta_{\mathcal{G}_0}\right)\to \left(W, \Theta_{\mathcal{F}_0}\right)$ is a holomorphic foliated map in the sense of \textnormal{Remark \ref{bcc1}} and $\Pi(\Theta_{\mathcal{G}_0})= \pi^* \Theta_{\mathcal{F}_0}$, where $\Pi:\Theta_M \to \pi^* \Theta_W $ is the natural map. In other words,  $\Theta_{\mathcal{G}_0}\cong \pi^* \Theta_{\mathcal{F}_0}$. 
\end{enumerate}
\end{definition}

\begin{remark}
The condition $(3)$ implies that the branch locus $C$ is non-singular and $(3)$ is actually a consequence of $(1)$ and $(2)$ $($see \cite{Wav68} \textnormal{Appendix p.957}$)$. We explain the condition $(6)$: we shall study deformations of a special kind of holomorphic maps between foliated complex manifolds in \textnormal{Part VII}. We recall \textnormal{Remark \ref{bcc1}} that $\pi:(M, \Theta_{\mathcal{G}_0})\to (W, \Theta_{\mathcal{F}_0})$ with $\Theta_{\mathcal{G}_0}$ and $\Theta_{\mathcal{F}_0}$ locally free is called a holomorphic foliated map if the exact sequence $0 \to \Theta_{\mathcal{F}_0} \to \Theta_W \to \frac{\Theta_W}{\Theta_{\mathcal{F}_0}}\to 0$ induces the exact sequence $0\to \pi^* \Theta_{\mathcal{F}_0} \to \pi^* \Theta_W \to \pi^*\left( \frac{\Theta_W}{\Theta_{\mathcal{F}_0}} \right) \to 0$, and $\Pi(\Theta_{\mathcal{G}_0})\subset \pi^* \Theta_{\mathcal{F}_0}$ where $\Pi: \Theta_M \to \pi^* \Theta_W$ is the natural map. When $\pi$ is a foliated $k$-sheeted cyclic covering, we have the exact sequence $ 0 \to \Theta_M \xrightarrow{\Pi} \pi^* \Theta_W \to \mathcal{N}_\pi\to 0$, i.e. $\Pi$ is injective. We claim that the condition $(6)$ implies $(5)$, i.e. $C$ is a nonsingular invariant hypersurface of $(W, \Theta_{\mathcal{F}_0})$. In fact,  let $\Theta_{\mathcal{F}_0}$ on a neighborhood $U$ of $W$ with the coordinate $\left(\xi_1,..., \xi_n \right)$ be generated by $T^1,..., T^p\in \Gamma(U, \Theta_W)$ where
\begin{align*}
T^\alpha= \sum_{\beta=1}^{n-1} T^{\alpha\beta} (\xi) \frac{\partial}{\partial \xi_\beta} + T_n^\alpha (\xi)\frac{\partial}{\partial \xi_n},\,\,\,\,\,\,\,\alpha=1,..., p
\end{align*}
Since $\pi^* \Theta_{\mathcal{F}_0}=\Pi(\Theta_{\mathcal{G}_0})$, we see that
\begin{align*}
\pi^* T^\alpha = \sum_{\beta=1}^{n-1} T^{\alpha \beta}\left(z_1,...,z_{n-1}, ( z_n)^k\right) \frac{\partial}{\partial \xi_\beta} + T_n^\alpha \left(z_1,..., z_{n-1}, z_n^k\right)\frac{\partial}{\partial \xi_n}
\end{align*} 
generates $\Pi(\Theta_{\mathcal{G}_0})$ on $\pi^{-1}(U)$, so that $\Theta_{\mathcal{G}_0}$ is generated by $Q^\alpha$
\begin{align*}
Q^\alpha = \sum_{\beta=1}^{n-1} T^{\alpha \beta}\left(z_1,..., z_{n-1}, (z_n)^k\right)\frac{\partial}{\partial z_\beta} + Q_n^\alpha(z_1,..., z_n)\frac{\partial}{\partial z_n}
\end{align*}
such that
\begin{align*}
Q_n^\alpha\left(z_1,..., z_n\right)k \left(z_n\right)^{k-1} = T_n^\alpha\left(z_1,..., z_{n-1}, \left(z_n\right)^k\right)
\end{align*}
This implies that $T_n^\alpha\left( z_1,..., z_{n-1}, (z_n)^k\right)$ has a factor $\left(z_n\right)^{k-1}$. When we expand $T_n^\alpha\left(\xi_1,..., \xi_n\right)$ in $\xi_n$ whose coefficients are holomorphic functions of $\xi_1,...,\xi_{n-1}$,  the term involving only $\xi_1,...,\xi_{n-1}$ is $0$, so that $T_n^\alpha$ is of the form $T_n^\alpha \left(\xi_1,..., \xi_n \right)= \xi_n G_n^\alpha \left(\xi_1,..., \xi_n \right)$. This implies that $\xi_n=0$ on $U$ is invariant with respect to $\Theta_{\mathcal{F}_0}$

\end{remark}

Before proving the following theorem \ref{bc1}, we define a notion of a foliated line bundle as an analogue of a notion of a Poisson line bundle (for the detail, we refer to Part V). We shall study deformations of foliated line bundles (and more generally deformations of foliated vector bundles) in Part V as an analogue of the author's previous work (\cite{Kim16}) on deformations of Poisson line bundles in the context of foliated deformations. Given a foliated line bundle, we will make the total space of the line bundle a foliated complex manifold, and we will also make the total space of the associated projective line bundle a foliated complex manifold.

\begin{definition}
Let $\left(W, \Theta_{\mathcal{F}_0}\right)$ be a foliated complex manifold with $\Theta_{\mathcal{F}_0}$ of locally free. A foliated connection on a line bundle $F$ on $\left(W, \Theta_{\mathcal{F}_0}\right)$ is a $\mathbb{C}$-linear bracket $\{-,-\}_F:\Theta_{\mathcal{F}_0} \times \mathcal{O}(F) \to \mathcal{O}(F)$ where $\mathcal{O}(F)$ is the sheaf on $W$ of germs of holomorphic sections of $F$, such that $\left\{fT, a \right\}_F= f \left\{T,a \right\}_F$ and $\left\{T, fa \right\}_F=\left[T,f \right]a+f\left\{T, a\right\}_F$ where $T\in \Theta_{\mathcal{F}_0}, f\in \mathcal{O}_W,a\in \mathcal{O}(F)$. Equivalently, a foliated connection is given by a homomorphism $\nabla_F^0 :\mathcal{O}(F)\to \mathscr{H}om_{\mathcal{O}_W}\left(\Theta_{\mathcal{F}_0}, \mathcal{O}(F) \right)$ which satisfies the identity $\nabla_F^0(fa)(T)= \left[T, f \right]a + f\nabla_F^0(a)(T)$ and the bracket is defined by the formula $\nabla_F^0(a)(T)= \left\{T,a \right\}_F$. A foliated connection $\nabla_F^0$ is called flat if it satisfies the Jacobi identity $\left\{T_1,\left\{T_2,a\right\}_F \right\}_F = \left\{ \left[T_1,T_2\right], a \right\}_F+ \left\{T_2, \left\{T_1,a \right\}_F \right\}_F$ for $T_1,T_2\in \Theta_\mathcal{F}, a\in \mathcal{O}(F)$.
A line bundle $F$ on $\left(W, \Theta_{\mathcal{F}_0}\right)$ equipped with a flat foliated connection $\nabla_F$ is called a foliated line bundle. We denote the foliated line bundle $\left(F, \nabla_F \right)$.
\end{definition}

Now let $\left( F, \nabla_F \right)$ be a foliated line bundle on a complex foliated manifold $\left(W, \Theta_{\mathcal{F}_0}\right)$ with $\Theta_{\mathcal{F}_0}$ locally free of rank $p$. Let $\left\{V_i \right\}$ be a covering of $W$ by polycylinders with coordinates $\left(\xi_i^1,..., \xi_i^n\right)$ on $V_i$ such that $F$ is trivialized and $\Gamma\left(V_i, \Theta_{\mathcal{F}_0}\right)$ is generated by $T_i^1,..., T_i^p\in \Gamma\left( V_i, \Theta_W \right)$ with
\begin{align*}
T_i^\alpha&=\sum_{\eta=1}^n T_i^{\alpha\eta} (\xi_i)\frac{\partial}{\partial \xi_i^\eta} \\
T_i^\alpha &= \sum_{\beta=1}^p r_{ij}^{\alpha\beta}(\xi_j) T_j^\beta \,\,\,\,\,\,\,\,\,\,\text{on}\,\,\,V_{ij}\\
\left[ T_i^\alpha, T_i^\beta \right] & =\sum_{\gamma=1}^p c_{i\alpha\beta}^\gamma(\xi_i) T_i^\gamma \,\,\,\,\,\,\,\,\,\text{on}\,\,\, V_i
\end{align*}

 Then $\left(F, \nabla_F\right)$ is represented by $\left( \left\{\phi_i \right\}, \left\{f_{ij} \right\} \right)$ where $\phi_i:\Gamma \left(V_i, \Theta_{\mathcal{F}_0}\right)\to \Gamma\left(V_i, \mathcal{O}_W \right)$ and $f_{ij}\in \Gamma\left(V_{ij}, \mathcal{O}_W^* \right)$ such that
\begin{align*}
\phi_j \left(T_i^\alpha \right)&- \phi_i\left(T_i^\alpha\right)+ \left[T_i^\alpha, \log f_{ij} \right]=0 \\
\left[T_i^\alpha, \phi_i\left(T_i^\beta\right) \right]&-\left[T_i^\beta, \phi_i\left(T_i^\alpha\right)\right]-\phi_i\left(\left[T_i^\alpha, T_i^\beta\right] \right) = 0 \\
&f_{ij} f_{jk}= f_{ik}
\end{align*}
Then $\left(F^k, \nabla_{F^k}\right)$ is naturally a foliated line bundle on $\left(W, \Theta_{\mathcal{F}_0}\right)$ represented by $\left( \left\{k\phi_i \right\}, \left\{f_{ij}^k\right\} \right)$ (for the detail, see Part V).

Next we will show that a foliated line bundle $\left(F, \nabla_F \right)$ on $\left(W, \Theta_{\mathcal{F}_0}\right)$ with $\Theta_{\mathcal{F}_0}$ locally free of rank $p$ as above defines a foliation $\Theta_{\tilde{\mathcal{F}}_0}$ on the total space $F$, so that $\left(F,  \Theta_{\tilde{\mathcal{F}}_0} \right)$ is a foliated complex manifold with $\Theta_{\tilde{\mathcal{F}}_0}$ locally free of rank $p$. Let the total space $F$ be described by
\begin{align*}
V_{ij}\times \mathbb{C} \subset V_j\times \mathbb{C} &\to V_{ij}\times \mathbb{C} \subset V_i\times \mathbb{C}\\
                             \left(\xi_j, t_j \right)= \left(\xi_j^1,...,\xi_j^n, t_j \right) &\mapsto \left(\xi_i, f_{ij}(\xi_j)t_j \right)= \left(\xi_i^1,...,\xi_i^n , t_i \right)
\end{align*}
We define $\tilde{T}_j^\beta \in \Gamma\left(V_j\times \mathbb{C}, \Theta_F \right)$ by
\begin{align*}
\tilde{T}_j^\beta= T_j^\beta+ t_j\phi_j \left(T_j^\beta \right) \frac{\partial}{\partial t_j} = \sum_{\eta=1}^p T_j^{\beta \eta}(\xi_j)\frac{\partial}{\partial \xi_j^\eta} + t_j\phi_j \left(T_j^\beta \right)\frac{\partial}{\partial t_j}
\end{align*}
From $t_i=f_{ij}(\xi_j) t_j$, we have $\frac{\partial}{\partial t_j}=f_{ij}(\xi_j)\frac{\partial}{\partial t_i}$ and $\frac{\partial}{\partial \xi_j^\alpha}=\sum_{\beta=1}^n \frac{\partial \xi_i^{\beta}}{\partial \xi_j^\alpha}\frac{\partial}{\partial \xi_i^\beta} + \frac{\partial f_{ij}(\xi_j)}{\partial \xi_j^\alpha}t_j\frac{\partial}{\partial t_i} = \sum_{\beta=1}^n \frac{\partial \xi_i^\beta}{\partial \xi_j^\alpha}\frac{\partial}{\partial \xi_i^\beta} + \frac{\partial \log f_{ij}(\xi_j)}{\partial \xi_j^\alpha} t_i\frac{\partial}{\partial t_i}$. Then we have
\begin{align*}
\sum_{\beta=1}^p r_{ij}^{\alpha\beta}(\xi_j)\tilde{T}_j^\beta & = \sum_{\beta=1}^p r_{ij}^{\alpha\beta}(\xi_j) \left( \sum_{\eta=1}^n T_j^{\beta \eta}(\xi_j)\frac{\partial}{\partial \xi_j^\eta} + t_j\phi_j\left(T_j^\beta\right)\frac{\partial}{\partial t_j}  \right)\\
&=\sum_{\beta=1}^p r_{ij}^{\alpha\beta}(\xi_j)\left(\sum_{\eta, \gamma=1}^n T_j^{\beta\eta}\left(\xi_j \right)\frac{\partial \xi_i^\gamma}{\partial \xi_j^\eta}\frac{\partial}{\partial \xi_i^\gamma}  +\sum_{\eta=1}^n T_j^{\beta\eta}(\xi_j)\frac{\partial \log f_{ij}}{\partial \xi_j^\eta} t_i\frac{\partial}{\partial t_i} + t_i\phi_j\left(T_j^\beta \right)\frac{\partial}{\partial t_i} \right)\\
&= T_i^\alpha + t_i \left[T_i^\alpha, \log f_{ij} \right] \frac{\partial}{\partial t_i} + t_i\phi_j(T_i^\alpha)\frac{\partial}{\partial t_i} = T_i^\alpha + t_i\phi_i(T_i^\alpha)\frac{\partial}{\partial t_i} = \tilde{T}_i^\alpha
\end{align*}
and we have
\begin{align*}
\left[\tilde{T}_i^\alpha , \tilde{T}_i^\beta\right]&=\left[ T_i^\alpha + t_i \phi_i\left(T_i^\alpha\right)\frac{\partial}{\partial t_i} , T_i^\beta +   t_i\phi_i\left(T_i^\beta\right)\frac{\partial}{\partial t_i}  \right] = \left[T_i^\alpha, T_i^\beta \right]+ t_i\left[T_i^\alpha, \phi_i\left(T_i^\beta\right)\right] \frac{\partial}{\partial t_i}-  t_i\left[ T_i^\beta, \phi_i\left(T_i^\alpha\right) \right]\frac{\partial}{\partial t_i}\\
&= \sum_{\gamma=1}^p c_{i\alpha\beta}^\gamma (\xi_i) T_i^\gamma + t_i \phi_i\left( \left[T_i^\alpha, T_i^\beta \right] \right) \frac{\partial}{\partial t_i}= \sum_{\gamma=1}^p c_{i\alpha\beta}^\gamma (\xi_i) \left( T_i^\gamma + t\phi_i\left(T_i^\gamma\right)\frac{\partial}{\partial t_i} \right) = \sum_{\gamma=1}^p c_{i\alpha\beta}^\gamma (\xi_i) \tilde{T}_i^\gamma
\end{align*}

Then the subsheaf $\Theta_{\tilde{\mathcal{F}}_0}$ of $\Theta_F$ generated by $\tilde{T}_i^1,..., \tilde{T}_i^p$ on $V_i\times \mathbb{C}$ is locally free of rank $p$ and defines a foliation on $F$, and the transition matrices of $\Theta_{\tilde{\mathcal{F}_0}}$ with respect to the covering $\{V_i\times \mathbb{C}\}$ is given by $\left\{\left(R_{ij}^T\right)^{-1}\right\}$ where $R_{ij}=\left( r_{ij}^{\alpha\beta}(\xi_j) \right)$.  Let $\pi:F\to W$ be the bundle map. Then we have $\pi^* \Theta_{\mathcal{F}_0} \cong \Theta_{\tilde{\mathcal{F}}_0}$ and $\pi: \left(F, \Theta_{\tilde{\mathcal{F}}_0} \right)\to \left(W, \Theta_{\mathcal{F}_0} \right)$ is a holomorphic foliated map. In fact, the natural homomrphism $\pi^* \Theta_{\mathcal{F}_0} \to \pi^*\Theta_W$ is injective, and $\Pi\left(\Theta_{\tilde{\mathcal{F}}_0}\right)\subset \pi^* \Theta_{\mathcal{F}_0} $ where $\Pi: \Theta_F \to \pi^* \Theta_W$ is the natural map.

Let $F^\sharp$ be the $\mathbb{P}_\mathbb{C}^1$-bundle over $W$ obtained from $F$. Let $F^\sharp$ be described by
\begin{center}
$\begin{CD}
\left(\xi_j,t_j\right)\in V_{ij}\times \mathbb{C} \subset V_j\times \mathbb{C}\subset V_j\times \mathbb{P}_\mathbb{C}^1 @>>> \left(\xi_i,t_i\right)=  \left(\xi_i,f_{ij}(\xi_j)t_j\right) \in V_{ij} \times \mathbb{C} \subset V_i\times \mathbb{C} \subset V_i \times \mathbb{P}_\mathbb{C}^1\\
@VVV @VVV\\
\left(\xi_j, t_j' \right)=\left(\xi_j, \frac{1}{t_j}\right)\in V_{ij}\times \mathbb{C} \subset V_j\times \mathbb{C} \subset V_j\times \mathbb{P}_\mathbb{C}^1 @>>> \left(\xi_i,t_i' \right)= \left(\xi_i, \frac{1}{f_{ij}(\xi_j)t_j} \right) \in V_{ij}\times \mathbb{C} \subset V_i\times \mathbb{C} \subset V_i\times \mathbb{P}_\mathbb{C}^1
\end{CD}$
\end{center}
We note that since $t_i= \frac{1}{t_i'}$, we have $t_i\frac{\partial}{\partial t_i}= - t_i' \frac{\partial}{\partial t_i'}$. Then we define on $(\xi_i,t')$ 
\begin{align*}
\tilde{T}_i'^\alpha= T_i^\alpha - t_i'\phi_i\left(T_i^\alpha \right)\frac{\partial}{\partial t_i'}
\end{align*}
Then $\tilde{T}_i^\alpha= \tilde{T}_i'^\alpha$ on $V_i\times \mathbb{C}^*$, so that $\tilde{T}_i^\alpha$ is naturally extended to a vector field on $V_i\times \mathbb{P}_\mathbb{C}^1$ and $\tilde{T}_i'^\alpha= \sum_{\beta=1}^p r_{ij}^{\alpha\beta}(\xi_j) \tilde{T}_j^\beta$ and $\tilde{T}_i^\alpha= \sum_{\beta=1}^p r_{ij}^{\alpha\beta}(\xi_j) \tilde{T}_j'^\beta$ and $\tilde{T}_i'^\alpha= \sum_{\beta=1}^q r_{ij}^{\alpha\beta}(\xi_j)\tilde{T}_j'^\beta$, and moreover we have
\begin{align*}
\left[\tilde{T}_i'^\alpha, \tilde{T}_i'^\beta\right]&= \left[ T_i^\alpha - t_i'\phi_i\left(T_i^\alpha\right)\frac{\partial}{\partial t_i'} , T_i^\beta - t_i'\phi_i \left(T_i^\beta\right)\frac{\partial}{\partial t_i'} \right]\\
                                                      & =\left[T_i^\alpha, T_i^\beta \right]- t_i' \left[T_i^\alpha, \phi_i \left(T_i^\beta\right)\right]\frac{\partial}{\partial t_i'} + t_i' \left[ T_i^\beta, \phi_i(T_i^\alpha) \right] \frac{\partial}{\partial t_i'} = \sum_{\gamma=1}^p c_{i\alpha\beta}^\gamma (\xi_i) T_i^\gamma - t_i'\phi_i \left( \left[T_i^\alpha, T_i^\beta\right] \right) \frac{\partial}{\partial t_i'}\\
                                                      &=\sum_{\gamma=1}^p c_{i\alpha\beta}^\gamma \left( T_i^\gamma - t_i'\phi_i\left(T_i^\gamma\right) \frac{\partial}{\partial t_i'} \right) = \sum_{\gamma=1}^p c_{i\alpha\beta}^\gamma \tilde{T}_i^\gamma
\end{align*}

Then the subsheaf $\Theta_{\tilde{\mathcal{F}}^\sharp_0}$ of $\Theta_{F^\sharp}$ generated by $\tilde{T}_i^\alpha$ on $V_i\times \mathbb{C}\subset V_i\times \mathbb{P}_\mathbb{C}^1 $ and $\tilde{T}_i'^\alpha$ on $V_i\times \mathbb{C}\subset V_i \times \mathbb{P}_\mathbb{C}^1$ is locally free of rank $p$ and defines a foliation on $F^\sharp$, and the transition matrices of $\Theta_{\tilde{\mathcal{F}}^\sharp_0}$ with respect to the covering $\{U_\alpha\}=\{U_{2i}, U_{2i+1}\}$ of $F^\sharp$ where $U_{2i}\cong U_{2i+1}\cong V_i\times \mathbb{C}$ and if coordinates in $U_\alpha$ are $(\xi_\alpha, u_\alpha)$, then $u_{2i}= f_{ij} u_{2j}$ in $U_{2i}\cap U_{2j}$ and $u_{2i}=\frac{1}{u_{2i+1}}$ is given by

\begin{equation} \label{cs2}
\phi_{\alpha\beta}^\sharp=\begin{cases}
\left(R_{ij}^T \left(\xi_j\right)\right)^{-1} \,\,\,\,\,\,\,\,\,\,\,\,\textnormal{if}\,\,\,\alpha=2i,\beta=2j \\
\textnormal{I}_p \,\,\,\,\,\,\,\,\,\,\,\,\,\,\,\,\, \,\,\,\,\,\,\,\,\,\,\,\,\,\,\,\,\,\,\,\,\,\,\textnormal{if}\,\,\, \alpha=2i, \beta = 2i+1
\end{cases}
\end{equation}
Let $\pi^\sharp:F^\sharp \to W$ be the bundle map. Then $\left(\pi^\sharp \right)^* \Theta_{\mathcal{F}_0}\cong \Theta_{\tilde{\mathcal{F}}_0^\sharp}$ and $\pi^\sharp: \left( F^\sharp, \Theta_{\tilde{\mathcal{F}}_0^\sharp} \right) \to (W, \Theta_{\mathcal{F}_0})$ is a holomorphic foliated map.

Let $s\in H^0\left(W, F^k\right)$ represented by $s_i=  f_{ij}^ks_j$ such that $\left[T_i^\alpha, s_i \right]= k\phi_i\left( T_i^\alpha\right)s_i$ for each $i$ and $\alpha=1,...,p$. Then we have
\begin{align*}
\left[ \tilde{T}_i^\alpha, t_i^k-s_i \right]= \left[ T_i^\alpha + t_i\phi_i\left( T_i^\alpha \right) \frac{\partial}{\partial t_i}, t_i^k - s_i   \right] = -\left[ T_i^\alpha, s_i \right]+ k t_i^k \phi_i\left( T_i^\alpha \right) = k \phi_i\left( T_i^\alpha \right)\left( t_i^k - s_i \right)
\end{align*}
This implies that $\left\{t_i^k -s_i=0 \right\}$ is invariant with respect to $\Theta_{\tilde{\mathcal{F}}_0}$. Moreover, we have
\begin{align*}
\left[ \tilde{T}_i'^\alpha, t_i'^k s_i-1  \right] = \left[ T_i^\alpha - t_i'\phi_i\left( T_i^\alpha \right) \frac{\partial}{\partial t_i'}, t_i'^k s_i -1  \right] = t_i'^k \left[ T_i^\alpha, s_i \right] - k t_i'^k \phi_i \left( T_i^\alpha \right) s_i = 0
\end{align*}
This implies that $\left\{ t_i^k - s_i=0, t_i'^k s_i-1 =0 \right\} $ is invariant with respect to $\Theta_{\tilde{\mathcal{F}}_0^\sharp}$.

\begin{theorem}[compare \cite{Wav68} Theorem1.2 p.929-931]  \label{bc1}
If $\left(M, \Theta_{\mathcal{G}_0}\right)$ is a foliated cyclic $k$-covering of a foliated complex manifold $\left(W, \Theta_{\mathcal{F}_0}\right)$ with $\Theta_{\mathcal{G}_0}$ and $\Theta_{\mathcal{F}_0}$ locally free, we can find a foliated line bundle $(F, \nabla_F)$ on $\left(W, \Theta_{\mathcal{F}_0}\right)$ so that, over a suitable fine covering $\left\{V_i \right\}$ of $W$, $M$ may be identified with the non-degenerate invariant complex submanifold of $\left(F, \Theta_{\tilde{\mathcal{F}}_0}\right)$ defined by the equation $\xi_i^k=\phi_i$, where $\xi_i$ is the fibre coordinate of $F$ over $V_i$ and $\phi_i$ is the local equation of $C$ in $V_i$.
\end{theorem}

\begin{proof}
We extend the argument of \cite{Wav68} Theorem1.2 p.929-931 to which we refer for the detail. We tried to maintain the notational consistency with \cite{Wav68}. We choose a covering $\left\{U_j, U_\lambda\right\}$ of $W$ by polycylinder having the following properties $(\pi: M \to W$ is the projection):
\begin{enumerate}
\item $U_j\cap C = \emptyset, U_\lambda \cap C \ne \emptyset$;
\item on $U_\lambda$, $C$ has the equation $\xi_\lambda^n =0$ where $\left(\xi_\lambda^1,..., \xi_\lambda^n\right)$ are local coordinates on $U_\lambda$;
\item $\pi^{-1} \left(U_j \right)= U_j^{(0)} \bigcup \cdots \bigcup U_j^{(k-1)}$ with $U_j^{(v)}\cap U_j^{(\mu)} = \emptyset$ for $v\ne \mu$;
\item $\pi^{-1}\left(U_\lambda\right)= U_\lambda^\sharp$ and coordinates in $U_\lambda^\sharp$ are $\left(z_\lambda^1,..., z_\lambda^n\right)$ and the map $\pi$ is given by
\begin{equation*}
\begin{cases}
\xi_\lambda^v= z_\lambda^v, v=1,..., n-1 \\
\xi_\lambda^n= \left(z_\lambda^n \right)^k
\end{cases}
\end{equation*}
\end{enumerate}
Let the coordinate changes on $U_j^{(v)}\cap U_\lambda^{\sharp}$ be 
\begin{align*}
z_\lambda^\alpha = g_{\lambda j}^{(v)\alpha}\left(z_j^{(v)} \right)= g_{\lambda j}^{(v)\alpha} \left( z_j^{(v)1},.., z_j^{(v)n} \right)
\end{align*}
where $z_j^{(v)}= \left(z_j^{(v)1},..., z_j^{(v)n}\right)$ are coordinates on $U_j^{(v)}$. Then we may assume that $U_j^{(v)}$ are numbered so that $g_{\lambda j}^{(v) n}= \rho^v g_{\lambda j}^{(0) n}$ with $\rho= \exp\left( \frac{2\pi i}{k} \right)$, and if $U_j^{(0)}\cap U_k^{(v)}\ne \emptyset$, then $U_j^{(s)} \cap U_k^{(v+s)} \ne \emptyset $ where by $v+s$ we understand that residue class of $v+s \mod k$. Then we can define a line bundle $F$ on $W$ by the following transition functions:
\begin{equation*}
\begin{cases}
f_{ij}= \rho^{-v} \,\,\,\,\, \textnormal{if} \,\,\,\,\, U_i^{(0)}\cap U_j^{(v)}\ne \emptyset \\
f_{\lambda j}= g_{\lambda j}^{(0) n}\\
f_{\lambda \mu}= \phi_{\lambda \mu}\left( \xi_\mu^1,..., \xi_\mu^{n-1}, \sqrt[k]{\xi_\mu^n}  \right)
\end{cases}
\end{equation*}
where $\phi_{\lambda \mu}(z_\mu)=\phi_{\lambda\mu}\left(z_\mu^1,..., z_\mu^n\right)$ is a non-vanishing holomorphic function on $U_\lambda^\sharp \cap U_\mu^\sharp\ne \emptyset$ such that $z_\lambda^n= \phi_{\lambda \mu}(z_\mu)\cdot z_\mu^n $ on $U_\lambda^\sharp \cap U_\mu^\sharp \ne \emptyset$ (for the detail, see \cite{Wav68} Lemma p.930).

We will give a foliated connection $\nabla_F$ on $F$ over $\left(W, \Theta_{\mathcal{F}_0}\right) $, so that $\left(F, \nabla_F\right)$ is a foliated line bundle. Let $\Theta_{\mathcal{F}_0}$ on $U_j$ be generated by
\begin{align*}
T_j^\alpha= \sum_{\beta=1}^{n} T_j^\alpha\left(\xi_j^1,..., \xi_j^n \right)\frac{\partial}{\partial \xi_j^\beta} , \,\,\,\,\,\,\,\,\alpha=1,..., p
\end{align*}
and $\Theta_{\mathcal{F}_0}$ on $U_\lambda$ be generated by
\begin{align*}
T_\lambda^\alpha= \sum_{\beta=1}^{n} T_\lambda^\alpha\left(\xi_\lambda^1,..., \xi_\lambda^n \right)\frac{\partial}{\partial \xi_\lambda^\beta} + \xi_\lambda^n G_\lambda^\alpha \left(\xi_\lambda^1,..., \xi_\lambda^n \right) \frac{\partial}{\partial \xi_\lambda^n},\,\,\,\,\,\,\,\,\alpha=1,..., p
\end{align*}
such that 
\begin{align*}
T_\sigma^\alpha &=\sum_{\beta=1}^p r_{\sigma \rho}^{\alpha\beta} T_\rho^\beta,\,\,\,\,\,\,\,\sigma,\rho= i,j,..., \lambda, \mu,...\\
\left[T_\rho^\alpha, T_\rho^\beta   \right]  & = \sum_{\gamma=1}^p c_{\rho \alpha\beta}^\gamma (\xi_\rho) T_\rho^\gamma,\,\,\,\,\,\,\,\,\,\rho= i,j,..., \lambda, \mu,....
\end{align*}
Then $\Theta_{G_0}$ on $U_j^{(v)}$ is generated by
\begin{align*}
Q_j^{(v)\alpha} =  \sum_{\beta=1}^{n} T_j^\alpha\left(z_j^{(v)1},..., z_j^{(v)n} \right)\frac{\partial}{\partial z_j^{(v)\beta} } , \,\,\,\,\,\,\,\,\alpha=1,..., p
\end{align*}
and $\Theta_{G_0}$ on $U_\lambda^\sharp$ is generated by
\begin{align*}
Q_\lambda^\alpha= \sum_{\beta=1}^{n-1} T_\lambda^{\alpha\beta}\left( z_\lambda^1,..., z_\lambda^{n-1}, \left( z_\lambda^n \right)^k  \right) \frac{\partial}{\partial z_\lambda^\beta}  + \frac{1}{k} z_\lambda^n G_\lambda^\alpha \left( z_\lambda^1,..., z_\lambda^{n-1}, \left( z_\lambda^n \right)^k \right) \frac{\partial}{\partial z_\lambda^n},\,\,\,\,\,\,\,\,\,\alpha=1,..., p
\end{align*}
such that
\begin{align*}
Q_i^{(v)\alpha} &= \sum_{\beta=1}^p r_{ij}^{\alpha\beta}\left(z_j^{(s)}\right) Q_j^{(s)\beta}\,\,\,\,\,\,\,\,\,\,\,\,\text{on}\,\,\,\,\, U_i^{(v)}\cap U_j^{(s)}\ne \emptyset\\
Q_\lambda^\alpha &= \sum_{\beta=1}^p r_{\lambda j}^{\alpha\beta}\left(z_j^{(v)} \right) Q_j^{(v)\beta}\,\,\,\,\,\,\, \text{on}\,\,\,\,U_\lambda^\sharp \cap U_j^{(v)}\ne \emptyset\\
Q_{\lambda }^\alpha &= \sum_{\beta=1}^p r_{\lambda \mu}^{\alpha\beta}\left(z_\mu^1,..., z_\mu^{n-1}, \left( z_\mu^n \right)^k \right) Q_\mu^\beta,\,\,\,\,\,\,\,\,\,\,\text{on}\,\,\,\,\,U_\lambda^\sharp \cap U_\mu^\sharp \ne \emptyset
\end{align*}

We define $\left\{\phi_j, \phi_\lambda\right\}\in C^0\left(\mathcal{U}, \mathscr{H}om_{\mathcal{O}_W}\left( \Theta_{\mathcal{F}_0}, \mathcal{O}_W \right) \right)$ where $\phi_j\in \Gamma\left(U_j, \mathscr{H}om_{\mathcal{O}_W}\left( \Theta_{\mathcal{F}_0}, \mathcal{O}_W \right) \right)$ such that
\begin{align*}
\phi_j:\Gamma\left(U_j, \Theta_{\mathcal{F}_0} \right) &\to \Gamma\left( U_j , \mathcal{O}_W \right) \\
 T_j^\alpha &\mapsto 0
\end{align*}
and $\phi_\lambda \in \Gamma\left( U_\lambda, \mathscr{H}om_{\mathcal{O}_W}\left( \Theta_{\mathcal{F}_0}, \mathcal{O}_W \right) \right)$ such that
\begin{align*}
\phi_\lambda:\Gamma\left(U_\lambda, \Theta_{\mathcal{F}_0} \right) &\to \Gamma\left( W, \mathcal{O}_W \right) \\
 T_\lambda^\alpha &\mapsto \frac{1}{k} G_\lambda^\alpha\left(\xi_\lambda^1,..., \xi_\lambda^n \right)
\end{align*}

We will show that $\left(\left\{\phi_j, \phi_\lambda \right\}, \left\{f_{ij},  f_{\lambda j}, f_{\lambda \mu} \right\} \right) \in C^0\left( \mathcal{U}, \mathscr{H}om_{\mathcal{O}_W}\left( \Theta_{\mathcal{F}_0}, \mathcal{O}_W \right) \right) \bigoplus C^1\left( \mathcal{U}, \mathcal{O}_W^* \right)$ satisfies the cocycle condition defining a foliated line bundle on $\left(W, \Theta_{\mathcal{F}_0}\right)$. It is sufficient to show that
\begin{align}
\phi_j\left(T_\lambda^\alpha\right)- \phi_\lambda\left(T_\lambda^\alpha\right) +\left[ T_\lambda^\alpha, \log g_{\lambda j}^{(0) n} \right]=0 \label{bc2}\\
\phi_\mu\left(T_\lambda^\alpha \right) - \phi_\lambda\left(T_\lambda^\alpha \right)  + \left[ T_\lambda^\alpha, \log \phi_{\lambda \mu}\left( \xi_\mu^1,..., \xi_\mu^{n-1}, \sqrt[k]{\xi_\mu^n} \right) \right]=0 \label{bc3}\\
\left[T_\lambda, \phi_\lambda\left( T_\lambda^\beta \right) \right] - \left[ T_\lambda^\beta, \phi_\lambda\left( T_\lambda^\alpha \right) \right]- \phi_\lambda\left(\left[ T_\lambda^\alpha, T_\lambda^\beta \right] \right)=0 \label{bc4}
\end{align}

Let us show $(\ref{bc2})$. We note that from $Q_\lambda^\alpha = \sum_{\beta=1}^p r_{\lambda j}^{\alpha\beta}\left(z_j^{(0)} \right) Q_j^{(0)\beta}$, and $z_\lambda^\alpha= g_{\lambda j}^{(0)\alpha}\left(z_j^{(0)} \right)$, we have
\begin{align*}
\sum_{\beta=1}^{n-1} T_\lambda^{\alpha\beta}\left( z_\lambda^1,..., z_\lambda^{n-1}, \left(  z_\lambda^n \right)^k \right) \frac{\partial}{\partial z_\lambda^\beta} + \frac{1}{k} z_\lambda^n  G_\lambda^\alpha\left( z_\lambda^1,..., z_\lambda^{n-1}, \left( z_\lambda^n \right)^k \right) \frac{\partial}{\partial z_\lambda^n} = \sum_{\gamma=1}^p \sum_{\eta=1}^n r_{\lambda j}^{\alpha \gamma}\left( z_j^{(0)} \right) T_j^{\gamma \eta}\left( z_j^{(0)} \right) \frac{\partial}{\partial z_j^{(0)\eta}}
\end{align*}

Then we have
\begin{align*}
\frac{1}{k} g_{\lambda j}^{(0)n} G_\lambda^\alpha\left(z_\lambda^1,..., z_\lambda^{n-1} , \left( z_\lambda^n  \right)^k    \right)= \sum_{\gamma=1}^p \sum_{\eta=1}^n r_{\lambda j}^{\alpha \gamma}\left(z_j^{(0)} \right) T_j^{\gamma \eta}\left(z_j^{(0)} \right)\frac{\partial g_{\lambda j}^{(0) n}}{\partial z_j^{(0)\eta}}
\end{align*}
By taking the inverse of $\pi$: $\left(\xi_j^1,..., \xi_j^n \right)\mapsto \left(z_j^{(0)1}.,..., z_j^{(0)n} \right)=\left(\xi_j^1,..., \xi_j^n \right)$ and $\left( \xi_\lambda^1, ..., \xi_\lambda^n \right)\mapsto \left(z_\lambda^1,..., z_\lambda^n \right)= \left( \xi_\lambda^1,..., \xi_\lambda^{n-1}, \sqrt[k]{\xi_\lambda^n} \right)$ on $U_\lambda \cap U_j$ which is well-defined, we have
\begin{align*}
\frac{1}{k} G_\lambda^\alpha\left( \xi_\lambda^1,..., \xi_\lambda^n \right) = \left[ T_\lambda^\alpha, \log g_{\lambda j}^{(0) n} \right]
\end{align*}
This implies $(\ref{bc2})$. 

Let us show $(\ref{bc3})$. From $Q_\lambda^\alpha= \sum_{\beta=1}^p r_{\lambda \mu}^{\alpha\gamma}\left(z_\mu^1,..., z_\mu^{n-1}, \left(  z_\mu^n \right)^k \right) Q_\mu^\gamma$, we have
{\small{\begin{align}
&\sum_{\beta=1}^{n-1} T_\lambda^{\alpha\beta}\left(z_\lambda^1,...,z_\lambda^{n-1}, \left(z_\lambda^n\right)^k \right) \frac{\partial}{\partial z_\lambda^\beta} + \frac{1}{k}z_\lambda^n G_\lambda^\alpha\left(z_\lambda^1,...,z_\lambda^{n-1}, \left(z_\lambda^n\right)^k \right) \frac{\partial}{\partial z_\lambda^n} \label{bc5}\\
&=\sum_{\gamma=1}^p r_{\lambda \mu}^{\alpha \gamma}\left(z_\mu^1,...,z_\mu^{n-1}, \left(z_\mu^n\right)^k \right) \left( \sum_{\eta=1}^{n-1} T_\mu^{\gamma\eta}\left(z_\mu^1,...,z_\mu^{n-1}, \left(z_\mu^n \right)^k \right) \frac{\partial}{\partial z_\mu^\eta} + \frac{1}{k}z_\mu^n G_\mu^\alpha\left(z_\mu^1,...,z_\mu^{n-1}, \left(z_\mu^n \right)^k \right) \frac{\partial}{\partial z_\mu^n}    \right) \notag
\end{align}}}
Since  $z_\lambda^n=\phi_{\lambda \mu}\left(z_\mu\right)\cdot z_\mu^n$ on $U_\lambda^\sharp \bigcap U_\mu^\sharp\ne \emptyset$, we have $\frac{\partial}{\partial z_\mu^\eta}= \sum_{\alpha=1}^{n-1} \frac{\partial z_\lambda^\alpha}{\partial z_\mu^\eta}\frac{\partial}{\partial z_\lambda^\alpha} +  \frac{ \partial \phi_{\lambda \mu}(z_\mu)}{\partial z_\mu^\eta} z_\mu^n \frac{\partial}{\partial z_\lambda^n } $ for $\eta=1,...,n-1$ and $\frac{\partial}{\partial z_\mu^n}=\sum_{\alpha=1}^{n-1}\frac{\partial z_\lambda^\alpha}{\partial z_\mu^n}\frac{\partial}{\partial z_\lambda^\alpha} + \left(\frac{\partial \phi_{\lambda \mu}(z_\mu)}{\partial z_\mu^n}z_\mu^n + \phi_{\lambda \mu}(z_\mu)\right)\frac{\partial}{\partial z_\lambda^n}$. Then from $(\ref{bc5})$ we have (by considering the coefficient of $\frac{\partial}{\partial z_\lambda^n}$)
\begin{align*}
&\frac{1}{k} z_\lambda^n G_\lambda^\alpha \left(z_\lambda^1,...,z_\lambda^{n-1}, \left(z_\lambda^n\right)^k \right)\\
&= \sum_{\gamma=1}^p r_{\lambda\mu}^{\alpha \gamma}\left(z_\mu^1,..., z_\mu^{n-1}, \left(z_\mu^n \right)^k\right) \left( \sum_{\eta=1}^{n-1} T_\mu^{\gamma \eta}\left(z_\mu^1,...,z_\mu^{n-1}, \left(z_\mu^n \right)^k \right) \frac{\partial \phi_{\lambda \mu}}{\partial z_\mu^\eta} z_\mu^n  + \frac{1}{k} z_\mu^n G_\mu^\alpha \left( z_\mu^1,...,z_\mu^{n-1}, \left(z_\mu^n \right)^k \right) \frac{\partial \phi_{\lambda \mu}}{\partial z_\mu^n} z_\mu^n \right)\\
&+ \sum_{\gamma=1}^p r_{\lambda \mu}^{\alpha \gamma}\left(z_\mu^1,..., z_\mu^{n-1}, \left(z_\mu^n\right)^k \right) \cdot  \frac{1}{k} z_\mu^n G_\mu^\alpha \left(z_\mu^1,..., z_\mu^{n-1}, \left(z_\mu^n \right)^k \right)\phi_{\lambda \mu}(z_\mu)
\end{align*}
By dividing $z_\lambda^n$, we have 
\begin{align}
&\frac{1}{k} G_\lambda^\alpha \left(z_\lambda^1,...,z_\lambda^{n-1}, \left(z_\lambda^n \right)^k \right) - \sum_{\gamma=1}^p r_{\lambda \mu}^{\alpha \gamma}\left(z_\mu^1,..., z_\mu^{n-1}, \left(z_\mu^n\right)^k \right) \cdot \frac{1}{k}  G_\mu^\alpha \left(z_\mu^1,..., z_\mu^{n-1}, \left(z_\mu^n\right)^k \right) \label{bc6}\\
&=\frac{1}{\phi_{\lambda \mu}} \sum_{\gamma=1}^p r_{\lambda\mu}^{\alpha \gamma}\left(z_\mu^1,..., z_\mu^{n-1}, \left(z_\mu^n \right)^k \right) \left( \sum_{\eta=1}^{n-1} T_\mu^{\gamma \eta}\left(z_\mu^1,...,z_\mu^{n-1},\left(z_\mu^n \right)^k \right) \frac{\partial \phi_{\lambda \mu}}{\partial z_\mu^\eta}  + \frac{1}{k} z_\mu^n G_\mu^\alpha \left( z_\mu^1,...,z_\mu^{n-1}, \left(z_\mu^n\right)^k  \right) \frac{\partial \phi_{\lambda \mu}}{\partial z_\mu^n} \right) \notag
\end{align}
We note that $\phi_{\lambda \mu}\left(z_\mu^1,..., z_\mu^n\right)= \phi_{\lambda \mu}\left(z_\mu^1,..., z_{\mu}^{n-1}, \rho^s \cdot z_\mu^n\right)$ for $s=1,...,{k-1}$ (see \cite{Wav68} p.930). By taking the derivative with respect to $z_\mu^n$, we have
\begin{align*}
&\frac{\partial \phi_{\lambda \mu}}{\partial z_\mu^n}\left(z_\mu^1,..., z_\mu^n\right)= \frac{\partial \phi_{\lambda \mu}}{\partial z_\mu^n}\left( z_\mu^1,..., z_\mu^n, \rho^s \cdot z_\mu^n \right)\cdot \rho^s\\
&\Longrightarrow  z_\mu^n\cdot \frac{\partial \phi_{\lambda \mu}}{\partial z_\mu^n}(z_\mu^1,..., z_\mu^n)= \frac{\partial \phi_{\lambda \mu}}{\partial z_\mu^n}\left( z_\mu^1,..., z_\mu^n, \rho^s \cdot z_\mu^n \right)\cdot  \rho^s z_\mu^n
\end{align*}
This implies that $\sqrt[k]{\xi_\mu^n}\cdot \frac{\partial \phi_{\lambda \mu}}{\partial z_\mu^n} \left(\xi_\mu^1,..., \xi_\mu^{n-1}, \sqrt[k]{\xi_\mu^n} \right)$ is a well-defined holomorphic function on $U_\lambda \cap U_\mu$. From $z_\lambda^n= \phi_{\lambda \mu}(z_\mu) z_\mu^n$, we have $\frac{\partial z_\lambda^n}{\partial z_\mu^n}= \phi_{\lambda \mu}(z_\mu)+ \frac{\partial \phi_{\lambda \mu}}{\partial z_\mu^n}\cdot z_\mu^n$, we see that $\frac{\partial z_\lambda^n}{\partial z_\mu^n}\left(\xi_\mu^1,..., \xi_\mu^{n-1}, \sqrt[k]{\xi_\mu^n}\right)$ is a well-defined holomorphic function. We note that
\begin{align*}
\frac{\partial \xi_\lambda^n}{\partial z_\mu^n}&= \sum_{\alpha=1}^{n-1} \frac{\partial \xi_\lambda^n}{\partial z_\lambda^\alpha}\cdot \frac{\partial z_\lambda^\alpha}{\partial z_\mu^n} + \frac{\partial \xi_\lambda^n}{ \partial z_\lambda^n}\frac{\partial z_\lambda^n}{\partial z_\mu^n} = k \left(z_\lambda^n\right)^{k-1} \frac{\partial z_\lambda^n}{\partial z_\mu^n} \\
\frac{\partial \xi_\lambda^n}{\partial z_\mu^n}&= \sum_{\alpha=1}^{n-1} \frac{\partial \xi_\lambda^n}{\partial \xi_\mu^n}\left(z_\mu^1,..., z_\mu^{n-1}, \left(z_\mu^n\right)^k\right) \frac{\partial \xi_\mu^\alpha}{\partial z_\mu^n} + \frac{\partial \xi_\lambda^n}{\partial \xi_\mu^n} \left(z_\mu^1,..., z_\mu^{n-1}, \left(z_\mu^n\right)^k\right)\frac{\partial \xi_\mu^n}{\partial z_\mu^n}= \frac{\partial \xi_\lambda^n}{\partial \xi_\mu^n} \left(z_\mu^1,..., z_\mu^{n-1}, \left(z_\mu^n\right)^k\right) k\left(z_\mu^n\right)^{k-1}
\end{align*}
This implies that
\begin{align*}
\frac{\partial \xi_\lambda^n}{\partial \xi_\mu^n}= \phi_{\lambda \mu}^{k-1}\left(\xi_\mu^1,..., \xi_\mu^{n-1}, \sqrt[k]{\xi_\mu^n} \right) \cdot \frac{\partial z_\lambda^n}{\partial z_\mu^n}\left(\xi_\mu^1,..., \xi_\mu^{n-1}, \sqrt[k]{\xi_\mu^n}\right)
\end{align*}
On the other hand, we note that from $z_\lambda^n=\phi_{\lambda \mu}(z_\mu)\cdot z_\mu$, we have $\left(z_\lambda^n\right)^k = \left(\phi_{\lambda\mu}(z_\mu) \right)^k \cdot  \left( z_\mu^n \right)^k$. By plugging $\left(\xi_\lambda^1,..., \xi_\lambda^{n-1}, \left(\xi_\lambda^n \right)^{\frac{1}{k}} \right)$ into $\left(z_\lambda^1,..., z_\lambda^n \right)$ and $\left(\xi_\mu^1,..., \xi_\mu^{n-1}, \left(\xi_\mu^n \right)^{\frac{1}{k}} \right)$ into $\left(z_\mu^1,..., z_\mu^n \right)$, we have $\xi_\lambda^n= \left( \phi_{\lambda \mu}\left(\xi_\mu^1,..., \xi_\mu^{n-1}, \sqrt[k]{\xi_\mu^n} \right) \right)^k \cdot \xi_\mu^n$ on $U_\lambda \cap U_\mu$. We set  $\Phi_{\lambda \mu}(\xi_\mu)= \phi_{\lambda \mu}\left(\xi_\mu^1,..., \xi_\mu^{n-1}, \sqrt[k]{\xi_\mu^n} \right) $. Then we have 
\begin{align*}
\frac{\partial \xi_\lambda^n}{\partial \xi_\mu^n}&= \frac{\partial (\Phi_{\lambda\mu}^k \xi_\mu^n )}{\partial \xi_\mu^n}= k \Phi_{\lambda \mu}^{k-1} \frac{\partial \Phi_{\lambda\mu}}{\partial \xi_\mu^n} \xi_\mu^n + \Phi_{\lambda\mu}^k\\
\phi_{\lambda\mu}^{k-1} \frac{\partial z_\lambda^n}{\partial z_\mu^n}&= \phi_{\lambda \mu}^{k-1} \frac{\partial (\phi_{\lambda \mu} z_\mu^n)}{\partial z_\mu^n}= \phi_{\lambda\mu}^{k-1} \frac{\partial \phi_{\lambda \mu}}{\partial z_\mu^n} z_\mu^n + \phi_{\lambda \mu}^k
\end{align*}
This implies that
\begin{align*}
\sqrt[k]{\xi_\mu^n}\cdot \frac{\partial \phi_{\lambda\mu}}{\partial z_\mu^n}\left(\xi_\mu^1,...,\xi_\mu^{n-1}, \sqrt[k]{\xi_\mu^n} \right)= k \xi_\mu^n\frac{\partial \phi_{\lambda \mu}\left(\xi_\mu^1,...,\xi_\mu^{n-1}, \sqrt[k]{\xi_\mu^n}\right) }{\partial \xi_\mu^n } = k \xi_\mu^n \frac{\partial \Phi_{\lambda\mu}}{\partial \xi_\mu^n}
\end{align*}
By plugging $\left(\xi_\lambda^1,..., \xi_\lambda^{n-1}, \left(\xi_\lambda^n \right)^{\frac{1}{k}} \right)$ into $\left(z_\lambda^1,..., z_\lambda^n \right)$, and $\left(\xi_\mu^1,..., \xi_\mu^{n-1}, \left(\xi_\mu^n \right)^{\frac{1}{k}} \right)$ into $\left(z_\mu^1,..., z_\mu^n \right)$ in $(\ref{bc6})$, we have
\begin{align*}
&\frac{1}{k} G_\lambda^\alpha \left(\xi_\lambda^1,...,\xi_\lambda^n \right) - \sum_{\gamma=1}^p r_{\lambda \mu}^{\alpha \gamma}\left(\xi_\mu \right) \frac{1}{k}  G_\mu^\alpha\left(\xi_\mu^1,..., \xi_\mu^n \right) \\
&=\frac{1}{\Phi_{\lambda\mu}(\xi_\mu)} \sum_{\gamma=1}^p r_{\lambda\mu}^{\alpha \gamma}\left(\xi_\mu \right) \left( \sum_{\eta=1}^{n-1} T_\mu^{\gamma \eta}\left(\xi_\mu^1,..., \xi_\mu^n\right) \frac{\partial  \Phi_{\lambda\mu}(\xi_\mu)}{\partial \xi_\mu^\eta}  + \xi_\mu^n G_\mu^\alpha \left( \xi_\mu^1,.. , \xi_\mu^n \right) \frac{\partial \Phi_{\lambda\mu}(\xi_\mu) }{\partial \xi_\mu^n} \right) \\
&= \left[T_\lambda^\alpha, \log \Phi_{\lambda\mu}\left(\xi_\mu\right) \right]
\end{align*}
This implies $(\ref{bc3})$.

Let us show $(\ref{bc4})$. Since  $\left[T_\lambda^\alpha, T_\lambda^\beta  \right] = \sum_{\gamma=1}^p c_{\lambda \alpha\beta}^\gamma(\xi_\lambda) T_\lambda^\gamma$,
\begin{align*}
\left[ \sum_{\eta=1}^{n-1} T_\lambda^{\alpha\eta}(\xi_\lambda) \frac{\partial}{\partial \xi_\lambda^\eta} + \xi_\lambda^n G_\lambda^\alpha(\xi_\lambda)\frac{\partial}{\partial \xi_\lambda^n},  \sum_{\delta=1}^{n-1} T_\lambda^{\beta\delta}(\xi_\lambda)\frac{\partial}{\partial \xi_\lambda^\delta} + \xi_\lambda^n G_\lambda^\beta(\xi_\lambda)\frac{\partial}{\partial \xi_\lambda^n}   \right] =\sum_{\gamma=1}^p c_{\lambda \alpha\beta}^\gamma (\xi_\lambda)\left(   \sum_{\eta=1}^{n-1} T_\lambda^{\gamma\eta}(\xi_\lambda) \frac{\partial}{\partial \xi_\lambda^\eta} + \xi_\lambda^n G_\lambda^\gamma (\xi_\lambda)\frac{\partial}{\partial \xi_\lambda^n}  \right)
\end{align*}
By considering the coefficient of $\frac{\partial}{\partial \xi_\lambda^n}$, we have
\begin{align*}
&\sum_{\eta=1}^{n-1} T_\lambda^{\alpha\eta}(\xi_\lambda)\xi_\lambda^n\frac{\partial G_\lambda^\beta}{\partial \xi_\lambda^\eta} \frac{\partial}{\partial \xi_\lambda^n} +\xi_\lambda^n G_\lambda^\alpha(\xi_\lambda)G_\lambda^\beta(\xi_\lambda)\frac{\partial}{\partial \xi_\lambda^n}+ \left(\xi_\lambda^n\right)^2 G_\lambda^\alpha \frac{\partial G_\lambda^\beta}{\partial \xi_\lambda^n}\frac{\partial}{\partial \xi_\lambda^n} \\
&- \sum_{\delta=1}^{n-1} T_\lambda^{\beta\delta}(\xi_\lambda)\xi_\lambda^n \frac{\partial G_\lambda^\alpha}{\partial \xi_\lambda^\delta}\frac{\partial}{\partial \xi_\lambda^n}- \xi_\lambda^n G_\lambda^\beta(\xi_\lambda) G_\lambda^\alpha(\xi_\lambda)\frac{\partial}{\partial \xi_\lambda^n} - \left( \xi_\lambda^n \right)^2 G_\lambda^\beta \frac{\partial G_\lambda^\alpha}{\partial \xi_\lambda^n}\frac{\partial}{\partial \xi_\lambda^n}\\
&= \sum_{\gamma=1}^p c_{\lambda \alpha\beta}^\gamma(\xi_\lambda) \xi_\lambda^n G_\lambda^\gamma(\xi_\lambda)\frac{\partial}{\partial \xi_\lambda^n}
\end{align*}
Then we have
\begin{align*}
\sum_{\eta=1}^{n-1} T_\lambda^{\alpha\eta} \frac{\partial G_\lambda^\beta}{\partial \xi_\lambda^\eta}+ \xi_\lambda^n G_\lambda^\alpha \frac{\partial G_\lambda^\beta}{\partial \xi_\lambda^n} - \sum_{\delta=1}^{n-1} T_\lambda^{\beta\delta} \frac{\partial G_\lambda^\alpha}{\partial \xi_\lambda^\delta}- \xi_\lambda^n G_\lambda^\beta \frac{\partial G_\lambda^\alpha}{\partial \xi_\lambda^n} = \sum_{\gamma=1}^p c_{\lambda\alpha\beta}^\gamma G_\lambda^\gamma \frac{\partial}{\partial \xi_\lambda^n}  \iff  \left[T_\lambda^ \alpha, G_\lambda^\beta \right]- \left[T_\lambda^\beta, G_\lambda^\alpha \right] - \sum_{\gamma=1}^p c_{\lambda \alpha\beta}^\gamma G_\lambda^\gamma=0
\end{align*}
This implies $(\ref{bc4})$. Hence there is a foliated connection $\nabla_F$ on $F$ over $\left(W, \Theta_{\mathcal{F}_0}\right)$, so that we can make $\left(F, \Theta_{\tilde{\mathcal{F}}_0}\right)$ as a  foliated complex manifold with $\Theta_{\tilde{\mathcal{F}}_0}$ locally free of rank $p$. We recall that $\Theta_{\tilde{\mathcal{F}}_0}$ is generated by 
\begin{align*}
\tilde{T}_j^\alpha&:= T_j^\alpha,\,\,\,\,\,\alpha=1,...,p\,\,\,\,\,\,\,\,\textnormal{on}\,\,\, U_j\times \mathbb{C}\,\,\,\,\,\textnormal{with coordinate}\,\,\,(\xi_j, t_j)\\
\tilde{T}_\lambda^\alpha&:= T_\lambda^\alpha + \frac{1}{k}t_\lambda G_\lambda^\alpha \frac{\partial}{\partial t_\lambda}, \,\,\,\,\,\alpha=1,..., p\,\,\,\,\,\,\,\,\textnormal{on}\,\,\, U_\lambda \times \mathbb{C}\,\,\,\,\,\textnormal{with coordinate}\,\,\,(\xi_\lambda, t_\lambda)
\end{align*}

Now we define a holomorphic map $\Psi:M\to F$ by
\begin{align*}
\Psi: U_j^{(v)} &\to U_j\times \mathbb{C}\\
\left(z_j^{(v)} \right) &\mapsto \left( z_j^{(v)} ; \rho^v \right)\\
\Psi:U_\lambda^\sharp &\to U_\lambda \times \mathbb{C}\\
(z_\lambda^1,..., z_\lambda^n) &\mapsto \left( z_\lambda^1,..., z_\lambda^{n-1} , \left(z_\lambda^n \right)^k ; z_\lambda^n \right) =\left( \xi_\lambda^1,..., \xi_\lambda^{n-1}, \xi_\lambda^n ; t_\lambda \right)
\end{align*}
Then we may identify $M$ with the image. The image has the equation $t_j^k=1$ over $U_j \times \mathbb{C}$ and $t_\lambda^k=\xi_\lambda^n$ over $U_\lambda\times \mathbb{C}$. We claim that the image is invariant with respect to $\Theta_{\tilde{\mathcal{F}}_0}$. In fact,
\begin{align*}
\left[ \tilde{T}_j^\alpha , t_j^k-1\right] &= 0\\
\left[ \tilde{T}_\lambda^\alpha, t_\lambda^k - \xi_\lambda^n \right] &= \left[   \sum_{\beta=1}^{n-1} T_\lambda^{\alpha\beta}(\xi_\lambda)\frac{\partial}{\partial \xi_\lambda^\beta}+ \xi_\lambda^n G_\lambda^\alpha ( \xi_\lambda) \frac{\partial}{\partial \xi_\lambda^n} + \frac{1}{k} t_\lambda G_\lambda^\alpha(\xi_\lambda) \frac{\partial}{\partial t_\lambda} , t_\lambda^k - \xi_\lambda^n    \right] \\
&= t_\lambda^k G_\lambda^\alpha(\xi_\lambda) - \xi_\lambda^n G_\lambda^\alpha(\xi_\lambda)= \left(t_\lambda^k - \xi_\lambda^n \right) G_\lambda^\alpha( \xi_\lambda)
\end{align*}

Next we will show that $\Psi$ induces a foliated isomorphism between $\left(M, \Theta_{\mathcal{G}_0}\right)$ and the image with the foliation induced from $\left(F, \Theta_{\tilde{\mathcal{F}}_0} \right)$, i.e. we claim that $\Psi: \left(M, \Theta_{\mathcal{G}_0}\right)\cong \left(F, \Theta_{\tilde{\mathcal{F}}_0} \right)|_{\Psi(M)}$. In fact,
\begin{align*}
Q_\lambda^\alpha= \sum_{\beta=1}^{n-1} T_\lambda^{\alpha\beta} & \left(z_\lambda^1,...,z_\lambda^{n-1}, \left( z_\lambda^n \right)^k \right) \frac{\partial}{\partial z_\lambda^\beta} + \frac{1}{k} z_\lambda^n G_\lambda^\alpha \left(z_\lambda^1,...,z_\lambda^{n-1}, \left(z_\lambda^n \right)^k \right)\frac{\partial}{\partial z_\lambda^n} \\
\Longrightarrow \tilde{\Psi}(Q_\lambda^\alpha) &= \sum_{\beta=1}^{n-1} T_\lambda^{\alpha\beta} \left(z_\lambda^1,...,z_\lambda^{n-1}, \left( z_\lambda^n \right)^k \right) \frac{\partial}{\partial \xi_\lambda^\beta} + \frac{1}{k} z_\lambda^n G_\lambda^\alpha \left(z_\lambda^1,...,z_\lambda^{n-1}, \left(z_\lambda^n \right)^k \right)\left( k (z_\lambda^n)^{k-1} \frac{\partial}{\partial \xi_\lambda^n} + \frac{\partial}{\partial t_\lambda} \right) \\
&=  \sum_{\beta=1}^{n-1} T_\lambda^{\alpha\beta} \left(z_\lambda^1,...,z_\lambda^{n-1}, \left( z_\lambda^n \right)^k \right) \frac{\partial}{\partial \xi_\lambda^\beta} + \left( z_\lambda^n \right)^k G_\lambda^\alpha \left(z_\lambda^1,...,z_\lambda^{n-1}, \left(z_\lambda^n \right)^k \right) \frac{\partial}{\partial \xi_\lambda^n} + \frac{1}{k} z_\lambda^n G_\lambda^\alpha \left( z_\lambda^1 , ..., z_\lambda^{n-1} , \left( z_\lambda^n \right)^k \right) \frac{\partial}{\partial t_\lambda}\\
&=\Psi^*\left( \tilde{T}_\lambda^\alpha \right)
\end{align*}
where $\tilde{\Psi}:\Theta_M \to \Psi^* \Theta_F$ is the natural map. On the other hand, since $ Q_j^{(v)\alpha} = \sum_{\beta=1}^n T_j^\alpha \left( z_j^{(v)1},..., z_j^{(v)n} \right) \frac{\partial}{\partial z_j^{(v)\beta}} $, 
\begin{align*}
\tilde{\Psi } \left( Q_j^{(v)\alpha} \right)  = \sum_{\beta=1}^n T_j^\alpha \left( z_j^{(v)1},..., z_j^{(v)n} \right) \frac{\partial}{\partial \xi_j^\beta} = \Psi^*\left( \tilde{T}_j^\alpha \right)
\end{align*}
which proves the claim.

\end{proof}

\begin{corollary}[compare \cite{Wav68} Corollary p.929] \label{cs1}
If $C\subset \left(W, \Theta_{\mathcal{F}_0}\right)$ is a nonsingular invariant hypersurface $($connected if $k\ne 2$$)$, then a necessary and sufficient condition for $C$ to be the invariant branch locus of a foliated cyclic $k$-sheeted covering of $\left(W, \Theta_{\mathcal{F}_0}\right)$ is that there is a foliated line bundle $\left(F, \nabla_F\right)$ over $\left(W, \Theta_{\mathcal{F}_0}\right)$ with $\left(F^k, \nabla_{F^k}\right) =\left[\left(C, \nabla_C \right)\right]$.
\end{corollary}

\begin{proof}

 let $C$ be defined by $\phi_i=0$ on $V_i$ for a covering $\mathcal{U}=\left\{V_i\right\}$ of $W$. Since $C$ is invariant, $\left[T_i^\alpha, \phi_i\right]= G_i^\alpha \cdot \phi_i$ for some $G_i^\alpha \in \Gamma\left(V_i, \mathcal{O}_W\right)$. Then $\left[\left(C, \nabla_C\right)\right]$ is defined by $\left( \left\{T_i^\alpha \mapsto G_i^\alpha\right\}, \left\{\frac{\phi_i}{\phi_j}\right\} \right)\in C^0\left( \mathcal{U}, \mathscr{H}om_{\mathcal{O}_W}\left(\Theta_{\mathcal{F}_0}, \mathcal{O}_W \right) \right)$ $\bigoplus C^1\left( \mathcal{U}, \mathcal{O}_W^* \right)$. Assume that there is a foliated line bundle $\left(F, \nabla_F\right)$ such that $\left(F^k, \nabla_{F^k}\right)=\left[\left( C, \nabla_C\right)\right]$. Then $\left(F, \nabla_F\right)$ is defined by $\left( \left\{ T_i^\alpha \to \frac{1}{k} G_\lambda  \right\} , \left\{ f_{ij} \right\} \right)$ with $f_{ij}^k= \frac{\phi_i}{\phi_j}$. We note that $\nabla_F$ makes $(F, \Theta_{\mathcal{G}})$ a foliated complex manifold as in discussion before Theorem \ref{bc1} and $\left\{t_i^k - \phi_i\right\}$ determine a nonsingular manifold $M$ which is a cyclic $k$-sheeted covering of $W$ having branch locus $C$. Moreover $t_i^k-\phi_i$ is invariant with respect to $\Theta_{\mathcal{G}}$. In fact, 
 \begin{align*}
 \left[ T_i^\alpha + t_i \frac{1}{k} G_i^\alpha \frac{\partial}{\partial t_i}, t_i^k- \phi_i  \right]= t_i^k G_i^\alpha - G_i^\alpha \cdot \phi_i= \left(t_i^k- \phi_i\right) G_i^\alpha
 \end{align*}
 Hence $\left(M, \Theta_\mathcal{G}|_M\right)$ is a foliated $k$-sheeted cyclic covering of $(W, \Theta_{\mathcal{F}_0})$.

Conversely, by Theorem \ref{bc1}, we can find a foliated line bundle $\left(F, \nabla_F\right)$ such that $M\subset \left(F, \Theta_{\tilde{\mathcal{F}}_0}\right)$ is an invariant hypersurface given by $t_j^k=1$ over $U_j\times \mathbb{C}$ and $t_\lambda^k= \xi_\lambda^n$ over $U_\lambda \times \mathbb{C}$ as in the proof of Theorem \ref{bc1}, and $[C,\nabla_C]$ is defined by $\left( \left\{ \left\{ T_j^\alpha\mapsto 0 \right\} , \left\{ T_\lambda^\alpha \mapsto G_\lambda^\alpha \right\}   \right\} ,  \left\{h_{\lambda \mu}= \frac{\xi_\lambda^n}{\xi_\mu^n } , h_{\lambda j} = \frac{\xi_j^n}{1} ,  h_{ij }=\frac{1}{1} \right\} \right)$, so that $\left[ \left(C, \nabla_C \right) \right]= \left(F^k, \nabla_{F^k}\right) $.

\end{proof}

If $F$ is a line bundle over $W$, we denote by $F^\sharp$ the $\mathbb{P}_\mathbb{C}^1$-bundle over $W$ obtained from $F$. We review the computation of cohomology groups on $F^\sharp$ and branched coverings contained in $F^\sharp$. Let $\{V_i\}$ be a covering of $W$ by polycylinders such that $L$ is trivialized with transition functions $\{f_{ij}\}$, we obtain a covering $\{U_\alpha\}=\{U_{2i}, U_{2i+1}\}$ of $F^\sharp$ by Stein manifolds in the following way: let $\pi:F^\sharp \to W$ denote projection $U_{2i}\cong U_{2i+1}\cong V_i\times \mathbb{C}, \pi^{-1}(V_i) = U_{2i} \bigcup U_{2i+1}$. If coordinates in $U_\alpha$ are $(z_\alpha;\xi_\alpha)$, then $\xi_{2i}= f_{ij}\cdot \xi_{2j}$ in $U_{2i}\bigcap U_{2j}$ and $\xi_{2i}= \frac{1}{\xi_{2i+1}}$ on $U_{2i}\bigcap U_{2i+1}$. The map $\pi$ is the projection $U_{2i}, U_{2i+1}\to V_i$.

\begin{proposition}[see \cite{Wav68} Proposition 2.1 p.933]\label{cs11}
Let $\mathcal{G}$ be the sheaf over $F^\sharp$ of germs of holomorphic sections of a vector bundle $G^\sharp$ whose transition matrices $\phi_{\alpha\beta}^\sharp$ with respect to the covering $\{U_\alpha\}$ have the form
\begin{equation*}
\phi_{\alpha\beta}^\sharp=\begin{cases}
\phi_{ij}\,\,\,\,\,\,\,\,\,\,\,\,\,\,\,\,\,\,\,\,\,\,\,\,\textnormal{if}\,\,\,\alpha=2i, \beta=2j\\
\frac{1}{\xi_{2i+1}^s}\cdot \textnormal{I}_r \,\,\,\,\,\,\,\,\,\,\textnormal{if}\,\,\,\alpha=2i, \beta= 2i+1
\end{cases}
\end{equation*}
where $s$ is a non-negative integer, and the $\phi_{ij}$ are transition matrices for a vector bundle $\Phi$ on $W$. Then
\begin{align*}
H^m \left( F^\sharp, \mathcal{G} \right) \cong  \bigoplus_{\lambda =0}^s H^m\left( W, \mathcal{O}\left( \Phi\otimes F^{-\lambda} \right) \right)
\end{align*}
\end{proposition}

\begin{proposition}[see \cite{Wav68} Proposition 2.2 p.935]\label{cs34}
Let $\mathcal{G}$ be the sheaf over $F^\sharp$ of germs of holomorphic sections of a line bundle $\Phi^\sharp$ whose transition functions $\phi_{\alpha\beta}^\sharp$ with respect to the covering $\{U_\alpha\}$ have the form
\begin{equation*}
\phi_{\alpha\beta}^\sharp=\begin{cases}
\phi_{ij}\,\,\,\,\,\,\,\,\,\,\,\,\,\,\,\,\,\,\,\,\,\,\,\,\textnormal{if}\,\,\,\alpha=2i, \beta=2j\\
\xi_{2i+1}^s  \,\,\,\,\,\,\,\,\,\,\textnormal{if}\,\,\,\alpha=2i, \beta= 2i+1
\end{cases}
\end{equation*}
where $s$ is a positive integer, and the $\phi_{ij}$ are transition matrices for a line bundle $\Phi$ on $W$. Then
\begin{align*}
H^m\left( F^\sharp, \mathcal{G} \right) &\cong \bigoplus_{\lambda=1}^{s-1} H^{m-1}\left( W, \mathcal{O}\left(\Phi\otimes F^\lambda  \right) \right) \\
 &\left( =0 \,\,\,\,\,\textnormal{for}\,\,\, s=1\right)
\end{align*}
\end{proposition}

If $M\subset F^\sharp$ is a nonsingular hypersurface, we denote by $[M]$ the line bundle over $F^\sharp$ defined by $M$ considered as a hypersurface in $F^\sharp$. 

\begin{remark}
Given a foliated $k$-cyclic covering $\pi: \left(M, \Theta_{\mathcal{G}_0} \right)\to \left( W, \Theta_{\mathcal{F}_0} \right)$ with $\Theta_{\mathcal{F}_0}$ locally free between compact foliated complex manifolds, $\left( M, \Theta_{\mathcal{F}_0}\right)$ is a non-degenerate invariant compact complex submanifold of $\left(F, \Theta_{\tilde{\mathcal{F}}_0} \right)$ defined by a foliated line bundle $\left(F, \nabla_F \right)$ over $\left(W, \Theta_{\mathcal{F}_0} \right)$ as in \textnormal{Theorem \ref{bc1}}. Let $\left(F^\sharp, \Theta_{\tilde{\mathcal{F}}_0^\sharp} \right)$ be the foliated $\mathbb{P}_\mathbb{C}^1$-bundle obtained from $\left(F, \Theta_{\tilde{\mathcal{F}}_0} \right)$ as in the discussion before \textnormal{Theorem \ref{bc1}}. Then we can find a polycylindrical covering $\left\{V_i\right\}$ of $W$ so that, with respect to the associated covering $\{U_\alpha\}$ of $F^\sharp$ as above, the bundle $[M]$ has transition functions
\begin{equation}\label{cs25}
\phi^\sharp_{\alpha\beta}= \begin{cases}
(f_{ij})^k \,\,\,\,\,\,\,\,\,\textnormal{if}\,\,\, \alpha=2i, \beta=2j\\
\frac{1}{\xi_{2i+1}^k}\,\,\,\,\,\,\,\,\,\,\,\textnormal{if}\,\,\, \alpha= 2i, \beta= 2i+1
\end{cases}
\end{equation}
We note that the local equation of $M$ on $U_{2i}$ is $\xi_{2i}^k- \phi_i=0$, where $\phi_i=0$ is the equation of $C$ in $V_i$, and the local equation of $M$ in $U_{2i+1}$ is $1- \xi_{2i+1}^k \phi_i=0$.
\end{remark}

\subsection{Foliated double covering of a foliated projective plane $\left(\mathbb{P}_\mathbb{C}^2, T_0  \right)$ of degree $2$ with $T_0 \in H^0\left(\mathbb{P}_\mathbb{C}^2, \Theta_{\mathbb{P}_\mathbb{C}^2}(1) \right)$}\label{cs5} \

We keep the notations in Appendix \ref{p1}. Consider a foliation $\Theta_{\mathcal{F}_0}$ on $\mathbb{P}_\mathbb{C}^2$ of degree $2$ defined by $T_0\in H^0\left( \mathbb{P}_\mathbb{C}^2, \Theta_{\mathbb{P}_\mathbb{C}^2}(1) \right)$
\begin{align*}
T_0=\xi_0\xi_2\frac{\partial}{\partial \xi_1}+ \xi_0 \xi_1\frac{\partial}{\partial \xi_2} - 2\xi_1\xi_2\frac{\partial}{\partial \xi_0}
\end{align*}
Then $\xi_0^2+ \xi_1^2+ \xi_2^2=0$ is a non-degenerate invariant complex submanifold $C$ of $\mathbb{P}_\mathbb{C}^2$, so that $\left(\mathcal{O}_{\mathbb{P}_\mathbb{C}^2}(2), \nabla_{\mathcal{F}_0}\right)$ is a foliated line bundle where $\nabla_{\mathcal{F}_0}$ is the foliated connection defined by $T_0$ and $\xi_0^2+\xi_1^2+ \xi_2^2=0$.  We will construct a foliated double covering of $\left( \mathbb{P}_\mathbb{C}^2, T \right)$ which has the branch locus $C$. By Corollary \ref{cs1} it is sufficient to find the foliated line bundle $\left( \mathcal{O}_{\mathbb{P}_\mathbb{C}^2}(1), \nabla_{\mathcal{O}_{\mathbb{P}_\mathbb{C}^2}(1)} \right)$ such that $\left( \mathcal{O}_{\mathbb{P}_\mathbb{C}^2}(1) , \nabla_{\mathcal{O}_{\mathbb{P}_\mathbb{C}^2}(1)}\right)^{\otimes 2} =\left( \mathcal{O}_{\mathbb{P}_\mathbb{C}^2}(2), \nabla_{\mathcal{O}_{\mathbb{P}_\mathbb{C}^2}(2)} \right)= \left( \mathcal{O}_{\mathbb{P}_\mathbb{C}^2}(2), \nabla_{\mathcal{F}_0} \right)$.

 We will explicitly describe the foliated line bundle $\left(\mathcal{O}_{\mathbb{P}_\mathbb{C}^2}(2), \nabla_{\mathcal{F}_0} \right)$ on $\left( \mathbb{P}_\mathbb{C}^2, \Theta_{\mathcal{F}_0}  \right)$ determined by $\xi_0^2+\xi_1^2+ \xi_2^2=0$ in terms of the covering $\mathcal{U} = \left\{ U_0, U_1, U_2 \right\}$. We note that $T_0$ on $U_0$ is given by
\begin{align*}
T_0 = X_2\left( 1+ 2 X_1^2 \right)\frac{\partial}{\partial X_1} + X_1\left(1 + 2 X_2^2 \right) \frac{\partial}{\partial X_2} 
\end{align*}
and $C$ on $U_0$ is given by $1+ X_1^2 + X_2^2=0$. Then we have
\begin{align*}
\left[ X_2\left( 1+ 2 X_1^2 \right)\frac{\partial}{\partial X_1} + X_1\left(1 + 2 X_2^2 \right) \frac{\partial}{\partial X_2}  , 1+ X_1^2 + X_2^2 \right] = 4 X_1 X_2 \left( 1+ X_1^2 + X_2^2 \right)
\end{align*}
We note that $T_0$ on $U_1$ is given by
\begin{align*}
T_0 = - Y_2\left( Y_0^2 + 2  \right)\frac{\partial}{\partial Y_0} - Y_0\left(  Y_2^2 -1 \right)\frac{\partial}{\partial Y_2}
\end{align*}
and $C$ on $U_1$ is given by $Y_0^2+ Y_2^2 +1=0$. Then we have
\begin{align*}
\left[ - Y_2\left( Y_0^2 + 2  \right)\frac{\partial}{\partial Y_0} - Y_0\left(  Y_2^2 -1 \right)\frac{\partial}{\partial Y_2} \right] = - 2 Y_0 Y_2\left( Y_0^2 + Y_2^2 +1 \right)
\end{align*}
We note that $T_0$ on $U_2$ is given by 
\begin{align*}
T_0 = - Z_1\left( 2 + Z_0^2 \right) \frac{\partial}{\partial Z_0} + Z_0 \left( 1 - Z_1^2\right)\frac{\partial}{\partial Z_1}
\end{align*}
and $C$ on $U_2$ is given by $Z_0^2 + Z_1^2 +1 =0$. Then we have
\begin{align*}
\left[ - Z_1\left( 2 + Z_0^2 \right) \frac{\partial}{\partial Z_0} + Z_0 \left( 1 - Z_1^2\right)\frac{\partial}{\partial Z_1}, Z_0^2 + Z_1^2 +1 \right] = - 2 Z_0Z_1\left( 1+ Z_0^2 + Z_1^2 \right)
\end{align*}
Then the foliated line bundle $\left( \mathcal{O}_{\mathbb{P}_\mathbb{C}^2}(2), \nabla_{\mathcal{F}_0} \right)$ is given by 
\begin{align} \label{css17}
\left(\left\{ \phi_0, \phi_1, \phi_2 \right\}, \left\{ g_{01}= X_1^2, g_{02}= X_2^2, g_{12}= Y_2^2  \right\} \right) &\in C^0\left( \mathcal{U}, \mathscr{H}om_{\mathcal{O}_{\mathbb{P}_\mathbb{C}^2}} \left( \mathcal{O}_{\mathbb{P}_\mathbb{C}^2}(-1), \mathcal{O}_{\mathbb{P}_\mathbb{C}^2} \right) \right)\bigoplus C^1\left( \mathcal{U}, \mathcal{O}_{\mathbb{P}_\mathbb{C}^2}^*  \right) \\
 &  =
 C^0\left( \mathcal{U},  \mathcal{O}_{\mathbb{P}_\mathbb{C}^2}(1) \right)\bigoplus C^1\left( \mathcal{U}, \mathcal{O}_{\mathbb{P}_\mathbb{C}^2}^* \right) \notag
\end{align}
where $\phi_0\left( T_0 \right)= 4 X_1X_2$ on $U_0$, $\phi_1 \left( T_0 \right)=- 2 Y_0 Y_2$ on $U_1$ and $\phi_2\left( T_0 \right)= - 2 Z_0Z_1$ on $U_2$. Hence if we take the foliated line bundle $\left( \mathcal{O}_{\mathbb{P}_\mathbb{C}^2}(1), \nabla_{\mathcal{O}_{\mathbb{P}_\mathbb{C}^2}(1)} \right)$ defined by
\begin{align*}
\left( \left\{ \phi_0', \phi_1' , \phi_2' \right\}, \left\{ f_{01}= X_1, f_{02}= X_2, f_{12}= Y_2  \right\} \right) \in C^0\left( \mathcal{U}, \mathcal{O}_{\mathbb{P}_\mathbb{C}^2}(1) \right) \bigoplus C^1\left( \mathcal{U}, \mathcal{O}_{\mathbb{P}_\mathbb{C}^2}^* \right)
\end{align*}
where $\phi_0'\left( T_0 \right)= 2 X_1X_2, \phi_1'\left( T_0  \right)= - Y_0Y_2$ and $\phi_2'\left( T_0  \right)=-  Z_0Z_1$, we have $\left( \mathcal{O}_{\mathbb{P}_\mathbb{C}^2}(1), \nabla_{\mathcal{O}_{\mathbb{P}_\mathbb{C}^2(1)}} \right)^{\otimes 2} = \left( \mathcal{O}_{\mathbb{P}_\mathbb{C}^2}, \nabla_{\mathcal{F}_0} \right)$.

Next we describe the associated foliated complex manifold $\left(Y=\mathcal{O}_{\mathbb{P}_\mathbb{C}^2}(1)^\sharp , \Theta_{\tilde{\mathcal{F}}_0^\sharp}\right) $ which is a $\mathbb{P}_\mathbb{C}^1$-bundle on $\mathbb{P}_\mathbb{C}^2$ induced from $\left(\mathcal{O}_{\mathbb{P}_\mathbb{C}^2}(1), \nabla_{\mathcal{O}_{\mathbb{P}_\mathbb{C}^2(1)}} \right)$ in terms of the covering $\mathcal{U}^\sharp= \left\{U_{00}, U_{01}, U_{10}, U_{11}, U_{20}, U_{21} \right\}$ of $Y=\mathcal{O}_{\mathbb{P}_\mathbb{C}^2}(1)^\sharp$ where $U_{i0}= U_i\times \mathbb{C}$ and $U_{i1}= U_i\times \mathbb{C}$ and $\left(x, t \right) \in U_{i0}$ with $x\in U_i$ is identified with $\left(x, t' \right)$ on $U_{i1}$ if $t'=\frac{1}{t}$ on $U_i \times \mathbb{C}^*$ for $i=0,1,2$ as in the previous subsection. The foliation $\Theta_{\tilde{\mathcal{F}}_0^\sharp}$ on $Y$ is given by $T\in H^0\left( Y, \Theta_Y \otimes \Theta_{\tilde{\mathcal{F}}_0^\sharp}^* \right)$, where $\Theta_{\tilde{\mathcal{F}}_0^\sharp}^*$ has the transition functions in terms of the covering $\mathcal{U}^\sharp$ from $(\ref{cs2})$
\begin{equation}\label{cs15}
\begin{cases}
\textnormal{transition functions of}\,\,\,\,\, \mathcal{O}_{\mathbb{P}_\mathbb{C}^2}(1) & \textnormal{between $U_{00}$ and $U_{10}$ and  $U_{20}$}\\
1 & \textnormal{between $U_{00}$ and $U_{01}$ }\\
1 & \textnormal{between $U_{10}$ and $U_{11}$ }\\
1 & \textnormal{between $U_{20}$ and $U_{21}$ }
\end{cases}
\end{equation}
We note that $T$ on $U_{00}$ is given by
\begin{align*}
T = X_2 \left(1+2X_1^2\right)\frac{\partial}{\partial X_1} + X_1\left(1+2X_2^2\right)\frac{\partial}{\partial X_2} +2 tX_1X_2\frac{\partial}{\partial t}
\end{align*}
and $T$ on $U_{01}$ is given by
\begin{align*}
T = X_2 \left(1+2X_1^2\right)\frac{\partial}{\partial X_1} + X_1\left(1+2X_2^2\right)\frac{\partial}{\partial X_2} -2 t'X_1X_2\frac{\partial}{\partial t'}
\end{align*}
We note that $T$ on $U_{10}$ is given by
\begin{align*}
 - Y_2\left( Y_0^2 + 2  \right)\frac{\partial}{\partial Y_0} - Y_0\left(  Y_2^2 -1 \right)\frac{\partial}{\partial Y_2} - t Y_0Y_2\frac{\partial}{\partial t}
\end{align*}
and $T$ on $U_{11}$ is given by
\begin{align*}
 - Y_2\left( Y_0^2 + 2  \right)\frac{\partial}{\partial Y_0} - Y_0\left(  Y_2^2 -1 \right)\frac{\partial}{\partial Y_2} + t' Y_0Y_2\frac{\partial}{\partial t'}
\end{align*}
We note that $T$ on $U_{20}$ is given by
\begin{align*}
- Z_1\left( 2 + Z_0^2 \right) \frac{\partial}{\partial Z_0} + Z_0 \left( 1 - Z_1^2\right)\frac{\partial}{\partial Z_1} - tZ_0 Z_1\frac{\partial}{\partial t}
\end{align*}
and $T$ on $U_{21}$ is given by
\begin{align*}
- Z_1\left( 2 + Z_0^2 \right) \frac{\partial}{\partial Z_0} + Z_0 \left( 1 - Z_1^2\right)\frac{\partial}{\partial Z_1} + t' Z_0 Z_1\frac{\partial}{\partial t'}
\end{align*}
The the foliated double covering $X$ is the non-degenerate invariant complex submanifold of $\left( Y, \Theta_{\tilde{\mathcal{F}}_0^\sharp} \right)$ defined by the equation $1+ X_1^2 + X_2^2- t^2=0$ on $U_{00}$, $t'^2\left( 1+ X_1^2 + X_2^2  \right)-1=0$ on $U_{01}$, $1+ Y_0^2+ Y_2^2 -t^2=0$ on $U_{10}$, $t'^2\left( 1+ Y_0^2 + Y_2^2 \right)-1=0$, $1+ Z_0^2+ Z_1^2 - t^2=0$ on $U_{20}$ and $t'^2\left( 1 + Z_0^2 + Z_1^2 \right)-1 =0$ on $U_{21}$. Let us denote the foliation on $X$ induced from $\left( Y, \Theta_{\tilde{\mathcal{F}}_0^\sharp} \right)$ by $\Theta_{\mathcal{F}_{X_0}}$. 

In the following subsections, we will study embedded deformations of $C$ in $\left( \mathbb{P}_\mathbb{C}^2, T_0 \right)$, and foliated deformations of $\left( Y, \Theta_{\tilde{\mathcal{F}}_0^\sharp} \right)$ in terms of tangent sheaves, and foliated deformations of $\left(X, \Theta_{\mathcal{F}_{X_0}} \right)$, and embedded foliated deformations of $\left( X, \Theta_{\mathcal{F}_{X_0}} \right)$ in $\left( Y, \Theta_{\tilde{\mathcal{F}}_0^\sharp} \right)$. We will show that they are all unobstructed in foliated deformations. By assuming the knowledge of simultaneous deformations of foliated complex manifolds and their non-degenerate invariant complex manifolds defined by locally free subsheaves of tangent sheaves, we will study simultaneous deformations of $C\subset \left(\mathbb{P}_\mathbb{C}^2, T \right)$, and simultaneous foliated deformations of $X \subset \left( Y, \Theta_{\tilde{\mathcal{F}}_0^\sharp} \right)$. We also show that the associated functors of Artin rings $\textnormal{FDef}_{C\subset \left( \mathbb{P}_\mathbb{C}^2, T_0 \right)}$ and $\textnormal{FDef}_{X\subset \left( Y, \Theta_{\tilde{\mathcal{F}}_0^\sharp} \right)}$ (see Part VI) are all unobstructed.

\subsection{Embedded deformations of the invariant submanifold $C$ in $\left( \mathbb{P}_\mathbb{C}^2, T_0 \right)$} \ 

We keep the notations in subsection \ref{cs5}. We study deformations of the invariant hypersurface $C$ of $\left( \mathbb{P}_\mathbb{C}^2, T_0 \right)$ (see section \ref{hs1}). We recall the foliated line bundle $\left( \mathcal{O}_{\mathbb{P}_{\mathbb{C}}^2}(2), \nabla_{\mathcal{F}_0} \right)$. Then we have the associated complex of sheaves (see Part V)
\begin{align*}
\mathcal{O}_{\mathbb{P}_\mathbb{C}^2}(2)^\bullet: \mathcal{O}_{\mathbb{P}_\mathbb{C}^2}(2)\xrightarrow{\nabla_{\mathcal{F}_0}} \mathscr{H}om_{\mathcal{O}_{\mathbb{P}_\mathbb{C}^2}}\left( \mathcal{O}_{\mathbb{P}_\mathbb{C}^2}(-1),  \mathcal{O}_{\mathbb{P}_\mathbb{C}^2}(2)\right)\cong \mathcal{O}_{\mathbb{P}_\mathbb{C}^2}(3) \to 0
\end{align*}
By restricting to $C$, we have from $(\ref{cs8})$
\begin{align*}
\mathcal{O}_C(2)^\bullet: \mathcal{O}_C(2)\xrightarrow{\nabla_{\mathcal{F}_0}|_C} \mathcal{O}_C(3) \to 0
\end{align*}
From the exact sequence $0 \to \mathcal{O}_{\mathbb{P}_\mathbb{C}^2}\to \mathcal{O}_{\mathbb{P}_\mathbb{C}^2}(2)\to \mathcal{O}_C(2) \to 0$ and $0\to \mathcal{O}_{\mathbb{P}_\mathbb{C}^2}(1)\to \mathcal{O}_{\mathbb{P}_\mathbb{C}^2}(3)\to \mathcal{O}_C(3) \to 0$, we see that $H^1\left( C, \mathcal{O}_C(2) \right)= H^1\left( C, \mathcal{O}_C (3)\right)=0$ and $H^0\left(C, \mathcal{O}_C \right)$ is given by 
\begin{align}\label{css21}
H^0\left(C, \mathcal{O}_C(2) \right)&\cong \left\{ a_1 + a_2 X_1 + a_3 X_2 + a_4 X_1^2 + a_5 X_1X_2+ a_6 X_2^2 | a_1,....,a_6\in \mathbb{C}\,\,\,\textnormal{and}\,\,\, 1+ X_1^2 + X_2^2=0 \right\} \\
& \cong \left\{a_2 X_1 + a_3 X_2 + a_4 X_1^2+ a_5 X_1 X_2 + a_6 X_2^2 |a_2,...,a_6\in \mathbb{C} \right\} \notag
\end{align}
and $H^0\left( C, \mathcal{O}_C(3) \right)$ is given by
\begin{align}\label{css23}
H^0\left( C, \mathcal{O}_C(3) \right) \cong &\left\{ b_0 + b_1 X_1 + b_2 X_2 + b_3 X_1^2+ b_4 X_1X_2 + b_5 X_2^2 + b_6 X_1^3+ b_7 X_1^2 X_2 + b_8 X_1 X_2^2+  b_9 X_2^3 \right\} \\
 &\textnormal{with}\,\,\, b_0,..., b_9\in \mathbb{C}, \,\,\, \textnormal{and}\,\,\, 1+ X_1^2 + X_2^2=0, X_1 + X_1^3+ X_2^2 X_1=0, X_2 + X_1^2 X_2 + X_2^3=0  \notag \\
  \cong & \left\{b_1X_1 + b_2 X_2+ b_3 X_1^2+ b_4 X_1X_2 + b_5 X_2^2 +  b_7 X_1^2 X_2 + b_8 X_1 X_2^2 | b_3,..., b_9\in \mathbb{C} \right\} \notag
\end{align}
We note that $\mathbb{H}^0\left(C, \mathcal{O}_C(2)^\bullet \right)= \textnormal{ker}\left( H^0\left( C, \mathcal{O}_C(2)  \right) \to H^0\left( C, \mathcal{O}_C(3) \right) \right)$ and we have
{\Small{\begin{align*}
&\left[ X_2 \left(1+2X_1^2\right)\frac{\partial}{\partial X_1} + X_1\left(1+2X_2^2\right)\frac{\partial}{\partial X_2} ,a_2 X_1 + a_3 X_2 + a_4 X_1^2 + a_5 X_1X_2+ a_6 X_2^2\right] - 4 X_1 X_2\left(a_2 X_1 + a_3 X_2 + a_4 X_1^2 + a_5 X_1X_2+ a_6 X_2^2 \right) \\
&= a_3X_1+ a_2X_2 +a_5X_1^2 + a_5X_2^2 + \left(2a_4+2a_6\right)X_1X_2 -2a_2X_1^2X_2 -2 a_3X_1X_2^2 =0\,\,\,\,\,\textnormal{on}\,\,\, C\iff  a_2=a_3=a_5=0, a_4+a_6=0
\end{align*}}}
Hence we have
\begin{align} \label{css16}
\mathbb{H}^0\left( C, \mathcal{O}_C(2)^\bullet\right) \cong \left\langle X_1^2 - X_2^2 \right\rangle
\end{align}
Then we note that $\xi_0^2+ (1- c)\xi_1^2+ (1+c)\xi_2^2=0$ defines an invariant submanifold of $\left(\mathbb{P}_\mathbb{C}^2, T_0\right)$ for sufficiently small $c\in B\subset \mathbb{C}$. Then the foliated analytic family of invariant hypersurfaces defined $\xi_0^2+ (1- c)\xi_1^2 + (1+c)\xi_2^2=0$ in $\left( \mathbb{P}_\mathbb{C}^2\times B, T_0 \right)$ has the bijective characteristic map $\rho_{d,0}$ at $c=0$ (see (\ref{cs9})). Hence the family of invariant submanifolds in $\left( \mathbb{P}_\mathbb{C}^2, T_0 \right)$ is complete at $c=0$ (see Part VI).

\subsection{Foliated deformations of $\left( Y, \Theta_{\tilde{\mathcal{F}}_0^\sharp} \right)$ in terms of tangent sheaves}\

We keep the notations in subsection \ref{cs5}. We study foliated deformations of $\left( Y, \Theta_{\tilde{\mathcal{F}}_0^\sharp} \right)$ in terms of tangent sheaves. We have the leaf complex
\begin{align*}
\Theta_{\tilde{\mathcal{F}}_0^\sharp}^\bullet : \Theta_Y \to \mathscr{H}om_{\mathcal{O}_Y}\left( \Theta_{\tilde{\mathcal{F}}_0^\sharp} , \frac{\Theta_Y}{\Theta_{\tilde{\mathcal{F}}_0^\sharp } }\right) \cong \Theta_{\tilde{\mathcal{F}}_0^\sharp}^* \otimes \frac{\Theta_Y}{\Theta_{\tilde{\mathcal{F}}_0^\sharp} } \to 0
\end{align*}
We will compute $\mathbb{H}^i\left(Y, \Theta_{\tilde{\mathcal{F}}_0^\sharp}^\bullet \right)$ for $i=0,1,2$. First we compute $H^i\left(Y, \Theta_Y \right)$ for $i=0,1,2$. We note that we have the exact sequence from \cite{Wav68} p.942
\begin{align}\label{cs10}
0\to \mathcal{B} \to \Theta_Y \to \mathcal{A}\to 0
\end{align}
where the sheaf $\mathcal{A}$ is defined by the transition maps in terms of the covering $\mathcal{U}^\sharp$
\begin{equation}\label{cs13}
\begin{cases}
\textnormal{transition maps of}\,\,\,\,\, \Theta_{\mathbb{P}_\mathbb{C}^2} & \textnormal{between $U_{00}$ and $U_{10}$ and  $U_{20}$}\\
\textnormal{Id} & \textnormal{between $U_{00}$ and $U_{01}$ }\\
\textnormal{Id} & \textnormal{between $U_{10}$ and $U_{11}$ }\\
\textnormal{Id} & \textnormal{between $U_{20}$ and $U_{21}$ }
\end{cases}
\end{equation}
and the sheaf $\mathcal{B}$ is defined by the transition functions in terms of the covering $\mathcal{U}^\sharp$
\begin{equation} \label{cs14}
\begin{cases}
\textnormal{transition functions of}\,\,\,\,\, \mathcal{O}_{\mathbb{P}_\mathbb{C}^2}(1) & \textnormal{between $U_{00}$ and $U_{10}$ and  $U_{20}$}\\
-t^2 & \textnormal{between $U_{00}$ and $U_{01}$ }\\
-t^2 & \textnormal{between $U_{10}$ and $U_{11}$ }\\
-t^2 & \textnormal{between $U_{20}$ and $U_{21}$ }
\end{cases}
\end{equation}
Then by Proposition \ref{cs11} we have $H^i\left( Y, \mathcal{A}  \right) \cong H^i\left(\mathbb{P}_\mathbb{C}^2 , \Theta_{\mathbb{P}_\mathbb{C}^2} \right)$ for $i=0,1,2$, so that we have
\begin{align*}
 \dim_\mathbb{C}H^0\left(Y, \mathcal{A} \right)=8\,\,\,\,\,\,\, \textnormal{and}\,\,\,\,\,\,\,\, H^1\left( Y, \mathcal{A}\right)=H^2\left( Y, \mathcal{A} \right)=0
 \end{align*}
 
 On the other hand, by Proposition \ref{cs11} we have $H^i\left(Y, \mathcal{B} \right)\cong H^i\left(\mathbb{P}_\mathbb{C}^2, \mathcal{O}_{\mathbb{P}_\mathbb{C}^2} (1) \right)\bigoplus H^i\left(\mathbb{P}_\mathbb{C}^2, \mathcal{O}_{\mathbb{P}_\mathbb{C}^2} \right)\bigoplus H^i\left( \mathbb{P}_\mathbb{C}^2, \mathcal{O}_{\mathbb{P}_\mathbb{C}^2}(-1) \right)$, so that we have
 \begin{align*}
 \dim_\mathbb{C} H^0\left( Y, \mathcal{B} \right)=4,\,\,\,\,\,H^1\left( Y, \mathcal{B} \right)=H^2\left( Y, \mathcal{B} \right)=0.
 \end{align*}
 Then from the long exact sequence from $(\ref{cs10})$, we infer that
 \begin{align}\label{cs18}
\dim_\mathbb{C} H^0\left(Y, \Theta_Y \right) = 12\,\,\,\,\,\,\,\,\textnormal{and}\,\,\,\,\,\,\, H^1\left( Y, \Theta_Y \right)= H^2\left( Y, \Theta_Y \right)=0
 \end{align}
 Next we compute $H^i\left( Y, \Theta_{\tilde{\mathcal{F}}_0^\sharp}\otimes \frac{\Theta_Y}{\Theta_{\tilde{\mathcal{F}}_0^\sharp} } \right)$ for $i=0,1$. The exact sequence $0\to \Theta_{\tilde{\mathcal{F}}_0^\sharp}\to \Theta_Y \to \frac{\Theta_Y}{\Theta_{\tilde{\mathcal{F}}_0^\sharp}}\to 0$ induces the exact sequence
 \begin{align}\label{cs12}
 0 \to \mathcal{O}_Y \to \Theta_{\tilde{\mathcal{F}}_0^\sharp}^* \otimes  \Theta_Y \to \Theta_{\tilde{\mathcal{F}}_0^\sharp}^*\otimes \frac{\Theta_Y}{\Theta_{\tilde{\mathcal{F}}_0^\sharp} } \to 0
 \end{align}
Since $H^i\left( Y, \mathcal{O}_Y \right)\cong H^i\left( \mathbb{P}_\mathbb{C}^2, \mathcal{O}_{\mathbb{P}_\mathbb{C}^2} \right)$ by Proposition \ref{cs11}, we see that from the long exact sequence from $(\ref{cs12})$, we have
\begin{align}\label{cs17}
\dim_\mathbb{C} H^0\left( Y,  \Theta_{\tilde{\mathcal{F}}_0^\sharp}^* \otimes \frac{\Theta_Y}{\Theta_{\tilde{\mathcal{F}}_0^\sharp}}\right)= \dim_{\mathbb{C}} H^0\left( Y, \Theta_{\tilde{\mathcal{F}}_0^\sharp}^*\otimes \Theta_Y \right)-1 \,\,\,\,\,\,\,\, \textnormal{and} \,\,\,\,\,\,\,\, H^1\left(Y, \Theta_{\tilde{\mathcal{F}}_0^\sharp}^*\otimes \frac{\Theta_Y}{\Theta_{\tilde{\mathcal{F}}_0^\sharp} } \right)\cong H^1\left( Y, \Theta_{\tilde{\mathcal{F}}_0^\sharp}^* \otimes \Theta_Y \right)
\end{align}
So we compute $H^i\left( Y, \Theta_{\tilde{\mathcal{F}}_0^\sharp}^*\otimes \Theta_Y \right)$ for $i=0,1$. We note that from $(\ref{cs10})$ we have the exact sequence
\begin{align}\label{cs16}
0\to  \Theta_{\tilde{\mathcal{F}}_0^\sharp}^*\otimes \mathcal{B} \to \Theta_{\tilde{\mathcal{F}}_0^\sharp}^*\otimes \Theta_Y \to \Theta_{\tilde{\mathcal{F}}_0^\sharp}^*\otimes \mathcal{A}\to 0
\end{align}
and we note that from $(\ref{cs13})$ and $(\ref{cs15})$, $\Theta_{\tilde{\mathcal{F}}_0^\sharp}^*\otimes \mathcal{A}$ is defined by the transition maps
\begin{equation} \label{cs23}
\begin{cases}
\textnormal{transition maps of}\,\,\,\,\, \Theta_{\mathbb{P}_\mathbb{C}^2}(1) & \textnormal{between $U_{00}$ and $U_{10}$ and  $U_{20}$}\\
\textnormal{Id} & \textnormal{between $U_{00}$ and $U_{01}$ }\\
\textnormal{Id} & \textnormal{between $U_{10}$ and $U_{11}$ }\\
\textnormal{Id} & \textnormal{between $U_{20}$ and $U_{21}$ }
\end{cases}
\end{equation}
and from $(\ref{cs14})$ and $(\ref{cs15})$, $\Theta_{\tilde{\mathcal{F}}_0^\sharp}^*\otimes \mathcal{B}$ is defined by the transition functions
\begin{equation} \label{cs24}
\begin{cases}
\textnormal{transition functions of}\,\,\,\,\, \mathcal{O}_{\mathbb{P}_\mathbb{C}^2}(2) & \textnormal{between $U_{00}$ and $U_{10}$ and  $U_{20}$}\\
-t^2 & \textnormal{between $U_{00}$ and $U_{01}$ }\\
-t^2 & \textnormal{between $U_{10}$ and $U_{11}$ }\\
-t^2 & \textnormal{between $U_{20}$ and $U_{21}$ }
\end{cases}
\end{equation}
Then by Proposition \ref{cs11} we have $H^i\left(Y, \Theta_{\tilde{\mathcal{F}}_0^\sharp}^* \otimes \mathcal{A} \right)= H^i\left( \mathbb{P}_\mathbb{C}^2, \Theta_{\mathbb{P}_\mathbb{C}^2}(1) \right)$ for $i=0,1$, so that we have
\begin{align*}
\dim_\mathbb{C} H^0\left( Y, \Theta_{\tilde{\mathcal{F}}_0^\sharp}^* \otimes \mathcal{A} \right) = 15,\,\,\,\,\,\,\,\,\,\,\,\,H^1\left( Y, \Theta_{\tilde{\mathcal{F}}_0^\sharp}^* \otimes \mathcal{A} \right)=0
\end{align*}
On the other hand, we have $H^i\left( Y , \Theta_{\tilde{\mathcal{F}}_0^\sharp}^*\otimes \mathcal{B} \right)= H^i\left(\mathbb{P}_\mathbb{C}^2, \mathcal{O}_{\mathbb{P}_\mathbb{C}^2}(2) \right) \bigoplus H^i\left( \mathbb{P}_\mathbb{C}^2, \mathcal{O}_{\mathbb{P}_\mathbb{C}^2}(1) \right) \bigoplus H^i\left( \mathbb{P}_\mathbb{C}^2, \mathcal{O}_{\mathbb{P}_\mathbb{C}^2} \right)$ for $i=0,1$, so that we have
\begin{align*}
\dim_\mathbb{C} H^0\left( Y, \Theta_{\tilde{\mathcal{F}}_0^\sharp}^*\otimes \mathcal{B} \right) = 10,\,\,\,\,\,\,\,\,H^1\left( Y, \Theta_{\tilde{\mathcal{F}}_0^\sharp}^*\otimes \mathcal{B} \right)=0
\end{align*}
Then from the long exact sequence from $(\ref{cs16})$, we have
\begin{align}\label{cs19}
\dim_\mathbb{C} H^0\left( Y, \Theta_{\tilde{\mathcal{F}}_0^\sharp}^*\otimes \Theta_Y \right) = 25,\,\,\,\,\,\,\,\,\,\,H^1\left( Y, \Theta_{\tilde{\mathcal{F}}_0^\sharp}^* \otimes \Theta_Y \right)=0
\end{align}
Hence from $(\ref{cs17})$, we have
\begin{align*}
\dim_\mathbb{C} H^0\left( Y, \Theta_{\tilde{\mathcal{F}}_0^\sharp}^*\otimes \frac{\Theta_Y}{\Theta_{\tilde{\mathcal{F}}_0^\sharp} } \right)= 24,\,\,\,\,\,\,\,\,\,\,\,H^1\left( Y, \Theta_{\tilde{\mathcal{F}}_0^\sharp}^* \otimes \frac{\Theta_Y}{\Theta_{\tilde{\mathcal{F}}_0^\sharp} } \right)= 0
\end{align*}
Then we have
\begin{align*}
\mathbb{H}^0\left( Y, \Theta_{\tilde{\mathcal{F}}_0^\sharp}^\bullet \right)&= \textnormal{ker}\left( H^0\left( Y, \Theta_Y \right)\to H^0\left( Y, \Theta_{\tilde{\mathcal{F}}_0^\sharp}^* \otimes \frac{\Theta_Y}{\Theta_{\tilde{\mathcal{F}}_0^\sharp}} \right) \right)\\
\mathbb{H}^1\left( Y,  \Theta_{\tilde{\mathcal{F}}_0^\sharp}^\bullet  \right)&=\textnormal{coker}\left( H^0\left( Y, \Theta_Y \right)\to H^0\left( Y, \Theta_{\tilde{\mathcal{F}}_0^\sharp}^* \otimes \frac{\Theta_Y}{\Theta_{\tilde{\mathcal{F}}_0^\sharp}} \right)\right)\\
\mathbb{H}^2\left( Y,  \Theta_{\tilde{\mathcal{F}}_0^\sharp}^\bullet  \right)&=0
\end{align*}
In particular, $\left( Y, \Theta_{\tilde{\mathcal{F}}_0^\sharp} \right)$ is unobstructed in foliated deformations in terms of tangent sheaves.

Now we will compute $H^0\left(Y, \Theta_Y \right)\to H^0\left( Y, \Theta_{\tilde{\mathcal{F}}_0^\sharp}^* \otimes \frac{\Theta_Y}{\Theta_{\tilde{\mathcal{F}}_0^\sharp}} \right)$. We need to describe $H^0\left(Y, \Theta_Y \right)$ and $H^0\left( Y, \Theta_{\tilde{\mathcal{F}}_0^\sharp}^* \otimes \frac{\Theta_Y}{\Theta_{\tilde{\mathcal{F}}_0^\sharp} } \right)$ explicitly.

\subsubsection{Explicit description of $H^0\left( Y, \Theta_Y \right)$ on $U_{00}$}\

We recall that $\dim_\mathbb{C} H^0\left( Y, \Theta_Y \right)=12$ from $(\ref{cs18})$. We explicitly describe $H^0\left(Y, \Theta_Y \right)$ on $U_{00}$. Any element in $H^0\left( Y, \Theta_Y \right)$ on $U_{00}$ is of the form
\begin{align*}
V_{00}=A\left( X_1,X_2, t \right)\frac{\partial}{\partial X_1} + B\left(X_1, X_2, t \right)\frac{\partial}{\partial X_2} + C\left(X_1,X_2, t\right) \frac{\partial}{\partial t}
\end{align*}
Then $V_{00}$ on $U_{00}$ is transformed on $U_{01}$ into
\begin{align*}
V_{00}=A\left( X_1,X_2, \frac{1}{t'} \right)\frac{\partial}{\partial X_1} + B\left(X_1, X_2, \frac{1}{t'} \right)\frac{\partial}{\partial X_2} -t'^2 C\left(X_1,X_2,  \frac{1}{t'} \right) \frac{\partial}{\partial t'}
\end{align*}
which should be holomorphic, so that $V_{00}$ is written as
\begin{align*}
V_{01}= A(X_1,X_2)\frac{\partial}{\partial X_1} + B(X_1,X_2)\frac{\partial}{\partial X_2} + \left(C_1\left(X_1,X_2 \right) +  tC_2( X_1,X_2) + t^2 C_3(X_1,X_2) \right)\frac{\partial}{\partial t}
\end{align*}
where $A(X_1,X_2)= \sum_{i,j\geq 0} a_{ij} X_1^i X_2^j, a_{ij}\in \mathbb{C}, B(X_1,X_2)= \sum_{i,j\geq 0} b_{ij} X_1^1 X_2^j, b_{ij}\in \mathbb{C}, C_1(X_1,X_2)= \sum_{i,j\geq 0} c_{ij} X_1^i X_2^j, c_{ij} \in \mathbb{C}$, $C_2(X_1,X_2)= \sum_{i,j\geq 0} f_{ij} X_1^i X_2^j, f_{ij}\in \mathbb{C}$, and $C_3(X_1,X_2)= \sum_{i,j\geq 0} g_{ij} X_1^i X_2^j, g_{ij}\in \mathbb{C}$. We note that since $(Y_0,Y_2,t')=\left( \frac{1}{X_1}, \frac{X_2}{X_1}, \frac{1}{X_1}t \right)$ on $U_{00}\cap U_{10}$, we have
\begin{align}\label{css10}
\frac{\partial}{\partial X_1} = - \frac{1}{X_1^2}\frac{\partial}{\partial Y_0}- \frac{X_2}{X_1^2}\frac{\partial}{\partial Y_2}- \frac{t}{X_1^2}\frac{\partial}{\partial t'}= - Y_0^2\frac{\partial}{\partial Y_0}- Y_0Y_2\frac{\partial}{\partial Y_2}- t'Y_0\frac{\partial}{\partial t},\,\,\,\,\,\,\,\,\,\,\,\frac{\partial}{\partial X_2} = Y_0 \frac{\partial}{\partial Y_2},\,\,\,\,\,\,\,\,\,\,\,\frac{\partial}{\partial t} = Y_0\frac{\partial}{\partial t'}
\end{align}
we see that $V_{00}$ on $U_{00}$ is transformed on $U_{10}$ into
{\small{\begin{align*}
V_{10}=&A\left( \frac{1}{Y_0},  \frac{Y_2}{Y_0} \right)\left( - Y_0^2\frac{\partial}{\partial Y_0} - Y_0Y_2 \frac{\partial}{\partial Y_2}-t'Y_0\frac{\partial}{\partial t'}\right) + B\left( \frac{1}{Y_0}, \frac{Y_2}{Y_0} \right)Y_0\frac{\partial}{\partial Y_2} \\
& + \left( C_1\left( \frac{1}{Y_0}, \frac{Y_2}{Y_0}\right)+ \frac{t'}{Y_0} C_2\left( \frac{1}{Y_0} , \frac{Y_2}{Y_0} \right) + \frac{t'^2}{Y_0^2} C_3 \left( \frac{1}{Y_0}, \frac{Y_2}{Y_0}\right)\right)Y_0\frac{\partial}{\partial t'}\\
&= - A\left( \frac{1}{Y_0},  \frac{Y_2}{Y_0} \right)Y_0^2\frac{\partial}{\partial Y_0} +\left( - A\left( \frac{1}{Y_0},  \frac{Y_2}{Y_0} \right) Y_0Y_2 + B\left( \frac{1}{Y_0}, \frac{Y_2}{Y_0}\right)Y_0 \right)\frac{\partial}{\partial Y_2} \\
&+ \left( Y_0C_1\left( \frac{1}{Y_0}, \frac{Y_2}{Y_0}\right)+ t' \left( C_2\left( \frac{1}{Y_0} , \frac{Y_2}{Y_0} \right) - Y_0 A\left( \frac{1}{Y_0}, \frac{Y_2}{Y_0} \right) \right) + \frac{t'^2}{Y_0} C_3 \left( \frac{1}{Y_0}, \frac{Y_2}{Y_0}\right) \right)\frac{\partial}{\partial t'}
\end{align*}}}
which should be holomorphic. We note that in order for
\begin{align*}
- A\left( \frac{1}{Y_0},  \frac{Y_2}{Y_0} \right) Y_0Y_2 + B\left( \frac{1}{Y_0}, \frac{Y_2}{Y_0}\right)Y_0 = \sum_{i,j\geq 0} \frac{-a_{ij} Y_2^{j+1} + b_{ij} Y_2^j }{Y_0^{i+j-1}} 
\end{align*}
to be holomorphic, we only need to consider $(i,j)=(0,0),(0,1), (1,0)$ which clearly defines a holomorphic function and $(0,2),(1,1),(2,0)$ in which case we have
\begin{align*}
\frac{- a_{20} Y_2 + b_{20} - a_{11} Y_2^2 + b_{11} Y_2 - a_{02} Y_2^3 + b_{02} Y_2^2}{Y_0} \Longrightarrow a_{02}=0, b_{20}=0, b_{11}= a_{20} , b_{02}= a_{11}
\end{align*}
Then we have
\begin{align*}
A(X_1,X_2)&= a_{00} + a_{10} X_1 + a_{01} X_2 + a_{20} X_1^2 + a_{11} X_1 X_2\\
B(X_1,X_2)&= b_{00} + b_{10} X_1 + b_{01} X_2 + a_{20} X_1 X_2 + a_{11} X_2^2
\end{align*}
We note that in order for
\begin{align*}
C_2\left( \frac{1}{Y_0}, \frac{Y_2}{Y_0} \right) - Y_0 A\left( \frac{1}{Y_0}, \frac{Y_2}{Y_0}   \right) = C_2\left( \frac{1}{Y_0}, \frac{Y_2}{Y_0} \right) - Y_0\left( a_{00}+ \frac{ a_{10} + a_{01} Y_2}{Y_0}  + \frac{a_{20}+ a_{11} Y_2}{Y_0^2} \right)
\end{align*}
to be holomorphic, we should have
\begin{align*}
\frac{f_{10} + f_{01}Y_2}{Y_0}  - \frac{a_{20}  + a_{11} Y_2}{Y_0} = 0 \Longrightarrow f_{10}= a_{20}, f_{01}= a_{11}
\end{align*}
Then we have 
\begin{align*}
C_1(X_1,X_2)&= c_{00} + c_{10} X_1 + c_{01} X_2 \\
C_2(X_1,X_2)&= f_{00} + a_{20} X_1 + a_{11} X_2\\
C_3(X_1,X_2)& =0
\end{align*}
Hence $H^0\left( Y, \Theta_Y \right)$ is given by 
\begin{align}
H^0\left( Y, \Theta_Y \right) = &\left\{ \left( A_0 + A_1 X_1 + A_2 X_2 + A_3 X_1X_2 + A_4 X_1^2\right)\frac{\partial}{\partial X_1} +\left( B_0 + B_1 X_1 + B_2 X_2 + A_3 X_2^2  + A_4 X_1 X_2 \right)\frac{\partial}{\partial X_2} \right. \label{css13}\\
&\left. +\left( \left( C_0 + C_1 X_1 + C_2 X_2 \right) + t \left(C_3 + A_3 X_2 + A_4 X_1 \right)\right)\frac{\partial}{\partial t} | A_0,A_1,A_2,A_3,A_4, B_0,B_1,B_2, C_0,C_1,C_2,C_3 \in \mathbb{C}  \right\} \notag
\end{align}

\subsubsection{Explicit description of $H^0\left( Y, \Theta_{\tilde{\mathcal{F}}_0^\sharp}^*\otimes \frac{\Theta_Y}{\Theta_{\tilde{\mathcal{F}}_0^\sharp}} \right)$ on $U_{00}$}\

First we explicitly describe $H^0\left( Y, \Theta_{\tilde{\mathcal{F}}_0^\sharp}^*\otimes \Theta_Y \right)$ on $U_{00}$, and then by using the exact sequence induced from $(\ref{cs12})$
\begin{align}\label{css11}
0\to H^0\left( Y, \mathcal{O}_Y \right) &\to H^0\left(Y, \Theta_{\tilde{\mathcal{F}}_0^\sharp}^* \otimes \Theta_Y \right)\to H^0\left( Y, \Theta_{\tilde{\mathcal{F}}_0^\sharp}^*\otimes \frac{\Theta_Y}{\Theta_{\tilde{\mathcal{F}}_0^\sharp} } \right) \to 0 \\
       1 & \mapsto T = X_2 \left(1+2X_1^2\right)\frac{\partial}{\partial X_1} + X_1\left(1+2X_2^2\right)\frac{\partial}{\partial X_2} +2 tX_1X_2\frac{\partial}{\partial t} \notag
\end{align}
we will deduce the explicit basis of $H^0\left( Y, \Theta_{\tilde{\mathcal{F}}_0^\sharp}^* \otimes \frac{\Theta_Y}{\Theta_{\tilde{\mathcal{F}}_0^\sharp} } \right)$. We recall that $\dim_\mathbb{C} H^0\left( Y, \Theta_{\tilde{\mathcal{F}}_0^\sharp}^*\otimes \Theta_Y \right)=25$ from $(\ref{cs19})$. Let us explicitly describe $H^0\left( Y, \Theta_{\tilde{\mathcal{F}}_0^\sharp}^* \otimes \Theta_Y \right)$ on $U_{00}$. In view of $(\ref{cs2})$ any element in $H^0\left( Y, \Theta_{\tilde{\mathcal{F}}_0^\sharp}^*\otimes \Theta_Y \right)$ is of the form
\begin{align*}
V_{00}= A(X_1,X_2)\frac{\partial}{\partial X_1}+ B(X_1,X_2)\frac{\partial}{\partial X_2} + \left(C_1(X_1,X_2, t)+ t C_2(X_1,X_2) + t^2 C_3(X_1,X_2) \right) \frac{\partial}{\partial t}
\end{align*}
where $A(X_1,X_2)= \sum_{i,j\geq 0} a_{ij} X_1^i X_2^j, a_{ij}\in \mathbb{C}, B(X_1,X_2)= \sum_{i,j\geq 0} b_{ij} X_1^1 X_2^j, b_{ij}\in \mathbb{C}, C_1(X_1,X_2)= \sum_{i,j\geq 0} c_{ij} X_1^i X_2^j, c_{ij} \in \mathbb{C}$, $C_2(X_1,X_2)= \sum_{i,j\geq 0} f_{ij} X_1^i X_2^j, f_{ij}\in \mathbb{C}$, and $C_3(X_1,X_2)= \sum_{i,j\geq 0} g_{ij} X_1^i X_2^j, g_{ij}\in \mathbb{C}$. We note that in view of $(\ref{cs15})$ and $(\ref{css10})$, we see that $V_{00}$ on $U_{00}$ is transformed on $U_{10}$ into
{\small{\begin{align*}
V_{10}=&Y_0\left(A\left( \frac{1}{Y_0},  \frac{Y_2}{Y_0} \right)\left( - Y_0^2\frac{\partial}{\partial Y_0} - Y_0Y_2 \frac{\partial}{\partial Y_2}-t'Y_0\frac{\partial}{\partial t'}\right) + B\left( \frac{1}{Y_0}, \frac{Y_2}{Y_0} \right)Y_0\frac{\partial}{\partial Y_2}\right)  \\
& + Y_0\left( \left( C_1\left( \frac{1}{Y_0}, \frac{Y_2}{Y_0}\right)+ \frac{t'}{Y_0} C_2\left( \frac{1}{Y_0} , \frac{Y_2}{Y_0} \right) + \frac{t'^2}{Y_0^2} C_3 \left( \frac{1}{Y_0}, \frac{Y_2}{Y_0}\right)\right)Y_0\frac{\partial}{\partial t'}\right)\\
&= - A\left( \frac{1}{Y_0},  \frac{Y_2}{Y_0} \right)Y_0^3\frac{\partial}{\partial Y_0} +\left( - A\left( \frac{1}{Y_0},  \frac{Y_2}{Y_0} \right) Y_0^2Y_2 + B\left( \frac{1}{Y_0}, \frac{Y_2}{Y_0}\right)Y_0^2 \right)\frac{\partial}{\partial Y_2} \\
&+ \left(  Y_0^2 C_1\left( \frac{1}{Y_0}, \frac{Y_2}{Y_0}\right)+ t'  \left(Y_0 C_2\left( \frac{1}{Y_0} , \frac{Y_2}{Y_0} \right) - A\left( \frac{1}{Y_0}, \frac{Y_2}{Y_0} \right)Y_0^2 \right)+ t'^2 C_3 \left( \frac{1}{Y_0}, \frac{Y_2}{Y_0}\right)\right)\frac{\partial}{\partial t'}
\end{align*}}}
which should be holomorphic. We note that in order for
\begin{align*}
- A\left( \frac{1}{Y_0},  \frac{Y_2}{Y_0} \right) Y_0^2Y_2 + B\left( \frac{1}{Y_0}, \frac{Y_2}{Y_0}\right)Y_0^2= \sum_{i,j\geq 0} \frac{  b_{ij} Y_2^j - a_{ij} Y_2^{j+1}}{Y_0^{i+j-2}}
\end{align*}
to be holomorphic, we only need to consider $(i,j)=(0,0),(0,1),(1,0),(0,2),(1,1),(2,0)$ which clearly defines a holomorphic function and $(0,3),(1,2),(2,1),(3,0)$ in which case we have
\begin{align*}
\frac{b_{30}- a_{30} Y_2 + b_{21}Y_2- a_{21}Y_2^2 + b_{12} Y_2^2 - a_{12} Y_2^3 + b_{03} Y_2^3 - a_{03} Y_2^4}{Y_0} \Longrightarrow b_{30}=0, b_{21}= a_{30}, b_{12}= a_{21},b_{03}= a_{12}, a_{03}=0
\end{align*}
Then we have
\begin{align*}
A(X_1,X_2)&= a_{00} + a_{10} X_1 + a_{01} X_2 + a_{20} X_1^2 + a_{11} X_1X_2 + a_{02} X_2^2 + a_{30} X_1^3 + a_{21} X_1^2 X_2+ a_{12} X_1 X_2^2  \\
B(X_1,X_2)&= b_{00} + b_{10} X_1 + b_{01} X_2 + b_{20} X_1^2 + b_{11} X_1 X_2 + b_{02} X_2^2 + a_{30}X_1^2 X_2 + a_{21} X_1 X_2 ^2 + a_{12} X_2^3
\end{align*}
We note that in order for
{\small{\begin{align*}
C_2\left( \frac{1}{Y_0}, \frac{Y_2}{Y_0} \right)Y_0 - A\left( \frac{1}{Y_0}, \frac{Y_2}{Y_0} \right) Y_0^2 = C_2\left( \frac{1}{Y_0}, \frac{Y_2}{Y_0}  \right)- Y_0^2\left( a_{00}+\frac{a_{10} + a_{01}Y_2 }{Y_0} + \frac{a_{20}+ a_{11} Y_2+ a_{02} Y_2^2}{Y_0^2} + \frac{a_{30}+ a_{21} Y_2+ a_{12} Y_2^2}{Y_0^3} \right)
\end{align*}}}
to be holomorphic, we should have
\begin{align*}
 \frac{ f_{20}+ f_{11} Y_2 + f_{02} Y_2^2 }{Y_0}  - \frac{a_{30} + a_{21} Y_2 + a_{12} Y_2^2}{Y_0} =0 \Longrightarrow f_{20}=a_{30}, f_{11}= a_{21} , f_{02}= a_{12}
\end{align*}
Then we have
\begin{align*}
C_1(X_1,X_2)&= c_{00} + c_{10}X_1 + c_{01} X_2 + c_{20}X_1^2 + c_{11} X_1 X_2 + c_{02} X_2^2 \\
C_2(X_1,X_2)&= f_{00} + f_{10} X_1 + f_{01} X_2 + a_{30}X_1^2 + a_{21}X_1X_2 + a_{12} X_2^2 \\
C_3(X_1,X_2) &= g_{00}
\end{align*}
Hence $H^0\left( Y, \Theta_{\tilde{\mathcal{F}}_0^\sharp}^* \otimes \Theta_Y \right)$ on $U_{00}$ is given by
{\small{\begin{align*}
H^0\left( Y, \Theta_{\tilde{\mathcal{F}}_0^\sharp}^* \otimes \Theta_Y \right) =\left\{ \left(a_0 + a_1 X_1 + a_2 X_2 + a_3 X_1^2 + a_4 X_1X_2 + a_5 X_2^2 + a_6 X_1^3 + a_7 X_1^2 X_2+ a_8 X_1 X_2^2 \right)\frac{\partial}{\partial X_1}\right. \\
\left. + \left( a_{9} + a_{10}X_1 + a_{11}X_2 + a_{12}X_1^2 + a_{13}X_1X_2 + a_{14}X_2^2 + a_6 X_1^2 X_2+ a_7 X_1 X_2^2 + a_8 X_2^3 \right)\frac{\partial}{\partial X_2}\right. \\
\left. + \left(  a_{15} + a_{16}X_1 + a_{17}X_2 + a_{18}X_1^2 + a_{19}X_1X_2 + a_{20}X_2^2 + \left( a_{21} + a_{22}X_1+ a_{23}X_2  + a_6 X_1^2+ a_7 X_1X_2 + a_8 X_2^2\right)t  + a_{24}t^2   \right) \frac{\partial}{\partial t}\right\}
\end{align*}}}
where $a_0,a_1,..., a_{24}\in \mathbb{C}$. Then from $(\ref{css11})$, since we have
\begin{align*}
a_7X_1^2 X_2 \frac{\partial}{\partial X_1} + a_7X_1X_2^2 \frac{\partial}{\partial X_2} +  a_7 t X_1 X_2 \frac{\partial}{\partial t} \equiv - \frac{1}{2} a_7 X_2 \frac{\partial}{\partial X_1} - \frac{1}{2} a_7 X_1 \frac{\partial}{\partial X_2} \,\,\,\,\,\,\,\textnormal{in}\,\,\,\,\,\, H^0\left(Y, \Theta_{\tilde{\mathcal{F}}_0^\sharp}^* \otimes \frac{\Theta_Y}{\Theta_{\tilde{\mathcal{F}}_0^\sharp}} \right),
\end{align*}
we see that $H^0\left( Y, \Theta_{\tilde{\mathcal{F}}_0^\sharp}^* \otimes \frac{\Theta_Y}{\Theta_{\tilde{\mathcal{F}}_0^\sharp}}  \right)$ on $U_{00}$ is given by
{\small{\begin{align}
H^0\left( Y, \Theta_{\tilde{\mathcal{F}}_0^\sharp}^* \otimes \frac{\Theta_Y}{\Theta_{\tilde{\mathcal{F}}_0^\sharp}} \right) =\left\{ \left(a_0 + a_1 X_1 + a_2 X_2 + a_3 X_1^2 + a_4 X_1X_2 + a_5 X_2^2 + a_6 X_1^3 + a_8 X_1 X_2^2 \right)\frac{\partial}{\partial X_1}\right. \label{css12} \\
\left. + \left( a_{9} + a_{10}X_1 + a_{11}X_2 + a_{12}X_1^2 + a_{13}X_1X_2 + a_{14}X_2^2 + a_6 X_1^2 X_2 + a_8 X_2^3 \right)\frac{\partial}{\partial X_2}\right. \notag \\
\left. + \left(  a_{15} + a_{16}X_1 + a_{17}X_2 + a_{18}X_1^2 + a_{19}X_1X_2 + a_{20}X_2^2 + \left( a_{21} + a_{22}X_1+ a_{23}X_2  + a_6 X_1^2 + a_8 X_2^2\right)t  + a_{24}t^2   \right) \frac{\partial}{\partial t}\right\} \notag
\end{align}}}
where $a_0,...,a_6,a_8,...,a_{24}\in \mathbb{C}$.

\subsubsection{Computation of $H^0\left( Y, \Theta_Y \right)\to H^0\left(Y, \Theta_{\tilde{\mathcal{F}}_0^\sharp}^* \otimes \frac{\Theta_Y}{\Theta_{\tilde{\mathcal{F}}_0^\sharp}} \right)$}\

By using $(\ref{css13})$ and $(\ref{css12})$ we compute $H^0\left( Y, \Theta_Y \right)\to H^0\left(Y, \Theta_{\tilde{\mathcal{F}}_0^\sharp}^* \otimes \frac{\Theta_Y}{\Theta_{\tilde{\mathcal{F}}_0^\sharp}} \right)$. The image is given by
{\small{\begin{align*}
&\left[ \left( A_0 + A_1 X_1 + A_2 X_2 + A_3 X_1X_2 + A_4 X_1^2\right)\frac{\partial}{\partial X_1},  X_2\left(1+2X_1^2\right)\frac{\partial}{\partial X_1} + X_1\left(1+2X_2^2\right)\frac{\partial}{\partial X_2} +2 tX_1X_2\frac{\partial}{\partial t}\right]\\
&+\left[ \left( B_0 + B_1 X_1 + B_2 X_2 + A_3 X_2^2  + A_4 X_1 X_2 \right)\frac{\partial}{\partial X_2} , X_2\left(1+2X_1^2\right)\frac{\partial}{\partial X_1} + X_1\left(1+2X_2^2\right)\frac{\partial}{\partial X_2} +2 tX_1X_2\frac{\partial}{\partial t} \right] \\
&+ \left[ \left( \left( C_0 + C_1 X_1 + C_2 X_2 \right) + t \left(C_3 + A_3 X_2 + A_4 X_1 \right)\right)\frac{\partial}{\partial t}, X_2\left(1+2X_1^2\right)\frac{\partial}{\partial X_1} + X_1\left(1+2X_2^2\right)\frac{\partial}{\partial X_2} +2 tX_1X_2\frac{\partial}{\partial t}  \right] \\
&= \left(B_0+\left(B_1 - A_2\right)X_1 + \left(B_2 - A_1 \right)X_2 + \left(4A_0-  A_4\right)X_1X_2 + \left(2A_1+2 B_2\right)X_1^2 X_2+ 2A_2 X_1X_2^2 \right)\frac{\partial}{\partial X_1}\\
&+\left(\left(2B_0 - A_3 \right) X_1^2+ 2B_1 X_1^3+ 2 A_3 X_1^2 X_2^2+2 A_4 X_1^3 X_2\right)\frac{\partial}{\partial X_1}\\
&+\left(A_0 + \left(A_1- B_2\right)X_1 + \left(A_2- B_1\right) X_2 + \left(4B_0-  A_3\right)X_1X_2+ 2 B_1 X_1^2 X_2+  \left(2 B_2+2 A_1\right) X_1 X_2^2\right)\frac{\partial}{\partial X_2}\\
&+ \left(\left(2A_0- A_4\right) X_2^2  + 2A_2X_2^3+ 2A_3X_1X_2^3+ 2A_4 X_1^2 X_2^2\right)\frac{\partial}{\partial X_2}\\
&+ \left(- C_2 X_1 - C_1 X_2+2C_0X_1X_2\right)\frac{\partial}{\partial t} \\
&+t\left(\left(2B_0- A_3\right)X_1+ \left( 2A_0 - A_4 \right)X_2 + 2B_1 X_1^2 + \left(2A_1+2 B_2\right)X_1X_2+ 2A_2X_2^2+ 2A_3X_1X_2^2+ 2A_4 X_1^2X_2\right)\frac{\partial}{\partial t}
\end{align*}}}

We note that from $\left(X_2+ 2 X_1^2 X_2 \right)\frac{\partial}{\partial X_1} + \left( X_1 + 2 X_1 X_2^2 \right)\frac{\partial}{\partial X_2} + 2 t X_1 X_2 \frac{\partial}{\partial t}\equiv 0$, we have
\begin{align*}
  2A_4X_1^3 X_2 \frac{\partial}{\partial X_1} + 2  A_4 X_1^2 X_2^2 \frac{\partial}{\partial X_2} +  2A_4 X_1^2 X_2 \frac{\partial}{\partial t} &\equiv - A_4 X_1 X_2 \frac{\partial}{\partial X_1} - A_4 X_1^2 \frac{\partial}{\partial X_2} \\
   2A_3X_1^2 X_2^2 \frac{\partial}{\partial X_1}  + 2 A_3 X_1 X_2^3 \frac{\partial}{\partial X_2} +   2 A_3 t X_1 X_2^2 \frac{\partial}{\partial t} &\equiv - A_3 X_2^2 \frac{\partial}{\partial X_1} - A_3 X_1 X_2 \frac{\partial}{\partial X_2} \\
2 \left(A_1+ B_2\right)X_1^2 X_2 \frac{\partial}{\partial X_1} + 2 \left(A_1+ B_2\right)     X_1 X_2^2 \frac{\partial}{\partial X_2} + 2\left(A_1+ B_2\right) t  X_1X_2 \frac{\partial}{\partial t} &\equiv - \left(A_1 + B_2 \right)X_2 \frac{\partial}{\partial X_1} - \left( A_1 + B_2 \right) X_1 \frac{\partial}{\partial X_2}
\end{align*}
Then the image of $H^0\left( Y, \Theta_Y \right)\to H^0\left( Y , \Theta_{\tilde{\mathcal{F}}_0^\sharp}^* \otimes \frac{\Theta_Y}{\Theta_{\tilde{\mathcal{F}}_0^\sharp}} \right)$ is given by
\begin{align}\label{css30}
\left(B_0+\left(B_1 - A_2\right)X_1 - 2 A_1 X_2 + \left(4A_0-  2A_4\right)X_1X_2 + 2A_2 X_1X_2^2 + \left(2B_0 - A_3 \right) X_1^2 - A_3 X_2^2 + 2B_1 X_1^3\right)\frac{\partial}{\partial X_1}\\
+ \left(A_0 - 2 B_2 X_1 + \left(A_2- B_1\right) X_2 + \left(4B_0-  2  A_3 \right)X_1X_2+ 2 B_1 X_1^2 X_2 + \left(2A_0- A_4\right) X_2^2 -A_4 X_1^2 + 2A_2X_2^3 \right)\frac{\partial}{\partial X_2} \notag \\
+ \left(- C_2 X_1 - C_1 X_2+2C_0X_1X_2 + t\left(\left(2B_0- A_3\right)X_1+ \left( 2A_0 - A_4 \right)X_2 + 2B_1 X_1^2 + 2A_2X_2^2 \right) \right)\frac{\partial}{\partial t}  \notag
\end{align}
Then from $(\ref{css12})$ we have $\dim_\mathbb{C} \textnormal{ker}\left( H^0\left( Y, \Theta_Y \right)\to H^0\left( Y, \Theta_{\tilde{\mathcal{F}}_0^\sharp}^* \otimes \frac{\Theta_Y}{\Theta_{\tilde{\mathcal{F}}_0^\sharp}} \right) \right)= 1$, so that we have
\begin{align*}
\dim_\mathbb{C} \mathbb{H}^0\left( Y , \Theta_{\tilde{\mathcal{F}}_0^\sharp}^\bullet \right) &= 1\\
\dim_\mathbb{C} \mathbb{H}^1\left( Y , \Theta_{\tilde{\mathcal{F}}_0^\sharp}^\bullet \right) & = 24 - 11 = 13
\end{align*}

\subsection{Foliated deformations of $\left(X, \Theta_{\mathcal{F}_{X_0}} \right)$}\

We keep the notations in subsection \ref{cs5}. We study foliated deformations of $\left( X, \Theta_{\mathcal{F}_{X_0}} \right)$. We have the leaf complex
\begin{align*}
\Theta_{\mathcal{F}_{X_0}}^\bullet : \Theta_X \to \mathscr{H}om_{\mathcal{O}_X}\left( \Theta_{\mathcal{F}_{X_0}}, \frac{\Theta_X}{\Theta_{\mathcal{F}_{X_0}}} \right) \cong \Theta_{\mathcal{F}_{X_0}}^* \otimes \frac{\Theta_X}{\Theta_{\mathcal{F}_{X_0}}} \to 0
\end{align*}
We note that the double covering $X$ of $\mathbb{P}_\mathbb{C}^2$ is isomorphic to quadratic surface $\xi_0^2+ \xi_1^2 + \xi_2^2+ \xi_3^2=0$ in $\mathbb{P}_\mathbb{C}^3$ with the homogeneous coordinates $\left[\xi_0: \xi_1: \xi_2: \xi_3 \right]$, so that we have
\begin{align}\label{css3}
\dim_\mathbb{C} H^0\left( X, \Theta_X \right)=6,\,\,\,\,\,\,\,\,\,\,\,\,H^1\left( X, \Theta_X \right)= H^2\left(X, \Theta_X \right)=0
\end{align}
This implies that we have
\begin{align}
\mathbb{H}^0\left(X, \Theta_{\mathcal{F}_{X_0}}^\bullet \right)&\cong \textnormal{ker}\left( H^0\left( X, \Theta_X \right)\to H^0\left( X, \Theta_{\mathcal{F}_{X_0}}^* \otimes \frac{\Theta_X}{\Theta_{\mathcal{F}_{X_0}}} \right) \right)\\
\mathbb{H}^1\left( X, \Theta_{\mathcal{F}_{X_0}}^\bullet \right)& \cong  \textnormal{coker} \left( H^0\left( X, \Theta_X \right)\to H^0\left( X, \Theta_{\mathcal{F}_{X_0}}^* \otimes \frac{\Theta_X}{\Theta_{\mathcal{F}_{X_0}}} \right) \right) \\
\mathbb{H}^2\left(X, \Theta_{\mathcal{F}_{X_0}}^\bullet \right)& \cong H^1\left(X, \Theta_{\mathcal{F}_{X_0}}^* \otimes \frac{\Theta_X}{\Theta_{\mathcal{F}_{X_0}} } \right) \label{css1}
\end{align}
We will compute that $H^i\left( X, \Theta_{\mathcal{F}_{X_0}}^*\otimes \frac{\Theta_X}{\Theta_{\mathcal{F}_{X_0}}} \right)$ for $i=0,1$. The exact sequence $0\to \Theta_{\mathcal{F}_{X_0}}\to \Theta_X \to \frac{\Theta_X}{\Theta_{\mathcal{F}_{X_0}}}\to 0$ induces the exact sequence
\begin{align*}
0\to \mathcal{O}_X \to \Theta_{\mathcal{F}_{X_0}}^* \otimes \Theta_X \to \Theta_{\mathcal{F}_{X_0}}^*\otimes \frac{\Theta_X}{\Theta_{\mathcal{F}_{X_0}}}\to 0
\end{align*}
from which we see that we have
\begin{align}\label{cs20}
H^1\left( X, \Theta_{\mathcal{F}_{X_0}}^* \otimes \Theta_X   \right)\cong H^1\left( X, \Theta_{\mathcal{F}_{X_0}}^*\otimes \frac{\Theta_X}{\Theta_{\mathcal{F}_{X_0}}} \right)
\end{align}
and we have
\begin{align}
 0\to H^0\left(X, \mathcal{O}_X \right)\cong \mathbb{C} &\to  H^0\left(X, \Theta_{\mathcal{F}_{X_0}}^*\otimes \Theta_X \right)\to H^0\left(X, \Theta_{\mathcal{F}_{X_0}}^*\otimes \frac{\Theta_X}{\Theta_{\mathcal{F}_{X_0}}} \right)\to 0 \label{cs21}\\
 1&\to  T|_X=    X_2 \left(1+2X_1^2\right)\frac{\partial}{\partial X_1} + X_1\left(1+2X_2^2\right)\frac{\partial}{\partial X_2} +2 tX_1X_2\frac{\partial}{\partial t}|_{1+X_1^2 + X_2^2 - t^2=0 } \notag
\end{align}
First we will compute $H^i\left( X, \Theta_{\mathcal{F}_{X_0}}^* \otimes \Theta_X \right)$ for $i=0,1$ and then by using $(\ref{cs20})$ and $(\ref{cs21})$ we will deduce $H^i\left( X, \Theta_{\mathcal{F}_{X_0}}^*\otimes \frac{\Theta_X}{\Theta_{\mathcal{F}_{X_0}}} \right)$ for $i=0,1$. We note that we have the exact sequence $0\to \Theta_X \to \Theta_Y|_X \to \mathcal{O}_Y(X)|_X\to 0$ induces
\begin{align}\label{cs22}
0\to \Theta_X \otimes \Theta_{\mathcal{F}_{X_0}}^* \to  \left(\Theta_Y\otimes \Theta_{\tilde{\mathcal{F}}_0^\sharp}^*\right) |_X \to \left(\mathcal{O}_Y(X)\otimes \Theta_{\tilde{\mathcal{F}}_0^\sharp}^*\right)|_X = \mathcal{N}_{X/Y}\otimes \Theta_{\mathcal{F}_{X_0}}^*  \to 0
\end{align}
By explicitly describing $H^i\left(X, \left(\Theta_{\tilde{\mathcal{F}}_0^\sharp}^* \otimes \Theta_Y\right) |_X  \right)$ and $H^i\left(X, \left( \mathcal{O}_Y(X)\otimes \Theta_{\tilde{\mathcal{F}}_0^\sharp}^* \right)|_X \right)$ for $i=0,1$, we will deduce $H^i\left( X, \Theta_{\mathcal{F}_{X_0}}^* \otimes \Theta_X \right)$ for $i=0,1$ from $(\ref{cs22})$. Let us compute $H^i\left( X, \left( \Theta_{\tilde{\mathcal{F}}_0^\sharp}^* \otimes \Theta_Y\right)|_X \right)$ for $i=0,1$. We note that we have the exact sequence
\begin{align} \label{cs29}
0 \to \Theta_{\tilde{\mathcal{F}}_0^\sharp}^* \otimes \Theta_Y \otimes \mathcal{O}_Y(-X) \to \Theta_{\tilde{\mathcal{F}}_0^\sharp}^* \otimes \Theta_Y \to \left( \Theta_{\tilde{\mathcal{F}}_0^\sharp}^* \otimes \Theta_Y  \right)|_X \to 0
\end{align}
We compute $H^i\left( Y, \Theta_{ \tilde{\mathcal{F}}_0^\sharp}^* \otimes \Theta_Y \otimes \mathcal{O}_Y(-X) \right)$ for $i=0,1$.  We note that by duality
\begin{align}\label{cs26}
H^i \left( Y, \Theta_{ \tilde{\mathcal{F}}_0^\sharp}^* \otimes \Theta_Y \otimes \mathcal{O}_Y(-X) \right) \cong H^{3-i}\left( Y, \Omega_Y^1\otimes \Theta_{\tilde{\mathcal{F}}_0^\sharp}\otimes \mathcal{O}_Y(X)\otimes K_Y\right)^*\,\,\,\,\,\,\,\,\textnormal{for}\,\,\, i=0,1,2
\end{align}
We note that from $(\ref{cs25})$ $\mathcal{O}_Y(X)$ is defined by the transition functions in terms of the covering $\mathcal{U}^\sharp$
\begin{equation} \label{cs28}
\begin{cases}
\textnormal{transition functions of}\,\,\,\,\, \mathcal{O}_{\mathbb{P}_\mathbb{C}^2}(2) & \textnormal{between $U_{00}$ and $U_{10}$ and  $U_{20}$}\\
t^2 & \textnormal{between $U_{00}$ and $U_{01}$ }\\
t^2 & \textnormal{between $U_{10}$ and $U_{11}$ }\\
t^2 & \textnormal{between $U_{20}$ and $U_{21}$ }
\end{cases}
\end{equation}

We note that we have the exact sequence from \cite{Wav68} p. 943
\begin{align*}
0\to \mathcal{C}\to \Omega_Y^1 \otimes \mathcal{O}_Y(X)\otimes K_Y \to \mathcal{D}\to 0 
\end{align*}
where the sheaf $\mathcal{C}$ is defined by the transition maps in terms of the covering $\mathcal{U}^\sharp$
\begin{equation}
\begin{cases}
\textnormal{transition maps of}\,\,\,\,\, \mathcal{O}_{\mathbb{P}_\mathbb{C}^2}(-3)\otimes \Omega_{\mathbb{P}_\mathbb{C}^2}^1 \otimes \mathcal{O}_{\mathbb{P}_\mathbb{C}^2}(-1)\otimes \mathcal{O}_{\mathbb{P}_\mathbb{C}^2}(2)  =\Omega_{\mathbb{P}_\mathbb{C}^2}^1(-2) & \textnormal{between $U_{00}$ and $U_{10}$ and  $U_{20}$}\\
-\textnormal{Id} & \textnormal{between $U_{00}$ and $U_{01}$ }\\
-\textnormal{Id} & \textnormal{between $U_{10}$ and $U_{11}$ }\\
-\textnormal{Id} & \textnormal{between $U_{20}$ and $U_{21}$ }
\end{cases}
\end{equation}
and the sheaf $\mathcal{D}$ is defined by the transition functions in terms of the covering $\mathcal{U}^\sharp$
\begin{equation}
\begin{cases}
\textnormal{transition functions of}\,\,\,\,\, \mathcal{O}_{\mathbb{P}_\mathbb{C}^2}(-3)& \textnormal{between $U_{00}$ and $U_{10}$ and  $U_{20}$}\\
\frac{1}{t^2} & \textnormal{between $U_{00}$ and $U_{01}$ }\\
\frac{1}{t^2} & \textnormal{between $U_{10}$ and $U_{11}$ }\\
\frac{1}{t^2} & \textnormal{between $U_{20}$ and $U_{21}$ }
\end{cases}
\end{equation}
Then from $(\ref{cs26})$ and $(\ref{cs15})$ we have the exact sequence
\begin{align} \label{cs27}
0\to \mathcal{C}\otimes \Theta_{\tilde{\mathcal{F}}_0^\sharp} \to \Omega_Y^1 \otimes \mathcal{O}_Y(X)\otimes K_Y \otimes  \Theta_{\tilde{\mathcal{F}}_0^\sharp} \to \mathcal{D} \otimes  \Theta_{\tilde{\mathcal{F}}_0^\sharp} \to 0  
\end{align}
where the sheaf $\mathcal{C}\otimes  \Theta_{\tilde{\mathcal{F}}_0^\sharp}$ is defined by the transition maps in terms of the covering $\mathcal{U}^\sharp$
\begin{equation}
\begin{cases}
\textnormal{transition maps of}\,\,\,\,\, \Omega_{\mathbb{P}_\mathbb{C}^2}^1(-3) & \textnormal{between $U_{00}$ and $U_{10}$ and  $U_{20}$}\\
-\textnormal{Id} & \textnormal{between $U_{00}$ and $U_{01}$ }\\
-\textnormal{Id} & \textnormal{between $U_{10}$ and $U_{11}$ }\\
-\textnormal{Id} & \textnormal{between $U_{20}$ and $U_{21}$ }
\end{cases}
\end{equation}
and the sheaf $\mathcal{D}\otimes  \Theta_{\tilde{\mathcal{F}}_0^\sharp}$ is defined by the transition functions in terms of the covering $\mathcal{U}^\sharp$
\begin{equation}
\begin{cases}
\textnormal{transition functions of}\,\,\,\,\, \mathcal{O}_{\mathbb{P}_\mathbb{C}^2}(-4)& \textnormal{between $U_{00}$ and $U_{10}$ and  $U_{20}$}\\
\frac{1}{t^2} & \textnormal{between $U_{00}$ and $U_{01}$ }\\
\frac{1}{t^2} & \textnormal{between $U_{10}$ and $U_{11}$ }\\
\frac{1}{t^2} & \textnormal{between $U_{20}$ and $U_{21}$ }
\end{cases}
\end{equation}
Then by Proposition \ref{cs11} we have $H^i\left(Y, \mathcal{C} \otimes \Theta_{\tilde{\mathcal{F}}_0^\sharp} \right)\cong H^i\left( \mathbb{P}_\mathbb{C}^2, \Omega_{\mathbb{P}_\mathbb{C}^2}^1(-3) \right)$ for $i=1,2,3$, so that we have
\begin{align*}
H^1\left( Y, \mathcal{C}\otimes \Theta_{\tilde{\mathcal{F}}_0^\sharp} \right)=H^3\left( Y ,  \mathcal{C}\otimes \Theta_{\tilde{\mathcal{F}}_0^\sharp} \right)=0,\,\,\,\,\,\,\,\,\,\dim_\mathbb{C} H^2\left( Y, \mathcal{C}\otimes \Theta_{\tilde{\mathcal{F}}_0^\sharp} \right)= \dim_\mathbb{C} H^2\left( \mathbb{P}_\mathbb{C}^2, \Omega_{\mathbb{P}_\mathbb{C}^2}^1(-3) \right) = 8
\end{align*}
On the other hand, by Proposition \ref{cs34} we have $H^i\left( Y, \mathcal{D}\otimes \Theta_{\tilde{\mathcal{F}}_0^\sharp} \right)\cong H^{i-1}\left( \mathbb{P}_\mathbb{C}^2, \mathcal{O}_{\mathbb{P}_\mathbb{C}^2}(-3) \right)$ for $i=1,2,3$, so that we have
\begin{align*}
H^i\left( Y, \mathcal{D}\otimes \Theta_{\tilde{\mathcal{F}}_0^\sharp} \right) =0\,\,\,\,\,\,\,\,\textnormal{for $i=1,2$},\,\,\,\,\,\,\, \dim_\mathbb{C} H^3\left( Y, \mathcal{D}\otimes \Theta_{\tilde{\mathcal{F}}_0^\sharp} \right)= 1
\end{align*}
Then from the long exact sequence from $(\ref{cs27})$ we see that
{\small{\begin{align*}
H^0\left( Y, \Omega_Y^1\otimes \mathcal{O}_Y(X)\otimes K_Y\otimes \Theta_{\tilde{\mathcal{F}}_0^\sharp} \right)=0,\dim_\mathbb{C} H^2 \left( Y, \Omega_Y^1\otimes \mathcal{O}_Y(X)\otimes K_Y\otimes \Theta_{\tilde{\mathcal{F}}_0^\sharp} \right)= 8, \dim_\mathbb{C} H^3\left( Y, \Omega_Y^1\otimes \mathcal{O}_Y(X)\otimes K_Y\otimes \Theta_{\tilde{\mathcal{F}}_0^\sharp} \right)= 1
\end{align*}}}
which implies that from $(\ref{cs26})$
\begin{align*}
\dim_\mathbb{C} H^0\left(Y, \Theta_{\tilde{\mathcal{F}}_0^\sharp}^* \otimes \Theta_Y \otimes \mathcal{O}_Y(-X) \right)= 1, \dim_\mathbb{C} H^1\left(Y, \Theta_{\tilde{\mathcal{F}}_0^\sharp}^* \otimes \Theta_Y \otimes \mathcal{O}_Y(-X) \right)= 8, H^2\left(Y, \Theta_{\tilde{\mathcal{F}}_0^\sharp}^* \otimes \Theta_Y \otimes \mathcal{O}_Y(-X \right)=0
\end{align*}
Then from $(\ref{cs19})$ and the long exact sequence from $(\ref{cs29})$ we have the exact sequence
{\small{\begin{align*}
0\to H^0\left(Y,  \Theta_{\tilde{\mathcal{F}}_0^\sharp}^* \otimes \Theta_Y \otimes \mathcal{O}_Y(-X) \right)\to  H^0\left( Y, \Theta_{\tilde{\mathcal{F}}_0^\sharp}^* \otimes \Theta_Y \right)\to H^0\left( X, \left( \Theta_{\tilde{\mathcal{F}}_0^\sharp}^* \otimes \Theta_Y \right)|_X \right) \to H^1\left( Y,  \Theta_{\tilde{\mathcal{F}}_0^\sharp}^* \otimes \Theta_Y \otimes \mathcal{O}_Y(-X)\right)\to 0
\end{align*}}}
which implies that
\begin{align}\label{cs32}
\dim_\mathbb{C} H^0\left( X, \left( \Theta_{\tilde{\mathcal{F}}_0^\sharp}^* \otimes \Theta_Y \right)|_X \right) = (25-1)+8 = 32
\end{align}
and we have
\begin{align}\label{cs33}
H^1\left( X, \left( \Theta_{\tilde{\mathcal{F}}_0^\sharp}^* \otimes \Theta_Y \right)|_X \right)=0
\end{align}

Next let us compute $H^i\left(X, \left(\mathcal{O}_Y(X)\otimes \Theta_{\tilde{\mathcal{F}}_0^\sharp}^* \right)|_X \right)$ for $i=0,1$. We note that we have the exact sequence
\begin{align}\label{cs30}
0\to \Theta_{\tilde{\mathcal{F}}_0^\sharp}^* \to \mathcal{O}_Y(X)\otimes \Theta_{\tilde{\mathcal{F}}_0^\sharp}^* \to \left( \mathcal{O}_Y(X) \otimes \Theta_{\tilde{\mathcal{F}}_0^\sharp}^* \right)|_X \to 0
\end{align}
We note that from $(\ref{cs15})$ and Proposition \ref{cs11}, we have $H^i\left( Y, \Theta_{\tilde{\mathcal{F}}_0^\sharp}^* \right)=H^i\left( \mathbb{P}_\mathbb{C}^2, \mathcal{O}_{\mathbb{P}_\mathbb{C}^2}(1)\right)$, so that
\begin{align*}
\dim_\mathbb{C} H^0\left( Y, \Theta_{\tilde{\mathcal{F}}_0^\sharp}^* \right)=3,\,\,\,\,\,\,\,\,H^1\left(Y, \Theta_{\tilde{\mathcal{F}}_0^\sharp}^* \right)= H^2\left(Y, \Theta_{\tilde{\mathcal{F}}_0^\sharp}^* \right) = 0
\end{align*}
We note that from $(\ref{cs28})$ and $(\ref{cs15})$ $\mathcal{O}_Y(X)\otimes \Theta_{\tilde{\mathcal{F}}_0^\sharp}^*$ is defined by the transition functions in terms of the covering $\mathcal{U}^\sharp$
\begin{equation}\label{cs35}
\begin{cases}
\textnormal{transition functions of}\,\,\,\,\, \mathcal{O}_{\mathbb{P}_\mathbb{C}^2}(3)& \textnormal{between $U_{00}$ and $U_{10}$ and  $U_{20}$}\\
t^2 & \textnormal{between $U_{00}$ and $U_{01}$ }\\
t^2 & \textnormal{between $U_{10}$ and $U_{11}$ }\\
t^2 & \textnormal{between $U_{20}$ and $U_{21}$ }
\end{cases}
\end{equation}
By Proposition \ref{cs11} we have $H^i\left( Y, \mathcal{O}_Y(X)\otimes \Theta_{\tilde{\mathcal{F}}_0^\sharp}^*  \right)= H^i\left(\mathbb{P}_\mathbb{C}^2, \mathcal{O}_{\mathbb{P}_\mathbb{C}^2}(3) \right)\bigoplus H^i\left( \mathbb{P}_\mathbb{C}^2, \mathcal{O}_{\mathbb{P}_\mathbb{C}^2}(2) \right)\bigoplus H^i\left( \mathbb{P}_\mathbb{C}^2, \mathcal{O}_{\mathbb{P}_\mathbb{C}^2}(1) \right)$ for $i=0,1$, so that we have
\begin{align*}
\dim_\mathbb{C} H^0\left(  Y, \mathcal{O}_Y(X)\otimes \Theta_{\tilde{\mathcal{F}}_0^\sharp}^* \right)= 19,\,\,\,\,\,\,\,\,\,\,H^1\left(Y, \mathcal{O}_Y(X)\otimes \Theta_{\tilde{\mathcal{F}}_0^\sharp}^* \right) = 0
\end{align*}
Then from the long exact sequence from $(\ref{cs30})$ we see that
\begin{align}\label{cs31}
\dim_\mathbb{C} H^0\left( X, \left(\mathcal{O}_Y(X)\otimes \Theta_{\tilde{\mathcal{F}}_0^\sharp}\right)|_X \right) = 19-3 = 16,\,\,\,\,\,\,\,\, H^1\left( X, \left(\mathcal{O}_Y(X)\otimes \Theta_{\tilde{\mathcal{F}}_0^\sharp}\right)|_X  \right) = 0
\end{align}

Hence from $(\ref{cs33}),(\ref{cs31})$ and the long exact sequence from $(\ref{cs22})$, we have the exact sequence
\begin{align}\label{cs39}
0\to H^0\left(X, \Theta_{\mathcal{F}_{X_0}}^* \otimes \Theta_X \right)\to H^0\left( X, \left( \Theta_{\tilde{\mathcal{F}}_0^\sharp}^* \otimes \Theta_Y \right)|_X \right) \to H^0\left( X, \left(\mathcal{O}_Y(X)\otimes \Theta_{\tilde{\mathcal{F}}_0^\sharp}^* \right)|_X \right)\to H^1\left( X, \Theta_{\mathcal{F}_{X_0}}^* \otimes \Theta_X \right)\to 0
\end{align}
We will compute $H^0\left(X, \Theta_{\mathcal{F}_{X_0}}^* \otimes \Theta_X \right)$ by computing 
\begin{align*}
\textnormal{ker}\left( H^0\left( X, \left( \Theta_{\tilde{\mathcal{F}}_0^\sharp}^* \otimes \Theta_Y \right)|_X \right) \to H^0\left( X, \left(\mathcal{O}_Y(X)\otimes \Theta_{\tilde{\mathcal{F}}_0^\sharp}^* \right)|_X \right)\right)
\end{align*}
Then we need to describe explicitly $H^0\left( X, \left( \Theta_{\tilde{\mathcal{F}}_0^\sharp}^* \otimes \Theta_Y \right)|_X \right)$ and $H^0\left( X, \left(\mathcal{O}_Y(X)\otimes \Theta_{\tilde{\mathcal{F}}_0^\sharp}^* \right)|_X  \right)$.

\subsubsection{Explicit description of $H^0\left( X, \left( \Theta_{\tilde{\mathcal{F}}_0^\sharp}^* \otimes \Theta_Y \right)|_X \right)$ on $U_{00}\cap X$}\

We recall from $(\ref{cs32})$ that $\dim_\mathbb{C} H^0\left( X, \left( \Theta_{\tilde{\mathcal{F}}_0^\sharp}^* \otimes \Theta_Y \right)|_X \right)=32$. We note that since $t^2=1+ X_1^2 + X_2^2$, we can write any element in $H^0\left( X, \left( \Theta_{\tilde{\mathcal{F}}_0^\sharp}^* \otimes \Theta_Y \right)|_X \right)$ on $U_{00}\cap X$ by
\begin{align*}
V_{00}=\left(A_0(X_1,X_2)+ t A_1(X_1,X_2) \right) \frac{\partial}{\partial X_1}+\left(B_0(X_1,X_2)+t B_1(X_1,X_2) \right) \frac{\partial}{\partial X_2} + \left( C_0(X_1,X_2)+ t C_1(X_1,X_2) \right)\frac{\partial}{\partial t}
\end{align*}
where $A_{0}(X_1,X_2)= \sum_{i,j\geq 0} a_{ij}X_1^i X_2^j, a_{ij}\in \mathbb{C}, A_1(X_1,X_2)=\sum_{i,j\geq 0} f_{ij} X_1^i X_2^j,f_{ij}\in \mathbb{C}, B_0(X_1,X_2)=\sum_{i,j\geq 0} b_{ij} X_1^iX_2^j, b_{ij}\in \mathbb{C}$, $B_1(X_1,X_2)= \sum_{i,j\geq 0} g_{ij}X_1^i X_2^j, g_{ij}\in \mathbb{C}, C_0(X_1,X_2)= \sum_{i,j\geq 0} h_{ij} X_1^i X_2^j, h_{ij}\in \mathbb{C},C_1(X_1,X_2)= \sum_{i,j\geq 0} c_{ij} X_1^i X_2^j, c_{ij}\in \mathbb{C}$.
We note that in view of $(\ref{cs15})$ $V_{00}$ on $U_{00}\cap X$ is transformed on $U_{01}\cap X$ into
\begin{align*}
V_{01}= \left(A_0(X_1,X_2)+ \frac{1}{t'} A_1(X_1,X_2) \right) \frac{\partial}{\partial X_1}+\left(B_0(X_1,X_2)+\frac{1}{t'} B_1(X_1,X_2) \right) \frac{\partial}{\partial X_2} + \left( - t'^2 C_0(X_1,X_2)- t' C_1(X_1,X_2) \right)\frac{\partial}{\partial t'}
\end{align*}
We note that since $t'^2\left(1+ X_1^2 + X_2^2 \right)=1$ on $U_{01}\cap X$, we see that $t'$ is a unit, so that $V_{01}$ is still holomorphic on $U_{01}\cap X$. We note that in view of $(\ref{cs15})$ and $(\ref{css10})$ $V_{00}$ on $U_{00}\cap X$ is transformed on $U_{10}\cap X$ into
\begin{align*}
V_{10}&=Y_0\left( A_0\left( \frac{1}{Y_0}, \frac{Y_2}{Y_0} \right) + \frac{t'}{Y_0} A_1\left( \frac{1}{Y_0}, \frac{Y_2}{Y_0} \right)\right)\left( - Y_0^2 \frac{\partial}{\partial Y_0}- Y_0 Y_2\frac{\partial}{\partial Y_2}- t' Y_0\frac{\partial}{\partial t'} \right)\\
&+Y_0\left( B_0\left( \frac{1}{Y_0}, \frac{Y_2}{Y_0} \right) + \frac{t'}{Y_0} B_1\left( \frac{1}{Y_0}, \frac{Y_2}{Y_0} \right)\right) Y_0\frac{\partial}{\partial Y_2} + Y_0\left( C_0\left( \frac{1}{Y_0}, \frac{Y_2}{Y_0} \right) + \frac{t'}{Y_0} C_1\left( \frac{1}{Y_0}, \frac{Y_2}{Y_0} \right) \right) Y_0 \frac{\partial}{\partial t'} \\
&= \left( - Y_0^3 A_0\left( \frac{1}{Y_0}, \frac{Y_2}{Y_0}\right) - t' Y_0^2 A_1\left( \frac{1}{Y_0}, \frac{Y_2}{Y_0} \right) \right) \frac{\partial}{\partial Y_0}\\
&+\left( - Y_0^2 Y_2 A_0\left( \frac{1}{Y_0}, \frac{Y_2}{Y_0} \right) + Y_0^2B_0\left( \frac{1}{Y_0}, \frac{Y_2}{Y_0} \right)+ t' \left( Y_0B_1\left( \frac{1}{Y_0}, \frac{Y_2}{Y_0} \right)  -Y_0 Y_2 A_1\left( \frac{1}{Y_0}, \frac{Y_2}{Y_0} \right) \right) \right) \frac{\partial}{\partial Y_2}\\
&+\left(Y_0^2 C_0\left( \frac{1}{Y_0}, \frac{Y_2}{Y_0} \right) + t' \left( Y_0 C_1\left( \frac{1}{Y_0}, \frac{Y_2}{Y_0} \right)  - Y_0^2 A_0\left( \frac{1}{Y_0}, \frac{Y_2}{Y_0}  \right) \right) - t'^2 Y_0 A_1\left( \frac{1}{Y_0}, \frac{Y_2}{Y_0} \right) \right) \frac{\partial}{\partial t'}
\end{align*}
Since $t'^2= 1+ Y_0^2 + Y_2^2$ on $U_{10}\cap X$, we have
\begin{align*}
V_{10}&=\left( - Y_0^3 A_0\left( \frac{1}{Y_0}, \frac{Y_2}{Y_0}\right) - t' Y_0^2 A_1\left( \frac{1}{Y_0}, \frac{Y_2}{Y_0} \right) \right) \frac{\partial}{\partial Y_0}\\
&+\left( - Y_0^2 Y_2 A_0\left( \frac{1}{Y_0}, \frac{Y_2}{Y_0} \right) + Y_0^2B_0\left( \frac{1}{Y_0}, \frac{Y_2}{Y_0} \right)+ t' \left( Y_0B_1\left( \frac{1}{Y_0}, \frac{Y_2}{Y_0} \right)  -Y_0 Y_2 A_1\left( \frac{1}{Y_0}, \frac{Y_2}{Y_0} \right) \right) \right) \frac{\partial}{\partial Y_2}\\
&+\left(Y_0^2 C_0\left( \frac{1}{Y_0}, \frac{Y_2}{Y_0} \right) - \left( Y_0 + Y_0^3 + Y_0 Y_2^2 \right)A_1\left( \frac{1}{Y_0}, \frac{Y_2}{Y_0} \right) + t' \left( Y_0 C_1\left( \frac{1}{Y_0}, \frac{Y_2}{Y_0} \right)  - Y_0^2 A_0\left( \frac{1}{Y_0}, \frac{Y_2}{Y_0}  \right) \right)  \right) \frac{\partial}{\partial t'}
\end{align*}
which should be holomorphic on $U_{10}\cap X$. We note that in order for
\begin{align*}
- Y_0^2 Y_2 A_0\left( \frac{1}{Y_0}, \frac{Y_2}{Y_0} \right) + Y_0^2 B_0\left( \frac{1}{Y_0}, \frac{Y_2}{Y_0} \right) =  \sum_{i,j\geq 0} \frac{- a_{ij} Y_2^{j+1}+ b_{ij} Y_2^j }{Y_0^{i+j-2}}
\end{align*}
to be holomorphic, we only need to consider $(i,j)=(0,0),(0,1),(1,0),(2,0),(1,1),(0,2)$ which clearly defines a holomorphic function on $U_{10}\cap X$ and $(0,3),(1,2),(2,1),(0,3)$ in which case we have
\begin{align*}
\frac{- a_{03} Y_2^4+ b_{03} Y_2^3- a_{12} Y_2^3 + b_{12} Y_2^2 - a_{21} Y_2^2 +b_{21} Y_2 -a_{30} Y_2 + b_{30}}{Y_0} \Longrightarrow a_{03}=0, b_{03}= a_{12} , b_{12} = a_{21}, b_{21}= a_{30}, b_{30}=0
\end{align*}
Then we have
\begin{align*}
A_0(X_1,X_2)&= a_{00} + a_{10} X_1 + a_{01} X_2 + a_{20} X_1^2 + a_{11} X_1 X_2 + a_{02}X_2^2 + a_{30}X_1^3 + a_{21} X_1^2 X_2+ a_{12} X_1 X_2^2\\
B_0(X_1,X_2)&=  b_{00} + b_{10} X_1+ b_{01} X_2 + b_{20} X_1^2 + b_{11} X_1 X_2 + b_{02} X_2^2  + a_{30}X_1^2 X_2 + a_{21} X_1 X_2^2 + a_{12} X_2^3
\end{align*}
We note that in order for
\begin{align*}
Y_0 B_1\left( \frac{1}{Y_0}, \frac{Y_2}{Y_0} \right)- Y_0 Y_2 A_1\left( \frac{1}{Y_0}, \frac{Y_2}{Y_0} \right) = \sum_{i,j\geq 0} \frac{g_{ij} Y_2^j - f_{ij} Y_2^{j+1}}{Y_0^{i+ j-1}}
\end{align*}
to be holomorphic, we only need to consider $(i,j)=(0,0),(0,1),(1,0)$ which clearly defines a holomorphic function on $U_{10}\cap X$ and $(0,2),(1,1),(2,0)$ in which case we have
\begin{align*}
 \frac{ g_{20} - f_{20} Y_2+ g_{11} Y_2- f_{11} Y_2^2 + g_{02} Y_2^2- f_{02} Y_2^3}{Y_0} \Longrightarrow g_{20}=0, f_{20}= g_{11}, f_{11}= g_{02}, f_{02}=0
\end{align*}
Then we have
\begin{align*}
A_1(X_1,X_2)&= f_{00} + f_{10} X_1 + f_{01} X_2 + f_{20}X_1^2 + f_{11} X_1 X_2 \\
B_1(X_1,X_2)&= g_{00} + g_{10} X_1 + g_{01} X_2 + f_{20}X_1X_2+ f_{11} X_2^2
\end{align*}
On the other hand, in order for
{\small{\begin{align*}
Y_0^2 C_0\left( \frac{1}{Y_0}, \frac{Y_2}{Y_0} \right) - \left( Y_0 + Y_0^3 + Y_0 Y_2^2 \right)A_1\left( \frac{1}{Y_0}, \frac{Y_2}{Y_0} \right) = Y_0^2 C_0\left( \frac{1}{Y_0}, \frac{Y_2}{Y_0} \right) - \left(Y_0 + Y_0^3+ Y_0 Y_2^2 \right)\left( f_{00} + \frac{f_{10}+f_{10} Y_2}{Y_0}  + \frac{f_{20}+ f_{11} Y_2 }{Y_0^2}\right)
\end{align*}}}
to be holomorphic, we should have
\begin{align*}
 \frac{  (h_{30}- f_{20}) + (h_{21}-f_{11})Y_2 + ( h_{12} - f_{20}) Y_2^2 + (h_{03} - f_{11}) Y_2^3 }{Y_0}=0 \Longrightarrow h_{30}=f_{20} , h_{21}= f_{11}, h_{12}= f_{20}, h_{03}= f_{11}
\end{align*}
Then we have
\begin{align*}
C_0(X_1,X_2)= h_{00} + h_{10} X_1 + h_{01} X_2  + h_{20}X_1^2 + h_{11} X_1 X_2 + h_{02} X_2^2+ f_{20}X_1^2+ f_{11} X_1^2 X_2 + f_{20}X_1 X_2^2 + f_{11} X_2^3
\end{align*}
We note that in order for
{\small{\begin{align*}
Y_0 C_1\left( \frac{1}{Y_0}, \frac{Y_2}{Y_0} \right) - Y_0^2 A_0 \left( \frac{1}{Y_0}, \frac{Y_2}{Y_0} \right) =Y_0 C_1\left( \frac{1}{Y_0}, \frac{Y_2}{Y_0} \right) - Y_0^2\left(a_{00}+ \frac{a_{10}+ a_{01} Y_2}{Y_0} +\frac{a_{20}+ a_{11} Y_2 + a_{02} Y_2^2}{Y_0^2} + \frac{a_{30}+ a_{21}Y_2+ a_{12} Y_2^2}{Y_0^3} \right) 
\end{align*}}}
to be holomorphic, we should have
\begin{align*}
 \frac{c_{02} Y_2^2 + c_{11}Y_0 + c_{20}}{Y_0} - \frac{ a_{12} Y_2^2 + a_{21} Y_2 + a_{30} }{Y_0}=0 \Longrightarrow c_{02}= a_{12}, c_{11}=a_{21}, c_{20}= a_{30}
\end{align*}
Then we have
\begin{align*}
C_1(X_1,X_2)= c_{00} + c_{10} X_1 + c_{01} X_2 + a_{30}X_1^2+ a_{21}X_1X_2+ a_{12} X_2^2
\end{align*}
Hence $H^0\left( X, \left( \Theta_{\tilde{\mathcal{F}}_0^\sharp}^* \otimes \Theta_Y \right)|_X \right)$ on $U_{00}\cap X$ is given by
{\small{\begin{align}\label{cs37}
H^0\left(X, \left( \Theta_{\tilde{\mathcal{F}}_0^\sharp}^* \otimes \Theta_Y \right)|_X \right)=\left\{\left( a_{00} + a_{10} X_1 + a_{01} X_2 + a_{20} X_1^2 + a_{11} X_1 X_2 + a_{02}X_2^2 + a_{30}X_1^3 + a_{21} X_1^2 X_2+ a_{12} X_1 X_2^2\right)\frac{\partial}{\partial X_1}\right.  \\
\left.+t \left( f_{00} + f_{10} X_1 + f_{01} X_2 + f_{20}X_1^2 + f_{11} X_1 X_2  \right) \frac{\partial}{\partial X_1}\right. \notag \\
\left.+ \left( b_{00} + b_{10} X_1+ b_{01} X_2 + b_{20} X_1^2 + b_{11} X_1 X_2 + b_{02} X_2^2  + a_{30}X_1^2 X_2 + a_{21} X_1 X_2^2 + a_{12} X_2^3 \right)\frac{\partial}{\partial X_2}\right. \notag \\
 \left.+ t \left( g_{00} + g_{10} X_1 + g_{01} X_2 + f_{20}X_1X_2+ f_{11} X_2^2 \right) \frac{\partial}{\partial X_2}\right. \notag\\
\left.+\left( h_{00} + h_{10} X_1 + h_{01} X_2 + h_{20}X_1^2 + h_{11} X_1 X_2 + h_{02} X_2^2+  f_{20}X_1^3+ f_{11} X_1^2 X_2+ f_{20}X_1 X_2^2 + f_{11} X_2^3 \right) \frac{\partial}{\partial t}\right.  \notag \\
\left.+t \left( c_{00} + c_{10} X_1 + c_{01} X_2+ a_{30} X_1^2 + a_{21} X_1 X_2+ a_{12} X_2^2  \right) \frac{\partial}{\partial t}\right\} \notag
\end{align}}}
where $a_{00},a_{10},a_{01}, a_{20},a_{11},a_{02},a_{30},a_{21},a_{12},f_{00},f_{10},f_{01},f_{20}, f_{11}, b_{00},b_{10},b_{01}, b_{20},b_{11},b_{02}, g_{00}, g_{10}, g_{01}, h_{00}, h_{10}$, $ h_{01}, h_{20}, h_{11}, h_{02}, c_{00}, c_{10}, c_{01} \in \mathbb{C}$.

\subsubsection{Explicit description of $H^0\left( X, \left(\mathcal{O}_Y(X)\otimes \Theta_{\tilde{\mathcal{F}}_0^\sharp}^* \right)|_X \right)$ on $U_{00}\cap X $}\

We recall from $(\ref{cs31})$ that $\dim_\mathbb{C} H^0\left( X, \left(\mathcal{O}_Y(X)\otimes \Theta_{\tilde{\mathcal{F}}_0^\sharp}^* \right)|_X \right)= 16$. We have the exact sequence
\begin{align}\label{cs36}
0\to H^0\left( Y, \Theta_{\tilde{\mathcal{F}}_0}^* \right)\xrightarrow{\cdot \left(1+ X_1^2+ X_2^2 - t^2 \right)} H^0\left( Y,  \mathcal{O}_Y(X)\otimes \Theta_{\tilde{\mathcal{F}}_0^\sharp}^* \right)\to H^0\left( X, \left(\mathcal{O}_Y(X)\otimes \Theta_{\tilde{\mathcal{F}}_0^\sharp}^* \right)|_X \right)\to 0
\end{align}
In view of $(\ref{cs15})$ $H^0\left(Y, \Theta_{\tilde{\mathcal{F}}_0^\sharp}^* \right)$ on $U_{00}$ is given by
\begin{align*}
H^0\left( Y, \Theta_{\tilde{\mathcal{F}}_0^\sharp} \right) =\left\{b_0+ b_1 X_1 + b_2 X_2 |b_0,b_1,b_2\in \mathbb{C} \right\}
\end{align*}
and in view of $(\ref{cs35})$ $H^0\left(Y, \mathcal{O}_Y(X)\otimes \Theta_{\tilde{\mathcal{F}}_0^\sharp}^* \right)$ on $U_{00}$ is given by
\begin{align*}
H^0\left(Y, \mathcal{O}_Y(X)\otimes \Theta_{\tilde{\mathcal{F}}_0^\sharp}^* \right)=\left\{ b_0 + b_1 X_1 + b_2 X_2+ b_3 X_1^2 + b_4 X_1 X_2 + b_5 X_2^2 + b_6 X_1^3+ b_7 X_1^2 X_2 + b_8 X_1 X_2^2 + b_9 X_2^3 \right.\\
\left. +\left(b_{10}+ b_{11}X_1 + b_{12} X_2 + b_{13}X_1^2 + b_{14}X_1X_2 + b_{15}X_2^2 \right)t +\left(b_{16} + b_{17}X_1 + b_{18}X_2  \right)t^2 \right\}
\end{align*}
where $b_0,b_1,..., b_{18}\in \mathbb{C}$. Then from $(\ref{cs36})$, we see that $H^0\left(X, \left( \mathcal{O}_Y(X)\otimes  \Theta_{\tilde{\mathcal{F}}_0^\sharp}^* \right)|_X \right)$ on $U_{00}\cap X$ is given by
\begin{align} \label{cs38}
H^0\left(X, \left( \mathcal{O}_Y(X)\otimes  \Theta_{\tilde{\mathcal{F}}_0^\sharp}^* \right)|_X\right)=&\left\{ b_0 + b_1 X_2 + b_2 X_2 + b_3 X_1^2 + b_4 X_1 X_2 + b_5 X_2^2 + b_6 X_1^3+ b_7 X_1^2 X_2 + b_8 X_1 X_2^2 + b_9 X_2^3 \right.\\
&\left. +\left(b_{10}+ b_{11}X_1 + b_{12} X_2 + b_{13}X_1^2 + b_{14}X_1X_2 + b_{15}X_2^2 \right)t  \right\} \notag
\end{align}
where $b_0,...,b_{15}\in \mathbb{C}$.

\subsubsection{Computation of $\textnormal{ker}\left( H^0\left(X,  \left( \Theta_{\tilde{\mathcal{F}}_0^\sharp} \otimes \Theta_Y \right)|_X \right)\to H^0\left( X, \left( \mathcal{O}_Y(X)\otimes \Theta_{\tilde{\mathcal{F}}_0^\sharp}^*  \right)|_X  \right) \right)$}\

By using $(\ref{cs37})$ and $(\ref{cs38})$, we will explicitly describe $H^0\left( X, \Theta_X \otimes \Theta_{\mathcal{F}_{X_0}} \right)$ on $U_{00}\cap X$ by computing 
\begin{align*}
\textnormal{ker}\left( H^0\left(X,  \left( \Theta_{\tilde{\mathcal{F}}_0^\sharp} \otimes \Theta_Y \right)|_X \right)\to H^0\left( X, \left( \mathcal{O}_Y(X)\otimes \Theta_{\tilde{\mathcal{F}}_0^\sharp}^*  \right)|_X  \right) \right)
\end{align*}
From $(\ref{cs37})$, we see that the image is given by
{\Small{\begin{align*}
&\left[ \left( a_{00} + a_{10} X_1 + a_{01} X_2 + a_{20} X_1^2 + a_{11} X_1 X_2 + a_{02}X_2^2 + a_{30}X_1^3 + a_{21} X_1^2 X_2+ a_{12} X_1 X_2^2\right)\frac{\partial}{\partial X_1} , X_1^2 + X_2^2 + 1- t^2\right]|_{t^2= 1+ X_1^2 + X_2^2}\\
&+\left[ t \left( f_{00} + f_{10} X_1 + f_{01} X_2 + f_{20}X_1^2 + f_{11} X_1 X_2  \right) \frac{\partial}{\partial X_1}, X_1^2 + X_2^2 +1 - t^2 \right]|_{t^2=1+ X_1^2 + X_2^2} \\
& +\left[  \left( b_{00} + b_{10} X_1+ b_{01} X_2 + b_{20} X_1^2 + b_{11} X_1 X_2 + b_{02} X_2^2  + a_{30}X_1^2 X_2 + a_{21} X_1 X_2^2 + a_{12} X_2^3 \right)\frac{\partial}{\partial X_2}, X_1^2 + X_2^2 +1 - t^2 \right]|_{t^2= 1+ X_1^2 + X_2^2} \\
&+\left[  t \left( g_{00} + g_{10} X_1 + g_{01} X_2 + f_{20}X_1X_2+ f_{11} X_2^2 \right) \frac{\partial}{\partial X_2}, X_1^2 + X_2^2+1 -t^2 \right]|_{t^2= 1+ X_1^2+ X_2^2} \\
&+\left[ \left( h_{00} + h_{10} X_1 + h_{01} X_2 + h_{20}X_1^2 + h_{11} X_1 X_2 + h_{02} X_2^2+  f_{20}X_1^3+ f_{11} X_1^2 X_2+ f_{20}X_1 X_2^2 + f_{11} X_2^3 \right) \frac{\partial}{\partial t}, X_1^2 + X_2^2 +1 - t^2 \right]|_{t^2= X_1^2 + X_2^2 +1} \\
&+\left[ t \left( c_{00} + c_{10} X_1 + c_{01} X_2+ a_{30} X_1^2 + a_{21} X_1 X_2+ a_{12} X_2^2  \right) \frac{\partial}{\partial t}, X_1^2 + X_2^2 +1 -  t^2 \right]|_{t^2= 1 + X_1^2 + X_2^2 }\\
&=- 2 c_{00}  + (2 a_{00}- 2 c_{10})X_1 + (2 b_{00}- 2 c_{01})X_2 + ( 2 a_{10}- 2 c_{00}- 2 a_{30})X_1^2 + ( 2 a_{01}+ 2 b_{10}- 2 a_{21}) X_1 X_2 \\
& + ( 2 b_{01} - 2 c_{00} - 2 a_{12})X_2^2 + (2 a_{20}- 2 c_{10})X_1^3 + ( 2 a_{11}+ 2 b_{20} - 2 c_{01}) X_1^2 X_2  + (2 a_{02}+ 2 b_{11}- 2 c_{10})X_1 X_2^2 + ( 2b_{02}- 2 c_{01})X_2^3\\
&+ t \left( - 2 h_{00} + (2 f_{00}- 2 h_{10}) X_1 + (2 g_{00}- 2 h_{01})X_2 + (2 f_{10}- 2 h_{20})X_1^2 + (2 f_{01}+ 2 g_{10}- 2 h_{11}) X_1X_2 + ( 2 g_{01}- 2 h_{02})X_2^2 \right)
\end{align*}}}
Then $\textnormal{ker}\left( H^0\left(X,  \left( \Theta_{\tilde{\mathcal{F}}_0^\sharp} \otimes \Theta_Y \right)|_X \right)\to H^0\left( X, \left( \mathcal{O}_Y(X)\otimes \Theta_{\tilde{\mathcal{F}}_0^\sharp}^*  \right)|_X  \right) \right)$ on $U_{00} \cap X$ is given by
\begin{align} \label{css8}
H^0\left( X, \Theta_{\mathcal{F}_{X_0}}^*\otimes \Theta_X \right)= \left\{ \left( a_{00} + a_{10} X_1 + a_{01} X_2 + a_{00} X_1^2 + a_{11} X_1 X_2 + a_{02}X_2^2 + a_{10}X_1^3 + a_{21} X_1^2 X_2+ a_{12} X_1 X_2^2\right)\frac{\partial}{\partial X_1}\right.\\
\left. +t \left( f_{00} + f_{10} X_1 + f_{01} X_2 + f_{20}X_1^2 + f_{11} X_1 X_2  \right) \frac{\partial}{\partial X_1} \right. \notag \\
\left. + \left( b_{00} + (a_{21}- a_{01} )X_1+ a_{12} X_2 + (b_{00}- a_{11}) X_1^2 + (a_{00}- a_{02}) X_1 X_2 + b_{00} X_2^2  + a_{10}X_1^2 X_2 + a_{21} X_1 X_2^2 + a_{12} X_2^3 \right)\frac{\partial}{\partial X_2}\right. \notag \\
 \left. + t \left( g_{00} + g_{10} X_1 + g_{01} X_2 + f_{20}X_1X_2+ f_{11} X_2^2 \right) \frac{\partial}{\partial X_2} \right. \notag \\
\left. +\left( f_{00} X_1 + g_{00} X_2 + f_{10}X_1^2 + (f_{01}+ g_{10}) X_1 X_2 + g_{01} X_2^2+  f_{20}X_1^3+ f_{11} X_1^2 X_2+ f_{20}X_1 X_2^2 + f_{11} X_2^3 \right) \frac{\partial}{\partial t} \right. \notag \\
\left. +t \left(  a_{00} X_1 + b_{00} X_2+ a_{10} X_1^2 + a_{21} X_1 X_2+ a_{12} X_2^2  \right) \frac{\partial}{\partial t} \right\} \notag
\end{align}
where $a_{00}, a_{10}, a_{01}, a_{11}, a_{02}, a_{21}, a_{12}, f_{00},f_{10}, f_{01}, f_{20}, f_{11},b_{00}, g_{00}, g_{10},g_{01}\in \mathbb{C}$, so that we have
\begin{align*}
\dim_\mathbb{C} H^0\left( X, \Theta_{\mathcal{F}_{X_0}}^* \otimes \Theta_X \right) = 16
\end{align*}
We note that from $(\ref{cs39})$, we infer that $H^1\left( X, \Theta_{\mathcal{F}_{X_0}}^* \otimes \Theta_X \right)=0$. Then from $(\ref{cs20})$, we see that $H^1\left(X, \Theta_{\mathcal{F}_{X_0}}^* \otimes \frac{\Theta_X}{ \Theta_{\mathcal{F}_{X_0}}} \right) =0$. Hence from $(\ref{css1})$ we have
\begin{align*}
\mathbb{H}^2\left( X, \Theta_{\mathcal{F}_{X_0}}^\bullet \right)=0
\end{align*}
In particular, $\left(X, \Theta_{\mathcal{F}_{X_0}} \right)$ is unobstructed in foliated deformations. We explicitly describe $H^0\left(X, \Theta_{\mathcal{F}_{X_0}}^* \otimes \frac{\Theta_X}{\Theta_{\mathcal{F}_{X_0}}} \right)$. From the exact sequence $(\ref{cs21})$ and $(\ref{css8})$, and
\begin{align*}
 a_{21} X_1^2X_2\frac{\partial}{\partial X_1}+ a_{21} X_1 X_2^2\frac{\partial}{\partial X_2} + a_{21}t X_1X_2\frac{\partial}{\partial t}\equiv -\frac{1}{2}a_{21}X_2\frac{\partial}{\partial X_1}- \frac{1}{2}a_{21}X_1\frac{\partial}{\partial X_2}\,\,\,\,\,\,\,\,\, \textnormal{in}\,\,\,\,\, H^0\left(X, \Theta_{\mathcal{F}_{X_0}}^*\otimes \frac{\Theta_X}{\Theta_{\mathcal{F}_{X_0}}^*} \right),
 \end{align*}
we see that $H^0\left( X, \Theta_{\mathcal{F}_{X_0}}^* \otimes \frac{\Theta_X}{\Theta_{\mathcal{F}_{X_0}}} \right)$ on $U_{00}\cap X$ is given by
\begin{align} \label{css9}
H^0\left( X, \Theta_{\mathcal{F}_{X_0}}^*\otimes  \frac{\Theta_X}{\Theta_{\mathcal{F}_{X_0}}} \right)= \left\{ \left( a_{00} + a_{10} X_1 + a_{01} X_2 + a_{00} X_1^2 + a_{11} X_1 X_2 + a_{02}X_2^2 + a_{10}X_1^3 + a_{12} X_1 X_2^2\right)\frac{\partial}{\partial X_1}\right.\\
\left. +t \left( f_{00} + f_{10} X_1 + f_{01} X_2 + f_{20}X_1^2 + f_{11} X_1 X_2  \right) \frac{\partial}{\partial X_1} \right. \notag \\
\left. + \left( b_{00} - a_{01} X_1+ a_{12} X_2 + (b_{00}- a_{11}) X_1^2 + (a_{00}- a_{02}) X_1 X_2 + b_{00} X_2^2  + a_{10}X_1^2 X_2  + a_{12} X_2^3 \right)\frac{\partial}{\partial X_2}\right. \notag \\
 \left. + t \left( g_{00} + g_{10} X_1 + g_{01} X_2 + f_{20}X_1X_2+ f_{11} X_2^2 \right) \frac{\partial}{\partial X_2} \right. \notag \\
\left. +\left( f_{00} X_1 + g_{00} X_2 + f_{10}X_1^2 + (f_{01}+ g_{10}) X_1 X_2 + g_{01} X_2^2+  f_{20}X_1^3+ f_{11} X_1^2 X_2+ f_{20}X_1 X_2^2 + f_{11} X_2^3 \right) \frac{\partial}{\partial t} \right. \notag \\
\left. +t \left(  a_{00} X_1 + b_{00} X_2+ a_{10} X_1^2 + a_{12} X_2^2  \right) \frac{\partial}{\partial t} \right\} \notag
\end{align}
where $a_{00}, a_{10}, a_{01}, a_{11}, a_{02}, a_{12}, f_{00},f_{10}, f_{01}, f_{20}, f_{11},b_{00}, g_{00}, g_{10},g_{01}\in \mathbb{C}$.

Now we will compute $\mathbb{H}^1\left(X, \Theta_{\mathcal{F}_{X_0}}^\bullet \right)=\textnormal{coker}\left(H^0\left( X, \Theta_X \right) \to H^0\left( X, \Theta_{\mathcal{F}_{X_0}}^* \otimes \frac{\Theta_X}{\Theta_{\mathcal{F}_{X_0}}} \right) \right)$. We will explicitly describe $H^0\left( X, \Theta_X \right)$ on $U_{00}\cap X$.

\subsubsection{Explicit description of $H^0\left(X, \Theta_X \right)$}\

We recall that $\dim_\mathbb{C} H^0\left(X, \Theta_X \right)=6$ from $(\ref{css3})$. We note that
\begin{align}
H^0\left( X, \Theta_X\right)= \textnormal{ker}\left( H^0\left(X, \Theta_Y|_X \right) \to H^0\left(X, \mathcal{O}_Y(X)|_X \right) \right) \label{css2}
\end{align}
By explicitly describing $H^0\left(X, \Theta_Y|_X \right)$ and $H^0\left(X,\mathcal{O}_Y(X)|_X \right)$, we will find $H^0\left( X, \Theta_X\right)$ by computing $(\ref{css2})$. First we note that from $(\ref{cs28})$ and Proposition \ref{cs11}, we have $\dim_\mathbb{C} H^0\left(Y, \mathcal{O}_Y(X) \right)=10$. From the exact sequence $0\to \mathcal{O}_Y\to \mathcal{O}_Y(X)\to \mathcal{O}_Y(X)|_X\to 0$, we see that $\dim_\mathbb{C} H^0\left( X, \mathcal{O}_Y(X)|_X\right)=10- 1= 9$. From the exact sequence $0\to \Theta_X \to \Theta_Y|_X\to \mathcal{O}_Y(X)|_X\to 0$, we see that $\dim_\mathbb{C} H^0\left( X, \Theta_Y|_X \right) = 6+ 9 =15$. 

We explicitly describe $H^0\left( X, \Theta_Y|_X \right)$ on $U_{00}\cap X$. We note that since $t^2=1+ X_1^2 + X_2^2$, we can write any element in $H^0\left( X, \Theta_Y|_X \right)$ on $U_{00}\cap X$ by
\begin{align*}
V_{00}=\left(A_0(X_1,X_2)+ t A_1(X_1,X_2) \right) \frac{\partial}{\partial X_1}+\left(B_0(X_1,X_2)+t B_1(X_1,X_2) \right) \frac{\partial}{\partial X_2} + \left( C_0(X_1,X_2)+ t C_1(X_1,X_2) \right)\frac{\partial}{\partial X_3}
\end{align*}
where $A_{0}(X_1,X_2)= \sum_{i,j\geq 0} a_{ij}X_1^i X_2^j, a_{ij}\in \mathbb{C}, A_1(X_1,X_2)=\sum_{i,j\geq 0} f_{ij} X_1^i X_2^j,f_{ij}\in \mathbb{C}, B_0(X_1,X_2)=\sum_{i,j\geq 0} b_{ij} X_1^iX_2^j, b_{ij}\in \mathbb{C}$, $B_1(X_1,X_2)= \sum_{i,j\geq 0} g_{ij}X_1^i X_2^j, g_{ij}\in \mathbb{C}, C_0(X_1,X_2)= \sum_{i,j\geq 0} h_{ij} X_1^i X_2^j, h_{ij}\in \mathbb{C},C_1(X_1,X_2)= \sum_{i,j\geq 0} c_{ij} X_1^i X_2^j, c_{ij}\in \mathbb{C}$. Then $V_{00}$ on $U_{00}\cap X$ is transformed on $U_{01}\cap X$ into
\begin{align*}
V_{01}= \left(A_0(X_1,X_2)+ \frac{1}{t'} A_1(X_1,X_2) \right) \frac{\partial}{\partial X_1}+\left(B_0(X_1,X_2)+\frac{1}{t'} B_1(X_1,X_2) \right) \frac{\partial}{\partial X_2} + \left( - t'^2 C_0(X_1,X_2)- t' C_1(X_1,X_2) \right)\frac{\partial}{\partial t'}
\end{align*}
We note that since $t'^2\left(1+ X_1^2 + X_2^2 \right)=1$ on $U_{01}\cap X$, we see that $t'$ is a unit, so that $V_{01}$ is still holomorphic on $U_{01}\cap X$. We note that in view of $(\ref{css10})$ $V_{00}$ on $U_{00}\cap X$ is transformed on $U_{10}\cap X$ into
\begin{align*}
V_{10}=&\left( A_0\left( \frac{1}{Y_0}, \frac{Y_2}{Y_0} \right) + \frac{t'}{Y_0} A_1\left( \frac{1}{Y_0}, \frac{Y_2}{Y_0} \right)\right)\left( - Y_0^2 \frac{\partial}{\partial Y_0}- Y_0 Y_2\frac{\partial}{\partial Y_2}- t' Y_0\frac{\partial}{\partial t'} \right)\\
&+\left( B_0\left( \frac{1}{Y_0}, \frac{Y_2}{Y_0} \right) + \frac{t'}{Y_0} B_1\left( \frac{1}{Y_0}, \frac{Y_2}{Y_0} \right)\right) Y_0\frac{\partial}{\partial Y_2} + \left( C_0\left( \frac{1}{Y_0}, \frac{Y_2}{Y_0} \right) + \frac{t'}{Y_0} C_1\left( \frac{1}{Y_0}, \frac{Y_2}{Y_0} \right) \right) Y_0 \frac{\partial}{\partial t'} \\
&=\left( - Y_0^2 A_0\left( \frac{1}{Y_0}, \frac{Y_2}{Y_0}\right) - t' Y_0 A_1\left( \frac{1}{Y_0}, \frac{Y_2}{Y_0} \right) \right) \frac{\partial}{\partial Y_0}\\
&+\left( - Y_0 Y_2 A_0\left( \frac{1}{Y_0}, \frac{Y_2}{Y_0} \right) + Y_0B_0\left( \frac{1}{Y_0}, \frac{Y_2}{Y_0} \right)+ t' \left( B_1\left( \frac{1}{Y_0}, \frac{Y_2}{Y_0} \right)  - Y_2 A_1\left( \frac{1}{Y_0}, \frac{Y_2}{Y_0} \right) \right) \right) \frac{\partial}{\partial Y_2}\\
&+\left(Y_0 C_0\left( \frac{1}{Y_0}, \frac{Y_2}{Y_0} \right) + t' \left(  C_1\left( \frac{1}{Y_0}, \frac{Y_2}{Y_0} \right)  - Y_0 A_0\left( \frac{1}{Y_0}, \frac{Y_2}{Y_0}  \right) \right) - t'^2  A_1\left( \frac{1}{Y_0}, \frac{Y_2}{Y_0} \right) \right) \frac{\partial}{\partial t'}
\end{align*}
Since $t'^2= 1+ Y_0^2 + Y_2^2$ on $U_{10}\cap X$, we have
\begin{align*}
V_{10}=&\left( - Y_0^2 A_0\left( \frac{1}{Y_0}, \frac{Y_2}{Y_0}\right) - t' Y_0 A_1\left( \frac{1}{Y_0}, \frac{Y_2}{Y_0} \right) \right) \frac{\partial}{\partial Y_0}\\
&+\left( - Y_0 Y_2 A_0\left( \frac{1}{Y_0}, \frac{Y_2}{Y_0} \right) + Y_0B_0\left( \frac{1}{Y_0}, \frac{Y_2}{Y_0} \right)+ t' \left( B_1\left( \frac{1}{Y_0}, \frac{Y_2}{Y_0} \right)  - Y_2 A_1\left( \frac{1}{Y_0}, \frac{Y_2}{Y_0} \right) \right) \right) \frac{\partial}{\partial Y_2}\\
&+\left(Y_0 C_0\left( \frac{1}{Y_0}, \frac{Y_2}{Y_0} \right) - \left( 1 + Y_0^2 + Y_2^2 \right)A_1\left( \frac{1}{Y_0}, \frac{Y_2}{Y_0} \right) + t' \left(  C_1\left( \frac{1}{Y_0}, \frac{Y_2}{Y_0} \right)  - Y_0A_0\left( \frac{1}{Y_0}, \frac{Y_2}{Y_0}  \right) \right)  \right) \frac{\partial}{\partial t'}
\end{align*}
which should be holomorphic on $U_{10}\cap X$. We note that in order for
\begin{align*}
 - Y_0 Y_2 A_0\left( \frac{1}{Y_0}, \frac{Y_2}{Y_0} \right) + Y_0B_0\left( \frac{1}{Y_0}, \frac{Y_2}{Y_0} \right) = \sum_{i,j\geq 0}  \frac{- a_{ij} Y_2^{j+1} + b_{ij} Y_2^j }{Y_0^{i+j-1}}
\end{align*}
to be holomorphic, we only need to consider $(i,j)=(0,0),(0,1),(1,0)$ which clearly defines a holomorphic function on $U_{10}\cap X$ and $(0,2),(1,1),(2,0)$ in which case we have
\begin{align*}
\frac{- a_{02} Y_2^3+ b_{02} Y_2^2 - a_{11} Y_2^2 + b_{11} Y_2 - a_{20}Y_2 + b_{20}}{Y_0} \Longrightarrow a_{02}=0, b_{02}= a_{11}, b_{11}= a_{20}, b_{20}=0
\end{align*}
Then we have
\begin{align*}
A_0(X_1,X_2) &= a_{00} + a_{10}X_1 + a_{01} X_2+ a_{20} X_1^2 + a_{11} X_1 X_2\\
B_0(X_1,X_2)&= b_{00} + b_{10}X_1 + b_{01} X_2 + a_{20}X_1X_2 + a_{11} X_2^2
\end{align*}
We note that in order for 
\begin{align*}
B_1\left( \frac{1}{Y_0}, \frac{Y_2}{Y_0} \right)  - Y_2 A_1\left( \frac{1}{Y_0}, \frac{Y_2}{Y_0} \right) = \sum_{i,j\geq 0} \frac{ g_{ij} Y_2^j - f_{ij} Y_2^{j+1} }{Y_0^{i+j}}
\end{align*}
to be holomorphic, we should have
\begin{align*}
\frac{g_{10}- f_{10}Y_2 + g_{01}Y_2 - f_{01} Y_2^2}{Y_0}=0 \Longrightarrow g_{10}=0, g_{01}= f_{10}, f_{01}=0
\end{align*}
Then we have
\begin{align*}
A_1(X_1,X_2)&= f_{00} + f_{10} X_1\\
B_1(X_1,X_2)&= g_{00}+ f_{10} X_2
\end{align*}
On the other hand, in order for
\begin{align*}
Y_0 C_0\left( \frac{1}{Y_0}, \frac{Y_2}{Y_0} \right) - \left( 1 + Y_0^2 + Y_2^2 \right)A_1\left( \frac{1}{Y_0}, \frac{Y_2}{Y_0} \right) = Y_0 C_0\left( \frac{1}{Y_0}, \frac{Y_2}{Y_0} \right) - \left( 1+ Y_0^2 + Y_2^2 \right)\left( f_{00}  + \frac{f_{10}}{Y_0} \right)
\end{align*}
to be holomorphic, we should have
\begin{align*}
\frac{h_{20} + h_{11} Y_2+ h_{02} Y_2^2- f_{10} - f_{10} Y_2^2}{Y_0} = 0 \Longrightarrow h_{20}= f_{10}, h_{11}=0, h_{02}= f_{10}
\end{align*}
Then we have
\begin{align*}
C_0(X_1,X_2)= h_{00} + h_{10} X_1 + h_{01} X_2+ f_{10} X_1^2 + f_{10} X_2^2
\end{align*}
We note that in order for
\begin{align*}
 C_1\left( \frac{1}{Y_0}, \frac{Y_2}{Y_0} \right)  - Y_0A_0\left( \frac{1}{Y_0}, \frac{Y_2}{Y_0}  \right)= C_1\left( \frac{1}{Y_0}, \frac{Y_2}{Y_0} \right) - Y_0\left( a_{00}+ \frac{a_{10} + a_{01} Y_2}{Y_0} + \frac{a_{20}+ a_{11} Y_2}{Y_0^2} \right)
\end{align*}
to be holomorphic, we should have
\begin{align*}
\frac{c_{10}+ c_{01} Y_2 - a_{20}- a_{11} Y_2}{Y_0} \Longrightarrow c_{10}= a_{20}, c_{01}= a_{11}
\end{align*}
Then we have
\begin{align*}
C_1(X_1,X_2)= c_{00} + a_{20} X_1 + a_{11} X_2
\end{align*}
Hence $H^0\left( X, \Theta_Y|_X \right)$ on $U_{00}\cap X$ is given by
\begin{align} \label{css7}
H^0\left( X, \Theta_Y|_X \right) &= \left\{\left( a_{00} + a_{10} X_1 + a_{01} X_2 + a_{20} X_1^2 + a_{11} X_1 X_2 + t \left( f_{00} + f_{10} X_1 \right)\right) \frac{\partial}{\partial X_1}\right. \\
&\left.+\left( b_{00} + b_{10} X_1 + b_{01} X_2 + a_{20}X_1X_2+ a_{11} X_2^2 + t \left(g_{00}+ f_{10} X_2 \right) \right) \frac{\partial}{\partial X_2} \right. \notag \\
&\left.+ \left( h_{00} + h_{10} X_1 + h_{01} X_2 + f_{10} X_1^2 + f_{10} X_2^2 + t \left( c_{00} + a_{20}X_1 + a_{11}X_2 \right) \right) \frac{\partial}{\partial t}\right\} \notag
\end{align}
where $a_{00}, a_{10}, a_{01}, a_{20}, a_{11},  f_{00}, f_{10}, b_{00}, b_{10}, b_{01},  g_{00}, h_{00}, h_{10}, h_{01}, c_{00} \in \mathbb{C}$.

Next we explicitly describe $H^0\left(X, \mathcal{O}_Y(X)|_X \right)=H^0\left(X, \mathcal{N}_{X/Y} \right)$. The exact sequence $0\to \mathcal{O}_Y\to \mathcal{O}_Y(X)\to \mathcal{O}_Y(X)|_X \to 0$ induces the exact sequence
\begin{align}\label{css5}
0\to H^0\left(Y, \mathcal{O}_Y \right) \cong \mathbb{C} \xrightarrow{\cdot \left(1 + X_1^2 + X_2^2 - t^2 \right)} H^0\left( Y,  \mathcal{O}_Y(X) \right)\to H^0\left( X, \mathcal{O}_Y(X)|_X \right)\to 0
\end{align}
In view of $(\ref{cs28})$ $H^0\left( Y, \mathcal{O}_Y(X) \right)$ on $U_{00}$ is given by
\begin{align*}
H^0\left(Y, \mathcal{O}_Y(X) \right)=\left\{a_0+ a_1 X_1 + a_2 X_2 + a_3 X_1^2 + a_4 X_1 X_2 + a_5 X_2^2 +t\left( b_0 + b_1 X_1 + b_2 X_2 \right)+ c_0 t^2\right\}
\end{align*}
where $a_0,...,a_5,b_0,b_1,b_2,c_0\in \mathbb{C}$. Then from $(\ref{css5})$, we see that $H^0\left( X, \mathcal{O}_Y(X) \right)$ on $U_{00}\cap X$ is given by
\begin{align}\label{css6}
H^0\left( X, \mathcal{O}_Y(X)|_X \right)=\left\{ a_0 + a_1 X_1 + a_2 X_2 + a_3 X_1^2 + a_4 X_1 X_2 + a_5 X_2^2+t\left(b_0 + b_1 X_1 + b_2 X_2 \right) \right\}
\end{align}
where $a_0,...,a_5,b_0,b_1,b_2\in \mathbb{C}$. 

By using $(\ref{css7})$ and $(\ref{css6})$ we will explicitly describe $H^0\left(X, \Theta_X \right)$ on $U_{00}\cap X$ by computing
\begin{align*}
\textnormal{ker}\left(H^0\left( X, \Theta_Y|_X \right)\to H^0\left( X, \mathcal{O}_Y(X)|_X \right) \right)
\end{align*}
From $(\ref{css7})$, we see that the image is given by
{\Small{\begin{align*}
&\left[ \left( a_{00} + a_{10} X_1 + a_{01} X_2 + a_{20} X_1^2 + a_{11} X_1 X_2 + t \left( f_{00} + f_{10} X_1 \right)\right) \frac{\partial}{\partial X_1}, X_1^2+ X_2^2 +1 - t^2 \right]|_{t^2= 1+ X_1^2 + X_2^2} \\
&+\left[ \left( b_{00} + b_{10} X_1 + b_{01} X_2 + a_{20}X_1X_2+ a_{11} X_2^2 + t \left(g_{00}+ f_{10} X_2 \right) \right) \frac{\partial}{\partial X_2} , X_1^2 + X_2^1 +1 - t^2 \right] |_{t^2= 1+ X_1^2 + X_2^2} \\
&+\left[  \left( h_{00} + h_{10} X_1 + h_{01} X_2 + f_{10} X_1^2 + f_{10} X_2^2 + t \left( c_{00} + a_{20}X_1 + a_{11}X_2 \right) \right) \frac{\partial}{\partial t} ,  X_1^2 + X_2^2 +1 - t^2 \right] |_{t^2= 1+ X_1^2 + X_2^2} \\
&= -2 c_{00} + \left(2a_{00}- 2a_{20}\right)X_1 + \left(2 b_{00}- 2 a_{11}\right)X_2 + \left(2 a_{10}- 2 c_{00}\right)X_1^2 + \left(2 a_{01} + 2 b_{10}\right)X_1X_2+ \left( 2 b_{01}- 2 c_{00}\right)X_2^2\\
&\,\,\,\,\,\,\,\,\,+t\left( -2 h_{00} + \left(2 f_{00}- 2 h_{10}\right) X_1+ \left(2 g_{00}- 2 h_{01}\right) X_2\right)
\end{align*}}}
Then $\textnormal{ker}\left(H^0\left( X, \Theta_Y|_X \right)\to H^0\left( X, \mathcal{O}_Y(X)|_X \right) \right)$ on $U_{00}\cap X$ is given by
\begin{align}
H^0\left(X, \Theta_X \right)=& \left\{ \left( a_{00}  + a_{01} X_2 + a_{00} X_1^2 + a_{11} X_1 X_2 + t \left( f_{00} + f_{10} X_1 \right)\right) \frac{\partial}{\partial X_1} \right. \\
&\left. +\left( a_{11} - a_{01} X_1  + a_{00}X_1X_2+ a_{11} X_2^2 + t \left(g_{00}+ f_{10} X_2 \right) \right) \frac{\partial}{\partial X_2}\right. \notag \\
&\left. + \left( f_{00} X_1 + g_{00} X_2 + f_{10} X_1^2 + f_{10} X_2^2 + t \left( a_{00}X_1 + a_{11}X_2 \right) \right) \frac{\partial}{\partial t} \right\} \notag
\end{align}
where $a_{00}, a_{01}, a_{11}, f_{00}, f_{10}, g_{00} \in \mathbb{C}$.

\subsubsection{Computation of $H^0\left( X, \Theta_X \right)\to H^0\left( X, \Theta_{\mathcal{F}_{X_0}}^* \otimes \frac{\Theta_X}{\Theta_{\mathcal{F}_{X_0}}} \right)$}\

By using $(\ref{css9})$ and $(\ref{css8})$, we compute $H^0\left( X, \Theta_X \right)\to H^0\left( X, \Theta_{\mathcal{F}_{X_0}}^* \otimes \frac{\Theta_X}{\Theta_{\mathcal{F}_{X_0}}} \right)$. The image is given by
{\small{\begin{align*}
&\left[ \left( a_{00}  + a_{01} X_2 + a_{00} X_1^2 + a_{11} X_1 X_2 + t \left( f_{00} + f_{10} X_1 \right)\right) \frac{\partial}{\partial X_1}, X_2\left( 1+ 2 X_1^2 \right)\frac{\partial}{\partial X_1} + X_1\left(1+2 X_2^2 \right) \frac{\partial}{\partial X_2} + 2 t X_1 X_2 \frac{\partial}{\partial t} \right]\\
&+ \left[ \left( a_{11} - a_{01} X_1  + a_{00}X_1X_2+ a_{11} X_2^2 + t \left(g_{00}+ f_{10} X_2 \right) \right) \frac{\partial}{\partial X_2} , X_2\left( 1+ 2 X_1^2 \right)\frac{\partial}{\partial X_1} + X_1\left(1+2 X_2^2 \right) \frac{\partial}{\partial X_2} + 2 t X_1 X_2 \frac{\partial}{\partial t} \right] \\
&+ \left[ \left( f_{00} X_1 + g_{00} X_2 + f_{10} X_1^2 + f_{10} X_2^2 + t \left( a_{00}X_1 + a_{11}X_2 \right) \right) \frac{\partial}{\partial t}, X_2\left( 1+ 2 X_1^2 \right)\frac{\partial}{\partial X_1} + X_1\left(1+2 X_2^2 \right) \frac{\partial}{\partial X_2} + 2 t X_1 X_2 \frac{\partial}{\partial t} \right]\\
&=\left(a_{11}- 2a_{01}X_1  + a_{11} X_1^2 + 3 a_{00}X_1 X_2 + 2 a_{01} X_1 X_2^2  - 2 a_{01} X_1^3  + 2a_{00}X_1^3X_2+ 2a_{11} X_1^2X_2^2\right) \frac{\partial}{\partial X_1} \\
 &+ t \left( g_{00}  +2g_{00}X_1^2+ 2f_{00}X_1X_2+ 2 f_{10} X_1^2 X_2 \right)\frac{\partial}{\partial X_1}\\
&+  \left(a_{00}  + 2 a_{01} X_2  + 3a_{11} X_1 X_2  +  a_{00} X_2^2 + 2a_{01} X_2^3 +2 a_{00} X_1^2X_2^2 + 2a_{11} X_1 X_2^3  -2 a_{01} X_1^2 X_2 \right)\frac{\partial}{\partial X_2} \\
&+ t  \left( f_{00}+2g_{00}X_1X_2 + 2f_{00} X_2^2+2 f_{10} X_1X_2^2 \right) \frac{\partial}{\partial X_2} \\
&+ \left( g_{00}X_1+ f_{00}X_2 + 2 g_{00} X_1^3 + 2 g_{00} X_1 X_2^2  + 2 f_{00}X_1^2X_2 +2 f_{00} X_2^3 + 2 f_{10} X_1^3 X_2 + 2f_{10} X_1 X_2^3 \right) \frac{\partial }{\partial t} \\
&+ t\left( a_{00}X_2 + a_{11} X_1 - 2 a_{01} X_1^2 + 2_{01} X_2^2  + 2 a_{00}X_1^2 X_2 + 2 a_{11}X_1 X_2^2 \right) \frac{\partial}{\partial t}
\end{align*}}}
We note that from $\left(X_2 + 2 X_1^2 X_2 \right)\frac{\partial}{\partial X_1} + \left( X_1 + 2 X_1 X_2^2\right)\frac{\partial}{\partial X_2} +2 tX_1 X_2 \frac{\partial}{\partial t} \equiv 0 $, we have
\begin{align*}
2 a_{00} X_1^3 X_2 \frac{\partial}{\partial X_1}+ 2a_{00} X_1^2 X_2^2\frac{\partial}{\partial X_2}+ 2a_{00} t X_1^2 X_2 \frac{\partial}{\partial t}  \equiv - a_{00} X_1 X_2 \frac{\partial}{\partial X_1} - a_{00} X_1^2 \frac{\partial}{\partial X_2}\\
2 a_{11} X_1^2 X_2^2 \frac{\partial}{\partial X_1}+  2a_{11} X_1 X_2^3 \frac{\partial}{\partial X_2}+ 2a_{11} t X_1 X_2^2 \frac{\partial}{\partial t}  \equiv - a_{11}  X_2^2 \frac{\partial}{\partial X_1} - a_{11} X_1X_2 \frac{\partial}{\partial X_2}
\end{align*}
and we have
{\small{\begin{align*}
&2 f_{10} X_1^2 X_2 \frac{\partial}{\partial X_1}  + 2 f_{10} X_1 X_2^2\frac{\partial}{\partial X_2} \equiv - f_{10} X_2 \frac{\partial}{\partial X_1} - f_{10} X_1 \frac{\partial}{\partial X_2}- 2 f_{10} t X_1X_2 \frac{\partial}{\partial t}\\
& \Longrightarrow 2 f_{10} tX_1^2 X_2 \frac{\partial}{\partial X_1}  + 2 f_{10} tX_1 X_2^2\frac{\partial}{\partial X_2} \equiv - f_{10}t X_2 \frac{\partial}{\partial X_1} - f_{10} tX_1 \frac{\partial}{\partial X_2}- 2 f_{10} \left(1+ X_1^2+ X_2^2 \right) X_1X_2 \frac{\partial}{\partial t}\\
&\Longrightarrow 2 f_{10} tX_1^2 X_2 \frac{\partial}{\partial X_1}  + 2 f_{10} tX_1 X_2^2\frac{\partial}{\partial X_2} \equiv - f_{10}t X_2 \frac{\partial}{\partial X_1} - f_{10} tX_1 \frac{\partial}{\partial X_2}+ \left(-2 f_{10} X_1 X_2-2 f_{10} X_1^3 X_2 - 2 f_{10} X_1 X_2^3\right) \frac{\partial}{\partial t}
\end{align*}}}
Then the image of $H^0\left( X, \Theta_X \right)\to H^0\left( X, \Theta_{\mathcal{F}_{X_0}}^* \otimes \frac{\Theta_X}{\Theta_{\mathcal{F}_{X_0}} } \right)$ is given by
\begin{align*}
\left(a_{11}- 2a_{01}X_1  + a_{11} X_1^2  + 2 a_{00}X_1 X_2 - a_{11} X_2^2 + 2 a_{01} X_1 X_2^2  - 2 a_{01} X_1^3  \right) \frac{\partial}{\partial X_1} \\
 + t \left( g_{00} - f_{10}X_2  +2g_{00}X_1^2+ 2f_{00}X_1X_2 \right)\frac{\partial}{\partial X_1}\\
  \left(a_{00}  + 2 a_{01} X_2  + 2a_{11} X_1 X_2 - a_{00}X_1^2  +  a_{00} X_2^2 + 2a_{01} X_2^3   -2 a_{01} X_1^2 X_2 \right)\frac{\partial}{\partial X_2} \\
 + t  \left( f_{00} - f_{10} X_1+2g_{00}X_1X_2 + 2f_{00} X_2^2 \right) \frac{\partial}{\partial X_2} \\
 \left( g_{00}X_1+ f_{00}X_2 -2 f_{10} X_1 X_2 + 2 g_{00} X_1^3 + 2 g_{00} X_1 X_2^2  + 2 f_{00}X_1^2X_2 +2 f_{00} X_2^3  \right) \frac{\partial }{\partial t} \\
+ t\left( a_{00}X_2 + a_{11} X_1 - 2 a_{01} X_1^2 + 2 a_{01} X_2^2  \right) \frac{\partial}{\partial t}
\end{align*}
Then from $(\ref{css9})$, we have $\dim_\mathbb{C} \textnormal{ker}\left( H^0\left( X, \Theta_X \right)\to H^0\left( X, \Theta_{\mathcal{F}_{X_0}}^* \otimes \frac{\Theta_X}{\Theta_{\mathcal{F}_{X_0}} } \right) \right)= 0 $, so that we have
\begin{align*}
\dim_\mathbb{C} \mathbb{H}^0 \left( X, \Theta_{\mathcal{F}_{X_0}}^\bullet \right) & =  0\\
\dim_\mathbb{C} \mathbb{H}^1 \left(X, \Theta_{\mathcal{F}_{X_0}}^\bullet \right) & = 15- 6 = 9
\end{align*}

\subsection{Embedded foliated deformations of the invariant submanifold $\left(X, \Theta_{\mathcal{F}_{X_0}} \right)$ in $\left( Y , \Theta_{\tilde{\mathcal{F}}_0^\sharp} \right)$}\

We keep the notations in subsection \ref{cs5}. We study deformations of the invariant submanifold $X$ of $\left( Y, \Theta_{\tilde{\mathcal{F}}_0^\sharp} \right)$ (see section \ref{hs1}). We note that $\mathcal{O}_Y(X)$ defines a foliated line bundle over $\left( Y, \Theta_{\tilde{\mathcal{F}}_0^\sharp}  \right)$ as follows: we note that we have on $U_{00}$ and $U_{01}$, respectively
\begin{align*}
&\left[ X_2 \left(1+2X_1^2\right)\frac{\partial}{\partial X_1} + X_1\left(1+2X_2^2\right)\frac{\partial}{\partial X_2} +2 tX_1X_2\frac{\partial}{\partial t},  1+ X_1^2+ X_2^2 - t^2 \right] = 4X_1X_2\left(1+ X_1^2 + X_2^2 - t^2 \right) \\
&\left[ X_2 \left(1+2X_1^2\right)\frac{\partial}{\partial X_1} + X_1\left(1+2X_2^2\right)\frac{\partial}{\partial X_2} +2 tX_1X_2\frac{\partial}{\partial t}, t'^2\left( 1+ X_1^2 + X_2^2 \right)-1  \right]=0
\end{align*}
and we have on $U_{10}$ and $U_{11}$, respectively
\begin{align*}
&\left[  - Y_2\left( Y_0^2 + 2  \right)\frac{\partial}{\partial Y_0} - Y_0\left(  Y_2^2 -1 \right)\frac{\partial}{\partial Y_2} - t Y_0Y_2\frac{\partial}{\partial t}, 1+ Y_0^2 + Y_2^2- t^2 \right]= -2 Y_0 Y_2 \left( 1+ Y_0^2 + Y_2^2 - t^2 \right) \\
&\left[  - Y_2\left( Y_0^2 + 2  \right)\frac{\partial}{\partial Y_0} - Y_0\left(  Y_2^2 -1 \right)\frac{\partial}{\partial Y_2} + t' Y_0Y_2\frac{\partial}{\partial t'},t'^2\left( 1+ Y_0^2 + Y_2^2 \right) - 1  \right]=  0 
\end{align*}
and we have on $U_{20}$ and $U_{21}$, respectively
\begin{align*}
&\left[ - Z_1\left( 2 + Z_0^2 \right) \frac{\partial}{\partial Z_0} + Z_0 \left( 1 - Z_1^2\right)\frac{\partial}{\partial Z_1} - tZ_0 Z_1\frac{\partial}{\partial t} , 1+ Z_0^2 + Z_1^2 - t^2 \right] = -2 Z_0Z_1\left( 1+ Z_0^2 + Z_1^2 - t^2 \right) \\
&\left[ - Z_1\left( 2 + Z_0^2 \right) \frac{\partial}{\partial Z_0} + Z_0 \left( 1 - Z_1^2\right)\frac{\partial}{\partial Z_1} + t'Z_0 Z_1\frac{\partial}{\partial t'} ,  t'^2\left(1+ Z_0^2 + Z_1^2 \right) - 1 \right] =  0
\end{align*}
Then we see that the foliated line bundle $\left( \mathcal{O}_X(Y), \nabla_{\tilde{\mathcal{F}}_0^\sharp} \right)$ is given by
\begin{align*}
\left( \left\{\phi_{00}, \phi_{01}, \phi_{10}, \phi_{11}, \phi_{20}, \phi_{21} \right\} , \left\{\textnormal{(\ref{cs28})} \right\} \right) \in C^0\left( \mathcal{U}^\sharp, \mathscr{H}om_{\mathcal{O}_Y}\left( \Theta_{\tilde{\mathcal{F}}_0^\sharp} , \mathcal{O}_Y \right) \right) \bigoplus C^1\left( \mathcal{U}^\sharp, \mathcal{O}_Y^* \right)
\end{align*}
where $\phi_{00}(T)= 4 X_1X_2,\phi_{01}(T)=0,\phi_{10}(T)= - 2 Y_0Y_2, \phi_{11}(T)= 0, \phi_{20}(T)=-2 Z_0Z_1, \phi_{21}(T)=0$. Then we have the associated complex of sheaves (see Part V)
\begin{align*}
\mathcal{O}_Y(X)^\bullet: \mathcal{O}_Y(X)\xrightarrow{\nabla_{\tilde{\mathcal{F}}_0^\sharp} } \Theta_{\tilde{\mathcal{F}}_0^\sharp}^* \otimes \mathcal{O}_Y(X) \to 0
\end{align*} 
By restricting to $X$, we have from $(\ref{cs8})$
\begin{align*}
\mathcal{O}_Y(X)|_X^\bullet : \mathcal{O}_Y(X)|_X \xrightarrow{\nabla_{\tilde{\mathcal{F}}_0^\sharp}|_X}  \left(\Theta_{\tilde{\mathcal{F}}_0^\sharp}^* \otimes \mathcal{O}_Y(X)\right)|_X \to 0
\end{align*}
By using $(\ref{css6})$ and $(\ref{cs38})$ we compute $\mathbb{H}^0\left( X, \mathcal{O}_Y(X)|_X^\bullet \right)$. Let us compute $H^0\left( X, \mathcal{O}_Y(X)|_X\right) \xrightarrow{\nabla_{\tilde{\mathcal{F}}_0^\sharp}|_X} H^0\left( X ,   \left(\Theta_{\tilde{\mathcal{F}}_0^\sharp}^* \otimes \mathcal{O}_Y(X)\right)|_X \right) $ on $U_{00}$. The image is given by
{\small{\begin{align*}
&\left[ \left( X_2+ 2 X_1^2 X_2 \right) \frac{\partial}{\partial X_1} +\left( X_1 + 2 X_1 X_2^2 \right) \frac{\partial}{\partial X_2} + 2 t X_1 X_2 \frac{\partial}{\partial t},  a_0 + a_1 X_1+ a_2 X_2 + a_3 X_1^2 + a_4 X_1 X_2+ a_5 X_2^2   + t \left(b_0+ b_1X_1+ b_2 X_2\right)  \right]\\
&- 4 X_1 X_2 \left(  a_0 + a_1 X_1+ a_2 X_2 + a_3 X_1^2 + a_4 X_1 X_2+ a_5 X_2^2  + t \left(b_0+ b_1X_1+ b_2 X_2\right) \right)\\
&=a_2 X_1 +a_1X_2  + a_4 X_1^2 +  \left(2 a_3 +2 a_5 - 4 a_0\right) X_1X_2+ a_4 X_2^2- 2 a_1X_1^2 X_2 -2 a_2 X_1X_2^2 + t \left( b_2 X_1 +b_1 X_2 - 2b_0X_1X_2\right) 
\end{align*}}}
From $(\ref{cs38})$ we see that $\textnormal{ker}\left( H^0\left( X, \mathcal{O}_Y(X)|_X\right) \xrightarrow{\nabla_{\tilde{\mathcal{F}}_0^\sharp}|_X} H^0\left( X ,   \left(\Theta_{\tilde{\mathcal{F}}_0^\sharp}^* \otimes \mathcal{O}_Y(X)\right)|_X \right) \right)$ is given by 
\begin{align}\label{css18}
\mathbb{H}^0\left(X, \mathcal{O}_Y(X)|_X^\bullet \right)=\left\{c_1+ c_2 X_1^2 + \left( 2c_1- c_2 \right)X_2^2 | c_1,c_2\in \mathbb{C} \right\}
\end{align}
Then we note that $(1- c_1)\xi_0^2+ (1-c_2)\xi_1^2+ (1- 2c_1+  c_2)\xi_2^2=0$ defines a family of non-degenerate invariant submanifold of $\left(\mathbb{P}_\mathbb{C}^2, T_0 \right)$ for sufficiently small $c_1,c_2\in \mathbb{C}$. This family has a surjective characteristic map at $(c_1,c_2)=(0,0)$ from $(\ref{css16})$, so that the family is complete at $(c_1,c_2)= (0,0)$. We note that $(1- c_1)\xi_0^2+ (1-c_2)\xi_1^2+ (1- 2c_1+  c_2)\xi_2^2=0$ defines the same foliated line bundle $\left( \mathcal{O}_{\mathbb{P}^2}(2),  \nabla_{\mathcal{F}_0} \right)$ as in $(\ref{css17})$ since we have
{\Small{\begin{align*}
&\left[ X_2\left( 1+ 2 X_1^2\right)\frac{\partial}{\partial X_1}+ X_1\left( 1+ 2 X_2^2 \right)\frac{\partial}{\partial X_2} , (1-c_1)+ (1-c_2)X_1^2 + (1-2c_1+c_2)X_2^2  \right] = 4X_1X_2\left( (1-c_1)+ (1-c_2)X_1^2 + (1-2c_1+c_2)X_2^2  \right) \\
&\left[ - Y_2\left( Y_0^2 + 2\right)\frac{\partial}{\partial Y_0}- Y_0\left(Y_2^2 - 1 \right) \frac{\partial}{\partial Y_2}, (1-c_1)Y_0^2+(1-c_2)+ (1-2c_1+c_2)Y_2^2  \right]=- 2 Y_0Y_2\left( (1-c_1)Y_0^2+(1-c_2)+ (1-2c_1+c_2)Y_2^2   \right)\\
&\left[ - Z_1\left( 2 + Z_0^2 \right)\frac{\partial}{\partial Z_0} + Z_0 \left( 1 - Z_1^2 \right)\frac{\partial}{\partial Z_1}, (1-c_1)Z_0^2+(1-c_2)Z_1^2+ 1-2c_1+c_2 \right]= - 2 Z_0Z_1\left( (1-c_1)Z_0^2+(1-c_2)Z_1^2+ 1-2c_1+c_2 \right)
\end{align*}}}
Then $t^2= (1-c_1)+ (1-c_2)X_1^2 + (1-2c_1+c_2) X_2^2$ on $U_{00}$, $t'^2\left((1-c_1)+ (1-c_2)X_1^2 + (1-2c_1+c_2) X_2^2 \right)=1$ on $U_{01}$, $t^2=(1-c_1)Y_0^2+(1-c_2)+ (1-2c_1+c_2)Y_2^2$ on $U_{10}$, $t'^2\left( (1-c_1)Y_0^2+(1-c_2)+ (1-2c_1+c_2)Y_2^2  \right)=1$ on $U_{11}$, $t^2= (1-c_1)Z_0^2+ (1-c_2)Z_1^2 + 1 - 2c_1 + c_2$ on $U_{21}$, and $t'^2 \left(   (1-c_1)Z_0^2+ (1-c_2)Z_1^2 + 1 - 2c_1 + c_2 \right)=1$ on $U_{21}$ defines foliated double coverings over $\left( \mathbb{P}_\mathbb{C}^2, T_0 \right)$, and they define a family of non-degenerate invariant complex submanifolds of $\left( Y=\mathcal{O}_{\mathbb{P}_\mathbb{C}^2}(1)^\sharp, \Theta_{\tilde{\mathcal{F}}_0^\sharp} \right)$ over a sufficiently small neighborhood of $(0,0)\in \mathbb{C}^2$. Moreover from $(\ref{css18})$, the characteristic map of the family is an isomorphism at $(c_1,c_2)=(0,0)$, so that the family is complete at $(c_1,c_2)=(0,0)$ by Theorem of completeness for deformations of non-degenerate invariant submanifold (see Part VI). Hence any embedded deformations of $X$ in $\left( Y, \Theta_{\tilde{\mathcal{F}}_0^\sharp} \right)$ is induced from the family. Hence we can conclude that any small embedded foliated deformations of $X$ in $\left(Y, \Theta_{\tilde{\mathcal{F}}_0^\sharp} \right)$ are foliated double coverings induced by some embeded deformations of $C$ in $\left( \mathbb{P}_\mathbb{C}^2, T_0 \right)$.

\subsection{Simultaneous deformations of the invariant submanifold $C\subset \left( \mathbb{P}_\mathbb{C}^2, T_0 \right)$}\

We keep the notations in subsection \ref{cs5}. We assume the knowledge of simultaneous deformations of foliated complex manifolds defined by locally free sheaves of tangent sheaves and their non-degenerate invariant complex submanifolds in Part VI. Simultaneous deformations of $C\subset \left( \mathbb{P}_\mathbb{C}^2, T_0 \right)$ are controlled by the following complex of sheaves (see Part VI)
\begin{align*}
\Theta_{\mathbb{P}_\mathbb{C}^2}^\bullet : \Theta_{\mathbb{P}_\mathbb{C}^2} \to \left(\Theta_{\mathcal{F}_0}^*\otimes \frac{\Theta_{\mathbb{P}_\mathbb{C}^2}}{\Theta_{\mathcal{F}_0}}\right) \bigoplus i_*\left( \mathcal{N}_{C/\mathbb{P}_\mathbb{C}^2} \right) \to i_*\left(\Theta_{\mathcal{F}_0}^* \otimes \mathcal{N}_{C/\mathbb{P}_\mathbb{C}^2} \right) \to 0
\end{align*}
Equivalently
\begin{align}\label{css20}
\Theta_{\mathbb{P}_\mathbb{C}^2}^\bullet: \Theta_{\mathbb{P}_\mathbb{C}^2}\to  \left(\frac{\Theta_{\mathbb{P}_\mathbb{C}^2}}{\Theta_{\mathcal{F}_0}}\otimes \mathcal{O}_{\mathbb{P}_\mathbb{C}^2}(1) \right)\bigoplus i_* \mathcal{O}_C(2) \to i_* \mathcal{O}_C(3) \to 0
\end{align}
where $i:C \hookrightarrow  \mathbb{P}_\mathbb{C}^2$ is the embedding. We will denote the $i$-th cohomology group of $(\ref{css20})$ by $\mathbb{H}^i\left( \mathbb{P}_\mathbb{C}^2, \Theta_{\mathbb{P}_\mathbb{C}^2}^\bullet  \right)$.\footnote{We remark that we use different notations for the cohomology groups in Part VI. We use $\mathbb{H}^i\left(\mathbb{P}_\mathbb{C}^2, \left( C, \mathbb{P}_\mathbb{C}^2, \Theta_{\mathcal{F}_0} \right)^\bullet \right)$ in Part VI.} Since higher cohomology groups of each term vanish, we can compute $\mathbb{H}^i\left(\mathbb{P}_\mathbb{C}^2, \Theta_{\mathbb{P}_\mathbb{C}^2} \right)$ by using the following complex
\begin{align*}
H^0\left(\mathbb{P}_\mathbb{C}^2, \Theta_{\mathbb{P}_\mathbb{C}^2} \right) \to H^0\left( \Theta_{\mathcal{F}_0}^*\otimes \frac{\Theta_{\mathbb{P}_\mathbb{C}^2}}{\Theta_{\mathcal{F}_0}} \right) \bigoplus H^0\left( C, \mathcal{O}_C(2) \right) \to H^0\left( C, \mathcal{O}_C(3) \right) \to 0
\end{align*}

We note that we have the description $H^0\left(\mathbb{P}_\mathbb{C}^2, \Theta_{\mathbb{P}_\mathbb{C}^2} \right)$ from $\textnormal{(\ref{ap1})}$ and the description $H^0\left( C, \mathcal{O}_C(2) \right)$ from $(\ref{css21})$ and the description $H^0\left( C,\mathcal{O}_C(3) \right)$ from $(\ref{css23})$. It remains to  explicitly describe  $H^0\left( \mathbb{P}_\mathbb{C}^2,   \Theta_{\mathcal{F}_0}^*\otimes \frac{\Theta_{\mathbb{P}_\mathbb{C}^2}}{\Theta_{\mathcal{F}_0}} \right)$. From the exact sequence $0\to \mathcal{O}_{\mathbb{P}_\mathbb{C}^2}\to \Theta_{\mathbb{P}_\mathbb{C}} (1) \to  \Theta_{\mathcal{F}_0}^* \otimes \frac{\Theta_{\mathbb{P}_\mathbb{C}^2}}{\Theta_{\mathcal{F}_0}}\to 0$, we have the exact sequence
\begin{align*}
0 \to H^0\left( \mathbb{P}_\mathbb{C}^2,  \mathcal{O}_{\mathbb{P}_\mathbb{C}^2} \right)\cong \mathbb{C} &\to H^0\left( \mathbb{P}_\mathbb{C}^2, \Theta_{\mathbb{P}_\mathbb{C}^2} (1) \right)\to H^0\left(  \mathbb{P}_\mathbb{C}^2,   \Theta_{\mathcal{F}_0}^*\otimes \frac{\Theta_{\mathbb{P}_\mathbb{C}^2}}{\Theta_{\mathcal{F}_0}} \right)\to 0 \\
 1 &\mapsto T_0=X_2\left(1 + 2 X_1^2 \right)\frac{\partial}{\partial X_1} + X_1\left( 1+ 2 X_2^2 \right)\frac{\partial}{\partial X_2}
\end{align*}
Since $X_1^2X_2\frac{\partial}{\partial X_1}+ X_1X_2^2\frac{\partial}{\partial X_2}\equiv - \frac{1}{2}X_2\frac{\partial}{\partial X_2} - \frac{1}{2} X_1\frac{\partial}{\partial X_2}$ in $H^0\left(  \mathbb{P}_\mathbb{C}^2,   \Theta_{\mathcal{F}_0}^*\otimes \frac{\Theta_{\mathbb{P}_\mathbb{C}^2}}{\Theta_{\mathcal{F}_0}}\right)$, we see that from $\textnormal{(\ref{p3})}$

\begin{align}\label{css24}
H^0\left(  \mathbb{P}_\mathbb{C}^2,   \Theta_{\mathcal{F}_0}^*\otimes \frac{\Theta_{\mathbb{P}_\mathbb{C}^2}}{\Theta_{\mathcal{F}_0}} \right) =&\left\{\left( t_1+ t_2 X_2+ t_3 X_1+ t_4 X_2^2+ t_5 X_1X_2+ t_6 X_1^2+ t_7 X_1X_2^2+ t_9 X_1^3 \right)\frac{\partial}{\partial X_1}\right.\\
&\left. + \left( t_{10} + t_{11}X_2+ t_{12}X_1+ t_{13} X_2^2 + t_{14}X_1X_2+ t_{15}X_1^2+ t_7 X_2^3+ t_9 X_1^2 X_2 \right) \frac{\partial}{\partial X_2} \right\} \notag
\end{align}
where $t_1,...,t_7, t_9,...,t_{15}\in \mathbb{C}$. 

First we compute $H^0\left( \mathbb{P}_\mathbb{C}^2 , \Theta_{\mathcal{F}_0}^*\otimes \frac{\Theta_{\mathbb{P}_\mathbb{C}^2}}{\Theta_{\mathcal{F}_0}} \right) \bigoplus H^0\left( C, \mathcal{O}_C(2) \right) \to H^0\left( C, \mathcal{O}_C(3) \right)$. From $(\ref{css24})$ and $(\ref{css21})$, the image is given by (see Part VI)
{\Small{\begin{align*}\\
&\left[ \left( t_1+ t_2 X_2+ t_3 X_1+ t_4 X_2^2+ t_5 X_1X_2+ t_6 X_1^2+ t_7 X_1X_2^2+ t_9 X_1^3 \right)\frac{\partial}{\partial X_1}, 1+ X_1^2 + X_2^2 \right]|_{1+ X_1^2 + X_2^2 =0}\\
&+\left[ \left( t_{10} + t_{11}X_2+ t_{12}X_1+ t_{13} X_2^2 + t_{14}X_1X_2+ t_{15}X_1^2+ t_7 X_2^3+ t_9 X_1^2 X_2 \right) \frac{\partial}{\partial X_2}, 1+ X_1^2 + X_2^2 \right]|_{1+ X_1^2 + X_2^2=0}\\
&+\left[ X_2 \left(1+2X_1^2\right)\frac{\partial}{\partial X_1} + X_1\left(1+2X_2^2\right)\frac{\partial}{\partial X_2} ,a_2 X_1 + a_3 X_2 + a_4 X_1^2 + a_5 X_1X_2+ a_6 X_2^2\right] - 4 X_1 X_2\left(a_2 X_1 + a_3 X_2 + a_4 X_1^2 + a_5 X_1X_2+ a_6 X_2^2 \right)\\
&=\left(2t_1- 2 t_6+ a_3\right)X_1+ \left(2t_2+2 t_{12}+ 2 a_4 + 2a_6\right) X_1X_2+ \left(2t_3-2 t_9+a_5\right) X_1^2+ \left(2t_4- 2 t_6+ 2 t_{14} - 2a_3\right) X_1X_2^2  \\
 &+ \left(2 t_{10} - 2 t_{13}+ a_2\right) X_2 + \left(2 t_{11}- 2 t_7 + a_5\right)X_2^2 + \left(2t_{15} - 2 t_{13} + 2 t_5 - 2 a_2\right) X_1^2X_2
\end{align*}}}
Then in view of $(\ref{css23})$ we see that $\textnormal{ker}\left( H^0\left( \Theta_{\mathcal{F}_0}^*\otimes \frac{\Theta_{\mathbb{P}_\mathbb{C}^2}}{\Theta_{\mathcal{F}_0}} \right) \bigoplus H^0\left( C, \mathcal{O}_C(2) \right) \to H^0\left( C, \mathcal{O}_C(3) \right) \right)$ is given by
\begin{align}\label{css27}
&\left( t_1+ t_2 X_2+ t_3 X_1+ t_4 X_2^2+ t_5 X_1X_2+ t_6 X_1^2+ t_7 X_1X_2^2+ t_9 X_1^3 \right)\frac{\partial}{\partial X_1}\\
&+ \left( t_{10} + \left(t_3+ t_7- t_9\right)X_2+ t_{12}X_1+ t_{13} X_2^2 + \left(- 2t_1-t_4+ 3 t_6\right)X_1X_2+ \left( 3 t_{13}- 2t_{10} - t_5\right) X_1^2+ t_7 X_2^3+ t_9 X_1^2 X_2 \right) \frac{\partial}{\partial X_2} \notag \\
&+ \left(2 t_{13} - 2 t_{10}\right)X_1 + \left(- 2t_1+ 2t_6\right)X_2 + a_4 X_1^2 + \left(-2t_3+ 2t_9\right)X_1X_2 + \left(- t_2 -t_{12} - a_4\right)X_2^2 \notag
\end{align}
where $t_1,t_2,t_3,t_4,t_5,t_6,t_7, t_9, t_{10}, t_{12}, t_{13}, a_4 \in \mathbb{C}$, so that 
\begin{align*}
\dim_\mathbb{C} \mathbb{H}^2\left( \mathbb{P}_\mathbb{C}^2, \Theta_{\mathbb{P}_\mathbb{C}^2}^\bullet \right)= \dim_\mathbb{C} \textnormal{coker}\left(  H^0\left( \Theta_{\mathcal{F}_0}^*\otimes \frac{\Theta_{\mathbb{P}_\mathbb{C}^2}}{\Theta_{\mathcal{F}_0}} \right) \bigoplus H^0\left( C, \mathcal{O}_C(2) \right) \to H^0\left( C, \mathcal{O}_C(3) \right)  \right) = 7- (19-12)=0
\end{align*}
Since $\mathbb{H}^2\left( \mathbb{P}_\mathbb{C}^2, \Theta_{\mathbb{P}_\mathbb{C}^2}^\bullet \right)$ is an obstruction space for the functor $\textnormal{FDef}_{ C\subset \left( \mathbb{P}_\mathbb{C}^2, T_0 \right)}$ of simultaneous deformations of over Artin rings (see Part VI), we see that the functor $\textnormal{FDef}_{ C\subset \left( \mathbb{P}_\mathbb{C}^2, T_0 \right)}$ is unobstructed.

Next we compute $H^0\left( \mathbb{P}_\mathbb{C}^2 , \Theta_{\mathbb{P}_\mathbb{C}^2} \right) \to H^0\left( \Theta_{\mathcal{F}_0}^*\otimes \frac{\Theta_{\mathbb{P}_\mathbb{C}^2}}{\Theta_{\mathcal{F}_0}} \right) \bigoplus H^0\left( C, \mathcal{O}_C(2) \right)$. From $\textnormal{(\ref{ap1})}$, the image is given by (see Part VI)
{\Small{\begin{align}
&\left[  \left( b_1 + b_2 X_1 + b_3 X_2 + b_7 X_1^2 + b_8 X_1X_2 \right) \frac{\partial}{\partial X_1}+ \left( b_4 + b_5 X_1 + b_6 X_2 + b_7 X_1X_2 + b_8 X_2^2 \right)\frac{\partial}{\partial X_2} , X_2\left(1 + 2 X_1^2 \right)\frac{\partial}{\partial X_1} + X_1\left( 1+ 2 X_2^2 \right)\frac{\partial}{\partial X_2} \right]\notag\\
&+\left[  \left( b_1 + b_2 X_1 + b_3 X_2 + b_7 X_1^2 + b_8 X_1X_2 \right) \frac{\partial}{\partial X_1}+ \left( b_4 + b_5 X_1 + b_6 X_2 + b_7 X_1X_2 + b_8 X_2^2 \right)\frac{\partial}{\partial X_2}, 1+ X_1^2 + X_2^2 \right]|_{1+X_1^2 + X_2^2=0} \notag\\
&\equiv \left( b_4 + \left(b_5- b_3\right) X_1 + \left(b_6 - b_2\right)X_2 + \left(2b_4 - b_8\right)X_1^2   + \left(4b_1 - 2b_7\right)X_1X_2 - b_8 X_2^2  + 2b_5 X_1^3 + \left(2b_6 +2b_2\right)X_1^2 X_2 + 2 b_3 X_1 X_2^2 \right)\frac{\partial}{\partial X_1}\label{css26}\\
&+\left( b_1 + \left(b_2- b_6 \right)X_1 + \left(b_3 - b_5\right)X_2  + \left(4 b_4 - 2 b_8\right)X_1X_2 - b_7 X_1^2 + \left(2b_1 - b_7\right) X_2^2 + 2b_5 X_1^2 X_2 + \left(2b_6 +2 b_2\right)X_1X_2^2  +2 b_3 X_2^3  \right)\frac{\partial}{\partial X_2} \notag \\
&+\left(2b_1- 2 b_7 \right)X_1 + \left(2b_4 - 2 b_8\right) X_2 +2 b_2 X_1^2  + \left(2b_3 + 2 b_5 \right)X_1X_2 +2 b_6 X_2^2\notag
\end{align}}}
Then in view of $(\ref{css24})$ and $(\ref{css21})$, we see that
\begin{align*}
\mathbb{H}^0\left( \mathbb{P}_\mathbb{C}^2, \Theta_{\mathbb{P}_\mathbb{C}^2}^\bullet \right)=\textnormal{ker}\left( H^0\left(\mathbb{P}_\mathbb{C}^2, \Theta_{\mathbb{P}_\mathbb{C}^2} \right)\to H^0\left( \Theta_{\mathcal{F}_0}^*\otimes \frac{\Theta_{\mathbb{P}_\mathbb{C}^2}}{\Theta_{\mathcal{F}_0}} \right) \bigoplus H^0\left( C, \mathcal{O}_C(2) \right)  \right) =0
\end{align*}
and so from $(\ref{css27})$ we have
\begin{align*}
\dim_\mathbb{C} \mathbb{H}^1\left( \mathbb{P}_\mathbb{C}^2, \Theta_{\mathbb{P}_\mathbb{C}^2}^\bullet \right) = 12- 8 = 4
\end{align*}
By comparing $(\ref{css27})$ and $(\ref{css26})$, we see that
\begin{align}
 \mathbb{H}^1\left( \mathbb{P}_\mathbb{C}^2, \Theta_{\mathbb{P}_\mathbb{C}^2}^\bullet \right)  &\cong \frac{\textnormal{ker}\left(  H^0\left( \Theta_{\mathcal{F}_0}^*\otimes \frac{\Theta_{\mathbb{P}_\mathbb{C}^2}}{\Theta_{\mathcal{F}_0}} \right) \bigoplus H^0\left( C, \mathcal{O}_C(2) \right) \to H^0\left( C, \mathcal{O}_C(3) \right) \right)}{\textnormal{im}\left( H^0\left( \mathbb{P}_\mathbb{C}^2, \Theta_{\mathbb{P}_\mathbb{C}^2} \right)\to  H^0\left( \Theta_{\mathcal{F}_0}^*\otimes \frac{\Theta_{\mathbb{P}_\mathbb{C}^2}}{\Theta_{\mathcal{F}_0}} \right) \bigoplus H^0\left( C, \mathcal{O}_C(2) \right) \right)} \notag \\
 &\cong \textnormal{Span}_\mathbb{C} \left\langle X_1 \frac{\partial}{\partial X_1} + X_2 \frac{\partial}{\partial X_2} - 2 X_1X_2, X_2^2 \frac{\partial}{\partial X_1} - X_1 X_2 \frac{\partial}{\partial X_2},X_1X_2 \frac{\partial}{\partial X_1} - X_1^2 \frac{\partial}{\partial X_2} , X_1^2 - X_2^2 \right\rangle  \label{css28}
\end{align}
Then we have simultaneous deformation $\mathcal{C}\subset \left(  \mathbb{P}_\mathbb{C}^2\times \textnormal{Spec}\left( A\right) , T \right)$ of $C\subset \left( \mathbb{P}_\mathbb{C}^2, T_0 \right)$ over $A=\frac{\mathbb{C}\left[t_1,t_2,t_3,t_4\right]}{ \left(t_1^2, t_2^2, t_3^2, t_4^2, t_1t_2,t_1t_3,t_1t_4,t_2t_3,t_2t_4,t_3t_4\right)}$ where $T$ is given by
\begin{align}\label{css29}
T=\left( X_2+ 2 X_1^2 X_2 - t_1 X_1 - t_3X_2^2 - t_4 X_1X_2  \right)\frac{\partial}{\partial X_1} + \left( X_2 + 2X_1 X_2^2 - t_1 X_2 + t_3 X_1X_2+ t_4 X_1^2 \right)\frac{\partial}{\partial X_2}
\end{align}
and $\mathcal{C}$ is given by
\begin{align}\label{ccs1}
\xi_0^2+ \xi_1^2 + \xi_2^2 +2t_1 \xi_1\xi_2 -t_2\left(\xi_1^2- \xi_2^2 \right) = 0
\end{align}

The author has not yet established theorem of existence and completeness for simultaneous deformations of compact foliated complex manifolds defined by locally free subsheaves of tangent sheaves, and their non-degenerate invariant complex submanifolds, though he believes that they are true. Under the assumption that theorem of existence and completeness holds, since $\mathbb{H}^2\left( \mathbb{P}_\mathbb{C}^2, \Theta_{\mathbb{P}_\mathbb{C}^2}^\bullet \right)=0$, we should have a foliated analytic family of simultaneous deformation $C\subset \left( \mathbb{P}_\mathbb{C}^2, T_0 \right)$ such that the characteristic map is an isomorphism onto $\mathbb{H}^1\left( \mathbb{P}_\mathbb{C}^2, \Theta_{\mathbb{P}_\mathbb{C}^2}^\bullet \right)$ given by $(\ref{css28})$ at the distinguished point. However the author have difficulties with directly constructing a foliated analytic family of simultaneous deformations of $C\subset \left( \mathbb{P}_\mathbb{C}^2, T_0 \right)$ extending $(\ref{css29})$ and $(\ref{ccs1})$.

\subsection{Simultaneous deformations of the invariant submanifold $\left( X, \Theta_{\mathcal{F}_{X_0}} \right)\subset \left(  Y , \Theta_{\tilde{\mathcal{F}}_0^\sharp} \right)$}\

We keep the notations in subsection \ref{cs5}. We assume the knowledge of simultaneous deformations of foliated complex manifolds defined by locally free subsheaves of tangent sheaves and their non-degenerate invariant complex submanifolds in Part VI. We note that simultaneous deformations of $X\subset \left( Y, \Theta_{\tilde{\mathcal{F}}_0^\sharp} \right)$ are controlled by the following complex of sheaves (see Part VI)
\begin{align}\label{css25}
\Theta_Y^\bullet : \Theta_Y \to  \left( \Theta_{\tilde{\mathcal{F}}_0^\sharp}^* \otimes \frac{\Theta_Y}{\Theta_{\tilde{\mathcal{F}}_0^\sharp}} \right) \bigoplus i_*\left(\mathcal{O}_Y(X)|_X  \right) \to i_*\left( \Theta_{\mathcal{F}_{X_0}}\otimes \mathcal{O}_Y(X)|_X \right) \to 0
\end{align}
where $i:X\hookrightarrow Y$ is the embedding. We will denote the $i$-th cohomology group of $(\ref{css25})$ by $\mathbb{H}^i\left( Y, \Theta_Y^\bullet \right)$.\footnote{We remark that we use different notations for the cohomology groups in Part VI. We use $\mathbb{H}^i\left(Y, \left( X, Y, \Theta_{\tilde{\mathcal{F}}_0^\sharp} \right)^\bullet \right)$ in Part VI.} Since higher cohomology groups of each term in $(\ref{css25})$ vanish, we can compute $\mathbb{H}^i\left( Y, \Theta_Y^\bullet \right)$ by the following complex
\begin{align*}
H^0\left( Y, \Theta_Y \right)\to H^0\left( Y, \Theta_{\tilde{\mathcal{F}}_0^\sharp}^* \otimes \frac{\Theta_Y}{\Theta_{\tilde{\mathcal{F}}_0^\sharp}} \right) \bigoplus H^0\left( X, \mathcal{O}_X(Y) \right) \to H^0\left( X, \left( \Theta_{\tilde{\mathcal{F}}_0^\sharp}^* \otimes \mathcal{O}_Y(X) \right)|_X \right) \to 0
\end{align*}
We note that we have the description $H^0\left( Y, \Theta_Y \right)$ from $(\ref{css13})$ and the description $H^0\left( Y, \Theta_{\tilde{\mathcal{F}}_0^\sharp}^* \otimes \frac{\Theta_Y}{\Theta_{\tilde{\mathcal{F}}_0^\sharp}} \right)$ from $(\ref{css12})$ and the description $H^0\left( X, \mathcal{O}_X(Y) \right)$ from $(\ref{css6})$ and the description $ H^0\left( X, \left( \Theta_{\tilde{\mathcal{F}}_0^\sharp}^* \otimes \mathcal{O}_Y(X) \right)|_X \right)$  from $(\ref{cs38})$. 

First we compute $H^0\left( Y, \Theta_{\tilde{\mathcal{F}}_0^\sharp}^* \otimes \frac{\Theta_Y}{\Theta_{\tilde{\mathcal{F}}_0^\sharp}} \right) \bigoplus H^0\left( X, \mathcal{O}_X(Y) \right) \to H^0\left( X, \left( \Theta_{\tilde{\mathcal{F}}_0^\sharp}^* \otimes \mathcal{O}_Y(X) \right)|_X \right)$. From $(\ref{css12})$ and $(\ref{css6})$, the image is given by (see Part VI)
{\Small{\begin{align*}
&\left[ \left(  b_0 + b_1 X_1 + b_2 X_2 + b_3 X_1^2 + b_4 X_1X_2 + b_5 X_2^2 + b_6 X_1^3 + b_7 X_1 X_2^2\right)\frac{\partial}{\partial X_1}, 1+ X_1^2 + X_2^2 - t^2 \right]|_{t^2=1+X_1^2+X_2^2}\\
&+\left[ \left( b_8 + b_9 X_1 + b_{10}X_2 + b_{11}X_1^2 + a_{12}X_1X_2 + b_{13}X_2^2 + b_6 X_1^2 X_2 + b_7 X_2^3 \right) \frac{\partial}{\partial X_2} , 1+ X_1^2 + X_2^2 - t^2\right]|_{t^2=1+ X_1^2 + X_2^2}\\
&+\left[ \left(  b_{14} + b_{15}X_1 + b_{16}X_2 + b_{17}X_1^2 + b_{18}X_1X_2 + b_{19}X_2^2 + \left( b_{20} + b_{21}X_1+ b_{22}X_2  + b_6 X_1^2 + b_7 X_2^2\right)t  + b_{23}t^2 \right) \frac{\partial}{\partial t} , 1 + X_1^2 + X_2^2 - t^2 \right]|_{ t^2 = 1 + X_1^2 + X_2^2} \\
&+\left[ \left( X_2+ 2 X_1^2 X_2 \right) \frac{\partial}{\partial X_1} +\left( X_1 + 2 X_1 X_2^2 \right) \frac{\partial}{\partial X_2} + 2 t X_1 X_2 \frac{\partial}{\partial t},  a_0 + a_1 X_1+ a_2 X_2 + a_3 X_1^2 + a_4 X_1 X_2+ a_5 X_2^2   + t \left(b_0+ b_1X_1+ b_2 X_2\right)  \right]\\
&- 4 X_1 X_2 \left(  a_0 + a_1 X_1+ a_2 X_2 + a_3 X_1^2 + a_4 X_1 X_2+ a_5 X_2^2  + t \left(b_0+ b_1X_1+ b_2 X_2\right) \right)\\
&= -2b_{20}+ \left( 2b_0 - 2 b_{21}+ a_2 \right)X_1+ \left(2 b_1 -2 b_6- 2 b_{20}+ a_4 \right) X_1^2 + \left(2b_2+ 2 b_9 + 2 a_3+ 2 a_5 - 4 a_0\right) X_1X_2 + \left( 2b_3 -2 b_{21}\right) X_1^3 \\
& + \left(2 b_4- 2 b_{22}+2 b_{11}- 2a_1\right) X_1^2X_2 + \left(2b_5 - 2 b_{21} + 2 b_{12} - 2 a_2\right) X_1 X_2^2   + \left(2b_8 - 2 b_{22}+ a_1\right) X_2 + \left( 2 b_{10} - 2 b_7 - 2 b_{20} + a_4\right) X_2^2   + \left( 2b_{13} - 2 b_{22}\right)X_2^3  \\
&+ t\left( -2b_{14}- 2 b_{23} + (d_2 -2b_{15})X_1+ \left(d_1 -2 b_{16}\right)X_2 - \left(2 b_{17}+ 2b_{23}\right)X_1^2 - \left(2 b_{18}+ 2 d_0\right)X_1X_2 -\left(2 b_{19}+ 2 b_{23}\right)X_2^2 \right)
\end{align*}}}
Then in view of $(\ref{cs38})$ we see that $\textnormal{ker}\left( H^0\left( Y, \Theta_{\tilde{\mathcal{F}}_0^\sharp}^* \otimes \frac{\Theta_Y}{\Theta_{\tilde{\mathcal{F}}_0^\sharp}} \right) \bigoplus H^0\left( X, \mathcal{O}_X(Y) \right) \to H^0\left( X, \left( \Theta_{\tilde{\mathcal{F}}_0^\sharp}^* \otimes \mathcal{O}_Y(X) \right)|_X \right) \right)$ is given by
{\small{\begin{align} \label{css32}
&\left( b_0 + b_1 X_1 + b_2 X_2 + b_3 X_1^2 + b_4 X_1X_2 + b_5 X_2^2 + b_6 X_1^3 + b_7 X_1 X_2^2 \right) \frac{\partial}{\partial X_1} \\
&+ \left( b_8 + b_9 X_1  + \left( b_7-b_6+ b_1 \right)X_2 + \left( 3 b_{13}- 2 b_8 - b_4 \right)X_1^2 + \left( 3 b_3- 2b_0 - b_5 \right)X_1X_2 + b_{13}X_2^2 + b_6 X_1^2 X_2 + b_7 X_2^3 \right)\frac{\partial}{\partial X_2}  \notag \\
& + \left( b_{17} + b_{15}X_1 + b_{16}X_2 + b_{17}X_1^2 + b_{18}X_1X_2 + b_{17}X_2^2 + \left( b_3X_1+ b_{13}X_2  + b_6 X_1^2 + b_7 X_2^2\right)t   - b_{17} t^2 \right)\frac{\partial}{\partial t} \notag \\
&+a_0+ \left(2b_{13}- 2 b_8\right)X_1 + \left(2b_3- 2b_0\right) X_2 + a_3 X_1^2 + \left(2b_6- 2b_1\right) X_1 X_2 + \left(2 a_0- a_3- b_2 - b_9\right)X_2^2 +t\left( - b_{18} + 2 b_{16}X_1 + 2b_{15} X_2 \right) \notag
\end{align}}}
where $b_0,b_1,b_2,b_3,b_4,b_5,b_6,b_7,b_8,b_9,b_{13},b_{15},b_{16},b_{17},b_{18},a_0,a_3 \in \mathbb{C}$, so that
{\small{\begin{align*}
\dim_\mathbb{C} \mathbb{H}^2\left( Y, \Theta_Y^\bullet \right)= \dim_\mathbb{C} \textnormal{coker}\left(H^0\left( Y, \Theta_{\tilde{\mathcal{F}}_0^\sharp}^* \otimes \frac{\Theta_Y}{\Theta_{\tilde{\mathcal{F}}_0^\sharp}} \right) \bigoplus H^0\left( X, \mathcal{O}_X(Y) \right) \to H^0\left( X, \left( \Theta_{\tilde{\mathcal{F}}_0^\sharp}^* \otimes \mathcal{O}_Y(X) \right)|_X \right) \right) = 16 -( 33-17 ) = 0
\end{align*}}}
Since $\mathbb{H}^2\left( Y, \Theta_Y^\bullet \right)$ is an obstruction space for the functor $\textnormal{FDef}_{X\subset \left( Y, \Theta_{\tilde{\mathcal{F}}_0^\sharp} \right)}$ of simultaneous deformations of $X\subset \left( Y, \Theta_{\tilde{\mathcal{F}}_0^\sharp} \right)$ over Artin rings (see Part VI), we see that the functor $\textnormal{FDef}_{X\subset \left( Y, \Theta_{\tilde{\mathcal{F}}_0^\sharp} \right)}$ is unobstructed.

Next we compute $H^0\left( Y, \Theta_Y \right)\to H^0\left( Y, \Theta_{\tilde{\mathcal{F}}_0^\sharp}^* \otimes \frac{\Theta_Y}{\Theta_{\tilde{\mathcal{F}}_0^\sharp}} \right) \bigoplus H^0\left( X, \mathcal{O}_X(Y) \right)$. From $(\ref{css13})$ and from $(\ref{css30})$, the image is given by (see Part VI)
{\Small{\begin{align}
&\left[ \left( A_0 + A_1 X_1 + A_2 X_2 + A_3 X_1X_2 + A_4 X_1^2\right)\frac{\partial}{\partial X_1},  X_2\left(1+2X_1^2\right)\frac{\partial}{\partial X_1} + X_1\left(1+2X_2^2\right)\frac{\partial}{\partial X_2} +2 tX_1X_2\frac{\partial}{\partial t}\right] \notag \\
&+\left[ \left( B_0 + B_1 X_1 + B_2 X_2 + A_3 X_2^2  + A_4 X_1 X_2 \right)\frac{\partial}{\partial X_2} , X_2\left(1+2X_1^2\right)\frac{\partial}{\partial X_1} + X_1\left(1+2X_2^2\right)\frac{\partial}{\partial X_2} +2 tX_1X_2\frac{\partial}{\partial t} \right] \notag \\
&+ \left[ \left( \left( C_0 + C_1 X_1 + C_2 X_2 \right) + t \left(C_3 + A_3 X_2 + A_4 X_1 \right)\right)\frac{\partial}{\partial t}, X_2\left(1+2X_1^2\right)\frac{\partial}{\partial X_1} + X_1\left(1+2X_2^2\right)\frac{\partial}{\partial X_2} +2 tX_1X_2\frac{\partial}{\partial t}  \right] \notag \\
&+ \left[  \left( A_0 + A_1 X_1 + A_2 X_2 + A_3 X_1X_2 + A_4 X_1^2\right)\frac{\partial}{\partial X_1} +  \left( B_0 + B_1 X_1 + B_2 X_2 + A_3 X_2^2  + A_4 X_1 X_2 \right)\frac{\partial}{\partial X_2} , 1+ X_1^2 + X_2^2 - t^2 \right]|_{t^2=1+ X_1^2 + X_2^2} \notag\\
&+ \left[  \left( \left( C_0 + C_1 X_1 + C_2 X_2 \right) + t \left(C_3 + A_3 X_2 + A_4 X_1 \right)\right)\frac{\partial}{\partial t}, 1+ X_1^2 + X_2^2 - t^2 \right]|_{t^2= 1+ X_1^2 + X_2^2} \notag \\
&\equiv \left(B_0+\left(B_1 - A_2\right)X_1 - 2 A_1 X_2 + \left(4A_0-  2A_4\right)X_1X_2 + 2A_2 X_1X_2^2 + \left(2B_0 - A_3 \right) X_1^2 - A_3 X_2^2 + 2B_1 X_1^3\right)\frac{\partial}{\partial X_1} \label{css31}\\
&+ \left(A_0 - 2 B_2 X_1 + \left(A_2- B_1\right) X_2 + \left(4B_0-  2  A_3 \right)X_1X_2+ 2 B_1 X_1^2 X_2 + \left(2A_0- A_4\right) X_2^2 -A_4 X_1^2 + 2A_2X_2^3 \right)\frac{\partial}{\partial X_2} \notag \\
&+ \left(- C_2 X_1 - C_1 X_2+2C_0X_1X_2 + t\left(\left(2B_0- A_3\right)X_1+ \left( 2A_0 - A_4 \right)X_2 + 2B_1 X_1^2 + 2A_2X_2^2 \right) \right)\frac{\partial}{\partial t}  \notag \\
&- 2C_3+(2A_0- 2 A_4)X_1 + (2 B_0 - 2 A_3) X_2  + (2A_1 - 2 C_3 ) X_1^2+ ( 2A_2 + 2 B_1) X_1X_2  +  (2B_2- 2 C_3) X_2^2+ t\left( - 2 C_0 - 2 C_1 X_1 - 2 C_2 X_2 \right)  \notag
\end{align}}}
Then in view of $(\ref{css12})$ and $(\ref{css6})$, we see that
\begin{align*}
\mathbb{H}^0\left( Y, \Theta_Y^\bullet \right) = \textnormal{ker}\left( H^0\left( Y, \Theta_Y \right)\to   H^0\left( Y, \Theta_{\tilde{\mathcal{F}}_0^\sharp}^* \otimes \frac{\Theta_Y}{\Theta_{\tilde{\mathcal{F}}_0^\sharp}} \right) \bigoplus H^0\left( X, \mathcal{O}_X(Y) \right) \right) = 0
\end{align*}
and so from $(\ref{css32})$ we have
\begin{align*}
\dim_\mathbb{C} \mathbb{H}^1\left( Y, \Theta_Y^\bullet \right)= 17 - 12 = 5
\end{align*}

By comparing $(\ref{css32})$ and $(\ref{css31})$, we see that
\begin{align}
\mathbb{H}^1\left( Y, \Theta_Y^\bullet \right) &\cong \frac{\textnormal{ker}\left(H^0\left( Y, \Theta_{\tilde{\mathcal{F}}_0^\sharp}^* \otimes \frac{\Theta_Y}{\Theta_{\tilde{\mathcal{F}}_0^\sharp}} \right) \bigoplus H^0\left( X, \mathcal{O}_X(Y) \right) \to H^0\left(X, \left( \Theta_{\tilde{\mathcal{F}}_0^\sharp}^* \otimes \mathcal{O}_Y(X) \right)|_X \right) \right)}{\textnormal{im}\left( H^0\left( Y, \Theta_Y \right)\to   H^0\left( Y, \Theta_{\tilde{\mathcal{F}}_0^\sharp}^* \otimes \frac{\Theta_Y}{\Theta_{\tilde{\mathcal{F}}_0^\sharp}} \right) \bigoplus H^0\left( X, \mathcal{O}_X(Y) \right)\right)} \notag \\
&\cong \textnormal{Span}_\mathbb{C}\left\langle X_1 \frac{\partial}{\partial X_1} + X_2\frac{\partial}{\partial X_2} - 2 X_1X_2, X_1X_2\frac{\partial}{\partial X_1}- X_1^2 \frac{\partial}{\partial X_2}, X_2^2\frac{\partial}{\partial X_1}- X_1X_2 \frac{\partial}{\partial X_2}, \right. \label{css35} \\
 &\left. \,\,\,\,\,\,\,\,\,\,\,\,\,\,\,\,\,\,\,\,\,\,\,\,\,\,\,\,\,\,\,\,\,\,\left( 1 + X_1^2 + X_2^2 - t^2 \right)\frac{\partial}{\partial t} , X_1^2 - X_2^2 \right\rangle \notag
\end{align}
Then we have simultaneous deformation $\mathcal{X}\subset \left(Y \times \textnormal{Spec}(A), \tilde{T} \right)$ of $X\subset\left( Y, \Theta_{\tilde{\mathcal{F}}_0^\sharp} \right)$ over $A=\frac{\mathbb{C}[c_1,c_2,x_3,c_4,c_5 ]}{\left(c_1^2,c_2^2,c_3^2,c_4^2,c_5^2 \right)}$ where $\tilde{T}$ is given by
{\Small{\begin{align} \label{css33}
\tilde{T}=  \left( X_2+ 2 X_1^2 X_2 - c_1 X_1 - c_3X_2^2 - c_4 X_1X_2  \right)\frac{\partial}{\partial X_1} + \left( X_1+ 2X_1 X_2^2 - c_1 X_2 + c_3 X_1X_2+ c_4 X_1^2 \right)\frac{\partial}{\partial X_2} +\left( 2t X_1X_2 +c_5\left(1+ X_1^2+ X_2^2-  t^2 \right) \right)\frac{\partial}{\partial t}
\end{align}}}
and $\mathcal{X}$ is given by
\begin{align} \label{css36}
1+ X_1^2 + X_2^2 -t^2 + 2 c_1 X_1X_2- c_2\left(X_1^2 - X_2^2 \right)= 0
\end{align}
The author has not yet established theorem of existence and completeness for simultaneous deformations of compact foliated complex manifolds defined by locally free subsheaves of tangent sheaves, and their non-degenerate invariant complex submanifolds, though he believes that they are true. Under the assumption that theorem of existence and completeness holds, since $\mathbb{H}^2\left( Y, \Theta_Y^\bullet \right)=0$, we should have a foliated analytic family of simultaneous deformation $X\subset \left( Y, \Theta_{\tilde{\mathcal{F}}_0^\sharp} \right)$ such that the characteristic map is an isomorphism onto $\mathbb{H}^1\left(Y,  \Theta_Y^\bullet \right)$ given by $(\ref{css35})$ at the distinguished point. However the author have difficulties with directly constructing a foliated analytic family of simultaneous deformations of $X\subset \left( Y, \Theta_{\tilde{\mathcal{F}}_0^\sharp} \right)$ extending $(\ref{css33})$ and $(\ref{css36})$. We also note that the coefficient of $\frac{\partial}{\partial t}$ of $\tilde{T}$ in $(\ref{css33})$ is not of the form $t\cdot g\left(X_1,  X_2,c_1,...,c_5\right)$ for some holomorphic function $g(X_1,X_2,c_1,...,c_5)$ which is the form in the construction of foliated cyclic coverings, so that we can conclude that simultaneous deformations of $X \subset \left( Y, \Theta_{\tilde{\mathcal{F}}_0^\sharp} \right)$ does not necessarily arise from the construction of foliated double coverings of simultaneous deformations of $ C \subset \left( \mathbb{P}_\mathbb{C}^2, T_0 \right)$.

\section{Foliations on $T\times \mathbb{P}_\mathbb{C}^1$ defined by Poisson structures, where $T$ is a complex torus with $\dim_\mathbb{C} T =2$} 

In this section we study foliated deformations of $\left(X=T\times \mathbb{P}_\mathbb{C}^1 ,\Lambda_{w_0}\right)$ in terms of cotangent sheaves and tangent sheaves, and unfoldings of $\left(X, \Lambda_{w_0} \right)$, where $T$ is a complex torus with $\dim_\mathbb{C} T=2$ and the (regular) foliations $\Lambda_{w_0} \in H^0\left( X, \Omega_X^1 \otimes \bigwedge^3 \Theta_X \right)$ are defined by holomorphic Poisson structures $\Lambda_0\in H^0\left( X, \bigwedge^2 \Theta_X \right)$ on $T\times \mathbb{P}_\mathbb{C}^1$. We shall study the relationship between Poisson deformations and foliated deformations in terms of cotangent sheaves of threefolds defined by holomorphic Poisson structures in general in Part VIII.

\subsection{Holomorphic Poisson structures on $X=T\times \mathbb{P}_\mathbb{C}^1$ and the associated codimension $1$-foliations on $X$} \

Let $T=\mathbb{C}^2/\mathcal{G}$ be a complex torus with dimension $2$, where $\mathcal{G}$ is the free abelian group generated by non-degenerate periods $\omega_j=\left( \omega_{j1}, \omega_{j2} \right),j=1,2,3,4$ and where, if $(z_1,z_2)$ are the coordinates of $\mathbb{C}^2$, $\omega_j\in \mathcal{G}$ operates on $\mathbb{C}^2$ by sending $z=(z_1,z_2)$ into $z+\omega_j = \left(z_1+ \omega_{j1}, z_2+ \omega_{j2}\right)$. Let $X=T\times \mathbb{P}_\mathbb{C}^1$ and $\xi$ is the inhomogeneous coordinates of $\mathbb{P}_\mathbb{C}^1$. Let $\pi_1: X\to T$ and $\pi_2 : X\to \mathbb{P}_\mathbb{C}^1$ be the natural projections. Then we have $\Theta_X= \pi_1^* \Theta_T \bigoplus \pi_2^* \Theta_{\mathbb{P}_\mathbb{C}^1}$, and $\bigwedge^2 \Theta_X= \pi_1^*\left( \bigwedge^2 \Theta_T \right) \bigoplus  \left( \pi_1^* \Theta_T \otimes \pi_2^* \Theta_{\mathbb{P}_\mathbb{C}^1} \right)$. Then by K\"unneth formula, we have
\begin{align*}
H^0\left(X, \bigwedge^2 \Theta_X \right)\cong H^0\left(T, \bigwedge^2 \Theta_T \right) \bigoplus \left(H^0\left( \mathbb{P}_\mathbb{C}^1, \Theta_{\mathbb{P}_\mathbb{C}^1} \right)\otimes H^0\left( T, \Theta_T \right)\right)
\end{align*}
Then $\dim_\mathbb{C} H^0\left( X, \bigwedge^2 \Theta_X \right)= 7$ and any element of $H^0\left( X, \bigwedge^2 \Theta_X \right)$ has an unique representation of the form
\begin{align}\label{ct1}
a\frac{\partial}{\partial z_1}\wedge \frac{\partial}{\partial z_2} + \left( b_0 + b_1 \xi + b_2 \xi^2 \right) \frac{\partial}{\partial z_2}\wedge \frac{\partial}{\partial \xi} + \left( c_0 + c_1 \xi + c_2 \xi^2  \right) \frac{\partial}{\partial \xi}\wedge \frac{\partial}{\partial z_1},\,\,\,\,\,\,\,a,b_0,b_1,b_2,c_0,c_1,c_2\in \mathbb{C}
\end{align}
We describe holomorphic Poisson structures on $X$. Let $\Lambda_0\in H^0\left(X, \bigwedge^2 \Theta_X \right)$ be an element as in the form $(\ref{ct1})$. We set
\begin{align*}
P_1:&=a \\
P_2:&= b_0 + b_1 \xi + b_2 \xi^2\\
P_3:&= c_0 + c_1 \xi + c_2 \xi^2
\end{align*}
Then $\Lambda_0$ defines a holomorphic Poisson structures on $X$ if and only if $\left[ \Lambda_0, \Lambda_0\right]=0$ where $[-.-]$ is the Schouten bracket if and only if
\begin{align*}
& P_1\left( \frac{\partial P_3}{\partial z_1} - \frac{\partial P_2}{\partial z_2} \right) + P_2 \left( \frac{\partial P_1}{\partial z_2} - \frac{\partial P_3}{\partial \xi} \right) + P_3\left( \frac{\partial P_2}{\partial \xi} - \frac{\partial P_1}{\partial z_1} \right) = 0 \\
& \iff b_1c_0- b_0c_1=0,\,\,\,\,\, b_2c_0- b_0 c_2=0,\,\,\,\,\, b_2c_1- b_1 c_2= 0
\end{align*}
Hence we can divide holomorphic Poisson structures on $X= T\times \mathbb{P}_\mathbb{C}^1$ into three classes
\begin{enumerate}
\item $\Lambda_0= D \frac{\partial}{\partial z_1}\wedge \frac{\partial}{\partial z_2}, D\in \mathbb{C}$.
\item $\Lambda_0= D\frac{\partial}{\partial z_1}\wedge \frac{\partial}{\partial z_2} + \left( A+ B \xi + C \xi^2 \right) \frac{\partial}{\partial z_2}\wedge \frac{\partial}{\partial \xi} + k\left( A + B \xi + C \xi^2 \right), D\in \mathbb{C},\left( A , B , C \right)\ne 0 \in \mathbb{C}^3, k\in \mathbb{C}$
\item $\Lambda_0= D \frac{\partial}{\partial z_1} \wedge \frac{\partial}{\partial z_2} + \left( A+ B \xi + C \xi^2 \right) \frac{\partial}{\partial \xi} \wedge \frac{\partial}{\partial z_1}, D \in \mathbb{C}, \left( A, B, C  \right)\ne 0 \in \mathbb{C}^3$
\end{enumerate}
Then $\Lambda_0\in H^0\left( X, \bigwedge^2 \Theta_X \right)$ with $\left[ \Lambda_0, \Lambda_0 \right]$ corresponds to $\omega_0\in H^0\left( X, \Omega_X^1 \otimes K_X^{-1} \right)$ with $\omega_0 \wedge d \omega_0=0$ where
\begin{enumerate}
\item $\omega_0= D\frac{\partial}{\partial \xi}, D\in C$
\item $\omega_0= Dd\xi + \left( A + B \xi + C \xi^2 \right) dz_1 +  k\left( A+ B \xi + C \xi^2 \right) dz_2, D\in \mathbb{C}, \left( A, B, C \right)\ne 0 \in \mathbb{C}^3, k\in \mathbb{C}$
\item $\omega_0= D d \xi + \left( A+ B \xi + C \xi^2 \right) dz_2, D\in \mathbb{C}, \left( A, B, C \right)\ne 0 \in \mathbb{C}^3$
\end{enumerate}
We note that if $D=0$, then singular locus of $\omega_0$ does not have codimension at least $2$, so that in order for $\omega_0$ to define a foliation on $X$, we should have $D\ne 0$. In the following we assume that $D\ne 0$, so that $\omega_0$ defines a regular foliation on $X$.

We shall study foliated deformations of $\left( X, \omega_0 \right)$ for each case. Since $\omega_0$ is regular, foliated deformations of $\left(X, \omega_0 \right)$ in terms of cotangent sheaves are equivalent to foliated deformation of $\left( X, \omega_0 \right)$ in terms of tangent sheaves. However we shall describe foliated deformations of $\left( X, \omega_0 \right)$ in terms of cotangent sheaves as well as in terms of tangent sheaves. We show that $\left(X, \omega_0 \right)$ is obstructed in foliated deformations for the case $(1)$, and $\left( X, \omega_0 \right)$ is unobstructed in foliated deformations for the case $(2)$ and $(3)$. We shall also study unfoldings of $\left( X, \omega_0 \right)$ for each case. We show that $\left( X, \omega_0 \right)$ is unobstructed in unfoldings for all cases. Since the case $(3)$ becomes the case (2) by changing $dz_1$ with $dz_2$, we only need to consider the case $(1)$ and $(2)$.

\subsection{Foliated deformations of $\left( X, \omega_0 \right)$ where $\omega_0\in H^0\left( X, \Omega_X^1 \otimes \bigwedge^3 \Theta_X \right)$ with $\omega_0 \wedge d\omega_0=0$ in terms of cotangent sheaves} \

 Let us consider the foliation on $X$ defined by $\omega_0  \in H^0\left(X, \Omega_X^1 \otimes \bigwedge^3 \Theta_X \right)$ with $\omega_0 \wedge d\omega_0=0$. Then we have the dual leaf complex from $(\ref{cp2})$
\begin{align*}
\mathcal{E}_{\bigwedge^3 \Theta_X} : \mathcal{E}_{ \bigwedge^3 \Theta_X}\xrightarrow{\tilde{E}_0} \Omega_X^1 \otimes \bigwedge^3 \Theta_X \xrightarrow{\tilde{E}_1} \bigwedge^3 \Theta_X \otimes \bigwedge^3 \Theta_X \otimes \bigwedge^3 \Omega_X^1 \to 0
\end{align*}
equivalently
\begin{align*}
\mathcal{E}_{\bigwedge^3 \Theta_X}^\bullet: \mathcal{E}_{\bigwedge^3 \Theta_X} \xrightarrow{\tilde{E}_0} \Omega_X^1\otimes \bigwedge^3 \Theta_X \xrightarrow{\tilde{E}_1} \bigwedge^3 \Theta_X  \to 0
\end{align*}

Let us consider the spectral sequence at $E_1$ from Dolbeault bicomplex associated to $\mathcal{E}_{\bigwedge^3 \Theta_X}^\bullet$:
\begin{equation} \label{ct7}
\begin{CD}
0 \\
@AAA \\
H^0\left(X, \bigwedge^3 \Theta_X \right)\cong H^0\left(X, \bigwedge^3 \Theta_X \otimes \bigwedge^3 \Theta_X \otimes \bigwedge^3 \Omega_X^1 \right) @.\,\,\,\,\cdots \\
@AAA @AAA \\
H^0\left(X, \Omega_X^1 \otimes  \bigwedge^3 \Theta_X \right) @. \,\,\,\,\,\,\,\,\,\,\,\, H^1\left(X , \Omega_X^1 \otimes \bigwedge^3 \Theta_X \right) \\
@AAA @AAA  \\
H^0\left( X, \mathcal{E}_{\bigwedge^3 \Theta_X} \right) @. \,\,\,\,\,\,\,\,\,\,\, H^1\left( X,  \mathcal{E}_{\bigwedge^3 \Theta_X} \right) @. \,\,\,\,\,\,\,\,\,H^2\left(X, \mathcal{E}_{\bigwedge^3 \Theta_X} \right)
\end{CD}
\end{equation}
Then we have
\begin{align} \label{ct12}
\mathbb{H}^1\left( X, \mathcal{E}_{\bigwedge^3 \Theta_X}^\bullet \right) \cong \frac{\textnormal{ker}\left( H^0\left(X, \Omega_X^1 \otimes \bigwedge^3 \Theta_X \right) \to H^0\left( X, \bigwedge^3 \Theta_X \right) \right)}{\textnormal{im}\left( H^0\left( X, \mathcal{E}_{\bigwedge^3 \Theta_X}  \right)\to H^0\left( X, \Omega_X^1 \otimes \bigwedge^3 \Theta_X \right) \right)} \bigoplus \textnormal{ker}\left(  H^1\left(X, \mathcal{E}_{\bigwedge^3 \Theta_X} \right) \to H^1\left( X, \Omega_X^1 \otimes \bigwedge^2 \Theta_X \right) \right)
\end{align}

We will describe each cohomology group in $(\ref{ct7})$. We note that $\dim_\mathbb{C} H^0\left(X, \Theta_X \right)=\dim_\mathbb{C} H^0\left(X, \pi_1^* T \bigoplus \pi_2^* \Theta_{\mathbb{P}_\mathbb{C}^1} \right)=5$ and any element $V$ of $H^0\left(X, \Theta_X \right)$ has an unique representation of the form
\begin{align}\label{ct8}
V:=a_1\frac{\partial}{\partial z_1} + a_2 \frac{\partial}{\partial z_1} +\left( a_3 + a_4 \xi + a_5 \xi^2 \right)\frac{\partial}{\partial \xi},\,\,\,\,\,a_1,...,a_5 \in \mathbb{C}
\end{align}

We note that $\dim_\mathbb{C} H^1\left(X, \Theta_X \right)= \dim_\mathbb{C} H^1\left( X, \pi_1^* T \bigoplus \pi^* \Theta_{\mathbb{P}_\mathbb{C}^1} \right)= 10$ and any element $\phi$ of $H^1\left( X, \Theta_X \right)$ has an unique representation of the form
\begin{align} \label{ct9}
\phi:=d \bar{z}_1\wedge \left( r_0 \frac{\partial}{\partial z_1} + r_1 \frac{\partial}{\partial z_2} + \left( a_0 + a_1 \xi + a_2 \xi^2 \right) \frac{\partial}{\partial \xi} \right) + d\bar{z}_2\wedge  \left( r_2 \frac{\partial}{\partial z_1} + r_3 \frac{\partial}{\partial z_2} + \left( b_0 + b_1 \xi + b_2 \xi^2   \right) \frac{\partial}{\partial \xi} \right) 
\end{align}
where $r_0,r_1,r_2,r_3, a_0,a_1,a_2,b_0,b_1,b_2\in \mathbb{C}$.

We note that $\bigwedge^3 \Theta_X = \bigwedge^3\left( \pi_1^* \Theta_T \bigoplus \pi_2^* \Theta_{\mathbb{P}_\mathbb{C}^1}  \right)=\bigwedge^2 \pi_1^*\Theta_T \otimes \pi_1^* \Theta_{\mathbb{P}_\mathbb{C}^1}\cong \pi_2^* \Theta_{\mathbb{P}_\mathbb{C}^3}$. Then $\dim_\mathbb{C} H^0\left(X, \bigwedge^3 \Theta_X \right)= 3$ and $H^0\left( X, \bigwedge^3 \Theta_X \otimes \bigwedge^3 \Theta_X \otimes \bigwedge^3 \Omega_X^1 \right)$ has an unique representation of the form
\begin{align}
\left(e_0 + e_1 \xi + e_2 \xi^2 \right) d\xi \wedge dz_1 \wedge dz_2\otimes \left( \frac{\partial}{\partial \xi} \wedge \frac{\partial}{\partial z_1} \wedge \frac{\partial}{\partial z_2}  \right) \otimes  \left( \frac{\partial}{\partial \xi} \wedge \frac{\partial}{\partial z_1} \wedge \frac{\partial}{\partial z_2} \right),\,\,\,\,\,  e_0, e_1 , e_2 \in \mathbb{C}
\end{align}

We note that $\Omega_X^1 \otimes \bigwedge^3 \Theta_X \cong \left( \pi_1^* \Omega_T^1 \bigoplus \pi_2^* \Omega_{\mathbb{P}_\mathbb{C}^1}^1 \right)\otimes \pi_2^* \Theta_{\mathbb{P}_\mathbb{C}^1}\cong \left( \pi_1^* \Omega_T^1\otimes \pi_2^* \Theta_{\mathbb{P}_\mathbb{C}^1} \right) \bigoplus \mathcal{O}_X$. Then we see that $\dim_\mathbb{C} H^0\left( X, \Omega_X^1 \otimes \bigwedge^3 \Theta_X \right)=7$ and any element $b$ of $H^0\left( X, \Omega_X^1 \otimes \bigwedge^3 \Theta_X \right)$ has an unique representation of the form
\begin{align} \label{ct14}
b=\left( a d\xi + \left( b_0 + b_1 \xi + b_2 \xi^2 \right) dz_1 + \left( c_0 + c_1 \xi + c_2 \xi^2 \right) dz_2 \right)\otimes \left( \frac{\partial}{\partial \xi} \otimes  \frac{\partial}{\partial z_1} \otimes \frac{\partial}{\partial z_2}\right)
\end{align}

We note that $\dim_\mathbb{C} H^1 \left(X, \Omega_X^1 \otimes \bigwedge^3 \Theta_X \right)=14$ and any element of $H^1 \left(X, \Omega_X^1 \otimes \bigwedge^3 \Theta_X \right)$ has an unique representation of the form
\begin{align} \label{ct15}
&d\bar{z}_1 \wedge \left( f_0 d\xi + \left(f_1+ f_2 \xi + f_3 \xi^2 \right) d z_1+ \left(f_4 + f_5 \xi + f_6 \xi^2 \right) dz_2 \right)\otimes \left( \frac{\partial}{\partial \xi} \wedge \frac{\partial}{\partial z_1} \wedge \frac{\partial}{\partial z_2} \right) \\
&+  d \bar{z}_2 \wedge \left( f_7 d \xi + \left( f_8 + f_9 \xi + f_{10}\xi^2 \right) dz_1+ \left( f_{11} + f_{12} \xi + f_{13} \xi^2 \right) dz_2  \right)  \otimes \left( \frac{\partial}{\partial \xi} \wedge \frac{\partial}{\partial z_1} \otimes \frac{\partial}{\partial z_2} \right) \notag
\end{align}
where $f_0,..., f_{13}\in \mathbb{C}$.

It remains to describe $H^0\left( X, \mathcal{E}_{\bigwedge^3 \Theta_X} \right)$ and $H^1\left( X, \mathcal{E}_{\bigwedge^3 \Theta_X} \right)$. First we note that
\begin{lemma}\label{ct10}
Let $M$ be a compact complex manifold with dimension $n$. Then 
\begin{align} \label{ctp1}
H^i\left(M, \mathcal{E}_{\bigwedge^n \Theta_M} \right) \cong H^i\left( M, \Theta_M \right) \bigoplus H^i\left( M, \mathcal{O}_M \right)
\end{align}
where $\mathcal{E}_{\bigwedge^n \Theta_M}$ is the Atiyah extension associated to the line bundle $\bigwedge^n \Theta_M\cong K_M^{-1}$ 
\end{lemma}

\begin{proof}
We will show that $\mathcal{E}_{\bigwedge^n \Theta_M} \cong \Theta_M \bigoplus \mathcal{O}_M$ as sheaves of abelian groups. Then this induces $(\ref{ctp1})$. In fact, let $\mathcal{U}=\left\{ U_i \right\}$ be a Stein open covering of $M$ with local coordinates $z_i=\left( z_i^1,..., z_i^n \right) $ on $U_i$ such that $\mathcal{E}_{\bigwedge^n \Theta_X}$ is defined in the following way: $\mathcal{E}_{\bigwedge^n \Theta_M}|_{U_j}$ is isomorphic to $\Theta_M|_{U_j}\bigoplus \mathcal{O}_M|_{U_j}$ and, $\left(\rho_j, a_j \right)\in \Theta_M|_{U_j}\bigoplus \mathcal{O}_M|_{U_j}$ is identified with $\left(\rho_k, a_k\right)\in \Theta_M|_{U_k}\bigoplus \mathcal{O}_M|_{U_k}$ on $U_j\cap U_k\ne \emptyset$ if and only if $\rho_j=\rho_k$ and $a_k-a_j=\left[ \rho_j, \log h_{jk}\right]$, where $h_{jk}$ is defined by $dz_j^1\wedge \cdots \wedge d z_j^n = h_{jk}dz_k^1\wedge \cdots \wedge dz_k^n$. We define locally on $U_k$
\begin{align}\label{pti1}
\Theta_M\bigoplus \mathcal{O}_M &\to  \mathcal{E}_{\bigwedge^n \Theta_M}:=\Theta_M|_{U_k} \bigoplus \mathcal{O}_M|_{U_k}\\
       \left(\rho= \sum_{\alpha=1}^n \rho_k^\alpha \frac{\partial}{\partial z_k^\alpha},  a\right) &\mapsto \left(\rho,  a- \sum_{\alpha=1}^n\frac{\partial \rho_k^\alpha}{\partial z_k^\alpha} \right) \notag
\end{align}
We claim that $(\ref{pti1})$ is well-defined. Let $z_j= f_{jk}(z_k)$ be the coordinate transformation from $U_k$ to $U_j$. Then $\rho= \sum_{\alpha, \beta=1}^n \rho_k^\alpha \frac{\partial f_{jk}^\beta}{\partial z_k^\alpha}\frac{\partial}{\partial z_j^\beta}$ on $U_j$, and we have $\rho= \sum_{\beta=1}^n \rho_j^\beta \frac{\partial}{\partial z_j^\beta}$ with $\rho_j^\beta= \sum_{\alpha=1}^n \rho_k^\alpha \frac{\partial f_{jk}^\beta}{\partial z_k^\alpha}$, so that $(\rho, a)$ on $U_j$ is mapped to
\begin{align*}
\left(\rho, a- \sum_{\alpha,\beta=1}^n \frac{\partial}{\partial z_j^\beta}\left(\rho_k^\alpha \frac{\partial f_{jk}^\beta}{\partial z_k^\alpha} \right)  \right)\,\,\,\,\,\,\,\,\,\,\textnormal{on}\,\,\,\,\,U_j
\end{align*}
Then we note that
\begin{align*}
&\left(a- \sum_{\alpha=1}^n \frac{\partial \rho_k^\alpha}{\partial z_k^\alpha}\right)-\left( a-  \sum_{\alpha,\beta=1}^n \frac{\partial}{\partial z_j^\beta}\left(\rho_k^\alpha \frac{\partial f_{jk}^\beta}{\partial z_k^\alpha} \right) \right) = - \sum_{\alpha,\beta=1}^n \frac{\partial f_{jk}^\beta}{\partial z_k^\alpha}\frac{\partial \rho_k^\alpha}{\partial z_j^\beta} + \sum_{\alpha,\beta=1}^n \frac{\partial \rho_k^\alpha}{\partial z_j^\beta} \frac{\partial f_{jk}^\beta}{\partial z_k^\alpha} + \sum_{\alpha,\beta=1}^n \rho_k^\alpha \frac{\partial}{\partial z_j^\beta}\left(\frac{\partial f_{jk}^\beta}{\partial z_k^\alpha} \right) \\
&=  \sum_{\alpha,\beta=1}^n \rho_k^\alpha \frac{\partial}{\partial z_j^\beta}\left(\frac{\partial f_{jk}^\beta}{\partial z_k^\alpha} \right)=\sum_{\alpha=1}^n \rho_k^\alpha \left(  \frac{1}{h_{jk}}\frac{\partial h_{jk}}{\partial z_k^\alpha}  \right) = \frac{1}{h_{jk}}\left[ \rho, h_{jk}  \right] = \left[ \rho, \log h_{jk} \right]
\end{align*}
which follows from $\sum_{\beta=1}^n \frac{\partial}{\partial z_j^\beta}\left( \frac{\partial f_{jk}^\beta}{\partial z_k^\alpha} \right) = \frac{1}{h_{jk}}\frac{\partial h_{jk}}{\partial z_k^\alpha}$. In fact, we set
\begin{align*}
H_{jk}(z_k):= \left(\begin{matrix} \frac{\partial f_{jk}^1}{\partial z_k^1} & \frac{\partial f_{jk}^2}{\partial z_k^1} & \cdots & \frac{\partial f_{jk}^n}{\partial z_k^1} \\ \frac{\partial f_{jk}^1}{\partial z_k^2} & \frac{\partial f_{jk}^2}{\partial z_k^2} & \cdots & \frac{\partial f_{jk}^n}{\partial z_k^2} \\ 
\cdots & \cdots & \cdots & \cdots \\ \frac{\partial f_{jk}^1}{\partial z_k^n} & \frac{\partial f_{jk}^2}{\partial z_k^n} & \cdots & \frac{\partial f_{jk}^n}{\partial z_k^n} \end{matrix}\right)
\end{align*}
Then we have $\det\left( H_{jk}(z_k) \right)= h_{jk}$ and we have
{\Small{\begin{align*}
 \left(\begin{matrix} \frac{\partial}{\partial z_k^\alpha}\left(\frac{\partial f_{jk}^1}{\partial z_k^1}\right) & \frac{\partial}{\partial z_k^\alpha}\left(\frac{\partial f_{jk}^2}{\partial z_k^1} \right) & \cdots & \frac{\partial}{\partial z_k^\alpha}\left(\frac{\partial f_{jk}^n}{\partial z_k^1}\right) \\ \frac{\partial}{\partial z_k^\alpha}\left(\frac{\partial f_{jk}^1}{\partial z_k^2} \right) & \frac{\partial}{\partial z_k^\alpha}\left(\frac{\partial f_{jk}^2}{\partial z_k^2} \right) & \cdots & \frac{\partial}{\partial z_k^\alpha}\left( \frac{\partial f_{jk}^n}{\partial z_k^2} \right) \\ 
\cdots & \cdots & \cdots & \cdots \\  \frac{\partial}{\partial z_k^\alpha} \left( \frac{\partial f_{jk}^1}{\partial z_k^n} \right) & \frac{\partial}{\partial z_k^\alpha}\left(\frac{\partial f_{jk}^2}{\partial z_k^n} \right) & \cdots & \frac{\partial}{\partial z_k^\alpha}\left(\frac{\partial f_{jk}^n}{\partial z_k^n} \right) \end{matrix}\right)=  \left(\begin{matrix} \frac{\partial f_{jk}^1}{\partial z_k^1} & \frac{\partial f_{jk}^2}{\partial z_k^1} & \cdots & \frac{\partial f_{jk}^n}{\partial z_k^1} \\ \frac{\partial f_{jk}^1}{\partial z_k^2} & \frac{\partial f_{jk}^2}{\partial z_k^2} & \cdots & \frac{\partial f_{jk}^n}{\partial z_k^2} \\ 
\cdots & \cdots & \cdots & \cdots \\ \frac{\partial f_{jk}^1}{\partial z_k^n} & \frac{\partial f_{jk}^2}{\partial z_k^n} & \cdots & \frac{\partial f_{jk}^n}{\partial z_k^n} \end{matrix}\right)\cdot \left(\begin{matrix} \frac{\partial}{\partial z_j^1}\left(\frac{\partial f_{jk}^1}{\partial z_k^\alpha}\right) & \frac{\partial}{\partial z_j^1}\left(\frac{\partial f_{jk}^2}{\partial z_k^\alpha} \right) & \cdots & \frac{\partial}{\partial z_j^1}\left(\frac{\partial f_{jk}^n}{\partial z_k^\alpha}\right) \\ \frac{\partial}{\partial z_j^2}\left(\frac{\partial f_{jk}^1}{\partial z_k^\alpha} \right) & \frac{\partial}{\partial z_j^2}\left(\frac{\partial f_{jk}^2}{\partial z_k^\alpha} \right) & \cdots & \frac{\partial}{\partial z_j^2}\left( \frac{\partial f_{jk}^n}{\partial z_k^\alpha} \right) \\ 
\cdots & \cdots & \cdots & \cdots \\  \frac{\partial}{\partial z_j^n} \left( \frac{\partial f_{jk}^1}{\partial z_k^\alpha} \right) & \frac{\partial}{\partial z_j^n}\left(\frac{\partial f_{jk}^2}{\partial z_k^\alpha} \right) & \cdots & \frac{\partial}{\partial z_j^n}\left(\frac{\partial f_{jk}^n}{\partial z_k^\alpha} \right) \end{matrix}\right)
\end{align*}}}
Then we note that
\begin{align*}
\frac{\partial h_{jk}}{\partial z_k^\alpha } =\frac{\partial}{\partial z_k^\alpha}\left(\det H_{jk} \right)=   h_{jk} \cdot \textnormal{Tr}\left( H_{jk}^{-1}  \cdot \frac{\partial H_{jk}}{\partial z_k^\alpha}      \right)= h_{jk} \left( \frac{\partial}{\partial z_j^1}\left(\frac{\partial f_{jk}^1}{\partial z_k^\alpha }\right) +\frac{\partial}{\partial z_j^2}\left( \frac{\partial f_{jk}^2}{\partial z_k^\alpha} \right) +\cdots + \frac{\partial}{\partial z_j^n}\left( \frac{\partial f_{jk}^n}{\partial z_k^\alpha} \right) \right) 
\end{align*}
Hence $(\ref{pti1})$ is well-defined with  the inverse $\left(\rho_j= \sum_{\alpha=1}^n \rho_j^\alpha \frac{\partial}{\partial z_j^\alpha} , a_j \right)\mapsto \left(\rho_j, a_j + \sum_{\alpha=1}^n \frac{\partial \rho_j^\alpha}{\partial z_j^\alpha}  \right)$.

\end{proof}

From the isomorphism $(\ref{pti1})$ in Lemma \ref{ct10} and $(\ref{ct8})$, we have an isomorphism
\begin{align*}
H^0\left( X, \Theta_X \right) \bigoplus H^0\left( X , \mathcal{O}_X \right) &\to H^0\left(X, \mathcal{E}_{\bigwedge^3 \Theta_X}  \right)\\
\left( a_1 \frac{\partial}{\partial z_1} + a_2 \frac{\partial}{\partial z_2} + \left( a_3+ a_4 \xi + a_5 \xi^2\right)\frac{\partial}{\partial \xi}  , a_6\right) &\mapsto \left( a_1 \frac{\partial}{\partial z_1} + a_2 \frac{\partial}{\partial z_2} + \left( a_3 + a_4 \xi + a_5 \xi^2 \right)\frac{\partial}{\partial \xi} , a_6- a_4 -2a_5 \xi \right)
\end{align*}
where $a_1,...,a_6\in \mathbb{C}$. Then we see that $\dim_\mathbb{C} H^0\left(X, \mathcal{E}_{\bigwedge^3 \Theta_X} \right)= 6$ and $H^0\left( X,  \mathcal{E}_{\bigwedge^3 \Theta_X} \right)$ has an unique representation of the form
\begin{align}\label{ct13}
\left( a_1 \frac{\partial}{\partial z_1} + a_2 \frac{\partial}{\partial z_2} + \left( a_3 + a_4 \xi + a_5 \xi^2 \right) \frac{\partial}{\partial \xi} , a_6- 2a_5 \xi\right),\,\,\,\,\,\,\,\,a_1,...,a_6 \in \mathbb{C}
\end{align}
On the other hand, from the isomorphism $(\ref{pti1})$ and $(\ref{ct9})$, we have an isomorphism
\begin{align*}
&H^1\left( X, \Theta_X \right) \bigoplus H^1\left( X , \mathcal{O}_X \right) \to H^1\left(X, \mathcal{E}_{\bigwedge^3 \Theta_X}  \right)\\
&\left(d \bar{z}_1\wedge \left( r_0 \frac{\partial}{\partial z_1} + r_1 \frac{\partial}{\partial z_2} + \left( a_0 + a_1 \xi + a_2 \xi^2 \right) \frac{\partial}{\partial \xi} \right) + d\bar{z}_2\wedge  \left( r_2 \frac{\partial}{\partial z_1} + r_3 \frac{\partial}{\partial z_2} + \left( b_0 + b_1 \xi + b_2 \xi^2   \right) \frac{\partial}{\partial \xi} \right)  , c_1d\bar{z}_1+ c_2d\bar{z}_2 \right)\\
&\mapsto \left(d \bar{z}_1\wedge \left( r_0 \frac{\partial}{\partial z_1} + r_1 \frac{\partial}{\partial z_2} + \left( a_0 + a_1 \xi + a_2 \xi^2 \right) \frac{\partial}{\partial \xi} \right) + d\bar{z}_2\wedge  \left( r_2 \frac{\partial}{\partial z_1} + r_3 \frac{\partial}{\partial z_2} + \left( b_0 + b_1 \xi + b_2 \xi^2   \right) \frac{\partial}{\partial \xi} \right)  , \right.\\
&\left. \,\,\,\,\,\,\,\,\, \left(c_1 - a_1- 2 a_2 \xi \right)d\bar{z}_1+ \left(c_2 - b_1 - 2b_2 \right) d\bar{z}_2 \right),\,\,\,\,\,\,\,\,r_0,r_1,r_2,r_3, a_0,a_1,a_2,b_1,b_2,b_3,c_1,c_2\in \mathbb{C}
\end{align*}
Then we see that $\dim_\mathbb{C} H^1\left( X,  \mathcal{E}_{\bigwedge^3 \Theta_X} \right)= 12$ and any element of $H^1\left( X, \mathcal{E}_{\bigwedge^3 \Theta_X} \right)$ has an unique representation of the form
{\Small{\begin{align}\label{ct11}
\left(d \bar{z}_1\wedge \left( r_0 \frac{\partial}{\partial z_1} + r_1 \frac{\partial}{\partial z_2} + \left( a_0 + a_1 \xi + a_2 \xi^2 \right) \frac{\partial}{\partial \xi} \right) + d\bar{z}_2\wedge  \left( r_2 \frac{\partial}{\partial z_1} + r_3 \frac{\partial}{\partial z_2} + \left( b_0 + b_1 \xi + b_2 \xi^2   \right) \frac{\partial}{\partial \xi} \right) , \left( a_3 - 2a_2 \xi \right)d\bar{z}_1 + \left( b_3 - 2b_2 \xi  \right) d\bar{z}_2 \right)
\end{align}}}
where $r_0,r_1,r_2,r_3,a_0,a_1,a_2,a_3,b_0,b_1,b_2,b_3\in \mathbb{C}$.

With this preparation, we will study foliated deformations of $\left(X, \omega_0 \right)$ in terms of cotangent sheaves where $\omega_0\in H^0\left(X, \Omega_X^1 \otimes \bigwedge^3 \Theta_X \right)$ with $\omega_0 \wedge d\omega_0 =0$.

\subsection{Foliated deformations of $\left( X, D d\xi \right), D\ne 0$ in terms of cotangent sheaves} \

Let us consider the foliation on $X$ defined by $\omega_0=D d\xi, D\ne 0$. We compute $\mathbb{H}^1\left( X, \mathcal{E}_{\bigwedge^3 \Theta_X}^\bullet  \right)$. From $(\ref{ct12})$, first we compute $H^0\left( X, \mathcal{E}_{\bigwedge^3 \Theta_X} \right) \to H^0\left( X, \Omega_X^1\otimes \bigwedge^3 \Theta_X \right)$. Then from $(\ref{ct13})$, the image is given by
\begin{align}\label{ct28}
\mathcal{L}_{a_1 \frac{\partial}{\partial z_1} + a_2 \frac{\partial}{\partial z_2} + \left(a_3 +a_4 \xi + a_5 \xi^2 \right) \frac{\partial}{\partial \xi}}\left( D d\xi \right) +  \left(a_6 - 2a_5\xi \right) \cdot D d\xi= D\left( a_4+ a_6\right) d\xi
\end{align}
Then from $(\ref{ct14})$ we have $\dim_\mathbb{C} \mathbb{H}^0\left(X, \mathcal{E}_{\bigwedge^3 \Theta_X}^\bullet \right)=\dim_\mathbb{C} \textnormal{ker}\left(H^0\left( X, \mathcal{E}_{\bigwedge^3 \Theta_X} \right) \to H^0\left( X, \Omega_X^1\otimes \bigwedge^3 \Theta_X \right) \right)=5$ and $\dim_\mathbb{C} \textnormal{im}\left(H^0\left( X, \mathcal{E}_{\bigwedge^3 \Theta_X} \right) \to H^0\left( X, \Omega_X^1\otimes \bigwedge^3 \Theta_X \right)  \right)=6-5=1$. Next we compute $H^0\left(X, \Omega_X^1 \otimes \bigwedge^3 \Theta_X \right) \to H^0\left( X, \bigwedge^3 \Theta_X \right)$. From $(\ref{ct14})$ the image is given by
\begin{align} \label{ct29}
b\wedge dw_0 + w_0 \wedge db=0,
\end{align}
so that $\dim_\mathbb{C} \textnormal{ker}\left( H^0\left(X, \Omega_X^1 \otimes \bigwedge^3 \Theta_X \right) \to H^0\left( X, \bigwedge^3 \Theta_X \right)  \right)= 7$. Then we have
\begin{align} \label{ct16}
\dim_\mathbb{C}\frac{\textnormal{ker}\left( H^0\left(X, \Omega_X^1 \otimes \bigwedge^3 \Theta_X \right) \to H^0\left( X, \bigwedge^3 \Theta_X \right) \right)}{\textnormal{im}\left(H^0\left( X , \mathcal{E}_{\bigwedge^3 \Theta_X}  \right)\to H^0\left( X, \Omega_X^1 \otimes \bigwedge^3 \Theta_X \right) \right)}= 7- 1 = 6
\end{align}
 On the other hand, let us compute $ H^1\left(X, \mathcal{E}_{\bigwedge^3 T_X} \right)\to H^1\left( X, \Omega_X^1\otimes \bigwedge^3 \Theta_X \right)$. From $(\ref{ct11})$ the image is given by
\begin{align}\label{cu1}
&d\bar{z}_1\wedge \mathcal{L}_{r_0 \frac{\partial}{\partial z_1} + r_1 \frac{\partial}{\partial z_2} +\left( a_0+ a_1\xi+ a_2 \xi^2 \right)\frac{\partial}{\partial \xi}  } \left( D d\xi \right) +d\bar{z}_2 \wedge \mathcal{L}_{r_2 \frac{\partial}{\partial z_1} + r_3 \frac{\partial}{\partial z_2} +(b_0+ b_1\xi + b_2\xi^2)\frac{\partial}{\partial \xi}  }\left(D d\xi \right) \\
&+ d\bar{z}_1\wedge \left( a_3 - 2a_2 \xi \right) \cdot \left(D d\xi \right)+d\bar{z}_2 \wedge \left(b_3- 2b_2\xi \right) \cdot \left( D d\xi\right) \notag \\
&=d\bar{z}_1\wedge D\left(a_1 + a_3 \right)d \xi + d\bar{z}_2\wedge D\left(b_1 + b_3 \right) d\xi, \notag
\end{align}
so that from $(\ref{ct15})$ we have $\dim_{\mathbb{C}} \textnormal{ker}\left( H^1\left( X, \mathcal{E}_{\bigwedge^3 \Theta_X} \right)\to H^1\left(X, \Omega_X^1 \otimes \bigwedge^3 \Theta_X  \right)  \right)=10$. Then we see that from $(\ref{ct12})$ and $(\ref{ct16})$
\begin{align*}
\dim_\mathbb{C} \mathbb{H}^1\left( X, \mathcal{E}_{\bigwedge^3 \Theta_X}^\bullet \right)= 10 + 6 = 16
\end{align*}

Moreover, we note that the natural map $\mathbb{H}^1\left(X, \mathcal{E}_{\bigwedge^3 \Theta_X}^\bullet \right)\to H^1\left(X, \Theta_X \right)$ is surjective. Hence since $X$ is obstructed in deformations of complex structures (see \cite{Kod58} p.436), $\left(X, \omega_0 \right)$ is obstructed in foliated deformations.

\subsection{Foliated deformations of $\left(X, D d\xi + \left( A+ B \xi + C \xi^2 \right)dz_1 + k\left( A+ B \xi + C \xi^2 \right) dz_2 \right), D\ne 0, \left( A, B, C \right) \ne 0 \in \mathbb{C}^3, k\in \mathbb{C}$ in terms of cotangent sheaves} \

Let us consider the foliation on $X$ defined by $\omega_0= D d\xi + \left( A+ B \xi + C \xi^2 \right)dz_1 + k\left( A+ B \xi + C \xi^2 \right) dz_2, D\ne 0, \left( A, B, C \right) \ne 0 \in \mathbb{C}^3, k\in \mathbb{C}$. We compute $\mathbb{H}^1\left(X , \mathcal{E}_{\bigwedge^3 \Theta_X}^\bullet \right)$. From $(\ref{ct12})$, first we compute $H^0\left( X, \mathcal{E}_{\bigwedge^3 \Theta_X} \right) \to H^0\left( X, \Omega_X^1\otimes \bigwedge^3 \Theta_X \right)$. Then from $(\ref{ct13})$, the image is given by
\begin{align}
&\mathcal{L}_{ a_1\frac{\partial}{\partial z_1} + a_2\frac{\partial}{\partial z_2} + \left( a_3+ a_4 \xi + a_5 \xi^2 \right) \frac{\partial}{\partial \xi}} \left( D d\xi+ \left(A+ B\xi+ C\xi^2\right)dz_1 + k\left(A+ B\xi+ C\xi^2\right)dz_2 \right) \notag \\
&+ \left( a_6-2a_5\xi  \right) \left( D d\xi+ \left(A+ B\xi+ C\xi^2\right)dz_1 + k\left(A+ B\xi+ C\xi^2\right) dz_2 \right) \notag \\
&= D(a_4+ a_6) d\xi +\left(  Ba_3+ A a_6 +\left( Ba_4+2 C a_3+ Ba_6- 2A a_5 \right) \xi +\left(  2 C a_4+ C a_6-  B a_5 \right) \xi^2 \right) dz_1 \label{cu7} \\
&+k \left(  Ba_3+ A a_6 +\left( Ba_4+2 C a_3+ Ba_6- 2A a_5 \right) \xi +\left( 2 C a_4+ C a_6-  B a_5 \right) \xi^2 \right) dz_2 \notag
\end{align}
Then from $(\ref{ct14})$ the kernel is given by
\begin{align*}
&a_4+ a_6=0\\
&B a_3+ Aa_6 =0\\
&Ba_4+ 2 Ca_3+ Ba_6- 2A a_5=0\\
&2Ca_4+ Ca_6- Ba_5=0
\end{align*}
Then $a_6= - a_4$, so that we have
\begin{align*}
Ba_3- A a_4=0\\
 Ca_3- Aa_5 =0 \\
Ca_4- Ba_5 =0
\end{align*}
This implies that $\dim_\mathbb{C} \mathbb{H}^0\left(X, \mathcal{E}_{\bigwedge^3 \Theta_X}^\bullet \right)=\dim_\mathbb{C} \textnormal{ker}\left(H^0\left( X, \mathcal{E}_{\bigwedge^3 \Theta_X} \right) \to H^0\left( X, \Omega_X^1\otimes \bigwedge^3 \Theta_X \right) \right)=3$ and $\dim_\mathbb{C} \textnormal{im}\left(H^0\left( X, \mathcal{E}_{\bigwedge^3 \Theta_X} \right) \to H^0\left( X, \Omega_X^1\otimes \bigwedge^3 \Theta_X \right)  \right)=6-3=3$. Next we compute $H^0\left(X, \Omega_X^1 \otimes \bigwedge^3 \Theta_X \right) \to H^0\left( X, \bigwedge^3 \Theta_X \right)$. we note that $d\omega_0= \left( B+ 2 C\xi\right) d\xi \wedge dz_1 + k \left(B + 2 C \xi  \right) d\xi \wedge dz_2$. Then from $(\ref{ct14})$ the image is given by
{\small{\begin{align}
&b\wedge d\omega_0 + \omega_0 \wedge db  \notag  \\
&=\left(D d\xi+ \left(A+ B\xi+ C\xi^2\right)dz_1 + k\left(A+ B\xi+ C\xi^2\right)dz_2 \right)\wedge \left( \left( b_1 + 2 b_2 \xi \right) d\xi\wedge dz_1 +\left( c_1 + 2 c_2\xi \right) d\xi\wedge dz_2\right) \notag  \\
&+\left(  a d\xi + \left(b_0+ b_1 \xi + b_2 \xi^2 \right) dz_1 + \left( c_0 + c_1 \xi + c_2 \xi^2  \right) dz_2 \right) \wedge  \left( \left(B+ 2C \xi \right)d\xi\wedge dz_1 + k\left( B + 2 C \xi\right) d\xi\wedge dz_2 \right) \notag \\
&=\left(\left(Bc_0- A c_1+ k \left(Ab_1- Bb_0\right) \right) + \left(2 Cc_0- 2Ac_2 + k\left(2Ab_2- 2Cb_0 \right) \right) \xi + \left( Cc_1- Bc_2 + k\left(Bb_2- Cb_1\right) \right) \xi^2 \right) d\xi \wedge dz_1 \wedge dz_2  \label{cu15}
\end{align}}}
Then the kernel is given by
\begin{align*}
B\left(c_0- kb_0\right)- A\left(c_1- kb_1\right)=0\\
C \left(c_0- kb_0 \right)- A \left(c_2- kb_2 \right)=0\\
C\left(c_1-kb_1\right)-B \left(c_2- kb_2 \right)=0
\end{align*}
Then we have
\begin{align}
c_0- kb_0&= t_1 A \Longrightarrow c_0= t_1 A+ kb_0 \notag\\
c_1- kb_1&= t_1 B \Longrightarrow c_1= t_1 B + kb_1 \label{cu16} \\
c_2- kb_2&= t_1C \Longrightarrow c_2= t_1 C + kb_2,\,\,\,\,\,\,\,t_1\in \mathbb{C} \notag
\end{align}
Hence $\dim_\mathbb{C} \textnormal{ker}\left( H^0\left(X, \Omega_X^1 \otimes \bigwedge^3 \Theta_X \right) \to H^0\left( X, \bigwedge^3 \Theta_X \right) \right)=5$ and any element of the kernel has an unique representation of the form
\begin{align}\label{cu6}
e_1 d\xi+ \left(e_2+ e_3 \xi+ e_4\xi^2 \right) dz_1+ \left(ke_2+ e_5 A + \left(ke_3+ e_5 B \right)\xi + \left(ke_4+ e_5 C \right)\xi^2 \right) dz_2,\,\,\,\,\,\,\,e_1,e_2,e_3,e_4,e_5\in \mathbb{C}
\end{align}
Then we have $\dim_\mathbb{C}\frac{\textnormal{ker}\left( H^0\left(X, \Omega_X^1 \otimes \bigwedge^3 \Theta_X \right) \to H^0\left( X, \bigwedge^3 \Theta_X \right) \right)}{\textnormal{im}\left(H^0\left( X , \mathcal{E}_{\bigwedge^3 \Theta_X}  \right)\to H^0\left( X, \Omega_X^1 \otimes \bigwedge^3 \Theta_X \right) \right)}= 5-3 =2$. We recall that from $(\ref{cu7})$ the image of $\tilde{E}_0$ is given by
\begin{align} \label{cu8}
&D(a_4+ a_6) d\xi +\left(  Ba_3+ A a_6 +\left( Ba_4+2 C a_3+ Ba_6- 2A a_5 \right) \xi +\left(  2 C a_4+ C a_6-  B a_5 \right) \xi^2 \right) dz_1 \\
&+k \left(  Ba_3+ A a_6 +\left( Ba_4+2 C a_3+ Ba_6- 2A a_5 \right) \xi +\left( 2 C a_4+ C a_6-  B a_5 \right) \xi^2 \right) dz_2 \notag
\end{align}
By setting $a_4=-a_6$, $(\ref{cu8})$ becomes
\begin{align*}
\left( Ba_3- Aa_4 + \left( 2Ca_3 - 2Aa_5\right) \xi + \left( Ca_4- B a_5\right) \xi^2 \right) dz_1 +k\left(B a_3- A a_4 +\left(2 C a_3 - 2 Aa_5 \right) \xi  +\left( Ca_4- B a_5\right) \xi^2 \right) dz_2
\end{align*}
Since $\textnormal{rank} \left[\begin{matrix}
B & -A & 0\\
2C & 0 & -2 A\\
0 & C & -B
\end{matrix}\right]=2$, there is $\left(F_0,F_1,F_2\right)\in \mathbb{C}^3$ such that $\left( F_1,F_2,F_3 \right)$ is not in the image of the matrix. We note that $\left( F_0+ F_1 \xi + F_2 \xi^2 \right)dz_1+ k\left( F_0 + F_1 \xi + F_2 \xi^2 \right)dz_2 $ is in the kernel of $\tilde{E}_1$.  It is clear that $\left(A + B \xi + C \xi^2 \right) dz_2$ is not in the image of $\tilde{E}_0$. Then we claim that
{\small{\begin{align}\label{cu9}
&\frac{\textnormal{ker}\left( H^0\left(X, \Omega_X^1 \otimes \bigwedge^3 \Theta_X \right) \to H^0\left( X, \bigwedge^3 \Theta_X \right) \right)}{\textnormal{im}\left(H^0\left( X , \mathcal{E}_{\bigwedge^3 \Theta_X}  \right)\to H^0\left( X, \Omega_X^1 \otimes \bigwedge^3 \Theta_X \right) \right)} \cong \textnormal{Span}_\mathbb{C} \left\langle \left( F_0+ F_1 \xi + F_2 \xi^2 \right)dz_1+ k\left( F_0 + F_1 \xi + F_2 \xi^2 \right)dz_2 ,  \left( A + B \xi + C\xi^2 \right) dz_2 \right\rangle
\end{align}}}
In fact, assume that for some $r_1,r_2\in \mathbb{C}$ we have
\begin{align*}
&r_1\left( \left( F_0+ F_1 \xi + F_2 \xi\right)dz_1+ k\left( F_0 + F_1 \xi + F_2 \xi^2 \right)dz_2\right) + r_2\left(\left( A + B \xi + C\xi^2 \right) dz_2 \right)\\
&= D(a_4+ a_6) d\xi +\left(  Ba_3+ A a_6 +\left( Ba_4+2 C a_3+ Ba_6- 2A a_5 \right) \xi +\left(  2 C a_4+ C a_6-  B a_5 \right) \xi^2 \right) dz_1 \\
&+k \left(  Ba_3+ A a_6 +\left( Ba_4+2 C a_3+ Ba_6- 2A a_5 \right) \xi +\left( 2 C a_4+ C a_6-  B a_5 \right) \xi^2 \right) dz_2 \\
&\in \textnormal{im}\left(H^0\left( X, \mathcal{E}_{\bigwedge^3 \Theta_X}\right)\to H^0\left(X, \Omega_X^1\otimes \bigwedge^3 \Theta_X \right) \right)
\end{align*}
Then $a_4+ a_6=0$, so that by our choice of $\left( F_1, F_2, F_3\right)$, we have $r_1=0$, and so $r_2=0$. This proves the claim.

On the other hand, let us compute $ H^1\left(X, \mathcal{E}_{\bigwedge^3 T_X} \right)\to H^1\left( X, \Omega_X^1\otimes \bigwedge^3 \Theta_X \right)$. From $(\ref{ct11})$ the image is given by
\begin{align}
&d\bar{z}_1\wedge \mathcal{L}_{r_0 \frac{\partial}{\partial z_1} + r_1 \frac{\partial}{\partial z_2} +\left( a_0+ a_1\xi+ a_2 \xi^2 \right)\frac{\partial}{\partial \xi}  } \left( D d\xi+(A+ B\xi+ C\xi^2)dz_1 + k(A+ B\xi+ C\xi^2)dz_2 \right) \notag \\
&d\bar{z}_2 \wedge \mathcal{L}_{r_2 \frac{\partial}{\partial z_1} + r_3 \frac{\partial}{\partial z_2} +(b_0+ b_1\xi + b_2\xi^2)\frac{\partial}{\partial \xi}  }\left(D d\xi+(A+ B\xi+ C\xi^2)dz_1 + k(A+ B\xi+ C\xi^2)dz_2 \right) \notag \\
&+ d\bar{z}_1\wedge \left( a_3 - 2a_2 \xi \right) \cdot \left(D d\xi+(A+ B\xi+ C\xi^2)dz_1 + k(A+ B\xi+ C\xi^2)dz_2 \right) \notag\\
&+d\bar{z}_2 \wedge \left(b_3- 2b_2\xi \right) \cdot \left( D d\xi+(A+ B\xi+ C\xi^2)dz_1 + k(A+ B\xi+ C\xi^2)dz_2\right)\notag\\
&=d\bar{z}_1\wedge \left( D(a_1+ a_3) d\xi +\left(  Ba_0+ A a_3 +\left( Ba_1+2 C a_0+ Ba_3- 2A a_2 \right) \xi +\left(  2 C a_1+ C a_3-  B a_2 \right) \xi^2 \right) dz_1\right) \label{cu18}\\
&+ d\bar{z}_1\wedge \left(k \left(  Ba_1+ A a_3 +\left( Ba_1+2 C a_0+ Ba_3- 2A a_2 \right) \xi +\left( 2 C a_1+ C a_3-  B a_2 \right) \xi^2 \right) dz_2 \right) \notag \\
&+ d\bar{z}_2 \wedge \left( D(b_1+ b_3) d\xi +\left(  Bb_0+ A b_3 +\left( Bb_1+2 C b_0+ Bb_3- 2A b_2 \right) \xi +\left(  2 C b_1+ C b_3-  B b_2 \right) \xi^2 \right) dz_1 \right) \notag \\
&+ d\bar{z}_2 \wedge \left( k \left(  Bb_0+ A b_3 +\left( Bb_1+2 C b_0+ Bb_3- 2A b_2 \right) \xi +\left( 2 C b_1+ C b_3-  B b_2 \right) \xi^2 \right) dz_2\right) \notag
\end{align}
Then we see that any element of $\textnormal{ker}\left( H^1\left(X, \mathcal{E}_{\bigwedge^3 T_X} \right)\to H^1\left( X, \Omega_X^1\otimes \bigwedge^3 \Theta_X \right)\right)$ has an unique representation of the form
\begin{align}
&\left( \left(t_3 \frac{\partial}{\partial z_1} + t_4 \frac{\partial}{\partial z_2} + t_7\left(A+ B \xi + C \xi^2 \right)\frac{\partial}{\partial \xi} \right) d\bar{z}_1 + \left(t_5 \frac{\partial}{\partial z_1} + t_6\frac{\partial}{\partial z_2} + t_8\left(A+ B \xi + C \xi^2 \right)\frac{\partial}{\partial \xi} \right) d\bar{z}_2, \right.\\
&\left. \,\,\,\,\,\,\,\,\,\left(- t_7 B- 2t_7 C \xi\right)d\bar{z}_1 + \left( - t_8 B - 2t_8 C \xi \right) d\bar{z}_2   \right) ,\,\,\,\,\,\,\,\, t_3,t_4,t_5,t_6,t_7,t_8\in \mathbb{C} \notag
\end{align}
Hence from $(\ref{ct12})$ and $(\ref{cu9})$ we have
\begin{align} \label{cu11}
\mathbb{H}^1\left( X, \mathcal{E}_{\bigwedge^3 \Theta_X}^\bullet \right) \cong &\left\{ \left(  t_2\left(F_0 + F_1 \xi + F_1 \xi^2 \right) dz_1 +  \left( kt_2F_0 + t_1A +\left(kt_2 F_1 + t_1 B \right) \xi + \left(kt_2 F_2+t_1C\right) \xi^2\right) dz_2  \right) \right.\\
&\left.  \left(t_3 \frac{\partial}{\partial z_1} + t_4 \frac{\partial}{\partial z_2} + t_7\left(A+ B \xi + C \xi^2 \right)\frac{\partial}{\partial \xi} \right) d\bar{z}_1 + \left(t_5 \frac{\partial}{\partial z_1} + t_6\frac{\partial}{\partial z_2} + t_8\left(A+ B \xi + C \xi^2 \right)\frac{\partial}{\partial \xi} \right) d\bar{z}_2, \right. \notag \\
&\left. \,\,\,\,\,\,\,\,\,\left(- t_7 B- 2t_7 C \xi\right)d\bar{z}_1 + \left( - t_8 B - 2t_8 C \xi \right) d\bar{z}_2   \,\,\,|\,\,\,\,\,\,\,\,t_1,t_2,t_3,t_4,t_5,t_6,t_7,t_8\in \mathbb{C} \right\} \notag
\end{align}

We will construct a foliated analytic family of deformations of $\left(X, \omega_0 \right)$ where $\omega_0=Dd\xi+ \left( A+ B \xi + C \xi^2 \right)dz_1+ k\left( A + B \xi + C\xi^2\right)dz_2,D\ne 0, \left(A,B,C\right)\ne 0\in \mathbb{C}^3, k\in \mathbb{C}$ such that the foliated Kodaira-Spencer map in terms of cotangent sheaves is an isomorphism at the distinguished point, so that $\left( X , \omega_0 \right)$ is unobstructed in foliated deformations. In fact, we recall the integrability conditions for a foliated analytic family in terms of cotangent sheaves from $(\ref{oni1})-(\ref{oni4})$. If we take
\begin{align*}
\left( Y(t) \left\{ \left( \varphi(t), Z(t) \right) \right\} \right) \in  A^{0,0}\left( X, \Omega_X^1 \otimes \bigwedge^3 \Theta_X \right) \bigoplus A^{0,1}\left( X, \mathcal{E}_{\bigwedge^3 \Theta_X} \right)
\end{align*}
where 
{\small{\begin{align*}
Y(t)=&Dd\xi+ \left( A+t_2F_0 +\left( B+  t_2F_1\right) \xi + \left( C+ t_2F_2\right) \xi^2 \right) dz_1 \\
&+ \left( k\left(A+ t_2  F_0\right) +  t_1\left( A  + t_2 F_0\right)+\left(k \left(B+t_2 F_1\right)+ t_1\left(B + t_2 F_1\right)  \right)\xi + \left( k \left(C+ t_2 F_2\right) + t_1 \left(C + t_2 F_2 \right) \right)\xi^2             \right)dz_2
\end{align*}}}
\begin{align*}
\varphi(t)=&d\bar{z}_1 \wedge \left(t_3\frac{\partial}{\partial z_1} + t_4 \frac{\partial}{\partial z_2} + t_7\left(A+ t_2 F_0 + \left( B+ t_2 F_1 \right)\xi  +\left(C+ t_2 F_2 \right)\xi^2\right) \frac{\partial}{\partial \xi} \right)  \\
& +d\bar{z}_2 \wedge \left( t_5 \frac{\partial}{\partial z_1} + t_6 \frac{\partial}{\partial z_2}  + t_8 \left( A+ t_2 F_0 + \left( B+ t_2 F_1 \right)\xi  +\left(C+ t_2 F_2 \right)\xi^2\right) \frac{\partial}{\partial \xi} \right) 
\end{align*}
\begin{align*}
Z(t)= \left( -t_7 \left(B + t_2 F_1\right)  -2 t_7  \left(C   + t_2 F_2 \right)\xi \right)d\bar{z}_1 +\left(- t_8 \left(B + t_2 F_1\right) - 2t_8 \left( C + t_2 F_2\right) \xi  \right)d\bar{z}_2
\end{align*}
Then we have $\varphi(0)=0,Y(0)=\omega_0$ and $Z(0)=0$ and the following equalities hold: 
\begin{align*}
\bar{\partial}\varphi(t)- \frac{1}{2}\left[ \varphi(t), \varphi(t) \right] = 0\\
\bar{\partial} Y(t) - \mathcal{L}_{\varphi(t)}\left( Y(t) \right)- Z(t) \wedge Y(t)=0\\
Y(t)\wedge \partial Y(t)=0\\
\bar{\partial}Z(t)- \left[ \varphi(t), Z(t) \right]= 0
\end{align*}
In other words, the integrability conditions $(\ref{oni1})-(\ref{oni4})$ are satisfied, so that there exists a foliated analytic family $\left(\mathcal{M}, \omega , B, \pi \right)$ where $\pi:\left( \mathcal{M}, \omega \right)\to B$ for a sufficiently small neighborhood $B$ of $0$ in $\mathbb{C}^8$ with coordinates $t=\left(t_1,t_2,t_3,t_4,t_5,t_6,t_7,t_8 \right)$ such that $\pi^{-1}(0)=\left(X, \omega_0 \right)$ and the foliated Kodaira-Spencer map from $(\ref{cu14})$ given by
\begin{align*}
\varphi: T_0(B) &\to \mathbb{H}^1\left( X, \mathcal{E}_{\bigwedge^3 \Theta_X}^\bullet \right)\\
\frac{\partial}{\partial t_1} &\mapsto \left( \left( A+ B \xi + C \xi^2\right) dz_2, \left(0 , 0 \right) \right)\\
\frac{\partial}{\partial t_2} & \mapsto \left( \left( F_0 + F_1 \xi + F_2 \xi^2 \right) dz_1 + k\left( F_0 + F_1 \xi + F_2  \xi^2  \right) dz_2, (0,0) \right)\\
\frac{\partial}{\partial t_3} &\mapsto \left( 0, \left( \frac{\partial}{\partial z_1} d\bar{z}_1, 0 \right) \right)\\
\frac{\partial}{\partial t_4} & \mapsto \left( 0 , \left(\frac{\partial}{\partial z_2}d\bar{z}_1, 0 \right) \right)\\
\frac{\partial}{\partial t_5} & \mapsto \left( 0 , \left(\frac{\partial}{\partial z_1}d\bar{z}_2, 0 \right) \right)\\
\frac{\partial}{\partial t_6} & \mapsto \left( 0 , \left(\frac{\partial}{\partial z_2}d\bar{z}_2, 0 \right) \right)\\
\frac{\partial}{\partial t_7} &\mapsto \left(  0, \left( \left(A+ B \xi + C \xi^2\right)\frac{\partial}{\partial \xi} d\bar{z}_1, \left( - B - 2 C \xi \right) d\bar{z}_1  \right)\right) \\
\frac{\partial}{\partial t_8} & \mapsto \left( 0, \left( \left( A+ B \xi + C \xi^2 \right)\frac{\partial}{\partial \xi} d\bar{z}_2,  \left(-B- 2C \xi\right) d\bar{z}_2\right) \right)
\end{align*}
is an isomorphism at $t=0$ by $(\ref{cu11})$, so that $\left( X, \omega_0 \right)$ is unobstructed in foliated deformations.

\subsection{Foliated deformations of $\left(X, \omega_0\right)$ in terms of tangent sheaves}\

Let us consider the foliation on $X$ defined by $Dd\xi \in H^0\left(X, \Omega_X^1 \otimes \bigwedge^3 \Theta_X \right)$. Then $\Theta_{\mathcal{F}_0}\cong \mathcal{O}_X\bigoplus \mathcal{O}_X$ is generated by $T_1=\frac{\partial}{\partial z_1}$ and $T_2=\frac{\partial}{\partial z_2}\in H^0\left(X, \Theta_X \right)$ with $\left[ T_1, T_2 \right]=0$ and we have the exact sequence
\begin{align*}
0\to \mathcal{O}_X \bigoplus \mathcal{O}_X &\to \Theta_X \to \frac{\Theta_X}{\mathcal{O}_X \bigoplus \mathcal{O}_X}\cong \bigwedge^3 \Theta_X\cong \pi_2^* \Theta_{\mathbb{P}_\mathbb{C}^1} \to 0 \\
(1,0)&\mapsto T_1= \frac{\partial}{\partial z_1}\\
(0,1) & \mapsto T_2=\frac{\partial}{\partial z_2}
\end{align*}
On the other hand, let us consider the foliation on $X$ defined by $\left(X, Dd\xi + \left( A+ B \xi + C\xi^2 \right)dz_1 + k\left( A + B\xi + C\xi^2 \right) dz_2 \right)$, $D\ne 0, \left(A, B, C \right)\ne 0 \in \mathbb{C}^3 , k\in \mathbb{C}$. Then $\Theta_{\mathcal{F}_0}\cong \mathcal{O}_X \bigoplus \mathcal{O}_X$ is generated by $T_1= D\frac{\partial}{\partial z_1} -\left( A+ B\xi + C \xi^2 \right) \frac{\partial}{\partial \xi}$ and $T_2= D \frac{\partial}{\partial z_2} - k \left( A+ B\xi + C \xi^2 \right)\frac{\partial}{\partial \xi}$ with $\left[ T_1, T_2 \right]$ and we have the exact sequence
\begin{align*}
0\to \mathcal{O}_X \bigoplus \mathcal{O}_X &\to \Theta_X \to \frac{ \Theta_X }{\mathcal{O}_X \bigoplus \mathcal{O}_X} \cong \bigwedge^3 \Theta_X \cong \pi_2^* \Theta_{\mathbb{P}_\mathbb{C}^1} \to 0 \\
\left(1, 0\right) &\mapsto T_1= D\frac{\partial}{\partial z_1} -\left( A+ B\xi + C \xi^2 \right) \frac{\partial}{\partial \xi}\\
(0,1) & \mapsto T_2= D \frac{\partial}{\partial z_2} - k \left( A+ B\xi + C \xi^2 \right)\frac{\partial}{\partial \xi}
\end{align*}

Then for $\left(X, \omega_0\right)$ we have the leaf complex (see Part II)
{\small{\begin{align} \label{ct6}
\mathcal{E}_{\mathcal{O}_X \bigoplus \mathcal{O}_X}^\bullet :\mathcal{E}_{\mathcal{O}_X \bigoplus \mathcal{O}_X} \cong \Theta_X \bigoplus \mathscr{H}om_{\mathcal{O}_X}\left( \mathcal{O}_X \bigoplus \mathcal{O}_X, \mathcal{O}_X\bigoplus \mathcal{O}_X \right)\xrightarrow{\tilde{D}_0} \mathscr{H}om_{\mathcal{O}_X}\left(\mathcal{O}_X \bigoplus \mathcal{O}_X, \Theta_X \right) \\
\xrightarrow{\tilde{D}_1} \mathscr{H}om_{\mathcal{O}_X}\left( \bigwedge^2 \left(\mathcal{O}_X\bigoplus \mathcal{O}_X \right) , \frac{\Theta_X}{\mathcal{O}_X \bigoplus \mathcal{O}_X} \right)\to 0 \notag
\end{align}}}
Let us consider the spectral sequence at $E_1$ from Dolbeault bicomplex associated to $\mathcal{E}_{\mathcal{O}_X \bigoplus \mathcal{O}_X}^\bullet$:
\begin{equation}\label{ct3}
\begin{CD}
0@. \,\,\,\,\,\,\,\,\, 0 \\ 
@AAA @AAA \\
H^0\left(X, \mathscr{H}om_{\mathcal{O}_X}\left( \mathcal{O}_X, \pi_2^*  \Theta_{\mathbb{P}_\mathbb{C}^1}  \right) \right) @.\,\,\,\,\,\,\,\,\,\,\,H^1\left( X, \mathscr{H}om_{\mathcal{O}_X}\left( \mathcal{O}_X, \pi_2^* \Theta_{\mathbb{P}_\mathbb{C}^1}  \right)\right) \\
@A\tilde{D}_1AA @A\tilde{D}_1AA \\
H^0\left(X, \mathscr{H}om_{\mathcal{O}_X}\left( \mathcal{O}_X \bigoplus \mathcal{O}_X, \Theta_X \right) \right) @.\,\,\,\,\,\,\,\,\,\,\,\, H^1\left( X, \mathscr{H}om_{\mathcal{O}_X}\left( \mathcal{O}_X \bigoplus \mathcal{O}_X, \Theta_X  \right) \right) @.\,\,\,\,\,\,\,\,\,\, \\
@A\tilde{D}_0 AA @A\tilde{D}_0 AA \\
H^0\left(X, \Theta_X \right) \bigoplus H^0\left(X, \mathscr{H}om_{\mathcal{O}_X}\left( \mathcal{O}_X \bigoplus\mathcal{O}_X, \mathcal{O}_X \bigoplus \mathcal{O}_X \right)\right) @. \,\,\,\,\,\,\,\,\,\,\,\,\,H^1\left(X, \Theta_X \right)\bigoplus  H^1 \left(X, \mathscr{H}om_{\mathcal{O}_X}\left( \mathcal{O}_X \bigoplus\mathcal{O}_X, \mathcal{O}_X \bigoplus \mathcal{O}_X \right)\right)@.
\end{CD}
\end{equation}
Then we have
{\small{\begin{align}\label{ct5}
\mathbb{H}^1\left( X, \mathcal{E}_{\mathcal{O}_X \bigoplus \mathcal{O}_X}^\bullet \right) &\cong  \frac{\textnormal{ker}\left(H^0\left( X, \mathscr{H}om_{\mathcal{O}_X}\left( \mathcal{O}_X \bigoplus \mathcal{O}_X, \Theta_X \right) \right) \to H^0\left( X , \mathscr{H}om_{\mathcal{O}_X}\left( \mathcal{O}_X, \pi_2^*  \Theta_{\mathbb{P}_\mathbb{C}^1}  \right) \right) \right)}{\textnormal{im}\left( H^0\left( X, \Theta_X \right)\bigoplus H^0\left( X, \mathscr{H}om_{\mathcal{O}_X}\left( \mathcal{O}_X \bigoplus \mathcal{O}_X, \mathcal{O}_X \bigoplus \mathcal{O}_X \right) \to   H^0\left( X, \mathscr{H}om_{\mathcal{O}_X}\left( \mathcal{O}_X \bigoplus \mathcal{O}_X, \Theta_X \right) \right) \right) \right)} \\
 &\bigoplus \textnormal{ker}\left( H^1\left( X, \Theta_X \right)\bigoplus H^1 \left( X, \mathscr{H}om_{\mathcal{O}_X}\left( \mathcal{O}_X \bigoplus \mathcal{O}_X, \mathcal{O}_X \bigoplus \mathcal{O}_X \right) \to   H^1 \left( X, \mathscr{H}om_{\mathcal{O}_X}\left( \mathcal{O}_X \bigoplus \mathcal{O}_X, \Theta_X \right) \right) \right)  \right) \notag
\end{align}}}

We will describe each cohomology group in $(\ref{ct5})$. We note that $\dim_\mathbb{C} H^0\left( X, \mathscr{H}om_{\mathcal{O}_X}\left( \mathcal{O}_X \bigoplus \mathcal{O}_X, \mathcal{O}_X, \mathcal{O}_X \right)\right)= \dim_\mathbb{C} H^0\left( X, \bigoplus^4 \mathcal{O}_X \right)= 4$ in $(\ref{ct3})$ and any element $\psi$ of $H^0\left( X, \mathscr{H}om_{\mathcal{O}_X}\left( \mathcal{O}_X \bigoplus \mathcal{O}_X, \mathcal{O}_X \bigoplus \mathcal{O}_X \right) \right)$ has an unique representation of the form
\begin{align}\label{ct18}
\psi(T_1)&= f_1 T_1 + f_2 T_2\\
\psi(T_2)&= f_3 T_1 + f_3 T_2,\,\,\,\,\,f_1,f_2,f_3,f_4 \in \mathbb{C} \notag
\end{align}

We note that $\dim_\mathbb{C} H^1\left(X, \mathscr{H}om_{\mathcal{O}_X}\left(  \mathcal{O}_X \bigoplus \mathcal{O}_{\mathcal{O}_X}, \mathcal{O}_X\bigoplus \mathcal{O}_X\right) \right)=\dim_\mathbb{C} H^1\left( X,\bigoplus^4 \mathcal{O}_X\right) = 8$ and any element $\psi$ of $H^1\left( X, \mathscr{H}om_{\mathcal{O}_X}\left(  \mathcal{O}_X \bigoplus \mathcal{O}_X,  \mathcal{O}_X\bigoplus \mathcal{O}_X \right)\right)$ has an unique representation of the form
\begin{align}\label{ct22}
\psi\left(T_1 \right)&=  d\bar{z}_1\wedge \left(g_1 T_1+ g_2 T_2 \right) + d\bar{z}_2 \wedge \left( g_3 T_1 + g_4 T_2 \right) \\
\psi\left(T_2 \right)&= d\bar{z}_1\wedge \left( g_5 T_1+ g_6 T_2 \right) + d\bar{z}_2 \wedge \left( g_7 T_2 + g_8 T_2 \right) ,\,\,\,\,\,g_1,...,g_8\in \mathbb{C} \notag
\end{align}

We note that $\dim_\mathbb{C} H^0\left(X, \mathscr{H}om_{\mathcal{O}_X}\left( \mathcal{O}_X \bigoplus \mathcal{O}_X, \Theta_X \right) \right)= \dim_\mathbb{C} H^0\left(X, \Theta_X \bigoplus \Theta_X \right) = 10$ in $(\ref{ct3})$ and any element $\rho$ of $H^0\left(X, \mathscr{H}om_{\mathcal{O}_X}\left( \mathcal{O}_X \bigoplus \mathcal{O}_X, \Theta_X \right) \right)$ has an unique representation in the form
\begin{align}\label{ct20}
\rho(T_1)&= e_1\frac{\partial}{\partial z_1} + e_2 \frac{\partial}{\partial z_2} + \left( c_1 + c_2 \xi + c_3 \xi^2 \right)\frac{\partial}{\partial \xi} \\
\rho(T_1)&= e_3\frac{\partial}{\partial z_1} + e_4 \frac{\partial}{\partial z_2} + \left( c_4 + c_5 \xi + c_6 \xi^2 \right)\frac{\partial}{\partial \xi},\,\,\,\,\, e_1,...,e_4,c_1,...,c_6 \in \mathbb{C} \notag
\end{align}

We note that $\dim_\mathbb{C} H^1\left(X, \mathscr{H}om_{\mathcal{O}_X}\left( \mathcal{O}_X \bigoplus \mathcal{O}_X, \Theta_X \right) \right)=\dim_\mathbb{C} H^1\left(X, \Theta_X \bigoplus \Theta_X \right) = 20$ and any element $\sigma$ of $H^1\left(X, \mathscr{H}om_{\mathcal{O}_X}\left( \mathcal{O}_X \bigoplus \mathcal{O}_X, \Theta_X \right) \right)$ has an unique representation of the form
\begin{align} \label{ct2}
\sigma\left( T_1 \right)&= d\bar{z}_1\wedge\left( h_1 \frac{\partial}{\partial z_1} + h_2 \frac{\partial}{\partial z_2}  + \left( h_3 + h_4 \xi + h_5 \xi^2 \right)\frac{\partial}{\partial \xi} \right) + d \bar{z}_2 \wedge \left( h_6 \frac{\partial}{\partial z_1} + h_7 \frac{\partial}{\partial z_2}  + \left( h_8 + h_9 \xi + h_{10} \xi^2 \right)\frac{\partial}{\partial \xi} \right) \\
\sigma\left( T_2 \right)&= d\bar{z}_1\wedge\left( h_{11} \frac{\partial}{\partial z_1} + h_{12} \frac{\partial}{\partial z_2}  + \left( h_{13} + h_{14} \xi + h_{15} \xi^2 \right)\frac{\partial}{\partial \xi} \right) + d \bar{z}_2 \wedge \left( h_{16} \frac{\partial}{\partial z_1} + h_{17} \frac{\partial}{\partial z_2}  + \left( h_{18} + h_{19} \xi + h_{20} \xi^2 \right)\frac{\partial}{\partial \xi} \right) \notag
\end{align}
where $h_1,...,h_{20}\in \mathbb{C}$.

We note that $\dim_\mathbb{C} H^0\left( X, \mathscr{H}om_{\mathcal{O}_X}\left( \mathcal{O}_X, \pi_2^* \Theta_{\mathbb{P}_\mathbb{C}^1} \right) \right)= \dim_\mathbb{C} H^0\left( X, \pi_2^* \Theta_{\mathbb{P}_\mathbb{C}^1}\right) = 3$ in $(\ref{ct3})$ and any element $\eta$ of $H^0\left(X , \mathscr{H}om_{\mathcal{O}_X}\left( \mathcal{O}_X, \pi_2^* \Theta_{\mathbb{P}_\mathbb{C}^1} \right) \right)$ has an unique representation in the form
\begin{align*}
\eta\left( T_1 \wedge T_2 \right) = \left( s_1 + s_2 \xi + s_3 \xi^2 \right)\frac{\partial}{\partial \xi},\,\,\,\,\,\,s_1,s_2,s_3 \in \mathbb{C}
\end{align*}

We recall the integrability conditions for a foliated analytic family of deformations of $\left(X, \mathcal{O}_X\bigoplus \mathcal{O}_X \right)$ in terms of tangent sheave from Part VI-1. We note that for a sufficiently small neighborhood $B$ of $0\in \mathbb{C}^m$ with coordinates $t=(t_1,...,t_m)$,  $\left(X(t) , \left\{\left( \varphi(t), U(t) \right) \right\}  \right) \in A^{0,0}\left( X, \mathscr{H}om_{\mathcal{O}_X}\left( \mathcal{O}_X \bigoplus \mathcal{O}_X, \Theta_X  \right) \right)\bigoplus A^{0,1}\left(X, \mathcal{E}_{\mathcal{O}_X \bigoplus \mathcal{O}_X}\right)$ $= A^{0,0}\left(X, \mathscr{H}om_{\mathcal{O}_X}\left(\mathcal{O}_X \bigoplus \mathcal{O}_X, \Theta_X \right) \right) \bigoplus A^{0,1}\left(X, \Theta_X \right)\bigoplus A^{0,1}\left(X, \mathscr{H}om_{\mathcal{O}_X}\left( \mathcal{O}_X \bigoplus \mathcal{O}_X, \mathcal{O}_X\bigoplus \mathcal{O}_X \right) \right)$ with $\varphi(t)=0, X(0)=\textnormal{Id}, U(0)=0$ defines a foliated analytic family $\left( \mathcal{M}, \Theta_\mathcal{F}, B, \pi \right)$ of deformations of $\pi^{-1}(0)=\left(X, \mathcal{O}_X \bigoplus \mathcal{O}_X \right)$ such that the foliated Kodaira-Spencer map at $t=0$ is given by
\begin{align} \label{cu33}
\varphi: T_0(B) &\to \mathbb{H}^1\left( X, \mathcal{E}_{\mathcal{O}_X \bigoplus \mathcal{O}_X}^\bullet \right)\\
\frac{\partial}{\partial t} &\mapsto \left( \frac{\partial X (t)}{\partial t}|_{t=0}, \left( \frac{\partial \varphi(t)}{\partial t}|_{t=0}, \frac{\partial U(t)}{\partial t}|_{t=0} \right)  \right) \notag
\end{align}
if the following integrability condition holds: if we set $U(t)\left(T_1\right)= U^{11}(t) T_1 + U^{12}(t) T_2$  and $U(t)\left(T_2\right)= U^{21}(t) T_1 + U^{22}(t) T_2$ for some $U^{\alpha\beta}(t)\in A^{0,1}\left( X, \mathcal{O}_X\right)$ for $\alpha,\beta=1,2$, then
\begin{align}
&\bar{\partial}\varphi(t) - \frac{1}{2}\left[ \varphi(t), \varphi(t) \right]=0 \label{cu31}\\
&\left[ X(t)\left( T_1 \right), X(t)\left( T_2 \right) \right] = K_1 (t) X(t)\left(T_1 \right)  + K_2 (t) X(t)\left(T_2 \right)\,\,\,\,\,\textnormal{for $K_1(t),K_2(t)\in A^{0,0}\left(X, \mathcal{O}_X \right)$}\\
&\bar{\partial} X(t)\left(T_1\right)-\left[ \varphi(t), X(t)\left(T_1\right)\right]+ U(t)\left(T_1\right)=0, \,\,\,\,\,\,\,\,\bar{\partial} X(t)\left(T_2\right)- \left[ \varphi(t), X(t)\left(T_2\right)\right] + U(t)\left(T_2\right)=0\\
&\bar{\partial} U^{\alpha\beta}(t) - \left[ \varphi(t), U^{\alpha\beta}(t)\right] + \sum_{\sigma=1}^2 U^{\alpha \sigma}(t) \wedge U^{\sigma \beta}(t)=0,\,\,\,\,\,\,\,\,\alpha,\beta=1,2 \label{cu32}
\end{align}

With this preparation, we will study foliated deformations of $\left( X, \omega_0 \right)$ in terms of tangent sheaves where $\omega_0 \in H^0 \left( X , \Omega_X^1 \otimes \bigwedge^3 \Theta_X \right)$ with $\omega_0 \wedge d \omega_0 = 0$

\subsection{Foliated deformations of $\left( X, Dd \xi\right),D\ne 0$ in terms of tangent sheaves}\

Let us consider the foliation on $X$ defined by $Dd\xi \in H^0\left(X, \Omega_X^1 \otimes \bigwedge^3 \Theta_X \right)$. Then $\Theta_{\mathcal{F}_0}\cong \mathcal{O}_X\bigoplus \mathcal{O}_X$ is generated by $T_1=\frac{\partial}{\partial z_1}$ and $T_2=\frac{\partial}{\partial z_2}\in H^0\left(X, \Theta_X \right)$ with $\left[ T_1, T_2 \right]=0$. We compute $\mathbb{H}^1\left( X, \mathcal{E}_{\mathcal{O}_X \bigoplus \mathcal{O}_X}^\bullet  \right)$. From $(\ref{ct5})$, first we compute $H^0\left( X, \Theta_X \right)\bigoplus H^0\left(X, \mathscr{H}om_{\mathcal{O}_X}\left( \mathcal{O}_X \bigoplus \mathcal{O}_X, \mathcal{O}_X\bigoplus \mathcal{O}_X \right) \right) \xrightarrow{\tilde{D}_0} H^0\left(X, \mathscr{H}om_{\mathcal{O}_X}\left(\mathcal{O}_X\bigoplus \mathcal{O}_X, \Theta_X \right) \right)$. Then from $(\ref{ct8})$ and $(\ref{ct18})$ the image of $\left( V , \psi \right)$ is given by $\rho=\tilde{D}_0 \left( V, \psi \right)$ where
\begin{align*}
\rho\left( T_1\right) =\left[ a_1\frac{\partial}{\partial z_1} + a_2 \frac{\partial}{\partial z_2} + \left( a_3 + a_4 \xi + a_5 \xi^2 \right)\frac{\partial}{\partial \xi} , \frac{\partial}{\partial z_1}\right] - \left( f_1 T_1 + f_2 T_2 \right) = -f_1 T_1 - f_2 T_2\\
\rho\left( T_2\right) =\left[ a_1\frac{\partial}{\partial z_1} + a_2 \frac{\partial}{\partial z_2} + \left( a_3 + a_4 \xi + a_5 \xi^2 \right)\frac{\partial}{\partial \xi} , \frac{\partial}{\partial z_2}\right] - \left( f_3 T_1 + f_4 T_2 \right) = -f_3T_1 - f_4T_2
\end{align*}
Then from $(\ref{ct20})$ we have $\dim_\mathbb{C} \mathbb{H}^0\left( X, \mathcal{E}_{\mathcal{O}_X \bigoplus \mathcal{O}_X}^\bullet \right)=\dim_\mathbb{C} \textnormal{ker}\left( \tilde{D}_1 \right)=5$ and so
{\small{\begin{align}\label{ct21}
\dim_\mathbb{C} \textnormal{im}\left( H^0\left( X, \Theta_X \right)\bigoplus H^0\left(X, \mathscr{H}om_{\mathcal{O}_X}\left( \mathcal{O}_X \bigoplus \mathcal{O}_X, \mathcal{O}_X\bigoplus \mathcal{O}_X \right) \right) \xrightarrow{\tilde{D}_0} H^0\left(X, \mathscr{H}om_{\mathcal{O}_X}\left(\mathcal{O}_X\bigoplus \mathcal{O}_X, \Theta_X \right) \right) \right)= 9 - 5= 4
\end{align}}}
 Next we compute $H^0\left( X, \mathscr{H}om_{\mathcal{O}_X} \left( \mathcal{O}_X \bigoplus \mathcal{O}_X , \Theta_X \right) \right) \xrightarrow{\tilde{D}_1} H^0\left( X, \mathscr{H}om_{\mathcal{O}_X} \left( \mathcal{O}_X, \pi_2^* \Theta_{\mathbb{P}_\mathbb{C}^1} \right)\right)$. We recall that the map $\tilde{D}_1$ is induced by the composition (see Part II)
{\small{\begin{align*}
H^0\left( X, \mathscr{H}om_{\mathcal{O}_X} \left( \mathcal{O}_X \bigoplus \mathcal{O}_X , \Theta_X \right) \right)& \xrightarrow{q}H^0\left( X, \mathscr{H}om_{\mathcal{O}_X}\left( \mathcal{O}_X \bigoplus \mathcal{O}_X, \frac{\Theta_X}{\mathcal{O}_X\bigoplus \mathcal{O}_X} \right) \right) \xrightarrow{D_1} H^0\left( X, \mathscr{H}om_{\mathcal{O}_X} \left( \mathcal{O}_X, \pi_2^* \Theta_{\mathbb{P}_\mathbb{C}^1} \right)\right)\\
\rho &\mapsto \bar{\rho}
\end{align*}}}
where $q$ is defined by the natural quotient $\Theta_X \to \frac{\Theta_X}{\mathcal{O}_X \bigoplus \mathcal{O}_X }\cong \pi_2^* \Theta_{\mathbb{P}_\mathbb{C}^1}$. Then since $\frac{\partial}{\partial z_1}\equiv 0$ and $\frac{\partial}{\partial z_2}\equiv 0$ in $\frac{\Theta_X}{\mathcal{O}_X\bigoplus \mathcal{O}_X}$. Then we have from $(\ref{ct20})$
\begin{align*}
\bar{\rho}\left( T_1\right) &\equiv \left( c_1 + c_2 \xi + c_3 \xi^2 \right) \frac{\partial}{\partial \xi} \\
\bar{\rho}\left( T_2 \right) &\equiv \left( c_4 + c_5 \xi + c_6 \xi^2 \right) \frac{\partial}{\partial \xi}
\end{align*}
Then we have
\begin{align*}
\tilde{D}_1\left(\rho \right)\left( T_1 \wedge T_2 \right)&= \left[ T_1, \bar{\rho}\left(T_2 \right)\right] - \left[ T_2, \bar{\rho}\left( T_1 \right) \right] - \bar{\rho}\left( \left[ T_1, T_2 \right] \right)\\
 &= \left[ \frac{\partial}{\partial z_1}, \left(c_4 + c_5 \xi + c_6 \xi^2 \right)\frac{\partial}{\partial \xi} \right) - \left[ \frac{\partial}{\partial z_2}, \left( c_1 + c_2 \xi + c_3 \xi^2  \right) \frac{\partial}{\partial \xi} \right] =0,
\end{align*}
so that $\dim_\mathbb{C} \textnormal{ker}\left( H^0\left( X, \mathscr{H}om_{\mathcal{O}_X} \left( \mathcal{O}_X \bigoplus \mathcal{O}_X , \Theta_X \right) \right)\to H^0\left( X, \mathscr{H}om_{\mathcal{O}_X} \left( \mathcal{O}_X, \pi_2^* \Theta_{\mathbb{P}_\mathbb{C}^1} \right)\right) \right)= 10 $. Then we have from $(\ref{ct21})$
\begin{align}\label{ct19}
\dim_\mathbb{C} \frac{\textnormal{ker}\left( H^0\left( X, \mathscr{H}om_{\mathcal{O}_X} \left( \mathcal{O}_X \bigoplus \mathcal{O}_X , \Theta_X \right) \right)\to H^0\left( X, \mathscr{H}om_{\mathcal{O}_X} \left( \mathcal{O}_X, \pi_2^* \Theta_{\mathbb{P}_\mathbb{C}^1} \right)\right) \right)}{\textnormal{im}\left( H^0\left( X, \Theta_X \right)\bigoplus H^0\left( X, \mathscr{H}om_{\mathcal{O}_X}\left( \mathcal{O}_X \bigoplus \mathcal{O}_X , \mathcal{O}_X \bigoplus \mathcal{O}_X \right)\right) \to H^0\left(X, \mathscr{H}om_{\mathcal{O}_X}\left( \mathcal{O}_X \bigoplus \mathcal{O}_X, \Theta_X  \right) \right) \right)} = 10 - 4 = 6
\end{align}

Let us compute $H^1\left( X, \Theta_X \right)\bigoplus H^1 \left( X, \mathscr{H}om_{\mathcal{O}_X}\left( \mathcal{O}_X \bigoplus \mathcal{O}_X, \mathcal{O}_X \bigoplus \mathcal{O}_X \right) \right) \xrightarrow{\tilde{D}_0}   H^1 \left( X, \mathscr{H}om_{\mathcal{O}_X}\left( \mathcal{O}_X \bigoplus \mathcal{O}_X, \Theta_X \right) \right)  $. Then the image of $\left( \phi, \psi \right)$ from $(\ref{ct9})$ and $(\ref{ct22})$ is given by $\sigma=\tilde{D}_0\left( \phi, \psi \right)$ where
\begin{align*}
\sigma(T_1)&=d\bar{z}_1\wedge \left[ r_0 \frac{\partial}{\partial z_1} + r_1 \frac{\partial}{\partial z_2} + \left( a_0 + a_1 \xi + a_2 \xi^2 \right) \frac{\partial}{\partial \xi }, \frac{\partial}{\partial z_1} \right] + d\bar{z}_2 \left[ r_2\frac{\partial}{\partial z_1} + r_3 \frac{\partial}{\partial z_2} + \left( b_0 + b_1 \xi + b_2 \xi^2 \right) \frac{\partial}{\partial \xi} , \frac{\partial}{\partial z_1} \right]  \\
&- d\bar{z}_1\wedge\left( g_1T_1 + g_2 T_2 \right) - d \bar{z}_2 \wedge \left( g_3 T_1 + g_4 T_2 \right)  \\
\sigma(T_2)&=d\bar{z}_1\wedge \left[ r_0 \frac{\partial}{\partial z_1} + r_1 \frac{\partial}{\partial z_2} + \left( a_0 + a_1 \xi + a_2 \xi^2 \right) \frac{\partial}{\partial \xi }, \frac{\partial}{\partial z_2} \right] + d\bar{z}_2 \left[ r_2\frac{\partial}{\partial z_1} + r_3 \frac{\partial}{\partial z_2} + \left( b_0 + b_1 \xi + b_2 \xi^2 \right) \frac{\partial}{\partial \xi} , \frac{\partial}{\partial z_2} \right] \\
&- d\bar{z}_1\wedge\left( g_5 T_1 + g_6 T_2 \right) - d \bar{z}_2 \wedge \left( g_7  T_2 + g_8 T_2  \right)  
\end{align*}
Then from $(\ref{ct2})$ we see that the kernel is $H^1\left(X, \Theta_X \right)$. Then we have from $(\ref{ct5})$ and $(\ref{ct19})$
\begin{align*}
\dim_\mathbb{C} \mathbb{H}^1\left(  X ,  \mathcal{E}_{\mathcal{O}_X \bigoplus \mathcal{O}_X}^\bullet \right) = 6 + 10 = 16
\end{align*}

Moreover, we note that the natural map $\mathbb{H}^1\left( X , \mathcal{E}_{\mathcal{O}_X \bigoplus \mathcal{O}_X}^\bullet \right)\to H^1\left( X, \Theta_X \right)$ is surjective. Hence since $X$ is obstructed in deformations of complex structures (see \cite{Kod58} p.436), $\left(X, \mathcal{O}_X \bigoplus \mathcal{O}_X\right)$ is obstructed in foliated deformations.

\subsection{Foliated deformations of $\left(X, Dd\xi + \left( A+ B \xi + C\xi^2 \right)dz_1 + k\left( A + B\xi  +  C\xi^2 \right) dz_2 \right), D\ne 0, \left(A, B, C \right)\ne 0 \in \mathbb{C}^2 , k\in \mathbb{C}$ in terms of tangent sheaves} \

Let us consider the foliation on $X$ defined by $\left(X, Dd\xi + \left( A+ B \xi + C\xi^2 \right)dz_1 + k\left( A + B\xi + C\xi^2 \right) dz_2 \right), D\ne 0, \left(A, B, C \right)\ne 0 \in \mathbb{C}^3 , k\in \mathbb{C}$. Then $\Theta_{\mathcal{F}_0}\cong \mathcal{O}_X \bigoplus \mathcal{O}_X$ is generated by $T_1= D\frac{\partial}{\partial z_1} -\left( A+ B\xi + C \xi^2 \right) \frac{\partial}{\partial \xi}$ and $T_2= D \frac{\partial}{\partial z_2} - k \left( A+ B\xi + C \xi^2 \right)\frac{\partial}{\partial \xi}\in H^0\left( X, \Theta_X \right)$ with $\left[ T_1, T_2 \right]=0$.
We compute $\mathbb{H}^1\left( X, \mathcal{E}_{\mathcal{O}_X \bigoplus \mathcal{O}_X}^\bullet  \right)$. From $(\ref{ct5})$, first we compute $H^0\left( X, \Theta_X \right)\bigoplus H^0\left(X, \mathscr{H}om_{\mathcal{O}_X}\left( \mathcal{O}_X \bigoplus \mathcal{O}_X, \mathcal{O}_X\bigoplus \mathcal{O}_X \right) \right) \xrightarrow{\tilde{D}_0} H^0\left(X, \mathscr{H}om_{\mathcal{O}_X}\left(\mathcal{O}_X\bigoplus \mathcal{O}_X, \Theta_X \right) \right)$. Then from $(\ref{ct8})$ and $(\ref{ct18})$ the image of $\left( V , \psi \right)$ is given by $\rho=\tilde{D}_0 \left( V, \psi \right)$ where
\begin{align}
\rho\left( T_1\right) &=\left[ a_1\frac{\partial}{\partial z_1} + a_2 \frac{\partial}{\partial z_2} + \left( a_3 + a_4 \xi + a_5 \xi^2 \right)\frac{\partial}{\partial \xi} ,  D\frac{\partial}{\partial z_1} -\left( A+ B\xi + C \xi^2 \right) \frac{\partial}{\partial \xi}\right] \notag \\
& - \left( f_1 \left( D\frac{\partial}{\partial z_1} -\left( A+ B\xi + C \xi^2 \right) \frac{\partial}{\partial \xi}\right)+ f_2 \left(  D\frac{\partial}{\partial z_2} -k\left( A+ B\xi + C \xi^2 \right) \frac{\partial}{\partial \xi} \right) \right) \notag \\
&= - f_1 D \frac{\partial}{\partial z_1} - f_2 D \frac{\partial}{\partial z_2} \label{cu24} \\
&  +\left( \left(f_1 + f_2k\right)A+ a_4 A  - a_3 B +\left(\left(f_1+ f_2k\right)B +2a_5 A- 2a_3 C \right) \xi + \left( \left(f_1+ f_2k \right)C+ a_5B - a_4 C\right) \xi^2  \right) \frac{\partial}{\partial \xi} \notag \\
\rho\left( T_2\right)& =\left[ a_1\frac{\partial}{\partial z_1} + a_2 \frac{\partial}{\partial z_2} + \left( a_3 + a_4 \xi + a_5 \xi^2 \right)\frac{\partial}{\partial \xi} ,  D\frac{\partial}{\partial z_2} - k\left( A+ B\xi + C \xi^2 \right) \frac{\partial}{\partial \xi}\right] \notag\\
& - \left( f_3 \left(  D\frac{\partial}{\partial z_1} -\left( A+ B\xi + C \xi^2 \right) \frac{\partial}{\partial \xi} \right)+ f_4 \left(  D\frac{\partial}{\partial z_2} - k \left( A+ B\xi + C \xi^2 \right) \frac{\partial}{\partial \xi} \right) \right) \notag\\
&= - f_3 D \frac{\partial}{\partial z_1} - f_4 D \frac{\partial}{\partial z_2} \label{cu25} \\
&  +\left( \left(f_3 + f_4k\right)A+ k\left( a_4 A  - a_3 B\right) +\left(\left(f_3+ f_4k\right)B + k\left(2a_5 A- 2a_3 C\right) \right) \xi + \left( \left(f_3+ f_4k \right)C+ k\left(a_5B - a_4 C\right)\right) \xi^2  \right) \frac{\partial}{\partial \xi} \notag
\end{align}
Then from $(\ref{ct20})$ we have $\dim_\mathbb{C} \mathbb{H}^0\left( X, \mathcal{E}_{\mathcal{O}_X \bigoplus \mathcal{O}_X}^\bullet \right)=\dim_\mathbb{C} \textnormal{ker}\left( \tilde{D}_1 \right)=3$ and so
{\small{\begin{align}
\dim_\mathbb{C} \textnormal{im}\left( H^0\left( X, \Theta_X \right)\bigoplus H^0\left(X, \mathscr{H}om_{\mathcal{O}_X}\left( \mathcal{O}_X \bigoplus \mathcal{O}_X, \mathcal{O}_X\bigoplus \mathcal{O}_X \right) \right) \xrightarrow{\tilde{D}_0} H^0\left(X, \mathscr{H}om_{\mathcal{O}_X}\left(\mathcal{O}_X\bigoplus \mathcal{O}_X, \Theta_X \right) \right) \right)= 9 - 3 = 6
\end{align}}}
 Next we compute $H^0\left( X, \mathscr{H}om_{\mathcal{O}_X} \left( \mathcal{O}_X \bigoplus \mathcal{O}_X , \Theta_X \right) \right) \xrightarrow{\tilde{D}_1} H^0\left( X, \mathscr{H}om_{\mathcal{O}_X} \left( \mathcal{O}_X, \pi_2^* \Theta_{\mathbb{P}_\mathbb{C}^1} \right)\right)$. We recall that the map $\tilde{D}_1$ is induced by the composition (see Part II)
{\small{\begin{align*}
H^0\left( X, \mathscr{H}om_{\mathcal{O}_X} \left( \mathcal{O}_X \bigoplus \mathcal{O}_X , \Theta_X \right) \right)& \xrightarrow{q}H^0\left( X, \mathscr{H}om_{\mathcal{O}_X}\left( \mathcal{O}_X \bigoplus \mathcal{O}_X, \frac{\Theta_X}{\mathcal{O}_X\bigoplus \mathcal{O}_X} \right) \right) \xrightarrow{D_1} H^0\left( X, \mathscr{H}om_{\mathcal{O}_X} \left( \mathcal{O}_X, \pi_2^* \Theta_{\mathbb{P}_\mathbb{C}^1} \right)\right)\\
\rho &\mapsto \bar{\rho}
\end{align*}}}
where $q$ is defined by the natural quotient $\Theta_X \to \frac{\Theta_X}{\mathcal{O}_X \bigoplus \mathcal{O}_X }\cong \pi_2^* \Theta_{\mathbb{P}_\mathbb{C}^1}$. Then since $\frac{\partial}{\partial z_1}\equiv \frac{1}{D}\left( A + B \xi + C \xi^2 \right)\frac{\partial}{\partial \xi}$ and $\frac{\partial}{\partial z_2}\equiv \frac{k}{D}\left( A + B \xi + C \xi^2 \right) \frac{\partial}{\partial \xi}$ in $\frac{\Theta_X}{\mathcal{O}_X\bigoplus \mathcal{O}_X}$. Then we have from $(\ref{ct20})$
\begin{align*}
\bar{\rho}\left(T_1\right)= e_1\frac{\partial}{\partial z_1} + e_2 \frac{\partial}{\partial z_2}+(c_1+c_2\xi+ c_3 \xi^2)\frac{\partial}{\partial \xi} \equiv \left( \frac{e_1 + ke_2}{D}\left( A+ B\xi + C \xi^2\right) +  c_1 + c_2 \xi + c_3\xi^2  \right) \frac{\partial}{\partial \xi}\\
\bar{\rho}\left(T_2\right)= e_3\frac{\partial}{\partial z_1} + e_4 \frac{\partial}{\partial z_2}+(c_4+c_5\xi+ c_6 \xi^2)\frac{\partial}{\partial \xi} \equiv \left(  \frac{e_3 + ke_4}{D}\left( A+ B\xi + C \xi^2\right) +  c_4 + c_5 \xi + c_6\xi^2  \right) \frac{\partial}{\partial \xi}
\end{align*}
Then we have
{\Small{\begin{align*}
&D_1\left(\bar{\rho}\right)\left(T_1 \wedge T_2  \right)=\left[ T_1, \bar{\rho}\left(T_2\right)\right]- \left[T_2, \bar{\rho}\left(T_1\right)\right]- \bar{\rho}\left( \left[T_1,T_2\right] \right) = \left[ T_1, \bar{\rho} \left(T_2\right) \right] - \left[T_2, \bar{\rho}\left(T_1\right)\right]\\
&= \left[ D\frac{\partial}{\partial z_1} - (A+ B\xi + C \xi^2)\frac{\partial}{\partial \xi} , \left( \frac{e_3 + ke_4}{D}\left( A+ B\xi + C \xi^2\right)+c_4 + c_5 \xi+ c_6 \xi^2\right) \frac{\partial}{\partial \xi} \right] \\
&- \left[ D\frac{\partial}{\partial z_2} - k(A+B + C\xi)\frac{\partial}{\partial \xi}, \left( \frac{e_1 + ke_2}{D}\left( A+ B\xi + C \xi^2\right) + c_1+ c_2\xi+ c_3 \xi^2 \right) \frac{\partial}{\partial \xi} \right]\\
&=\left( B \left(c_4-kc_1\right) - A\left(c_5-kc_2\right) +2\left(C \left(c_4-kc_1\right)- A \left(c_6-kc_3\right) \right)\xi + 2 \left( C\left(c_5- kc_2\right) -B \left(c_6-kc_3\right) \right) \xi^2 \right)\frac{\partial}{\partial \xi}=0\\
\end{align*}}}
Then the kernel is given by
\begin{align*}
B\left(c_4- kc_1\right)- A\left(c_5- kc_2\right)=0\\
C \left(c_4- kc_1 \right)- A \left(c_6- kc_3 \right)=0\\
C\left(c_5-kc_2\right)-B \left(c_6- k c_3\right)=0
\end{align*}
Then we have
\begin{align}
c_4- kc_1&= t_1 A \Longrightarrow c_4= t_1 A+ kc_1 \notag\\
c_5- kc_2&= t_1 B \Longrightarrow c_5= t_1 B + kc_2 \\
c_6- kc_3&= t_1C \Longrightarrow c_6= t_1 C + kc_3,\,\,\,\,\,\,\,t_1\in \mathbb{C} \notag
\end{align}
Hence $\dim_\mathbb{C} \textnormal{ker}\left(H^0\left( X, \mathscr{H}om_{\mathcal{O}_X} \left( \mathcal{O}_X \bigoplus \mathcal{O}_X , \Theta_X \right) \right) \xrightarrow{\tilde{D}_1} H^0\left( X, \mathscr{H}om_{\mathcal{O}_X} \left( \mathcal{O}_X, \pi_2^* \Theta_{\mathbb{P}_\mathbb{C}^1} \right)\right) \right)=8$ and any element $\rho$ of the kernel has an unique representation of the form
\begin{align}
\rho\left(T_1\right)&=e_1\frac{\partial}{\partial z_1}+ e_2 \frac{\partial}{\partial z_2} + \left(c_1+ c_2\xi+ c_3\xi^2 \right)\frac{\partial}{\partial \xi}\\
\rho\left(T_1\right)&=e_3\frac{\partial}{\partial z_1}+ e_4 \frac{\partial}{\partial z_2} + \left(t_1A+ kc_1+ \left(t_1B+ kc_2\right)\xi+ \left(t_1C+ kc_3\right) \xi^2 \right)\frac{\partial}{\partial \xi},\,\,\,\,\,e_1,e_2,e_3,e_4, c_1,c_2,c_3, t_1 \in \mathbb{C} \notag
\end{align}
Then we have $\dim_\mathbb{C} \frac{\textnormal{ker}\left( H^0\left( X, \mathscr{H}om_{\mathcal{O}_X} \left( \mathcal{O}_X \bigoplus \mathcal{O}_X , \Theta_X \right) \right) \xrightarrow{\tilde{D}_1} H^0\left( X, \mathscr{H}om_{\mathcal{O}_X} \left( \mathcal{O}_X, \pi_2^* \Theta_{\mathbb{P}_\mathbb{C}^1} \right)\right) \right)}{\textnormal{im}\left( H^0\left( X, \Theta_X \right)\bigoplus H^0\left(X, \mathscr{H}om_{\mathcal{O}_X}\left( \mathcal{O}_X \bigoplus \mathcal{O}_X, \mathcal{O}_X\bigoplus \mathcal{O}_X \right) \right) \xrightarrow{\tilde{D}_0} H^0\left(X, \mathscr{H}om_{\mathcal{O}_X}\left(\mathcal{O}_X\bigoplus \mathcal{O}_X, \Theta_X \right) \right)\right)}= 8 -6 =2$. We recall that from $(\ref{cu24})$ and $(\ref{cu25})$ the image of $\tilde{D}_0$ is given by
\begin{align}
\rho\left( T_1\right) &= - f_1 D \frac{\partial}{\partial z_1} - f_2 D \frac{\partial}{\partial z_2} \label{cu26} \\
&  +\left( \left(f_1 + f_2k\right)A+ a_4 A  - a_3 B +\left(\left(f_1+ f_2k\right)B +2a_5 A- 2a_3 C \right) \xi + \left( \left(f_1+ f_2k \right)C+ a_5B - a_4 C\right) \xi^2  \right) \frac{\partial}{\partial \xi} \notag  \\
\rho\left( T_2\right)&= - f_3 D \frac{\partial}{\partial z_1} - f_4 D \frac{\partial}{\partial z_2} \label{cu27}  \\
&  +\left( \left(f_3 + f_4k\right)A+ k\left(a_4 A  - a_3 B\right) +\left(\left(f_3+ f_4k\right)B + k\left(2a_5 A- 2a_3 C \right) \right) \xi + \left( \left(f_3+ f_4k \right)C+  k\left( a_5B - a_4 C \right) \right) \xi^2  \right) \frac{\partial}{\partial \xi} \notag
\end{align}
By setting $f_1=f_2=f_3=f_4=0$, $(\ref{cu26})$ and $(\ref{cu27})$ become
\begin{align}
\rho\left( T_1\right) &= \left(  a_4 A  - a_3 B + \left(2a_5 A- 2a_3 C \right) \xi + \left( a_5B - a_4 C\right) \xi^2  \right) \frac{\partial}{\partial \xi}  \label{uc28}  \\
\rho\left( T_2\right)&= \left(  k\left(a_4 A  - a_3 B\right) + k\left( 2a_5 A- 2a_3 C \right) \xi + k\left( a_5B - a_4 C\right) \xi^2  \right) \frac{\partial}{\partial \xi} 
\end{align}
Since $\textnormal{rank} \left[\begin{matrix}
-B & A & 0\\
-2C & 0 & 2 A\\
0 & -C & B
\end{matrix}\right]=2$, there is $\left(F_0,F_1,F_2\right)\in \mathbb{C}^3$ such that $\left( F_1,F_2,F_3 \right)$ is not in the image of the matrix. We note that $\rho_1\in H^0\left(  X, \mathscr{H}om_{\mathcal{O}_X}\left( \mathcal{O}_X \bigoplus \mathcal{O}_X, \Theta_X \right)\right)$ defined by $\rho_1\left( T_1\right) = \left( F_0 + F_1 \xi + F_2 \xi^2\right)\frac{\partial}{\partial \xi}$ and $\rho_1\left(T_2  \right) = k \left( F_0 + F_1 \xi + F_2 \xi^2\right)\frac{\partial}{\partial \xi}$ is in the kernel of $\tilde{D}_1$. It is clear that $\rho_2\in H^0\left(  X, \mathscr{H}om_{\mathcal{O}_X}\left( \mathcal{O}_X \bigoplus \mathcal{O}_X, \Theta_X \right)\right)$ defined by $\rho_2 \left( T_1 \right) = 0$ and $\rho_2\left( T_2 \right)=\left( A + B \xi + C\xi^2 \right) \frac{\partial}{\partial \xi}$ is not in the image of $\tilde{D}_0$. Then we claim that
\begin{align}\label{cu30}
\frac{\textnormal{ker}\left( H^0\left( X, \mathscr{H}om_{\mathcal{O}_X} \left( \mathcal{O}_X \bigoplus \mathcal{O}_X , \Theta_X \right) \right) \xrightarrow{\tilde{D}_1} H^0\left( X, \mathscr{H}om_{\mathcal{O}_X} \left( \mathcal{O}_X, \pi_2^* \Theta_{\mathbb{P}_\mathbb{C}^1} \right)\right) \right)}{\textnormal{im}\left( H^0\left( X, \Theta_X \right)\bigoplus H^0\left(X, \mathscr{H}om_{\mathcal{O}_X}\left( \mathcal{O}_X \bigoplus \mathcal{O}_X, \mathcal{O}_X\bigoplus \mathcal{O}_X \right) \right) \xrightarrow{\tilde{D}_0} H^0\left(X, \mathscr{H}om_{\mathcal{O}_X}\left(\mathcal{O}_X\bigoplus \mathcal{O}_X, \Theta_X \right) \right)\right)} \cong \textnormal{Span}_\mathbb{C}\left\langle \rho_1, \rho_2 \right\rangle
\end{align}
In fact, assume that $r_1\rho_1+ r_2 \rho_2\equiv 0$ for some $r_1,r_2\in \mathbb{C}$, i.e. $r_1\rho_1 + r_2 \rho_2 \in \textnormal{im}\left( \tilde{D}_0 \right)$. Then $r_1\rho(T_1)= r_1\left( F_0 + F_1 \xi + F_2 \xi^2 \right) \frac{\partial}{\partial \xi}$ is of the form $(\ref{uc28})$, so that by our choice of $\left(F_0,F_1,F_2 \right)$, $r_1$ should be $0$. Then $r_2(T_2)= r_2\left( A+ B\xi + C \xi^2\right)\frac{\partial}{\partial \xi}=0$, so that $r_2=0$. This proves the claim. On the other hand, let us compute $H^1\left( X, \Theta_X \right)\bigoplus H^1\left(X, \mathscr{H}om_{\mathcal{O}_X}\left( \mathcal{O}_X \bigoplus \mathcal{O}_X, \mathcal{O}_X\bigoplus \mathcal{O}_X \right) \right) \xrightarrow{\tilde{D}_0} H^1\left(X, \mathscr{H}om_{\mathcal{O}_X}\left(\mathcal{O}_X\bigoplus \mathcal{O}_X, \Theta_X \right) \right)$. From $(\ref{ct9})$ and $(\ref{ct22})$, the image of $\sigma=\tilde{D}_0 \left( \phi, \psi\right)$ is given by
{\tiny{\begin{align*}
\sigma\left( T_1 \right)&= d\bar{z}_1\wedge \left[ r_0 \frac{\partial}{\partial z_1} + r_1 \frac{\partial}{\partial z_2} + \left( a_0 + a_1 \xi + a_2 \xi^2 \right)\frac{\partial}{\partial \xi}, D\frac{\partial}{\partial z_1} -\left( A+ B \xi + C \xi^2 \right) \frac{\partial}{\partial \xi} \right] \\
&- d\bar{z}_1\left( g_1\left( D \frac{\partial}{\partial z_1} - \left( A + B \xi + C \xi^2 \right) \frac{\partial}{\partial \xi} \right) + g_2\left( D \frac{\partial}{\partial z_2} - k \left( A + B \xi + C \xi^2 \right)\frac{\partial}{\partial \xi} \right) \right)\\
&+ d\bar{z}_2\wedge \left[ r_2\frac{\partial}{\partial z_1} + r_3 \frac{\partial}{\partial z_2} + \left( b_0 + b_1 \xi + b_2 \xi^2 \right)\frac{\partial}{\partial \xi}, D\frac{\partial}{\partial z_1} -\left( A+ B \xi + C \xi^2 \right) \frac{\partial}{\partial \xi} \right] \\
&- d\bar{z}_2\left( g_3\left( D \frac{\partial}{\partial z_1} - \left( A + B \xi + C \xi^2 \right) \frac{\partial}{\partial \xi} \right) + g_4\left( D \frac{\partial}{\partial z_2} - k \left( A + B \xi + C \xi^2 \right)\frac{\partial}{\partial \xi} \right) \right)\\
&= d\bar{z}_1 \wedge \left(   - g_1 D \frac{\partial}{\partial z_1} - g_2 D \frac{\partial}{\partial z_2} + \left( \left(g_1 + g_2k\right)A+ a_1 A  - a_0 B +\left(\left(g_1+ g_2k\right)B +2a_2 A- 2a_0 C \right) \xi + \left( \left(g_1+ g_2k \right)C+ a_2B - a_1 C\right) \xi^2  \right) \frac{\partial}{\partial \xi}      \right)\\
&+d\bar{z}_2 \wedge \left(   - g_3 D \frac{\partial}{\partial z_1} - g_4 D \frac{\partial}{\partial z_2} + \left( \left(g_3 + g_4k\right)A+ b_1 A  - b_0 B +\left(\left(g_3+ g_4k\right)B +2b_2 A- 2b_0 C \right) \xi + \left( \left(g_3+ g_4k \right)C+ b_2B - b_1 C\right) \xi^2  \right) \frac{\partial}{\partial \xi}      \right)
\end{align*}}}
{\tiny{\begin{align*}
\sigma\left( T_2 \right)&= d\bar{z}_1\wedge \left[ r_0 \frac{\partial}{\partial z_1} + r_1 \frac{\partial}{\partial z_2} + \left( a_0 + a_1 \xi + a_2 \xi^2 \right)\frac{\partial}{\partial \xi}, D\frac{\partial}{\partial z_2} - k \left( A+ B \xi + C \xi^2 \right) \frac{\partial}{\partial \xi} \right] \\
&- d\bar{z}_1\left( g_5\left( D \frac{\partial}{\partial z_1} - \left( A + B \xi + C \xi^2 \right) \frac{\partial}{\partial \xi} \right) + g_6\left( D \frac{\partial}{\partial z_2} - k \left( A + B \xi + C \xi^2 \right)\frac{\partial}{\partial \xi} \right) \right)\\
&+ d\bar{z}_2\wedge \left[ r_2\frac{\partial}{\partial z_1} + r_3 \frac{\partial}{\partial z_2} + \left( b_0 + b_1 \xi + b_2 \xi^2 \right)\frac{\partial}{\partial \xi}, D\frac{\partial}{\partial z_2} -k\left( A+ B \xi + C \xi^2 \right) \frac{\partial}{\partial \xi} \right] \\
&- d\bar{z}_2\left( g_7\left( D \frac{\partial}{\partial z_1} - \left( A + B \xi + C \xi^2 \right) \frac{\partial}{\partial \xi} \right) + g_8\left( D \frac{\partial}{\partial z_2} - k \left( A + B \xi + C \xi^2 \right)\frac{\partial}{\partial \xi} \right) \right) \\
&= d\bar{z}_1 \wedge \left(   - g_5 D \frac{\partial}{\partial z_1} - g_6 D \frac{\partial}{\partial z_2} + \left( \left(g_5 + g_6k\right)A+ k\left(a_1 A  - a_0 B\right) +\left(\left(g_5+ g_6k\right)B + k\left(2a_2 A- 2a_0 C\right) \right) \xi + \left( \left(g_5+ g_6k \right)C+ k\left( a_2B - a_1 C\right) \right) \xi^2  \right) \frac{\partial}{\partial \xi}      \right)\\
&+d\bar{z}_2 \wedge \left(   - g_7 D \frac{\partial}{\partial z_1} - g_8 D \frac{\partial}{\partial z_2} + \left( \left(g_7 + g_8k\right)A+ k\left(b_1 A  - b_0 B\right) +\left(\left(g_7+ g_8k\right)B +k\left(2b_2 A- 2b_0 C\right) \right) \xi + \left( \left(g_3+ g_4k \right)C+ k\left( b_2B - b_1 C \right) \right) \xi^2  \right) \frac{\partial}{\partial \xi}      \right)
\end{align*}}}
Then the kernel is given by $g_1= \cdots = g_8=0$ and $a_1A- a_0B=0, a_2A- a_0 C=0, a_2B- a_1 C=0$ and $b_1A- b_0B=0, b_2A- b_0C=0,b_2B- b_1C =0$, so that any element $\left(\phi,\psi \right)$ of $\textnormal{ker}\left( \tilde{D}_0 \right)$ has a unique representation of the form
\begin{align*}
\phi&= d\bar{z}_1 \wedge \left( t_3 \frac{\partial}{\partial z_1} + t_4 \frac{\partial}{\partial z_2} + t_7\left( A + B \xi + C \xi^2 \right) \frac{\partial}{\partial \xi}   \right) + d \bar{z}_2\wedge \left( t_4 \frac{\partial}{\partial z_1} + t_5 \frac{\partial}{\partial z_2} + t_8\left( A + B \xi + C \xi^2 \right) \frac{\partial}{\partial \xi} \right)\\
\psi&=0
\end{align*}
Then from $(\ref{ct5})$ and from $(\ref{cu30})$ we have
\begin{align} \label{cu34}
\mathbb{H}^1\left( X, \mathcal{E}_{\mathcal{O}_X \bigoplus \mathcal{O}_X}^\bullet \right) \cong &\left\{ \left(  t_2\rho_1 + t_1 \rho_2  \right) \right.\\
&\left.  \left(t_3 \frac{\partial}{\partial z_1} + t_4 \frac{\partial}{\partial z_2} + t_7\left(A+ B \xi + C \xi^2 \right)\frac{\partial}{\partial \xi} \right) d\bar{z}_1 + \left(t_5 \frac{\partial}{\partial z_1} + t_6\frac{\partial}{\partial z_2} + t_8\left(A+ B \xi + C \xi^2 \right)\frac{\partial}{\partial \xi} \right) d\bar{z}_2, \right. \notag \\
&\left.  \,\,\,|\,\,t_1,t_2,t_3,t_4,t_5,t_6,t_7,t_8\in \mathbb{C} \right\} \notag
\end{align}

We will construct a foliated analytic family of deformations of $\left(X, \omega_0 \right)$ where $\omega_0=Dd\xi+ \left( A+ B \xi + C \xi^2 \right)dz_1+ k\left( A + B \xi + C\xi^2\right)dz_2,D\ne 0, \left(A,B,C\right)\ne 0\in \mathbb{C}^3, k\in \mathbb{C}$ in terms of tangent sheaves such that the foliated Kodaira-Spencer map in terms of tangent sheaves is an isomorphism at the distinguished point, so that $\left( X , \omega_0 \right)$ is unobstructed in foliated deformations. In fact, we recall the integrability conditions for a foliated analytic family in terms of tangent sheaves from $(\ref{cu31})-(\ref{cu32})$. If we take
{\small{\begin{align*}
\left( X(t), \left( \varphi(t), U(t)=0 \right) \right) \in A^{0,0}\left( X, \mathscr{H}om_{\mathcal{O}_X}\left( \mathcal{O}_X \bigoplus \mathcal{O}_X , \Theta_X \right)\right)\bigoplus A^{0,1}\left( X , \Theta_X \right) \bigoplus A^{01}\left( X, \mathscr{H}om_{\mathcal{O}_X}\left(\mathcal{O}_X \bigoplus \mathcal{O}_X, \mathcal{O}_X\bigoplus \mathcal{O}_X \right) \right)
\end{align*}}}
where
{\small{\begin{align*}
X(t)\left( T_1 \right)& =T_1(t):= D\frac{\partial}{\partial z_1} - (A + t_2 F_0 + (B+ t_2 F_1) \xi + (C+ t_2 F_2)\xi^2)\frac{\partial}{\partial \xi} \\
X(t)\left( T_2 \right)&=T_2(t):= D\frac{\partial}{\partial z_2}-  \left( k(A+ t_2 F_0) + t_1( A + t_2 F_0)  +\left( k(B+ t_2 F_1)+ t_1 (B+ t_2 F_1) \right)\xi + \left( k(C+t_2F_2) + t_1 ( C + t_2 F_2) \right)\xi^2\right)\frac{\partial}{\partial \xi}\\
\varphi(t)=&d\bar{z}_1 \wedge \left(t_3\frac{\partial}{\partial z_1} + t_4 \frac{\partial}{\partial z_2} + t_7\left(A+ t_2 F_0 + \left( B+ t_2 F_1 \right)\xi  +\left(C+ t_2 F_2 \right)\xi^2\right) \frac{\partial}{\partial \xi} \right)  \\
& +d\bar{z}_2 \wedge \left( t_5 \frac{\partial}{\partial z_1} + t_6 \frac{\partial}{\partial z_2}  + t_8 \left( A+ t_2 F_0 + \left( B+ t_2 F_1 \right)\xi  +\left(C+ t_2 F_2 \right)\xi^2\right) \frac{\partial}{\partial \xi} \right) 
\end{align*}}}
Then we have $\varphi(0)=0$, $X(0)= \textnormal{Id}$ and the following equalities hold:
\begin{align*}
&\bar{\partial} \varphi(t) - \frac{1}{2} \left[ \varphi(t) , \varphi(t) \right] = 0\\
&\left[ T_1(t), T_2(t) \right] =0\\
&\bar{\partial} T_1(t)- \left[ \varphi(t), T_1(t)\right] = 0,\,\,\,\,\,\bar{\partial} T_2(t)- \left[ \varphi(t), T_2(t) \right] =0
\end{align*}
In other words, the integrability conditions $(\ref{cu31})-(\ref{cu32})$ are satisfied, so that there exists a foliated analytic family $\left(\mathcal{M}, \Theta_\mathcal{F} , B, \pi \right)$ where $\pi:\left( \mathcal{M}, \Theta_{\mathcal{F}} \right)\to B$ for a sufficiently small neighborhood $B$ of $0$ in $\mathbb{C}^8$ with coordinates $t=\left(t_1,t_2,t_3,t_4,t_5,t_6,t_7,t_8 \right)$ such that $\pi^{-1}(0)=\left(X, \omega_0 \right)$ and the foliated Kodaira-Spencer map in terms of tangent sheaves from $(\ref{cu33})$ given by
\begin{align*}
\varphi: T_0(B) &\to \mathbb{H}^1\left( X, \mathcal{E}_{ \mathcal{O}_X \bigoplus \mathcal{O}_X }^\bullet \right)\\
\frac{\partial}{\partial t_1} &\mapsto \left( \rho, \left(0 , 0 \right) \right),\,\,\,\textnormal{where}\,\,\,\rho(T_1)=0, \rho(T_2)= - \left( A+ B\xi + C \xi^2\right)\frac{\partial}{\partial \xi} \\
\frac{\partial}{\partial t_2} & \mapsto \left( \rho, (0,0) \right),\,\,\, \textnormal{where}\,\,\, \rho(T_1)=- \left( F_0 + F_1 \xi + F_2 \xi^2 \right) \frac{\partial}{\partial \xi}, \rho(T_2)= -k\left( F_0 + F_1 \xi + F_2 \xi^2 \right)\frac{\partial}{\partial \xi}\\
\frac{\partial}{\partial t_3} &\mapsto \left( 0, \left( \frac{\partial}{\partial z_1} d\bar{z}_1, 0 \right) \right)\\
\frac{\partial}{\partial t_4} & \mapsto \left( 0 , \left(\frac{\partial}{\partial z_2}d\bar{z}_1, 0 \right) \right)\\
\frac{\partial}{\partial t_5} & \mapsto \left( 0 , \left(\frac{\partial}{\partial z_1}d\bar{z}_2, 0 \right) \right)\\
\frac{\partial}{\partial t_6} & \mapsto \left( 0 , \left(\frac{\partial}{\partial z_2}d\bar{z}_2, 0 \right) \right)\\
\frac{\partial}{\partial t_7} &\mapsto \left(  0, \left( \left(A+ B \xi + C \xi^2\right)\frac{\partial}{\partial \xi} d\bar{z}_1, 0 \right)\right) \\
\frac{\partial}{\partial t_8} & \mapsto \left( 0, \left( \left( A+ B \xi + C \xi^2 \right)\frac{\partial}{\partial \xi} d\bar{z}_2,  0 \right) \right)
\end{align*}
is an isomorphism at $t=0$ by $(\ref{cu34})$, so that $\left( X, \omega_0 \right)$ is unobstructed in foliated deformations in terms of tangent sheaves.

\subsection{Unfoldings of $\left( X, \omega_0 \right)$ where $\omega_0\in H^0\left( X, \Omega_X^1 \otimes \bigwedge^3 \Theta_X \right)$ with $\omega_0 \wedge d\omega_0=0$}\

Let us consider the foliation $\omega_0$ on $X$ where $\omega_0\in H^0\left( X, \Omega_X^1 \otimes \bigwedge^3 \Theta_X \right)$ with $\omega_0 \wedge d\omega_0 = 0 $. Then we have the unfolded dual leaf complex from $(\ref{cp3})$
\begin{align*}
\tilde{\mathcal{E}}_{\bigwedge^3 \Theta_X}^\bullet : \mathcal{E}_{\bigwedge^3 \Theta_X} \xrightarrow{\tilde{G}_0} \left( \Omega_X^1 \otimes \bigwedge^3 \Theta_X \right) \bigoplus \bigwedge^3 \Theta_X \xrightarrow{\tilde{G}_1} \left( \bigwedge^3 \Theta_X \right) \bigoplus \left( \bigwedge^3 \Theta_X \otimes \tilde{\mathcal{S}}^1 \right) \xrightarrow{G_2} \bigwedge^3 \Theta_X \to 0
\end{align*}
We note that $\bigwedge^3 \Theta_X$ in the third term and in the fourth term is precisely
\begin{align*}
\tilde{\mathcal{S}}^2\otimes \bigwedge^3 \Theta_X \cong \bigwedge^3 \Theta_X \otimes \bigwedge^3 \Omega_X^1 \otimes \bigwedge^3 \Theta_X \cong \bigwedge^3 \Theta_X
\end{align*}
and we have
\begin{align*}
\bigwedge^3 \Theta_X \otimes \tilde{S}^1 \subset \bigwedge^3 \Theta_X \otimes \bigwedge^3 \Theta_X \otimes \bigwedge^2 \Omega_X^1 
\end{align*}
We note that since $\bigwedge^3 \Theta_X\cong \pi_2^* \Theta_{\mathbb{P}_\mathbb{C}^1}$, we have
\begin{align*}
\bigwedge^3 \Theta_X \otimes \bigwedge^3 \Theta_X\otimes \bigwedge^2 \Omega_X^1 &\cong  \pi_2^* \left( \Theta_{\mathbb{P}_\mathbb{C}^1}\otimes \Theta_{\mathbb{P}_\mathbb{C}^1} \right) \otimes \left( \bigwedge^2 \Omega_X^1 \right) = \pi_2^* \mathcal{O}_{\mathbb{P}_\mathbb{C}^1}(4) \otimes\left(\pi_1^* \left(\bigwedge^2 \Omega_T^1 \right) \bigoplus \left(\pi_1^*\Omega_T^1 \otimes \pi_2^* \Omega_{\mathbb{P}_\mathbb{C}^1}^1 \right) \right)\\
& \cong \left(  \pi_2^* \mathcal{O}_{\mathbb{P}_\mathbb{C}^1}(4) \otimes \pi_1^*\left(\bigwedge^2 \Omega_T^1 \right) \right) \bigoplus \left(\pi_1^* \Omega_T^1 \otimes \pi_2^* \Omega_{\mathbb{P}_\mathbb{C}^1}^1(4)\right)
\end{align*}
Since $H^0\left( \mathbb{P}_\mathbb{C}^1 ,   \mathcal{O}_{\mathbb{P}^1}(4) \right)=5$ and $H^0\left( \mathbb{P}_\mathbb{C}^1 , \Omega_{\mathbb{P}_\mathbb{C}^1}^1(4) \right)=3$, we see that $H^0\left( X,  \left( \bigwedge^3 \Theta_X \right)^{\otimes 2} \otimes \bigwedge^2 \Omega_X^1\right)=5+ 6 = 11$  by K\"unneth formula. Then $H^0\left( X, \left( \bigwedge^3 \Theta_X \right)^{\otimes 2} \otimes \bigwedge^2 \Omega_X^1 \right)$ has an unique representation of the form

\begin{align} \label{cu2}
\left( \left( a_0  + a_1 \xi + a_2 \xi^2 \right) d\xi \wedge dz_1+ \left( b_0 + b_1 \xi + b_2 \xi^2 \right) d\xi \wedge dz_2 + \left( c_0 + c_1 \xi + c_2 \xi^2 + c_3 \xi^3 + c_4 \xi^4 \right) dz_1 \wedge dz_2\right) \\
\otimes \left( \frac{\partial}{\partial \xi}\wedge \frac{\partial}{\partial z_1} \wedge \frac{\partial}{\partial z_2} \right) \otimes \left( \frac{\partial}{\partial \xi} \wedge \frac{\partial}{\partial z_1} \wedge \frac{\partial}{\partial z_2} \right) \notag
\end{align}
where $a_0,a_1,a_2, b_0,b_1,b_2,c_0,c_1,c_2,c_3,c_4\in \mathbb{C}$.

Let us consider the spectral sequence at $E_1$ from Dolbeualt bicomplex associated to $\tilde{\mathcal{E}}_{\bigwedge^3 \Theta_X}^\bullet$:
\begin{center}
$\begin{CD}
0\\
@AAA \\
H^0\left( X, \bigwedge^3 \Theta_X \right) \\
@AG_2 AA \\
\binom{H^0\left(X, \bigwedge^3 \Theta_X \right) \bigoplus H^0\left( X, \bigwedge^3 \Theta_X \otimes \tilde{S}^1 \right)}{\subset H^0\left(X, \bigwedge^3 \Theta_X\right) \bigoplus H^0\left( X, \left(\bigwedge^3 \Theta_X \right)^{\otimes 2} \otimes \bigwedge^2 \Omega_X^1 \right)} @.\,\,\,\,\cdots \\
@A\tilde{G}_1AA @A\tilde{G}_1AA \\
H^0\left(X, \Omega_X^1 \otimes  \bigwedge^3 \Theta_X \right) \bigoplus H^0\left(X, \bigwedge^3 \Theta_X \right) @. \,\,\,\,\,\,\,\,\,\,\,\, H^1\left(X , \Omega_X^1 \otimes \bigwedge^3 \Theta_X \right) \bigoplus H^1\left(  X, \bigwedge^3 \Theta_X \right) \\
@A\tilde{G}_0 AA @A\tilde{G}_0 AA  \\
H^0\left( X, \mathcal{E}_{\bigwedge^3 \Theta_X} \right) @. \,\,\,\,\,\,\,\,\,\,\, H^1\left( X,  \mathcal{E}_{\bigwedge^3 \Theta_X} \right) @. \,\,\,\,\,\,\,\,\,H^2\left(X, \mathcal{E}_{\bigwedge^3 \Theta_X} \right)
\end{CD}$
\end{center}
Then we have
{\small{\begin{align} \label{ct23}
\mathbb{H}^1\left( X, \tilde{ \mathcal{E}}_{\bigwedge^3 \Theta_X}^\bullet \right) \cong & \frac{\textnormal{ker}\left( H^0\left( X, \Omega_X^1\otimes \bigwedge^3 \Theta_X \right) \bigoplus H^0\left( X, \bigwedge^3 \Theta_X  \right) \xrightarrow{\tilde{G}_1} H^0\left( X, \bigwedge^3 \Theta_X \right) \bigoplus  H^0\left(X, \left( \bigwedge^3 \Theta_X \right)^{\otimes 2} \otimes \bigwedge^2 \Omega_X^1 \right) \right)}{\textnormal{im}\left( H^0\left(X, \mathcal{E}_{\bigwedge^3 \Theta_X} \right) \xrightarrow{\tilde{G}_0} H^0\left( X, \Omega_X^1\otimes \bigwedge^3 \Theta_X \right) \bigoplus H^0\left(X, \bigwedge^3 \Theta_X \right) \right) } \\
&\bigoplus \textnormal{ker}\left(H^1\left( X, \mathcal{E}_{\bigwedge^3 \Theta_X} \right) \xrightarrow{\tilde{G}_0}  H^1\left(X, \Omega_X^1 \otimes \bigwedge^3 \Theta_X \right) \bigoplus H^1\left( X, \bigwedge^3 \Theta_X \right) \right) \notag
\end{align}}}

It remains to describe the cohomology groups $H^0\left(X, \bigwedge^3 \Theta_X \right)$ and $H^1\left(X, \bigwedge^3 \Theta_X \right)$ in the second term. Since $\bigwedge^3 \Theta_X \cong \pi^* \Theta_{\mathbb{P}_\mathbb{C}^1}$, we see that $\dim_\mathbb{C} H^0\left(X, \bigwedge^3 \Theta_X \right)=3$ and we have
\begin{align} \label{ct27}
H^0\left( X, \bigwedge^3 \Theta_X \right)=\left\{ f = g_0 + g_1 \xi + g_2 \xi^2 | g_0,g_1,g_2\in \mathbb{C} \right\}
\end{align}
On the other hand, $\dim_\mathbb{C}H^1 \left( X, \bigwedge^3 \Theta_X \right)=6$ and we have
\begin{align} \label{ct25}
H^1 \left( X , \bigwedge^3 \Theta_X \right) =\left\{ \left( h_1+h_2 \xi + h_3 \xi^2\right) d\bar{z}_1 + \left(  h_4  + h_5 \xi + h_6 \xi^2 \right) d\bar{z}_2 | h_1,...,h_6 \in \mathbb{C} \right\}
\end{align}

\subsection{Unfoldings of $\left(X, D d \xi \right),D\ne 0$}\

Let us  consider the foliation on $X$ defined by $\omega= D d\xi, D\ne 0$. We compute $\mathbb{H}^1\left( X, \tilde{\mathcal{E}}_{\bigwedge^3 \Theta_X}^\bullet \right)$. From $(\ref{ct23})$, first we compute $H^0\left(X, \mathcal{E}_{\bigwedge^3 \Theta_X}\right)\xrightarrow{\tilde{G}_0} H^0\left(  X, \Omega_X^1 \otimes \bigwedge^3 \Theta_X \right)\bigoplus H^0\left( X, \bigwedge^3 \Theta_X \right)$. We note that from $(\ref{ct13})$ the image of $H^0\left( X, \mathcal{E}_{\bigwedge^3 \Theta_X} \right)\to H^0\left( X, \Omega_X^1 \otimes \bigwedge^3 \Theta_X \right)$ is given by (from $(\ref{ct28})$)
\begin{align*}
D(a_4+a_6)d\xi
\end{align*}
and the image of $H^0\left( X, \mathcal{E}_{\bigwedge^3 \Theta_X} \right)\to H^0\left( X, \bigwedge^3 \Theta_X \right)$ is given by
\begin{align*}
i_{a_1\frac{\partial}{\partial z_1} \wedge + a_2 \frac{\partial}{\partial z_2} + \left( a_3 + a_4 \xi + a_5 \xi^2 \right) \frac{\partial}{\partial \xi}} \left( D d\xi \right)= D\left(a_3+ a_4 \xi + a_5 \xi^2 \right)
\end{align*}
Then from $(\ref{ct14})$ and $(\ref{ct27})$, we have 
\begin{align*}
\dim_\mathbb{C} \mathbb{H}^0\left( X, \tilde{\mathcal{E}}_{\bigwedge^3 \Theta_X}^\bullet \right)= \dim_\mathbb{C} \textnormal{ker}\left( H^0\left(X, \mathcal{E}_{\bigwedge^3 \Theta_X}\right)\xrightarrow{\tilde{G}_0} H^0\left(  X, \Omega_X^1 \otimes \bigwedge^3 \Theta_X \right)\bigoplus H^0\left( X, \bigwedge^3 \Theta_X \right) \right)=2
\end{align*}
and so $\dim_\mathbb{C} \textnormal{im}\left( H^0\left(X, \mathcal{E}_{\bigwedge^3 \Theta_X}\right)\xrightarrow{\tilde{G}_0} H^0\left(  X, \Omega_X^1 \otimes \bigwedge^3 \Theta_X \right)\bigoplus H^0\left( X, \bigwedge^3 \Theta_X \right) \right)= 6 - 2 =4$.

We compute $H^0\left( X, \Omega_X^1 \otimes \bigwedge^3 \Theta_X  \right) \bigoplus H^0\left(X, \bigwedge^3 \Theta_X \right) \xrightarrow{\tilde{G}_1} H^0\left( X, \bigwedge^3 \Theta_X \right) \bigoplus H^0\left(X, \left(\bigwedge^3 \Theta_X \right)^{\otimes 2} \otimes \bigwedge^2 \Omega_X^1  \right)$. From $(\ref{ct14})$ and $(\ref{ct27})$ the image of $(b, f)$ under $\tilde{G}_1$ is given by (from $(\ref{ct29})$)
\begin{align*}
\left( \omega_0 \wedge db + b \wedge d\omega_0=0,  f d\omega_0 + \omega_0 \wedge df - \omega_0 \wedge b \right)
\end{align*}
and we note that
\begin{align*}
&f d\omega_0 + \omega_0 \wedge df - \omega_0 \wedge b= \omega_0 \wedge df- \omega_0 \wedge b \\
&= D d\xi \wedge d \left(g_0 + g_1 \xi + g_2\xi^2 \right) - D d\xi\wedge \left( a d\xi+ \left( b_0 + b_1 \xi + b_2 \xi^2 \right) dz_1 + \left( c_0 + c_1 \xi + c_2 \xi^2 \right) dz_2 \right)\\
&=- D\left( b_0+ b_1 \xi + b_2 \xi^2 \right) d\xi \wedge dz_1 - D\left(c_0 + c_1 \xi + c_2 \xi^2 \right) d\xi \wedge dz_2 = 0 \iff b_0=b_1=b_2=c_0=c_1=c_2=0
\end{align*}
Hence from $(\ref{cu2})$ we have
\begin{align*}
\dim_\mathbb{C} \textnormal{ker}\left( H^0\left( X, \Omega_X^1 \otimes \bigwedge^3 \Theta_X  \right) \bigoplus H^0\left(X, \bigwedge^3 \Theta_X \right) \to H^0\left( X, \bigwedge^3 \Theta_X \right) \bigoplus H^0\left(X, \left(\bigwedge^3 \Theta_X \right)^{\otimes 2} \otimes \bigwedge^2 \Omega_X^1  \right) \right)=4
\end{align*}
 and we have
\begin{align} \label{ct24}
\dim_\mathbb{C} \frac{\textnormal{ker}\left( H^0\left( X, \Omega_X^1 \otimes \bigwedge^3 \Theta_X  \right) \bigoplus H^0\left(X, \bigwedge^3 \Theta_X \right) \xrightarrow{\tilde{G}_0} H^0\left( X, \bigwedge^3 \Theta_X \right) \bigoplus H^0\left(X, \left(\bigwedge^3 \Theta_X \right)^{\otimes 2} \otimes \bigwedge^2 \Omega_X^1  \right) \right)}{\textnormal{im}\left( H^0\left(X, \mathcal{E}_{\bigwedge^3 \Theta_X}\right)\xrightarrow{\tilde{G}_0} H^0\left(  X, \Omega_X^1 \otimes \bigwedge^3 \Theta_X \right)\bigoplus H^0\left( X, \bigwedge^3 \Theta_X \right) \right)} = 4 - 4 = 0
\end{align}

On the other hand, we compute $H^1\left( X, \mathcal{E}_{\bigwedge^3 \Theta_X} \right)\to H^1\left( X, \Omega_X^1 \otimes \bigwedge^3 \Theta_X \right) \bigoplus H^1\left(X, \bigwedge^3 \Theta_X \right)$. From $(\ref{ct11})$ the image of $H^1\left( X, \mathcal{E}_{\bigwedge^3 \Theta_X} \right)\to H^1\left(X, \Omega_X^1 \otimes \bigwedge^3 \Theta_X \right)$ is given by (from $(\ref{cu1})$)
\begin{align*}
d\bar{z}_1\wedge D(a_1+ a_3)d\xi + d\bar{z}_2 \wedge D(b_1+b_3) d\xi
\end{align*}
and the image of $H^1 \left( X, \mathcal{E}_{\bigwedge^3 \Theta_X} \right)\to H^1\left( X, \bigwedge^3 \Theta_X \right)$ is given by
\begin{align*}
&d\bar{z}_1\wedge i_{r_0 \frac{\partial}{\partial z_1} + r_1 \frac{\partial}{\partial z_2} + \left( a_0 + a_1 \xi + a_2 \xi^2 \right) \frac{\partial}{\partial \xi} }\left(D \xi \right) + d\bar{z}_2 \wedge i_{r_2 \frac{\partial}{\partial z_1} + r_3 \frac{\partial}{\partial z_2} + \left( b_0 + b_1 \xi + b_2 \xi^2 \right) \frac{\partial}{\partial \xi}}\left( D\xi\right)\\
&=D \left(a_0+ a_1 \xi  + a_2 \xi^2\right)d\bar{z}_1 + D\left( b_0 + b_1 \xi + b_2 \xi^2\right) d\bar{z}_2= 0 \iff a_0=a_1=a_2=b_0=b_1=b_2=0
\end{align*}
Then from $(\ref{ct15})$ and $(\ref{ct25})$ $\dim_\mathbb{C} \textnormal{ker}\left(H^1\left( X, \mathcal{E}_{\bigwedge^3 \Theta_X} \right)\to H^1\left( X, \Omega_X^1 \otimes \bigwedge^3 \Theta_X \right) \bigoplus H^1\left(X, \bigwedge^3 \Theta_X \right) \right)=4$ and the kernel is given by
\begin{align*}
\left(\phi=\left(t_1 \frac{\partial}{\partial z_1} + t_2 \frac{\partial}{\partial z_2} \right)d \bar{z}_1 + \left( t_3 \frac{\partial}{\partial z_1} + t_4 \frac{\partial}{\partial z_2} \right)  d\bar{z}_2, 0 \right),\,\,\,\,\,\,\,t_1,t_2,t_3,t_4 \in \mathbb{C}
\end{align*}
Then from $(\ref{ct23})$ and from $(\ref{ct24})$ we have
\begin{align}\label{cu10}
\mathbb{H}^1\left( X, \tilde{\mathcal{E}}_{\bigwedge^3 \Theta_X}^\bullet \right)\cong \left\{ \left(t_1 \frac{\partial}{\partial z_1} + t_2 \frac{\partial}{\partial z_2} \right)d \bar{z}_1 + \left( t_3 \frac{\partial}{\partial z_1} + t_4 \frac{\partial}{\partial z_2} \right)  d\bar{z}_2 | t_1,t_2,t_3, t_4 \in \mathbb{C} \right\}
\end{align}

We will construct an unfolded analytic family of unfoldings of $\left(X, \omega_0 \right)$ where $\omega= D d\xi, D\ne 0$ such that the unfolded Kodaira-Spencer map is an isomorphism at the distinguished point, so that $\left( X , \omega_0 \right)$ is unobstructed in unfoldings. In fact, we recall the integrability conditions for an unfolding from $(\ref{uti4})-(\ref{uti8})$. If we take
\begin{align*}
\left( Y(t)=0, Q^\eta(t)=0, \left\{ \left( \varphi(t), Z(t)=0 \right) \right\} \right) \in  A^{0,0}\left( X, \Omega_X^1 \otimes \bigwedge^3 \Theta_X \right) \bigoplus A^{0,0}\left( X, \bigwedge^3 \Theta_X \right)\bigoplus A^{0,1}\left( X, \mathcal{E}_{\bigwedge^3 \Theta_X} \right),\eta=1,2,3,4
\end{align*}
where $\varphi(t) := \left( t_1 \frac{\partial}{\partial z_1} + t_2 \frac{\partial}{\partial z_2} \right)d \bar{z}_1 + \left( t_3 \frac{\partial}{\partial z_1} + t_4 \frac{\partial}{\partial z_2}\right)\in A^{0,1}\left(X, \Theta_X \right)$. Then we have $\bar{\partial}\varphi(t)- \frac{1}{2}\left[ \varphi(t), \varphi(t) \right]=0$ and so the integrability conditions $(\ref{uti4})-(\ref{uti8})$ are satisfied, so that there exists an unfolded analytic family $\left(\mathcal{M}, \omega , B, \pi \right)$ where $\pi:\left( \mathcal{M}, \omega \right)\to B$ for a sufficiently small neighborhood $B$ of $0$ in $\mathbb{C}^4$ with coordinates $t=\left(t_1,t_2,t_3,t_4 \right)$ such that $\pi^{-1}(0)=\left(X, \omega_0 \right)$ and  the unfolded Kodaira-Spencer map from $(\ref{cu5})$ given by
\begin{align*}
\varphi:T_0\left(B \right) &\to \mathbb{H}^1\left( X, \tilde{\mathcal{E}}_{\bigwedge^3 \Theta_X}^\bullet \right) \\
\frac{\partial}{\partial t_1}&\mapsto \left( 0, 0, \left( \frac{\partial}{\partial z_1} d\bar{z}_1   , 0 \right)  \right)\\
\frac{\partial}{\partial t_2}&\mapsto \left( 0, 0, \left( \frac{\partial}{\partial z_2} d\bar{z}_1   , 0 \right)  \right)\\
\frac{\partial}{\partial t_3}&\mapsto \left( 0, 0, \left( \frac{\partial}{\partial z_1} d\bar{z}_2   , 0 \right)  \right)\\
\frac{\partial}{\partial t_4}&\mapsto \left( 0, 0, \left( \frac{\partial}{\partial z_2} d\bar{z}_2   , 0 \right)  \right)
\end{align*}
is an isomorphism at $t=0$ by $(\ref{cu10})$, so that $\left( X, \omega_0= Dd\xi \right),D\ne 0$ is unobstructed in unfoldings.

\subsection{Unfoldings of $\left(X, Dd\xi + \left( A+ B \xi + C\xi^2 \right)dz_1 + k\left( A + B\xi + C \xi^2 \right) dz_2 \right), D\ne 0, \left(A, B, C \right)\ne 0 \in \mathbb{C}^3 , k\in \mathbb{C}$}\

Let us consider the foliation on $X$ defined by $\omega_0=Dd\xi + \left( A+ B \xi + C\xi^2 \right)dz_1 + k\left( A + B\xi + C \xi^2 \right) dz_2,D\ne 0, \left(A, B, C \right)\ne 0 \in \mathbb{C}^3 , k\in \mathbb{C}$. We compute $\mathbb{H}^1\left( X, \tilde{\mathcal{E}}_{\bigwedge^3 \Theta_X}^\bullet \right)$. From $(\ref{ct23})$, first we compute $H^0\left(X, \mathcal{E}_{\bigwedge^3 \Theta_X}\right)\xrightarrow{\tilde{G}_0} H^0\left(  X, \Omega_X^1 \otimes \bigwedge^3 \Theta_X \right)\bigoplus H^0\left( X, \bigwedge^3 \Theta_X \right)$. We note that from $(\ref{ct13})$ the image of $H^0\left( X, \mathcal{E}_{\bigwedge^3 \Theta_X} \right)\to H^0\left( X, \Omega_X^1 \otimes \bigwedge^3 \Theta_X \right)$ is given by (from $(\ref{cu7})$)
\begin{align*}
&D(a_4+ a_6) d\xi +\left(  Ba_3+ A a_6 +\left( Ba_4+2 C a_3+ Ba_6- 2A a_5 \right) \xi +\left(  2 C a_4+ C a_6-  B a_5 \right) \xi^2 \right) dz_1 \\
&+k \left(  Ba_3+ A a_6 +\left( Ba_4+2 C a_3+ Ba_6- 2A a_5 \right) \xi +\left( 2 C a_4+ C a_6-  B a_5 \right) \xi^2 \right) dz_2 = 0 \\
&\iff a_6= - a_4, Ba_3- Aa_4=0 , Ca_3- Aa_5=0, Ca_4- B a_5 =0
\end{align*}
and the image of $H^0\left( X, \mathcal{E}_{\bigwedge^3 \Theta_X} \right)\to H^0\left( X, \bigwedge^3 \Theta_X \right)$ is given by
\begin{align*}
&i_{a_1\frac{\partial}{\partial z_1} \wedge + a_2 \frac{\partial}{\partial z_2} + \left( a_3 + a_4 \xi + a_5 \xi^2 \right) \frac{\partial}{\partial \xi}} \left( Dd\xi + \left( A+ B \xi + C\xi^2 \right)dz_1 + k\left( A + B\xi + C \xi^2 \right) dz_2 \right)\\
&=(a_1+ka_2)A+a_3D +\left( (a_1+ ka_2)B+ a_4D \right) \xi + \left( (a_1+k a_2) C + a_5 D  \right)\xi^2 =0\\
&\iff a_3= - \frac{a_1 + ka_2}{D}A , a_4= - \frac{a_1 + ka_2}{D}B, a_5= - \frac{a_1 + ka_2}{D}C,a_6= \frac{a_1+ k a_2}{D} B
\end{align*}
Then from $(\ref{ct14})$ and $(\ref{ct27})$, we have 
\begin{align*}
\dim_\mathbb{C} \mathbb{H}^0\left( X, \tilde{\mathcal{E}}_{\bigwedge^3 \Theta_X}^\bullet \right)= \dim_\mathbb{C} \textnormal{ker}\left( H^0\left(X, \mathcal{E}_{\bigwedge^3 \Theta_X}\right)\xrightarrow{\tilde{G}_0} H^0\left(  X, \Omega_X^1 \otimes \bigwedge^3 \Theta_X \right)\bigoplus H^0\left( X, \bigwedge^3 \Theta_X \right) \right)=2
\end{align*}
and so $\dim_\mathbb{C} \textnormal{im}\left( H^0\left(X, \mathcal{E}_{\bigwedge^3 \Theta_X}\right)\xrightarrow{\tilde{G}_0} H^0\left(  X, \Omega_X^1 \otimes \bigwedge^3 \Theta_X \right)\bigoplus H^0\left( X, \bigwedge^3 \Theta_X \right) \right)= 6 - 2 =4$.

We compute $H^0\left( X, \Omega_X^1 \otimes \bigwedge^3 \Theta_X  \right) \bigoplus H^0\left(X, \bigwedge^3 \Theta_X \right) \xrightarrow{\tilde{G}_1} H^0\left( X, \bigwedge^3 \Theta_X \right) \bigoplus H^0\left(X, \left(\bigwedge^3 \Theta_X \right)^{\otimes 2} \otimes \bigwedge^2 \Omega_X^1  \right)$. From $(\ref{ct14})$ and $(\ref{ct27})$ the image of $(b, f)$ under $\tilde{G}_1$ is given by 
\begin{align*}
\left( \omega_0 \wedge db + b \wedge d\omega_0,  f d\omega_0 + \omega_0 \wedge df - \omega_0 \wedge b \right)
\end{align*}
where we have from $(\ref{cu15})$  and $(\ref{cu16})$
{\small{\begin{align}
&b\wedge d\omega_0 + \omega_0 \wedge db \notag \\
&=\left(\left(Bc_0- A c_1+ k \left(Ab_1- Bb_0\right) \right) + \left(2 Cc_0- 2Ac_2 + k\left(2Ab_2- 2Cb_0 \right) \right) \xi + \left( Cc_1- Bc_2 + k\left(Bb_2- Cb_1\right) \right) \xi^2 \right) d\xi \wedge dz_1 \wedge dz_2=0 \notag \\
&\iff   c_0- kb_0= t_1 A, c_1-kb_1=t_1B , c_2-kb_2= t_1C,\,\,\,\,\,t_1\in \mathbb{C} \label{cu17}            
\end{align}}}
and 
\begin{align*}
&f d\omega_0 + \omega_0 \wedge df - \omega_0 \wedge b \\
&=  \left(g_0 + g_1\xi+ g_2 \xi^2\right)d \omega_0 + \omega_0 \wedge d\left(g_0 + g_1 \xi+ g_2 \xi^2 \right)  - \omega_0 \wedge \left( a d\xi + \left(b_0+ b_1 \xi + b_2 \xi^2 \right) dz_1 + \left( c_0 + c_1 \xi + c_2 \xi^2  \right) dz_2 \right)\\
&=\left( \left( - Ag_1 + B g_0 + Aa - Db_0 \right) +\left(-2 Ag_2+ 2 C g_0+ Ba - Db_1 \right)\xi +\left( - B g_2 + Cg_1 + Ca - Db_2 \right) \xi^2 \right) d\xi \wedge dz_1\\
&+\left(k\left(- Ag_1+ B g_0 + Aa \right)- Dc_0  +\left(  k\left(- 2 A g_2+ 2 C g_0 + Ba\right) - Dc_1 \right)\xi + \left( k \left(- Bg_2 + Cg_1+ Ca \right)- Dc_2 \right) \xi^2  \right) d\xi \wedge dz_2\\
&-\left(  Ac_0+ Bc_0\xi + Cc_0\xi^2 + Ac_1\xi+ Bc_1\xi^2+ Cc_1\xi^3+ Ac_2\xi^2+ Bc_2\xi^3+ Cc_2\xi^4 \right) dz_1\wedge dz_2\\
&+ k \left(  Ab_0+ Bb_0\xi + Cb_0\xi^2 + Ab_1\xi+ Bb_1\xi^2+ Cb_1\xi^3+ Ab_2\xi^2+ Bb_2\xi^3+ Cb_2\xi^4 \right) dz_1\wedge dz_2
\end{align*}
From $(\ref{cu17})$, if we set $b= a+ \left(b_0 + b_1 \xi+ b_2 \xi^2 \right)dz_1+ \left( kb_0 + t_1 A + \left( kb_1 + t_1 B \right)\xi + \left( kb_2 + t_1 C \right) \xi^2 \right)$, then
\begin{align*}
&f d\omega_0 + \omega_0 \wedge df - \omega_0 \wedge b\\
&=\left( \left( - Ag_1 + B g_0 + Aa - D b_0  \right) +\left(-2 Ag_2+ 2 C g_0+ Ba - Db_1  \right)\xi +\left( - B g_2 + Cg_1 + Ca - Db_2  \right) \xi^2 \right) d\xi \wedge dz_1\\
&+\left(k\left(- Ag_1+ B g_0 + Aa - Db_0 \right) +k\left( - 2 A g_2+ 2 C g_0 + Ba - Db_1  \right)\xi +k \left(- Bg_2 + Cg_1+ Ca - Db_2  \right) \xi^2  \right) d\xi \wedge dz_2\\
&- Dt_1\left( A+ B \xi + C\xi^2 \right)d\xi \wedge dz_2 -t_1\left( A + B \xi + C\xi^2 \right)^2 dz_1 \wedge dz_2
\end{align*}
Then $\textnormal{ker}\left( H^0\left( X, \Omega_X^1 \otimes \bigwedge^3 \Theta_X  \right) \bigoplus H^0\left(X, \bigwedge^3 \Theta_X \right) \xrightarrow{\tilde{G}_1} H^0\left( X, \bigwedge^3 \Theta_X \right) \bigoplus H^0\left(X, \left(\bigwedge^3 \Theta_X \right)^{\otimes 2} \otimes \bigwedge^2 \Omega_X^1  \right) \right)$ is given by
\begin{align*}
&\,\,\,\,\,\,\,\, t_1=0\\
&-Ag_1+ Bg_0+ Aa- Db_0=0\\
&-2Ag_2+ 2Cg_0+ Ba- Db_1=0\\
&-Bg_2+ Cg_1+ Ca- Db_2 = 0 
\end{align*}
Since $D\ne 0$, the rank of the matrix from the last three equations is $3$, so that we have
\begin{align*}
\dim_\mathbb{C} \textnormal{ker}\left( H^0\left( X, \Omega_X^1 \otimes \bigwedge^3 \Theta_X  \right) \bigoplus H^0\left(X, \bigwedge^3 \Theta_X \right) \xrightarrow{\tilde{G}_1} H^0\left( X, \bigwedge^3 \Theta_X \right) \bigoplus H^0\left(X, \left(\bigwedge^3 \Theta_X \right)^{\otimes 2} \otimes \bigwedge^2 \Omega_X^1  \right) \right)= 4
\end{align*}
and so we have
\begin{align}
\dim_\mathbb{C} \frac{\textnormal{ker}\left( H^0\left( X, \Omega_X^1 \otimes \bigwedge^3 \Theta_X  \right) \bigoplus H^0\left(X, \bigwedge^3 \Theta_X \right) \xrightarrow{\tilde{G}_0} H^0\left( X, \bigwedge^3 \Theta_X \right) \bigoplus H^0\left(X, \left(\bigwedge^3 \Theta_X \right)^{\otimes 2} \otimes \bigwedge^2 \Omega_X^1  \right) \right)}{\textnormal{im}\left( H^0\left(X, \mathcal{E}_{\bigwedge^3 \Theta_X}\right)\xrightarrow{\tilde{G}_0} H^0\left(  X, \Omega_X^1 \otimes \bigwedge^3 \Theta_X \right)\bigoplus H^0\left( X, \bigwedge^3 \Theta_X \right) \right)} = 4 - 4 = 0
\end{align}

On the other hand, we compute $H^1\left( X, \mathcal{E}_{\bigwedge^3 \Theta_X} \right)\to H^1\left( X, \Omega_X^1 \otimes \bigwedge^3 \Theta_X \right) \bigoplus H^1\left(X, \bigwedge^3 \Theta_X \right)$. From $(\ref{ct11})$ the image of $H^1\left( X, \mathcal{E}_{\bigwedge^3 \Theta_X} \right)\to H^1\left(X, \Omega_X^1 \otimes \bigwedge^3 \Theta_X \right)$ is given by (from $(\ref{cu18})$)
\begin{align*}
&d\bar{z}_1\wedge \left( D(a_1+ a_3) d\xi +\left(  Ba_0+ A a_3 +\left( Ba_1+2 C a_0+ Ba_3- 2A a_2 \right) \xi +\left(  2 C a_1+ C a_3-  B a_2 \right) \xi^2 \right) dz_1\right) \\
&+ d\bar{z}_1\wedge \left(k \left(  Ba_1+ A a_3 +\left( Ba_1+2 C a_0+ Ba_3- 2A a_2 \right) \xi +\left( 2 C a_1+ C a_3-  B a_2 \right) \xi^2 \right) dz_2 \right)\\
&+ d\bar{z}_2 \wedge \left( D(b_1+ b_3) d\xi +\left(  Bb_0+ A b_3 +\left( Bb_1+2 C b_0+ Bb_3- 2A b_2 \right) \xi +\left(  2 C b_1+ C b_3-  B b_2 \right) \xi^2 \right) dz_1 \right) \\
&+ d\bar{z}_2 \wedge \left( k \left(  Bb_0+ A b_3 +\left( Bb_1+2 C b_0+ Bb_3- 2A b_2 \right) \xi +\left( 2 C b_1+ C b_3-  B b_2 \right) \xi^2 \right) dz_2\right) \\
\end{align*}
and the image of $H^1 \left( X, \mathcal{E}_{\bigwedge^3 \Theta_X} \right)\to H^1\left( X, \bigwedge^3 \Theta_X \right)$ is given by
\begin{align*}
&d\bar{z}_1\wedge i_{r_0 \frac{\partial}{\partial z_1} + r_1 \frac{\partial}{\partial z_2} + \left( a_0 + a_1 \xi + a_2 \xi^2 \right) \frac{\partial}{\partial \xi} }\left(D \xi + \left( A+ B \xi + C \xi^2 \right) dz_1  + k\left( A+ B \xi + C \xi^2 \right) dz_2 \right) \\
& + d\bar{z}_2 \wedge i_{r_2 \frac{\partial}{\partial z_1} + r_3 \frac{\partial}{\partial z_2} + \left( b_0 + b_1 \xi + b_2 \xi^2 \right) \frac{\partial}{\partial \xi}}\left( D\xi + \left( A+ B \xi + C \xi^2 \right) dz_1  + k\left( A+ B \xi + C \xi^2 \right) dz_2\right)\\
&= \left((r_0+kr_1)A+a_0D +\left( (r_0+ kr_1)B+ a_1D \right) \xi + \left( (r_0+k r_1) C + a_2 D  \right)\xi^2 \right) d\bar{z}_1\\
&+  \left((r_2 +kr_3)A+b_0 D +\left( (r_2+ kr_3)B+ b_1 D \right) \xi + \left( (r_2+k r_3) C + b_2 D  \right)\xi^2 \right) d\bar{z}_2
\end{align*}
Then from $(\ref{ct15})$ and $(\ref{ct25})$ $\dim_\mathbb{C} \textnormal{ker}\left(H^1\left( X, \mathcal{E}_{\bigwedge^3 \Theta_X} \right)\to H^1\left( X, \Omega_X^1 \otimes \bigwedge^3 \Theta_X \right) \bigoplus H^1\left(X, \bigwedge^3 \Theta_X \right) \right)=4$ and the kernel is given by
\begin{align}
&\left( \left(t_1 \frac{\partial}{\partial z_1} + t_2 \frac{\partial}{\partial z_2} -\frac{t_1+ kt_2 }{D}\left(A+ B \xi + C \xi^2 \right)\frac{\partial}{\partial \xi} \right) d\bar{z}_1 + \left(t_3 \frac{\partial}{\partial z_1} + t_4\frac{\partial}{\partial z_2} - \frac{t_3 + kt_4 }{D}\left(A+ B \xi + C \xi^2 \right)\frac{\partial}{\partial \xi} \right) d\bar{z}_2, \right.\\
&\left. \,\,\,\,\,\,\,\,\,\left(\frac{t_1+ kt_2}{D}B + 2\frac{t_1+ kt_2}{D}C \xi\right)d\bar{z}_1 + \left( \frac{ t_3+ kt_4}{D} B + 2\frac{t_3 + k t_4}{D} C \xi \right) d\bar{z}_2   \right) ,\,\,\,\,\,\,\,\, 
t_1, t_2 , t_3 , t_4  \in \mathbb{C} \notag
\end{align}
Then from $(\ref{ct23})$ and from $(\ref{ct24})$ we have
\begin{align} \label{cu21}
&\mathbb{H}^1\left( X, \tilde{\mathcal{E}}_{\bigwedge^3 \Theta_X}^\bullet \right)\\
&\cong \left\{ \left( \left(t_1 \frac{\partial}{\partial z_1} + t_2 \frac{\partial}{\partial z_2} -\frac{t_1+ kt_2 }{D}\left(A+ B \xi + C \xi^2 \right)\frac{\partial}{\partial \xi} \right) d\bar{z}_1 + \left(t_3 \frac{\partial}{\partial z_1} + t_4\frac{\partial}{\partial z_2} - \frac{t_3 + kt_4 }{D}\left(A+ B \xi + C \xi^2 \right)\frac{\partial}{\partial \xi} \right) d\bar{z}_2, \right. \right. \notag \\
&\left.\left. \,\,\,\,\,\,\,\,\,\left(\frac{t_1+ kt_2}{D}B + 2\frac{t_1+ kt_2}{D}C \xi\right)d\bar{z}_1 + \left( \frac{ t_3+ kt_4}{D} B + 2\frac{t_3 + k t_4}{D} C \xi \right) d\bar{z}_2   \right) \,\,\,| \,\,\,t_1, t_2 , t_3 , t_4  \in \mathbb{C} \right\} \notag
\end{align}
We will construct an unfolded analytic family of unfoldings of $\left(X, \omega_0 \right)$ where $\omega= D d\xi, D\ne 0$ such that the unfolded Kodaira-Spencer map is an isomorphism at the distinguished point, so that $\left( X , \omega_0 \right)$ is unobstructed in unfoldings. In fact, we recall the integrability conditions for an unfolding from $(\ref{uti4})-(\ref{uti8})$. If we take
\begin{align*}
\left( Y(t)=\omega_0, Q^\eta(t)=0, \left\{ \left( \varphi(t), Z(t) \right) \right\} \right) \in  A^{0,0}\left( X, \Omega_X^1 \otimes \bigwedge^3 \Theta_X \right) \bigoplus A^{0,0}\left( X, \bigwedge^3 \Theta_X \right)\bigoplus A^{0,1}\left( X, \mathcal{E}_{\bigwedge^3 \Theta_X} \right),\eta=1,2,3,4
\end{align*}
where $\varphi(t):= \left(t_1 \frac{\partial}{\partial z_1} + t_2 \frac{\partial}{\partial z_2} -\frac{t_1+ kt_2 }{D}\left(A+ B \xi + C \xi^2 \right)\frac{\partial}{\partial \xi} \right) d\bar{z}_1 + \left(t_3 \frac{\partial}{\partial z_1} + t_4\frac{\partial}{\partial z_2} - \frac{t_3 + kt_4 }{D}\left(A+ B \xi + C \xi^2 \right)\frac{\partial}{\partial \xi} \right) d\bar{z}_2$ and $Z(t):=\left(\frac{t_1+ kt_2}{D}B + 2\frac{t_1+ kt_2}{D}C \xi\right)d\bar{z}_1 + \left( \frac{ t_3+ kt_4}{D} B + 2\frac{t_3 + k t_4}{D} C \xi \right) d\bar{z}_2$. Then we have $\varphi(0)=0$ and $Z(0)=0$ and the following equalities hold:
\begin{align*}
&\bar{\partial}\varphi(t) - \frac{1}{2}\left[ \varphi(t), \varphi(t) \right]=0\\
&\bar{\partial} \omega_0 - \mathcal{L}_{\varphi(t)}\left( \omega_0 \right)- Z(t)\wedge \omega_0 =0\\
&\bar{\partial} Q^\eta - \left[ \varphi(t), Q^\eta\right]- i_{\frac{\partial \varphi(t)}{\partial t_\eta}}\left(\omega_0 \right) - Z(t) Q^\eta= - i_{\frac{\partial \varphi(t)}{\partial t_\eta}}\left( \omega_0 \right)=0,\,\,\,\,\eta=1,2,3,4\\
&\bar{\partial} Z(t) - \left[ \varphi(t), Z(t) \right] = 0
\end{align*}
In other words, the integrability conditions $(\ref{uti4})-(\ref{uti8})$ are satisfied, so that there exists an unfolded analytic family $\left(\mathcal{M}, \omega , B, \pi \right)$ where $\pi:\left( \mathcal{M}, \omega \right)\to B$ for a sufficiently small neighborhood $B$ of $0$ in $\mathbb{C}^4$ with coordinates $t=\left(t_1,t_2,t_3,t_4\right)$ such that $\pi^{-1}(0)=\left(X, \omega_0 \right)$ and  the unfolded Kodaira-Spencer map from $(\ref{cu5})$ given by
\begin{align*}
\varphi:T_0\left(B \right) &\to \mathbb{H}^1\left( X, \tilde{\mathcal{E}}_{\bigwedge^3 \Theta_X}^\bullet \right) \\
\frac{\partial}{\partial t_1}&\mapsto \left(  0,0, d\bar{z}_1\left( \frac{\partial}{\partial z_1} -\frac{1}{D}\left( A + B \xi + C \xi^2 \right) \frac{\partial}{\partial \xi} \right)  ,d\bar{z}_1\frac{1}{D}\left( B+ 2 C \xi\right)              \right)\\
\frac{\partial}{\partial t_2}&\mapsto \left(  0,0, d\bar{z}_1\left( \frac{\partial}{\partial z_2} -\frac{k}{D}\left( A + B \xi + C \xi^2 \right) \frac{\partial}{\partial \xi} \right)  ,d\bar{z}_1\frac{k}{D}\left( B+ 2 C \xi\right)              \right) \\
\frac{\partial}{\partial t_3}&\mapsto \left(  0,0, d\bar{z}_2\left( \frac{\partial}{\partial z_1} -\frac{1}{D}\left( A + B \xi + C \xi^2 \right) \frac{\partial}{\partial \xi} \right)  ,d\bar{z}_2 \frac{1}{D}\left( B+ 2 C \xi\right)              \right)\\
\frac{\partial}{\partial t_4}&\mapsto \left(  0,0, d\bar{z}_2\left( \frac{\partial}{\partial z_2} -\frac{k}{D}\left( A + B \xi + C \xi^2 \right) \frac{\partial}{\partial \xi} \right)  ,d\bar{z}_2\frac{k}{D}\left( B+ 2 C \xi\right)              \right)
\end{align*}
is an isomorphism at $t=0$ by $(\ref{cu21})$, so that $\left( X, \omega_0 \right)$ is unobstructed in unfoldings.

\section{Foliations on a complex torus $T$ with $\dim_\mathbb{C} T= 3$ defined by Poisson structures}

In this section we study foliated deformations of $\left(T, w_0 \right)$ in terms of cotangent sheaves and tangent sheaves, and unfoldings of $\left( T, w_0 \right)$, where $T$ is a complex torus with $\dim_\mathbb{C} T =3$ and the (regular) foliation $w_0 \in H^0\left( T, \Omega_T^1\otimes K_T^{-1} \right)\cong H^0\left( T, \Omega_T^1\right)$ are defined by holomorphic Poisson structures $\Lambda_0\in H^0\left(  T, \bigwedge^2 \Theta_T \right)$ on $T$.

\subsection{Holomorphic Poisson structures on a 3-dimensional complex torus $T$ and the associated codimension $1$-foliations on $T$}\

Let $T=\mathbb{C}^3/\mathcal{G}$ be a complex torus with dimension $3$ where $\mathcal{G}$ is the free abelian group generated by non-degenerate periods $\omega_j=\left( \omega_{j1}, \omega_{j2}, \omega_{j3}\right),j=1,2,3,4,5,6$ and where, if $(z_1,z_2,z_3)$ are the coordinates of $\mathbb{C}^3$, $\omega_j\in \mathcal{G}$ operates on $\mathbb{C}^3$ by sending $z=(z_1,z_2,z_3)$ into $z+\omega_j= \left(z_1 + \omega_{j1}, z_2 + \omega_{j2}, z_3 + \omega_{j3} \right)$. Then any element of $H^0\left( T, \bigwedge^2 \Theta_T \right)$ has an unique representation of the form $\Lambda_0=a_1\frac{\partial}{\partial z_1}\wedge \frac{\partial}{\partial z_3} + a_2 \frac{\partial}{\partial z_3}\wedge \frac{\partial}{\partial z_1} + a_3 \frac{\partial}{\partial z_1} \wedge \frac{\partial}{\partial z_2}$ where $\left(a_1,a_2,a_3\right)\in \mathbb{C}^3$ which defines holomorphic Poisson structures on $T$. When $(a_1,a_2,a_3)\ne 0\in \mathbb{C}^3$, $\Lambda_0$ corresponds to a regular foliation $\omega_0=a_1dz_1 + a_2 dz_2 + a_3 dz_3\in H^0\left( T, \Omega_T^1\otimes K_T^{-1} \right)= H^0\left(T, \Omega_T^1\otimes \mathcal{O}_T \right)= H^0\left( T, \Omega_T^1 \right)$ on $T$. In the following we will assume that $a_3\ne 0$. We shall study foliated deformations and unfoldings of $\left(T, \omega_0 \right)$. Since $\omega_0$ is regular, foliated deformations of $\left( T, \omega_0 \right)$ in terms of cotangent sheaves are equivalent to foliated deformations of $\left( T, \omega_0 \right)$ in terms of tangent sheaves. However we shall describe foliated deformations of $\left( T, \omega_0 \right)$ in terms of cotangent sheaves as well as tangent sheaves. We show that $\left(  T, \omega_0 \right)$ is unobstructed in foliated deformations and unfoldings.

\subsection{Foliated deformations of $\left( T, \omega_0= a_1dz_1 + a_2dz_2 + a_3 dz_3 \right),a_3\ne 0$ where $\omega_0\in H^0\left( T, \Omega_T^1 \right)$ in terms of cotangent sheaves}\

Let us consider the foliation on $T$ defined by $\omega_0= a_1dz_1 + a_2 dz_2 + a_3 dz_3\in H^0\left(T, \Omega_T^1\otimes K_T^{-1} \right)= H^0\left(T, \Omega_T^1\otimes \mathcal{O}_T \right)$. Then we have the dual leaf complex from $(\ref{cp2})$
\begin{align*}
\mathcal{E}_{\mathcal{O}_T}^\bullet : \mathcal{E}_{\mathcal{O}_T}= \Theta_T \bigoplus \mathcal{O}_T \to \Omega_T^1\otimes \mathcal{O}_T \to \bigwedge^{1+2}\Omega_T \otimes \mathcal{O}_T \to 0
\end{align*}
equivalently
\begin{align*}
\mathcal{E}_{\mathcal{O}_T}^\bullet : \Theta_T \bigoplus \mathcal{O}_T \to \Omega_T^1 \to \mathcal{O}_T \to 0
\end{align*}
Let us consider the spectral sequence at $E_1$ from Dolbeault bicomplex associated to $\mathcal{E}_{\mathcal{O}_T}^\bullet$:
\begin{equation}\label{com1}
\begin{CD}
0 \\
@AAA \\
H^0\left( T, \mathcal{O}_T \right) @.\,\,\,\cdots \\
@AAA @AAA \\
H^0\left( T, \Omega_T^1 \right) @.\,\,\,\,\,\,\,H^1\left( T, \Omega_T^1 \right) \\
@AAA @AAA \\
H^0\left( T, \Theta_T \right)\bigoplus H^0\left( T, \mathcal{O}_T \right) @. \,\,\,\,\,\,\,\,\,\, H^1\left( T, \Theta_T \right) \bigoplus H^1\left( T, \mathcal{O}_T \right) @.\,\,\,\,\,\,\,\,\,\, H^2\left( T, \Theta_T \right) \bigoplus H^2\left( T, \mathcal{O}_T \right)
\end{CD}
\end{equation}
Then we have
\begin{align} \label{com2}
\mathbb{H}^1\left( T, \mathcal{E}_{\mathcal{O}_T}^\bullet \right) \cong \frac{\textnormal{ker}\left(H^0\left( T, \Omega_T^1 \right)\to H^0\left( T, \mathcal{O}_T \right) \right)}{\textnormal{im}\left( H^0\left( T, \Theta_T \right) \bigoplus H^0\left( T, \mathcal{O}_T \right) \to H^0\left( T, \Omega_T^1 \right)\right)} \bigoplus \textnormal{ker}\left( H^1\left( T, \Theta_T \right) \bigoplus H^1\left( T, \mathcal{O}_T \right) \to H^1\left(T, \Omega_T^1  \right) \right)
\end{align}
We describe each cohomology group in $(\ref{com1})$. We note that $\dim_\mathbb{C} H^0\left( T, \Theta_T \right)=3$ and any element of $H^0\left(T, \Theta_T \right)$ has an unique representation of the form
\begin{align} \label{com3}
b_1 \frac{\partial}{\partial z_1} + b_2 \frac{\partial}{\partial z_2} + b_3 \frac{\partial}{\partial z_3},\,\,\,\,\, b_1 , b_2 , b_3 \in \mathbb{C}
\end{align}
We note that $\dim_\mathbb{C} H^0\left(T, \Omega_T^1 \right)=3$ and any element of $H^0\left( T, \Omega_T^1 \right)$ has an unique representation of the form
\begin{align}\label{com5}
c_1 dz_1 + c_2 dz_2 + c_3 dz_3,\,\,\,\,\, c_1, c_2 , c_3 \in \mathbb{C}
\end{align}
We note that $\dim_\mathbb{C} H^0\left( T, \mathcal{O}_T\right)=1$ and any element of $H^0\left( T, \mathcal{O}_T\right)$ is given by
\begin{align}\label{com6}
H^0\left( T, \mathcal{O}_T \right)= \left\{ h\,\,\,|\,\,\, h\in \mathbb{C} \right\}
\end{align}
We note that $\dim_\mathbb{C} H^1\left( T, \Theta_T \right)=9$ and any element of $H^1\left( T, \Theta_T \right)$ has au unique representation of the form
\begin{align}\label{com7}
d\bar{z}_1\wedge \left( e_1 \frac{\partial}{\partial z_1} + e_2 \frac{\partial}{\partial z_3} + e_3 \frac{\partial}{\partial z_3} \right) + d\bar{z}_2 \wedge \left( e_4 \frac{\partial}{\partial z_1} + e_5 \frac{\partial}{\partial z_3} + e_6 \frac{\partial}{\partial z_3}\right) + d\bar{z}_3 \wedge \left( e_7 \frac{\partial}{\partial z_1} + e_8 \frac{\partial}{\partial z_3} + e_9 \frac{\partial}{\partial z_3} \right),\,\,\,e_1,...,e_9\in \mathbb{C}
\end{align}
We note that $\dim_\mathbb{C} \mathbb{H}^1\left(T, \Omega_T^1 \right)=9$ and any element of $H^1\left( T, \Omega_T^1 \right)$ has an unique representation of the form
\begin{align}\label{com8}
d\bar{z}_1 \wedge\left( g_1 dz_1 + g_2 dz_2 + g_3 dz_3 \right) + d\bar{z}_2\wedge \left( g_4 dz_1 + g_5 dz_2 + g_6 dz_3 \right) + d \bar{z}_3 \wedge \left( g_7 dz_1 + g_8 dz_2 + g_9 dz_3\right),\,\,\, g_1,...,g_9 \in \mathbb{C}
\end{align}
We note that $\dim_\mathbb{C} H^1\left( T, \mathcal{O}_T \right)=3$ and any element of $H^1\left( T, \mathcal{O}_T \right)$ has an unique representation of the form
\begin{align}\label{com10}
f_1 d\bar{z}_1 + f_2 d\bar{z}_2 + f_3 d\bar{z}_3,\,\,\,f_1,f_2,f_3\in \mathbb{C}
\end{align}

Now we compute $\mathbb{H}^1\left(T, \mathcal{E}_{\mathcal{O}_T}^\bullet \right)$. From $(\ref{com2})$, first we compute $H^0\left(T, \Theta_T \right) \bigoplus H^0\left( T, \mathcal{O}_T \right) \to H^0\left( T, \Omega_T^1 \right)$. Then from $(\ref{com3})$ and $(\ref{com6})$, the image is given by 
\begin{align} \label{com9}
\mathcal{L}_{b_1\frac{\partial}{\partial z_1} + b_2 \frac{\partial}{\partial z_2} + b_3 \frac{\partial}{\partial z_3}}\left( a_1 dz_1 + a_2 dz_2 + a_3 dz_3\right)+ h\cdot \left( a_1 dz_1 + a_2 dz_2 + a_3 dz_3  \right)= ha_1 dz_1 + ha_2 dz_2 + h a_3 dz_3
\end{align}
Then from $(\ref{com5})$ and $a_3\ne 0$, we have $\dim_\mathbb{C} \mathbb{H}^0\left( T, \mathcal{E}_{\mathcal{O}_T}^\bullet \right)= \dim_\mathbb{C} \textnormal{ker}\left(H^0\left( T, \Theta_T \right)\bigoplus H^0\left( T, \mathcal{O}_T \right)\to H^0\left( T , \Omega_T^1 \right) \right)= 3$ and $\dim_\mathbb{C} \textnormal{im}\left(H^0\left( T, \Theta_T \right) \bigoplus H^0\left( T, \mathcal{O}_T \right) \to H^0\left( T, \Omega_T^1 \right) \right)= 4-3 =1$. Next we compute $H^0\left( T, \Omega_T^1\right) \to H^0\left(T, \mathcal{O}_T \right)\cong H^0\left( T, \bigwedge^3 \Omega_T^1\right)\cong \mathbb{C}$. From $(\ref{com5})$, the image is given by
\begin{align*}
\left( c_1 dz_1 + c_2 dz_2 + c_3 dz_3 \right) \wedge d\omega_0+ \omega_0 \wedge d\left(c_1  dz_1+ c_2 dz_2 + c_3 dz_3 \right)=0
\end{align*}
Then the kernel is $H^0\left( T, \Omega_T^1\right)$. Then from $(\ref{com9})$, we see that since $a_3\ne 0$, 
\begin{align}\label{com11}
\frac{\textnormal{ker}\left( H^0\left( T, \Omega_T^1 \right) \to H^0\left( T, \mathcal{O}_T \right) \right) }{\textnormal{im}\left( H^0\left(T, \Theta_T \right)\bigoplus H^0\left( T, \mathcal{O}_T \right) \to H^0\left( T, \Omega_T^1 \right) \right)}\cong \textnormal{Span}_\mathbb{C}\left\langle  dz_1 ,  dz_2 \right\rangle
\end{align}

On the other hand, let us compute $H^1\left( T, \Theta_T \right)\bigoplus H^1\left( T, \mathcal{O}_T \right) \to H^1\left( T, \Omega_T^1 \right)$. From $(\ref{com7})$ and $(\ref{com10})$ the image is given by
{\small{\begin{align*}
&d\bar{z}_1\wedge \left(\mathcal{L}_{e_1\frac{\partial}{\partial z_1} + e_2 \frac{\partial}{\partial z_2} + e_3 \frac{\partial}{\partial z_3}}\left( \omega_0 \right) + f_1 \omega_0 \right) + d\bar{z}_2 \wedge \left(\mathcal{L}_{e_4\frac{\partial}{\partial z_1} + e_5 \frac{\partial}{\partial z_2} + e_6 \frac{\partial}{\partial z_3}}\left( \omega_0 \right) + f_2 \omega_0 \right) + d\bar{z}_3 \wedge \left(\mathcal{L}_{e_7 \frac{\partial}{\partial z_1} + e_8 \frac{\partial}{\partial z_2} + e_9 \frac{\partial}{\partial z_3}}\left( \omega_0 \right) + f_3 \omega_0 \right)\\
&= d\bar{z}_1\wedge f_1 \omega_0 + d\bar{z}_2 \wedge f_2 \omega+0 + d\bar{z}_3 \wedge f_3 \omega_0\\
&= d\bar{z}_1 \wedge \left( f_1a_1 dz_1 + f_1a_2 dz_2 + f_1 a_3 dz_3 \right)+ d\bar{z}_2 \wedge \left( f_2 a_1 dz_1 + f_2 a_2 dz_2 + f_2 a_3 dz_3 \right)  + d\bar{z}_3 \wedge \left( f_3 a_1 dz_1 + f_3a_2 dz_2 + f_3 a_3 dz_3 \right) =0\\
&\iff f_1=f_2=f_3=0
\end{align*}}}
Then since $a_3\ne 0$, we see that 
\begin{align} \label{com29}
\textnormal{ker}\left( H^1\left( T, \Theta_T \right) \bigoplus H^1\left( T, \mathcal{O}_T \right) \to H^1\left( T, \Omega_T^1 \right) \right)= H^1\left( T, \Theta_T \right). 
\end{align}
Hence from $(\ref{com2})$ and $(\ref{com11})$ and $(\ref{com7})$ we have
\begin{align}
\mathbb{H}^1\left( T, \mathcal{E}_{\mathcal{O}_T}^\bullet \right) &\cong \left\{ t_{10}dz_1  + t_{11} dz_2 , \right. \label{com12} \\
&\left. d\bar{z}_1\wedge \left( t_1 \frac{\partial}{\partial z_1} + t_2 \frac{\partial}{\partial z_3} + t_3 \frac{\partial}{\partial z_3} \right) + d\bar{z}_2 \wedge \left( t_4 \frac{\partial}{\partial z_1} + t_5 \frac{\partial}{\partial z_3} + t_6 \frac{\partial}{\partial z_3}\right) + d\bar{z}_3 \wedge \left( t_7 \frac{\partial}{\partial z_1} + t_8 \frac{\partial}{\partial z_3} + t_9 \frac{\partial}{\partial z_3} \right) \right. \notag \\
&\left. | t_1,..., t_{11} \in \mathbb{C} \right\} \notag
\end{align}

We will construct a foliated analytic family of deformations of $\left( T, \omega_0 \right)$, where $\omega_0 = a_1 dz_1 + a_2 dz_2 + a_3 dz_3, a_3\ne 0$ such that the foliated Kodaira-Spencer map in terms of cotangent sheaves is an isomorphism at the distinguished point, so that $\left( T, \omega_0 \right)$ is unobstructed in foliated deformations. In fact, we recall the integrability conditions for a foliated analytic family in terms of cotangent sheaves from $(\ref{oni1})-(\ref{oni4})$. If we take
\begin{align*}
\left(Y(t), \left(\varphi(t), Z(t) \right) \right) \in A^{0,0}\left( T, \Omega_T^1 \right) \bigoplus A^{0,1}\left( T, \Theta_T \right) \bigoplus A^{0,1}\left( T, \mathcal{O}_T \right)
\end{align*}
where 
\begin{align*}
Y(t) & = \left( a_1 + t_{10}   \right) dz_1 + \left( a_2 + t_{11}  \right) dz_2 + a_3 dz_3\\
\varphi(t) &= d\bar{z}_1\wedge \left( t_1 \frac{\partial}{\partial z_1} + t_2 \frac{\partial}{\partial z_3} + t_3 \frac{\partial}{\partial z_3} \right) + d\bar{z}_2 \wedge \left( t_4 \frac{\partial}{\partial z_1} + t_5 \frac{\partial}{\partial z_3} + t_6 \frac{\partial}{\partial z_3}\right) + d\bar{z}_3 \wedge \left( t_7 \frac{\partial}{\partial z_1} + t_8 \frac{\partial}{\partial z_3} + t_9 \frac{\partial}{\partial z_3} \right)\\
Z(t)& = 0 
\end{align*}
Then we have $\varphi(0)=0,Y(0)=\omega_0$ and $Z(0)=0$ and the following equalities hold:
\begin{align*}
\bar{\partial}\varphi(t)- \frac{1}{2}\left[ \varphi(t), \varphi(t) \right] = 0\\
\bar{\partial} Y(t) - \mathcal{L}_{\varphi(t)}\left( Y(t) \right)- Z(t) \wedge Y(t)=0\\
Y(t)\wedge \partial Y(t)=0\\
\bar{\partial}Z(t)- \left[ \varphi(t), Z(t) \right]= 0
\end{align*}
In other words, the integrability conditions $(\ref{oni1})-(\ref{oni4})$ are satisfied, so that there exists a foliated analytic family $\left(\mathcal{M}, \omega , B, \pi \right)$ where $\pi:\left( \mathcal{M}, \omega \right)\to B$ for a sufficiently small neighborhood $B$ of $0$ in $\mathbb{C}^{11}$ with coordinates $t=\left(t_1,t_2,t_3,t_4,t_5,t_6,t_7,t_8, t_9, t_{10}, t_{11} , t_{12} \right)$ such that $\pi^{-1}(0)=\left( T, \omega_0 \right)$ and the foliated Kodaira-Spencer map from $(\ref{cu14})$ given by
\begin{align*}
\varphi_0: T_0 (B) &\to \mathbb{H}^1\left( T, \mathcal{E}_{\mathcal{O}_T}^\bullet \right) \\
 \frac{\partial}{\partial t_1}  &\mapsto \left( 0, \left( \frac{\partial}{\partial z_1} d\bar{z}_1,0 \right) \right), \frac{\partial}{\partial t_2} \mapsto \left( 0, \left( \frac{\partial}{\partial z_2} d\bar{z}_1 ,0\right) \right), \frac{\partial}{\partial t_3} \mapsto \left( 0, \left( \frac{\partial}{\partial z_3} d\bar{z}_1 ,0 \right) \right)\\
  \frac{\partial}{\partial t_4}  &\mapsto \left( 0, \left( \frac{\partial}{\partial z_1} d\bar{z}_2 , 0 \right) \right), \frac{\partial}{\partial t_5} \mapsto \left( 0, \left( \frac{\partial}{\partial z_2} d\bar{z}_2 , 0 \right)  \right), \frac{\partial}{\partial t_6} \mapsto \left( 0, \left( \frac{\partial}{\partial z_3} d\bar{z}_2 , 0 \right) \right)\\
   \frac{\partial}{\partial t_7}  &\mapsto \left( 0, \left( \frac{\partial}{\partial z_1} d\bar{z}_3 , 0 \right) \right), \frac{\partial}{\partial t_8} \mapsto \left( 0, \left( \frac{\partial}{\partial z_2} d\bar{z}_3 , 0  \right) \right), \frac{\partial}{\partial t_9} \mapsto \left( 0, \left( \frac{\partial}{\partial z_3} d\bar{z}_3 , 0  \right) \right)\\
   \frac{\partial}{\partial t_{10} } & \mapsto \left(dz_1 , (0,0) \right), \frac{\partial}{\partial t_{11}} \mapsto \left( dz_2 , (0,0) \right)
\end{align*}
is an isomorphism at $t=0$ by $(\ref{com12})$, so that $\left( T, \omega_0 \right)$ is unobstructed in foliated deformations.

\subsection{Foliated deformations of $\left(T, \omega_0= a_1 dz_1 + a_2 dz_2 + a_3 dz_3 \right),a_3\ne 0$ in terms of tangent sheaves} \

Let us consider the foliation on $T$ defined by $a_1 dz_1 + a_2 dz_2 + a_3 dz_3\in H^0\left( T, \Omega_T^1 \right), a_3\ne 0$. Then $\Theta_{\mathcal{F}_0} \cong \mathcal{O}_T \bigoplus \mathcal{O}_T$ is generated by $T_1= a_3 \frac{\partial}{\partial z_1} - a_1 \frac{\partial}{\partial z_3}$ and $T_2 = a_3 \frac{\partial}{\partial z_2} - a_2 \frac{\partial}{\partial z_3}$ with $\left[ T_1, T_2 \right]=0$ and we have the exact sequence
\begin{align*}
0 \to \mathcal{O}_T \bigoplus \mathcal{O}_T &\to \Theta_T \to \frac{\Theta_T}{\mathcal{O}_T \bigoplus \mathcal{O}_T }\cong  \mathcal{O}_T \to 0 \\
(1,0) &\mapsto T_1= a_3 \frac{\partial}{\partial z_1} - a_1 \frac{\partial}{\partial z_3}\\
(0,1) & \mapsto T_2 = a_3 \frac{\partial}{\partial z_2} - a_2 \frac{\partial}{\partial z_3}
\end{align*} 
Then we have the leaf complex (see Part II)
\begin{align}
\mathcal{E}_{\mathcal{O}_T\bigoplus \mathcal{O}_T}^\bullet: \mathcal{E}_{\mathcal{O}_T \bigoplus \mathcal{O}_T}\cong \Theta_T \bigoplus \mathscr{H}om_{\mathcal{O}_T}\left( \mathcal{O}_T \bigoplus \mathcal{O}_T, \mathcal{O}_T \bigoplus \mathcal{O}_T \right) \to \mathscr{H}om_{\mathcal{O}_T}\left( \mathcal{O}_T \bigoplus \mathcal{O}_T, \Theta_T \right) \to  \\
\to \mathscr{H}om_{\mathcal{O}_T}\left(\bigwedge^2 \left( \mathcal{O}_T \bigoplus \mathcal{O}_T \right), \frac{\Theta_T}{\mathcal{O}_T\bigoplus \mathcal{O}_T} \right)\cong \mathcal{O}_T \to 0 \notag
\end{align}
Let us consider the spectral sequence at $E_1$ from Dolbeault bicomplex associated to $\mathcal{E}_{\mathcal{O}_T \bigoplus \mathcal{O}_T}^\bullet$:
\begin{equation} \label{com13}
\begin{CD}
0@. \,\,\, 0 \\ 
@AAA @AAA \\
H^0\left(T, \mathcal{O}_T \right) @.\,\,\,\,\,\,\,\,\,\,\,H^1\left( T, \mathcal{O}_T  \right) \\
@A\tilde{D}_1AA @A\tilde{D}_1AA \\
H^0\left( T, \mathscr{H}om_{\mathcal{O}_T }\left( \mathcal{O}_T \bigoplus \mathcal{O}_T , \Theta_T \right) \right) @.\,\,\,\,\,\,\,\,\,\,\,\, H^1\left( T , \mathscr{H}om_{\mathcal{O}_T}\left( \mathcal{O}_T \bigoplus \mathcal{O}_T , \Theta_T  \right) \right) @.\,\,\,\,\,\,\,\,\,\, \\
@A\tilde{D}_0 AA @A\tilde{D}_0 AA \\
H^0\left(T, \Theta_T \right) \bigoplus H^0\left(T, \mathscr{H}om_{\mathcal{O}_T}\left( \mathcal{O}_T \bigoplus\mathcal{O}_T, \mathcal{O}_T \bigoplus \mathcal{O}_T \right)\right) @. \,\,\,\,\,\,\,\,\,\,\,\,\,H^1\left(T, \Theta_T \right)\bigoplus  H^1 \left(T, \mathscr{H}om_{\mathcal{O}_T}\left( \mathcal{O}_T \bigoplus\mathcal{O}_T, \mathcal{O}_T \bigoplus \mathcal{O}_T \right)\right)@.
\end{CD}
\end{equation}
Then we have
{\small{\begin{align} \label{com18}
\mathbb{H}^1\left( T, \mathcal{E}_{\mathcal{O}_T \bigoplus \mathcal{O}_T}^\bullet \right) &\cong  \frac{\textnormal{ker}\left(H^0\left( T, \mathscr{H}om_{\mathcal{O}_T}\left( \mathcal{O}_T \bigoplus \mathcal{O}_T, \Theta_T \right) \right) \to H^0\left( T , \mathcal{O}_T \right) \right)}{\textnormal{im}\left( H^0\left( T , \Theta_T \right)\bigoplus H^0\left( T, \mathscr{H}om_{\mathcal{O}_T}\left( \mathcal{O}_T \bigoplus \mathcal{O}_T, \mathcal{O}_T \bigoplus \mathcal{O}_T \right) \to   H^0\left( T , \mathscr{H}om_{\mathcal{O}_T}\left( \mathcal{O}_T \bigoplus \mathcal{O}_T, \Theta_T \right) \right) \right) \right)} \\
 &\bigoplus \textnormal{ker}\left( H^1\left( T , \Theta_T \right)\bigoplus H^1 \left( T , \mathscr{H}om_{\mathcal{O}_T }\left( \mathcal{O}_T \bigoplus \mathcal{O}_T , \mathcal{O}_T \bigoplus \mathcal{O}_T \right) \to   H^1 \left( T , \mathscr{H}om_{\mathcal{O}_T}\left( \mathcal{O}_T \bigoplus \mathcal{O}_T , \Theta_T \right) \right) \right)  \right) \notag
\end{align}}}

We will describe each cohomology group in $(\ref{com13})$. We note that $\dim_\mathbb{C} H^0\left( T, \mathscr{H}om_{\mathcal{O}_T}\left( \mathcal{O}_T \bigoplus \mathcal{O}_T, \mathcal{O}_T, \mathcal{O}_T \right)\right)= \dim_\mathbb{C} H^0\left( T, \bigoplus^4 \mathcal{O}_T \right)= 4$ in $(\ref{com13})$ and any element $\psi$ of $H^0\left( T , \mathscr{H}om_{\mathcal{O}_T }\left( \mathcal{O}_T \bigoplus \mathcal{O}_T , \mathcal{O}_T \bigoplus \mathcal{O}_T \right) \right)$ has an unique representation of the form
\begin{align}\label{com14}
\psi(T_1)&= f_1 T_1 + f_2 T_2\\
\psi(T_2)&= f_3 T_1 + f_4 T_2,\,\,\,\,\,f_1,f_2,f_3,f_4 \in \mathbb{C} \notag
\end{align}

We note that $\dim_\mathbb{C} H^1\left( T , \mathscr{H}om_{\mathcal{O}_T}\left(  \mathcal{O}_T \bigoplus \mathcal{O}_{\mathcal{O}_T}, \mathcal{O}_T \bigoplus \mathcal{O}_T \right) \right)=\dim_\mathbb{C} H^1\left( T,\bigoplus^4 \mathcal{O}_T\right) = 12$ and any element $\psi$ of $H^1 \left( X, \mathscr{H}om_{\mathcal{O}_T}\left(  \mathcal{O}_T \bigoplus \mathcal{O}_T,  \mathcal{O}_T\bigoplus \mathcal{O}_T \right)\right)$ has an unique representation of the form
\begin{align}\label{com15}
\psi\left(T_1 \right)&=  d\bar{z}_1\wedge \left(g_1 T_1+ g_2 T_2 \right) + d\bar{z}_2 \wedge \left( g_3 T_1 + g_4 T_2 \right) + d\bar{z}_3 \wedge \left(g_5 T_1 + g_6 T_2 \right) \\
\psi\left(T_2 \right)&= d\bar{z}_1\wedge \left( g_7 T_1+ g_8 T_2 \right) + d\bar{z}_2 \wedge \left( g_9 T_2 + g_{10} T_2 \right) + d\bar{z}_3 \wedge \left( g_{11} T_1 + g_{12} T_2 \right) ,\,\,\,\,\,g_1,...,g_{12}\in \mathbb{C} \notag
\end{align}

We note that $\dim_\mathbb{C} H^0\left(T, \mathscr{H}om_{\mathcal{O}_T }\left( \mathcal{O}_T \bigoplus \mathcal{O}_T , \Theta_T \right) \right)= \dim_\mathbb{C} H^0\left( T , \Theta_T \bigoplus \Theta_T \right) = 6$ in $(\ref{com13})$ and any element $\rho$ of $H^0\left(T, \mathscr{H}om_{\mathcal{O}_T}\left( \mathcal{O}_T \bigoplus \mathcal{O}_T , \Theta_T \right) \right)$ has an unique representation in the form
\begin{align}\label{com16}
\rho(T_1)&= h_1 \frac{\partial}{\partial z_1} + h_2 \frac{\partial}{\partial z_2}+  h_3 \frac{\partial}{\partial z_3},\,\,\,\,\,\,\,\,\,
\rho(T_1)= h_4 \frac{\partial}{\partial z_1} + h_5 \frac{\partial}{\partial z_2} + h_6 \frac{\partial}{\partial z_3} ,\,\,\,\,\, h_1,...,h_6 \in \mathbb{C} 
\end{align}

We note that $\dim_\mathbb{C} H^1\left(T, \mathscr{H}om_{\mathcal{O}_T}\left( \mathcal{O}_T \bigoplus \mathcal{O}_T, \Theta_T \right) \right)=\dim_\mathbb{C} H^1\left(T, \Theta_T \bigoplus \Theta_T \right) = 18$ and any element $\sigma$ of $H^1\left(T, \mathscr{H}om_{\mathcal{O}_T}\left( \mathcal{O}_T \bigoplus \mathcal{O}_T, \Theta_T \right) \right)$ has an unique representation of the form
\begin{align} \label{com17}
\sigma\left( T_1 \right)&= d\bar{z}_1\wedge\left( l_1 \frac{\partial}{\partial z_1} + l_2 \frac{\partial}{\partial z_2}  + l_3 \frac{\partial}{\partial z_3} \right) + d \bar{z}_2 \wedge \left( l_4 \frac{\partial}{\partial z_1} + l_5 \frac{\partial}{\partial z_2}  + l_6 \frac{\partial}{\partial z_3} \right) + d \bar{z}_3 \wedge \left( l_7 \frac{\partial}{\partial z_1} + l_8 \frac{\partial}{\partial z_2} + l_9 \frac{\partial}{\partial z_3} \right) \\
\sigma\left( T_2 \right)&= d\bar{z}_1\wedge\left( l_{10} \frac{\partial}{\partial z_1} + l_{11} \frac{\partial}{\partial z_2}  + l_{12} \frac{\partial}{\partial z_3} \right) + d \bar{z}_2 \wedge \left( l_{13} \frac{\partial}{\partial z_1} + l_{14} \frac{\partial}{\partial z_2}  + l_{15} \frac{\partial}{\partial z_3} \right)  + d \bar{z}_3 \wedge \left( l_{16} \frac{\partial}{\partial z_1} +  l_{17} \frac{\partial}{\partial z_2} + l_{18} \frac{\partial}{\partial z_3} \right) \notag
\end{align}
where $l_1,..,l_{18} \in \mathbb{C}$.

 We compute $\mathbb{H}^1\left(  T, \mathcal{E}_{\mathcal{O}_T \bigoplus \mathcal{O}_T}^\bullet \right)$. From $(\ref{com18})$, we compute $H^0\left( T, \Theta_T \right)\bigoplus H^0\left(T, \mathscr{H}om_{\mathcal{O}_T}\left( \mathcal{O}_T \bigoplus \mathcal{O}_T, \mathcal{O}_T \bigoplus \mathcal{O}_T \right) \right) \xrightarrow{\tilde{D}_0} H^0\left(T, \mathscr{H}om_{\mathcal{O}_T}\left(\mathcal{O}_T \bigoplus \mathcal{O}_T , \Theta_T \right) \right)$ first. Then from $(\ref{com3})$ and $(\ref{com14})$, the image is given by $\rho$, where
 \begin{align}\label{com21}
 \rho(T_1)&= \left[ b_1 \frac{\partial}{\partial z_1} + b_2 \frac{\partial}{\partial z_2} + b_3 \frac{\partial }{\partial z_3},  a_3 \frac{\partial}{\partial z_1} - a_1 \frac{\partial}{\partial z_3} \right] - \left(f_1 T_1 + f_2 T_2 \right) = - f_1 T_1 - f_2 T_2 = 0 \iff f_1= f_2=0\\
 \rho(T_2)&=  \left[ b_1 \frac{\partial}{\partial z_1} + b_2 \frac{\partial}{\partial z_2} + b_3 \frac{\partial }{\partial z_3},  a_3 \frac{\partial}{\partial z_2} - a_2 \frac{\partial}{\partial z_3} \right] - \left(f_3 T_1 + f_4 T_2 \right) = - f_3 T_1 - f_4 T_2 = \iff f_3= f_4 = 0  \notag 
 \end{align}
 Then we have $\dim_\mathbb{C} \mathbb{H}^0\left( T, \mathcal{E}_{\mathcal{O}_T \bigoplus \mathcal{O}_T}^\bullet \right)= \dim_\mathbb{C} H^0\left( T, \Theta_T \right)=3$, and so
{\small{\begin{align} \label{com20}
\dim_\mathbb{C} \textnormal{im}\left( H^0\left( T, \Theta_T \right)\bigoplus H^0\left( T, \mathscr{H}om_{\mathcal{O}_T}\left( \mathcal{O}_T \bigoplus \mathcal{O}_T, \mathcal{O}_T \bigoplus \mathcal{O}_T \right) \right) \xrightarrow{\tilde{D}_0} H^0\left( T, \mathscr{H}om_{\mathcal{O}_T}\left( \mathcal{O}_T \bigoplus \mathcal{O}_T, \Theta_T \right) \right) \right)  = 7-3 = 4
 \end{align}}}
  Next we compute $H^0\left( T, \mathscr{H}om_{\mathcal{O}_T} \left( \mathcal{O}_T \bigoplus \mathcal{O}_T , \Theta_T \right) \right) \xrightarrow{\tilde{D}_1} H^0\left( T, \mathscr{H}om_{\mathcal{O}_T} \left( \mathcal{O}_T, \mathcal{O}_T \right)\right)$. We note that for $\rho\in H^0\left( T, \mathscr{H}om_{\mathcal{O}_T}\left( \mathcal{O}_T \bigoplus \mathcal{O}_T, \Theta_T \right) \right)$ given by $(\ref{com16})$, we have
  \begin{align*}
 \left[ T_1, \rho\left( T_2 \right)\right] - \left[ T_2, \rho\left( T_1 \right)\right] - \rho\left( \left[ T_1, T_2\right] \right) =0
  \end{align*}
 In other words, $\tilde{D}_1$ is the zero map, so that $\textnormal{ker}\left( \tilde{D}_1 \right)= H^0\left(T, \mathscr{H}om_{\mathcal{O}_T}\left( \mathcal{O}_T \bigoplus \mathcal{O}_T, \Theta_T \right) \right)$. Then we have from $(\ref{com20})$
 \begin{align*}
 \dim_\mathbb{C} \frac{ H^0\left(T, \mathscr{H}om_{\mathcal{O}_T}\left(\mathcal{O}_T \bigoplus \mathcal{O}_T, \Theta_T \right) \right)=\textnormal{ker}\left( H^0\left( T, \mathscr{H}om_{\mathcal{O}_T}\left( \mathcal{O}_T \bigoplus \mathcal{O}_T , \Theta_T \right) \right) \to H^0\left( T, \mathcal{O}_T \right) \right)}{\textnormal{im}\left( H^0\left( T, \Theta_T \right) \bigoplus H^0\left( T, \mathscr{H}om_{\mathcal{O}_T}\left( \mathcal{O}_T \bigoplus \mathcal{O}_T, \mathcal{O}_T \bigoplus \mathcal{O}_T \right) \right) \to H^0\left( T, \mathscr{H}om_{\mathcal{O}_T}\left( \mathcal{O}_T \bigoplus \mathcal{O}_T, \Theta_T \right) \right)  \right)}= 6 - 4 = 2
 \end{align*}
 and from $(\ref{com21})$ we can check that $ \frac{H^0\left( T, \mathscr{H}om_{\mathcal{O}_T}\left( \mathcal{O}_T \bigoplus \mathcal{O}_T , \Theta_T \right) \right) }{\textnormal{im}\left( H^0\left( T, \Theta_T \right) \bigoplus H^0\left( T, \mathscr{H}om_{\mathcal{O}_T}\left( \mathcal{O}_T \bigoplus \mathcal{O}_T, \mathcal{O}_T \bigoplus \mathcal{O}_T \right) \right) \to H^0\left( T, \mathscr{H}om_{\mathcal{O}_T}\left( \mathcal{O}_T \bigoplus \mathcal{O}_T, \Theta_T \right) \right)  \right)}$ is generated by $\phi, \psi\in H^0\left( T, \mathscr{H}om_{\mathcal{O}_T}\left( \mathcal{O}_T \bigoplus \mathcal{O}_T, \Theta_T \right) \right)$, where
 \begin{align}\label{com23}
 \phi(T_1)&=   -  \frac{\partial}{\partial z_3},\,\,\,\,\,\,\,\,\phi(T_2)= 0 \\
 \psi(T_1)&= 0,\,\,\,\,\,\,\,\,\,\psi(T_2)= -  \frac{\partial}{\partial z_3} \notag
 \end{align}
 
 Let us compute $H^1\left( T, \Theta_T \right)\bigoplus H^1\left( T, \mathscr{H}om_{\mathcal{O}_T}\left( \mathcal{O}_T \bigoplus \mathcal{O}_T, \mathcal{O}_T \bigoplus \mathcal{O}_T \right)\right)\xrightarrow{\tilde{D}_0} H^1\left( T, \mathscr{H}om_{\mathcal{O}_T}\left( \mathcal{O}_T \bigoplus \mathcal{O}_T, \Theta_T \right) \right)$. Then from $(\ref{com7})$ and $(\ref{com15})$ the image is given by $\sigma$, where
{\small{ \begin{align*}
 \sigma\left(T_1 \right)&=d\bar{z}_1\wedge \left( \left[e_1 \frac{\partial}{\partial z_1} + e_2 \frac{\partial}{\partial z_2} + e_3 \frac{\partial}{\partial z_3}, T_1 \right] \right)+ d\bar{z}_2\wedge \left( \left[e_4 \frac{\partial}{\partial z_1} + e_5 \frac{\partial}{\partial z_2} + e_6 \frac{\partial}{\partial z_3}, T_1  \right] \right) + d\bar{z}_3\wedge \left( \left[e_7 \frac{\partial}{\partial z_1} + e_8 \frac{\partial}{\partial z_2} + e_9 \frac{\partial}{\partial z_3},  T_1 \right] \right) \\
 &-d\bar{z}_1\wedge \left( g_1 T_1 + g_2 T_2 \right) - d\bar{z}_2\wedge \left( g_3 T_1 + g_4 T_2 \right) - d \bar{z}_3 \wedge \left( g_5 T_1 + g_6 T_2 \right) = 0 \iff g_1=g_2=g_3=g_4=g_5=g_6=0\\
  \sigma\left(T_2 \right)&=d\bar{z}_1\wedge \left( \left[e_1 \frac{\partial}{\partial z_1} + e_2 \frac{\partial}{\partial z_2} + e_3 \frac{\partial}{\partial z_3}, T_2 \right] \right)+ d\bar{z}_2\wedge \left( \left[e_4 \frac{\partial}{\partial z_1} + e_5 \frac{\partial}{\partial z_2} + e_6 \frac{\partial}{\partial z_3}, T_2  \right] \right) + d\bar{z}_3\wedge \left( \left[e_7 \frac{\partial}{\partial z_1} + e_8 \frac{\partial}{\partial z_2} + e_9 \frac{\partial}{\partial z_3},  T_2 \right] \right) \\
 &-d\bar{z}_1\wedge \left( g_7 T_1 + g_8 T_2 \right) - d\bar{z}_2\wedge \left( g_9 T_1 + g_{10} T_2 \right) - d \bar{z}_3 \wedge \left( g_{11} T_1 + g_{12} T_2 \right) =0 \iff g_7=g_8=g_9=g_{10}=g_{11}=g_{12}=0
 \end{align*}}}
Then from $(\ref{com17})$ we see that the kernel is $H^1\left( T, \Theta_T \right)$. Then from $(\ref{com18})$ and $(\ref{com23})$ we have
{\small{\begin{align} \label{com25}
\mathbb{H}^1\left( T, \mathcal{E}_{\mathcal{O}_T \bigoplus \mathcal{O}_T}^\bullet \right) \cong &\left\{t_{10} \phi + t_{11} \psi   , \right. \\
& d\bar{z}_1\wedge\left(t_1\frac{\partial}{\partial z_1} + t_2 \frac{\partial}{\partial z_2} + t_3 \frac{\partial}{\partial z_3}  \right) +  d\bar{z}_2\wedge\left(t_4\frac{\partial}{\partial z_1} + t_5 \frac{\partial}{\partial z_2} + t_6 \frac{\partial}{\partial z_3}  \right) + d\bar{z}_3 \wedge\left(t_7\frac{\partial}{\partial z_1} + t_{8} \frac{\partial}{\partial z_2} + t_{9} \frac{\partial}{\partial z_3}  \right) \notag\\
&\left. | t_1,...,t_{11}\in \mathbb{C} \right\} \notag 
\end{align}}}

We will construct a foliated analytic family of deformations of $\left(X, \omega_0 \right)$ where $\omega_0= a_1dz_1 + a_2 dz_2 + a_3 dz_3 ,a_3\ne 0$ in terms of tangent sheaves such that the foliated Kodaira-Spencer map in terms of tangent sheaves is an isomorphism at the distinguished point, so that $\left( X , \omega_0 \right)$ is unobstructed in foliated deformations. In fact, we recall the integrability conditions for a foliated analytic family in terms of tangent sheaves from $(\ref{cu31})-(\ref{cu32})$. If we take
{\small{\begin{align*}
\left( X(t), \left( \varphi(t), U(t)=0 \right) \right) \in A^{0,0}\left( T, \mathscr{H}om_{\mathcal{O}_T}\left( \mathcal{O}_T \bigoplus \mathcal{O}_T , \Theta_T \right)\right)\bigoplus A^{0,1}\left( T , \Theta_T \right) \bigoplus A^{01}\left( T, \mathscr{H}om_{\mathcal{O}_T}\left(\mathcal{O}_T \bigoplus \mathcal{O}_T, \mathcal{O}_T \bigoplus \mathcal{O}_T \right) \right)
\end{align*}}}
where
\begin{align*}
X(t)\left(T_1\right)&= T_1(t):= a_3 \frac{\partial}{\partial z_1} - \left( a_1 +  t_{10} \right) \frac{\partial}{\partial z_3}\\
X(t)\left(T_2\right)&= T_2(t):= a_3 \frac{\partial}{\partial z_2} - \left( a_2 +  t_{11} \right) \frac{\partial}{\partial z_3}\\
\varphi(t)&:= d\bar{z}_1\wedge\left(t_1 \frac{\partial}{\partial z_1} + t_2 \frac{\partial}{\partial z_2} + t_3 \frac{\partial}{\partial z_3}  \right) +  d\bar{z}_2\wedge\left(t_4 \frac{\partial}{\partial z_1} + t_5 \frac{\partial}{\partial z_2} + t_6 \frac{\partial}{\partial z_3}  \right) + d\bar{z}_3 \wedge\left(t_7 \frac{\partial}{\partial z_1} + t_{8} \frac{\partial}{\partial z_2} + t_{9} \frac{\partial}{\partial z_3}  \right)
\end{align*}
Then we have $\varphi(0)=0, X(0)=\textnormal{Id}$ and the integrability conditions $(\ref{cu31})-(\ref{cu32})$ are satisfied, so that there exists a foliated analytic family $\left(\mathcal{M}, \Theta_\mathcal{F} , B, \pi \right)$ where $\pi:\left( \mathcal{M}, \Theta_{\mathcal{F}} \right)\to B$ for a sufficiently small neighborhood $B$ of $0$ in $\mathbb{C}^{11}$ with coordinates $t=\left(t_1,t_2,t_3,t_4,t_5,t_6,t_7,t_8, t_9, t_{10}, t_{11} \right)$ such that $\pi^{-1}(0)=\left( T , \omega_0 \right)$ and the foliated Kodaira-Spencer map in terms of tangent sheaves from $(\ref{cu33})$ given by
\begin{align*}
\varphi: T_0(B) &\to \mathbb{H}^1\left( T, \mathcal{E}_{\mathcal{O}_T \bigoplus \mathcal{O}_T}^\bullet \right) \\
 \frac{\partial}{\partial t_1}  &\mapsto \left( 0, \left( \frac{\partial}{\partial z_1} d\bar{z}_1,0 \right) \right), \frac{\partial}{\partial t_2} \mapsto \left( 0, \left( \frac{\partial}{\partial z_2} d\bar{z}_1 ,0\right) \right), \frac{\partial}{\partial t_3} \mapsto \left( 0, \left( \frac{\partial}{\partial z_3} d\bar{z}_1 ,0 \right) \right)\\
  \frac{\partial}{\partial t_4}  &\mapsto \left( 0, \left( \frac{\partial}{\partial z_1} d\bar{z}_2 , 0 \right) \right), \frac{\partial}{\partial t_5} \mapsto \left( 0, \left( \frac{\partial}{\partial z_2} d\bar{z}_2 , 0 \right)  \right), \frac{\partial}{\partial t_6} \mapsto \left( 0, \left( \frac{\partial}{\partial z_3} d\bar{z}_2 , 0 \right) \right)\\
   \frac{\partial}{\partial t_7}  &\mapsto \left( 0, \left( \frac{\partial}{\partial z_1} d\bar{z}_3 , 0 \right) \right), \frac{\partial}{\partial t_8} \mapsto \left( 0, \left( \frac{\partial}{\partial z_2} d\bar{z}_3 , 0  \right) \right), \frac{\partial}{\partial t_9} \mapsto \left( 0, \left( \frac{\partial}{\partial z_3} d\bar{z}_3 , 0  \right) \right)\\
\frac{\partial}{\partial t_{10} } &\mapsto \left( \phi, (0,0) \right),\,\,\,\,\,\,\textnormal{where}\,\,\,\phi(T_1)= -  \frac{\partial}{\partial z_3},\,\,\,\phi(T_2)= 0 \\
\frac{\partial}{\partial t_{11} } &\mapsto \left( \psi, (0,0) \right),\,\,\,\,\,\,\textnormal{where}\,\,\,\psi(T_1)= 0 ,\,\,\,\psi(T_2)=  -  \frac{\partial}{\partial z_3}
\end{align*}
is an isomorphism at $t=0$ by $(\ref{com25})$, so that $\left( T, \omega_0 \right)$ is unobstructed in foliated deformations in terms of tangent sheaves.

\subsection{Unfoldings of $\left(T, \omega_0 \right)$, where $\omega_0= a_1dz_1 + a_2 dz_2 + a_3 dz_3, a_3\ne 0$}\

Let us consider the foliation $\omega_0$ on $T$ where $\omega_0= a_1 dz_1 + a_2 dz_2 + a_3 dz_3, a_3 \ne 0 \in H^0\left( T, \Omega_T^1 \right)$. Then we have the unfolded dual leaf complex from $(\ref{cp3})$
\begin{align*}
\tilde{\mathcal{E}}_{O_T}^\bullet : \Theta_T \bigoplus \mathcal{O}_T \xrightarrow{\tilde{G}_0}  \Omega_T^1 \bigoplus \mathcal{O}_T \xrightarrow{\tilde{G}_1} \mathcal{O}_T \bigoplus \mathscr{H}om_{\mathcal{O}_T} \left( \mathcal{N}_{\mathcal{F}_0}^*, \tilde{\mathcal{S}}^1 \right) \left(\subset \mathcal{O}_T \bigoplus \bigwedge^{2}\Omega_T^1\right) \xrightarrow{G_2} \mathcal{O}_T \to 0
\end{align*}

Let us consider the spectral sequence at $E_1$ from Dolbeault bicomplex associated to $\tilde{\mathcal{E}}_{\mathcal{O}_T}^\bullet$:
\begin{center}
$\begin{CD}
0 \\
@AAA \\
H^0\left( T, \mathcal{O}_T \right) \\
@AG_2AA \\
H^0\left(T, \mathcal{O}_T \right)\bigoplus H^0\left(T, \tilde{\mathcal{S}}^1 \right) \subset H^0\left( T, \mathcal{O}_T \right) \bigoplus H^0\left(T, \bigwedge^2 \Omega_T^1 \right)\\
@A\tilde{G}_1AA \\
H^0\left( T, \Omega_T^1 \right) \bigoplus H^0\left( T, \mathcal{O}_T \right) @. \,\,\,\,\,\,\,\,\,\,\,\,\,\, H^1\left( T, \Omega_T^1 \right) \bigoplus H^1\left( T, \mathcal{O}_T \right)  \\
@A\tilde{G}_0 AA @A\tilde{G}_0AA\\
H^0\left(T, \Theta_T \right)\bigoplus H^0\left(T, \mathcal{O}_T \right) @.\,\,\,\,\,\,\,\,\,\,\,\,\, H^1\left( T, \Theta_T \right) \bigoplus H^1\left(T, \mathcal{O}_T \right)
\end{CD}$
\end{center}
Then we have
\begin{align}\label{com26}
\mathbb{H}^1\left( T, \tilde{\mathcal{E}}_{\mathcal{O}_T}^\bullet \right) \cong & \frac{\textnormal{ker}\left( H^0\left( T, \Omega_T^1 \right)\bigoplus H^0\left( T, \mathcal{O}_T \right) \xrightarrow{\tilde{G}_1} H^0\left( T, \mathcal{O}_T \right) \bigoplus H^0\left(T, \bigwedge^2 \Omega_T^1 \right) \right)}{\textnormal{im}\left(H^0\left( T, \Theta_T \right) \bigoplus H^0\left( T, \mathcal{O}_T \right) \xrightarrow{\tilde{G}_0} H^0\left( T, \Omega_T^1 \right) \bigoplus H^0\left( T, \mathcal{O}_T \right) \right)}\\
&\bigoplus \textnormal{ker}\left( H^1\left( T, \Theta_T \right) \bigoplus H^1\left( T, \mathcal{O}_T \right) \xrightarrow{\tilde{G}_0} H^1\left( T, \Omega_T^1 \right) \bigoplus H^1\left(T, \mathcal{O}_T \right) \right) \notag
\end{align}

We compute $\mathbb{H}^1\left(T, \tilde{\mathcal{E}}_{\mathcal{O}_T}^\bullet \right)$. From $(\ref{com26})$, first we compute $H^0\left( T, \Theta_T \right) \bigoplus H^0\left( T, \mathcal{O}_T \right) \xrightarrow{\tilde{G}_0} H^0\left( T, \Omega_T^1 \right) \bigoplus H^0\left( T, \mathcal{O}_T \right) $. We note that from $(\ref{com9})$ the image of $H^0\left( T, \Theta_T \right) \bigoplus H^0\left( T, \mathcal{O}_T \right) \to H^0\left( T, \Omega_T^1 \right)  $ is given by 
\begin{align*} 
\mathcal{L}_{b_1\frac{\partial}{\partial z_1} + b_2 \frac{\partial}{\partial z_2} + b_3 \frac{\partial}{\partial z_3}}\left( a_1 dz_1 + a_2 dz_2 + a_3 dz_3\right)+ h\cdot \left( a_1 dz_1 + a_2 dz_2 + a_3 dz_3  \right)= ha_1 dz_1 + ha_2 dz_2 + h a_3 dz_3=0 \iff h=0
\end{align*}
and the image of $H^0\left( T, \Theta_T \right)\bigoplus H^0\left(T, \mathcal{O}_T \right) \to H^0\left(T, \mathcal{O}_T \right)$ is given by
\begin{align*}
i_{b_1 \frac{\partial}{\partial z_1} + b_2 \frac{\partial}{\partial z_2} + b_3 \frac{\partial}{\partial z_3}}\left( a_1dz_1 + a_2 dz_2 + a_3 dz_3 \right)= a_1b_1 + a_2 b_2 + a_3 b_3
\end{align*}
Then we see that $\dim_\mathbb{C} \mathbb{H}^0\left( T, \tilde{\mathcal{E}}_{\mathcal{O}_T}^\bullet \right)= 2$, so that we have 
\begin{align}\label{com27}
\dim_\mathbb{C} \textnormal{im}\left( H^0\left( T, \Theta_T \right) \bigoplus H^0\left( T, \mathcal{O}_T \right) \xrightarrow{\tilde{G}_0} H^0\left( T, \Omega_T^1 \right) \bigoplus H^0\left( T, \mathcal{O}_T \right) \right) = 4-2 = 2
\end{align}
We compute $H^0\left( T, \Omega_T^1  \right) \bigoplus H^0\left( T, \mathcal{O}_T \right) \xrightarrow{\tilde{G}_1} H^0\left( T, \mathcal{O}_T \right) \bigoplus H^0\left( T,  \bigwedge^2 \Omega_T^1 \right)$. For $c=c_1dz_1+ c_2 dz_2 + c_3 dz_3\in H^0\left( T, \Omega_T^1 \right)$ and $h\in H^0\left(T, \mathcal{O}_T \right)$, the image of $(c,h)$ under $\tilde{G}_1$ is given by
\begin{align*}
&\left(\omega_0 \wedge dc + c \wedge d\omega_0 ,  - \omega_0 \wedge c + h \cdot d\omega_0 + \omega_0 \wedge dh \right) \\
&= \left( 0, - ( a_1 c_2 - c_1 a_2) dz_1 \wedge dz_2 - \left(a_2 c_3 -  c_2 a_3 \right) dz_2 \wedge dz_3 - \left(a_3c_1 - c_3 a_1 \right) dz_3 \wedge dz_1 \right)
\end{align*}
This implies that 
\begin{align}\label{com28}
\dim_\mathbb{C} \textnormal{ker}\left( H^0\left( T, \Omega_T^1  \right) \bigoplus H^0\left( T, \mathcal{O}_T \right) \xrightarrow{\tilde{G}_1} H^0\left( T, \mathcal{O}_T \right) \bigoplus H^0\left( T,  \bigwedge^2 \Omega_T^1 \right) \right) = 2
\end{align}
Then from $(\ref{com27})$ and $(\ref{com28})$ we have
\begin{align} \label{com31}
\frac{\textnormal{ker}\left( H^0\left( T, \Omega_T^1  \right) \bigoplus H^0\left( T, \mathcal{O}_T \right) \xrightarrow{\tilde{G}_1} H^0\left( T, \mathcal{O}_T \right) \bigoplus H^0\left( T,  \bigwedge^2 \Omega_T^1 \right) \right)}{\textnormal{im}\left( H^0\left( T, \Theta_T \right) \bigoplus H^0\left( T, \mathcal{O}_T \right) \xrightarrow{\tilde{G}_0} H^0\left( T, \Omega_T^1 \right) \bigoplus H^0\left( T, \mathcal{O}_T \right) \right)} = 2- 2 = 0
\end{align}

On the other hand, we compute $H^1\left( T, \Theta_T \right) \bigoplus H^1\left( T, \mathcal{O}_T \right) \xrightarrow{\tilde{G}_0} H^1\left( T, \Omega_T^1 \right) \bigoplus H^1\left( T, \mathcal{O}_T \right)$. First we note that from $(\ref{com29})$, we have $\textnormal{ker}\left( H^1\left( T, \Theta_T \right)\bigoplus H^1\left( T, \mathcal{O}_T \right) \to H^1\left( T ,  \Omega_T^1 \right) \right)= H^1\left(T, \Theta_T \right)$, and we have
\begin{align*}
&d\bar{z}_1 \wedge i_{e_1 \frac{\partial}{\partial z_1} + e_2 \frac{\partial}{\partial z_2} + e_3 \frac{\partial}{\partial z_3}}\left( a_1 dz_1 + a_2 dz_2 + a_3 dz_3\right) + d\bar{z}_2 \wedge i_{e_4 \frac{\partial}{\partial z_1} + e_5 \frac{\partial}{\partial z_2} + e_6 \frac{\partial}{\partial z_3}} \left( \omega_0 \right) + d\bar{z}_3 \wedge i_{e_7 \frac{\partial}{\partial z_1} + e_8 \frac{\partial}{\partial z_2} + e_9 \frac{\partial}{\partial z_3}} \left(\omega_0 \right)\\
&= \left(e_1a_1 + e_2 a_2 + e_3 a_3 \right)d\bar{z}_1 + \left( e_4 a_1 + e_5 a_2 + e_6 a_3 \right) d\bar{z}_2  + \left( e_7 a_1 + e_8 a_2 + a_9 a_3 \right) d\bar{z}_3
\end{align*}
Hence $\textnormal{ker}\left(H^1\left( T , \Theta_T \right) \bigoplus H^1 \left( T, \mathcal{O}_T  \right) \xrightarrow{\tilde{G}_0} H^1\left( T, \Omega_T^1 \right) \bigoplus H^1\left(  T, \mathcal{O}_T \right) \right)$ is given by
\begin{align} \label{com30}
&t_1\left(a_3 \frac{\partial}{\partial z_1} -a_1 \frac{\partial}{\partial z_3}\right) d \bar{z}_1 + t_2 \left(  a_3 \frac{\partial}{\partial z_2} -a_2 \frac{\partial}{\partial z_3} \right) d\bar{z}_1 +t_3\left(a_3 \frac{\partial}{\partial z_1} -a_1 \frac{\partial}{\partial z_3}\right) d \bar{z}_2 + t_4 \left(  a_3 \frac{\partial}{\partial z_2} -a_2 \frac{\partial}{\partial z_3} \right) d\bar{z}_2 \\
&+t_5\left( a_3\frac{\partial}{\partial z_1} - a_1 \frac{\partial}{\partial z_3}\right) d \bar{z}_3 + t_6 \left(  a_3 \frac{\partial}{\partial z_2} - a_2 \frac{\partial}{\partial z_3} \right) d\bar{z}_3, \,\,\,\,\,t_1,t_2,t_3,t_4,t_5,t_6 \in \mathbb{C} \notag
\end{align}
Then from $(\ref{com26})$ and $(\ref{com31})$ and $(\ref{com30})$ we have
\begin{align} \label{comp1}
&\mathbb{H}^1\left( T, \tilde{\mathcal{E}}_{\mathcal{O}_T}^\bullet \right) \\
&\cong \left\{ t_1\left(a_3 \frac{\partial}{\partial z_1} -a_1 \frac{\partial}{\partial z_3}\right) d \bar{z}_1 + t_2 \left(  a_3 \frac{\partial}{\partial z_2} -a_2 \frac{\partial}{\partial z_3} \right) d\bar{z}_1 +t_3\left(a_3 \frac{\partial}{\partial z_1} -a_1 \frac{\partial}{\partial z_3}\right) d \bar{z}_2 + t_4 \left(  a_3 \frac{\partial}{\partial z_2} -a_2 \frac{\partial}{\partial z_3} \right) d\bar{z}_2 \right. \notag \\
 &\,\,\,\,\,\,\,\,\,\,\, \left. +t_5\left( a_3\frac{\partial}{\partial z_1} - a_1 \frac{\partial}{\partial z_3}\right) d \bar{z}_3 + t_6 \left(  a_3 \frac{\partial}{\partial z_2} - a_2 \frac{\partial}{\partial z_3} \right) d\bar{z}_3 \,\,\, | \,\,\, t_1,t_2,t_3,t_4,t_5,t_6 \in \mathbb{C} \right\} \notag
\end{align}
We will construct an unfolded analytic family of unfoldings of $\left(T, \omega_0 \right)$ where $\omega= a_1 dz_1 + a_2 dz_2 + a_3 dz_3 , a_3 \ne 0$ such that the unfolded Kodaira-Spencer map is an isomorphism at the distinguished point, so that $\left( T , \omega_0 \right)$ is unobstructed in unfoldings. In fact, we recall the integrability conditions for an unfolding from $(\ref{uti4})-(\ref{uti8})$. If we take
{\small{\begin{align*}
\left( Y(t)=\omega_0, Q^\eta(t)=0, \left\{ \left( \varphi(t), Z(t)=0 \right) \right\} \right) \in  A^{0,0}\left( T, \Omega_T^1 \right) \bigoplus A^{0,0}\left( T, \mathcal{O}_T \right)\bigoplus A^{0,1}\left( T, \Theta_T \right)\bigoplus A^{0,1}\left( T, \mathcal{O}_T \right),\eta=1,2,3,4,5,6
\end{align*}}}
where $\varphi(t)= t_1\left(a_3 \frac{\partial}{\partial z_1} -a_1 \frac{\partial}{\partial z_3}\right) d \bar{z}_1 + t_2 \left(  a_3 \frac{\partial}{\partial z_2} -a_2 \frac{\partial}{\partial z_3} \right) d\bar{z}_1 +t_3\left(a_3 \frac{\partial}{\partial z_1} -a_1 \frac{\partial}{\partial z_3}\right) d \bar{z}_2 + t_4 \left(  a_3 \frac{\partial}{\partial z_2} -a_2 \frac{\partial}{\partial z_3} \right) d\bar{z}_2 +t_5\left( a_3\frac{\partial}{\partial z_1} - a_1 \frac{\partial}{\partial z_3}\right) d \bar{z}_3 + t_6 \left(  a_3 \frac{\partial}{\partial z_2} - a_2 \frac{\partial}{\partial z_3} \right) d\bar{z}_3$. Then we have $\varphi(0)=0$ and $Z(0)=0$ and the integrability conditions $(\ref{uti4})-(\ref{uti8})$ are satisfied, so that there exists an unfolded analytic family $\left(\mathcal{M}, \omega , B, \pi \right)$ where $\pi:\left( \mathcal{M}, \omega \right)\to B$ for a sufficiently small neighborhood $B$ of $0$ in $\mathbb{C}^6$ with coordinates $t=\left(t_1,t_2,t_3,t_4, t_5, t_6\right)$ such that $\pi^{-1}(0)=\left(X, \omega_0 \right)$ and  the unfolded Kodaira-Spencer map from $(\ref{cu5})$ given by
\begin{align*}
\varphi : T_0(B) & \to \mathbb{H}^1\left( T , \tilde{\mathcal{E}}_{\mathcal{O}_T}^\bullet \right) \\
\frac{\partial}{\partial t_1} & \mapsto \left( 0,0, \left( \left( a_3 \frac{\partial}{\partial z_1} - a_1 \frac{\partial}{\partial z_3} \right) d\bar{z}_1 , 0\right) \right), \frac{\partial}{\partial t_2}  \mapsto \left( 0,0, \left( \left( a_3 \frac{\partial}{\partial z_2} - a_2 \frac{\partial}{\partial z_3} \right) d\bar{z}_1 , 0\right) \right) \\
\frac{\partial}{\partial t_3} & \mapsto \left( 0,0, \left( \left( a_3 \frac{\partial}{\partial z_1} - a_1 \frac{\partial}{\partial z_3} \right) d\bar{z}_2 , 0\right) \right), \frac{\partial}{\partial t_4}  \mapsto \left( 0,0, \left( \left( a_3 \frac{\partial}{\partial z_2} - a_2 \frac{\partial}{\partial z_3} \right) d\bar{z}_2 , 0\right) \right) \\
\frac{\partial}{\partial t_5} & \mapsto \left( 0,0, \left( \left( a_3 \frac{\partial}{\partial z_1} - a_1 \frac{\partial}{\partial z_3} \right) d\bar{z}_3 , 0\right) \right), \frac{\partial}{\partial t_6}  \mapsto \left( 0,0, \left( \left( a_3 \frac{\partial}{\partial z_2} - a_2 \frac{\partial}{\partial z_3} \right) d\bar{z}_3 , 0\right) \right) 
\end{align*}
is an isomorphism at $t=0$ by $(\ref{comp1})$, so that $\left( T, \omega_0 \right)$ is unobstructed in unfoldings.

\appendix

\section{Description of $H^0 \left( \mathbb{P}_\mathbb{C}^2, \mathcal{E}_{\mathcal{O}_{\mathbb{P}_\mathbb{C}^2}}(m) \right),H^0\left(\mathbb{P}_\mathbb{C}^2, \Theta_{\mathbb{P}_\mathbb{C}^2}(1) \right)$ and $H^0\left( \mathbb{P}_\mathbb{C}^2, \Omega_{\mathbb{P}_\mathbb{C}^2}^1(4) \right)$} \label{p1}

Let $\left[ \xi_0 : \xi_1 : \xi_2 \right]$ be the homogeneous coordinates on $\mathbb{P}_\mathbb{C}^2$. Let $\mathcal{U}=\left\{ U_0, U_1, U_2 \right\}$ be an affine open covering of $\mathbb{P}_\mathbb{C}^2$, where $U_0=\left\{\left( X_1, X_2 \right)\in \mathbb{C}^2 | \xi_0 \ne 0, X_1= \frac{\xi_1}{\xi_0}, X_2= \frac{\xi_2}{\xi_0} \right\}$ and $U_1=\left\{ Y_0, Y_2 | \xi_1\ne 0, Y_0=\frac{ \xi_0}{\xi_1}, Y_2= \frac{\xi_2}{\xi_1} \right\}$ and $U_2=\left\{\left( Z_0, Z_1 \right) \in \mathbb{C}^2  |\xi_2\ne 0, Z_0=\frac{\xi_0}{\xi_2}, Z_1= \frac{\xi_0}{\xi_2} \right\}$. Let us consider the line bundle $\mathcal{O}_{\mathbb{P}_{\mathbb{C}}^2 }(m)$ on $\mathbb{P}_\mathbb{C}^2$ for some $m \in \mathbb{Z}$. In terms of the covering $\mathcal{U}$ of $\mathbb{P}_\mathbb{C}^2$, $\mathcal{O}_{\mathbb{P}_\mathbb{C}^2}(m)$ is described by the following transition functions
\begin{align*}
\left\{ f_{01}= X_1^m, f_{02}= X_2^m , f_{12}= Y_2^m \right\} \in C^0\left( \mathcal{U}, \mathcal{O}_{\mathbb{P}_\mathbb{C}^2}^* \right)
\end{align*}
where a local section $s_i\in \Gamma\left(U_i, \mathcal{O}_{\mathbb{P}_\mathbb{C}^2} \right)$ of $\mathcal{O}_{\mathbb{P}_\mathbb{C}^2}(m)$ is identified with $s_j\in \Gamma\left( U_j, \mathcal{O}_{\mathbb{P}_\mathbb{C}^2} \right)$ if $s_i= f_{ij}a_j$ on $U_i\cap U_j$.

We note that $\dim_\mathbb{C} H^0\left( \mathbb{P}_\mathbb{C}^2, \Theta_{\mathbb{P}_\mathbb{C}^2} \right)= 8$ and any element $T\in H^0\left(\mathbb{P}_\mathbb{C}^2, \Theta_{\mathbb{P}_\mathbb{C}^2} \right)$ on $U_0$ is of the form (see \cite{Kod05} p.223)
\begin{align}\label{ap1}
\left( b_1 + b_2 X_1 + b_3 X_2 + b_7 X_1^2 + b_8 X_1X_2 \right) \frac{\partial}{\partial X_1} + \left( b_4 + b_5 X_1 + b_6 X_2 + b_7 X_1X_2 + b_8 X_2^2 \right)\frac{\partial}{\partial X_2}
\end{align}
We note that from the Atiyah extension
\begin{align*}
0\to \mathcal{O}_{\mathbb{P}_{\mathbb{C}}^2} \to \mathcal{E}_{\mathcal{O}_{\mathbb{P}_\mathbb{C}^2}(m)} \to \Theta_{\mathbb{P}_\mathbb{C}^2}\to 0,
\end{align*}
we see that $\dim_\mathbb{C} H^0\left( \mathbb{P}_\mathbb{C}^2, \mathcal{E}_{\mathcal{O}_{\mathbb{P}_\mathbb{C}^2}(m)} \right)= 9$. We explicitly describe $ H^0\left( \mathbb{P}_\mathbb{C}^2, \mathcal{E}_{\mathcal{O}_{\mathbb{P}_\mathbb{C}^2}(m)} \right)$ on $U_0$. Let $f+ T\in H^0\left( \mathbb{P}_\mathbb{C}^2, \mathcal{E}_{\mathcal{O}_{\mathbb{P}_\mathbb{C}^2}}(m) \right)$ on $U_0$, where $f=\sum_{i,j \geq 0} c_{ij}X_1^i X_2^j, c_{ij}\in \mathbb{C}$ is a holomorphic function on $U_0$, and $T\in H^0\left( \mathbb{P}_\mathbb{C}^2, \Theta_{\mathbb{P}_\mathbb{C}^2} \right)$ is given by $\textnormal{(\ref{ap1})}$. By the definition of $\mathcal{E}_{\mathcal{O}_{\mathbb{P}_\mathbb{C}^2}}(m)$, $f+ T$ on $U_0$ is identified with $f+ \left[ T, \log X_1^m \right] + T $ on $U_1$, so that from  $X_1= \frac{1}{Y_0}$ and $X_2=\frac{Y_2}{Y_0}$ on $U_0\cap U_1$, we see that
\begin{align*}
&\sum_{i,j \geq 0} c_{ij} X_1^i X_2^j +\left[ \left( b_1 + b_2 X_1 + b_3 X_2 + b_7 X_1^2 + b_8 X_1X_2 \right) \frac{\partial}{\partial X_1} + \left( b_4 + b_5 X_1 + b_6 X_2 + b_7 X_1X_2 + b_8 X_2^2 \right)\frac{\partial}{\partial X_2}, \log X_1^m \right]\\
&= \sum_{i,j\geq 0} c_{ij} Y_0^{-i-j} Y_2^j + m\left( b_1Y_0 + b_2 + b_3 Y_2 + b_7 \frac{1}{Y_0} + b_8 \frac{Y_2}{Y_0} \right)
\end{align*}
should be holomorphic on $U_1$, i.e. $c_{10}=- m b_7$ and $c_{01}=-mb_8$, so that $f= c_{00}- mb_7 X_1 - mb_8 X_2$. Hence we have, on $U_0$
\begin{align}\label{p2}
&H^0\left( \mathbb{P}_\mathbb{C}^2, \mathcal{E}_{\mathcal{O}_{\mathbb{P}_\mathbb{C}^2}(m)} \right)\\
&= \left\{ f, T_1\frac{\partial}{\partial X_1} + T_2 \frac{\partial}{\partial X_2} \right\} \notag \\
&=\left\{ a- mb_7X_1- mb_8 X_2, \left( b_1 + b_2 X_1 + b_3 X_2 + b_7 X_1^2 + b_8 X_1X_2 \right) \frac{\partial}{\partial X_1} + \left( b_4 + b_5 X_1 + b_6 X_2 + b_7 X_1X_2 + b_8 X_2^2 \right)\frac{\partial}{\partial X_2}   \right\} \notag
\end{align}
where $a,b_1,...,b_8\in \mathbb{C}$.

\subsection{Description of $H^0\left( \mathbb{P}_\mathbb{C}^2, \Theta_{\mathbb{P}_\mathbb{C}^2} (1) \right)$}\

We note that from the Euler sequence $0\to \mathcal{O}_{\mathbb{P}_\mathbb{C}^2}\to \bigoplus^3 \mathcal{O}_{\mathbb{P}_\mathbb{C}^2}(1) \to \Theta_{\mathbb{P}_\mathbb{C}^2}\to 0$, we have $0 \to \mathcal{O}_{\mathbb{P}_\mathbb{C}^2}(1)\to \bigoplus^3 \mathcal{O}_{\mathbb{P}_\mathbb{C}^2}(2) \to \Theta_{\mathbb{P}_\mathbb{C}^2}(1) \to 0$, we see that $\dim_\mathbb{C} H^0\left( \mathbb{P}_\mathbb{C}^2, \Theta_{\mathbb{P}_\mathbb{C}^2}(1) \right)= 15
$. We explicitly describe $H^0\left( \mathbb{P}_\mathbb{C}^2, \Theta_{\mathbb{P}_\mathbb{C}^2}(1) \right)$ on $U_0$. Let
\begin{align*}
T_0= v_0^1\left( X_1, X_2\right) \frac{\partial}{\partial X_1} + v_0^2\left( X_1, X_2 \right) \frac{\partial}{\partial X_2} \in H^0\left(  \mathbb{P}_\mathbb{C}^2, \Theta_{\mathbb{P}_\mathbb{C}^2}(1)\right) \,\,\,\textnormal{on}\,\,\, U_0
\end{align*}
be identified with 
\begin{align*}
T_1= v_1^0\left( Y_0, Y_2 \right) \frac{\partial}{\partial Y_0} + v_1^2 \left( Y_0 , Y_2 \right) \frac{\partial}{\partial Y_2} \in H^0\left( \mathbb{P}_\mathbb{C}^2 , \Theta_{\mathbb{P}_\mathbb{C}^3} (1) \right)\,\,\, \textnormal{on}\,\,\, U_1
\end{align*}
so that we have 
\begin{align*}
T_0= X_1 T_1 = \frac{1}{Y_0} T_1
\end{align*}
Since $Y_0= \frac{1}{X_1}, Y_2 = \frac{X_2}{X_1}$, we have $\frac{\partial}{\partial X_1}= - \frac{1}{X_1^2} \frac{\partial}{\partial Y_0} - \frac{X_2}{X_1^2} \frac{\partial}{\partial Y_2}, \frac{\partial}{\partial X_2} = \frac{1}{X_1} \frac{\partial}{\partial Y_2}$. Then we have
\begin{align*}
T_0= v_0^1\left( \frac{1}{Y_0}, \frac{Y_2}{Y_0} \right)\left( - Y_0^2 \frac{\partial}{\partial Y_0} - Y_0 Y_2 \frac{\partial}{\partial Y_2} \right) + Y_0 v_0^2\left( \frac{1}{Y_0}, \frac{Y_2}{Y_0} \right)\frac{\partial}{\partial Y_2}= \frac{1}{Y_0}\left(  v_1^0\left( Y_0, Y_2 \right) \frac{\partial}{\partial Y_0} + v_1^2 \left( Y_0 , Y_2 \right) \frac{\partial}{\partial Y_2} \right)
\end{align*}
Then by letting $v_0^1\left( X_1, X_2 \right) = \sum_{h,k\geq 0} a_{hk} X_1^h X_2^k, a_{hk}\in \mathbb{C}$ and $v_0^2\left( X_1, X_2 \right) = \sum_{h,k\geq 0} b_{hk} X_1^h X_2^k, b_{hk} \in \mathbb{C}$, we get
\begin{align*}
v_1^0 \left(Y_0, Y_2 \right) & = - Y_0^3 v_0^1\left( \frac{1}{Y_0}, \frac{Y_2}{Y_0} \right) = \sum_{h,k,r\geq 0} \frac{- a_{hk} Y_2^k}{Y_0^{k+r-3}} \\
v_1^2\left( Y_0, Y_2 \right) & =  - Y_0^2 Y_2 v_0^1\left( \frac{1}{Y_0}, \frac{Y_2}{Y_0} \right) + Y_0^2 v_0^2 \left( \frac{1}{Y_0}, \frac{Y_2}{Y_0} \right) = \sum_{h,k\geq 0} \frac{- a_{hk} Y_2^{k+1} + b_{hk} Y_2^k }{Y_0^{h+k-2}}
\end{align*}
which should be holomorphic on $U_1$. In order for $v_1^2\left( Y_0, Y_2\right)$ to be holomorphic on $U_1$, we only need to consider the case
\begin{align*}
(h,k)= (0,0),(0,1),(1,0),(0,2),(1,1),(2,0),(0,3),(1,2),(2,1),(3,0)
\end{align*}
The case $(0,0),(0,1),(1,0),(0,2),(1,1),(2,0)$ for $v_1^2\left( Y_0, Y_2\right)$ clearly defines a holomorphic function on $U_1$. Let us consider the case $(0,3),(1,2),(2,1),(3,0)$ which is
\begin{align*}
&\frac{\left(- a_{03} Y_2^4+ b_{03}Y_2^3 \right) + \left(- a_{12} Y_2^3 + b_{12} Y_2^2 \right) +\left( - a_{21} Y_2^2 + b_{21}Y_2 \right) + \left( - a_{30} Y_2 + b_{30} \right)}{Y_0}\\
& = \frac{b_{30}+ \left( b_{21} - a_{30}\right) Y_2 +\left( b_{12}- a_{21} \right) Y_2^2 + \left( b_{03}- a_{12}\right) Y_2^3 - a_{03}Y_2^4 }{Y_0}
\end{align*}
In order for $v_1^2\left( Y_0, Y_2 \right)$ to be holomorphic on $U_1$, we should have $b_{30}=0, b_{21}= a_{30}, b_{12}= a_{21}, b_{03}= a_{12}, a_{03}=0$. On the other hand, in order for $v_1^0\left( Y_0, Y_2\right)$ to be holomorphic on $U_1$, we should have $a_{hk}=0$ except for
\begin{align*}
(h,k)= (0,0),(0,1),(1,0),(0,2),(1,1),(2,0),(0,3),(1,2),(2,1),(3,0)
\end{align*}
Then we set
\begin{align*}
&t_1=a_{00},\,\,\,\,t_2=a_{01},\,\,\,\,\,t_3=a_{10},\,\,\,\,\,t_4=a_{02},\,\,\,\,\,t_5=a_{11},\,\,\,\,\,t_6=a_{20}\\
&t_7=a_{12}= b_{03},\,\,\,\,\,t_8=a_{21}= b_{12},\,\,\,\,\,t_9=a_{30}= b_{21},\,\,\,\,\,t_{10}=b_{00},\,\,\,\,\,t_{11}=b_{01},\,\,\,\,\,t_{12}=b_{10}\\
&t_{13}=b_{02},\,\,\,\,\,t_{14}=b_{11},\,\,\,\,\,t_{15}=b_{20}
\end{align*}
Then $H^0\left( \mathbb{P}_\mathbb{C}^2, \Theta_{\mathbb{P}_\mathbb{C}^2}(1) \right)$ is given by
\begin{align} \label{p3}
&\left( t_1+ t_2 X_2+ t_3 X_1+ t_4 X_2^2+ t_5 X_1X_2+ t_6 X_1^2+ t_7 X_1X_2^2+ t_8 X_1^2 X_2+ t_9 X_1^3 \right)\frac{\partial}{\partial X_1}\\
&+ \left( t_{10} + t_{11}X_2+ t_{12}X_1+ t_{13} X_2^2 + t_{14}X_1X_2+ t_{15}X_1^2+ t_7 X_2^3+ t_8 X_1X_2^2+ t_9 X_1^2 X_2 \right) \frac{\partial}{\partial X_2} \notag
\end{align}

\subsection{Description of $H^0\left( \mathbb{P}_\mathbb{C}^2, \Omega_{\mathbb{P}_\mathbb{C}^2}^1(4) \right)$}\

We have $\dim_\mathbb{C} H^0\left( \mathbb{P}_\mathbb{C}^2, \Omega_{\mathbb{P}_\mathbb{C}^2}^1(4) \right)= 15$ by Bott formula (see \cite{OSS88} p.4). We explicitly describe $H^0\left( \mathbb{P}_\mathbb{C}^2 , \Omega_{\mathbb{P}_\mathbb{C}^2}^1(4) \right)$ on $U_0$. Let
\begin{align*}
v_0= v_0^1\left( X_1, X_2 \right) dX_1 + v_0^2\left(X_1, X_2 \right) dX_2 \in H^0\left( \mathbb{P}_\mathbb{C}^2, \Omega_{\mathbb{P}_\mathbb{C}^2}^1(4) \right)\,\,\,\textnormal{on}\,\,\, U_0
\end{align*}
be identified with
\begin{align*}
v_1= v_1^0\left( Y_0, Y_2 \right) dY_0 + v_1^2 \left( Y_0, Y_2 \right) dY_2 \in H^0\left( \mathbb{P}_\mathbb{C}^2, \Omega_{\mathbb{P}_\mathbb{C}^2}^1(4) \right)\,\,\,\textnormal{on}\,\,\, U_1
\end{align*}
so that we have 
\begin{align*}
v_0= X_1^4 v_1 = \frac{1}{Y_0^4} v_1
\end{align*}
Since $X_1= \frac{1}{Y_0}, X_2=\frac{Y_2}{Y_0}$, we have $dX_1= - \frac{1}{Y_0^2} dY_0$ and $dX_2= - \frac{Y_2}{Y_0^2} dY_0+ \frac{1}{Y_0} dY_2$, and so
{\small{\begin{align*}
v_0&= - \frac{1}{Y_0^2} v_0^1\left( \frac{1}{Y_0}, \frac{Y_2}{Y_0} \right) dY_0 + v_0^2\left( \frac{1}{Y_0}, \frac{Y_2}{Y_0} \right)\left( - \frac{Y_2}{Y_0^2} dY_0+ \frac{1}{Y_0} dY_2 \right) \\
&= \left( - \frac{1}{Y_0^2}v_0^1\left(  \frac{1}{Y_0}, \frac{Y_2}{Y_0} \right) - \frac{Y_2}{Y_0^2} v_0^2\left( \frac{1}{Y_0}, \frac{Y_2}{Y_0} \right) \right) dY_0 + \frac{1}{Y_0} v_0^2\left( \frac{1}{Y_0}, \frac{Y_2}{Y_0} \right) dY_2 = \frac{1}{Y_0^4} v_1^0\left( Y_0, Y_2\right) dY_0 + \frac{1}{Y_0^4} v_1^2\left( Y_0, Y_2 \right) dY_2
\end{align*}}}
Hence we have
\begin{align*}
v_1^0\left( Y_0 , Y_2 \right)& = - Y_0^2 v_0^1\left( \frac{1}{Y_0}, \frac{Y_2}{Y_0} \right)- Y_0^2 Y_2 v_0^2\left( \frac{1}{Y_0}, \frac{Y_2}{Y_0} \right)\\
v_1^2\left( Y_0, Y_2 \right)&= Y_0^3 v_0^2\left(\frac{1}{Y_0}, \frac{Y_2}{Y_0} \right)
\end{align*}
By letting $v_0^1\left(X_1, X_2 \right)= \sum_{h,k\geq 0} a_{hk}X_1^h X_2^k, a_{hk}\in \mathbb{C}$ and $v_0^2\left( X_1, X_2 \right)= \sum_{h,k\geq 0} b_{hk} X_1^h X_2^k, b_{hk}\in \mathbb{C}$, we get
\begin{align*}
v_1^0\left( Y_0 , Y_2 \right)&= \sum_{h,k\geq 0} \frac{- a_{hk} Y_2^k - b_{hk} Y_2^{k+1}}{Y_0^{h+k -2}}\\
v_1^2\left( Y_0, Y_2 \right)&=\sum_{h,k\geq 0} \frac{b_{hk}Y_2^k}{ Y_0^{h+k-3}}
\end{align*}
In order for $v_1^0\left( Y_0, Y_2\right)$ to be holomorphic on $U_1$, we only need to consider the case
\begin{align*}
(h,k)= (0,0),(0,1),(1,0),(0,2),(1,1),(2,0),(0,3),(1,2),(2,1),(3,0)
\end{align*}
The case $(0,0),(0,1),(1,0),(0,2),(1,1),(2,0)$ for $v_1^0\left( Y_0, Y_2\right)$ clearly defines a holomorphic function on $U_1$. Let us consider the case $(0,3),(1,2),(2,1),(3,0)$ which is
\begin{align*}
&\frac{\left(a_{03}Y_2^3+b_{03}Y_2^4 \right)+ \left(a_{12}Y_2^2+b_{12}Y_2^3\right)+ \left(a_{21}Y_2+b_{21}Y_2^2 \right)+ \left(a_{30}+b_{30}Y_2\right)}{Y_0}\\
&=\frac{a_{30}+ \left(a_{21}+b_{30} \right)Y_2+ \left( a_{12}+b_{21}\right)Y_2^2+(a_{03}+b_{12})Y_2^3+b_{03}Y_2^4}{Y_0}
\end{align*}
In order for $v_1^0\left( Y_0, Y_2 \right)$ to be holomorphic on $U_1$, we should have $a_{30}=0, a_{21}= - b_{30}, a_{12}= - b_{21}, a_{03}= - b_{12}, b_{03}=0$. On the other hand, in order for $v_1^2\left( Y_0, Y_2 \right)$ to be holomorphic on $U_1$, we should have $b_{hk}=0$ except for
\begin{align*}
(h,k)= (0,0),(0,1),(1,0),(0,2),(1,1),(2,0),(0,3),(1,2),(2,1),(3,0)
\end{align*}
Then we set
\begin{align*}
&t_1=b_{00},\,\,\,\,t_2=b_{01},\,\,\,\,\,t_3=b_{10},\,\,\,\,\,t_4=b_{02},\,\,\,\,\,t_5=b_{11},\,\,\,\,\,t_6=b_{20}\\
&t_7=b_{12}= - a_{03},\,\,\,\,\,t_8=b_{21}= - a_{12},\,\,\,\,\,t_9=b_{30}= - a_{21},\,\,\,\,\,t_{10}=a_{00},\,\,\,\,\,t_{11}=a_{01},\,\,\,\,\,t_{12}=a_{10}\\
&t_{13}=a_{02},\,\,\,\,\,t_{14}=a_{11},\,\,\,\,\,t_{15}=a_{20}
\end{align*}
Then $H^0\left( \mathbb{P}_\mathbb{C}^2 , \Omega_{\mathbb{P}_\mathbb{C}^2}^1 ( 4) \right)$ is given by
\begin{align} \label{pa1}
b=&\left(t_{10}+t_{11}X_2+t_{12}X_1+t_{13}X_2^2+t_{14}X_1X_2+t_{15}X_1^2-t_9X_1^2X_2-t_8X_1X_2^2-t_7X_2^3 \right)dX_1\\
&+\left( t_1+t_2X_2+t_3X_1+t_4X_2^2+t_5X_1X_2+t_6X_1^2+t_7X_1X_2^2+t_8X_1^2X_2+t_9X_1^3\right)dX_2 \notag
\end{align}

\section{Description of $H^0\left(\mathbb{P}_\mathbb{C}^3, \Theta_{\mathbb{P}_\mathbb{C}^3} \right),H^0\left( \mathbb{P}_\mathbb{C}^3, \mathcal{E}_{\mathcal{O}_{\mathbb{P}_\mathbb{C}^3}(m)} \right), H^0\left( \mathbb{P}_\mathbb{C}^3, \Theta_{\mathbb{P}_\mathbb{C}^3}(1) \right),H^0\left( \mathbb{P}_\mathbb{C}^3, \Theta_{\mathbb{P}_\mathbb{C}^3}(-1) \right)$}\label{cp7}

Let $\left[ \xi_0: \xi_1: \xi_2 : \xi_3 \right]$ be the homogenous coordinate on $\mathbb{P}_\mathbb{C}^3$. Let $\mathcal{U}=\left\{ U_0, U_1, U_2, U_3 \right\}$ be an affine open covering of $\mathbb{P}_\mathbb{C}^3$ where 
\begin{align*}
U_0 &=\left\{\left( X_1 , X_2, X_3\right) |\xi_0\ne 0, X_1=\frac{\xi_1}{\xi_0}, X_2=\frac{\xi_2}{\xi_0}, X_3 = \frac{\xi_3}{\xi_0} \right\}\\
U_1&=\left\{\left( Y_0 , Y_2, Y_3\right) |\xi_1\ne 0, Y_0=\frac{\xi_0}{\xi_1}, Y_2=\frac{\xi_2}{\xi_1}, Y_3 = \frac{\xi_3}{\xi_1} \right\}\\
U_2&=\left\{\left( Z_0 , Z_1, Z_3\right) |\xi_2\ne 0, Z_0=\frac{\xi_0}{\xi_2}, Z_1=\frac{\xi_1}{\xi_2}, Z_3 = \frac{\xi_3}{\xi_2} \right\}\\
U_3&=\left\{\left( W_0 , W_1, W_2\right) |\xi_3\ne 0, W_0=\frac{\xi_0}{\xi_3}, W_1=\frac{\xi_1}{\xi_3}, W_2 = \frac{\xi_2}{\xi_3} \right\}
\end{align*}

\subsection{Description of $H^0\left( \mathbb{P}_\mathbb{C}^3 , \Theta_{\mathbb{P}_\mathbb{C}^3} \right)$ } \

We note that from the Euler sequence
\begin{align*}
0\to \mathcal{O}_{\mathbb{P}_\mathbb{C}^3} \to \bigoplus^4 \mathcal{O}_{\mathbb{P}_\mathbb{C}^3}(1) \to \Theta_{\mathbb{P}_\mathbb{C}^3} \to 0
\end{align*}
we see that $\dim_\mathbb{C} H^0\left( \mathbb{P}_\mathbb{C}^3, \Theta_{\mathbb{P}_\mathbb{C}^3} \right)=15$. We explicitly describe $H^0\left( \mathbb{P}_\mathbb{C}^3 , \Theta_{\mathbb{P}_\mathbb{C}^3}  \right)$ on $U_0$. Let 
\begin{align*}
T= v_0^1\left(X_1,X_2,X_3 \right)\frac{\partial}{\partial X_1} + v_0^2\left(X_1,X_2,X_3 \right)\frac{\partial}{\partial X_2} + v_0^3 \left( X_1,X_2, X_3 \right) \frac{\partial}{\partial X_3}\in H^0\left(\mathbb{P}_\mathbb{C}^3, \Theta_{\mathbb{P}_\mathbb{C}^3} \right)\,\,\,\,\,\textnormal{on}\,\,\, U_0
\end{align*}
where
\begin{align*}
v_0^1\left( X_1, X_2, X_3 \right)= \sum_{h,k,r\geq 0} a_{hkr} X_1^h X_2^k X_3^r,\,\,\,\,\,\,\, a_{hkr}\in \mathbb{C}\\
v_0^2\left( X_1, X_2, X_3 \right)= \sum_{h,k,r\geq 0} b_{hkr} X_1^h X_2^k X_3^r,\,\,\,\,\,\,\, b_{hkr}\in \mathbb{C}\\
v_0^3\left( X_1, X_2, X_3 \right)= \sum_{h,k,r\geq 0} c_{hkr} X_1^h X_2^k X_3^r,\,\,\,\,\,\,\, c_{hkr}\in \mathbb{C}
\end{align*}
We note that $Y_0= \frac{1}{X_1}, Y_2=\frac{X_2}{X_1}, Y_3= \frac{X_3}{X_1}$. Then we have $\frac{\partial}{\partial X_1}= - \frac{1}{X_1^2}\frac{\partial}{\partial Y_0}- \frac{X_2}{X_1^2}\frac{\partial}{\partial Y_2}- \frac{X_3}{X_1^2}\frac{\partial}{\partial Y_3}, \frac{\partial}{\partial X_2}= \frac{1}{X_1}\frac{\partial}{\partial Y_2}, \frac{\partial}{\partial X_3}= \frac{1}{X_1}\frac{\partial}{\partial Y_3}$. Then $T$ on $U_0$ is transformed on $U_1$ into 
\begin{align*}
&v_1^0\left(Y_0,Y_2, Y_3 \right)\frac{\partial}{\partial Y_0} + v_1^2\left(Y_0,Y_2, Y_3 \right)\frac{\partial}{\partial Y_2}+ v_1^3\left(Y_0,Y_2, Y_3 \right)\frac{\partial}{\partial Y_3} \\
&= v_0^1\left(\frac{1}{Y_0}, \frac{Y_2}{Y_0}, \frac{Y_3}{Y_0} \right) \left(- Y_0^2 \frac{\partial}{\partial Y_0} - Y_0 Y_2 \frac{\partial}{\partial Y_2} - Y_0 Y_3 \frac{\partial}{\partial Y_3} \right)  + v_0^2\left(\frac{1}{Y_0}, \frac{Y_2}{Y_0}, \frac{Y_3}{Y_0} \right) Y_0\frac{\partial}{\partial Y_2} + v_0^3\left( \frac{1}{Y_0}, \frac{Y_2}{Y_0}, \frac{Y_3}{Y_0} \right) Y_0 \frac{\partial}{\partial Y_3}\\
\end{align*}
Then we have
\begin{align*}
v_1^0\left(Y_0, Y_2, Y_3 \right) &= - Y_0^2 v_0^1\left(\frac{1}{Y_0}, \frac{Y_2}{Y_0}, \frac{Y_3}{Y_0} \right) =   \sum_{h,k,r \geq 0} \frac{- a_{hkr} Y_2^k Y_3^r}{Y_0^{h+k+r-2}}\\
v_1^2\left( Y_0, Y_2, Y_3 \right)& = - Y_0Y_2v_0^1\left(\frac{1}{Y_0}, \frac{Y_2}{Y_0}, \frac{Y_3}{Y_0} \right) + Y_0 v_0^2\left( \frac{1}{Y_0}, \frac{Y_2}{Y_0}, \frac{Y_3}{Y_0} \right) = \sum_{h,k,r \geq 0} \frac{- a_{hkr} Y_2^{k+1} Y_3 ^r + b_{hkr} Y_2^k Y_3^r}{ Y_0^{h+k+r -1}}\\
v_1^3\left( Y_0, Y_2, Y_3 \right)&= - Y_0 Y_3 v_0^1\left(\frac{1}{Y_0}, \frac{Y_2}{Y_0}, \frac{Y_3}{Y_0} \right) + Y_0 v_0^3\left(\frac{1}{Y_0}, \frac{Y_2}{Y_0}, \frac{Y_3}{Y_0} \right) = \sum_{h,k,r \geq 0} \frac{ - a_{hkr} Y_2^k Y_3^{r+1}+ c_{hkr} Y_2^k Y_3^r}{Y_0^{h+k+r-1}}
\end{align*}
which should be holomorphic on $U_1$. We only need to consider the case
\begin{align*}
(h,k,r)=&\left( 0,0,0\right)\\
       &(0,0,1),(0,1,0),(1,0,0)\\
       &(0,0,2),(0,1,1),(0,2,0), (1,0,1),(1,1,0),(2,0,0)
\end{align*}
For $v_1^0\left( Y_0, Y_2, Y_3 \right)$, all the case $(h,k,r)$ in the above defines a holomorphic function $U_1$, and for $v_1^2\left( Y_0, Y_2, Y_3\right)$, the case $(0,0,0)$, $(0,0,1)$,$(0,1,0)$,$(1,0,0)$ clearly defines a holomorphic function on $U_1$. Let us consider the remaining case $(0,0,2),(0,1,1),(0,2,0), (1,0,1),(1,1,0),(2,0,0)$ for $v_1^2\left(Y_0, Y_2, Y_3 \right)$ which is
\begin{align*}
&  \frac{\left(- a_{002}Y_2Y_3^2 + b_{002} Y_3^2\right) +\left(- a_{011} Y_2^2Y_3 + b_{011} Y_2Y_3 \right) +\left( - a_{020} Y_2^3 + b_{020}Y_2^2    \right) + \left(- a_{101} Y_2Y_3 + b_{101} Y_3 \right)}{Y_0} \\
&+ \frac{\left(- a_{110} Y_2^2+ b_{110} Y_2 \right) + \left(- a_{200} Y_2 + b_{200}  \right)}{Y_0}\\
&= \frac{ b_{200}+ \left(b_{110}- a_{200} \right)Y_2 +  b_{101} Y_3 + \left( b_{020}- a_{110} \right)Y_2^2 + \left( b_{011}- a_{101} \right) Y_2Y_3 + b_{002} Y_3^2 - a_{002} Y_2Y_3^2 - a_{011} Y_2^2 Y_3 - a_{020} Y_2^3}{Y_0}
\end{align*}
In order for $v_1^2\left( Y_0, Y_2 , Y_3\right)$ to be holomorphic on $U_1$, we should have $b_{200}=0, b_{110}= a_{200}, b_{101}=0, b_{020}= a_{110}, b_{011}= a_{101}, b_{002}=0, a_{002}=0, a_{011}=0, a_{020}=0$. On the other hand, $v_1^3\left( Y_0, Y_2, Y_3\right)$, the case $(0,0,0)$, $(0,0,1)$,$(0,1,0)$,$(1,0,0)$ clearly defines a holomorphic function on $U_1$. Let us consider the remaining case $(0,0,2),(0,1,1),(0,2,0), (1,0,1),(1,1,0),(2,0,0)$ for $v_1^3\left(Y_0, Y_2, Y_3 \right)$ which is
\begin{align*}
&\frac{\left(- a_{002} Y_3^3 + c_{002}Y_3^2 \right)  + \left( - a_{011} Y_2 Y_3^2+ c_{011} Y_2Y_3\right) + \left( - a_{020} Y_2^2 Y_3+ c_{020} Y_2^2 \right)  + \left( - a_{101} Y_3^2+ c_{101} Y_3 \right)}{Y_0} \\
&+ \frac{ \left( - a_{110} Y_2Y_3+ c_{110} Y_2 \right) + \left( - a_{200} Y_3 + c_{200} \right) }{Y_0}\\
&= c_{200}+ c_{110} Y_2 + \left( c_{101} - a_{200} \right) Y_3 + \left( c_{011} - a_{110} \right) Y_2Y_3 + c_{020} Y_2^2  + \left( c_{002} - a_{101} \right) Y_3^2 - a_{011} Y_2 Y_3^2 - a_{020} Y_2^2 Y_3 - a_{002} Y_3^3
\end{align*}
In order for $v_1^3\left( Y_0, Y_2 , Y_3\right)$ to be holomorphic on $U_1$, we should have $c_{200}=0, c_{110}=0, c_{101}= a_{200}, c_{011}= a_{110}, c_{020}=0, c_{002}= a_{101}, a_{011}=0, a_{020}=0, a_{002}=0$. Then all in all, we have
\begin{align*}
&a_{000}= a_1, a_{001} = a_2, a_{010}= a_3, a_{100}= a_4, a_{101}= a_5, a_{110}= a_6 , a_{200}=a_7\\
&b_{000}= a_8, b_{001}= a_9, b_{010}=a_{10}, b_{100}= a_{11}, b_{002}=0, b_{011}=a_5, b_{020}= a_6, b_{101}=0, b_{110}= a_7, b_{200}= 0\\
&c_{000}= a_{12}, c_{001}= a_{13}, c_{010}= a_{14}, c_{100}= a_{15}, c_{002}= a_5, c_{011}= a_6, c_{020}=0, c_{101}= a_7, c_{110}=0 , c_{200}= 0
\end{align*}

Hence we have
\begin{align} \label{ap2}
H^0\left( \mathbb{P}_\mathbb{C}^3,  \Theta_{\mathbb{P}_\mathbb{C}^3}  \right) = &\left\{  \left( a_1 + a_2 X_3+ a_3 X_2 + a_4 X_1 + a_5 X_1X_3 + a_6 X_1 X_2 + a_7 X_1^2 \right) \frac{\partial}{\partial X_1}\right. \\
 &\left.+\left( a_8 + a_9 X_3 + a_{10}X_2 + a_{11}X_1 + a_5 X_2X_3 + a_6 X_2^2 + a_7 X_1X_2\right) \frac{\partial}{\partial X_2} \right. \notag \\ 
 &\left.  +\left( a_{12}+ a_{13} X_3 + a_{14} X_2 + a_{15}X_1 + a_5 X_3^2 + a_6 X_2 X_3  + a_7 X_1X_3 \right)  \frac{\partial}{\partial X_3}, a_1,..., a_{15}\in \mathbb{C} \right\} \notag
\end{align}

\subsection{Description of $H^0\left( \mathbb{P}_\mathbb{C}^3, \Theta_{\mathbb{P}_\mathbb{C}^3}(1) \right)$}\

From the Euler sequence, we have $0\to \mathcal{O}_{\mathbb{P}_\mathbb{C}^3}(1)\to \bigoplus^4 \mathcal{O}_{\mathbb{P}_\mathbb{C}^3}(2) \to \Theta_{\mathbb{P}_\mathbb{C}^3}(1) \to 0$, so that we see that $\dim_\mathbb{C} H^0\left( \mathbb{P}_\mathbb{C}^3, \Theta_{\mathbb{P}_\mathbb{C}^3}(1) \right)= 36$. We explicitly describe $H^0\left( \mathbb{P}_\mathbb{C}^3, \Theta_{\mathbb{P}_\mathbb{C}^3}(1) \right)$ on $U_0$. Let
\begin{align*}
T_0=v_0^1\left( X_1,X_2, X_3 \right) \frac{\partial}{\partial X_1}+ v_0^2 \left( X_1, X_2, X_3 \right) \frac{\partial}{\partial X_2} + v_0^3 \left( X_1, X_2, X_3 \right) \frac{\partial}{\partial X_3}\in H^0\left( \mathbb{P}_\mathbb{C}^3, \Theta_{\mathbb{P}_\mathbb{C}^3}(1) \right)\,\,\,\,\,\textnormal{on} \,\,\, U_0
\end{align*}
be identified with
\begin{align*}
T_1= v_1^0\left( Y_0, Y_2, Y_3 \right) \frac{\partial}{\partial Y_0} + v_1^2\left( Y_0, Y_2, Y_3\right)\frac{\partial}{\partial Y_2} + v_1^3\left(  Y_0, Y_2, Y_3\right) \frac{\partial}{\partial Y_3} \in H^0\left(\mathbb{P}_\mathbb{C}^3, \Theta_{\mathbb{P}_\mathbb{C}^3}(1) \right)\,\,\,\,\,\textnormal{on}\,\,\, U_1
\end{align*}
so that we have
\begin{align*}
v_0= X_1 v_1 = \frac{1}{Y_0} v_1
\end{align*}
Then we have
\begin{align*}
v_0&= v_0^1\left(\frac{1}{Y_0}, \frac{Y_2}{Y_0}, \frac{Y_3}{Y_0} \right) \left(- Y_0^2 \frac{\partial}{\partial Y_0} - Y_0 Y_2 \frac{\partial}{\partial Y_2} - Y_0 Y_3 \frac{\partial}{\partial Y_3} \right)  + v_0^2\left(\frac{1}{Y_0}, \frac{Y_2}{Y_0}, \frac{Y_3}{Y_0} \right) Y_0\frac{\partial}{\partial Y_2} + v_0^3\left( \frac{1}{Y_0}, \frac{Y_2}{Y_0}, \frac{Y_3}{Y_0} \right) Y_0 \frac{\partial}{\partial Y_3}\\
&= \frac{1}{Y_0}\left(v_1^0\left(Y_0,Y_2, Y_3 \right)\frac{\partial}{\partial Y_0} + v_1^2\left(Y_0,Y_2, Y_3 \right)\frac{\partial}{\partial Y_2}+ v_1^3\left(Y_0,Y_2, Y_3 \right)\frac{\partial}{\partial Y_3} \right)
\end{align*}
Then by letting $v_0^1\left( X_1, X_2, X_3 \right)= \sum_{h,k,r \geq 0} a_{hkr}X_1^h X_2^k X_3^r, v_0^2\left( X_1, X_2, X_3 \right)= \sum_{h,kr\geq 0} b_{hkr} X_1^h X_2^k X_3^r$, $v_0^3\left(X_1,X_2,X_3 \right)=\sum_{h,k,r\geq 0} c_{hkr}X_1^h X_2^k X_3^r, a_{hkr}, b_{hkr}, c_{hkr}\in \mathbb{C}$, we have
\begin{align*}
v_1^0\left(Y_0, Y_2, Y_3 \right) &= - Y_0^3 v_0^1\left(\frac{1}{Y_0}, \frac{Y_2}{Y_0}, \frac{Y_3}{Y_0} \right) =   \sum_{h,k,r \geq 0} \frac{- a_{hkr} Y_2^k Y_3^r}{Y_0^{h+k+r-3}}\\
v_1^2\left( Y_0, Y_2, Y_3 \right)& = - Y_0^2Y_2v_0^1\left(\frac{1}{Y_0}, \frac{Y_2}{Y_0}, \frac{Y_3}{Y_0} \right) + Y_0^2 v_0^2\left( \frac{1}{Y_0}, \frac{Y_2}{Y_0}, \frac{Y_3}{Y_0} \right) = \sum_{h,k,r \geq 0} \frac{- a_{hkr} Y_2^{k+1} Y_3 ^r + b_{hkr} Y_2^k Y_3^r}{ Y_0^{h+k+r -2}}\\
v_1^3\left( Y_0, Y_2, Y_3 \right)&= - Y_0^2 Y_3 v_0^1\left(\frac{1}{Y_0}, \frac{Y_2}{Y_0}, \frac{Y_3}{Y_0} \right) + Y_0^2 v_0^3\left(\frac{1}{Y_0}, \frac{Y_2}{Y_0}, \frac{Y_3}{Y_0} \right) = \sum_{h,k,r \geq 0} \frac{ - a_{hkr} Y_2^k Y_3^{r+1}+ c_{hkr} Y_2^k Y_3^r}{Y_0^{h+k+r-2}}
\end{align*}
which should be holomorphic on $U_1$. We only need to consider the case
\begin{align*}
(h,k,r)=&\left( 0,0,0\right)\\
       &(0,0,1),(0,1,0),(1,0,0)\\
       &(0,0,2),(0,1,1),(0,2,0), (1,0,1),(1,1,0),(2,0,0)\\
&(0,0,3),(0,1,2),(0,2,1),(0,3,0),(1,0,2),(1,1,1),(1,2,0),(2,0,1),(2,1,0),(3,0,0)
\end{align*}
For $v_1^0\left( Y_0, Y_2, Y_3 \right)$, all the case $(h,k,r)$ in the above defines a holomorphic function on $U_1$, and for $v_1^2\left( Y_0 , Y_2, Y_3 \right)$, the case $(0,0,0),(0,0,1),(0,1,0),(1,0,0), (0,0,2),(0,1,1),(0,2,0), (1,0,1),(1,1,0),(2,0,0)$ defines a holomorphic function on $U_1$. Let us consider the remaining case which is
{\small{\begin{align*}
&\frac{\left( - a_{003} Y_2 Y_3^3 + b_{003} Y_3^3\right) + \left( - a_{012}Y_2^2 Y_3^2 + b_{012} Y_2 Y_3^2 \right) +\left( - a_{021}Y_2^3 Y_3+ b_{021}Y_2^2Y_3 \right)+\left(- a_{030} Y_2^4 + b_{030} Y_2^3 \right)+\left(- a_{102} Y_2 Y_3^2+ b_{102} Y_3^2 \right)}{Y_0}\\
&+ \frac{\left( - a_{111} Y_2^2 Y_3 + b_{111} Y_2Y_3 \right)+\left(- a_{120} Y_2^3 + b_{120} Y_2^2 \right) +\left( - a_{201} Y_2 Y_3 + b_{201} Y_3 \right) + \left( - a_{210} Y_2^2 + b_{210} Y_2\right) + \left( - a_{300} Y_2 + b_{300} \right) }{Y_0}\\
&= \frac{b_{300}+ \left( b_{210}- a_{300} \right) Y_2 + b_{201} Y_3 + \left( b_{120}- a_{210} \right) Y_2^2 + b_{102} Y_3^2 + \left( b_{111} - a_{201} \right) Y_2 Y_3 +\left( b_{030} - a_{120} \right) Y_2^3 + b_{003}Y_3^3 }{Y_0}\\
&+\frac{  \left( b_{021} - a_{111}\right) Y_2^2 Y_3+ \left( b_{012}- a_{102} \right) Y_2 Y_3^2 - a_{030} Y_2^4 - a_{021} Y_2^3 Y_3 - a_{012} Y_2^2 Y_3^2 - a_{003} Y_2 Y_3^3 }{Y_0}
\end{align*}}}
In order for $v_1^2\left( Y_0, Y_2 , Y_3 \right)$ to be holomorphic on $U_1$, we should have $b_{300}=0, b_{210} = a_{300}, b_{201}=0, b_{120}= a_{210}, b_{102}=0, b_{111}= a_{201}, b_{030}= a_{120}, b_{003}=0, b_{021}= a_{111}, b_{012}= a_{102}, a_{030}=0, a_{021}=0, a_{012}=0, a_{003}=0$. On the other hand, for $v_1^3\left( Y_0 , Y_2, Y_3\right)$, the case $(0,0,0),(0,0,1),(0,1,0),(1,0,0), (0,0,2),(0,1,1),(0,2,0)$,  $(1,0,1),(1,1,0),(2,0,0)$ defines a holomorphic function on $U_1$. Let us consider the remaining case which is
{\small{\begin{align*}
&\frac{ \left( - a_{003} Y_3^4+ c_{003} Y_3^3 \right) + \left(  - a_{012} Y_2 Y_3^3+ c_{012} Y_2 Y_3^2\right) + \left( - a_{021} Y_2^2 Y_3^2+ c_{021} Y_2^2 Y_3\right) +\left( - a_{030} Y_2^3 Y_3 + c_{030} Y_2^3 \right) + \left( - a_{102} Y_3^3 + c_{102} Y_3^2 \right) }{Y_0} \\
& + \frac{\left( - a_{111} Y_2 Y_3^2 + c_{111} Y_2 Y_3 \right) + \left( -a_{120} Y_2^2 Y_3 + c_{120} Y_2^2\right) + \left( - a_{201} Y_3^2 + c_{201} Y_3 \right) + \left( - a_{210} Y_2Y_3 + c_{210} Y_2 \right) + \left( - a_{300} Y_3 + c_{300} \right) }{Y_0} \\
&=\frac{c_{300} + c_{210} Y_2 + \left( c_{201} - a_{300} \right) Y_3 + c_{120} Y_2^2 + \left( c_{102} - a_{201} \right) Y_3^2 + \left( c_{111} - a_{210}  \right) Y_2 Y_3 + c_{030} Y_2^3 + \left( c_{003} - a_{102} \right) Y_3^3  }{Y_0}\\
&+ \frac{  \left( c_{021} - a_{120} \right) Y_2^2 Y_3 + \left( c_{012} - a_{111} \right)Y_2 Y_3^2 - a_{003} Y_3^4 - a_{030} Y_2^3 Y_3 - a_{021} Y_2^2 Y_3^2 - a_{012} Y_2 Y_3^3 }{Y_0}
\end{align*}}}
In order for $v_1^3\left( Y_0, Y_2, Y_3 \right)$ to be holomorphic on $U_1$, we should have $c_{300}=0, c_{210}=0, c_{201}= a_{300}, c_{120}=0, c_{102}= a_{201}, c_{111}= a_{210}, c_{030}=0, c_{003} = a_{102}, c_{021}= a_{120}, c_{012}= a_{111}, a_{003}=0, a_{030}=0, a_{021}=0, a_{012}=0$. Then all in all, we set
\begin{align*}
&t_1= a_{000}, \,\,\,\,\,t_2= a_{001},\,\,\,\,\, t_3= a_{010},\,\,\,\,\,t_4= a_{100},\,\,\,\,\, t_5 = a_{002},\,\,\,\,\, t_6= a_{011},\,\,\,\,\, t_7= a_{020},\,\,\,\,\, t_8= a_{101}\\
&t_9=a_{101},\,\,\,\,\, t_{10}= a_{200},\,\,\,\,\,t_{11}=a_{102}= b_{012}= c_{003},\,\,\,\,\,t_{12}= a_{111}= b_{021}= c_{012},\,\,\,\,\,t_{13}= a_{120}= b_{030}= c_{021},\,\,\,\,\,\\
&t_{14}= a_{201}= b_{111}= c_{102},\,\,\,\,\, t_{15}= a_{210} = b_{120}= c_{111} ,\,\,\,\,\,t_{16}= a_{300} = b_{210} = c_{201}\\
& t_{17}= b_{000},\,\,\,\,\,t_{18}= b_{001},\,\,\,\,\,t_{19}=b_{010},\,\,\,\,\,t_{20}= b_{100},\,\,\,\,\, t_{21}= b_{002},\,\,\,\,\,t_{22}= b_{011},\,\,\,\,\, t_{23}= b_{020},\,\,\,\,\,t_{24}= b_{101},\,\,\,\,\,t_{25}= b_{110},\,\,\,\,\,t_{26}= b_{200} \\
&t_{27} = c_{000},\,\,\,\,\, t_{28}= c_{001},\,\,\,\,\,t_{29}= c_{010},\,\,\,\,\, t_{30}= c_{100},\,\,\,\,\,t_{31}= c_{002},\,\,\,\,\, t_{32}= c_{011},\,\,\,\,\, t_{33}= c_{020},\,\,\,\,\, t_{34}= c_{101},\,\,\,\,\, t_{35} = c_{110},\,\,\,\,\, t_{36}= c_{200}
\end{align*}

Hence we have
\begin{align} \label{j9}
H^0\left( \mathbb{P}_\mathbb{C}^3, \Theta_{\mathbb{P}_\mathbb{C}^3}(1) \right)=&\left\{\left(t_1+t_2X_3+t_3X_2+t_4X_1+t_5X_3^2+t_6X_2X_3+t_7X_2^2+t_8X_1X_3+t_9X_1X_2 + t_{10}X_1^2\right.\right.\\
&\left. \,\,\,\,\,\,\,\,\,\, +t_{11}X_1X_3^2+t_{12}X_1X_2X_3+t_{13}X_1X_2^2+t_{14}X_1^2X_3+t_{15}X_1^2X_2+t_{16}X_1^3 \right)\frac{\partial}{\partial X_1} \notag \\
&+\left(t_{17}+t_{18}X_3+t_{19}X_2 +t_{20}X_1+t_{21}X_3^2+t_{22}X_2X_3 +t_{23}X_2^2+t_{24}X_1X_3 +t_{25}X_1X_2+t_{26}X_1^2\right. \notag \\
&\left.\,\,\,\,\,\,\,\,\,\,\,\,\,\,+t_{11}X_2X_3^2+t_{12}X_2^2X_3+t_{13}X_2^3 +t_{14}X_1X_2X_3+t_{15}X_1X_2^2+t_{16}X_1^2X_2\right)\frac{\partial}{\partial X_2} \notag \\
&+\left(t_{27}+t_{28}X_3+t_{29}X_2+t_{30}X_1 +t_{31}X_3^2+t_{32}X_2X_3 +t_{33}X_2^2 +t_{34}X_1X_3+t_{35}X_1X_2 +t_{36}X_1^2\right. \notag \\
&\,\,\,\,\,\,\,\,\,\,\,\,\,+\left.\left. t_{11}X_3^3+t_{12}X_2X_3^2+t_{13}X_2^2X_3 +t_{14}X_1X_3^2+t_{15}X_1X_2X_3+t_{16}X_1^2X_3\right)\frac{\partial}{\partial X_3}, t_1,..., t_{36} \in \mathbb{C} \right\} \notag
\end{align}

\subsection{Description of $H^0\left(  \mathbb{P}_\mathbb{C}^3, \Theta_{\mathbb{P}_\mathbb{C}^3}(-1)\right)$}\

From the Euler sequence, we have $0\to \mathcal{O}_{\mathbb{P}_\mathbb{C}^3}(-1) \to \bigoplus^4 \mathcal{O}_{\mathbb{P}_\mathbb{C}^3} \to \Theta_{\mathbb{P}_\mathbb{C}^3}(-1) \to 0$, so that we see that $\dim_\mathbb{C} H^0\left( \mathbb{P}_\mathbb{C}^3, \Theta_{\mathbb{P}_\mathbb{C}^3}(-1) \right)= 4$. We explicitly describe $H^0\left( \mathbb{P}_\mathbb{C}^3 , \Theta_{\mathbb{P}_\mathbb{C}^3}(-1) \right)$ on $U_0$. Let
\begin{align*}
T_0=v_0^1\left( X_1,X_2, X_3 \right) \frac{\partial}{\partial X_1}+ v_0^2 \left( X_1, X_2, X_3 \right) \frac{\partial}{\partial X_2} + v_0^3 \left( X_1, X_2, X_3 \right) \frac{\partial}{\partial X_3}\in H^0\left( \mathbb{P}_\mathbb{C}^3, \Theta_{\mathbb{P}_\mathbb{C}^3}(-1) \right)\,\,\,\,\,\textnormal{on} \,\,\, U_0
\end{align*}
be identified with
\begin{align*}
T_1= v_1^0\left( Y_0, Y_2, Y_3 \right) \frac{\partial}{\partial Y_0} + v_1^2\left( Y_0, Y_2, Y_3\right)\frac{\partial}{\partial Y_2} + v_1^3\left(  Y_0, Y_2, Y_3\right) \frac{\partial}{\partial Y_3} \in H^0\left(\mathbb{P}_\mathbb{C}^3, \Theta_{\mathbb{P}_\mathbb{C}^3}(-1) \right)\,\,\,\,\,\textnormal{on}\,\,\, U_1
\end{align*}
so that we have
\begin{align*}
v_0= X_1^{-1} v_1 = Y_0 v_1
\end{align*}
Then we have
\begin{align*}
v_0&= v_0^1\left(\frac{1}{Y_0}, \frac{Y_2}{Y_0}, \frac{Y_3}{Y_0} \right) \left(- Y_0^2 \frac{\partial}{\partial Y_0} - Y_0 Y_2 \frac{\partial}{\partial Y_2} - Y_0 Y_3 \frac{\partial}{\partial Y_3} \right)  + v_0^2\left(\frac{1}{Y_0}, \frac{Y_2}{Y_0}, \frac{Y_3}{Y_0} \right) Y_0\frac{\partial}{\partial Y_2} + v_0^3\left( \frac{1}{Y_0}, \frac{Y_2}{Y_0}, \frac{Y_3}{Y_0} \right) Y_0 \frac{\partial}{\partial Y_3}\\
&=Y_0 \left(v_1^0\left(Y_0,Y_2, Y_3 \right)\frac{\partial}{\partial Y_0} + v_1^2\left(Y_0,Y_2, Y_3 \right)\frac{\partial}{\partial Y_2}+ v_1^3\left(Y_0,Y_2, Y_3 \right)\frac{\partial}{\partial Y_3} \right)
\end{align*}
Then by letting $v_0^1\left( X_1, X_2, X_3 \right)= \sum_{h,k,r \geq 0} a_{hkr}X_1^h X_2^k X_3^r, v_0^2\left( X_1, X_2, X_3 \right)= \sum_{h,kr\geq 0} b_{hkr} X_1^h X_2^k X_3^r$, $v_0^3\left(X_1,X_2,X_3 \right)=\sum_{h,k,r\geq 0} c_{hkr}X_1^h X_2^k X_3^r, a_{hkr}, b_{hkr}, c_{hkr}\in \mathbb{C}$, we have
\begin{align*}
v_1^0\left(Y_0, Y_2, Y_3 \right) &= - Y_0 v_0^1\left(\frac{1}{Y_0}, \frac{Y_2}{Y_0}, \frac{Y_3}{Y_0} \right) =   \sum_{h,k,r \geq 0} \frac{- a_{hkr} Y_2^k Y_3^r}{Y_0^{h+k+r-1}}\\
v_1^2\left( Y_0, Y_2, Y_3 \right)& = - Y_2v_0^1\left(\frac{1}{Y_0}, \frac{Y_2}{Y_0}, \frac{Y_3}{Y_0} \right) + v_0^2\left( \frac{1}{Y_0}, \frac{Y_2}{Y_0}, \frac{Y_3}{Y_0} \right) = \sum_{h,k,r \geq 0} \frac{- a_{hkr} Y_2^{k+1} Y_3 ^r + b_{hkr} Y_2^k Y_3^r}{ Y_0^{h+k+r }}\\
v_1^3\left( Y_0, Y_2, Y_3 \right)&= - Y_3 v_0^1\left(\frac{1}{Y_0}, \frac{Y_2}{Y_0}, \frac{Y_3}{Y_0} \right) + v_0^3\left(\frac{1}{Y_0}, \frac{Y_2}{Y_0}, \frac{Y_3}{Y_0} \right) = \sum_{h,k,r \geq 0} \frac{ - a_{hkr} Y_2^k Y_3^{r+1}+ c_{hkr} Y_2^k Y_3^r}{Y_0^{h+k+r}}
\end{align*}
which should be holomorphic on $U_1$. We only need to consider the case $(h,k,r) = (0,0,0), (0,0,1),(0,1,0),(1,0,0)$. For $v_1^2\left( Y_0, Y_2, Y_3 \right)$, let us consider the rational part $\frac{\left( - a_{001} Y_2Y_3 + b_{001} Y_3 \right) + \left( - a_{010} Y_2^2 + b_{010} Y_2 \right) + \left( - a_{100} Y_2 + b_{100} \right) }{Y_0}$, so that we have $b_{100}=0, b_{010}= a_{100}, b_{001}=0, a_{001}=0, a_{010}=0$. On the other hand, for $v_1^3\left( Y_0, Y_2, Y_3 \right)$, let us consider the ration part $\frac{ \left( - a_{001} Y_3^2 + c_{001} Y_3\right) +\left( - a_{010} Y_2Y_3 + c_{010} Y_2 \right) + \left(  - a_{100} Y_3 + c_{100} \right) }{Y_0}$, so that we have $c_{100}=0, c_{010}= 0, c_{001}= a_{100}, a_{001}=0, a_{010}=0$. Then we set $t_1= a_{000}, t_2= a_{100}= b_{010}= c_{001}, t_3= b_{000}, t_4= c_{000}$. Then we have
\begin{align}\label{j10}
H^0\left( \mathbb{P}_{\mathbb{C}^3}, \Theta_{\mathbb{P}_\mathbb{C}^3}(-1) \right)= \left\{  \left( t_1 + t_2 X_1 \right) \frac{\partial}{\partial X_1} + \left( t_3 + t_2 X_2 \right) \frac{\partial}{\partial X_2} + \left( t_4 + t_2 X_3 \right) \frac{\partial}{\partial X_3}, t_1,...,t_4 \in \mathbb{C} \right\}
\end{align}

\subsection{Description of $H^0\left( \mathbb{P}_\mathbb{C}^3, \mathcal{E}_{\mathcal{O}_{\mathbb{P}_\mathbb{C}^3}(m)} \right)$ } \

Let us consider the line bundle $\mathcal{O}_{\mathbb{P}_\mathbb{C}^3}(m)$ on $\mathbb{P}_\mathbb{C}^3$ for some $m \in \mathbb{Z}$. In terms of the covering $\mathcal{U}$ of $\mathbb{P}_\mathbb{C}^3$, $\mathcal{O}_{\mathbb{P}_\mathbb{C}^3}(m)$ is described by the following transition functions
\begin{align*}
\left\{ f_{01}= X_1^m , f_{02}= X_2^m , f_{03}= X_3^m, f_{12}= Y_2^m , f_{13}= Y_3^m , f_{23}= Z_3^m  \right\} \in C^0\left( \mathcal{U}, \mathcal{O}_{\mathbb{P}_\mathbb{C}^3}^* \right)
\end{align*}
where a local section $s_i\in \Gamma\left(U_i, \mathcal{O}_{\mathbb{P}_\mathbb{C}^3} \right)$ of $\mathcal{O}_{\mathbb{P}_\mathbb{C}^3}(m)$ is identified with $s_j\in \Gamma\left( U_j, \mathcal{O}_{\mathbb{P}_\mathbb{C}^3} \right)$ if $s_i= f_{ij}a_j$ on $U_i\cap U_j$.

We note that from the Atiyah extension
\begin{align*}
0\to \mathcal{O}_{\mathbb{P}_\mathbb{C}^3} \to \mathcal{E}_{\mathcal{O}_{\mathbb{P}_\mathbb{C}^3}(m)} \to \Theta_{\mathbb{P}_\mathbb{C}^3} \to 0 
\end{align*}
we see that $\dim_\mathbb{C} H^0\left( \mathbb{P}_\mathbb{C}^3 , \mathcal{E}_{\mathcal{O}_{\mathbb{P}_\mathbb{C}^3} (m) } \right)= 16$. We explicitly describe $H^0\left( \mathbb{P}_\mathbb{C}^3, \mathcal{E}_{\mathcal{O}_{\mathbb{P}_\mathbb{C}^3}(m)} \right)$ on $U_0$. Let $f+T\in H^0\left(X, \mathcal{E}_{\mathcal{O}_{\mathbb{P}_\mathbb{C}^3} (m)}\right)$ on $U_0$, where $f=\sum_{i,j,k\geq  0} c_{ijk}X_1^iX_2^jX_3^k,c_{ijk}\in \mathbb{C}$ is a holomorphic function on $U_0$, and $T\in H^0\left(\mathbb{P}_{\mathbb{C}}^3 , \Theta_{\mathbb{P}_\mathbb{C}^3} \right)$ is given by $\textnormal{(\ref{ap2})}$.  By the definition of $\mathcal{E}_{\mathcal{O}_{\mathbb{P}_\mathbb{C}^3}(m)}$, $f+T$ on $U_0$ is identified with $f+ \left[T,\log X_1^m \right]+T$ on $U_1$, so that
 \begin{align*}
 &\sum_{i,j,k\geq 0} c_{ijk}X_1^iX_2^jX_3^k+\left[ T, \log X_1^m\right]\\
 &=\sum_{i,j,k\geq 0} c_{ijk}X_1^iX_2^j X_3^k+ m \left( \frac{a_1}{X_1}+\frac{a_2X_3}{X_1} + \frac{a_3 X_2}{X_1} +a_4+a_5X_3+a_6X_2+a_7X_1\right)\\
 &=\sum_{i,k,k\geq 0} c_{ijk}Y_0^{-i-j-k}Y_2^j Y_3^k+m\left( a_1Y_0+ a_2Y_3 + a_3 Y_2 +a_4+a_5\frac{Y_3}{Y_0}+a_6\frac{Y_2}{Y_0}+a_7\frac{1}{Y_0} \right)
 \end{align*}
 should be holomorphic on $U_1$, i.e. $c_{001}=-ma_5$ and $c_{010}=-ma_6$ and $c_{100}=-ma_7$, so that 
 \begin{align*}
 f=c_{000}-ma_7 X_1-ma_6 X_2 - ma_5 X_3
 \end{align*}
 Then we have
 \begin{align}\label{cp8}
 H^0\left( \mathbb{P}_\mathbb{C}^3 ,\mathcal{E}_{\mathcal{O}_{\mathbb{P}_\mathbb{C}^3} (m)}\right)=&\left\{f+ T_1\frac{\partial}{\partial X_1}+ T_2 \frac{\partial}{\partial X_2} + T_3 \frac{\partial}{\partial X_3}\right.\\
 &\left.=a_{16}-ma_7X_1-ma_6X_2-ma_5X_3 \right. \notag \\
 &\left.+ \left(a_1+a_2X_3+a_3X_2+a_4X_1+a_5X_1X_3+a_6X_1X_2+a_7X_1^2\right)\frac{\partial}{\partial X_1}\right.\notag\\
 &\left.+\left(a_8+a_9X_3+a_{10}X_2 +a_{11}X_1+a_5X_2X_3+a_6X_2^2+a_7X_1X_2 \right)\frac{\partial}{\partial X_2} \right. \notag \\
 &\left.+\left(a_{12}+a_{13}X_3+a_{14}X_2+a_{15}X_1+a_5X_3^2+a_6X_2X_3+a_7X_1X_3 \right)\frac{\partial}{\partial X_3}, a_1,...,a_{16}\in \mathbb{C} \right\} \notag
 \end{align}
 
 \section{Description of $H^0\left( \mathbb{P}_\mathbb{C}^3, \Omega_{\mathbb{P}_\mathbb{C}^3}^1(4) \right)$ and $H^0\left( \mathbb{P}_\mathbb{C}^3 , \bigwedge^2 \Omega_{\mathbb{P}_\mathbb{C}^3}^1(8) \right)$ and $H^0\left( \mathbb{P}_\mathbb{C}^3, \mathcal{O}_{\mathbb{P}_\mathbb{C}^3}(4) \right)$}\label{ja1}
 
 We note that $\dim_\mathbb{C} H^0\left( \mathbb{P}_\mathbb{C}^3, \Omega_{\mathbb{P}_\mathbb{C}^3}^1(4)\right)= \frac{6 \cdot 5 \cdot 5 \cdot 3}{2}=45 $ by Bott formula (see \cite{OSS88} p.4). We explicitly describe $H^0\left( \mathbb{P}_\mathbb{C}^3, \Omega_{\mathbb{P}_\mathbb{C}^3}^1(4) \right)$ on $U_0$. Let
 \begin{align*}
v_0= v_0^1\left( X_1, X_2, X_3 \right) dX_1 + v_0^2 \left( X_1, X_2, X_3\right)dX_2 + v_0^3 \left( X_1, X_2, X_3 \right) dX_3 \in H^0\left( \mathbb{P}_\mathbb{C}^3 , \Omega_{\mathbb{P}_\mathbb{C}^3}^1 (4)  \right) \,\,\,\,\,\textnormal{on}\,\,\, U_0
 \end{align*}
 be identified with
 \begin{align*}
 v_1= v_1^0\left(Y_0, Y_2, Y_3 \right) dY_0 + v_1^2\left( Y_0, Y_2, Y_3 \right) dY_2 + v_1^3 \left( Y_0, Y_2, Y_3 \right) dY_3 \in H^0\left( \mathbb{P}_\mathbb{C}^3, \Omega_{\mathbb{P}_\mathbb{C}^3}^1(4) \right)\,\,\,\,\,\textnormal{on}\,\,\, U_1
 \end{align*}
 so that we have
 \begin{align*}
 v_0= X_1^4 v_1= \frac{1}{Y_0^4} v_1
 \end{align*}
From $X_1= \frac{1}{Y_0}, X_2= \frac{Y_2}{Y_0}, X_3=\frac{Y_3}{Y_0}$, we have
\begin{align*}
dX_1=-\frac{1}{Y_0^2}dY_0,\,\,\,\,\, dX_2=-\frac{Y_2}{Y_0^2}dY_0+\frac{1}{Y_0}dY_2,\,\,\,\,\,dX_3=-\frac{Y_3}{Y_0^2}dY_0+\frac{1}{Y_0}dY_3
\end{align*}
Then we have
{\Small{\begin{align*}
v_0&=-\frac{1}{Y_0^2}v_0^1\left(\frac{1}{Y_0},\frac{Y_2}{Y_0},\frac{Y_3}{Y_0} \right) dY_0 +v_0^2\left(\frac{1}{Y_0},\frac{Y_2}{Y_0},\frac{Y_3}{Y_0} \right)\left(-\frac{Y_2}{Y_0^2}dY_0+\frac{1}{Y_0} dY_2 \right) +v_0^3\left(\frac{1}{Y_0},\frac{Y_2}{Y_0},\frac{Y_3}{Y_0} \right)\left(-\frac{Y_3}{Y_0^2}dY_0+\frac{1}{Y_0}dY_3 \right)\\
&=\left(-\frac{1}{Y_0^2}v_0^1\left(\frac{1}{Y_0},\frac{Y_2}{Y_0},\frac{Y_3}{Y_0} \right)-\frac{Y_2}{Y_0^2}v_0^2\left(\frac{1}{Y_0},\frac{Y_2}{Y_0},\frac{Y_3}{Y_0} \right)-\frac{Y_3}{Y_0^2}v_0^3\left(\frac{1}{Y_0},\frac{Y_2}{Y_0},\frac{Y_3}{Y_0} \right) \right)dY_0+\frac{1}{Y_0}v_0^2\left(\frac{1}{Y_0},\frac{Y_2}{Y_0},\frac{Y_3}{Y_0} \right)dY_2+\frac{1}{Y_0}v_0^3\left( \frac{1}{Y_0},\frac{Y_2}{Y_0},\frac{Y_3}{Y_0}\right)dY_3\\
&=\frac{1}{Y_0^4}v_1^0\left( Y_0, Y_2, Y_3 \right)dY_0+\frac{1}{Y_0^4}v_1^2\left(Y_0,Y_2,Y_3 \right)dY_2 + \frac{1}{Y_0^4}v_1^3\left( Y_0, Y_2, Y_3 \right)
\end{align*}}}
Hence we have
\begin{align*}
v_1^0\left(Y_0,Y_2,Y_3\right)&=-Y_0^2v_0^1\left(\frac{1}{Y_0},\frac{Y_2}{Y_0},\frac{Y_3}{Y_0} \right)-Y_0^2 Y_2v_0^2\left(\frac{1}{Y_0},\frac{Y_2}{Y_0},\frac{Y_3}{Y_0} \right)-Y_0^2Y_3v_0^3\left(\frac{1}{Y_0},\frac{Y_2}{Y_0},\frac{Y_3}{Y_0} \right)\\
v_1^2\left(Y_0,Y_2,Y_3\right)&=Y_0^3v_0^2\left( \frac{1}{Y_0},\frac{Y_2}{Y_0},\frac{Y_3}{Y_0}\right)\\
v_1^3 \left( Y_0,Y_2,Y_3 \right)&= Y_0^3v_0^3\left(\frac{1}{Y_0},\frac{Y_2}{Y_0},\frac{Y_3}{Y_0} \right)
\end{align*}
 By letting
\begin{align*}
v_0^1\left(X_1,X_2,X_3\right)&=\sum_{h,k,r\geq 0} a_{hkr}X_1^h X_2^k X_3^r,\,\,\,\,\,a_{hkr}\in \mathbb{C}\\
v_0^2\left(X_1,X_2,X_3\right)&=\sum_{h,k,r\geq 0} b_{hkr}X_1^h X_2^k X_3^r,\,\,\,\,\,b_{hkr}\in \mathbb{C}\\
v_0^2\left(X_1,X_2,X_3\right)&=\sum_{h,k,r\geq 0} c_{hkr}X_1^h X_2^k X_3^r,\,\,\,\,\, c_{hkr}\in \mathbb{C}
\end{align*}
we get
\begin{align*}
v_1^0\left(Y_0,Y_2,Y_3 \right)&=\sum_{h,k,r\geq 0}\frac{-a_{hkr}Y_2^kY_3^r -b_{hkr} Y_2^{k+1}Y_3^r-c_{hkr}Y_2^k Y_3^{r+1}}{Y_0^{h+k+r-2}} \\
v_1^2\left(Y_0, Y_2, Y_3 \right)&=\sum_{h,k,r\geq 0} b_{hkr}\frac{Y_2^k Y_3^r}{Y_0^{h+k+r-3}}\\
v_1^3\left(Y_0, Y_2, Y_3 \right)&=\sum_{h,k,r\geq 0} c_{hkr}\frac{Y_2^k Y_3^r}{Y_0^{h+k+r-3}}
\end{align*}
In order for $v_1^0\left( Y_0, Y_2, Y_3\right)$ to be holomorphic on $U_1$, we only need to consider the case
\begin{align*}
(h,k,r)=&(0,0,0),\\
             &(0,0,1),(0,1,0),(1,0,0),\\
             &(0,0,2),(0,1,1),(0,2,0),(1,0,1),(1,1,0),(2,0,0) \\
             &(0,0,3),(0,1,2),(0,2,1),(0,3,0),(1,0,2),(1,1,1),(1,2,0),(2,0,1),(2,1,0),(3,0,0)
\end{align*}
The case $(0,0,0),(0,0,1),(0,1,0), (1,0,0),(0,0,2),(0,1,1),(0,2,0),(1,0,1),(1,1,0),(2,0,0)$ for $v_1^0\left( Y_0,Y_2,Y_3 \right)$ clearly defines a holomorphic function on $U_1$. Let us consider the remaining case which is
\begin{align*}
&\frac{\left( a_{003}Y_3^3+b_{003}Y_2 Y_3^3+c_{003}Y_3^4\right) + \left( a_{012}Y_2 Y_3^2+b_{012} Y_2^2 Y_3^2+c_{012}Y_2 Y_3^3\right)    +  \left( a_{021}Y_2^2 Y_3+b_{021} Y_2^3 Y_3+c_{021} Y_2^2 Y_3^2\right)}{Y_0}\\
+&\frac{\left(a_{030}Y_2^3+b_{030}Y_2^4+c_{030} Y_2^3 Y_3 \right)+ \left(  a_{102}Y_3^2+b_{102} Y_2 Y_3^2+c_{102} Y_3^3   \right) +     \left(a_{111} Y_2 Y_3+b_{111} Y_2^2  Y_3+c_{111} Y_2 Y_3^2     \right)}{Y_0}\\
+&\frac{\left(a_{120}Y_2^2+b_{120}Y_2^3+c_{120}Y_2^2Y_3        \right) +   \left( a_{201}Y_3+b_{201}Y_2 Y_3+c_{201} Y_3^2    \right) +     \left( a_{210} Y_2+b_{210}Y_2^2+c_{210}Y_2 Y_3     \right)}{Y_0}\\
+&\frac{\left( a_{300}+b_{300}Y_2+c_{300}Y_3     \right)}{Y_0}\\
=& \frac{a_{300}+\left( a_{210}+b_{300}\right)Y_2+\left(a_{201}+c_{300} \right)Y_3 + \left(a_{120}+b_{210} \right) Y_2^2 + \left( a_{102}+c_{201}\right)Y_3^2 + \left(a_{111}+b_{201}+c_{210} \right)Y_2 Y_3}{Y_0} \\
+&\frac{\left( a_{030}+b_{120}\right)Y_2^3 +\left( a_{003}+c_{102}     \right)Y_3^3 +\left(a_{021}+b_{111}+c_{120}     \right)Y_2^2 Y_3 +\left(a_{012}+b_{102}+c_{111} \right)Y_2 Y_3^2}{Y_0}\\
+&\frac{b_{030}Y_2^4+c_{003}Y_3^4 +\left(b_{003}+c_{012} \right)Y_2 Y_3^3 + \left(b_{012}+c_{021}   \right)Y_2^2 Y_3^2 +\left(b_{021}+ c_{030} \right)Y_2^3Y_3}{Y_0}
\end{align*}
In order for $v_1^0\left( Y_0, Y_2, Y_3\right)$ to be holomorphic on $U_1$, we should have $a_{300}=0, a_{210}=-b_{300}, a_{201}=-c_{300},a_{120}=-b_{210}, a_{102}=-c_{201},a_{111}=-b_{201}-c_{210},a_{030}=-b_{120},a_{003}=-c_{102}, a_{021}=-b_{111}-c_{120}, a_{012}=-b_{102}-c_{111}, b_{030}=0, c_{003}=0, b_{003}=-c_{012},b_{012}=-c_{021},b_{021}=-c_{030}$.

On the other hand, in order for $v_1^2 \left( Y_0, Y_2, Y_3 \right)$ and $v_1^3\left( Y_0, Y_2, Y_3 \right)$ to be holomorphic on $U_1$, we should have $b_{hkr}=c_{hkr}=0$ except for
\begin{align*}
(h,k,r)=&(0,0,0)\\
           &(0,0,1),(0,1,0),(1,0,0)\\
           &(0,0,2),(0,1,1),(0,2,0),(1,0,1),(1,1,0),(2,0,0)\\
           &(0,0,3),(0,1,2),(0,2,1),(0,3,0),(1,0,2),(1,1,1),(1,2,0),(2,0,1),(2,1,0),(3,0,0)
\end{align*}
Then we set
\begin{align*}
&t_1=b_{000}\\
&t_2=b_{001},\,\,\,t_3=b_{010},\,\,\,t_4=b_{100}\\
&t_5=b_{002},\,\,\,t_6=b_{011},\,\,\,t_7=b_{020},\,\,\,t_8=b_{101},\,\,\,t_9=b_{110},\,\,\,t_{10}=b_{200}\\
&t_{11}=b_{003},\,\,\,t_{12}=b_{012},\,\,\,t_{13}=b_{021},\,\,\,t_{14}=b_{102},\,\,\,t_{15}=b_{111},\,\,\,t_{16}=b_{120},\,\,\,t_{17}=b_{201},\,\,\,t_{18}=b_{210},\,\,\,t_{19}=b_{300}
\end{align*}

\begin{align*}
&t_{20}=c_{000}\\
&t_{21}=c_{001},\,\,\,t_{22}=c_{010},\,\,\,t_{23}=c_{100}\\
&t_{24}=c_{002},\,\,\,t_{25}=c_{011},\,\,\,t_{26}=c_{020},\,\,\,t_{27}=c_{101},\,\,\,t_{28}=c_{110},\,\,\,t_{29}=c_{200}\\
&c_{012}=-t_{11},\,\,\,c_{021}=-t_{12},\,\,\,c_{030}=-t_{13},\,\,\,t_{30}=c_{102},\,\,\,t_{31}=c_{111},\,\,\,t_{32}=c_{120},\,\,\,t_{33}=c_{201},\,\,\,t_{34}=c_{210},\,\,\,t_{35}=c_{300}
\end{align*}
\begin{align*}
&t_{36}=a_{000} \\
&t_{37}=a_{001},\,\,\,t_{38}=a_{010},\,\,\,t_{39}=a_{100}\\
&t_{40}=a_{002},\,\,\,t_{41}=a_{011},\,\,\,t_{42}=a_{020},\,\,\,t_{43}=a_{101},\,\,\,t_{44}=a_{110},\,\,\,t_{45}=a_{200}
\end{align*}
We note that
\begin{align*}
a_{300}&=0\\
a_{210}&=-b_{300}=-t_{19}\\
a_{201}&=-c_{300}=-t_{35}\\
a_{120}&=-b_{210}=-t_{18}\\
a_{102}&=-c_{201}=-t_{33}\\
a_{111}&=-b_{201}-c_{210}=-t_{17}-t_{34}\\
a_{030}&=-b_{120}=-t_{16}\\
a_{003}&=-c_{102}=-t_{30}\\
a_{021}&=-b_{111}-c_{120}=-t_{15}-t_{32}\\
a_{012}&=-b_{102}-c_{111}=-t_{14}-t_{31}\\
b_{030}&=0\\
c_{003}&=0\\
b_{003}&=-c_{012} = t_{11}\\
b_{012}&=-c_{021}= t_{12}\\
b_{021}&=-c_{030}= t_{13}
\end{align*}
Hence $H^0\left( \mathbb{P}_\mathbb{C}^3, \Omega_{\mathbb{P}_\mathbb{C}^3}^1(4) \right)$ on $U_0$ is given by 
{\Small{\begin{align}\label{cm1}
&b=\left(t_{36}+t_{37}X_3+t_{38}X_2+t_{39}X_1+t_{40}X_3^2+t_{41}X_2 X_3+t_{42}X_2^2+t_{43}X_1X_3+t_{44}X_1X_2+t_{45}X_1^2 \right) dX_1\\
+&\left(-t_{19}X_1^2X_2-t_{35}X_1^2X_3-t_{18}X_1X_2^2-t_{33}X_1X_3^2-(t_{17}+t_{34})X_1X_2X_3-t_{16}X_2^3-t_{30}X_3^3-(t_{15}+t_{32})X_2^2X_3-(t_{14}+t_{31})X_2X_3^2 \right)dX_1\notag \\ 
+&\left(t_1+t_2X_3+t_3X_2+t_4X_1+t_5X_3^2+t_6X_2X_3+t_7X_2^2 +t_8X_1X_3+t_9X_1X_2+t_{10}X_1^2 +t_{11}X_3^3+t_{12}X_2X_3^2 \right)dX_2 \notag\\
+&\left( t_{13}X_2^2X_3 +t_{14}X_1X_3^2+t_{15}X_1X_2X_3+t_{16}X_1X_2^2+t_{17}X_1^2X_3+t_{18}X_1^2X_2+t_{19}X_1^3   \right)dX_2 \notag\\
+&\left(t_{20}+t_{21}X_3+t_{22}X_2+t_{23}X_1+t_{24}X_3^2+t_{25}X_2X_3+t_{26}X_2^2+t_{27}X_1X_3+t_{28}X_1X_2+t_{29}X_1^2 -t_{11}X_2X_3^2-t_{12}X_2^2X_3-t_{13}X_2^3 \right)dX_3 \notag \\
+&\left(t_{30}X_1X_3^2+t_{31}X_1X_2X_3+t_{32}X_1X_2^2+t_{33}X_1^2X_3+t_{34}X_1^2X_2+t_{35}X_1^3 \right)dX_3 \notag
\end{align}}}

We compute
\begin{align*}
&d b=(t_{38}+t_{41}X_3+2t_{42}X_2+t_{44}X_1)dX_2\wedge dX_1+(t_{37}+2t_{40}X_3+t_{41}X_2+t_{43}X_1)dX_3\wedge dX_1\\
&+(-t_{19}X_1^2-2t_{18}X_1X_2-(t_{17}+t_{34})X_1X_3-3t_{16}X_2^2-2(t_{15}+t_{32})X_2X_3-(t_{14}+t_{31})X_3^2)dX_2\wedge dX_1\\
&+\left(-t_{35}X_1^2-2t_{33}X_1X_3-(t_{17}+t_{34})X_1X_2-3t_{30}X_3^2 -(t_{15}+t_{32})X_2^2-2(t_{14}+t_{31})X_2X_3 \right)dX_3\wedge dX_1\\
&+\left( t_4+t_8X_3+t_9X_2+2t_{10}X_1 + t_{14}X_3^2+t_{15}X_2X_3+t_{16}X_2^2+2t_{17}X_1X_3+2t_{18}X_1X_2+3t_{19}X_1^2\right)dX_1\wedge dX_2\\
&+\left( t_2+2t_5X_3+t_6X_2+t_8X_1+3t_{11}X_3^2+2t_{12}X_2X_3+t_{13}X_2^2+2t_{14}X_1X_3+t_{15}X_1X_2+t_{17}X_1^2\right)dX_3\wedge dX_2\\
&+\left(t_{23}+t_{27}X_3+t_{28}X_2+2t_{29}X_1 +t_{30}X_3^2+t_{31}X_2X_3+t_{32}X_2^2+2t_{33}X_1X_3+2t_{34}X_1X_2+3t_{35}X_1^2\right)dX_1\wedge dX_3\\
&+\left( t_{22} + t_{25}X_3+2t_{26}X_2+t_{28}X_1-t_{11}X_3^2-2t_{12}X_2X_3-3t_{13}X_2^2 + t_{31}X_1X_3+2t_{32}X_1X_2+ t_{34}X_1^2 \right)dX_2\wedge dX_3
\end{align*}

Consider the coefficient of $dX_2\wedge dX_3$ of $d b$:
\begin{align*}
&- t_2-2t_5X_3-t_6X_2-t_8X_1-3t_{11}X_3^2-2t_{12}X_2X_3-t_{13}X_2^2-2t_{14}X_1X_3-t_{15}X_1X_2-t_{17}X_1^2\\
&+ t_{22} + t_{25}X_3+2t_{26}X_2+t_{28}X_1-t_{11}X_3^2-2t_{12}X_2X_3-3t_{13}X_2^2 + t_{31}X_1X_3+2t_{32}X_1X_2+ t_{34}X_1^2 \\
&=(t_{22}-t_2)+(t_{28}-t_8)X_1+(2t_{26}-t_6)X_2+(t_{25}-2t_5)X_3+(t_{34}-t_{17})X_1^2+(-3t_{13}-t_{13})X_2^2+(-4t_{11})X_3^2\\
&+( 2t_{32}-t_{15})X_1X_2 + (t_{31}-2t_{14})X_1X_3+ (-4t_{12} )X_2X_3
\end{align*}

Consider the coefficient of $dX_3\wedge dX_1$ of $d b$:
\begin{align*}
&t_{37}+2t_{40}X_3+t_{41}X_2+t_{43}X_1\\
&-t_{35}X_1^2-2t_{33}X_1X_3-(t_{17}+t_{34})X_1X_2-3t_{30}X_3^2 -(t_{15}+t_{32})X_2^2-2(t_{14}+t_{31})X_2X_3\\
&-t_{23}-t_{27}X_3-t_{28}X_2-2t_{29}X_1 -t_{30}X_3^2-t_{31}X_2X_3-t_{32}X_2^2-2t_{33}X_1X_3-2t_{34}X_1X_2-3t_{35}X_1^2\\
&=(t_{37}-t_{23})+(t_{43}-2t_{29})X_1 +(t_{41}-t_{28})X_2+(2t_{40}-t_{27})X_3-4t_{35}X_1^2-(t_{15}+2t_{32})X_2^2-4t_{30}X_3^2\\
&-(t_{17}+3t_{34} )X_1X_2 -4t_{33}X_1X_3-(2t_{14}+3t_{31})X_2X_3
\end{align*}

Consider the coefficient of $dX_1\wedge dX_2$ of $d b$:
{\small{\begin{align*}
&-t_{38}-t_{41}X_3-2t_{42}X_2-t_{44}X_1 \\
&+t_{19}X_1^2+2t_{18}X_1X_2+(t_{17}+t_{34})X_1X_3+3t_{16}X_2^2+2(t_{15}+t_{32})X_2X_3+(t_{14}+t_{31})X_3^2\\
&+t_4+t_8X_3+t_9X_2+2t_{10}X_1 + t_{14}X_3^2+t_{15}X_2X_3+t_{16}X_2^2+2t_{17}X_1X_3+2t_{18}X_1X_2+3t_{19}X_1^2\\
&=(t_4-t_{38})+(2t_{10}-t_{44})X_1+(t_9-2t_{42})X_2+(t_8-t_{41})X_3+4t_{19}X_1^2+4t_{16}X_2^2+(2t_{14}+t_{31})X_3^2\\
&+4t_{18}X_1X_2 + \left(3t_{17}+t_{34}\right)X_1X_3+(3t_{15}+2t_{32})X_2X_3
\end{align*}}}

All in all, we have
{\Small{\begin{align}\label{cm2}
db&= \left(t_{22}-t_2+(t_{28}-t_8)X_1+(2t_{26}-t_6)X_2+(t_{25}-2t_5)X_3+(t_{34}-t_{17})X_1^2-4t_{13}X_2^2-4t_{11}X_3^2 \right.\\
&\left. \,\,\,\,\,\,\,\,\,\,+( 2t_{32}-t_{15})X_1X_2 + (t_{31}-2t_{14})X_1X_3-4t_{12} X_2X_3 \right )dX_2\wedge dX_3 \notag\\
&+\left(t_{37}-t_{23}+(t_{43}-2t_{29})X_1 +(t_{41}-t_{28})X_2+(2t_{40}-t_{27})X_3-4t_{35}X_1^2-(t_{15}+2t_{32})X_2^2-4t_{30}X_3^2 \right. \notag\\
&\left.\,\,\,\,\,\,\,\,\,\,-(t_{17}+3t_{34} )X_1X_2 -4t_{33}X_1X_3-(2t_{14}+3t_{31})X_2X_3 \right )dX_3\wedge dX_1 \notag \\
&+ \left(t_4-t_{38}+(2t_{10}-t_{44})X_1+(t_9-2t_{42})X_2+(t_8-t_{41})X_3+4t_{19}X_1^2+4t_{16}X_2^2+(2t_{14}+t_{31})X_3^2 \right. \notag \\
&\left. \,\,\,\,\,\,\,\,\,\,+4t_{18}X_1X_2 + (3t_{17}+t_{34})X_1X_3+(3t_{15}+2t_{32})X_2X_3 \right)dX_1\wedge dX_2 \notag
\end{align}}}

\subsection{Description of $H^0\left(\mathbb{P}_\mathbb{C}^3, \bigwedge^2 \Omega_{\mathbb{P}_\mathbb{C}^3}^1 (8) \right)$}\

We note that $\dim_\mathbb{C} H^0\left(\mathbb{P}_\mathbb{C}^3, \bigwedge^2 \Omega_{\mathbb{P}_\mathbb{C}^3}^1 (8) \right)=189$ by Bott formula (see \cite{OSS88} p.4). We explicitly describe $H^0\left( \mathbb{P}_\mathbb{C}^3, \bigwedge^2 \Omega_{\mathbb{P}_\mathbb{C}^3}^1(8) \right)$ on $U_0$. Let
\begin{align*}
v_0= A_0\left(X_1,X_2,X_3 \right) dX_1 \wedge dX_2  + B_0 \left(X_1,X_2,X_3 \right) dX_2\wedge dX_3 + C_0\left( X_1,X_2,X_3 \right)dX_3\wedge dX_1
\end{align*}
be identified with
\begin{align*}
v_1= A_1\left( Y_0,Y_2,Y_3\right) dY_0\wedge dY_2 + B_1\left(Y_0,Y_2,Y_3\right)dY_2\wedge dY_3 + C_1\left(Y_0, Y_2, Y_3 \right) dY_3\wedge dY_0
\end{align*}
so that we have
\begin{align*}
v_0= X_1^8 v_1= \frac{1}{Y_0^8}v_1
\end{align*}
We note that 
{\small{\begin{align*}
dX_1\wedge dX_2=- \frac{1}{Y_0^3}dY_0\wedge dY_2,\,\,\,\,\,\,\,\, dX_2\wedge dX_3= \frac{Y_2}{Y_0^3} dY_3 \wedge dY_0 + \frac{Y_3}{Y_0^3}d Y_0\wedge dY_2 + \frac{1}{Y_0^2} dY_2\wedge dY_3,\,\,\,\,\,\,\,\,\, dX_3\wedge dX_1 = - \frac{1}{Y_0^3} dY_3\wedge dY_0
\end{align*}}}
Then we have
{\small{\begin{align*}
v_0&=\left( - \frac{1}{Y_0^3}A_0\left( \frac{1}{Y_0}, \frac{Y_2}{Y_0}, \frac{Y_3}{Y_0} \right) + \frac{Y_3}{Y_0^3}B_0\left( \frac{1}{Y_0}, \frac{Y_2}{Y_0}, \frac{Y_3}{Y_0} \right) \right) dY_0\wedge dY_2+ \frac{1}{Y_0^2}B_0\left( \frac{1}{Y_0} , \frac{Y_2}{Y_0}, \frac{Y_3}{Y_0}\right) dY_2 \wedge dY_3 \\
&+ \left(  - \frac{1}{Y_0^3}C_0\left( \frac{1}{Y_0}, \frac{Y_2}{Y_0}, \frac{Y_3}{Y_0} \right) + \frac{Y_2}{Y_0^3}B_0\left( \frac{1}{Y_0}, \frac{Y_2}{Y_0}, \frac{Y_3}{Y_0} \right) \right)dY_3 \wedge dY_0 = \frac{1}{Y_0^8} A_1\left(Y_0, Y_2, Y_3 \right) + \frac{1}{Y_0^8}B_1\left( Y_0, Y_2, Y_3\right) + \frac{1}{Y_0^8} C_1\left( Y_0, Y_2, Y_3 \right)
\end{align*}}}
Letting $A_0\left(X_1,X_2,X_3\right)=\sum_{h,k,r\geq 0} a_{hkr}X_1^h X_2^k X_3^r, B_0\left(X_1,X_2,X_3\right)=\sum_{h,k,r\geq 0} b_{hkr}X_1^h X_2^k X_3^r,  C_0\left(X_1,X_2,X_3\right)=\sum_{h,k,r\geq 0} c_{hkr}X_1^h X_2^k X_3^r, a_{hkr}, b_{hkr}, c_{hkr}\in \mathbb{C}$, we have
\begin{align*}
A_1\left( Y_0, Y_2, Y_3 \right)&= \sum_{h,k,r \geq 0} \frac{- a_{hkr} Y_2^k Y_3^r + b_{hkr} Y_2^k Y_3^{r+1}}{Y_0^{h+k+r -5}}\\
B_1\left( Y_0, Y_2, Y_3 \right)&= \sum_{h,k,r \geq 0} \frac{b _{hkr} Y_2^k Y_3^r}{Y_0^{h+k+r-6}}\\
C_1\left( Y_0, Y_2, Y_3 \right)&= \sum_{h,k,r\geq 0} \frac{-c _{hkr}Y_2^k Y_3^r + b_{hkr} Y_2^{k+1} Y_3^r}{Y_0^{h+k+r-5}}
\end{align*} 
In order for $A_1\left( Y_0, Y_2, Y_3 \right)$ to be holomorphic on $U_1$, we only need to consider
\begin{align*}
(h,k,r)=&(0,0,0),\\
             &(0,0,1),(0,1,0),(1,0,0),\\
             &(0,0,2),(0,1,1),(0,2,0),(1,0,1),(1,1,0),(2,0,0)  \\
             &(0,0,3),(0,1,2),(0,2,1),(0,3,0),(1,0,2), (1,1,1),(1,2,0), (2,0,1),(2,1,0),(3,0,0)\\
             &(0,0,4),(0,1,3),(0,2,2),(0,3,1),(0,4,0),(1,0,3),(1,1,2),(1,2,1),(1,3,0),(2,0,2),(2,1,1),(2,2,0),\\
             &(3,0,1),(3,1,0),(4,0,0)\\
             &(0,0,5),(0,1,4),(0,2,3),(0,3,2),(0,4,1),(0,5,0),(1,0,4),(1,1,3),(1,2,2),(1,3,1),(1,4,0),(2,0,3),\\
             &(2,1,2),(2,2,1),(2,3,0),(3,0,2),(3,1,1),(3,2,0),(4,0,1),(4,1,0),(5,0,0)\\
             &(0,0,6),(0,1,5),(0,2,4),(0,3,3),(0,4,2),(0,5,1),(0,6,0),(1,0,5),(1,1,4),(1,2,3),(1,3,2),(1,4,1),\\
             &(1,5,0),(2,0,4),(2,1,3),(2,2,2),(2,3,1),(2,4,0),(3,0,3),(3,1,2),(3,2,1),(3,3,0),(4,0,2),(4,1,1),\\
             &(4,2,0),(5,0,1),(5,1,0),(6,0,0)
\end{align*}
The case $(0,0,0)-(5,0,0)$ for $A_1\left(Y_0, Y_2, Y_3 \right)$ clearly define s holomorphic function on $U_1$. Let us consider the remaining case which is
{\Small{\begin{align*}
&\frac{\left(- a_{006} Y_3^6+ b_{006} Y_3^7 \right) +\left( - a_{015}Y_2Y_3^5 + b_{015}Y_2Y_3^6 \right) + \left( - a_{024} Y_2^2Y_3^4+ b_{024} Y_2^2 Y_3^5\right) + \left(- a_{033} Y_2^3 Y_3^3 + b_{033}Y_2^3Y_3^4 \right)+\left(- a_{042} Y_2^4 Y_3^2 + b_{042}Y_2^4Y_3^3 \right)}{Y_0}\\
& + \frac{\left(- a_{051} Y_2^5 Y_3 + b_{051}Y_2^5Y_3^2 \right) + \left(- a_{060} Y_2^6 + b_{060}Y_2^6Y_3 \right) + \left(- a_{105}  Y_3^5 + b_{105}Y_3^6 \right) + \left(- a_{114} Y_2 Y_3^4 + b_{114}Y_2Y_3^5 \right) + \left(- a_{123} Y_2^2 Y_3^3 + b_{123}Y_2^2Y_3^4 \right)}{Y_0}\\
& + \frac{\left(- a_{132} Y_2^3 Y_3^2 + b_{132}Y_2^3Y_3^3 \right) + \left(- a_{141} Y_2^4 Y_3 + b_{141}Y_2^4Y_3^2 \right) + \left(- a_{150}  Y_2^5 + b_{150}Y_2^5Y_3 \right) + \left(- a_{204} Y_3^4 + b_{204}Y_3^5 \right) + \left(- a_{213} Y_2 Y_3^3 + b_{213}Y_2Y_3^4 \right)}{Y_0}\\
& + \frac{\left(- a_{222} Y_2^2 Y_3^2 + b_{222}Y_2^2Y_3^3 \right) + \left(- a_{231} Y_2^3 Y_3 + b_{231}Y_2^3Y_3^2 \right) + \left(- a_{240}  Y_2^4 + b_{240}Y_2^4Y_3\right) + \left(- a_{303}  Y_3^3 + b_{303}Y_3^4 \right) + \left(- a_{312} Y_2 Y_3^2 + b_{312}Y_2Y_3^3 \right)}{Y_0}\\
& + \frac{\left(- a_{321} Y_2^2 Y_3 + b_{321}Y_2^2Y_3^2 \right) + \left(- a_{330} Y_2^3 + b_{330}Y_2^3Y_3 \right) + \left(- a_{402}  Y_3^2 + b_{402}Y_3^3 \right) + \left(- a_{411} Y_2 Y_3 + b_{411}Y_2Y_3^2 \right) + \left(- a_{420} Y_2^2 + b_{420}Y_2^2Y_3 \right)}{Y_0}\\
& + \frac{\left(- a_{501}  Y_3 + b_{501}Y_3^2 \right) + \left(- a_{510} Y_2 + b_{510}Y_2Y_3 \right) + \left(- a_{600}   + b_{600} Y_3 \right) }{Y_0}\\
&= \frac{- a_{600} - a_{510} Y_2 + (b_{600}- a_{501}) Y_3 -a_{420} Y_2^2 +(b_{510}- a_{411}) Y_2 Y_3 +(b_{501}- a_{402}) Y_3^2 + (b_{402}- a_{303}) Y_3^3 + (b_{411}- a_{312}) Y_2Y_3^2}{Y_0}\\
&+ \frac{(b_{420}- a_{321})Y_2^2 Y_3 - a_{330} Y_2^3 - a_{240} Y_2^4 + ( b_{330} - a_{231}) Y_2^3 Y_3 + (b_{321}- a_{222} ) Y_2^2 Y_3^2 + (b_{312} - a_{213}) Y_2 Y_3^3 + ( b_{303}- a_{204}) Y_3^4 - a_{150} Y_2^5}{Y_0}\\
&+\frac{ (b_{240} - a_{141}) Y_2^4 Y_3 + ( b_{231} - a_{132}) Y_2^3 Y_3^2 + (b_{222} - a_{123}) Y_2^2 Y_3^3 + (b_{213} - a_{114})Y_2 Y_3^4 + ( b_{204} - a_{105}) Y_3^5 - a_{060}Y_2^6  + ( b_{150} - a_{051} ) Y_2^5 Y_3}{Y_0}  \\
&+ \frac{(b_{141} - a_{042}) Y_2^4 Y_3^2 + (b_{132} - a_{033}) Y_2^3 Y_3^3 + ( b_{123} - a_{024}) Y_2^2 Y_3^4 + (  b_{114} - a_{015}) Y_2 Y_3^5 + ( b_{105} - a_{006} ) Y_3^6 + b_{006}Y_3^7+ b_{015}Y_2Y_3^6 + b_{024}Y_2^2 Y_3^5}{Y_0}\\
&+ \frac{b_{033} Y_2^3 Y_3^4+ b_{042}Y_2^4 Y_3^3+ b_{051} Y_2^5 Y_3^2+ b_{060}Y_2^6Y_3}{Y_0}
\end{align*}}}
Then we should have $a_{600}=0, a_{510}=0, b_{600}= a_{501}, a_{420}=0,  b_{510}= a_{411}, b_{501}= a_{402}, b_{402}= a_{303}, b_{411}=a_{312}, b_{420}=a_{321}, a_{330}=0, a_{240}=0,  b_{330} = a_{231}, b_{321}= a_{222}, b_{312}= a_{213}, b_{303}= a_{204}, a_{150}=0, b_{240}= a_{141}, b_{231} = a_{132} , b_{222}= a_{123}, b_{213}= a_{114}, b_{204}= a_{105}, a_{060}=0,b_{150}= a_{051}, b_{141}= a_{042}, b_{132}= a_{033}, b_{123}= a_{024}, b_{114}= a_{015}, b_{105}= a_{006}, b_{006}=0, b_{015}=0, b_{024}=0, b_{033}=0, b_{042}=0, b_{051}= 0, b_{060}=0 $. On the other hand, similarly for $C_1\left( Y_0, Y_2, Y_3 \right)$ to be holomorphic on $U_1$, we should have
{\Small{\begin{align*}
&\frac{\left(- c_{006}Y_3^6+ b_{006} Y_2 Y_3^6 \right) + \left( - c_{015} Y_2 Y_3^5 + b_{015} Y_2^2 Y_3^5 \right) + \left( - c_{024} Y_2^2Y_3^4 + b_{024}Y_2^3Y_3^4 \right) +\left(- c_{033} Y_2^3Y_3^3 + b_{033} Y_2^4Y_3^3  \right) + \left( - c_{042}X_2^4 X_3^2+ b_{042}X_2^5 X_3^2 \right) }{Y_0}\\
&+ \frac{\left(- c_{051} Y_2^5 Y_3 + b_{051} Y_2^6 Y_3\right)+ \left( - c_{060}X_2^6 + b_{060} Y_2^7 \right) + \left(- c_{105} Y_3^5 + b_{105} Y_2 Y_3^5 \right) + \left( - c_{114} Y_2 Y_3^4 + b_{114} Y_2^2 Y_3^4 \right) +\left( - c_{123} Y_2^2 Y_3^3 + b_{123}Y_2^3 Y_3^3 \right)  }{Y_0}\\
&+ \frac{ \left( - c_{132} Y_2^3 Y_3^2 + b_{132} Y_2^4 Y_3^2\right)  + \left( - c_{141} Y_2^4 Y_3 + b_{141} Y_2^5 Y_3 \right) +\left( - c_{150}Y_2^5 + b_{150}Y_2^6 \right) +\left( - c_{024} Y_2^2 Y_3^4 + b_{024} Y_2^3 Y_3^4 \right) +\left( - c_{033} Y_2^3 Y_3^3+ b_{033} Y_2^4 Y_3^3 \right)}{Y_0} \\
&+ \frac{\left(- c_{042} Y_2^4 Y_3^2+ b_{042} Y_2^5 Y_3^2 \right) +\left( - c_{051} Y_2^5 Y_3 + b_{051} Y_2^6 Y_3 \right) +\left(- c_{060} Y_2^6 + b_{060} Y_2^7 \right) +\left( - c_{105} Y_3^5+ b_{105} Y_1 Y_3^5 \right) +\left( - c_{114} Y_2 Y_3^4 + b_{114} Y_2^2 Y_3^4 \right) }{Y_0} \\
&+ \frac{\left( - c_{123} Y_2^2 Y_3^3 + b_{123} Y_2^3 Y_3^3 \right) +\left( - c_{132} Y_2^3 Y_3^2 + b_{132} Y_2^4 Y_3^2 \right) +\left( - c_{141} Y_2^4 Y_3 + b_{141} Y_2^5 Y_3 \right) +\left(- c_{150} Y_2^5 + b_{150} Y_2^6 \right) +\left( - c_{204} Y_3^4 + b_{204} Y_2 Y_3^4\right)}{Y_0}\\
&+ \frac{ \left( - c_{213} Y_2 Y_3^3+ b_{213} Y_2^2 Y_3^3\right) +\left( - c_{222} Y_2^2 Y_3^2 + b_{222} Y_2^3 Y_3^2 \right) +\left( - c_{231} Y_2^3 Y_3 + b_{231} Y_2^4 Y_3 \right) +\left( - c_{240} Y_2^4 + b_{240} Y_2^5 \right) +\left(- c_{303} Y_3^3 + b_{303} Y_2 Y_3^3 \right) }{Y_0} \\
&+ \frac{ \left(- c_{312} Y_2 Y_3^2 + b_{312} Y_2^2 Y_3^2\right)+\left(- c_{321} Y_2^2 Y_3 + b_{321} Y_2^3 Y_3 \right)+\left(- c_{330} Y_2^3 + b_{330} Y_2^4 \right)+\left( - c_{402} Y_3^2 + b_{402} Y_2 Y_3^2 \right) +\left( - c_{411} Y_2 Y_3 + b_{411} Y_2^2 Y_3 \right) }{Y_0}\\
&+ \frac{ \left( - c_{420} Y_2^2 + b_{420} Y_2^3\right) + \left(- c_{501} Y_3 + b_{501} Y_2 Y_3 \right)+\left(- c_{510} Y_2 + b_{510} Y_2^2 \right) + \left( - c_{600} + b_{600} Y_2 \right)}{Y_0} \\
&= \frac{- c_{600} + (b_{600}- c_{510}) Y_2 - c_{501} Y_3 + (b_{510}- c_{420})Y_2^2 - c_{402} Y_3^2 + (b_{501}- c_{411}) Y_2Y_3 + (b_{402}- c_{312}) Y_2Y_3^2+ (b_{411}- c_{ 321}) Y_2^2 Y_3 }{Y_0}\\
&+ \frac{( b_{420} - c_{330}) Y_2^3  - c_{303}  Y_3^3 + (b_{330} - c_{240} ) Y_2^4 + (b_{321} - c_{231} ) Y_2^3 Y_3 + ( b_{312} - c_{222}) Y_2^2 Y_3^2 + ( b_{303} - c_{213} ) Y_2 Y_3^3 - c_{204} Y_3^4 + ( b_{240} - c_{150}) Y_2^5 }{Y_0}\\
&+ \frac{( b_{231}- c_{141} ) Y_2^4 Y_3 + ( b_{222} - c_{132}) Y_2^3 Y_3^2 + ( b_{213} - c_{123} ) Y_2^2 Y_3^3 + ( b_{204}  - c_{114}) Y_2 Y_3^4 - c_{105} Y_3^5 + ( b_{150} - c_{060} ) Y_2^6 + ( b_{141} - c_{051} ) Y_2^5 Y_3 }{Y_0}\\
&+ \frac{(b_{132}  - c_{042}) Y_2^4 Y_3^2 + ( b_{123} - c_{033} ) Y_2^3 Y_3^3 + ( b_{114} - c_{024}) Y_2^2 Y_3^4 + (  b_{105} - c_{015} ) Y_2 Y_3^5 - c_{006} Y_3^6 + b_{006}Y_2Y_3^6 + b_{015}Y_2^2 Y_3^5 + b_{024} Y_2^3 Y_3^4  }{Y_0}\\
&+ \frac{b_{033} Y_2^4 Y_3^3+ b_{042} X_2^5X_3^2 + b_{051} Y_2^6 Y_3 + b_{060} Y_2^7 }{Y_0}
\end{align*}}}
Then we should have $c_{600}=0, b_{600}= c_{510}, c_{501}=0, b_{510}= c_{420}, c_{402}=0, b_{501}= c_{411}, b_{402}= c_{312}, b_{411}= c_{321}, b_{420}= c_{330}, c_{303}=0, b_{330}= c_{240}, b_{321}= c_{231}, b_{312}= c_{222}, b_{303}= c_{213}, c_{204}=0, b_{240}= c_{150}, b_{231}= c_{141}, b_{222}=c_{132}, b_{213}= c_{123}, b_{204}= c_{114}, c_{105}=0, b_{150}= c_{060}, b_{141}= c_{051}, b_{132}= c_{042}, b_{123}=c_{033}, b_{114}= c_{024}, b_{105}= c_{015}, c_{006}=0, b_{006}=0, b_{015}=0, b_{024}=0, b_{033}=0, b_{042}= 0 , b_{051}=0, b_{060}=0$.

Then $H^0\left( \mathbb{P}_\mathbb{C}^3 , \bigwedge^2 \Omega_{\mathbb{P}_\mathbb{C}^3}^1(8) \right)$ on $U_0$ is given by
{\Small{\begin{align} \label{ja8}
&\left( a_{000} + a_{001}X_3+ a_{010}X_2 + a_{100}X_1 + a_{002} X_3^2 + a_{011}X_2X_3+ a_{020}X_2^2 + a_{101}X_1X_3+ a_{110}X_1X_2+ a_{200}X_1^2 + a_{003}X_3^3+ a_{012}X_2X_3^2 \right) dX_1\wedge dX_2\\
&+ \left( a_{021} X_2^2 X_3 + a_{030}X_2^3+ a_{102}X_1 X_3^2+ a_{111} X_1X_2X_3+ a_{120} X_1 X_2^2+ a_{201}X_1^2 X_3+ a_{210}X_1^2 X_2 + a_{300}X_1^3 + a_{004}X_3^4 + a_{013}X_2X_3^3\right) dX_1 \wedge dX_2 \notag \\
&+\left( a_{022}X_2^2 X_3^2+ a_{031}X_2^3 X_3 + a_{040}X_2^4 + a_{103}X_1 X_3^3+ a_{112}X_1X_2X_3^2  + a_{121}X_1X_2^2 X_3 + a_{130}X_1X_2^3+ a_{202}X_1^2 X_3^2+ a_{211}X_1^2 X_2X_3\right) dX_1 \wedge dX_2 \notag \\
&+ \left(  a_{220}X_1^2 X_2^2 + a_{301}X_1^3 X_3 + a_{310}X_1^3 X_2 + a_{400}X_1^4 + a_{005}X_3^5+ a_{014}X_2 X_3^4+ a_{023}X_2^2 X_3^3 + a_{032}X_2^3 X_3^2 + a_{041}X_2^4 X_3+ a_{050}X_2^5\right) dX_1 \wedge dX_2 \notag\\
&+ \left(a_{104}X_1 X_3^4+  a_{113}X_1X_2X_3^3 + a_{122}X_1 X_2^2 X_3^2 + a_{131}X_1 X_2^3 X_3 + a_{140}X_1X_2^4+ a_{203}X_1^2 X_3^3+ a_{212}X_1^2 X_2 X_3^2 + a_{221}X_1^2 X_2^2 X_3 \right) dX_1 \wedge dX_2 \notag \\
&+ \left( a_{230}X_1^2 X_2^3 + a_{302}X_1^3 X_3^2 + a_{311}X_1^3 X_2X_3 + a_{320}X_1^3 X_2^2+ a_{401}X_1^4 X_3+ a_{410}X_1^4 X_2+ a_{500}X_1^5 + b_{105}X_3^6+ b_{114}X_2X_3^5+ b_{123}X_2^2 X_3^4 \right) dX_1 \wedge dX_2 \notag \\
&+ \left( b_{132} X_2^3 X_3^3 + b_{141}X_2^4X_3^2+ b_{150}X_2^5 X_3+ b_{204} X_1 X_3^5+ b_{213}X_1X_2 X_3^4+ b_{222}X_1X_2^2X_3^3+ b_{231}X_1X_2^3 X_3^2 + b_{240}X_1X_2^4X_3 \right) dX_1 \wedge dX_2 \notag \\
&+ \left(  b_{303}X_1^2 X_3^4+ b_{312}X_1^2 X_2 X_3^3+ b_{321}X_1^2 X_2^2 X_3^2 + b_{330}X_1^2 X_2^3 X_3+ b_{402}X_1^3 X_3^3+ b_{411}X_1^3 X_2 X_3^2 \right) dX_1 \wedge dX_2 \notag  \\
&+ \left( b_{420}X_1^3 X_2^2 X_3+ b_{501}X_1^4 X_3^2 + b_{510}X_1^4 X_2X_3  + b_{600}X_1^5 X_3  \right) dX_1 \wedge dX_2  \notag
\end{align}
\begin{align*}
& + \left( b_{000} + b_{001}X_3+ b_{010}X_2 + b_{100}X_1 + b_{002} X_3^2 + b_{011}X_2X_3+ b_{020}X_2^2 + b_{101}X_1X_3+ b_{110}X_1X_2+ b_{200}X_1^2 + b_{003}X_3^3+ b_{012}X_2X_3^2 \right) dX_2\wedge dX_3\\
&+ \left( b_{021} X_2^2 X_3 + b_{030}X_2^3+ b_{102}X_1 X_3^2+ b_{111} X_1X_2X_3+ b_{120} X_1 X_2^2+ b_{201}X_1^2 X_3+ b_{210}X_1^2 X_2 + b_{300}X_1^3 + b_{004}X_3^4 + b_{013}X_2X_3^3\right) dX_2 \wedge dX_3 \\
&+\left( b_{022}X_2^2 X_3^2+ b_{031}X_2^3 X_3 + b_{040}X_2^4 + b_{103}X_1 X_3^3+ b_{112}X_1X_2X_3^2  + b_{121}X_1X_2^2 X_3 + b_{130}X_1X_2^3+ b_{202}X_1^2 X_3^2+ b_{211}X_1^2 X_2X_3\right) dX_2 \wedge dX_3 \\
&+ \left(  b_{220}X_1^2 X_2^2 + b_{301}X_1^3 X_3 + b_{310}X_1^3 X_2 + b_{400}X_1^4 + b_{005}X_3^5+ b_{014}X_2 X_3^4+ b_{023}X_2^2 X_3^3 + b_{032}X_2^3 X_3^2 + b_{041}X_2^4 X_3+ b_{050}X_2^5\right) dX_2 \wedge dX_3 \\
&+ \left(b_{104}X_1 X_3^4+  b_{113}X_1X_2X_3^3 + b_{122}X_1 X_2^2 X_3^2 + b_{131}X_1 X_2^3 X_3 + b_{140}X_1X_2^4+ b_{203}X_1^2 X_3^3+ b_{212}X_1^2 X_2 X_3^2 + b_{221}X_1^2 X_2^2 X_3 \right) dX_2 \wedge dX_3 \\
&+ \left( b_{230}X_1^2 X_2^3 + b_{302}X_1^3 X_3^2 + b_{311}X_1^3 X_2X_3 + b_{320}X_1^3 X_2^2+ b_{401}X_1^4 X_3+ b_{410}X_1^4 X_2+ b_{500}X_1^5 \right) dX_2 \wedge dX_3 \\
&+ \left(  b_{105}X_1 X_3^5+ b_{114}X_1X_2 X_3^4+ b_{123}X_1X_2^2X_3^3+ b_{132}X_1X_2^3 X_3^2 + b_{141} X_1X_2^4X_3 \right) dX_2 \wedge dX_3 \\
&+ \left( b_{150} X_1 X_2^5+ b_{204}X_1^2 X_3^4+ b_{213}X_1^2 X_2 X_3^3+ b_{222} X_1^2 X_2^2 X_3^2 + b_{231}X_1^2 X_2^3 X_3+ b_{240}X_1^2 X_2^4+ b_{303}X_1^3 X_3^3+ b_{312}X_1^3 X_2 X_3^2 \right) dX_2 \wedge dX_3 \\
&+ \left( b_{321}X_1^3 X_2^2 X_3 + b_{330}X_1^3 X_2^3+ b_{402}X_1^4 X_3^2 + b_{411}X_1^4 X_2X_3 + b_{420} X_1^4 X_2^2 + b_{501}X_1^5 X_3+ b_{510}X_1^5 X_2 + b_{600}X_1^6 \right) dX_2 \wedge dX_3 
\end{align*}

\begin{align*}
&+\left( c_{000} + c_{001}X_3+ c_{010}X_2 + c_{100}X_1 + c_{002} X_3^2 + c_{011}X_2X_3+ c_{020}X_2^2 + c_{101}X_1X_3+ c_{110}X_1X_2+ c_{200}X_1^2 + c_{003}X_3^3+ c_{012}X_2X_3^2 \right) dX_3\wedge dX_1\\
&+ \left( c_{021} X_2^2 X_3 + c_{030}X_2^3+ c_{102}X_1 X_3^2+ c_{111} X_1X_2X_3+ c_{120} X_1 X_2^2+ c_{201}X_1^2 X_3+ c_{210}X_1^2 X_2 + c_{300}X_1^3 + c_{004}X_3^4 + c_{013}X_2X_3^3\right) dX_3 \wedge dX_1 \\
&+\left( c_{022}X_2^2 X_3^2+ c_{031}X_2^3 X_3 + c_{040}X_2^4 + c_{103}X_1 X_3^3+ c_{112}X_1X_2X_3^2  + c_{121}X_1X_2^2 X_3 + c_{130}X_1X_2^3+ c_{202}X_1^2 X_3^2+ c_{211}X_1^2 X_2X_3\right) dX_3 \wedge dX_1 \\
&+ \left(  c_{220}X_1^2 X_2^2 + c_{301}X_1^3 X_3 + c_{310}X_1^3 X_2 + c_{400}X_1^4 + c_{005}X_3^5+ c_{014}X_2 X_3^4+ c_{023}X_2^2 X_3^3 + c_{032}X_2^3 X_3^2 + c_{041}X_2^4 X_3+ c_{050}X_2^5\right) dX_3 \wedge dX_1\\
&+ \left( c_{104}X_1 X_3^4+  c_{113}X_1X_2X_3^3 + c_{122}X_1 X_2^2 X_3^2 + c_{131}X_1 X_2^3 X_3 + c_{140}X_1X_2^4+ c_{203}X_1^2 X_3^3+ c_{212}X_1^2 X_2 X_3^2 + c_{221}X_1^2 X_2^2 X_3 \right) dX_3 \wedge dX_1 \\
&+ \left( c_{230}X_1^2 X_2^3 + c_{302}X_1^3 X_3^2 + c_{311}X_1^3 X_2X_3 + c_{320}X_1^3 X_2^2+ c_{401}X_1^4 X_3+ c_{410}X_1^4 X_2+ c_{500}X_1^5+ b_{105}X_2X_3^5+ b_{114}X_2^2 X_3^4 \right) dX_3 \wedge dX_1 \\
&+ \left( b_{123}X_2^3 X_3^3 + b_{132}X_2^4X_3^2+ b_{141}X_2^5 X_3+ b_{150}X_2^6+ b_{204}X_1X_2 X_3^4+ b_{213}X_1X_2^2X_3^3+ b_{222}X_1X_2^3 X_3^2 + b_{231}X_1X_2^4X_3 \right) dX_3 \wedge dX_1 \\
&+ \left( b_{240} X_1 X_2^5+ b_{303}X_1^2 X_2 X_3^3+ b_{312}X_1^2 X_2^2 X_3^2 + b_{321}X_1^2 X_2^3 X_3+ b_{330}X_1^2 X_2^4+ b_{402}X_1^3 X_2 X_3^2 \right) dX_3 \wedge dX_1 \\
&+ \left( b_{411}X_1^3 X_2^2 X_3 + b_{420}X_1^3 X_2^3+ b_{501}X_1^4 X_2X_3 + b_{510} X_1^4 X_2^2+ b_{600}X_1^5 X_2  \right) dX_3 \wedge dX_1 
\end{align*}}}

\subsubsection{\textnormal{Description of $H^0\left( \mathbb{P}_\mathbb{C}^3, \mathcal{O}_{\mathbb{P}_\mathbb{C}^3}(4) \right)$}}\

We note that $\dim_\mathbb{C} H^0\left( \mathbb{P}_\mathbb{C}^3, \mathcal{O}_{\mathbb{P}_\mathbb{C}^3}(4) \right)=35$, and any element $c$ of $H^0\left( \mathbb{P}_\mathbb{C}^3, \mathcal{O}_{\mathbb{P}_\mathbb{C}^3}(4) \right)$ on $U_0$ is given by
{\small{\begin{align}\label{ja6}
c=&c_1+ c_2X_1+ c_3 X_2+ c_4X_3 +c_5 X_1^2+ c_6 X_1X_2+ c_7 X_1X_3+ c_8 X_2^2 + c_9 X_2X_3+ c_{10}X_3^2 + c_{11}X_1^3 + c_{12}X_1^2 X_2 + c_{13}X_1^2 X_3+ c_{14}X_1X_2^2 \\
+& c_{15}X_1X_2X_3 + c_{16}X_1X_3^2 + c_{17}X_2^3+ c_{18}X_2^2 X_3+ c_{19}X_2X_3^2 + c_{20}X_3^3 +c_{21}X_1^4 + c_{22}X_2^4 + c_{23} X_3^4 + c_{24} X_1^3 X_2 + c_{25}X_1^3 X_3 + c_{26} X_1^2 X_2^2  \notag \\
+& c_{27} X_1^2 X_2 X_3 + c_{28}X_1^2 X_3^2 + c_{29}X_1 X_2^3 + c_{30}X_1X_2^2 X_3+ c_{31} X_1X_2 X_3^2 + c_{32}X_1X_3^3 + c_{33}X_2^3X_3 + c_{34}X_2^2 X_3^2 + c_{35} X_2 X_3^3 \notag
\end{align}}}
for constants $c_1,c_2, c_3,..., c_{35}\in \mathbb{C}$. Then we have
{\small{\begin{align}\label{ja7}
dc&= \left( c_2 + 2c_5X_1+ c_6X_2 +c_7X_3 + 3 c_{11} X_1^2 + c_{14}X_2^2 + c_{16}X_3^2 + 2 c_{12}X_1X_2 + 2 c_{13}X_1X_3 + c_{15}X_2 X_3\right)dX_1  \\
&+\left( 4c_{21}X_1^3+ c_{29}X_2^3 + c_{32}X_3^3+ 3 c_{24}X_1^2 X_2 + 3 c_{25}X_1^2 X_3 + 2 c_{26}X_1X_2^2 + 2 c_{27}X_1X_2X_3 + 2 c_{28}X_1X_3^2 + c_{30}X_2^2 X_3 + c_{31}X_2X_3^2 \right)dX_1 \notag \\
&+\left( c_3 + c_6X_1 + 2c_8X_2 + c_9 X_3+ c_{12}X_1^2 + 3 c_{17}X_2^2 + c_{19}X_3^2 + 2 c_{14} X_1 X_2 + c_{15}X_1X_3 +  2c_{18}X_2X_3\right)dX_2 \notag \\
&+ \left(c_{24}X_1^3+4c_{22}X_2^3 + c_{35}X_3^3 +  2 c_{26}X_1^2 X_2 + c_{27}X_1^2 X_3 + 3 c_{29}X_1X_2^2 + 2 c_{30}X_1X_2 X_3 + c_{31}X_1X_3^2 + 3 c_{33}X_2^2 X_3 + 2 c_{34}X_2X_3^2 \right)dX_2  \notag \\
&+\left( c_4 + c_7X_1 + c_9X_2 + 2 c_{10} X_3 + c_{13} X_1^2 + c_{18}X_2^2 + 3 c_{20}X_3^2  + c_{15}X_1X_2+ 2 c_{16}X_1X_3 + 2c_{19}X_2X_3\right)dX_3 \notag \\
&+ \left(c_{25}X_1^3+c_{33}X_2^3+4 c_{23}X_3^3 + c_{27}X_1^2 X_2 + 2 c_{28}X_1^2 X_3  + c_{30}X_1X_2^2 + 2c_{31} X_1X_2X_3 + 3 c_{32}X_1X_3^2 + 2 c_{34}X_2^2 X_3 + 3 c_{35}X_2X_3^2\right)dX_3 \notag
\end{align}}}

\section{Poisson cyclic coverings}

We can construct `Poisson cyclic coverings' in a similar way to construct foliated cyclic coverings in section \ref{c1}. Let $\left( M , \Lambda_0 \right)$ be a compact holomorphic Poisson manifold, where $\Lambda_0\in H^0\left(M, \bigwedge^2 \Theta_M \right)$ with $\left[ \Lambda_0 , \Lambda_0 \right]=0$. Our presentation is based on \cite{Wav68} to which we refer for the detail. We will use the notion of a Poisson line bundle on a holomorphic Poisson manifold (see \cite{Kim16}).

\begin{definition}[compare \cite{Wav68}]
Let $\left(M, \pi, W\right)$ be a $k$-sheeted branched covering of $W$, where $\pi:M\to W$. Let $C$ denote the branch locus. Moreover let $\left(M,\Pi_0\right)$ be a holomorphic Poisson manifold with $\Pi_0\in H^0\left(M, \bigwedge^2 \Theta_M \right)$, and $\left(W, \Lambda_0\right)$ be a holomorphic Poisson manifold with $\Lambda_0\in H^0\left( W, \bigwedge^2 \Theta_W \right)$.
\begin{enumerate}
\item for each $x\in C, \pi^{-1}(x)$ consists of one point.
\item the group of covering transformations of $M- \pi^{-1}(C)$ over $W- C$ is cyclic of order $k$.
\item for each $p\in C$, we can find a neighborhood $U$ with coordinates $\left(\xi_1,..., \xi_n\right)$ on $U$ and $\left(z_1,..., z_n\right)$ on $\pi^{-1}(U)$ such that the map is given by $\xi_1=z_i(1\leq i \leq n-1), \xi_n=\left(z_n\right)^k$.
\item if $k\ne 2$, $C$ is connected.
\item $C$ is a holomorphic Poisson submanifold of $\left(W, \Lambda_0\right)$.
\item $\pi:\left(M, \Pi_0\right)\to \left(W, \Lambda_0\right)$ is a holomorphic Poisson map, i.e. $\pi_{*} \left( \Pi_{0,x} \right) = \Lambda_{0,\pi(x)}$ for any $x\in M$ under the natural map $\pi_{* }:\bigwedge^2 \Theta_{M, x} \to \bigwedge^2 \Theta_{W, \pi(x)}$.
\end{enumerate}
\end{definition}

\begin{remark}
The condition $(3)$ implies that the branch locus $C$ is non-singular and $(3)$ is actually a consequence of $(1)$ and $(2)$ $($see \cite{Wav68} \textnormal{Appendix p.957}$)$. The condition $(6)$ implies the condition $(5)$. In fact, let $\Lambda_0$ on a neighborhood $U$ of $W$ with the coordinate $\left(\xi_1,..., \xi_n \right)$ be written as
\begin{align*}
\Lambda = \sum_{\alpha, \beta=1}^{n-1} \Lambda_{\alpha\beta}(\xi_1,..., \xi_n)\frac{\partial}{\partial \xi_\alpha} \wedge \frac{\partial}{\partial \xi_\beta}  + \sum_{\gamma=1}^{n-1} \Lambda_{\gamma n}(\xi_1,..., \xi_n)\frac{\partial}{\partial \xi_\gamma}\wedge \frac{\partial}{\partial \xi_n}
\end{align*}
and let $\Pi_0$ on $\pi^{-1}(U)$ be written as
\begin{align*}
\Pi_0= \sum_{\alpha,\beta=1}^{n-1} \Pi_{\alpha\beta}\left( z_1,..., z_n \right)\frac{\partial}{\partial z_\alpha} \wedge \frac{\partial}{\partial z_\beta} + \sum_{\gamma=1}^{n-1} \Pi_{\gamma n}\left(z_1,..., z_n\right)\frac{\partial}{\partial z_\gamma}\wedge \frac{\partial}{\partial z_n}
\end{align*}
Since $\pi$ is a holomorphic Poisson map, we have
\begin{align*}
\Lambda_{\gamma n}\left(z_1,..., z_{n-1}, (z_n)^k \right)  = \Pi_{\gamma n}\left(z_1,..., z_n\right)\cdot k \left( z_n\right)^{k-1}
\end{align*}
This implies that $\Lambda_{\gamma n}\left(z_1,..., z_{n-1}, \left(z_n \right)^k \right)$ has a factor $\left(  z_n \right)^{k-1}$. If $\Lambda_{\gamma n}\left(\xi_1,..., \xi_n\right)$ is expended in powers of $\xi_n$ whose coefficients are holomorphic functions of $\xi_1,..., \xi_{n-1}$, the term involving only $\xi_1,.., \xi_{n-1}$ is $0$, so that $\Lambda_{\gamma n}$ is of the form $\Lambda_{\gamma n}(\xi_1,..., \xi_n)= \xi_n G_{\gamma n}\left(\xi_1,..., \xi_n\right)$. This implies that $\xi_n=0$ on $U$ defines a Poisson submanifold.

\end{remark} 

Before proving the following theorem \ref{pcc1}, we recall a notion of Poisson line bundle on a holomorphic Poisson manifold (see \cite{Kim16}).  Given a Poisson line bundle on a holomorphic Poisson manifold $\left(W, \Lambda_0\right)$, we will make the total space of the line bundle a holomorphic Poisson manifold, and we will also make the total space of the associated projective line bundle a holomorphic Poisson manifold.

\begin{definition}
Let $\left( W, \Lambda_0 \right)$ be a holomorphic Poisson manifold with the Poisson bracket $\left\{ -, - \right\}$ induced by $\Lambda_0\in H^0\left( W, \bigwedge^2 \Theta_W \right)$. A Poisson connection on a line bundle $F$ on $\left(W, \Lambda_0 \right)$ is a $\mathbb{C}$-linear bracket $\left\{ - , - \right\}_F: \mathcal{O}_W \times \mathcal{O}\left( F \right)\to \mathcal{O}\left( F \right)$ where $\mathcal{O} \left(F \right)$ is the sheaf on $W$ of germs of holomorphic sections of $F$, such that $\left\{ fg, s \right\}_F = f\left\{ g, s\right\}_F + g \left\{ f,g s\right\}_F$ and $\left\{ f, gs\right\}_F=  \left\{ f, g \right\} s  + g \left\{f, s \right\}_F$ where $f,g \in \mathcal{O}_W, s\in \mathcal{O}\left( F \right)$. Equivalently, a Poisson connection is given by a homomorphism $\nabla_F:\mathcal{O}\left( F \right)\to \mathscr{H}om_{\mathcal{O}_W}\left(\Omega_{W}^1, \mathcal{O}\left( F \right) \right)=\textnormal{Der}\left( \mathcal{O}_W, \mathcal{O}\left(  F \right) \right) $ which satisfies the identity $\nabla_F(gs)= - \left[ \Lambda_0, g\right]  \otimes s + g \cdot \nabla_F(s) $ and the bracket is defined by the formula $\nabla_F (s)f=\left\{ f, s\right\}_F$. A Poisson connection $v$ is called flat if it satisfies the Jacobi identity $\left\{ f, \left\{ g, s \right\}_F \right\}_F = \left\{ \left\{ f, g \right\}, s \right\}_F + \left\{ g, \left\{ f, g \right\}_F \right\}_F$. A line bundle $F$ on $\left( W, \Lambda_0\right)$ equipped with a flat Poisson connection $\nabla_F$ is called a Poisson line bundle. We denote the Poisson line bundle by $\left( F , \nabla_F \right)$.
\end{definition}

Now let $\left(F, \nabla_F\right)$ be a Poisson line bundle on a holomorphic Poisson manifold $\left(W, \Lambda_0\right)$. Let $\left\{V_i \right\}$ be a covering of $W$ by polycylinders with coordinates $\left(\xi_i^1,..., \xi_i^n\right)$ such that $F$ is trivialized and $\Lambda_0$ is of the form on $V_i$:
\begin{align*}
\Lambda_0= \Lambda_i= \sum_{\alpha,\beta=1}^n \Lambda_{\alpha\beta}^i(\xi_i)\frac{\partial}{\partial \xi_i^\alpha}\wedge \frac{\partial}{\partial \xi_i^\beta}
\end{align*}
and $\left(F, \nabla_F\right)$ is represented by $\left( \left\{T_i\right\}, \left\{f_{ij}\right\} \right)$ where $T_i\in \Gamma\left( V_i, \Theta_W\right)$ and $f_{ij}\in \Gamma\left( V_{ij},  \mathcal{O}_W^* \right)$ such that
\begin{align*}
\left[ \Lambda_0, T_i\right]&=0 \\
T_j- T_i + \left[ \Lambda_0, \log f_{ij} \right]&=0 \\
f_{ij} f_{jk}&= f_{ik}
\end{align*}
Then $\left(F^k, \nabla_{F^k}\right)$ is naturally a Poisson line bundle on $\left(W, \Lambda_0\right)$ represented by $\left( \left\{ k T_i\right\}, \left\{ f_{ij}^k \right\} \right)$ (for the detail, see \cite{Kim16}).

Next we will show that a Poisson line bundle $\left(F, \nabla_F\right)$ on $\left(W, \Lambda_0\right)$ as above defines a Poisson structure $\tilde{\Lambda}_0$ on the total space $F$, so that $\left(F,\tilde{\Lambda}_0\right)$ is a holomorphic Poisson manifold. Let the total space $F$ be described by
\begin{align*}
V_{ij}\times \mathbb{C} \subset V_j\times \mathbb{C} &\to V_{ij}\times \mathbb{C} \subset V_i\times \mathbb{C}\\
                             \left(\xi_j, t_j \right)= \left(\xi_j^1,...,\xi_j^n, t_j \right) &\mapsto \left(\xi_i, f_{ij}\left(\xi_j\right)t_j \right)= \left(\xi_i^1,...,\xi_i^n , t_i \right)
\end{align*}
We define $\tilde{\Lambda}_j \in \Gamma\left( V_j\times \mathbb{C}, \bigwedge^2 \Theta_L \right)$ by
\begin{align*}
\tilde{\Lambda}_j:= \Lambda_j -  t_j T_j \wedge \frac{\partial}{\partial t_j}
\end{align*}
From $t_i= f_{ij}\left(\xi_j\right) t_j$, we have $\frac{\partial}{\partial t_j}= f_{ij}\left(\xi_j\right)\frac{\partial}{\partial t_i}$ and $\frac{\partial}{\partial \xi_j^\alpha}= \sum_{\beta=1}^n \frac{\partial \xi_i^\beta}{\partial \xi_j^\alpha}\frac{\partial}{\partial \xi_i^\beta} + \frac{\partial \log f_{ij}(\xi_j)}{\partial \xi_j^\alpha} t_i \frac{\partial}{\partial t_i}$. Then we have
{\small{\begin{align*}
\tilde{\Lambda}_j&= \Lambda_j - t_j T_j \wedge \frac{\partial}{\partial t_j} = \sum_{\alpha, \beta=1}^n \Lambda_{\alpha\beta}^j (\xi_j)\frac{\partial}{\partial \xi_j^\alpha}\wedge \frac{\partial}{\partial \xi_j^\beta} - t_j \left(T_i- \left[ \Lambda_0, \log f_{ij}  \right] \right) \wedge \left( f_{ij} \frac{\partial}{\partial t_i} \right)  \\
&= \sum_{\delta, \eta,\alpha,\beta=1}^n \Lambda_{\alpha\beta}^j(\xi_j)\frac{\partial \xi_i^\delta}{\partial \xi_j^\alpha}\frac{\partial \xi_i^\eta}{\partial \xi_j^\beta}\frac{\partial}{\partial \xi_i^\delta}\wedge \frac{\partial}{\partial \xi_i^\eta} + \sum_{\alpha, \beta, \delta=1}^n \Lambda_{\alpha\beta}^j(\xi_j) \frac{\partial \xi_i^\delta}{\partial \xi_j^\alpha} \frac{\partial \log f_{ij}}{\partial \xi_j^\beta} t_i  \frac{\partial}{\partial \xi_i^\delta}\wedge \frac{\partial}{\partial t_i} - \sum_{\alpha, \beta, \eta=1}^n \Lambda_{\alpha\beta}^j(\xi_j) \frac{\partial \xi_i^\eta}{\partial \xi_j^\beta} \frac{\partial \log f_{ij}}{\partial \xi_j^\alpha} t_i \frac{\partial}{\partial \xi_i^\eta} \wedge \frac{\partial}{\partial t_i} \\
&- t_i T_i \wedge \frac{\partial}{\partial t_i} + t_i \left[ \Lambda_0, \log f_{ij} \right]\wedge \frac{\partial}{\partial t_i} = \Lambda_i - t_i \left[ \Lambda_0, \log f_{ij} \right] \wedge \frac{\partial}{\partial t_i} - t_i T_i \wedge \frac{\partial}{\partial t_i} + t_i \left[ \Lambda_0, \log f_{ij} \right] \wedge \frac{\partial}{\partial t_i} = \Lambda_i - t_i T_i \wedge \frac{\partial}{\partial t_i} = \tilde{\Lambda}_i
\end{align*}}}
and we have
\begin{align*}
\left[ \tilde{\Lambda}_i, \tilde{\Lambda}_i \right]= \left[\Lambda_i - t_i T_i \wedge \frac{\partial}{\partial t_i}, \Lambda_i - t_i T_i \wedge \frac{\partial}{\partial t_i} \right] = 0
\end{align*}
Then $\left\{ \tilde{\Lambda}_i \right\}$ defines a Poisson structure $\tilde{\Lambda}_0$ on $F$. Let $\pi: F\to W$ be the bundle map. Then $\pi : \left(F, \tilde{\Lambda}_0 \right)\to \left(W, \Lambda_0 \right)$ is a Poisson map.

Let $F^\sharp$ be the $\mathbb{P}_\mathbb{C}^1$-bundle over $W$ obtained from $F$. Let $F^\sharp$ be described by
\begin{center}
$\begin{CD}
\left(\xi_j,t_j\right)\in V_{ij}\times \mathbb{C} \subset V_j\times \mathbb{C}\subset V_j\times \mathbb{P}_\mathbb{C}^1 @>>> \left(\xi_i,t_i\right)= \left(\xi_i,f_{ij}(\xi_j)t_j\right) \in V_{ij} \times \mathbb{C} \subset V_i\times \mathbb{C} \subset V_i \times \mathbb{P}_\mathbb{C}^1\\
@VVV @VVV\\
\left(\xi_j, t_j' \right)=\left(\xi_j, \frac{1}{t_j}\right)\in V_{ij}\times \mathbb{C} \subset V_j\times \mathbb{C} \subset V_j\times \mathbb{P}_\mathbb{C}^1 @>>> \left(\xi_i,t_i' \right)= \left(\xi_i, \frac{1}{f_{ij}(\xi_j)t_j} \right) \in V_{ij}\times \mathbb{C} \subset V_i\times \mathbb{C} \subset V_i\times \mathbb{P}_\mathbb{C}^1
\end{CD}$
\end{center}
We note that since $t_i= \frac{1}{t_i'}$, we have $t_i\frac{\partial}{\partial t_i}= - t_i'\frac{\partial}{\partial t_i'}$. Then we define on $\left(\xi_i, t_i'\right)$
\begin{align*}
\tilde{\Lambda}_i'= \Lambda_i + t_i' T_i \wedge \frac{\partial}{\partial t_i'}
\end{align*}
Then $\tilde{\Lambda}_i= \tilde{\Lambda}_i'$ on $V_i\times \mathbb{C}^*$, so that $\tilde{\Lambda}_i$ is naturally extended to a bivector field on $V_i\times \mathbb{P}_\mathbb{C}^1$, and $\left[ \tilde{\Lambda}_i', \tilde{\Lambda}_i'\right]=0$. Hence $\left\{\tilde{\Lambda}_i, \tilde{\Lambda}_i' \right\}$ defines a Poisson structure $\tilde{\Lambda}_0^\sharp$ on $F^\sharp$. Let $\pi^\sharp: F^\sharp \to W$ be the bundle map. Then $\pi^\sharp: \left( F^\sharp, \tilde{\Lambda}_0^\sharp \right) \to \left( W, \Lambda_0 \right)$ is a Poisson map.

Let $s  \in H^0\left( W, F^k \right)$ represented by $s_i= f_{ij}^k s_j$ such that $\left[ \Lambda_0 , s_i \right] = k s_i T_i$ for each $i$ and $s$ defines a smooth divisor on $W$. We note that $0=\left[\left[ \Lambda_0, s_i\right], s_i \right]= k s_i \left[ T_i, s_i \right] $, so that $\left[T_i, s_i \right]=0$. Then we have
\begin{align*}
\left[ \tilde{\Lambda}_i, t_i^k -s_i \right] =\left[ \Lambda_i - t_i T_i \wedge \frac{\partial}{\partial t_i}, t_i^k - s_i  \right] = - \left[ \Lambda_i, s_i \right] + k t_i^k T_i  + t_i \left[ T_i, s_i \right] \frac{\partial}{\partial t_i} = - ks_i T_i + k t_i^k T_i= k\left(t_i^k - s_i \right) T_i
\end{align*} 
This implies that $\left\{t_i^k - s_i =0 \right\}$ defines a Poisson submanifold of $\left(F, \tilde{\Lambda}_0 \right)$. Moreover, we have
\begin{align*}
\left[ \tilde{\Lambda}_i', t_i'^k s_i -1 \right] = \left[ \Lambda_i + t_i' T_i \wedge \frac{\partial}{\partial t_i'}, t_i'^k s_i -1 \right] = t_i'^k \left[\Lambda_i,  s_i \right] + t_i'^{k+1} \left[ T_i, s_i \right]  \frac{\partial}{\partial t_i} - k t_i'^k s_i T_i = 0
\end{align*}
This implies that $\left\{ t_i^k-s_i=0, t_i'^k -1 =0 \right\}$ is a Poisson submanifold of $\left( F^\sharp , \tilde{\Lambda}_0^\sharp \right)$. 

\begin{theorem}[compare \cite{Wav68} Theorem 1.2 p.929-931]\label{pcc1}
If $\left(M, \Pi_0\right)$ is a Poisson $k$-sheeted cyclic covering of $\left(W, \Lambda_0\right)$, we can find a Poisson line bundle $\left(F, \nabla_F\right)$ on $\left(W, \Lambda_0\right)$ so that, over suitable fine covering $\left\{V_i\right\}$ of $W$, $\left(M, \Pi_0\right)$ may be identified with the Poisson submanifold of $\left(F, \Lambda_\nabla\right)$ defined by the equation $\xi_i^k=\phi_i$, where $\xi_i$ is the fibre coordinate of $F$ over $V_i$ and $\phi_i$ is the local equation of $C$ in $V_i$.
\end{theorem}

\begin{proof}
We extend the argument of \cite{Wav68} Theorem1.2 p.929-931 to which we refer for the detail. We tried to maintain the notational consistency with \cite{Wav68}. We choose a covering $\left\{U_j, U_\lambda\right\}$ of $W$ by polycylinder having the following properties $\left(\pi: M \to W\right)$ is the projection):
\begin{enumerate}
\item $U_j\cap C = \emptyset, U_\lambda \cap C \ne \emptyset$;
\item on $U_\lambda$, $C$ has the equation $\xi_\lambda^n =0$ where $\left(\xi_\lambda^1,..., \xi_\lambda^n\right)$ are local coordinates on $U_\lambda$;
\item $\pi^{-1}\left(U_j\right)= U_j^{(0)} \bigcup \cdots \bigcup U_j^{(k-1)}$ with $U_j^{(v)}\cap U_j^{(\mu)} = \emptyset$ for $v\ne \mu$;
\item $\pi^{-1}\left(U_\lambda\right)= U_\lambda^\sharp$ and coordinates in $U_\lambda^\sharp$ are $(z_\lambda^1,..., z_\lambda^n)$ and the map $\pi$ is given by
\begin{equation*}
\begin{cases}
\xi_\lambda^v= z_\lambda^v, v=1,..., n-1 \\
\xi_\lambda^n= \left(z_\lambda^n\right)^k
\end{cases}
\end{equation*}
\end{enumerate}
Let the coordinate changes on $U_j^{(v)}\cap U_\lambda^{\sharp}$ be 
\begin{align*}
z_\lambda^\alpha = g_{\lambda j}^{(v)\alpha}\left(z_j^{(v)} \right)= g_{\lambda j}^{(v)\alpha} \left( z_j^{(v)1},.., z_j^{(v)n} \right)
\end{align*}
where $z_j^{(v)}= \left(z_j^{(v)1},..., z_j^{(v)n}\right)$ are coordinates on $U_j^{(v)}$. Then we may assume that $U_j^{(v)}$ are numbered so that $g_{\lambda j}^{(v) n}= \rho^v g_{\lambda j}^{(0) n}$ with $\rho= \exp\left( \frac{2\pi i}{k} \right)$, and if $U_j^{(0)}\cap U_k^{(v)}\ne \emptyset$, then $U_j^{(s)} \cap U_k^{(v+s)}$ where by $v+s$ we understand that residue class of $v+s \mod k$. Then we can define a line bundle $F$ on $W$ by the following transition functions:
\begin{equation*}
\begin{cases}
f_{ij}= \rho^{-v} \,\,\,\,\, \textnormal{if} \,\,\,\,\, U_i^{(0)}\cap U_j^{(v)}\ne \emptyset \\
f_{\lambda j}= g_{\lambda j}^{(0) n}\\
f_{\lambda \mu}= \phi_{\lambda \mu}\left( \xi_\mu^1,..., \xi_\mu^{n-1}, \sqrt[k]{\xi_\mu^n}  \right)
\end{cases}
\end{equation*}
where $\phi_{\lambda \mu}\left(z_\mu\right)=\phi_{\lambda\mu}\left(z_\mu^1,..., z_\mu^n\right)$ is a non-vanishing holomorphic function on $U_\lambda^\sharp \cap U_\mu^\sharp\ne \emptyset$ such that $z_\lambda^n= \phi_{\lambda \mu}(z_\mu)\cdot z_\mu^n $ on $U_\lambda^\sharp \cap U_\mu^\sharp \ne \emptyset$ (for the detail, see \cite{Wav68} Lemma p.930).

We will give a Poisson connection $\nabla_F$ on $F$ over $\left(W, \Lambda_0\right)$, so that $\left(L, \nabla_L\right)$ is a Poisson line bundle. Let $\Lambda_0$ on $U_j$ be written as
\begin{align*}
\Lambda_0:=\Lambda_j=\sum_{\alpha,\beta=1}^n \Lambda_{\alpha\beta}^j(\xi_j)\frac{\partial}{\partial \xi_j^\alpha}\wedge \frac{\partial}{\partial \xi_j^\beta}
\end{align*}
On the other hand, since $\xi_\lambda^n=0$ on $U_\lambda$ is a Poisson submanifold of $\left(W, \Lambda_0\right)$, we can write $\Lambda_0$ on $U_\lambda$ as
\begin{align*}
\Lambda_0= \Lambda_\lambda=\sum_{\alpha, \beta=1}^{n-1} \Lambda_{\alpha\beta}^\lambda(\xi_\lambda) \frac{\partial}{\partial \xi_\lambda^\alpha} \wedge \frac{\partial}{\partial \xi_\lambda^\beta} + \sum_{\gamma=1}^{n-1} \xi_\lambda^n G_\lambda^\gamma(\xi_\lambda)\frac{\partial}{\partial \xi_\lambda^\gamma} \wedge \frac{\partial}{\partial \xi_\lambda^n}
\end{align*}
Then $\Pi_0$ on $U_j^{(v)}$ is of the form
\begin{align*}
\Pi_0:=\Pi_j^{(v)}=\sum_{\alpha,\beta=1}^n \Lambda_{\alpha\beta}^j \left( z_j^{(v)} \right)\frac{\partial}{\partial z_j^{(v)\alpha}}\wedge \frac{\partial}{\partial z_j^{(v)\beta}}
\end{align*}
and $\Pi_0$ on $U_\lambda^{\sharp}$ is of the form
\begin{align*}
\Pi_0:=\Pi_\lambda=\sum_{\alpha,\beta=1}^{n-1}  \Lambda_{\alpha\beta}^\lambda\left( z_\lambda^1,..., z_\lambda^{n-1}, \left( z_\lambda^n \right)^k \right) \frac{\partial}{\partial z_\lambda^\alpha} \wedge \frac{\partial}{\partial z_\lambda^\beta} + \frac{1}{k} \sum_{\gamma=1}^{n-1} z_\lambda^n G_\lambda^\gamma\left( z_\lambda^1,..., z_\lambda^{n-1}, \left( z_\lambda^n \right)^k \right) \frac{\partial}{\partial z_\lambda^\gamma} \wedge \frac{\partial}{\partial z_\lambda^n}
\end{align*}

We define $\left\{T_j, T_\lambda \right\} \in C^0\left(\mathcal{U}, \Theta_W \right)$ where $T_j=0 \in \Gamma\left(U_j, \Theta_W\right)$ and $T_\lambda = - \sum_{\gamma=1}^{n-1} \frac{1}{k} G_\lambda^\gamma\left( \xi_\lambda^1,..., \xi_\lambda^n \right)\frac{\partial}{\partial \xi_\lambda^\gamma}\in \Gamma(U_\lambda, \Theta_W)$. We will show that $\left(\left\{ T_j, T_\lambda \right\}, \left\{ f_{ij}, f_{\lambda j}, f_{\lambda \mu} \right\}  \right) \in C^0\left( \mathcal{U}, \Theta_W \right) \bigoplus C^1\left(\mathcal{U}, \mathcal{O}_W^* \right)$ satisies the cocycle condition defining a Poisson line bundle on $\left(W, \Lambda_0\right)$. It is sufficient to show that
\begin{align}
&\left[ \Lambda_0, T_\lambda \right]=0 \label{pcc2} \\
&T_j- T_\lambda  +\left[ \Lambda_0, \log g_{\lambda j}^{(v) n}\right] =0 \label{pcc3} \\
&T_\mu - T_\lambda+ \left[\Lambda_0, \log \phi_{\lambda \mu}\left( \xi_\mu^1,..., \xi_\mu^{n-1}, \sqrt[k]{\xi_\mu^n} \right)   \right]=0 \label{pcc4}
\end{align}
Let us show $\textnormal{(\ref{pcc2})}$. We can write $\Lambda_0= \Lambda'- k\xi_\lambda^n T_\lambda \wedge \frac{\partial}{\partial \xi_\lambda^n}$, where $\Lambda'=\sum_{\alpha,\beta=1}^{n-1} \Lambda_{\alpha\beta}^\lambda (\xi_\lambda )\frac{\partial}{\partial \xi_\lambda^\alpha} \wedge \frac{\partial}{\partial \xi_\lambda^\beta}$. Then
\begin{align*}
\left[ \Lambda_0, T_\lambda \right]= \left[ \Lambda' - k \xi_\lambda^n T_\lambda \wedge \frac{\partial}{\partial \xi_\lambda^n}, T_\lambda \right] = \left[ \Lambda', T_\lambda \right] + k\xi_\lambda^n T_\lambda \wedge \frac{\partial T_\lambda }{\partial \xi_\lambda^n}
\end{align*}

Using $\left[\Lambda_0, \Lambda_0 \right]=0$, we consider the coefficient of $\frac{\partial}{\partial \xi_\lambda^\alpha}\wedge \frac{\partial}{\partial \xi_\lambda^\beta}\wedge \frac{\partial}{\partial \xi_\lambda^n}$. Then
\begin{align*}
0=\left[\Lambda' -  k\xi_\lambda^n T_\lambda \wedge \frac{\partial}{\partial \xi_\lambda^n}, \Lambda' -  k\xi_\lambda^n T_\lambda \wedge \frac{\partial}{\partial \xi_\lambda^n}    \right] = - 2\left[ \Lambda', k\xi_\lambda^n T_\lambda \right] \wedge \frac{\partial}{\partial \xi_\lambda^n} - 2 k^2 \left( \xi_\lambda^n \right)^2 T_\lambda \wedge \frac{\partial T_\lambda}{\partial \xi_\lambda^n} \wedge \frac{\partial}{\partial \xi_\lambda^n} + \cdots, 
\end{align*}
so that we have $2k\xi_\lambda^n\left( \left[ \Lambda',  T_\lambda \right] + k \xi_\lambda^n T_\lambda \wedge \frac{\partial T_\lambda}{\partial \xi_\lambda^n}  \right)=0$ which implies $\textnormal{(\ref{pcc2})}$.

Let us show $\textnormal{(\ref{pcc3})}$. since $z_\lambda^n= g_{\lambda j}^{(0)n}\left(z_j^{(0)}\right)$ on $U_j^{(0)}\cap U_\lambda^\sharp\ne \emptyset$, we have 
\begin{align*}
&\left[ \sum_{\alpha,\beta=1}^{n-1} \Lambda_{\alpha\beta}^j (z_j^{(0)})\frac{\partial}{\partial z_j^{(0)\alpha}}\wedge \frac{\partial}{\partial z_j^{(0)\beta}} + \sum_{\gamma=1}^{n-1} \Lambda_\gamma^j \left(z_j^{(0)}\right)\frac{\partial}{\partial z_j^{(0)\gamma} } \wedge \frac{\partial}{\partial z_j^{(0)n}}, g_{\lambda j}^{(0)n}\left(z_j^{(0)}\right) \right]\\
&= \left[ \sum_{\alpha,\beta=1}^n \Lambda_{\alpha\beta}^\lambda \left(z_\lambda^1,...,z_\lambda^{n-1}, \left(z_\lambda^n\right)^k\right) \frac{\partial}{\partial z_\lambda^\alpha}\wedge \frac{\partial}{\partial z_\lambda^\beta} + \sum_{\gamma=1}^{n-1} \frac{1}{k} z_\lambda^n G_\lambda^\gamma \left(z_\lambda^1,...,z_\lambda^{n-1}, \left(z_\lambda^n \right)^k \right) \frac{\partial}{\partial z_\lambda^\gamma} \wedge \frac{\partial}{\partial z_\lambda^n} , z_\lambda^n \right]\\
&=-g_{\lambda j}^{(0)n} \left(z_j^{(0)} \right) \sum_{\gamma=1}^{n-1} \frac{1}{k} G_\lambda^\gamma \left(z_\lambda^1,...,z_\lambda^{n-1}, \left(z_\lambda^n\right)^k \right)\frac{\partial}{\partial z_\lambda^\gamma}
\end{align*}
By taking the inverse of $\pi:\left(\xi_j^1,..., \xi_j^n \right)\mapsto \left(z_j^{(0)1},..., z_j^{(0)n} \right)$ and $\left( \xi_\lambda^1,..., \xi_\lambda^n \right)\mapsto \left(z_\lambda^1,..., z_\lambda^n \right)=\left( \xi_\lambda^1,..., \xi_\lambda^{n-1}, \sqrt[k]{\xi_\lambda^n}  \right)$ on $U_\lambda \cap U_j$ which is well-defined, we have $\left[ \Lambda_0, \log g_{\lambda j}^{(0)n} \right]= - \sum_{\gamma=1}^{n-1} \frac{1}{k} G_\lambda^\gamma\left(\xi_\lambda^1,..., \xi_\lambda^n \right)\frac{\partial}{\partial \xi_\lambda^\gamma}$ which implies $\textnormal{(\ref{pcc3})}$.

Let us show $\textnormal{(\ref{pcc4})}$. On $U_\lambda^\sharp \cap U_\mu^\sharp\ne \emptyset$, we have $z_\lambda^n=\phi_{\lambda\mu}\left(z_\mu\right) z_\mu^n$, so that we have $\frac{\partial}{\partial z_\mu^\alpha}= \sum_{\beta=1}^{n-1}\frac{\partial z_\lambda^\beta}{\partial z_\mu^\alpha} \frac{\partial}{\partial z_\lambda^\beta} + \frac{\partial \phi_{\lambda \mu}}{\partial z_\mu^\alpha} z_\mu^n \frac{\partial}{\partial z_\lambda^n}$ for $\alpha=1,..., n-1$, and $\frac{\partial}{\partial z_\mu^n}=\sum_{\alpha=1}^{n-1} \frac{\partial z_\lambda^\alpha}{\partial z_\mu^n}\frac{\partial}{\partial z_\lambda^\alpha}+\left( \frac{\partial \phi_{\lambda \mu}}{\partial z_\mu^n} z_\mu^n + \phi_{\lambda \mu} \right)\frac{\partial}{\partial z_\lambda^n}$. Then we have
\begin{align*}
&\Pi_0= \sum_{\alpha,\beta=1}^{n-1} \Lambda_{\alpha\beta}^\lambda \left(z_\lambda^1,...,z_\lambda^{n-1}, \left(z_\lambda^n \right)^k \right) \frac{\partial}{\partial z_\lambda^\alpha}\wedge \frac{\partial}{\partial z_\lambda^\beta} + \sum_{\gamma=1}^{n-1} \frac{1}{k} z_\lambda^n G_\lambda^\gamma \left(z_\lambda^1,...,z_\lambda^{n-1}, \left(z_\lambda^n\right)^k\right) \frac{\partial}{\partial z_\lambda^\gamma} \wedge \frac{\partial}{\partial z_\lambda^n}\\
&= \sum_{\alpha,\beta=1}^{n-1} \Lambda_{\alpha\beta}^\mu \left(z_\mu^1,...,z_\mu^{n-1}, \left(z_\mu^n\right)^k \right) \frac{\partial}{\partial z_\mu^\alpha}\wedge \frac{\partial}{\partial z_\mu^\beta} + \sum_{\gamma=1}^{n-1} \frac{1}{k} z_\mu^n G_\mu^\gamma\left(z_\mu^1,...,z_\mu^{n-1}, \left(z_\mu^n\right)^k\right) \frac{\partial}{\partial z_\mu^\gamma} \wedge \frac{\partial}{\partial z_\mu^n}\\
&=z_\mu^n \sum_{\alpha,\beta ,\gamma =1}^{n-1} \Lambda_{\alpha\beta}^\mu\left(z_\mu^1,...,z_\mu^{n-1}, \left(z_\mu^n \right)^k \right)\frac{\partial \phi_{\lambda \mu}}{\partial z_\mu^\alpha}\frac{\partial z_\lambda^\gamma}{\partial z_\mu^\beta} \frac{\partial}{\partial z_\lambda^n}\wedge\frac{\partial}{\partial z_\lambda^\gamma}+ z_\mu^n  \sum_{\alpha,\beta , \gamma =1}^{n-1} \Lambda_{\alpha\beta}^\mu\left(z_\mu^1,...,z_\mu^{n-1}, \left(z_\mu^n \right)^k \right) \frac{\partial z_\lambda^\gamma}{\partial z_\mu^\alpha}\frac{\partial \phi_{\lambda \mu}}{\partial z_\mu^\beta}\frac{\partial}{\partial z_\lambda^\gamma}\wedge\frac{\partial}{\partial z_\lambda^n}\\
&+ \sum_{\gamma,\alpha=1}^{n-1} \frac{1}{k} z_\mu^n G_\mu^\alpha \left(z_\mu^1,...,z_\mu^{n-1}, \left(z_\mu^n\right)^k \right)\frac{\partial z_\lambda^\gamma}{\partial z_\mu^\alpha} \left( \frac{\partial \phi_{\lambda \mu}}{\partial z_\mu^n} z_\mu^n+\phi_{\lambda \mu}  \right) \frac{\partial}{\partial z_\lambda^\gamma} \wedge \frac{\partial}{\partial z_\lambda^n}\\
& + \sum_{\gamma,\alpha=1}^{n-1} \frac{1}{k} z_\mu^n G_\mu^\alpha \left(z_\mu^1,...,z_\mu^{n-1}, \left(z_\mu^n\right)^k \right) \frac{\partial \phi_{\lambda\mu}}{\partial z_\mu^\alpha}z_\mu^n\frac{\partial z_\lambda^\gamma}{\partial z_\mu^n} \frac{\partial}{\partial z_\lambda^n} \wedge \frac{\partial}{\partial z_\lambda^\gamma}+ \cdots
\end{align*}
By considering the coefficient of $\frac{\partial}{\partial z_\lambda^\gamma}\wedge \frac{\partial}{\partial z_\lambda^n}$ for $\gamma=1,..., n-1$, we have
\begin{align}\label{pcc6}
&\frac{1}{k} z_\lambda^n G_\lambda^\gamma \left(z_\lambda^1,..., z_\lambda^{n-1}, \left(z_\lambda^n \right)^k \right)\\
&= - z_\mu^n \sum_{\alpha,\beta=1}^{n-1} \Lambda_{\alpha\beta}^{\mu }\left(z_\mu^1,..., z_\mu^{n-1}, \left(z_\mu^n \right)^k \right)\frac{\partial \phi_{\lambda \mu}}{\partial z_\mu^\alpha} \frac{\partial z_\lambda^\gamma}{\partial z_\mu^\beta} + z_\mu^n \sum_{\alpha,\beta=1}^{n-1} \Lambda_{\alpha\beta}^\mu \left(z_\mu^1,...,z_\mu^{n-1}, \left(z_\mu^n \right)^k \right) \frac{\partial z_\lambda^\gamma}{\partial z_\mu^\alpha} \frac{\partial \phi_{\lambda \mu}}{\partial z_\mu^\beta} \notag \\
&+\sum_{\alpha=1}^{n-1} \frac{1}{k} z_\mu^n G_\mu^\alpha\left( z_\mu^1,...,z_\mu^{n-1}, \left(z_\mu^n\right)^k \right) \frac{\partial z_\lambda^\gamma}{\partial z_\mu^\alpha}\frac{\partial \phi_{\lambda\mu}}{\partial z_\mu^n} z_\mu^n + \sum_{\alpha=1}^{n-1} \frac{1}{k} z_\lambda^n G_\mu^\alpha\left(z_\mu^1,..., z_\mu^{n-1}, \left(z_\mu^n\right)^k \right) \frac{\partial z_\lambda^\gamma}{\partial z_\mu^\alpha} \notag \\
&- \sum_{\alpha=1}^{n-1} \frac{1}{k} z_\mu^n G_\mu^\alpha \left(z_\mu^1,...,z_\mu^{n-1}, \left(z_\mu^n\right)^k \right) \frac{\partial \phi_{\lambda\mu}}{\partial z_\mu^\alpha}z_\mu^n\frac{\partial z_\lambda^\gamma}{\partial z_\mu^n} \notag
\end{align}
By dividing $\textnormal{(\ref{pcc6})}$ by $z_\lambda^n$, we obtain
\begin{align}\label{pcc7}
&\frac{1}{k}  G_\lambda^\gamma \left(z_\lambda^1,..., z_\lambda^{n-1}, \left(z_\lambda^n \right)^k \right)\\
&= - \frac{1}{\phi_{\lambda \mu}} \sum_{\alpha,\beta=1}^{n-1} \Lambda_{\alpha\beta}^{\mu }\left(z_\mu^1,..., z_\mu^{n-1}, \left(z_\mu^n \right)^k \right)\frac{\partial \phi_{\lambda \mu}}{\partial z_\mu^\alpha} \frac{\partial z_\lambda^\gamma}{\partial z_\mu^\beta} + \frac{1}{\phi_{\lambda\mu}} \sum_{\alpha,\beta=1}^{n-1} \Lambda_{\alpha\beta}^\mu \left(z_\mu^1,...,z_\mu^{n-1}, \left(z_\mu^n \right)^k \right) \frac{\partial z_\lambda^\gamma}{\partial z_\mu^\alpha} \frac{\partial \phi_{\lambda \mu}}{\partial z_\mu^\beta} \notag \\
&+\sum_{\alpha=1}^{n-1} \frac{1}{k}  G_\mu^\alpha\left( z_\mu^1,...,z_\mu^{n-1}, \left(z_\mu^n\right)^k \right) \frac{\partial z_\lambda^\gamma}{\partial z_\mu^\alpha}\frac{1}{\phi_{\lambda\mu}}\frac{\partial \phi_{\lambda\mu}}{\partial z_\mu^n} z_\mu^n + \sum_{\alpha=1}^{n-1} \frac{1}{k}  G_\mu^\alpha\left(z_\mu^1,..., z_\mu^{n-1}, \left(z_\mu^n\right)^k \right) \frac{\partial z_\lambda^\gamma}{\partial z_\mu^\alpha} \notag \\
&- \sum_{\alpha=1}^{n-1} \frac{1}{k}  G_\mu^\alpha \left(z_\mu^1,...,z_\mu^{n-1}, \left(z_\mu^n\right)^k \right)\frac{1}{\phi_{\lambda\mu}} \frac{\partial \phi_{\lambda\mu}}{\partial z_\mu^\alpha}z_\mu^n\frac{\partial z_\lambda^\gamma}{\partial z_\mu^n} \notag
\end{align}

On the other hand, we note that $\phi_{\lambda_\mu}\left(z_\mu^1,..., z_\mu^n\right)= \phi_{\lambda \mu}\left(z_\mu^1,..., z_{\mu}^{n-1}, \rho^s \cdot z_\mu^n\right)$ for $s=1,...,{k-1}$ (see \cite{Wav68} p.930). By taking the derivative with respect to $z_\mu^n$, we have
\begin{align*}
&\frac{\partial \phi_{\lambda \mu}}{\partial z_\mu^n}\left(z_\mu^1,..., z_\mu^n\right)= \frac{\partial \phi_{\lambda \mu}}{\partial z_\mu^n}\left( z_\mu^1,..., z_\mu^n, \rho^s \cdot z_\mu^n \right)\cdot \rho^s\\
&\Longrightarrow  z_\mu^n\cdot \frac{\partial \phi_{\lambda \mu}}{\partial z_\mu^n}\left(z_\mu^1,..., z_\mu^n\right)= \frac{\partial \phi_{\lambda \mu}}{\partial z_\mu^n}\left( z_\mu^1,..., z_\mu^n, \rho^s \cdot z_\mu^n \right)\cdot  \rho^s z_\mu^n
\end{align*}
This implies that $\sqrt[k]{\xi_\mu^n}\cdot \frac{\partial \phi_{\lambda \mu}}{\partial z_\mu^n} \left(\xi_\mu^1,..., \xi_\mu^{n-1}, \sqrt[k]{\xi_\mu^n} \right)$ is a well-defined holomorphic function on $U_\lambda \cap U_\mu$. From $z_\lambda^n= \phi_{\lambda \mu}\left(z_\mu\right) z_\mu^n$, we have $\frac{\partial z_\lambda^n}{\partial z_\mu^n}= \phi_{\lambda \mu}(z_\mu)+ \frac{\partial \phi_{\lambda \mu}}{\partial z_\mu^n}\cdot z_\mu^n$, we see that $\frac{\partial z_\lambda^n}{\partial z_\mu^n}\left(\xi_\mu^1,..., \xi_\mu^{n-1}, \sqrt[k]{\xi_\mu^n}\right)$ is a well-defined holomorphic function. We note that
\begin{align*}
\frac{\partial \xi_\lambda^n}{\partial z_\mu^n}&= \sum_{\alpha=1}^{n-1} \frac{\partial \xi_\lambda^n}{\partial z_\lambda^\alpha}\cdot \frac{\partial z_\lambda^\alpha}{\partial z_\mu^n} + \frac{\partial \xi_\lambda^n}{ \partial z_\lambda^n}\frac{\partial z_\lambda^n}{\partial z_\mu^n} = k \left(z_\lambda^n\right)^{k-1} \frac{\partial z_\lambda^n}{\partial z_\mu^n} \\
\frac{\partial \xi_\lambda^n}{\partial z_\mu^n}&= \sum_{\alpha=1}^{n-1} \frac{\partial \xi_\lambda^n}{\partial \xi_\mu^n}(z_\mu^1,..., z_\mu^{n-1}, (z_\mu^n)^k) \frac{\partial \xi_\mu^\alpha}{\partial z_\mu^n} + \frac{\partial \xi_\lambda^n}{\partial \xi_\mu^n} \left(z_\mu^1,..., z_\mu^{n-1}, \left(z_\mu^n\right)^k\right)\frac{\partial \xi_\mu^n}{\partial z_\mu^n}= \frac{\partial \xi_\lambda^n}{\partial \xi_\mu^n} \left(z_\mu^1,..., z_\mu^{n-1}, \left(z_\mu^n\right)^k\right) k\left(z_\mu^n\right)^{k-1}
\end{align*}
This implies that
\begin{align*}
\frac{\partial \xi_\lambda^n}{\partial \xi_\mu^n}= \phi_{\lambda \mu}^{k-1}\left(\xi_\mu^1,..., \xi_\mu^{n-1}, \sqrt[k]{\xi_\mu^n} \right) \cdot \frac{\partial z_\lambda^n}{\partial z_\mu^n}\left(\xi_\mu^1,..., \xi_\mu^{n-1}, \sqrt[k]{\xi_\mu^n}\right)
\end{align*}

On the other hand, we note that from $z_\lambda^n=\phi_{\lambda \mu}\left(z_\mu\right)\cdot z_\mu$, we have $\left(z_\lambda^n\right)^k = \left(\phi_{\lambda\mu}\left(z_\mu\right) \right)^k \cdot  \left( z_\mu^n \right)^k$. By plugging $\left(\xi_\lambda^1,..., \xi_\lambda^{n-1}, \left(\xi_\lambda^n \right)^{\frac{1}{k}} \right)$ into $\left(z_\lambda^1,..., z_\lambda^n \right)$ and $\left(\xi_\mu^1,..., \xi_\mu^{n-1}, \left(\xi_\mu^n \right)^{\frac{1}{k}} \right)$ into $\left(z_\mu^1,..., z_\mu^n \right)$, we have $\xi_\lambda^n= \left( \phi_{\lambda \mu}\left(\xi_\mu^1,..., \xi_\mu^{n-1}, \sqrt[k]{\xi_\mu^n} \right) \right)^k \cdot \xi_\mu^n$ on $U_\lambda \cap U_\mu$. We set  $\Phi_{\lambda \mu}(\xi_\mu)= \phi_{\lambda \mu}\left(\xi_\mu^1,..., \xi_\mu^{n-1}, \sqrt[k]{\xi_\mu^n} \right) $. Then we have 
\begin{align*}
\frac{\partial \xi_\lambda^n}{\partial \xi_\mu^n}&= \frac{\partial \left(\Phi_{\lambda\mu}^k \xi_\mu^n \right)}{\partial \xi_\mu^n}= k \Phi_{\lambda \mu}^{k-1} \frac{\partial \Phi_{\lambda\mu}}{\partial \xi_\mu^n} \xi_\mu^n + \Phi_{\lambda\mu}^k\\
\phi_{\lambda\mu}^{k-1} \frac{\partial z_\lambda^n}{\partial z_\mu^n}&= \phi_{\lambda \mu}^{k-1} \frac{\partial \left(\phi_{\lambda \mu} z_\mu^n \right)}{\partial z_\mu^n}= \phi_{\lambda\mu}^{k-1} \frac{\partial \phi_{\lambda \mu}}{\partial z_\mu^n} z_\mu^n + \phi_{\lambda \mu}^k
\end{align*}
This implies that
\begin{align*}
\sqrt[k]{\xi_\mu^n}\cdot \frac{\partial \phi_{\lambda\mu}}{\partial z_\mu^n}\left(\xi_\mu^1,...,\xi_\mu^{n-1}, \sqrt[k]{\xi_\mu^n} \right)= k \xi_\mu^n\frac{\partial \phi_{\lambda \mu}\left(\xi_\mu^1,...,\xi_\mu^{n-1}, \sqrt[k]{\xi_\mu^n}\right) }{\partial \xi_\mu^n } = k \xi_\mu^n \frac{\partial \Phi_{\lambda\mu}}{\partial \xi_\mu^n}
\end{align*}

We note that
\begin{align*}
\frac{\partial}{\partial \xi_\mu^\gamma}&= \sum_{\alpha=1}^{n-1}\frac{\partial \xi_\lambda^\alpha}{\partial \xi_\mu^\gamma}\frac{\partial}{\partial \xi_\lambda^\alpha}+ k\Phi_{\lambda\mu}^{k-1}\frac{\partial \Phi_{\lambda\mu}}{\partial \xi_\mu^\gamma} \xi_\mu^n\frac{\partial}{\partial \xi_\lambda^n} = \sum_{\alpha=1}^{n-1}\frac{\partial \xi_\mu^\alpha}{\partial \xi_\mu^\gamma}\frac{\partial}{\partial \xi_\lambda^\alpha} + k\frac{\partial \log \Phi_{\lambda\mu}}{\partial \xi_\mu^\gamma} \xi_\lambda^n \frac{\partial}{\partial \xi_\lambda^n},\,\,\,\,\,\,\gamma=1,..., n-1\\
\frac{\partial}{\partial \xi_\mu^n}&= \sum_{\alpha=1}^{n-1}\frac{\partial \xi_\lambda^\alpha}{\partial \xi_\mu^n}\frac{\partial}{\partial \xi_\lambda^\alpha}+ \left(\Phi_{\lambda\mu}\right)^k\frac{\partial}{\partial \xi_\lambda^n} + k\frac{\partial \log \Phi_{\lambda \mu}}{\partial \xi_\mu^n} \xi_\lambda^n \frac{\partial}{\partial \xi_\lambda^n}
\end{align*}
and we note that
\begin{align*}
\frac{\partial \xi_\lambda^\gamma}{\partial z_\mu^n}&= \sum_{\alpha=1}^n \frac{\partial \xi_\lambda^\gamma}{\partial z_\lambda^\alpha}\frac{\partial z_\lambda^\alpha}{\partial z_\mu^n} = \frac{\partial z_\lambda^\gamma}{\partial z_\mu^n}\\
\frac{\partial \xi_\lambda^\gamma}{\partial z_\mu^n} & = \sum_{\alpha=1}^n \frac{\partial \xi_\lambda^\gamma}{\partial \xi_\mu^\alpha} \frac{\partial \xi_\mu^\alpha}{\partial z_\mu^n} = \frac{\partial \xi_\lambda^\gamma}{\partial \xi_\mu^n} \frac{\partial \xi_\mu^n}{\partial z_\mu^n}= \frac{\partial \xi_\lambda^\gamma}{\partial \xi_\mu^n} \left(z_\mu^1,..., z_\mu^{n-1}, \left(z_\mu^n\right)^k \right) \cdot k \left(z_\mu^n \right)^{k-1}
\end{align*}
Then
\begin{align*}
z_\mu^n \frac{\partial z_\lambda^\gamma}{\partial z_\mu^n} = \frac{\partial \xi_\lambda^\gamma}{\partial \xi_\mu^n} \left(z_\mu^1,..., z_\mu^{n-1}, \left(z_\mu^n\right)^k \right) \cdot k \left(z_\mu^n \right)^k
\end{align*}
By plugging $\left(\xi_\mu^1,..., \xi_\mu^{n-1} , \sqrt[k]{\xi_\mu^n}\right)$ into $\left(z_\mu^1,..., z_\mu^n \right)$, we see that $\sqrt[k]{\xi_\mu^n}\frac{\partial z_\lambda^\gamma}{\partial z_\mu^n}\left(\xi_\mu^1,..., \xi_\mu^{n-1}, \sqrt[k]{\xi_\mu^n}\right)$ is a well-defined holomorphic function on $U_\lambda \cap U_\mu$, and
\begin{align*}
\sqrt[k]{\xi_\mu^n}\frac{\partial z_\lambda^\gamma}{\partial z_\mu^n} \left(\xi_\mu^1,..., \xi_\mu^{n-1}, \sqrt[k]{\xi_\mu^n}\right) = k \xi_\mu^n \frac{\partial \xi_\lambda^\gamma}{\partial \xi_\mu^n}
\end{align*}

Then by plugging $\left(\xi_\lambda^1,..., \xi_\lambda^{n-1}, \sqrt[k]{\xi_\mu^n} \right)$ into $\left(z_\lambda^1,..., z_\lambda^n \right)$ and $\left(\xi_\mu^1,..., \xi_\mu^{n-1},  \sqrt[k]{\xi_\mu^n} \right)$ into $\left(z_\mu^1,..., z_\mu^n \right)$ in $\textnormal{(\ref{pcc7})}$, we have
\begin{align*}
&\frac{1}{k} G_\lambda^\gamma\left(\xi_\lambda^1,..., \xi_\lambda^n \right)\\
&= - \frac{1}{\Phi_{\lambda \mu}} \sum_{\alpha,\beta=1}^{n-1} \Lambda_{\alpha\beta}^{\mu }\left(\xi_\mu^1,...,  \xi_\mu^n \right)\frac{\partial \Phi_{\lambda \mu}}{\partial \xi_\mu^\alpha}\frac{\partial \xi_\lambda^\gamma}{\partial \xi_\mu^\beta} + \frac{1}{\Phi_{\lambda\mu}} \sum_{\alpha,\beta=1}^{n-1} \Lambda_{\alpha\beta}^\mu \left(\xi_\mu^1,..., \xi_\mu^n \right) \frac{\partial \xi_\lambda^\gamma}{\partial \xi_\mu^\alpha} \frac{\partial \Phi_{\lambda \mu}}{\partial \xi_\mu^\beta} \\
&+\sum_{\alpha=1}^{n-1} \frac{1}{k}  G_\mu^\alpha\left( \xi_\mu^1,..., \xi_\mu^n \right) \frac{\partial \xi_\lambda^\gamma}{\partial \xi_\mu^\alpha}\frac{1}{\Phi_{\lambda\mu}} k \xi_\mu^n \frac{\partial \Phi_{\lambda\mu}}{\partial \xi_\mu^n} + \sum_{\alpha=1}^{n-1} \frac{1}{k}  G_\mu^\alpha\left(\xi_\mu^1,..., \xi_\mu^n  \right) \frac{\partial \xi_\lambda^\gamma}{\partial \xi_\mu^\alpha}\\
&- \sum_{\alpha=1}^{n-1} \frac{1}{k}  G_\mu^\alpha \left(\xi_\mu^1,...,\xi_\mu^n \right)\frac{1}{\Phi_{\lambda\mu}} \frac{\partial \Phi_{\lambda\mu}}{\partial \xi_\mu^\alpha} k \xi_\mu^n \frac{\partial \xi_\lambda^\gamma}{\partial \xi_\mu^n}
\end{align*}
Then we have
{\small{\begin{align*}
&\sum_{\gamma=1}^{n-1}\frac{1}{k} G_\lambda^\gamma\left(\xi_\lambda^1,..., \xi_\lambda^n \right) \frac{\partial}{\partial \xi_\lambda^\gamma}\\
&= -  \sum_{\alpha,\beta,\gamma=1}^{n-1} \Lambda_{\alpha\beta}^{\mu }\left(\xi_\mu^1,...,  \xi_\mu^n \right)\frac{\partial \log \Phi_{\lambda \mu}}{\partial \xi_\mu^\alpha}\frac{\partial \xi_\lambda^\gamma}{\partial \xi_\mu^\beta}\frac{\partial}{\partial \xi_\lambda^\gamma} + \sum_{\alpha,\beta,\gamma=1}^{n-1} \Lambda_{\alpha\beta}^\mu \left(\xi_\mu^1,..., \xi_\mu^n \right)  \frac{\partial \log \Phi_{\lambda \mu}}{\partial \xi_\mu^\beta} \frac{\partial \xi_\lambda^\gamma}{\partial \xi_\mu^\alpha}\frac{\partial}{\partial \xi_\lambda^\gamma} \\
&+\sum_{\alpha,\gamma=1}^{n-1}  G_\mu^\alpha\left( \xi_\mu^1,..., \xi_\mu^n \right)  \xi_\mu^n \frac{\partial \log \Phi_{\lambda\mu}}{\partial \xi_\mu^n} \frac{\partial \xi_\lambda^\gamma}{\partial \xi_\mu^\alpha}\frac{\partial}{\partial \xi_\lambda^\gamma} + \sum_{\alpha,\gamma=1}^{n-1} \frac{1}{k}  G_\mu^\alpha\left(\xi_\mu^1,..., \xi_\mu^n  \right) \frac{\partial \xi_\lambda^\gamma}{\partial \xi_\mu^\alpha} \frac{\partial}{\partial \xi_\lambda^\gamma}- \sum_{\alpha,\gamma=1}^{n-1}   G_\mu^\alpha \left(\xi_\mu^1,...,\xi_\mu^n \right) \frac{\partial \log \Phi_{\lambda\mu}}{\partial \xi_\mu^\alpha}  \xi_\mu^n \frac{\partial \xi_\lambda^\gamma}{\partial \xi_\mu^n} \frac{\partial}{\partial \xi_\lambda^\gamma}\\
&= -  \sum_{\alpha,\beta=1}^{n-1} \Lambda_{\alpha\beta}^{\mu }\left(\xi_\mu^1,...,  \xi_\mu^n \right)\frac{\partial \log \Phi_{\lambda \mu}}{\partial \xi_\mu^\alpha} \left( \frac{\partial }{\partial \xi_\mu^\beta} - k \frac{\partial \log \Phi_{\lambda\mu}}{\partial \xi_\mu^\beta} \xi_\lambda^n \frac{\partial}{\partial \xi_\lambda^n} \right) + \sum_{\alpha,\beta=1}^{n-1} \Lambda_{\alpha\beta}^\mu \left(\xi_\mu^1,..., \xi_\mu^n \right)  \frac{\partial \log \Phi_{\lambda \mu}}{\partial \xi_\mu^\beta}\left( \frac{\partial}{\partial \xi_\mu^\alpha}- k \frac{\partial \log \Phi_{\lambda\mu}}{\partial \xi_\mu^\alpha} \xi_\lambda^n \frac{\partial}{\partial \xi_\lambda^n}    \right)\\
&+\sum_{\alpha=1}^{n-1}  G_\mu^\alpha\left( \xi_\mu^1,..., \xi_\mu^n \right) \xi_\mu^n \frac{\partial \log \Phi_{\lambda\mu}}{\partial \xi_\mu^n} \left( \frac{\partial}{\partial \xi_\mu^\alpha}- k \frac{\partial \log \Phi_{\lambda\mu}}{\partial \xi_\mu^\alpha} \xi_\lambda^n \frac{\partial}{\partial \xi_\lambda^n}    \right) + \sum_{\alpha=1}^{n-1} \frac{1}{k}  G_\mu^\alpha\left(\xi_\mu^1,..., \xi_\mu^n  \right) \left( \frac{\partial}{\partial \xi_\mu^\alpha}- k \frac{\partial \log \Phi_{\lambda\mu}}{\partial \xi_\mu^\alpha} \xi_\lambda^n \frac{\partial}{\partial \xi_\lambda^n}    \right)\\
&- \sum_{\alpha=1}^{n-1}   G_\mu^\alpha \left(\xi_\mu^1,...,\xi_\mu^n \right)\frac{\partial \log \Phi_{\lambda\mu}}{\partial \xi_\mu^\alpha}  \xi_\mu^n \left( \frac{\partial}{\partial \xi_\mu^n} - \left(\Phi_{\lambda \mu} \right)^k \frac{\partial}{\partial \xi_\lambda^n} - k \frac{\partial \log \Phi_{\lambda \mu}}{\partial \xi_\mu^n} \xi_\lambda^n \frac{\partial}{\partial \xi_\lambda^n}  \right) \\
&=  -  \sum_{\alpha,\beta=1}^{n-1} \Lambda_{\alpha\beta}^{\mu }\left(\xi_\mu^1,...,  \xi_\mu^n \right)\frac{\partial \log \Phi_{\lambda \mu}}{\partial \xi_\mu^\alpha}  \frac{\partial }{\partial \xi_\mu^\beta}  + \sum_{\alpha,\beta=1}^{n-1} \Lambda_{\alpha\beta}^\mu \left(\xi_\mu^1,..., \xi_\mu^n \right)  \frac{\partial \log \Phi_{\lambda \mu}}{\partial \xi_\mu^\beta} \frac{\partial}{\partial \xi_\mu^\alpha} \\
& + \sum_{\alpha=1}^{n-1} \xi_\mu^n G_\mu^\alpha\left(\xi_\mu^1,..., \xi_\mu^n\right)\frac{\partial \log \Phi_{\lambda \mu}}{\partial \xi_\mu^n}\frac{\partial}{\partial \xi_\mu^\alpha} - \sum_{\alpha=1}^{n-1} \xi_\mu^n G_\mu^\alpha\left( \xi_\mu^1,..., \xi_\mu^n  \right) \frac{\partial \log \Phi_{\lambda \mu}}{\partial \xi_\mu^\alpha}\frac{\partial}{\partial \xi_\mu^n}  + \sum_{\alpha=1}^{n-1} \frac{1}{k} G_\mu^\alpha \left( \xi_\mu^1,..., \xi_\mu^n \right) \frac{\partial}{\partial \xi_\mu^\alpha} \\
&=- \left[\Lambda_0, \log \Phi_{\lambda \mu} \right] + \sum_{\alpha=1}^{n-1} \frac{1}{k} G_\mu^\alpha \left( \xi_\mu^1,..., \xi_\mu^n \right) \frac{\partial}{\partial \xi_\mu^\alpha} \iff - T_\lambda = - \left[ \Lambda_0, \log \Phi_{\lambda \mu} \right] - T_\mu
\end{align*}}}
which shows $\textnormal{(\ref{pcc4})}$. Hence there is a Poisson connection $\nabla_F$ on $F$ over $\left(W, \Lambda_0\right)$, so that we can make $\left(F, \tilde{\Lambda}_0\right)$ as holomorphic Poisson manifold. We recall that $\Lambda_{\nabla_F}$ is of the form
\begin{align*}
\tilde{\Lambda}_0&=\tilde{\Lambda}_j= \Lambda_j\,\,\,\,\,\,\,\textnormal{on}\,\,\, U_j\times \mathbb{C}\\
\tilde{\Lambda}_0&=\tilde{\Lambda}_\lambda = \Lambda_\lambda - t_\lambda T_\lambda \wedge \frac{\partial}{\partial t_\lambda}\,\,\,\,\,\,\,\,\textnormal{on}\,\,\,U_\lambda \times \mathbb{C}
\end{align*}

Now we define a holomorphic map $\Psi:M\to F$ by
\begin{align*}
\Psi: U_j^{(v)} &\to U_j\times \mathbb{C}\\
\left(z_j^{(v)} \right) &\mapsto \left( z_j^{(v)} ; \rho^v \right)\\
\Psi:U_\lambda^\sharp &\to U_\lambda \times \mathbb{C}\\
\left(z_\lambda^1,..., z_\lambda^n\right) &\mapsto \left( z_\lambda^1,..., z_\lambda^{n-1} , \left(z_\lambda^n \right)^k ; z_\lambda^n \right) =\left( \xi_\lambda^1,..., \xi_\lambda^{n-1}, \xi_\lambda^n ; t_\lambda \right)
\end{align*}
Then we may identify $M$ with the image. The image has the equation $t_j^k=1$ over $U_j\times \mathbb{C}$ and $t_\lambda^k=\xi_\lambda^n$ over $U_\lambda \times \mathbb{C}$. We calim that the image is a Poisson submanifold of $(F, \tilde{\Lambda}_0)$. In fact, 
\begin{align*}
\left[ \tilde{\Lambda}_0, t_j^k-1 \right]&=\left[\Lambda_j, t_j^k-1 \right] =0\\
\left[ \tilde{\Lambda}_0, t_\lambda^k - \xi_\lambda^n \right]&=\left[ \sum_{\alpha,\beta=1}^{n-1} \Lambda_{\alpha\beta}^\lambda \left( \xi_\lambda\right)\frac{\partial}{\partial \xi_\lambda^\alpha}\wedge \frac{\partial}{\partial \xi_\lambda^\beta} + \sum_{\gamma=1}^{n-1} \xi_\lambda^n G_\lambda^\gamma \left(\xi_\lambda\right)\frac{\partial}{\partial \xi_\lambda^\gamma}\wedge \frac{\partial}{\partial \xi_\lambda^n} + t_\lambda \sum_{\gamma=1}^{n-1}\frac{1}{k} G_\lambda^\gamma \left(\xi_\lambda\right)\frac{\partial}{\partial \xi_\lambda^\gamma}\wedge \frac{\partial}{\partial t_\lambda} , t_\lambda^k - \xi_\lambda^n \right]\\
&= \sum_{\gamma=1}^{n-1} \xi_\lambda^n G_\lambda^\gamma \left(\xi_\lambda\right)\frac{\partial}{\partial \xi_\lambda^\gamma} - t_\lambda^k \sum_{\gamma=1}^{n-1} G_\lambda^\gamma \left(\xi_\lambda\right) \frac{\partial}{\partial \xi_\lambda^\gamma} =  \left( \xi_\lambda^n - t_\lambda^k \right) \sum_{\gamma=1}^{n-1} G_\lambda^\gamma \left(\xi_\lambda\right) \frac{\partial}{\partial \xi_\lambda^\gamma}
\end{align*} 

Next we will show that $\Psi$ induces a Poisson isomorphism between $\left(M, \Pi_0\right)$ and the image with the Poisson structure induced from $\left(F, \tilde{\Lambda}_0\right)$. It is sufficient to show that $\Psi$ is a Poisson map. In fact,
{\small{\begin{align*}
&\Pi_\lambda= \sum_{\alpha,\beta=1}^{n-1}\Lambda_{\alpha\beta}^\lambda \left(z_\lambda^1,..., z_\lambda^{n-1}, \left(  z_\lambda^n \right)^k \right) \frac{\partial}{\partial z_\lambda^\alpha}\wedge \frac{\partial}{\partial z_\lambda^\beta} + \frac{1}{k}\sum_{\gamma=1}^{n-1} z_\lambda^n G_\lambda ^\gamma \left(z_\lambda^1,..., z_\lambda^{n-1}, \left( z_\lambda^n \right)^k \right) \frac{\partial}{\partial z_\lambda^\gamma}\wedge \frac{\partial}{\partial z_\lambda^n}\\
& \Longrightarrow \tilde{\Psi} (\Pi_\lambda)= \sum_{\alpha,\beta=1}^{n-1}\Lambda_{\alpha\beta}^\lambda \left(z_\lambda^1,..., z_\lambda^{n-1}, \left(  z_\lambda^n \right)^k \right) \frac{\partial}{\partial \xi_\lambda^\alpha}\wedge \frac{\partial}{\partial \xi_\lambda^\beta} + \frac{1}{k}\sum_{\gamma=1}^{n-1} z_\lambda^n G_\lambda ^\gamma \left(z_\lambda^1,..., z_\lambda^{n-1}, \left( z_\lambda^n \right)^k \right) \frac{\partial}{\partial \xi_\lambda^\gamma}\wedge\left( k \left( z_\lambda^n \right)^{k-1} \frac{\partial}{\partial \xi_\lambda^n} + \frac{\partial}{\partial t_\lambda}    \right)\\
&= \sum_{\alpha,\beta=1}^{n-1}\Lambda_{\alpha\beta}^\lambda \left(z_\lambda^1,..., z_\lambda^{n-1}, \left(  z_\lambda^n \right)^k \right) \frac{\partial}{\partial \xi_\lambda^\alpha}\wedge \frac{\partial}{\partial \xi_\lambda^\beta} + \sum_{\gamma=1}^{n-1}\left( z_\lambda^n \right)^k G_\lambda^\gamma \left(z_\lambda^1,..., z_\lambda^{n-1}, \left( z_\lambda^n \right)^k \right) \frac{\partial}{\partial \xi_\lambda^\gamma} \wedge \frac{\partial}{\partial \xi_\lambda^n} \\
& + \frac{1}{k}\sum_{\gamma=1}^{n-1} z_\lambda^n G_\lambda ^\gamma \left(z_\lambda^1,..., z_\lambda^{n-1}, \left( z_\lambda^n \right)^k \right) \frac{\partial}{\partial \xi_\lambda^\gamma} \wedge \frac{\partial}{\partial t_\lambda} = \Psi^* \left( \tilde{\Lambda}_\lambda \right)
\end{align*}}}
where $\tilde{\Psi}: \bigwedge^2 \Theta_M \to \bigwedge^2 \Psi^* \Theta_F$ is the natural map. On the other hand, 
\begin{align*}
\tilde{\Psi}\left(\Lambda_j^{(v)} \right)= \sum_{\alpha,\beta=1}^n \Lambda_{\alpha\beta}^j \left( z_j^{(v)}   \right) \frac{\partial}{\partial \xi_j^\alpha}\wedge \frac{\partial}{\partial \xi_j^\beta} = \Psi^* \left( \tilde{\Lambda}_j \right)
\end{align*}
which proves the claim. 

\end{proof}

\begin{corollary}[compare \cite{Wav68} Corollary p.929]
If $C\subset \left(W, \Lambda_0\right)$ is a nonsingular Poisson hypersurface $($connected if $k\ne 2$$)$, then a necessary and sufficient condition for $C$ to be the Poisson branch locus of a Poisson $k$-sheeted cyclic covering of $\left(W, \Lambda_0\right)$ is that there is a Poisson line bundle $\left(F, \nabla_F\right)$ over $\left(W, \Lambda_0\right)$ with $\left(F^k, \nabla_{F^k}\right)=\left[ \left(C, \nabla_C\right)\right]$
\end{corollary}

\begin{proof}
Let $C$ be defined by $\phi_i$ on $V_i$ for a covering $W= \left\{V_i \right\}$. Since $C$ is a Poisson hypersurface, $\left[ \Lambda_0, \phi_i  \right] = \phi_i T_i $ for some $T_i \in \Gamma\left( V_i, \Theta_W \right)$. Then $\left[(C, \nabla_C) \right]$ is defined by $\left( \{T_i \} ,\left\{ \frac{\phi_i}{\phi_j} \right\} \right)$. Assume that there is a Poisson line bundle $(F, \nabla_F)$ over $(W, \Lambda_0)$ such that $\left( F^k, \nabla_{F^k} \right) = \left[ \left(C, \nabla_C \right) \right]$. Then $(F, \nabla_F)$ is defined by $\left( \left\{ \frac{1}{k} T_i \right\} , \left\{ f_{ij}\right\}   \right)$ with $f_{ij}^k = \frac{\phi_i}{\phi_j}$. We note that $\nabla_F$ makes $\left(F, \tilde{\Lambda}_0\right)$ a holomorphic Poisson manifold as in the discussion before Theorem \ref{pcc1} and $\left\{t_i^k- \phi_i\right\}$ determine a nonsingular manifold $M$ which is a cyclic $k$-sheeted covering of $W$ having branch locus $C$. Moreover $\left\{ t_i^k- \phi_i\right\}$ defines a Poisson submanifold of $\left(F, \tilde{\Lambda}_0\right)$. In fact, 
\begin{align*}
\left[\Lambda_i - t_i \frac{1}{k} T_i\wedge \frac{\partial}{\partial t_i}, t_i^k- \phi_i   \right] = - \phi_i T_i + t_i^k T_i=  \left(t_i^k- \phi_i \right)T_i
\end{align*}
Hence $\left(M, \tilde{\Lambda}_0|_M\right)$ is a Poisson $k$-sheeted cyclic covering of $\left(W, \Lambda_0\right)$.

Conversely, by Theorem \ref{pcc1}, we can find a Poisson line bundle $\left(F, \nabla_F\right)$ such that $M\subset \left( F, \tilde{\Lambda}_0 \right)$ is a Poisson hypersurface given by $t_j^k=1$ over $U_j\times \mathbb{C}$ and $t_\lambda^k= \xi_\lambda^n$ over $U_\lambda \times \mathbb{C} $ as in the proof of Theorem \ref{pcc1}, and $\left[\left(C, \nabla_C\right) \right]$ is defined by $\left( \left\{ T_j=0  \right\} , \left\{ T_\lambda = - \sum_{\alpha=1}^{n-1} G_\lambda^\alpha \left(\xi_\lambda\right)\frac{\partial}{\partial \xi_\lambda^\alpha}\right\} , \left\{h_{\lambda \mu}= \frac{\xi_\lambda^n}{\xi_\mu^n}, h_{\lambda j}= \frac{\xi_\lambda^n}{1}, h_{ij}=\frac{1}{1}  \right\}   \right)$, so that $\left( F^k, \nabla_{F^k} \right) = \left[ \left(C, \nabla_C\right)  \right]$.
\end{proof}

\bibliographystyle{amsalpha}
\bibliography{References-Rev9}

\providecommand{\bysame}{\leavevmode\hbox to3em{\hrulefill}\thinspace}
\providecommand{\MR}{\relax\ifhmode\unskip\space\fi MR }
\providecommand{\MRhref}[2]{%
  \href{http://www.ams.org/mathscinet-getitem?mr=#1}{#2}
}
\providecommand{\href}[2]{#2}
\begin{thebibliography}{CLN96}

\bibitem[CLN96]{CLN96}
D.~Cerveau and A.~Lins~Neto, \emph{Irreducible components of the space of
  holomorphic foliations of degree two in {$\bold C{\rm P}(n)$}, {$n\geq 3$}},
  Ann. of Math. (2) \textbf{143} (1996), no.~3, 577--612. \MR{1394970}

\bibitem[Har10]{Har10}
Robin Hartshorne, \emph{Deformation theory}, Graduate Texts in Mathematics,
  vol. 257, Springer, New York, 2010. \MR{2583634 (2011c:14023)}

\bibitem[Kim14]{Kim16}
Chunghoon Kim, \emph{Deformations of nonsingular {P}oisson varieties and
  {P}oisson invertible sheaves}, preprint (2014).

\bibitem[Kod05]{Kod05}
Kunihiko Kodaira, \emph{Complex manifolds and deformation of complex
  structures}, english ed., Classics in Mathematics, Springer-Verlag, Berlin,
  2005, Translated from the 1981 Japanese original by Kazuo Akao. \MR{2109686
  (2005h:32030)}

\bibitem[KS53]{KS53}
K.~Kodaira and D.~C. Spencer, \emph{On a theorem of {L}efschetz and the lemma
  of {E}nriques-{S}everi-{Z}ariski}, Proc. Nat. Acad. Sci. U.S.A. \textbf{39}
  (1953), 1273--1278. \MR{68286}

\bibitem[KS58]{Kod58}
\bysame, \emph{On deformations of complex analytic structures. {I}, {II}}, Ann.
  of Math. (2) \textbf{67} (1958), 328--466. \MR{0112154 (22 \#3009)}

\bibitem[LNS20]{LS20}
Alcides Lins~Neto and Bruno Sc\'{a}rdua, \emph{Complex algebraic foliations},
  De Gruyter Expositions in Mathematics, vol.~67, De Gruyter, Berlin, [2020]
  \copyright 2020. \MR{4201904}

\bibitem[Mal98]{Mal98}
Daniel Mall, \emph{On holomorphic and transversely holomorphic foliations on
  {H}opf surfaces}, J. Reine Angew. Math. \textbf{501} (1998), 41--69.
  \MR{1637833}

\bibitem[MK06]{MK06}
James Morrow and Kunihiko Kodaira, \emph{Complex manifolds}, AMS Chelsea
  Publishing, Providence, RI, 2006, Reprint of the 1971 edition with errata.
  \MR{2214741}

\bibitem[OSS11]{OSS88}
Christian Okonek, Michael Schneider, and Heinz Spindler, \emph{Vector bundles
  on complex projective spaces}, 1988 ed., Modern Birkh\"{a}user Classics,
  Birkh\"{a}user/Springer Basel AG, Basel, 2011, With an appendix by S. I.
  Gelfand. \MR{2815674}

\bibitem[Ser06]{Ser06}
Edoardo Sernesi, \emph{Deformations of algebraic schemes}, Grundlehren der
  Mathematischen Wissenschaften [Fundamental Principles of Mathematical
  Sciences], vol. 334, Springer-Verlag, Berlin, 2006. \MR{2247603
  (2008e:14011)}

\bibitem[Wav68]{Wav68}
John~J. Wavrik, \emph{Deformations of {B}anach [branched] coverings of complex
  manifolds}, Amer. J. Math. \textbf{90} (1968), 926--960. \MR{233384}

\end{thebibliography}

\end{document}